\documentclass[12pt]{article}
\usepackage{srcltx}
\usepackage{eurosym}
\usepackage{mathtools}
\usepackage{amsmath}
\usepackage{amsfonts}
\usepackage{amssymb}
\usepackage{amsthm}
\usepackage{graphicx}
\usepackage{mathrsfs}
\usepackage{xcolor}
\usepackage{exscale}
\usepackage{latexsym}
\usepackage{authblk}

\usepackage{bbm}
\usepackage[colorlinks,plainpages=true,pdfpagelabels,hypertexnames=true,colorlinks=true,pdfstartview=FitV,linkcolor=blue,citecolor=red,urlcolor=black]{hyperref}
\PassOptionsToPackage{unicode}{hyperref}
\PassOptionsToPackage{naturalnames}{hyperref}
\usepackage{enumerate}
\usepackage[shortlabels]{enumitem}
\usepackage{bookmark}
\usepackage{wasysym}
\usepackage{esint}
\usepackage[ddmmyyyy]{datetime}
\usepackage[margin=1in]{geometry}
\makeatletter
\g@addto@macro\normalsize{%
	\setlength\abovedisplayskip{4pt}
	\setlength\belowdisplayskip{4pt}
	\setlength\abovedisplayshortskip{4pt}
	\setlength\belowdisplayshortskip{4pt}
}
\numberwithin{equation}{section}
\everymath{\displaystyle}
\usepackage[capitalize,nameinlink]{cleveref}
\crefname{section}{Section}{Sections}
\crefname{subsection}{Subsection}{Subsections}
\crefname{condition}{Condition}{Conditions}
\crefname{hypothesis}{Hypothesis}{Conditions}
\crefname{assumption}{Assumption}{Assumptions}
\crefname{lemma}{Lemma}{Lemmas}
\crefname{definition}{Definition}{Definitions}

\crefformat{equation}{\textup{#2(#1)#3}}
\crefrangeformat{equation}{\textup{#3(#1)#4--#5(#2)#6}}
\crefmultiformat{equation}{\textup{#2(#1)#3}}{ and \textup{#2(#1)#3}}
{, \textup{#2(#1)#3}}{, and \textup{#2(#1)#3}}
\crefrangemultiformat{equation}{\textup{#3(#1)#4--#5(#2)#6}}%
{ and \textup{#3(#1)#4--#5(#2)#6}}{, \textup{#3(#1)#4--#5(#2)#6}}%
{, and \textup{#3(#1)#4--#5(#2)#6}}

\Crefformat{equation}{#2Equation~\textup{(#1)}#3}
\Crefrangeformat{equation}{Equations~\textup{#3(#1)#4--#5(#2)#6}}
\Crefmultiformat{equation}{Equations~\textup{#2(#1)#3}}{ and \textup{#2(#1)#3}}
{, \textup{#2(#1)#3}}{, and \textup{#2(#1)#3}}
\Crefrangemultiformat{equation}{Equations~\textup{#3(#1)#4--#5(#2)#6}}%
{ and \textup{#3(#1)#4--#5(#2)#6}}{, \textup{#3(#1)#4--#5(#2)#6}}%
{, and \textup{#3(#1)#4--#5(#2)#6}}

\crefdefaultlabelformat{#2\textup{#1}#3}
\newtheorem{theorem} {Theorem}[section]
\newtheorem{proposition} [theorem]{Proposition}

\newtheorem{lemma}[theorem]{Lemma}
\newtheorem{corollary}[theorem]{Corollary}

\newtheorem{counter example}[theorem]{Counter Example}
\newtheorem{remark}[theorem] {Remark}

\newtheorem{claim}[theorem] {Claim}

\def\CC{{\rm \kern.24em \vrule width.02em height1.4ex depth-.05ex \kern-.26emC}}

\def\TagOnRight

\def\AA{{it I} \hskip-3pt{\tt A}}

\def\QQ{\rlap {\raise 0.4ex \hbox{$\scriptscriptstyle |$}} {\hskip -0.1em Q}}

\makeatletter
\newcommand{\vo}{\vec{o}\@ifnextchar{^}{\,}{}}
\makeatother
\def\YYint#1#2#3{{\setbox0=\hbox{$#1{#2#3}{\iint}$}
		\vcenter{\hbox{$#2#3$}}\kern-.50\wd0}}
\def\XXint#1#2#3{{\setbox0=\hbox{$#1{#2#3}{\int}$}
		\vcenter{\hbox{$#2#3$}}\kern-.50\wd0}}

\makeatletter
\def\namedlabel#1#2{\begingroup
	\def\@currentlabel{#2}%
	\label{#1}\endgroup
}
\makeatother
\makeatletter
\newcommand{\rmh}[1]{\mathpalette{\raisem@th{#1}}}
\newcommand{\raisem@th}[3]{\hspace*{-1pt}\raisebox{#1}{$#2#3$}}
\makeatother

\newcounter{desccount}
\newcommand{\descitem}[2]{\item[#1]\refstepcounter{desccount}\label{#2}}
\newcommand{\descref}[2]{\hyperref[#1]{\textnormal{\textcolor{black}{}\textcolor{blue}{ #2}\textcolor{black}{}}}}
\newcommand{\dref}[2]{\hyperref[#1]{\textcolor{black}{(}\textcolor{blue}{\bf #2}\textcolor{black}{)}}}
\newcommand{\be} {\begin{eqnarray}}
	\newcommand{\ee} {\end{eqnarray}}
\newcommand{\Bea} {\begin{eqnarray*}}
	\newcommand{\Eea} {\end{eqnarray*}}
\newcommand{\pa} {\partial}

\newcommand{\al} {\alpha}
\newcommand{\rr}{\rightarrow}

\newcommand{\B} {\beta}
\newcommand{\de} {\delta}

\newcommand{\p}  {\prime}
\newcommand{\e}  {\varepsilon}
\newcommand{\De} {\Delta}

\newcommand{\la} {\lambda}
\newcommand{\si} {\sigma}

\newcommand{\La} {\Lambda}
\newcommand{\f}{\infty}
\newcommand{\R}{\mathbb{R}}
\newcommand{\eps} {\epsilon}

\newcommand{\noi} {\noindent}

\newcommand{\va} {\varphi}
\newcommand{\ga}{\gamma}
\newcommand{\ep}{\varepsilon}

\newcommand{\sig}{\sigma}
\newcommand{\ta}{\tau}

\newcommand{\norm}[1]{\left|\hspace{-0.2mm}\left| #1 \right|\hspace{-0.2mm}\right|}
\newcommand{\abs}[1]{\left| #1\right|}

\newcounter{whitney}
\refstepcounter{whitney}

\newcounter{ineqcounter}
\refstepcounter{ineqcounter}
\makeatletter
\def\ps@pprintTitle{%
	\let\@oddhead\@empty
	\let\@evenhead\@empty
	\def\@oddfoot{}%
	\let\@evenfoot\@oddfoot}
\makeatother
\usepackage[doublespacing]{setspace}
\usepackage[titletoc,toc,page]{appendix}

\makeatletter
\newcommand{\refcheckize}[1]{%
	\expandafter\let\csname @@\string#1\endcsname#1%
	\expandafter\DeclareRobustCommand\csname relax\string#1\endcsname[1]{%
		\csname @@\string#1\endcsname{##1}\wrtusdrf{##1}}%
	\expandafter\let\expandafter#1\csname relax\string#1\endcsname
}
\makeatother
\refcheckize{\cref}
\refcheckize{\Cref}


\makeatletter
\newcommand{\mainsectionstyle}{%
	\renewcommand{\@secnumfont}{\bfseries}
	\renewcommand\section{\@startsection{section}{2}%
		\z@{.5\linespacing\@plus.7\linespacing}{-.5em}%
		{\normalfont\bfseries}}%
}
\makeatother
\usepackage{pgf,tikz}
\usetikzlibrary{arrows}
\usetikzlibrary{decorations.pathreplacing}
\usepackage{xpatch}
\xpatchcmd{\MaketitleBox}{\hrule}{}{}{}% remove first horizontal rule (above abstract)
\xpatchcmd{\MaketitleBox}{\hrule}{}{}{}% remoce second horizonral rule (below keywords)

\date{}
\usepackage{scalerel}
\makeatletter

\title{Vanishing viscosity limit for $n\times n$ strictly hyperbolic system of conservation laws in 1-d with nonlinear viscosity: Part-I Uniform BV estimates}

\author[1,a]{Boris Haspot}
\author[2,b,c]{Animesh Jana}
\affil[a]{\footnotesize	 Universit\'e Paris Dauphine, PSL Research University, Ceremade (UMR CNRS 7534), Place du Mar\' echal De Lattre De Tassigny 75775 Paris cedex 16, France.}
\affil[b]{Harish-Chandra Research Institute, A CI of Homi Bhabha National Institute, Chhatnag Road, Jhunsi, Prayagraj 211019, India.}
\affil[c]{\footnotesize Istituto di Matematica Applicata e Tecnologie Informatiche ``Enrico Magenes'' - Consiglio Nazionale delle Ricerche, Via Ferrata, 5/a - 27100 Pavia, Italy.}
\affil[1]{\em \footnotesize	 haspot@ceremade.dauphine.fr}
\affil[2]{\em \footnotesize	animeshjana@hri.res.in}

\allowdisplaybreaks

\usepackage{accents}
\newlength{\dhatheight}

	\DeclarePairedDelimiter\floor{\lfloor}{\rfloor}

\begin{document}
\maketitle

\begin{abstract}
	We consider  the following parabolic approximation for hyperbolic system of conservation laws in 1-D with non-singular viscosity matrix $B(u)$ and $A(u)$ strictly hyperbolic,
	\[u^\e_t+A(u^\e)u^\e_x=\varepsilon(B(u^\e)u^\e_x)_x.\] We prove global in time uniform $BV$ bound for solution to this parabolic system when $\varepsilon>0$ provided that the initial data is small in $BV$ and the matrix $A(u)$ and $B(u)$ commutate. Moreover, in the case where the system is conservative, we show that the sequence 
$(u^\e)_{\e>0}$  admits a limit 
$u$, which is the unique global weak solution to the limiting strictly hyperbolic system. We provide a concrete application of this result in the study of the visco-dispersive limit of the Navier-Stokes-Korteweg system.
\end{abstract}
\tableofcontents
\section{Introduction}
In this article, we consider the following parabolic approximation for strictly hyperbolic system of conservation laws in 1-D with invertible viscosity matrix $B(u)$
\begin{equation}
	\begin{cases}
		\begin{aligned}
			&u^\e_t+A(u^\e)u^\e_x=\e(B(u^\e)u^\e_x)_x \quad \quad\,\mbox{ for }t>0,x\in\R,%\label{eqn-parabolic}
			\\
			&u^\e(0,x)=\bar{u}(x)\hspace{3,2cm}\mbox{for }x\in\R,\label{eqn-parabolic}%\label{eqn:initial-data}
		\end{aligned}
	\end{cases}
\end{equation}
where $u^\e:[0,\f)\times\R\rr\R^n$ and $A,B$ are $n\times n$ matrices satisfying the following conditions for some open $\mathcal{U}\subset\R^n$.
\begin{description}
	\descitem{($\mathcal{H}_A$.)}{A1} The matrix $A(u)$ is smoothly depending on $u$  and has $n$ distinct real eigenvalues $\la_1(u)<\cdots<\la_n(u)$ for $u\in\mathcal{U}$ such that it exists $c_0>0$ satisfying:
	\begin{equation}
		0<c_0\leq \sup_{u\in \mathcal{U}}\max_{i\ne j} (|\la_i(u)-\la_j(u)|).
		\label{stricthyperbo}
	\end{equation}
	% when $\abs{u-u^*}\leq \rho$ for some $u^*\in\R^n$ and $\rho>0$.
	\descitem{($\mathcal{H}_B$.)}{B} The matrix $B(u)$ is smoothly depending on $u$ and  $B(u)$ has $n$ positive regular real eigenvalues $\{\mu_i(u)\}_{i=1}^{n}$ satisfying $\mu_i(u)\geq c_1$ for $u\in\mathcal{U}$ and for some $c_1>0$. In addition these eigenvalues are associated to the $n$ following eigenvectors $(r_1(u),\cdots,r_n(u))$ with for any $i\in\{1,\cdots,n\}$, $\|r_i(u)\|=1$. 
	
\end{description}
For a general system, proving that we can select the physical solution of a strictly hyperbolic system via a process of vanishing viscosity limit has been during longtime a challenging open problem. One of the main difficulty consists in
obtaining uniform total variation bound with respect to $t$ and $\ep$ for a solution $u^\e$ of \eqref{eqn-parabolic}. In their seminal paper \cite{BB-vv-lim-ann-math}, Bianchini and Bressan solved this problem for $B=I$ by developing a new impressive $L^1$ theory. However since in many physical systems,
the viscosity depends generally on the macroscopic variables, it seems important to extend the results of \cite{BB-vv-lim-ann-math} to general invertible matrix $B(u)$. It corresponds in particular to the Problem 2 that proposes Bressan in his survey paper on open problems for hyperbolic conservation laws in one dimension \cite{survey}.

Our goal in this paper is to establish that the solution $u^\e$ of  \eqref{eqn-parabolic} satisfies global in time uniform total variation estimate and the bound is independent of $\e$ provided that the initial data $\bar{u}$ is sufficiently small in $TV$ and that $A(u)$ and $B(u)$ commutate.
 We observe by rescaling the coordinate as follows $(t,x)\mapsto (\frac{t}{\e},\frac{x}{\e})$, it is sufficient to study the following system
\begin{equation}\label{eqn-main}
\begin{cases}
\begin{aligned}
	&u_t+A(u)u_x=(B(u)u_x)_x\\
	&u(0,\cdot)=\bar{u}(\e\cdot).
	\end{aligned}
	\end{cases}
\end{equation}
Indeed this change of coordinate preserves the total variation norm of the initial data. Our problem is then restricted to prove global in time $TV$ bound for the system \eqref{eqn-main}, it will enable us to obtain uniform $TV$ bound in $\e$ for the solution $u^\e$ of system \eqref{eqn-parabolic}. We recall that it is one of
the main step in order to obtain some results of vanishing viscosity limit to \eqref{eqn-parabolic}, indeed by using compactness argument of Helly type,  we expect  that up to a subsequence $(u^\e)_{\e>0}$ converges to a limit $u$ which % the limit $\lim_{\eps\rightarrow 0}u_\eps=u$ (with $u_\eps$ solution of \eqref{eqn-parabolic}) 
satisfies the following hyperbolic system in a  sense to precise,
\begin{equation}\label{eqn-hyperbolic}
	u_t+A(u)u_x=0.
\end{equation}
We should point out that this system is a priori not conservative, it implies that we do not expect existence of global weak solution for such system. However we can define a solution to such system associated to a vanishing viscosity process provided that the solution $u^\e$ of  \eqref{eqn-parabolic} have a limit, in other words that the sequence $(u^\e)_{\e>0}$ has a unique point of accumulation. It is exactly the result that Bianchini and Bressan have proved in \cite{BB-vv-lim-ann-math} when $B(u)=I$. When the system is conservative which means $A(u)=Df(u)$ with $f$ regular then we can show that $u$ %the limit $\lim_{\eps\rightarrow 0}u_\eps=u$ 
is a global weak solution of \eqref{eqn-hyperbolic}
and this limit belongs to a class of uniqueness defined by Bressan and De Lellis in \cite{BDL}. It implies in particular that $(u^\e)_{\e>0}$ has a unique accumulation point and then admits a limit which is the unique solution of \eqref{eqn-hyperbolic}.

Before coming back on these different issues, we would like to recall some important features of the conservation laws. Hyperbolic conservation laws appear in different areas of science and engineering fields such as fluid mechanics, traffic flows etc. We refer to the book of Dafermos \cite{Dafermos} for more on applications of conservation laws.

One of the important characteristic of hyperbolic conservation laws is that even from a smooth initial data, discontinuity can appear in the solution in finite time. That is why we consider solution in some functional spaces which allows to have some discontinuities. For scalar conservation laws, well-posedness for BV initial data has been established by Lax \cite{Lax} and Oleinik \cite{Oleinik} in one space dimension. In multi-D,  Kruzkov \cite{Kruzkov} proved the well-posedness for $L^\f$ initial data. For the system of hyperbolic conservation laws, Glimm \cite{Glimm} made the breakthrough by proving existence of BV solution in 1-D for initial data with small total variation provided that the field are genuinely non linear or linearly degenerate. It has been extended by Bressan and Bianchini in \cite{BB-vv-lim-ann-math} for the general case where the matrix $A(u)$ is assumed only strictly hyperbolic.

For a general class of hyperbolic conservation laws, well-posedness is studied in \cite{BrCo-semi,Bre-Cra-Pic, Bre-Mar-Vai}. For large data, existence of BV solutions are limited to special systems (see \cite{Bressan-Goatin-2000,Nishida-Smoller}). For a general system, obtaining total variation bound for large initial data is still an open question (see \cite{Bre-Chen-Zhang} for more explanations on how difficult in this question). On the other hand, the uniqueness of BV solution has been studied by Bressan-Liu-Yang \cite{Bre-Liu-Yang} ( where the authors proved by a simpler way as in \cite{BrCo-semi,Bre-Cra-Pic} the $L^1$ contraction principle)
and Bressan-LeFloch \cite{Bressan-LeFloch} for $n\times n$ hyperbolic system for initial data with small total variation (see also the article by Bressan and Goatin \cite{Bressan-Goatin-1999}). We also refer to the book of Bressan \cite{Bressan-book}.
Uniqueness for large initial data has been established in \cite{Bressan-Colombo-uniq} for $2\times 2$ systems. Later, the requirement on solutions for belonging to the class of uniqueness has been improved in \cite{BDL} by removing the so called ``tame condition".
By ``a-contraction" method, stability of BV solution has been proved in \cite{Chen-Krupa-Vas} in the regime of entropy solutions satisfying strong trace condition. The uniqueness problem becomes again more challenging for \eqref{eqn-hyperbolic} in the non-conservative case since we lose the notion of global weak solution. For Temple class $A(u)$ (that is when $r_i(u)\cdot Dr_i(u)=0$ for all $i$) and for general hyperbolic matrix $A(u)$, stability of vanishing viscosity solution has been proved in \cite{BB-temple} and \cite{BB-vv-lim-ann-math} respectively. We refer to \cite{Bianchini-arma2003} where Riemann solver has been constructed for non-conservative cases. 

%This is why the natural space where we study the 

The viscosity approximation plays an important role in the study of hyperbolic conservation laws. It not only provides the existence of entropy solutions, but it allows also to select the physically relevant solutions via a vanishing viscosity process. For scalar conservation laws, vanishing viscosity limit has been established by Oleinik \cite{Oleinik} in 1-D and later Kruzkov \cite{Kruzkov} proved it in multi space dimension for $L^\f$ initial data by using essentially maximum principle.

For $n\times n$  hyperbolic system of conservation laws in one space dimension, the problem is again more challenging since we cannot directly obtained such results by maximum principle.
% vanishing viscosity limit satisfies the entropy condition and it is unique. We also refer to \cite{} for detailed discussion on scalar case. 
For $2\times 2$ system, the vanishing viscosity limit was studied by DiPerna \cite{DiPerna} via compensated compactness. This study has been extended by Chen and Perepelitsa \cite{Chen-Per} where the authors consider  the vanishing viscosity limit for the compressible Navier-Stokes equations. More recently Chen, Kang and Vasseur in \cite{Chen-Kang-Vas} proves the convergence when the viscosity $\e>0$ goes to $0$ of the solution of the compressible Navier-Stokes equations to a unique solution of the Euler system in one dimension resolving a long outstanding open problem (even if the Navier-Stokes system is a $2\times 2$ system this question is particularly delicate since the matrix $B(u)$ in this framework is degenerate).
For $n\times n$ hyperbolic system the vanishing viscosity has been studied by Goodman and Xin \cite{GX}. They established the viscosity limit for Riemann problem solutions. For Temple system it has been studied first by Serre \cite{Serre-1}. Later, Bianchini and Bressan \cite{BB-temple} established a global existence of the parabolic-hyperbolic equation \eqref{eqn-parabolic} when $A$ belongs to the Temple class and $B(u)=I$ and study the vanishing viscosity limit of such solution even when the system is not conservative. In \cite{HJ-temple-class}, we generalize their approach to the case of non-constant viscosity matrix $B$ satisfying an assumption of commutativity $AB=BA$. On the other hand, when $A(u)$ is a $2\times2$ triangular strictly hyperbolic matrix vanishing viscosity limit has been obtained for $B(u)=I$ and non-constant $B(u)$ satisfying $AB=BA$ in \cite{BB-triangular} and \cite{HJ-triangular} respectively (see also \cite{Spinolo}).

For a general system, proving that we can select the physical solution of a strictly hyperbolic system via a process of vanishing viscosity limit has been during longtime a challenging open problem. As we mentioned previously, one of the main difficulty consists in
obtaining uniform total variation bound with respect to $t$ and $\ep$ for a solution $u^\e$ of \eqref{eqn-parabolic}. In \cite{BB-vv-lim-ann-math}, Bianchini and Bressan solved this problem for $B=I$ by developing a completely new 
$L^1$  theory incorporating, in particular, the stability of small viscous shocks. One of their main ideas is to decompose the gradient of  solution $u_x$ in a basis of travelling waves selected by a center manifold techniques as $u_x=\sum_iv_i\tilde{r}_i$ and $u_t=\sum_i(w_i-\la_i^*v_i)\tilde{r}_i$. We point out that this decomposition makes sense only for long time. From a certain point of view, it is a matter of waiting for the parabolic effects on the solution, since it will regularize and its derivatives will become sufficiently small for sufficiently large times.
By this way, it enables to diagonalize the system satisfied by $u_x, u_u$, more precisely the components $v_i, w_i$ satisfy some equations of the form
$$
\begin{aligned}
	&v_{i,t}+(\tilde{\la}_i v_i)_x-v_{i,xx}=\phi_i,\\
	&w_{i,t}+(\tilde{\la}_i w_i)_x-w_{i,xx}=\psi_i,
\end{aligned}$$
where the remainder terms $\phi_i, \psi_i$ can be controlled by quantities of the form
\begin{equation*}
	\sum\limits_{i\neq j}|v_iv_j|,\quad (w_iv_{i,x}-w_{i,x}v_i),\quad v_i^2\left(\frac{w_i}{v_i}\right)_x^2\mathbbm{1}_{\{|w_i|/|v_i|\leq\de_1\}}\mbox{ and }v_{i,x}^2 \mathbbm{1}_{\{|w_i|/|v_i|>\de_1\}}.
\end{equation*}
In order to control the total variation of $u_x$ by maximum principle it is sufficient now to estimate the $L^{1}_{t,x}$ norm of $\phi_i$ and $\psi_i$ (we mention in particular that getting $L^1$ norm in time on $\phi_i$, $\psi_i$ requires in a natural way to take into account the solutions of type travelling waves, it explains why Bianchini and Bressan operate a decomposition as gradient of travelling waves).
By their method of transversal interaction estimate \cite{BB-temple}, they established the $L^1_tL^1_x$ estimate for the first types of the terms. For terms of second and third types that appear due to interaction of same family viscous waves, they used the Lyapunov functional \cite{BB-Lyapunov} related to shortening of curve. Last type of the terms are also proved to be small in $L^1_tL^1_x$ by using length estimate and energy estimate. We also refer to the work of Christoforou \cite{Christoforou} in the more general framework of hyperbolic systems of balance laws always with $B(u)=I$.

In our case, that is when $B$ is depending $u$, several difficulties appear. First and foremost, we remark the importance of $AB=BA$ condition, which says that $A$ and $B$ share the same set of eigenvectors. This property guarantees the coupling arising due to viscosity matrix $B$ is quadratic, by which we mean that in the parabolic equation of the $v_i$, the terms involving $v_{j,xx}$ do appear in  the form of $w_{j,xx}v_j-v_{j,xx}w_j$. In this paper we are going to show that terms $w_{j,xx}v_j-v_{j,xx}w_j$ are small in $L^1_tL^1_x$. We remark that without the condition $AB=BA$, it can appear $ v_{j,xx}$ (with $j\neq i$) in the advection diffusion equation of $v_i$. In addition it is interesting to point out that this condition of commutativity is reminiscent of the Chuey-Conley-Smoller theory for positively invariant regions of nonlinear diffusion equations (see \cite{CCS}).
Although in our case the matrix $B$ is non-singular, we would also like to comment that the assumption $AB=BA$ prevents the parabolic-hyperbolic equation \eqref{eqn-parabolic} to satisfy the Kawashima-Shizuta condition \cite{Kawa-Shi} when $B$ is singular. 

For non-constant $B(u)$, the major difficulty is to get an $L^1_tL^1_x$ estimate for $w_{i,xx}v_i-v_{i,xx}w_i$. To obtain an $L^1_tL^1_x$ estimate of this term, one way to approach is to invoke the shorten estimate of Bianchini and Bressan \cite{BB-vv-lim-ann-math} two times with the pairs $(v_i,w_{i,x})$ and $(v_{i,x},w_i)$. The key hurdle one will face in this approach is to get the bounds for forcing parts of the advection-diffusion equations satisfied by $w_{i,x}$ and $v_{i,x}$. There are two issues appear if we wish to calculate directly the advection-diffusion equations of $w_{i,x}$ and $v_{i,x}$: (1) It can be checked that one need to have a bound for terms of the form $|v_{i,x}|/|v_i|$ on the set $\{|w_i|/|v_i|\leq \de_1\}$. (2) There appear terms of the form $w_{i,xxx}v_i-v_{i,xxx}w_i$ which prevents to end the loop of estimates. 

This is a challenging situation to deal with the issue (1) since for a general initial data if we decompose the gradient $u_x$ as in Bianchini and Bressan \cite{BB-vv-lim-ann-math}, $v_{i,x}$ can only be estimated as $v_{i,x}=O(1)(|w_i|+|v_i|)+O(1)\sum\limits_{j\neq i}\de_0|v_j|$. Note that we can have desired estimate only if for each $j\ne i$, the quantity $v_j$ is bounded in terms of $v_i$. To ensure this, we modify the basis $\tilde{r}_i$ in the decomposition of $u_x$ and $u_t$ such a way that when $v_i$ is dominant field (that is, $|v_i|^2\geq |v_j|$ for all $j\neq i$) then $\tilde{r}_i$ is the direction of travelling wave and in this case $\tilde{r}_j,j\neq i$ reduces to the direction of $j$-th eigenvector of $A$ if $|v_j|\lesssim |v_i|^2$. Furthermore, we need a cutoff function as well to make $\tilde{r}_i$ as eigenvector of $A$ when $|w_i|/|v_i|$ is large enough (it is important to point out that it is conform with a solution of type travelling wave). 
These choices are in fact suggested by the coupling formula between $\mu_i v_{i,x}$ and $(\tilde{\la}_i-\la_i^*)v_i+w_i$ as in \eqref{eqn-v-i-x-1} which plays a crucial role. For $B(u)=I$ this has been observed in \cite{BB-vv-lim-ann-math} where $v_{i,x}$ is estimated by $w_i,v_i$ and $v_j (j\neq i)$. In our case, this can be reduced to $v_{i,x}=O(1)|v_i|$ in the set when $|w_i|\leq \de_1|v_i|$ and $\{|v_i|>v_j^2;j\neq i\}$.% Furthermore, we also calculate and estimate $v_{i,xx}$ and $v_{i,xxx}$ by using the formula \eqref{eqn-v-i-x-1} in requirement of our analysis. 

To summarize the properties of our new basis, we would like to mention that whenever $i$-th component is dominant and $|w_i|/|v_i|$ is small then we use direction of travelling wave otherwise we use the eigenvector of $A(u)$. 

Even with this new basis, we still face the issue (2) as directly differentiation of the advection-diffusion equation of $v_i$ and $w_{i}$ gives terms of the form $w_{i,xxx}v_i-v_{i,xxx}w_i$. To resolve this, we introduce a special variables $z_i$ and $\hat{z}_i$ as below
\begin{align}
	z_i&:=\mu_iv_{i,x}-(\tilde{\la}_i-\la_i^*)v_i+\sum\limits_{j\neq i}a_{ij}\left(w_{j,x}-\frac{w_j}{v_j}v_{j,x}\right),\label{intro-definition-z-i}\\
	\hat{z}_i&:=\mu_iw_{i,x}-(\tilde{\la}_i-\la_i^*)w_i+\sum\limits_{j\neq i}\hat{a}_{ij}\left(w_{j,x}-\frac{w_j}{v_j}v_{j,x}\right).\label{intro-definition-z-hat-i}
\end{align}
where $a_{ij}$ and $\hat{a}_{ij}$ are chosen in a suitable way in order to cancel out the extra terms $w_{i,xxx}v_i-v_{i,xxx}w_i$. One can think of these quantities as effective fluxes for $v_i$ and $w_i$. Now we observe that combining the shorten estimates corresponding to the pairs $(z_i,w_i)$ and $(\hat{z}_i,v_i)$ we get the $L^1_tL^1_x$ estimate of $w_{i,xx}v_i-v_{i,xx}w_i$. 

To derive the PDEs satisfied by the variables $z_i$ and $\hat{z}_i$ we do take differentiation of advection-diffusion equation of $w_{i}$ and $v_i$. Although due to the choices of $z_i$ and $\hat{z}_i$ there is no term of the form $w_{i,xxx}v_i-v_{i,xxx}w_i$ but it does have terms like $v_{i,xx}v_j$ for $j\neq i$. We can estimate such error terms by transversal estimate of Bianchini and Bressan \cite{BB-temple,BB-vv-lim-ann-math} by using the interaction between $z_i$ and $v_j$. We would also like to mention that $w_{j,xxx}, v_{j,xxx}$ do appear in the forcing part of the equation of $z_i$ and $\hat{z}_i$ but they appear in the following form
\begin{equation}\label{eqn-3rd-order-error}
	\sum\limits_{j\neq i}|v_{i,xxx}v_iv_j|\mbox{ and }	\sum\limits_{j\neq i}|w_{i,xxx}v_iv_j|.
\end{equation}
We estimate such term by using the equation, non transversal interaction and  bootstrap argument. Indeed we are not in measure to directly estimating the $L^\infty$ norm of the third derivatives terms
$v_{i,xxx}v_i$. We refer to section \ref{section:Lambda-2} for more details. 

Apart from the above-mentioned key technical challenges, there are a few other new estimates and observations that we would like to discuss. We purposefully use a cutoff function to work with only in the region $v_i\neq 0$. We note that $v_{i,xx}$ and $w_{i,xx}$ becomes singular on the set $\{v_i=w_i=0\}$. This is a challenge which is also present in the analysis of Bianchini-Bressan  \cite{BB-vv-lim-ann-math} in the case of $B(u)=I$. They resolved this by taking a smooth approximation of $v_i$ and $w_i$. In our paper, we try to deal with this problem in a slightly different manner. We incorporate a cut-off in the basis function $\tilde{r}_i$ such that on the set $\{v_i=w_i=0\}$, $\tilde{r}_i$ becomes the $i$-th eigenvector of $A(u)$. This enables us to ensure the regularity of $w_i$ and $v_i$ on whole $\R$. 

In our analysis, the coupling formula between $\mu_i v_{i,x}$ and $(\tilde{\la}_i-\la_i^*)v_i+w_i$ as in \eqref{eqn-v-i-x-1} has a key role as in \cite{BB-vv-lim-ann-math}, however in our setting, we also need to differentiate the coupling formula in order to estimate, in particular, the following terms $v_{i,xx}$ and $v_{i,xxx}$ by using the formula \eqref{eqn-v-i-x-1}. 

On the other hand, since we need to deduce the PDEs for $z_i$ and $\hat{z}_i$, the variables $v_i$ and $w_i$ are needed to have bound on the second order derivatives of $v_i$ and $w_i$. In other words, we need to have higher order derivatives of $u_x$ and $u_t$. For non-constant viscosity matrix, the regularity of $u_x$ is established by using parabolic regularizing effects as in \cite{HJ-temple-class} by introducing suitable change of unknown in order to recover the case $B(u)=I$. We follow the line of arguments as in \cite{HJ-temple-class} but for higher order derivatives of $u_x$ which makes it more involved.

\section{Main result}
We state our main result of this article which concerns about global existence of smooth solutions to \eqref{eqn-main} for small BV initial data.
\begin{theorem}\label{theorem:BV-estimate}
	Consider the Cauchy problem hyperbolic system with viscosity,
	\begin{equation}\label{eqn-thm-1}
		u_t+A(u)u_x=(B(u)u_x)_x,\quad u(0,x)=\bar{u}(x).
	\end{equation}
	We assume that the drift $A$ satisfies \descref{A1}{($\mathcal{H}_A$.)} and viscosity matrix $B$ verifies \descref{B}{($\mathcal{H}_B$.)}. Furthermore, we assume that
	\begin{equation}
		A(u)B(u)=B(u)A(u)\mbox{ for all }u\in \mathcal{U}.
	\end{equation} There exist $L_1,L_2>0$ and $\de_0>0$ such that the following holds. If $\bar{u}$ satisfies
	\begin{equation}\label{condition-data-thm-1}
		TV(\bar{u})\leq\de_0\mbox{ and }\lim\limits_{x\rr-\f}\bar{u}(x)\in \mathcal{K},
	\end{equation}
	for some compact set $\mathcal{K}\subset\mathcal{U}$ then there exists unique solution $u$ to the Cauchy problem \eqref{eqn-thm-1} and it satisfies the following properties
	\begin{align}
		TV(u(t))&\leq L_1TV(\bar{u}),\label{thm-1:BV-bound}\\
		\norm{u(t)-u(s)}_{L^1}&\leq L_2\left(\abs{t-s}+ \abs{\sqrt{t}-\sqrt{s}}\right).\label{L1-cont}
	\end{align} 
	%\begin{equation}\label{thm-1:BV-bound}
	%		TV(u(t))\leq L_1TV(\bar{u}).
	%	\end{equation} 
	\end{theorem}
Next, let us consider \eqref{eqn-parabolic} in the conservative case, that is when $A(u)=Df(u)$ for some smooth function $f:\R^n\rr\R^n$,
	\begin{equation}
		u^\eps_t+(f(u^\eps))_x=\eps(B(u)u^\eps_x)_x\mbox{ and }u^\eps(0,x)=\bar{u}(x).
	\end{equation}  
	We have the following result as an application of Theorem \ref{theorem:BV-estimate}.
	\begin{corollary}\label{corollary:vv-limit}
		Consider the Cauchy problem for hyperbolic system with viscosity,
		\begin{equation}\label{eqn-coro-vis-1}
			u^\eps_t+(f(u^\eps))_x=\eps(B(u)u^\eps_x)_x\mbox{ and }u^\eps(0,x)=\bar{u}(x).
		\end{equation}  
		We assume that the drift $A(u)=Df(u)$ satisfies \descref{A1}{($\mathcal{H}_A$.)} and viscosity matrix $B$ verifies \descref{B}{($\mathcal{H}_B$.)}. Furthermore, we assume that
		\begin{equation}
			A(u)B(u)=B(u)A(u)\mbox{ for all }u\in \mathcal{U}.
		\end{equation} 
		There exist $L_1,L_2>0$ and $\de_0>0$ such that the following holds. If $\bar{u}$ satisfies
		\begin{equation}
			TV(\bar{u})\leq\de_0\mbox{ and }\lim\limits_{x\rr-\f}\bar{u}(x)\in \mathcal{K},
		\end{equation}
		for some compact set $\mathcal{K}\subset\mathcal{U}$ then for every $\eps>0$ there exists unique solution $u^\eps$ to the Cauchy problem \eqref{eqn-coro-vis-1} and it satisfies the following properties
		\begin{align}
			TV(u^\eps(t))&\leq L_1TV(\bar{u}),\label{coro-1:BV-bound}\\
			\norm{u^\eps(t)-u^\eps(s)}_{L^1}&\leq L_2\left(\abs{t-s}+ \sqrt{\eps}\abs{\sqrt{t}-\sqrt{s}}\right).\label{coor-L1-cont}
		\end{align} 
		Furthermore, as $\eps\rr0$, $u^{\eps}(t,\cdot)$ converges in $L^1_{loc}(\R)$ to $u^\f(t,\cdot)$ for each $t>0$, which is the unique global weak solution to the following Cauchy problem for the hyperbolic system of conservation laws
		\begin{equation}\label{eqn-hyperbolic1}
			u^\f_t+f(u^\f)_x=0,\mbox{ for }(t,x)\in\R_+\times \R\mbox{ and }u^\f(0,x)=\bar{u}(x), \mbox{ for }x\in\R,
		\end{equation}
		and satisfying the {\it Liu admissibility condition} for shock. Furthermore, $u^\f$ satisfies the following bounds
		%there exist a subsequence $\eps_k\rr0$ such that $u^{\eps_k}$ converges to a weak solution $u^\f$ to the following 
		%with $u^\f(0,\cdot)=\bar{u}$ and the following BV bound and time continuity hold for $u^\f$,
		\begin{align*}
			TV(u^\f(t))&\leq L_1 TV(\bar{u}),\\
			\norm{u^\f(t)-u^\f(s)}_{L^1}&\leq L_2 \abs{t-s}.
		\end{align*}
		%Moreover, $u^\f$ is the unique weak solution from initial data $\bar{u}$ satisfying the {\it Liu admissibility condition}.
		%		
		%		
		%		 as $\eps\rr0+$ up to a subsequence $u^\eps$ converges to
		%		
		%		Let $\{u^\eps\}_{\eps>0}$ be solving \eqref{viscous-eqn-conservative} with $u^\eps(0,\cdot)=\bar{u}(\cdot)$ with $A=Df$ and $B$ as in Theorem \ref{theorem:BV-estimate}. Assume that the initial data $\bar{u}$ satisfies the hypothesis \eqref{condition-data-thm-1}. Then $u^\ep$ converges to the unique weak solution to \eqref{eqn-hyperbolic1} satisfying Liu condition.
	\end{corollary}
We would now like to provide a concrete physical example to which the Corollary \ref{corollary:vv-limit}
can be applied. Let us consider the following Navier-Stokes Korteweg system (see \cite{Germain,Has} and the references therein):
\begin{equation}
\begin{cases}
\begin{aligned}
&\rho_t+(\rho u)_x=0,\\
&(\rho u)_t+(\rho u^2)_x-(\mu(\rho) u_x)+(P(\rho))_x=c(K[\rho])_x,\\
&(\rho,u)(0,\cdot)=(\rho_{0},u_{0}),
\end{aligned}
\end{cases}
\label{3systemea}
\end{equation}
with $c>0$.
Here $u=u(t,x)\in\R$ stands for the velocity, $\rho=\rho(t,x)\in\R^{+}$ is the density, $\mu(\rho)>0$ for all $\rho>0$ is the viscosity coefficient (a smooth function in $\rho>0$) and $P(\rho)$ (a smooth function in $\rho$) is the pressure term . The Korteweg system (\ref{3systemea}) is largely used to model the time evolution of a mixture of two or more
compressible fluids with different densities.
The capillary tensor $(K[\rho])_x$  has been introduced by Korteweg in order to describe the variation of density at the interfaces between the two fluids or between two phases, generally a mixture liquid-vapor. 
In this
case the two fluids are not anymore separated by a sharp interface as it is the case for the compressible Navier-Stokes system but rather
by a thin layer where the density, although passing continuously from one
fluid to another, can have large variation. This type of models are
called \textit{Diffusive Interface Models}. The Korteweg tensor reads as:
$$K[\rho]=\rho\kappa(\rho)\rho_{xx}+\frac{1}{2}(\rho\kappa'(\rho)-\kappa(\rho))(\rho_x)^2,$$
with $\kappa(\rho)>0$ for all $\rho>0$ the capillary coefficient a smooth function in $\rho>0$.
We observe that:
$$c(K[\rho])_x=c (\mu_1(\rho)[\va_1(\rho)]_{xx})_x,$$
 with $\mu_1(\rho)=\sqrt{\rho^3\kappa(\rho)}$ and $\va_1'(\rho)=\frac{\mu_1(\rho)}{\rho^2}$.
We observe in particular that system \eqref{3systemea} satisfies the simplified formulation
$$
\begin{aligned}
&\rho V_t+\rho U\cdot V_x-(A V_x)_x+F_x=0,
\end{aligned}
$$
with 
$$
V=\begin{pmatrix}
		[\va_1(\rho)]_x\\
		u
	\end{pmatrix},\;U=\begin{pmatrix}
		u&0\\
		0&u
	\end{pmatrix},\;
A=
\begin{pmatrix}
		0&-\mu_1(\rho)\\
		c\mu_1(\rho)&
		\mu(\rho)
	\end{pmatrix}\;\;\mbox{and}\;\;F=\begin{pmatrix}
		0\\
		P(\rho)
	\end{pmatrix}.$$
	It is then natural to diagonalize the previous system in order to reformulate the system \eqref{3systemea} in an extended version, to do this we introduce the following new effective velocity
	generalizing the approach of \cite{Has}, we set
\begin{equation}
v_i=u+\frac{\mu(\rho)}{\rho^2}r_i(\rho)\rho_x,
\label{effectivevi}
\end{equation}
with $i\in\{1,2\}$. In addition we have
$$
\begin{aligned}
r_1(\rho)&=\frac{1}{2}\Big(1+\sqrt{1-4c\frac{\mu_1^2(\rho)}{\mu^2(\rho)}}\Big)\;\;\;\mbox{if}\;\;\;4c\frac{\mu_1^2(\rho)}{\mu^2(\rho)}\leq 1\\
&=\frac{1}{2}\Big(1+i\sqrt{4c\frac{\mu_1^2(\rho)}{\mu^2(\rho)}-1}\Big)\;\;\;\mbox{if}\;\;\;4c\frac{\mu_1^2(\rho)}{\mu^2(\rho)}\geq 1,
\end{aligned}
$$
$$\begin{aligned}
r_2(\rho)&=\frac{1}{2}\Big(1-\sqrt{1-4c\frac{\mu_1^2(\rho)}{\mu^2(\rho)}}\Big)\;\;\;\mbox{if}\;\;\;4c\frac{\mu_1^2(\rho)}{\mu^2(\rho)}\leq 1\\
&=\frac{1}{2}\Big(1-i\sqrt{4c\frac{\mu_1^2(\rho)}{\mu^2(\rho)}-1}\Big)\;\;\;\mbox{if}\;\;\;4c\frac{\mu_1^2(\rho)}{\mu^2(\rho)}\geq 1,
\end{aligned}
$$
which are roots of the equation  $x^2-x+c\frac{\mu_1^2(\rho)}{\mu^2(\rho)}=x(x-1)+c\frac{\mu_1^2(\rho)}{\mu^2(\rho)}=0$.
Setting $h(\rho)=\frac{\mu_1(\rho)}{\mu(\rho)}$ we can rewrite the system \eqref{3systemea} as follows
$$
\begin{cases}
\begin{aligned}
&\rho_t+(\rho u)_x=0,\\
&(\rho v_1)_t+(\rho u v_1)_x-(\mu(\rho)r_2(\rho) v_{1,x})_x+(P(\rho))_x-\frac{c}{2}\biggl(\rho_x^2\frac{\mu^2(\rho)}{\rho^2}[h^2(\rho)]'(1-4c h^2(\rho))^{-\frac{1}{2}}\biggl)_x=0,\\
&(\rho v_2)_t+(\rho u  v_2)_x-(\mu(\rho)r_1(\rho) v_{2,x})_x+(P(\rho))_x+\frac{c}{2}\biggl(\rho_x^2\frac{\mu^2(\rho)}{\rho^2}[h^2(\rho)]'(1-4c h^2(\rho))^{-\frac{1}{2}}\biggl)_x=0.
\end{aligned}
\end{cases}
$$
We can rewrite the previous system in conservative form as follows
$$
\begin{cases}
\begin{aligned}
&\rho_t+(\rho u)_x=0,\\
&(\rho v_1)_t+(\rho  v_1 v_2)_x-(\frac{\mu(\rho)}{\rho} r_2(\rho) (\rho v_{1})_x)_x+(P(\rho))_x-\frac{c}{2}\biggl(\rho_x^2\frac{\mu^2(\rho)}{\rho^2}[h^2(\rho)]'(1-4c h^2(\rho))^{-\frac{1}{2}}\biggl)_x=0,\\
&(\rho v_2)_t+(\rho v_1  v_2)_x-(\frac{\mu(\rho)}{\rho}r_1(\rho) (\rho v_{2})_x)_x+(P(\rho))_x+\frac{c}{2}\biggl(\rho_x^2\frac{\mu^2(\rho)}{\rho^2}[h^2(\rho)]'(1-4c h^2(\rho))^{-\frac{1}{2}}\biggl)_x=0.
\end{aligned}
\end{cases}
$$
We can now rewrite this system as follows with $
V=\begin{pmatrix}
		\rho\\
		\rho v_1\\
		\rho v_2
	\end{pmatrix}$,
\begin{equation}
\begin{aligned}
&V_t+A(V)V_x-(B(V) V_x)_x+F_x=0,
\end{aligned}
\label{cleimpoy}
\end{equation}
and
$$
\begin{aligned}
&A(V)=\begin{pmatrix}
		0&\theta&1-\theta\\
		-v_1v_2+P'(\rho)&v_2&v_1\\
		-v_1v_2+P'(\rho)&v_2&v_1\\
	\end{pmatrix},\;F=\begin{pmatrix}
		0\\
		-\frac{c}{2}\biggl(\rho_x^2\frac{\mu^2(\rho)}{\rho^2}[h^2(\rho)]'(1-4c h^2(\rho))^{-\frac{1}{2}}\biggl)\\
		\frac{c}{2}\biggl(\rho_x^2\frac{\mu^2(\rho)}{\rho^2}[h^2(\rho)]'(1-4c h^2(\rho))^{-\frac{1}{2}}\biggl)
	\end{pmatrix},
	\end{aligned}
	$$
	$$
	\begin{aligned}
	&B(V)=
\begin{pmatrix}
		\frac{\mu(\rho)}{\rho}(\theta r_1(\rho)+(1-\theta)r_2(\rho))&0&0\\
		0&\frac{\mu(\rho)}{\rho}r_2(\rho)&0\\
		0&0&\frac{\mu(\rho)}{\rho}r_1(\rho)
	\end{pmatrix},
	\end{aligned}$$
	and $\theta\in[0,1]$.
%\frac{1}{2}(1-\sqrt{1-4c\frac{\mu_1^2(\rho)}{\mu^2(\rho)}}).$$
%with:
%$$r_1(\rho)=\frac{1}{2}(1+\sqrt{1-4c\frac{\mu_1^2(\rho)}{\mu^2(\rho)}})\;\;\;\mbox{and}\;\;\;
%r_2(\rho)=\frac{1}{2}(1-\sqrt{1-4c\frac{\mu_1^2(\rho)}{\mu^2(\rho)}}).$$
We wish now to consider the vanishing viscosity capillary process,
we therefore replace the viscosity coefficient $\mu(\rho)$ by $\mu^\e(\rho)=\e\mu(\rho)$ and the coefficient $c$ by $c(\e)=\e^2\delta(\e)$ with $c(\e)\rightarrow_{\e\rightarrow 0}0$. 
We observe then that when $h(\rho)=\frac{\mu_1(\rho)}{\mu(\rho)}=\delta_1$ for any $\rho>0$ we can rewrite the system \eqref{cleimpoy} as follows with $
V=\begin{pmatrix}
		\rho^\e\\
		\rho^\e v_1^\e\\
		\rho^\e v_2^\e
	\end{pmatrix}$,
\begin{equation}
\begin{aligned}
&V^\e_t+A(V^\e)V^\e_x-\e(B^\e(V^\e) V^\e_x)_x=0,
\end{aligned}
%\label{cleimpoy}
\end{equation}
and
	$$
	\begin{aligned}
	&B_\e(V)=
\begin{pmatrix}
		\frac{\mu(\rho)}{\rho}(\theta r_{1,\e}(\rho)+(1-\theta)r_{2,\e}(\rho))&0&0\\
		0&\frac{\mu(\rho)}{\rho}r_{2,\e}(\rho)&0\\
		0&0&\frac{\mu(\rho)}{\rho}r_{1,\e}(\rho)
	\end{pmatrix},
	\end{aligned}$$
and
$$
\begin{aligned}
r_{1,\e}(\rho)&=\frac{1}{2}\big(1+\sqrt{1-4\delta(\e)\delta_1^2 }\big)\;\;\;\mbox{if}\;\;\;4\delta(\e)\delta_1^2\leq 1\\
&=\frac{1}{2}\big(1+i\sqrt{4\delta(\e)\delta_1^2-1}\big)\;\;\;\mbox{if}\;\;\;4\delta(\e)\delta_1^2\geq 1,
\end{aligned}
$$
$$\begin{aligned}
r_{2,\e}(\rho)&=\frac{1}{2}\big(1-\sqrt{1-4\delta(\e)\delta_1^2}\big)\;\;\;\mbox{if}\;\;\;4\delta(\e)\delta_1^2\leq 1\\
&=\frac{1}{2}\big(1-i\sqrt{4\delta(\e)\delta_1^2-1}\big)\;\;\;\mbox{if}\;\;\;4\delta(\e)\delta_1^2\geq 1.
\end{aligned}
$$
It implies that when $\sqrt{\rho^3\kappa(\rho)}=\delta_1\mu(\rho)$ and $\delta(\e)=c_1>0$ for any $\e>0$ then the Navier-Stokes Korteweg system can be rewritten in the form of a system of the type \eqref{eqn-parabolic}. In addition we consider initial data $V_{0,\e}=(\rho_{0,\e},\rho_{0,\e}v_{1,0,\e},\rho_{0,\e}v_{2,0,\e})$ with:
$$v_{1,0,\e}=u_{0}+\e\frac{\mu(\rho)}{\rho^2}r_{1,\e}(\rho_{0,\e})(\rho_{0,\e})_x\;\;\mbox{and}\;\;
v_{2,0,\e}=u_{0}+\e\frac{\mu(\rho)}{\rho^2}r_{2,\e}(\rho_{0,\e})(\rho_{0,\e})_x.$$
Furthermore we observe that $A(V)$ has the following characteristic polynomial
\[
\chi_A(\lambda)
=
-\lambda
\left(
\lambda^2-(v_1+v_2)\lambda+v_1v_2-P'(\rho)
\right).
\]
The spectrum of $A$ corresponds to
\[
\operatorname{Sp}(A)
=
\left\{
0,\,
\frac{v_1+v_2+\sqrt{(v_1-v_2)^2+4P'(\rho)}}{2},
\,
\frac{v_1+v_2-\sqrt{(v_1-v_2)^2+4P'(\rho)}}{2}
\right\}.
\]
In particular $A$ is strictly hyperbolic if $P'(\rho)>0$ and $v_1v_2\ne P'(\rho)$.
In addition when $\delta(\e)=\frac{1}{4\delta_1^2}$, we observe that $A(V)$ and $B_\e(V)=B(V)$ commutate which corresponds to our framework (indeed in this case $B_\e$ does not depend on $\e>0$). In addition we observe that in this case we have $v_1^\e=v_2^\e=v^\e$, we consider then only a $2\times2$ system.  It implies that we can apply the Corollary \ref{corollary:vv-limit} provided that $(V_{0,\e})_{\e>0}$ is an adapted family of initial data as in \cite{Chen-Kang-Vas}.
More precisely, we assume that for any $\e>0$
\begin{align}
&TV(V_{0,\e})\leq \delta_0,\;\lim_{x\rightarrow-\infty}u_{0}(x)=\bar{u}, \;\lim_{x\rightarrow-\infty}\rho_{0,\e}(x)=\bar{\rho},  \;\lim_{x\rightarrow-\infty}(\rho_{0,\e})_x(x)=0,\\
&\e TV((\rho_{0,\e})_x)\rightarrow_{\e\rightarrow 0}0,\\
&\|V_{0,\e}-V_0\|_{L^1}\rightarrow_{\e\rightarrow 0}0.
\end{align}
with $P'(\bar{\rho})\ne 0$, $\bar{\rho}>0$. Under these conditions we observe that the family of initial data $(V_{0,\e})_{\e>0}$ are living in a region where the matrix $A$ is strictly hyperbolic for $\e$ and $\delta_0$ sufficiently small. In addition from the Corollary \ref{corollary:vv-limit}, it implies that $V^\e(t)$ converges to $V(t)$ in $L^1_{loc}$ for any $t>0$ with  $V=(\rho,m)$ unique global weak solution for small $BV$ initial data of the following system% with $\theta\in[0,1]$
%\begin{equation}
%\begin{cases}
%\begin{aligned}
%&\rho_t+\theta m_{1,x}+(1-\theta)m_{2,x}=0,\\
%&m_{1,t}+(\frac{m_1m_2}{\rho})_x+(P(\rho))_x=0,\\
%&m_{2,t}+(\frac{m_1m_2}{\rho})_x+(P(\rho))_x=0,\\
%&(\rho(0,\cdot),m_{1}(0,\cdot),m_2(0,\cdot))=(\rho_0,\rho_0u_0,\rho_0u_0).
%\end{aligned}
%\end{cases}
%\label{Eulernew}
%\end{equation}
%Let us consider the unique small $BV$ global weak solution of the Euler system $(\rho',m')$ (in the sense of Bressan-De Lellis \cite{BDL}) with
\begin{equation}
\begin{cases}
\begin{aligned}
&\rho_t+m_x=0,\\
&m_{t}+(\frac{m^2}{\rho})_x+(P(\rho))_x=0,\\
&(\rho(0,\cdot), m(0,\cdot))=(\rho_0,\rho_0u_0).
\end{aligned}
\end{cases}
\label{Euler}
\end{equation}
%We observe then that $(\rho',m',m')$ is solution of \eqref{Eulernew}, and then by uniqueness we deduce that $\rho=\rho'$, $m_1=m'$, $m_2=m'$.
 In particular it proves that $(\rho^\e(t), \rho^\e v^\e(t))_{\e>0}$ converges strongly in $L^1_{loc}$ to $(\rho(t),m(t))$ for each $t>0$ with $(\rho,m)$ unique solution of the Euler system.
\section{Outline of the proof}

In this article, we focus on getting uniform total variation bound for solutions to \eqref{eqn-main}. As in \cite{BB-vv-lim-ann-math}, we can obtain such result by proceeding in two steps: (i) in Step 1 we prove the parabolic regularizing effects of solution to \eqref{eqn-main} up to a  time $\hat{t}$ and (ii) in Step 2 we establish the uniform global $BV$ bound by obtaining several interaction estimates. 

The parabolic regularizing with non-constant viscosity matrix $B(u)$ has been studied in \cite{HJ-temple-class} but up to first order. We can obtain higher order regularity in a similar way (see Proposition \ref{prop:parabolic}). For sake of completeness, we provide a proof in the Appendix \ref{appendix:higher-regularity}. For $2\times 2$ triangular system, we \cite{HJ-triangular} have established higher order regularity by using the explicit form of the system. By using a similar argument as in \cite{BB-vv-lim-ann-math}, we obtain that the solution has small total variation up to the time $\hat{t}$. 

Now we move to the most important part of this paper, that is to obtain the uniform total variation estimate for all time $t>0$. For a Temple class $A(u)$, we can decompose $u_x$ in the basis of eigenvector $\{r_{i}(u)\}_{i=1}^n$ and get the uniform total variation bound by proving the interaction estimates (see \cite{BB-temple,HJ-temple-class}). This does not work for a general $A(u)$ as we can not able to establish the $L^1_tL^1_x$ estimate of source terms (see the section 3 by Bianchini and Bressan \cite{BB-vv-lim-ann-math}). This problem can be resolved by taking the decomposition of $u_x$ with respect to the basis of viscous travelling waves. This is one of the crucial point of \cite{BB-vv-lim-ann-math}. We have mentioned in the introduction that decomposing only in the directions of viscous travelling waves does not work in our case. To explain in details let $u_x=\sum_{i}v_i\tilde{r}_i$ and $u_t=\sum_{i}(w_i-\la_i^*v_i)\tilde{r}_i$. It can be checked that $v_i$ and $w_i$ satisfy the following PDE
\begin{align*}
	\pa_tv_{i}+\pa_x(\tilde{\la}_iv_i)-\pa_x(\mu_i\pa_xv_i)&=\sum\limits_{j\neq i}(\mu_j-\mu_i)q_{ij}(t,x)(w_{j,xx}v_j-w_jv_{j,xx})+\phi_i,\\
	\pa_tw_{i}+\pa_x(\tilde{\la}_iw_i)-\pa_x(\mu_i\pa_xw_i)&=\sum\limits_{j\neq i}(\mu_j-\mu_i)\hat{q}_{ij}(t,x)(w_{j,xx}v_j-w_jv_{j,xx})+\psi_i,
\end{align*}
for some bounded functions $q_{ij}$ and $\hat{q}_{ij}$. To estimate $w_{i,xx}v_i-w_iv_{i,xx}$, we need to take derivative of the PDEs which demand a bound on the terms like $|v_{i,x}|/|v_i|$. To avoid these type of problems to occur, we choose the following basis: $\{\tilde{r}_i(u,\bar{v}_i^\e\xi_i,\si_i)\}$, where $\bar{v}_i^\e, \xi_i$ and $\si_i$ are defined as follows
\begin{equation*}
	\bar{v}_i^\e=v_i\chi^\e_i\prod\limits_{j\neq i}\eta\left(\frac{v_j^2}{v_i}\right),\, \chi^\e_i=\chi\left(\frac{v^N_i}{\e}\right)\mbox{ and }
	\si_i=\la_i^*-\theta\left(\frac{w_i}{v_i}\right),
\end{equation*} 
for some appropriately chosen cut-off functions $\eta,\chi$ and $\theta$. Note that due to the assumption $AB=BA$ and choice of these cut-off functions we have $\tilde{r}_i(u,\bar{v}_i^\e\xi_i,\si_i)=r_i(u)$ when $v_j^2>|v_i|$ for some $j\neq i$ or $|w_i|>\de_1|v_i|$ or $|v_i|^N>\e$. This choice of new basis helps us to estimate  terms like $|v_{i,x}|/|v_i|$ as in our analysis these type of term only appears when $v_j^2<|v_i|, |w_i|\leq \de_1|v_i|$ due to the cut-off functions. 

We first derive the PDEs of $v_i,w_i$ by explicit calculation in section \ref{section:phi-i} and \ref{section:psi-i}. While deriving the equations we note that the forcing part consists of the following terms: (i) $v_jv_i,\,v_{j,x}v_i,\,w_jv_i$ and $w_{j}v_{i,x}$, $(j\neq i)$ terms that appear due to interaction between different family viscous waves, (ii) $w_{i,x}v_i-w_iv_{i,x}$, it appears due to change in speed (linear), (iii) $v_i^2\left(\frac{w_i}{v_i}\right)_x^2\mathbbm{1}_{\{|w_i|/|v_i|\leq \de_1\}}$, it occurs due to change in speed (quadratic), (iv) $(w_{i,x}^2+v_{i,x}^2)\mathbbm{1}_{\{|w_i|/|v_i|> \de_1\}}$ due to wrong speed and finally (v) $w_{i,xx}v_i-w_iv_{i,xx}$ which arises due to change in speed (second order) and different eigenvalues of $B(u)$. We note that terms of the types (i) to (iv) we estimate by appropriate adaptation of estimates of \cite{BB-vv-lim-ann-math}. 

To establish the $L^1_tL^1_x$ estimate of $w_{i,xx}v_i-w_iv_{i,xx}$ we introduce two new variables $z_i,\hat{z}_i$ defined as in \eqref{intro-definition-z-i}--\eqref{intro-definition-z-hat-i}. In sections \ref{section1zi} and \ref{section1hatzi} respectively, we will deduce that $z_i$ and$ \hat{z}_i$ satisfy the following advection-diffusion equations
\begin{align*}
	z_{i,t}+(\tilde{\la}_{i}z_{i})_x-(\mu_iz_{i,x})_x=\Phi_i,\\
	\hat{z}_{i,t}+(\tilde{\la}_{i}\hat{z}_{i})_x-(\mu_i\hat{z}_{i,x})_x=\Psi_i,
\end{align*}
for some forcing terms $\Phi_i$ and $\Psi_i$. We now focus on deriving the expression of $\Phi_i$ and $\Psi_i$. Note that it contains the derivative of forcing parts of equations of $v_i$ and $w_i$, that is,  $\phi_{i,x}$ and $\psi_{i,x}$ (see \eqref{eqfonda1} and \eqref{eqfonda2} in sections \ref{section1zi} and \ref{section1hatzi} respectively). We carefully analyze $\phi_{i,x}$ and $\psi_{i,x}$ in section \ref{sectionphii}. We have discussed in details how to estimate the source part. One of the crucial terms that appears in $\phi_{i,x}$ is the following 
\begin{equation*}
	(v_{j,x}-\mu_j^{-1}(\tilde{\la}_i-\la_j^*+\theta_j)v_j)_{xx}(B-\mu_jI_n)[\tilde{r}_{j}+v_j\xi_j\pa_{v_j}\bar{v}_j\tilde{r}_{j,v}].
\end{equation*}
To deal with this type of terms we derive the second order derivative of coupling equation \eqref{eqn-v-i-x-1} which is a complex term to deal with. In section \ref{sectionphii} we discuss this with detailed explanation. We use Lemma \ref{lemme9.11} to prove that we can bound $v_j	(v_{j,x}-\mu_j^{-1}(\tilde{\la}_i-\la_j^*+\theta_j)v_j)_{xx}$ by the terms of types (i) to (v) and additionally terms consisting third order derivatives as in \eqref{eqn-3rd-order-error}.

%To obtain the estimates of $(v_{j,x}-\mu_j^{-1}(\tilde{\la}_i-\la_j^*+\theta_j)v_j)_{xx}$, we need to derive estimates of derivatives of $\tilde{r}_i, \psi_{ik}(=r_k(u)\cdot \tilde{r}_i)$. See Lemma \ref{lemme9.11} in section \ref{section:directional}. Due to the choices of our basis $\tilde{r}_i$ we can bound  

Furthermore, we note that $\phi_{i,x}$ and $\psi_{i,x}$ can be estimated in terms of the interaction terms as before and additionally, we have $v_{j,xx}v_i, w_{j,xx}v_i$ with $j\neq i$ and $v_{j,xxx}v_jv_i$ to bound. By using transversal estimate corresponding to $z_j$ and $v_i$ we can have $L^1_tL^1_x$ bounds of $v_{j,xx}v_i$. For the therm involving third order derivative of $v_j$ or $w_j$ we use a bootstrap argument along with the regularity estimate of $v_{i,tx}$ and $w_{i,tx}$ (see section \ref{section:Lambda-2}).

\section{Parabolic estimates}\label{section:parabolic}
In this section, we prove that for initial data with small total variation up to a small time $\hat{t}>0$, the solution $u$ to \eqref{eqn-main} gain regularity from the equation due to the parabolic assumptions ($\mathcal{H}_B$) on $B(u)$. Before, we jump into the main regularizing results of the section, let us first state the following elementary lemma regarding the change of variable.
\begin{lemma}\label{lemma:transformation}
	Let $n\geq 1, p\geq1$, $u\in C^4(\R,\R^n)$ and $d_i:\R^n\rr[c_0,\f),1\leq i\leq n$ are $C^4$ functions for some $c_0>0$. Consider $\mathcal{T}:L^p(\R,\R^n)\rr L^p(\R,\R^n)$ defined as follows
	\begin{equation}
		\mathcal{T}(f)_i(x)=f_i(X_i(x))\mbox{ where }X_i(x)=\int\limits_{0}^{x}\frac{1}{\sqrt{d_i(u(z))}}\,dz.
	\end{equation}
	Then $\mathcal{T}$ is well-defined and it satisfies
	\begin{equation}\label{inequality-T-1}
		\norm{\mathcal{T}(f)}_{L^p}\leq m^{\frac{1}{2p}}\norm{f}_{L^p}\leq c_0^{-\frac{1}{2p}}m^\frac{1}{2p}\norm{\mathcal{T}(f)}_{L^p}\mbox{ where }m=\sup\limits_{1\leq i\leq n}\norm{d_i(u(\cdot))}_{L^\f}.
	\end{equation} 
	In addition $\mathcal{T}$ is invertible with and $ \mathcal{T}^{-1}$ is defined as follows for any $f\in L^p(\R,\R^n)$
	\begin{equation}
		\mathcal{T}^{-1}(f)_i(x)=f_i(X_i^{-1}(x)),
	\end{equation}
	where  for any $i\in\{1,\cdots,n\}$, $X_i^{-1}$ is the inverse of the function $X_i$ with $(X_i^{-1})'(X_i(x))=\sqrt{d_i(u(x))}$. We have in particular for any $f\in L^p(\R,\R^n)$
	\begin{equation}\label{inequality-T-1bis}
		\norm{\mathcal{T}^{-1}(f)}_{L^p}\leq c_0^{-\frac{1}{2p}}\norm{f}_{L^p}.
	\end{equation} 
	Furthermore, for $f\in C^2$ we have
	\begin{align}
		&\norm{(\mathcal{T}(f))^\p}_{L^p}\leq c_0^{-\frac{p-1}{2p}} \norm{f^\p}_{L^p}\leq c_0^{-\frac{p-1}{2p}}m^\frac{p-1}{2p}\norm{\mathcal{T}(f)^\p}_{L^p},	\label{inequality-T-2}\\
		&\norm{(\mathcal{T}(f))^{\p\p}}_{L^p}\leq 2c_0^{-\frac{2p-1}{2p}}\norm{f^{\p\p}}_{L^p}+ m_1c_0^{-\frac{3p-1}{2p}}\norm{u_x}_{L^\f}\norm{f^\p}_{L^p},	\label{inequality-T-3a}\\
		& \norm{(\mathcal{T}(f))^{\p\p\p}}_{L^p}\leq 5^{1-1/p}\frac{3}{2}  m_1 \norm{u_x}_{L^\infty} c_o^{-\frac{4p-1}{2p}}\|f''\|_{L^p}+5^{1-1/p}(c_0)^{-\frac{3p-1}{2p}}\|f'''\|_{L^p}
		\nonumber \\
		&+5^{1-1/p}\frac{3}{4} \norm{u_x}_{L^\infty}^{2}m_1^{2}(c_0)^{-\frac{5p-1}{2p}}\norm{f'}_{L^p}  +\frac{5^{1-1/p}}{2}\norm{u_{xx}}_{L^\infty} m_1 (c_0)^{-\frac{3p-1}{2p}} \norm{f'}_{L^p} \nonumber\\
		&+ \frac{5^{1-1/p}}{2}\norm{u_{x}}_{L^\infty}^{2} m_1 (c_o)^{-\frac{3p-1}{2p}} \norm{f'}_{L^p}\label{4.7}\\
		&\norm{f^{\p\p}}_{L^p}\leq 2m^{\frac{2p-1}{2p}}\norm{(\mathcal{T}(f))^{\p\p}}_{L^p}+ m_1c_0^{-\frac{1}{2}}\norm{u_x}_{L^\f}\norm{f^\p}_{L^p},\label{inequality-T-3b}\\
		&\norm{f'''}_{L^p}
		\leq 5^{1-1/p}m^{\frac{3p-1}{2p}} \norm{(\mathcal{T}(f))'''}+5^{1-1/p}\frac{3}{2}c_0^{-\frac{1}{2}}m_1\|u_x\|_{L^\infty} \norm{ f'' }_{L^p} \\
		&+5^{1-1/p}\norm{ f'}_{L^p}\big(\frac{3}{4} c_0^{-1}\|u_x\|_{L^\infty}^{2}m_1^{2}+\frac{1}{2} \|u_{xx}\|_{L^\infty}m_1
		+\frac{1}{2} \|u_{x}\|_{L^\infty}^{2}m_1
		\big).\label{4.7bis}
	\end{align}
	where $m_1:=max(\sup\limits_{1\leq i\leq n}\norm{Dd_i(u(\cdot))}_{L^\f}, \sup\limits_{1\leq i\leq n}\norm{D^2d_i(u(\cdot))}_{L^\f})$. 
	%We have in addition:
\end{lemma}
\begin{remark}
	In the sequel we will assume that for any $i\in\{1,\cdots,n\}$, $d_i=\mu_i$ with $\mu_i$ defined on $\mathcal{U}$. 
\end{remark}
We search now to construct  $BV$ solution $u(t,x)$ of \eqref{eqn-main} satisfying
$$u^*=\lim_{x\rightarrow-\infty}u(t,x),$$
with  $u^*\in \mathcal{K}$ clearly independent of time and with $\mathcal{K}$ defined in the Theorem \ref{theorem:BV-estimate}.
We consider $G$ as the fundamental solution of the following parabolic equation
\begin{equation}
	w_t+A_2^*w_x=w_{xx},
\end{equation}
where $A_2^*=A_1 B_1^{-1/2}(u^*)$. Furthermore $A_1(u)$ and $B_1(u)$ are defined as follows $A_1(u)=\mbox{diag}(\lambda_1(u),\cdots,\lambda_n(u))=P(u)A(u)P(u)^{-1}$ and $B_1(u)=\mbox{diag}(\mu_1(u),\cdots,\mu_n(u))$ with $B_1(u)=P(u)B(u)P(u)^{-1}$. In the sequel it will be convenient to use "$\bullet$" to denote a directional derivative, so that $z\bullet A(u)=DA(u)\cdot z$. We will also use the notation $z\cdot DA(u)$ to define the differential of $A$ in $u$ applied to the vector $z$.\\
The function $G$ satisfies the following estimates
\begin{equation}\label{def:kappa}
	\norm{G(t,\cdot)}_{L^1}\leq \kappa_G,\quad	\norm{G_x(t,\cdot)}_{L^1}\leq \frac{\kappa_G}{\sqrt{t}},
\end{equation}
for some constant $\kappa_G>0$.
To have the parabolic estimates for \eqref{eqn-main}, we define the following constants with $\bar{\delta}>0$ sufficiently small such that $B(u^*,\bar{\delta})$ is included in $\mathcal{U}$,
\begin{align}
	&\kappa_1=2\max\{c_0^{-\frac{1}{2}},m^\frac{1}{2},m_1\},\nonumber\\
	&\kappa_A=\sup\limits_{\abs{u-u^*}\leq \bar{\de}}\max\limits_{k=1,2,3,4}\left\{\abs{A(u)},\abs{D^kA(u)}, |D^k (A_1^* B_1^{-1/2})(u)|\right\},\nonumber\\
	&\kappa_{B_1}=\sup\limits_{\abs{u-u^*}\leq \bar{\de}}\max\limits_{k=1,2,3,4, i=1,2,3} \left\{\abs{B_1(u)},\abs{B_1^{-\frac{i}{2}}(u)},\abs{D^{k} B_1(u)}, \abs{D^{k} B^{-\frac{i}{2}}_1(u)} \right\},\nonumber\\
	&\kappa_B=\sup\limits_{\abs{u-u^*}\leq \bar{\de}}\max\limits_{k=1,2,3,4,1\leq i\leq n}\left\{\abs{B(u)},\abs{B^{-1}(u)}, \abs{D^{k} B(u)},(\mu_i+\mu_i^{-1}),\abs{D^k\mu_i}\right\},\nonumber\\
	&\kappa_P=\sup\limits_{\abs{u-u^*}\leq \bar{\de}}\max\limits_{k=1,2,3,4}\left\{\abs{P(u)},\abs{P^{-1}(u)},\abs{D^{k} P(u)},\abs{D^{k} P^{-1}(u)}\right\}. \label{def:kappa-P}
\end{align}
where $m,m_1$ are as in Lemma \ref{lemma:transformation} and $c_0=c_1$ as in ($\mathcal{H}_B)$. We have in particular
\begin{align*}
	&m=\sup\limits_{1\leq i\leq n}\norm{\mu_i(u(\cdot))}_{L^\f(B(u^*,\bar{\delta}))}, \\
	&m_1=\max(\sup\limits_{1\leq i\leq n}\norm{D\mu_i(u(\cdot))}_{L^\f(B(u^*,\bar{\delta}))}, \sup\limits_{1\leq i\leq n}\norm{D^2\mu_i(u(\cdot))}_{L^\f(B(u^*,\bar{\delta}))}).
\end{align*}
We set $\kappa:=\max\{\kappa_A,\kappa_{B_1},\kappa_B,\kappa_G,\kappa_P,\kappa_1\}$. Now, we state the main parabolic regularizing effect result for $u$ provided that $u_x$ is small in $L^1$ up to a time $\hat{t}$.
\begin{proposition}\label{prop:parabolic}
	Let $u$ be a solution to the equation \eqref{eqn-main} satisfying 
	\begin{equation}\label{assumption:u-L1}
		\norm{u_x(t,\cdot)}_{L^1}\leq \de_0\mbox{ for all }t\in[0,\hat{t}]\mbox{ where }\hat{t}:=\left(\frac{1}{2^{55}\kappa^{42}\de_0}\right)^2,
	\end{equation}
	for some $0<\de_0<\frac{1}{2^{56}\kappa^{42}}$ and $\kappa$ is defined as above. Then we have for any $t\in]0,\hat{t}]$
	\begin{equation}\label{estimate:parabolic-1}
		\norm{\pa_x^{(k)}u_{x}(t,\cdot)}_{L^1}\leq \frac{2^{a_k}\kappa^{b_k}\de_0}{t^\frac{k}{2}}\mbox{ for }k=1,2,3,4,
	\end{equation}
	the $a_k$ sequence is defined by the recurrence relation $a_k=a_{k-1}+1+\floor*{k/2}$ with $a_0=0$. In addition we have $b_1=4$, $b_2=6$, $b_3=10$, $b_4=10$. Furthermore, we have for any $t\in]0,\hat{t}]$
		\begin{equation}\label{estimate:parabolic-2}
			\norm{u_{xxxxx}(t,\cdot)}_{L^\f}\leq \frac{2^{11}\kappa^{11}\de_0}{t^{5/2}}.
	\end{equation}
\end{proposition}
Note that $a_1=1,a_2=3,a_3=5$ and $a_4=8$. Proof of Proposition \ref{prop:parabolic} is similar to the \cite[Proposition 3.2]{HJ-temple-class}, but here we need higher order estimates. For completeness we decide to give a rigorous proof in Appendix \ref{appendix:proof-prarabolic}.

Next we show that the assumption \eqref{assumption:u-L1}, as in Proposition \ref{prop:parabolic} holds true for sufficiently small initial data. 
\begin{proposition}\label{prop:local-existence}
	Let $u$ be solution of \eqref{eqn-main} such that $TV(\bar{u})\leq \frac{\de_0}{4\kappa}$.  Then $u$ is well-defined on $[0,\hat{t}]$ where $\hat{t}$ is defined as in \eqref{assumption:u-L1}. Moreover, we have
	\begin{equation}\label{local-estimate-u-h}
		\norm{u_x(t)}_{L^1}\leq \frac{\de_0}{2}\mbox{ for }t\in[0,\hat{t}].
	\end{equation}
\end{proposition}
Proof of the proposition follows from estimates \eqref{estimate:parabolic-1} with the argument as in \cite[Proposition 2.3]{BB-vv-lim-ann-math}, \cite[Proposition 3.6]{HJ-temple-class}. We omit the proof here.
\begin{corollary}\label{coro2.2}
	Let $T>\hat{t}$. Assume that the solution $u$ of  \eqref{eqn-main} satisfies on $[0,T]$
	$$\|u_{x}(t,\cdot)\|_{L^1}\leq \delta_0,$$
	then for all $t\in[\hat{t},T]$ we have
	\begin{equation}
		\norm{\pa_x^{k}u_{x}(t,\cdot)}_{L^1}=O(1)\de_0^{k+1}\mbox{ for }k=1,2,3,4.
	\end{equation}
 Furthermore, we have for all  $t\in[\hat{t},T]$ 
		\begin{equation}\label{estimate:parabolic-2a}
			\norm{u_{xxxxx}(t,\cdot)}_{L^\f}=O(1)\delta_0^6.
	\end{equation}
\end{corollary}
\begin{proof} It suffices to apply Proposition \eqref{prop:parabolic} on $[t-\hat{t},t]$.
\end{proof}

%\bibitem{BB-vv-lim-ann-math}

\section{Viscous travelling wave}
For a Temple class system $A(u)$ we can prove the uniform total variation estimate by decomposing $u_x$ in the basis of $\{r_i(u)\}_{i=1}^n$. Without Temple class assumption, that is for a general $A(u)$ satisfying \descref{A1}{($\mathcal{H}_A$.)}, we need to decompose $u_x$ in the basis of travelling waves. We refer to \cite{BB-vv-lim-ann-math} for the explanation why it is necessary. In this section, we construct the basis $\{\tilde{r}_i(u,v_i,\si_i)\}$ of travelling waves. Although the construction of $\tilde{r}_i$ is similar to the one in Bianchini and Bressan \cite{BB-vv-lim-ann-math} we deduce a few more new identities and estimates which will be crucially used in our article.

Travelling wave $u(t,x)=U(x-\sigma t)$ of \eqref{eqn-main} satisfy the following equation
\begin{equation*}
	(A(U)-\si \mbox{I}_n)U^\p=(B(U)U^\p)^\p.
\end{equation*}
We can write
\begin{equation}
	U^{\p\p}=B^{-1}(U)(A(U)-\si \mbox{I}_n)U^\p-B^{-1}(U)(U^\p\cdot DB(U))U^\p.
	\label{equaU}
\end{equation}
We consider the following system of ODE
\begin{equation}\label{ODE-1}
	\begin{cases}
		\dot{u}&=v,\\
		\dot{v}&=B^{-1}(u)(A(u)-\si \mbox{I}_n)v-B^{-1}(u)(v\cdot DB(u))v,\\
		\dot{\si}&=0,
	\end{cases}
\end{equation}
where we define $F$ as follows:
$$F(u,v,\si)=\begin{pmatrix}
		v\\
		B^{-1}(u)(A(u)-\si \mbox{I}_n)v-B^{-1}(u)(v\cdot DB(u))v\\
		0
		\end{pmatrix}.$$
We note that $P^*_i:=(u^*,0,\la_{i}(u^*))$ are equilibrium points for $1\leq i\leq n$. We linearize near the point $P^*_i$ and get
\begin{equation}
	\begin{cases}
		\dot{u}&=v,\\
		\dot{v}&=B^{-1}(u^*)(A(u^*)-\la_i(u^*)\mbox{I}_n)v,\\
		\dot{\si}&=0.
	\end{cases}
\end{equation} 
Let $\{r_i(u)\}_{1\leq i\leq n}$ and $\{l_i(u)\}_{1\leq i\leq n}$ be the sets of right and left eigenvector of $A(u)$ of norm one. We denote $r_i^*=r_i(u^*),l_i^*=l_i(u^*)$. We define $V_i,1\leq i\leq n$ as follows
\begin{equation}
	v=\sum\limits_{j}V_jr_j^*,\quad V_j:=l_j^*\cdot v.
\end{equation}
Then we have
\begin{equation}
	Z=DF(u^*,0,\la_i(u^*))=\begin{pmatrix}
		O_n&I_n& 0\\
		O_n&B^{-1}(u^*)(A(u^*)-\la_i(u^*)I_n)&0\\
		0&0&0
	\end{pmatrix}.
\end{equation}
Subsequently, we get
\begin{equation}
	Z^2=\begin{pmatrix}
		O_n&B^{-1}(u^*)(A(u^*)-\la_i(u^*)I_n)& 0\\
		O_n&(B^{-1}(u^*)(A(u^*)-\la_i(u^*)I_n))^2&0\\
		0&0&0
	\end{pmatrix}.
\end{equation}
Now since $B(u^*)$ and $A(u^*)$ commutate, we deduce that\begin{equation}
	Z^2=\begin{pmatrix}
		O_n&B^{-1}(u^*)(A(u^*)-\la_i(u^*)I_n)& 0\\
		O_n&(B^{-1}(u^*))^2 (A(u^*)-\la_i(u^*)I_n)^2&0\\
		0&0&0
	\end{pmatrix}.
\end{equation}
Therefore the center subspace will look like
\begin{equation}
	\mathcal{N}_i:=\left\{(u,v,\si)\in\R^n\times\R^n\times\R;\,V_j=0,j\neq i\right\}.
\end{equation}
Note that $\mbox{dim}(\mathcal{N}_i)=n+2$, by Center Manifold Theorem \cite{Vanderbauwhede}, there exists a smooth manifold $\mathcal{M}_i\subset \R^{2n+1}$ which is tangent to $\mathcal{N}_i$ at $P^*_i$. Furthermore, $\mathcal{M}_i$ has dimension $n+2$ and is locally invariant under the flow of \eqref{ODE-1}. We can write
\begin{equation}
	V_j=\varphi_j(u,V_i,\si),\quad j\neq i.
\end{equation}
We further consider 
\begin{equation}
	\mathcal{D}_i:=\left\{\abs{u-u^*}<\e,\,\abs{V_i}<\e,\,\abs{\si-\la_i(u^*)}<\e\right\}.
\end{equation} 
Since $\mathcal{M}_i$ is tangent to $\mathcal{N}_i$ we have
\begin{equation}\label{quadratic-varphi-1}
	\varphi_j(u,V_i,\si)=O(1)\left(\abs{u-u^*}^2+\abs{V_i}^2+\abs{\si-\la_i(u^*)}^2\right).
\end{equation}
Note that equilibrium points $(u,0,\si)$ with $\abs{u-u^*}<\e,\abs{\si-\la_i(u^*)}<\e$ lie in $\mathcal{M}_i$ we have
\begin{equation}\label{equilibrium-1}
	\varphi_j(u,0,\si)=0\mbox{ for all }j\neq i.
\end{equation}
Hence, we may write
\begin{equation}
	\varphi_j(u,V_i,\si)=\psi_j(u,V_i,\si) V_i,
\end{equation}
for some $\psi_j$.  Now, we would like to make a change of coordinates $V_i\mapsto \tilde{V}_i$ as follows
\begin{equation}
	\tilde{V}_i=l_i(u)\cdot v=V_i l_i(u)\cdot (r_i^*+\sum_{j\ne i}\psi_j(u,V_i,\sigma)r_j^*)
\end{equation}
It implies that $\tilde{V}_i=\zeta_{ii}(u,V_i,\si)V_i$ with $\zeta_{ii}(u,V_i,\si)=l_i(u)\cdot (r_i^*+\sum_{j\ne i}\psi_j(u,V_i,\sigma)r_j^*)$
a $C^2$ function such that $\zeta_{ii}(u^*,0,\si)=1$. If we consider the function $f_{(u,\sigma)}(V_i)=\zeta_{ii}(u,V_i,\si)V_i$, we observe that $f_{(u,\sigma)}'(V_i)=\zeta_{ii}(u,V_i,\si)+V_i\zeta_{ii,V_i}(u,V_i,\si)\ne 0$ in ${\cal D}_i$. It implies that $f_{(u,\sigma)}$ is locally invertible and we can make the change of coordinates  $V_i\mapsto \tilde{V}_i$.
%For $k\neq i$ we have  $\tilde{V}_k=\tilde{\phi}_k(u,\tilde{V}_i,\si)$ for some smooth function $\tilde{\phi}$. From \eqref{equilibrium-1} we deduce that $\tilde{V}_k=\tilde{\psi}_k(u,\tilde{V}_i,\si)\tilde{V}_i$. 
Therefore, for any point $(u,v,\si)\in\mathcal{M}_i$ we can write
\begin{equation}
	\begin{aligned}
		v&=\sum\limits_{k}l_k(u)\cdot v\, r_k(u)=\sum\limits_{k}V_i l_k(u)\cdot (r_i^*+\sum_{j\ne i}\psi_{j}(u,V_i,\sigma)r_{j}^*)r_{k}(u)\\
		&=\tilde{V}_ir_i(u)+\sum\limits_{k\ne i}\tilde{V}_i \frac{1}{\xi_{ii}(u,V_i,\sigma)}l_k(u)\cdot (r_i^*+\sum_{j\ne i}\psi_{j}(u,V_i,\sigma)r_{j}^*)\,r_{k}(u)\\
		&=\tilde{V}_i\left(r_i(u)+\sum\limits_{k\neq i}\tilde{\psi}_k(u,\tilde{V}_i,\si)r_k(u)\right)=\tilde{V}_i \tilde{r_i} (u,\tilde{V}_i,\si),
	\end{aligned}
\end{equation}
with $\tilde{\psi}_k(u,\tilde{V}_i,\si)=\frac{1}{\xi_{ii}(u,V_i,\sigma)}l_k(u)\cdot (r_i^*+\sum_{j\ne i}\psi_{j}(u,V_i,\sigma)r_{j}^*)$ and $\tilde{r_i} (u,\tilde{V}_i,\si)=r_i(u)+\sum\limits_{k\neq i}\tilde{\psi}_k(u,\tilde{V}_i,\si)r_k(u)$.
%Now, we would like to make a change of coordinates $V_k\mapsto \tilde{V}_k$ as follows
%\begin{equation}
%	\tilde{V}_k=v\cdot r_k(u). 
%\end{equation}
%We observe that $\tilde{V}_i=\zeta_{ii}(u,V_i,\si)V_i$ for some $C^2$ function $\zeta_{ii}$ with $\zeta_{ii}(u^*,0,\si)=1$. For $k\neq i$ we have  $\tilde{V}_k=\tilde{\phi}_k(u,\tilde{V}_i,\si)$ for some smooth function $\tilde{\phi}$. From \eqref{equilibrium-1} we deduce that $\tilde{V}_k=\tilde{\psi}_k(u,\tilde{V}_i,\si)\tilde{V}_i$. Therefore, for any point $(u,v,\si)\in\mathcal{M}_i$ we can write
%\begin{equation}
%	v=\sum\limits_{k}\tilde{V}_kr_k(u)=\tilde{V}_i\left(r_i(u)+\sum\limits_{j\neq i}\tilde{\psi}_j(u,\tilde{V}_i,\si)r_j(u)\right)=:\tilde{V}_i\tilde{r}_i(u,\tilde{V}_i,\si).
%\end{equation}
From \eqref{quadratic-varphi-1} it follows that
\begin{equation}
	\tilde{r}_i(u,\tilde{V}_i,\si)\rr r_i^*\mbox{ as }(u,\tilde{V}_i,\si)\rr(u^*,0,\la_i(u^*)).
\end{equation}

Therefore, we may write $	v=\tilde{V}_i \tilde{r}_i$,
\begin{equation}
	\mathcal{M}_i=\left\{(u,v,\si_i);\,\,v=\tilde{V}_i \tilde{r}_i(u,\tilde{V}_i,\si_i)\right\}\mbox{ for }1\leq i\leq n,
\end{equation}
provided $(u,v_i,\si_i)\in\R^n\times\R^n\times\R$ are in a sufficiently small neighbourhood of $(u^*,0,\la_i(u^*))$. Next we derive few more estimates on $\tilde{\psi}_k,\,k\neq i$. 
\begin{claim}\label{claim-3.1}
	In a neighbourhood of $(u^*,0,\la_i^*)$ we have
	\begin{equation}
		\tilde{\psi}_{k}(u,\tilde{V}_i,\si)=O(1)\abs{\tilde{V}_i}.
		\label{superimpo}
	\end{equation}
	Furthermore, we have
	\begin{equation}
		\tilde{\psi}_{k,u}(u,\tilde{V}_i,\si),\tilde{\psi}_{k,\si}(u,\tilde{V}_i,\si),\tilde{\psi}_{k,u\si}(u,\tilde{V}_i,\si),\tilde{\psi}_{k,\si\si}(u,\tilde{V}_i,\si), \tilde{\psi}_{k,u\si\si}(u,\tilde{V}_i,\si)=O(1)\abs{\tilde{V}_i}.
		\label{5.19}
	\end{equation}
\end{claim}
\begin{proof}[Proof of Claim \ref{claim-3.1}:]
Let us consider a travelling wave $u(t,x)=U(x-\sigma t)$ of \eqref{eqn-main} contained in the center manifold $\mathcal{M}_i$ such that:
$$U'=\tilde{V}_i \tilde{r}_i(u,\tilde{V}_i,\si).$$
	We observe that
	\begin{equation}
		U^{\p\p}=\tilde{V}_i^{\p}[\tilde{r}_i+\tilde{V}_i\tilde{r}_{i,v}]+\tilde{V}^2_i\tilde{r}_{i,u}\tilde{r}_i.
	\end{equation}
	Then we have using \eqref{equaU}
	\begin{equation}
		\tilde{V}_i^{\p}B(u)[\tilde{r}_i+\tilde{V}_i\tilde{r}_{i,v}]+\tilde{V}^2_iB\tilde{r}_{i,u}\tilde{r}_i+\tilde{V}_i^2\tilde{r}_i\cdot DB(u)\tilde{r}_i=\tilde{V}_i(A(u)-\si)\tilde{r}_i.
		\label{3.19}
	\end{equation}
	We deduce then that
	\begin{align}
		\tilde{V}_i^{\p}=\frac{1}{\langle B(u)[\tilde{r}_i+\tilde{V}_i\tilde{r}_{i,v}],\tilde{r}_i \rangle }(\tilde{\la}_i-\si\|\tilde{r}_i\|^2)\tilde{V}_i+O(1)\abs{\tilde{V}_i}^2,
	\end{align}
	with $\tilde{\lambda}_i=\langle A(U)\tilde{r}_i\,\tilde{r}_i\rangle$. It gives then
	\begin{equation}
		\frac{1}{\langle B(u)[\tilde{r}_i+\tilde{V}_i\tilde{r}_{i,v}],\tilde{r}_i \rangle }(\tilde{\la}_i-\si\|\tilde{r}_i\|^2)\tilde{V}_iB(u)[\tilde{r}_i+\tilde{V}_i\tilde{r}_{i,v}]+O(1)\abs{\tilde{V}_i}^2=\tilde{V}_i(A(u)-\si)\tilde{r}_i.
	\end{equation}
	Dividing by $\tilde{V}_i$ and considering the limit $\tilde{V}_i\rr0$, we obtain
	\begin{equation}\label{eqn:vi=0}
		\frac{(\tilde{\la}_i(u,0,\si)-\si\|\tilde{r_i}(u,0,\sigma)\|^2)}{\langle B(u)\tilde{r}_i(u,0,\si),\tilde{r}_i(u,0,\si) \rangle }B(u)\tilde{r}_i(u,0,\si)= (A(u)-\si)\tilde{r}_i(u,0,\si).
	\end{equation}
	Set $b(u,0,\si)=\langle B(u)\tilde{r}_i(u,0,\si),\tilde{r}_i(u,0,\si) \rangle>0$ (indeed we use the fact that $\tilde{r}_i(u,0,\si)$ is close from $r_i^*$ and the assumption  $\mathcal{H}_B$ on $B(u)$) we get
	\begin{align*}
		B(u)\tilde{r}_i(u,0,\si)&=\left(Br_i(u)+\sum\limits_{j\neq i}\tilde{\psi}_j(u,0,\si)Br_j(u)\right)\\
		&=\left(\alpha_i(u)r_i(u)+\sum\limits_{j\neq i}\tilde{\psi}_j(u,0,\si)\alpha_j(u)r_j(u)\right),\\
		A(u)\tilde{r}_i(u,0,\si)&=\left(Ar_i(u)+\sum\limits_{j\neq i}\tilde{\psi}_j(u,0,\si)Ar_j(u)\right)\\
		&=\left(\la_i(u)r_i(u)+\sum\limits_{j\neq i}\tilde{\psi}_j(u,0,\si)\la_j(u)r_j(u)\right).
	\end{align*}
	Applying above identities on \eqref{eqn:vi=0} we obtain
	\begin{align*}
		&\frac{(\tilde{\la}_i(u,0,\si)-\si\|\tilde{r_i}(u,0,\sigma)\|^2)}{b(u,0,\si)}\alpha_i(u)r_i(u)+	\frac{(\tilde{\la}_i(u,0,\si)-\si\|\tilde{r_i}(u,0,\sigma)\|^2)}{b(u,0,\si)}\\
		&\times\sum\limits_{j\neq i}\tilde{\psi}_j(u,0,\si)\alpha_j(u)r_j(u)=(\la_i(u)-\si)r_i(u)+\sum\limits_{j\neq i}\tilde{\psi}_j(u,0,\si)(\la_j(u)-\si)r_j(u).
	\end{align*}
	Therefore, we have 
	\begin{align*}
		\frac{(\tilde{\la}_i(u,0,\si)-\si\|\tilde{r_i}(u,0,\sigma)\|^2)}{b(u,0,\si)}\alpha_i(u)&=\la_i(u)-\si,\\
		\frac{(\tilde{\la}_i(u,0,\si)-\si\|\tilde{r_i}(u,0,\sigma)\|^2)}{b(u,0,\si)}\alpha_j(u)&=\la_j(u)-\si\,\,\mbox{ if }\,\,\tilde{\psi}_j(u,0,\si)\neq0.
	\end{align*}
	This implies  since $\tilde{\lambda}_i(u,0,\si)=\langle A(u)\tilde{r}_i(u,0,\si), \tilde{r}_i(u,0,\si)\rangle$ and using the fact that $\tilde{r}_i(u,0,\si)$ and $\sigma$ are respectively close from $r_i^*$ and $\lambda_i^*$ we can find $\e_0$ small enough such  that  we have using the assumption $\mathcal{H}_A$,  $\mathcal{H}_B$ respectively on $A(u)$ and $B(u)$
	\begin{equation*}
		c_0\leq \abs{\la_i(u)-\la_j(u)}=\frac{\abs{\tilde{\la}_i(u,0,\si)-\si \|\tilde{r_i}(u,0,\sigma)\|^2}}{b(u,0,\si)}\abs{\alpha_i(u)-\alpha_j(u)}\leq \e_0.
	\end{equation*}
	This is a contradiction. Hence, $\tilde{\psi}_j(u,0,\si)=0$ for all $j\neq i$.
\end{proof}
Then we may write 
\begin{equation*}
	\tilde{\psi}_{j}(u,v_i,\si_i)=v_i \phi_j(u,v_i,\si_i)\mbox{ for }j\neq i.
\end{equation*}
We would like to write $\tilde{r}_{i}(u,v_i,\si_i)$ as follows
\begin{equation}
	\tilde{r}_i(u,v_i,\si_i)=r_i(u)+\sum\limits_{j\neq i}\psi_{ij}(u,v_i,\si_i)r_j(u).
	\label{5.36}
\end{equation}
\begin{remark}
We observe in particular using \eqref{superimpo} that for any $i\in\{1,\cdots,n\}$
\begin{equation}
\tilde{r}_i(u,0,\si_i)=r_i(u).
\label{observationimp}
\end{equation}
It implies in particular that
\begin{equation}
\begin{aligned}
&\tilde{r}_i(u,v_i,\si_i)-r_i(u)=O(1)v_i,\;\tilde{r}_{i,\si}(u,v_i,\si_i)=O(1)v_i,\\
&\tilde{r}_{i,u\sig}(u,v_i,\si_i), \tilde{r}_{i,u\sig \sig}(u,v_i,\si_i)=O(1)v_i,\;\;\;\;\;\;\;\;\,\tilde{r}_{i,\si \si}(u,v_i,\si_i)=O(1)v_i.
\end{aligned}
\label{observationimp1}
\end{equation}
\end{remark}
On the manifold $\mathcal{M}_i$ we have as in \eqref{3.19}
\begin{equation}
	v_{i,x}B(u)\tilde{r}_i+v_{i}B(u)\tilde{r}_{i,x}+v_i^2\tilde{r}_{i}\cdot DB(u)\tilde{r}_i=v_{i}(A(u)-\si_i)\tilde{r}_i.
\end{equation}
Then we have
\begin{equation}\label{identity-1}
	v_{i,x}B[\tilde{r}_i+v_i\tilde{r}_{i,v}]+v^2_iB\tilde{r}_{i,u}\tilde{r}_i+v_i^2\tilde{r}_i\cdot DB(u)\tilde{r}_i=v_i(A(u)-\si_i)\tilde{r}_i.
\end{equation}
First, observe that taking the product with $l_i(u)$ we have then (indeed we observe that $\langle B (u)\tilde{r}_i,l_i\rangle= \mu_i(u)$ and $\langle B\tilde{r}_{i,v},l_i\rangle=0$)
$$v_{i,x}\mu_i(u)=v_i(\lambda_i(u)-\sigma_i-v_i\langle B\tilde{r}_{i,u}\tilde{r}_i+\tilde{r}_i\cdot DB(u)\tilde{r}_i ,l_i(u)\rangle).$$
We set 
\begin{equation}
\tilde{\la}_i:= \la_i-v_i
\langle B\tilde{r}_{i,u}\tilde{r}_i+\tilde{r}_i\cdot DB(u)\tilde{r}_i ,l_i(u)\rangle.
\label{deflatilde}
\end{equation}
 Therefore, we have%v_i\langle B\tilde{r}_{i,u}\tilde{r}_i+\tilde{r}_i\cdot DB(u)\tilde{r}_i,r_i\rangle$. Therefore, we have
\begin{equation}
	v_{i,x}=\mu_i^{-1}(\tilde{\la}_i-\si_i)v_{i}.
\end{equation}
From \eqref{identity-1} we obtain
\begin{equation}
	\mu_i^{-1}(\tilde{\la}_i-\si_i)v_{i}B[\tilde{r}_i+v_i\tilde{r}_{i,v}]+v^2_iB\tilde{r}_{i,u}\tilde{r}_i+v_i^2\tilde{r}_i\cdot DB(u)\tilde{r}_i=v_i(A(u)-\si_i)\tilde{r}_i,
\end{equation}
or equivalently,
\begin{equation}\label{identity-2}
	v_iB\tilde{r}_{i,u}\tilde{r}_i+v_i\tilde{r}_i\cdot DB(u)\tilde{r}_i-A(u)\tilde{r}_i=-\si_i\tilde{r}_i-\mu_i^{-1}(\tilde{\la}_i-\si_i)B[\tilde{r}_i+v_i\tilde{r}_{i,v}].
\end{equation}

\begin{lemma}
\label{lemme.laisig} We observe that
	\begin{equation}
		\tilde{\lambda}_{i,u},\tilde{\lambda}_{i,v}=O(1), \tilde{\la}_{i,\sig v}=O(1)|v_i|\;\mbox{and}\;\;\tilde{\lambda}_{i,\sig}, \tilde{\lambda}_{i,\sig\sig}, \tilde{\lambda}_{i,u\sig}=O(1)|v_i|^2.
		\label{tildelamb}
	\end{equation}
\end{lemma}

\begin{proof} We recall that for any $i\in\{1,\cdots,n\}$ we have
$$
	\tilde{r}_i(u,v_i,\si_i)=r_i(u)+\sum\limits_{j\neq i}\tilde{\psi}_{ij}(u,v_i,\si_i)r_j(u).
$$
Applying \eqref{5.19}, we deduce that
\begin{equation}
\tilde{r}_{i,\si}(u,v_i,\si_i), \tilde{r}_{i,\si\si}(u,v_i,\si_i), \tilde{r}_{i,u\si}(u,v_i,\si_i),  \tilde{r}_{i,u\si\si}(u,v_i,\si_i)=O(1)|v_i|.
\label{5.43}
\end{equation}
We recall now that
 $$\tilde{\la}_i(u,v_i,\si_i)= \la_i(u)-v_i
\langle B(u)\tilde{r}_{i,u}(u,v_i,\si_i)\tilde{r}_i(u,v_i,\si_i)+\tilde{r}_i(u,v_i,\si_i)\cdot DB(u)\tilde{r}_i (u,v_i,\si_i),l_i(u)\rangle.$$
Applying \eqref{5.43}, we can conclude the proof.
\end{proof}

\section{Gradient decomposition}\label{sec:Gradient-Decomp}
Let us assume that $u:\R^+\times\R\rightarrow \R^n$ is a smooth function with small total variation and is a solution of \eqref{eqn-main}, then it is important to try to decompose the gradient $u_x$ of the solution in a suitable basis of vectors $\tilde{r}_i(u,v_i,\si_i)$, $i\in\{ 1,\cdots,n\}$. As in \cite{BB-vv-lim-ann-math}, a key point of our analysis consist in decomposing $u_x$ in a basis of gradient of travelling wave $(\tilde{r}_i(u,v_i,\si_i))_{1\leq i\leq n}$. Since we search to write $u_x$ on the following form
\begin{equation}
	u_x=\sum_i v_i\tilde{r}_i(u,\bar{v}_i,\si_i),
\end{equation}
we start by fixing the vector $\tilde{r}_i$, and to do this we need to define the wave strengths $v_i$ and $\si_i$ in terms of $u$, $u_x$ and $u_{xx}$. Let us consider in a first time the particular case of a $j$ viscous travelling wave $u(t,x)=U(x-\si_jt)$ such that following our previous analysis we get
\begin{equation}
u_x=v_j \tilde{r}_j(u,v_j,\si_j).
\end{equation}
We can note that $u_t=-\si_j u_x$ where $\si_j$ is the speed of the wave. Now, we can write
\begin{equation}
	u_t=(B(u)u_x)_x-A(u)u_x=\omega_j\tilde{r}_j(u,v_j,\si_j)\mbox{ for some function }w_j.
\end{equation}
Therefore, we have $\si_j=-\frac{w_j}{v_j}$. For a general function $u$ we search then to decompose the vectors $u_x$, $u_t$ as follows
\begin{align}
	u_x&=\sum\limits_{i}v_i\tilde{r}_i(u,v_i,\si_i),\label{formula:u_x}\\
	u_t&=\sum\limits_{i}\omega_i\tilde{r}_i(u,v_i,\si_i),\mbox{ where }\si_i=-\frac{\omega_i}{v_i}.\label{formula:u_t}
\end{align}
Here $u_t$ is defined by $u_t=(B(u)u_x)_x-A(u)u_x.$ Unfortunately we have seen previously that we are able to define the vector $\tilde{r}_i(u,v_i,\si_i)$ provided that $\si_i$ belongs to a neighborhood of $\la_i(u^*)$ (here $u^*$ corresponds to $\lim_{x\rightarrow-\infty}u(t,x)=u^*$), it is not the case when the ratio $-\frac{w_i}{v_i}.$ is too large. It is in particular the case when $v_i$ is sufficiently small which can appear for example when $u_x$ is sufficiently small. In order to overcome this difficulty, as in \cite{BB-vv-lim-ann-math} we introduce a cutoff odd function $\theta:\R\rr[-2\de_1,2\de_1]$ such that
\begin{equation}
	\theta(s)=\left\{\begin{array}{rl}
		s&\mbox{ if }\abs{s}\leq \de_1,\\
		0&\mbox{ if }\abs{s}\geq 3\de_1,
	\end{array}\right.\quad \abs{\theta^\p}\leq 1\mbox{ and  }\abs{\theta^{\p\p}}\leq \frac{4}{\de_1}.
	\label{choixtheta}
\end{equation}
We consider a new variable $w_i$ defined as $\omega_i=w_i-\la^*_iv_i$. In our analysis, $\si_i$ is required to coincide with $-\omega_i/v_i$ whenever the ratio is in a sufficiently small neighborhood of $\la_i^*$. Let $\eta$ and $\xi$ be $C^3$ functions satisfying the following conditions
\begin{align}
	\eta(s)&:=\left\{\begin{array}{cl}
		1&\mbox{ for }0\leq \abs{s}\leq 3/4,\\
		0\leq \eta\leq 1&\mbox{ for }3/4\leq \abs{s}\leq 4/5,\\
		0&\mbox{ for }\abs{s}\geq 4/5,
	\end{array}\right.\\
	\xi(s)&:=\left\{\begin{array}{cl}
		1&\mbox{ for }0\leq s\leq \frac{\de_1}{2},\\
		0\leq \xi\leq 1&\mbox{ for }\frac{\de_1}{2}\leq s\leq \de_1,\\
		0&\mbox{ for }s\geq \de_1.
	\end{array}\right.
\end{align} 
Now, we consider $\xi_i=\xi(\frac{w_i}{v_i})$ and $\eta_{ij},\eta_i$ are defined as follows
\begin{equation}
	\eta_{ij}:=\eta\left(\frac{v_j^2}{v_i}\right)\mbox{ and }\eta_i:=\prod\limits_{j\neq i}\eta\left(\frac{v_j^2}{v_i}\right).
\end{equation}
\noi\textbf{Cut-off function $\chi^\ep_i$:}
Let $\chi:\R\rr[0,\f)$ be a $C^4$ function such that
\begin{equation}
	\chi(s)=0\mbox{ for }s\in[-1,1]\mbox{ and }\chi(s)=1\mbox{ for } s\in (-\f,-2]\cup[2,+\f).
\end{equation}
Now, we consider $\chi^\ep_i=\chi(\frac{v_i^{2N}}{\ep})$ for some $\ep>0$ and a positive integer $N>0$. We define $\bar{v}_i$ as follows
\begin{equation}
	\bar{v}_i^\e=v_i\chi^\ep_i\prod\limits_{j\neq i}\eta\left(\frac{v_j^2}{v_i}\right)=v_i\chi^\ep_i\prod\limits_{j\neq i}\eta_{ij}=v_i\chi^\ep_i\eta_i=\chi^\ep_i\hat{v}_i,
\end{equation} 
with $\tilde{v}_i=v_i\eta_i$.
We would like now to decompose
\begin{align}
	u_x&=\sum\limits_{i}v_i\tilde{r}_i(u,\bar{v}_i^\e\xi_i,\si_i),                            \label{eqn-v-i}\\
	u_t&=\sum\limits_{i}(w_i-\la_i^*v_i)\tilde{r}_i(u,\bar{v}_i^\e\xi_i,\si_i),           \label{eqn-w-i}
\end{align}
where 
\begin{equation*}
	u_t=(B(u)u_x)_x-A(u)u_x,\quad \si_i=\la_i^*-\theta\left(\frac{w_i}{v_i}\right).
\end{equation*}
\begin{remark}
We remark that $\si_i$ is not well defined when $v_i=w_i=0$. However from \eqref{observationimp}, we have in this case $\tilde{r}_i=r_i(u)$ regardless the value of $\si_i$.
\end{remark}
\begin{remark}
It is important to mention that our decomposition depends on $\e>0$. In our analysis we will obtain uniform $BV$ estimates in $\e>0$ for the solution $u$ of \eqref{eqn-main} on some time interval of the form $[0,C(\frac{1}{\e})]$ with $C(\frac{1}{\e})$ going to $+\infty$ when $\e\rightarrow 0$. It will be sufficient to prove that we control the $BV$ norm of $u$ for all time.\\
In the sequel for simplicity in the notations, we will adopt a slight abuse of notations inasmuch as we should write $v_i^\e$, $\bar{v}_i^\e$, $w_i^\e$ and $\si_i^\e$ in \eqref{eqn-v-i} and  \eqref{eqn-w-i}. In addition when it will be clear, we will note even $\bar{v}_i^\e$ by simply $\bar{v}_i$.
\end{remark}
\begin{remark}
	We would like to point out that it is important to choose $\xi$ such that the support of $\xi$ is included in the region where $\theta$ is equal to $1$. Indeed it will ensure that when $\xi\left(\frac{w_i}{v_i}\right)\ne 0$ then we have $\theta\left(\frac{w_i}{v_i}\right)=\frac{w_i}{v_i}$. It will simplify in particular the computation later.
\end{remark}
We would like now to make some comments on this choice of gradient decomposition.
\begin{itemize}
\item It is important to remark that if $u(t,x)=U(x-\si t)$ is a travelling wave such that $u_x=v_i\tilde{r}_i(u,v_i,\si)$ then we observe that $(u_x,u_t)$ satisfies the previous decomposition with $v_i$ and $w_i=(\la_i^*-\si)v_i$ provided that $\si$ is sufficiently close from $\la_i^*$ in terms of $\delta_1$ and that $|v_i^{2N}|\geq  2\e$.
\item We observe that when it exists $j\ne i$ such that $v_j^2\geq \frac{4}{5}|v_i|$ then the vector $\tilde{r}_i(u,v_i,\si_i)$ is not activated in the sense that $\tilde{r}_i(u,\xi_i \bar{v}_i,\si_i)=\tilde{r}_i(u,0,\si_i)=r_i(u)$. We recover here the decomposition with a classical eigenvector $r_i(u)$ of $A(u)$. 
\item Similarly we observe that if $|v_i|^{2N}\leq \e$ then we have  $\tilde{r}_i(u,\xi_i \bar{v}_i,\si_i)=\tilde{r}_i(u,0,\si_i)=r_i(u)$ this is due to the fact that $\chi^\ep_i=\chi(\frac{v_i^{2N}}{\ep})=0$, this approximation in $\e$ is necessary in order to prove in the sequel that the functions $v_i$ and $w_i$ are regular, we need in particular to avoid the zone $\{v_i=0\}$.
\item We insist on the fact that compared with \cite{BB-vv-lim-ann-math}, our gradient decomposition is significantly more complicated since in our analysis we have to deal with new interaction terms of the form $w_{i,xx}v_i-v_{i,xx}w_i$ which requires more precise estimates on the reminder terms $\phi_i$.
\end{itemize}
We need now to analyse $\hat{v}_i=v_i\eta_i$ and $\bar{v}_i^\e$. First we have the following Lemma.
\begin{lemma}\label{estimimpo1bis}
	We have setting $\varrho_i=\prod_{j\ne i}\mathbbm{1}_{\{0\leq \frac{v_j^2}{|v_i|}\leq\frac{4}{5}\}}$
	\begin{equation*}
		\pa_{v_i}\hat{v}_i=O(1)\varrho_i,\quad\pa_{v_i v_i}\hat{v}_i=O(1)\frac{1}{\abs{v_i}}\varrho_i,\quad
		\pa_{v_i v_i v_i}\hat{v}_i=O(1)\frac{1}{\abs{v_i}^2}\varrho_i.
	\end{equation*}
	
	If $i\ne j$, $k\ne i$ and $l\ne i$, we have
	\begin{align*}
		&\pa_{v_j}\hat{v}_i=O(1)\abs{v_j}\varrho_i ,\,
		\pa_{v_j v_k}\hat{v}_i=O(1)\varrho_i,\, \pa_{v_i v_j}\hat{v}_i=O(1)\frac{\abs{v_j}}{\abs{v_i}} \varrho_i,\\
		%& \pa_{w_j}\bar{v}_i=O(1)\abs{w_j} \varrho_i,\, \pa_{w_j v_k}\bar{v}_i,\pa_{w_j w_k}\bar{v}_i=O(1)\varrho_i,\, \pa_{v_i w_j}\bar{v}_i=O(1)\frac{\abs{w_j}}{\abs{v_i}}\varrho_i,\\
		&\pa_{v_j v_iv_i}\hat{v}_i=O(1)\frac{\abs{v_j}}{v^2_i}\varrho_i,%\,\pa_{w_j v_iv_i}\bar{v}_i=O(1)\frac{\abs{w_j}}{v^2_i}\varrho_i
		\,\pa_{v_j v_kv_i}\hat{v}_i=O(1)\frac{\abs{v_jv_k}}{v_i^2}\varrho_i
		%&\pa_{w_k v_j v_i}\bar{v}_i=O(1)\frac{\abs{v_jw_k}}{v_i^2}\varrho_i
		,\,\pa_{v_j v_kv_l}\hat{v}_i=O(1)\frac{\abs{v_j}}{\abs{v_i}}\varrho_i.%,\,\pa_{w_k v_j v_l}\bar{v}_i=O(1)\frac{\abs{v_j}}{\abs{v_i}}\varrho_i.
	\end{align*}
\end{lemma}
\begin{proof}
	We have
	\begin{align*}
		\pa_{v_i}\hat{v_i}=\prod_{k\ne i}\eta(\frac{v_k^2}{v_i})-\sum_{k\ne i}\left[\prod_{l\ne k,i}\eta(\frac{v_l^2}{v_i})\right]\eta'(\frac{v_k^2}{v_i}) \frac{v_k^2}{v_i}=O(1)\varrho_i.
	\end{align*}
	Similarly we obtain for $i\ne j$
	\begin{align*}
		&\pa_{v_j}\hat{v_i}=2v_j\left[\prod_{k\ne i,j}\eta(\frac{v_k^2}{v_i})\right] \eta'(\frac{v_j^2}{v_i})=O(1)|v_j|\varrho_i.
		%&\pa_{w_j}\tilde{v_i}=2w_j\left[\prod_{k\ne i,j}\eta(\frac{v_k^2+w_k^2}{v_i})\right]  \eta'(\frac{v_j^2+w_j^2}{v_i} ).
	\end{align*}
	Let us consider now the second derivative which gives
	\begin{align*}
		\pa_{v_iv_i}\hat{v_i}&=\sum_{k\ne i}\left[\prod_{l\ne k,i}\eta(\frac{v_l^2}{v_i})\right]\eta''(\frac{v_k^2}{v_i}) \frac{1}{v_i}(\frac{v_k^2}{v_i})^2+\sum_{k\ne i}\eta'(\frac{v_k^2}{v_i}) \frac{v_k^2}{v_i}\sum_{l\ne k,i}\left[\prod_{l'\ne l,k,i}
		\eta(\frac{v_{l'}^2}{v_i})\right] \eta'(\frac{v_l^2}{v_i})\frac{v_l^2}{v^2_i}\\
		&=O(1)\frac{1}{|v_i|}\varrho_i.
	\end{align*}
	Similarly we obtain:
	\begin{align*}
		&\pa_{v_k v_j}\hat{v_i}=4\frac{v_jv_k}{v_i}\left[\prod_{l\ne i,j,k}\eta(\frac{v_l^2}{v_i})\right] \eta'(\frac{v_j^2}{v_i})\eta'(\frac{v_k^2}{v_i})=O(1)\varrho_i.
		%&\pa_{w_k w_j}\bar{v_i}=4\frac{w_j w_k}{v_i}\left[\prod_{k\ne i,j}\eta(\frac{v_k^2+w_k^2}{v_i})\right]  \eta'(\frac{v_j^2+w_j^2}{v_i} )\eta'(\frac{v_k^2+w_k^2}{v_i})=O(1).
	\end{align*}
	We proceed similarly for the other terms. We omit explicit calculation of the terms. This completes the proof Lemma \ref{estimimpo1}.
\end{proof}

To analyze the derivatives of $\eta_i,\eta_{ij},\chi_i^\ep$ we introduce the following notations. %We set
\begin{equation}
	\varkappa^\ep_i=\mathbbm{1}_{\left\{\ep\leq \abs{v_i}^{2N}\leq 2\ep\right\}},\,\,\rho^\ep_i=\mathbbm{1}_{\left\{ \abs{v_i}^{2N}\geq \ep \right\}}\prod\limits_{k\neq i}\mathbbm{1}_{\left\{0\leq\frac{\abs{v_k}^2}{\abs{v_i}}\leq \frac{4}{5}\right\}},
\end{equation}
and
\begin{align*}
	\Delta^\ep_{ik}&=\mathbbm{1}_{\left\{ \abs{v_i}^{2N}\geq \ep\right\}}\mathbbm{1}_{\left\{\frac{3}{4}\leq\frac{\abs{v_k}^2}{\abs{v_i}}\leq \frac{4}{5}\right\}}\prod\limits_{l\neq i,k}\mathbbm{1}_{\left\{0\leq\frac{\abs{v_l}^2}{\abs{v_i}}\leq \frac{4}{5}\right\}},\quad \Delta^\ep_i=\sum\limits_{k\neq i}\Delta^\ep_{ik},\\
	\Box^\ep_{ijk} &=\mathbbm{1}_{\left\{ \abs{v_i}^{2N}\geq \ep\right\}}\mathbbm{1}_{\left\{\frac{3}{4}\leq\frac{\abs{v_j}^2}{\abs{v_i}}\leq \frac{4}{5}\right\}}\mathbbm{1}_{\left\{\frac{3}{4}\leq\frac{\abs{v_k}^2}{\abs{v_i}}\leq \frac{4}{5}\right\}}\prod\limits_{l\neq i,j,k}\mathbbm{1}_{\left\{0\leq\frac{\abs{v_l}^2}{\abs{v_i}}\leq \frac{4}{5}\right\}},\\
	\Box^\ep_{ij}&=\sum\limits_{k\neq i,j}\Box^\ep_{ijk}\mbox{ and }\Box^\ep_i=\sum\limits_{k\neq i}\Box^\ep_{ik}.
\end{align*}
Furthermore, we set
\begin{equation}
	\mathfrak{A}_i:=\mathbbm{1}_{\left\{0\leq \abs{\frac{w_i}{v_i}}\leq 3\de_1\right\}}.
\end{equation}
\begin{lemma}\label{estimimpo1}
	For $1\leq i\leq n$ we have
	\begin{align}
		&\pa_{v_i}\eta_i=O(1)\frac{\Delta_i^\e}{|v_i|}, 
		\pa_{v_i}\bar{v}_i^\e=
		\chi_i^\ep\eta_i+O(1)\left(\Delta^\ep_i+\varkappa^\ep_i\rho^\ep_i\right),
			\pa_{v_iv_i}\bar{v}_i=O(1)\left(\Delta^\ep_i+\varkappa^\ep_i\rho^\ep_i\right)\frac{1}{\abs{v_i}},
			\\
		&\pa_{v_iv_iv_i}\bar{v}_i^\e=
		\left(
		\Delta^\ep_i+\varkappa^\ep_i\rho^\ep_i\right)\frac{1}{\abs{v_i}^2}.
	\end{align}
	Furthermore, for $j\neq i$, $k\neq i,j$ and $l\neq i,j,k$ we have
	\begin{align}
		&\pa_{v_j}\eta_i=O(1)\frac{|v_j|}{|v_i|}\Delta^{\e}_{ij},,\,\,\quad\quad\quad\pa_{v_j}\bar{v}_i^\e=O(1)\abs{v_j} \Delta^\ep_{ij},\,\,\quad	\pa_{v_iv_j}\bar{v}_i^\e={O}(1) \Delta^\ep_{ij}\frac{\abs{v_j}}{\abs{v_i}},\\
		 &\pa_{v_jv_j}\bar{v}_i^\e=O(1) \Delta^\ep_{ij},\,\,\quad\quad\quad\quad \pa_{v_jv_k}\bar{v}_i={O}(1) \Box^\ep_{ijk},\\
		 &\pa_{v_iv_iv_j}\bar{v}_i^\e=O(1) \Delta^\ep_{ij}\frac{\abs{v_j}}{\abs{v_i}^2},\quad\quad \pa_{v_iv_jv_j}\bar{v}_i^\e=O(1) \Delta^\ep_{ij}\frac{1}{\abs{v_i}},\\
		& \pa_{v_jv_jv_j}\bar{v}_i^\e=O(1) \Delta^\ep_{ij}\frac{\abs{v_j}}{\abs{v_i}},	\quad\quad \pa_{v_jv_jv_k}\bar{v}_i^\e=O(1) \Box^\ep_{ijk}\frac{\abs{v_k}}{\abs{v_i}},\\
		&\pa_{v_iv_jv_k}\bar{v}_i^\e={O}(1)\Box^\ep_{ijk}\frac{1}{\abs{v_i}}, \quad \pa_{v_jv_kv_l}\bar{v}_i^\e=O(1) \frac{\abs{v_jv_kv_l}}{\abs{v_i}^2}\Box^\ep_{ijk}\Box^\ep_{ikl}.
	\end{align}
\end{lemma}
\begin{proof}
First we observe that
\begin{align*}
&\pa_{v_i}\eta_i=-\frac{1}{v_i}\sum_{j\ne i}\eta'_{ij}\frac{v_j^2}{v_i} \prod\limits_{l\neq i, j}\eta'_{il}=O(1)\frac{\Delta_{i}^\e}{|v_i|}.
\end{align*}
	From the definition of $\bar{v}_i^\e$, we deduce that
	\begin{align*}
		\pa_{v_i}\bar{v}_i^\e&=%\frac{2N (v_i)^{2N}}{\e}(\chi')_i^\e\eta_i+\chi^\e_i \pa_{v_i}\hat{v}_i=O(1)\rho_i^\e.\\
		-\chi^\ep_i \sum\limits_{k\neq i}\frac{(v_k)^2}{v_i}\eta^{\p}_{ik} \prod\limits_{l\neq i, k}\eta_{il}+\frac{2N(v_i)^{2N}}{\ep}(\chi')^{\ep}_i\eta_i+\chi_i^\ep\eta_i\\
		&=\chi_i^\ep\eta_i+{O}(1)\left(\Delta^\ep_i+\varkappa^\ep_i\rho^\ep_i\right),\\
		\pa_{v_iv_i}\bar{v}_i^\e&=%=\chi^\e_i \pa_{v_iv_i}\hat{v}_i+4N (\chi')_i^\e \frac{v_i^{2N-1}}{\e} \pa_{v_i}\hat{v}_i+\hat{v}_i \big(2N(2N-1)\frac{(v_i)^{2N-2}}{\e}(\chi')_i^\e+4N^2(\frac{ (v_i)^{2N-1}}{\e})^2(\chi'')_i^\e\big),\\
				\chi^\ep_i\sum\limits_{k\neq i}\sum\limits_{m\neq i,k}\frac{(v_k)^2(v_m)^2}{(v_i)^3}\eta^{\p}_{ik} \eta^{\p}_{im} \prod\limits_{l\neq i, k,m}\eta_{il}+\chi^\ep_i \sum\limits_{k\neq i}\frac{v_k^4}{v^3_i}\eta^{\p\p}_{ik} \prod\limits_{l\neq i, k}\eta_{il}\\
&-\frac{4N(v_i)^{2N-1}}{\ep}(\chi')^{\ep}_i\sum\limits_{k\neq i}\frac{(v_k)^2}{v_i}\eta^{\p}_{ik} \prod\limits_{l\neq i, k}\eta_{il}\\
&+\left(\frac{2N(2N+1)(v_i)^{2N-1}}{\ep}(\chi'){\ep}_i+\frac{4N^2(v_i)^{4N-1}}{\ep^2}(\chi'')^{\ep}_i\right)\eta_{i}\\
		&={O}(1)\left(\Delta^\ep_i+\varkappa^\ep_i\rho^\ep_i\right)\frac{1}{\abs{v_i}},
	\end{align*}
	and also we calculate
	\begin{align*}
		\pa_{v_iv_iv_i}\bar{v}_i^\e&=-\chi^\ep_i\sum\limits_{k\neq i}\sum\limits_{m\neq i,k}\sum\limits_{p\neq i,k,m}\frac{(v_k)^2(v_m)^2(v_p)^2}{(v_i)^5}\eta^{\p}_{ik} \eta^{\p}_{im}\eta^{\p}_{ip} \prod\limits_{l\neq i, k,m,p}\eta_{il}\\
		&-3\chi^\ep_i\sum\limits_{k\neq i}\sum\limits_{m\neq i,k}\frac{(v_k)^4(v_m)^2}{(v_i)^5}\eta^{\p\p}_{ik} \eta^{\p}_{im} \prod\limits_{l\neq i, k,m}\eta_{il}
		-3\chi^\ep_i\sum\limits_{k\neq i}\sum\limits_{m\neq i,k}\frac{(v_k)^2(v_m)^2}{(v_i)^4}\eta^{\p}_{ik} \eta^{\p}_{im} \prod\limits_{l\neq i, k,m}\eta_{il}\\
		&-3\chi^\ep_i \sum\limits_{k\neq i}\frac{v_k^4}{v^4_i}\eta^{\p\p}_{ik} \prod\limits_{l\neq i, k}\eta_{il}-\chi^\ep_i \sum\limits_{k\neq i}\frac{v_k^6}{v^{5}_i}\eta^{\p\p\p}_{ik} \prod\limits_{l\neq i, k}\eta_{il}\\
		&+\frac{6N(v_i)^{2N-1}}{\ep}(\chi')^{\ep}_i\left(\sum\limits_{k\neq i}\sum\limits_{m\neq i,k}\frac{(v_k)^2(v_m)^2}{(v_i)^3}\eta^{\p}_{ik} \eta^{\p}_{im} \prod\limits_{l\neq i, k,m}\eta_{il}+\sum\limits_{k\neq i}\frac{v_k^4}{v^3_i}\eta^{\p\p}_{ik} \prod\limits_{l\neq i, k}\eta_{il}\right)\\
		&-\left(\frac{4N(2N-2) (v_i)^{2N-2}}{\ep}(\chi')^{\ep}_i+\frac{8N^2(v_i)^{2N-2}}{\ep^2}(\chi^{\p\p})^{\ep}_i\right)\sum\limits_{k\neq i}\frac{(v_k)^2}{v_i}\eta^{\p}_{ik} \prod\limits_{l\neq i, k}\eta_{il}\\
	&+\left(\frac{2N(4N^2-1)(v_i)^{2N-2}}{\ep}\chi^{\ep,\p}_i+\frac{24N^3(v_i)^{4N-2}}{\ep^2}\chi^{\ep,\p\p}_i+\frac{8N^3(v_i)^{6N-2}}{\ep^3}\chi^{\ep,\p\p\p}_i\right)\eta_{i}\\
		&{\color{black}{-}}\left(\frac{2N(2N+1)(v_i)^{2N-1}}{\ep}\chi^{\ep,\p}_i+\frac{4N^2(v_i)^{4N-1}}{\ep^2}\chi^{\ep,\p\p}_i\right)\sum\limits_{k\neq i}\frac{(v_k)^2}{v_i^2}\eta^{\p}_{ik} \prod\limits_{l\neq i, k}\eta_{il}\\
		&={O}(1)\left(
		\Delta^\ep_i+\varkappa^\ep_i\rho^\ep_i\right)\frac{1}{\abs{v_i}^2}.
	\end{align*}
	Now, for $j\neq i$ we get
	\begin{equation*}
		\pa_{v_j}\bar{v}^\ep_i=\chi^\ep_{i} 2v_j \eta^{\p}_{ij}\prod\limits_{l\neq i, j}\eta_{il}={O}(1)\abs{v_j} \Delta^\ep_{ij}.
	\end{equation*}
	Similarly,
	\begin{align*}
		\pa_{v_iv_j}\bar{v}_i^\e&=-\chi^\ep_i 2v_j\frac{(v_j)^2}{(v_i)^2} \eta^{\p\p}_{ij}  \prod\limits_{l\neq i, j}\eta_{il}-\chi^\ep_{i}2v_j \eta^{\p}_{ij} \sum\limits_{k\neq j,i}\frac{(v_k)^2}{(v_i)^2}\eta^{\p}_{ik} \prod\limits_{l\neq i, j,k}\eta_{il}\\
		&+\frac{4N(v_i)^{2N-1}v_j}{\ep}\chi^{\ep,\p}_i \eta^{\p}_{ij}  \prod\limits_{l\neq i, j}\eta_{il}=O(1) \Delta^\ep_{ij}\frac{\abs{v_j}}{\abs{v_i}},\\
		\pa_{v_jv_j}\bar{v}_i^\e&=\chi^\ep_{i} 2 \eta^{\p}_{ij}\prod\limits_{l\neq i, j}\eta_{il}+\chi^\ep_{i} \frac{4v^2_j}{v_i} \eta^{\p\p}_{ij}\prod\limits_{l\neq i, j}\eta_{il}=O(1) \Delta^\ep_{ij},\\
		\pa_{v_jv_k}\bar{v}_i^\e&=-\chi^\ep_i \frac{4v_jv_k}{v_i}\eta^{\p}_{ij}\eta^{\p}_{ik}  \prod\limits_{l\neq i, j,k}\eta_{il}=O(1) \Box^\ep_{ijk}.
	\end{align*}
	Furthermore, we can calculate
		\begin{align*}
		\pa_{v_iv_iv_j}\bar{v}_i^\e&=\chi^\ep_i 4\frac{(v_j)^3}{(v_i)^3} \eta^{\p\p}_{ij}  \prod\limits_{l\neq i, j}\eta_{il}+\chi^\ep_i 2\frac{(v_j)^5}{(v_i)^4} \eta^{\p\p\p}_{ij}  \prod\limits_{l\neq i, j}\eta_{il}+\chi^\ep_i 4v_j\frac{(v_j)^2}{(v_i)^2} \eta^{\p\p}_{ij}  \sum\limits_{k\neq i,j}\frac{v_k^2}{v_i^2}\eta_{ik}^\ep\prod\limits_{l\neq i, j,k}\eta_{il}\\
		&
		+4\chi^\ep_{i}v_j \eta^{\p}_{ij} \sum\limits_{k\neq j,i}\frac{(v_k)^2}{(v_i)^3}\eta^{\p}_{ik} \prod\limits_{l\neq i, j,k}\eta_{il}\\
		&+\chi^\ep_{i}2v_j \eta^{\p}_{ij} \sum\limits_{k\neq j,i}\sum\limits_{m\neq i,j,k}\frac{v_k^2v_m^2}{(v_i)^4}\eta^{\p}_{ik} \eta^{\p}_{im}\prod\limits_{l\neq i, j,k,m}\eta_{il}+2\chi^\ep_{i}v_j \eta^{\p}_{ij} \sum\limits_{k\neq j,i}\frac{(v_k)^4}{(v_i)^4}\eta^{\p\p}_{ik} \prod\limits_{l\neq i, j,k}\eta_{il}\\
		&{\color{black}{-}}\frac{4N(v_i)^{2N-1}v_j}{\ep}(\chi')^{\ep}_i \left(\frac{(v_j)^2}{(v_i)^2} \eta^{\p\p}_{ij}  \prod\limits_{l\neq i, j}\eta_{il}+2 \eta^{\p}_{ij} \sum\limits_{k\neq j,i}\frac{(v_k)^2}{(v_i)^2}\eta^{\p}_{ik} \prod\limits_{l\neq i, j,k}\eta_{il}\right)\\
		&+\left(\frac{4N(2N-1)(v_i)^{2N-2}v_j}{\ep}(\chi')^{\ep}_i +\frac{8N^2(v_i)^{4N-2}v_j}{\ep^2}(\chi'')^{\ep}_i  \right) \eta^{\p}_{ij} \prod\limits_{l\neq i, j}\eta_{il}\\
		&=O(1) \Delta^\ep_{ij}\frac{\abs{v_j}}{\abs{v_i}^2},\\
		%----------------------------------------------------------
		\pa_{v_iv_jv_j}\bar{v}_i^\e&=-\chi^\ep_i \frac{6v_j^2}{(v_i)^2} \eta^{\p\p}_{ij}  \prod\limits_{l\neq i, j}\eta_{il}-4\chi^\ep_i \frac{(v_j)^4}{(v_i)^3} \eta^{\p\p\p}_{ij}  \prod\limits_{l\neq i, j}\eta_{il}-2\chi^\ep_{i} \eta^{\p}_{ij} \sum\limits_{k\neq j,i}\frac{(v_k)^2}{(v_i)^2}\eta^{\p}_{ik} \prod\limits_{l\neq i, j,k}\eta_{il}\\
		&-4\chi^\e_i\frac{v_j^2}{v_i}\eta_{ij}^{\p\p}
		 \sum\limits_{k\neq j,i}\frac{(v_k)^2}{(v_i)^2}\eta^{\p}_{ik} \prod\limits_{l\neq i, j,k}\eta_{il}+\frac{4N(v_i)^{2N-1}}{\ep}(\chi')^{\ep}_i \eta^{\p}_{ij}  \prod\limits_{l\neq i, j}\eta_{il}\\
		 &+\frac{8N(v_i)^{2N-1}}{\ep }\frac{v_j^2}{v_i}(\chi')^{\ep}_i \eta^{\p}_{ij}  \prod\limits_{l\neq i, j}\eta_{il}=O(1) \Delta^\ep_{ij}\frac{1}{\abs{v_i}},\\
		%----------------------------------------------------------
		\pa_{v_jv_jv_j}\bar{v}_i^\e&=\chi^\ep_{i} \frac{12v_j}{v_i} \eta^{\p\p}_{ij}\prod\limits_{l\neq i, j}\eta_{il}+\chi^\ep_{i} \frac{8v^3_j}{v_i^2} \eta^{\p\p\p}_{ij}\prod\limits_{l\neq i, j}\eta_{il}=O(1) \Delta^\ep_{ij}\frac{\abs{v_j}}{\abs{v_i}},\\
		\pa_{v_iv_jv_k}\bar{v}_i^\e&=\chi^\ep_i \frac{4v_jv_k}{v^2_i}\eta^{\p}_{ij}\eta^{\p}_{ik}  \prod\limits_{l\neq i, j,k}\eta_{il}+\chi^\ep_i \frac{4v_j^3v_k}{v^3_i}\eta^{\p\p}_{ij}\eta^{\p}_{ik}  \prod\limits_{l\neq i, j,k}\eta_{il}+\chi^\ep_i \frac{4v_jv_k^3}{v^3_i}\eta^{\p}_{ij}\eta^{\p\p}_{ik}  \prod\limits_{l\neq i, j,k}\eta_{il}\\
		&+\chi^\ep_i \frac{4v_jv_k}{v_i}\eta^{\p}_{ij}\eta^{\p}_{ik}  \sum\limits_{m\neq i,j,k}\frac{v_m^2}{v_i^2}\eta_{im}^\p \prod\limits_{l\neq i, j,k,m}\eta_{il}-\frac{8 N (v_i)^{2N-1}}{\ep}(\chi')^{\ep}_i \frac{v_jv_k}{v_i}\eta^{\p}_{ij}\eta^{\p}_{ik}  \prod\limits_{l\neq i, j,k}\eta_{il}\\
		&={O}(1)\Box^\ep_{ijk}\frac{1}{\abs{v_i}},\\
		%----------------------------------------------------------
		\pa_{v_jv_jv_k}\bar{v}_i^\e&=-\chi^\ep_i \frac{4v_k}{v_i}\eta^{\p}_{ij}\eta^{\p}_{ik}  \prod\limits_{l\neq i, j,k}\eta_{il}-\chi^\ep_i \frac{8v_j^2v_k}{v_i^2}\eta^{\p\p}_{ij}\eta^{\p}_{ik}  \prod\limits_{l\neq i, j,k}\eta_{il}=O(1) \Box^\ep_{ijk}\frac{\abs{v_k}}{\abs{v_i}}\\
		%----------------------------------------------------------
		\pa_{v_jv_kv_l}\bar{v}_i^\e&=-\chi^\ep_i \frac{4v_jv_kv_l}{v_i^2}\eta^{\p}_{ij}\eta^{\p}_{ik}\eta^{\p}_{il}  \prod\limits_{m\neq i, j,k,l}\eta_{im}=O(1) \frac{\abs{v_jv_kv_l}}{\abs{v_i}^2}\Box^\ep_{ijk}\Box^\ep_{ikl}.
	\end{align*}
	This completes the proof of Lemma \ref{estimimpo1}.
\end{proof}

\begin{lemma}
\label{lemme6.61}
	 There exists $\hat{\de}>0$ such that the following holds true. For $\abs{u-u^*},\abs{u_x},\abs{u_{xx}}\leq\de_0$, the system of $2n$ equations \eqref{eqn-v-i}--\eqref{eqn-w-i} has a unique solution $(v,w)=(v_1,\cdots,v_n,w_1,\cdots,w_n)$ for any $\e>0$. Furthermore for any $\e>0$, the map $(u,u_x,u_{xx})\mapsto (v,w)$ is smooth for $\e>0$ sufficiently small
	on a neighborhood $V$ independent on $\e$ of the point $(u^*,0,0)$. 
\end{lemma}

\begin{proof}
	%\noi\underline{Existence:}
	 We consider $\La$ defined as follows,
	\begin{align}
		\La(u,v,w)&=\sum\limits_{i=1}^{n}\La_i(u,v_i,w_i),\\
		\La_i(u,v_i,w_i)&=\begin{pmatrix}
			v_i\tilde{r}_i(u,\xi_i\bar{v}_i,\la^*_i-\theta(w_i/v_i))\\
			(w_i-\la^*_iv_i)\tilde{r}_i(u,\xi_i\bar{v}_i,\la^*_i-\theta(w_i/v_i)).
		\end{pmatrix}.
	\end{align}
	First we observe that if $\R^n$ is endowed with the $\|\cdot\|_\infty$ norm then we have for $v\in B(0,\e^{\frac{1}{2 N}})$ and any $i\in\{1,\cdots,n\}$:
\begin{align}
		\La_i(u,v_i,w_i)&=\begin{pmatrix}
			v_i r_i(u)\\
			(w_i-\la^*_iv_i)r_i(u)
		\end{pmatrix}.
	\end{align}	
	This is due to the fact that in this case 
	we have $\bar{v_i}=0$ and we can apply \eqref{observationimp}.
	We calculate 
	\begin{equation*}
		\frac{\pa\La}{\pa(v_i,w_i)}=\begin{pmatrix}
			\tilde{r}_i&0\\
			-\la_i^*\tilde{r}_i&\tilde{r}_i
		\end{pmatrix}+\begin{pmatrix}
			Q^{i}_{11}&Q^{i}_{12}\\
			Q^{i}_{21}&Q^{i}_{22}
		\end{pmatrix},
	\end{equation*}
	where $Q_{lk}^i$ are defined as follows
	\begin{align*}
		Q_{11}^i&:=\left(\xi_iv_i\pa_{v_i}\bar{v}_i\tilde{r}_{i,v}+\sum\limits_{j\neq i}\xi_jv_j\pa_{v_i}\bar{v}_j\tilde{r}_{j,v}\right)+\frac{w_i}{v_i}[-\xi_i^\p\bar{v}_i\tilde{r}_{i,v}+\theta^\p_i\tilde{r}_{i,\si}], \\
		Q_{12}^i&:=[\xi_i^\p\bar{v}_i\tilde{r}_{i,v}-\theta^\p_i\tilde{r}_{i,\si}], \\
		Q_{21}^i&:=\left(\xi_i(w_i-\la_i^*v_i)\pa_{v_i}\bar{v}_i\tilde{r}_{i,v}+\sum\limits_{j\neq i}\xi_j(w_j-\la_j^*v_j)\pa_{v_i}\bar{v}_j\tilde{r}_{j,v}\right)+\frac{(w_i-\la_i^*v_i)w_i}{v^2_i}[-\xi_i^\p\bar{v}_i\tilde{r}_{i,v}+\theta^\p_i\tilde{r}_{i,\si}],\\
		Q_{22}^i&:=(\frac{w_i}{v_i}-\la_i^*)[\xi_i^\p\bar{v}_i\tilde{r}_{i,v}-\theta^\p_i\tilde{r}_{i,\si}].%\\&(\la_i^*-(w_i/v_i))\theta_i^\p\tilde{r}_{i,\si}.
	\end{align*}
	 We can note again that $Q_{11}^i, Q_{12}^i, Q_{21}^i, Q_{22}^i=0$
	 when $\|v\|_{\infty}<\e^{\frac{1}{2N}}$. This is due to the fact that $\bar{v}_i=0$ in this case with the fact that $\tilde{r}_{i,\si}=O(1)\xi_i\bar{v}_i$ and $\pa_{v_i}\bar{v}_i=O(1)\rho_i^\e$, $\pa_{v_i}\bar{v}_j=O(1)|v_j|\rho_i^\e$ from the Lemma \ref{estimimpo1} and \eqref{observationimp1}.
	Then we may write
	\begin{equation*}
		\frac{\pa\La}{\pa(v_i,w_i)}
		=\begin{pmatrix}
			r_i^*&0\\
			-\la_i^*{r}^*_i& {r}^*_i
		\end{pmatrix}+\begin{pmatrix}
			\tilde{r}_i-r_i^*&0\\
			-\la_i^*(\tilde{r}_i-{r}^*_i)& \tilde{r}_i-{r}^*_i
		\end{pmatrix}+\begin{pmatrix}
			Q^{i}_{11}&Q^{i}_{12}\\
			Q^{i}_{21}&Q^{i}_{22}
		\end{pmatrix}.
	\end{equation*}
	Denote
	\begin{equation}
		\frac{\pa\La}{\pa(v,w)}=\mathcal{B}_0(u,v,w)+\mathcal{B}_1(u,v,w)+\mathcal{B}_2(u,v,w).
	\end{equation}
	We observe that using again  \eqref{observationimp1} and Lemma  \ref{estimimpo1} that
	\begin{equation}
		\abs{\mathcal{B}_1(u,v,w)}+\abs{\mathcal{B}_2(u,v,w)}=O(1)\left(\abs{v}+\abs{w}+\abs{u-u^*}\right).
	\end{equation}
	Since $\{r_i^*\}$ forms a basis, $B_0$ is uniformly invertible. Hence, in a small enough neighbourhood of $(u^*,0,0)$, $	\frac{\pa\La}{\pa(v,w)}$ is uniformly invertible in $\e$ and $u$. In fact we note that for any sequence $(v_n,w_n)\rr0$ the limit $\lim\limits_{n\rr\f}	\frac{\pa\La}{\pa(v,w)}(u,v_n,w_n)$ is uniformly invertible provided that $\abs{u-u^*}$ is sufficiently small and for any $\e>0$  small enough. Hence, by generalized implicit function theorem \cite{Vanderbauwhede}, we obtain that $\La^{-1}$ exists and is regular in a neighbourhood of $(0,0)$ independent on $\e$. Therefore for a given $(u,u_x,u_{xx})$ there exist unique value $(v,w)$ such that
	\begin{equation}\label{La:v-w:u-x-u-xx}
		\La(u,v,w)=(u_x,B(u)u_{xx}+u_x\bullet B(u) u_x-A(u)u_x).
	\end{equation}
	This completes the proof of Lemma \ref{lemme6.61}.
\end{proof}
\begin{remark} We would like here explain why we choose to introduce the variable $\bar{v}_i^\e=\chi_i^\e\tilde{v}_i$ for $\e>0$ whereas the choice $\hat{v}_i$ seems more natural. Indeed this last choice enables to rewrite a travelling wave with the corresponding decomposition as $U'=v_i\tilde{r}_i(u,\xi_i\tilde{v}_i,\si_i)$ whereas it is not the case with $\bar{v}_i^\e$ in particular if $|v_i|^{2N}\leq \e$.From this point our approximation seems less efficient.
	
First by using the previous Lemma, we observe that we can decompose the vectors $(u_x,u_t)$ as follows:
\begin{align}
	u_x&=\sum\limits_{i}v_i\tilde{r}_i(u,\bar{v}_i^\e\xi_i,\si_i),         \\%                   \label{eqn-v-i}\\
	u_t&=\sum\limits_{i}(w_i-\la_i^*v_i)\tilde{r}_i(u,\bar{v}_i^\e\xi_i,\si_i),    %       \label{eqn-w-i}
\end{align}
In addition for $\e>0$ small enough  we have $(v,w)=\La^{-1}(u,u_x,(B(u)u_x)_x-A(u)u_x)$ provided that $(u,u_x,(B(u)u_x)_x-A(u)u_x)$ is living in $V$.
We deduce that if $u$ is the solution of \eqref{eqn-main} then provided that $(u(t,\cdot),u_x(t,\cdot),(B(u)u_x)_x(t,\cdot)-A(u)u_x)(t,\cdot))$ is living in $V$ for $t$ large enough (it will be the case), then we can ensure that $(v,w)$ are regular in $t,x$ since $\La^{-1}$ is regular and $u$ as solution of a parabolic equation.
\\
In contrast if we replace $\bar{v}_i^\e$ by $\hat{v}_i$ we can check by similar argument as in 
\cite{BB-vv-lim-ann-math} that $\La^{-1}$ is only continuously differentiable which not allow to ensure that $(v,w)$ is $C^2$ for example. In fact $\La^{-1}$ is smooth outside the $n$ manifolds ${\cal N}_i=\{v_i=w_i=0\}$. It explains our approximation $\bar{v}_i^\e=\chi_i^\e\tilde{v}_i$ of $\hat{v}_i$ which avoids the zone $v_i=0$.
\end{remark}

\begin{remark}In addition we observe that for $(v,w)$ we have the following estimate which is independent on $\e>0$
	\begin{equation}\label{estimate:La-u}
		\frac{\pa\La}{\pa u}(v,w)=O(1)(\abs{v}+\abs{w}).
	\end{equation}
\end{remark}

\begin{lemma}\label{lemme6.5}
	Suppose that $u$ satisfies the assumptions of Proposition \ref{prop:parabolic} and further we assume that \eqref{assumption:u-L1} holds for larger interval $[0,T]$ for some $T>\hat{t}$. Then for $t\in[\hat{t},T]$ the decomposition \eqref{eqn-v-i}--\eqref{eqn-w-i} is well defined and $v_i,w_i$ satisfy the following bounds which are independent on $\e>0$
	\begin{align}
		\norm{v_i(t)}_{L^1},\,\norm{w_i(t)}_{L^1}&=O(1)\de_0,\\
		\norm{v_{i,t}(t)}_{L^1},\,\norm{w_{i,t}(t)}_{L^1},\,	\norm{v_i(t)}_{L^\f},\,\norm{w_i(t)}_{L^\f},\,\norm{v_{i,x}(t)}_{L^1},\,\norm{w_{i,x}(t)}_{L^1}&=O(1)\de^2_0,\\
		\norm{v_{i,t}(t)}_{L^\f},\,\norm{w_{i,t}(t)}_{L^\f},\;\norm{v_{i,x}(t)}_{L^\f},\,\norm{w_{i,x}(t)}_{L^\f}&=O(1)\de^3_0.
		%\norm{v_{i,t}(t)}_{L^\f},\,\norm{w_{i,t}(t)}_{L^\f}&=O(1)\de^4_0.
	\end{align}  
\end{lemma}
\begin{proof} First by the Lemma \ref{lemme6.6} the map $(u,u_x,u_{xx})\rightarrow (v,w)$ is well defined for any $\e>0$ sufficiently small on a neighborhood $V$ of $(u^*,0,0)$ independent on $\e$, from Corollary \ref{coro2.2} we deduce that $(u(t,\cdot),u_x(t,\cdot),u_{xx}(t,\cdot))$ belongs to $V$ for $\hat{t}\leq t\leq T$ provided that $\delta_0>0$ is sufficiently small, it particular it guarantees that the decomposition \eqref{eqn-v-i}-\eqref{eqn-w-i} is well defined.
	From \eqref{La:v-w:u-x-u-xx} and using invertibility of $\La$ and the fact that $\La^{-1}(u,\cdot)$ is uniformly Lipschitz in $\e$, $u$ we obtain
	\begin{equation}
		\abs{v_i},\,\abs{w_i}=O(1)\left(\abs{u_x}+\abs{u_{xx}}+\abs{u_x}^2\right).
	\end{equation}
	From Corollary \ref{coro2.2}, we deduce the estimates on $\|(v_i,w_i)(t,\cdot)\|_{L^1}$ and 
	 $\|(v_i,w_i)(t,\cdot)\|_{L^\infty}$.
	Differentiating \eqref{La:v-w:u-x-u-xx} we obtain
	\begin{align}
		&\frac{\pa\La}{\pa u}u_x+\frac{\pa\La}{\pa(v,w)}(v_x,w_x)\\
		&=\Big(u_{xx},B(u)u_{xxx}+u_x\bullet B(u)u_{xx}+u_{xx}\bullet B(u)u_{x}+u_x\otimes u_x: D^2B(u) u_x\\
		&\quad\quad\quad\quad\quad-A(u)u_{xx}-u_x\bullet A(u)u_x\Big).
	\end{align}
	Using \eqref{estimate:La-u} and uniform invertibility of $\frac{\pa\La}{\pa(v,w)}$ we obtain
	\begin{equation}
		(v_{x},w_{x})=O(1)\left(\abs{u_{xx}}+\abs{u_{xxx}}+\abs{u_x}\abs{u_{xx}}+\abs{u_x}^3+\abs{u_x}^2+\abs{u_x}(\abs{v}+\abs{w})\right).
	\end{equation}
	Applying again the Corollary \ref{coro2.2} and the previous estimates on $v_i,w_i$, we deduce the estimates on  $\|(v_i,w_i)_x (t,\cdot)\|_{L^1}$ and 
	 $\|(v_i,w_i)_x(t,\cdot)\|_{L^\infty}$.
	
%\end{proof}

%\begin{proof}
%	From \eqref{La:v-w:u-x-u-xx} and using invertibility of $\La$ we obtain
%	\begin{equation}
%		\abs{v_1},\,\abs{w_i}=\mathcal{O}(1)\left(\abs{u_x}+\abs{u_{xx}}+\kappa_{B,1}\abs{u_x}^2\right).
%	\end{equation}
%	Differentiating \eqref{La:v-w:u-x-u-xx} with respect to $x$ we obtain
%	\begin{align}
%		&\frac{\pa\La}{\pa u}u_x+\frac{\pa\La}{\pa(v,w)}(v_x,w_x)\\
%		&=\Big(u_{xx},B(u)u_{xxx}+2u_x\bullet B(u)u_{xx}+(u_x\otimes u_x):D^2B(u)  u_x\\
%		&\quad\quad\quad\quad\quad-A(u)u_{xx}-u_x\bullet A(u)u_x\Big).
%	\end{align}
%	Using \eqref{estimate:La-u} and uniform invertibility of $\frac{\pa\La}{\pa(v,w)}$ we obtain
%	\begin{equation}
%		(v_{x},w_{x})=\mathcal{O}(1)\left(\abs{u_{xx}}+\abs{u_{xxx}}+\abs{u_x}\abs{u_{xx}}+\abs{u_x}^3+\abs{u_x}^2+\abs{u_x}(\abs{v}+\abs{w})\right).
%	\end{equation}
	Similarly, differentiating \eqref{La:v-w:u-x-u-xx} with respect to $t$ it gives 
	\begin{align*}
		&\frac{\pa\La}{\pa u}u_t+\frac{\pa\La}{\pa(v,w)}(v_t,w_t)\\
		&=\Big(u_{xt},B(u)u_{txx}+u_t\bullet B(u)u_{xx}+u_x\bullet B(u)u_{tx}+u_{tx}\bullet B(u) u_x+(u_t\otimes u_x):D^2B(u)  u_x\\
		&\hspace{5cm}-A(u)u_{tx}-u_t\bullet A(u)u_x\Big).
	\end{align*}
	Using \eqref{estimate:La-u} and uniform invertibility of $\frac{\pa\La}{\pa(v,w)}$ it yields
	\begin{align*}
		(v_{t},w_{t})&=O(1)\Big(\abs{u_{xt}}+\abs{u_{txx}}+\abs{u_t}\abs{u_{xx}}+\abs{u_x}\abs{u_{tx}}+\abs{u_t}\abs{u_x}^2\\
		&\hspace{5cm}+\abs{u_x}\abs{u_t}+\abs{u_t}(\abs{v}+\abs{w})\Big).
	\end{align*}
	Applying again the Corollary \ref{coro2.2} and the previous estimates on $v_i,w_i$, we deduce the estimates on  $\|(v_{i,t},w_{i,t}) (t,\cdot)\|_{L^1}$ and 
	 $\|(v_{i,t},w_{i,t})(t,\cdot)\|_{L^\infty}$.
\end{proof}
\begin{lemma}\label{lemme6.6}
	Suppose that $u$ satisfies the assumptions of Proposition \ref{prop:parabolic} and further we assume that \eqref{assumption:u-L1} holds for larger interval $[0,T]$ for some $T>\hat{t}$. Then for $t\in[\hat{t},T]$, the variables $v_i,w_i$ satisfy the following bounds
	\begin{align}
		\norm{ v_{i,xx}(t)}_{L^1},\,\norm{w_{i,xx}(t)}_{L^1}&=O(1)\de_0^3,\quad \norm{ v_{i,xx}(t)}_{L^\f},\,\norm{w_{i,xx}(t)}_{L^\f}=O(1)\de_0^4,	\label{estimate-v-w-xx}\\
		\norm{ v_{i,tx}(t)}_{L^1},\,\norm{w_{i,tx}(t)}_{L^1}&=O(1)\de^{3}_0,\quad	\norm{ v_{i,tx}(t)}_{L^\f},\,\norm{w_{i,tx}(t)}_{L^\f}=O(1)\de^{4}_0, \label{estimate-v-w-tx}
	\end{align}  
	for $1\leq i\leq n$. 
\end{lemma}

\begin{proof}[Sketch of the proof:]
	Let us define the vectors $r^{VV}_i$, $r^{VW}_i$, $r^{WW}_i$ and $r^{WV}_i$ as below.
	\begin{align}
		r_i^{VV}&= \tilde{r}_i+v_i\pa_{v_i}\bar{v}_i\xi_i\tilde{r}_{i,v}-\xi_i^\p\bar{v}_i \frac{w_i}{v_i}\tilde{r}_{i,v}+ \theta_i^\p\frac{w_i}{v_i} \tilde{r}_{i,\si}+ \sum\limits_{j\neq i} v_j\pa_{v_i}\bar{v}_j\xi_j\tilde{r}_{j,v},\nonumber\\
		r^{VW}_{i}&=\xi_i^\p\bar{v}_i \tilde{r}_{i,v}-\theta_i^\p \tilde{r}_{i,\si},\\
		r_i^{WW}&=\tilde{r}_i+\xi_i^\p\bar{v}_i\left(\frac{w_i}{v_i}-\la_i^*\right) \tilde{r}_{i,v}-\theta_i^\p\left(\frac{w_i}{v_i}-\la_i^*\right) \tilde{r}_{i,\si},\nonumber\\
		r^{WV}_{i}&=-\la_i^*\tilde{r}_i+(w_i-\la_i^*v_i)\pa_{v_i}\bar{v}_i\xi_i\tilde{r}_{i,v}-\frac{w_i}{v_i}\left(\frac{w_i}{v_i}-\la_i^*\right) [\xi_i^\p\bar{v}_i\tilde{r}_{i,v}-\theta_i^\p\tilde{r}_{i,\si}]\nonumber\\
		&+\sum\limits_{j\neq i} v_j\pa_{v_i}\bar{v}_j \xi_j\left(\frac{w_j}{v_j}-\la_j^*\right)\tilde{r}_{j,v}.
		\label{defriVVs}
	\end{align}
	Consider the matrices $	{\cal Q}^{VV}, {\cal Q}^{VW}, {\cal Q}^{WW},  {\cal Q}^{WV}$ defined as
	\begin{align*}
		{\cal Q}^{VV} &=\left[r^{VV}_1 r^{VV}_2\cdots r^{VV}_n\right] ,\quad\quad {\cal Q}^{VW} =\left[r^{VW}_1 r^{VW}_2\cdots r^{VW}_n\right] ,\\
		{\cal Q}^{WW}&=\left[r^{WW}_1 r^{WW}_2\cdots r^{WW}_n\right] ,\quad\quad {\cal Q}^{WV}=\left[r^{WV}_1 r^{WV}_2\cdots r^{WV}_n\right].
	\end{align*}
	By taking twice derivatives in $x$ on $u_x$ and $u_t$ in \eqref{eqn-v-i} and \eqref{eqn-w-i} respectively, we can obtain
	\begin{align}
		\begin{pmatrix}
			u_{xxx}\\
			u_{txx}
		\end{pmatrix}&=\begin{pmatrix}
			{\cal Q}^{VV} & {\cal Q}^{VW}\\
			{\cal Q}^{WV} & {\cal Q}^{WW}
		\end{pmatrix}\begin{pmatrix}
			v_{i,xx}\\
			w_{i,xx}
		\end{pmatrix}+\mbox{lower order terms}.
	\end{align}
	The lower order terms are in quadratic form and since ${\cal Q}^{VW}$ is small and ${\cal Q}^{VV}$, ${\cal Q}^{WW}$ are invertible for small $(v,w)$, we can then invert the matrix $\begin{pmatrix}
		{\cal Q}^{VV} & {\cal Q}^{VW}\\
		{\cal Q}^{WV} & {\cal Q}^{WW}
	\end{pmatrix}$.	 Hence we can conclude
	\begin{equation}
		\norm{v_{i,xx}}_{L^1},\norm{w_{i,xx}}_{L^1}=O(1)\de_0^3\mbox{ and }\norm{v_{i,xx}}_{L^\f},\norm{w_{i,xx}}_{L^\f}=\mathcal{O}(1)\de_0^4.
	\end{equation}
	Proof of the estimate \eqref{estimate-v-w-tx} follows from a similar argument. We omit here. For detailed proof we refer to Appendix \ref{appendix:higher-regularity}. 
\end{proof}

\subsection{Some coupling identities}
\label{section6.1}
We now consider a global smooth solution $u=u(t,x)$ of \eqref{eqn-main} satisfying $\lim_{x\rightarrow-\infty}u(t,x)=u^*$ for $t\geq 0$ and
which satisfies on the time interval $[\hat{t},T]$
\begin{equation}
\|u_x(t)\|_{L^1}\leq\delta_0,
\label{supercle}
\end{equation}
for $T>\hat{t}.$
In view of Lemma \ref{lemme6.5} 
the decomposition  \eqref{eqn-v-i}-\eqref{eqn-w-i} is well defined, furthermore from the Lemma \ref{lemme6.61} we deduce that the component $v_i$, $w_i$ are smooth for $(t,x)\in[\hat{t},T]\times\R$. We wish now to establish a coupling relation between $w_i$ and $\mu_i v_{i,x}-(\tilde{\la}_i-\la_i^*)v_i$. We recall that $\tilde{\la}_i$ is defined in \eqref{deflatilde}.\\
We can now derive from  \eqref{eqn-v-i}- \eqref{eqn-w-i}
\begin{align*}
	&B(u)u_{xx}+u_x\bullet B(u)u_x-A(u)u_x\\
	&=\sum\limits_{i}v_{i,x}B(u)\tilde{r}_{i}+\sum\limits_{i}v_{i}B(u)\tilde{r}_{i,x}+\sum\limits_{ij}v_iv_j\tilde{r}_i\cdot DB(u)\tilde{r}_j-\sum\limits_{i}v_{i}A(u)\tilde{r}_{i}\\
	&=\sum\limits_{i}v_{i,x}B(u)\tilde{r}_{i}+\sum\limits_{i}v_{i,x}v_iB(u)\xi_i\pa_{v_i}\bar{v}_i\tilde{r}_{i,v}+\sum\limits_{i}\sum\limits_{j\neq i}v_{j,x}v_i\xi_i\pa_{v_j}\bar{v}_iB(u)\tilde{r}_{i,v}\\
	&+\sum\limits_{i}v_i\left[v_iB(u)\tilde{r}_{i,u}\tilde{r}_i+v_i\tilde{r}_i\cdot DB(u)\tilde{r}_i-A(u)\tilde{r}_i\right]\\
	&+\sum\limits_{i}v_i\left(\frac{w_i}{v_i}\right)_xB(u)[-\theta_i^\p\tilde{r}_{i,\si}+\bar{v}_i\xi_i^\p\tilde{r}_{i,v}]+\sum\limits_{i\neq j}v_iv_j\tilde{r}_i\cdot DB(u)\tilde{r}_j+\sum\limits_{i\neq j}v_iv_jB(u)\tilde{r}_{i,u}\tilde{r}_j,
\end{align*}
or equivalently,
\begin{align*}
	&B(u)u_{xx}+u_x\bullet B(u)u_x-A(u)u_x=\sum\limits_{i}v_{i,x}B(u)\tilde{r}_{i}+\sum\limits_{i}v_{i,x}v_iB(u)\xi_i\pa_{v_i}\bar{v}_i\tilde{r}_{i,v}\\
	&+\sum\limits_{i}\sum\limits_{j\neq i}v_{j,x}v_i\xi_i\pa_{v_j}\bar{v}_iB(u)\tilde{r}_{i,v}+\sum\limits_{i}v_i\left[\xi_i\bar{v}_iB(u)\tilde{r}_{i,u}\tilde{r}_i+\xi_i\bar{v}_i\tilde{r}_i\cdot DB(u)\tilde{r}_i-A(u)\tilde{r}_i\right]\\
	&+\sum\limits_{i}(v_i-\xi_i\bar{v}_i)v_i\left[B(u)\tilde{r}_{i,u}\tilde{r}_i+\tilde{r}_i\cdot DB(u)\tilde{r}_i\right]+\sum\limits_{i}v_i\left(\frac{w_i}{v_i}\right)_xB(u)[-\theta_i^\p\tilde{r}_{i,\si}+\bar{v}_i\xi_i^\p\tilde{r}_{i,v}]\\
	&+\sum\limits_{i\neq j}v_iv_j\tilde{r}_i\cdot DB(u)\tilde{r}_j+\sum\limits_{i\neq j}v_iv_jB(u)\tilde{r}_{i,u}\tilde{r}_j.
\end{align*}
From the key equality \eqref{identity-2}, we deduce that
\begin{equation*}
	\xi_i\bar{v}_iB\tilde{r}_{i,u}\tilde{r}_i+\xi_i\bar{v}_i\tilde{r}_i\cdot DB(u)\tilde{r}_i-A(u)\tilde{r}_i=-\si_i\tilde{r}_i-\mu_i^{-1}(\tilde{\la}_i-\si_i)B[\tilde{r}_i+\xi_i\bar{v}_i\tilde{r}_{i,v}].
\end{equation*}
We get then
\begin{align*}
	&B(u)u_{xx}+u_x\bullet B(u)u_x-A(u)u_x=\sum\limits_{i}v_{i,x}B(u)\tilde{r}_{i}+\sum\limits_{i}v_{i,x}v_iB(u)\xi_i\pa_{v_i}\bar{v}_i\tilde{r}_{i,v}\\
	&+\sum\limits_{i}\sum\limits_{j\neq i}v_{j,x}v_i\xi_i\pa_{v_j}\bar{v}_iB(u)\tilde{r}_{i,v}-\sum\limits_{i}v_i\si_{i}\tilde{r}_i-\sum\limits_{i}v_i\mu_i^{-1}(\tilde{\la}_i-\si_i)B[\tilde{r}_i+\xi_i\bar{v}_i\tilde{r}_{i,v}]\\
	&+\sum\limits_{i}(v_i-\xi_i\bar{v}_i)v_i\left[B(u)\tilde{r}_{i,u}\tilde{r}_i+\tilde{r}_i\cdot DB(u)\tilde{r}_i\right]+\sum\limits_{i}v_i\left(\frac{w_i}{v_i}\right)_xB(u)[-\theta_i^\p\tilde{r}_{i,\si}+\bar{v}_i\xi_i^\p\tilde{r}_{i,v}]\\
	&+\sum\limits_{i\neq j}v_iv_j\tilde{r}_i\cdot DB(u)\tilde{r}_j+\sum\limits_{i\neq j}v_iv_jB(u)\tilde{r}_{i,u}\tilde{r}_j,
\end{align*}
or equivalently,
\begin{align}
	&B(u)u_{xx}+u_x\bullet B(u)u_x-A(u)u_x=\sum\limits_{i}(v_{i,x}-\mu_i^{-1}(\tilde{\la}_i-\si_i)v_i)B(u)[\tilde{r}_{i}+(v_i\xi_i\pa_{v_i}\bar{v}_i )\tilde{r}_{i,v}]\nonumber\\
	&+\sum\limits_{i}\sum\limits_{j\neq i}v_{j,x}v_i\xi_i\pa_{v_j}\bar{v}_iB(u)\tilde{r}_{i,v}+\sum\limits_{i}v_i\left(\frac{w_i}{v_i}\right)_xB(u)[-\theta_i^\p\tilde{r}_{i,\si}+\bar{v}_i\xi_i^\p\tilde{r}_{i,v}]\nonumber\\
	&+\sum\limits_{i}\mu_i^{-1}(\tilde{\la}_i-\si_i)v_i\xi_i[(\pa_{v_i}\bar{v}_i)v_i-\bar{v}_i]B\tilde{r}_{i,v}-\sum\limits_{i}v_i\si_i\tilde{r}_i\nonumber\\
	&+\sum\limits_{i}(v_i-\xi_i\bar{v}_i)v_i\left[B(u)\tilde{r}_{i,u}\tilde{r}_i+\tilde{r}_i\cdot DB(u)\tilde{r}_i\right]+\sum\limits_{i\neq j}v_iv_j\tilde{r}_i\cdot DB(u)\tilde{r}_j+\sum\limits_{i\neq j}v_iv_jB(u)\tilde{r}_{i,u}\tilde{r}_j.\label{superimpu}
\end{align}
Since $u_t=B(u)u_{xx}+u_x\bullet B(u)u_x-A(u)u_x$ we obtain
\begin{align}
	&\sum\limits_{i}(w_i-\la_i^*v_i)\tilde{r}_i=\sum\limits_{i}(v_{i,x}-\mu_i^{-1}(\tilde{\la}_i-\si_i)v_i)B(u)[\tilde{r}_{i}+v_i\xi_i\pa_{v_i}\bar{v}_i \tilde{r}_{i,v}]\nonumber\\
	&+\sum\limits_{i}\sum\limits_{j\neq i}\xi_iv_{j,x}v_i\pa_{v_j}\bar{v}_iB(u)\tilde{r}_{i,v}+\sum\limits_{i}v_i\left(\frac{w_i}{v_i}\right)_xB(u)[-\theta_i^\p\tilde{r}_{i,\si}+\bar{v}_i\xi_i^\p\tilde{r}_{i,v}]\nonumber\\
	&+\sum\limits_{i}\mu_i^{-1}(\tilde{\la}_i-\si_i)v_i\xi_i[(\pa_{v_i}\bar{v}_i)v_i-\bar{v}_i]B\tilde{r}_{i,v}-\sum\limits_{i}v_i\si_i\tilde{r}_i\nonumber\\
	&+\sum\limits_{i}(v_i-\xi_i\bar{v}_i)v_i\left[B(u)\tilde{r}_{i,u}\tilde{r}_i+\tilde{r}_i\cdot DB(u)\tilde{r}_i\right]+\sum\limits_{i\neq j}v_iv_j\tilde{r}_i\cdot DB(u)\tilde{r}_j+\sum\limits_{i\neq j}v_iv_jB(u)\tilde{r}_{i,u}\tilde{r}_j.\label{ntych}
\end{align}
Note that from the definition of $\tilde{r}_i $ in \eqref{5.36} we deduce that
\begin{align*}
	&\langle l_i, \tilde{r}_i \rangle =1,\quad \langle l_i, \tilde{r}_j \rangle =\psi_{ji}(u,\xi_j\bar{v}_j,\si_j),\quad  \langle l_i, B\tilde{r}_i \rangle =\mu_i(u),\quad \langle l_i, B\tilde{r}_j \rangle =\mu_i(u)\psi_{ji}(u,\xi_j\bar{v}_j,\si_j),\\
	&\langle l_i, B(u)\tilde{r}_{i,v}\rangle=0,  \langle l_i, B(u)\tilde{r}_{j,v}\rangle=\mu_i(u)\psi_{ji,v}(u,\xi_j\bar{v}_j,\si_j),\\
	&\langle l_i, B(u)[-\theta_i^\p\tilde{r}_{i,\si}+\bar{v}_i\xi_i^\p\tilde{r}_{i,v}]\rangle =0, \langle l_i, B(u)[-\theta_j^\p\tilde{r}_{j,\si}+\bar{v}_j\xi_j^\p\tilde{r}_{j,v}]\rangle =\mu_i(u)[-\theta_j^\p\psi_{ji,\si}+\bar{v}_j\xi_j^\p\psi_{ji,v}]
\end{align*}
Therefore taking the scalar product of \eqref{ntych} with $l_i$ we get
\begin{align*}
	&w_i-\la_i^*v_i+\sum\limits_{j\neq i}\psi_{ji}(u,\xi_j\bar{v}_j,\si_j)(w_j-\la_j^*v_j)=(\mu_i v_{i,x}-(\tilde{\la}_i-\si_i)v_i)\\
	&+\sum\limits_{j\neq i}(v_{j,x}-\mu_j^{-1}(\tilde{\la}_j-\si_j)v_j)\mu_i[\psi_{ji}+v_j\xi_j\pa_{v_j}\bar{v}_j \psi_{ji,v}]+\sum\limits_{j\neq i}\sum\limits_{k\neq j}\xi_jv_{k,x}v_j\pa_{v_k}\bar{v}_j\mu_i \psi_{ji,v}\\
	&+\sum\limits_{j\neq i}v_j\left(\frac{w_j}{v_j}\right)_x\mu_i [-\theta_j^\p\psi_{ji,\si}+\bar{v}_j\xi_j^\p\psi_{ji,v}]\\
	&+\sum\limits_{j\neq i}\mu_j^{-1}(\tilde{\la}_j-\si_j)v_j\xi_j[(\pa_{v_j}\bar{v}_j)v_j-\bar{v}_j]\mu_i\psi_{ji,v}-v_i\si_i-\sum\limits_{j\neq i}v_j\si_j\psi_{ji}\\
	&+\sum\limits_{k}v_k^2(1-\xi_k\chi_k^\ep\prod\limits_{j\neq k}\eta_j)\langle l_i, \left[B(u)\tilde{r}_{k,u}\tilde{r}_k+\tilde{r}_k\cdot DB(u)\tilde{r}_k\right]\rangle \\
	&+\sum\limits_{k\neq j}v_kv_j\langle l_i,\tilde{r}_k\cdot DB(u)\tilde{r}_j+B(u)\tilde{r}_{k,u}\tilde{r}_j\rangle,
\end{align*}
or equivalently, 
\begin{align}
	&\mu_i v_{i,x}-(\tilde{\la}_i-\la_i^*)v_i-w_i =\sum\limits_{j\neq i}\psi_{ji}(u,\xi_j\bar{v}_j,\si_j)(w_j-\la_j^*v_j)\nonumber\\
	&-\sum\limits_{j\neq i}(v_{j,x}-\mu_j^{-1}(\tilde{\la}_j-\si_j)v_j)\mu_i[\psi_{ji}(u,\xi_j\bar{v}_j,\si_j)+v_j\xi_j\pa_{v_j}\bar{v}_j \psi_{ji,v}]\nonumber\\
	&-\sum\limits_{j\neq i}\sum\limits_{k\neq j}\xi_jv_{k,x}v_j\pa_{v_k}\bar{v}_j\mu_i \psi_{ji,v}-\sum\limits_{j\neq i}v_j\left(\frac{w_j}{v_j}\right)_x\mu_i [-\theta_j^\p\psi_{ji,\si}+\bar{v}_j\xi_j^\p\psi_{ji,v}] \nonumber\\
	&-\sum\limits_{j\neq i}\mu_j^{-1}(\tilde{\la}_j-\si_j)v_j\xi_j[(\pa_{v_j}\bar{v}_j)v_j-\bar{v}_j]\mu_i\psi_{ji,v}+\sum\limits_{j\neq i}v_j\si_j\psi_{ji} \nonumber\\
	&-\sum\limits_{k}v_k^2(1-\xi_k\chi_k^\ep\eta_k)\langle l_i, \left[B(u)\tilde{r}_{k,u}\tilde{r}_k+\tilde{r}_k\cdot DB(u)\tilde{r}_k\right]\rangle \nonumber\\
	&-\sum\limits_{k\neq j}v_kv_j\langle l_i,\tilde{r}_k\cdot DB(u)\tilde{r}_j+B(u)\tilde{r}_{k,u}\tilde{r}_j\rangle.         \label{eqn-v-i-x-1a}
\end{align}
We denote $\tilde{v}_j=
\pa_{v_j}\bar{v}_j v_j-\bar{v}_j$. Note that
\begin{equation}\label{form-of-tilde-v-i}
\tilde{v}_j=2N\frac{v_j^{2N}}{\e}(\chi')^\e_j\hat{v}_j-\chi^\e_j\sum_{l\ne j}\prod_{k\ne l,j}\eta\left(\frac{v_k^2}{v_j}\right)\eta^\p\left(\frac{v_l^2}{v_j}\right)v_l^2.	
\end{equation}
We observe since $\psi_{ji}(u,\xi_j\bar{v}_j,\si_j)+v_j\xi_j\pa_{v_j}\bar{v}_j \psi_{ji,v}=O(1)|\xi_j v_j \rho_j^\e|$, $\psi_{ji,\sig}=O(1)\xi_j \bar{v}_j$,  $\psi_{ji}=O(1)\xi_j \bar{v}_j$ and because $\mbox{supp}\xi$ is included in $\{x\in\R,\theta(x)=1\}$, that we have
\begin{align}
	&\mu_i v_{i,x}-(\tilde{\la}_i-\la_i^*)v_i-w_i  \nonumber\\
	&=-\sum\limits_{j\neq i}(\mu_j v_{j,x}-(\tilde{\la}_j-\la_j^*)v_j-w_j)\frac{\mu_i}{\mu_j}[\psi_{ji}(u,\xi_j\bar{v}_j,\si_j)+v_j\xi_j\pa_{v_j}\bar{v}_j \psi_{ji,v}]\nonumber\\
	&-\sum\limits_{j\neq i}\sum\limits_{k\neq j}\xi_jv_{k,x}v_j\pa_{v_k}\bar{v}_j\mu_i \psi_{ji,v}-\sum\limits_{j\neq i}v_j\left(\frac{w_j}{v_j}\right)_x\mu_i [-\psi_{ji,\si}+\bar{v}_j\xi_j^\p\psi_{ji,v}] \nonumber\\
	&-\sum\limits_{j\neq i}\mu_j^{-1}(\tilde{\la}_j-\si_j)v_j\xi_j \tilde{v}_j\mu_i\psi_{ji,v} -\sum\limits_{k}v_k^2(1-\xi_k\chi_k^\ep\eta_k)\langle l_i, \left[B(u)\tilde{r}_{k,u}\tilde{r}_k+\tilde{r}_k\cdot DB(u)\tilde{r}_k\right]\rangle \nonumber\\
	&-\sum\limits_{j}\sum\limits_{k\neq j}v_kv_j\langle l_i,\tilde{r}_k\cdot DB(u)\tilde{r}_j+B(u)\tilde{r}_{k,u}\tilde{r}_j\rangle.         \label{eqn-v-i-x-1}
\end{align}
Note that from the Lemma \ref{estimimpo1} and \eqref{form-of-tilde-v-i}, we have for $k\ne j$
\begin{equation}
	\abs{\pa_{v_k}\bar{v}_j}=O(1)\abs{v_k}  \Delta^\ep_{jk}
	%\mathbbm{1}_{\mbox{\tiny supp}(\eta_{jk}^\p)}
	\;\mbox{and}\;\;
	\tilde{v}_j=O(1)\sum_{l\ne j}v_l^2  \Delta^\ep_{jl}+O(1) (\chi')^\e_j\hat{v}_j.
	\label{6.tech}
\end{equation}
Hence we get the following lemma from \eqref{eqn-v-i-x-1} and \eqref{6.tech}.
\begin{lemma}
	\label{lemme6.8}
	For $1\leq i \leq n$ we have
	\begin{align}
		\mu_i v_{i,x}&=(\tilde{\la}_i-\la_i^*)v_i+w_i+O(1)\sum\limits_{j\neq i}\xi_j\rho_j^\e \abs{(\mu_j v_{j,x}-(\tilde{\la}_j-\la_j^*)v_j-w_j)v_j} \nonumber\\
		&+O(1)\sum\limits_{j\neq i}\sum\limits_{k\neq j}\left[\abs{v_kv_{k,x}v_j}+\abs{v^2_kv_j}\right]\xi_j\Delta^\e_{jk}%\mathbbm{1}_{\mbox{\tiny supp}(\eta_{jk}^\p)}\nonumber\\
		+O(1)\sum\limits_{j\neq i}v^2_j\left(\frac{w_j}{v_j}\right)_x\mathfrak{A}_j\rho_j^\e
		\nonumber\\
		&+O(1)\sum\limits_{k}\sum\limits_{l\neq k}\abs{v_kv_l}+O(1)\sum\limits_{k}v_k^2(1-\xi_k\chi_k^\ep\eta_k)+O(1)\sum\limits_{j\neq i}v_j\xi_j (\chi')^\e_j\hat{v}_j.
	\end{align}
	It gives also
	\begin{align}
		\mu_i v_{i,x}&=(\tilde{\la}_i-\la_i^*)v_i+w_i+O(1)\sum\limits_{j\neq i}\xi_j\rho_j^\e \abs{(\mu_j v_{j,x}-(\tilde{\la}_j-\la_j^*)v_j-w_j)v_j} \nonumber\\
		&+O(1)\sum\limits_{j\neq i}\sum\limits_{k\neq j}\left[\abs{v_kv_{k,x}v_j}+\abs{v^2_kv_j}\right]\xi_j\Delta^\e_{jk}%\mathbbm{1}_{\mbox{\tiny supp}(\eta_{jk}^\p)}\nonumber\\
		+O(1)\sum\limits_{j\neq i}v^2_j\left(\frac{w_j}{v_j}\right)_x\mathfrak{A}_j\rho_j^\e\nonumber\\
		&+O(1)\sum\limits_{k}\sum\limits_{l\neq k}\abs{v_kv_l}+O(1)\sum\limits_{k}v_k^2(1-\xi_k\chi_k^\ep\eta_k)+O(1)\e^{\frac{1}{N}}.
	\end{align}
\end{lemma}
From the Lemma \ref{lemme6.5} we know that $v_j=O(1)\de_0^2$ for small enough $\de_0>0$ then we can obtain
\begin{align}
	&\mu_i v_{i,x}-(\tilde{\la}_i-\la_i^*)v_i-w_i\nonumber\\
	&=O(1)\sum\limits_{j=1}^{n}\sum\limits_{k\neq j}\left[\abs{v_kv_{k,x}v_j}+\abs{v^2_kv_j}\right]\xi_j \Delta^\e_{jk} +O(1)\sum\limits_{j=1}^{n}v^2_j\left(\frac{w_j}{v_j}\right)_x\mathfrak{A}_j\rho^\e_j\nonumber\\
	&+ O(1)\sum\limits_{k=1}^{n}\sum\limits_{l\neq k}\abs{v_kv_l}+O(1)\sum\limits_{k=1}^{n}v_k^2(1-\xi_k\chi_k^\ep\eta_k)+O(1)\sum\limits_{j}v_j\xi_j (\chi')^\e_j\hat{v}_j. \label{6.54}
\end{align}
From the definition \eqref{deflatilde}
\begin{equation}
\tilde{\la}_i= \la_i-v_i\langle B\tilde{r}_{i,u}\tilde{r}_i+\tilde{r}_i\cdot DB(u)\tilde{r}_i ,l_i(u)\rangle,
\end{equation}
we deduce from Lemma \ref{lemme6.5} and \eqref{supercle} that
\begin{equation}
\tilde{\la}_i-\la_i^*=O(1)\delta_0. 
\label{estimla}
\end{equation}
We introduce now the following terms that we will describe later
\begin{align}
	&\La_i^1=\sum\limits_{j\neq i}\left(\abs{v_j}+\abs{v_{j,x}}\right)( \abs{v_i}+\abs{v_{i,x}}+\abs{v_{i,xx}}+|w_i|+|w_{i,x}|+|w_{i,xx}|)\nonumber\\
	&+\sum\limits_{j\neq i}\left(\abs{w_j}+\abs{w_{j,x}}\right)( \abs{v_i}+\abs{v_{i,x}}+\abs{v_{i,xx}}+|w_i|+|w_{i,x}|+|w_{i,xx}|)\nonumber\\
	&+\sum\limits_{j\neq i}\left(\abs{v_{j,xx}}+\abs{w_{j,xx}}\right)( \abs{v_i}+\abs{v_{i,x}}+\abs{v_{i,xx}}+|w_i|+|w_{i,x}|+|w_{i,xx}|)\nonumber\\
	&\La_i^2=\sum\limits_{j\neq i}(\abs{v_{j,xxx}}+\abs{w_{j,xxx}})(\abs{w_j}+\abs{v_j})\abs{v_i},\label{def:La_i-1}\\
	&\La_i^3= \abs{v_i}\left(\frac{w_i}{v_i}\right)_x^2\mathfrak{A}_i  \mathbbm{1}_{\{v_i^{2N}\geq\e\}},\nonumber\\
	&\La_i^4=\abs{w_{i,x}v_i-w_iv_{i,x}},\;\La_i^5=\abs{w_{i,xx}v_i-w_iv_{i,xx}},\nonumber\\
	&\La_i^6=v_{i,x}^2\mathbbm{1}_{\{|\frac{w_i}{v_i}|\geq \frac{\de_1}{2}\}},\;\La_i^{6,1}=w_{i,x}^2\mathbbm{1}_{\{|\frac{w_i}{v_i}|\geq \frac{\de_1}{2}\}}.\nonumber
\end{align} 
We observe using \eqref{6.54} that if $\abs{\frac{w_i}{v_i}}\leq 3\delta_1$ then
\begin{equation}
v_{i,x}=O(1)v_i+O(1)\sum_j(\La_j^1+\La_j^4)+O(1)\sum_{j\ne i}|v_j|^2.
\label{6.54bis}
\end{equation}
Similarly we deduce that if $\abs{\frac{w_i}{v_i}}\geq\frac{\delta_1}{2}$ then for $\delta_0$ sufficiently small in term of $\delta_1$ we get from \eqref{estimla} and  \eqref{6.54}
\begin{equation}\label{6.45}
	w_i, v_i=O(1)v_{i,x}+O(1)\sum_j(\Lambda_j^1+\Lambda_j^4)
	+O(1)\sum\limits_{k\ne i}^{n}v_k^2(1-\xi_k\chi_k^\ep\eta_k)+O(1)\sum\limits_{j\ne i}v_j\xi_j (\chi')^\e_j\hat{v}_j.
\end{equation}
It implies also that
\begin{align}
&\mu_i v_{i,x}-(\tilde{\la}_i-\la_i^*)v_i-w_i\nonumber\\
&=O(1)\sum_j(\La_j^1+\La_j^4)+O(1)\sum\limits_{k=1}^{n}v_k^2(1-\eta_k \chi_k^\ep)\mathbbm{1}_{\{|\frac{w_k}{v_k}|\leq\frac{\delta_1}{2}\}}+O(1)\sum\limits_{k=1}^{n}v_k^2\mathbbm{1}_{\{|\frac{w_k}{v_k}|\geq\frac{\delta_1}{2}\}}\nonumber\\
&+O(1)\sum\limits_{j}v_j\xi_j (\chi')^\e_j\hat{v}_j.\label{ngtech}
\end{align}
Applying \eqref{6.45}, we deduce that
\begin{align}
v_k^2\mathbbm{1}_{\{|\frac{w_k}{v_k}|\geq\frac{\delta_1}{2}\}}&=O(1)\mathbbm{1}_{\{|\frac{w_k}{v_k}|\geq\frac{\delta_1}{2}\}}v_{k,x} v_k+O(1)\sum_l(\La_l^1+\La_l^4)%+O(1)\sum\limits_{j}v_k v_j\xi_j (\chi')^\e_j\hat{v}_j
\nonumber\\
&=O(1)\mathbbm{1}_{\{|\frac{w_k}{v_k}|\geq\frac{\delta_1}{2}\}}v_{k,x}^2+O(1)\sum_l(\La_l^1+\La_l^4).%+O(1)\sum\limits_{j}v_k v_j\xi_j (\chi')^\e_j\hat{v}_j
%\nonumber\\
%&+O(1)\sum\limits_{j}v_{k,x} v_j\xi_j (\chi')^\e_j\hat{v}_j.
\label{ngtech1}
\end{align}
Combining \eqref{ngtech}, \eqref{ngtech1} and Lemma \ref{lemme6.5}, we deduce that
\begin{align}
&\mu_i v_{i,x}-(\tilde{\la}_i-\la_i^*)v_i-w_i\nonumber\\
&=O(1)\sum_j(\La_j^1+\La_j^4+\La_j^6)+O(1)\sum\limits_{k=1}^{n}v_k^2(1-\eta_k \chi_k^\ep)\mathbbm{1}_{\{|\frac{w_k}{v_k}|\leq\frac{\delta_1}{2}\}}+O(1)\sum\limits_{j}v_j\xi_j (\chi')^\e_j\hat{v}_j.\label{ngtech2}
\end{align}
We have now $v_k^2(1-\eta_k \chi_k^\ep)=v_k^2(1-\eta_k)+\eta_k v_k^2(1- \chi_k^\ep)$ and using the fact that $1-\eta_k\ne 0$ if it exists $j\ne k$ such that $v_j^2\geq \frac{3}{4}|v_k|$ (it implies in particular that $0\leq 1-\eta_k \leq \sum_{j\ne k}\mathbbm{1}_{\{|\frac{v_j^2}{v_k}|\geq \frac{3}{4}\}}$) we deduce using Lemma \ref{lemme6.5}
\begin{equation}
O(1)\sum\limits_{k=1}^{n}v_k^2(1-\eta_k \chi_k^\ep)=O(1)\sum\limits_{k=1}^{n} \eta_k v_k^2(1- \chi_k^\ep)+O(1)\sum_j \La^1_j.\label{ngtech3}
\end{equation}
Combining \eqref{ngtech2} and \eqref{ngtech3}, we get
\begin{align}
&\mu_i v_{i,x}-(\tilde{\la}_i-\la_i^*)v_i-w_i\nonumber\\
&=O(1)\sum_j(\La_j^1+\La_j^4+\La_j^6)%+O(1)\sum\limits_{k=1}^{n}v_k^2(1-\eta_k \chi_k^\ep)\mathbbm{1}_{\{|\frac{w_k}{v_k}|\leq\frac{\delta_1}{2}\}}
+O(1)\sum\limits_{j}v_j\xi_j (\chi')^\e_j\hat{v}_j+O(1)\sum\limits_{k=1}^{n} \eta_k v_k^2(1- \chi_k^\ep).\label{ngtech4}
\end{align}
From \eqref{6.54} and Lemmas \ref{lemme6.5}, \ref{lemme6.6}, we deduce also that
\begin{align}
	&\mathfrak{A}_i \rho_i v_{i,x}, \mathfrak{A}_i \rho_i (\mu_i v_{i,x}-(\tilde{\la}_i-\la_i^*)v_i-w_i)\nonumber\\
	&=O(1)v_i \mathfrak{A}_i \rho_i +O(1)\sum\limits_{j\ne i}v^2_j\left(\frac{w_j}{v_j}\right)_x \mathfrak{A}_j\rho^\e_j \mathfrak{A}_i \rho_i +O(1)\sum\limits_{j}v_j\xi_j (\chi')^\e_j\hat{v}_j \mathfrak{A}_i \rho_i \nonumber\\
	&=O(1)v_i \mathfrak{A}_i \rho_i +O(1)\sum\limits_{j\ne i}v^2_j\left(\frac{w_j}{v_j}\right)_x \mathfrak{A}_j\rho^\e_j \mathfrak{A}_i \rho_i=O(1)v_i \mathfrak{A}_i \rho_i +O(1)\sum_j\La_j^4 ,\nonumber\\
	&=O(1)\sqrt{|v_i|} \mathfrak{A}_i \rho_i . \label{ngtech5}
\end{align}
Similarly from \eqref{ngtech4} and \eqref{estimla} we get
\begin{align}\label{ngtech6}
	\mathbbm{1}_{\{|\frac{w_i}{v_i}|\geq\frac{\delta_1}{2}\}} v_i=&O(1)v_{i,x}\mathbbm{1}_{\{|\frac{w_i}{v_i}|\geq\frac{\delta_1}{2}\}} 
	+O(1)\sum_j(\La_j^1+\La_j^4+\La_j^6)%+O(1)\sum\limits_{k=1}^{n}v_k^2(1-\eta_k \chi_k^\ep)\mathbbm{1}_{\{|\frac{w_k}{v_k}|\leq\frac{\delta_1}{2}\}}
+O(1)\sum\limits_{j}v_j\xi_j (\chi')^\e_j\hat{v}_j\nonumber\\
&+O(1)\sum\limits_{k=1}^{n} \eta_k v_k^2(1- \chi_k^\ep).
\end{align}
\section{Estimation of the remainder terms}
Under the same assumption \eqref{supercle} as in the previous section, we are going to prove that $v_i$ and $w_i$ satisfy the following equations when $t\in[ \hat{t},T]$ for $T>\hat{t}$
\begin{align}
	&v_{i,t}+(\tilde{\la}_iv_i)_x-(\mu_iv_{i,x})_x=\widetilde{\phi}_i,\label{eqn-v-i-1}\\
	&w_{i,t}+(\tilde{\la}_iw_i)_x-(\mu_iw_{i,x})_x=\widetilde{\psi}_i,\label{eqn-w-i-1}
\end{align}
for some remainder terms $\widetilde{\phi}_i$ and $\widetilde{\psi}_i$, $1\leq i\leq n$. To show the uniform $BV$ estimates of $u_x$ on $[\hat{t},T]$ it is enough to prove the bound of the remainder terms $\widetilde{\phi}_i$ and $\widetilde{\psi}_i$ in $L^1([\hat{t},T],L^1(\R))$. We will see later that we can choose $T=+\infty$, furthermore
we can write  $\widetilde{\phi}_i$ and $\widetilde{\psi}_i$ as follows
\begin{equation}
\begin{aligned}
	&\widetilde{\phi}_i=\phi_i+\sum\limits_{j\neq i}(\mu_i-\mu_j)b_{ij}\left(w_{j,x}-\frac{w_j}{v_j}v_{j,x}\right)_x,\\
	&\widetilde{\psi}_i=\psi_i+\sum\limits_{j\neq i}(\mu_i-\mu_j)\hat{b}_{ij}\left(w_{j,x}-\frac{w_j}{v_j}v_{j,x}\right)_x,\label{deftilde}
\end{aligned}
\end{equation}
for some appropriate functions $b_{ij}, \hat{b}_{ij}$ and $\phi_i,\psi_i$ that we will express explicitly in the sequel. Our main aim is to estimate the terms $\widetilde{\phi_i}$ and $\widetilde{\psi_i}$, we will see that these terms can be bound by the terms defined in \eqref{def:La_i-1}.

Note that in $\widetilde{\phi_i}$ and $\widetilde{\psi_i}$ the second order terms of the type $\La_i^5=w_{i,xx}v_i-w_iv_{i,xx}$ appears due the fact that eigenvalues of $B(u)$ can be different. These terms does not appear in the analysis of Bianchini and Bressan \cite{BB-vv-lim-ann-math}. In order to deal with this term, we will introduce new unknowns  some effective flux $z_i$ and $\hat{z}_i$ that we will estimate also in $L^{1}([\hat{t},T],L^1(\R))$ norm.

\subsection{Calculation for the remainder term $\phi_i$}\label{section:phi-i}
From \eqref{superimpu} we have
\begin{align*}
	&B(u)u_{xx}+u_x\bullet B(u)u_x-A(u)u_x=\sum\limits_{i}(v_{i,x}-\mu_i^{-1}(\tilde{\la}_i-\si_i)v_i)B(u)[\tilde{r}_{i}+v_i\xi_i\pa_{v_i}\bar{v}_i\tilde{r}_{i,v}]\\
	&+\sum\limits_{i}\sum\limits_{j\neq i}v_{j,x}v_i\xi_i\pa_{v_j}\bar{v}_iB(u)\tilde{r}_{i,v}+\sum\limits_{i}\left(w_{i,x}-\frac{w_i}{v_i}v_{i,x}\right)B(u)[-\theta_i^\p\tilde{r}_{i,\si}+\bar{v}_i\xi_i^\p\tilde{r}_{i,v}]\\
	&+\sum\limits_{i}\mu_i^{-1}(\tilde{\la}_i-\si_i)v_i\xi_i[(\pa_{v_i}\bar{v}_i)v_i-\bar{v}_i]B\tilde{r}_{i,v}-\sum\limits_{i}v_i\si_i\tilde{r}_i\\
	&+\sum\limits_{i}(v_i-\xi_i\bar{v}_i)v_i\left[B(u)\tilde{r}_{i,u}\tilde{r}_i+\tilde{r}_i\cdot DB(u)\tilde{r}_i\right]+\sum\limits_{i\neq j}v_iv_j\tilde{r}_i\cdot DB(u)\tilde{r}_j+\sum\limits_{i\neq j}v_iv_jB(u)\tilde{r}_{i,u}\tilde{r}_j,
\end{align*}
with $\bar{v}_i=\chi_i^\e\hat{v}_i$.
After a rearrangement, we can write
\begin{align*}
	&B(u)u_{xx}+u_x\bullet B(u)u_x-A(u)u_x\\
	&=\sum\limits_{i}(\mu_iv_{i,x}-\tilde{\la}_iv_i)[\tilde{r}_{i}+(v_i\xi_i\pa_{v_i}\bar{v}_i-(w_i/v_i)\xi^\p_i\bar{v}_i) \tilde{r}_{i,v}+(w_i/v_i)\theta^\p_i\tilde{r}_{i,\si}]\\
	&+\sum\limits_{i}(\mu_iw_{i,x}-\tilde{\la}_iw_i)(\xi_i^\p\bar{v}_i\tilde{r}_{i,v}-\theta^\p_i\tilde{r}_{i,\si})-\sum\limits_{i}v_i\si_i\tilde{r}_i\\
	&+\sum\limits_{i}\si_iv_i[\tilde{r}_{i}+v_i\xi_i\pa_{v_i}\bar{v}_i \tilde{r}_{i,v}]+\sum\limits_{i}\sum\limits_{j\neq i}v_{j,x}v_i\xi_i\pa_{v_j}\bar{v}_iB(u)\tilde{r}_{i,v}\\
	&+\sum\limits_{i}(v_{i,x}-\mu_i^{-1}\tilde{\la}_iv_i)(B-\mu_iI_n)[\tilde{r}_{i}+(v_i\xi_i\pa_{v_i}\bar{v}_i-(w_i/v_i)\xi^\p_i\bar{v}_i)  \tilde{r}_{i,v}+(w_i/v_i)\theta^\p_i\tilde{r}_{i,\si}]\\
	&+\sum\limits_{i}(w_{i,x}-\mu_i^{-1}\tilde{\la}_iw_i)(B-\mu_iI_n)[\xi_i^\p\bar{v}_i\tilde{r}_{i,v}-\theta^\p_i\tilde{r}_{i,\si}]\\
	&+\sum\limits_{i}\mu_i^{-1}\si_iv_i(B-\mu_iI_n)[\tilde{r}_{i}+v_i\xi_i\pa_{v_i}\bar{v}_i \tilde{r}_{i,v}]+
	\sum\limits_{i}\mu_i^{-1}(\tilde{\la}_i-\si_i)v_i\xi_i[(\pa_{v_i}\bar{v}_i)v_i-\bar{v}_i]B\tilde{r}_{i,v}\\
	&+\sum\limits_{i}(v_i-\xi_i\bar{v}_i)v_i\left[B(u)\tilde{r}_{i,u}\tilde{r}_i+\tilde{r}_i\cdot DB(u)\tilde{r}_i\right]+\sum\limits_{i\neq j}v_iv_j\tilde{r}_i\cdot DB(u)\tilde{r}_j+\sum\limits_{i\neq j}v_iv_jB(u)\tilde{r}_{i,u}\tilde{r}_j.	
\end{align*}
Thus, we write 
\begin{align}
	&B(u)u_{xx}+u_x\bullet B(u)u_x-A(u)u_x=\sum_{i=1}^9 J_i,\label{egaliteimp}
\end{align}
where $\{J_i\}_{i=1}^{9}$ are defined as follows
\begin{align*}
	J_1&=\sum\limits_{i}(\mu_iv_{i,x}-\tilde{\la}_iv_i)[\tilde{r}_{i}+(v_i\xi_i\pa_{v_i}\bar{v}_i-(w_i/v_i)\xi^\p_i\bar{v}_i)\tilde{r}_{i,v}+(w_i/v_i)\theta^\p_i\tilde{r}_{i,\si}],\\
	J_2&=\sum\limits_{i}(\mu_iw_{i,x}-\tilde{\la}_iw_i)[\xi_i^\p\bar{v}_i\tilde{r}_{i,v}-\theta^\p_i\tilde{r}_{i,\si}],\\
	J_3&=\sum\limits_{i}(\la_i^*-\theta_i)v^2_i\xi_i\pa_{v_i}\bar{v}_i \tilde{r}_{i,v},\\
	J_4	&=\sum\limits_{i}\sum\limits_{j\neq i}v_{j,x}v_i\xi_i\pa_{v_j}\bar{v}_iB(u)\tilde{r}_{i,v},\\
	J_5&=\sum\limits_{i}(v_{i,x}-\mu_i^{-1}(\tilde{\la}_i-\la_i^*+\theta_i)v_i)(B-\mu_iI_n)[\tilde{r}_{i}+v_i\xi_i\pa_{v_i}\bar{v}_i \tilde{r}_{i,v}],\\
	J_6&=\sum\limits_{i}(w_{i,x}-(w_i/v_i)v_{i,x})(B-\mu_iI_n)[\xi_i^\p\bar{v}_i\tilde{r}_{i,v}-\theta^\p_i\tilde{r}_{i,\si}],\\
	J_7&=\sum\limits_{i}\mu_i^{-1}(\tilde{\la}_i-\si_i)v_i\xi_i[(\pa_{v_i}\bar{v}_i)v_i-\bar{v}_i]B\tilde{r}_{i,v},\\
	J_8&=\sum\limits_{i}(v_i-\xi_i\bar{v}_i)v_i\left[B(u)\tilde{r}_{i,u}\tilde{r}_i+\tilde{r}_i\cdot DB(u)\tilde{r}_i\right],\\
	J_{9}&=\sum\limits_{i\neq j}v_iv_j\tilde{r}_i\cdot DB(u)\tilde{r}_j+\sum\limits_{i\neq j}v_iv_jB(u)\tilde{r}_{i,u}\tilde{r}_j.
\end{align*}
Next we calculate derivatives of each $J_i$'s with respect to $x$. 
\begin{align*}
	&J_{1,x}=\sum\limits_{i}(\mu_iv_{i,x}-\tilde{\la}_iv_i)_x\Big[\tilde{r}_{i}+(v_i\xi_i\pa_{v_i}\bar{v}_i-(w_i/v_i)\xi_i^\p\bar{v}_i) \tilde{r}_{i,v}+(w_i/v_i)\theta^\p_i\tilde{r}_{i,\si}\Big]\\
	&+\sum\limits_{i}\tilde{r}_{i,u}\tilde{r}_i\left[(\mu_iv_{i,x}-\tilde{\la}_iv_i)v_i\right] +\sum\limits_{i}\sum\limits_{j\neq i}\tilde{r}_{i,u}\tilde{r}_j\left[(\mu_iv_{i,x}-\tilde{\la}_iv_i)v_j\right] \\
	&+\sum\limits_{i}\tilde{r}_{i,v}(2\xi_i\pa_{v_i}\bar{v}_i+v_i\xi_i\pa_{v_iv_i}\bar{v}_i-\xi_i^\p\pa_{v_i}\bar{v}_i \frac{w_i}{v_i})\left[(\mu_iv_{i,x}-\tilde{\la}_iv_i)v_{i,x}\right]\\
	&+\sum\limits_{i}\tilde{r}_{i,v}\left[v_i\xi_i^\p\pa_{v_i}\bar{v}_i-\xi^{\p\p}_i\frac{w_i}{v_i}\bar{v}_i\right](\mu_iv_{i,x}-\tilde{\la}_iv_i)\left(\frac{w_i}{v_i}\right)_x\\
	&+\sum\limits_{i}\sum\limits_{j\neq i}\tilde{r}_{i,v}(\xi_i\pa_{v_j}\bar{v}_i+v_i\xi_i\pa_{v_iv_j}\bar{v}_i-(w_i/v_i)\xi_i^\p\pa_{v_j}\bar{v}_i)\left[(\mu_iv_{i,x}-\tilde{\la}_iv_i)v_{j,x}\right] \\ 
	&-\sum\limits_{i}\tilde{r}_{i,\si}\left[(\mu_iv_{i,x}-\tilde{\la}_iv_i)\theta_i^\p\left(\frac{w_i}{v_i}\right)_x\right] \\
	&+\sum\limits_{i}\tilde{r}_{i,uv}\tilde{r}_i (v_i\xi_i\pa_{v_i}\bar{v}_i-(w_i/v_i)\xi_i^\p\bar{v}_i)\left[(\mu_iv_{i,x}-\tilde{\la}_iv_i)v_i\right] \\
	&+\sum\limits_{i}\sum\limits_{j\neq i}\tilde{r}_{i,uv}\tilde{r}_j(v_i\xi_i\pa_{v_i}\bar{v}_i-(w_i/v_i)\xi_i^\p\bar{v}_i)\left[(\mu_iv_{i,x}-\tilde{\la}_iv_i)v_j\right] \\
	&+\sum\limits_{i}\tilde{r}_{i,vv}(\xi_i\pa_{v_i}\bar{v}_i)(v_i\xi_i\pa_{v_i}\bar{v}_i-(w_i/v_i)\xi_i^\p\bar{v}_i)\left[(\mu_iv_{i,x}-\tilde{\la}_iv_i)v_{i,x}\right]\\
	&+\sum\limits_{i}\tilde{r}_{i,vv}\xi_i^\p\bar{v}_i\left(\frac{w_i}{v_i}\right)_x(v_i\xi_i\pa_{v_i}\bar{v}_i-(w_i/v_i)\xi_i^\p\bar{v}_i)\left[(\mu_iv_{i,x}-\tilde{\la}_iv_i)\right]\\
	& +\sum\limits_{i}\sum\limits_{j\neq i}\tilde{r}_{i,vv}\xi_i\pa_{v_j}\bar{v}_i(v_i\xi_i\pa_{v_i}\bar{v}_i-(w_i/v_i)\xi_i^\p\bar{v}_i)\left[(\mu_iv_{i,x}-\tilde{\la}_iv_i)v_{j,x}\right] \\
	&-\sum\limits_{i}\tilde{r}_{i,v\si}(v_i\xi_i\pa_{v_i}\bar{v}_i-(w_i/v_i)\xi_i^\p\bar{v}_i)\left[(\mu_iv_{i,x}-\tilde{\la}_iv_i)\theta_i^\p \left(\frac{w_i}{v_i}\right)_x\right] \\
	&+\sum\limits_{i}\tilde{r}_{i,u\si}\tilde{r}_i\left[(\mu_iv_{i,x}-\tilde{\la}_iv_i)w_i\theta^\p_i\right] +\sum\limits_{i}\sum\limits_{j\neq i}\tilde{r}_{i,u\si}\tilde{r}_j\left[(\mu_iv_{i,x}-\tilde{\la}_iv_i)v_j(w_i/v_i)\theta^\p_i\right] \\
	&+\sum\limits_{i}\tilde{r}_{i,v\si}\xi_i\pa_{v_i}\bar{v}_i\left[(\mu_iv_{i,x}-\tilde{\la}_iv_i)v_{i,x}(w_i/v_i)\theta^\p_i\right]+\sum\limits_{i}\sum\limits_{j\neq i}\tilde{r}_{i,v\si}\xi_i\pa_{v_j}\bar{v}_i\left[(\mu_iv_{i,x}-\tilde{\la}_iv_i)v_{j,x}\frac{w_i}{v_i}\theta^\p_i\right] \\
	&+\sum\limits_{i}\tilde{r}_{i,v\si}\xi^\p_i\bar{v}_i\left(\frac{w_i}{v_i}\right)_x\left[(\mu_iv_{i,x}-\tilde{\la}_iv_i)(w_i/v_i)\theta^\p_i\right]-\sum\limits_{i}\tilde{r}_{i,\si\si}\left[(\mu_iv_{i,x}-\tilde{\la}_iv_i)(\theta_i^\p)^2\left(\frac{w_i}{v_i}\right)_x(w_i/v_i)\right] \\
	&+\sum\limits_{i}\tilde{r}_{i,\si}(\mu_iv_{i,x}-\tilde{\la}_iv_i)\left[\theta_i^\p\left(\frac{w_i}{v_i}\right)_x+\theta_i^{\p\p}\frac{w_i}{v_i}\left(\frac{w_i}{v_i}\right)_x\right]=\sum\limits_{k=1}^{20}\sum\limits_{i}\al_{i}^{1,k}.
\end{align*}
Taking derivative of $J_2$ we obtain,
\begin{align*}
	&J_{2,x}=\sum\limits_{i}(\mu_iw_{i,x}- \tilde{\la}_iw_i)_x\left[\xi_i^\p\bar{v}_i\tilde{r}_{i,v}-\theta^\p_i\tilde{r}_{i,\si}\right]+\sum\limits_{i}(\mu_iw_{i,x}- \tilde{\la}_iw_i)v_i\xi_i^\p\bar{v}_i\tilde{r}_{i,uv}\tilde{r}_i\\
	&+\sum\limits_{i}\sum\limits_{j\neq i}(\mu_iw_{i,x}- \tilde{\la}_iw_i)v_j\xi_i^\p\bar{v}_i\tilde{r}_{i,uv}\tilde{r}_j+\sum\limits_{i}(\mu_iw_{i,x}- \tilde{\la}_iw_i)v_{i,x}\xi_i\xi_i^\p\bar{v}_i\pa_{v_i}\bar{v}_i\tilde{r}_{i,vv}\\
	&+\sum\limits_{i}\sum\limits_{j\neq i}(\mu_iw_{i,x}- \tilde{\la}_iw_i)v_{j,x}\xi_i\xi_i^\p\bar{v}_i\pa_{v_j}\bar{v}_i\tilde{r}_{i,vv}
	+\sum\limits_{i}(\mu_iw_{i,x}- \tilde{\la}_iw_i)(\xi_i^\p)^2\left(\frac{w_i}{v_i}\right)_x\bar{v}^2_i\tilde{r}_{i,vv}\\
	&-\sum\limits_{i}(\mu_iw_{i,x}- \tilde{\la}_iw_i)\xi_i^\p \bar{v}_{i}\left(\frac{w_i}{v_i}\right)_x\theta^\p_i\tilde{r}_{i,v\si}+\sum\limits_{i}(\mu_iw_{i,x}- \tilde{\la}_iw_i)v_{i,x}\xi_i^\p\pa_{v_i}\bar{v}_i\tilde{r}_{i,v}\\
	&+\sum\limits_{i}\sum\limits_{j\neq i}(\mu_iw_{i,x}- \tilde{\la}_iw_i)v_{j,x}\xi_i^\p\pa_{v_j}\bar{v}_i\tilde{r}_{i,v}+\sum\limits_{i}(\mu_iw_{i,x}- \tilde{\la}_iw_i)\xi_i^{\p\p}\left(\frac{w_i}{v_i}\right)_x\bar{v}_i\tilde{r}_{i,v}\\
	&-\sum\limits_{i}(\mu_iw_{i,x}- \tilde{\la}_iw_i)v_i\theta^\p_i\tilde{r}_{i,u\si}\tilde{r}_i-\sum\limits_{i}\sum\limits_{j\neq i}(\mu_iw_{i,x}- \tilde{\la}_iw_i)v_j\theta^\p_i\tilde{r}_{i,u\si}\tilde{r}_j\\
	&-\sum\limits_{i}(\mu_iw_{i,x}- \tilde{\la}_iw_i)\xi_iv_{i,x}\pa_{v_i}\bar{v}_i\theta^\p_i\tilde{r}_{i,v\si}-\sum\limits_{i}\sum\limits_{j\neq i}(\mu_iw_{i,x}- \tilde{\la}_iw_i)\xi_iv_{j,x}\pa_{v_j}\bar{v}_i\theta^\p_i\tilde{r}_{i,v\si}\\
	&-\sum\limits_{i}(\mu_iw_{i,x}- \tilde{\la}_iw_i)\xi_i^\p \bar{v}_{i}\left(\frac{w_i}{v_i}\right)_x\theta^\p_i\tilde{r}_{i,v\si}+\sum\limits_{i}(\mu_iw_{i,x}- \tilde{\la}_iw_i)\left(\frac{w_i}{v_i}\right)_x(\theta^\p_i)^2\tilde{r}_{i,\si\si}\\
	&-\sum\limits_{i}(\mu_iw_{i,x}- \tilde{\la}_iw_i)\left(\frac{w_i}{v_i}\right)_x\theta^{\p\p}_i\tilde{r}_{i,\si}=:\sum\limits_{k=1}^{17}\sum\limits_{i}\al_{i}^{2,k}.
\end{align*}
Next we calculate the derivative of $J_3$,
\begin{align*}
	&J_{3,x}=2\sum\limits_{i}(\la_i^*-\theta_i)\xi_iv_iv_{i,x} \pa_{v_i}\bar{v}_i \tilde{r}_{i,v}-\sum\limits_{i}\theta_i^\p\left(\frac{w_i}{v_i}\right)_x\xi_iv^2_i \pa_{v_i}\bar{v}_i \tilde{r}_{i,v}\\
	&+\sum\limits_{i}(\la_i^*-\theta_i)\xi_iv_{i,x}v^2_i \pa_{v_iv_i}\bar{v}_i \tilde{r}_{i,v}+\sum\limits_{i}\sum\limits_{j\neq i}(\la_i^*-\theta_i)\xi_iv^2_i v_{j,x}\pa_{v_iv_j}\bar{v}_i \tilde{r}_{i,v}\\
	%&+\sum\limits_{i}\sum\limits_{j\neq i}(\la_i^*-\theta_i)\xi_iv^2_i w_{j,x}\pa_{v_iw_j}\bar{v}_i \tilde{r}_{i,v}\\
	&+\sum\limits_{i}(\la_i^*-\theta_i)\xi_iv_i^3\pa_{v_i}\bar{v}_i \tilde{r}_{i,uv}\tilde{r}_i+\sum\limits_{i}\sum\limits_{j\neq i}(\la_i^*-\theta_i)\xi_iv_i^2v_j\pa_{v_i}\bar{v}_i \tilde{r}_{i,uv}\tilde{r}_j\\
	&+\sum\limits_{i}(\la_i^*-\theta_i)\xi^2_iv_i^2v_{i,x}(\pa_{v_i}\bar{v}_i)^2\tilde{r}_{i,vv}+\sum\limits_{i}\sum\limits_{j\neq i}(\la_i^*-\theta_i)\xi^2_iv_i^2v_{j,x}\pa_{v_i}\bar{v}_i\pa_{v_j}\bar{v}_i\tilde{r}_{i,vv}\\
	%&+\sum\limits_{i}\sum\limits_{j\neq i}(\la_i^*-\theta_i)\xi^2_iv_i^2w_{j,x}\pa_{v_i}\bar{v}_i\pa_{w_j}\bar{v}_i\tilde{r}_{i,vv}\\
	&+\sum\limits_{i}(\la_i^*-\theta_i)\xi_i\xi^\p_iv_i^2\left(\frac{w_i}{v_i}\right)_x\bar{v}_i \pa_{v_i}\bar{v}_i \tilde{r}_{i,vv}-\sum\limits_{i}\theta_i^\p(\la_i^*-\theta_i)\xi_iv_i^2\left(\frac{w_i}{v_i}\right)_x\pa_{v_i}\bar{v}_i\tilde{r}_{i,v\si}\\
	&+\sum\limits_{i}(\la_i^*-\theta_i)\xi^\p_iv^2_i \left(\frac{w_i}{v_i}\right)_x\pa_{v_i}\bar{v}_i \tilde{r}_{i,v}=\sum\limits_{k=1}^{11}\sum\limits_{i}\al_{i}^{3,k}.
\end{align*}
For the term $J_4$ we get
\begin{align*}
	&J_{4,x}=\sum\limits_{i}\sum\limits_{j\neq i}B(u)\tilde{r}_{i,v}\pa_{v_j}\bar{v}_i\left[v_{j,xx}v_i\xi_i+v_{j,x}v_{i,x}\xi_i+v_{j,x}v_i\xi^\p_i\left(\frac{w_i}{v_i}\right)_x\right]\\
	&+\sum\limits_{i}\sum\limits_{j\neq i}\sum\limits_{k}B(u)\tilde{r}_{i,v}v_i\pa_{v_jv_k}\bar{v}_i\big[\xi_iv_{j,x}v_{k,x}\big]%+\sum\limits_{i}\sum\limits_{j\neq i}\sum\limits_{k\neq i}B(u)\tilde{r}_{i,v}v_i\pa_{v_jw_k}\bar{v}_i\big[\xi_iv_{j,x}w_{k,x}\big]\\
	+\sum\limits_{i}\sum\limits_{j\neq i}B(u)\tilde{r}_{i,uv}\tilde{r}_i\pa_{v_j}\bar{v}_i\big[v_{j,x}v^2_i\xi_i\big]\\
	&+\sum\limits_{i}\sum\limits_{j\neq i}\sum\limits_{k\neq i}B(u)\tilde{r}_{i,uv}\tilde{r}_k v_i\pa_{v_j}\bar{v}_i\big[v_{j,x}v_k\xi_i\big]
	+\sum\limits_{i}\sum\limits_{j\neq i}\sum\limits_{k}B(u)\tilde{r}_{i,vv}\pa_{v_j}\bar{v}_i\pa_{v_k}\bar{v}_i\big[v_{k,x}v_{j,x}v_i\xi_i^2\big]\\
	%&+\sum\limits_{i}\sum\limits_{j\neq i}\sum\limits_{k\neq i}B(u)\tilde{r}_{i,vv}\pa_{v_j}\bar{v}_i\pa_{w_k}\bar{v}_i\big[w_{k,x}v_{j,x}v_i\xi_i^2\big]\\
	&+\sum\limits_{i}\sum\limits_{j\neq i}B(u)\tilde{r}_{i,vv}\pa_{v_j}\bar{v}_i\left[v_{j,x}v_i\xi_i\xi_i^\p\left(\frac{w_i}{v_i}\right)_x\bar{v}_i\right]-\sum\limits_{i}\sum\limits_{j\neq i} B(u)\tilde{r}_{i,v\si}\pa_{v_j}\bar{v}_i\left[v_{j,x}v_i\xi_i\left(\frac{w_i}{v_i}\right)_x\theta^\p_i\right]\\
	&+\sum\limits_{i}\sum\limits_{j\neq i}v_{j,x}v_i\xi_i\pa_{v_j}\bar{v}_iu_x\cdot DB(u)\tilde{r}_{i,v}=:\sum\limits_{l=1}^{8}\sum\limits_{i}\al_{i}^{4,l}.
\end{align*}
%For the term $J_5$ we get
%\begin{align*}
%	J_{5,x}&=\sum\limits_{i}\sum\limits_{j\neq i}B(u)\tilde{r}_{i,v}\pa_{w_j}\bar{v}_i\left[w_{j,xx}v_i\xi_i+w_{j,x}v_{i,x}\xi_i+w_{j,x}v_i\xi^\p_i\left(\frac{w_i}{v_i}\right)_x\right]\\
%	&+\sum\limits_{i}\sum\limits_{j\neq i}\sum\limits_{k}B(u)\tilde{r}_{i,v}v_i\pa_{w_jv_k}\bar{v}_i\big[\xi_iw_{j,x}v_{k,x}\big]\\
%	&+\sum\limits_{i}\sum\limits_{j\neq i}\sum\limits_{k\neq i}B(u)\tilde{r}_{i,v}v_i\pa_{w_jw_k}\bar{v}_i\big[\xi_iw_{j,x}w_{k,x}\big]\\
%	&+\sum\limits_{i}\sum\limits_{j\neq i}B(u)\tilde{r}_{i,uv}\tilde{r}_i\pa_{w_j}\bar{v}_i\big[w_{j,x}v^2_i\xi_i\big]\\
%	&+\sum\limits_{i}\sum\limits_{j\neq i}\sum\limits_{k\neq i}B(u)\tilde{r}_{i,uv}\tilde{r}_k v_i\pa_{w_j}\bar{v}_i\big[w_{j,x}v_k\xi_i\big]\\
%	&+\sum\limits_{i}\sum\limits_{j\neq i}\sum\limits_{k}B(u)\tilde{r}_{i,vv}\pa_{w_j}\bar{v}_i\pa_{v_k}\bar{v}_i\big[v_{k,x}w_{j,x}v_i\xi_i^2\big]\\
%	&+\sum\limits_{i}\sum\limits_{j\neq i}\sum\limits_{k\neq i}B(u)\tilde{r}_{i,vv}\pa_{w_j}\bar{v}_i\pa_{w_k}\bar{v}_i\big[w_{k,x}w_{j,x}v_i\xi_i^2\big]\\
%	&+\sum\limits_{i}\sum\limits_{j\neq i}B(u)\tilde{r}_{i,vv}\pa_{w_j}\bar{v}_i\left[w_{j,x}v_i\xi_i\xi^\p_{i}\left(\frac{w_i}{v_i}\right)_x\bar{v}_i\right]\\
%	&-\sum\limits_{i}\sum\limits_{j\neq i} B(u)\tilde{r}_{i,v\si}\pa_{w_j}\bar{v}_i\left[w_{j,x}v_i\xi_i\left(\frac{w_i}{v_i}\right)_x\theta^\p_i\right]\\
%	&+\sum\limits_{i}\sum\limits_{j\neq i}w_{j,x}v_i\xi_i\pa_{w_j}\bar{v}_iu_x\cdot DB(u)\tilde{r}_{i,v},\\
%	&=:\sum\limits_{l=1}^{10}\sum\limits_{i}\al_{i}^{5,l}.
%\end{align*}
Next we calculate the derivative of term $J_5$.
\begin{align*}
	&J_{5,x}=\sum\limits_{i}(v_{i,x}-\mu_i^{-1}(\tilde{\la}_i-\la_i^*+\theta_i)v_i)_x(B-\mu_iI_n)[\tilde{r}_{i}+v_i\xi_i\pa_{v_i}\bar{v}_i\tilde{r}_{i,v}]\\
	&+\sum\limits_{i}(v_{i,x}-\mu_i^{-1}(\tilde{\la}_i-\la^*_i+\theta_i)v_i)v_i\tilde{r}_i\cdot D(B-\mu_iI_n)[\tilde{r}_{i}+v_i\xi_i\pa_{v_i}\bar{v}_i\tilde{r}_{i,v}]\\
	&+\sum\limits_{i}\sum\limits_{j\neq i}(v_{i,x}-\mu_i^{-1}(\tilde{\la}_i-\la^*_i+\theta_i)v_i)v_j\tilde{r}_j\cdot D(B-\mu_iI_n)[\tilde{r}_{i}+v_i\xi_i\pa_{v_i}\bar{v}_i\tilde{r}_{i,v}]\\
	&+\sum\limits_{i}(v_{i,x}-\mu_i^{-1}(\tilde{\la}_i-\la_i^*+\theta_i)v_i)v_{i}(B-\mu_iI_n)[\tilde{r}_{i,u}\tilde{r}_i+v_i\xi_i\pa_{v_i}\bar{v}_i\tilde{r}_{i,uv}\tilde{r}_{i}]\\
	&+\sum\limits_{i}\sum\limits_{j\neq i}(v_{i,x}-\mu_i^{-1}(\tilde{\la}_i-\la_i^*+\theta_i)v_i)v_{j}(B-\mu_iI_n)[\tilde{r}_{i,u}\tilde{r}_j+v_i\xi_i\pa_{v_i}\bar{v}_i\tilde{r}_{i,uv}\tilde{r}_{j}]\\
	&+\sum\limits_{i}(v_{i,x}-\mu_i^{-1}(\tilde{\la}_i-\la_i^*+\theta_i)v_i)v_{i,x}\xi_i(B-\mu_iI_n)[(2\pa_{v_i}\bar{v}_i+v_i\pa_{v_iv_i}\bar{v}_i)\tilde{r}_{i,v}+\xi_iv_i(\pa_{v_i}\bar{v}_i)^2\tilde{r}_{i,vv}]\\
	&+\sum\limits_{i}\sum\limits_{j\neq i}(v_{i,x}-\mu_i^{-1}(\tilde{\la}_i-\la_i^*+\theta_i)v_i)v_{j,x}\xi_i(B-\mu_iI_n)[(\pa_{v_j}\bar{v}_i+v_i\pa_{v_iv_j}\bar{v}_i)\tilde{r}_{i,v}+\xi_iv_i\pa_{v_i}\bar{v}_i\pa_{v_j}\bar{v}_i\tilde{r}_{i,vv}]\\
	%&+\sum\limits_{i}\sum\limits_{j\neq i}(v_{i,x}-\mu_i^{-1}(\tilde{\la}_i-\la_i^*+\theta_i)v_i)w_{j,x}\xi_i(B-\mu_iI_n)[(\pa_{w_j}\bar{v}_i+v_i\pa_{v_iw_j}\bar{v}_i)\tilde{r}_{i,v}+\xi_iv_i\pa_{v_i}\bar{v}_i\pa_{w_j}\bar{v}_i\tilde{r}_{i,vv}]\\
	&+\sum\limits_{i}(v_{i,x}-\mu_i^{-1}(\tilde{\la}_i-\la_i^*+\theta_i)v_i)\xi^\p_i\left(\frac{w_i}{v_i}\right)_x(B-\mu_iI_n)[(v_i\pa_{v_i}\bar{v}_i+\bar{v}_i)\tilde{r}_{i,v}+\xi_iv_i\pa_{v_i}\bar{v}_i\bar{v}_i\tilde{r}_{i,vv}]\\
	&-\sum\limits_{i}(v_{i,x}-\mu_i^{-1}(\tilde{\la}_i-\la_i^*+\theta_i)v_i)\theta_i^\p\left(\frac{w_i}{v_i}\right)_x(B-\mu_iI_n)[\tilde{r}_{i,\si}+\xi_iv_i\pa_{v_i}\bar{v}_i\tilde{r}_{i,v\si}]\\
	&=:\sum\limits_{l=1}^{9}\sum\limits_{i}\al_{i}^{5,l}.
\end{align*}
Similarly, we have
\begin{align*}
	J_{6,x}&=\sum\limits_{i}(w_{i,x}-(w_i/v_i)v_{i,x})_x(B-\mu_iI_n)[\xi_i^\p\bar{v}_i\tilde{r}_{i,v}-\theta^\p_i\tilde{r}_{i,\si}]\\
	&+\sum\limits_{i}(w_{i,x}-(w_i/v_i)v_{i,x})\left(\frac{w_i}{v_i}\right)_x(B-\mu_iI_n)[\xi^{\p\p}_i\bar{v}_i\tilde{r}_{i,v}-\theta^{\p\p}_i\tilde{r}_{i,\si}]\\
	&+\sum\limits_{i}(w_{i,x}v_i-w_iv_{i,x})\tilde{r}_i\cdot D(B-\mu_iI_n)[\xi_i^\p\bar{v}_i\tilde{r}_{i,v}-\theta^\p_i\tilde{r}_{i,\si}]\\
	&+\sum\limits_{i}\sum\limits_{j\neq i}(w_{i,x}-(w_i/v_i)v_{i,x})v_j\tilde{r}_j\cdot D(B-\mu_iI_n)[\xi_i^\p\bar{v}_i\tilde{r}_{i,v}-\theta^\p_i\tilde{r}_{i,\si}]\\
	&+\sum\limits_{i}(w_{i,x}v_i-w_iv_{i,x}) (B-\mu_iI_n)[\xi_i^\p\bar{v}_i\tilde{r}_{i,uv}\tilde{r}_i-\theta^\p_i\tilde{r}_{i,u\si}\tilde{r}_i]\\
	&+\sum\limits_{i}\sum\limits_{j\neq i}(w_{i,x}-(w_i/v_i)v_{i,x})v_j(B-\mu_iI_n)[\xi_i^\p\bar{v}_i\tilde{r}_{i,uv}\tilde{r}_j-\theta^\p_i\tilde{r}_{i,u\si}\tilde{r}_j]\\
	&+\sum\limits_{i}(w_{i,x}-(w_i/v_i)v_{i,x})v_{i,x}\xi_i\pa_{v_i}\bar{v}_i(B-\mu_iI_n)[\xi_i^\p\bar{v}_i\tilde{r}_{i,vv}-\theta^\p_i\tilde{r}_{i,v\si}]\\
	&+\sum\limits_{i}\sum\limits_{j\neq i}(w_{i,x}-(w_i/v_i)v_{i,x})v_{j,x}\xi_i\pa_{v_j}\bar{v}_i\theta^\p_i(B-\mu_iI_n)[\xi_i^\p\bar{v}_i\tilde{r}_{i,vv}-\theta^\p_i\tilde{r}_{i,v\si}]\\
	%&+\sum\limits_{i}\sum\limits_{j\neq i}(w_{i,x}-(w_i/v_i)v_{i,x})w_{j,x}\xi_i\pa_{w_j}\bar{v}_i\theta^\p_i(B-\mu_iI_n)[\xi_i^\p\bar{v}_i\tilde{r}_{i,vv}-\theta^\p_i\tilde{r}_{i,v\si}]\\
	&+\sum\limits_{i}(w_{i,x}-(w_i/v_i)v_{i,x})\left(\frac{w_i}{v_i}\right)_x\xi_i^\p\bar{v}_i(B-\mu_iI_n)[\xi_i^\p\bar{v}_i\tilde{r}_{i,vv}-\theta^\p_i\tilde{r}_{i,v\si}]\\
	&-\sum\limits_{i}(w_{i,x}-(w_i/v_i)v_{i,x})\theta^\p_i\left(\frac{w_i}{v_i}\right)_x(B-\mu_iI_n)[\xi_i^\p\bar{v}_i\tilde{r}_{i,v\si}-\theta^\p_i\tilde{r}_{i,\si\si}]\\
	&+\sum\limits_{i}(w_{i,x}-(w_i/v_i)v_{i,x})v_{i,x}(B-\mu_iI_n)\xi_i^\p\pa_{v_i}\bar{v}_i\tilde{r}_{i,v}\\
	&+\sum\limits_{i}\sum\limits_{j\neq i}(w_{i,x}-(w_i/v_i)v_{i,x})v_{j,x}(B-\mu_iI_n)\xi_i^\p\pa_{v_j}\bar{v}_i\tilde{r}_{i,v}
	%&+\sum\limits_{i}\sum\limits_{j\neq i}(w_{i,x}-(w_i/v_i)v_{i,x})w_{j,x}(B-\mu_iI_n)\xi_i^\p\pa_{w_j}\bar{v}_i\tilde{r}_{i,v}\\
	=:\sum\limits_{l=1}^{12}\sum\limits_{i}\al_{i}^{6,l}.
\end{align*}
For the term $J_7$ and $J_8$ we get
\begin{align*}
	&J_{7,x}=\sum\limits_{i}\sum\limits_{k}v_k\tilde{r}_k\cdot (D_u\mu_i^{-1})(\tilde{\la}_i-\si_i)v_i\xi_i[(\pa_{v_i}\bar{v}_i)v_i-\bar{v}_i]B\tilde{r}_{i,v}\\
	&+\sum\limits_{i}\mu_i^{-1}(\tilde{\la}_i-\si_i)\xi^\p_i\left(\frac{w_i}{v_i}\right)_xv_i[(\pa_{v_i}\bar{v}_i)v_i-\bar{v}_i]B\tilde{r}_{i,v}+\sum\limits_{i}\mu_i^{-1}(\tilde{\la}_i-\si_i)\xi_i v_{i,x}[(\pa_{v_i}\bar{v}_i)v_i-\bar{v}_i]B\tilde{r}_{i,v}\\
	&+\sum\limits_{i}\sum\limits_{k}\mu_i^{-1}\tilde{r}_k\cdot D \tilde{\la}_iv_kv_i \xi_i[(\pa_{v_i}\bar{v}_i)v_i-\bar{v}_i]B\tilde{r}_{i,v}+\sum\limits_{i}\mu_i^{-1}\tilde{\la}_{i,v}\pa_{v_i}\bar{v}_iv_{i,x}v_i\xi_i^2[(\pa_{v_i}\bar{v}_i)v_i-\bar{v}_i]B\tilde{r}_{i,v}\\
	&+\sum\limits_{i}\sum\limits_{j\neq i}\mu_i^{-1}\tilde{\la}_{i,v}\pa_{v_j}\bar{v}_iv_{j,x}v_i \xi_i^2[(\pa_{v_i}\bar{v}_i)v_i-\bar{v}_i]B\tilde{r}_{i,v}
	+\sum\limits_{i}\mu_i^{-1}\tilde{\la}_{i,v}\xi_i^\p\left(\frac{w_i}{v_i}\right)_x\bar{v}_i v_i\xi_i[(\pa_{v_i}\bar{v}_i)v_i-\bar{v}_i]B\tilde{r}_{i,v}\\
	&-\sum\limits_{i}\mu_i^{-1}\tilde{\la}_{i,\si}\theta_i^\p\left(\frac{w_i}{v_i}\right)_xv_i\xi_i[(\pa_{v_i}\bar{v}_i)v_i-\bar{v}_i]B\tilde{r}_{i,v}+\sum\limits_{i}\mu_i^{-1}\theta_i^\p\left(\frac{w_i}{v_i}\right)_xv_i\xi_i[(\pa_{v_i}\bar{v}_i)v_i-\bar{v}_i]B\tilde{r}_{i,v}\\
	&+\sum\limits_{i}\mu_i^{-1}(\tilde{\la}_i-\si_i)v^2_iv_{i,x}\xi_i(\pa_{v_iv_i}\bar{v}_i)B\tilde{r}_{i,v}+\sum\limits_{i}\sum\limits_{j\neq i}\mu_i^{-1}(\tilde{\la}_i-\si_i)v_iv_{j,x}\xi_i[(\pa_{v_iv_j}\bar{v}_i)v_i-\pa_{v_j}\bar{v}_i]B\tilde{r}_{i,v}\\
	&+\sum\limits_{i}\mu_i^{-1}(\tilde{\la}_i-\si_i)v_i\xi_i[(\pa_{v_i}\bar{v}_i)v_i-\bar{v}_i]u_x\cdot DB\tilde{r}_{i,v}\\
	&+\sum\limits_{i}\sum\limits_{k}\mu_i^{-1}(\tilde{\la}_i-\si_i)v_i\xi_i v_k[(\pa_{v_i}\bar{v}_i)v_i-\bar{v}_i]B\tilde{r}_{i,uv}\tilde{r}_k\\
	&+\sum\limits_{i} \mu_i^{-1}(\tilde{\la}_i-\si_i)v_iv_{i,x}\xi_i^2\pa_{v_i}\bar{v}_i[(\pa_{v_i}\bar{v}_i)v_i-\bar{v}_i]B\tilde{r}_{i,vv}\\
	&+\sum\limits_{i}\sum\limits_{j\neq i}\mu_i^{-1}(\tilde{\la}_i-\si_i)v_iv_{j,x}\xi_i^2\pa_{v_j}\bar{v}_i[(\pa_{v_i}\bar{v}_i)v_i-\bar{v}_i]B\tilde{r}_{i,vv}\\
	&+\sum\limits_{i}\mu_i^{-1}(\tilde{\la}_i-\si_i)v_i\xi_i\xi_i^\p\bar{v}_i\left(\frac{w_i}{v_i}\right)_x[(\pa_{v_i}\bar{v}_i)v_i-\bar{v}_i]B\tilde{r}_{i,vv}\\
	&-\sum\limits_{i}\mu_i^{-1}(\tilde{\la}_i-\si_i)v_i\theta_i^\p\left(\frac{w_i}{v_i}\right)_x\xi_i[(\pa_{v_i}\bar{v}_i)v_i-\bar{v}_i]B\tilde{r}_{i,v\si}=:\sum\limits_{l=1}^{17}\sum\limits_{i}\al_{i}^{7,l},
\end{align*}
and 
\begin{align*}
	J_{8,x}&=\sum\limits_{i}(2v_i-\xi_i\bar{v}_i-\xi_iv_i\pa_{v_i}\bar{v}_i)v_{i,x}\left[B(u)\tilde{r}_{i,u}\tilde{r}_i+\tilde{r}_i\cdot DB(u)\tilde{r}_i\right]\\
	&-\sum\limits_{i}\xi^\p_iv_i\bar{v}_i\left(\frac{w_i}{v_i}\right)_x\left[B(u)\tilde{r}_{i,u}\tilde{r}_i+\tilde{r}_i\cdot DB(u)\tilde{r}_i\right]\\
	&-\sum\limits_{i}\sum\limits_{j\neq i}\xi_iv_iv_{j,x}\pa_{v_j}\bar{v}_i\left[B(u)\tilde{r}_{i,u}\tilde{r}_i+\tilde{r}_i\cdot DB(u)\tilde{r}_i\right]\\
	&+\sum\limits_{i}(1-\eta_i\xi_i \chi_i^\e)v^2_i\left[B(u)\tilde{r}_{i,u}\tilde{r}_i+\tilde{r}_i\cdot DB(u)\tilde{r}_i\right]_x=:\sum\limits_{l=1}^{4}\sum\limits_{i}\al_{i}^{8,l}.
\end{align*}
Finally, for the term $J_{9}$ we calculate
\begin{align*}
	J_{9,x}&=\sum\limits_{i}\sum\limits_{j\neq i}(v_{i,x}v_j+v_iv_{j,x})\left[\tilde{r}_i\cdot DB(u)\tilde{r}_j+B(u)\tilde{r}_{i,u}\tilde{r}_j\right]\\
	&+\sum\limits_{i}\sum\limits_{j\neq i}v_iv_j\left[\tilde{r}_i\cdot DB(u)\tilde{r}_j+B(u)\tilde{r}_{i,u}\tilde{r}_j\right]_x=:\sum\limits_{l=1}^{2}\sum\limits_{i}\al_{i}^{9,l}.
\end{align*}
Combining $J_{1,x}, J_{2,x}$ and $J_{3,x}$ we obtain (see the Appendix \ref{sectionE} for more details on the computation)
\begin{align}
	&J_{1,x}+J_{2,x}+J_{3,x}%+\sum\limits_{i}\al_{i}^{5,1}
	-\sum\limits_{i}(w_i-\la_i^*v_i)v_{i}\tilde{r}_{i,u}\tilde{r}_i-\sum\limits_{i}\sum\limits_{j\neq i}(w_j-\la_j^*v_j)v_{i}\tilde{r}_{i,u}\tilde{r}_j\nonumber\\
	&=\sum\limits_{i}(\mu_iv_{i,x}-\tilde{\la}_iv_i)_x \left[\tilde{r}_{i}+(v_i\xi_i\pa_{v_i}\bar{v}_i-(w_i/v_i)\xi_i^\p\bar{v}_i) \tilde{r}_{i,v}+(w_i/v_i)\theta^\p_i\tilde{r}_{i,\si}\right]\nonumber\\
	&+\sum\limits_{i}(\mu_iw_{i,x}- \tilde{\la}_iw_i)_x\left[\xi_i^\p\bar{v}_i\tilde{r}_{i,v}-\theta^\p_i\tilde{r}_{i,\si}\right]+\sum\limits_{i}\sum\limits_{j\neq i}\tilde{r}_{i,u}\tilde{r}_j\left[(\mu_iv_{i,x}-\tilde{\la}_iv_i)v_j-(w_j-\la_j^*v_j)v_{i}\right] \nonumber\\
	&+\sum\limits_{i}\tilde{r}_{i,u}\tilde{r}_i\left[(\mu_iv_{i,x}-(\tilde{\la}_i-\la_i^*)v_i-w_i)v_i\right] \nonumber\\
	&+\sum\limits_{i}\tilde{r}_{i,v}(2\xi_i\pa_{v_i}\bar{v}_i+v_i\xi_i\pa_{v_iv_i}\bar{v}_i)\left[(\mu_iv_{i,x}-(\tilde{\la}_i-\la_i^*)v_i-\theta_i v_i)v_{i,x}\right]\nonumber\\
	&+\sum\limits_{i}\tilde{r}_{i,v}\mu_i v_{i,x}\xi^\p_i\left(w_{i,x}-\frac{w_i}{v_i}v_{i,x}\right) \pa_{v_i}\bar{v}_i+\sum\limits_{i}\tilde{r}_{i,v}v_i\xi_i^\p\pa_{v_i}\bar{v}_i(\mu_iv_{i,x}-(\tilde{\la}_i-\la_i^*)v_i-\theta_i v_i)\left(\frac{w_i}{v_i}\right)_x\nonumber\\
	&+\sum\limits_{i}\tilde{r}_{i,v}\mu_i\xi^{\p\p}_i\left(w_{i,x}-\frac{w_i}{v_i}v_{i,x}\right) \bar{v}_i\left(\frac{w_i}{v_i}\right)_x-\sum_i v_i^2\left(\frac{w_i}{v_i}\right)_x\pa_{v_i}\bar{v}_i\tilde{r}_{i,v}\theta'_i\xi_i   \nonumber\\
	&+\sum\limits_{i}\sum\limits_{j\neq i}\tilde{r}_{i,v}(\xi_i\pa_{v_j}\bar{v}_i)\left[(\mu_iv_{i,x}-\tilde{\la}_iv_i)v_{j,x}\right]
	+\sum\limits_{i}\sum\limits_{j\neq i}\tilde{r}_{i,v}v_i\xi_i\pa_{v_iv_j}\bar{v}_i\left[(\mu_iv_{i,x}-(\tilde{\la}_i-\la_i^*)v_i-\theta_i v_i)v_{j,x}\right]\nonumber\\
	& +\sum\limits_{i}\sum\limits_{j\neq i}\tilde{r}_{i,v}\mu_i\xi_i^\p\pa_{v_j}\bar{v}_iv_i\left(\frac{w_i}{v_i}\right)_xv_{j,x}
	+\sum\limits_{i}\tilde{r}_{i,uv}\tilde{r}_i v_i\xi_i\pa_{v_i}\bar{v}_i\left[(\mu_iv_{i,x}-(\tilde{\la}_i-\la_i^*)v_i- \theta_i v_i)v_i\right]\nonumber\\
	&+\sum\limits_{i}\tilde{r}_{i,uv}\tilde{r}_i \mu_i\xi_i^\p\bar{v}_i(w_{i,x}v_i-w_iv_{i,x})+\sum\limits_{i}\sum\limits_{j\neq i}\tilde{r}_{i,uv}\tilde{r}_j(\xi_i\pa_{v_i}\bar{v}_i)v_jv_i(\mu_iv_{i,x}-(\tilde{\la}_i-\la_i^*)v_i-\theta_i v_i)\nonumber\\
	&+\sum\limits_{i}\sum\limits_{j\neq i}\tilde{r}_{i,uv}\tilde{r}_j(w_{i,x}-(w_i/v_i)v_{i,x})\xi_i^\p\bar{v}_iv_j\mu_i+\sum\limits_{i}\tilde{r}_{i,vv}(\xi_i\pa_{v_i}\bar{v}_i)^2v_iv_{i,x}(\mu_iv_{i,x}-(\tilde{\la}_i-\la_i^*)v_i-\theta_i v_i)\nonumber\\
	&+\sum\limits_{i}\tilde{r}_{i,vv}(\xi_i\pa_{v_i}\bar{v}_i)\xi_i^\p\bar{v}_i(w_{i,x}-(w_i/v_i)v_{i,x})v_{i,x}\mu_i\nonumber\\
	&+\sum\limits_{i}\tilde{r}_{i,vv}\xi_i^\p\bar{v}_i\left(\frac{w_i}{v_i}\right)_x(v_i\xi_i\pa_{v_i}\bar{v}_i)\left[(\mu_iv_{i,x}-(\tilde{\la}_i-\la_i^*)v_i-\theta_i v_i)\right]\nonumber\\
	&+\sum\limits_{i}\tilde{r}_{i,vv}(\xi_i^\p)^2\mu_i\bar{v}_i^{2}\left(\frac{w_i}{v_i}\right)_x(w_{i,x}-(w_i/v_i)v_{i,x})\nonumber\\
	&+\sum\limits_{i}\sum\limits_{j\neq i}\tilde{r}_{i,vv}\xi_i^2\pa_{v_j}\bar{v}_i\pa_{v_i}\bar{v}_iv_i(\mu_iv_{i,x}-(\tilde{\la}_i-\la_i^*)v_i-\theta_i v_i)v_{j,x}\nonumber\\
	&+\sum\limits_{i}\sum\limits_{j\neq i}\tilde{r}_{i,vv}\mu_i\xi_i\xi^\p_i\pa_{v_j}\bar{v}_i\bar{v}_i\left(w_{i,x}-\frac{w_i}{v_i}v_{i,x}\right)v_{j,x}\nonumber\\
	&-\sum\limits_{i}\tilde{r}_{i,v\si}v_i\xi_i\pa_{v_i}\bar{v}_i(\mu_iv_{i,x}-(\tilde{\la}_i-\la_i^*)v_i-\theta_i v_i)\theta_i^\p \left(\frac{w_i}{v_i}\right)_x\nonumber\\
	&-2\sum\limits_{i}\tilde{r}_{i,v\si}\xi_i^\p\bar{v}_i\mu_i(w_{i,x}-(w_i/v_i)v_{i,x}) \theta_i^\p \left(\frac{w_i}{v_i}\right)_x-\sum\limits_{i}\sum\limits_{k}\tilde{r}_{i,v\si}\xi_i\pa_{v_k}\bar{v}_iv_{k,x}\theta^\p_i\mu_i\left[w_{i,x}-\frac{w_i}{v_i}v_{i,x}\right]\nonumber\\
	&+\sum\limits_{i}\tilde{r}_{i,u\si}\tilde{r}_i \theta^\p_i \mu_i(w_iv_{i,x}-w_{i,x}v_i)
	-\sum\limits_{i}\sum\limits_{j\neq i}\tilde{r}_{i,u\si}\tilde{r}_j\mu_i \theta^\p_iv_j\left[w_{i,x}-(w_i/v_i)v_{i,x}\right]\nonumber\\
	&+\sum\limits_{i}\tilde{r}_{i,\si\si} \mu_i(\theta_i^\p)^2\left(\frac{w_i}{v_i}\right)_x\left[w_{i,x}-(w_i/v_i)v_{i,x}\right]-\sum\limits_{i}\tilde{r}_{i,\si}\mu_i\theta_i^{\p\p}\left(\frac{w_i}{v_i}\right)_x\left[w_{i,x}-(w_i/v_i)v_{i,x}\right]\label{equaJ1x}\\
	&=:\sum\limits_{i}(\mu_iv_{i,x}-\tilde{\la}_iv_i)_x\left[\tilde{r}_{i}+(v_i\xi_i\pa_{v_i}\bar{v}_i-(w_i/v_i)\xi_i^\p\bar{v}_i) \tilde{r}_{i,v}+(w_i/v_i)\theta^\p_i\tilde{r}_{i,\si}\right]\nonumber\\
	&+\sum\limits_{i}(\mu_iw_{i,x}- \tilde{\la}_iw_i)_x\left[\xi_i^\p\bar{v}_i\tilde{r}_{i,v}-\theta^\p_i\tilde{r}_{i,\si}\right]+\sum\limits_{l=1}^{27}\sum\limits_{i}\widetilde{\al}_{i}^{1,l}.\label{defalpha1}
\end{align}
Now, we set
\begin{equation}
\begin{aligned}
	\hat{r}_i&:=\tilde{r}_{i}+(v_i\xi_i\pa_{v_i}\bar{v}_i-(w_i/v_i)\xi_i^\p\bar{v}_i) \tilde{r}_{i,v}+(w_i/v_i)\theta^\p_i\tilde{r}_{i,\si},\\
	r_i^\dagger&:=\xi_i^\p\bar{v}_i\tilde{r}_{i,v}-\theta^\p_i\tilde{r}_{i,\si}.
	\label{defAo}
\end{aligned}
\end{equation}
Therefore, we get
\begin{align}
	&J_{1,x}+J_{2,x}+J_{3,x}-\sum\limits_{i}(w_i-\la_i^*v_i)v_{i}\tilde{r}_{i,u}\tilde{r}_i-\sum\limits_{i}\sum\limits_{j\neq i}(w_j-\la_j^*v_j)v_{i}\tilde{r}_{i,u}\tilde{r}_j\nonumber\\
	&=\sum\limits_{i}(\mu_iv_{i,x}-\tilde{\la}_iv_i)_x\hat{r}_i+\sum\limits_{i}(\mu_iw_{i,x}- \tilde{\la}_iw_i)_xr_i^\dagger+\sum\limits_{l=1}^{27}\sum\limits_{i}\widetilde{\al}_{i}^{1,l}.
	\label{J1J2J3}
\end{align}
Taking derivative in time on \eqref{eqn-v-i} we get
\begin{align*}
	u_{tx}&=\sum\limits_{i}v_{i,t}\tilde{r}_i+\sum\limits_{i}v_{i}\tilde{r}_{i,t}\\
	&=\sum\limits_{i}v_{i,t}\tilde{r}_i+\sum\limits_{ij}(w_j-\la_j^*v_j)v_{i}\tilde{r}_{i,u}\tilde{r}_j+\sum\limits_{i}\left(v_{i,t}\pa_{v_i}\bar{v_i}+\sum\limits_{j\neq i}v_{j,t}\pa_{v_j}\bar{v_i}\right)\xi_iv_{i}\tilde{r}_{i,v}\\
	&+\sum\limits_{i}\left(w_{i,t}-\frac{w_i}{v_i}v_{i,t}\right)[\xi_i^\p\bar{v}_i\tilde{r}_{i,v}-\theta_i^\p\tilde{r}_{i,\si}].
\end{align*}
After a rearrangement of terms we obtain
\begin{align*}
	u_{tx}&=\sum\limits_{i}v_{i,t}\left[\tilde{r}_i+(\xi_iv_i\pa_{v_i}\bar{v}_i-(w_i/v_i)\xi_i^\p\bar{v}_i)\tilde{r}_{i,v}+\theta_i^\p\frac{w_i}{v_i}\tilde{r}_{i,\si}\right]+\sum\limits_{i}w_{i,t}[\xi^\p\bar{v}_i\tilde{r}_{i,v}-\theta_i^\p\tilde{r}_{i,\si}]\\
	&+\sum\limits_{i}(w_i-\la_i^*v_i)v_{i}\tilde{r}_{i,u}\tilde{r}_i+\sum\limits_{i}\sum\limits_{j\neq i}v_{i,t}v_j\xi_j\pa_{v_i}\bar{v_j}\tilde{r}_{j,v}	+\sum\limits_{i}\sum\limits_{j\neq i}(w_j-\la_j^*v_j)v_{i}\tilde{r}_{i,u}\tilde{r}_j.
\end{align*}
Hence, we may write
\begin{align}
	u_{tx}
	&=\sum\limits_{i}v_{i,t}\hat{r}_i+\sum\limits_{i}w_{i,t}r^\dagger_i+\sum\limits_{i}(w_i-\la_i^*v_i)v_{i}\tilde{r}_{i,u}\tilde{r}_i+\sum\limits_{i}\sum\limits_{j\neq i}(w_j-\la_j^*v_j)v_{i}\tilde{r}_{i,u}\tilde{r}_j\nonumber\\
	&+\sum\limits_{i}\sum\limits_{j\neq i}v_{i,t}v_j\xi_j\pa_{v_i}\bar{v_j}\tilde{r}_{j,v}	\nonumber\\
	&=\sum\limits_{i}v_{i,t}\hat{r}_i+\sum\limits_{i}w_{i,t}r^\dagger_i+\sum\limits_{i}(w_i-\la_i^*v_i)v_{i}\tilde{r}_{i,u}\tilde{r}_i \nonumber\\
	&+\sum\limits_{i}\sum\limits_{j\neq i}(w_j-\la_j^*v_j)v_{i}\tilde{r}_{i,u}\tilde{r}_j+\sum_i\alpha_i^{10,1}.\label{egalitesuperimp}
\end{align}
Since $u_{tx}=[B(u)u_{xx}+u_x\bullet B(u)u_x-A(u)u_x]_x$ we have by combining \eqref{egaliteimp}, \eqref{J1J2J3} and \eqref{egalitesuperimp}
\begin{align}
	&\sum\limits_{i}(v_{i,t}+(\tilde{\la}_iv_i)_x-(\mu_iv_{i,x})_x)\hat{r}_i+\sum\limits_{i}(w_{i,t}+(\tilde{\la}_iw_i)_x-(\mu_iw_{i,x})_x)r^\dagger_i \nonumber\\
	&=\sum\limits_{l=1}^{27}\sum\limits_{i}\widetilde{\al}_{i}^{1,l}+\sum\limits_{k=4}^{10}\sum\limits_{l=1}^{N_k}\sum\limits_{i}\al^{k,l}_i,	\label{cal-phi-1}
\end{align}
where $N_4=8,\,N_5=9,\,N_6=12,\, N_7=17,\,N_8=4,\,N_9=2$ and  $N_{10}=1$. From $\al_i^{6,1}$ we can write by rearranging terms and using \eqref{5.36}
\begin{align*}
	\sum\limits_{i}\al_i^{6,1}
	&=\sum\limits_{i}(w_{i,x}-(w_i/v_i)v_{i,x})_x(B-\mu_iI_n)[\xi_i^\p\bar{v}_i\tilde{r}_{i,v}-\theta^\p_i\tilde{r}_{i,\si}]\\
	&=\sum\limits_{i}(w_{i,x}-(w_i/v_i)v_{i,x})_x\sum\limits_{j\neq i}(\mu_j-\mu_i)[\xi_i^\p\bar{v}_i\psi_{ij,v} -\theta^\p_i\psi_{ij,\si} ]r_j\\
	&=\sum\limits_{i} \sum\limits_{j\neq i}(\mu_i-\mu_j)\left[\xi_j^\p\bar{v}_j\psi_{ji,v} -\theta^\p_j\psi_{ji,\si} \right]\left(w_{j,x}-\frac{w_j}{v_j}v_{j,x}\right)_xr_i\\
	&=\sum\limits_{i} \sum\limits_{j\neq i}(\mu_i-\mu_j)b_{ij}\left(w_{j,x}-\frac{w_j}{v_j}v_{j,x}\right)_xr_i,
\end{align*}
where $b_{ij}$ is defined as follows
\begin{equation}\label{def:b-ija}
	b_{ij}:= \xi_j^\p\bar{v}_j\psi_{ji,v} -\theta^\p_j\psi_{ji,\si}.
\end{equation}
It is important to recall that $\theta^\p=1$ on $\mbox{supp}\xi$, since $\psi_{ji,\si}=O(1)\xi_j \bar{v}_j$ we deduce that $\theta^\p_j\psi_{ji,\si}=\psi_{ji,\si}$. We have then
\begin{equation}\label{def:b-ij}
	b_{ij}:= \xi_j^\p\bar{v}_j\psi_{ji,v} -\psi_{ji,\si}.
\end{equation}
We can write
\begin{align}
	\sum\limits_{i}\al_i^{6,1}
	&=\sum\limits_{i}\left(\sum\limits_{j\neq i}(\mu_i-\mu_j)b_{ij}(w_{j,x}-(w_j/v_j)v_{j,x})_x\right)\hat{r}_i\\
	&+\sum\limits_{i}\left(\sum\limits_{j\neq i}(\mu_i-\mu_j)b_{ij}(w_{j,x}-(w_j/v_j)v_{j,x})_x\right)(r_i-\hat{r}_i)\\
	&=\sum\limits_{i}\left(\sum\limits_{j\neq i}(\mu_i-\mu_j)b_{ij}(w_{j,x}-(w_j/v_j)v_{j,x})_x\right)\hat{r}_i+\sum\limits_{i}\sum\limits_{l=1}^{3}\al_i^{6,1,l},
\end{align}
where $\al_i^{6,1,l}$ are defined as follows
\begin{align*}
	\al^{6,1,1}_i&=-\left(\sum\limits_{j\neq i}(\mu_i-\mu_j)b_{ij}(w_{j,x}-(w_j/v_j)v_{j,x})_x\right)\sum\limits_{k\neq i}\psi_{ik} r_k(u),\\
	\al^{6,1,2}_i&=-\left(\sum\limits_{j\neq i}(\mu_i-\mu_j)b_{ij}(w_{j,x}-(w_j/v_j)v_{j,x})_x\right)\sum\limits_{k\neq i}(v_i\xi_i\pa_{v_i}\bar{v}_i-(w_i/v_i)\xi_i^\p\bar{v}_i)\psi_{ik,v} r_k,\\
	\al^{6,1,3}_i&=-\left(\sum\limits_{j\neq i}(\mu_i-\mu_j)b_{ij}(w_{j,x}-(w_j/v_j)v_{j,x})_x\right)\sum\limits_{k\neq i}(w_i/v_i)\theta^\p_i\psi_{ik,\si} r_k.
\end{align*}
Therefore, from \eqref{cal-phi-1} we may write
\begin{align}
	&\sum\limits_{i}\left(v_{i,t}+(\tilde{\la}_iv_i)_x-(\mu_iv_{i,x})_x-\sum\limits_{j\neq i}(\mu_i-\mu_j)b_{ij}(w_{j,x}-(w_j/v_j)v_{j,x})_x\right)\hat{r}_i \nonumber\\
	&+\sum\limits_{i}\left(w_{i,t}+(\tilde{\la}_iw_i)_x-(\mu_iw_{i,x})_x-\sum\limits_{j\neq i}(\mu_i-\mu_j)\hat{b}_{ij}(w_{j,x}-(w_j/v_j)v_{j,x})_x\right)r^\dagger_i \nonumber\\
	&=\sum\limits_{l=1}^{27}\sum\limits_{i}\widetilde{\al}_{i}^{1,l}+\sum\limits_{k=4}^{5}\sum\limits_{l=1}^{N_k}\sum\limits_{i}\al^{k,l}_i+\sum_{l=2}^{N_6}\sum_i\al_i^{6,l}+\sum\limits_{i}\sum\limits_{l=1}^{3}\al_i^{6,1,l}+\sum\limits_{k=7}^{10}\sum\limits_{l=1}^{N_k}\sum\limits_{i}\al^{k,l}_i+\sum\limits_{i}\al_{i}^{11},
	\label{eqn-remainder-1}
\end{align}
where $\hat{b}_{ij}$ will be defined later (see \eqref{def:hat-b-ij}) and $\al_i^{11}$ is as follows,
\begin{equation}\label{def:al-11}
	\al_i^{11}=-\left(\sum\limits_{j\neq i}(\mu_i-\mu_j)\hat{b}_{ij}(w_{j,x}-(w_j/v_j)v_{j,x})_x\right)r^\dagger_i.
\end{equation}

%----------------------------------------------------------------------------------------------------------------------------------------------

\subsection{Calculation for the remainder term $\psi_i$}\label{section:psi-i}
In this subsection, we calculate the term $\mathcal{R}_i^w$. Similar to previous section, we derive 
from \eqref{eqn-v-i}- \eqref{eqn-w-i}
\begin{align*}
	&B(u)u_{tx}+u_t\bullet B(u)u_x-A(u)u_t\\
	&=\sum\limits_{i}(w_{i,x}-\la_i^*v_{i,x})B(u)\tilde{r}_{i}+\sum\limits_{i}(w_{i}-\la_i^*v_{i})B(u)\tilde{r}_{i,x}+\sum\limits_{ij}(w_{i}-\la_i^*v_{i})v_j\tilde{r}_i\cdot DB(u)\tilde{r}_j\\
	&-\sum\limits_{i}(w_{i}-\la_i^*v_{i})A(u)\tilde{r}_{i}\\
	&=\sum\limits_{i}(w_{i,x}-\la_i^*v_{i,x})B(u)\tilde{r}_{i}+\sum\limits_{i}v_{i,x}\xi_i(w_{i}-\la_i^*v_{i})B(u)\pa_{v_i}\bar{v}_i\tilde{r}_{i,v}\\
	&+\sum\limits_{i}\sum\limits_{j\neq i}\xi_iv_{j,x}(w_{i}-\la_i^*v_{i})\pa_{v_j}\bar{v}_iB\tilde{r}_{i,v}+\sum\limits_{i}(w_{i}-\la_i^*v_{i})\left(\frac{w_i}{v_i}\right)_xB(u)[\xi^\p_i\bar{v}_i\tilde{r}_{i,v}-\theta_i^\p\tilde{r}_{i,\si}]\\
	&+\sum\limits_{i}(w_{i}-\la_i^*v_{i})\left[v_iB(u)\tilde{r}_{i,u}\tilde{r}_i+v_i\tilde{r}_i\cdot DB(u)\tilde{r}_i-A(u)\tilde{r}_i\right]\\
	&+\sum\limits_{i\neq j}(w_{i}-\la_i^*v_{i})v_j\tilde{r}_i\cdot DB(u)\tilde{r}_j+\sum\limits_{i\neq j}(w_{i}-\la_i^*v_{i})v_jB(u)\tilde{r}_{i,u}\tilde{r}_j,
\end{align*}
or equivalently,
\begin{align*}
	&B(u)u_{tx}+u_t\bullet B(u)u_x-A(u)u_t\\
	&=\sum\limits_{i}(w_{i,x}-\la_i^*v_{i,x})B(u)\tilde{r}_{i}+\sum\limits_{i}v_{i,x}\xi_i(w_{i}-\la_i^*v_{i})B(u)\pa_{v_i}\bar{v}_i\tilde{r}_{i,v}\\
	&+\sum\limits_{i}\sum\limits_{j\neq i}v_{j,x}\xi_i(w_{i}-\la_i^*v_{i})\pa_{v_j}\bar{v}_iB(u)\tilde{r}_{i,v}+\sum\limits_{i}(w_{i}-\la_i^*v_{i})\left(\frac{w_i}{v_i}\right)_xB(u)[\xi^\p_i\bar{v}_i\tilde{r}_{i,v}-\theta_i^\p\tilde{r}_{i,\si}]\\
	&+\sum\limits_{i}(w_{i}-\la_i^*v_{i})\left[\xi_i\bar{v}_iB(u)\tilde{r}_{i,u}\tilde{r}_i+\xi_i\bar{v}_i\tilde{r}_i\cdot DB(u)\tilde{r}_i-A(u)\tilde{r}_i\right]\\
	&+\sum\limits_{i}(v_i-\xi_i\bar{v}_i)(w_{i}-\la_i^*v_{i})\left[B(u)\tilde{r}_{i,u}\tilde{r}_i+\tilde{r}_i\cdot DB(u)\tilde{r}_i\right]\\
	&+\sum\limits_{i\neq j}(w_{i}-\la_i^*v_{i})v_j\tilde{r}_i\cdot DB(u)\tilde{r}_j+\sum\limits_{i\neq j}(w_{i}-\la_i^*v_{i})v_jB(u)\tilde{r}_{i,u}\tilde{r}_j.
\end{align*}
By using \eqref{identity-2} we get
\begin{align*}
	&B(u)u_{tx}+u_t\bullet B(u)u_x-A(u)u_t\\
	&=\sum\limits_{i}(w_{i,x}-\la_i^*v_{i,x})B(u)\tilde{r}_{i}+\sum\limits_{i}v_{i,x}(w_{i}-\la_i^*v_{i})\xi_iB(u)\pa_{v_i}\bar{v}_i\tilde{r}_{i,v}\\
	&+\sum\limits_{i}\sum\limits_{j\neq i}\xi_iv_{j,x}(w_{i}-\la_i^*v_{i})\pa_{v_j}\bar{v}_iB(u)\tilde{r}_{i,v}\\
	&+\sum\limits_{i}(w_{i}-\la_i^*v_{i})\left(\frac{w_i}{v_i}\right)_xB(u)[\xi^\p_i\bar{v}_i\tilde{r}_{i,v}-\theta_i^\p\tilde{r}_{i,\si}]\\
	&-\sum\limits_{i}(w_{i}-\la_i^*v_{i})\si_i\tilde{r}_i-\sum\limits_{i}(w_{i}-\la_i^*v_{i})\mu_i^{-1}(\tilde{\la}_i-\si_i)B[\tilde{r}_i+\xi_i\bar{v}_i\tilde{r}_{i,v}]\\
	&+\sum\limits_{i}(v_i-\xi_i\bar{v}_i)(w_{i}-\la_i^*v_{i})\left[B(u)\tilde{r}_{i,u}\tilde{r}_i+\tilde{r}_i\cdot DB(u)\tilde{r}_i\right]\\
	&+\sum\limits_{i\neq j}(w_{i}-\la_i^*v_{i})v_j\tilde{r}_i\cdot DB(u)\tilde{r}_j+\sum\limits_{i\neq j}(w_{i}-\la_i^*v_{i})v_jB(u)\tilde{r}_{i,u}\tilde{r}_j.
\end{align*}
After rearrangement,
\begin{align*}
	&B(u)u_{tx}+u_t\bullet B(u)u_x-A(u)u_t\\
	&=\sum\limits_{i}(w_{i,x}-\mu_i^{-1}(\tilde{\la}_i-\la^*_i+\theta_i)w_i)B(u)[\tilde{r}_i+(w_i/v_i)\xi_i^\p\bar{v}_i\tilde{r}_{i,v}-(w_i/v_i)\theta_i^\p\tilde{r}_{i,\si}]\\
	&+\sum\limits_{i}(v_{i,x}-\mu_i^{-1}(\tilde{\la}_i-\la^*_i+\theta_i)v_i)B(u)[(w_i\xi_i(\pa_{v_i}\bar{v}_i)-(w_i/v_i)^2\xi_i^\p\bar{v}_i)\tilde{r}_{i,v}+(w_i/v_i)^2\theta_i^\p\tilde{r}_{i,\si}]\\
	&-\sum\limits_{i}\la_i^*(v_{i,x}-\mu_i^{-1}(\tilde{\la}_i-\la^*_i+\theta_i^*)v_i)B(u)[\tilde{r}_{i}+(\xi_iv_i\pa_{v_i}\bar{v}_i-(w_i/v_i)\xi_i^\p\bar{v}_i)\tilde{r}_{i,v}+(w_i/v_i)\theta_i^\p\tilde{r}_{i,\si}]\\
	&-\sum\limits_{i}\la_i^*(w_{i,x}-\mu_i^{-1}(\tilde{\la}_i-\la^*_i+\theta_i)w_i)B(u)[\xi_i^\p\bar{v}_i\tilde{r}_{i,v}-\theta_i^\p\tilde{r}_{i,\si}]\\
	&+\sum\limits_{i}\sum\limits_{j\neq i}\xi_iv_{j,x}(w_{i}-\la_i^*v_{i})\pa_{v_j}\bar{v}_iB(u)\tilde{r}_{i,v}-\sum\limits_{i}w_{i}(\la_i^*-\theta_i)\tilde{r}_i+\sum\limits_{i}\la_i^*v_{i}(\la_i^*-\theta_i)\tilde{r}_i\\
	&+\sum\limits_{i}(w_{i}-\la_i^*v_{i})\mu_i^{-1}(\tilde{\la}_i-\si_i)\xi_i(v_i\pa_{v_i}\bar{v}_i-\bar{v}_i)B\tilde{r}_{i,v}\\
	&+\sum\limits_{i}(v_i-\xi_i\bar{v}_i)(w_{i}-\la_i^*v_{i})\left[B(u)\tilde{r}_{i,u}\tilde{r}_i+\tilde{r}_i\cdot DB(u)\tilde{r}_i\right]\\
	&+\sum\limits_{i\neq j}(w_{i}-\la_i^*v_{i})v_j\tilde{r}_i\cdot DB(u)\tilde{r}_j+\sum\limits_{i\neq j}(w_{i}-\la_i^*v_{i})v_jB(u)\tilde{r}_{i,u}\tilde{r}_j.
\end{align*}
Furthermore, we have
\begin{align*}
	&B(u)u_{tx}+u_t\bullet B(u)u_x-A(u)u_t\\
	&=\sum\limits_{i}(\mu_iw_{i,x}-(\tilde{\la}_i-\la^*_i+\theta_i)w_i)[\tilde{r}_i+(w_i/v_i)\xi_i^\p\bar{v}_i\tilde{r}_{i,v}-(w_i/v_i)\theta_i^\p\tilde{r}_{i,\si}]\\
	&+\sum\limits_{i}(\mu_iv_{i,x}-(\tilde{\la}_i-\la^*_i+\theta_i)v_i)[(w_i\xi_i(\pa_{v_i}\bar{v}_i)-(w_i/v_i)^2\xi_i^\p\bar{v}_i)\tilde{r}_{i,v}+(w_i/v_i)^2\theta_i^\p\tilde{r}_{i,\si}]\\
	&-\sum\limits_{i}\la_i^*(\mu_iv_{i,x}-(\tilde{\la}_i-\la^*_i+\theta_i)v_i)[\tilde{r}_{i}+(\xi_iv_i\pa_{v_i}\bar{v}_i-(w_i/v_i)\xi_i^\p\bar{v}_i)\tilde{r}_{i,v}+(w_i/v_i)\theta_i^\p\tilde{r}_{i,\si}]\\
	&-\sum\limits_{i}\la_i^*(\mu_iw_{i,x}-(\tilde{\la}_i-\la^*_i+\theta_i)w_i)[\xi_i^\p\bar{v}_i\tilde{r}_{i,v}-\theta_i^\p\tilde{r}_{i,\si}]\\
	&+\sum\limits_{i}\sum\limits_{j\neq i}\xi_iv_{j,x}(w_{i}-\la_i^*v_{i})\pa_{v_j}\bar{v}_iB(u)\tilde{r}_{i,v}-\sum\limits_{i}w_{i}(\la_i^*-\theta_i)\tilde{r}_i+\sum\limits_{i}\la_i^*v_{i}(\la_i^*-\theta_i)\tilde{r}_i\\
	&+\sum\limits_{i}(w_{i,x}-\mu_i^{-1}(\tilde{\la}_i-\la^*_i+\theta_i)w_i)(B-\mu_iI_n)[\tilde{r}_i+(w_i/v_i)\xi_i^\p\bar{v}_i\tilde{r}_{i,v}-(w_i/v_i)\theta_i^\p\tilde{r}_{i,\si}]\\
	&+\sum\limits_{i}(v_{i,x}-\mu_i^{-1}(\tilde{\la}_i-\la^*_i+\theta_i)v_i)(B-\mu_iI_n)[(w_i\xi_i(\pa_{v_i}\bar{v}_i)-(w_i/v_i)^2\xi_i^\p\bar{v}_i)\tilde{r}_{i,v}+(w_i/v_i)^2\theta_i^\p\tilde{r}_{i,\si}]\\
	&-\sum\limits_{i}\la_i^*(v_{i,x}-\mu_i^{-1}(\tilde{\la}_i-\la^*_i+\theta_i^*)v_i)(B-\mu_iI_n)[\tilde{r}_{i}+(\xi_iv_i\pa_{v_i}\bar{v}_i-(w_i/v_i)\xi_i^\p\bar{v}_i)\tilde{r}_{i,v}+(w_i/v_i)\theta_i^\p\tilde{r}_{i,\si}]\\
	&-\sum\limits_{i}\la_i^*(w_{i,x}-\mu_i^{-1}(\tilde{\la}_i-\la^*_i+\theta_i)w_i)(B-\mu_iI_n)[\xi_i^\p\bar{v}_i\tilde{r}_{i,v}-\theta_i^\p\tilde{r}_{i,\si}]\\
	&+\sum\limits_{i}(w_{i}-\la_i^*v_{i})\mu_i^{-1}(\tilde{\la}_i-\si_i)\xi_i(v_i\pa_{v_i}\bar{v}_i-\bar{v}_i)B\tilde{r}_{i,v}\\
	&+\sum\limits_{i}(v_i-\xi_i\bar{v}_i)(w_{i}-\la_i^*v_{i})\left[B(u)\tilde{r}_{i,u}\tilde{r}_i+\tilde{r}_i\cdot DB(u)\tilde{r}_i\right]\\
	&+\sum\limits_{i\neq j}(w_{i}-\la_i^*v_{i})v_j\tilde{r}_i\cdot DB(u)\tilde{r}_j+\sum\limits_{i\neq j}(w_{i}-\la_i^*v_{i})v_jB(u)\tilde{r}_{i,u}\tilde{r}_j.
\end{align*}
Thus, we may write
\begin{equation}\label{eqn-B-u-tx}
	B(u)u_{tx}+u_t\bullet B(u)u_x-A(u)u_t=\sum\limits_{j=1}^{15}K_{j}
\end{equation}
where $\{K_i\}_{i=1}^{15}$ are defined as follows
\begin{align*}
	K_{1}&=\sum\limits_{i}(\mu_iw_{i,x}-\tilde{\la}_iw_i)[\tilde{r}_i+(w_i/v_i)\xi_i^\p\bar{v}_i\tilde{r}_{i,v}-(w_i/v_i)\theta_i^\p\tilde{r}_{i,\si}],\\
	K_{2}&=\sum\limits_{i}(\mu_iv_{i,x}-\tilde{\la}_iv_i)[(w_i\xi_i(\pa_{v_i}\bar{v}_i)-(w_i/v_i)^2\xi_i^\p\bar{v}_i)\tilde{r}_{i,v}+(w_i/v_i)^2\theta_i^\p\tilde{r}_{i,\si}],\\
	K_{3}&=-\sum\limits_{i}\la_i^*(\mu_iv_{i,x}-\tilde{\la}_iv_i)[\tilde{r}_{i}+(\xi_iv_i\pa_{v_i}\bar{v}_i-(w_i/v_i)\xi_i^\p\bar{v}_i)\tilde{r}_{i,v}+(w_i/v_i)\theta_i^\p\tilde{r}_{i,\si}],\\
	K_{4}&=-\sum\limits_{i}\la_i^*(\mu_iw_{i,x}-\tilde{\la}_iw_i)[\xi_i^\p\bar{v}_i\tilde{r}_{i,v}-\theta_i^\p\tilde{r}_{i,\si}],\\
	K_{5}&=\sum\limits_{i}w_{i}v_i(\la_i^*-\theta_i)\xi_i\pa_{v_i}\bar{v}_i\tilde{r}_{i,v},\\
	K_{6}&=-\sum\limits_{i}\la_i^*v^2_{i}\xi_i\pa_{v_i}\bar{v}_i(\la_i^*-\theta_i)\tilde{r}_{i,v}\\
	K_{7}&=\sum\limits_{i}\sum\limits_{j\neq i}\xi_iv_{j,x}(w_{i}-\la_i^*v_{i})\pa_{v_j}\bar{v}_iB(u)\tilde{r}_{i,v},\\
	K_{8}&=\sum\limits_{i}(w_{i,x}-\mu_i^{-1}(\tilde{\la}_i-\la^*_i+\theta_i)w_i)(B-\mu_iI_n)\tilde{r}_i,\\
	K_{9}&=\sum\limits_{i}(v_{i,x}-\mu_i^{-1}(\tilde{\la}_i-\la^*_i+\theta_i)v_i)(B-\mu_iI_n)[w_i\xi_i(\pa_{v_i}\bar{v}_i)\tilde{r}_{i,v}],\\
	K_{10}&=-\sum\limits_{i}\la_i^*(v_{i,x}-\mu_i^{-1}(\tilde{\la}_i-\la^*_i+\theta_i^*)v_i)(B-\mu_iI_n)[\tilde{r}_{i}+(\xi_iv_i\pa_{v_i}\bar{v}_i)\tilde{r}_{i,v}],\\
	K_{11}&=\sum\limits_{i}(w_{i,x}-(w_i/v_i)v_{i,x})\frac{w_i}{v_i}(B-\mu_iI_n)[\xi^\p_i\bar{v}_i\tilde{r}_{i,v}-\theta_i^\p\tilde{r}_{i,\si}],\\
	K_{12}&=-\sum\limits_{i}\la_i^*(w_{i,x}-(w_i/v_i)v_{i,x})(B-\mu_iI_n)[\xi^\p_i\bar{v}_i\tilde{r}_{i,v}-\theta_i^\p\tilde{r}_{i,\si}],\\
	K_{13}&=\sum\limits_{i}(w_{i}-\la_i^*v_{i})\mu_i^{-1}(\tilde{\la}_i-\si_i)\xi_i(v_i\pa_{v_i}\bar{v}_i-\bar{v}_i)B\tilde{r}_{i,v},\\
	K_{14}&=\sum\limits_{i}(v_i-\xi_i\bar{v}_i)(w_{i}-\la_i^*v_{i})\left[B(u)\tilde{r}_{i,u}\tilde{r}_i+\tilde{r}_i\cdot DB(u)\tilde{r}_i\right],\\
	K_{15}&=\sum\limits_{i\neq j}(w_{i}-\la_i^*v_{i})v_j\left[\tilde{r}_i\cdot DB(u)\tilde{r}_j+B(u)\tilde{r}_{i,u}\tilde{r}_j\right].
\end{align*}
Indeed we can observe that
\begin{align*}
&K_8+K_{9}+K_{11}=\sum\limits_{i}(w_{i,x}-\mu_i^{-1}(\tilde{\la}_i-\la^*_i+\theta_i)w_i)(B-\mu_iI_n)[\tilde{r}_i+(w_i/v_i)\xi_i^\p\bar{v}_i\tilde{r}_{i,v}-(w_i/v_i)\theta_i^\p\tilde{r}_{i,\si}]\\
	&+\sum\limits_{i}(v_{i,x}-\mu_i^{-1}(\tilde{\la}_i-\la^*_i+\theta_i)v_i)(B-\mu_iI_n)[(w_i\xi_i(\pa_{v_i}\bar{v}_i)-(w_i/v_i)^2\xi_i^\p\bar{v}_i)\tilde{r}_{i,v}+(w_i/v_i)^2\theta_i^\p\tilde{r}_{i,\si}],\\[2mm]
	&K_{10}+K_{12}=-\sum\limits_{i}\la_i^*(w_{i,x}-\mu_i^{-1}(\tilde{\la}_i-\la^*_i+\theta_i)w_i)(B-\mu_iI_n)[\xi_i^\p\bar{v}_i\tilde{r}_{i,v}-\theta_i^\p\tilde{r}_{i,\si}]\\
	&-\sum\limits_{i}\la_i^*(v_{i,x}-\mu_i^{-1}(\tilde{\la}_i-\la^*_i+\theta_i^*)v_i)(B-\mu_iI_n)[\tilde{r}_{i}+(\xi_iv_i\pa_{v_i}\bar{v}_i-(w_i/v_i)\xi_i^\p\bar{v}_i)\tilde{r}_{i,v}+(w_i/v_i)\theta_i^\p\tilde{r}_{i,\si}].\\[2mm]
	&K_1+K_2+K_5=\sum\limits_{i}(\mu_iv_{i,x}-(\tilde{\la}_i-\la^*_i+\theta_i)v_i)[(w_i\xi_i(\pa_{v_i}\bar{v}_i)-(w_i/v_i)^2\xi_i^\p\bar{v}_i)\tilde{r}_{i,v}+(w_i/v_i)^2\theta_i^\p\tilde{r}_{i,\si}]\\
	&+\sum\limits_{i}(\mu_iw_{i,x}-(\tilde{\la}_i-\la^*_i+\theta_i)w_i)[\tilde{r}_i+(w_i/v_i)\xi_i^\p\bar{v}_i\tilde{r}_{i,v}-(w_i/v_i)\theta_i^\p\tilde{r}_{i,\si}]-\sum\limits_{i}w_{i}(\la_i^*-\theta_i)\tilde{r}_i.\\[2mm]
	&K_3+K_4+K_6=-\sum\limits_{i}\la_i^*(\mu_iw_{i,x}-(\tilde{\la}_i-\la^*_i+\theta_i)w_i)[\xi_i^\p\bar{v}_i\tilde{r}_{i,v}-\theta_i^\p\tilde{r}_{i,\si}]+\sum\limits_{i}\la_i^*v_{i}(\la_i^*-\theta_i)\tilde{r}_i \\
	&-\sum\limits_{i}\la_i^*(\mu_iv_{i,x}-(\tilde{\la}_i-\la^*_i+\theta_i)v_i)[\tilde{r}_{i}+(\xi_iv_i\pa_{v_i}\bar{v}_i-(w_i/v_i)\xi_i^\p\bar{v}_i)\tilde{r}_{i,v}+(w_i/v_i)\theta_i^\p\tilde{r}_{i,\si}].
	%&+\sum\limits_{i}\la_i^*v_{i}(\la_i^*-\theta_i)\tilde{r}_i\\\
	\end{align*}
Next we calculate derivatives of each $K_j$'s with respect to $x$. 
\begin{align}
	&K_{1,x}=\sum\limits_{i}(\mu_iw_{i,x}-\tilde{\la}_iw_i)_x[\tilde{r}_i+(w_i/v_i)\xi_i^\p\bar{v}_i\tilde{r}_{i,v}-(w_i/v_i)\theta_i^\p\tilde{r}_{i,\si}]\nonumber\\
	&+\sum\limits_{i}(\mu_iw_{i,x}-\tilde{\la}_iw_i)v_i\tilde{r}_{i,u}\tilde{r}_i+\sum\limits_{i}\sum\limits_{j\neq i}(\mu_iw_{i,x}-\tilde{\la}_iw_i)v_j\tilde{r}_{i,u}\tilde{r}_j\nonumber\\
	&+\sum\limits_{i}(\mu_iw_{i,x}-\tilde{\la}_iw_i)v_{i,x}\xi_i(\pa_{v_i}\bar{v}_i)\tilde{r}_{i,v}+\sum\limits_{i}\sum\limits_{j\neq i}(\mu_iw_{i,x}-\tilde{\la}_iw_i)v_{j,x}\xi_i(\pa_{v_j}\bar{v}_i)\tilde{r}_{i,v}\nonumber\\
	&+2\sum\limits_{i}(\mu_iw_{i,x}-\tilde{\la}_iw_i)\xi^\p_i\bar{v}_i\left(\frac{w_i}{v_i}\right)_x\tilde{r}_{i,v}+\sum\limits_{i}(\mu_iw_{i,x}-\tilde{\la}_iw_i)\xi^{\p\p}_i\bar{v}_i\frac{w_i}{v_i}\left(\frac{w_i}{v_i}\right)_x\tilde{r}_{i,v}\nonumber\\
	&+\sum\limits_{i}(\mu_iw_{i,x}-\tilde{\la}_iw_i)v_{i,x}\xi^\p_i(\pa_{v_i}\bar{v}_i)\frac{w_i}{v_i}\tilde{r}_{i,v}+\sum\limits_{i}\sum\limits_{j\neq i}(\mu_iw_{i,x}-\tilde{\la}_iw_i)v_{j,x}\xi_i^\p(\pa_{v_j}\bar{v}_i)\frac{w_i}{v_i}\tilde{r}_{i,v}\nonumber\\
	&-2\sum\limits_{i}(\mu_iw_{i,x}-\tilde{\la}_iw_i)\theta_i^\p\left(\frac{w_i}{v_i}\right)_x\tilde{r}_{i,\si}-\sum\limits_{i}(\mu_iw_{i,x}-\tilde{\la}_iw_i)\frac{w_i}{v_i}\left(\frac{w_i}{v_i}\right)_x\theta_i^{\p\p} \tilde{r}_{i,\si}\nonumber\\
	&+\sum\limits_{i}(\mu_iw_{i,x}-\tilde{\la}_iw_i)v_i\frac{w_i}{v_i}\xi_i^\p\bar{v}_i\tilde{r}_{i,uv}\tilde{r}_i+\sum\limits_{i}\sum\limits_{j\neq i}(\mu_iw_{i,x}-\tilde{\la}_iw_i)v_j\frac{w_i}{v_i}\xi_i^\p\bar{v}_i\tilde{r}_{i,uv}\tilde{r}_{j}\nonumber\\
	&+\sum\limits_{i}(\mu_iw_{i,x}-\tilde{\la}_iw_i)v_{i,x}\frac{w_i}{v_i}\xi_i^\p\xi_i\pa_{v_i}\bar{v}_i\bar{v}_i\tilde{r}_{i,vv}+\sum\limits_{i}\sum\limits_{j\neq i}(\mu_iw_{i,x}-\tilde{\la}_iw_i)v_{j,x}\frac{w_i}{v_i}\xi_i^\p\xi_i\bar{v}_i(\pa_{v_j}\bar{v}_i)\tilde{r}_{i,vv}\nonumber\\
	&+\sum\limits_{i}(\mu_iw_{i,x}-\tilde{\la}_iw_i)\frac{w_i}{v_i}(\xi_i^\p)^2\bar{v}_i^{2}\left(\frac{w_i}{v_i}\right)_x\tilde{r}_{i,vv}-\sum\limits_{i}(\mu_iw_{i,x}-\tilde{\la}_iw_i)\frac{w_i}{v_i}(\xi_i^\p\bar{v}_i)\theta_i^\p\left(\frac{w_i}{v_i}\right)_x\tilde{r}_{i,v\si}\nonumber\\
	&-\sum\limits_{i}(\mu_iw_{i,x}-\tilde{\la}_iw_i) w_i\theta_i^\p \tilde{r}_{i,u\si}\tilde{r}_i-\sum\limits_{i}\sum\limits_{j\neq i}(\mu_iw_{i,x}-\tilde{\la}_iw_i) (w_i/v_i)v_j\theta_i^\p \tilde{r}_{i,u\si}\tilde{r}_j\nonumber\\
	&-\sum\limits_{i}(\mu_iw_{i,x}-\tilde{\la}_iw_i) \xi_i\frac{w_i}{v_i}v_{i,x}\theta_i^\p \pa_{v_i}\bar{v}_i\tilde{r}_{i,v\si}-\sum\limits_{i}\sum\limits_{j\neq i}(\mu_iw_{i,x}-\tilde{\la}_iw_i)\xi_i \frac{w_i}{v_i}v_{j,x}\theta_i^\p \pa_{v_j}\bar{v}_i\tilde{r}_{i,v\si}\nonumber\\
	&-\sum\limits_{i}(\mu_iw_{i,x}-\tilde{\la}_iw_i) \xi_i^\p\bar{v}_i\frac{w_i}{v_i}\left(\frac{w_i}{v_i}\right)_x\theta_i^\p \tilde{r}_{i,v\si} +\sum\limits_{i}(\mu_iw_{i,x}-\tilde{\la}_iw_i) \frac{w_i}{v_i}\left(\frac{w_i}{v_i}\right)_x(\theta_i^\p)^2 \tilde{r}_{i,\si\si}\nonumber\\
	&=\sum_{l=1}^{23}\sum_i\beta_i^{1,l}. \label{eqn-K-1}
\end{align}
Next we calculate the derivative of term $K_2$.
\begin{align}
	&K_{2,x}=\sum\limits_{i}(\mu_iv_{i,x}-\tilde{\la}_iv_i)_x[(w_i\xi_i(\pa_{v_i}\bar{v}_i)-(w_i/v_i)^2\xi_i^\p\bar{v}_i)\tilde{r}_{i,v}+(w_i/v_i)^2\theta_i^\p\tilde{r}_{i,\si}]\nonumber\\
	&+\sum\limits_{i}(\mu_iv_{i,x}-\tilde{\la}_iv_i)w_{i,x}\xi_i\pa_{v_i}\bar{v}_i\tilde{r}_{i,v}+\sum\limits_{i}(\mu_iv_{i,x}-\tilde{\la}_iv_i)\xi_iw_{i}v_{i,x}\pa_{v_iv_i}\bar{v}_i\tilde{r}_{i,v}\nonumber\\
	&+\sum\limits_{i}\sum\limits_{j\neq i}(\mu_iv_{i,x}-\tilde{\la}_iv_i)\xi_iw_{i}v_{j,x}\pa_{v_jv_i}\bar{v}_i\tilde{r}_{i,v}
	+\sum\limits_{i}(\mu_iv_{i,x}-\tilde{\la}_iv_i)\xi_iw_{i}v_i\pa_{v_i}\bar{v}_i\tilde{r}_{i,uv}\tilde{r}_i\nonumber\\
	&+\sum\limits_{i}\sum\limits_{j\neq i}(\mu_iv_{i,x}-\tilde{\la}_iv_i)\xi_iw_{i}v_j\pa_{v_i}\bar{v}_i\tilde{r}_{i,uv}\tilde{r}_j+\sum\limits_{i}(\mu_iv_{i,x}-\tilde{\la}_iv_i)\xi_i^2w_{i}v_{i,x}(\pa_{v_i}\bar{v}_i)^2\tilde{r}_{i,vv}\nonumber\\
	&+\sum\limits_{i}\sum\limits_{j\neq i}(\mu_iv_{i,x}-\tilde{\la}_iv_i)v_{j,x}w_{i}\xi_i^2\pa_{v_i}\bar{v}_i\pa_{v_j}\bar{v}_i\tilde{r}_{i,vv}-\sum\limits_{i}(\mu_iv_{i,x}-\tilde{\la}_iv_i)w_{i}\xi_i\pa_{v_i}\bar{v}_i\theta_i^\p\left(\frac{w_i}{v_i}\right)_x\tilde{r}_{i,v\si}\nonumber\\
	&+\sum\limits_{i} \xi_i^\p\bar{v}_i(\mu_iv_{i,x}-\tilde{\la}_iv_i) w_{i}\xi_i\pa_{v_i}\bar{v}_i\left(\frac{w_i}{v_i}\right)_x\tilde{r}_{i,vv}+\sum\limits_{i}(\mu_iv_{i,x}-\tilde{\la}_iv_i)w_i\xi^\p_i\left(\frac{w_i}{v_i}\right)_x(\pa_{v_i}\bar{v}_i)\tilde{r}_{i,v} \nonumber\\
	&-2\sum\limits_{i}(\mu_iv_{i,x}-\tilde{\la}_iv_i) \frac{w_i}{v_i}\left(\frac{w_i}{v_i}\right)_x\xi_i^\p\bar{v}_i\tilde{r}_{i,v}-\sum\limits_{i}(\mu_iv_{i,x}-\tilde{\la}_iv_i) \frac{w^2_i}{v^2_i}\xi_i^\p v_{i,x}\pa_{v_i}\bar{v}_i\tilde{r}_{i,v}\nonumber\\
	&-\sum\limits_{i}\sum\limits_{j\neq i}(\mu_iv_{i,x}-\tilde{\la}_iv_i) \frac{w^2_i}{v^2_i}v_{j,x}\xi_i^\p\pa_{v_j}\bar{v}_i\tilde{r}_{i,v}
	-\sum\limits_{i}(\mu_iv_{i,x}-\tilde{\la}_iv_i) \frac{w^2_i}{v^2_i}\xi_i^{\p\p}\bar{v}_i\left(\frac{w_i}{v_i}\right)_x\tilde{r}_{i,v}\nonumber\\
	&-\sum\limits_{i}(\mu_iv_{i,x}-\tilde{\la}_iv_i) \frac{w^2_i}{v^2_i}\xi_i^\p \bar{v}_i v_{i}\tilde{r}_{i,uv}\tilde{r}_i-\sum\limits_{i}\sum\limits_{j\neq i}(\mu_iv_{i,x}-\tilde{\la}_iv_i) \frac{w^2_i}{v^2_i}\xi_i^\p \bar{v}_i v_{j}\tilde{r}_{i,uv}\tilde{r}_j\nonumber\\
	&-\sum\limits_{i}(\mu_iv_{i,x}-\tilde{\la}_iv_i) \frac{w^2_i}{v^2_i}\bar{v}_i\xi_i^\p \xi_iv_{i,x}\pa_{v_i}\bar{v}_i\tilde{r}_{i,vv}-\sum\limits_{i}\sum\limits_{j\neq i}(\mu_iv_{i,x}-\tilde{\la}_iv_i)v_{j,x} \frac{w^2_i}{v^2_i}\bar{v}_i\xi_i^\p\xi_i\pa_{v_j}\bar{v}_i\tilde{r}_{i,vv}\nonumber\\
	&-\sum\limits_{i}(\mu_iv_{i,x}-\tilde{\la}_iv_i) \frac{w^2_i}{v^2_i}(\bar{v}_i\xi_i^\p)^2
			\left(\frac{w_i}{v_i}\right)_x
	\tilde{r}_{i,vv}+\sum\limits_{i}(\mu_iv_{i,x}-\tilde{\la}_iv_i) \frac{w^2_i}{v^2_i}\bar{v}_i\xi_i^\p  \theta_i^\p\left(\frac{w_i}{v_i}\right)_x\tilde{r}_{i,v\si}\nonumber\\
	&+2\sum\limits_{i}(\mu_iv_{i,x}-\tilde{\la}_iv_i)\frac{w_i}{v_i}\left(\frac{w_i}{v_i}\right)_x\theta_i^\p\tilde{r}_{i,\si}+\sum\limits_{i}(\mu_iv_{i,x}-\tilde{\la}_iv_i)(w_i/v_i)^2\theta^{\p\p}_i\left(\frac{w_i}{v_i}\right)_x\tilde{r}_{i,\si}\nonumber\\
	&+\sum\limits_{i}(\mu_iv_{i,x}-\tilde{\la}_iv_i)w_i(w_i/v_i)\theta^\p_i\tilde{r}_{i,u\si}\tilde{r}_i+\sum\limits_{i}\sum\limits_{j\neq i}(\mu_iv_{i,x}-\tilde{\la}_iv_i)(w_i/v_i)^2v_j\theta^\p_i\tilde{r}_{i,u\si}\tilde{r}_j \nonumber\\
	&+\sum\limits_{i}(\mu_iv_{i,x}-\tilde{\la}_iv_i)v_{i,x}\xi_i(w_i/v_i)^2\pa_{v_i}\bar{v}_i\theta^\p_i\tilde{r}_{i,v\si}+\sum\limits_{i}\sum\limits_{j\neq i}(\mu_iv_{i,x}-\tilde{\la}_iv_i)\xi_iv_{j,x}(w_i/v_i)^2\pa_{v_j}\bar{v}_i\theta^\p_i\tilde{r}_{i,v\si} \nonumber\\
	&+\sum\limits_{i}(\mu_iv_{i,x}-\tilde{\la}_iv_i)\xi_i^\p\left(\frac{w_i}{v_i}\right)_{x}(w_i/v_i)^2\bar{v}_i\theta^\p_i\tilde{r}_{i,v\si}-\sum\limits_{i}(\mu_iv_{i,x}-\tilde{\la}_iv_i)(w_i/v_i)^2(\theta^\p_i)^2\left(\frac{w_i}{v_i}\right)_x\tilde{r}_{i,\si\si}\nonumber\\
	&=\sum_{l=1}^{29}\sum_i\beta_i^{2,l} . \label{eqn-K-2}
\end{align}
We remark now that
\begin{equation}\label{K3K4}
K_{3,x}=-\sum\limits_{k=1}^{20}\sum\limits_{i}\la_i^*\al_{i}^{1,k}\;\;\mbox{and}\;\;K_{4,x}=-\sum\limits_{k=1}^{17}\sum\limits_{i}\la_i^*\al_{i}^{2,k}.
\end{equation}
For the term $K_5$ we get,
\begin{align}
	K_{5,x}&=-\sum\limits_{i}\theta_i^\p (\pa_{v_i}\bar{v}_i)v_iw_i\left(\frac{w_i}{v_i}\right)_x\xi_i\tilde{r}_{i,v}+\sum\limits_{i}(\la_i^*-\theta_i)\xi_iv_{i,x}(\pa_{v_i}\bar{v}_i+v_i\pa_{v_iv_i}\bar{v}_i)w_i\tilde{r}_{i,v}\nonumber\\
	&+\sum\limits_{i}(\la_i^*-\theta_i)v_{i}\xi_i^\p\left(\frac{w_i}{v_i}\right)_x(\pa_{v_i}\bar{v}_i)w_{i}\tilde{r}_{i,v}+\sum\limits_{i}(\la_i^*-\theta_i)v_{i}\xi_i(\pa_{v_i}\bar{v}_i)w_{i,x}\tilde{r}_{i,v}\nonumber\\
	&+\sum\limits_{i}\sum\limits_{j\neq i}(\la_i^*-\theta_i)\xi_iv_{j,x}w_iv_i\pa_{v_iv_j}\bar{v}_i\tilde{r}_{i,v}+\sum\limits_{i}(\la_i^*-\theta_i)v^2_{i}\xi_i(\pa_{v_i}\bar{v}_i)w_{i}\tilde{r}_{i,uv}\tilde{r}_i \nonumber\\
	&+\sum\limits_{i}\sum\limits_{j\neq i}(\la_i^*-\theta_i)\xi_i v_{i}w_{i}v_j(\pa_{v_i}\bar{v}_i)\tilde{r}_{i,uv}\tilde{r}_j+\sum\limits_{i}(\la_i^*-\theta_i)v_{i}\xi_i(\pa_{v_i}\bar{v}_i)w_{i}\bar{v}_i\xi^\p_i\left(\frac{w_i}{v_i}\right)_x\tilde{r}_{i,vv} \nonumber\\
	&+\sum\limits_{i}(\la_i^*-\theta_i)\xi_i^2v_{i,x}v_{i}w_{i}(\pa_{v_i}\bar{v}_i)^2\tilde{r}_{i,vv}+\sum\limits_{i}\sum\limits_{j\neq i}(\la_i^*-\theta_i)\xi_i^2v_{j,x}v_{i}w_{i}(\pa_{v_i}\bar{v}_i)(\pa_{v_j}\bar{v}_i)\tilde{r}_{i,vv}\nonumber\\
	&-\sum\limits_{i}(\la_i^*-\theta_i)\xi_iv_{i}w_{i}(\pa_{v_i}\bar{v}_i)\theta^\p_i\left(\frac{w_i}{v_i}\right)_x\tilde{r}_{i,v\si}\nonumber\\
	&=\sum_{l=1}^{11}\sum_i\beta_i^{5,l}. \label{eqn-K-5}
\end{align}
We observe again that 
\begin{equation}\label{K6}
K_{6,x}=-\sum\limits_{k=1}^{11}\sum\limits_{i}\la_i^*\al_i^{3,k}.
\end{equation}
For the term $K_7$ we observe that we can rewrite these terms in function of $J_4$
\begin{equation*}
	K_{7,x}= \sum\limits_{i\ne j}\left(\frac{w_i}{v_i}\right)_x\xi_iv_{j,x}v_i\pa_{v_j}\bar{v}_iB(u)\tilde{r}_{i,v}+
	\sum\limits_{l=1}^{8}\sum\limits_{i}\left(\frac{w_i}{v_i}-\la_i^*\right)\al_{i}^{4,l}.
\end{equation*}
Taking derivative of $K_8$ we have
\begin{align*}
	K_{8,x}&=\sum\limits_{i}(w_{i,x}-\mu_i^{-1}(\tilde{\la}_i-\la_i^*+\theta_i)w_i)_x(B-\mu_iI_n)\tilde{r}_{i}\\
	&+\sum\limits_{i}(w_{i,x}-\mu_i^{-1}(\tilde{\la}_i-\la^*_i+\theta_i)w_i)v_i\tilde{r}_i\cdot D(B-\mu_iI_n)\tilde{r}_{i}\\
	&+\sum\limits_{i}\sum\limits_{j\neq i}(w_{i,x}-\mu_i^{-1}(\tilde{\la}_i-\la^*_i+\theta_i)w_i)v_j\tilde{r}_j\cdot D(B-\mu_iI_n)\tilde{r}_{i}\\
	&+\sum\limits_{i}(w_{i,x}-\mu_i^{-1}(\tilde{\la}_i-\la_i^*+\theta_i)w_i)v_{i}(B-\mu_iI_n)\tilde{r}_{i,u}\tilde{r}_i\\
	&+\sum\limits_{i}\sum\limits_{j\neq i}(w_{i,x}-\mu_i^{-1}(\tilde{\la}_i-\la_i^*+\theta_i)w_i)v_{j}(B-\mu_iI_n)\tilde{r}_{i,u}\tilde{r}_j\\
	&+\sum\limits_{i}(w_{i,x}-\mu_i^{-1}(\tilde{\la}_i-\la_i^*+\theta_i)w_i)v_{i,x}\xi_i(B-\mu_iI_n)\pa_{v_i}\bar{v}_i\tilde{r}_{i,v}\\
	&+\sum\limits_{i}\sum\limits_{j\neq i}(w_{i,x}-\mu_i^{-1}(\tilde{\la}_i-\la_i^*+\theta_i)w_i)v_{j,x}\xi_i(B-\mu_iI_n)(\pa_{v_j}\bar{v}_i)\tilde{r}_{i,v}\\
	&+\sum\limits_{i}(w_{i,x}-\mu_i^{-1}(\tilde{\la}_i-\la_i^*+\theta_i)w_i)\bar{v}_i\xi^\p_i\left(\frac{w_i}{v_i}\right)_x(B-\mu_iI_n)\tilde{r}_{i,v}\\
	&-\sum\limits_{i}(w_{i,x}-\mu_i^{-1}(\tilde{\la}_i-\la_i^*+\theta_i)w_i)\theta_i^\p\left(\frac{w_i}{v_i}\right)_x(B-\mu_iI_n)\tilde{r}_{i,\si}\\
	&=:\sum\limits_{l=1}^{9}\sum\limits_{i}\B^{8,l}_i.
\end{align*}
Similarly, we obtain
\begin{align*}
	K_{9,x}&=\sum\limits_{i}(v_{i,x}-\mu_i^{-1}(\tilde{\la}_i-\la_i^*+\theta_i)v_i)_x(B-\mu_iI_n)[w_i\xi_i\pa_{v_i}\bar{v}_i\tilde{r}_{i,v}]\\
	&+\sum\limits_{i}(v_{i,x}-\mu_i^{-1}(\tilde{\la}_i-\la^*_i+\theta_i)v_i)v_i\tilde{r}_i\cdot D(B-\mu_iI_n)[w_i\xi_i\pa_{v_i}\bar{v}_i\tilde{r}_{i,v}]\\
	&+\sum\limits_{i}\sum\limits_{j\neq i}(v_{i,x}-\mu_i^{-1}(\tilde{\la}_i-\la^*_i+\theta_i)v_i)v_j\tilde{r}_j\cdot D(B-\mu_iI_n)[w_i\xi_i\pa_{v_i}\bar{v}_i\tilde{r}_{i,v}]\\
	&+\sum\limits_{i}(v_{i,x}-\mu_i^{-1}(\tilde{\la}_i-\la_i^*+\theta_i)v_i)v_{i}(B-\mu_iI_n)[w_i\xi_i\pa_{v_i}\bar{v}_i\tilde{r}_{i,uv}\tilde{r}_{i}]\\
	&+\sum\limits_{i}\sum\limits_{j\neq i}(v_{i,x}-\mu_i^{-1}(\tilde{\la}_i-\la_i^*+\theta_i)v_i)v_{j}(B-\mu_iI_n)[w_i\xi_i\pa_{v_i}\bar{v}_i\tilde{r}_{i,uv}\tilde{r}_{j}]\\
	&+\sum\limits_{i}(v_{i,x}-\mu_i^{-1}(\tilde{\la}_i-\la_i^*+\theta_i)v_i)v_{i,x}\xi_i(B-\mu_iI_n)[(w_i\pa_{v_iv_i}\bar{v}_i)\tilde{r}_{i,v}+w_i(\pa_{v_i}\bar{v}_i)^2\tilde{r}_{i,vv}]\\
	&+\sum\limits_{i}\sum\limits_{j\neq i}(v_{i,x}-\mu_i^{-1}(\tilde{\la}_i-\la_i^*+\theta_i)v_i)v_{j,x}\xi_i(B-\mu_iI_n)[(w_i\pa_{v_iv_j}\bar{v}_i)\tilde{r}_{i,v}+w_i\pa_{v_i}\bar{v}_i\pa_{v_j}\bar{v}_i\tilde{r}_{i,vv}]\\
	%&+\sum\limits_{i}\sum\limits_{j\neq i}(v_{i,x}-\mu_i^{-1}(\tilde{\la}_i-\la_i^*+\theta_i)v_i)w_{j,x}\xi_i(B-\mu_iI_n)[(w_i\pa_{v_iw_j}\bar{v}_i)\tilde{r}_{i,v}+w_i\pa_{v_i}\bar{v}_i\pa_{w_j}\bar{v}_i\tilde{r}_{i,vv}]\\
	&+\sum\limits_{i}(v_{i,x}-\mu_i^{-1}(\tilde{\la}_i-\la_i^*+\theta_i)v_i)\xi^\p_i\left(\frac{w_i}{v_i}\right)_x(B-\mu_iI_n)[(w_i\pa_{v_i}\bar{v}_i)\tilde{r}_{i,v}+\xi_iw_i\pa_{v_i}\bar{v}_i\bar{v}_i\tilde{r}_{i,vv}]\\
	&-\sum\limits_{i}(v_{i,x}-\mu_i^{-1}(\tilde{\la}_i-\la_i^*+\theta_i)v_i)\theta_i^\p\left(\frac{w_i}{v_i}\right)_x(B-\mu_iI_n)[\xi_iw_i\pa_{v_i}\bar{v}_i\tilde{r}_{i,v\si}]\\
	&+\sum\limits_{i}(v_{i,x}-\mu_i^{-1}(\tilde{\la}_i-\la_i^*+\theta_i)v_i)(B-\mu_iI_n)[\xi_iw_{i,x}\pa_{v_i}\bar{v}_i\tilde{r}_{i,v}]\\
	&=:\sum\limits_{l=1}^{10}\sum\limits_{i}\B^{9,l}_i.
\end{align*}
For the term $K_{10},K_{11},K_{12},K_{13}$ and $K_{14}$ we observe again that we can rewrite this terms in function of $J_5$, $J_6$, $J_7$ and $J_8$
\begin{align*}
	K_{10,x}&=-\sum\limits_{k=1}^{9}\sum\limits_{i}\la_i^*\alpha_{i}^{5,k},\\
	K_{11,x}&=\sum\limits_{i}\left(\frac{w_i}{v_i}\right)_x(w_{i,x}-(w_i/v_i)v_{i,x})(B-\mu_iI_n)[\xi_i^\p\bar{v}_i\tilde{r}_{i,v}-\theta^\p_i\tilde{r}_{i,\si}]\\
	&+\sum\limits_{k=1}^{12}\sum\limits_{i}\frac{w_i}{v_i}\alpha_{i}^{6,k},\\
	K_{12,x}&=-\sum\limits_{k=1}^{12}\sum\limits_{i}\la_i^*\al_{i}^{6,k},\\
	K_{13,x}&=\sum\limits_{i}\mu_i^{-1}(\tilde{\la}_i-\si_i)\left(\frac{w_i}{v_i}\right)_xv_i\xi_i[(\pa_{v_i}\bar{v}_i)v_i-\bar{v}_i]B\tilde{r}_{i,v}+\sum\limits_{k=1}^{17}\sum\limits_{i}\left(\frac{w_i}{v_i}-\la_i^*\right)\alpha_{i}^{7,k},\\
	K_{14,x}&=\sum\limits_{i}\left(\frac{w_i}{v_i}\right)_x(v_i-\xi_i\bar{v}_i)v_i\left[B(u)\tilde{r}_{i,u}\tilde{r}_i+\tilde{r}_i\cdot DB(u)\tilde{r}_i\right]+\sum\limits_{k=1}^{4}\sum\limits_{i}\left(\frac{w_i}{v_i}-\la_i^*\right)\al_{i}^{8,k}.
\end{align*}
Finally we calculate
\begin{align*}
	K_{15,x}&=\sum\limits_{i}\sum\limits_{j\neq i}((w_{i,x}-\la_i^*v_{i,x})v_j+(w_{i}-\la_i^*v_{i})v_{j,x})\left[\tilde{r}_i\cdot DB(u)\tilde{r}_j+B(u)\tilde{r}_{i,u}\tilde{r}_j\right]\\
	&+\sum\limits_{i}\sum\limits_{j\neq i}(w_{i}-\la_i^*v_{i})v_j\left[\tilde{r}_i\cdot DB(u)\tilde{r}_j+B(u)\tilde{r}_{i,u}\tilde{r}_j\right]_x=:\sum\limits_{l=1}^{2}\sum\limits_{i}\B^{15,l}_i.
\end{align*}
Taking derivative in time on \eqref{eqn-w-i} we get
\begin{align*}
	&u_{tt}=\sum\limits_{i}(w_i-\la_i^*v_i)_{t}\tilde{r}_i+\sum\limits_{i}(w_i-\la_i^*v_i)\tilde{r}_{i,t}\\
	&=\sum\limits_{i}(w_{i,t}-\la_i^*v_{i,t})\tilde{r}_i+\sum\limits_{ij}(w_j-\la_j^*v_j)(w_i-\la_i^*v_i)\tilde{r}_{i,u}\tilde{r}_j\\
	&+\sum\limits_{i}\left(v_{i,t}\xi_i\pa_{v_i}\bar{v_i}+\sum\limits_{j\neq i}v_{j,t}\xi_i\pa_{v_j}\bar{v_i}\right)(w_i-\la_i^*v_i)\tilde{r}_{i,v}\\
	&+\sum\limits_{i}\left(w_{i,t}-\frac{w_i}{v_i}v_{i,t}\right)\left(\frac{w_i}{v_i}-\la_i^*\right)[\xi_i^\p\bar{v}_i\tilde{r}_{i,v}-\theta_i^\p\tilde{r}_{i,\si}]\\
	&=\sum\limits_{i}w_{i,t}\left[\tilde{r}_{i}+\frac{w_i}{v_i}\xi_i^\p\bar{v}_i\tilde{r}_{i,v}-\frac{w_i}{v_i}\theta_i^\p\tilde{r}_{i,\si}\right]+\sum\limits_{i}v_{i,t}\left[\left(\xi_iw_i\pa_{v_i}\bar{v}_i-\xi^\p_i\bar{v}_i\frac{w_i^2}{v_i^2}\right)\tilde{r}_{i,v}+\frac{w^2_i}{v^2_i}\theta_i^\p\tilde{r}_{i,\si}\right]\\
	&-\sum\limits_{i}\la_i^*v_{i,t}\left[\tilde{r}_i+\left(v_i \xi_i\pa_{v_i}\bar{v}_i-\xi^\p_i\bar{v}_i\frac{w_i}{v_i}\right)\tilde{r}_{i,v}+\theta_i^\p\frac{w_i}{v_i}\tilde{r}_{i,\si}\right]-\sum\limits_{i}\la_i^*w_{i,t}[\xi_i^\p\bar{v}_i\tilde{r}_{i,v}-\theta_i^\p\tilde{r}_{i,\si}]\\
	&+\sum\limits_{i}\left( \sum\limits_{j\neq i}v_{j,t}\xi_i\pa_{v_j}\bar{v_i} \right)(w_i-\la_i^*v_i)\tilde{r}_{i,v}+\sum\limits_{i}(w_i-\la_i^*v_i)w_i\tilde{r}_{i,u}\tilde{r}_i\\
	&-\sum\limits_{i}\la_i^*v_i(w_i-\la_i^*v_i)\tilde{r}_{i,u}\tilde{r}_i+\sum\limits_{i}\sum\limits_{j\neq i}(w_j-\la_j^*v_j)(w_i-\la_i^*v_i)\tilde{r}_{i,u}\tilde{r}_j.
\end{align*}
Due to the following change of summation
\begin{align*}
	&\sum\limits_{i}\left( \sum\limits_{j\neq i}v_{j,t}\xi_i\pa_{v_j}\bar{v_i}\right)(w_i-\la_i^*v_i)\tilde{r}_{i,v}=\sum\limits_{i}\left( \sum\limits_{j\neq i}v_{i,t}\xi_j\pa_{v_i}\bar{v_j}\right)(w_j-\la_j^*v_j)\tilde{r}_{j,v},
\end{align*}
we deduce then
\begin{align}
	&u_{tt}=\sum\limits_{i}w_{i,t}\left[\tilde{r}_{i}+\frac{w_i}{v_i}\xi_i^\p\bar{v}_i\tilde{r}_{i,v}-\frac{w_i}{v_i}\theta_i^\p\tilde{r}_{i,\si}\right] \nonumber\\
	&+\sum\limits_{i}v_{i,t}\left[\left(w_i\xi_i\pa_{v_i}\bar{v}_i-\xi^\p_i\bar{v}_i\frac{w_i^2}{v_i^2}\right)\tilde{r}_{i,v}+\frac{w^2_i}{v^2_i}\theta_i^\p\tilde{r}_{i,\si}+\sum\limits_{j\neq i}(w_j-\la_j^*v_j)\xi_j(\pa_{v_i}\bar{v}_j)\tilde{r}_{j,v}\right]  \nonumber\\
	&-\sum\limits_{i}\la_i^*v_{i,t}\left[\tilde{r}_i+\left(v_i\xi_i\pa_{v_i}\bar{v}_i-\xi^\p_i\bar{v}_i\frac{w_i}{v_i}\right)\tilde{r}_{i,v}+\theta_i^\p\frac{w_i}{v_i}\tilde{r}_{i,\si}\right] \nonumber\\
	&-\sum\limits_{i}\la_i^*w_{i,t}[\xi_i^\p\bar{v}_i\tilde{r}_{i,v}-\theta_i^\p\tilde{r}_{i,\si}]+\sum\limits_{i}(w_i-\la_i^*v_i)w_i\tilde{r}_{i,u}\tilde{r}_i-\sum\limits_{i}\la_i^*v_i(w_i-\la_i^*v_i)\tilde{r}_{i,u}\tilde{r}_i \nonumber\\
	&+\sum\limits_{i}\sum\limits_{j\neq i}(w_j-\la_j^*v_j)(w_i-\la_i^*v_i)\tilde{r}_{i,u}\tilde{r}_j.\label{eqn-u-tt}
\end{align}

%------------------------------------------------------------------------------------------------------------------------------------

%------------------------------------------------------------------------------------------------------------------------------------
%------------------------------------------------------------------------------------------------------------------------------------
%------------------------------------------------------------------------------------------------------------------------------------
From \eqref{eqn-K-1}, \eqref{eqn-K-2} and \eqref{eqn-K-5} we get (see the Appendix \ref{sectionE} for more details on the computation)
\begin{align}
	&K_{1,x}+K_{2,x}+K_{5,x}-\sum\limits_{i}(w_i-\la_i^*v_i)w_i\tilde{r}_{i,u}\tilde{r}_i-\sum\limits_{i}\sum\limits_{j\neq i}(w_j-\la_j^*v_j)w_i\tilde{r}_{i,u}\tilde{r}_j\nonumber\\
	&=\sum\limits_{i}(\mu_iw_{i,x}-\tilde{\la}_iw_i)_x\left[\tilde{r}_{i}+\frac{w_i}{v_i}\xi_i^\p\bar{v}_i\tilde{r}_{i,v}-\frac{w_i}{v_i}\theta_i^\p\tilde{r}_{i,\si}\right]\nonumber\\
	&+\sum\limits_{i}(\mu_iv_{i,x}-\tilde{\la}_iv_i)_x\left[\left(w_i\xi_i\pa_{v_i}\bar{v}_i-\xi^\p_i\bar{v}_i\frac{w_i^2}{v_i^2}\right)\tilde{r}_{i,v}+\frac{w^2_i}{v^2_i}\theta_i^\p\tilde{r}_{i,\si}\right]\nonumber\\
	&+\sum\limits_{i}\sum\limits_{j\neq i}[(\mu_iw_{i,x}-\tilde{\la}_iw_i)v_j-(w_j-\la_j^*v_j)w_i]\tilde{r}_{i,u}\tilde{r}_j\nonumber\\
	&+\sum\limits_{i}[\mu_i(w_{i,x}v_i-w_iv_{i,x})+(\mu_iv_{i,x}-(\tilde{\la}_i-\la_i^*)v_i-w_i)w_i]\tilde{r}_{i,u}\tilde{r}_i\nonumber\\
	&+\sum\limits_{i}(\mu_iw_{i,x}-(\tilde{\la}_i-\la_i^*+\theta_i)w_i)v_{i,x}\xi_i(\pa_{v_i}\bar{v}_i)\tilde{r}_{i,v}+\sum\limits_{i}\sum\limits_{j\neq i}(\mu_iw_{i,x}-\tilde{\la}_iw_i )v_{j,x}\xi_i(\pa_{v_j}\bar{v}_i)\tilde{r}_{i,v}\nonumber\\
	&+2\sum\limits_{i}\mu_i(w_{i,x}-(w_i/v_i)v_{i,x})\xi^\p_i\bar{v}_i\left(\frac{w_i}{v_i}\right)_x\tilde{r}_{i,v}+\sum\limits_{i}\mu_i(w_{i,x}-(w_i/v_i)v_{i,x})\xi^{\p\p}_i\bar{v}_i\frac{w_i}{v_i}\left(\frac{w_i}{v_i}\right)_x\tilde{r}_{i,v}\nonumber\\
	&+\sum\limits_{i}\mu_i(w_{i,x}-(w_i/v_i)v_{i,x})v_{i,x}\xi^\p_i(\pa_{v_i}\bar{v}_i)\frac{w_i}{v_i}\tilde{r}_{i,v}+\sum\limits_{i}\sum\limits_{j\neq i}\mu_i(w_{i,x}-(w_i/v_i)v_{i,x})v_{j,x}\xi_i^\p(\pa_{v_j}\bar{v}_i)\frac{w_i}{v_i}\tilde{r}_{i,v}\nonumber\\
	&+\sum_i(\mu_i v_{i,x}-v_i(\tilde{\lambda}_i-\lambda_i^*)-\theta_i v_i)\xi_i w_iv_{i,x}\pa_{v_iv_i}\bar{v}_i\tilde{r}_{i,v}\nonumber\\
	&+\sum_{i\ne j}(\mu_i v_{i,x}-v_i(\tilde{\lambda}_i-\lambda_i^*)-\theta_i v_i) \xi_iw_iv_{j,x}\pa_{v_jv_i}\bar{v}_i\tilde{r}_{i,v}\nonumber\\
	&+\sum_{i}(\mu_i v_{i,x}-v_i(\tilde{\lambda}_i-\lambda_i^*)-\theta_i v_i) \xi_iw_{i,x}\pa_{v_i}\bar{v}_i\tilde{r}_{i,v}\nonumber\\
	&+\sum_{i}(\mu_i v_{i,x}-v_i(\tilde{\lambda}_i-\lambda_i^*)-\theta_i v_i) w_i\xi_i^\p\left(\frac{w_i}{v_i}\right)_x\pa_{v_i}\bar{v}_i\tilde{r}_{i,v}-\sum_{i}\theta^\p_i\pa_{v_i}\bar{v}_iv_iw_i\left(\frac{w_i}{v_i}\right)_x\xi_i \tilde{r}_{i,v}\nonumber\\
	&-2\sum\limits_{i}\mu_i(w_{i,x}-(w_i/v_i)v_{i,x})\theta_i^\p\left(\frac{w_i}{v_i}\right)_x\tilde{r}_{i,\si}-\sum\limits_{i}\mu_i(w_{i,x}-(w_i/v_i)v_{i,x})\frac{w_i}{v_i}\left(\frac{w_i}{v_i}\right)_x\theta_i^{\p\p} \tilde{r}_{i,\si}\nonumber\\
	&+\sum_{i}\xi_i w_i v_i\pa_{v_i}\bar{v}_i (\mu_i v_{i,x}-v_i(\tilde{\lambda}_i-\lambda_i^*)-\theta_i v_i) \tilde{r}_{i,uv}\tilde{r}_i\nonumber\\
	&+\sum_{i\ne j}\xi_i w_i v_j\pa_{v_i}\bar{v}_i (\mu_i v_{i,x}-v_i(\tilde{\lambda}_i-\lambda_i^*)-\theta_i v_i) \tilde{r}_{i,uv}\tilde{r}_j\nonumber\\
	&+\sum\limits_{i}\mu_i(w_{i,x}-(w_i/v_i)v_{i,x})w_i \xi_i^\p\bar{v}_i\tilde{r}_{i,uv}\tilde{r}_i+\sum\limits_{i}\sum\limits_{j\neq i}\mu_i(w_{i,x}-(w_i/v_i)v_{i,x})v_j\frac{w_i}{v_i}\xi_i^\p\bar{v}_i\tilde{r}_{i,uv}\tilde{r}_i\nonumber\\
	&+\sum_i \xi_i^2w_iv_{i,x}(\pa_{v_i}\bar{v}_i)^2 (\mu_iv_{i,x}-v_i(\tilde{\lambda}_i-\lambda_i^*)-\theta_i v_i)\tilde{r}_{i,vv}\nonumber\\
	&+\sum_{i\ne j}\xi_i^2v_{j,x}w_i \pa_{v_i}\bar{v}_i\pa_{v_j}\bar{v}_i(\mu_iv_{i,x}-v_i(\tilde{\lambda}_i-\lambda_i^*)-\theta_i v_i)\tilde{r}_{i,vv}\nonumber\\
	&+\sum_{i}\xi_i' \bar{v}_i w_i \xi_i \pa_{v_i}\bar{v}_i\left(\frac{w_i}{v_i}\right)_x(\mu_iv_{i,x}-v_i(\tilde{\lambda}_i-\lambda_i^*)-\theta_i v_i)\tilde{r}_{i,vv}\nonumber\\
	&+\sum\limits_{i}\mu_i(w_{i,x}-(w_i/v_i)v_{i,x})v_{i,x}\frac{w_i}{v_i}\xi_i^\p\xi_i\pa_{v_i}\bar{v}_i\bar{v}_i\tilde{r}_{i,vv}\nonumber\\
	&+\sum\limits_{i}\sum\limits_{j\neq i}\mu_i(w_{i,x}-(w_i/v_i)v_{i,x})v_{j,x}\frac{w_i}{v_i}\xi_i^\p\xi_i\bar{v}_i\pa_{v_j}\bar{v}_i\tilde{r}_{i,vv}\nonumber\\
	&+\sum\limits_{i}\mu_i(w_{i,x}-(w_i/v_i)v_{i,x})\frac{w_i}{v_i}(\bar{v}_i\xi_i^\p)^2\left(\frac{w_i}{v_i}\right)_x\tilde{r}_{i,vv}\nonumber\\
	&-\sum\limits_{i}\mu_i(w_{i,x}-(w_i/v_i)v_{i,x})\frac{w_i}{v_i}(\bar{v}_i\xi_i^\p)\theta_i^\p\left(\frac{w_i}{v_i}\right)_x\tilde{r}_{i,v\si}\nonumber\\
	&-\sum\limits_{i}\mu_i(w_{i,x}-(w_i/v_i)v_{i,x}) \xi_i\frac{w_i}{v_i}v_{i,x}\theta_i^\p \pa_{v_i}\bar{v}_i\tilde{r}_{i,v\si}\nonumber\\
	&-\sum\limits_{i}\sum\limits_{j\neq i}\mu_i(w_{i,x}-(w_i/v_i)v_{i,x})\xi_i \frac{w_i}{v_i}v_{j,x}\theta_i^\p \pa_{v_j}\bar{v}_i\tilde{r}_{i,v\si}\nonumber\\
	&-\sum\limits_{i}\mu_i(w_{i,x}-(w_i/v_i)v_{i,x}) \xi_i^\p\frac{w_i}{v_i}\left(\frac{w_i}{v_i}\right)_x\theta_i^\p \bar{v}_i\tilde{r}_{i,v\si}\nonumber\\
	&-\sum_i w_i\xi_i\pa_{v_i}\bar{v}_i\theta'_i\left(\frac{w_i}{v_i}\right)_x(\mu_iv_{i,x}-v_i(\tilde{\lambda}_i-\lambda_i^*)-\theta_i v_i)\tilde{r}_{i,v\si}\nonumber\\
	&-\sum\limits_{i}\mu_i(w_{i,x}-(w_i/v_i)v_{i,x}) w_i\theta_i^\p \tilde{r}_{i,u\si}\tilde{r}_i-\sum\limits_{i}\sum\limits_{j\neq i}\mu_i(w_{i,x}-(w_i/v_i)v_{i,x}) (w_i/v_i)v_j\theta_i^\p \tilde{r}_{i,u\si}\tilde{r}_j\nonumber\\
	&+\sum\limits_{i}\mu_i(w_{i,x}-(w_i/v_i)v_{i,x}) \frac{w_i}{v_i}\left(\frac{w_i}{v_i}\right)_x(\theta_i^\p)^2 \tilde{r}_{i,\si\si}\nonumber\\
	&=:\sum\limits_{i}(\mu_iw_{i,x}-\tilde{\la}_iw_i)_x\left[\tilde{r}_{i}+\frac{w_i}{v_i}\xi_i^\p\bar{v}_i\tilde{r}_{i,v}-\frac{w_i}{v_i}\theta_i^\p\tilde{r}_{i,\si}\right]\nonumber\\
	&+\sum\limits_{i}(\mu_iv_{i,x}-\tilde{\la}_iv_i)_x\left[\left(w_i\xi_i\pa_{v_i}\bar{v}_i-\xi^\p_i\bar{v}_i\frac{w_i^2}{v_i^2}\right)\tilde{r}_{i,v}+\frac{w^2_i}{v^2_i}\theta_i^\p\tilde{r}_{i,\si}\right]+\sum\limits_{l=1}^{33}\sum\limits_{i}\widetilde{\beta}_{i}^{1,l}.\label{equaK1x}
\end{align}
From \eqref{K3K4}, \eqref{K6} and \eqref{J1J2J3} we deduce that
\begin{align*}
	&K_{3,x}+K_{4,x}+K_{6,x}+\sum\limits_{i}\la_i^*(w_i-\la_i^*v_i)v_{i}\tilde{r}_{i,u}\tilde{r}_i+\sum\limits_{i}\sum\limits_{j\neq i}\la_i^*(w_j-\la_j^*v_j)v_{i}\tilde{r}_{i,u}\tilde{r}_j\\
	&=-\sum\limits_{i}\la_i^*(\mu_iv_{i,x}-\tilde{\la}_iv_i)_x\hat{r}_i-\sum\limits_{i}\la_i^*(\mu_iw_{i,x}- \tilde{\la}_iw_i)_xr_i^\dagger-\sum\limits_{l=1}^{{\color{black}{27}}}\sum\limits_{i}\la_i^*\widetilde{\al}_{i}^{1,l}.
	\end{align*}
We set
\begin{align}
	r^\ddagger_{i}&=\tilde{r}_{i}+\frac{w_i}{v_i}\xi_i^\p\bar{v}_i\tilde{r}_{i,v}-\frac{w_i}{v_i}\theta_i^\p\tilde{r}_{i,\si},\nonumber\\
	r_i^\#&=\left(w_i\pa_{v_i}\bar{v}_i\xi_i-\xi^\p_i\bar{v}_i\frac{w_i^2}{v_i^2}\right)\tilde{r}_{i,v}+\frac{w^2_i}{v^2_i}\theta_i^\p\tilde{r}_{i,\si},\nonumber\\
	r_i^\clubsuit&=-\la_i^*\hat{r}_i=-\la_i^*\left[\tilde{r}_{i}+(v_i\xi_i\pa_{v_i}\bar{v}_i-(w_i/v_i)\xi_i^\p\bar{v}_i) \tilde{r}_{i,v}+(w_i/v_i)\theta^\p_i\tilde{r}_{i,\si}\right],\nonumber\\
	r_i^\diamondsuit&=-\la_i^*r_i^\dagger=-\la_i^*\left[\xi_i^\p\bar{v}_i\tilde{r}_{i,v}-\theta^\p_i\tilde{r}_{i,\si}\right].\label{defAo1}
\end{align}
Therefore, we can write
\begin{align}
	&K_{1,x}+K_{2,x}+K_{5,x}-\sum\limits_{i}(w_i-\la_i^*v_i)w_i\tilde{r}_{i,u}\tilde{r}_i-\sum\limits_{i}\sum\limits_{j\neq i}(w_j-\la_j^*v_j)w_i\tilde{r}_{i,u}\tilde{r}_j \nonumber\\
	&=\sum\limits_{i}(\mu_iw_{i,x}-\tilde{\la}_iw_i)_xr_i^\ddagger+\sum\limits_{i}(\mu_iv_{i,x}-\tilde{\la}_iv_i)_xr_i^\#+\sum\limits_{l=1}^{33}\sum\limits_{i}\widetilde{\B}_{i}^{1,l},\label{eqn-K-1-2-5}
\end{align}
and 
\begin{align}
	&K_{3,x}+K_{4,x}+K_{6,x}+\sum\limits_{i}\la_i^*(w_i-\la_i^*v_i)v_{i}\tilde{r}_{i,u}\tilde{r}_i+\sum\limits_{i}\sum\limits_{j\neq i}\la_i^*(w_j-\la_j^*v_j)v_{i}\tilde{r}_{i,u}\tilde{r}_j \nonumber\\
	&=\sum\limits_{i} (\mu_iv_{i,x}-\tilde{\la}_iv_i)_xr_i^\clubsuit+\sum\limits_{i} (\mu_iw_{i,x}- \tilde{\la}_iw_i)_xr_i^\diamondsuit-\sum\limits_{l=1}^{27}\sum\limits_{i}\la_i^*\widetilde{\al}_{i}^{1,l}. \label{eqn-K-3-4-6}
\end{align}
From \eqref{eqn-u-tt} we have
\begin{align}
	u_{tt}&=\sum\limits_{i}w_{i,t}r^\ddagger_{i}
	+\sum\limits_{i}v_{i,t}r_i^\#+\sum\limits_{i}v_{i,t}\sum\limits_{j\neq i}(w_j-\la_j^*v_j)\xi_j (\pa_{v_i}\bar{v}_j)\tilde{r}_{j,v}  \nonumber\\
	&+\sum\limits_{i}v_{i,t}r_i^\clubsuit +\sum\limits_{i}w_{i,t}r_i^\diamondsuit +\sum\limits_{i}(w_i-\la_i^*v_i)w_i\tilde{r}_{i,u}\tilde{r}_i-\sum\limits_{i}\la_i^*v_i(w_i-\la_i^*v_i)\tilde{r}_{i,u}\tilde{r}_i \nonumber\\
	&+\sum\limits_{i}\sum\limits_{j\neq i}(w_j-\la_j^*v_j)(w_i-\la_i^*v_i)\tilde{r}_{i,u}\tilde{r}_j.\label{eqn-u-tt1}
\end{align}
Observe that
\begin{equation}\label{eqn-ux-A-comm}
	u_x\cdot DA u_t-u_t\cdot DA u_x=\sum\limits_{j\neq i}(w_i-\la_i^*v_i)v_j\left[\tilde{r}_j\cdot DA\tilde{r}_i-\tilde{r}_i\cdot DA\tilde{r}_j\right].
\end{equation}
Since 
\begin{equation}
	u_{tt}=\left[Bu_{tx}+u_t\cdot DBu_x-Au_t\right]_x+u_x\cdot DA u_t-u_t\cdot DA u_x,
\end{equation}
now combining \eqref{eqn-B-u-tx}, \eqref{eqn-K-1-2-5}, \eqref{eqn-K-3-4-6}, \eqref{eqn-ux-A-comm} and \eqref{eqn-u-tt1} we can obtain
\begin{align}
	&\sum\limits_{i}(w_{i,t}+(\tilde{\la}_iw_i)_x-(\mu_iw_{i,x})_x)r_i^\ddagger+\sum\limits_{i}(v_{i,t}+(\tilde{\la}_iv_i)_x-(\mu_iv_{i,x})_x)r^\#_i     \nonumber\\	
	&+\sum\limits_{i}(v_{i,t}+(\tilde{\la}_iv_i)_x-(\mu_iv_{i,x})_x)r_i^\clubsuit+\sum\limits_{i} (w_{i,t}+(\tilde{\la}_iw_i)_x-(\mu_iw_{i,x})_x)r_i^\diamondsuit      \nonumber\\
	&=\sum\limits_{l=1}^{33}\sum\limits_{i} \widetilde{\B}_i^{1,l}-\sum\limits_{l=1}^{27}\sum\limits_{i} \la_i^*\widetilde{\al}_i^{1,l}+\sum\limits_{j=7}^{15}K_{j,x} +\sum\limits_{j\neq i}(w_i-\la_i^*v_i)v_j\left[\tilde{r}_j\cdot DA\tilde{r}_i-\tilde{r}_i\cdot DA\tilde{r}_j\right]  \nonumber\\
	&-\sum\limits_{i}v_{i,t}\sum\limits_{j\neq i}(w_j-\la_j^*v_j)\xi_j(\pa_{v_i}\bar{v}_j)\tilde{r}_{j,v}.\label{eqn-K-all-1}
\end{align}
Due to the choice of the support of $\xi_i$ and $\theta_i$, it follows since $\psi_{ij}(u,\xi_i\bar{v}_i,\sigma_i)=O(1)\xi_i\bar{v}_i$ that $\theta_i v_i=w_i$ when this quantity is multiplied by  $\psi_{ij}(u,\xi_i\bar{v}_i,\sigma_i)$ we deduce then that 
\begin{align*}
	&\sum\limits_{i}\B^{8,1}_i=\sum\limits_{i}(w_{i,x}-\mu_i^{-1}(\tilde{\la}_i-\la_i^*+\theta_i)w_i)_x(B-\mu_iI_n)\tilde{r}_{i}\\
	&=\sum\limits_{i}((\theta_iv_i)_{x}-\mu_i^{-1}(\tilde{\la}_i-\la_i^*+\theta_i)\theta_iv_i)_x\sum\limits_{j\neq i}(\mu_j-\mu_i)\psi_{ij} r_{j}\\
	&=\sum\limits_{i}\theta_i^\p\left(w_{i,x}-\frac{w_i}{v_i}v_{i,x}\right)_x\sum\limits_{j\neq i}(\mu_j-\mu_i)\psi_{ij} r_{j}+\sum\limits_{i}\theta_i^{\p\p}v_i\left(\frac{w_i}{v_i}\right)_x^2\sum\limits_{j\neq i}(\mu_j-\mu_i)\psi_{ij}  {r}_{j}\\
	&+\sum\limits_{i}(\theta_i(v_{i,x}-\mu_i^{-1}(\tilde{\la}_i-\la_i^*+\theta_i)v_i))_x\sum\limits_{j\neq i}(\mu_j-\mu_i)\psi_{ij}  {r}_{j}.
\end{align*}
Interchanging the indices in the first term we obtain combining the fact that $\psi_{i,j}=O(1)\xi_i\bar{v}_i$ and  $\theta''_i\xi_i=0$ due to the fact that $\theta''(x)=0$ when $x\in\mbox{supp}\xi$
\begin{align}
	\sum\limits_{i}\B^{8,1}_i&=\sum\limits_{i}\sum\limits_{j\neq i}\left[\theta_j^\p\left(w_{j,x}-\frac{w_j}{v_j}v_{j,x}\right)_x(\mu_i-\mu_j)\psi_{ji}\right] {r}_{i} \nonumber\\
	%&+\sum\limits_{i}\theta_i^{\p\p}v_i\left(\frac{w_i}{v_i}\right)_x^2\sum\limits_{j\neq i}(\mu_j-\mu_i)\psi_{ij} r_{j} \nonumber\\
	&+\sum\limits_{i}\theta_i(v_{i,x}-\mu_i^{-1}(\tilde{\la}_i-\la_i^*+\theta_i)v_i)_x\sum\limits_{j\neq i}(\mu_j-\mu_i)\psi_{ij}  r_{j}  \nonumber\\
	&+\sum\limits_{i}\theta_i^\p\left(\frac{w_i}{v_i}\right)_x(v_{i,x}-\mu_i^{-1}(\tilde{\la}_i-\la_i^*+\theta_i)v_i)\sum\limits_{j\neq i}(\mu_j-\mu_i)\psi_{ij}  r_{j} \nonumber\\
	&=:\sum\limits_{i}\sum\limits_{j\neq i}\left[\theta_j^\p\left(w_{j,x}-\frac{w_j}{v_j}v_{j,x}\right)_x(\mu_i-\mu_j)\psi_{ji}\right] r_{i}+\sum\limits_{k=1}^{2}\sum\limits_{i}\B^{8,1,k}_i.\label{eqn-B-9-1}
\end{align}
Next we rearrange the first term in $K_{12,x}$ which corresponds to $-\sum_i\la_i^*\al_{i}^{6,1}$,
\begin{align}
	-\sum\limits_{i}\la_i^*\al_{i}^{6,1}
	&=-\sum\limits_{i}\la_i^*\left(w_{i,x}-\frac{w_i}{v_i}v_{i,x}\right)_x(B-\mu_iI_n)[\xi_i^\p\bar{v}_i\tilde{r}_{i,v}-\theta^\p_i\tilde{r}_{i,\si}]  \nonumber\\
	&=-\sum\limits_{i}\la_i^*\left(w_{i,x}-\frac{w_i}{v_i}v_{i,x}\right)_x\sum\limits_{j\neq i}(\mu_j-\mu_i)\Big[\xi_i^\p\bar{v}_i\psi_{ij,v}-\theta^\p_i\psi_{ij,\si}\Big]r_j \nonumber\\
	&=-\sum\limits_{i}\sum\limits_{j\neq i}\left(w_{j,x}-\frac{w_j}{v_j}v_{j,x}\right)_x\la_j^*(\mu_i-\mu_j)\Big[\xi_j^\p\bar{v}_j\psi_{ji,v}-\theta^\p_j\psi_{ji,\si}\Big]r_i. \label{beta-11-1-a}
\end{align}
Recall from \eqref{def:b-ij} that $	b_{ij}= \xi_j^\p\bar{v}_j\psi_{ji,v} -\theta^\p_j\psi_{ji,\si}$. We rewrite \eqref{beta-11-1-a} as
\begin{align}
	-\sum\limits_{i}\la_i^*\al_{i}^{6,1}&=-\sum\limits_{i}\sum\limits_{j\neq i}\left(w_{j,x}-\frac{w_j}{v_j}v_{j,x}\right)_x(\mu_i-\mu_j)b_{ij}\la_i^*r_i  \nonumber\\
	&+\sum\limits_{i}\sum\limits_{j\neq i}\left(w_{j,x}-\frac{w_j}{v_j}v_{j,x}\right)_x(\la_i^*-\la_j^*)(\mu_i-\mu_j)\Big[\xi_j^\p\bar{v}_j\psi_{ji,v}-\theta^\p_j\psi_{ji,\si}\Big]r_i . \label{eqn-B-11-1}
\end{align}
Similarly we have considering the first term $\frac{w_i}{v_i}\al_{i}^{6,1}$ in $K_{11,x}$
\begin{align}
	\sum\limits_{i}\frac{w_i}{v_i}\al_{i}^{6,1}
	&=\sum\limits_{i}\left(w_{i,x}-\frac{w_i}{v_i}v_{i,x}\right)_x\frac{w_i}{v_i}(B-\mu_iI_n)[\xi_i^\p\bar{v}_i\tilde{r}_{i,v}-\theta^\p_i\tilde{r}_{i,\si}]  \nonumber\\
	&=\sum\limits_{i}\sum\limits_{j\neq i}\left(w_{j,x}-\frac{w_j}{v_j}v_{j,x}\right)_x\frac{w_j}{v_j}(\mu_i-\mu_j)\Big[\xi_j^\p\bar{v}_j\psi_{ji,v}-\theta^\p_j\psi_{ji,\si}\Big]r_i. \label{beta-11-1-ah}
\end{align}
Now, we can combine \eqref{eqn-B-9-1}, \eqref{eqn-B-11-1} and \eqref{beta-11-1-ah}as follows
\begin{align*}
	&\sum\limits_{i}\B^{8,1}_i-\sum\limits_{i}\la_i^*\al_{i}^{6,1}+\sum\limits_{i}\frac{w_i}{v_i}\al_{i}^{6,1}\\
	&=\sum\limits_{i}\sum\limits_{j\neq i}\left[\theta_j^\p\left(w_{j,x}-\frac{w_j}{v_j}v_{j,x}\right)_x(\mu_i-\mu_j)\psi_{ji}(u,\bar{v}_j,\si_j)\right] {r}_{i}\\
	&+\sum\limits_{i}\sum\limits_{j\neq i}\left(w_{j,x}-\frac{w_j}{v_j}v_{j,x}\right)_x\left[(\la_i^*-\la_j^*)+\frac{w_j}{v_j}\right](\mu_i-\mu_j)\Big[\xi_j^\p\bar{v}_j\psi_{ji,v}-\theta^\p_j\psi_{ji,\si}\Big]r_i\\
	&-\sum\limits_{i}\sum\limits_{j\neq i}\left(w_{j,x}-\frac{w_j}{v_j}v_{j,x}\right)_x(\mu_i-\mu_j)b_{ij}\la_i^*r_i+\sum\limits_{k=1}^{2}\sum\limits_{i}\B^{8,1,k}_i.
\end{align*}
Then we set
\begin{align}
	\hat{b}_{ij}&=\theta_j^\p  \psi_{ji} +\left[(\la_i^*-\la_j^*)+\frac{w_j}{v_j}\right]\big[\xi_j^\p\bar{v}_j\psi_{ji,v}-\theta^\p_j\psi_{ji,\si}\big]. \label{def:hat-b-ija}
\end{align}
As previously since $\theta^\p=1$ on $\mbox{supp}\xi$, we deduce that
\begin{align}
	\hat{b}_{ij}&=\psi_{ji}+\left[(\la_i^*-\la_j^*)+\frac{w_j}{v_j}\right]\big[\xi_j^\p\bar{v}_j\psi_{ji,v}-\psi_{ji,\si}\big]. \label{def:hat-b-ij}
\end{align}
Hence, we may write
\begin{align}
	&\sum\limits_{i}\B^{8,1}_i+\sum\limits_{i}\la_i^*\al_{i}^{6,1}+\sum\limits_{i}\frac{w_i}{v_i}\al_{i}^{6,1}\nonumber\\
	&=\sum\limits_{i}\sum\limits_{j\neq i}\left(w_{j,x}-\frac{w_j}{v_j}v_{j,x}\right)_x(\mu_i-\mu_j)\hat{b}_{ij}r_i -\sum\limits_{i}\sum\limits_{j\neq i}\left(w_{j,x}-\frac{w_j}{v_j}v_{j,x}\right)_x(\mu_i-\mu_j)b_{ij}\la_i^*r_i \nonumber\\
	&+\sum\limits_{k=1}^{3}\sum\limits_{i}\B^{9,1,k}_i \nonumber\\
	&=\sum\limits_{i}\sum\limits_{j\neq i}\left(w_{j,x}-\frac{w_j}{v_j}v_{j,x}\right)_x(\mu_i-\mu_j)\hat{b}_{ij}r^\ddagger_i +\sum\limits_{i}\sum\limits_{j\neq i}\left(w_{j,x}-\frac{w_j}{v_j}v_{j,x}\right)_x(\mu_i-\mu_j)b_{ij}r_i^\#  \nonumber\\
	&+\sum\limits_{i}\sum\limits_{j\neq i}\left(w_{j,x}-\frac{w_j}{v_j}v_{j,x}\right)_x(\mu_i-\mu_j)b_{ij}r_i^\clubsuit +\sum\limits_{i}\sum\limits_{j\neq i}\left(w_{j,x}-\frac{w_j}{v_j}v_{j,x}\right)_x(\mu_i-\mu_j)\hat{b}_{ij}r_i^\diamondsuit  \nonumber\\
	&+\sum\limits_{l=1}^{5}\sum\limits_{i}\B^{16,1,l}_i+\sum\limits_{k=1}^{3}\sum\limits_{i}\B^{8,1,k}_i,  \label{eqn-B-8-9-11-all}
\end{align}
where $\B^{16,1,l}_i$ are defined as follows
\begin{align*}
	\B^{16,1,1}&=-\sum\limits_{j\neq i}\left(w_{j,x}-\frac{w_j}{v_j}v_{j,x}\right)_x(\mu_i-\mu_j)\hat{b}_{ij}    \sum\limits_{k\neq i}\psi_{ik}r_k,\\
	\B^{16,1,2}&=-\sum\limits_{j\neq i}\left(w_{j,x}-\frac{w_j}{v_j}v_{j,x}\right)_x(\mu_i-\mu_j)\hat{b}_{ij}  (\frac{w_i}{v_i}-\la_i^*)\xi_i^\p\bar{v}_i\tilde{r}_{i,v},\\
	\B^{16,1,3}&=\sum\limits_{j\neq i}\left(w_{j,x}-\frac{w_j}{v_j}v_{j,x}\right)_x(\mu_i-\mu_j)\hat{b}_{ij}  (\frac{w_i}{v_i}-\la_i^*)\theta_i^\p\tilde{r}_{i,\si},\\
	\B^{16,1,4}&=\sum\limits_{j\neq i}\left(w_{j,x}-\frac{w_j}{v_j}v_{j,x}\right)_x(\mu_i-\mu_j){b}_{ij}   \la_i^* \sum\limits_{k\neq i}\psi_{ik}r_k,\\
	\B^{16,1,5}&=\sum\limits_{j\neq i}\left(w_{j,x}-\frac{w_j}{v_j}v_{j,x}\right)_x(\mu_i-\mu_j){b}_{ij}   (\la_i^*-\frac{w_i}{v_i}) \Big((v_i\xi_i\pa_{v_i}\bar{v}_i-\frac{w_i}{v_i}\xi_i^\p\bar{v}_i) \tilde{r}_{i,v}\\
	&+\frac{w_i}{v_i}\theta^\p_i\tilde{r}_{i,\si}\Big).
\end{align*}
Plugging \eqref{eqn-B-8-9-11-all} in \eqref{eqn-K-all-1} we obtain
\begin{align}
	&\sum\limits_{i}\left(w_{i,t}+(\tilde{\la}_iw_i)_x-(\mu_iw_{i,x})_x-\sum\limits_{j\neq i}\hat{b}_{ij}(\mu_i-\mu_j)\left(w_{j,x}-\frac{w_j}{v_j}v_{j,x}\right)_x\right)r_i^\ddagger  \nonumber\\
	&+\sum\limits_{i}\left (v_{i,t}+(\tilde{\la}_iv_i)_x-(\mu_iv_{i,x})_x-\sum\limits_{j\neq i}b_{ij}(\mu_i-\mu_j)\left(w_{j,x}-\frac{w_j}{v_j}v_{j,x}\right)_x\right)r^\#_i   \nonumber\\	
	&+\sum\limits_{i}\left(v_{i,t}+(\tilde{\la}_iv_i)_x-(\mu_iv_{i,x})_x-\sum\limits_{j\neq i}b_{ij}(\mu_i-\mu_j)\left(w_{j,x}-\frac{w_j}{v_j}v_{j,x}\right)_x\right)r_i^\clubsuit   \nonumber\\
	&+\sum\limits_{i} \left (w_{i,t}+(\tilde{\la}_iw_i)_x-(\mu_iw_{i,x})_x-\sum\limits_{j\neq i}\hat{b}_{ij}(\mu_i-\mu_j)\left(w_{j,x}-\frac{w_j}{v_j}v_{j,x}\right)_x\right)r_i^\diamondsuit \nonumber\\
	&=\sum\limits_{l=1}^{33}\sum\limits_{i} \widetilde{\B}_i^{1,l}-\sum\limits_{l=1}^{27}\sum\limits_{i} \la_i^*\widetilde{\al}_i^{1,l}+K_{7,x}+\sum\limits_{l=2}^{9}\sum\limits_{i}\B^{8,l}_i+ \sum\limits_{j=9}^{10}K_{j,x}  \nonumber\\
	&+\sum\limits_{k=2}^{12}\sum\limits_{i}(\frac{w_i}{v_i}-\la_i^*)\alpha_{i}^{6,k}
	+\sum\limits_{j=13}^{15}K_{j,x}+\sum\limits_{l=1}^{5}\sum\limits_{i}\B^{16,1,l}_i  \nonumber\\
	&+\sum\limits_{j\neq i} (w_i-\la_i^*v_i)v_j\left[\tilde{r}_j\cdot DA\tilde{r}_i-\tilde{r}_i\cdot DA\tilde{r}_j\right]  \nonumber\\
	&-\sum\limits_{i}v_{i,t}\sum\limits_{j\neq i}(w_j-\la_j^*v_j)\xi_j(\pa_{v_i}\bar{v}_j)\tilde{r}_{j,v}+\sum\limits_{k=1}^{2}\sum\limits_{i}\B^{8,1,k}\nonumber\\
	&+\sum\limits_{i}\left(\frac{w_i}{v_i}\right)_x(w_{i,x}-(w_i/v_i)v_{i,x})(B-\mu_iI_n)[\xi_i^\p\bar{v}_i\tilde{r}_{i,v}-\theta^\p_i\tilde{r}_{i,\si}].
	\label{eqn-remainder-2}
\end{align}
with
\begin{align}
\sum\limits_{k=1}^{2}\sum\limits_{i}\B^{8,1,k}& =
	\sum\limits_{i}\theta_i(v_{i,x}-\mu_i^{-1}(\tilde{\la}_i-\la_i^*+\theta_i)v_i)_x\sum\limits_{j\neq i}(\mu_j-\mu_i)\psi_{ij}  r_{j}  \nonumber\\
	&+\sum\limits_{i}\theta_i^\p\left(\frac{w_i}{v_i}\right)_x(v_{i,x}-\mu_i^{-1}(\tilde{\la}_i-\la_i^*+\theta_i)v_i)\sum\limits_{j\neq i}(\mu_j-\mu_i)\psi_{ij}  r_{j} .
\end{align}

%---------------------------------------

\subsection{Derivation of the remainder terms}
We set 
\begin{align}
	\phi_i&:=v_{i,t}+(\tilde{\la}_iv_i)_x-(\mu_iv_{i,x})_x{\color{black}{-}}\sum\limits_{j\neq i}(\mu_i-\mu_j)b_{ij}\left(w_{j,x}-\frac{w_j}{v_j}v_{j,x}\right)_x,\label{7.51}\\
	\psi_i&:=w_{i,t}+(\tilde{\la}_iw_i)_x-(\mu_iw_{i,x})_x{\color{black}{-}}\sum\limits_{j\neq i}(\mu_i-\mu_j)\hat{b}_{ij}\left(w_{j,x}-\frac{w_j}{v_j}v_{j,x}\right)_x.\label{7.52}
\end{align}
From \eqref{eqn-remainder-1} and  \eqref{eqn-remainder-2}  we get under the assumption of section \ref{section6.1} for $t\in[\hat{t},T]$
\begin{align}
	\sum\limits_{i}\phi_i\hat{r}_i +\sum\limits_{i}\psi_ir^\dagger_i&=\mathcal{J}_1,\\
	\sum\limits_{i}\psi_ir_i^\ddagger +\sum\limits_{i}\phi_ir_i^\#+ \sum\limits_{i}\phi_ir_i^\clubsuit +\sum\limits_{i}\psi_ir_i^\diamondsuit &=\mathcal{J}_2,
\end{align}
where $\mathcal{J}_1$ and $\mathcal{J}_2$ are defined as follows 
\begin{align}
	\mathcal{J}_1&=\sum\limits_{l=1}^{27}\sum\limits_{i}\widetilde{\al}_{i}^{1,l}+\sum\limits_{k=4}^{5}\sum\limits_{l=1}^{N_k}\sum\limits_{i}\al^{k,l}_i+\sum_{l=2}^{N_6}\sum_i\al_i^{6,l}+\sum\limits_{i}\sum\limits_{l=1}^{3}\al_i^{6,1,l}+\sum\limits_{k=7}^{10}\sum\limits_{l=1}^{N_k}\sum\limits_{i}\al^{k,l}_i+\sum\limits_{i}\al_{i}^{11},
	\label{J1}\\
	\mathcal{J}_2&=\sum\limits_{l=1}^{33}\sum\limits_{i} \widetilde{\B}_i^{1,l}-\sum\limits_{l=1}^{27}\sum\limits_{i} \la_i^*\widetilde{\al}_i^{1,l}+K_{7,x}+\sum\limits_{l=2}^{9}\sum\limits_{i}\B^{8,l}_i+ \sum\limits_{j=9}^{10}K_{j,x} +\sum\limits_{k=2}^{12}\sum\limits_{i}(\frac{w_i}{v_i}-\la_i^*)\alpha_{i}^{6,k} \nonumber\\
	&
	+\sum\limits_{j=13}^{15}K_{j,x}+\sum\limits_{l=1}^{5}\sum\limits_{i}\B^{16,1,l}_i +\sum\limits_{k=1}^{3}\sum\limits_{i}\B^{8,1,k} +\sum\limits_{j\neq i} (w_i-\la_i^*v_i)v_j\left[\tilde{r}_j\cdot DA\tilde{r}_i-\tilde{r}_i\cdot DA\tilde{r}_j\right]  \nonumber\\
	&-\sum\limits_{i}v_{i,t}\sum\limits_{j\neq i}(w_j-\la_j^*v_j)\xi_j(\pa_{v_i}\bar{v}_j)\tilde{r}_{j,v}
	+\sum\limits_{i}v_i \left(\frac{w_i}{v_i}\right)^2_x(B-\mu_iI_n)[\xi_i^\p\bar{v}_i\tilde{r}_{i,v}-\theta^\p_i\tilde{r}_{i,\si}],\label{J2}
\end{align}
with
\begin{align}
\sum\limits_{k=1}^{2}\sum\limits_{i}\B^{8,1,k}& =
	\sum\limits_{i}\theta_i(v_{i,x}-\mu_i^{-1}(\tilde{\la}_i-\la_i^*+\theta_i)v_i)_x\sum\limits_{j\neq i}(\mu_j-\mu_i)\psi_{ij}  r_{j}  \nonumber\\
	&+\sum\limits_{i}\theta_i^\p\left(\frac{w_i}{v_i}\right)_x(v_{i,x}-\mu_i^{-1}(\tilde{\la}_i-\la_i^*+\theta_i)v_i)\sum\limits_{j\neq i}(\mu_j-\mu_i)\psi_{ij}  r_{j} .
\end{align}
Consider the $2n \times 2n$ matrix $\mathcal{A}=\begin{pmatrix}
		\mathcal{A}_{11}&\mathcal{A}_{12}\\
		\mathcal{A}_{21}&\mathcal{A}_{22}
	\end{pmatrix}$ and $2n$ vector $\mathcal{V}$ defined as follows
\begin{align*}
	&\mathcal{A}_{11}=\left(\hat{r}_1,\cdots,\hat{r}_n\right),\quad \mathcal{A}_{12}=\left({r}^\dagger_1,\cdots,{r}^\dagger_n\right),\quad \mathcal{A}_{21}=\left({r}^\#_1+r_1^\clubsuit,\cdots,{r}^\#_n+r_n^\clubsuit\right),\\
	& \mathcal{A}_{22}=\left({r}^\ddagger_1+r^\diamondsuit_1,\cdots ,r^\ddagger_n+r^\diamondsuit_n\right),\quad\mathcal{V}=(\phi_1,\cdots,\phi_n,\psi_1,\cdots,\psi_n)^T.
\end{align*}
Then we have $\mathcal{A}\mathcal{V}=\begin{pmatrix}
	\mathcal{J}_1\\
	\mathcal{J}_2
\end{pmatrix}$. Note that $\mathcal{A}_{11}$ and $\mathcal{A}_{22}$ are invertible matrices, indeed using Lemmas \ref{lemme6.5} and the fact that $\psi_{ij}, \tilde{r}_{i,\sig}=O(1)\xi_i\bar{v}_i$ we get for any $i\in\{1,\cdot,n\}$ and $t\in[\hat{t},T]$
\begin{align*}
\hat{r}_i, r^\ddagger_{i}+r^\diamondsuit_i=r_i+O(1)|v_i|=r_i+O(1)\delta_0^2.
\end{align*}
Similarly we observe that $\mathcal{A}_{12}=O(1)\de_0$. This implies that for small enough $\de_0$, the matrix $\mathcal{A}$ is invertible. Hence, to estimate $\phi_i$ and $\psi_i$ it is enough to prove some bounds on $\mathcal{J}_1$ and $\mathcal{J}_2$. Indeed since ${\mathcal{A}}^{-1}=O(1)$ we deduce that for any $i\in\{1,\cdots,n\}$ we have
\begin{equation}
\phi_i,\psi_i=O(1)(|{\cal J}_1|+|{\cal J}_2|).
\label{supercle0}
\end{equation}
We claim that all the terms in $\phi_i$ and $\psi_i$ can be bounded as follows
\begin{align}
\phi_i,\psi_i=O(1)\sum_j (\La_j^1+\delta_0^2 \La_j^3+\La_j^4+\La_j^5+\La_j^6+\La_j^{6,1})+R_\e,
\label{boundphi}
\end{align}
with a remainder term $R_\e$ depending on $\e>0$ that we will precise later satisfying for any $T_1>\hat{t}$, $\int_{\hat{t}}^{T_1}\int_\R|R_\e(s,x)| ds dx=O(1)\delta_0^2$ for $\e>0$ sufficiently small in terms of $\delta_0$ and $T_1-\hat{t}$.

%------------------------------------------------------------------------------------------------------------------------------------------------
%------------------------------------------------------------------------------------------------------------------------------------------------
%------------------------------------------------------------------------------------------------------------------------------------------------
%------------------------------------------------------------------------------------------------------------------------------------------------

%------------------------------------------------------------------------------------------------------------------------------------------------------------------------------------

\section{Derivation of equations for effective fluxes}
In this section we assume the same assumptions as the section \ref{section6.1}.
We define 
\begin{equation}
	\mathcal{E}_i=-\sum\limits_{j\neq i}(\mu_i-\mu_j)b_{ij}\left(w_{j,x}-\frac{w_j}{v_j}v_{j,x}\right)_x,\quad \mathcal{F}_i=-\sum\limits_{j\neq i}(\mu_i-\mu_j)\hat{b}_{ij}\left(w_{j,x}-\frac{w_j}{v_j}v_{j,x}\right)_x.
	\label{defEjFj}
\end{equation}
We also denote that from \eqref{7.51} and \eqref{7.52} we have
\begin{align}
	v_{i,t}+(\tilde{\la}_iv_i)_x-(\mu_iv_{i,x})_x+\mathcal{E}_i&=\phi_i,\label{eqvi}\\
	w_{i,t}+(\tilde{\la}_iw_i)_x-(\mu_iw_{i,x})_x+\mathcal{F}_i&=\psi_i.\label{eqwi}
\end{align}
We formally define 
\begin{align}
	z_i&:=\mu_iv_{i,x}-(\tilde{\la}_i-\la_i^*)v_i+\sum\limits_{j\neq i}a_{ij}\left(w_{j,x}-\frac{w_j}{v_j}v_{j,x}\right),\label{defzi}\\
	\hat{z}_i&:=\mu_iw_{i,x}-(\tilde{\la}_i-\la_i^*)w_i+\sum\limits_{j\neq i}\hat{a}_{ij}\left(w_{j,x}-\frac{w_j}{v_j}v_{j,x}\right),\label{defzihat}
\end{align}
where $a_{ij}$ and $\hat{a}_{ij}$ will be chosen later. Now, it follows,
\begin{align}
	z_{i,x}&=(\mu_iv_{i,x})_x-((\tilde{\la}_i-\la_i^*)v_i)_x+\sum\limits_{j\neq i}a_{ij}\left(w_{j,x}-\frac{w_j}{v_j}v_{j,x}\right)_x+\sum\limits_{j\neq i}\left(w_{j,x}-\frac{w_j}{v_j}v_{j,x}\right)a_{ij,x}\nonumber\\
	&=v_{i,t}-\phi_i+\mathcal{E}_i+\la_i^*v_{i,x}+\sum\limits_{j\neq i}a_{ij}\left(w_{j,x}-\frac{w_j}{v_j}v_{j,x}\right)_x+\sum\limits_{j\neq i}\left(w_{j,x}-\frac{w_j}{v_j}v_{j,x}\right)a_{ij,x},\label{cle1}\\
	z_{i,t}&=(\mu_iv_{i,x})_t-(\tilde{\la}_iv_i)_t+\la_i^* v_{i,t}+\sum\limits_{j\neq i}a_{ij}\left(w_{j,x}-\frac{w_j}{v_j}v_{j,x}\right)_t+\sum\limits_{j\neq i}\left(w_{j,x}-\frac{w_j}{v_j}v_{j,x}\right)a_{ij,t}.\label{cle2}
\end{align}
Similarly, we get
\begin{align*}
	\hat{z}_{i,x}&=(\mu_iw_{i,x})_x-((\tilde{\la}_i-\la_i^*)w_i)_x+\sum\limits_{j\neq i}\hat{a}_{ij}\left(w_{j,x}-\frac{w_j}{v_j}v_{j,x}\right)_x+\sum\limits_{j\neq i}\left(w_{j,x}-\frac{w_j}{v_j}v_{j,x}\right)\hat{a}_{ij,x}\\
	&=w_{i,t}-\psi_i+\mathcal{F}_i+\la_i^* w_{i,x}+\sum\limits_{j\neq i}\hat{a}_{ij}\left(w_{j,x}-\frac{w_j}{v_j}v_{j,x}\right)_x+\sum\limits_{j\neq i}\left(w_{j,x}-\frac{w_j}{v_j}v_{j,x}\right)\hat{a}_{ij,x}.
\end{align*}
and
\begin{equation}\label{eqn-hat-z-i-t}
	\hat{z}_{i,t}=(\mu_iw_{i,x})_t-(\tilde{\la}_iw_i)_t+\la_i^* w_{i,t}+\sum\limits_{j\neq i}\hat{a}_{ij}\left(w_{j,x}-\frac{w_j}{v_j}v_{j,x}\right)_t+\sum\limits_{j\neq i}\left(w_{j,x}-\frac{w_j}{v_j}v_{j,x}\right)\hat{a}_{ij,t}.
\end{equation}
\subsection{Derivation of equations for effective fluxes for \eqref{eqn-v-i-1}}\label{section1zi}
Now we first establish the parabolic equation for $z_i$. To do this we observe that using  \eqref{eqn-w-i} and \eqref{cle1} we get
\begin{align*}
	&(\mu_iv_{i,x})_t-(\tilde{\la}_iv_i)_t+\la_i^* v_{i,t}\\
	&=\mu_iv_{i,tx}-(\tilde{\la}_{i}-\la_i^*)v_{i,t}+\sum\limits_{j}(v_{i,x}\mu_{i,u}-v_i\tilde{\la}_{i,u})\tilde{r}_j(w_j-\la_j^*v_j)\\
	&-\tilde{\la}_{i,v}v_{i,t}\xi_i\pa_{v_i}\bar{v}_iv_{i}-\sum\limits_{j\neq i}\tilde{\la}_{i,v}v_{j,t}v_{i}\xi_i\pa_{v_j}\bar{v}_i-(\tilde{\la}_{i,v}\bar{v}_i\xi_i^\p-\tilde{\la}_{i,\si}\theta_i^\p)v_{i}\left(\frac{w_i}{v_i}\right)_t\\
	&=(\mu_iz_{i,x}-\tilde{\la}_{i}z_{i})_x+\mu_i\phi_{i,x}-\mu_i\mathcal{E}_{i,x}-\sum\limits_{j\neq i}\mu_ia_{ij}\left(w_{j,x}-\frac{w_j}{v_j}v_{j,x}\right)_{xx}\\
	&-2\sum\limits_{j\neq i}\mu_i\left(w_{j,x}-\frac{w_j}{v_j}v_{j,x}\right)_x a_{ij,x}-\sum\limits_{j\neq i}\mu_i\left(w_{j,x}-\frac{w_j}{v_j}v_{j,x}\right)a_{ij,xx}-\mu_{i,x}(z_{i,x}-\la_i^* v_{i,x})+\tilde{\lambda}_{i,x}(z_i-\la_i^* v_{i})\\
	&-(\tilde{\la}_i-\la_i^*)\phi_i+(\tilde{\la}_i-\la_i^*)\mathcal{E}_i+\tilde{\la}_i\sum\limits_{j\neq i}a_{ij}\left(w_{j,x}-\frac{w_j}{v_j}v_{j,x}\right)_{x}+\tilde{\la}_i\sum\limits_{j\neq i}a_{ij,x}\left(w_{j,x}-\frac{w_j}{v_j}v_{j,x}\right)\\
	&+\sum\limits_{j}(v_{i,x}\mu_{i,u}-v_i\tilde{\la}_{i,u})\tilde{r}_j(w_j-\la_j^*v_j)
	-\tilde{\la}_{i,v}v_{i,t}\xi_i\pa_{v_i}\bar{v}_iv_{i}-\sum\limits_{j\neq i}\tilde{\la}_{i,v}v_{j,t}v_{i}\xi_i\pa_{v_j}\bar{v}_i\\
	%&-\sum\limits_{j\neq i}\tilde{\la}_{i,v}w_{j,t}v_{i}\pa_{w_j}\bar{v}_i\xi_i\\
	&-(\tilde{\la}_{i,v}\bar{v}_i\xi_i^\p-\tilde{\la}_{i,\si}\theta_i^\p)v_{i}\left(\frac{w_i}{v_i}\right)_t.
\end{align*}
Hence, we get
\begin{align*}
	&(\mu_iv_{i,x})_t-(\tilde{\la}_iv_i)_t+\la_i^* v_{i,t}=(\mu_iz_{i,x}-\tilde{\la}_{i}z_{i})_x+\mu_i\phi_{i,x}-\mu_i\mathcal{E}_{i,x}-\sum\limits_{j\neq i}\mu_ia_{ij}\left(w_{j,x}-\frac{w_j}{v_j}v_{j,x}\right)_{xx}\\
	&-2\sum\limits_{j\neq i}\mu_i\left(w_{j,x}-\frac{w_j}{v_j}v_{j,x}\right)_x a_{ij,x}-\sum\limits_{j\neq i}\mu_i\left(w_{j,x}-\frac{w_j}{v_j}v_{j,x}\right)a_{ij,xx}\\
	&-(\tilde{\la}_i-\la_i^*)\phi_i+(\tilde{\la}_i-\la_i^*)\mathcal{E}_i+\tilde{\la}_i\sum\limits_{j\neq i}a_{ij}\left(w_{j,x}-\frac{w_j}{v_j}v_{j,x}\right)_{x}+\tilde{\la}_i\sum\limits_{j\neq i}a_{ij,x}\left(w_{j,x}-\frac{w_j}{v_j}v_{j,x}\right)\\
	&+\sum\limits_{j}(v_{i,x}(w_j-\la_j^*v_j)-v_j (z_{i,x}-\la_i^*v_{i,x}))\mu_{i,u}\tilde{r}_j-\sum_j(v_i(w_j-\lambda_j^*v_j)-(z_i-\la_i^*v_i) v_j) \tilde{\la}_{i,u}\tilde{r}_j\\
	&+
	\tilde{\la}_{i,v}\xi_i\pa_{v_i}\bar{v}_i (v_{i,x}(z_i-\la_i^* v_{i})
	-v_{i}v_{i,t})+\sum\limits_{j\neq i}\tilde{\la}_{i,v}((z_i-\la_i^* v_{i}) v_{j,x}-v_{j,t}v_{i})\xi_i\pa_{v_j}\bar{v}_i\\
	%&+\sum\limits_{j\neq i}\tilde{\la}_{i,v}((z_i-\la_i^* v_{i}) w_{j,x}-w_{j,t}v_{i})\pa_{w_j}\bar{v}_i\xi_i\\
	&+(\tilde{\la}_{i,v}\xi_i'\bar{v}_i -\tilde{\la}_{i,\sig}\theta'_i)\left[\left(\frac{w_i}{v_i}\right)_x(z_i-\la_i^* v_{i})-\left(\frac{w_i}{v_i}\right)_t v_i\right].
\end{align*}
Note that from \eqref{cle2}, we deduce that
\begin{align*}
	z_{i,t}&=(\mu_iz_{i,x}-\tilde{\la}_{i}z_{i})_x+\mu_i\phi_{i,x}\\
	&-\mu_i\mathcal{E}_{i,x}+\sum\limits_{j\neq i}a_{ij}\left[\left(w_{j,x}-\frac{w_j}{v_j}v_{j,x}\right)_t-\mu_i \left(w_{j,x}-\frac{w_j}{v_j}v_{j,x}\right)_{xx}\right]\\
	&+\sum\limits_{j\neq i}(a_{ij,t}-\mu_i a_{ij,xx})\left(w_{j,x}-\frac{w_j}{v_j}v_{j,x}\right)-2\sum\limits_{j\neq i}\mu_i\left(w_{j,x}-\frac{w_j}{v_j}v_{j,x}\right)_x a_{ij,x}\\
	&-(\tilde{\la}_i-\la_i^*)\phi_i+(\tilde{\la}_i-\la_i^*)\mathcal{E}_i\\
	&+\tilde{\la}_i\sum\limits_{j\neq i}a_{ij}\left(w_{j,x}-\frac{w_j}{v_j}v_{j,x}\right)_{x}+\tilde{\la}_i\sum\limits_{j\neq i}a_{ij,x}\left(w_{j,x}-\frac{w_j}{v_j}v_{j,x}\right)\\
	&+\sum\limits_{j}(v_{i,x}(w_j-\la_j^*v_j)-v_j (z_{i,x}-\la_i^*v_{i,x}))\mu_{i,u}\tilde{r}_j\\
	&-\sum_j(v_i(w_j-\lambda_j^*v_j)-(z_i -\la_i^*v_{i})v_j) \tilde{\la}_{i,u}\tilde{r}_j\\
	&+
	\tilde{\la}_{i,v}\xi_i\pa_{v_i}\bar{v}_i (v_{i,x}(z_i-\la_i^*v_{i})-v_{i}v_{i,t})\\
	&+\sum\limits_{j\neq i}\tilde{\la}_{i,v}((z_i -\la_i^*v_{i})v_{j,x}-v_{j,t}v_{i})\xi_i\pa_{v_j}\bar{v}_i\\
	%&+\sum\limits_{j\neq i}\tilde{\la}_{i,v}((z_i-\la_i^*v_{i}) w_{j,x}-w_{j,t}v_{i})\pa_{w_j}\bar{v}_i\xi_i\\
	&+(\tilde{\la}_{i,v}\xi_i'\bar{v}_i-\tilde{\lambda}_{i,\sig}\theta^\p_i)\left[\left(\frac{w_i}{v_i}\right)_x(z_i-\la_i^*v_{i})-\left(\frac{w_i}{v_i}\right)_t v_i\right]\\
	&=(\mu_iz_{i,x}-\tilde{\la}_{i}z_{i})_x+\mu_i\phi_{i,x}+\sum\limits_{l=1}^{9}\mathcal{T}_l,
\end{align*}
where $\mathcal{T}_l$ are defined as follows
\begin{align}
\mathcal{T}_1&=-\mu_i\mathcal{E}_{i,x}+\sum\limits_{j\neq i}a_{ij}\left[\left(w_{j,x}-\frac{w_j}{v_j}v_{j,x}\right)_t-\mu_i \left(w_{j,x}-\frac{w_j}{v_j}v_{j,x}\right)_{xx}\right],\nonumber\\
\mathcal{T}_2&=\sum\limits_{j\neq i}(a_{ij,t}-\mu_i a_{ij,xx})\left(w_{j,x}-\frac{w_j}{v_j}v_{j,x}\right)-2\sum\limits_{j\neq i}\mu_i\left(w_{j,x}-\frac{w_j}{v_j}v_{j,x}\right)_x a_{ij,x},\nonumber\\
\mathcal{T}_3&=-(\tilde{\la}_i-\la_i^*)\phi_i+(\tilde{\la}_i-\la_i^*)\mathcal{E}_i,\nonumber\\
\mathcal{T}_4&=\tilde{\la}_i\sum\limits_{j\neq i}a_{ij}\left(w_{j,x}-\frac{w_j}{v_j}v_{j,x}\right)_{x}+\tilde{\la}_i\sum\limits_{j\neq i}a_{ij,x}\left(w_{j,x}-\frac{w_j}{v_j}v_{j,x}\right),\nonumber\\
\mathcal{T}_5&=\sum\limits_{j}(v_{i,x}(w_j-\la_j^*v_j)-v_j (z_{i,x}-\la_i^*v_{i,x}))\mu_{i,u}\tilde{r}_j,\nonumber\\
\mathcal{T}_6&=-\sum_j(v_i(w_j-\lambda_j^*v_j)-(z_i -\la_i^*v_{i})v_j) \tilde{\la}_{i,u}\tilde{r}_j,\nonumber\\
\mathcal{T}_7&=\tilde{\la}_{i,v}\xi_i\pa_{v_i}\bar{v}_i (v_{i,x}(z_i-\la_i^*v_{i})-v_{i}v_{i,t}),\nonumber\\
\mathcal{T}_8&=\sum\limits_{j\neq i}\tilde{\la}_{i,v}((z_i -\la_i^*v_{i})v_{j,x}-v_{j,t}v_{i})\xi_i\pa_{v_j}\bar{v}_i,\nonumber\\
%\mathcal{T}_9&=\sum\limits_{j\neq i}\tilde{\la}_{i,v}((z_i-\la_i^*v_{i}) w_{j,x}-w_{j,t}v_{i})\pa_{w_j}\bar{v}_i\xi_i,\\
\mathcal{T}_{9}&=(\tilde{\la}_{i,v}\xi_i'\bar{v}_i-\tilde{\lambda}_{i,\sig}\theta^\p_i)\left[\left(\frac{w_i}{v_i}\right)_x(z_i-\la_i^*v_{i})-\left(\frac{w_i}{v_i}\right)_t v_i\right].\label{defglobT}
\end{align}
We calculate using \eqref{eqvi}, \eqref{eqwi} and
the fact that $\left(\frac{w_j}{v_j}\right)_{xx}=\frac{w_{j,xx}}{v_j}-\frac{v_{j,xx}w_j}{v_j^2}-2\frac{v_{j,x}}{v_j}\left(\frac{w_j}{v_j}\right)_x$
\begin{align*}
	&\left(w_{j,x}-\frac{w_j}{v_j}v_{j,x}\right)_t-\mu_i\left(w_{j,x}-\frac{w_j}{v_j}v_{j,x}\right)_{xx}\\
	&=w_{j,tx}-\mu_i w_{j,xxx}-\frac{w_j}{v_j}\left(v_{j,tx}-\mu_iv_{j,xxx}\right)-v_{j,x}\left[\left(\frac{w_j}{v_j}\right)_t-\mu_i\left(\frac{w_j}{v_j}\right)_{xx}\right]+2\mu_i\left(\frac{w_j}{v_j}\right)_x v_{j,xx}\\
	&=w_{j,tx}-\mu_i w_{j,xxx}-\frac{w_j}{v_j}\left(v_{j,tx}-\mu_iv_{j,xxx}\right)-\frac{v_{j,x}}{v_j}\left(w_{j,t}-\mu_iw_{j,xx}\right)\\
	&+\frac{w_jv_{j,x}}{v^2_j}\left(v_{j,t}-\mu_iv_{j,xx}\right)-\frac{2\mu_iv_{j,x}^2}{v_j}\left(\frac{w_j}{v_j}\right)_x+2\mu_i\left(\frac{w_j}{v_j}\right)_x v_{j,xx}\\
	&=\left(\mu_{j,x}w_{j,x}-\tilde{\la}_{j,x}w_j-\tilde{\la}_jw_{j,x}\right)_x+\mu_{j,x}w_{j,xx}+(\mu_j-\mu_i)w_{j,xxx}+\psi_{j,x}-\mathcal{F}_{j,x}\\
	&-\frac{w_j}{v_j}\left(\mu_{j,x}v_{j,x}-\tilde{\la}_{j,x}v_j-\tilde{\la}_jv_{j,x}\right)_x-\frac{w_j}{v_j}\mu_{j,x}v_{j,xx}-\frac{w_j}{v_j}(\mu_j-\mu_i)v_{j,xxx}-\frac{w_j}{v_j}\phi_{j,x}+\frac{w_j}{v_j}\mathcal{E}_{j,x}\\
	&-\frac{v_{j,x}}{v_j}\left(\mu_{j,x}w_{j,x}-\tilde{\la}_{j,x}w_j-\tilde{\la}_jw_{j,x}+(\mu_j-\mu_i)w_{j,xx}+\psi_j-\mathcal{F}_j\right)\\
	&+\frac{w_jv_{j,x}}{v^2_j}\left(\mu_{j,x}v_{j,x}-\tilde{\la}_{j,x}v_j-\tilde{\la}_jv_{j,x}+(\mu_j-\mu_i)v_{j,xx}+\phi_j-\mathcal{E}_j\right)-\frac{2\mu_iv_{j,x}^2}{v_j}\left(\frac{w_j}{v_j}\right)_x\\
	&+2\mu_i\left(\frac{w_j}{v_j}\right)_x v_{j,xx}.
\end{align*}
Then,
\begin{align}
	&\left(w_{j,x}-\frac{w_j}{v_j}v_{j,x}\right)_t-\mu_i\left(w_{j,x}-\frac{w_j}{v_j}v_{j,x}\right)_{xx}\nonumber\\
	&=\mu_{j,xx}v_j\left(\frac{w_j}{v_j}\right)_x+(2\mu_{j,x}-\tilde{\la}_{j})\left(w_{j,xx}-\frac{w_j}{v_j}v_{j,xx}\right)-2\tilde{\la}_{j,x}v_j\left(\frac{w_j}{v_j}\right)_x\nonumber\\
	&+(\mu_j-\mu_i)\left(w_{j,xxx}-\frac{w_j}{v_j}v_{j,xxx}\right)+\psi_{j,x}-\frac{w_j}{v_j}\phi_{j,x}+\frac{w_j}{v_j}\mathcal{E}_{j,x}-\mathcal{F}_{j,x}-(\mu_{j,x}-\tilde{\la}_j)v_{j,x}\left(\frac{w_j}{v_j}\right)_x\nonumber\\
	&-\frac{v_{j,x}}{v_j}(\mu_j-\mu_i)\left(w_{j,xx}-\frac{w_j}{v_j}v_{j,xx}\right)-\frac{v_{j,x}}{v_j}\left(\psi_j-\mathcal{F}_j\right)+\frac{w_jv_{j,x}}{v^2_j}\left(\phi_j-\mathcal{E}_j\right)-\frac{2\mu_iv_{j,x}^2}{v_j}\left(\frac{w_j}{v_j}\right)_x\nonumber\\
	&+2\mu_i\left(\frac{w_j}{v_j}\right)_x v_{j,xx}.\label{z-cal-1-1}
\end{align}
From the definition of $\mathcal{E}_i$ we have since $\left(\frac{w_j}{v_j}\right)_{xx}=\frac{1}{v_j}(w_{j,xx}-\frac{w_j}{v_j}v_{j,xx})-2\frac{v_{j,x}}{v_j}(\frac{w_j}{v_j})_x$
\begin{align}
	\mathcal{E}_{i,x}&=-\sum\limits_{j\neq i}(\mu_i-\mu_j)b_{ij}\left(w_{j,x}-\frac{w_j}{v_j}v_{j,x}\right)_{xx}-\sum\limits_{j\neq i}((\mu_i-\mu_j)b_{ij})_x\left(w_{j,x}-\frac{w_j}{v_j}v_{j,x}\right)_{x}\nonumber\\
	&=-\sum\limits_{j\neq i}(\mu_i-\mu_j)b_{ij}\left(w_{j,xxx}-\frac{w_j}{v_j}v_{j,xxx}\right)+\sum\limits_{j\neq i}2(\mu_i-\mu_j)b_{ij}\left(\frac{w_j}{v_j}\right)_xv_{j,xx}\nonumber\\
	&+\sum\limits_{j\neq i}(\mu_i-\mu_j)b_{ij}\frac{v_{j,x}}{v_j}\left(w_{j,xx}-\frac{w_j}{v_j}v_{j,xx}\right)-\sum\limits_{j\neq i}2(\mu_i-\mu_j)b_{ij}\frac{v^2_{j,x}}{v_j}\left(\frac{w_j}{v_j}\right)_x\nonumber\\
	&-\sum\limits_{j\neq i}(\mu_i-\mu_j)b_{ij,x}\left(w_{j,x}-\frac{w_j}{v_j}v_{j,x}\right)_{x}-\sum\limits_{j\neq i}(\mu_i-\mu_j)_xb_{ij}\left(w_{j,x}-\frac{w_j}{v_j}v_{j,x}\right)_{x}.\nonumber
\end{align}
Therefore,
\begin{align}
	&-\mu_i\mathcal{E}_{i,x}+\sum\limits_{j\neq i}a_{ij}(\mu_j-\mu_i)\left(w_{j,xxx}-\frac{w_j}{v_j}v_{j,xxx}\right)\nonumber\\
	&=\sum\limits_{j\neq i}(\mu_i-\mu_j)(-a_{ij}+\mu_ib_{ij})\left(w_{j,xxx}-\frac{w_j}{v_j}v_{j,xxx}\right)-\mu_i\sum\limits_{j\neq i}2(\mu_i-\mu_j)b_{ij}\left(\frac{w_j}{v_j}\right)_xv_{j,xx}\nonumber\\
	&-\mu_i\sum\limits_{j\neq i}(\mu_i-\mu_j)b_{ij}\frac{v_{j,x}}{v_j}\left(w_{j,xx}-\frac{w_j}{v_j}v_{j,xx}\right)+\mu_i\sum\limits_{j\neq i}2(\mu_i-\mu_j)b_{ij}\frac{v^2_{j,x}}{v_j}\left(\frac{w_j}{v_j}\right)_x\nonumber\\
	&+\mu_i\sum\limits_{j\neq i}(\mu_i-\mu_j)b_{ij,x}\left(w_{j,x}-\frac{w_j}{v_j}v_{j,x}\right)_{x}+\mu_i\sum\limits_{j\neq i}(\mu_i-\mu_j)_xb_{ij}\left(w_{j,x}-\frac{w_j}{v_j}v_{j,x}\right)_{x}.\label{eqn-E-ab}
\end{align}
Subsequently, it follows from \eqref{z-cal-1-1}
\begin{align*}
	\mathcal{T}_1&=-\mu_i\mathcal{E}_{i,x}+\sum\limits_{j\neq i}a_{ij}\left[\left(w_{j,x}-\frac{w_j}{v_j}v_{j,x}\right)_t-\mu_i \left(w_{j,x}-\frac{w_j}{v_j}v_{j,x}\right)_{xx}\right],\\
	&=-\mu_i\mathcal{E}_{i,x}+\sum\limits_{j\neq i}a_{ij}(\mu_j-\mu_i)\left(w_{j,xxx}-\frac{w_j}{v_j}v_{j,xxx}\right)\\
	&+\sum\limits_{j\neq i}a_{ij}\mu_{j,xx}v_j\left(\frac{w_j}{v_j}\right)_x+\sum\limits_{j\neq i}a_{ij}(2\mu_{j,x}-\tilde{\la}_{j})\left(w_{j,xx}-\frac{w_j}{v_j}v_{j,xx}\right)-\sum\limits_{j\neq i}2 a_{ij} \tilde{\la}_{j,x}v_j\left(\frac{w_j}{v_j}\right)_x\\
	&+\sum\limits_{j\neq i} a_{ij}\big(\psi_{j,x}-\frac{w_j}{v_j}\phi_{j,x}+\frac{w_j}{v_j}\mathcal{E}_{j,x}-\mathcal{F}_{j,x}\big)-\sum\limits_{j\neq i} a_{ij} (\mu_{j,x}-\tilde{\la}_j)v_{j,x}\left(\frac{w_j}{v_j}\right)_x\\
	&-\sum\limits_{j\neq i}a_{ij} \frac{v_{j,x}}{v_j}(\mu_j-\mu_i)\left(w_{j,xx}-\frac{w_j}{v_j}v_{j,xx}\right)-\sum\limits_{j\neq i}a_{ij} \frac{v_{j,x}}{v_j}\left(\psi_j-\mathcal{F}_j\right)+\sum\limits_{j\neq i}a_{ij} \frac{w_jv_{j,x}}{v^2_j}\left(\phi_j-\mathcal{E}_j\right)\\
	&-\sum\limits_{j\neq i}\frac{2 a_{ij} \mu_iv_{j,x}^2}{v_j}\left(\frac{w_j}{v_j}\right)_x+\sum\limits_{j\neq i}2a_{ij} \mu_i\left(\frac{w_j}{v_j}\right)_x v_{j,xx}.
\end{align*}
By using \eqref{eqn-E-ab} we have
\begin{align*}
	\mathcal{T}_1&=\sum\limits_{j\neq i}(\mu_i-\mu_j)(-a_{ij}+\mu_ib_{ij})\left(w_{j,xxx}-\frac{w_j}{v_j}v_{j,xxx}\right)-\mu_i\sum\limits_{j\neq i}2(\mu_i-\mu_j)b_{ij}\left(\frac{w_j}{v_j}\right)_xv_{j,xx}\nonumber\\
	&-\mu_i\sum\limits_{j\neq i}(\mu_i-\mu_j)b_{ij}\frac{v_{j,x}}{v_j}\left(w_{j,xx}-\frac{w_j}{v_j}v_{j,xx}\right)+\mu_i\sum\limits_{j\neq i}2(\mu_i-\mu_j)b_{ij}\frac{v^2_{j,x}}{v_j}\left(\frac{w_j}{v_j}\right)_x\nonumber\\
	&+\mu_i\sum\limits_{j\neq i}(\mu_i-\mu_j)b_{ij,x}\left(w_{j,x}-\frac{w_j}{v_j}v_{j,x}\right)_{x}+\mu_i\sum\limits_{j\neq i}(\mu_i-\mu_j)_xb_{ij}\left(w_{j,x}-\frac{w_j}{v_j}v_{j,x}\right)_{x}\\
	&+\sum\limits_{j\neq i}a_{ij}\mu_{j,xx}v_j\left(\frac{w_j}{v_j}\right)_x+\sum\limits_{j\neq i}a_{ij}(2\mu_{j,x}-\tilde{\la}_{j})\left(w_{j,xx}-\frac{w_j}{v_j}v_{j,xx}\right)-\sum\limits_{j\neq i}2a_{ij}\tilde{\la}_{j,x}v_j\left(\frac{w_j}{v_j}\right)_x\\
	&+\sum\limits_{j\neq i}a_{ij}\big(\psi_{j,x}-\frac{w_j}{v_j}\phi_{j,x}+\frac{w_j}{v_j}\mathcal{E}_{j,x}-\mathcal{F}_{j,x}\big)-\sum\limits_{j\neq i}a_{ij}(\mu_{j,x}-\tilde{\la}_j)v_{j,x}\left(\frac{w_j}{v_j}\right)_x\\
	&-\sum\limits_{j\neq i} a_{ij}\frac{v_{j,x}}{v_j}(\mu_j-\mu_i)\left(w_{j,xx}-\frac{w_j}{v_j}v_{j,xx}\right)-\sum\limits_{j\neq i}a_{ij} \frac{v_{j,x}}{v_j}\left(\psi_j-\mathcal{F}_j\right)+\sum\limits_{j\neq i} a_{ij} \frac{w_jv_{j,x}}{v^2_j}\left(\phi_j-\mathcal{E}_j\right)\\
	&-\sum\limits_{j\neq i}\frac{2 a_{ij} \mu_iv_{j,x}^2}{v_j}\left(\frac{w_j}{v_j}\right)_x+\sum\limits_{j\neq i}2a_{ij}\mu_i\left(\frac{w_j}{v_j}\right)_x v_{j,xx}.
\end{align*}
Now, we may choose $a_{ij}=\mu_i b_{ij}$ to obtain
\begin{align}
	\mathcal{T}_1&=-\mu_i\sum\limits_{j\neq i}(\mu_i-\mu_j)b_{ij}\frac{v_{j,x}}{v_j}\left(w_{j,xx}-\frac{w_j}{v_j}v_{j,xx}\right)+\mu_i\sum\limits_{j\neq i}2(\mu_i-\mu_j)b_{ij}\frac{v^2_{j,x}}{v_j}\left(\frac{w_j}{v_j}\right)_x\nonumber\\
	&+\mu_i\sum\limits_{j\neq i}(\mu_i-\mu_j)b_{ij,x}\left(w_{j,x}-\frac{w_j}{v_j}v_{j,x}\right)_{x}+\mu_i\sum\limits_{j\neq i}(\mu_i-\mu_j)_xb_{ij}\left(w_{j,x}-\frac{w_j}{v_j}v_{j,x}\right)_{x}\nonumber\\
	&+\sum\limits_{j\neq i}a_{ij}\mu_{j,xx}v_j\left(\frac{w_j}{v_j}\right)_x+\sum\limits_{j\neq i}a_{ij}(2\mu_{j,x}-\tilde{\la}_{j})\left(w_{j,xx}-\frac{w_j}{v_j}v_{j,xx}\right)-\sum\limits_{j\neq i}2a_{ij}\tilde{\la}_{j,x}v_j\left(\frac{w_j}{v_j}\right)_x\nonumber\\
	&+\sum\limits_{j\neq i}a_{ij}\big(\psi_{j,x}-\frac{w_j}{v_j}\phi_{j,x}+\frac{w_j}{v_j}\mathcal{E}_{j,x}-\mathcal{F}_{j,x}\big)-\sum\limits_{j\neq i}a_{ij}(\mu_{j,x}-\tilde{\la}_j)v_{j,x}\left(\frac{w_j}{v_j}\right)_x\nonumber\\
	&-\sum\limits_{j\neq i}a_{ij}\frac{v_{j,x}}{v_j}(\mu_j-\mu_i)\left(w_{j,xx}-\frac{w_j}{v_j}v_{j,xx}\right)-\sum\limits_{j\neq i}a_{ij} \frac{v_{j,x}}{v_j}\left(\psi_j-\mathcal{F}_j\right)+\sum\limits_{j\neq i}a_{ij} \frac{w_jv_{j,x}}{v^2_j}\left(\phi_j-\mathcal{E}_j\right)\nonumber\\
	&-\sum\limits_{j\neq i}\frac{2 a_{ij} \mu_iv_{j,x}^2}{v_j}\left(\frac{w_j}{v_j}\right)_x+\sum\limits_{j\neq i}2a_{ij}\mu_i\left(\frac{w_j}{v_j}\right)_x v_{j,xx}.
	\label{T1}
\end{align}
Furthermore using \eqref{eqn-v-i-x-1}, we recall that
\begin{equation} \label{eqn-v-i-x-1a-1}
	w_i-\la_i^*v_i=\mu_i v_{i,x}-\tilde{\la}_iv_i +{\cal A}_i=z_i-\la_i^*v_i-\sum_{i\ne j}a_{ij}\left(w_{j x}-\frac{w_j}{v_j}v_{jx}\right)+{\cal A}_i,     
\end{equation}
with
\begin{align}
	&{\cal A}_i=\sum\limits_{j\neq i}(\mu_j v_{j,x}-(\tilde{\la}_j-\la_j^*)v_j-w_j)\frac{\mu_i}{\mu_j}[\psi_{ji}+v_j\xi_j\pa_{v_j}\bar{v}_j \psi_{ji,v}]\nonumber\\
	&+\sum\limits_{j\neq i}\sum\limits_{k\neq j}\xi_jv_{k,x}v_j\pa_{v_k}\bar{v}_j\mu_i \psi_{ji,v}+\sum\limits_{j\neq i}v_j\left(\frac{w_j}{v_j}\right)_x\mu_i [-\psi_{ji,\si}+\bar{v}_j\xi_j^\p\psi_{ji,v}]  \nonumber\\
	&+\sum\limits_{j\neq i}\mu_j^{-1}(\tilde{\la}_j-\si_j)v_j\xi_j \tilde{v}_j\mu_i\psi_{ji,v}+\sum\limits_{k}v_k^2(1-\xi_k\chi_k^\ep \eta_k)\langle l_i, \left[B(u)\tilde{r}_{k,u}\tilde{r}_k+\tilde{r}_k\cdot DB(u)\tilde{r}_k\right]\rangle \nonumber\\
	&+\sum\limits_{k\neq j}v_kv_j\langle l_i,\tilde{r}_k\cdot DB(u)\tilde{r}_j+B(u)\tilde{r}_{k,u}\tilde{r}_j\rangle. \label{Ai}
\end{align} 
We observe using \eqref{eqn-v-i-x-1a-1} that
\begin{align*}
	v_{i,x}(w_i-\la_i^*v_i)-v_i(z_{i,x}-\la_i^*v_{i,x})&=\mu_iv_{i,x}v_{i,x}-\tilde{\la}_iv_iv_{i,x}+{\cal A}_i v_{i,x}+\la_i^*v_i v_{i,x}-v_iz_{i,x}\\
	&=v_{i,x}z_i-v_{i}z_{i,x}-v_{i,x}\sum\limits_{j\neq i}a_{ij}\left(w_{j,x}-\frac{w_j}{v_j}v_{j,x}\right)+{\cal A}_iv_{i,x}.
\end{align*}
As a result we have the following simplification of $\mathcal{T}_5$,
\begin{align}
	\mathcal{T}_5&=(v_{i,x}(w_i-\la_i^*v_i)-v_i (z_{i,x}-\la_i^*v_{i,x}))\mu_{i,u}\tilde{r}_i+\sum\limits_{j\neq i}(v_{i,x}(w_j-\la_j^*v_j)-v_j (z_{i,x}-\la_i^*v_{i,x}))\mu_{i,u}\tilde{r}_j\nonumber\\
	&=\biggl[(v_{i,x}z_i-v_{i}z_{i,x})-\sum\limits_{j\neq i}v_{i,x}a_{ij}\left(w_{j,x}-\frac{w_j}{v_j}v_{j,x}\right)+{\cal A}_iv_{i,x}\big]\mu_{i,u}\tilde{r}_i\nonumber\\
	&\hspace{5cm}+\sum\limits_{j\neq i}(v_{i,x}(w_j-\la_j^*v_j)-v_j (z_{i,x}-\la_i^*v_{i,x}))\mu_{i,u}\tilde{r}_j.	\label{T5n}
\end{align}
Similarly, it follows from \eqref{eqn-v-i-x-1a-1}
\begin{align*}
	&(v_i(w_i-\la_i^*v_i)-(z_i- \la_i^*v_i)v_i)=
	{\cal A}_i v_i-v_{i}\sum\limits_{j\neq i}a_{ij}\left(w_{j,x}-\frac{w_j}{v_j}v_{j,x}\right).
\end{align*}
Then, 
\begin{align}
	&\mathcal{T}_6=-(v_i(w_i-\lambda_i^*v_i)-(z_i -\la_i^*v_{i})v_i) \tilde{\la}_{i,u}\tilde{r}_i-\sum_{j\neq i}(v_i(w_j-\lambda_j^*v_j)-(z_i -\la_i^*v_{i})v_j) \tilde{\la}_{i,u}\tilde{r}_j\nonumber\\
	&=\big(-{\cal A}_i v_i+v_{i}\sum\limits_{j\neq i}a_{ij}\left(w_{j,x}-\frac{w_j}{v_j}v_{j,x}\right)\big)\tilde{\la}_{i,u}\tilde{r}_i-\sum_{j\neq i}(v_i(w_j-\lambda_j^*v_j)-(z_i -\la_i^*v_{i})v_j) \tilde{\la}_{i,u}\tilde{r}_j.
	\label{T6n}
\end{align}
We also observe that using \eqref{cle1}
\begin{align*}
	&(v_{i,x}(z_i-\la_i^*v_{i})-v_{i,t}v_{i})\nonumber\\
	&=v_{i,x}z_i-z_{i,x}v_i-\phi_iv_{i}+\mathcal{E}_iv_{i}+\sum\limits_{j\neq i}a_{ij}v_{i}\left(w_{j,x}-\frac{w_j}{v_j}v_{j,x}\right)_x+\sum\limits_{j\neq i}\left(w_{j,x}-\frac{w_j}{v_j}v_{j,x}\right)v_{i}a_{ij,x}\nonumber\\
	&=v_{i,x}z_i-z_{i,x}v_i-\phi_iv_{i}+\mathcal{E}_iv_{i}+\sum\limits_{j\neq i}a_{ij}v_{i}\left(w_{j,xx}-\frac{w_j}{v_j}v_{j,xx}\right)\nonumber\\
	&-\sum\limits_{j\neq i}a_{ij}v_{i}\left( \frac{w_j}{v_j}\right)_xv_{j,x}+\sum\limits_{j\neq i}\left(w_{j,x}-\frac{w_j}{v_j}v_{j,x}\right)v_{i}a_{ij,x}.
\end{align*}
It yields then
\begin{align}
	\mathcal{T}_7&=\tilde{\la}_{i,v}\xi_i\pa_{v_i}\bar{v}_i (v_{i,x}(z_i-\la_i^*v_{i})-v_{i}v_{i,t})\nonumber\\
	&=(\tilde{\la}_{i,v}\xi_i\pa_{v_i}\bar{v}_i )(v_{i,x}z_i-z_{i,x}v_i)-\phi_iv_{i}(\tilde{\la}_{i,v}\xi_i\pa_{v_i}\bar{v}_i )+\mathcal{E}_iv_{i}(\tilde{\la}_{i,v}\xi_i\pa_{v_i}\bar{v}_i )\nonumber\\
	&+(\tilde{\la}_{i,v}\xi_i\pa_{v_i}\bar{v}_i )\sum\limits_{j\neq i}a_{ij}v_{i}\left(w_{j,xx}-\frac{w_j}{v_j}v_{j,xx}\right)-(\tilde{\la}_{i,v}\xi_i\pa_{v_i}\bar{v}_i )\sum\limits_{j\neq i}a_{ij}v_{i}\left( \frac{w_j}{v_j}\right)_xv_{j,x}\nonumber\\
	&+(\tilde{\la}_{i,v}\xi_i\pa_{v_i}\bar{v}_i )\sum\limits_{j\neq i}\left(w_{j,x}-\frac{w_j}{v_j}v_{j,x}\right)v_{i}a_{ij,x}.\label{T7n}
\end{align}
%Observe that
%\begin{align*}
%	&(w_{j,x}(z_i-\la_i^*v_{i})-w_{j,t}v_{i})
%	=w_{j,x}z_i-\hat{z}_{j,x}v_{i}-\psi_jv_{i}+\mathcal{F}_jv_{i}+\sum\limits_{k\neq j}\hat{a}_{jk}v_{i}\left(w_{k,xx}-\frac{w_k}{v_k}v_{k,xx}\right)\\
%	&-(\la_i^*-\la_j^*)v_iw_{j,x}-\sum\limits_{k\neq j}\hat{a}_{jk}\left(\frac{w_k}{v_k}\right)_xv_{k,x}v_{i}+\sum\limits_{k\neq j}\left(w_{k,x}-\frac{w_k}{v_k}v_{k,x}\right)v_i\hat{a}_{jk,x}.
%\end{align*}
Furthermore again using \eqref{cle1}, we get
\begin{align}
	\mathcal{T}_8&
	=\sum\limits_{j\neq i}\tilde{\la}_{i,v}(v_{j,x}z_i-z_{j,x}v_{i})\xi_i\pa_{v_j}\bar{v}_i-\sum\limits_{j\neq i}\tilde{\la}_{i,v}\xi_iv_{i}(\pa_{v_j}\bar{v}_i)[\phi_j -\mathcal{E}_j]\nonumber\\
	&+\sum\limits_{j\neq i}\tilde{\la}_{i,v}v_{i}\xi_i(\pa_{v_j}\bar{v}_i)\sum\limits_{k\neq j}a_{jk}\left(w_{k,xx}-\frac{w_k}{v_k}v_{k,xx}-\left(\frac{w_k}{v_k}\right)_xv_{k,x}\right)\nonumber\\
	&+\sum\limits_{j\neq i}\tilde{\la}_{i,v}v_{i}\xi_i(\pa_{v_j}\bar{v}_i)\sum\limits_{k\neq j}\left(w_{k,x}-\frac{w_k}{v_k}v_{k,x}\right)a_{jk,x}-\sum\limits_{j\neq i}\tilde{\la}_{i,v}\xi_i\pa_{v_j}\bar{v}_i(\la_i^*-\la_j^*)v_iv_{j,x}. \label{T8n}
%	,\\
%	\mathcal{T}_9&=
%	\sum\limits_{j\neq i}\tilde{\la}_{i,v}(w_{j,x}z_i-\hat{z}_{j,x}v_{i})\xi_i\pa_{w_j}\bar{v}_i-\sum\limits_{j\neq i}\tilde{\la}_{i,v}v_{i}\xi_i(\pa_{w_j}\bar{v}_i)[\psi_j-\mathcal{F}_j]\\
%	&+\sum\limits_{j\neq i}\tilde{\la}_{i,v}v_{i}\xi_i(\pa_{w_j}\bar{v}_i)\sum\limits_{k\neq j}\hat{a}_{jk}\left(w_{k,xx}-\frac{w_k}{v_k}v_{k,xx}-\left(\frac{w_k}{v_k}\right)_xv_{k,x}\right)\\
%	&+\sum\limits_{j\neq i}\tilde{\la}_{i,v}v_{i}\xi_i(\pa_{w_j}\bar{v}_i)\sum\limits_{k\neq j}\left(w_{k,x}-\frac{w_k}{v_k}v_{k,x}\right)\hat{a}_{jk,x}-\sum\limits_{j\neq i}\tilde{\la}_{i,v}\xi_i\pa_{w_j}\bar{v}_i(\la_i^*-\la_j^*)v_iw_{j,x}.
\end{align}
Similarly applying \eqref{cle1} and \eqref{cle2}, we obtain
\begin{align*}
	&(z_{i}-\la_i^*v_i)\left(\frac{w_i}{v_i}\right)_x-v_{i}\left(\frac{w_i}{v_i}\right)_t\\
	&=\frac{1}{v_i}\left[(z_{i}-\la_i^* v_i)w_{i,x}-\frac{(z_{i}-\la_i^*v_i)w_iv_{i,x}}{v_i}-v_iw_{i,t}+ w_iv_{i,t}\right]\\
	&=\frac{1}{v_i}\left[z_{i}w_{i,x}-\frac{z_{i}w_iv_{i,x}}{v_i}-v_i\hat{z}_{i,x}+ w_iz_{i,x}\right]-\psi_i+\mathcal{F}_i+\frac{w_i}{v_i}(\phi_i-\mathcal{E}_i)\\
	&+ \sum\limits_{j\neq i}\hat{a}_{ij}\left(w_{j,x}-\frac{w_j}{v_j}v_{j,x}\right)_x-\frac{w_i}{v_i}\sum\limits_{j\neq i}a_{ij}\left(w_{j,x}-\frac{w_j}{v_j}v_{j,x}\right)_x\\
	&+ \sum\limits_{j\neq i}\hat{a}_{ij,x}\left(w_{j,x}-\frac{w_j}{v_j}v_{j,x}\right)-\frac{w_i}{v_i}\sum\limits_{j\neq i}a_{ij,x}\left(w_{j,x}-\frac{w_j}{v_j}v_{j,x}\right)\\
	&=\frac{1}{v_i}\left[\frac{w_i}{v_i}(z_{i,x}v_i-z_iv_{i,x})+\mu_i v_{i,x}w_{i,x}-(\tilde{\la}_i-\la_i^*)v_iw_{i,x}-v_i\hat{z}_{i,x}\right]\\
	&-\psi_i+\mathcal{F}_i+\frac{w_i}{v_i}(\phi_i-\mathcal{E}_i)+\sum\limits_{j\neq i}a_{ij}\frac{w_{i,x}}{v_i}\left(w_{j,x}-\frac{w_j}{v_j}v_{j,x}\right)\\
	&+  \sum\limits_{j\neq i}\hat{a}_{ij}\left(w_{j,x}-\frac{w_j}{v_j}v_{j,x}\right)_x-\frac{w_i}{v_i}\sum\limits_{j\neq i}a_{ij}\left(w_{j,x}-\frac{w_j}{v_j}v_{j,x}\right)_x\\
	&+ \sum\limits_{j\neq i}\hat{a}_{ij,x}\left(w_{j,x}-\frac{w_j}{v_j}v_{j,x}\right)-\frac{w_i}{v_i}\sum\limits_{j\neq i}a_{ij,x}\left(w_{j,x}-\frac{w_j}{v_j}v_{j,x}\right).
\end{align*}
We can obtain from \eqref{defzihat}
\begin{align*}
	&\mu_i v_{i,x}w_{i,x}-(\tilde{\la}_i-\la_i^*)v_iw_{i,x}-v_i\hat{z}_{i,x}\\
	&=v_{i,x}(\mu_i w_{i,x}-(\tilde{\la}_i-\la_i^*)w_i)-v_i\hat{z}_{i,x}+(\tilde{\la}_i-\la_i^*)[v_{i,x}w_i-v_iw_{i,x}]\\
	&=v_{i,x}\hat{z}_i-v_i\hat{z}_{i,x}+(\tilde{\la}_i-\la_i^*)[v_{i,x}w_i-v_iw_{i,x}]-\sum\limits_{j\neq i}\hat{a}_{ij}v_{i,x}\left(w_{j,x}-\frac{w_j}{v_j}v_{j,x}\right).
\end{align*}
It implies that
\begin{align*}
	&(z_{i}-\la_i^*v_i)\left(\frac{w_i}{v_i}\right)_x-v_{i}\left(\frac{w_i}{v_i}\right)_t\\
	&=\frac{1}{v_i}\left[\frac{w_i}{v_i}(z_{i,x}v_i-z_iv_{i,x})+v_{i,x}\hat{z}_i-v_i\hat{z}_{i,x}+(\tilde{\la}_i-\la_i^*)(v_{i,x}w_i-v_iw_{i,x})\right]\\
	&-\psi_i+\mathcal{F}_i+\frac{w_i}{v_i}(\phi_i-\mathcal{E}_i)+\sum\limits_{j\neq i}a_{ij}\frac{w_{i,x}}{v_i}\left(w_{j,x}-\frac{w_j}{v_j}v_{j,x}\right)-\sum\limits_{j\neq i}\hat{a}_{ij}\frac{v_{i,x}}{v_i}\left(w_{j,x}-\frac{w_j}{v_j}v_{j,x}\right)\\
	&+  \sum\limits_{j\neq i}\hat{a}_{ij}\left(w_{j,x}-\frac{w_j}{v_j}v_{j,x}\right)_x-\frac{w_i}{v_i}\sum\limits_{j\neq i}a_{ij}\left(w_{j,x}-\frac{w_j}{v_j}v_{j,x}\right)_x\\
	&+ \sum\limits_{j\neq i}\hat{a}_{ij,x}\left(w_{j,x}-\frac{w_j}{v_j}v_{j,x}\right)-\frac{w_i}{v_i}\sum\limits_{j\neq i}a_{ij,x}\left(w_{j,x}-\frac{w_j}{v_j}v_{j,x}\right).
\end{align*}
Therefore, 
\begin{align}
	\mathcal{T}_{9}&=(\tilde{\la}_{i,v}\xi_i^\p\bar{v}_i-\tilde{\lambda}_{i,\sig}\theta^\p_i)\left[\left(\frac{w_i}{v_i}\right)_x(z_i-\la_i^*v_{i})-\left(\frac{w_i}{v_i}\right)_t v_i\right]\nonumber\\
	&=(\tilde{\la}_{i,v}\xi_i^\p\bar{v}_i-\tilde{\lambda}_{i,\sig}\theta^\p_i)\frac{1}{v_i}\left[\frac{w_i}{v_i}(z_{i,x}v_i-z_iv_{i,x})+v_{i,x}\hat{z}_i-v_i\hat{z}_{i,x}+(\tilde{\la}_i-\la_i^*)(v_{i,x}w_i-v_iw_{i,x})\right]\nonumber\\
	&+(\tilde{\la}_{i,v}\xi_i^\p\bar{v}_i-\tilde{\lambda}_{i,\sig}\theta^\p_i)\left[-\psi_i+\mathcal{F}_i+\frac{w_i}{v_i}(\phi_i-\mathcal{E}_i)+\sum\limits_{j\neq i}a_{ij}\frac{w_{i,x}}{v_i}\left(w_{j,x}-\frac{w_j}{v_j}v_{j,x}\right)\right]\nonumber\\
	&-(\tilde{\la}_{i,v}\xi_i^\p\bar{v}_i-\tilde{\lambda}_{i,\sig}\theta^\p_i)\sum\limits_{j\neq i}\hat{a}_{ij}\frac{v_{i,x}}{v_i}\left(w_{j,x}-\frac{w_j}{v_j}v_{j,x}\right)\nonumber\\
	&+ (\tilde{\la}_{i,v}\xi_i^\p\bar{v}_i-\tilde{\lambda}_{i,\sig}\theta^\p_i)\left[\sum\limits_{j\neq i}\hat{a}_{ij}\left(w_{j,x}-\frac{w_j}{v_j}v_{j,x}\right)_x-\frac{w_i}{v_i}\sum\limits_{j\neq i}a_{ij}\left(w_{j,x}-\frac{w_j}{v_j}v_{j,x}\right)_x\right]\nonumber\\
	&+ (\tilde{\la}_{i,v}\xi_i^\p\bar{v}_i-\tilde{\lambda}_{i,\sig}\theta^\p_i)\left[\sum\limits_{j\neq i}\hat{a}_{ij,x}\left(w_{j,x}-\frac{w_j}{v_j}v_{j,x}\right)-\frac{w_i}{v_i}\sum\limits_{j\neq i}a_{ij,x}\left(w_{j,x}-\frac{w_j}{v_j}v_{j,x}\right)\right].\label{T9n}
\end{align}
Summarizing, the above calculations we get that
\begin{equation}\label{eqfonda1}
	z_{i,t}+(\tilde{\la}_{i}z_{i})_x-(\mu_iz_{i,x})_x=\mu_i\phi_{i,x}+\sum\limits_{l=1}^{9}\mathcal{T}_l=:\Phi_i.
\end{equation}
\subsection{Derivation of equations for effective fluxes for \eqref{eqn-w-i-1}}\label{section1hatzi}
Next we deduce the parabolic equation satisfies by $\hat{z}_i$. Using the definition of $\hat{z}_i$ and \eqref{cle2}, we get
\begin{align*}
	&(\mu_iw_{i,x})_t-(\tilde{\la}_iw_i)_t+\la_i^*w_{i,t}\\
	&=\mu_iw_{i,tx}-(\tilde{\la}_{i}-\la_i^*)w_{i,t}+\sum\limits_{j}(w_{i,x}\mu_{i,u}-w_i\tilde{\la}_{i,u})\tilde{r}_j(w_j-\la_j^*v_j)\\
	&-\tilde{\la}_{i,v}v_{i,t}\xi_i\pa_{v_i}\bar{v}_iw_{i}-\sum\limits_{j\neq i}\tilde{\la}_{i,v}v_{j,t}\xi_iw_{i}\pa_{v_j}\bar{v}_i
	-(\tilde{\la}_{i,v}\bar{v}_i\xi_i^\p-\tilde{\la}_{i,\si}\theta_i^\p)w_{i}\left(\frac{w_i}{v_i}\right)_t\\
	&=(\mu_i\hat{z}_{i,x})_x-(\tilde{\lambda}_i \hat{z}_i)_x-\mu_{i,x}(\hat{z}_{i,x}-\la_i^*w_{i,x})+\tilde{\la}_{i,x}(\hat{z}_i-\la_i^*w_i)
	+\mu_i\psi_{i,x}-\mu_i\mathcal{F}_{i,x}
	\\
	&-\sum\limits_{j\neq i}\hat{a}_{ij}\mu_i\left(w_{j,x}-\frac{w_j}{v_j}v_{j,x}\right)_{xx}-2\mu_i\sum\limits_{j\neq i}\left(w_{j,x}-\frac{w_j}{v_j}v_{j,x}\right)_x \hat{a}_{ij,x}\\
	&-\sum\limits_{j\neq i}\mu_i\left(w_{j,x}-\frac{w_j}{v_j}v_{j,x}\right)\hat{a}_{ij,xx}-(\tilde{\la}_i-\la_i^*)\psi_i+(\tilde{\la}_i-\la_i^*)\mathcal{F}_i\\
	&+\tilde{\la}_i\sum\limits_{j\neq i}\hat{a}_{ij}\left(w_{j,x}-\frac{w_j}{v_j}v_{j,x}\right)_{x}+\tilde{\la}_i\sum\limits_{j\neq i}\hat{a}_{ij,x}\left(w_{j,x}-\frac{w_j}{v_j}v_{j,x}\right)\\
	&+\sum\limits_{j}(w_{i,x}\mu_{i,u}-w_i\tilde{\la}_{i,u})\tilde{r}_j(w_j-\la_j^*v_j)-\tilde{\la}_{i,v}v_{i,t}\xi_i\pa_{v_i}\bar{v}_iw_{i}-\sum\limits_{j\neq i}\tilde{\la}_{i,v}v_{j,t}\xi_iw_{i}\pa_{v_j}\bar{v}_i\\
	&-(\tilde{\la}_{i,v}\bar{v}_i\xi_i^\p-\tilde{\la}_{i,\si}\theta_i^\p)w_{i}\left(\frac{w_i}{v_i}\right)_t.
\end{align*}	
With further rearrangement of terms, we obtain
\begin{align*}
	&(\mu_iw_{i,x})_t-(\tilde{\la}_iw_i)_t+\la_i^*w_{i,t}=(\mu_i\hat{z}_{i,x})_x-(\tilde{\lambda}_i \hat{z}_i)_x-(\tilde{\la}_i-\la_i^*)\psi_i+(\tilde{\la}_i-\la_i^*)\mathcal{F}_i\\
	&+\mu_i\psi_{i,x}-\mu_i\mathcal{F}_{i,x}-\sum\limits_{j\neq i}\hat{a}_{ij}\mu_i\left(w_{j,x}-\frac{w_j}{v_j}v_{j,x}\right)_{xx}\\
	&-2\mu_i\sum\limits_{j\neq i}\left(w_{j,x}-\frac{w_j}{v_j}v_{j,x}\right)_x \hat{a}_{ij,x}-\sum\limits_{j\neq i}\mu_i\left(w_{j,x}-\frac{w_j}{v_j}v_{j,x}\right)\hat{a}_{ij,xx}\\
	&+\tilde{\la}_i\sum\limits_{j\neq i}\hat{a}_{ij}\left(w_{j,x}-\frac{w_j}{v_j}v_{j,x}\right)_{x}+\tilde{\la}_i\sum\limits_{j\neq i}\hat{a}_{ij,x}\left(w_{j,x}-\frac{w_j}{v_j}v_{j,x}\right)\\
	&
	+\sum\limits_{j}(w_{i,x}(w_j-\la_j^*v_j)-v_j(\hat{z}_{i,x}-\la_i^*w_{i,x}))\mu_{i,u}\tilde{r}_j\\
	&-\sum\limits_{j}(w_i(w_j-\la_j^*v_j)-(\hat{z}_i-\la_i^*w_i)v_j)\tilde{\la}_{i,u}\tilde{r}_j+\tilde{\la}_{i,v}(v_{i,x}(\hat{z}_i-\la_i^*w_i)-v_{i,t}w_{i})\xi_i\pa_{v_i}\bar{v}_i\\
	&+\sum\limits_{j\neq i}\tilde{\la}_{i,v}(v_{j,x}(\hat{z}_i-\la_i^*w_i)-v_{j,t}w_{i})\xi_i\pa_{v_j}\bar{v}_i\\
	&+(\tilde{\la}_{i,v}\bar{v}_i\xi_i^\p-\tilde{\la}_{i,\si}\theta_i^\p)\left[(\hat{z}_i-\la_i^*w_i)\left(\frac{w_i}{v_i}\right)_t-w_{i}\left(\frac{w_i}{v_i}\right)_t\right].
\end{align*}
From \eqref{eqn-hat-z-i-t} we deduce that
\begin{align}
	\hat{z}_{i,t}&=(\mu_i\hat{z}_{i,x}-\tilde{\la}_{i}\hat{z}_{i})_x-(\tilde{\la}_i-\la_i^*)\psi_i+(\tilde{\la}_i-\la_i^*)\mathcal{F}_i+\mu_i\psi_{i,x}-\mu_i\mathcal{F}_{i,x}
	\nonumber\\
	&+\sum\limits_{j\neq i}\hat{a}_{ij}\left[ \left(w_{j,x}-\frac{w_j}{v_j}v_{j,x}\right)_t-\mu_i\left(w_{j,x}-\frac{w_j}{v_j}v_{j,x}\right)_{xx}\right]\nonumber\\
	&-2\mu_i\sum\limits_{j\neq i}\left(w_{j,x}-\frac{w_j}{v_j}v_{j,x}\right)_x \hat{a}_{ij,x}+\sum\limits_{j\neq i}\left(w_{j,x}-\frac{w_j}{v_j}v_{j,x}\right)(\hat{a}_{ij,t}-\mu_i\hat{a}_{ij,xx})\nonumber\\
	&+\tilde{\la}_i\sum\limits_{j\neq i}\hat{a}_{ij}\left(w_{j,x}-\frac{w_j}{v_j}v_{j,x}\right)_{x}+\tilde{\la}_i\sum\limits_{j\neq i}\hat{a}_{ij,x}\left(w_{j,x}-\frac{w_j}{v_j}v_{j,x}\right)\nonumber\\
	&
	+\sum\limits_{j}(w_{i,x}(w_j-\la_j^*v_j)-v_j(\hat{z}_{i,x}-\la_i^*w_{i,x}))\mu_{i,u}\tilde{r}_j\nonumber\\
	&-\sum\limits_{j}(w_i(w_j-\la_j^*v_j)-(\hat{z}_i-\la_i^*w_i)v_j)\tilde{\la}_{i,u}\tilde{r}_j\nonumber\\
	&+\tilde{\la}_{i,v}(v_{i,x}(\hat{z}_i-\la_i^*w_i)-v_{i,t}w_{i})\xi_i\pa_{v_i}\bar{v}_i \nonumber\\
	&+\sum\limits_{j\neq i}\tilde{\la}_{i,v}(v_{j,x}(\hat{z}_i-\la_i^*w_i)-v_{j,t}w_{i})\xi_i\pa_{v_j}\bar{v}_i\nonumber\\
	&+(\tilde{\la}_{i,v}\bar{v}_i\xi_i^\p-\tilde{\la}_{i,\si}\theta_i^\p)\left[(\hat{z}_i-\la_i^*w_i)\left(\frac{w_i}{v_i}\right)_t-w_{i}\left(\frac{w_i}{v_i}\right)_t\right]\nonumber\\
		&=(\mu_i\hat{z}_{i,x}-\tilde{\lambda}_i \hat{z}_i)_x+\mu_i\psi_{i,x}+\sum\limits_{l=1}^{9}\widehat{\mathcal{T}}_l,
\end{align}
where $\widehat{\mathcal{T}}_l$ are defined as follows
\begin{align}
\widehat{\mathcal{T}}_1&=   -\mu_i\mathcal{F}_{i,x}+\sum\limits_{j\neq i}\hat{a}_{ij}\left[ \left(w_{j,x}-\frac{w_j}{v_j}v_{j,x}\right)_t-\mu_i\left(w_{j,x}-\frac{w_j}{v_j}v_{j,x}\right)_{xx}\right],\nonumber\\
\widehat{\mathcal{T}}_2&=-(\tilde{\la}_i-\la_i^*)\psi_i+(\tilde{\la}_i-\la_i^*)\mathcal{F}_i,\nonumber\\
\widehat{\mathcal{T}}_3&=	-2\mu_i\sum\limits_{j\neq i}\left(w_{j,x}-\frac{w_j}{v_j}v_{j,x}\right)_x \hat{a}_{ij,x}+\sum\limits_{j\neq i}\left(w_{j,x}-\frac{w_j}{v_j}v_{j,x}\right)(\hat{a}_{ij,t}-\mu_i\hat{a}_{ij,xx}),\nonumber\\
\widehat{\mathcal{T}}_4&=\tilde{\la}_i\sum\limits_{j\neq i}\hat{a}_{ij}\left(w_{j,x}-\frac{w_j}{v_j}v_{j,x}\right)_{x}+\tilde{\la}_i\sum\limits_{j\neq i}\hat{a}_{ij,x}\left(w_{j,x}-\frac{w_j}{v_j}v_{j,x}\right),\nonumber\\
\widehat{\mathcal{T}}_5&=\sum\limits_{j}(w_{i,x}(w_j-\la_j^*v_j)-v_j(\hat{z}_{i,x}-\la_i^*w_{i,x}))\mu_{i,u}\tilde{r}_j,\nonumber\\
\widehat{\mathcal{T}}_6&=-\sum\limits_{j}(w_i(w_j-\la_j^*v_j)-(\hat{z}_i-\la_i^*w_i)v_j)\tilde{\la}_{i,u}\tilde{r}_j,\nonumber\\
\widehat{\mathcal{T}}_7&=\tilde{\la}_{i,v}(v_{i,x}(\hat{z}_i-\la_i^*w_i)-v_{i,t}w_{i})\xi_i\pa_{v_i}\bar{v}_i,\nonumber\\
\widehat{\mathcal{T}}_8&=\sum\limits_{j\neq i}\tilde{\la}_{i,v}(v_{j,x}(\hat{z}_i-\la_i^*w_i)-v_{j,t}w_{i})\xi_i\pa_{v_j}\bar{v}_i,\nonumber\\
\widehat{\mathcal{T}}_{9}&=(\tilde{\la}_{i,v}\bar{v}_i\xi_i^\p-\tilde{\la}_{i,\si}\theta_i^\p)\left[(\hat{z}_i-\la_i^*w_i)\left(\frac{w_i}{v_i}\right)_t-w_{i}\left(\frac{w_i}{v_i}\right)_t\right].\label{Tihat}
\end{align}
By using \eqref{z-cal-1-1} we obtain
\begin{align}
	&\sum\limits_{j\neq i}\hat{a}_{ij}\left[\left(w_{j,x}-\frac{w_j}{v_j}v_{j,x}\right)_t-\mu_i\left(w_{j,x}-\frac{w_j}{v_j}v_{j,x}\right)_{xx}\right]\nonumber\\
	&=\sum\limits_{j\neq i}\hat{a}_{ij}(\mu_j-\mu_i)\left(w_{j,xxx}-\frac{w_j}{v_j}v_{j,xxx}\right)\nonumber\\
	&+\sum\limits_{j\neq i}\hat{a}_{ij}\left[\mu_{j,xx}v_j\left(\frac{w_j}{v_j}\right)_x+(2\mu_{j,x}-\tilde{\la}_{j})\left(w_{j,xx}-\frac{w_j}{v_j}v_{j,xx}\right)-2\tilde{\la}_{j,x}v_j\left(\frac{w_j}{v_j}\right)_x\right]\nonumber\\
	&+\sum\limits_{j\neq i}\hat{a}_{ij}\left[\psi_{j,x}-\frac{w_j}{v_j}\phi_{j,x}+\frac{w_j}{v_j}\mathcal{E}_{j,x}-\mathcal{F}_{j,x}-(\mu_{j,x}-\tilde{\la}_j)v_{j,x}\left(\frac{w_j}{v_j}\right)_x\right]\nonumber\\
	&+\sum\limits_{j\neq i}\hat{a}_{ij}\left[-\frac{v_{j,x}}{v_j}(\mu_j-\mu_i)\left(w_{j,xx}-\frac{w_j}{v_j}v_{j,xx}\right)-\frac{v_{j,x}}{v_j}\left(\psi_j-\mathcal{F}_j\right)+\frac{w_jv_{j,x}}{v^2_j}\left(\phi_j-\mathcal{E}_j\right)\right]\nonumber\\
	&+\sum\limits_{j\neq i}\hat{a}_{ij}\left[2\mu_i\left(\frac{w_j}{v_j}\right)_x v_{j,xx}-\frac{2\mu_iv_{j,x}^2}{v_j}\left(\frac{w_j}{v_j}\right)_x\right].\label{T1bis}
\end{align}
From the definition of $\mathcal{F}_i$ and the fact that  $\left(\frac{w_j}{v_j}\right)_{xx}=\frac{w_{j,xx}}{v_j}-\frac{v_{j,xx}w_j}{v_j^2}-2\frac{v_{j,x}}{v_j}\left(\frac{w_j}{v_j}\right)_x$ we have
\begin{align*}
	\mathcal{F}_{i,x}&=-\sum\limits_{j\neq i}(\mu_i-\mu_j)\hat{b}_{ij}\left(w_{j,x}-\frac{w_j}{v_j}v_{j,x}\right)_{xx}-\sum\limits_{j\neq i}[(\mu_i-\mu_j)\hat{b}_{ij}]_x\left(w_{j,x}-\frac{w_j}{v_j}v_{j,x}\right)_{x}\\
	&=-\sum\limits_{j\neq i}(\mu_i-\mu_j)\hat{b}_{ij}\left(w_{j,xxx}-\frac{w_j}{v_j}v_{j,xxx}\right)+\sum\limits_{j\neq i}2(\mu_i-\mu_j)\hat{b}_{ij}\left(\frac{w_j}{v_j}\right)_xv_{j,xx}\\
	&+\sum\limits_{j\neq i}(\mu_i-\mu_j)\hat{b}_{ij}\frac{v_{j,x}}{v_j}\left(w_{j,xx}-\frac{w_j}{v_j}v_{j,xx}\right)-\sum\limits_{j\neq i}2(\mu_i-\mu_j)\hat{b}_{ij}\frac{v^2_{j,x}}{v_j}\left(\frac{w_j}{v_j}\right)_x\\
	&-\sum\limits_{j\neq i}(\mu_i-\mu_j)\hat{b}_{ij,x}\left(w_{j,x}-\frac{w_j}{v_j}v_{j,x}\right)_{x}-\sum\limits_{j\neq i}(\mu_i-\mu_j)_x\hat{b}_{ij}\left(w_{j,x}-\frac{w_j}{v_j}v_{j,x}\right)_{x}.
\end{align*}
Therefore,
\begin{align}
	&-\mu_i\mathcal{F}_{i,x}+\sum\limits_{j\neq i}\hat{a}_{ij}(\mu_j-\mu_i)\left(w_{j,xxx}-\frac{w_j}{v_j}v_{j,xxx}\right)\nonumber\\
	&=\sum\limits_{j\neq i}(\mu_j-\mu_i)(\hat{a}_{ij}-\mu_i\hat{b}_{ij})\left(w_{j,xxx}-\frac{w_j}{v_j}v_{j,xxx}\right)-\mu_i\sum\limits_{j\neq i}2(\mu_i-\mu_j)\hat{b}_{ij}\left(\frac{w_j}{v_j}\right)_xv_{j,xx}\nonumber\\
	&-\mu_i\sum\limits_{j\neq i}(\mu_i-\mu_j)\hat{b}_{ij}\frac{v_{j,x}}{v_j}\left(w_{j,xx}-\frac{w_j}{v_j}v_{j,xx}\right)+\mu_i\sum\limits_{j\neq i}2(\mu_i-\mu_j)\hat{b}_{ij}\frac{v^2_{j,x}}{v_j}\left(\frac{w_j}{v_j}\right)_x\nonumber\\
	&+\mu_i\sum\limits_{j\neq i}(\mu_i-\mu_j)\hat{b}_{ij,x}\left(w_{j,x}-\frac{w_j}{v_j}v_{j,x}\right)_{x}+\mu_i\sum\limits_{j\neq i}(\mu_i-\mu_j)_x\hat{b}_{ij}\left(w_{j,x}-\frac{w_j}{v_j}v_{j,x}\right)_{x}.\label{T1bisa}
\end{align}
Taking now $\hat{a}_{ij}=\mu_i\hat{b}_{ij}$, it gives from \eqref{T1}, \eqref{T1bis} and \eqref{T1bisa}
\begin{align}
	\widehat{\mathcal{T}}_1
	&=-\mu_i\sum\limits_{j\neq i}2(\mu_i-\mu_j)\hat{b}_{ij}\left(\frac{w_j}{v_j}\right)_xv_{j,xx}\nonumber\\
	&-\mu_i\sum\limits_{j\neq i}(\mu_i-\mu_j)\hat{b}_{ij}\frac{v_{j,x}}{v_j}\left(w_{j,xx}-\frac{w_j}{v_j}v_{j,xx}\right)+\mu_i\sum\limits_{j\neq i}2(\mu_i-\mu_j)\hat{b}_{ij}\frac{v^2_{j,x}}{v_j}\left(\frac{w_j}{v_j}\right)_x\nonumber\\
	&+\mu_i\sum\limits_{j\neq i}(\mu_i-\mu_j)\hat{b}_{ij,x}\left(w_{j,x}-\frac{w_j}{v_j}v_{j,x}\right)_{x}+\mu_i\sum\limits_{j\neq i}(\mu_i-\mu_j)_x\hat{b}_{ij}\left(w_{j,x}-\frac{w_j}{v_j}v_{j,x}\right)_{x}\nonumber\\
	&+\sum\limits_{j\neq i}\hat{a}_{ij}\left[\mu_{j,xx}v_j\left(\frac{w_j}{v_j}\right)_x+(2\mu_{j,x}-\tilde{\la}_{j})\left(w_{j,xx}-\frac{w_j}{v_j}v_{j,xx}\right)-2\tilde{\la}_{j,x}v_j\left(\frac{w_j}{v_j}\right)_x\right]\nonumber\\
	&+\sum\limits_{j\neq i}\hat{a}_{ij}\left[\psi_{j,x}-\frac{w_j}{v_j}\phi_{j,x}+\frac{w_j}{v_j}\mathcal{E}_{j,x}-\mathcal{F}_{j,x}-(\mu_{j,x}-\tilde{\la}_j)v_{j,x}\left(\frac{w_j}{v_j}\right)_x\right]\nonumber\\
	&+\sum\limits_{j\neq i}\hat{a}_{ij}\left[-\frac{v_{j,x}}{v_j}(\mu_j-\mu_i)\left(w_{j,xx}-\frac{w_j}{v_j}v_{j,xx}\right)-\frac{v_{j,x}}{v_j}\left(\psi_j-\mathcal{F}_j\right)+\frac{w_jv_{j,x}}{v^2_j}\left(\phi_j-\mathcal{E}_j\right)\right]\nonumber\\
	&+\sum\limits_{j\neq i}\hat{a}_{ij}\left[2\mu_i\left(\frac{w_j}{v_j}\right)_x v_{j,xx}-\frac{2\mu_iv_{j,x}^2}{v_j}\left(\frac{w_j}{v_j}\right)_x\right].\label{hatT1}
\end{align}
Furthermore from \eqref{eqn-v-i-x-1a-1}  and the definition of $\hat{z}_i$ we get
\begin{align*}
	&w_{i,x}(w_i-\la_i^*v_i)-v_i(\hat{z}_{i,x}-\la_i^*w_{i,x})\\
	&=v_{i,x}\hat{z}_i-v_{i}\hat{z}_{i,x}+(\tilde{\la}_{i}-\la_i^*)(w_iv_{i,x}-v_iw_{i,x})-v_{i,x}\sum\limits_{j\neq i}\hat{a}_{ij}\left(w_{j,x}-\frac{w_j}{v_j}v_{j,x}\right)+{\cal A}_iw_{i,x}.
\end{align*}
Now, we can simplify $\widehat{\mathcal{T}}_5$ as follows
\begin{align}
	&\widehat{\mathcal{T}}_5=(w_{i,x}(w_i-\la_i^*v_i)-v_i(\hat{z}_{i,x}-\la_i^*w_{i,x}))\mu_{i,u}\tilde{r}_i+\sum\limits_{j\neq i}(w_{i,x}(w_j-\la_j^*v_j)-v_j(\hat{z}_{i,x}-\la_i^*w_{i,x}))\mu_{i,u}\tilde{r}_j\nonumber\\
&=\left[(v_{i,x}\hat{z}_i-v_{i}\hat{z}_{i,x})+(\tilde{\la}_{i}-\la_i^*)(w_iv_{i,x}-v_iw_{i,x})\right]\mu_{i,u}\tilde{r}_i\nonumber\\
&+\left[-v_{i,x}\sum\limits_{j\neq i}\hat{a}_{ij}\left(w_{j,x}-\frac{w_j}{v_j}v_{j,x}\right)+{\cal A}_iw_{i,x}\right]\mu_{i,u}\tilde{r}_i\nonumber\\
&+\sum\limits_{j\neq i}(w_{i,x}(w_j-\la_j^*v_j)-v_j(\hat{z}_{i,x}-\la_i^*w_{i,x}))\mu_{i,u}\tilde{r}_j.\label{That5}
\end{align}
Similarly we observe that using \eqref{eqn-v-i-x-1a-1}  and the definition of $\hat{z}_i$
\begin{equation*}
	w_i(w_i-\la_i^*v_i)-(\hat{z}_i-\la_i^*w_i)v_i=\mu_i (v_{i,x}w_i-v_iw_{i,x})+{\cal A}_iw_{i}-v_{i}\sum\limits_{j\neq i}\hat{a}_{ij}\left(w_{j,x}-\frac{w_j}{v_j}v_{j,x}\right).
\end{equation*}
Subsequently, it follows
\begin{align}
	\widehat{\mathcal{T}}_6&=-(w_i(w_i-\la_i^*v_i)-(\hat{z}_i-\la_i^*w_i)v_i)\tilde{\la}_{i,u}\tilde{r}_i\nonumber\\
	&-\sum\limits_{j\neq i}(w_i(w_j-\la_j^*v_j)-(\hat{z}_i-\la_i^*w_i)v_j)\tilde{\la}_{i,u}\tilde{r}_j\nonumber\\
	&=-\left[\mu_i (v_{i,x}w_i-v_iw_{i,x})+{\cal A}_iw_{i}-v_{i}\sum\limits_{j\neq i}\hat{a}_{ij}\left(w_{j,x}-\frac{w_j}{v_j}v_{j,x}\right)\right]\tilde{\la}_{i,u}\tilde{r}_i\nonumber\\
	&-\sum\limits_{j\neq i}(w_i(w_j-\la_j^*v_j)-(\hat{z}_i-\la_i^*w_i)v_j)\tilde{\la}_{i,u}\tilde{r}_j.\label{That6}
\end{align}
We also observe that using \eqref{cle1} and the definitions of $z_i, \hat{z}_i$
\begin{align*}
	&(v_{i,x}(\hat{z}_i-\la_i^*w_i)-v_{i,t}w_{i})\\
	&=(\mu_iw_{i,x}v_{i,x}-\tilde{\la}_iw_iv_{i,x}-z_{i,x}w_i+\la_i^*v_{i,x}w_i)+\sum\limits_{j\neq i}\hat{a}_{ij}v_{i,x}\left(w_{j,x}-\frac{w_j}{v_j}v_{j,x}\right)\\
	&-\phi_iw_{i}+\mathcal{E}_iw_{i}+\sum\limits_{j\neq i}a_{ij}w_{i}\left(w_{j,x}-\frac{w_j}{v_j}v_{j,x}\right)_x+\sum\limits_{j\neq i}\left(w_{j,x}-\frac{w_j}{v_j}v_{j,x}\right)w_{i}a_{ij,x}\\
	&=(z_iw_{i,x}-z_{i,x}w_i)+(\tilde{\la}_i-\la_i^*)(w_{i,x}v_i-w_iv_{i,x})\\
	&-\sum\limits_{j\neq i}a_{ij}w_{i,x}\left(w_{j,x}-\frac{w_j}{v_j}v_{j,x}\right)+\sum\limits_{j\neq i}\hat{a}_{ij}v_{i,x}\left(w_{j,x}-\frac{w_j}{v_j}v_{j,x}\right)\\
	&-\phi_iw_{i}+\mathcal{E}_iw_{i}+\sum\limits_{j\neq i}a_{ij}w_{i}\left(w_{j,x}-\frac{w_j}{v_j}v_{j,x}\right)_x+\sum\limits_{j\neq i}\left(w_{j,x}-\frac{w_j}{v_j}v_{j,x}\right)w_{i}a_{ij,x}.
\end{align*}
Then we have
\begin{align}
	\widehat{\mathcal{T}}_7&=(v_{i,x}(\hat{z}_i-\la_i^*w_i)-v_{i,t}w_{i})\tilde{\la}_{i,v}\xi_i\pa_{v_i}\bar{v}_i\nonumber\\
	&=\tilde{\la}_{i,v}\xi_i\pa_{v_i}\bar{v}_i\left[(z_iw_{i,x}-z_{i,x}w_i)+(\tilde{\la}_i-\la_i^*)(w_{i,x}v_i-w_iv_{i,x})\right]\nonumber\\
	&-\sum\limits_{j\neq i}a_{ij}w_{i,x}\left(w_{j,x}-\frac{w_j}{v_j}v_{j,x}\right)\tilde{\la}_{i,v}\xi_i\pa_{v_i}\bar{v}_i+\sum\limits_{j\neq i}\hat{a}_{ij}v_{i,x}\left(w_{j,x}-\frac{w_j}{v_j}v_{j,x}\right)\tilde{\la}_{i,v}\xi_i\pa_{v_i}\bar{v}_i\nonumber\\
	&-\phi_iw_{i}\tilde{\la}_{i,v}\xi_i\pa_{v_i}\bar{v}_i+\mathcal{E}_iw_{i}\tilde{\la}_{i,v}\xi_i\pa_{v_i}\bar{v}_i+\sum\limits_{j\neq i}a_{ij}w_{i}\left(w_{j,x}-\frac{w_j}{v_j}v_{j,x}\right)_x\tilde{\la}_{i,v}\xi_i\pa_{v_i}\bar{v}_i\nonumber\\
	&+\sum\limits_{j\neq i}\left(w_{j,x}-\frac{w_j}{v_j}v_{j,x}\right)w_{i}a_{ij,x}\tilde{\la}_{i,v}\xi_i\pa_{v_i}\bar{v}_i.\label{That7}
\end{align}
Furthermore from \eqref{cle1},
\begin{align}
	\widehat{\mathcal{T}}_8&=\sum\limits_{j\neq i}\tilde{\la}_{i,v}\xi_i\pa_{v_j}\bar{v}_i(v_{j,x}(\hat{z}_i-\la_i^*w_i)-v_{j,t}w_{i})\nonumber\\
	&=\sum\limits_{j\neq i}\tilde{\la}_{i,v}\xi_i\pa_{v_j}\bar{v}_i(v_{j,x}(\hat{z}_i-\la_i^*w_i)-(z_{j,x}-\la_j^*v_{j,x})w_{i})+\sum\limits_{j\neq i}\tilde{\la}_{i,v}\xi_iw_i\pa_{v_j}\bar{v}_i(-\phi_j+\mathcal{E}_j)\nonumber\\
	&+\sum\limits_{j\neq i}\tilde{\la}_{i,v}\xi_iw_i\pa_{v_j}\bar{v}_i\sum\limits_{k\neq j}a_{jk}\left(w_{k,xx}-\frac{w_k}{v_k}v_{k,xx}-\left(\frac{w_k}{v_k}\right)_xv_{k,x}\right)\nonumber\\
	&+\sum\limits_{j\neq i}\tilde{\la}_{i,v}\xi_iw_i\pa_{v_j}\bar{v}_i\sum\limits_{k\neq j}\left(w_{k,x}-\frac{w_k}{v_k}v_{k,x}\right)a_{jk,x}.\label{That8}
\end{align}
%---------------------Related to derivative with respect to w_i-----------------------------------
%Similarly,
%\begin{align*}
%	(w_{j,x}\hat{z}_i-w_{j,t}w_{i})
%	&=w_{j,x}(\hat{z}_i-\la_i^*w_i)-(\hat{z}_{j,x}-\la_j^*w_{j,x})w_{i}\\
%	&-\psi_jw_{i}+\mathcal{F}_jw_{i}+\sum\limits_{k\neq j}\hat{a}_{jk}w_{i}\left(w_{k,xx}-\frac{w_k}{v_k}v_{k,xx}\right)\\
%	&-\sum\limits_{k\neq j}\hat{a}_{jk}\left(\frac{w_k}{v_k}\right)_xw_{i}v_{k,x}+\sum\limits_{k\neq j}\left(w_{k,x}-\frac{w_k}{v_k}v_{k,x}\right)w_{i}\hat{a}_{kj,x},
%\end{align*}
%subsequently, we get
%\begin{align*}
%	&\sum\limits_{j\neq i}\tilde{\la}_{i,v}\xi_i\pa_{w_j}\bar{v}_i(w_{j,x}(\hat{z}_i-\la_i^*w_i)-w_{j,t}w_{i})\\
%	&=\sum\limits_{j\neq i}\tilde{\la}_{i,v}\xi_i\pa_{w_j}\bar{v}_i(w_{j,x}(\hat{z}_i-\la_i^*w_i)-(\hat{z}_{j,x}-\la_j^*w_{j,x})w_{i}-\psi_jw_{i}+\mathcal{F}_jw_{i})\\
%	&+\sum\limits_{j\neq i}\tilde{\la}_{i,v}\xi_i\pa_{w_j}\bar{v}_i\sum\limits_{k\neq j}\hat{a}_{jk}w_{i}\left(w_{k,xx}-\frac{w_k}{v_k}v_{k,xx}-\left(\frac{w_k}{v_k}\right)_xv_{k,x}\right)\\
%	&+\sum\limits_{j\neq i}\tilde{\la}_{i,v}\xi_i\pa_{w_j}\bar{v}_i\sum\limits_{k\neq j}\left(w_{k,x}-\frac{w_k}{v_k}v_{k,x}\right)w_{i}\hat{a}_{kj,x}.
%\end{align*}
Again, it follows using \eqref{cle1}, \eqref{cle2}
\begin{align*}
	&\left[(\hat{z}_{i}-\la_i^*w_i)\left(\frac{w_i}{v_i}\right)_x-w_{i}\left(\frac{w_i}{v_i}\right)_t\right]\\
	&=\frac{1}{v_i}\left[(\hat{z}_{i}-\la_i^*w_i)w_{i,x}-\frac{(\hat{z}_{i}-\la_i^*w_i)w_iv_{i,x}}{v_i}-w_iw_{i,t}+\frac{w_i^2v_{i,t}}{v_i}\right]\\
	&=\frac{1}{v_i}\left[\hat{z}_{i}w_{i,x}-w_i\hat{z}_{i,x}-\frac{\hat{z}_{i}w_iv_{i,x}}{v_i}+\frac{z_{i,x}w_i^2}{v_i}\right]\\
	&-\frac{w_i}{v_i}\psi_i+\frac{w_i}{v_i}\mathcal{F}_i+\frac{w^2_i}{v^2_i}\phi_i-\frac{w^2_i}{v^2_i}\mathcal{E}_i\\
	&+\frac{w_i}{v_i}\sum\limits_{j\neq i}\hat{a}_{ij}\left(w_{j,x}-\frac{w_j}{v_j}v_{j,x}\right)_x+\frac{w_i}{v_i}\sum\limits_{j\neq i}\left(w_{j,x}-\frac{w_j}{v_j}v_{j,x}\right)\hat{a}_{ij,x}\\
	&-\frac{w^2_i}{v^2_i}\sum\limits_{j\neq i} {a}_{ij}\left(w_{j,x}-\frac{w_j}{v_j}v_{j,x}\right)_x-\frac{w^2_i}{v^2_i}\sum\limits_{j\neq i}\left(w_{j,x}-\frac{w_j}{v_j}v_{j,x}\right){a}_{ij,x}.
\end{align*}
Note that from the definition of $z_i$ and $\hat{z}_i$ we get
\begin{align*}
	&-\frac{\hat{z}_{i}w_iv_{i,x}}{v_i}+\frac{z_{i,x}w_i^2}{v_i}=\frac{w_i}{v_i}\left( z_{i,x}w_i-z_iw_{i,x}+z_iw_{i,x}-\hat{z}_iv_{i,x} \right)\\
	&=\frac{w_i}{v_i}\left( z_{i,x}w_i-z_iw_{i,x} \right)+\frac{w_i}{v_i}\left[\mu_i v_{i,x}w_{i,x}-(\tilde{\la}_{i}-\la_i^*)w_{i,x}v_i-\mu_iw_{i,x}v_{i,x}+(\tilde{\la}_{i}-\la_i^*)w_iv_{i,x}\right]\\
	&+\sum\limits_{j\neq i}w_{i,x}\frac{w_i}{v_i}\left(w_{j,x}-\frac{w_j}{v_j}v_{j,x}\right){a}_{ij}-\sum\limits_{j\neq i}\frac{w_i}{v_i}\left(w_{j,x}-\frac{w_j}{v_j}v_{j,x}\right)v_{i,x}\hat{a}_{ij}.
\end{align*}
Therefore,
\begin{align}
\widehat{\mathcal{T}}_9&=(\tilde{\la}_{i,v}\bar{v}_i\xi_i^\p-\tilde{\la}_{i,\si}\theta_i^\p)\left[(\hat{z}_{i}-\la_i^*w_i)\left(\frac{w_i}{v_i}\right)_x-w_{i}\left(\frac{w_i}{v_i}\right)_t\right]\nonumber\\
	&=\frac{(\tilde{\la}_{i,v}\bar{v}_i\xi_i^\p-\tilde{\la}_{i,\si}\theta_i^\p)}{v_i}\left[\hat{z}_{i}w_{i,x}-w_i\hat{z}_{i,x}+\frac{w_i}{v_i}\left( z_{i,x}w_i-z_iw_{i,x} \right)+\frac{(\tilde{\la}_{i}-\la_i^*)w_i}{v_i}(v_{i,x}w_i-v_iw_{i,x})\right]\nonumber\\
	&+(\tilde{\la}_{i,v}\bar{v}_i\xi_i^\p-\tilde{\la}_{i,\si}\theta_i^\p)\left[-\frac{w_i}{v_i}\psi_i+\frac{w_i}{v_i}\mathcal{F}_i+\frac{w^2_i}{v^2_i}\phi_i-\frac{w^2_i}{v^2_i}\mathcal{E}_i\right]\nonumber\\
	&+(\tilde{\la}_{i,v}\bar{v}_i\xi_i^\p-\tilde{\la}_{i,\si}\theta_i^\p)\left[\frac{w_i}{v_i}\sum\limits_{j\neq i}\hat{a}_{ij}\left(w_{j,x}-\frac{w_j}{v_j}v_{j,x}\right)_x+\frac{w_i}{v_i}\sum\limits_{j\neq i}\left(w_{j,x}-\frac{w_j}{v_j}v_{j,x}\right)\hat{a}_{ij,x}\right]\nonumber\\
	&+(\tilde{\la}_{i,v}\bar{v}_i\xi_i^\p-\tilde{\la}_{i,\si}\theta_i^\p)\left[-\frac{w^2_i}{v^2_i}\sum\limits_{j\neq i} {a}_{ij}\left(w_{j,x}-\frac{w_j}{v_j}v_{j,x}\right)_x-\frac{w^2_i}{v^2_i}\sum\limits_{j\neq i}\left(w_{j,x}-\frac{w_j}{v_j}v_{j,x}\right){a}_{ij,x}\right]\nonumber\\
	&+(\tilde{\la}_{i,v}\bar{v}_i\xi_i^\p-\tilde{\la}_{i,\si}\theta_i^\p)\left[\sum\limits_{j\neq i}\left(w_{j,x}-\frac{w_j}{v_j}v_{j,x}\right)\frac{w_{i,x}w_i}{v_i^2}{a}_{ij}-\sum\limits_{j\neq i}\left(w_{j,x}-\frac{w_j}{v_j}v_{j,x}\right)\frac{w_i v_{i,x}}{v^2_i}\hat{a}_{ij}\right].\label{That9}
\end{align}
Finally, we write
\begin{equation}\label{eqfonda2}
	\hat{z}_{i,t}+(\tilde{\la}_{i}\hat{z}_{i})_x-(\mu_i\hat{z}_{i,x})_x=\mu_i\psi_{i,x}+\sum\limits_{l=1}^{9}\widehat{\mathcal{T}}_l=:\Psi_i.
\end{equation}
From \eqref{eqfonda1} and \eqref{eqfonda2} it follows
\begin{align}
	z_{i,t}+(\tilde{\la}_iz_i)_x-(\mu_iz_{i,x})_x&=\Phi_i,\\
	\hat{z}_{i,t}+(\tilde{\la}_i\hat{z}_i)_x-(\mu_i\hat{z}_{i,x})_x&=\Psi_i.
\end{align}
We easily observe that in $\Phi_i$, $\Psi_i$ it appears some new terms of the form
\begin{align}
	&\La_i^7=\sum_{j\ne i}\big((|z_i|+|\hat{z}_i|)(|v_j|+|w_j|+|v_{j,x}|+|w_{j,x}|+|z_j|+|\hat{z}_j|)\nonumber\\
	&\hspace{1cm}+(|z_{i,x}|+|\hat{z}_{i,x}|)(|v_j|+|w_j|+|z_j|+|\hat{z}_j|)\big),\label{Lai7}\\
	&\La_i^8=|z_iw_{i,x}-w_i z_{i,x}|+|z_iv_{i,x}-v_i z_{i,x}|+|\hat{z}_iw_{i,x}-w_i \hat{z}_{i,x}|+|\hat{z}_iv_{i,x}-v_i \hat{z}_{i,x}|\label{Lai8}
\end{align}
We claim that $\Phi_i,\Psi_i$ can be bounded as 
\begin{equation}
	\Phi_i,\Psi_i=O(1)\sum_j(\La_j^1+\La_j^2+\delta_0^2\La_j^3+\La_j^4+\La_j^5+\La_j^{6}+\La_j^{6,1}+\La_j^7+\La_j^8)+\widetilde{R}_\e
	\label{8.32}
\end{equation}
with a remainder term $\widetilde{R}_\e$ depending on $\e>0$ that we will precise later satisfying for any $T_1>\hat{t}$, $\int_{\hat{t}}^{T_1}\int_\R|\widetilde{R}_\e(s,x)| ds dx=O(1)\delta_0^2$ for $\e>0$ sufficiently small in terms of $\delta_0$ and $T_1-\hat{t}$.
The proof of this result will be the object of the section \ref{sectionphii}. It is also important to point out that in $\Phi_i$ and $\Psi_i$ appears some terms of the form $\phi_{i,x}$ and $\psi_{i,x}$. It will require to estimate this type of term, we recall from the definition of ${\mathcal{A}}$, ${\cal J}_1$, ${\cal J}_2$ and ${\cal V}$ with ${\cal V}=(\phi_1,\cdots,\phi_n,\psi_1,\cdots,\psi_n)^T$ that
$$\mathcal{A} {\cal V}=\begin{pmatrix}
	\mathcal{J}_1\\
	\mathcal{J}_2
\end{pmatrix}$$
with $\mathcal{A} $ invertible and $\mathcal{A}^{-1}=O(1)$.
Since we must estimate ${\cal V}_x=(\phi_{1,x},\cdots,\phi_{n,x},\psi_{1,x},\cdots,\psi_{n,x})^T$, we get in particular
\begin{equation}
	{\cal V}_x=\mathcal{A}^{-1}\begin{pmatrix}
		\mathcal{J}_{1,x}\\
		\mathcal{J}_{2,x}
	\end{pmatrix}-\mathcal{A}^{-1}\mathcal{A}_x\mathcal{A}^{-1}\begin{pmatrix}
		\mathcal{J}_1\\
		\mathcal{J}_2
	\end{pmatrix}
\end{equation}
Using the Lemma \ref{lemme9.9}, we will see in particular that $\mathcal{A}_x=O(1)$ which enables to claim since  $\mathcal{A}^{-1}=O(1)$ that for any $i\in\{1,\cdots,n\}$
\begin{equation}
	\phi_{i,x},\psi_{i,x}=O(1)(|\mathcal{J}_{1,x}|+|\mathcal{J}_{2,x}|+|\mathcal{J}_{1}|+|\mathcal{J}_{2}|).
	\label{supercle1}
\end{equation}

\section{Global BV estimates}
\label{globalBVn}
The goal of this section is to prove that the system \eqref{eqn-main}  admits a global solution $u$ satisfying for any $t>0$
\begin{equation}
\|u_x(t,\cdot)\|_{L^1}\leq\delta_0,
\end{equation}
provided that the initial data satisfies $\mbox{TV}\bar{u}\leq \widetilde{\delta}_0$ for $\widetilde{\delta}_0$ small enough with in addition $\lim_{x\rightarrow_\infty}\bar{u}(x)=u^*$ with $u^*\in K$.
 We have seen in the proposition \ref{prop:local-existence}  that system \eqref{eqn-main} admits a unique solution defined on $[0,\hat{t}]$ satisfying for any $t\in]0,\hat{t}]$
\begin{equation}
	\|u_x(t,\cdot)\|_{L^1}\leq\frac{\delta_0}{8C}.
	\label{7.1}
\end{equation}
with $C>1$ defined as in \eqref{bonC} provided that $\|u_0\|_{TV}\leq \frac{\delta_0}{16C\kappa}$. We observe that the constant $\kappa$ depends on the matrix $A_1^*B_1^{-1/2}(u^*)$, but it is clear that this constant $\kappa$ remains uniformly bounded in $u^*$ when $u^*$ lives in a compact set $K\subset\R^n$. We mention that this solution $u$ can be extended belong the interval time $[0,\hat{t}]$ provided that the total variation of $u$ remains small by using again the proposition  \ref{prop:local-existence}. Assume here that the solution $u$ is defined on $[0,\widetilde{T}[$ inasmuch as  $u_x(t,\cdot)\in L^1$ for any $t\in]0,\widetilde{T}[$ with $\widetilde{T}\in\bar{\R}$. Let us define $T$ the maximal time for which the solution $u$ satisfies
$$\forall t\in[0,T],\;\|u_x(t,\cdot)\|_{L^1}\leq \delta_0,$$
with $T\in\bar{\R}$. We mention that we can define this time $T$ by continuity of the application $t\rightarrow \|u_x(t,\cdot)\|_{L^1}$, this is a consequence of Corollary \ref{coro2.2}. It implies in particular that
\begin{equation}
\|u_x(T,\cdot)\|_{L^1}=\delta_0.
\label{superimpot}
\end{equation}
 Let us assume by absurd that $T<+\infty$.
Using now the Corollary \ref{coro2.2} we observe that all the quantities $|u(t,\cdot)-u^*|$,  $|u_{x}(t,\cdot)|$, $|u_{xx}(t,\cdot)|$ are small for $t\in[\hat{t},T]$ provided that $\delta_0$ is sufficiently small. In particular it enables us to apply on $u_x, u_t$ the decomposition  in $v_i, w_i$ \eqref{eqn-v-i}-\eqref{eqn-v-i} for $t\in[\hat{t},T]$ with $\e>0$ that we are going to fix later. We recall in particular that $v_i$, $w_i$ are arbitrary regular on $[\hat{t},T]\times\R$. In particular due to the Lemmas \ref{lemme6.5}, \ref{lemme6.6} we deduce $L^1$ and $L^\infty$ estimate for each $v_{i}(t,\cdot)$, $v_{i,x}(t,\cdot)$, $v_{i,xx}(t,\cdot)$, $w_{i}(t,\cdot)$, $w_{i,x}(t,\cdot)$, $w_{i,xx}(t,\cdot)$ with $t\in[\hat{t},T]$. Let us fix now
 $\e>0$, to do this we take $\e>0$ sufficiently small depending on $\delta_0$ and $T-\hat{t}$ such that the estimates \eqref{12.2bis}, \eqref{12.11}, \eqref{12.11bis}, \eqref{12.11bis1}, \eqref{12.11bis2}, \eqref{12.11bis3}, \eqref{12.11bis4}, \eqref{12.11bis5}, \eqref{12.11bis6} are satisfied. In addition we assume that $\e^{\frac{1}{N}}(T-\hat{t})\leq 1$.%if for $T\geq\hat{t}$, $u$ satisfies
%\begin{equation}\label{7.2}
%	\|u_x(t,\cdot)\|_{L^1}\leq \delta_0\mbox{ for }t\in[\hat{t},T],
%\end{equation}
%then we can define the components $v_i$, $w_i$ for $t\in[\hat{t},T]$ because we know that $\|u_{x}(t,\cdot)\|_{L^\infty}$ and $\| u_{xx}(t,\cdot)\|_{L^\infty}$ are respectively of order $\delta_0^2$ and $\delta_0^3$. To prove the uniform total variation bound for all $t\in(0,+\f)$, we fix $T\in (\hat{t},+\f)$ and then we show that $\|u_x(t,\cdot)\|_{L^1}\leq \delta_0$ for all $t\in[0,T]$ for small enough $\de_0$. Since the choice of $\de_0$ does not depend on $T$, the uniform bound holds for all $T>\hat{t}$. We fix $\ep=\frac{\de_2}{1+T}$ for some $\de_2>0$ which we choose later. 
We define now the time
$T^*$ as follows:
\begin{equation}\label{def:T*}
	T^*=\sup\Big\{t\in[\hat{t},T];\, \sum_{i=1}^n\int_{\hat{t}}^t\int_{\R}	\mathcal{G}_i dx ds\leq \frac{\delta_0}{4},\Big\},
\end{equation}
where $\mathcal{G}_i$ is defined as follows
\begin{equation}
	\mathcal{G}_i:=|\widetilde{\phi}_i(s,x)|+|\widetilde{\psi}_i(s,x)|+|\Phi_i(s,x)|+|\Psi_i(s,x)|,
\end{equation}
with $\widetilde{\phi}_i$ and $\widetilde{\phi}_i$ defined in \eqref{deftilde} and $\Phi_i$, $\Psi_i$ in \eqref{eqfonda1}, \eqref{eqfonda2}. We are going to prove that $T^*=T$ and that necessary we must have $T=+\infty$. Before to do this, let us explain why $\mathcal{G}_i$ belongs to $L^1([\hat{t},T],L^1(\R))$. First combining \eqref{deftilde}, \eqref{boundphi},  \eqref{defEjFja} and the results of the section \ref{sectionphii}, we deduce that for any $t\in[\hat{t},T]$
$$|\widetilde{\phi}_i(t,x)|+|\widetilde{\psi}_i(t,x)|=O(1)\sum_j(\La_j^1+\delta_0^2\La_j^3+\La_j^4+\La_j^5+\La_j^6+\La_j^{6,1})(t,x)+R_\e(t,x),$$
with $\int^T_{\hat{t}}\int_{\R}|R_\e(s,x)| ds dx=O(1)\delta_0^2,$
due to our choice on $\e$. Due to the definition of the different terms $\La_j^l$  (see \eqref{def:La_i-1}) and the Lemmas \ref{lemme6.5}, \ref{lemme6.6}, we deduce that $\widetilde{\phi}_i$ and $\widetilde{\psi}_i$ belong to $L^1([\hat{t},T],L^1(\R))$. We point out in particular that $\La_j^3$ is in  $L^1([\hat{t},T],L^1(\R))$ due to the fact that $\La_j^3\ne 0$ provided that $v_i^{2N}\geq \e$. Similarly combining \eqref{supercle1} and the results of the sections in the sections \ref{sectionphii}, \ref{Phii} we deduce that for any $t\in[\hat{t},T]$ we have
\begin{align}
&|\Phi_i(t,x)|+|\Psi_i(t,x)|=&O(1)\sum_j(\La_j^1+\La_j^2+\delta_0^2\La_j^3+\La_j^4+\La_j^5+\La_j^6+\La_j^{6,1}+\La_j^7+\La_j^8)(t,x)\nonumber\\
&+\widetilde{R}_\e(t,x),\label{mimptech}
\end{align}
with $\int^T_{\hat{t}}\int_{\R}|\widetilde{R}_\e(s,x)| ds dx=O(1)\delta_0^2,$
due to our choice on $\e$. We wish now to show that $\Phi_i$ and $\Psi_i$ are in  $L^1([\hat{t},T],L^1(\R))$. From \eqref{mimptech} and what we have seen previously, it suffices to prove that $\La_j^2,\La_j^7,\La_j^8$ are in $L^1([\hat{t},T],L^1(\R))$. From \eqref{defzi}, \eqref{defzihat}, it gives
\begin{align*}
&z_{i,x}=\big(\mu_iv_{i,x}-(\tilde{\la}_i-\la_i^*)v_i)_x +\sum\limits_{j\neq i}\big[a_{ij}\big(w_{j,xx}-\frac{w_j}{v_j}v_{j,xx}-\left(\frac{w_j}{v_j}\right)_xv_{j,x}\big) +a_{ij,x}v_j\left(\frac{w_j}{v_j}\right)_x  \big],\\
	&\hat{z}_{i,x}=\big(\mu_iw_{i,x}-(\tilde{\la}_i-\la_i^*)w_i)_x +\sum\limits_{j\neq i}\big[\hat{a}_{ij}\big(w_{j,xx}-\frac{w_j}{v_j}v_{j,xx}-\left(\frac{w_j}{v_j}\right)_xv_{j,x}\big) +\hat{a}_{ij,x}v_j\left(\frac{w_j}{v_j}\right)_x  \big]
\end{align*}
We deduce in particular from Lemmas 
\ref{lemme6.5}, \ref{lemme6.6}, \ref{lemme11.3}, \ref{lemme9.6} and the fact that $a_{i,j}, \hat{a}_{ij}=O(1)\mathfrak{A}_j \bar{v}_j$ that $z_i, z_{i,x}, \hat{z}_i, \hat{z}_{i,x}$ are uniformly bounded in  $L^1(\R)\cap L^\infty(\R)$ on $[\hat{t},T]$, it implies from
 \eqref{Lai7}, \eqref{Lai8} that $\La_i^7, \La_i^8$ are in $L^1([\hat{t},T],L^1(\R))$. Using now Lemmas \ref{lemme6.5}, \ref{lemmeLa2} and the fact for any $j\in\{1,\cdots,n\}$ $\La_j^3$ is in  $L^1([\hat{t},T],L^1(\R))$, we deduce that for any $j\in\{1,\cdots,n\}$  $\La_j^2$ belongs to  $L^1([\hat{t},T],L^1(\R))$. This implies finally that $\Phi_i$ and $\Psi_i$ are in  $L^1([\hat{t},T],L^1(\R))$. It means in particular that $T^*$ is well defined.
 
Let us assume now by absurd that $T^*<T$, it implies then that 
\begin{equation}
\sum_{i=1}^n\int_{\hat{t}}^{T^*}\int_{\R}\left(|\widetilde{\phi}_i(s,x)|+|\widetilde{\psi}_i(s,x)|+|\Phi_i(s,x)|+|\Psi_i(s,x)|\right) dx ds= \frac{\delta_0}{4}.
\label{contradiction}
\end{equation}
In order to obtain a contradiction it suffices to show that
\begin{equation}\label{estimate:remainder-1}
	\sum_{i=1}^n\int_{\hat{t}}^{T^*}\int_{\R}\left(|\widetilde{\phi}_i(s,x)|+|\widetilde{\psi}_i(s,x)|+|\Phi_i(s,x)|+|\Psi_i(s,x)|\right) dx ds=O(1)\delta_0^{2},
\end{equation}
with $\delta_0$ sufficiently small.
To do this we are going to study the $L^1_{t,x}$ of each terms of $\widetilde{\phi}_i$, $\widetilde{\psi}_i$, $\Phi_i$ and $\Psi_i$. From the claims \eqref{boundphi} and \eqref{8.32} we have to estimate the $L^1([\hat{t},T^*]\times\R)$ norm of $\sum_{l=1,l\ne 3}^8\sum_j\Lambda_j^l$ and $\delta_0^2\sum_j\Lambda_j^3$, it is what we are going to do in the next subsections. We mention that the claims \eqref{boundphi} and \eqref{8.32} are proved in the sections \ref{sectionphii}, \ref{Phii}. Let us assume now that $T^*=T$, in order to conclude the proof of the Theorem \ref{theorem:BV-estimate}, we have to show that $T<+\infty$ is absurd. 
Applying the maximum principle to the equation \eqref{7.51}, we deduce that for any $t\in[\hat{t},T]$ 
\begin{align}
	\|u_x(t,\cdot)\|_{L^1}&\leq \sum_{i=1}^n\|\tilde{r}_i(u,\xi_i\bar{v}_i,\sigma_i)\|_{L^\infty([\hat{t},T^*],L^\infty(\R))}\|v_i(t,\cdot)\|_{L^1}\nonumber\\
	&\leq \max_{1\leq i\leq n}\|\tilde{r}_i(u,\xi_i\bar{v}_i,\sigma_i)\|_{L^\infty([\hat{t},T^*],L^\infty(\R))}\sum_{i=1}^n
	\left(\|v_i(\hat{t},\cdot)\|_{L^1}+\int^t_{\hat{t}}\int_{\R}|\phi_i(s,x)| dxds\right).\label{paraboesti}
\end{align}
Since for any $i\in\{1,\cdots,n\}$ we have $r_{i}(u)=1$, we obtain from the Lemmas \ref{lemme6.5}, \ref{lemme6.6}, the definition of $\tilde{r}_i(u,\xi_i\bar{v}_i,\sigma_i)=r_i(u)+\sum_{j\ne i}\psi_{ij}(u,\xi_i\bar{v}_i,\sigma_i)r_j(u)$ and  the fact that $\psi_{i,j}(u,\xi_i\bar{v}_i^\e,\sigma_i)=O(1)\xi_i\bar{v}_i^\e=O(1)\delta_0^2$ that for $\delta_0$ sufficiently small
\begin{align*}
	&\|\tilde{r}_i(u,\xi_i\bar{v}_i,\sigma_i)\|_{L^\infty([\hat{t},T^*],L^\infty(\R))}\leq 2.
\end{align*}
Since $u_x=\sum_{i=1}^n v_i \tilde{r}_i(u,\xi_i\bar{v}_i,\sigma_i)$, we deduce that for any $i\in\{1,\cdots,n\}$
\begin{align*}
	&\langle u_x,l_i(u)\rangle=v_i(t,x)+\sum_{j\ne i}\psi_{ji} v_j(t,x)
\end{align*}
It implies again since $\psi_{j,i}(u,\xi_i\bar{v}_j,\sigma_j)=O(1)\delta_0^2$ for $t\in[\hat{t},T^*]$ that there exists $C>0$ such that:
\begin{align}
	&\sum_{i=1}^n\|v_i(\hat{t},\cdot)\|\leq C\|u_x(\hat{t},\cdot)\|\label{bonC}
\end{align}
It implies in particular that for $t\in[\hat{t},T]$
\begin{align*}
	&\|u_x(t,\cdot)\|_{L^1}\leq 2C\|u_x(\hat{t},\cdot)\|_{L^1}+\frac{\delta_0}{2}\leq \frac{3\delta_0}{4}<\delta_0.
\end{align*}
It contradicts the assumption \eqref{superimpot} which implies finally that $T=+\infty$, it shows then that the total variation remains strictly less than $\delta_0$ for all $t\in[\hat{t},+\infty[$ and that the solution is globally defined.\\
It remains now to verify that \eqref{estimate:remainder-1} is satisfied. It will be the object of the next subsections.

\subsection{Transversal wave interactions}\label{section:transversal-wave-interaction}
We wish in a first time to prove that
\begin{align}
	&\int_{\hat{t}}^{T^*}\int_{\R}\sum\limits_{j\neq i}\big(|v_jv_i|+|v_j w_i|+|v_j v_{i,x}|+|w_j w_i|+|w_j w_{i,x}|+|w_jv_{i,x}|\big)\,dxds=O(1)\delta_0^2,\label{transvers1}\\
	&\int_{\hat{t}}^{T^*}\int_{\R}\sum_{j\ne i}(|z_i|+|\hat{z}_i|)(|v_j|+|w_j|+|v_{j,x}|+|w_{j,x}|+|z_j|+|\hat{z}_j|)\,dxds=O(1)\delta_0^2,\label{transvers2}\\
	&\int_{\hat{t}}^{T^*}\int_{\R}\sum_{j\ne i}(|z_{i,x}|+|\hat{z}_{i,x}|)(|v_j|+|w_j|+|z_j|+|\hat{z}_j|) \,dxds=O(1)\delta_0^2.\label{transvers2.1}
\end{align}
We notice in particular that $\eqref{transvers2}$ and \eqref{transvers2.1} imply in particular for any $i\in\{1,\cdots,n\}$ that
\begin{align}
&\int^{T^*}_{\hat{t}}\int_\R|\La_i^7(s,x)| ds dx=O(1)\delta_0^2.
\label{La7estim}
\end{align}
Let us recall for the moment the following Lemmas which are issue of \cite{BB-vv-lim-ann-math},\cite{HJ-temple-class} ( see section \ref{sec:proof-of-transversal} for a proof).
\begin{lemma}\label{lemma:transversal-1}
	Let $z,z^\#$ be solutions of the two independent scalar equations,
	\begin{align}
		&z_t+(\la(t,x)z)_x-(\mu z_{x})_x=\varphi(t,x),\label{eqn-z-1}\\
		&z^\#_t+(\la^\#(t,x)z^\#)_x-(\mu^\# z^\#_{x})_x=\varphi^\#(t,x),\label{eqn-z-2}
	\end{align}
	which is valid for $t\in[0,T]$. We assume that
	\begin{equation}\label{relqtion-la-12}
		\inf\limits_{t,x}\la^\#(t,x)-\sup\limits_{t,x}\la(t,x)\geq c>0,
	\end{equation}
	and $\norm{\mu_x,\mu^\#_x}_{L^\f}\leq \e_0<c/2$. Then we have
	\begin{equation}\label{est:transversal-1}
		\int\limits_{0}^{T}\int\limits_{\R}\abs{z(t,x)}\abs{z^\#(t,x)}\,dxdt\leq \frac{1}{c_0}E_1E_2,
	\end{equation}
	where $E_1,E_2$ are defined as follows
	\begin{align}
		E_1&:=\int\limits_{\R}\abs{z(0,x)}\,dx+\int\limits_{0}^{T}\int\limits_{\R}\abs{\varphi(t,x)}\,dxdt,\\
		E_2&:=\int\limits_{\R}\abs{z^\#(0,x)}\,dx+\int\limits_{0}^{T}\int\limits_{\R}\abs{\varphi^\#(t,x)}\,dxdt.
	\end{align}
\end{lemma}
Applying the Lemma \ref{lemma:transversal-1} enables us to estimate the terms  $|v_iv_j|$, $|v_j w_i|$, $|w_iw_j|$, $|z_i v_j|$, $|z_i w_j|$, $|\hat{z}_i v_j|$, $|\hat{z}_i w_j|$, $|z_iz_j|$, $|z_i\hat{z}_j|$, $|\hat{z}_i\hat{z}_j$ in \eqref{transvers1} and \eqref{transvers2}. Indeed from \eqref{stricthyperbo},  \eqref{deflatilde} and Lemma \ref{lemme6.5} we deduce that for $\delta_0$ sufficiently small we have \eqref{relqtion-la-12} which is satisfied for $\widetilde{\la}_i$ and $\widetilde{\la}_j$ when $i\ne j$ on $[\hat{t},T^*]\times\R$. Furthermore from Lemma \ref{lemme6.5} we observe that the norm $\|(\mu_i(u))_x\|_{L^\infty}$, $\|(\mu_j(u))_x\|_{L^\infty}$ can be chosen arbitrary small in term of $\delta_0$  on $[\hat{t},T^*]\times\R$. From the definition of $T^*$ and using Lemma \ref{lemme6.5}, \eqref{defzi}, \eqref{defzihat} (see \eqref{def:T*}), we have in addition that for any $i\in\{1,\cdots, n\}$
\begin{align}
&\int_{\hat{t}}^{T^*}\big(|\widetilde{\phi}_i(s,x)|+|\widetilde{\psi}_i(s,x)|+|\Phi_i(s,x)|+|\Psi_i(s,x)|\big) dx ds\leq\frac{\delta_0}{4}\label{defrT}\\
&\|v_i(\hat{t},\cdot)\|_{L1}, \|w_i(\hat{t},\cdot)\|_{L1}, \|z_i(\hat{t},\cdot)\|_{L1}, \|\hat{z}_i(\hat{t},\cdot)\|_{L1}=O(1)\delta_0,\label{defrT1}
\end{align}
which implies that the product $E_1E_2$ is of order $O(1)\delta_0^2$.
\begin{lemma}\label{lemma:transversal-2}
	Let $z,z^\#$ be solutions of \eqref{eqn-z-1}, \eqref{eqn-z-2} respectively and we assume that \eqref{relqtion-la-12} holds along with the following estimates
	\begin{align}
		&\int\limits_{0}^{T}\int\limits_{\R}\abs{\varphi(t,x)}dxdt\leq \de_0,\quad\quad 	\int\limits_{0}^{T}\int\limits_{\R}\abs{\varphi^\#(t,x)}dxdt\leq \de_0,\\
		&\norm{z(t)}_{L^1},\,\norm{z^\#(t)}_{L^1}\leq\de_0,\quad\quad \norm{z_x(t)}_{L^1},\,\norm{z^\#(t)}_{L^\f}\leq C_*\de^2_0,\\
		&\norm{\la_x(t)}_{L^\f},\,\norm{\la_x(t)}_{L^1}\leq C_*\de_0,\quad\quad \lim\limits_{x\rr-\f}\la(t,x)=0,\label{condsimpli}
	\end{align}
	for all $t\in[0,T]$. Then we have
	\begin{equation}\label{est:transversal-2}
		\int\limits_{0}^{T}\int\limits_{\R}\abs{z_x(t,x)}\abs{z^\#(t,x)}\,dxdt=O(1)\de^2_0.
	\end{equation}
\end{lemma}
Applying the Lemma \ref{lemma:transversal-2}  on $[\hat{t},T^*]$ provides the required estimates on the terms  $|v_{i,x}v_j|$, $|w_{i,x} w_j|$, $|w_j v_{i,x}|$, $|z_i v_{j,x}|$, $|z_i w_{j,x}|$, $|\hat{z}_i v_{j,x}|$, $|\hat{z}_i w_{j,x}|$, $|v_{j}z_{i,x}|$,  $|v_{j}\hat{z}_{i,x}|$,  $|w_{j}z_{i,x}|$,  $|w_{j}\hat{z}_{i,x}|$, $|z_i z_{j,x}|$, $|z_i \hat{z}_{j,x}|$, $|\hat{z}_i \hat{z}_{j,x}|$  in \eqref{transvers1}, \eqref{transvers2} and \eqref{transvers2.1}. We mention that we can apply this Lemma again due to \eqref{defzi}, \eqref{defzihat}, \eqref{defrT}, \eqref{defrT1}, Lemmas \ref{lemme6.5}, \ref{lemme6.6}, \ref{lemme11.3}. Let us note that the simplifying condition $\lim_{w\rightarrow-\infty}\lambda(t,x)$ can be easily extended provided that we use a new space coordinate $x'=x-\la_j^* t$.

\subsection{Functionals related to shortening curves}\label{sec:shorten}
We recall now the following Lemma related to shortening curves (see \cite{BB-vv-lim-ann-math}, \cite{HJ-triangular}).
%\begin{align}
%	v_t+(\la v)_x-(\mu v_x)_x&=\phi,\\
%	w_t+(\la w)_x-(\mu w_x)_x&=\psi,
%\end{align}
%We can deduce
%\begin{equation}
%	w_{t}+\la w_x-\mu w_{xx}=\psi-\la_x w+\mu_x w_x.
%\end{equation}
%Define
%\begin{equation}
%	\ga(t,x)=\left(\int\limits_{-\f}^{x}v(t,y)\,dy,\int\limits_{-\f}^{x}w(t,y)\,dy\right).
%\end{equation}
%We get
%\begin{equation}
%	\ga_t+\la \ga_x-\mu \ga_{xx}=\Phi(t,x):=\left(\int\limits_{-\f}^{x}\phi(t,y)\,dy,\int\limits_{-\f}^{x}\psi_1(t,y)\,dy\right),
%\end{equation}
%where $\psi_1:=(\psi-\la_x w+\mu_x w_x)_x$. 
\begin{lemma}\label{lemma-sc-1}
	Let $\zeta_1,\zeta_2$ be solutions of the following equations on $[0,T]$ for some $\varphi_1$ and $\varphi_2$ 
	\begin{equation}
		\begin{aligned}
			&\zeta_{1,t}+(\la \zeta_1)_x-(\alpha_1\zeta_{1,x})_x=\varphi_1,\\
			&\zeta_{2,t}+(\la \zeta_2)_x-(\alpha_1\zeta_{2,x})_x=\varphi_2,
		\end{aligned}
		\label{8.6}
	\end{equation}
	with $\alpha_1(t,x)\geq c>0$.
	For each $t$, we assume that $x\mapsto \zeta_1(t,x),\,x\mapsto \zeta_2(t,x)$ and $x\mapsto \la(t,x)$ are $C^{1,1}$. Then, we have
	\begin{equation}
		\begin{aligned}
			\frac{d}{dt}\mathcal{A}(t)\leq& -\int\limits_{\R}\alpha_1(t,x)\abs{\zeta_{1,x}(t,x)\zeta_2(t,x)-\zeta_1(t,x)\zeta_{2,x}(t,x)}\,dx\\
			&+\norm{\zeta_1(t)}_{L^1}\norm{\varphi_2(t)}_{L^1}+\norm{\zeta_2(t)}_{L^1}\norm{\varphi_1(t)}_{L^1},
		\end{aligned}
		\label{8.11}
	\end{equation}
	where $\mathcal{A}$ is defined as below
	\begin{equation}
		\mathcal{A}(t)=\frac{1}{2}\int\int\limits_{x<y}\abs{\zeta_1(t,x)\zeta_2(t,y)-\zeta_1(t,y)\zeta_2(t,x)}\,dxdy.
	\end{equation}
\end{lemma}
Applying Lemma \ref{lemma-sc-1} to $\zeta_1=v_i$ and $\zeta_2=w_i$, we deduce from Lemma \ref{lemme6.5}, \descref{B}{$(\mathcal{H}_B)$} and \eqref{defrT}
\begin{align}
	&c_1 \int_{\hat{t}}^{T^*}\int_{\R}|v_{i,x}w_i-w_{i,x}v_i|(s,x) ds dx\leq \frac{1}{2}\int\int\limits_{x<y}\abs{v_i(\hat{t},x)w_i(\hat{t},y)-v_i(\hat{t},y)w_i(\hat{t},x)}\,dxdy\nonumber\\
	&+\|v_i\|_{L^\infty([\hat{t},T^*],L^1)}\int_{\hat{t}}^{T^*}|\psi_i(s,x)| ds dx+\|w_i\|_{L^\infty([\hat{t},T^*],L^1)}\int_{\hat{t}}^{T^*}|\phi_i(s,x)| ds dx=O(1)\delta_0^2.\label{La4}
\end{align}
Similarly applying Lemma \ref{lemma-sc-1} respectively to $\zeta_1=z_i$ and $\zeta_2=w_i$, $\zeta_1=z_i$ and $\zeta_2=v_i$, $\zeta_1=\hat{z}_i$ and $\zeta_2=w_i$ and $\zeta_1=\hat{z}_i$ and $\zeta_2=v_i$ enables us to prove that
\begin{align}
	\int_{\hat{t}}^{T^*}\int_{\R}\big( |z_iw_{i,x}-w_i z_{i,x}|+|z_iv_{i,x}-v_i z_{i,x}|\big)\, dx ds
	&=O(1)\delta_0^2,\label{La4bis-1}\\
		\int_{\hat{t}}^{T^*}\int_{\R}\big( |\hat{z}_iw_{i,x}-w_i \hat{z}_{i,x}|+|\hat{z}_iv_{i,x}-v_i \hat{z}_{i,x}|\big)\, dx ds&=O(1)\delta_0^2.	\label{La4bis-2}
\end{align}
It implies in particular from \eqref{def:La_i-1} and \eqref{Lai8} that
\begin{align}
	\int_{\hat{t}}^{T^*}\int_{\R}\big( |\La^4(s,x)|+  |\La^8(s,x)| \big) ds dx=O(1)\delta_0^2.\label{La8estim}
	\end{align}
\begin{lemma}\label{lemme9.4}
	Under the assumptions of Lemma \ref{lemma-sc-1} with $\zeta_{1}=v_i$, $\zeta_{2}=w_i$, $\lambda=\tilde{\lambda}_i$, $\varphi_1=\phi_i, \varphi_2=\psi_i$, $\alpha_i=\mu_i$ with $i\in\{1,\cdots,n\} $ we have for any $\e'>0$ and any $t\in[\hat{t},T]$ %assume that $\gamma_x(t,x)\ne 0$ for every $x$ with 
	%$$\gamma(t,x)=\big(\int^x_{-\infty} v_i(t,y)dy,\int^x_{-\infty} w_i(t,y)dy\big).$$
	%Then
	\begin{equation} \label{utile}
		\frac{1}{(1+9\delta_1^2)^{\frac{3}{2}}}\int^t_{\hat{t}}\int_{\left\{|\frac{w_i}{v_i}|\leq 3\delta_1\right\}\cap \{v_i^{2N}\geq \e'\}}|\mu_iv_i| \abs{\left(\frac{w_i}{v_i}\right)_x}^2 dx\leq \mathcal{L}(\hat{t})+ \int_{\hat{t}}^t\big(\|\phi_i(s)\|_{L^1}+\|\psi_i(s)\|_{L^1}\big) ds,
	\end{equation}
	with
		\begin{equation}\label{def:L(t)}
		\mathcal{L}(t):=\int\sqrt{v_i^2(t,x)+w_i^2(t,x)}\,dx.
	\end{equation}
%	$${\cal L}(t)=\int_{\R}\sqrt{v_i^2(t)+w_i^2(t)} dx.$$
\end{lemma}
\begin{proof}
Setting $\gamma(t,x)=\big(\int^x_{-\infty} v_i(t,y)dy,\int^x_{-\infty} w_i(t,y)dy\big)$ with $t\in[\hat{t},T]$,
%	We set
%	\begin{equation}
%		\mathcal{L}(t):=\int\sqrt{v_i^2(t,x)+w_i^2(t,x)}\,dx.
%	\end{equation}
	we get
	\begin{equation}
		\ga_x=(v_i,w_i)\mbox{ and }\ga_{tx}+(\widetilde{\la}_i \ga_x)_x-(\mu_i\ga_{xx})_x=(\phi_i,\psi_i).
	\end{equation}
	We calculate
	\begin{align}
		\abs{\ga_{xx}}^2\abs{\ga_{x}}^2-\langle \ga_x,\ga_{xx}\rangle^2
		&=(v_{i,x}^2+w_{i,x}^2)(v_i^2+w_i^2)-(v_iv_{i,x}+w_{i}w_{i,x})^2\nonumber\\
		&=(v_iw_{i,x}-v_{i,x}w_i)^2=v_i^4\abs{\left(\frac{w_i}{v_i}\right)_x}^2.\label{techpract}
	\end{align}
	We set $\mathcal{L}^a(t):=\int\sqrt{ae^{-x^2}+v_i^2(t,x)+w_i^2(t,x)}\,dx$ for $a>0$. Then we calculate using integration by parts
	\begin{align*}
		&\frac{d}{dt}\mathcal{L}^a(t)=\int\frac{\langle\ga_x,\ga_{tx}\rangle}{\sqrt{ae^{-x^2}+\langle \ga_x,\ga_x\rangle}}\,dx\\
		&=\int\left[\frac{\langle\ga_x,(\mu_i\ga_{xx})_{x}\rangle}{\sqrt{a e^{-x^2}+\langle \ga_x,\ga_x\rangle}}-\frac{\langle\ga_x,(\widetilde{\la}_i \ga_x)_x\rangle}{\sqrt{ae^{-x^2}+\langle \ga_x,\ga_x\rangle}}+\frac{\langle\ga_x,(\phi_i,\psi_i)\rangle}{\sqrt{ae^{-x^2}+\langle \ga_x,\ga_x\rangle}}\right]\,dx\\
		&=\int\left[\left(\mu_i \big(\sqrt{ae^{-x^2}+\langle \ga_x,\ga_x\rangle}\big)_x \right)_{x}-\left(\widetilde{\la}_i\sqrt{a e^{-x^2}+\langle \ga_x,\ga_x\rangle}\right)_x+a\frac{-\widetilde{\la}_i xe^{-x^2}+ \mu_i e^{-x^2}(1-2x^2)}{\sqrt{ae^{-x^2}+\langle \ga_x,\ga_x\rangle}}\right]\,dx\\
		&+\int\left[-\mu_i \frac{\abs{\ga_{xx}}^2-\langle\ga_x/\sqrt{ae^{-x^2}+\langle \ga_x,\ga_x\rangle},\ga_{xx}\rangle^2}{\sqrt{ae^{-x^2}+\langle \ga_x,\ga_x\rangle}}+\frac{\langle\ga_x,(\phi_i,\psi_i)\rangle}{\sqrt{ae^{-x^2}+\langle \ga_x,\ga_x\rangle}}\right]\,dx\\
		&+\int\mu_i \frac{a^2x^2e^{-2x^2}-2axe^{-x^2}\langle\ga_x,\ga_{xx}\rangle}{\left(\sqrt{ae^{-x^2}+\langle \ga_x,\ga_x\rangle}\right)^{3/2}}\,dx+\int\frac{a e^{-x^2}( \mu_{i,x} x+\widetilde{\la}_{i,x})}{\sqrt{ae^{-x^2}+\langle \ga_x,\ga_x\rangle}}\,dx\\[2mm]
		&=\int\left[a\frac{-\widetilde{\la}_i xe^{-x^2}+ \mu_i e^{-x^2}(1-2x^2)}{\sqrt{ae^{-x^2}+\langle \ga_x,\ga_x\rangle}}\right]\,dx-\int \mu_i  \frac{a^2 x^2 e^{-2 x^2}}{\left(\sqrt{ae^{-x^2}+\langle \ga_x,\ga_x\rangle}\right)^{3/2}} \\
		&+\int\left[-\mu_i \frac{\abs{\ga_{xx}}^2-\langle\ga_x/\sqrt{ae^{-x^2}+\langle \ga_x,\ga_x\rangle},\ga_{xx}\rangle^2}{\sqrt{ae^{-x^2}+\langle \ga_x,\ga_x\rangle}}+\frac{\langle\ga_x,(\phi_i,\psi_i)\rangle}{\sqrt{ae^{-x^2}+\langle \ga_x,\ga_x\rangle}}\right]\,dx\\
		&+\int\mu_i ax e^{-x^2} \frac{2a x e^{-x^2}-2\langle\ga_x,\ga_{xx}\rangle}{\left(\sqrt{ae^{-x^2}+\langle \ga_x,\ga_x\rangle}\right)^{3/2}}\,dx+\int\frac{a e^{-x^2}( \mu_{i,x} x+\widetilde{\la}_{i,x})}{\sqrt{ae^{-x^2}+\langle \ga_x,\ga_x\rangle}}\,dx
		%&\leq -\int\frac{\abs{v}\abs{(w/v)_x}^2}{(1+(w/v)^2)^{3/2}}+\int\frac{v\phi_i+w\psi_i}{\sqrt{v^2+w^2}}\,dx.
	\end{align*}
	By integrations by parts, we obtain
	\begin{align*}
		&\frac{d}{dt}\mathcal{L}^a(t)=\int\left[a\frac{-\widetilde{\la}_i xe^{-x^2}+ \mu_i e^{-x^2}(1-2x^2)}{\sqrt{ae^{-x^2}+\langle \ga_x,\ga_x\rangle}}\right]\,dx-\int \mu_i  \frac{a^2 x^2 e^{-2 x^2}}{\left(\sqrt{ae^{-x^2}+\langle \ga_x,\ga_x\rangle}\right)^{3/2}} \\
		&+\int\left[-\mu_i \frac{\abs{\ga_{xx}}^2-\langle\ga_x/\sqrt{ae^{-x^2}+\langle \ga_x,\ga_x\rangle},\ga_{xx}\rangle^2}{\sqrt{ae^{-x^2}+\langle \ga_x,\ga_x\rangle}}+\frac{\langle\ga_x,(\phi_i,\psi_i)\rangle}{\sqrt{ae^{-x^2}+\langle \ga_x,\ga_x\rangle}}\right]\,dx\\
		&-2\int \big[\mu_{i,x}axe^{-x^2}+\mu_i a e^{-x^2}(1-2x^2)\big] \frac{1}{\sqrt{ae^{-x^2}+\langle \ga_x,\ga_x\rangle}}\,dx+\int\frac{a e^{-x^2}( \mu_{i,x} x+\widetilde{\la}_{i,x})}{\sqrt{ae^{-x^2}+\langle \ga_x,\ga_x\rangle}}\,dx
		%&\leq -\int\frac{\abs{v}\abs{(w/v)_x}^2}{(1+(w/v)^2)^{3/2}}+\int\frac{v\phi_i+w\psi_i}{\sqrt{v^2+w^2}}\,dx.
	\end{align*}
	Applying \eqref{techpract}, we deduce now since from Lemmas \ref{lemme6.5}, \ref{lemme11.3} $\mu_i,\mu_{i,x}, \widetilde{\la}_i, \widetilde{\la}_{i,x}$ are bounded on $[\hat{t},T]\times\R$  that
	\begin{align*}
		&\frac{d}{dt}\mathcal{L}^a(t)\leq \int_{\{v_i^{2N}\geq\e'\}}\left[-\mu_i \frac{\abs{\ga_{xx}}^2-\langle\ga_x/\sqrt{ae^{-x^2}+\langle \ga_x,\ga_x\rangle},\ga_{xx}\rangle^2}{\sqrt{ae^{-x^2}+\langle \ga_x,\ga_x\rangle}}\right]dx+\|\phi_i(t)\|_{L^1}+\|\psi_i(t)\|_{L^1}\\
		&+O(\sqrt{a}),\\
		&\leq \int_{\{v_i^{2N}\geq\e'\}}\left[-\mu_i \frac{\abs{\ga_{xx}}^2-\langle\frac{\ga_x}{|\ga_x|} ,\ga_{xx}\rangle^2}{\sqrt{ae^{-x^2}+\langle \ga_x,\ga_x\rangle}}\right]dx- \int_{\{v_i^{2N}\geq\e'\}}\frac{\mu_i ae^{-x^2}}{|\gamma_x|^2(ae^{-x^2}+|\gamma_x|^2)}dx		
		+\|\phi_i(t)\|_{L^1}\\
		&+\|\psi_i(t)\|_{L^1}+O(\sqrt{a})\\
		&\leq -\int _{\{v_i^{2N}\geq\e'\}} \frac{\abs{v_i}\abs{(w_i/v_i)_x}^2}{(1+(w_i/v_i)^2)\sqrt{ae^{-x^2}/v_i^2+1+(w_i/v_i)^2}}+\|\phi_i(t)\|_{L^1}+\|\psi_i(t)\|_{L^1}+O(\sqrt{a}).
		%\\
		%&-2\int \big[\mu_{i,x}axe^{-x^2}+\mu_i a e^{-x^2}(1-2x^2)\big] \frac{1}{\sqrt{ae^{-x^2}+\langle \ga_x,\ga_x\rangle}}\,dx+\int\frac{a e^{-x^2}( \mu_{i,x} x+\widetilde{\la}_{i,x})}{\sqrt{ae^{-x^2}+\langle \ga_x,\ga_x\rangle}}\,dx
		%&\leq -\int\frac{\abs{v}\abs{(w/v)_x}^2}{(1+(w/v)^2)^{3/2}}+\int\frac{v\phi_i+w\psi_i}{\sqrt{v^2+w^2}}\,dx.
	\end{align*}
	Integrating the previous estimate on $[\hat{t},t]$ and applying dominated convergence Theorem as $a\rightarrow 0^+$ we obtain the desired estimate \eqref{utile}.
	This completes the proof of Lemma \ref{lemme9.4}.
\end{proof}
From the Lemmas \ref{lemme9.4}, \ref{lemme6.5}, \eqref{eqn-v-i-1}, \eqref{eqn-w-i-1} \eqref{contradiction}, ($\mathcal{H}_B$.) we show that for $i\in\{1,\cdots,n\}$
\begin{equation}\label{La3}
	\int_{\hat{t}}^{T^*}\int_{\R} |\La_i^3(s,x)|ds dx=\int_{\hat{t}}^{T^*}\int_{\left\{|\frac{w_i}{v_i}|\leq 3\delta_1\right\}\cap \{v_i^{2N}\geq\e\}}|v_i| \abs{\left(\frac{w_i}{v_i}\right)_x}^2\, dxds=O(1)\delta_0.
\end{equation}
\subsection{Energy estimates}\label{sec:energy}

Consider a cut-off function $\widetilde{\eta}:\R\rightarrow [0,1]$ such that:
\begin{equation}
	\begin{aligned}
		&\widetilde{\eta}(s)=\begin{cases}
			\begin{aligned}
				&0\;\;\mbox{if}\;|s|\leq\frac{3}{9}\delta_1,\\[2mm]
				&1\;\;\mbox{if}\;|s|\geq\frac{3}{8}\delta_1,
			\end{aligned}
		\end{cases}
	\end{aligned}
\end{equation}
with $\widetilde{\eta}$ smooth even function. We consider $\bar{\eta}(s)=\widetilde{\eta}(|s|-\de_1/24)$.
\begin{lemma} 
			We have
	\begin{align}
		&\int_{\hat{t}}^{T^*}\int_{\R}\mathbbm{1}_{\left\{|\frac{w_i}{v_i}|\geq \frac{\delta_1}{2}\right\}}v_{i,x}^2 ds dx=O(1)\delta_0^2.
		\label{La6}
	\end{align}
\end{lemma}
\begin{proof}
	We denote $\widetilde{\eta}_i=\widetilde{\eta}(w_i/v_i)$. Multiplying the equation \eqref{7.51} by $\widetilde{\eta}_i v_i$ and integrating by parts, we obtain as in \cite{BB-vv-lim-ann-math}
	\begin{align*}
		\int_{\R}\widetilde{\eta}_i v_i\widetilde{\phi}_i dx=\frac{d}{dt}\int_{\R}\frac{1}{2}\widetilde{\eta}_i v_i^2 dx
		+\int_{\R}\Big[\tilde{\la}_{i,x}\widetilde{\eta}_i\frac{ v_i^2}{2}-\frac{1}{2}v_i^2\big(\widetilde{\eta}_{i,t} +\tilde{\la}_i \widetilde{\eta}_{i,x}-(\mu_i\widetilde{\eta}_{i,x})_x\big)+\mu_i \widetilde{\eta}_i v_{i,x}^2\quad\quad&\\
		+2\mu_iv_iv_{i,x}\widetilde{\eta}_{i,x} \Big]\,dx.&
	\end{align*}
	It implies that
	\begin{align}
		\int_{\R}\mu_i \widetilde{\eta}_i v_{i,x}^2dx=&
		-\frac{d}{dt}\int_{\R}\frac{1}{2}\widetilde{\eta}_i v_i^2 dx+\int_{\R}\widetilde{\eta}_i v_i\widetilde{\phi}_i dx\nonumber\\
		&+\int_{\R}\big(\frac{1}{2}v_i^2\Big(\widetilde{\eta}_{i,t} +\tilde{\la}_i \widetilde{\eta}_{i,x}-(\mu_i\widetilde{\eta}_{i,x})_x\big)-\tilde{\la}_{i,x}\widetilde{\eta}_i\frac{ v_i^2}{2}-2\mu_iv_iv_{i,x}\widetilde{\eta}_{i,x} \Big)dx.\label{energ1}
	\end{align}
	Direct computations give:
	\begin{align*}
		&(\mu_i \eta_{i,x})_x=\mu_i\widetilde{\eta}''_i\left(\frac{w_i}{v_i}\right)_x^2+\widetilde{\eta}'_i\left[\frac{1}{v_i}\big((\mu_i w_{i,x})_x-\frac{w_i}{v_i}(\mu_i v_{i,x})_x\big)-2\mu_i\frac{v_{i,x}}{v_i}\left(\frac{w_i}{v_i}\right)_x\right].
	\end{align*}
	We have then
	\begin{align}
		&\widetilde{\eta}_{i,t} +\tilde{\la}_i \widetilde{\eta}_{i,x}-(\mu_i\widetilde{\eta}_{i,x})_x=
		\widetilde{\eta}'_{i}(\frac{\tilde{\psi}_i}{v_i}-\frac{w_i}{v_i^2}\tilde{\phi}_i)+2\mu_i\widetilde{\eta}'_{i}\frac{v_{i,x}}{v_i}\left(\frac{w_i}{v_i}\right)_x-\mu_i\widetilde{\eta}''_i\left(\frac{w_i}{v_i}\right)_x^2.
		\label{energ2}
	\end{align}
	Integrating by parts, it yields
	\begin{align}
		|\int_{\R}\tilde{\la}_{i,x}\widetilde{\eta}_i\frac{ v_i^2}{2}dx|=|\int_{\R}(\tilde{\la}_i-\la_i^*)\left[\widetilde{\eta}'_{i}\left(\frac{w_i}{v_i}\right)_x\frac{v_i^2}{2}+\widetilde{\eta}_iv_iv_{i,x}\right]dx|
	\end{align}
	Now using the Lemma \ref{lemme6.8}, the fact that $|\tilde{\la}_i-\la_i^*|=O(1)\delta_0<<\delta_1<1$ and that $\|v_i(t,\cdot)\|_{L^\infty}=O(1)\delta^2_0, \|w_i(t,\cdot)\|_{L^\infty}=O(1)\delta_0^2$ (for $t\in[\hat{t},T^*]$)
	we deduce that $\widetilde{\eta}_i v_i=O(1) (v_{i,x}+\sum_{j\ne i}|v_j|)$. It yields for $\delta_0$ sufficiently small
	\begin{align}
		|\tilde{\la}_{i,x}\widetilde{\eta}_i\frac{ v_i^2}{2}dx|\leq \int_{\R}|w_{i,x}v_i-v_{i,x}w_i| dx+\frac{1}{2}\int_{\R}\mu_i \widetilde{\eta}_iv^2_{i,x}dx+\int_{\R}\sum_{j\ne i}|v_{i,x}v_j| dx.
		\label{energ3}
	\end{align}
	Combining \eqref{energ1}, \eqref{energ2} and \eqref{energ3} we deduce that
	\begin{align}
		\frac{1}{2}\int_{\R}\mu_i \widetilde{\eta}_iv^2_{i,x}dx \leq&-\frac{d}{dt}\int_{\R}\frac{1}{2}\widetilde{\eta}_i v_i^2 dx
		+\frac{1}{2}\int_{\R}|\widetilde{\eta}'_i|(|v_i\tilde{\psi}_i|+|w_i\tilde{\phi_i}|)dx\nonumber\\
		&+\int_{\R}|\mu_i\widetilde{\eta}'_iv_i v_{i,x}\left(\frac{w_i}{v_i}\right)_x| dx+\frac{1}{2}\int_{\R} \mu_i|\widetilde{\eta}''_iv_i^2\left(\frac{w_i}{v_i}\right)_x^2| dx\nonumber\\
		&+ \int_{\R}|w_{i,x}v_i-v_{i,x}w_i| dx+\int_{\R}\sum_{j\ne i}|v_{i,x}v_j| dx\nonumber\\
		&+2\int \mu_i|v_iv_{i,x}\widetilde{\eta}_{i,x} | dx+\int_{\R}|\widetilde{\eta}_i v_i\widetilde{\phi}_i |dx.\label{9.16}
	\end{align}
	From the definition of $\widetilde{\eta}$ we remark that $\widetilde{\eta}^\p_i\ne 0$ provided that $\abs{\frac{w_i}{v_i}}\leq \frac{3}{8}\delta_1$, it turns out that using the Lemma \ref{lemme6.8} we deduce that
	\begin{align}
		|\mu_i v_iv_{i,x}\widetilde{\eta}_{i,x} |&=|\mu_i v_iv_{i,x}\widetilde{\eta}'_{i}\left(\frac{w_i}{v_i}\right)_x|\\
		&\leq O(1)|w_{i,x}v_i-v_{i,x}w_i|+O(1)\sum_{j\ne i}(|v_j w_{i,x}|+|v_jv_{i,x}|).
		\label{9.17}
	\end{align}
	Combining \eqref{9.16}, \eqref{9.17}, \eqref{transvers1}, \eqref{La4}, \eqref{contradiction} and Lemma \ref{lemme6.6} we get for any $t\in[\hat{t},T^*]$
	\begin{align*}
		\int_{\hat{t}}^t\int_{\R} \widetilde{\eta}_iv^2_{i,x}dxds &\leq O(1)\int_{\R} \widetilde{\eta}_i v_i^2(\hat{t},x) dx
		+O(1)\int_{\hat{t}}^t \int_{\R}(|v_i\tilde{\psi}_i|+|w_i\tilde{\phi_i}|+|v_i\tilde{\phi}_i|)dx ds\\
		&+O(1)\int_{\hat{t}}^t \int_{\R}|w_{i,x}v_i-v_{i,x}w_i| dxds+O(1)\int_{\hat{t}}^t \int_{\R}\sum_{j\ne i}(|v_{i,x}v_j|+|v_jw_{i,x}|)dxds
		\nonumber\\
		&+\int_{\hat{t}}^t \int_{\{|\frac{w_i}{v_i}|\leq \delta_1 \}}|v_i\left(\frac{w_i}{v_i}\right)_x |^2d xds\\
		&\leq O(1)\delta_0^2.
	\end{align*}
	It achieves the proof due to the definition of $ \widetilde{\eta}_i$.
\end{proof}
\begin{lemma} We have
\begin{align}
&\int_{\hat{t}}^{T^*}\int_{\R}\mathbbm{1}_{\{|\frac{w_i}{v_i}|\geq \frac{\delta_1}{2}\}}w_{i,x}^2 ds dx=O(1)\delta_0^2.
\label{La61}
\end{align}
\end{lemma}
\begin{proof}
	We denote $\bar{\eta}_i=\bar{\eta}(w_i/v_i)$. Multiplying the equation \eqref{7.52} by $\bar{\eta}_i w_i$ and integrating by parts, we obtain
as in \cite{BB-vv-lim-ann-math}
\begin{align*} 
\int_{\R}\bar{\eta}_i w_i\widetilde{\psi}_i dx=&\frac{d}{dt}\int_{\R}\frac{1}{2}\bar{\eta}_i w_i^2 dx
+\int_{\R}\big(\tilde{\la}_{i,x}\bar{\eta}_i\frac{ w_i^2}{2}-\frac{1}{2}w_i^2\big(\bar{\eta}_{i,t} +\tilde{\la}_i \bar{\eta}_{i,x}-(\mu_i\bar{\eta}_{i,x})_x\big)+\mu_i \bar{\eta}_i w_{i,x}^2\\
&+2\mu_i w_iw_{i,x}\bar{\eta}_{i,x} \big)dx
\end{align*}
It implies that
\begin{align}
\int_{\R}\mu_i \bar{\eta}_i w_{i,x}^2dx=&
-\frac{d}{dt}\int_{\R}\frac{1}{2}\bar{\eta}_i w_i^2 dx
+\int_{\R}\big(\frac{1}{2}w_i^2\big(\bar{\eta}_{i,t} +\tilde{\la}_i \bar{\eta}_{i,x}-(\mu_i\bar{\eta}_{i,x})_x\big)-\tilde{\la}_{i,x}\bar{\eta}_i\frac{ w_i^2}{2}-2\mu_iw_iw_{i,x}\bar{\eta}_{i,x} \big)dx\nonumber\\
&+\int_{\R}\bar{\eta}_i w_i\widetilde{\psi}_i dx\label{lenerg1}
\end{align}
%Direct computations give:
%\begin{align*}
%&(\mu_i \eta_{i,x})_x=\mu_i\widetilde{\eta}''_i\left(\frac{w_i}{v_i}\right)_x^2+\widetilde{\eta}'_i\left[\frac{1}{v_i}\big((\mu_i w_{i,x})_x-\frac{w_i}{v_i}(\mu_i v_{i,x})_x\big)-2\mu_i\frac{v_{i,x}}{v_i}\left(\frac{w_i}{v_i}\right)_x\right].
%\end{align*}
%We have then
%\begin{align}
%&\widetilde{\eta}_{i,t} +\tilde{\la}_i \widetilde{\eta}_{i,x}-(\mu_i\widetilde{\eta}_{i,x})_x=
%\widetilde{\eta}'_{i}(\frac{\tilde{\psi}_i}{v_i}-\frac{w_i}{v_i^2}\tilde{\phi}_i)+2\mu_i\widetilde{\eta}'_{i}\frac{v_{i,x}}{v_i}\left(\frac{w_i}{v_i}\right)_x-\mu_i\widetilde{\eta}''_i\left(\frac{w_i}{v_i}\right)_x^2.
%\label{energ2}
%\end{align}
Integrating by parts, it yields
\begin{align}
|\int_{\R}\tilde{\la}_{i,x}\bar{\eta}_i\frac{ w_i^2}{2}dx|=|\int_{\R}(\tilde{\la}_i-\la_i^*)(\bar{\eta}'_{i}\left(\frac{w_i}{v_i}\right)_x\frac{w_i^2}{2}+\bar{\eta}_iw_iw_{i,x})dx|
\end{align}
Now using the Lemma \ref{lemme6.8} and the fact  that $\|v_i(t,\cdot)\|_{L^\infty}=O(1)\delta^2_0, \|w_i(t,\cdot)\|_{L^\infty}=O(1)\delta_0^2$ (for $t\in[\hat{t},T^*]$)
we deduce that 
\begin{align*}
&|\bar{\eta}_i w_i w_{i,x}|=O(1) \bar{\eta}_i (|v_{i,x}w_{i,x}|+\sum_{j\ne i}|v_j w_{i,x}|)=O(1)\bar{\eta}_i (v_{i,x}^2+w_{i,x}^2+\sum_{j\ne i}|v_j w_{i,x}|).
\end{align*}
It yields since  $|\tilde{\la}_i-\la_i^*|=O(1)\delta_0<<\delta_1<1$ that for $\delta_0$ sufficiently small and using the fact that $\bar{\eta}_i'\frac{w_i^2}{v_i^2}=O(1)$ we get
\begin{align}
|\int_\R\tilde{\la}_{i,x}\bar{\eta}_i\frac{ w_i^2}{2}dx|\leq \int_{\R}|w_{i,x}v_i-v_{i,x}w_i| dx+\frac{1}{2}\int_{\R}\mu_i (\bar{\eta}_iv^2_{i,x}+\bar{\eta}_iw^2_{i,x})dx+\int_{\R}\sum_{j\ne i}|w_{i,x}v_j| dx.
\label{lenerg3}
\end{align}
Combining \eqref{lenerg1}, \eqref{energ2} (that we adapt to the unknown $\bar{\eta}_i$) and \eqref{lenerg3} we deduce that using the fact that $|\bar{\eta}'_i\frac{w_i}{v_i}|\leq |\bar{\eta}'_i|\delta_1$
\begin{align}
&\frac{1}{2}\int_{\R}\mu_i \bar{\eta}_iw^2_{i,x}dx \leq-\frac{d}{dt}\int_{\R}\frac{1}{2}\bar{\eta}_i w_i^2 dx
+\frac{1}{2}\delta_1^2\int_{\R}|\bar{\eta}'_i|(|v_i\tilde{\psi}_i|+|w_i\tilde{\phi_i}|)dx+\delta_1^2\int_{\R}|\mu_i\bar{\eta}'_iv_i v_{i,x}\left(\frac{w_i}{v_i}\right)_x| dx\nonumber\\
&+\frac{1}{2}\delta_1^2\int_{\R} \mu_i|\bar{\eta}''_iv_i^2\left(\frac{w_i}{v_i}\right)_x^2| dx+\frac{1}{2}\int_{\R}\mu_i \bar{\eta}_iv^2_{i,x}dx
+ \int_{\R}|w_{i,x}v_i-v_{i,x}w_i| dx+\int_{\R}\sum_{j\ne i}|w_{i,x}v_j| dx\nonumber\\
&+2\int \mu_i|w_iw_{i,x}\bar{\eta}_{i,x} | dx
+\int_{\R}|\bar{\eta}_i v_i\widetilde{\phi}_i |dx.\label{l9.16}
\end{align}
We observe now since $\bar{\eta}'_i\ne 0$ if $\frac{3}{8}\delta_1<|\frac{w_i}{v_i}|<\frac{10 \delta_1}{24}$ and that $\widetilde{\eta}_i=1$ when $\bar{\eta}'_i\ne 0$ that we have chosen $\delta_1$ sufficiently small
$$|\bar{\eta}_i' |\,|v_i\left(\frac{w_i}{v_i}\right)_x|^2\geq\frac{1}{2} |\bar{\eta}_i' |w_{i,x}^2-\frac{1}{2}\widetilde{\eta}_i v_{i,x}^2.$$
It implies that
\begin{align}
|w_iw_{i,x}\bar{\eta}_{i,x}|&=O(1) |v_i\left(\frac{w_i}{v_i}\right)_x|^2\mathbbm{1}_{\{|\frac{w_i}{v_i}|\leq 3 \delta_1\}}+O(1)|\bar{\eta}_i'| w_{i,x}^2\nonumber\\
&=O(1) |v_i\left(\frac{w_i}{v_i}\right)_x|^2\mathbbm{1}_{\{|\frac{w_i}{v_i}|\leq 3 \delta_1\}}+O(1)\widetilde{\eta}_i v_{i,x}^2.
\label{9.60}
\end{align}
%From the definition of $\widetilde{\eta}$ we remark that $\widetilde{\eta}'_i\ne 0$ provided that $|\frac{w_i}{v_i}|\leq \frac{3}{8}\delta_1$, it turns out that using the Lemma \ref{lemme6.8} we deduce using the fact that $|\bar{eta}'_i\frac{w_i}{v_i}|\leq |\bar{eta}'_i|\delta_1$
%\begin{align}
%|\mu_i v_iv_{i,x}\widetilde{\eta}_{i,x} |&=|\mu_i v_iv_{i,x}\widetilde{\eta}'_{i}\left(\frac{w_i}{v_i}\right)_x\\
%&\leq O(1)|w_{i,x}v_i-v_{i,x}w_i|+O(1)\sum_{j\ne i}(|v_j w_{i,x}|+|v_jv_{i,x}|).
%\label{9.17}
%\end{align}
Combining \eqref{l9.16}, \eqref{9.60}, \eqref{transvers1}, \eqref{La4}, \eqref{La3}, \eqref{La6} and Lemmas \ref{lemme6.5}, \ref{lemme6.6} we get for any $t\in[\hat{t},T^*]$
\begin{align}
&\int_{\hat{t}}^t\int_{\R} \bar{\eta}_iw^2_{i,x}dxds \leq O(1)\delta_0^2.
\end{align}
It achieves the proof due to the definition of $ \bar{\eta}_i$.
\end{proof}
\subsection{Estimate of terms $\La_k^5=w_{k,xx}v_k-v_{k,xx}w_k$}\label{sec:Lambda-5}
Combining \eqref{defzihat} and direct computations gives
\begin{align}
	&\hat{z}_{k,x}v_k-v_{k,x}\hat{z}_k=\mu_{k}w_{k,xx}v_k+\mu_{k,x}w_{k,x}v_k-\tilde{\la}_{k,x}w_kv_k- (\tilde{\la}_k-\la_k^*)(w_{k,x}v_k-v_{k,x}w_k)\nonumber\\
	&+\sum_{j\ne k} v_k (\hat{a}_{k,j}(w_{j,x}-\frac{w_j}{v_j}v_{j,x}))_x-\mu_k v_{k,x}w_{k,x}-\sum_{j\ne k} v_{k,x} \hat{a}_{k,j}(w_{j,x}-\frac{w_j}{v_j}v_{j,x}).\label{restum1}
\end{align}
From  \eqref{eqn-v-i-x-1a-1} we get
\begin{align*}
	&\mu_k v_{k,x}w_{k,x}=(\tilde{\la}_k-\la_k^*)w_{k,x}v_k+w_k w_{k,x}-{\cal A}_k w_{k,x}\\
	&w_{k,x}=\mu_{k,x}v_{k,x}+\mu_kv_{k,xx}-(\tilde{\la}_k-\la_k^*)v_{k,x}-\tilde{\la}_{k,x}v_k+{\cal A}_{k,x}.
\end{align*}
It implies that
\begin{align}
	\mu_k v_{k,x}w_{k,x}=&(\tilde{\la}_k-\la_k^*)(w_{k,x}v_k-v_{k,x}w_k)+
	\mu_{k,x}v_{k,x}w_k+\mu_kv_{k,xx}w_k-\tilde{\la}_{k,x}v_kw_k\nonumber\\
	&+{\cal A}_{k,x}w_k-{\cal A}_k w_{k,x}.\label{restum2}
\end{align}
Finally we get from \eqref{restum1}, \eqref{restum2}
\begin{align*}
	&\hat{z}_{k,x}v_k-v_{k,x}\hat{z}_k\\
	&=\mu_{k}(w_{k,xx}v_k- v_{k,xx}w_k) -2(\tilde{\la}_k-\la_k^*)(w_{k,x}v_k-v_{k,x}w_k)+\mu_{k,x}(w_{k,x}v_k-v_{k,x}w_k)\\
	&+\sum_{j\ne k} v_k \hat{a}_{k,j, x}(w_{j,x}-\frac{w_j}{v_j}v_{j,x})+\sum_{j\ne k} v_k \hat{a}_{k,j}(w_{j,xx}-\frac{w_j}{v_j}v_{j,xx})-\sum_{j\ne k} v_k \hat{a}_{k,j}\left(\frac{w_j}{v_j}\right)_x v_{j,x}
	\\
	&-\sum_{j\ne k} v_{k,x} \hat{a}_{k,j}(w_{j,x}-\frac{w_j}{v_j}v_{j,x})+{\cal A}_k w_{k,x}-{\cal A}_{k,x}w_k.
\end{align*}
Now using \descref{B}{$(\mathcal{H}_B)$} , the fact that $\hat{a}_{k,j}=O(1)|\xi_j \bar{v}_j|$ and  $v_k=O(1)\delta_0^2$ from Lemma \ref{lemme6.5} , we deduce  that on $[\hat{t},T^*]$
\begin{align*}
	&\sum_k |w_{k,xx}v_k- v_{k,xx}w_k|
	=O(1)\biggl(
	\sum_k |\hat{z}_{k,x}v_k-v_{k,x}\hat{z}_k|
	+\sum_k |w_{k,x}v_k-v_{k,x}w_k|\\
	&+\sum_k \sum_{j\ne k} |v_k \hat{a}_{k,j, x}(w_{j,x}-\frac{w_j}{v_j}v_{j,x})|
	+\sum_{j\ne k} |v_k \hat{a}_{k,j}\left(\frac{w_j}{v_j}\right)_x v_{j,x}|
	+\sum_k \sum_{j\ne k}| v_{k,x} \hat{a}_{k,j}(w_{j,x}-\frac{w_j}{v_j}v_{j,x})|\\
	&+\sum_k|{\cal A}_k w_{k,x}-{\cal A}_{k,x}w_k|\biggl).
\end{align*}
From the Lemmas \ref{lemme6.5}, \ref{lemme6.6}, \ref{lemme9.6} we know that $ \hat{a}_{k,j, x}=O(1)\mathfrak{A}_j\rho_j^\e$ on $[\hat{t},T^*]$, we deduce then using Lemmas \ref{lemme6.5}, \ref{lemme6.6}, \eqref{transvers1}, \eqref{La8estim}%, \eqref{La4}, \eqref{La4bis-1},  \eqref{La4bis-2}, 
and the fact that  $\hat{a}_{k,j}=O(1)\mathfrak{A}_j \bar{v}_j$ that
\begin{align}
	&\sum_k\int_{\hat{t}}^{T^*}\int_{\R} |w_{k,xx}v_k- v_{k,xx}w_k|(s,x) ds dx=O(1)\delta_0^2
	+\sum_i\int_{\hat{t}}^{T^*}\int_{\R}|{\cal A}_i w_{i,x}-{\cal A}_{i,x}w_i|dx ds.\label{endtech4}
\end{align}
We recall now from \eqref{Ai} that
\begin{align*}
	&{\cal A}_i=\sum\limits_{j\neq i}(\mu_j v_{j,x}-(\tilde{\la}_j-\la_j^*)v_j-w_j)\frac{\mu_i}{\mu_j}[\psi_{ji}+v_j\xi_j\pa_{v_j}\bar{v}_j \psi_{ji,v}]\nonumber\\
	&+\sum\limits_{j\neq i}\sum\limits_{k\neq j}\xi_jv_{k,x}v_j\pa_{v_k}\bar{v}_j\mu_i \psi_{ji,v}+\sum\limits_{j\neq i}v_j\left(\frac{w_j}{v_j}\right)_x\mu_i [-\psi_{ji,\si}+\bar{v}_j\xi_j^\p\psi_{ji,v}]  \nonumber\\
	&+\sum\limits_{j\neq i}\mu_j^{-1}(\tilde{\la}_j-\si_j)v_j\xi_j \tilde{v}_j\mu_i\psi_{ji,v}+\sum\limits_{k}v_k^2(1-\xi_k\chi_k^\ep \eta_k)\langle l_i, \left[B(u)\tilde{r}_{k,u}\tilde{r}_k+\tilde{r}_k\cdot DB(u)\tilde{r}_k\right]\rangle \nonumber\\
	&+\sum\limits_{k\neq j}v_kv_j\langle l_i,\tilde{r}_k\cdot DB(u)\tilde{r}_j+B(u)\tilde{r}_{k,u}\tilde{r}_j\rangle.% \\[2mm]
	%&=\sum\limits_{j\neq i}(\mu_j v_{j,x}-(\tilde{\la}_j-\la_j^*)v_j-w_j)\frac{\mu_i}{\mu_j}[\psi_{ji}(u,\xi_j\bar{v}_j,\si_j)+v_j\xi_j\pa_{v_j}\bar{v}_j \psi_{ji,v}]\nonumber\\
	%&+\sum\limits_{j\neq i}\sum\limits_{k\neq j}\xi_jv_{k,x}v_j\pa_{v_k}\bar{v}_j\mu_i \psi_{ji,v}-\sum\limits_{j\neq i}\sum\limits_{k\neq j}\xi_jw_{k,x}v_j\pa_{w_k}\bar{v}_j\mu_i \psi_{ji,v}\nonumber\\
	%&+\sum\limits_{j\neq i}v_j\left(\frac{w_j}{v_j}\right)_x\mu_i [-\psi_{ji,\si}+\bar{v}_j\xi_j^\p\psi_{ji,v}] \nonumber\\
	%&+\sum\limits_{j\neq i}\mu_j^{-1}(\tilde{\la}_j-\si_j)v_j\xi_j \tilde{v}_j\mu_i\psi_{ji,v} \nonumber\\
	%&+\sum\limits_{k}v_k^2(1-\xi_k\eta_k)\langle l_i, \left[B(u)\tilde{r}_{k,u}\tilde{r}_k+\tilde{r}_k\cdot DB(u)\tilde{r}_k\right]\rangle \nonumber\\
	%&+\sum\limits_{k\neq j}v_kv_j\langle l_i,\tilde{r}_k\cdot DB(u)\tilde{r}_j+B(u)\tilde{r}_{k,u}\tilde{r}_j\rangle. 
\end{align*} 
%\begin{align*}
%	&{\cal A}_i=\sum\limits_{j\neq i}(\mu_j v_{j,x}-(\tilde{\la}_j-\la_j^*)v_j-w_j)\frac{\mu_i}{\mu_j}[\psi_{ji}(u,\xi_j\bar{v}_j,\si_j)+v_j\xi_j\pa_{v_j}\bar{v}_j \psi_{ji,v}]\nonumber\\
%	&+\sum\limits_{j\neq i}\sum\limits_{k\neq j}\xi_jv_{k,x}v_j\pa_{v_k}\bar{v}_j\mu_i \psi_{ji,v}-\sum\limits_{j\neq i}\sum\limits_{k\neq j}\xi_jw_{k,x}v_j\pa_{w_k}\bar{v}_j\mu_i \psi_{ji,v}\nonumber\\
%	&+\sum\limits_{j\neq i}v_j\left(\frac{w_j}{v_j}\right)_x\mu_i [-\psi_{ji,\si}+\bar{v}_j\xi_j^\p\psi_{ji,v}] +\sum\limits_{j\neq i}\mu_j^{-1}(\tilde{\la}_j-\si_j)v_j\xi_j \tilde{v}_j\mu_i\psi_{ji,v} \nonumber\\
%	&+\sum\limits_{k}v_k^2(1-\xi_k\eta_k)\langle l_i, \left[B(u)\tilde{r}_{k,u}\tilde{r}_k+\tilde{r}_k\cdot DB(u)\tilde{r}_k\right]\rangle \nonumber\\
%	&+\sum\limits_{k\neq j}v_kv_j\langle l_i,\tilde{r}_k\cdot DB(u)\tilde{r}_j+B(u)\tilde{r}_{k,u}\tilde{r}_j\rangle.
%\end{align*} 
Using Lemmas \ref{lemme6.5}, \ref{lemme6.6}, \eqref{6.54},  the fact that $\psi_{ji,\sig}=O(1)\mathfrak{A}_j \bar{v}_j$ and again \eqref{transvers1}, \eqref{La8estim} we obtain that
\begin{align}
&\int_{\hat{t}}^{T^*}\int_{\R}|{\cal A}_i w_{i,x}|dx ds\leq O(1)\delta_0^2+O(1)\sum_i\int_{\hat{t}}^{T^*}\int_{\R}|v_i^2(1-\chi_i^\e \xi_i\eta_i) w_{i,x}|dx ds\nonumber\\
&\leq O(1)\delta_0^2+O(1)\sum_i\int_{\hat{t}}^{T^*}\int_{\R}|v_i^2w_{i,x}|\big((1-\eta_i)\mathbbm{1}_{\{|\frac{w_i}{v_i}|\leq\frac{\delta_1}{2}\}}+\mathbbm{1}_{\{|\frac{w_i}{v_i}|\geq\frac{\delta_1}{2}\}} +\mathbbm{1}_{\{v_i^{2N}\leq 2\e \}}\big)dx ds\nonumber\\
&\leq O(1)\delta_0^2+\delta_0^2(T^*-\hat{t})\e^{\frac{1}{N}}=O(1)\delta_0^2.\label{endtech1}
\end{align}
We have in particular used the fact that $\e>0$ satisfies $\e^{\frac{1}{N}}(T-\hat{t})\leq 1$. Let us deal now with the term $\int_{\hat{t}}^{T^*}\int_{\R}|{\cal A}_{i,x}w_i|dx ds$. First we observe that using Lemmas \ref{estimimpo1}, \ref{lemme6.5}, \ref{lemme6.6},  \ref{lemme11.3}, \ref{lemme11.3a}, \ref{lemme11.4} and the fact that $\psi_{ji}, \psi_{ji,\sig}=O(1)\mathfrak{A}_j \bar{v}_j$, we have for $t\in[\hat{t},T^*]$
%\begin{align*}
%	&\sum_k\int_{\hat{t}}^{T^*}\int_{\R} |w_{k,xx}v_k- v_{k,xx}w_k|(s,x) ds dx=O(1)\delta_0^2
%	+\sum_i\int_{\hat{t}}^{T^*}\int_{\R}|{\cal A}_i w_{i,x}-{\cal A}_{i,x}w_i|dx ds.
%\end{align*}
%We recall now from \eqref{Ai} that
\begin{align}
	&{\cal A}_{i,x}=O(1)\sum_{j\ne i}(|v_j|+|v_{j,x}|+|w_j|+|w_{j,x}|)+2v_i v_{i,x}(1-\xi_i\chi_i^\ep \eta_i)\langle l_i, \left[B(u)\tilde{r}_{i,u}\tilde{r}_i+\tilde{r}_i\cdot DB(u)\tilde{r}_i\right]\rangle \nonumber\\
	&-v_i^2 \big[\xi_i'\left(\frac{w_i}{v_i}\right)_x \chi_i^\ep \eta_i+ 2N\frac{v_i^{2N-1}}{\e}v_{i,x}\xi_i(\chi')_i^\e \eta_i+\xi_i\chi_i^\e\pa_x(\eta_i)\big]\langle l_i, \left[B(u)\tilde{r}_{i,u}\tilde{r}_i+\tilde{r}_i\cdot DB(u)\tilde{r}_i\right]\rangle \nonumber\\
	&+O(1)v_i^2(1-\xi_i\chi_i^\ep \eta_i)+O(1)\sum\limits_{k\neq j}[v_kv_j+v_{k,x}v_j+v_{k}v_{j,x}].\label{endtech2}% \\[2mm]
	%&=\sum\limits_{j\neq i}(\mu_j v_{j,x}-(\tilde{\la}_j-\la_j^*)v_j-w_j)\frac{\mu_i}{\mu_j}[\psi_{ji}(u,\xi_j\bar{v}_j,\si_j)+v_j\xi_j\pa_{v_j}\bar{v}_j \psi_{ji,v}]\nonumber\\
	%&+\sum\limits_{j\neq i}\sum\limits_{k\neq j}\xi_jv_{k,x}v_j\pa_{v_k}\bar{v}_j\mu_i \psi_{ji,v}-\sum\limits_{j\neq i}\sum\limits_{k\neq j}\xi_jw_{k,x}v_j\pa_{w_k}\bar{v}_j\mu_i \psi_{ji,v}\nonumber\\
	%&+\sum\limits_{j\neq i}v_j\left(\frac{w_j}{v_j}\right)_x\mu_i [-\psi_{ji,\si}+\bar{v}_j\xi_j^\p\psi_{ji,v}] \nonumber\\
	%&+\sum\limits_{j\neq i}\mu_j^{-1}(\tilde{\la}_j-\si_j)v_j\xi_j \tilde{v}_j\mu_i\psi_{ji,v} \nonumber\\
	%&+\sum\limits_{k}v_k^2(1-\xi_k\eta_k)\langle l_i, \left[B(u)\tilde{r}_{k,u}\tilde{r}_k+\tilde{r}_k\cdot DB(u)\tilde{r}_k\right]\rangle \nonumber\\
	%&+\sum\limits_{k\neq j}v_kv_j\langle l_i,\tilde{r}_k\cdot DB(u)\tilde{r}_j+B(u)\tilde{r}_{k,u}\tilde{r}_j\rangle. 
\end{align} 
Using \eqref{endtech2}, \eqref{transvers1} and proceeding as in the term \eqref{endtech1}, we deduce that
\begin{align}
&\int_{\hat{t}}^{T^*}\int_{\R}|{\cal A}_{i,x} w_{i}|dx ds=O(1)\delta_0^2.\label{endtech3}
\end{align}
Finally combining \eqref{endtech1}, \eqref{endtech3} and \eqref{endtech4} we obtain
%\begin{align*}
%	\sum_k\int_{\hat{t}}^{T^*}\int_{\R} |w_{k,xx}v_k- v_{k,xx}w_k|\,ds dx&=O(1)\delta_0^2\\
%&+O(1)\sum_i\int_{\hat{t}}^{T^*}\int_{\R}\big(|v_i^2(1-\xi_i\eta_i) w_{i,x}|+|{\cal A}_{i,x}w_i|\big)dx ds.
%\end{align*}
%We observe now that using Corollary \ref{corollary-regularity} and Lemma \ref{lemme6.8}
%\begin{align*}
%	&|w_{i,x} |v_i^2(1-\xi_i\eta_i)|=O(1)\left(|w_{i,x}||v_i|^2\mathbbm{1}_{\left\{|\frac{w_i}{v_i}|\geq\frac{\delta_1}{2}\right\}}+\sum_{k\ne i}\mathbbm{1}_{\left\{|\frac{v_k^2+w_k^2}{|v_i|}|\geq\frac{3}{4}\right\}}|w_{i,x}||v_i|^2\right)\\
%	&=O(1)\big(|v_{i,x}|+\sum_{k\ne i}(|v_k|+|w_k|)+\La_{i}^4)^2||w_{i,x}|\mathbbm{1}_{\left\{|\frac{w_i}{v_i}|\geq\frac{\delta_1}{2}\right\}}+\sum_{k\ne i}\mathbbm{1}_{\left\{|\frac{v_k^2+w_k^2}{|v_i|}|\geq\frac{3}{4}\right\}}|w_{i,x}||v_i|(v_k^2+w_k^2)\big).
%\end{align*}
%Using now Corollary \ref{corollary-regularity} and \eqref{transvers1}, \eqref{La4} we deduce that
%\begin{align*}
%	\sum_i\int_{\hat{t}}^{T^*}\int_{\R}|v_i^2(1-\xi_i\eta_i) w_{i,x}|dxds=O(1)\delta_0^2.
%\end{align*}
%Proceeding in a similar way for dealing with the term $\sum_i\int_{\hat{t}}^{T^*}\int_{\R}|{\cal A}_{i,x}w_i| dx ds$ we obtain finally
\begin{align}
	&\sum_k\int_{\hat{t}}^{T^*}\int_{\R} |w_{k,xx}v_k- v_{k,xx}w_k|(s,x) ds dx=O(1)\delta_0^2.
	\label{La5}
\end{align}
\subsection{Transversal terms of high order}\label{sec:higher-order}
We wish now to estimate the transversal terms of high order, in other words we are going to prove that
\begin{align}
	& \int\limits_{\hat{t}}^{T^*}\int\limits_{\R}\sum_{j\ne i}\Big(|v_i|(|v_{j,xx}|+|w_{j,xx}|)+|v_{i,x}|(|v_{j,x}|+|v_{j,xx}|+|w_{j,x}|+|w_{j,xx}|)\nonumber\\
	&\hspace{2cm}+|w_i|(|v_{j,xx}|+|w_{j,xx}|)+|w_{i,x}|(|w_{j,x}|+|v_{j,xx}|+|w_{j,xx}|)\Big) dx ds=O(1)\delta_0^2.\label{transvers3a}
\end{align}
From \eqref{transvers2}, we know that
\begin{align*}
	& \int\limits_{\hat{t}}^{T^*}\int\limits_{\R}\sum_{j\ne i}|z_i v_{j,x}| dx ds=O(1)\delta_0^2,\,\int\limits_{\hat{t}}^{T^*}\int\limits_{\R}\sum_{j\ne i}|\hat{z}_i v_{j,x}| dx ds=O(1)\delta_0^2, .
\end{align*}
In particular from \eqref{defzi}, \eqref{defzihat} it implies that
\begin{align}
	& \int\limits_{\hat{t}}^{T^*}\int\limits_{\R}\sum_{j\ne i}|\big(\mu_iv_{i,x}-(\tilde{\lambda}_i-\la_i^*)v_i+\sum_{k\ne i}a_{ik}(w_{k,x}-\frac{w_k}{v_k}v_{k,x}) \big)v_{j,x}| dx ds=O(1)\delta_0^2,\nonumber\\
	& \int\limits_{\hat{t}}^{T^*}\int\limits_{\R}\sum_{j\ne i}|\big(\mu_iw_{i,x}-(\tilde{\lambda}_i-\la_i^*)w_i+\sum_{k\ne i}\hat{a}_{ik}(w_{k,x}-\frac{w_k}{v_k}v_{k,x}) \big)v_{j,x}| dx ds=O(1)\delta_0^2.
	\label{9.62}
\end{align}
Combining \eqref{9.62}, \descref{B}{$(\mathcal{H}_B)$}, \eqref{transvers1}, \eqref{La8estim}, Lemma \ref{lemme6.5} and the fact that $a_{ik}, \hat{a}_{ik}=O(1)\mathfrak{A}_k\bar{v}_k$ we deduce that
\begin{align}
	&\sum_{j\ne i} \int\limits_{\hat{t}}^{T^*}\int\limits_{\R}|v_{i,x}v_{j,x}| dx ds=O(1)\delta_0^2,\, \sum_{j\ne i} \int\limits_{\hat{t}}^{T^*}\int\limits_{\R}|w_{i,x}v_{j,x}| dx ds=O(1)\delta_0^2 .
	\label{9.59}
\end{align}
Similarly we have from \eqref{transvers2}
\begin{align}
	& \int\limits_{\hat{t}}^{T^*}\int\limits_{\R}\sum_{j\ne i}|z_{j,x} v_{i}| dx ds=O(1)\delta_0^2,  \,\int\limits_{\hat{t}}^{T^*}\int\limits_{\R}\sum_{j\ne i}|\hat{z}_{j,x} v_{i}| dx ds=O(1)\delta_0^2 .\nonumber
\end{align}
It implies from \eqref{defzi}, \eqref{defzihat} that
\begin{align}
	& \int\limits_{\hat{t}}^{T^*}\int\limits_{\R}\sum_{j\ne i}|v_i||\mu_j v_{j,xx}+\mu_{j,x}v_{j,x}-\tilde{\lambda}_{j,x}v_j-(\tilde{\la}_j-\la_i^*)v_{j,x}+\sum_{k\ne j}\big(a_{j,k}(w_{k,x}-\frac{w_k}{v_k}v_{k,x})\big)_x| dx ds,\nonumber\\
	& \int\limits_{\hat{t}}^{T^*}\int\limits_{\R}\sum_{j\ne i}|v_i||\mu_j w_{j,xx}+\mu_{j,x}w_{j,x}-\tilde{\lambda}_{j,x}w_j-(\tilde{\la}_j-\la_i^*)w_{j,x}+\sum_{k\ne j}\big(\hat{a}_{j,k}(w_{k,x}-\frac{w_k}{v_k}v_{k,x})\big)_x| dx ds\nonumber\\
	&=O(1)\delta_0^2.
\end{align}
Using Lemmas \ref{lemme6.5}, \ref{lemme6.6}, \ref{lemme9.6} we know that $a_{jk,x}, \hat{a}_{jk,x}=O(1)\mathfrak{A}_k\rho_k^\e$, applying now Lemmas  \ref{lemme6.5}, \ref{lemme6.6}, \ref{lemme9.6}, \ref{lemme11.3}, \eqref{transvers1},  \descref{B}{$(\mathcal{H}_B)$}, \eqref{La5}, \eqref{La8estim} and the fact that $a_{jk},\hat{a}_{jk}=O(1)\mathfrak{A}_k\bar{v}_k$  it yields
\begin{align}
	&\sum_{j\ne i} \int\limits_{\hat{t}}^{T^*}\int\limits_{\R}|v_{i}v_{j,xx}| dx ds,\,\sum_{j\ne i} \int\limits_{\hat{t}}^{T^*}\int\limits_{\R}|v_{i}w_{j,xx}| dx ds =O(1)\delta_0^2.
	\label{9.64}
\end{align}
Similarly from \eqref{transvers2.1}, \eqref{defzi}, \eqref{defzihat} we have
\begin{align}
	& \int\limits_{\hat{t}}^{T^*}\int\limits_{\R}\sum_{j\ne i}|z_i z_{j,x}| dx ds=
	\int\limits_{\hat{t}}^{T^*}\int\limits_{\R}\sum_{j\ne i}|\mu_i v_{i,x}-(\tilde{\lambda}_i-\la_i^*)v_i+\sum_{k\ne i}a_{i,k}(w_{k,x}-\frac{w_k}{v_k}v_{k,x})| \nonumber\\
	&\times|\mu_j v_{j,xx}+\mu_{j,x}v_{j,x}-\tilde{\lambda}_{j,x}v_j-(\tilde{\la}_j-\la_j^*)v_{j,x}+\sum_{k\ne j}\big(a_{j,k}(w_{k,x}-\frac{w_k}{v_k}v_{k,x})\big)_x| dx ds=O(1)\delta_0^2,\nonumber\\
	& \int\limits_{\hat{t}}^{T^*}\int\limits_{\R}\sum_{j\ne i}|\hat{z}_i \hat{z}_{j,x}| dx ds=
	\int\limits_{\hat{t}}^{T^*}\int\limits_{\R}\sum_{j\ne i}|\mu_i w_{i,x}-(\tilde{\lambda}_i-\la_i^*)w_i+\sum_{k\ne j}\hat{a}_{j,k}(w_{k,x}-\frac{w_k}{v_k}v_{k,x})| \nonumber\\
	&\times|\mu_j w_{j,xx}+\mu_{j,x}w_{j,x}-\tilde{\lambda}_{j,x}w_j-(\tilde{\la}_j-\la_j^*)w_{j,x}+\sum_{k\ne j}\big(\hat{a}_{j,k}(w_{k,x}-\frac{w_k}{v_k}v_{k,x})\big)_x| dx ds=O(1)\delta_0^2. \label{9.66}
\end{align}
Again combining Lemmas \ref{lemme6.5}, \ref{lemme6.6}, \eqref{9.66}, \eqref{9.59}, \eqref{9.64}, \eqref{transvers1}, \eqref{9.59}, \eqref{La5}, \eqref{La8estim}, \descref{B}{$(\mathcal{H}_B)$}, \eqref{ngtech5}  and the fact that $a_{ik}, \hat{a}_{ik},a_{jk}, \hat{a}_{jk}=O(1)\mathfrak{A}_k\rho_k^\e v_k$, $a_{jk,x}, \hat{a}_{jk,x}=O(1)\mathfrak{A}_k\rho_k^\e$ (see Lemma \ref{lemme9.6}) we get
\begin{equation}
	\sum_{j\ne i} \int\limits_{\hat{t}}^{T^*}\int\limits_{\R}|v_{i,x}v_{j,xx}| dx ds=O(1)\delta_0^2,\, \sum_{j\ne i} \int\limits_{\hat{t}}^{T^*}\int\limits_{\R}|w_{i,x}w_{j,xx}| dx ds=O(1)\delta_0^2. \label{9.66a}
\end{equation}
We proceed in a similar manner for estimating all the other terms in $\La_i^1$ (see \eqref{def:La_i-1}). We omit here the explicit proofs for the other terms and combining \eqref{transvers1}, \eqref{9.59}, \eqref{9.64}, \eqref{9.66a} we conclude that
\begin{align}
	\int_{\hat{t}}^{T^*}\int_{\R}|\La^1(s,x)| ds dx=O(1)\delta_0^2.\label{transvers3}
	\end{align}
\subsection{Estimate of terms involving third order derivatives $\La_i^2$}\label{section:Lambda-2}
We wish now to estimate
 the terms involving third order derivatives $w_{j,xxx}$ and $v_{j,xxx}$ which are of the following form
\begin{equation}
	\La_j^2=\sum_{k\ne j}(\abs{w_{j,xxx}}+\abs{v_{j,xxx}})(\abs{w_j}+\abs{v_j})\abs{v_k},
\end{equation}
\begin{lemma}
\label{lemmeLa2}
We have for any $t\in[\hat{t},T]$
\begin{equation}
\sum_j( |v_{j,xxx}|+|w_{j,xxx}|)=O(1)+O(1)\sum_k\delta_0^2\La_k^3+\widetilde{R}_\e,
\end{equation}
with
$ \int^T_{\hat{t}}\int_{\R}|\widetilde{R}_\e(s,x)| ds dx=O(1)\delta_0^2$.
%\label{estimwjxxx3}
\end{lemma}
\begin{proof}
 From the equation \eqref{eqn-w-i-1} and \eqref{deftilde} we have for any $t\in[\hat{t},T]$
\begin{align*}
	&w_{j,xxx}=%\left(\frac{1}{\mu_j}\right)_x\biggl[w_{j,t}+\tilde{\la}_{j,x}w_j+\tilde{\la}_j w_{j,x}-\mu_{j,x}w_{j,x}-\psi_j\\
	%&-\sum\limits_{m\neq j}(\mu_{m}-\mu_{j})\hat{b}_{jm}\left(w_{m,xx}-v_{m,x}\left(\frac{w_m}{v_m}\right)_x-\frac{w_m}{v_m}v_{m,xx}\right)\biggl]\\
	\frac{1}{\mu_j}\biggl[w_{j,tx}+\tilde{\la}_{j,xx}w_j+2\tilde{\la}_{j,x}w_{j,x}+\tilde{\la}_jw_{j,xx}-\mu_{j,xx}w_{j,x}-2\mu_{j,x}w_{j,xx}\biggl]\\
	&-\frac{1}{\mu_j}\psi_{j,x}-\frac{1}{\mu_j}\sum\limits_{m\neq j}(\mu_j-\mu_m)\hat{b}_{jm,x}\left(w_{m,xx}-v_{m,x}\left(\frac{w_m}{v_m}\right)_x-\frac{w_m}{v_m}v_{m,xx}\right)\\
	&-\frac{1}{\mu_j}\sum\limits_{m\neq j}(\mu_{j,x}-\mu_{m,x})\hat{b}_{jm}\left(w_{m,xx}-v_{m,x}\left(\frac{w_m}{v_m}\right)_x-\frac{w_m}{v_m}v_{m,xx}\right)\\
	&-\frac{1}{\mu_j}\sum\limits_{m\neq j}(\mu_j-\mu_m)\hat{b}_{jm}\left(w_{m,xxx}-2 v_{m,xx}\left(\frac{w_m}{v_m}\right)_{x}-\frac{w_m}{v_m}v_{m,xxx}\right)\\
	&+\frac{1}{\mu_j}\sum\limits_{m\neq j}(\mu_j-\mu_m)\hat{b}_{jm}v_{m,x}\left(\frac{w_m}{v_m}\right)_{xx}.
\end{align*}
From the Lemmas \ref{lemme6.5}, \ref{lemme6.6}, \ref{lemme11.6} we get for any $t\in[\hat{t},T]$
\begin{align*}
		\tilde{\la}_{j,xx}&=O(1)+O(1)\mathfrak{A}_j \rho_j^\e \biggl(|\left(\frac{w_j}{v_j}\right)_{xx}v_j|+|v_j\left(\frac{w_j}{v_j}\right)_{x}^2|+| \left(\frac{w_j}{v_j}\right)_{x}v_{j,x}|\\
		&+| \left(\frac{w_j}{v_j}\right)_{x}|\sum_{k\ne j}|v_kv_{k,x}|+\frac{v_{j,x}^2}{|v_j|}+|v_{j,x}|\sum_{k\ne j}|\frac{v_kv_{k,x}}{v_j}|%+\sum_{k\ne j}|v_kv_{k,xx}|+\sum_{k\ne j,l\ne j}|v_{k,x}v_{l,x}|
		\biggl).
	\end{align*}
%\begin{align*}
%&\tilde{\la}_{j,xx}=O(1)+O(1)\mathfrak{A}_j\biggl(|\left(\frac{w_j}{v_j}\right)_{xx}v_j|+|v_j\left(\frac{w_j}{v_j}\right)_{x}^2|+|\left(\frac{w_j}{v_j}\right)_{x}|+|v_{j,x}|\\
%		&+\sum_{k\ne j}|v_kv_{k,x}|+\frac{v_{j,x}^2}{|v_j|}+|v_{j,x}|\sum_{k\ne j}|\frac{v_kv_{k,x}}{v_j}|+|v_{j,xx}|+\sum_{k\ne j}|v_kv_{k,xx}|+\sum_{k\ne j,l\ne j}|v_{k,x}v_{l,x}|\biggl).
%		\end{align*}
%\begin{align*}
%\tilde{\la}_{j,xx}=O(1)+\mathbbm{1}_{\{|\frac{w_j}{v_j}|\leq 3\delta_1\}}\biggl(|\left(\frac{w_j}{v_j}\right)_{xx}v_j|+|v_j\left(\frac{w_j}{v_j}\right)_{x}^2|+|\left(\frac{w_j}{v_j}\right)_{x}|
%		+\frac{v_{j,x}^2}{|v_j|}&\\
%		+|v_{j,x}|\sum_{k\ne j}(|\frac{v_kv_{k,x}}{v_j}|+|\frac{w_kw_{k,x}}{v_j}|)\biggl).&
%	\end{align*}
It implies in particular again from Lemmas \ref{lemme6.5}, \ref{lemme6.6} that for any $t\in[\hat{t},T]$
\begin{align}
&w_j\tilde{\la}_{j,xx}=O(1).\label{10.51}
\end{align}
Using \eqref{10.51}, \eqref{ngtech5},  the Lemmas \ref{lemme6.5}, \ref{lemme6.6}, \ref{lemme11.3} and the fact that $\hat{b}_{jm}=O(1)\mathfrak{A}_m\bar{v}_m$ we deduce that
\begin{align}
	&w_{j,xxx}=O(1)+O(1)|\psi_{j,x}|-\frac{1}{\mu_j}\sum\limits_{m\neq j}(\mu_m-\mu_j)\hat{b}_{jm,x}\left(w_{m,xx}-v_{m,x}\left(\frac{w_m}{v_m}\right)_x-\frac{w_m}{v_m}v_{m,xx}\right)\nonumber\\
	&+O(1)\sum\limits_{m\neq j}(|v_mw_{m,xxx}|+|v_m v_{m,xxx}|)+O(1)\sum\limits_{m\neq j}\mathfrak{A}_m\rho_m^\e ||v_{m,x}^2\left(\frac{w_m}{v_m}\right)_{x}|\nonumber\\
	&=O(1)+O(1)|\psi_{j,x}|-\frac{1}{\mu_j}\sum\limits_{m\neq j}(\mu_m-\mu_j)\hat{b}_{jm,x}\left(w_{m,xx}-v_{m,x}\left(\frac{w_m}{v_m}\right)_x-\frac{w_m}{v_m}v_{m,xx}\right)\nonumber\\
	&+O(1)\sum\limits_{m\neq j}(|v_mw_{m,xxx}|+|v_m v_{m,xxx}|).\label{wjxxx}
\end{align}
Applying now Lemmas \ref{lemme6.5}, \ref{lemme6.6}, \ref{lemme9.6} and \eqref{ngtech5} we get
\begin{align*}
&-\frac{1}{\mu_j}\sum\limits_{m\neq j}(\mu_m-\mu_j)\hat{b}_{jm,x}\left(w_{m,xx}-v_{m,x}\left(\frac{w_m}{v_m}\right)_x-\frac{w_m}{v_m}v_{m,xx}\right)\\
&=O(1)+O(1)\sum_{m\ne j}\mathfrak{A}_m\rho_m^\e \big(|v_m v_{m,x}\left(\frac{w_m}{v_m}\right)_x^2|+|v_{m,x}^2\left(\frac{w_m}{v_m}\right)_x|+|v_{m,x}\left(\frac{w_m}{v_m}\right)_x|\sum_{k\ne m}|v_kv_{k,x}|\big)\\
&=O(1)+O(1)\sum_k\delta_0^2\La_k^3.
\end{align*}
From \eqref{wjxxx}, we have
\begin{align}
&w_{j,xxx}=O(1)+O(1)|\psi_{j,x}|+O(1)\sum_k\delta_0^2\La_k^3+O(1)\sum\limits_{m\neq j}(|v_mw_{m,xxx}|+|v_m v_{m,xxx}|).
\label{estimwjxxx}
	\end{align}
We recall that in the section \ref{sectionphii}, we are going to prove that
\begin{align}
|\psi_{j,x}|=O(1)\sum_l(\La_k^1+\La_k^2+\delta_0^2\La_k^3+\La_k^4+\La_k^5+\La_k^6+\La_k^{6,1})+\widetilde{R}_{\e},
\label{estimwjxxx1}
\end{align}
with $\int^T_{\hat{t}}\int_{\R}|\widetilde{R}_\e(s,x)| ds dx|=O(1)\delta_0^2$.
Combining \eqref{estimwjxxx},  \eqref{estimwjxxx1} and Lemmas \ref{lemme6.5}, \ref{lemme6.6} we obtain that
\begin{align}
&w_{j,xxx}=O(1)+O(1)\sum_k(\La_k^2+\delta_0^2\La_k^3)+O(1)\sum\limits_{m\neq j}(|v_mw_{m,xxx}|+|v_m v_{m,xxx}|)+\widetilde{R}_{\e},.
\label{estimwjxxx2}
\end{align}
In a similar way we can show that for any $j\in\{1,\cdot,n\}$
\begin{align}
&v_{j,xxx}=O(1)+O(1)\sum_k(\La_k^2+\delta_0^2\La_k^3)+O(1)\sum\limits_{m\neq j}(|v_mw_{m,xxx}|+|v_m v_{m,xxx}|)+\widetilde{R}_{\e},.
\label{estimwjxxx3}
\end{align}
From \eqref{estimwjxxx2}, \eqref{estimwjxxx3} and using the fact that 
$\|v_m(t,\cdot)\|_{L^\infty}, \|w_m(t,\cdot)\|_{L^\infty}=O(1)\delta_0^2$ for any $t\in[\hat{t},T]$, $m\in\{1,\cdot,n\}$
(see the Lemma \ref{lemme6.5}), we deduce by a bootstrap argument provided that $\delta_0$ is sufficiently small that for any $t\in[\hat{t},T]$ we have
\begin{equation}
\sum_j( |v_{j,xxx}|+|w_{j,xxx}|)=O(1)+O(1)\sum_k\delta_0^2\La_k^3+\widetilde{R}_{\e},
%\label{estimwjxxx3}
\end{equation}
with $\int^T_{\hat{t}}\int_{\R}|\widetilde{R}_\e(s,x)| ds dx|=O(1)\delta_0^2$.
\end{proof}
The previous Lemma implies in particular that
\begin{align}
&\int^{T^*}_{\hat{t}}\int_{\R}\La_j^2(s,x) ds dx=\int^{T^*}_{\hat{t}}\int_{\R}\sum_{k\ne j}(\abs{w_{j,xxx}}+\abs{v_{j,xxx}})(\abs{w_j}+\abs{v_j})\abs{v_k}(s,x) ds dx,\nonumber\\
&=O(1)\int^{T^*}_{\hat{t}}\sum_{k\ne j} \int_{\R}(\abs{w_j}+\abs{v_j})\abs{v_k}(s,x) ds dx+O(1)\delta_0^2\sum_k \int^{T^*}_{\hat{t}}\int_{\R}\La_k^3(s,x) ds dx\nonumber\\
&+\int^{T^*}_{\hat{t}}\int_{\R}|\widetilde{R}_\e(s,x)| ds dx.
\end{align}
From Lemma \ref{lemmeLa2}, \eqref{transvers1} and \eqref{La3} we deduce that for any $j\in\{1,\cdots,n\}$
\begin{align}
&\int^{T^*}_{\hat{t}}\int_{\R}\sum_{k\ne j}(\abs{w_{j,xxx}}+\abs{v_{j,xxx}})(\abs{w_j}+\abs{v_j})\abs{v_k}(s,x) ds dx=O(1)\delta_0^2,
\end{align}
which allows to estimate the terms $\La_k^2$.
\subsection{Proof of Theorem \ref{theorem:BV-estimate}}
From the analysis in section \ref{section:transversal-wave-interaction}, \ref{sec:shorten}, \ref{sec:Lambda-5}, \ref{sec:energy} \ref{sec:higher-order}, \ref{section:Lambda-2} we have established \eqref{estimate:remainder-1}, that is,
\begin{equation}\label{pf-thm-1-est-1}
	\sum_{i=1}^n\int_{\hat{t}}^{T^*}\int_{\R}\left(|\widetilde{\phi}_i(s,x)|+|\widetilde{\psi}_i(s,x)|+|\Phi_i(s,x)|+|\Psi_i(s,x)|\right) dx ds=O(1)\delta_0^{2}.
\end{equation}
From the definition of $T^*$  in \eqref{def:T*} we deduce that $T^*=T$ and we have seen at the beginning of the section \ref{globalBVn} that in this case we have necessary $T=+\infty$. It implies that the solution is global and $\norm{u_x(t)}_{L^1(\R)}\leq\de_0$ holds for all $t\in [0,+\infty)$. This proves the estimate \eqref{thm-1:BV-bound} in Theorem \ref{theorem:BV-estimate}.

Next we show the estimate \eqref{L1-cont} by a similar argument as in \cite{BB-vv-lim-ann-math,HJ-temple-class,HJ-triangular}. Without loss of generality we assume $t>s>0$. By using \eqref{eqn-thm-1} we have
\begin{align*}
	\norm{u(t)-u(s)}_{L^1(\R)}&\leq \int\limits_{s}^{t}\int\limits_{\R}\abs{\pa_tu(\tau,x)}\,dxd\tau\\
	&\leq \int\limits_{s}^{t}\int\limits_{\R} \left(|A(u)u_x(\tau,x)|+|B(u)u_{xx}(\tau,x)|+|DB(u)||u_x(\tau,x)|^2\right)\,dxd\tau\\
	&\leq \widetilde{C}_1\int\limits_{s}^{t}\left(\norm{u_x(\tau)}_{L^1}+\norm{u_{xx}(\tau)}_{L^1}+\norm{u_x(\tau)}_{L^1} \norm{u_{xx}(\tau)}_{L^1}\right)\,dxd\tau,
\end{align*}
for some constant $\widetilde{C}_1>0$ independent of $t,s$. Now, if $s>\hat{t}$ then by Corollary \ref{coro2.2} we get $\norm{u_x(\tau)}_{L^1},\norm{u_{xx}(\tau)}_{L^1}=O(1)\de_0$ and we get the estimate \eqref{L1-cont}. When $t<\hat{t}$, by using Proposition \ref{prop:parabolic} and \ref{prop:local-existence} we have
\begin{align*}
	\norm{u(t)-u(s)}_{L^1(\R)}&\leq \widetilde{C}_2\int\limits_{s}^{t}\left(\de_0+\frac{\de_0}{\sqrt{\tau}}\right)\,dxd\tau\\
	&\leq  \widetilde{C}_2\de_0\left(|t-s|+|\sqrt{t}-\sqrt{s}|\right),
\end{align*}
for some constant $\widetilde{C}_2>0$. Now, if we consider $0<s<\hat{t}<t$, by a similar argument we get
\begin{equation*}
		\norm{u(t)-u(s)}_{L^1(\R)}\leq \widetilde{C}_3\de_0\left(|t-\hat{t}|+|\hat{t}-s|+|\sqrt{\hat{t}}-\sqrt{s}|\right)\leq  \widetilde{C}_3\de_0\left(|t-s|+|\sqrt{t}-\sqrt{s}|\right),
\end{equation*}
for some constant $\widetilde{C}_3>0$. Therefore, the time continuity estimate \eqref{L1-cont} follows. This completes the proof of Theorem \ref{theorem:BV-estimate}.

\subsection{Proof of Corollary \ref{corollary:vv-limit}}\label{subsection:pf-of-coro}

In this subsection, we show the vanishing viscosity limit for conservative case, that is when $A(u)=Df(u)$ for some smooth function $f:\R^n\rr\R^n$. The proof follows from a similar argument as in \cite{BB-vv-lim-ann-math,HJ-temple-class,HJ-triangular}. For sake of completeness we present a brief proof of Corollary \ref{corollary:vv-limit}.

\begin{proof}[Proof of Corollary \ref{corollary:vv-limit}:]  First, observe that $u^\eps(t,x)=u(t/\eps,x/\eps)$ where $u$ solves the following problem with fix viscosity but scaled initial data, 
	\begin{equation}\label{viscous-eqn-conservative}
		u_t+(f(u))_x=(B(u)u_x)_x\mbox{ and }u(0,x)=\bar{u}(\eps x).
	\end{equation}
	From the scaling of $\bar{u}$ we note that
	\begin{equation*}
		TV(\bar{u}(\eps \cdot ))= TV(\bar{u}(\cdot)),\quad		\norm{\bar{u}(\eps\cdot )}_{L^1}=\frac{1}{\eps}	\norm{\bar{u}(\cdot )}_{L^1}.
	\end{equation*}
	By using Theorem \ref{theorem:BV-estimate} we get
	\begin{align*}
		TV(u(t))&\leq L_1 TV(\bar{u}),\\
		\norm{u(t)-u(s)}_{L^1}&\leq L_2\left(\abs{t-s}+\abs{\sqrt{t}-\sqrt{s}}\right).
	\end{align*}
	Observe that
	\begin{equation}
		TV(u^\eps(t, \cdot ))= TV(u(t/\eps,\cdot))\mbox{ and }\norm{u^\eps(t)-u^\eps(s)}_{L^1}=\ep \norm{u(t/\eps)-u(s/\eps)}_{L^1}.
	\end{equation}Therefore, it yields
	\begin{align}
		TV(u^\eps(t))&\leq L_1 TV(\bar{u}),\label{est-1}\\
		\norm{u^\eps(t)-u^\eps(s)}_{L^1}
		&\leq L_2\left(\abs{t-s}+\sqrt{\e}\abs{\sqrt{t}-\sqrt{s}}\right).     \label{est-3}
	\end{align}
	The convergence of $u^\eps$ as $\eps\rr0$ follows from a standard argument with an application of Helly's theorem and the $L^1$ continuity \eqref{est-3}. Indeed, due to the uniform TV estimate \eqref{est-1} by using Helly's theorem we can pass to a the limit (up to a subsequence) for a countable dense set $\{t_n\}$ and then applying $L^1$ continuity we can define the limit function at all time $t>0$. We set
	\begin{equation}
		L^1_{loc}-\lim\limits_{k\rr\f}u^{\eps_k}(t,\cdot)=u^\f(t,\cdot).
	\end{equation}
	In particular up to a subsequence (denoting by same notation) $\{u^{\eps_k}\}_{k}$ converges almost everywhere to $u^\f$, it implies in particular that for any function $\varphi$ in $C^\infty_0(\R^+\times\R)$ by  applying dominated convergence, we get
	\begin{align}
		\int_{\R^+}\int_\R\varphi_x(s,x)f(u^{\eps_k}(s,x)) \,ds dx\overset{k\rightarrow\infty}{\longrightarrow}\int_{\R^+}\int_\R\varphi_x(s,x)f(u^\f(s,x))\, ds dx.
	\end{align}
	Similarly we have for any $T>0$
	\begin{equation}
		\eps\int^T_0 \int_\R |B(u^\eps)u^\eps_x| dx ds\leq \sup_{\{|v-u^*|\leq 2\delta_0\}}|B(v)|L_1\delta_0\eps T,
	\end{equation}
	then $\eps B(u^\e)u^\eps_x$ converges to zero in $L^1_{loc, t,x}$. It implies that $u^\f$ is a global weak solution of the following hyperbolic system of conservation laws. 
	\begin{equation}
		u^\f_t+f(u^\f)_x=0.
	\end{equation}
	By a similar argument as in \cite{BB-vv-lim-ann-math,Bressan-book}, we can show that there exists a measure zero set $\mathcal{N}\subset\R_+$ such that for $\tau\notin \mathcal{N}$  and $\xi\in\R$
	\begin{equation*}
		u_-=\lim\limits_{x\rr\xi-}u^\f(\tau,x)\mbox{ and }	u_+=\lim\limits_{x\rr\xi+}u^\f(\tau,x).
	\end{equation*}
	Then $u$ satisfies
	\begin{align*}
		\lim\limits_{r\rr0+}\frac{1}{r^2}\int\limits_{-r}^{r}\int\limits_{-\kappa r}^{\kappa r}\abs{u^\f(\tau+t,\xi+x)-U(t,x)}dxdt&=0,\\
		\lim\limits_{r\rr0+}\frac{1}{r}\int\limits_{-\kappa r}^{\kappa r}\abs{u^\f(\tau+t,\xi+x)-U(r,x)}dx&=0,
	\end{align*}
	for every $\kappa>0$ where $U$ is defined as follows
	\begin{equation}\label{def:U-RIemann}
		U(t,x)=\left\{\begin{array}{rl}
			u_-&\mbox{ for }x<\la t,\\
			u_+&\mbox{ for }x>\la t,
		\end{array}\right.
	\end{equation}
	for some $\la\in\R$. Furthermore, we can show
	\begin{equation*}
		\la (u_+-u_-)=f(u_+)-f(u_-).
	\end{equation*}
	Due to the condition $AB=BA$, we note that the matrix $\mathcal{S}(z)=-z^2B(u^*)-zi A(u^*)$ satisfies $Re(\zeta_k(z))\leq -c_2z^2$ for all $z\in\R$ and for some $c_2>0$ where $\{\zeta_j(z)\}_{j=1}^n$ are eigenvalues of $\mathcal{S}(z)$ and $i^2=-1$. Then by \cite[Theorem 2.1, page 235]{Majda-Pego} it follows that $B$ is {\em strictly stable} in a neighbourhood of $u^*$. Hence, from \cite[Corollary 2, page 233]{Majda-Pego}, we know that $U(t,x)$ (defined in \eqref{def:U-RIemann}) is {\em Liu admissible} solution if and only if the function $U$ is vanishing viscosity limit of \eqref{eqn-coro-vis-1} with $\bar{u}$ as Riemann data defined as follows
	\begin{equation}\label{def:U}
		\bar{u}(x)=\left\{\begin{array}{rl}
			u_-&\mbox{ for }x<0,\\
			u_+&\mbox{ for }x>0,
		\end{array}\right.
	\end{equation}
	with $|u_\pm-u^*|<\de_1$ for some small enough $\de_1>0$. Therefore, $u^\f$ satisfies {\em Liu condition} at discontinuities. We can invoke \cite{BDL} to show uniqueness of the limit $u^\f$. Therefore, we established the following result on vanishing viscosity limit for conservative case with possibly non-constant $B(u)$.
\end{proof}

\section{Some basic estimates}
Since we are going to use it many times let us mention here the calculation of $\left(\frac{w_i}{v_i}\right)_{xx}$ and $\left(\frac{w_i}{v_i}\right)_{xxx}$ when $v_i\ne 0$.
\begin{align}
	\left(\frac{w_i}{v_i}\right)_{xx}&=\frac{1}{v_i}\left(w_{i,xx}-\frac{w_i}{v_i}v_{i,xx}\right)-\frac{2v_{i,x}}{v_i}\left(\frac{w_i}{v_i}\right)_{x}, \label{identity:wi-vi-xx}\\
	\left(\frac{w_i}{v_i}\right)_{xxx}&=\frac{1}{v_i}\left(w_{i,xxx}-\frac{w_i}{v_i}v_{i,xxx}\right)-\frac{3v_{i,x}}{v_i^2}\left(w_{i,xx}-\frac{w_i}{v_i}v_{i,xx}\right)-\frac{3v_{i,xx}}{v_i}\left(\frac{w_i}{v_i}\right)_x\nonumber\\
		&+\frac{6v_{i,x}^2}{v_i^2}\left(\frac{w_i}{v_i}\right)_x.
	%-\frac{v_{i,x}}{v_i^2}\left(w_{i,xx}-\frac{w_i}{v_i}v_{i,xx}\right)
	%+\frac{1}{v_i}\left(w_{i,xxx}-\frac{w_i}{v_i}v_{i,xxx}\right)-\frac{1}{v_i}\left(\frac{w_i}{v_i}\right)_x v_{i,xx} \nonumber\\
	%&-\frac{2v_{i,x}}{v_i}\left(\frac{w_i}{v_i}\right)_{xx}-2(\frac{v_{i,xx}}{v_i}-\frac{v_{i,x}^2}{v_i^2})\left(\frac{w_i}{v_i}\right)_{x}.
	\label{identity:wi-vi-xxx}
\end{align}
In this section we are under the same assumption as the section \ref{section6.1}, in particular all the estimates concerned time $t\in[\hat{t},T]$ such that the decomposition  \eqref{eqn-v-i}- \eqref{eqn-w-i} is satisfied.

\subsection{First order estimate}
\begin{lemma}	\label{lemme9.2}
	Let $\rho_i^\ep,\De_{ij}^\ep$ and $\Box^\ep_{ijk}$ as in section \ref{sec:Gradient-Decomp}. The following estimate holds true.
	\begin{align}
		\bar{v}^\e_{i,x}&=O(1)\left(\abs{v_{i,x}}\rho^\ep_i+\sum_{j\ne i}\abs{v_j v_{j,x}}\De^\ep_{ij}\right),\,	\bar{v}^\e_{i,t}=O(1)\left(\abs{v_{i,t}}\rho^\ep_i+\sum_{j\ne i}\abs{v_j v_{j,t}}\De^\ep_{ij}\right),\\
		\bar{v}^\e_{i,xx}&=O(1)\left(\frac{\abs{v_{i,x}}^2}{\abs{v_i}}+\abs{v_{i,xx}}\right)\rho_i^\e+O(1)\left(\sum\limits_{j\neq i}(\abs{v_{j,xx}}+\abs{\frac{v_{i,x}v_{j,x}}{v_i}})|v_j|\De^\ep_{ij}+\sum\limits_{j,k\neq i}\abs{v_{j,x}v_{k,x}}\Box_{ijk}^\ep\right).
	\end{align}
\end{lemma}
\begin{proof} A direct computation gives using the Lemma \ref{estimimpo1}
	\begin{align*}
		\bar{v}^\e_{i,x}&=\pa_{v_i}\bar{v}_i^\e v_{i,x}+\sum_{i\ne j}\pa_{v_j}\bar{v}_i^\e v_{j,x}\\
		&=O(1)\left(\abs{v_{i,x}}\rho^\ep_i+\sum_{j\ne i}\abs{v_j v_{j,x}}\De^\ep_{ij}\right).
	\end{align*}
	Similarly we have
	\begin{align*}
		\bar{v}^\e_{i,xx}&=\pa_{v_iv_i}\bar{v}_i ^\e (v_{i,x})^2+\pa_{v_i}\bar{v}_i^\e v_{i,xx}+2\sum_{i\ne j}(\pa_{v_i v_j}\bar{v}_i^\e) v_{j,x}v_{i,x}\\
		&+\sum_{j,k\ne i}(\pa_{v_k v_j}\bar{v}_i^\e) v_{j,x}v_{k,x}+\sum_{i\ne j}(\pa_{v_j}\bar{v}_i^\e) v_{j,xx}\\
		&=O(1)\left(\frac{\abs{v_{i,x}}^2}{\abs{v_i}}+\abs{v_{i,xx}}\right)\rho_i^\e+O(1)\left(\sum\limits_{j\neq i}(\abs{v_jv_{j,xx}}+\abs{\frac{v_j}{v_i}v_{i,x}v_{j,x}})\De^\ep_{ij}+\sum\limits_{j,k\neq i}\abs{v_{j,x}v_{k,x}}\Box_{ijk}^\ep\right).
	\end{align*}
	We have also from Lemma \ref{estimimpo1}
\begin{equation*}
		\bar{v}_{i,t}^\e=\pa_{v_i}\bar{v}_i ^\e v_{i,t}+\sum_{i\ne j}\pa_{v_j}\bar{v}_i^\e v_{j,t}=O(1)\left(\abs{v_{i,t}}\rho^\ep_i+\sum_{j\ne i}\abs{v_j v_{j,t}}\De^\ep_{ij}\right).
	\end{equation*}
This completes the proof. 
\end{proof}

\begin{lemma}\label{lemme11.3}
	For each $1\leq k\leq n$, we have
	\begin{align*}
		&\pa_x(\tilde{r}_{k}),\pa_x(\tilde{r}_{k,u})=O(1)\left(\abs{v_{k,x}}\rho_k^\e\mathfrak{A}_k+\sum_{j}\abs{v_j}\right)+O(1)%\left(\sum\limits_{j\neq k}\abs{v_{j,x}v_j}\rho^\e_k\mathfrak{A}_k+
		|v_k\left(\frac{w_k}{v_k}\right)_x|\rho^\e_k \mathfrak{A}_k%\right)
		,\\
		&\pa_x(\tilde{r}_{k,v}), \pa_x(\tilde{r}_{k,vv}),	\pa_x(\tilde{r}_{k,v\si})=O(1)\left(\abs{v_{k,x}}\rho_k^\e\mathfrak{A}_k+\sum_{j}\abs{v_j}%+\sum\limits_{j\neq k} \abs{v_{j,x}v_j}\rho_k^\e\mathfrak{A}_k
		\right)+O(1)\abs{\left(\frac{w_k}{v_k}\right)_x}\mathfrak{A}_k,\\%\label{estimate-derivative-r-k-x}\\
		%&	\pa_x(\tilde{r}_{i,v}),\pa_x(\tilde{r}_{i,vv}),	\pa_x(\tilde{r}_{i,v\si})=O(1)\left(|v_i|+\sum\limits_{k\ne i}^{n}\abs{v_k}+\abs{\left(\frac{w_i}{v_i}\right)_x}\mathfrak{A}_i+ |v_{i,x}|\rho_i\mathfrak{A}_i\right),\\%\label{estimate:r-i-v}\\
	&	\pa_x(\tilde{r}_{k,\si}),	\pa_x(\tilde{r}_{k,\si\si})=O(1)\left(\sum_j |v_j v_k|+\sum\limits_{j\ne k}^{n}\abs{v_jv_{j,x}}+\abs{v_k\left(\frac{w_k}{v_k}\right)_x}+  |v_{k,x}|\right)\rho_k^\e\mathfrak{A}_k.\label{estimate:r-i-si}\\
&\tilde{r}_{k,t}, \pa_t(\tilde{r}_{k,\sig})=O(1)\left(|v_{k,t}|\rho_k^\e\mathfrak{A}_k+ \sum\limits_{j}(|w_j|+|v_j|)+\sum\limits_{j\neq k}|v_{j,t}v_j|\rho_k^\e\mathfrak{A}_k+|v_k\left(\frac{w_k}{v_k}\right)_t| \rho_k^\e\mathfrak{A}_k \right),% \label{derivative-r-k-t}
		\\
		&\pa_t(\tilde{r}_{k,v})=O(1)\left(|v_{k,t}|\rho_k^\e\mathfrak{A}_k+ \sum\limits_{j}(|w_j|+|v_j|)+\sum\limits_{j\neq k}|v_{j,t}v_j|\rho_k^\e\mathfrak{A}_k+|\left(\frac{w_k}{v_k}\right)_t\mathfrak{A}_k| \right), %\label{derivative-r-k-t1}
		\\
		&\tilde{\la}_{j,x},  \pa_x(\tilde{\la}_{j,v})%\pa_x(\tilde{\la}_{j,u}\cdot \tilde{r}_k),  
		=O(1)\left(\sum_k |v_k|+(|v_j\left(\frac{w_j}{v_j}\right)_x|+|v_{j,x}|) \mathfrak{A}_j\rho_j^\e\right),\\%+O(1)\sum_{k\ne j}|v_k v_{k,x}|\rho_j^\e\mathfrak{A}_j
		&\pa_x(\tilde{\la}_{j,\sig})=\mathfrak{A}_j\rho_j^\e |v_j|\left(\sum_k |v_k|+|\left(\frac{w_j}{v_j}\right)_x|+|v_{j,x}|\right).
		\end{align*}
\end{lemma}
\begin{proof}
		Taking derivative of $\tilde{r}_k$ we get recalling that $\theta_k^\p \tilde{r}_{k,\si}= \tilde{r}_{k,\si}$ (indeed from \eqref{observationimp1} we have  $\tilde{r}_{k,\si}=O(1)\xi_k \bar{v}_k$ and $\mbox{supp}\xi$ is included in $\{x,\theta'(x)=1\}$)
	\begin{align}
		\tilde{r}_{k,x}&=\sum\limits_{j}v_j\tilde{r}_{k,u}\tilde{r}_j+v_{k,x}\xi_k\pa_{v_k}\bar{v}_k\tilde{r}_{k,v}+\sum\limits_{j\neq k}v_{j,x}\xi_k\pa_{v_j}\bar{v}_k\tilde{r}_{k,v}
		\\
		&+\xi_k^\p\left(\frac{w_k}{v_k}\right)_x\bar{v}_k\tilde{r}_{k,v}-\left(\frac{w_k}{v_k}\right)_x\tilde{r}_{k,\si}. \label{derivative-r-k-x}
	\end{align}
	We obtain the desired estimate from Lemmas \ref{estimimpo1}, \ref{lemme6.5}, \ref{lemme6.6} and since $\tilde{r}_{k,\sig}=O(1)\xi_k\bar{v}_k$.
	We can proceed similarly for $\tilde{r}_{k,u}$, $\tilde{r}_{k,v}$, $ \tilde{r}_{k,vv}$, $\tilde{r}_{k,v\si}$ provided that we observe from \eqref{observationimp1} that $\tilde{r}_{k,u\sig}=O(1)\xi_k\bar{v}_k$.
	Concerning the terms $\tilde{r}_{k,\sig}, \tilde{r}_{k,\sig \sig}$ we simply exploit in addition the fact that $\tilde{r}_{k,\sig u}, \tilde{r}_{k,\sig \sig u}=O(1)\xi_k\bar{v}_k$ from \eqref{observationimp1}. \\% Therefore using the Lemma \ref{estimimpo1}, the estimate \eqref{estimate-derivative-r-k-x} follows.
	 In a similar way, we have
	\begin{align}
		\tilde{r}_{k,t}&=\sum\limits_{j}(w_j-\la_j^*v_j)\tilde{r}_{k,u}\tilde{r}_j+v_{k,t}\xi_k\pa_{v_k}\bar{v}_k\tilde{r}_{k,v}+\sum\limits_{j\neq k}v_{j,t}\xi_k\pa_{v_j}\bar{v}_k\tilde{r}_{k,v}
		\nonumber\\
		&+\xi_k^\p\left(\frac{w_k}{v_k}\right)_t\bar{v}_k\tilde{r}_{k,v}-\theta'_k\left(\frac{w_k}{v_k}\right)_t \tilde{r}_{k,\si}. \label{derivative-r-k-t}
	\end{align}
Using the fact $\tilde{r}_{k,\si}=O(1)\xi_k\bar{v}_k$ and the Lemma \ref{estimimpo1} we deduce the desired estimate. We proceed similarly for $\pa_t(\tilde{r}_{k,\sig})$ and $\pa_t(\tilde{r}_{k,v})$.
	Finally we have
	\begin{align}
		&\tilde{\la}_{j,x}=\sum_k v_k\tilde{r}_k\cdot \tilde{\la}_{j,u}+\left(\frac{w_j}{v_j}\right)_x(\xi'_j\tilde{\la}_{j,v}\bar{v}_j-\tilde{\la}_{j,\sig}\theta'_j)+\tilde{\la}_{j,v}\xi_j\left[\pa_{v_j}\bar{v}_jv_{j,x}+\sum_{k\ne j}\pa_{v_k}\bar{v}_j v_{k,x}\right].\label{lajx}\\
		&\pa_x(\tilde{\la}_{j,\sig})=\sum_k v_k\tilde{r}_k\cdot \tilde{\la}_{j,\sig u}+\left(\frac{w_j}{v_j}\right)_x(\xi'_j\tilde{\la}_{j,\sig v}\bar{v}_j-\tilde{\la}_{j,\sig \sig}\theta'_j)+\tilde{\la}_{j,\sig v}\xi_j\left[\pa_{v_j}\bar{v}_jv_{j,x}+\sum_{k\ne j}\pa_{v_k}\bar{v}_j v_{k,x}\right].\\
		&\pa_x(\tilde{\la}_{j,v})=\sum_k v_k\tilde{r}_k\cdot \tilde{\la}_{j,v u}+\left(\frac{w_j}{v_j}\right)_x(\xi'_j\tilde{\la}_{j,v v}\bar{v}_j-\tilde{\la}_{j,v \sig}\theta'_j)+\tilde{\la}_{j,v v}\xi_j\left[\pa_{v_j}\bar{v}_jv_{j,x}+\sum_{k\ne j}\pa_{v_k}\bar{v}_j v_{k,x}\right].
	\end{align}
	We conclude using \eqref{tildelamb} and Lemmas \ref{estimimpo1}, \ref{lemme6.5}.

\end{proof}

\begin{lemma}\label{lemme11.3a}
For $j\neq i$, the following estimates are true.
	\begin{align*}
		\pa_x\psi_{ji}, \pa_x\psi_{ji,\sig}&=O(1)\biggl[\sum_{k\ne j}(|v_k v_j|+|v_k v_{k,x}|)+\abs{v_j \left(\frac{w_j}{v_j}\right)_x}+|v_{j,x}|+ |v_j|^2\biggl]\mathfrak{A}_j\rho_j^\e ,\\
		\pa_x\psi_{ji, v}&=O(1)\sum_{k}|v_k |+O(1)\biggl[\abs{\left(\frac{w_j}{v_j}\right)_x}%+\sum_{k\ne j}|v_k v_{k,x}|\rho_j^\e 
		+|v_{j,x}|\rho_j^\e\biggl] \mathfrak{A}_j,\\
		\pa_t\psi_{ji,\sig}, \pa_t\psi_{ji}&=O(1)\biggl[\sum_k(|v_jw_k|+|v_kv_j|)+|w_{j,xx}-v_{j,xx}\frac{w_j}{v_j}|+|v_j\left(\frac{w_j}{v_j}\right)_x|\\
		&\hspace{2cm}+|\tilde{\psi}_j|+|\tilde{\phi}_j|+|v_{j,xx}|+|v_{j,x}|+|v_j|+\sum_{k\ne j}|v_{k,t}v_k|\biggl] \mathfrak{A}_j\rho_j^\e\\
		\pa_t\psi_{ji,v}&=O(1)\sum_k(|v_k|+|w_k|)+O(1)\biggl[|\frac{1}{v_j}||w_{j,xx}-v_{j,xx}\frac{w_j}{v_j}|\\
		&+|\left(\frac{w_j}{v_j}\right)_x|+|\frac{1}{v_j}|(|\tilde{\psi}_j|+|\tilde{\phi}_j|)+\rho_j^\e |v_{j,xx}|+\rho_j^\e|v_{j,x}|+\rho_j^\e|v_j|+\rho_j^\e|\tilde{\phi}_j|\\
		&+\rho_j^\e \sum_{k\ne j}|v_{k,t}v_k|\biggl] \mathfrak{A}_j.
	\end{align*}

\end{lemma}

\begin{proof}[Proof of Lemma \ref{lemme11.3a}:]
		In order to simplify the notation, we skip in the sequel all the exponent $\e>0$ for the terms of the type $\bar{v}_i^{\e}$, $\pa_{v_j}\bar{v}_i^\e$. We have from  \eqref{eqn-v-i}
	\begin{align}
		\pa_x\psi_{ji}&=\sum_k\psi_{ji,u}\tilde{r}_k v_k+\left(\frac{w_j}{v_j}\right)_x(\xi'_j\bar{v}_j\psi_{ji,v}-\theta^\p_j \psi_{ji,\sig})+ \psi_{ji,v}\xi_j (\pa_{v_j}\bar{v}_j)v_{j,x}\nonumber\\
		&+\sum_{k\ne j} \psi_{ji,v}\xi_j (\pa_{v_k}\bar{v}_j)v_{k,x},\nonumber\\
		\pa_x\psi_{ji,v}&=\sum_k\psi_{ji,uv}\tilde{r}_k v_k+\left(\frac{w_j}{v_j}\right)_x(\xi'_j\bar{v}_j\psi_{ji,vv}-\theta'_j \psi_{ji,v\sig})+ \psi_{ji,vv}\xi_j (\pa_{v_j}\bar{v}_j)v_{j,x}\nonumber\\
		&+\sum_{k\ne j} \psi_{ji,vv}\xi_j (\pa_{v_k}\bar{v}_j)v_{k,x},\nonumber\\
		\pa_x\psi_{ji,\si}&=\sum_k\psi_{ji,u\si}\tilde{r}_k v_k+\left(\frac{w_j}{v_j}\right)_x(\xi^\p_j\bar{v}_j \psi_{ji,v\si}-\theta'_j \psi_{ji,\si\si})+ \psi_{ji,v\sig}\xi_j (\pa_{v_j}\bar{v}_j)v_{j,x}\nonumber\\
		&+\sum_{k\ne j} \psi_{ji,v\si}\xi_j (\pa_{v_k}\bar{v}_j)v_{k,x}.\label{eqn-derivative-psi-x-sig}
	\end{align}
	We deduce the results by using the estimate \eqref{5.19} and Lemma \ref{estimimpo1}, \ref{lemme6.5}, \ref{lemme6.6}.
	Next we have
	\begin{align*}
	&\pa_t\psi_{ji}=\sum_k(w_k-\la_k^* v_k)\psi_{ji, u}\tilde{r}_k+\left(\frac{w_j}{v_j}\right)_t(\psi_{ji, v}\xi'_j\bar{v}_j-\theta'_j\psi_{ji, \sig})+\psi_{ji, v}\xi_j\pa_t( \bar{v}_j)\\
		&\pa_t\psi_{ji,\sig}=\sum_k(w_k-\la_k^* v_k)\psi_{ji,\sig u}\tilde{r}_k+\left(\frac{w_j}{v_j}\right)_t(\psi_{ji,\sig v}\xi'_j\bar{v}_j-\theta'_j\psi_{ji,\sig \sig})+\psi_{ji,\sig v}\xi_j\pa_t( \bar{v}_j)\\
		&\pa_t\psi_{ji,v}=\sum_k(w_k-\la_k^* v_k)\psi_{ji,v u}\tilde{r}_k+\left(\frac{w_j}{v_j}\right)_t(\psi_{ji,v v}\xi'_j\bar{v}_j-\theta'_j\psi_{ji,v \sig})+\psi_{ji,v v}\xi_j\pa_t( \bar{v}_j).
	\end{align*}
	From \eqref{eqn-v-i-1}-\eqref{eqn-w-i-1} we obtain
	\begin{equation}
		v_j \left(\frac{w_j}{v_j}\right)_t=\mu_j(w_{j,xx}-v_{j,xx}\frac{w_j}{v_j})+(\mu_{j,x}-\tilde{\la}_j)v_j\left(\frac{w_j}{v_j}\right)_x+\tilde{\psi}_j-\frac{w_j}{v_j}\tilde{\phi}_j,
		\label{formuleimpo}
	\end{equation}
	we deduce again using \eqref{5.19},  \eqref{eqn-v-i-1} and Lemmas \ref{estimimpo1},  \ref{lemme9.2}, \ref{lemme11.3}
	\begin{align*}
		\pa_t \psi_{ji},\pa_t\psi_{ji,\sig}&=O(1)\bigg[\sum_k(|v_jw_k|+|v_kv_j|) +|w_{j,xx}-v_{j,xx}\frac{w_j}{v_j}|+|v_j\left(\frac{w_j}{v_j}\right)_x|\\
		&+|\tilde{\psi}_j|+|\tilde{\phi}_j|+|v_{j,xx}|+|v_{j,x}|+|v_j|+\sum_{k\ne j}|v_{k,t}v_k|\bigg] \rho_j \mathfrak{A}_j\\
		|\pa_t\psi_{ji,v}|&=O(1)\sum_k(|v_k|+|w_k|)+O(1)\biggl[|\frac{1}{v_j}||w_{j,xx}-v_{j,xx}\frac{w_j}{v_j}|\\
		&+|\left(\frac{w_j}{v_j}\right)_x|+|\frac{1}{v_j}|(|\tilde{\psi}_j|+|\tilde{\phi}_j|)+\rho_j |v_{j,xx}|+\rho_j|v_{j,x}|+\rho_j |v_j|\\
		&+\rho_j |\tilde{\phi}_j|+\rho_j \sum_{k\ne j}|v_{k,t}v_k|\biggl]  \mathfrak{A}_j.
	\end{align*}
	This completes the proof of Lemma \ref{lemme11.3a}.
\end{proof}
\begin{lemma}
\label{lemme11.4}
For $1\leq i,j\leq n$ the following estimates hold true.
\begin{align}
	\pa_x(\tilde{r}_{i,u}\tilde{r}_j)
	&=O(1)\biggl(\sum\limits_{k=1}^{n}\abs{v_k}+|v_{i,x}|\rho_i^\e\mathfrak{A}_i+|v_{j,x}|\rho_j^\e\mathfrak{A}_j
	+\abs{v_i\left(\frac{w_i}{v_i}\right)_x}\rho_i^\e\mathfrak{A}_i+\abs{v_j\left(\frac{w_j}{v_j}\right)_x}\rho_j^\e\mathfrak{A}_j%+\sum_{k\ne j}|v_kv_{k,x}|\rho_j\mathfrak{A}_j+\sum_{k\ne i}|v_kv_{k,x}|\rho_j\mathfrak{A}_j
	\biggl),\label{11.10u}\\
    \pa_x(\tilde{r}_{i,u\si}\tilde{r}_j)
	&=O(1)\biggl(\sum\limits_{k=1}^{n}\abs{v_iv_k}+|v_{i,x}|+|v_{j,x}v_i|\rho_j^\e\mathfrak{A}_j
	+\abs{v_i\left(\frac{w_i}{v_i}\right)_x}+\abs{v_iv_j\left(\frac{w_j}{v_j}\right)_x}\rho_j^\e\mathfrak{A}_j\nonumber\\
	&\hspace{1cm}%\sum_{k\ne j}|v_iv_kv_{k,x}|\rho_j^\e\mathfrak{A}_j
	+\sum_{k\ne i}|v_kv_{k,x}|\biggl)\rho_i^\e\mathfrak{A}_i,\label{11.10-u-sig-1}\\
	\pa_x(\tilde{r}_{i,uv}\tilde{r}_j)
	&=O(1)\biggl(\sum\limits_{k=1}^{n}\abs{v_k}+|v_{i,x}|\rho_i^\e\mathfrak{A}_i+|v_{j,x}|\rho_j^\e\mathfrak{A}_j
	+\abs{\left(\frac{w_i}{v_i}\right)_x}\mathfrak{A}_i+\abs{v_j\left(\frac{w_j}{v_j}\right)_x}\rho_j^\e\mathfrak{A}_j
	\biggl),\label{11.11-uv1}\\
		\pa_x(\tilde{r}_{i,u}\tilde{r}_{j,v})
	&=O(1)\biggl(\sum\limits_{k=1}^{n}\abs{v_k}+|v_{i,x}|\rho_i^\e\mathfrak{A}_i+|v_{j,x}|\rho_j^\e\mathfrak{A}_j
	+\abs{v_i\left(\frac{w_i}{v_i}\right)_x}\rho_i^\e\mathfrak{A}_i+\abs{\left(\frac{w_j}{v_j}\right)_x}\mathfrak{A}_j%+\sum_{k\ne j}|v_kv_{k,x}|\rho_j\mathfrak{A}_j+\sum_{k\ne i}|v_kv_{k,x}|\rho_j\mathfrak{A}_j
	\biggl),\label{11.11-uv2}\\
	\pa_x(\tilde{r}_{i,u}\tilde{r}_{j,\si})
	&=O(1)\biggl(\sum\limits_{k=1}^{n}\abs{v_jv_k}+|v_jv_{i,x}|\rho_i^\e\mathfrak{A}_i+|v_{j,x}|
	+\abs{v_iv_j\left(\frac{w_i}{v_i}\right)_x}\rho_i^\e\mathfrak{A}_i+\abs{v_j\left(\frac{w_j}{v_j}\right)_x}\nonumber\\
	&\hspace{1cm}+\sum_{k\ne j}|v_kv_{k,x}|%+\sum_{k\ne i}|v_jv_kv_{k,x}|\rho_i\mathfrak{A}_i
	\biggl)\rho_j\mathfrak{A}_j^\e.\label{11.10-u-sig-2}\\
	\pa_t(\tilde{r}_{i,u}\tilde{r}_j)
		&=O(1)\biggl(\sum\limits_{k=1}^{n}(\abs{v_k}+|w_k|)+|v_{i,t}| \rho_i^\e\mathfrak{A}_i+|v_{j,t}|
		\rho_j^\e\mathfrak{A}_j+\abs{v_i\left(\frac{w_i}{v_i}\right)_t}\rho_i^\e\mathfrak{A}_i\nonumber\\
		&+\abs{v_j\left(\frac{w_j}{v_j}\right)_t}\rho_j^\e\mathfrak{A}_j
		%\sum_{k\ne j}|w_kw_{k,t}|\al_j+\sum_{k\ne i}|w_kw_{k,t}|\al_i
		\biggl)
		.\label{riut}
\end{align}
\end{lemma}
\begin{proof}
	Taking derivative of $\tilde{r}_{i,u}\tilde{r}_j$, we have
	\begin{align}
		\pa_x(\tilde{r}_{i,u}\tilde{r}_j)&=\sum\limits_{k=1}^{n}v_k[\tilde{r}_{i,uu}:\tilde{r}_{j}\otimes\tilde{r}_k+\tilde{r}_{i,u}\tilde{r}_{j,u}\tilde{r}_k]+(\xi_i \bar{v}_{i})_x\tilde{r}_{i,uv}\tilde{r}_j+(\xi_j \bar{v}_{j})_x\tilde{r}_{i,u}\tilde{r}_{j,v}\nonumber\\
		&-\left(\frac{w_i}{v_i}\right)_x\theta_i^\p\tilde{r}_{i,u\si}\tilde{r}_j-\left(\frac{w_j}{v_j}\right)_x\theta_j^\p\tilde{r}_{i,u}\tilde{r}_{j,\si}\nonumber\\
		&=\sum\limits_{k=1}^{n}v_k[\tilde{r}_{i,uu}:\tilde{r}_{j}\otimes\tilde{r}_k+\tilde{r}_{i,u}\tilde{r}_{j,u}\tilde{r}_k]+\left(\frac{w_i}{v_i}\right)_x\left[\xi^\p_i\bar{v}_i\tilde{r}_{i,uv}\tilde{r}_j-\theta_i^\p\tilde{r}_{i,u\si}\tilde{r}_j\right]\nonumber\\
		&+\left(\frac{w_j}{v_j}\right)_x\left[\xi^\p_j\bar{v}_j\tilde{r}_{i,u}\tilde{r}_{j,v}-\theta_j^\p\tilde{r}_{i,u}\tilde{r}_{j,\si}\right]+v_{i,x}\xi_i\pa_{v_i}\bar{v}_i\tilde{r}_{i,uv}\tilde{r}_j+v_{j,x}\xi_j\pa_{v_j}\bar{v}_j\tilde{r}_{i,u}\tilde{r}_{j,v}\nonumber\\
		&+\sum\limits_{k\neq i}\left( v_{k,x}\pa_{v_k}\bar{v}_i\right)\xi_i\tilde{r}_{i,uv}\tilde{r}_j+\sum\limits_{k\neq j}\left( v_{k,x}\pa_{v_k}\bar{v}_j\right)\xi_j\tilde{r}_{i,u}\tilde{r}_{j,v}.\label{rix}
	\end{align}
	By using the fact that $\psi_{j,i}=O(1)\xi_j\bar{v}_j$ and the Lemmas \ref{lemme6.5}, \ref{lemme6.6} we deduce the estimates \eqref{11.10u}. 
	We proceed similarly for the term $\tilde{r}_{i,u\si}\tilde{r}_j$ where we use the fact that $\tilde{r}_{i,u\si}, \tilde{r}_{i,uu\si}, \tilde{r}_{i,u\si\si}=O(1)\xi_i\bar{v}_i$
	. We treat the other terms in the same manner.\\
	 Let us deal now with the term $\pa_t(\tilde{r}_{i,u}\tilde{r}_j)$, we have in particular 
		\begin{align}
		\pa_t(\tilde{r}_{i,u}\tilde{r}_j)&=\sum\limits_{k=1}^{n}(w_k-\la_k^*v_k)[\tilde{r}_{i,uu}:\tilde{r}_j\otimes\tilde{r}_k+\tilde{r}_{i,u}\tilde{r}_{j,u}\tilde{r}_k]+\left(\frac{w_i}{v_i}\right)_t\left[\xi^\p_i\bar{v}_i\tilde{r}_{i,uv}\tilde{r}_j-\theta_i^\p\tilde{r}_{i,u\si}\tilde{r}_j\right]\nonumber\\
		&+\left(\frac{w_j}{v_j}\right)_t\left[\xi^\p_j\bar{v}_j\tilde{r}_{i,u}\tilde{r}_{j,v}-\theta_j^\p\tilde{r}_{i,u}\tilde{r}_{j,\si}\right]+v_{i,t}\xi_i\pa_{v_i}\bar{v}_i\tilde{r}_{i,uv}\tilde{r}_j+v_{j,t}\xi_j\pa_{v_j}\bar{v}_j\tilde{r}_{i,u}\tilde{r}_{j,v}\nonumber\\
		&+\sum\limits_{k\neq i}v_{k,t}\pa_{v_k}\bar{v}_i\xi_i\tilde{r}_{i,uv}\tilde{r}_j+\sum\limits_{k\neq j}v_{k,t}\pa_{v_k}\bar{v}_j\xi_j\tilde{r}_{i,u}\tilde{r}_{j,v}.\label{rit}
	\end{align}
Using again the fact that $\tilde{r}_{i,u \sig}=O(1)\xi_i\bar{v}_ij$ and $\tilde{r}_{j,\sig}=O(1)\xi_j\bar{v}_j$ we deduce again the result from Lemmas \ref{estimimpo1}, \ref{lemme6.5}.
\end{proof}
\subsection{Second order estimate}
\begin{lemma}\label{lemme11.5}
	For $1\leq i,j\leq n$ with $i\neq j$ the following estimates are true.
	\begin{align}
		\pa_{xx}\psi_{ji},	\pa_{xx}\psi_{ji,\si}&=O(1)\left(\sum\limits_{k}(\abs{v_kv_j}+|v_{k,x}v_j|+\abs{v_kv_{j,x}}+ \abs{v_k w_{j,x}})
		\right)\rho_j^\e\mathfrak{A}_j \nonumber\\
		&+O(1) \left(\abs{v_{j,xx}}+\sum\limits_{k\neq j}\abs{v_kv_{k,xx}}+\sum\limits_{l\neq j} \frac{\abs{v_l v_{l,x}v_{j,x}}}{\abs{v_j}}+\frac{v_{j,x}^2}{|v_j|} \right)\rho_j^{\e}\mathfrak{A}_j \nonumber\\
		&+ O(1)\left(\sum\limits_{k}\sum\limits_{l\neq j}\abs{v_kv_{l,x}}+\sum\limits_{k\ne j,l\ne j}\abs{v_{k,x}v_{l,x}}+\abs{v_j\left(\frac{w_j}{v_j}\right)_{xx}}+\abs{v_j\left(\frac{w_j}{v_j}\right)_{x}^2}\right)\rho_j^\e\mathfrak{A}_j \nonumber \\
		&+O(1)|\left(\frac{w_j}{v_j}\right)_x|\left(\sum\limits_{k\neq j}\abs{v_{k,x}v_k}+\abs{v_{j,x}}+\abs{v_j} \right)\rho_j^\e\mathfrak{A}_j,\label{estimate-psi-xx-sigma}\\
			\pa_{xx}\psi_{ji,v}
			&=O(1)\sum_k (|v_k|+|v_{k,x}|)+O(1) |\left(\frac{w_j}{v_j}\right)_x|\biggl(\sum_k|v_k|+|v_{j,x}|\rho_j^\e\biggl)\mathfrak{A}_j\nonumber\\
			&+O(1)\biggl[ \left(\frac{w_j}{v_j}\right)^2_x+\abs{\left(\frac{w_j}{v_j}\right)_{xx}}
			\biggl]\mathfrak{A}_j\nonumber\\
			&+O(1)\biggl(\frac{v_{j,x}^2}{|v_j|}+|v_{j,xx}|+\sum_{k\ne j} |\frac{v_k}{v_j}| |v_{k,x}v_{j,x}|+\sum_{k\ne j}|v_{k}v_{k,xx}| \biggl)\rho_j^\e \mathfrak{A}_j.
		\label{estimate-psi-xx-v}
	\end{align}
\end{lemma}
\begin{proof}
		We proceed to prove the estimate \eqref{estimate-psi-xx-sigma} for $\pa_{xx}\psi_{ij,\si}$. From \eqref{eqn-derivative-psi-x-sig} we deduce that %similarly for $\psi_{ji,\sig}$. We deduce that 
		\begin{align}
			\pa_{xx}\psi_{ji,\sig}&=\sum_{k,l}\psi_{ji,uu\si}(\tilde{r}_k\otimes\tilde{r}_l) v_lv_k+\sum_k\psi_{ji,u\si}[\tilde{r}_{k,x} v_k+\tilde{r}_kv_{k,x}]+\sum_k\psi_{ji,uv\si}\tilde{r}_k \xi_j(\pa_{v_j}\bar{v}_j)v_{j,x}v_k\nonumber\\
			&+\sum_{k}\sum\limits_{l\neq j}\psi_{ji,uv\si}\tilde{r}_k \xi_j(v_{l,x}\pa_{v_l}\bar{v}_j)v_k
			+\sum_k(\xi_j^\p\bar{v}_j\psi_{ji,uv\si}-\theta_j^\p\psi_{ji,u\si\si})\tilde{r}_k v_k\left(\frac{w_j}{v_j}\right)_x\nonumber\\
			&+\left(\frac{w_j}{v_j}\right)_{xx}(\xi^\p_j\bar{v}_j \psi_{ji,v\si}-\theta^\p_j \psi_{ji,\sig\si})+\left(\frac{w_j}{v_j}\right)^2_x(\xi^{\p\p}_j\bar{v}_j \psi_{ji,v\si}-\theta^{\p\p}_j \psi_{ji,\sig\si})\nonumber\\
			&+\left(\frac{w_j}{v_j}\right)_x\xi^\p_j\pa_x(\bar{v}_j) \psi_{ji,v\si}+\left(\frac{w_j}{v_j}\right)_x\sum\limits_{k}v_k(\xi^\p_j\bar{v}_j \psi_{ji,uv\si}-\theta^\p_j \psi_{ji,u\sig\si})\tilde{r}_k\nonumber\\
			&+\left(\frac{w_j}{v_j}\right)^2_x(\xi^\p_j\bar{v}_j \psi_{ji,vv\si}-\theta^\p_j \psi_{ji,v\sig\si})\xi^\p_j\bar{v}_j+\left(\frac{w_j}{v_j}\right)_x\sum\limits_{k}(\xi^\p_j\bar{v}_j \psi_{ji,vv\si}-\theta^\p_j \psi_{ji,v\sig\si})\xi_j\pa_{v_k}\bar{v}_jv_{k,x}\nonumber\\
			&-\left(\frac{w_j}{v_j}\right)^2_x(\xi^\p_j\bar{v}_j \psi_{ji,v\si\si}-\theta'_j \psi_{ji,\si\sig\si})\theta_j^\p
			+\sum_k\sum\limits_{l} \psi_{ji,uv\si}\tilde{r}_l\xi_j (\pa_{v_k}\bar{v}_j)v_lv_{k,x}\nonumber\\
			&+\sum_k \left(\frac{w_j}{v_j}\right)_x\psi_{ji,vv\si}\xi_j^\p{\bar{v}_j}\xi_j (\pa_{v_k}\bar{v}_j)v_{k,x}
			+\sum_k\sum\limits_{l} \psi_{ji,vv\si}\xi_j^2(\pa_{v_l}\bar{v}_j) (\pa_{v_k}\bar{v}_j)v_{l,x}v_{k,x}\nonumber\\
			&-\sum_k \psi_{ji,v\si\si}\xi_j (\pa_{v_k}\bar{v}_j)v_{k,x}\theta_j^\p\left(\frac{w_j}{v_j}\right)_x+\sum_k \psi_{ji,v\si}\xi_j^\p\left(\frac{w_j}{v_j}\right)_x (\pa_{v_k}\bar{v}_j)v_{k,x}\nonumber\\
			&+\sum_k\sum\limits_{l} \psi_{ji,v\si}\xi_j (\pa_{v_kv_l}\bar{v}_j)v_{l,x}v_{k,x}
			+\sum_k \psi_{ji,v\si}\xi_j (\pa_{v_k}\bar{v}_j)v_{k,xx}.\label{psijixx}
		\end{align}
		By using Lemmas \ref{estimimpo1}, \ref{lemme11.3} and the fact that $\psi_{ji,u\si},\psi_{ji,u u \si}, \psi_{ji,u \si \si}=O(1)\bar{v}_j \mathfrak{A}_j$ we obtain
		\begin{align*}
			&\sum_{k,l}\psi_{ji,uu\si}(\tilde{r}_k\otimes\tilde{r}_l) v_lv_k+\sum_k\psi_{ji,u\si}[\tilde{r}_{k,x} v_k+\tilde{r}_kv_{k,x}]+\sum_k\psi_{ji,uv\si}\tilde{r}_k \xi_j(\pa_{v_j}\bar{v}_j)v_{j,x}v_k\\
			&=O(1) \left(\sum\limits_{k,l}\abs{v_j v_kv_l}+\sum\limits_{k}(\abs{v_kv_j}+|v_{k,x}v_j|+\abs{v_kv_{j,x}})\right)\rho_j^\e\mathfrak{A}_j,\\
			&\sum_{k}\sum\limits_{l\neq j}\psi_{ji,uv\si}\tilde{r}_k \xi_j(v_{l,x}\pa_{v_l}\bar{v}_j)v_k
			+\sum_k(\xi_j^\p\bar{v}_j\psi_{ji,uv\si}-\theta_j^\p\psi_{ji,u\si\si})\tilde{r}_k v_k\left(\frac{w_j}{v_j}\right)_x\\
			&=O(1) \left(\sum\limits_{k}\sum\limits_{l\neq j}\abs{v_k}\abs{v_{l,x}}+\sum\limits_{k}\abs{v_k}(\abs{v_{j,x}}+\abs{w_{j,x}})\right)\rho_j^\e\mathfrak{A}_j.
		\end{align*}
		Using now the fact that $\psi_{ji,\sig\si}=O(1)\bar{v}_j \mathfrak{A}_j$ we get
		\begin{align*}
			&\left(\frac{w_j}{v_j}\right)_{xx}(\xi^\p_j\bar{v}_j \psi_{ji,v\si}-\theta^\p_j \psi_{ji,\sig\si})+\left(\frac{w_j}{v_j}\right)^2_x(\xi^{\p\p}_j\bar{v}_j \psi_{ji,v\si}-\theta^{\p\p}_j \psi_{ji,\sig\si})\\
			&=O(1)\left(\abs{\bar{v}_j\left(\frac{w_j}{v_j}\right)_{xx}}+\abs{\bar{v}_j\left(\frac{w_j}{v_j}\right)_{x}^2}\right)\rho_j^\e\mathfrak{A}_j,\\
			&\left(\frac{w_j}{v_j}\right)_x\xi^\p_j\pa_x(\bar{v}_j) \psi_{ji,v\si}+\left(\frac{w_j}{v_j}\right)_x\sum\limits_{k}v_k(\xi^\p_j\bar{v}_j \psi_{ji,uv\si}-\theta^\p_j \psi_{ji,u\sig\si})\tilde{r}_k\\
			&=O(1)|\left(\frac{w_j}{v_j}\right)_x|\left(\sum\limits_{k\neq j}\abs{v_{k,x}v_k}+\abs{v_{j,x}}+\abs{v_j}\right)\rho_j^\e\mathfrak{A}_j.
		\end{align*}
		Similarly we have
		\begin{align*}
			&\left(\frac{w_j}{v_j}\right)^2_x(\xi^\p_j\bar{v}_j \psi_{ji,vv\si}-\theta^\p_j \psi_{ji,v\sig\si})\xi^\p_j\bar{v}_j+\left(\frac{w_j}{v_j}\right)_x\sum\limits_{k}(\xi^\p_j\bar{v}_j \psi_{ji,vv\si}-\theta^\p_j \psi_{ji,v\sig\si})\xi_jv_{k,x}\pa_{v_k}\bar{v}_j\\
			&-\left(\frac{w_j}{v_j}\right)^2_x(\xi^\p_j\bar{v}_j \psi_{ji,v\si\si}-\theta^\p_j \psi_{ji,\si\sig\si})\theta_j^\p\\
			&=O(1) \left(\left(\frac{w_j}{v_j}\right)^2_x|\bar{v}_j|+|\left(\frac{w_j}{v_j}\right)_x|\sum_{k\ne j} |v_{k,x}v_k| \right)\rho_j^\e\mathfrak{A}_j,\\
			&\sum_k\sum\limits_{l} \psi_{ji,uv\si}\tilde{r}_l\xi_j (\pa_{v_k}\bar{v}_j)v_lv_{k,x}+\sum_k \left(\frac{w_j}{v_j}\right)_x\psi_{ji,vv\si}\xi_j^\p{\bar{v}_j}\xi_j (\pa_{v_k}\bar{v}_j)v_{k,x}\\
			&=O(1)\left(\sum\limits_{k,l}\abs{v_{k,x}v_l}+\big (\abs{v_{j,x}\bar{v}_j}+|\bar{v_j}|\sum\limits_{k\neq j}\abs{v_kv_{k,x}}\big)\abs{\left(\frac{w_j}{v_j}\right)_x} \right)\rho_j^\e\mathfrak{A}_j.
		\end{align*}
		We have again
		\begin{align*}
			&\sum_k\sum\limits_{l} \psi_{ji,vv\si}\xi_j^2(\pa_{v_l}\bar{v}_j) (\pa_{v_k}\bar{v}_j)v_{l,x}v_{k,x}=O(1)\rho_j^\e\mathfrak{A}_j\sum\limits_{k,l}\abs{v_{k,x}v_{l,x}},\\
			&-\sum_k \psi_{ji,v\si\si}\xi_j (\pa_{v_k}\bar{v}_j)v_{k,x}\theta_j^\p\left(\frac{w_j}{v_j}\right)_x+\sum_k \psi_{ji,v\si}\xi_j^\p\left(\frac{w_j}{v_j}\right)_x (\pa_{v_k}\bar{v}_j)v_{k,x}\\
			&=O(1)\left( \sum\limits_{k\ne j}\abs{v_kv_{k,x}}\abs{\left(\frac{w_j}{v_j}\right)_x}+ \abs{v_{j,x}\left(\frac{w_j}{v_j}\right)_x}\right)\rho_j^\e\mathfrak{A}_j.
		\end{align*}
		Similarly it yields using Lemma \ref{estimimpo1}
		\begin{align*}
			&\sum_k\sum\limits_{l} \psi_{ji,v\si}\xi_j (\pa_{v_kv_l}\bar{v}_j)v_{l,x}v_{k,x}+\sum_k \psi_{ji,v\si}\xi_j (\pa_{v_k}\bar{v}_j)v_{k,xx}\\
			&=O(1) \left(\abs{v_{j,xx}}+\sum\limits_{k\neq j}\abs{v_kv_{k,xx}}+\sum\limits_{l\neq j}\frac{\abs{v_l v_{l,x}v_{j,x}}}{\abs{v_j}}+\frac{v_{j,x}^2}{|v_j|}+\sum\limits_{k\neq j,l\neq j}\abs{v_{k,x}v_{l,x}} \right)\rho_j^\e\mathfrak{A}_j .
		\end{align*}
		Combining all the above estimates we obtain
		\begin{align*}
			\pa_{xx}\psi_{ji,\si}&=O(1)\left(\sum\limits_{k}(\abs{v_kv_j}+|v_{k,x}v_j|+\abs{v_kv_{j,x}}+ \abs{v_k w_{j,x}})
			\right)\rho_j^\e\mathfrak{A}_j \\
			&+O(1) \left(\abs{v_{j,xx}}+\sum\limits_{k\neq j}\abs{v_kv_{k,xx}}+\sum\limits_{l\neq j} \frac{\abs{v_l v_{l,x}v_{j,x}}}{\abs{v_j}}+\frac{v_{j,x}^2}{|v_j|} \right)\rho_j^\e\mathfrak{A}_j \\
			&+ O(1)\left(\sum\limits_{k}\sum\limits_{l\neq j}\abs{v_kv_{l,x}}+\sum\limits_{k\ne j,l\ne j}\abs{v_{k,x}v_{l,x}} +\abs{v_j\left(\frac{w_j}{v_j}\right)_{xx}}+\abs{v_j\left(\frac{w_j}{v_j}\right)_{x}^2}\right)\rho_j^\e\mathfrak{A}_j  \\
			&+O(1)|\left(\frac{w_j}{v_j}\right)_x|\left(\sum\limits_{k\neq j}\abs{v_{k,x}v_k}+\abs{v_{j,x}}+\abs{v_j}\right)\rho_j^\e\mathfrak{A}_j.
		\end{align*}
		We obtain a similar result for $\pa_{xx}\psi_{ji}$. Next we calculate $\pa_{xx}\psi_{ji,v}$ as in \eqref{psijixx}
		\begin{align*}
			\pa_{xx}\psi_{ji,v}&=\sum_{k,l}\psi_{ji,uuv}(\tilde{r}_k\otimes\tilde{r}_l) v_lv_k+\sum_k\psi_{ji,uv}[\tilde{r}_{k,x} v_k+\tilde{r}_kv_{k,x}]+\sum_k\psi_{ji,uvv}\tilde{r}_k \xi_j(\pa_{v_j}\bar{v}_j)v_{j,x}v_k\\
			&+\sum_{k}\sum\limits_{l\neq j}\psi_{ji,uvv}\tilde{r}_k \xi_j(v_{l,x}\pa_{v_l}\bar{v}_j)v_k
			+\sum_k(\xi_j^\p\bar{v}_j\psi_{ji,uvv}-\theta_j^\p\psi_{ji,uv\si})\tilde{r}_k v_k\left(\frac{w_j}{v_j}\right)_x\\
			&+\left(\frac{w_j}{v_j}\right)_{xx}(\xi^\p_j\bar{v}_j \psi_{ji,vv}-\theta^\p_j \psi_{ji,v\si})+\left(\frac{w_j}{v_j}\right)^2_x(\xi^{\p\p}_j\bar{v}_j \psi_{ji,vv}-\theta^{\p\p}_j \psi_{ji,v\si})\\
			&+\left(\frac{w_j}{v_j}\right)_x\xi^\p_j\pa_x(\bar{v}_j) \psi_{ji,vv}+\left(\frac{w_j}{v_j}\right)_x\sum\limits_{k}v_k(\xi^\p_j\bar{v}_j \psi_{ji,uvv}-\theta^\p_j \psi_{ji,uv\si})\tilde{r}_k\\
			&+\left(\frac{w_j}{v_j}\right)^2_x(\xi^\p_j\bar{v}_j \psi_{ji,vvv}-\theta^\p_j \psi_{ji,vv\si})\xi^\p_j\bar{v}_j+\left(\frac{w_j}{v_j}\right)_x\sum\limits_{k}(\xi^\p_j\bar{v}_j \psi_{ji,vvv}-\theta^\p_j \psi_{ji,vv\si})\xi_j\pa_{v_k}\bar{v}_j {v_{k,x}}\\
			&-\left(\frac{w_j}{v_j}\right)^2_x(\xi^\p_j\bar{v}_j \psi_{ji,vv\si}-\theta'_j \psi_{ji,v\sig\si})\theta_j^\p\\
			&+\sum_k\sum\limits_{l} \psi_{ji,uvv}\tilde{r}_l\xi_j (\pa_{v_k}\bar{v}_j)v_lv_{k,x}+\sum_k {\bar{v}_j}\left(\frac{w_j}{v_j}\right)_x\psi_{ji,vvv}\xi_j^\p\xi_j (\pa_{v_k}\bar{v}_j)v_{k,x}\\
			&+\sum_k\sum\limits_{l} \psi_{ji,vvv}\xi_j^2(\pa_{v_l}\bar{v}_j) (\pa_{v_k}\bar{v}_j)v_{l,x}v_{k,x}\\
			&-\sum_k \psi_{ji,vv\si}\xi_j (\pa_{v_k}\bar{v}_j)v_{k,x}\theta_j^\p\left(\frac{w_j}{v_j}\right)_x+\sum_k \psi_{ji,vv}\xi_j^\p\left(\frac{w_j}{v_j}\right)_x (\pa_{v_k}\bar{v}_j)v_{k,x}\\
			&+\sum_k\sum\limits_{l} \psi_{ji,vv}\xi_j (\pa_{v_kv_l}\bar{v}_j)v_{l,x}v_{k,x}
			+\sum_k \psi_{ji,vv}\xi_j (\pa_{v_k}\bar{v}_j)v_{k,xx}.
		\end{align*}
		We have now
		\begin{align*}
			&\sum_{k,l}\psi_{ji,uuv}(\tilde{r}_k\otimes\tilde{r}_l) v_lv_k+\sum_k\psi_{ji,uv}[\tilde{r}_{k,x} v_k+\tilde{r}_kv_{k,x}]+\sum_k\psi_{ji,uvv}\tilde{r}_k \xi_j(\pa_{v_j}\bar{v}_j)v_{j,x}v_k\\
			&=O(1)\sum_{k}(|v_k|+|v_{k,x}|),\\
			&\sum_{k}\sum\limits_{l\neq j}\psi_{ji,uvv}\tilde{r}_k \xi_j(v_{l,x}\pa_{v_l}\bar{v}_j)v_k
			+\sum_k(\xi_j^\p\bar{v}_j\psi_{ji,uvv}-\theta_j^\p\psi_{ji,uv\si})\tilde{r}_k v_k\left(\frac{w_j}{v_j}\right)_x\\
			&=O(1)\sum_k |v_k|+O(1)\sum_k |v_k\left(\frac{w_j}{v_j}\right)_x|\mathfrak{A}_j.
		\end{align*}
		Similarly we have
		\begin{align*}
			&\left(\frac{w_j}{v_j}\right)_{xx}(\xi^\p_j\bar{v}_j \psi_{ji,vv}-\theta^\p_j \psi_{ji,v\si})+\left(\frac{w_j}{v_j}\right)^2_x(\xi^{\p\p}_j\bar{v}_j \psi_{ji,vv}-\theta^{\p\p}_j \psi_{ji,v\si})\\
			&=O(1)\biggl(\abs{\left(\frac{w_j}{v_j}\right)_{xx}}+\abs{\left(\frac{w_j}{v_j}\right)_{x}^2}\biggl)\mathfrak{A}_j,\\
			&\left(\frac{w_j}{v_j}\right)_x\xi^\p_j\pa_x(\bar{v}_j) \psi_{ji,vv}+\left(\frac{w_j}{v_j}\right)_x\sum\limits_{k}v_k(\xi^\p_j\bar{v}_j \psi_{ji,uvv}-\theta^\p_j \psi_{ji,uv\si})\tilde{r}_k\\
			&=O(1) |\left(\frac{w_j}{v_j}\right)_x|\biggl(\sum_k|v_k|+|v_{j,x}|\rho_j^\e\biggl)\mathfrak{A}_j.
		\end{align*}
		Next we get
		\begin{align*}
			&\left(\frac{w_j}{v_j}\right)^2_x(\xi^\p_j\bar{v}_j \psi_{ji,vvv}-\theta'_j \psi_{ji,vv\si})\xi^\p_j\bar{v}_j+\left(\frac{w_j}{v_j}\right)_x\sum\limits_{k}(\xi^\p_j\bar{v}_j \psi_{ji,vvv}-\theta^\p_j \psi_{ji,vv\si})\xi_j\pa_{v_k}\bar{v}_j v_{k,x}\\
			&-\left(\frac{w_j}{v_j}\right)^2_x(\xi^\p_j\bar{v}_j \psi_{ji,vv\si}-\theta'_j \psi_{ji,v\sig\si})\theta_j^\p\\
			&=O(1)\biggl[ \left(\frac{w_j}{v_j}\right)^2_x
			+|\left(\frac{w_j}{v_j}\right)_x|\big(|v_{j,x}|+\sum_{k\ne j}|v_{k}v_{k,x}|\big)\rho_j^\e \biggl]\mathfrak{A}_j.
		\end{align*}
		We have now
		\begin{align*}
			&\sum_k\sum\limits_{l} \psi_{ji,uvv}\tilde{r}_l\xi_j (\pa_{v_k}\bar{v}_j)v_lv_{k,x}+\sum_k {\bar{v}_j}\left(\frac{w_j}{v_j}\right)_x\psi_{ji,vvv}\xi_j^\p\xi_j (\pa_{v_k}\bar{v}_j)v_{k,x}\\
			&+\sum_k\sum\limits_{l} \psi_{ji,vvv}\xi_j^2(\pa_{v_l}\bar{v}_j) (\pa_{v_k}\bar{v}_j)v_{l,x}v_{k,x}\\
			&=O(1)\sum_k |v_{k,x}|+O(1)|v_j\left(\frac{w_j}{v_j}\right)_x| \biggl(|v_{j,x}|+\sum_{k\ne j}|v_{k}v_{k,x}|\biggl)\rho_j^\e\mathfrak{A}_j.
		\end{align*}
		Next we have
		\begin{align*}
			&-\sum_k \psi_{ji,vv\si}\xi_j (\pa_{v_k}\bar{v}_j)v_{k,x}\theta_j^\p\left(\frac{w_j}{v_j}\right)_x+\sum_k \psi_{ji,vv}\xi_j^\p\left(\frac{w_j}{v_j}\right)_x (\pa_{v_k}\bar{v}_j)v_{k,x}\\
			&+\sum_k\sum\limits_{l} \psi_{ji,vv}\xi_j (\pa_{v_kv_l}\bar{v}_j)v_{l,x}v_{k,x}+\sum_k \psi_{ji,vv}\xi_j (\pa_{v_k}\bar{v}_j)v_{k,xx}\\
			&=O(1)\sum_k|v_{k,x}|+O(1)|\left(\frac{w_j}{v_j}\right)_x| \biggl(|v_{j,x}|+\sum_{k\ne j}|v_{k}v_{k,x}|\biggl)\rho_j^\e\mathfrak{A}_j\\
			&+O(1)\biggl(\frac{v_{j,x}^2}{|v_j|}+|v_{j,xx}|+\sum_{k\ne j} |\frac{v_k}{v_j}| |v_{k,x}v_{j,x}|+\sum_{k\ne j}|v_{k}v_{k,xx}| \biggl)\rho_j^\e\mathfrak{A}_j.
		\end{align*}
		Combining all the previous estimate, we get from Lemma \ref{lemme6.5}, \ref{lemme6.6}
		\begin{align*}
			\pa_{xx}\psi_{ji,v}
			&=O(1)\sum_k (|v_k|+|v_{k,x}|)+O(1) |\left(\frac{w_j}{v_j}\right)_x|\biggl(\sum_k|v_k|+|v_{j,x}|\rho_j^\e\biggl)\mathfrak{A}_j\\
			&+O(1)\biggl[ \left(\frac{w_j}{v_j}\right)^2_x+\abs{\left(\frac{w_j}{v_j}\right)_{xx}}
			%+|\left(\frac{w_j}{v_j}\right)_x|\big(|v_{j,x}|+\sum_{k\ne j}|v_{k}v_{k,x}|\big)\rho_j 
			\biggl]\mathfrak{A}_j\\
			%&+O(1)\sum_k|v_{k,x}|+O(1)|\left(\frac{w_j}{v_j}\right)_x| \biggl(|v_{j,x}|+\sum_{k\ne j}|v_{k}v_{k,x}|\biggl)\rho_j^\e\mathfrak{A}_j\\
			&+O(1)\biggl(\frac{v_{j,x}^2}{|v_j|}+|v_{j,xx}|+\sum_{k\ne j} |\frac{v_k}{v_j}| |v_{k,x}v_{j,x}|+\sum_{k\ne j}|v_{k}v_{k,xx}| \biggl)\rho_j^\e \mathfrak{A}_j.
		\end{align*}
		This completes the proof of estimate \ref{estimate-psi-xx-v}.
\end{proof}

\begin{lemma}
\label{lemme11.6}
	For $1\leq i,j,k\leq n$ the following bounds hold.
	\begin{align*}
		&\tilde{r}_{k,xx}=O(1)\biggl(\sum\limits_{j}|v_{j,x}|+\sum\limits_{j,l}|v_lv_{j}|+\sum\limits_{j}|v_j^2
		\left(\frac{w_j}{v_j}\right)_x|\rho_j^\e\mathfrak{A}_j\biggl)\\
		&
		+O(1)\rho_k^\e\mathfrak{A}_k\bigg( |\left(\frac{w_k}{v_k}\right)_x v_k| \sum\limits_{j}|v_{j}|+|v_{k,x}| \sum\limits_{j}|v_{j}|
		+|v_{k,xx}|+|v_{k,x}\left(\frac{w_k}{v_k}\right)_x|+ \frac{v_{k,x}^2}{|v_k|}\\
		&+\sum_{l\ne k}|\frac{v_lv_{k,x}v_{l,x}}{v_k}|+| \left(\frac{w_k}{v_k}\right)_x|\sum\limits_{j\neq k}|v_{j,x}v_j|
		+\sum_{j\ne k}|v_{j,xx}v_j|+\left(\frac{w_k}{v_k}\right)_x^2|v_k|	+| \left(\frac{w_k}{v_k}\right)_{xx} v_k|\biggl),\\
			&\tilde{\la}_{j,xx}=O(1)\sum_k(|v_k|+|v_{k,x}|)+O(1)\mathfrak{A}_j \rho_j^\e \biggl(|\left(\frac{w_j}{v_j}\right)_{xx}v_j|+|v_j\left(\frac{w_j}{v_j}\right)_{x}^2|+|v_j \left(\frac{w_j}{v_j}\right)_{x}|\\
		&+| \left(\frac{w_j}{v_j}\right)_{x}v_{j,x}|+| \left(\frac{w_j}{v_j}\right)_{x}|\sum_{k\ne j}|v_kv_{k,x}|+\frac{v_{j,x}^2}{|v_j|}+|v_{j,x}|\sum_{k\ne j}|\frac{v_kv_{k,x}}{v_j}|+|v_{j,xx}|%+\sum_{k\ne j}|v_kv_{k,xx}|+\sum_{k\ne j,l\ne j}|v_{k,x}v_{l,x}|
		\biggl).
	\end{align*}
	Moreover, we have
	\begin{align*}
		&\pa_{xx}(\tilde{r}_{i,u}\tilde{r}_j)
		=O(1)\biggl(\sum\limits_{k=1}|v_{k,x}|+\sum_{k,l} |v_k v_l|+\sum_k\rho_k^\e \mathfrak{A}_k|v_k\bar{v}_k\left(\frac{w_k}{v_k}\right)_x|
		\biggl)\\
		&+O(1)\rho_i^\e\mathfrak{A}_i\biggl( |v_i\left(\frac{w_i}{v_i}\right)_x| \sum_k|v_k|+\sum_{k\ne i}|v_kv_{k,xx}|
		+|v_i\left(\frac{w_i}{v_i}\right)_{xx}|
		+  |v_i\left(\frac{w_i}{v_i}\right)^2_x|\\
		&\hspace{2cm}+|v_{i,x}\left(\frac{w_i}{v_i}\right)_x|+ \sum_{k\ne i}|\left(\frac{w_i}{v_i}\right)_x||v_kv_{k,x}|+\rho_j^\e\mathfrak{A}_j|v_iv_{j,x} \left(\frac{w_i}{v_i}\right)_x|\\
		&\hspace{2cm}+|\left(\frac{w_i}{v_i}\right)_x v_i\left(\frac{w_j}{v_j}\right)_x|\rho_j^{\e}\mathfrak{A}_j+|v_{i,xx}|+\frac{v_{i,x}^2}{|v_i|}+\sum_{k\ne i}|\frac{v_{k,x}v_kv_{i,x}}{v_i}|\\
		&\hspace{2cm}+\sum\limits_{k=1}^{n}\abs{v_{i,x}v_k}+\rho_j^\e\mathfrak{A}_j|v_{j,x}v_{i,x}|+\abs{v_{i,x}v_j\left(\frac{w_j}{v_j}\right)_x}\rho_j^\e\mathfrak{A}_j\biggl)\\
		&
		+O(1)\rho_j^\e\mathfrak{A}_j \biggl(|v_j\left(\frac{w_j}{v_j}\right)_x| \sum_k|v_k|+\sum_{k\ne j}|v_kv_{k,xx}|+|v_j \left(\frac{w_j}{v_j}\right)_{xx}|+  |v_j\left(\frac{w_j}{v_j}\right)^2_x|+|v_{j,x}\left(\frac{w_j}{v_j}\right)_x|\\
		&\hspace{1.5cm}+|\left(\frac{w_j}{v_j}\right)_x|\sum_{k\ne j}|v_kv_{k,x}|+\rho_i^{\e}\mathfrak{A}_i |\left(\frac{w_j}{v_j}\right)_x v_j v_{i,x}|+\abs{\left(\frac{w_j}{v_j}\right)_x v_j\left(\frac{w_i}{v_i}\right)_x}\rho_i^\e\mathfrak{A}_i\\
		&\hspace{1cm}+|v_{j,xx}|+\frac{v_{j,x}^2}{|v_j|}
		+\sum_{k\ne j}|\frac{v_{k,x}v_kv_{j,x}}{v_j}|
		+\sum\limits_{k=1}^{n}\abs{v_{j,x}v_k}+\rho_i^\e\mathfrak{A}_i |v_{j,x}v_{i,x}|
		+\abs{v_{j,x}v_i\left(\frac{w_i}{v_i}\right)_x}\rho_i^\e\mathfrak{A}_i\biggl).
	\end{align*}
\end{lemma}

\begin{proof}
	
	\noi\underline{Estimate of $\tilde{r}_{k,xx}$:} 
	Taking derivative of $\tilde{r}_{k,x}$ in \eqref{derivative-r-k-x} we get %recalling that $\theta_k^\p \tilde{r}_{k,\si}= \tilde{r}_{k,\si}$
	\begin{align}
		\tilde{r}_{k,xx}&=\sum\limits_{j}(v_j\tilde{r}_{k,u}\tilde{r}_j)_x+(v_{k,x}\xi_k\pa_{v_k}\bar{v}_k\tilde{r}_{k,v})_x+\sum\limits_{j\neq k}(v_{j,x}\xi_k\pa_{v_j}\bar{v}_k\tilde{r}_{k,v})_x\\
		&+\left(\xi_k^\p\left(\frac{w_k}{v_k}\right)_x\bar{v}_k\tilde{r}_{k,v}\right)_x-\left(\left(\frac{w_k}{v_k}\right)_x\tilde{r}_{k,\si}\right)_x. 
	\end{align}
   Calculating the terms explicitly, we deduce from Lemmas \ref{lemme6.5}, \ref{lemme6.6}, \ref{estimimpo1} and the fact that $\tilde{r}_{k,u\sig}=O(1)\xi_k\bar{v}_k$, $\tilde{r}_{j,\sig}=O(1)\xi_j\bar{v}_j$, $\tilde{r}_{k,\sig},\tilde{r}_{k,\sig\sig},\tilde{r}_{k,u\sig}=O(1)\xi_k\bar{v}_k$ that
	\begin{align*}
		&\sum\limits_{j}(v_j\tilde{r}_{k,u}\tilde{r}_j)_x=\sum\limits_{j}v_{j,x}\tilde{r}_{k,u}\tilde{r}_j+\sum\limits_{j,l}v_lv_{j}\tilde{r}_{k,uu}:\tilde{r}_j\otimes\tilde{r}_l+\sum\limits_{j}\xi'_k \left(\frac{w_k}{v_k}\right)_x \bar{v}_kv_{j}\tilde{r}_{k,uv}\tilde{r}_j\\
		&+\sum\limits_{j}\sum_l \xi_k v_{l,x}\pa_{v_l}\bar{v}_k v_{j}\tilde{r}_{k,uv}\tilde{r}_j-\sum\limits_{j}\theta'_k \left(\frac{w_k}{v_k}\right)_x v_j\tilde{r}_{k,u \sig}\tilde{r}_j\\
		&+\sum\limits_{j,l}v_lv_j\tilde{r}_{k,u}\tilde{r}_{j,u}\tilde{r}_l+\sum\limits_{j}\xi'_j \left(\frac{w_j}{v_j}\right)_x \bar{v}_j v_j\tilde{r}_{k,u}\tilde{r}_{j,v}+\sum\limits_{j, l}\xi_j\pa_{v_l}\bar{v}_j v_{l,x}v_j\tilde{r}_{k,u}\tilde{r}_{j,v}\\
		&-\sum\limits_{j}\theta'_j \left(\frac{w_j}{v_j}\right)_x v_j\tilde{r}_{k,u}\tilde{r}_{j,\sig},\\
		&=O(1)\biggl(\sum_j|v_{j,x}|+\sum_{j,l}|v_j v_l|+\sum_j|\left(\frac{w_k}{v_k}\right)_x v_kv_j|\mathfrak{A}_k\rho_k^\e|+\sum_j\sum_{l\ne k}|v_{l,x}v_l v_j|\mathfrak{A}_k\rho_k^\e\\
		&+\sum_j|v_{k,x}v_j|\mathfrak{A}_k\rho_k^\e+\sum_jv_j^2|\left(\frac{w_j}{v_j}\right)_x|\mathfrak{A}_j\rho_j^\e+\sum_j\sum_{l\ne j}|v_lv_{l,x}v_j|\mathfrak{A}_j\rho_j^\e+\sum_j|v_jv_{j,x}|\mathfrak{A}_j\rho_j^\e\biggl).
		\\[2mm]
		&(v_{k,x}\xi_k\pa_{v_k}\bar{v}_k\tilde{r}_{k,v})_x
		=v_{k,xx}\xi_k\pa_{v_k}\bar{v}_k\tilde{r}_{k,v}+v_{k,x}\xi'_k\left(\frac{w_k}{v_k}\right)_x \pa_{v_k}\bar{v}_k\tilde{r}_{k,v}+\sum_l v_{k,x}\xi_k\pa_{v_l v_k}\bar{v}_k v_{l,x}\tilde{r}_{k,v}\\
		&+
		\sum_l v_l v_{k,x}\xi_k\pa_{v_k}\bar{v}_k\tilde{r}_{k,uv}\tilde{r}_l+\xi'_k\left(\frac{w_k}{v_k}\right)_x \bar{v}_k v_{k,x}\xi_k\pa_{v_k}\bar{v}_k\tilde{r}_{k,vv}\\
		&
		+\sum_l v_{k,x}\xi_k\pa_{v_k}\bar{v}_k \pa_{v_l}\bar{v}_k v_{l,x}\tilde{r}_{k,vv}-\theta'_k \left(\frac{w_k}{v_k}\right)_x v_{k,x}\xi_k\pa_{v_k}\bar{v}_k\tilde{r}_{k,v \sig},\\
		&=O(1)\mathfrak{A}_k\rho_k^\e\biggl(|v_{k,xx}|+|v_{k,x}\left(\frac{w_k}{v_k}\right)_x|+\frac{v_{k,x}^2}{|v_k|}+\sum_{l\ne k}\frac{|v_lv_{l,x}v_{k,x}|}{|v_k|}+\sum_l|v_lv_{k,x}|\biggl).
		\\[2mm]
		&\sum\limits_{j\neq k}(v_{j,x}\xi_k\pa_{v_j}\bar{v}_k\tilde{r}_{k,v})_x=\sum\limits_{j\neq k}v_{j,xx}\xi_k\pa_{v_j}\bar{v}_k\tilde{r}_{k,v}+\sum\limits_{j\neq k}v_{j,x}\xi'_k\left(\frac{w_k}{v_k}\right)_x\pa_{v_j}\bar{v}_k\tilde{r}_{k,v}\\
		&+\sum\limits_{j\neq k}\sum_l v_{j,x}\xi_k\pa_{v_j v_l}\bar{v}_kv_{l,x} \tilde{r}_{k,v}+\sum\limits_{j\neq k}\sum_l v_l v_{j,x}\xi_k\pa_{v_j}\bar{v}_k\tilde{r}_{k,uv}\tilde{r}_l\\
		&+\sum\limits_{j\neq k}\xi'_k\left(\frac{w_k}{v_k}\right)_x\bar{v}_k v_{j,x}\xi_k\pa_{v_j}\bar{v}_k\tilde{r}_{k,vv}+\sum\limits_{j\neq k}\sum_l \xi_k \pa_{v_l}\bar{v}_k v_{l,x} v_{j,x}\xi_k\pa_{v_j}\bar{v}_k\tilde{r}_{k,vv}\\
		&-\sum\limits_{j\neq k}\theta'_k\left(\frac{w_k}{v_k}\right)_x v_{j,x}\xi_k\pa_{v_j}\bar{v}_k\tilde{r}_{k,v \sig },\\
		&=O(1)\mathfrak{A}_k\rho_k^\e\biggl(\sum_{j\ne k}|v_jv_{j,xx}|+\sum_{j\ne k}|\left(\frac{w_k}{v_k}\right)_x||v_{j,x}v_j|+\sum_{j\ne k, l\ne k}|v_{j,x}v_{l,x}|+\sum_{j\ne k}\frac{|v_{j,x}v_jv_{k,x}|}{|v_k|}\\
		&+\sum_{j\ne k}\sum_l |v_lv_jv_{j,x}|+\sum_{j\ne k}|v_{k,x}v_{j,x}v_j|\biggl).
		\\[2mm]
		&\left(\xi_k^\p\left(\frac{w_k}{v_k}\right)_x\bar{v}_k\tilde{r}_{k,v}\right)_x=\xi_k''\left(\frac{w_k}{v_k}\right)_x^2\bar{v}_k\tilde{r}_{k,v}+\xi_k^\p\left(\frac{w_k}{v_k}\right)_{xx}\bar{v}_k\tilde{r}_{k,v}+\sum_l \xi_k^\p \left(\frac{w_k}{v_k}\right)_x \pa_{v_l}\bar{v}_kv_{l,x}\tilde{r}_{k,v}\\
		&+\sum_l v_l\xi_k^\p\left(\frac{w_k}{v_k}\right)_x\bar{v}_k\tilde{r}_{k,uv}\tilde{r}_l+\big(\xi'_k\left(\frac{w_k}{v_k}\right)_x\bar{v}_k\big)^2\tilde{r}_{k,vv}+\sum_l \xi_k\pa_{v_l}\bar{v}_k v_{l,x}\xi_k^\p\left(\frac{w_k}{v_k}\right)_x\bar{v}_k\tilde{r}_{k,vv}\\
		&-\theta'_k \xi_k^\p\left(\frac{w_k}{v_k}\right)_x^2 \bar{v}_k\tilde{r}_{k,v\sig},\\
		&=O(1)\mathfrak{A}_k\rho_k^\e\biggl(|v_k\left(\frac{w_k}{v_k}\right)_x^2|+| \left(\frac{w_k}{v_k}\right)_{xx} v_k|+\sum_{l\ne k}|\left(\frac{w_k}{v_k}\right)_x v_lv_{l,x}|+|\left(\frac{w_k}{v_k}\right)_x v_{k,x}|\\
		&+|\left(\frac{w_k}{v_k}\right)_x v_{k}|\sum_l|v_l|\biggl).\\[2mm]
		&-\left(\left(\frac{w_k}{v_k}\right)_x\tilde{r}_{k,\si}\right)_x=-\left(\frac{w_k}{v_k}\right)_{xx}\tilde{r}_{k,\si}-\sum_l v_l \left(\frac{w_k}{v_k}\right)_x\tilde{r}_{k,u\si}\tilde{r}_l-\xi'_k \left(\frac{w_k}{v_k}\right)_x^2\bar{v}_k \tilde{r}_{k,v\si}\\
		&-\sum_l \xi_k\pa_{v_l}\bar{v}_k v_{l,x}\left(\frac{w_k}{v_k}\right)_x\tilde{r}_{k,\sig v}+\theta'_k \left(\frac{w_k}{v_k}\right)_x^2\tilde{r}_{k,\si \si},\\
		&=O(1)\mathfrak{A}_k\rho_k^\e\biggl(|v_k\left(\frac{w_k}{v_k}\right)_x^2|+| \left(\frac{w_k}{v_k}\right)_{xx} v_k|+\sum_{l}|\left(\frac{w_k}{v_k}\right)_x v_l v_k|+\sum_{l\ne k}|\left(\frac{w_k}{v_k}\right)_x v_lv_{l,x}|\\
		&+|\left(\frac{w_k}{v_k}\right)_x v_{k,x}|\biggl).
	\end{align*}
	We deduce from Lemmas \ref{lemme6.5}, \ref{lemme6.6}, \ref{estimimpo1} and the fact that $\tilde{r}_{k,u\sig}=O(1)\xi_k\bar{v}_k$, $\tilde{r}_{j,\sig}=O(1)\xi_j\bar{v}_j$, $\tilde{r}_{k,\sig},\tilde{r}_{k,\sig\sig},\tilde{r}_{k,u\sig}=O(1)\xi_k\bar{v}_k$ that
	\begin{align*}
		&\tilde{r}_{k,xx}=O(1)\biggl(\sum\limits_{j}|v_{j,x}|+\sum\limits_{j,l}|v_lv_{j}|+\sum\limits_{j}|v_j^2
		\left(\frac{w_j}{v_j}\right)_x|\rho_j^\e\mathfrak{A}_j\biggl)\\
		&
		+O(1)\rho_k^\e\mathfrak{A}_k\bigg( |\left(\frac{w_k}{v_k}\right)_x v_k| \sum\limits_{j}|v_{j}|+|v_{k,x}| \sum\limits_{j}|v_{j}|
		+|v_{k,xx}|+|v_{k,x}\left(\frac{w_k}{v_k}\right)_x|+ \frac{v_{k,x}^2}{|v_k|}\\
		&+\sum_{l\ne k}|\frac{v_lv_{k,x}v_{l,x}}{v_k}|+| \left(\frac{w_k}{v_k}\right)_x|\sum\limits_{j\neq k}|v_{j,x}v_j|
		+\sum_{j\ne k}|v_{j,xx}v_j|+\left(\frac{w_k}{v_k}\right)_x^2|v_k|	+| \left(\frac{w_k}{v_k}\right)_{xx} v_k|\biggl).
	\end{align*}
	
	\noi\underline{Estimate of $\pa_{xx}(\tilde{r}_{i,u}\tilde{r}_j)$:} Taking now the derivative of $\pa_x(\tilde{r}_{i,u}\tilde{r}_j)$ in \eqref{rix}, we have
	\begin{align}
		\pa_{xx}(\tilde{r}_{i,u}\tilde{r}_j)
		&=\sum\limits_{k=1}^{n}\left(v_k[\tilde{r}_{i,uu}:\tilde{r}_{j}\otimes\tilde{r}_k+\tilde{r}_{i,u}\tilde{r}_{j,u}\tilde{r}_k]\right)_x+\left(\left(\frac{w_i}{v_i}\right)_x\left[\xi^\p_i\bar{v}_i\tilde{r}_{i,uv}\tilde{r}_j-\theta_i^\p\tilde{r}_{i,u\si}\tilde{r}_j\right]\right)_x\nonumber\\
		&+\left(\left(\frac{w_j}{v_j}\right)_x\left[\xi^\p_j\bar{v}_j\tilde{r}_{i,u}\tilde{r}_{j,v}-\theta_j^\p\tilde{r}_{i,u}\tilde{r}_{j,\si}\right]\right)_x+(v_{i,x}\xi_i\pa_{v_i}\bar{v}_i\tilde{r}_{i,uv}\tilde{r}_j)_x+(v_{j,x}\xi_j\pa_{v_j}\bar{v}_j\tilde{r}_{i,u}\tilde{r}_{j,v})_x\nonumber\\
		&+\sum\limits_{k\neq i}(\left( v_{k,x}\pa_{v_k}\bar{v}_i\right)\xi_i\tilde{r}_{i,uv}\tilde{r}_j)_x+\sum\limits_{k\neq j}(\left( v_{k,x}\pa_{v_k}\bar{v}_j\right)\xi_j\tilde{r}_{i,u}\tilde{r}_{j,v})_x.\label{11.26}
	\end{align}
	Calculating first three terms we have using the fact that $\tilde{r}_{i,uu\sig},\tilde{r}_{i,u\sig}=O(1)\xi_i\bar{v}_i$, $\tilde{r}_{j,u\sig},\tilde{r}_{j,\sig}=O(1)\xi_j\bar{v}_j$, $\tilde{r}_{k,\sig}=O(1)\xi_k\bar{v}_k
	$

	\begin{align}
		&\sum\limits_{k=1}^{n}\left(v_k[\tilde{r}_{i,uu}:\tilde{r}_{j}\otimes\tilde{r}_k+\tilde{r}_{i,u}\tilde{r}_{j,u}\tilde{r}_k]\right)_x\nonumber\\
		&
		=\sum\limits_{k=1}^{n}v_{k,x}[\tilde{r}_{i,uu}:\tilde{r}_{j}\otimes\tilde{r}_k+\tilde{r}_{i,u}\tilde{r}_{j,u}\tilde{r}_k]+\sum_{k,l} v_k v_l \tilde{r}_{i,uuu}:\tilde{r}_{j}\otimes\tilde{r}_k\otimes \tilde{r}_l\nonumber\\
		&+\sum_k v_{k}[\tilde{r}_{i,uuv}:\tilde{r}_{j}\otimes\tilde{r}_k](\xi_i \bar{v}_i)_x-
		\sum_k v_{k}[\tilde{r}_{i,uu \sig}:\tilde{r}_{j}\otimes\tilde{r}_k]\theta'_i\left(\frac{w_i}{v_i}\right)_x\nonumber\\
		&+\sum_{k,l} v_lv_k \tilde{r}_{i,uu}:\tilde{r}_{j,u}\tilde{r}_k\otimes \tilde{r}_l+\sum_kv_k\tilde{r}_{i, uu}:\tilde{r}_{jv}\otimes\tilde{r}_k (\xi_j\bar{v}_j)_x-\sum_kv_k\tilde{r}_{i ,uu}:\tilde{r}_{j\sig}\otimes\tilde{r}_k \theta'_j\left(\frac{w_j}{v_j}\right)_x\nonumber\\
		&+\sum_{k,l} v_lv_k \tilde{r}_{i,uu}:\tilde{r}_{j}\otimes \tilde{r}_{k,u}\tilde{r}_l+\sum_kv_k\tilde{r}_{i uu}:\tilde{r}_{j}\otimes\tilde{r}_{k,v} (\xi_k\bar{v}_k)_x-\sum_kv_k\tilde{r}_{i uu}:\tilde{r}_{j}\otimes\tilde{r}_{k,\sig} \theta'_k\left(\frac{w_k}{v_k}\right)_x\nonumber\\
		&+\sum_{k,l} v_lv_k\tilde{r}_{i,uu}:\tilde{r}_{j,u}\tilde{r}_k\otimes\tilde{r}_l+\sum_k v_k \tilde{r}_{i,uv}\tilde{r}_{j,u}\tilde{r}_k(\xi_i\bar{v}_i)_x-\sum_k v_k \tilde{r}_{i,u\sig}\tilde{r}_{j,u}\tilde{r}_k\theta'_i\left(\frac{w_i}{v_i}\right)_x\nonumber\\
		&+\sum_{k,l}v_lv_{k}\tilde{r}_{i,u}\tilde{r}_{j,uu}:\tilde{r}_l\otimes\tilde{r}_k+\sum_{k}v_{k}\tilde{r}_{i,u}\tilde{r}_{j,uv}\tilde{r}_k(\xi_j\bar{v}_j)_x-\sum_{k}v_{k}\tilde{r}_{i,u}\tilde{r}_{j,u\sig}\tilde{r}_k \theta^\p_j\left(\frac{w_j}{v_j}\right)_x\nonumber\\
		&+\sum_{k,l}v_lv_{k}\tilde{r}_{i,u}\tilde{r}_{j,u}\tilde{r}_{k,u}\tilde{r}_l+\sum_{k}v_{k}\tilde{r}_{i,u}\tilde{r}_{ju}\tilde{r}_{kv}(\xi_k\bar{v}_k)_x-\sum_{k}v_{k}\tilde{r}_{i,u}\tilde{r}_{j,u}\tilde{r}_{k,\sig} \theta^\p_k\left(\frac{w_k}{v_k}\right)_x,\nonumber\\
		&=O(1)\biggl(\sum_k|v_{k,x}|+\sum_{k,l}|v_kv_l|+\sum_k\rho_i^\e\mathfrak{A}_i|v_k|(|v_i\left(\frac{w_i}{v_i}\right)_x|+|v_{i,x}|+\sum_{l\ne i}|v_lv_{l,x}|)\nonumber\\
		&%+\sum_k\rho_i^\e\mathfrak{A}_i |v_kv_i\left(\frac{w_i}{v_i}\right)_x|
		+\sum_k \rho_j^\e\mathfrak{A}_j|v_k|(|v_j\left(\frac{w_j}{v_j}\right)_x|+|v_{j,x}|+\sum_{l\ne j}|v_lv_{l,x}|)\nonumber\\
		&
		+\sum_k \rho_k^\e\mathfrak{A}_k|v_k|(|v_k\left(\frac{w_k}{v_k}\right)_x|+|v_{k,x}|+\sum_{l\ne k}|v_lv_{l,x}|)\biggl),\label{11.27}
	\end{align}
	and
	\begin{align}
		\left(\left(\frac{w_i}{v_i}\right)_x\left[\xi^\p_i\bar{v}_i\tilde{r}_{i,uv}\tilde{r}_j-\theta_i^\p\tilde{r}_{i,u\si}\tilde{r}_j\right]\right)_x
		&=\left(\frac{w_i}{v_i}\right)_{xx}\left[\xi^\p_i\bar{v}_i\tilde{r}_{i,uv}\tilde{r}_j-\theta_i^\p\tilde{r}_{i,u\si}\tilde{r}_j\right]\nonumber\\
		&+\left(\frac{w_i}{v_i}\right)_{x}\left[\xi^\p_i\bar{v}_i\tilde{r}_{i,uv}\tilde{r}_j-\theta_i^\p\tilde{r}_{i,u\si}\tilde{r}_j\right]_x,\label{11.28}\\
		\left(\left(\frac{w_j}{v_j}\right)_x\left[\xi^\p_j\bar{v}_j\tilde{r}_{i,u}\tilde{r}_{j,v}-\theta_j^\p\tilde{r}_{i,u}\tilde{r}_{j,\si}\right]\right)_x
		&=\left(\frac{w_j}{v_j}\right)_{xx}\left[\xi^\p_j\bar{v}_j\tilde{r}_{i,u}\tilde{r}_{j,v}-\theta_j^\p\tilde{r}_{i,u}\tilde{r}_{j,\si}\right]\nonumber\\
		&+\left(\frac{w_j}{v_j}\right)_{x}\left[\xi^\p_j\bar{v}_j\tilde{r}_{i,u}\tilde{r}_{j,v}-\theta_j^\p\tilde{r}_{i,u}\tilde{r}_{j,\si}\right]_x.\label{11.29}
	\end{align}
	First we observe using Lemmas \ref{lemme6.5}, \ref{lemme6.6}, \ref{estimimpo1}, \ref{lemme9.2},  \eqref{11.11-uv1}, \eqref{11.10-u-sig-1} and the fact that $\tilde{r}_{i,u\sig}=O(1)\xi_i\bar{v}_i$,
	\begin{align}
		&\left(\frac{w_i}{v_i}\right)_{x}\left[\xi^\p_i\bar{v}_i\tilde{r}_{i,uv}\tilde{r}_j-\theta_i^\p\tilde{r}_{i,u\si}\tilde{r}_j\right]_x\nonumber\\
		&=O(1) \rho_i^\e\mathfrak{A}_i |\left(\frac{w_i}{v_i}\right)_x|   \biggl( |v_i\left(\frac{w_i}{v_i}\right)_x|+|v_{i,x}|+ \sum_{k\ne i}|v_kv_{k,x}|+\sum\limits_{k=1}^{n} \abs{v_iv_k}\nonumber\\
		&\hspace{3cm}+ \rho^\e_j\mathfrak{A}_j|v_iv_{j,x}| +\abs{v_iv_j\left(\frac{w_j}{v_j}\right)_x}\rho_j^\e\mathfrak{A}_j\biggl).
		\label{vtech1}
	\end{align}
	Similarly we have using \eqref{11.11-uv2}, \eqref{11.10-u-sig-2} and Lemmas \ref{lemme6.5}, \ref{lemme6.6}
	\begin{align}
		&\left(\frac{w_j}{v_j}\right)_{x}\left[\xi^\p_j\bar{v}_j\tilde{r}_{i,u}\tilde{r}_{j,v}-\theta_j^\p\tilde{r}_{i,u}\tilde{r}_{j,\si}\right]_x\nonumber\\
		&=O(1)\rho_j^\e\mathfrak{A}_j |\left(\frac{w_j}{v_j}\right)_x| \biggl(|v_j\left(\frac{w_j}{v_j}\right)_x|+|v_{j,x}|+\sum_{k\ne j}|v_kv_{k,x}|+\sum\limits_{k=1}^{n}\abs{v_jv_k}\nonumber\\
		&\hspace{5cm}+\rho_i^\e\mathfrak{A}_i |v_j v_{i,x}|+\abs{v_jv_i\left(\frac{w_i}{v_i}\right)_x}\rho_i^\e\mathfrak{A}_i\biggl).\label{vtech2}
	\end{align}
	We have now from Lemmas \ref{estimimpo1} and \eqref{11.11-uv1}
	\begin{align}
		(v_{i,x}\xi_i\pa_{v_i}\bar{v}_i\tilde{r}_{i,uv}\tilde{r}_j)_x&=\rho_i^\e\mathfrak{A}_i\biggl(|v_{i,xx}|+|v_{i,x}\left(\frac{w_i}{v_i}\right)_x|+\frac{v_{i,x}^2}{|v_i|}
		+\sum_{k\ne i}|\frac{v_{k,x}v_kv_{i,x}}{v_i}|
		+\sum\limits_{k=1}^{n}\abs{v_{i,x}v_k}\nonumber\\\
		&\hspace{2cm}+\rho_j^\e\mathfrak{A}_j|v_{j,x}v_{i,x}|+\abs{v_{i,x}v_j\left(\frac{w_j}{v_j}\right)_x}\rho_j^\e\mathfrak{A}_j\biggl).\label{11.22}
	\end{align}
	In a similar way, we obtain using Lemma \ref{estimimpo1} and \eqref{11.11-uv2}
	\begin{align}
		(v_{j,x}\xi_j\pa_{v_j}\bar{v}_j\tilde{r}_{i,u}\tilde{r}_{j,v})_x&=\rho^\e_j\mathfrak{A}_j\biggl(|v_{j,xx}|+|v_{j,x}\left(\frac{w_j}{v_j}\right)_x|+\frac{v_{j,x}^2}{|v_j|}
		+\sum_{k\ne j}|\frac{v_{k,x}v_kv_{j,x}}{v_j}|+\sum\limits_{k=1}^{n}\abs{v_{j,x}v_k}\nonumber\\
		&\hspace{3cm}+\rho_i^\e\mathfrak{A}_i|v_{j,x}v_{i,x}|+\abs{v_{j,x}v_i\left(\frac{w_i}{v_i}\right)_x}\rho^\e_i\mathfrak{A}_i\biggl).\label{11.23}
	\end{align}
	Similarly we have using Lemma \ref{estimimpo1} and \eqref{11.11-uv1}
	\begin{align}
		&\sum\limits_{k\neq i}[\left( v_{k,x}\pa_{v_k}\bar{v}_i\right)\xi_i\tilde{r}_{i,uv}\tilde{r}_j]_x\nonumber\\
		&=O(1)\rho_i^\e\mathfrak{A}_i\biggl(\sum_{k\ne i}|v_{k,xx}v_k|+\sum_{k\ne i}|\frac{v_{k,x}v_kv_{i,x}}{v_i}|+\sum_{k\ne i}\sum_{l\ne i}|v_{k,x}v_{l,x}|+|\left(\frac{w_i}{v_i}\right)_x|\sum_{k\ne i}|v_{k,x}v_k|\biggl)
		\nonumber\\
		&+O(1)\rho_i^\e\mathfrak{A}_i \sum_{l\ne i}|v_lv_{l,x}|
		\biggl(\sum\limits_{k=1}^{n}\abs{v_k}+\rho_j^\e\mathfrak{A}_j|v_{j,x}|+\abs{v_j\left(\frac{w_j}{v_j}\right)_x}\rho_j^\e\mathfrak{A}_j\biggl).
		\label{11.24}
	\end{align}
	Similarly we get from  \eqref{11.11-uv2}
	\begin{align}
		&\sum\limits_{k\neq j}[\left( v_{k,x}\pa_{v_k}\bar{v}_j\right)\xi_j\tilde{r}_{i,u}\tilde{r}_{j,v}]_x\nonumber\\
		&=O(1)\rho_j^\e\mathfrak{A}_j\biggl(\sum_{k\ne j} |v_{k,xx}v_k|+\sum_{k\ne j}|\frac{v_{k,x}v_kv_{j,x}}{v_j}|
		+\sum_{k\ne j}\sum_{l\ne j}|v_{k,x}v_{l,x}|+|\left(\frac{w_j}{v_j}\right)_x|\sum_{k\ne j}|v_{k,x}v_k|\biggl)\nonumber\\
		&+O(1)\rho_j^\e\mathfrak{A}_j\sum_{k\ne j}|v_kv_{k,x}|\biggl(\sum\limits_{k=1}^{n}\abs{v_k}+\rho_i^\e\mathfrak{A}_i|v_{i,x}|+\abs{v_i\left(\frac{w_i}{v_i}\right)_x}\rho_i^\e\mathfrak{A}_i\biggl).
		\label{11.25}
	\end{align}
	We deduce then using Lemmas \ref{lemme6.5}, \ref{lemme6.6}, \eqref{11.26}, \eqref{11.27},  \eqref{11.28},  \eqref{11.29}, \eqref{vtech1}, \eqref{vtech2}, \eqref{11.22}, \eqref{11.23}, \eqref{11.24}, \eqref{11.25}	and the fact that $\tilde{r}_{i,u\sig}=O(1)\xi_i\bar{v}_i$, $\tilde{r}_{j,\sig}=O(1)\xi_j\bar{v}_j$
	\begin{align*}
		&\pa_{xx}(\tilde{r}_{i,u}\tilde{r}_j)
		=O(1)\biggl(\sum\limits_{k=1}|v_{k,x}|+\sum_{k,l} |v_k v_l|+\sum_k\rho_k^\e \mathfrak{A}_k|v_k\bar{v}_k\left(\frac{w_k}{v_k}\right)_x|
		\biggl)\\
		&+O(1)\rho_i^\e\mathfrak{A}_i\biggl( |v_i\left(\frac{w_i}{v_i}\right)_x| \sum_k|v_k|+\sum_{k\ne i}|v_kv_{k,xx}|
		+|v_i\left(\frac{w_i}{v_i}\right)_{xx}|
		+  |v_i\left(\frac{w_i}{v_i}\right)^2_x|\\
		&\hspace{2cm}+|v_{i,x}\left(\frac{w_i}{v_i}\right)_x|+ \sum_{k\ne i}|\left(\frac{w_i}{v_i}\right)_x||v_kv_{k,x}|+\rho_j^\e\mathfrak{A}_j|v_iv_{j,x} \left(\frac{w_i}{v_i}\right)_x|\\
		&\hspace{2cm}+|\left(\frac{w_i}{v_i}\right)_x v_i\left(\frac{w_j}{v_j}\right)_x|\rho_j^{\e}\mathfrak{A}_j+|v_{i,xx}|+\frac{v_{i,x}^2}{|v_i|}+\sum_{k\ne i}|\frac{v_{k,x}v_kv_{i,x}}{v_i}|\\
		&\hspace{2cm}+\sum\limits_{k=1}^{n}\abs{v_{i,x}v_k}+\rho_j^\e\mathfrak{A}_j|v_{j,x}v_{i,x}|+\abs{v_{i,x}v_j\left(\frac{w_j}{v_j}\right)_x}\rho_j^\e\mathfrak{A}_j\biggl)\\
		&
		+O(1)\rho_j^\e\mathfrak{A}_j \biggl(|v_j\left(\frac{w_j}{v_j}\right)_x| \sum_k|v_k|+\sum_{k\ne j}|v_kv_{k,xx}|+|v_j \left(\frac{w_j}{v_j}\right)_{xx}|+  |v_j\left(\frac{w_j}{v_j}\right)^2_x|+|v_{j,x}\left(\frac{w_j}{v_j}\right)_x|\\
		&\hspace{1.5cm}+|\left(\frac{w_j}{v_j}\right)_x|\sum_{k\ne j}|v_kv_{k,x}|+\rho_i^{\e}\mathfrak{A}_i |\left(\frac{w_j}{v_j}\right)_x v_j v_{i,x}|+\abs{\left(\frac{w_j}{v_j}\right)_x v_j\left(\frac{w_i}{v_i}\right)_x}\rho_i^\e\mathfrak{A}_i\\
		&\hspace{1cm}+|v_{j,xx}|+\frac{v_{j,x}^2}{|v_j|}
		+\sum_{k\ne j}|\frac{v_{k,x}v_kv_{j,x}}{v_j}|
		+\sum\limits_{k=1}^{n}\abs{v_{j,x}v_k}+\rho_i^\e\mathfrak{A}_i |v_{j,x}v_{i,x}|
		+\abs{v_{j,x}v_i\left(\frac{w_i}{v_i}\right)_x}\rho_i^\e\mathfrak{A}_i\biggl).
	\end{align*}
	
	\noi\underline{Estimate of $\tilde{\la}_{j,xx}$:} Taking the derivative in $x$ of $\tilde{\la}_{j,x}$ in \eqref{lajx} we have
	\begin{align*}
		\tilde{\la}_{j,xx}&=\sum_k v_{k,x}\tilde{r}_k\cdot \tilde{\la}_{j,u}+\sum_k v_{k}\tilde{r}_{k,x}\cdot \tilde{\la}_{j,u}+\sum_{k,l} v_{k}v_l \tilde{r}_k\otimes \tilde{r}_l: \tilde{\la}_{j,u u}\\
		&+\sum_k v_{k}\left(\frac{w_j}{v_j}\right)_x \tilde{r}_{k}\cdot (\xi'_j\bar{v}_j \tilde{\la}_{j,u v}-\theta'_j \tilde{\la}_{j,u \sig})+\sum_k v_{k} \xi_j \tilde{r}_{k}\cdot \tilde{\la}_{j,u v}(\pa_{v_j}\bar{v}_j v_{j,x}+\sum_{l\ne j}\pa_{v_l}\bar{v}_j v_{l,x})\\
		&+\left(\frac{w_j}{v_j}\right)_{xx}(\xi'_j\tilde{\la}_{j,v}\bar{v}_j-\tilde{\la}_{j,\sig}\theta'_j)
		+\left(\frac{w_j}{v_j}\right)_{x}^2(\xi''_j\tilde{\la}_{j,v}\bar{v}_j-\tilde{\la}_{j,\sig}\theta''_j)\\
		&
		+\left(\frac{w_j}{v_j}\right)_x(\xi'_j\pa_x(\tilde{\la}_{j,v})\bar{v}_j-\pa_x(\tilde{\la}_{j,\sig})\theta'_j)
		+\left(\frac{w_j}{v_j}\right)_x\xi'_j\tilde{\la}_{j,v}\pa_x(\bar{v}_j)+(\pa_x(\tilde{\la}_{j,v})\xi_j\\
		&+\tilde{\la}_{j,v}\xi'_j\left(\frac{w_j}{v_j}\right)_x)\left[\pa_{v_j}\bar{v}_jv_{j,x}+\sum_{k\ne j}(\pa_{v_k}\bar{v}_j v_{k,x})\right]\\
		&+\tilde{\la}_{j, v}\xi_j\biggl[\pa_x(\pa_{v_j}\bar{v}_j)v_{j,x}+\pa_{v_j}\bar{v}_j v_{j,xx}+\sum_{k\ne j}(\pa_{v_k}\bar{v}_j v_{k,xx})+\sum_{k\ne j}\pa_x(\pa_{v_k}\bar{v}_j )v_{k,x}\biggl].
	\end{align*}
	We deduce from \eqref{tildelamb} and Lemmas \ref{estimimpo1}, \ref{lemme6.5}, \ref{lemme6.6}, \ref{lemme11.3}
	\begin{align*}
		\tilde{\la}_{j,xx}&=O(1)\sum_k(|v_k|+|v_{k,x}|)+O(1)\mathfrak{A}_j \rho_j^\e \biggl(|\left(\frac{w_j}{v_j}\right)_{xx}v_j|+|v_j\left(\frac{w_j}{v_j}\right)_{x}^2|+|v_j \left(\frac{w_j}{v_j}\right)_{x}|\\
		&+| \left(\frac{w_j}{v_j}\right)_{x}v_{j,x}|+| \left(\frac{w_j}{v_j}\right)_{x}|\sum_{k\ne j}|v_kv_{k,x}|+\frac{v_{j,x}^2}{|v_j|}+|v_{j,x}|\sum_{k\ne j}|\frac{v_kv_{k,x}}{v_j}|+|v_{j,xx}|%+\sum_{k\ne j}|v_kv_{k,xx}|+\sum_{k\ne j,l\ne j}|v_{k,x}v_{l,x}|
		\biggl).
	\end{align*}
	
	This completes the proof of Lemma \ref{lemme11.6}.
\end{proof}

\begin{lemma}\label{lemme9.6}
	For $1\leq i,j\leq n$, it follows
\begin{align}
	&a_{ij,x},b_{ij,x}, \hat{b}_{ij,x}, \hat{a}_{ij,x}
	=O(1)\left[v_j^2+|v_j \left(\frac{w_j}{v_j}\right)_x|+|v_{j,x}|+\sum_{k\ne j}(\abs{v_kv_j}+|v_kv_{k,x}|)\right]\mathfrak{A}_j\rho_j^\e ,\label{aijx}\\
	&a_{ij,t}, b_{ij,t},\hat{b}_{ij,t}, \hat{a}_{ij,t}=O(1)\biggl[\sum_k\big(|w_kv_j|+|v_kv_j|\big)+|w_{j,xx}-v_{j,xx}\frac{w_j}{v_j}|\nonumber\\
	&+|v_j\left(\frac{w_j}{v_j}\right)_x|+|\tilde{\psi}_j|+|\tilde{\phi}_j|+|v_{j,xx}|+|v_{j,x}|+|v_j|+\sum_{k\ne j}|v_{k,t}v_k|\biggl]\mathfrak{A}_j\rho_j^\e ,\label{aijt}\\[2mm]
	&a_{ij,xx},\hat{a}_{ij,xx} =O(1)\rho_j^\e\mathfrak{A}_j \biggl(|v_j|+|v_j \left(\frac{w_j}{v_j}\right)_x|+\abs{v_j\left(\frac{w_j}{v_j}\right)_{x}^2}+\abs{v_j\left(\frac{w_j}{v_j}\right)_{xx}}+|v_{j,x} \left(\frac{w_j}{v_j}\right)_x| \nonumber\\
&|v_{j,x}|+\abs{v_{j,xx}}
+\sum\limits_{k} \abs{v_k w_{j,x}}+\sum\limits_{k\neq j}\abs{v_kv_{k,xx}}+\sum\limits_{l\neq j} \frac{\abs{v_l v_{l,x}v_{j,x}}}{\abs{v_j}}+\frac{v_{j,x}^2}{|v_j|} +\sum\limits_{k}\sum\limits_{l\neq j}\abs{v_kv_{l,x}} \nonumber\\
&+\sum\limits_{k\ne j,l\ne j}\abs{v_{k,x}v_{l,x}}+|\left(\frac{w_j}{v_j}\right)_x|\sum\limits_{k\neq j}\abs{v_{k,x}v_k}+\sum_{k\ne j} |v_kv_{k,x}|\biggl).\label{aijxx}
\end{align}
\end{lemma}
\begin{proof}Using the fact that $a_{i,j}=\mu_i b_{i,j}$ and from \eqref{def:b-ij} we have
\begin{align*}
	a_{ij,x}&=\sum_k v_k\mu_{i,u}\tilde{r}_k(\xi'_j\bar{v}_j\psi_{ji,v}-\psi_{ji,\sig})+ \mu_i \xi_j^{\p\p}\left(\frac{w_j}{v_j}\right)_x\bar{v}_j\psi_{ji,v}
	+\mu_i \xi_j^\p \bar{v}_{j,x}\psi_{ji,v}\nonumber\\
	&+\mu_i \xi_j^\p\bar{v}_j\pa_x( \psi_{ji,v})-\mu_i \pa_x( \psi_{ji,\sig}),\\
	a_{ij,t}&=\sum_k (w_k-\la_k^*v_k)\mu_{i,u}\tilde{r}_k(\xi'_j\bar{v}_j\psi_{ji,v}-\psi_{ji,\sig})+\mu_i \xi_j^{\p\p}\left(\frac{w_j}{v_j}\right)_t\bar{v}_j\psi_{ji,v}+\mu_i \xi_j^\p \bar{v}_{j,t}\psi_{ji,v}\nonumber\\
	&+\mu_i \xi_j^\p\bar{v}_j\pa_t (\psi_{ji,v})-\mu_i\pa_t(\psi_{ji,\sig}).
\end{align*}
%\begin{equation}
%	v_j \left(\frac{w_j}{v_j}\right)_t=\mu_j(w_{j,xx}-v_{j,xx}\frac{w_j}{v_j})+(\mu_{j,x}+\tilde{\la}_j)v_j\left(\frac{w_j}{v_j}\right)_x+\tilde{\psi}_j-\frac{w_j}{v_j}\tilde{\phi}_j,
%	\label{formuleimpo}
%\end{equation}
Using the Lemmas \ref{lemme9.2}, \ref{lemme11.3a} the fact that $\psi_{ji,\sig}=O(1) v_j\rho_j^\e\mathfrak{A}_j$ and  \eqref{formuleimpo} we deduce that
\begin{align}
	&a_{ij,x}=O(1)\left[\sum_k|v_jv_k|+|v_j \left(\frac{w_j}{v_j}\right)_x|+|v_{j,x}|+\sum_{k\ne j}|v_kv_{k,x}|\right] \rho_j^\e\mathfrak{A}_j ,\nonumber\\
	&a_{ij,t}=O(1)\biggl[\sum_k\big(|w_kv_j|+|v_kv_j|\big)+|w_{j,xx}-v_{j,xx}\frac{w_j}{v_j}|+|v_j\left(\frac{w_j}{v_j}\right)_x|\nonumber\\
	&+|\tilde{\psi}_j|+|\tilde{\phi}_j|+|v_{j,xx}|+|v_{j,x}|+|v_j|+\sum_{k\ne j}|v_{k,t}v_k|\biggl]\rho_j^\e\mathfrak{A}_j .\label{5aij}
\end{align}
	Similarly we obtain the result for $b_{ij,x}$, $b_{ij,t}$. Let us consider now $\hat{b}_{ij,x}$ and $\hat{b}_{ij,t}$, we have using the fact that $\hat{b}_{ij}=\psi_{ji}+[\la_i^*-\la_j^*+\frac{w_j}{v_j}]b_{ij}$ (see \eqref{def:hat-b-ij})
\begin{align*}
	\hat{b}_{ij,x}=&\pa_x(\psi_{ji})+[\la_i^*-\la_j^*+\frac{w_j}{v_j}]\pa_x(b_{ij})+\left(\frac{w_j}{v_j}\right)_x b_{ij},\\
	\hat{b}_{ij,t}=&\pa_t(\psi_{ji})+[\la_i^*-\la_j^*+\frac{w_j}{v_j}]\pa_t(b_{ij})+\left(\frac{w_j}{v_j}\right)_t b_{ij}.
	%\sum\limits_{k\neq i,j}\pa_x(\psi_{ki,v})(w_k-\la_k^*v_k)\xi_k\pa_{w_j}\bar{v}_i+\sum\limits_{k\neq i,j}\psi_{ki,v}(w_{k,x}-\la_k^*v_{k,x})\xi_k\pa_{w_j}\bar{v}_i\\
	%&+\sum\limits_{k\neq i,j}\psi_{ki,v}(w_k-\la_k^*v_k)\xi_k^\p\left(\frac{w_k}{v_k}\right)_x\pa_{w_j}\bar{v}_i+\sum\limits_{k\neq i,j}\psi_{ki,v}(w_k-\la_k^*v_k)\xi_k\pa_x(\pa_{w_j}\bar{v}_i)\\
	%&+\theta_j^{\p\p}\left(\frac{w_j}{v_j}\right)_x  \psi_{ji} +\theta_j^{\p}  \pa_x(\psi_{ji})+(\la_i^*-\la_j^*)\left(\frac{w_j}{v_j}\right)_x\big[\xi_j^{\p\p}\bar{v}_j\psi_{ji,v}-\theta^{\p\p}_j\psi_{ji,\si}\big] \\
	%&+(\la_i^*-\la_j^*)\big[\xi_j^\p\pa_x(\bar{v}_j)\psi_{ji,v}+\xi_j^\p\bar{v}_j\pa_x(\psi_{ji,v})-\theta^\p_j\pa_x(\psi_{ji,\si})\big]\\
	%&+\sum\limits_{k\neq i,j} \la_i^*\pa_x(\psi_{ki,v})  v_k\xi_k\pa_{w_j}\bar{v}_i+\sum\limits_{k\neq i,j} \la_i^*\psi_{ki,v}  v_{k,x}\xi_k\pa_{w_j}\bar{v}_i\\
	%&+\sum\limits_{k\neq i,j} \la_i^*\psi_{ki,v}  v_k\xi_k^\p\left(\frac{w_k}{v_k}\right)\pa_{w_j}\bar{v}_i+\sum\limits_{k\neq i,j} \la_i^*\psi_{ki,v}  v_k\xi_k\pa_x(\pa_{w_j}\bar{v}_i). 
\end{align*}
We deduce now the result by using Lemma \ref{lemme11.3a} and \eqref{5aij},  \eqref{formuleimpo}. We proceed similarly for $\hat{a}_{ij,x}$ and $\hat{a}_{ij,t}$ using the fact that $\hat{a}_{ij}=\mu_i\hat{b}_{ij}$.
We observe now from \eqref{aijx} that
\begin{align*}
	&a_{ij,xx}=\sum_k \big[(v_{k,x}\mu_{i,u}+\sum_l v_{k}v_l\tilde{r}_l\otimes \tilde{r}_k: \mu_{i,uu})\tilde{r}_k+ v_{k}\mu_{i,u}\pa_x(\tilde{r}_k)\big]     (\xi'_j\bar{v}_j\psi_{ji,v}-\psi_{ji,\sig})\\%+\sum_{k,l} v_{k}v_l\tilde{r}_l\otimes \tilde{r}_k: \mu_{i,uu}(\xi'_j\bar{v}_j\psi_{ji,v}
	%-\psi_{ji,\sig})\\
	%+\sum_k v_{k}\mu_{i,u}\pa_x(\tilde{r}_k)(\xi'_j\bar{v}_j\psi_{ji,v}-\psi_{ji,\sig})
	&+\sum_k v_{k}\mu_{i,u}\tilde{r}_k\big[(\xi''_j\left(\frac{w_j}{v_j}\right)_x\bar{v}_j+\xi'_j \bar{v}_{j,x})\psi_{ji,v}+\xi'_j \bar{v}_j\pa_x(\psi_{ji,v})-\pa_x(\psi_{ji,\sig})\big]+\mu_{i,x} \xi_j^{\p\p}\left(\frac{w_j}{v_j}\right)_x\bar{v}_j\psi_{ji,v}\\
	&+\mu_{i} \big[(\xi'''_j\left(\frac{w_j}{v_j}\right)^2_x+\xi''_j\left(\frac{w_j}{v_j}\right)_{xx}) \bar{v}_j\psi_{ji,v}+2\xi''_j\left(\frac{w_j}{v_j}\right)_{x}(\bar{v}_{jx}\psi_{ji,v}+\bar{v}_j\pa_x(\psi_{ji,v})) \big]\\
	%&-\mu_{i} \xi''_j\left(\frac{w_j}{v_j}\right)_{xx}\bar{v}_j\psi_{ji,v}-
	%2\mu_{i} \xi''_j\left(\frac{w_j}{v_j}\right)_{x}\bar{v}_{jx}\psi_{ji,v}-2\mu_{i} \xi''_j\left(\frac{w_j}{v_j}\right)_{x}\bar{v}_j\pa_x(\psi_{ji,v})\\
	&+(\mu_{i,x}\bar{v}_{j,x}+\mu_{i} \bar{v}_{j,xx})    \xi_j^\p \psi_{ji,v}
	+2\mu_{i} \xi_j^\p \bar{v}_{j,x}\pa_x(\psi_{ji,v})+(\mu_{i,x}\pa_x( \psi_{ji,v})+\mu_{i} \pa_{xx}( \psi_{ji,v}))  \xi_j^\p\bar{v}_j\\\
	&+\mu_{i,x} \pa_x( \psi_{ji,\sig})+\mu_{i} \pa_{xx}( \psi_{ji,\sig}).		%&+\sum\limits_{k\neq i,j}\pa_x(\psi_{ki,v})  v_k\xi_k\pa_{w_j}\bar{v}_i+\sum\limits_{k\neq i,j}\psi_{ki,v}  v_{k,x}\xi_k\pa_{w_j}\bar{v}_i\\
	%&+\sum\limits_{k\neq i,j}\psi_{ki,v}  v_k\xi_k^\p\left(\frac{w_k}{v_k}\right)_x\pa_{w_j}\bar{v}_i+\sum\limits_{k\neq i,j}\psi_{ki,v}  v_k\xi_k\pa_x(\pa_{w_j}\bar{v}_i).
\end{align*}
We have now using the Lemmas \ref{estimimpo1}, \ref{lemme6.5}, \ref{lemme6.6}, \ref{lemme9.2}, \ref{lemme11.3}, \ref{lemme11.3a}, \ref{lemme11.5}
and the fact that $\bar{v}_j\pa_x(\psi_{ji,v})=O(1)\bar{v}_j+O(1)\bar{v}_j\left(\frac{w_j}{v_j}\right)_x$, $\psi_{ji,\sig}=O(1)\mathfrak{A}_j\bar{v}_j$
 \begin{align}
& a_{ij,xx}=O(1)\rho_j^\e\mathfrak{A}_j \biggl(|v_j|+|v_j \left(\frac{w_j}{v_j}\right)_x|+\abs{v_j\left(\frac{w_j}{v_j}\right)_{x}^2}+\abs{v_j\left(\frac{w_j}{v_j}\right)_{xx}}+|v_{j,x} \left(\frac{w_j}{v_j}\right)_x|\nonumber\\
&|v_{j,x}|+\abs{v_{j,xx}}
+\sum\limits_{k} \abs{v_k w_{j,x}}+\sum\limits_{k\neq j}\abs{v_kv_{k,xx}}+\sum\limits_{l\neq j} \frac{\abs{v_l v_{l,x}v_{j,x}}}{\abs{v_j}}+\frac{v_{j,x}^2}{|v_j|} +\sum\limits_{k}\sum\limits_{l\neq j}\abs{v_kv_{l,x}}\nonumber\\
&+\sum\limits_{k\ne j,l\ne j}\abs{v_{k,x}v_{l,x}}+|\left(\frac{w_j}{v_j}\right)_x|\sum\limits_{k\neq j}\abs{v_{k,x}v_k}+\sum_{k\ne j} |v_kv_{k,x}|\biggl).\label{aijxxb}
 \end{align}
%\begin{align*}
%	a_{ij,x}&=\sum_k v_k\mu_{i,u}\tilde{r}_k(\xi'_j\bar{v}_j\psi_{ji,v}-\psi_{ji,\sig})+ \mu_i \xi_j^{\p\p}\left(\frac{w_j}{v_j}\right)_x\bar{v}_j\psi_{ji,v}
%	+\mu_i \xi_j^\p \bar{v}_{j,x}\psi_{ji,v}\nonumber\\
	%&+\mu_i \xi_j^\p\bar{v}_j\pa_x( \psi_{ji,v})-\mu_i \pa_x( \psi_{ji,\sig}).
%\end{align*}
We can now estimate $\hat{a}_{ij,xx}$, using the fact that $\hat{a}_{ij}=\mu_i \hat{b}_{ij}=\mu_i\psi_{ji}+[\la_i^*-\la_j^*+\frac{w_j}{v_j}]a_{ij}$, we have then:
\begin{align*}
	&\hat{a}_{ij,xx}=\mu_{i,xx}\psi_{ji}+2\mu_{i,x}\psi_{ji,x}+\mu_i\psi_{ji,xx}+[\la_i^*-\la_j^*+\frac{w_j}{v_j}]a_{ij,xx}+2\left(\frac{w_j}{v_j}\right)_x a_{ij,x}\\
	&+\left(\frac{w_j}{v_j}\right)_{xx}a_{ij}.
\end{align*}
From \eqref{aijx}, \eqref{aijxxb}, Lemmas \ref{lemme11.3a}, \ref{lemme11.5} and the fact that $a_{ij}=O(1)\mathfrak{A}_j \bar{v}_j$, $\mu_{i,xx}, \mu_{i,x}=O(1)$ we obtain the desired estimate \eqref{aijxx} for $\hat{a}_{ij,xx}$.
\end{proof}

\subsection{Derivative of directional error terms}
\begin{lemma}\label{lemme9.7}
	The following estimate is true.
	\begin{align}
		&\left(\mu_i v_{i,x}-(\tilde{\la}_i-\la_i^*)v_{i}-w_{i} \right)_x\nonumber\\
		&=O(1)\sum\limits_{j\neq i}\left(\mu_j v_{j,x}-(\tilde{\la}_j-\la_j^*)v_j-w_j\right)_x v_j\xi_j	\nonumber\\
		&+O(1)\sum_{j\ne i}(|v_j|+|v_{j,x}|+|w_{j,x}|)(|v_j|+|w_{j,x}|+|v_{j,x}|+\sum_{k} |v_k||v_j|)\nonumber\\
		&+O(1)\left(\sum\limits_{j\neq i}\sum\limits_{k\neq j} (|w_{j,x}|+|v_{j,x}| +|v_j|)|v_{k,x} v_k|+\sum\limits_{j\neq i}\sum\limits_{k\neq j} |v_j v_{k,xx}v_k|\right)\nonumber\\
		%&+O(1)\left(\sum\limits_{j\neq i}\sum\limits_{k\neq j} |v_j w_{k,xx}w_k|+\sum\limits_{j\neq i}\sum\limits_{k\neq j}\sum_{l\ne j} |v_j |(|v_{k,x}|+|w_{k,x}|)|(|v_{l,x}|+|w_{l,x}|)\right)\nonumber\\
		&+O(1)\left(\sum_{i\ne j}|w_{j,xx}v_j-v_{j,xx}w_j|
		+\sum_k|v_kv_{k,x}||1-\chi_k^\e\xi_k\eta_k|\right)\nonumber\\
		&+O(1)\left(\sum_k v_k^2\abs{\left(\frac{w_k}{v_k}\right)_x}+\sum_k \sum_{l\ne k}|v_l|(|v_l||v_{k,x}|+|v_kv_{l,x}|) 
		+O(1)\sum\limits_{k}\sum\limits_{l\ne k}v_k^2|v_l| \right)\nonumber\\
		&+O(1)\sum\limits_{k}|v_k|^3(1-\chi_k^\e\xi_k\eta_k)+O(1)\sum_k \varkappa^\ep_k\rho_k^\e\mathfrak{A}_k|v_k v_{k,x}|\nonumber\\
		&+O(1)\left(\sum_{k\ne j}(|v_{k,x}v_j|+|v_k v_{j,x}|+|v_kv_j w_{k,x}|+|v_kv_jv_{j,x}|+|v_kv_jw_{j,x}|)+\sum_{j\ne k}\sum_l|v_kv_jv_l|
		\right).
		\label{paraestim1}
	\end{align}
	It implies in particular that
	\begin{align}
		&\left(\mu_i v_{i,x}-(\tilde{\la}_i-\la_i^*)v_{i}-w_{i} \right)_x\nonumber\\
		&=O(1)\sum\limits_{j\neq i}\left(\mu_j v_{j,x}-(\tilde{\la}_j-\la_j^*)v_j-w_j\right)_x \bar{v}_j\xi_j	\nonumber\\
		&+O(1)\sum_{j\ne i}\Big(|v_j|+|v_{j,x}|+|w_{j,x}|\Big)\left(|v_j|+|w_{j,x}|+|v_{j,x}|+\sum_{k} |v_k||v_j|\right)\nonumber\\
		&+O(1)\sum_k(\La_k^1+\La_k^4+\La_k^5+\La_k^6)+
		O(1)\big(\sum_k|v_kv_{k,x}|\eta_k(1-\chi^\e_k)+\sum\limits_{k}|v_k|^3\eta_k (1-\chi^\e_k)\big)
		\nonumber\\
		&+O(1)\sum_k \varkappa^\ep_k\rho_k^\e\mathfrak{A}_k|v_k v_{k,x}|.
	\label{estimate-v-i-xx-1aa}
	\end{align}
	Furthermore, we have
	\begin{align}
		&v_{i,xx}=O(1)(|v_{i,x}|+|w_{i,x}|+|v_i|)+O(1)\sum\limits_{j\neq i}\left(\mu_j v_{j,x}-(\tilde{\la}_j-\la_j^*)v_j-w_j\right)_x \bar{v}_j\xi_j	\nonumber\\
		&+O(1)\sum_{j\ne i}\Big(|v_j|+|v_{j,x}|+|w_{j,x}|\Big)\left(|v_j|+|w_{j,x}|+|v_{j,x}|+\sum_{k} |v_k||v_j|\right)\nonumber\\
		&+O(1)\sum_k(\La_k^1+\La_k^4+\La_k^5+\La_k^6)+
		O(1)\big(\sum_k|v_kv_{k,x}|\eta_k(1-\chi^\e_k)+\sum\limits_{k}|v_k|^3\eta_k (1-\chi^\e_k)\big)
		\nonumber\\
		&+O(1)\sum_k \varkappa^\ep_k\rho_k^\e\mathfrak{A}_k|v_k v_{k,x}|.
	\label{estimate-v-i-xx-1aabis}
	\end{align}
\end{lemma}
\begin{proof}
	Taking derivative of \eqref{eqn-v-i-x-1} we get
	\begin{equation}
		\mu_i v_{i,xx}+\mu_{i,x}v_{i,x}-\tilde{\la}_{i,x}v_i-(\tilde{\la}_i-\la_i^*)v_{i,x}-w_{i,x}=\sum\limits_{s=2}^{7}\mathcal{I}_s,
		\label{9.14}
	\end{equation}
	where $\mathcal{I}_s$ are defined as follows% using the fact that $\theta_j^\p\psi_{ji,\si}=\psi_{ji,\si}$,
	\begin{align*}
		\mathcal{I}_2&=-\pa_x\left(\sum\limits_{j\neq i}(\mu_j v_{j,x}-(\tilde{\la}_j-\la_j^*)v_j-w_j)\frac{\mu_i}{\mu_j}[\psi_{ji} +v_j\xi_j\pa_{v_j}\bar{v}_j \psi_{ji,v}]\right),\\
		\mathcal{I}_3&=-\pa_x\left(\sum\limits_{j\neq i}\sum\limits_{k\neq j}\xi_jv_{k,x}v_j\pa_{v_k}\bar{v}_j\mu_i \psi_{ji,v}\right),\\
		\mathcal{I}_4&=-\pa_x\left(\sum\limits_{j\neq i}v_j\left(\frac{w_j}{v_j}\right)_x\mu_i [-\psi_{ji,\si}+\bar{v}_j\xi_j^\p\psi_{ji,v}]\right),\\
		\mathcal{I}_5&=-\pa_x\left(\sum\limits_{j\neq i}\mu_j^{-1}(\tilde{\la}_j-\si_j)v_j\xi_j\tilde{v}_j\mu_i\psi_{ji,v}\right),\\
		\mathcal{I}_6&=-\pa_x\left(\sum\limits_{k}v_k^2(1-\chi_\e^k\xi_k\eta_k)\langle l_i, \left[B(u)\tilde{r}_{k,u}\tilde{r}_k+\tilde{r}_k\cdot DB(u)\tilde{r}_k\right]\rangle \right),\\
		\mathcal{I}_7&=-\pa_x\left(\sum\limits_{k\neq j}v_kv_j\langle l_i,\tilde{r}_k\cdot DB(u)\tilde{r}_j+B(u)\tilde{r}_{k,u}\tilde{r}_j\rangle\right).
	\end{align*}
	For convenience, we set $\mathcal{I}_1:=\mu_i v_{i,xx}+\mu_{i,x}v_{i,x}-\tilde{\la}_{i,x}v_i-(\tilde{\la}_i-\la_i^*)v_{i,x}-w_{i,x}$. Next we estimate each $\mathcal{I}_s$ for $2\leq s\leq 7$. By using \eqref{eqn-v-i-x-1} we get
	\begin{align*}
		\mathcal{I}_2&=-\left(\sum\limits_{j\neq i}\left(\mu_j v_{j,x}-(\tilde{\la}_j-\la_j^*)v_j-w_j\right)_x\frac{\mu_i}{\mu_j}[\psi_{ji} +v_j\xi_j\pa_{v_j}\bar{v}_j \psi_{ji,v}]\right)\\
		&-\left(\sum\limits_{j\neq i}\mathcal{H}_j \pa_x(\frac{\mu_i}{\mu_j})[\psi_{ji} +v_j\xi_j\pa_{v_j}\bar{v}_j \psi_{ji,v}]\right)\\
		&-\left(\sum\limits_{j\neq i} (\mu_j v_{j,x}-(\tilde{\la}_j-\la_j^*)v_j-w_j)\frac{\mu_i}{\mu_j}[\pa_x \psi_{ji} +v_{j,x}\xi_j\pa_{v_j}\bar{v}_j \psi_{ji,v}+v_{j}\xi_j^\p\left(\frac{w_j}{v_j}\right)_x \pa_{v_j}\bar{v}_j \psi_{ji,v}]\right)\\
		&-\left(\sum\limits_{j\neq i} \mathcal{H}_j \frac{\mu_i}{\mu_j}[v_{j}\xi_j\pa_x(\pa_{v_j}\bar{v}_j) \psi_{ji,v}+v_{j}\xi_j\pa_{v_j}\bar{v}_j \pa_x(\psi_{ji,v})]\right)\\
		&=:\sum\limits_{l=1}^{4}\mathcal{I}_{2,l}.%\nonumber\\
	\end{align*}
	Using the Lemma \ref{estimimpo1}, we have $\pa_x(\pa_{v_j}\bar{v}_j)=O(1)\rho_j^\e\big(\sum_{k\ne j}|\frac{v_k}{v_j}v_{k,x}|+|\frac{v_{j,x}}{v_j}|\big)$ which gives by applying \eqref{superimpo} and Lemmas \ref{estimimpo1}, \ref{lemme11.3a} and the fact that $|w_j|\mathfrak{A}_j=O(1)|v_j|\mathfrak{A}_j$
	\begin{align}
		\mathcal{I}_2&=O(1)\sum\limits_{j\neq i}\left(\mu_j v_{j,x}-(\tilde{\la}_j-\la_j^*)v_j-w_j\right)_x\bar{v_j}\xi_j\nonumber\\
		&+O(1)\sum_{j\ne i}\sum_{k\ne j}|v_k v_{k,x}|(|v_j|+|v_{j,x}|) \mathfrak{A}_j\rho_j^\e\nonumber\\
		&+O(1)\sum_{j\ne i}(|v_j|+|v_{j,x}|)(|w_{j,x}|+|v_{j,x}|+\sum_{k} |v_k||v_j|)\mathfrak{A}_j\rho_j^\e\label{I2a}
	\end{align}
	Similarly, we obtain using again Lemmas \ref{lemme11.3a}, \ref{estimimpo1}
	\begin{align}
		\mathcal{I}_3&=-\sum\limits_{j\neq i}\sum\limits_{k\neq j}\xi_j^\p\left(\frac{w_j}{v_j}\right)_xv_{k,x}v_j\pa_{v_k}\bar{v}_j\mu_i \psi_{ji,v}-\sum\limits_{j\neq i}\sum\limits_{k\neq j}\xi_jv_{k,xx}v_j\pa_{v_k}\bar{v}_j\mu_i \psi_{ji,v} \nonumber\\
		&-\sum\limits_{j\neq i}\sum\limits_{k\neq j}\xi_jv_{k,x}v_{j,x}\pa_{v_k}\bar{v}_j\mu_i \psi_{ji,v}-\sum\limits_{j\neq i}\sum\limits_{k\neq j}\xi_jv_{k,x}v_j\pa_{v_k}\bar{v}_j\mu_{i,x} \psi_{ji,v}\nonumber\\
		&-\sum\limits_{j\neq i}\sum\limits_{k\neq j}\xi_jv_{k,x}v_jv_{j,x}\pa_{v_kv_j}\bar{v}_j\mu_i \psi_{ji,v}-\sum\limits_{j\neq i}\sum\limits_{k\neq j}\sum\limits_{l\neq j}\xi_jv_{k,x}v_jv_{l,x}\pa_{v_lv_k}\bar{v}_j\mu_i \psi_{ji,v}\nonumber\\
		&%-\sum\limits_{j\neq i}\sum\limits_{k\neq j}\sum\limits_{l\neq j}\xi_jv_{k,x}v_jw_{l,x}\pa_{w_lv_k}\bar{v}_j\mu_i \psi_{ji,v}
		-\sum\limits_{j\neq i}\sum\limits_{k\neq j}\xi_jv_{k,x}v_j\pa_{v_k}\bar{v}_j\mu_i \pa_x(\psi_{ji,v})\nonumber\\
		%&-\sum\limits_{j\neq i}\sum\limits_{k\neq j}\xi^\p_j\left(\frac{w_j}{v_j}\right)_xw_{k,x}v_j\pa_{w_k}\bar{v}_j\mu_i \psi_{ji,v}-\sum\limits_{j\neq i}\sum\limits_{k\neq j}\xi_jw_{k,xx}v_j\pa_{w_k}\bar{v}_j\mu_i \psi_{ji,v}\nonumber\\
		%&-\sum\limits_{j\neq i}\sum\limits_{k\neq j}\xi_jw_{k,x}v_{j,x}\pa_{w_k}\bar{v}_j\mu_i \psi_{ji,v}-\sum\limits_{j\neq i}\sum\limits_{k\neq j}\xi_jw_{k,x}v_jv_{j,x}\pa_{w_kv_j}\bar{v}_j\mu_i \psi_{ji,v}\nonumber\\
		%&-\sum\limits_{j\neq i}\sum\limits_{k\neq j}\sum\limits_{l\neq j}\xi_jw_{k,x}v_jv_{l,x}\pa_{v_lw_k}\bar{v}_j\mu_i \psi_{ji,v}-\sum\limits_{j\neq i}\sum\limits_{k\neq j}\sum\limits_{l\neq j}\xi_jw_{k,x}w_{l,x}v_j\pa_{w_kw_l}\bar{v}_j\mu_i \psi_{ji,v}\nonumber\\
		%&-\sum\limits_{j\neq i}\sum\limits_{k\neq j}\xi_jw_{k,x}v_j\pa_{w_k}\bar{v}_j\mu_{i,x} \psi_{ji,v}-\sum\limits_{j\neq i}\sum\limits_{k\neq j}\xi_jw_{k,x}v_j\pa_{w_k}\bar{v}_j\mu_i \pa_x(\psi_{ji,v})\nonumber\\
		&=:\sum\limits_{l=1}^{7}\mathcal{I}_{3,l}\nonumber\\
		&=O(1)\left(\sum\limits_{j\neq i}\sum\limits_{k\neq j} (|w_{j,x}|+|v_{j,x}| +|v_j|)|v_{k,x} v_k|\mathfrak{A}_j\rho_j^\e
		+\sum\limits_{j\neq i}\sum\limits_{k\neq j} |v_j v_{k,xx}v_k|\mathfrak{A}_j\rho_j^\e\right)%\nonumber\\
		%&+O(1)\left(\sum\limits_{j\neq i}\sum\limits_{k\neq j} |v_j w_{k,xx}w_k|+\sum\limits_{j\neq i}\sum\limits_{k\neq j}\sum_{l\ne j} |v_j |(|v_{k,x}|+|w_{k,x}|)|(|v_{l,x}|+|w_{l,x}|)\right)
		.\label{I3a}
	\end{align}
	Next we calculate $	\mathcal{I}_4$ 
	\begin{align*}
		\mathcal{I}_4
		&=-\sum\limits_{j\neq i}v_{j,x}\left(\frac{w_j}{v_j}\right)_x\mu_i [-\psi_{ji,\si}+\bar{v}_j\xi_j^\p\psi_{ji,v}]-\sum\limits_{j\neq i}v_j\left(\frac{w_j}{v_j}\right)_{xx}\mu_i [-\psi_{ji,\si}+\bar{v}_j\xi_j^\p\psi_{ji,v}]\nonumber\\
		&-\sum\limits_{j\neq i}v_j\left(\frac{w_j}{v_j}\right)^2_x\mu_i [\bar{v}_j\xi_j^{\p\p}\psi_{ji,v}]-\sum\limits_{j\neq i}v_j\left(\frac{w_j}{v_j}\right)_x\mu_i v_{j,x}\pa_{v_j}\bar{v}_j\xi_j^\p\psi_{ji,v}\nonumber\\
		&-\sum\limits_{j\neq i}v_j\left(\frac{w_j}{v_j}\right)_x\mu_i\xi_j^\p\psi_{ji,v}\sum\limits_{l\neq j}v_{l,x}\pa_{v_l}\bar{v}_j-\sum\limits_{j\neq i}v_j\left(\frac{w_j}{v_j}\right)_x\mu_{i,x} [-\psi_{ji,\si}+\bar{v}_j\xi_j^\p\psi_{ji,v}]\nonumber\\
		&-\sum\limits_{j\neq i}v_j\left(\frac{w_j}{v_j}\right)_x\mu_i [-\pa_x(\psi_{ji,\si})+\bar{v}_j\xi_j^\p\pa_x(\psi_{ji,v})]\nonumber\\
		&=:\sum\limits_{l=1}^{7}\mathcal{I}_{4,l}.
	\end{align*}
	Note that
	\begin{equation*}
		v_j\left(\frac{w_j}{v_j}\right)_{xx}=\frac{1}{v_j}(w_{j,xx}v_j-w_jv_{j,xx})-\frac{2v_{j,x}}{v_j}\left(w_{j,x}-\frac{w_j}{v_j}v_{j,x}\right).
	\end{equation*}
	Hence, using again Lemmas \ref{estimimpo1} and \ref{lemme11.3a} we get
	\begin{align}
		\mathcal{I}_4&=O(1)\sum_{j\ne i}(|v_{j,x}|+|w_{j,x}|)^2\mathfrak{A}_j\rho_j^\e
		+O(1)\sum_{j\ne i}\sum_{j\ne k}(|w_{j,x}|+|v_{j,x}|)|v_{k,x}v_k||\mathfrak{A}_j\rho_j^\e\nonumber\\
		&+O(1)\sum_{j\ne i}\sum_{k}(|w_{j,x}|+|v_{j,x}|)|v_jv_k|\mathfrak{A}_j\rho_j^\e+\sum_{j\ne i}|w_{j,xx}v_j-v_{j,xx}w_j|\mathfrak{A}_j\rho_j^\e.\label{I4a}
	\end{align}
	First we observe that
	 $$\tilde{v}_j=2N\frac{v_j^{2N}}{\e}\chi'(\frac{v_j^{2N}}{\e})\hat{v}_j-\chi_j^\e \sum_{l\ne j}\prod_{k\ne l,j}\eta\left(\frac{v_k^2}{v_j}\right)\eta^\p\left(\frac{v_l^2}{v_j}\right)v_l^2.$$ We deduce from Lemma \ref{lemme9.2} that $\tilde{v}_j=O(1)v_j\rho_j^\e$ and $\pa_{v_j}\tilde{v_j}=O(1)\rho_j^\e$ whereas  $\pa_{v_k}\tilde{v_j}=O(1)|v_k|\rho_j^\e$  if $k\ne j$. Using Lemma \ref{lemme11.3a}
	we can estimate $I_5$ as follows
	\begin{align}
		\mathcal{I}_5
		&=-\sum\limits_{j\neq i}\pa_x(\mu_j^{-1})(\tilde{\la}_j-\si_j)v_j\xi_j\tilde{v}_j\mu_i\psi_{ji,v}-\sum\limits_{j\neq i}\mu_j^{-1}\left(\tilde{\la}_{j,x}+\theta_j^\p\left(\frac{w_j}{v_j}\right)_x\right)v_j\xi_j\tilde{v}_j\mu_i\psi_{ji,v}\nonumber\\
		&-\sum\limits_{j\neq i}\mu_j^{-1}(\tilde{\la}_j-\si_j)v_{j,x}\xi_j \tilde{v}_j\mu_i\psi_{ji,v}-\sum\limits_{j\neq i}\mu_j^{-1}(\tilde{\la}_j-\si_j)v_j\xi_j^\p\left(\frac{w_j}{v_j}\right)_x\tilde{v}_j\mu_i\psi_{ji,v}\nonumber\\
		&-\sum\limits_{j\neq i}\mu_j^{-1}(\tilde{\la}_j-\si_j)v_j\xi_jv_{j,x}\pa_{v_j} \tilde{v}_j\mu_i\psi_{ji,v}-\sum\limits_{j\neq i}\mu_j^{-1}(\tilde{\la}_j-\si_j)v_j\xi_j\sum\limits_{l\neq j}v_{l,x}\pa_{v_l} \tilde{v}_j\mu_i\psi_{ji,v}\nonumber\\
		%&-\sum\limits_{j\neq i}\mu_j^{-1}(\tilde{\la}_j-\si_j)v_j\xi_j\sum\limits_{l\neq j}w_{l,x}\pa_{w_l} \tilde{v}_j\mu_i\psi_{ji,v}\nonumber\\
		&-\sum\limits_{j\neq i}\mu_j^{-1}(\tilde{\la}_j-\si_j)v_j\xi_j\tilde{v}_j[\mu_{i,x}\psi_{ji,v}+\mu_i\pa_x\psi_{ji,v}]\nonumber\\
		&=:\sum\limits_{l=1}^{7}\mathcal{I}_{5,l}\nonumber\\
		&=O(1)\left(\sum\limits_{j\neq i}|v_j|^2\mathfrak{A}_j\rho_j^\e		
		+\sum\limits_{j\neq i}(\abs{v_jw_{j,x}}+\abs{v_jv_{j,x}})\mathfrak{A}_j\rho_j^\e+\sum\limits_{j\neq i}\sum\limits_{k\neq j}|v_j||v_{k,x}v_k|\mathfrak{A}_j\rho_j^\e\right).\label{I5a}
	\end{align}
	We observe that $\chi_k^\e v_k^2\pa_{v_k}\eta_k=O(1)\sum_{l\ne k}v_l^2 \rho_k^\e$ and  $\chi_k^\e v_k^2\pa_{v_l}\eta_k=O(1)|v_lv_k|\rho_k^\e$. Using now Lemmas \ref{lemme11.3}, \ref{lemme11.4}, \ref{lemme6.5}, \ref{lemme6.6} we get
	\begin{align}
		\mathcal{I}_6
		&=-\sum\limits_{k}2v_kv_{k,x}(1-\chi_k^\e\xi_k\eta_k)\langle l_i, \left[B(u)\tilde{r}_{k,u}\tilde{r}_k+\tilde{r}_k\cdot DB(u)\tilde{r}_k\right]\rangle\nonumber \\
		&+\sum\limits_{k}\frac{2N v_k^{2N+1}}{\e}v_{k,x} (\chi')_k^\e \xi_k\eta_k\langle l_i, \left[B(u)\tilde{r}_{k,u}\tilde{r}_k+\tilde{r}_k\cdot DB(u)\tilde{r}_k\right]\rangle \\
		&+\sum\limits_{k}v_k^2\chi_k^\e \xi^\p_k\eta_k\left(\frac{w_k}{v_k}\right)_x\langle l_i, \left[B(u)\tilde{r}_{k,u}\tilde{r}_k+\tilde{r}_k\cdot DB(u)\tilde{r}_k\right]\rangle \nonumber\\
		&+\sum\limits_{k}v_k^2\chi_k^\e \xi_k\left[v_{k,x}\pa_{v_k}\eta_k+\sum\limits_{l\neq k}v_{l,x}\pa_{v_l}\eta_k\right]\langle l_i, \left[B(u)\tilde{r}_{k,u}\tilde{r}_k+\tilde{r}_k\cdot DB(u)\tilde{r}_k\right]\rangle \nonumber\\
		&-\sum\limits_{k}\sum\limits_{l}v_k^2v_l(1-\chi_k^\e\xi_k\eta_k)\langle l_i, \left[\tilde{r}_l\cdot DB\tilde{r}_{k,u}\tilde{r}_k+\tilde{r}_k\otimes \tilde{r}_l: D^2B\tilde{r}_k\right]\rangle \nonumber\\
		&-\sum\limits_{k}\sum\limits_{l}v_k^2v_l(1-\chi_k^\e\xi_k\eta_k)\langle \tilde{r}_l\cdot Dl_i, \left[B(u)\tilde{r}_{k,u}\tilde{r}_k+\tilde{r}_k\cdot DB(u)\tilde{r}_k\right]\rangle\nonumber\\
		&-\sum\limits_{k}v_k^2(1-\chi_k^\e\xi_k\eta_k)\langle l_i, \left[B(u)\pa_x(\tilde{r}_{k,u}\tilde{r}_{k})+\tilde{r}_k\cdot DB(u)\tilde{r}_{k,x}+\tilde{r}_{k,x}\cdot DB(u)\tilde{r}_k\right]\rangle \nonumber\\
		&=:\sum\limits_{l=1}^{7}\mathcal{I}_{6,l}\nonumber\\
		&=O(1)\left(\sum_k |v_kv_{k,x}|(1-\chi_k^\e\xi_k\eta_k)+\sum_k v_k^2|\left(\frac{w_k}{v_k}\right)_x|\mathfrak{A}_k\rho_k^\e \right)\nonumber\\
		&+O(1)\sum_k \sum_{l\ne k}|v_l|(|v_l||v_{k,x}|+|v_kv_{l,x}|) \rho_k^\e\mathfrak{A}_k
		+O(1)\sum\limits_{k}\sum\limits_{l\ne k}v_k^2|v_l|\nonumber\\
		&+O(1)\sum\limits_{k}|v_k|^3(1-\chi_k^\e\xi_k\eta_k)+O(1)\sum_k \varkappa^\ep_k\rho_k^\e\mathfrak{A}_k|v_kv_{k,x}|.
		\label{I6a}
	\end{align}
	In a similar way we have using Lemmas \ref{lemme11.3}, \ref{lemme11.4}
	\begin{align}
		\mathcal{I}_7
		&=-\sum\limits_{k\neq j}(v_{k,x}v_j+v_kv_{j,x})\langle l_i,\tilde{r}_k\cdot DB(u)\tilde{r}_j+B(u)\tilde{r}_{k,u}\tilde{r}_j\rangle\nonumber\\
		&-\sum\limits_{k\neq j}\sum\limits_{l}v_kv_jv_l\langle \tilde{r}_l\cdot Dl_i,\tilde{r}_k\cdot DB(u)\tilde{r}_j+B(u)\tilde{r}_{k,u}\tilde{r}_j\rangle\nonumber\\
		&-\sum\limits_{k\neq j}\sum\limits_{l}v_kv_jv_l\langle l_i,\tilde{r}_k\otimes\tilde{r}_l: D^2B\tilde{r}_j+\tilde{r}_l\cdot DB\tilde{r}_{k,u}\tilde{r}_j\rangle\nonumber\\
		&-\sum\limits_{k\neq j}v_kv_j\langle l_i,\tilde{r}_{k,x}\cdot DB\tilde{r}_j+\tilde{r}_k\cdot DB\tilde{r}_{j,x}+B(u)\pa_x(\tilde{r}_{k,u}\tilde{r}_j)\rangle\nonumber\\
		&=:\sum\limits_{l=1}^{4}\mathcal{I}_{7,l}\nonumber\\
		&=O(1)\left(\sum_{k\ne j}(|v_{k,x}v_j|+|v_k v_{j,x}|+|v_kv_j w_{k,x}|+|v_kv_jv_{j,x}|+|v_kv_jw_{j,x}|)+\sum_{j\ne k}\sum_l|v_kv_jv_l|
		\right).\label{I7a}
	\end{align}
	Combining \eqref{I2a}, \eqref{I3a}, \eqref{I4a}, \eqref{I5a}, \eqref{I6a} and \eqref{I7a} we deduce \eqref{paraestim1}
	\begin{align}
		&\left(\mu_i v_{i,x}-(\tilde{\la}_i-\la_i^*)v_{i}-w_{i} \right)_x\nonumber\\
		&=O(1)\sum\limits_{j\neq i}\left(\mu_j v_{j,x}-(\tilde{\la}_j-\la_j^*)v_j-w_j\right)_x \bar{v}_j\xi_j	\nonumber\\
		&+O(1)\sum_{j\ne i}(|v_j|+|v_{j,x}|+|w_{j,x}|)(|v_j|+|w_{j,x}|+|v_{j,x}|+\sum_{k} |v_k||v_j|)\nonumber\\
		&+O(1)\left(\sum\limits_{j\neq i}\sum\limits_{k\neq j} (|w_{j,x}|+|v_{j,x}| +|v_j|)|v_{k,x} v_k|+\sum\limits_{j\neq i}\sum\limits_{k\neq j} |v_j v_{k,xx}v_k|\right)\nonumber\\
		&+O(1)\left(\sum_{i\ne j}|w_{j,xx}v_j-v_{j,xx}w_j|
		+\sum_k|v_kv_{k,x}||1-\chi_k^\e\xi_k\eta_k|\right)\nonumber\\
		&+O(1)\left(\sum_k v_k^2\abs{\left(\frac{w_k}{v_k}\right)_x}+\sum_k \sum_{l\ne k}|v_l|(|v_l||v_{k,x}|+|v_kv_{l,x}|) 
		+O(1)\sum\limits_{k}\sum\limits_{l\ne k}v_k^2|v_l| \right)\nonumber\\
		&+O(1)\sum\limits_{k}|v_k|^3(1-\chi_k^\e\xi_k\eta_k)+O(1)\sum_k \varkappa^\ep_k\rho_k^\e\mathfrak{A}_k|v_kv_{k,x}|\nonumber\\
		&+O(1)\left(\sum_{k\ne j}(|v_{k,x}v_j|+|v_k v_{j,x}|+|v_kv_j w_{k,x}|+|v_kv_jv_{j,x}|+|v_kv_jw_{j,x}|)+\sum_{j\ne k}\sum_l|v_kv_jv_l|
		\right).\nonumber
	\end{align}
	Now from \eqref{6.45}, we have	 using Lemmas \ref{lemme6.5}, \ref{lemme6.6}
	\begin{align*}
	&\sum_k|v_kv_{k,x}||1-\chi_k^\e\xi_k\eta_k|
		+\sum\limits_{k}|v_k|^3(1-\chi_k^\e\xi_k\eta_k)=O(1)\sum_k|v_kv_{k,x}||1-\chi_k^\e\eta_k|
		\\
		&+O(1)\sum\limits_{k}|v_k|^3(1-\chi_k^\e\eta_k)
		+O(1)\sum_k|v_kv_{k,x}|\mathbbm{1}_{\{|\frac{w_k}{v_k}|\geq\frac{\delta_1}{2}\}}+
		O(1)\sum\limits_{k}|v_k|^3 \mathbbm{1}_{\{|\frac{w_k}{v_k}|\geq\frac{\delta_1}{2}\}}\\
		&=O(1)\big(\sum_k|v_kv_{k,x}||1-\eta_k|+\sum_k|v_kv_{k,x}|\eta_k(1-\chi^\e_k)+\sum\limits_{k}|v_k|^3(1-\eta_k)\\
		&+\sum\limits_{k}|v_k|^3\eta_k (1-\chi^\e_k)\big)+O(1)\sum_l(\La_l^1+\La_l^4+\La_l^6)\\
		&=O(1)\sum_l(\La_l^1+\La_l^4+\La_l^6)+
		O(1)\big(\sum_k|v_kv_{k,x}|\eta_k(1-\chi^\e_k)+\sum\limits_{k}|v_k|^3\eta_k (1-\chi^\e_k)\big).
		\end{align*}
		Combining the previous estimate with \eqref{paraestim1} we deduce \eqref{estimate-v-i-xx-1aabis}.This completes the proof of Lemma \ref{lemme9.7}. 
\end{proof}
In the sequel we set for any $j\in\{1,\cdot,n\}$
$$\mathcal{H}_j:=\mu_j v_{j,x}-(\tilde{\la}_j-\la_j^*)v_j-w_j.$$
\begin{lemma}
\label{lemme9.11}
The following estimate is true.
\begin{align}
		&(\mu_i v_{i,x}-(\tilde{\la}_i-\la_i^*)v_i-w_i)_{xx}\nonumber\\
		&=O(1)\biggl(
		\sum_{k\ne j}(|v_{k,xx}v_j|+|v_{k,x}v_{j,x}|+|v_{k,x}v_j|+|v_kv_j|+|w_{k,x}v_j|)\nonumber\\
		&+
		\sum_{i\ne j}(|v_jv_{j,xxx}|+\abs{w_{j,xxx}v_j}+\abs{v_{j,xxx}w_j})
		+\sum_{j\ne i, k\ne j}|v_{k,xxx}v_k v_j|\mathfrak{A}_j\rho_j^\e
		\nonumber\\
		&+\sum_{i\ne j}(|v_j|+|w_j|+|v_{j,x}|+|w_{j,x}|+|v_{j,xx}|+|w_{j,xx}|+|v_j|\left(\frac{w_j}{v_j}\right)_{x}^2 \mathfrak{A}_j\rho_j^\e )\nonumber\\
		&+(|v_{i,x}^2|+|v_iv_{i,x}|+|v_iv_{i,xx}|+|v_i^2|)(1-\chi_i^\e\xi_i\eta_i)+\mathfrak{A}_i \Delta_i^\e\big(
		|v_i v_{i,x}|+v_{i,x}^2+|v_iv_{i,xx}|\big)
		\nonumber\\
		&+\mathfrak{A}_i\rho_i^\e\big(|v_iv_{i,x}\left(\frac{w_i}{v_i}\right)_x|
		+|v_i^2\left(\frac{w_i}{v_i}\right)^2_x|+|v_i^2\left(\frac{w_i}{v_i}\right)_x|+|w_{i,xx}v_i-v_{i,xx}w_i|
		\big)
		\nonumber\\
		&+(|v_i |(|v_{i,x}|+|v_{i,xx}|)+|v_{i,x}|(|v_{i,x}|+|w_{i,x}|))\varkappa^\ep_i\mathfrak{A}_i\rho_i^\e\biggl).\label{9.14primefin}
\end{align}
	It implies in particular that
	\begin{align}
		&(\mu_i v_{i,x}-(\tilde{\la}_i-\la_i^*)v_i-w_i)_{xx}\nonumber\\
		&=O(1)\biggl(
		\sum_j(\La_j^1+\La_j^2+\delta_0^2\La_j^3+\La_j^4+\La_j^5+\La_j^6)+
		%\\[5mm]
		\sum_{i\ne j}(|v_jv_{j,xxx}|+\abs{w_{j,xxx}v_j}+\abs{v_{j,xxx}w_j})
		\nonumber\\
		&+\sum_{i\ne j}(|v_j|+|w_j|+|v_{j,x}|+|w_{j,x}|+|v_{j,xx}|+|w_{j,xx}|+|v_j|\left(\frac{w_j}{v_j}\right)_{x}^2 \mathfrak{A}_j\rho_j^\e )\nonumber\\
		&+
(|v_iv_{i,x}|+|v_iv_{i,xx}|+|v_i^2|)\mathbbm{1}_{\{v_i^{2N}\leq 2\e\}}
+ |v_{i,x} v_{i,xx}|\mathbbm{1}_{\{|\frac{w_i}{v_i}|\geq\frac{\delta_1}{2}\}}\nonumber\\
		&+(|v_i |(|v_{i,x}|+|v_{i,xx}|)+|v_{i,x}|(|v_{i,x}|+|w_{i,x}|))\varkappa^\ep_i\mathfrak{A}_i\rho_i^\e\biggl).\label{9.15primefin}
\end{align}
\end{lemma}
\begin{proof}
From \eqref{9.14}, we deduce that
\begin{equation}
		(\mu_i v_{i,x}-(\tilde{\la}_i-\la_i^*)v_i-w_i)_{xx}=\sum\limits_{s=2}^{7}\mathcal{I}_{s,x}.
		\label{9.14prime}
	\end{equation}
We are now going to estimate each term $\mathcal{I}_{s,x}$ with $2\leq s\leq 7$ by calculating explicitly all these terms.

	\noi\underline{{Estimate of $\mathcal{I}_{2,x}$:}} 
	First we get
		\begin{align}
		\mathcal{I}_{2,1,x}&=-\sum\limits_{j\neq i}\mathcal{H}_{j,xx} \frac{\mu_i}{\mu_j}[\psi_{ji} +v_j\xi_j\pa_{v_j}\bar{v}_j \psi_{ji,v}]-\sum\limits_{j\neq i}\mathcal{H}_{j,x} \pa_x(\frac{\mu_i}{\mu_j})[\psi_{ji} +v_j\xi_j\pa_{v_j}\bar{v}_j \psi_{ji,v}]\nonumber\\
		&-\sum\limits_{j\neq i}\mathcal{H}_{j,x} \frac{\mu_i}{\mu_j}[\pa_x(\psi_{ji} )+v_j\xi_j\pa_{v_j}\bar{v}_j \pa_x(\psi_{ji,v})]\nonumber\\
		&-\sum\limits_{j\neq i}\mathcal{H}_{j,x} \frac{\mu_i}{\mu_j}[v_{j,x}\xi_j\pa_{v_j}\bar{v}_j \psi_{ji,v}+v_{j}\xi_j'\left(\frac{w_j}{v_j}\right)_x\pa_{v_j}\bar{v}_j \psi_{ji,v}+ v_{j}\xi_j\pa_x(\pa_{v_j}\bar{v}_j) \psi_{ji,v} ]\label{I21x}
		\end{align}
A direct computation gives
\begin{align*}
\mathcal{H}_{j,xx}=&\mu_jv_{j,xxx}+2\mu_{j,x}v_{j,xx}+\mu_{j,xx}v_{j,x}-\tilde{\la}_{j,xx}v_j-2\tilde{\la}_{j,x}v_{j,x}-(\tilde{\la}_j-\la_j^*)v_{j,xx}-w_{j,xx}.
\end{align*}
From Lemma \ref{lemme11.6}, \ref{lemme6.5}, \ref{lemme6.6}, we observe that $v_j\tilde{\la}_{j,xx}=O(1)$ it implies from Lemmas \ref{lemme6.5}, \ref{lemme6.6} and the fact that $\mu_{j,xx},\mu_{j,x}=O(1)$ that
\begin{align}
&\widetilde{\mathcal{H}}_{j,xx}=O(1)+O(1)|v_{j,xxx}|.\label{Hijxx}
\end{align}
Now combining \eqref{I21x}, \eqref{Hijxx} , the fact that $\psi_{ji}=O(1)\xi_j\bar{v}_j, \pa_x(\frac{\mu_i}{\mu_j})=O(1)$ and the Lemmas \ref{lemme6.5}, \ref{lemme6.6}, \ref{estimimpo1}, \ref{lemme11.3a} we deduce that	
\begin{align*}
		\mathcal{I}_{2,1,x}
		&=O(1)\sum_{i\ne j}(|v_jv_{j,xxx}|+|v_j|+|v_{j,x}|+|w_{j,x}|+|v_{j,xx}|)
		\end{align*}
	Next we have from Lemmas \ref{lemme6.5}, \ref{lemme6.6}, \ref{estimimpo1}, \ref{lemme11.3}, \ref{lemme11.3a}
	\begin{align*}
	\mathcal{I}_{2,2,x}&=-\left(\sum\limits_{j\neq i} \mathcal{H}_{j,x}  \pa_x(\frac{\mu_i}{\mu_j})[\psi_{ji} +v_j\xi_j\pa_{v_j}\bar{v}_j \psi_{ji,v}]\right)-\left(\sum\limits_{j\neq i} \mathcal{H}_{j}  \pa_{xx}(\frac{\mu_i}{\mu_j})[\psi_{ji} +v_j\xi_j\pa_{v_j}\bar{v}_j \psi_{ji,v}]\right)\\
	&-\left(\sum\limits_{j\neq i} \mathcal{H}_{j}  \pa_{x}(\frac{\mu_i}{\mu_j})[\pa_x(\psi_{ji} )+v_j\xi_j\pa_{v_j}\bar{v}_j\pa_x( \psi_{ji,v})]\right)-\left(\sum\limits_{j\neq i} \mathcal{H}_{j}  \pa_{x}(\frac{\mu_i}{\mu_j})[v_{j,x}\xi_j\pa_{v_j}\bar{v}_j \psi_{ji,v}]\right)\\
	&-\left(\sum\limits_{j\neq i} \mathcal{H}_{j}  \pa_{x}(\frac{\mu_i}{\mu_j})[v_j\xi_j'\left(\frac{w_j}{v_j}\right)_x\pa_{v_j}\bar{v}_j \psi_{ji,v}+v_j\xi_j\pa_x(\pa_{v_j}\bar{v}_j )\psi_{ji,v}]\right)\\
	&=O(1)\sum_{i\ne j}(|v_j|+|w_j|+|v_{j,x}|+|w_{j,x}|+|v_{j,xx}|).
	\end{align*}
From \eqref{ngtech5} and Lemmas \ref{estimimpo1}, \ref{lemme6.5}, \ref{lemme6.6}, \ref{lemme11.5} we deduce that
\begin{align}
&\left(\sum\limits_{j\neq i} \mathcal{H}_{j} \frac{\mu_i}{\mu_j}\pa_{xx} \psi_{ji} \right)=\label{techj68}\\
&O(1)\biggl(\sum_{i\ne j}(|v_j|+|w_j|+|v_{j,x}|)\mathfrak{A}_j\rho_j^\e+\sum_{i\ne j} \mathcal{H}_{j} \left(\sum\limits_{l\neq j}\frac{\abs{v_l v_{l,x}v_{j,x}}}{\abs{v_j}}+\frac{v_{j,x}^2}{|v_j|} \right)\mathfrak{A}_j\rho_j^\e \nonumber\\
		&+\sum_{i\ne j} \mathcal{H}_{j}\left(|v_j\left(\frac{w_j}{v_j}\right)_{xx}|+|v_j\left(\frac{w_j}{v_j}\right)_{x}^2|\right)\mathfrak{A}_j\rho_j^\e +\sum_{i\ne j}  \mathcal{H}_{j}|\left(\frac{w_j}{v_j}\right)_x|\left(\sum\limits_{k\neq j}\abs{v_{k,x}v_k}+\abs{v_{j,x}}\right)\mathfrak{A}_j\rho_j^\e\biggl) \nonumber\\
		&=O(1)\sum_{i\ne j}(|v_j|+|w_j|+|v_{j,x}|)\big(1+|v_j|\left(\frac{w_j}{v_j}\right)_{x}^2  \big)\mathfrak{A}_j\rho_j^\e.\nonumber
	\end{align}	
	We obtain similarly from \eqref{techj68} and Lemmas \ref{estimimpo1}, \ref{lemme6.5}, \ref{lemme6.6}, \ref{lemme11.3a}
	\begin{align}
	&-\left(\sum\limits_{j\neq i} \widetilde{\mathcal{H}}_{j} \frac{\mu_i}{\mu_j}[\pa_{xx} \psi_{ji} +v_{j,x}\xi_j\pa_{v_j}\bar{v}_j \pa_x(\psi_{ji,v})+v_{j}\xi_j^\p\left(\frac{w_j}{v_j}\right)_x \pa_{v_j}\bar{v}_j \pa_x(\psi_{ji,v})]\right)\nonumber\\
	&=O(1)\sum_{i\ne j}(|v_j|+|w_j|+|v_{j,x}|)\big(1+|v_j|\left(\frac{w_j}{v_j}\right)_{x}^2  \big)\mathfrak{A}_j\rho_j^\e  .\label{techj681}
	\end{align}
	From \eqref{techj681}, \eqref{ngtech5}  and Lemmas \ref{lemme6.5}, \ref{lemme6.6}, \ref{estimimpo1}, \ref{lemme11.3a} we finally get
	\begin{align*}
	&\mathcal{I}_{2,3,x}=-\left(\sum\limits_{j\neq i} \mathcal{H}_{j,x} \frac{\mu_i}{\mu_j}[\pa_x \psi_{ji} +v_{j,x}\xi_j\pa_{v_j}\bar{v}_j \psi_{ji,v}+v_{j}\xi_j^\p\left(\frac{w_j}{v_j}\right)_x \pa_{v_j}\bar{v}_j \psi_{ji,v}]\right)\\
	&-\left(\sum\limits_{j\neq i} \mathcal{H}_{j}\pa_x( \frac{\mu_i}{\mu_j})[\pa_x \psi_{ji} +v_{j,x}\xi_j\pa_{v_j}\bar{v}_j \psi_{ji,v}+v_{j}\xi_j^\p\left(\frac{w_j}{v_j}\right)_x \pa_{v_j}\bar{v}_j \psi_{ji,v}]\right)\\
	&-\left(\sum\limits_{j\neq i} \mathcal{H}_{j} \frac{\mu_i}{\mu_j}[\pa_{xx} \psi_{ji} +v_{j,x}\xi_j\pa_{v_j}\bar{v}_j \pa_x(\psi_{ji,v})+v_{j}\xi_j^\p\left(\frac{w_j}{v_j}\right)_x \pa_{v_j}\bar{v}_j \pa_x(\psi_{ji,v})]\right)\\
	&-\left(\sum\limits_{j\neq i} \mathcal{H}_{j} \frac{\mu_i}{\mu_j}[v_{j,xx}\xi_j\pa_{v_j}\bar{v}_j \psi_{ji,v}+v_{j,x}\xi_j'\left(\frac{w_j}{v_j}\right)_x \pa_{v_j}\bar{v}_j \psi_{ji,v}+v_{j,x}\xi_j\pa_x(\pa_{v_j}\bar{v}_j )\psi_{ji,v}]\right)\\
	&-\left(\sum\limits_{j\neq i} \mathcal{H}_{j} \frac{\mu_i}{\mu_j}\big[(v_{j,x}\xi_j^\p+v_j\xi_j''\left(\frac{w_j}{v_j}\right)_x)\left(\frac{w_j}{v_j}\right)_x \pa_{v_j}\bar{v}_j \psi_{ji,v}
	+v_{j}\xi_j^\p\big(\left(\frac{w_j}{v_j}\right)_{xx} \pa_{v_j}\bar{v}_j +\left(\frac{w_j}{v_j}\right)_{x} \pa_x(\pa_{v_j}\bar{v}_j)\big)\psi_{ji,v}\big]\right)\\
	&=O(1)\sum_{i\ne j}(|v_j|+|w_j|+|v_{j,x}|+|w_{j,x}|+|v_{j,xx}|)\big(1+|v_j|\left(\frac{w_j}{v_j}\right)_{x}^2  \big)\mathfrak{A}_j\rho_j^\e .
	\end{align*}
Let us consider now the term $\mathcal{I}_{2,4,x}$. First we have using again Lemmas  \ref{estimimpo1}, \ref{lemme6.5}, \ref{lemme6.6},\ref{lemme11.3a}, \ref{lemme11.5} and \eqref{ngtech5}
\begin{align}
&-\left(\sum\limits_{j\neq i} \widetilde{\mathcal{H}}_{j}\frac{\mu_i}{\mu_j}[v_{j}\xi_j\pa_x(\pa_{v_j}\bar{v}_j)\pa_x \psi_{ji,v}+v_{j}\xi_j\pa_{v_j}\bar{v}_j \pa_{xx}(\psi_{ji,v})]\right)\nonumber\\
&=O(1)\sum_{j\ne i}\biggl((|v_j|+|v_{j,x}|+|w_j|)\rho_j^\e\mathfrak{A}_j+
\rho_j^\e\mathfrak{A}_j\widetilde{\mathcal{H}}_{j} \big(|v_j\left(\frac{w_j}{v_j}\right)_{xx}|+|v_j|\left(\frac{w_j}{v_j}\right)_{x}^2\big)\nonumber\\
&+ \widetilde{\mathcal{H}}_{j}|\left(\frac{w_j}{v_j}\right)_x|\rho_j^\e\mathfrak{A}_j\big(|v_{j,x}|+\sum_{k\ne j}|v_kv_{k,x}|\biggl)\nonumber\\
&=O(1)\sum_{j\ne i}\big(|v_j|+|v_{j,x}|+|w_j|\big)\big(1+
|v_j|\left(\frac{w_j}{v_j}\right)_{x}^2\big)\rho_j^\e\mathfrak{A}_j.\label{techh181}
\end{align}
We recall now using Lemma \ref{estimimpo1}, \ref{lemme11.3}
\begin{align}
&v_j \pa_{xx}(\pa_{v_j}\bar{v}_j)=v_j\pa_{v_jv_jv_j}\bar{v}_j(v_{j,x})^2+v_j\pa_{v_jv_j}\bar{v}_j v_{j,xx}+2v_j\sum_{k\ne j}\pa_{v_k v_jv_j}\bar{v}_j v_{k,x}v_{j,x}\nonumber\\
&+v_j\sum_{k\ne j}\pa_{v_k v_j}\bar{v}_j v_{k,xx}+v_j\sum_{k\ne j}\sum_{l\ne j}\pa_{v_k v_l v_j}\bar{v}_j v_{k,x}v_{l,x}
\nonumber\\
%&+2v_j\sum_{k\ne j}\sum_{l\ne j}\pa_{v_k w_l v_j}\bar{v}_j v_{k,x}w_{l,x}\nonumber\\
&=O(1)\rho_j^\e\biggl(\frac{v_{j,x}^2}{|v_j|}+|v_{j,xx}|+\sum_{k\ne j}(|\frac{v_kv_{k,x}v_{j,x}}{v_j}|+|v_kv_{k,xx}|)+\sum_{k\ne j}\sum_{l\ne j}|v_{k,x}v_{l,x}|\biggl).\label{techh18}
\end{align}
Similarly we have using Lemmas \ref{lemme6.5}, \ref{lemme6.6}, \ref{lemme6.8}, \ref{estimimpo1} and \eqref{techh18}, \eqref{ngtech5}
\begin{align}
&-\left(\sum\limits_{j\neq i} \widetilde{\mathcal{H}}_{j}\frac{\mu_i}{\mu_j}\big[\big(v_{j,x}\xi_j+v_j\xi_j'\left(\frac{w_j}{v_j}\right)_x\big)\pa_x(\pa_{v_j}\bar{v}_j) \psi_{ji,v}+v_{j}\xi_j\pa_{xx}(\pa_{v_j}\bar{v}_j) \psi_{ji,v}\right)\nonumber\\
&=O(1)\sum_{j\ne i}\big(|v_{j,x}|+|v_j|+|w_j|\big))\big(1+
|v_j|\left(\frac{w_j}{v_j}\right)_{x}^2\big)\rho_j^\e\mathfrak{A}_j.\label{techh182}
\end{align}
Combining now \eqref{techh181}, \eqref{techh182}, \eqref{ngtech5} and using Lemmas \ref{estimimpo1}, \ref{lemme11.3a}
\begin{align*}
&\mathcal{I}_{2,4,x}=-\left(\sum\limits_{j\neq i} \widetilde{\mathcal{H}}_{j,x}\frac{\mu_i}{\mu_j}[v_{j}\xi_j\pa_x(\pa_{v_j}\bar{v}_j) \psi_{ji,v}+v_{j}\xi_j\pa_{v_j}\bar{v}_j \pa_x(\psi_{ji,v})]\right)\\
&-\left(\sum\limits_{j\neq i} \widetilde{\mathcal{H}}_{j}\pa_x(\frac{\mu_i}{\mu_j})[v_{j}\xi_j\pa_x(\pa_{v_j}\bar{v}_j) \psi_{ji,v}+v_{j}\xi_j\pa_{v_j}\bar{v}_j \pa_x(\psi_{ji,v})]\right)\\
&-\left(\sum\limits_{j\neq i} \widetilde{\mathcal{H}}_{j}\frac{\mu_i}{\mu_j}[v_{j}\xi_j\pa_x(\pa_{v_j}\bar{v}_j)\pa_x \psi_{ji,v}+v_{j}\xi_j\pa_{v_j}\bar{v}_j \pa_{xx}(\psi_{ji,v})]\right)\\
&-\left(\sum\limits_{j\neq i} \widetilde{\mathcal{H}}_{j}\frac{\mu_i}{\mu_j}\big[\big(v_{j,x}\xi_j+v_j\xi_j'\left(\frac{w_j}{v_j}\right)_x\big)\pa_x(\pa_{v_j}\bar{v}_j) \psi_{ji,v}+v_{j}\xi_j\pa_{xx}(\pa_{v_j}\bar{v}_j) \psi_{ji,v}\right)\\
&-\left(\sum\limits_{j\neq i} \widetilde{\mathcal{H}}_{j}\frac{\mu_i}{\mu_j}\big[\big(v_{j,x}\xi_j+v_j\xi_j'\left(\frac{w_j}{v_j}\right)_x\big)\pa_{v_j}\bar{v}_j \pa_x(\psi_{ji,v})+v_j\xi_j \pa_x(\pa_{v_j}\bar{v}_j)\pa_x(\psi_{ji,v})\big]\right)\\
&=O(1)\sum_{j\ne i}\big(|v_{j,xx}|+|v_j|+|v_{j,x}|+|w_j|+|w_{j,x}|\big)\big(1+
|v_j|\left(\frac{w_j}{v_j}\right)_{x}^2\big)\rho_j^\e\mathfrak{A}_j.
\end{align*}
		\noi\underline{{Estimate of $\mathcal{I}_{3,x}$:}}
		We have using the Lemmas \ref{lemme6.5}, \ref{lemme6.6},\ref{lemme11.3a}, \ref{estimimpo1} and \eqref{ngtech5}
	\begin{align*}
		\mathcal{I}_{3,1,x}&=-\sum\limits_{j\neq i}\sum\limits_{k\neq j}\xi_j^{\p\p}\left(\frac{w_j}{v_j}\right)^2_xv_{k,x}v_j\pa_{v_k}\bar{v}_j\mu_i \psi_{ji,v}-\sum\limits_{j\neq i}\sum\limits_{k\neq j}\xi_j^\p\left(\frac{w_j}{v_j}\right)_{xx}v_{k,x}v_j\pa_{v_k}\bar{v}_j\mu_i \psi_{ji,v}\\
		&-\sum\limits_{j\neq i}\sum\limits_{k\neq j}\xi_j^\p\left(\frac{w_j}{v_j}\right)_xv_{k,xx}v_j\pa_{v_k}\bar{v}_j\mu_i \psi_{ji,v}-\sum\limits_{j\neq i}\sum\limits_{k\neq j}\xi_j^\p\left(\frac{w_j}{v_j}\right)_xv_{k,x}v_{j,x}\pa_{v_k}\bar{v}_j\mu_i \psi_{ji,v}\\
		&-\sum\limits_{j\neq i}\sum\limits_{k\neq j}\xi_j^\p\left(\frac{w_j}{v_j}\right)_xv_{k,x}v_j\pa_x(\pa_{v_k}\bar{v}_j)\mu_i \psi_{ji,v}\\
		&-\sum\limits_{j\neq i}\sum\limits_{k\neq j}\xi_j^\p\left(\frac{w_j}{v_j}\right)_xv_{k,x}v_j\pa_{v_k}\bar{v}_j[\mu_{i,x} \psi_{ji,v}+\mu_i\pa_x(\psi_{ji,v})]\\
		&=O(1)\sum\limits_{j\neq i}\big(|v_{j,xx}|+|w_{j,xx}|+|v_{j,x}|+|w_{j,x}|+ |v_j\left(\frac{w_j}{v_j}\right)_x^2|\big)\mathfrak{A}_j\rho_j^\e.
		\end{align*}
			Using Lemmas \ref{lemme6.5}, \ref{lemme6.6}, \ref{estimimpo1}, \ref{lemme11.3a} and \eqref{ngtech5} we obtain
		\begin{align*}
		\mathcal{I}_{3,2,x}&=-\sum\limits_{j\neq i}\sum\limits_{k\neq j}\xi^\p_j\left(\frac{w_j}{v_j}\right)_x v_{k,xx}v_j\pa_{v_k}\bar{v}_j\mu_i \psi_{ji,v}
		-\sum\limits_{j\neq i}\sum\limits_{k\neq j}\xi_jv_{k,xxx}v_j\pa_{v_k}\bar{v}_j\mu_i \psi_{ji,v}\\
		&-\sum\limits_{j\neq i}\sum\limits_{k\neq j}\xi_jv_{k,xx}v_{j,x}\pa_{v_k}\bar{v}_j\mu_i \psi_{ji,v}-\sum\limits_{j\neq i}\sum\limits_{k\neq j}\xi_jv_{k,xx}v_j\pa_x(\pa_{v_k}\bar{v}_j)\mu_i \psi_{ji,v}\\
		&-\sum\limits_{j\neq i}\sum\limits_{k\neq j}\xi_jv_{k,xx}v_j\pa_{v_k}\bar{v}_j[\mu_{i,x} \psi_{ji,v}+\mu_i\pa_x(\psi_{ji,v})]\\
		&=O(1)\biggl(\sum\limits_{j\neq i}\left( \abs{w_{j,x}}+\abs{v_{j,x}}+|v_j|\right)\mathfrak{A}_j\rho_j^\e+\sum_{j\ne i, k\ne j}|v_{k,xxx}v_k v_j|\mathfrak{A}_j\rho_j^\e\biggl),
		\end{align*}
		\begin{align*}
		\mathcal{I}_{3,3,x}&=-\sum\limits_{j\neq i}\sum\limits_{k\neq j}\xi_j^\p\left(\frac{w_j}{v_j}\right)_xv_{k,x}v_{j,x}\pa_{v_k}\bar{v}_j\mu_i \psi_{ji,v}-\sum\limits_{j\neq i}\sum\limits_{k\neq j}\xi_jv_{k,xx}v_{j,x}\pa_{v_k}\bar{v}_j\mu_i \psi_{ji,v}\\
		&-\sum\limits_{j\neq i}\sum\limits_{k\neq j}\xi_jv_{k,x}v_{j,xx}\pa_{v_k}\bar{v}_j\mu_i \psi_{ji,v}-\sum\limits_{j\neq i}\sum\limits_{k\neq j}\xi_jv_{k,x}v_{j,x}\pa_x(\pa_{v_k}\bar{v}_j)\mu_i \psi_{ji,v}\\
		&-\sum\limits_{j\neq i}\sum\limits_{k\neq j}\xi_jv_{k,x}v_{j,x}\pa_{v_k}\bar{v}_j[\mu_{i,x} \psi_{ji,v}+\mu_i\pa_x(\psi_{ji,v})]\\
		&=O(1)\sum\limits_{j\neq i}\left( \abs{w_{j,x}}+\abs{v_{j,x}}+|v_{j,xx}|\right)\mathfrak{A}_j\rho_j^\e,
		\end{align*}
		\begin{align*}
		\mathcal{I}_{3,4,x}&=-\sum\limits_{j\neq i}\sum\limits_{k\neq j}\xi_j^\p\left(\frac{w_j}{v_j}\right)_xv_{k,x}v_j\pa_{v_k}\bar{v}_j\mu_{i,x} \psi_{ji,v}-\sum\limits_{j\neq i}\sum\limits_{k\neq j}\xi_jv_{k,xx}v_j\pa_{v_k}\bar{v}_j\mu_{i,x} \psi_{ji,v}\\
		&-\sum\limits_{j\neq i}\sum\limits_{k\neq j}\xi_jv_{k,x}v_{j,x}\pa_{v_k}\bar{v}_j\mu_{i,x} \psi_{ji,v}-\sum\limits_{j\neq i}\sum\limits_{k\neq j}\xi_jv_{k,x}v_j\pa_x(\pa_{v_k}\bar{v}_j)\mu_{i,x} \psi_{ji,v}\\
		&-\sum\limits_{j\neq i}\sum\limits_{k\neq j}\xi_jv_{k,x}v_j\pa_{v_k}\bar{v}_j[\mu_{i,xx} \psi_{ji,v}+\mu_{i,x}\pa_x(\psi_{ji,v})]\\
		&=O(1)\sum\limits_{j\neq i}\left( \abs{w_{j,x}}+\abs{v_{j,x}}+\abs{v_j}\right)\mathfrak{A}_j\rho_j^\e,
		\end{align*}
		\begin{align*}
		\mathcal{I}_{3,5,x}&=-\sum\limits_{j\neq i}\sum\limits_{k\neq j}\xi_j^\p\left(\frac{w_j}{v_j}\right)_xv_{k,x}v_jv_{j,x}\pa_{v_kv_j}\bar{v}_j\mu_i \psi_{ji,v}-\sum\limits_{j\neq i}\sum\limits_{k\neq j}\xi_jv_{k,xx}v_jv_{j,x}\pa_{v_kv_j}\bar{v}_j\mu_i \psi_{ji,v}\\
		&-\sum\limits_{j\neq i}\sum\limits_{k\neq j}\xi_jv_{k,x}v^2_{j,x}\pa_{v_kv_j}\bar{v}_j\mu_i \psi_{ji,v}-\sum\limits_{j\neq i}\sum\limits_{k\neq j}\xi_jv_{k,x}v_jv_{j,xx}\pa_{v_kv_j}\bar{v}_j\mu_i \psi_{ji,v}\\
		&-\sum\limits_{j\neq i}\sum\limits_{k\neq j}\xi_jv_{k,x}v_jv_{j,x}\pa_x(\pa_{v_kv_j}\bar{v}_j)\mu_i \psi_{ji,v}-\sum\limits_{j\neq i}\sum\limits_{k\neq j}\xi_jv_{k,x}v_jv_{j,x}\pa_{v_kv_j}\bar{v}_j[\mu_{i,x} \psi_{ji,v}+\mu_i\pa_x(\psi_{ji,v})]\\
		&=O(1)\sum\limits_{j\neq i}\left(|v_{j,xx}|+ \abs{w_{j,x}}+\abs{v_{j,x}}\right)\mathfrak{A}_j\rho_j^\e.
		\end{align*}
		We have now using Lemmas \ref{lemme6.5}, \ref{lemme6.6}, \ref{estimimpo1} and \eqref{6.54}
		\begin{align}
		&-\sum\limits_{j\neq i}\sum\limits_{k\neq j}\sum\limits_{l\neq j}\xi_jv_{k,x}v_jv_{l,x}\pa_x(\pa_{v_lv_k}\bar{v}_j)\mu_i \psi_{ji,v}\nonumber\\
		&=O(1)\sum_{j\ne i}\sum\limits_{k\neq j}\mathfrak{A}_j\rho_j^\e |v_{k,x}v_{l,x}|\big(|v_{j,x}|+\sum_{p\ne j}|v_lv_{p,x}|\big)\nonumber\\
		&=
		O(1)\sum_{j\ne i} |v_{j,x}|\mathfrak{A}_j\rho_j^\e +O(1)\sum_{j\ne i}\sum\limits_{k\neq j}\mathfrak{A}_j\rho_j^\e  (|v_k|+|w_k|+|v_j|)
		|v_{l,x}|\sum_{p\ne j}|v_lv_{p,x}|\nonumber\\
		&+O(1)\sum_{j\ne i}\sum\limits_{k\neq j}\mathfrak{A}_j\rho_j^\e  \sum_{p'\ne j}|v_{p'}^2\left(\frac{w_{p'}}{v_{p'}}\right)_x| \mathfrak{A}_{p'}\rho_{p'}^\e
		|v_{l,x}|\sum_{p\ne j}|v_lv_{p,x}|\nonumber\\
		&=O(1)\sum_j(|v_{j,x}|+|v_j|)\mathfrak{A}_j\rho_j^\e.\label{techhh1}
		\end{align}
		Combining \eqref{techhh1} and Lemmas \ref{estimimpo1}, \ref{lemme11.3} we obtain
		\begin{align*}
		&\mathcal{I}_{3,6,x}=-\sum\limits_{j\neq i}\sum\limits_{k\neq j}\sum\limits_{l\neq j}\xi_j^\p\left(\frac{w_j}{v_j}\right)_xv_{k,x}v_jv_{l,x}\pa_{v_lv_k}\bar{v}_j\mu_i \psi_{ji,v}-\sum\limits_{j\neq i}\sum\limits_{k\neq j}\sum\limits_{l\neq j}\xi_jv_{k,xx}v_jv_{l,x}\pa_{v_lv_k}\bar{v}_j\mu_i \psi_{ji,v}\\
		&-\sum\limits_{j\neq i}\sum\limits_{k\neq j}\sum\limits_{l\neq j}\xi_jv_{k,x}v_{j,x}v_{l,x}\pa_{v_lv_k}\bar{v}_j\mu_i \psi_{ji,v}-\sum\limits_{j\neq i}\sum\limits_{k\neq j}\sum\limits_{l\neq j}\xi_jv_{k,x}v_jv_{l,xx}\pa_{v_lv_k}\bar{v}_j\mu_i \psi_{ji,v}\\
		&-\sum\limits_{j\neq i}\sum\limits_{k\neq j}\sum\limits_{l\neq j}\xi_jv_{k,x}v_jv_{l,x}\pa_x(\pa_{v_lv_k}\bar{v}_j)\mu_i \psi_{ji,v}-\sum\limits_{j\neq i}\sum\limits_{k\neq j}\sum\limits_{l\neq j}\xi_jv_{k,x}v_jv_{l,x}\pa_{v_lv_k}\bar{v}_j[\mu_{i,x} \psi_{ji,v}+\mu_i\pa_x(\psi_{ji,v})]\\
		&=O(1)\sum\limits_{j\neq i}\left( \abs{w_{j,x}}+\abs{v_{j,x}}+|v_j|\right)\mathfrak{A}_j\rho_j^\e.
		\end{align*}
		Applying Lemmas \ref{lemme6.5}, \ref{lemme6.6}, \ref{lemme11.3a}, \ref{lemme11.5} and \eqref{ngtech5} we get
	\begin{align*}
		\mathcal{I}_{3,7,x}&=-\sum\limits_{j\neq i}\sum\limits_{k\neq j}\xi_j^\p\left(\frac{w_j}{v_j}\right)_xv_{k,x}v_j\pa_{v_k}\bar{v}_j\mu_i \pa_x(\psi_{ji,v})-\sum\limits_{j\neq i}\sum\limits_{k\neq j}\xi_jv_{k,xx}v_j\pa_{v_k}\bar{v}_j\mu_i \pa_x(\psi_{ji,v})\\
		&-\sum\limits_{j\neq i}\sum\limits_{k\neq j}\xi_jv_{k,x}v_{j,x}\pa_{v_k}\bar{v}_j\mu_i \pa_x(\psi_{ji,v})-\sum\limits_{j\neq i}\sum\limits_{k\neq j}\xi_jv_{k,x}v_j\pa_x(\pa_{v_k}\bar{v}_j)\mu_i \pa_x(\psi_{ji,v})\\
		&-\sum\limits_{j\neq i}\sum\limits_{k\neq j}\xi_jv_{k,x}v_j\pa_{v_k}\bar{v}_j[\mu_{i,x} \pa_x\psi_{ji,v}+\mu_i\pa_{xx}(\psi_{ji,v})]\\
		&=O(1)\sum\limits_{j\neq i}\left(|v_{j,xx}|+|w_{j,xx}|+ \abs{w_{j,x}}+\abs{v_{j,x}}+|v_j|+|v_j|\left(\frac{w_j}{v_j}\right)_x^2\right)\mathfrak{A}_j\rho_j^\e.
		\end{align*}
		
	\noi\underline{{Estimate of $\mathcal{I}_{4,x}$:}}
	Next we have applying Lemmas \ref{estimimpo1}, \ref{lemme11.3a}, \eqref{ngtech5}  and the fact that $\psi_{ji}, \psi_{ji,\sig}=O(1)\xi_j \bar{v}_j$
	\begin{align*}
		\mathcal{I}_{4,1,x}&=-\sum\limits_{j\neq i}v_{j,xx}\left(\frac{w_j}{v_j}\right)_x\mu_i [-\psi_{ji,\si}+\bar{v}_j\xi_j^\p\psi_{ji,v}]-\sum\limits_{j\neq i}v_{j,x}\left(\frac{w_j}{v_j}\right)_{xx}\mu_i [-\psi_{ji,\si}+\bar{v}_j\xi_j^\p\psi_{ji,v}]\\
		&-\sum\limits_{j\neq i}v_{j,x}\left(\frac{w_j}{v_j}\right)_x\mu_{i,x} [-\psi_{ji,\si}+\bar{v}_j\xi_j^\p\psi_{ji,v}]-\sum\limits_{j\neq i}v_{j,x}\left(\frac{w_j}{v_j}\right)^2_x\mu_i [\bar{v}_j\xi_j^{\p\p}\psi_{ji,v}]\\
		&-\sum\limits_{j\neq i}v_{j,x}\left(\frac{w_j}{v_j}\right)_x\mu_i [-\pa_x(\psi_{ji,\si})+\bar{v}_j\xi_j^\p\pa_x(\psi_{ji,v})+\pa_x(\bar{v}_j)\xi_j^\p\psi_{ji,v}]\\
		&=O(1)\sum\limits_{j\neq i}\left(|v_{j,x}|+|w_{j,x}|+\abs{v_j}\left(\frac{w_j}{v_j}\right)_x^2\right)\mathfrak{A}_j\rho_j^\e.
	\end{align*}
	Then we have the following estimate using Lemmas \ref{estimimpo1}, \ref{lemme6.5}, \ref{lemme6.6}, \ref{lemme11.3a}, \eqref{ngtech5},  \eqref{identity:wi-vi-xxx} and the fact that $\psi_{ji}, \psi_{ji,\sig}=O(1)\xi_j \bar{v}_j$
	\begin{align*}
		&\mathcal{I}_{4,2,x}=-\sum\limits_{j\neq i}v_{j,x}\left(\frac{w_j}{v_j}\right)_{xx}\mu_i [-\psi_{ji,\si}+\bar{v}_j\xi_j^\p\psi_{ji,v}]-\sum\limits_{j\neq i}v_j\left(\frac{w_j}{v_j}\right)_{xxx}\mu_i [-\psi_{ji,\si}+\bar{v}_j\xi_j^\p\psi_{ji,v}]\\
		&-\sum\limits_{j\neq i}v_j\left(\frac{w_j}{v_j}\right)_{xx}\left(\frac{w_j}{v_j}\right)_{x}\mu_i [\bar{v}_j\xi_j^{\p\p}\psi_{ji,v}]-\sum\limits_{j\neq i}v_j\left(\frac{w_j}{v_j}\right)_{xx}\mu_{i,x} [-\psi_{ji,\si}+\bar{v}_j\xi_j^\p\psi_{ji,v}]\\
		&-\sum\limits_{j\neq i}v_j\left(\frac{w_j}{v_j}\right)_{xx}\mu_i [-\pa_x(\psi_{ji,\si})+\pa_x(\bar{v}_j)\xi_j^\p\psi_{ji,v}+\bar{v}_j\xi_j^\p\pa_x(\psi_{ji,v})]\\
		&=O(1)\sum\limits_{j\neq i}\biggl(\abs{w_{j,xxx}v_j}+\abs{v_{j,xxx}w_j}+\abs{v_{j,xx}}+|w_{j,xx}|+\abs{w_{j,x}}+\abs{v_{j,x}}\\
		&+|v_j|+ \abs{v_j}\left(\frac{w_j}{v_j}\right)_x^2 \biggl)\mathfrak{A}_j\rho_j^\e.
		\end{align*}
		Applying Lemma \ref{lemme6.5}, \ref{lemme6.6}, \ref{estimimpo1}, \ref{lemme11.3a}
		\begin{align*}
		\mathcal{I}_{4,3,x}&=-\sum\limits_{j\neq i}v_{j,x}\left(\frac{w_j}{v_j}\right)^2_x\mu_i [\bar{v}_j\xi_j^{\p\p}\psi_{ji,v}]-\sum\limits_{j\neq i}2v_j\left(\frac{w_j}{v_j}\right)_x\left(\frac{w_j}{v_j}\right)_{xx}\mu_i [\bar{v}_j\xi_j^{\p\p}\psi_{ji,v}]\\
		&-\sum\limits_{j\neq i}v_j\left(\frac{w_j}{v_j}\right)^2_x\mu_{i,x} [\bar{v}_j\xi_j^{\p\p}\psi_{ji,v}]-\sum\limits_{j\neq i}v_j\left(\frac{w_j}{v_j}\right)^3_x\mu_i [\bar{v}_j\xi_j^{\p\p\p}\psi_{ji,v}]\\
		&-\sum\limits_{j\neq i}v_j\left(\frac{w_j}{v_j}\right)^2_x\mu_i [\bar{v}_j\xi_j^{\p\p}\pa_x(\psi_{ji,v})+\pa_x(\bar{v}_j)\xi_j^{\p\p}\psi_{ji,v}]\\
		&=O(1)\sum\limits_{j\neq i}\left(\abs{w_{j,x}}+\abs{v_{j,x}}+|v_j \left(\frac{w_j}{v_j}\right)_x^2| \right)\mathfrak{A}_j\rho_j^\e.
		\end{align*}
		Applying Lemmas \ref{estimimpo1}, \ref{lemme6.5}, \ref{lemme6.6}, \ref{lemme11.3a} and \eqref{ngtech5}
		\begin{align*}
		\mathcal{I}_{4,4,x}&=-\sum\limits_{j\neq i}v_{j,x}\left(\frac{w_j}{v_j}\right)_x\mu_i v_{j,x}\pa_{v_j}\bar{v}_j\xi_j^\p\psi_{ji,v}-\sum\limits_{j\neq i}v_j\left(\frac{w_j}{v_j}\right)_{xx}\mu_i v_{j,x}\pa_{v_j}\bar{v}_j\xi_j^\p\psi_{ji,v}\\
		&-\sum\limits_{j\neq i}v_j\left(\frac{w_j}{v_j}\right)_x\mu_{i,x} v_{j,x}\pa_{v_j}\bar{v}_j\xi_j^\p\psi_{ji,v}-\sum\limits_{j\neq i}v_j\left(\frac{w_j}{v_j}\right)_x\mu_i v_{j,xx}\pa_{v_j}\bar{v}_j\xi_j^\p\psi_{ji,v}\\
		&-\sum\limits_{j\neq i}v_j\left(\frac{w_j}{v_j}\right)_x\mu_i v_{j,x}\pa_x(\pa_{v_j}\bar{v}_j)\xi_j^\p\psi_{ji,v}-\sum\limits_{j\neq i}v_j\left(\frac{w_j}{v_j}\right)^2_x\mu_i v_{j,x}\pa_{v_j}\bar{v}_j\xi_j^{\p\p}\psi_{ji,v}\\
		&-\sum\limits_{j\neq i}v_j\left(\frac{w_j}{v_j}\right)_x\mu_i v_{j,x}\pa_{v_j}\bar{v}_j\xi_j^\p\pa_x(\psi_{ji,v})\\
		&=O(1)\sum\limits_{j\neq i}\left(\abs{w_{j,x}}+\abs{v_{j,x}}+\abs{v_j}\left(\frac{w_j}{v_j}\right)_x^2\right)\mathfrak{A}_j\rho_j^\e.
		\end{align*}
		Applying again Lemmas \ref{lemme6.5}, \ref{lemme6.6}, \ref{estimimpo1}, \ref{lemme11.3a}
		and \eqref{ngtech5} we get
		\begin{align*}
		\mathcal{I}_{4,5,x}&=-\sum\limits_{j\neq i}v_{j,x}\left(\frac{w_j}{v_j}\right)_x\mu_i\xi_j^\p\psi_{ji,v}\sum\limits_{l\neq j}
		v_{l,x}\pa_{v_l}\bar{v}_j-\sum\limits_{j\neq i}v_j\left(\frac{w_j}{v_j}\right)_{xx}\mu_i\xi_j^\p\psi_{ji,v}\sum\limits_{l\neq j}v_{l,x}\pa_{v_l}\bar{v}_j\\
		&-\sum\limits_{j\neq i}v_j\left(\frac{w_j}{v_j}\right)_x\mu_{i,x}\xi_j^\p\psi_{ji,v}\sum\limits_{l\neq j}v_{l,x}\pa_{v_l}\bar{v}_j-\sum\limits_{j\neq i}v_j\left(\frac{w_j}{v_j}\right)^2_x\mu_i\xi_j^{\p\p}\psi_{ji,v}\sum\limits_{l\neq j}v_{l,x}\pa_{v_l}\bar{v}_j\\
		&-\sum\limits_{j\neq i}v_j\left(\frac{w_j}{v_j}\right)_x\mu_i\xi_j^\p\pa_x(\psi_{ji,v})\sum\limits_{l\neq j}v_{l,x}\pa_{v_l}\bar{v}_j-\sum\limits_{j\neq i}v_j\left(\frac{w_j}{v_j}\right)_x\mu_i\xi_j^\p\psi_{ji,v}\sum\limits_{l\neq j}v_{l,xx}\pa_{v_l}\bar{v}_j\\
		&-\sum\limits_{j\neq i}v_j\left(\frac{w_j}{v_j}\right)_x\mu_i\xi_j^\p\psi_{ji,v}\sum\limits_{l\neq j}v_{l,x}\pa_x(\pa_{v_l}\bar{v}_j)\\
		&=O(1)\sum\limits_{j\neq i}\left(|w_{j,xx}|+|v_{j,xx}|+\abs{w_{j,x}}+\abs{v_{j,x}}+|v_j\left(\frac{w_j}{v_j}\right)_x^2|\right)\mathfrak{A}_j\rho_j^\e.
		\end{align*}
		Similarly using the fact that $\psi_{ji}, \psi_{ji,\sig}=O(1)\xi_j\bar{v}_j$ we have from Lemmas \ref{lemme11.3a}, \ref{estimimpo1}
		\begin{align*}
		\mathcal{I}_{4,6,x}&=-\sum\limits_{j\neq i}v_{j,x}\left(\frac{w_j}{v_j}\right)_x\mu_{i,x} [-\psi_{ji,\si}+\bar{v}_j\xi_j^\p\psi_{ji,v}]-\sum\limits_{j\neq i}v_j\left(\frac{w_j}{v_j}\right)_{xx}\mu_{i,x} [-\psi_{ji,\si}+\bar{v}_j\xi_j^\p\psi_{ji,v}]\\
		&-\sum\limits_{j\neq i}v_j\left(\frac{w_j}{v_j}\right)_x\mu_{i,xx} [-\psi_{ji,\si}+\bar{v}_j\xi_j^\p\psi_{ji,v}]-\sum\limits_{j\neq i}v_j\left(\frac{w_j}{v_j}\right)^2_x\mu_{i,x} [\bar{v}_j\xi_j^{\p\p}\psi_{ji,v}]\\
		&-\sum\limits_{j\neq i}v_j\left(\frac{w_j}{v_j}\right)_x\mu_{i,x} [-\pa_x(\psi_{ji,\si})+\bar{v}_j\xi_j^\p\pa_x(\psi_{ji,v})+\pa_x(\bar{v}_j)\xi_j^\p\psi_{ji,v}]\\
		&=O(1)\sum\limits_{j\neq i}\left(|v_j|+\abs{w_{j,x}}+\abs{v_{j,x}}\right)\mathfrak{A}_j\rho_j^\e.
		\end{align*}
		Applying Lemmas \ref{estimimpo1}, \ref{lemme6.5}, \ref{lemme6.6}, \ref{lemme11.3a}, \ref{lemme11.5} and \eqref{ngtech5} we obtain
		\begin{align*}
		\mathcal{I}_{4,7,x}&=-\sum\limits_{j\neq i}(v_{j,x}\left(\frac{w_j}{v_j}\right)_x+v_j\left(\frac{w_j}{v_j}\right)_{xx})
		\mu_i [-\pa_x(\psi_{ji,\si})+\bar{v}_j\xi_j^\p\pa_x(\psi_{ji,v})]\\
		&-\sum\limits_{j\neq i}v_j\left(\frac{w_j}{v_j}\right)_x\mu_{i,x} [-\pa_x(\psi_{ji,\si})+\bar{v}_j\xi_j^\p\pa_x(\psi_{ji,v})]-\sum\limits_{j\neq i}v_j\left(\frac{w_j}{v_j}\right)^2_x\mu_i [\bar{v}_j\xi_j^{\p\p}\pa_x(\psi_{ji,v})]\\
		&-\sum\limits_{j\neq i}v_j\left(\frac{w_j}{v_j}\right)_x\mu_i [-\pa_{xx}(\psi_{ji,\si})+\bar{v}_j\xi_j^\p\pa_{xx}(\psi_{ji,v})+\pa_x(\bar{v}_j)\xi_j^\p\pa_x(\psi_{ji,v})]\\
		&=O(1)\sum\limits_{j\neq i}\left(|w_{j,xx}|+|v_{j,xx}|+\abs{w_{j,x}}+\abs{v_{j,x}}+|v_j\left(\frac{w_j}{v_j}\right)_x^2|\right)\mathfrak{A}_j\rho_j^\e.
	\end{align*}
	\noi\underline{{Estimate of $\mathcal{I}_{5,x}$:}}
	Now, we estimate the terms of $\mathcal{I}_{5,k,x}$ for $1\leq k\leq 7$. Using the Lemmas \ref{lemme6.5}, \ref{lemme6.6}, \ref{estimimpo1}, \ref{lemme11.3}, \ref{lemme11.3a} with the fact that $\tilde{v}_j =\pa_{v_j}\bar{v}_j v_j-\bar{v}_j$, it yields
	\begin{align*}
		&\mathcal{I}_{5,1,x}=-\sum\limits_{j\neq i}\pa_{xx}(\mu_j^{-1})(\tilde{\la}_j-\si_j)v_j\xi_j\tilde{v}_j\mu_i\psi_{ji,v}-\sum\limits_{j\neq i}\pa_x(\mu_j^{-1})\left(\tilde{\la}_{j,x}+\theta_j^\p\left(\frac{w_j}{v_j}\right)_x\right)v_j\xi_j  \tilde{v}_j\mu_i\psi_{ji,v}\\
		&-\sum\limits_{j\neq i}\pa_x(\mu_j^{-1})(\tilde{\la}_j-\si_j)v_{j,x}\xi_j \tilde{v}_j\mu_i\psi_{ji,v}-\sum\limits_{j\neq i}\pa_x(\mu_j^{-1})(\tilde{\la}_j-\si_j)v_j\xi_j^\p\left(\frac{w_j}{v_j}\right)_x \tilde{v}_j\mu_i\psi_{ji,v}\\
		&-\sum\limits_{j\neq i}\pa_x(\mu_j^{-1})(\tilde{\la}_j-\si_j)v_j\xi_j\pa_x\tilde{v}_j\mu_i\psi_{ji,v}-\sum\limits_{j\neq i}\pa_x(\mu_j^{-1})(\tilde{\la}_j-\si_j)v_j\xi_j \tilde{v}_j[\mu_{i,x}\psi_{ji,v}+\mu_i\pa_x(\psi_{ji,v})]\\
		&=O(1)\sum\limits_{j\neq i}|v_j|\mathfrak{A}_j \rho_j^\e.
		\end{align*}
		Applying now the Lemmas \ref{lemme6.5}, \ref{lemme6.6}, \ref{estimimpo1}, \ref{lemme11.3a}, \ref{lemme11.6}
		\begin{align*}
		\mathcal{I}_{5,2,x}&=-\sum\limits_{j\neq i}\pa_x(\mu_j^{-1})\left(\tilde{\la}_{j,x}+\theta_j^\p\left(\frac{w_j}{v_j}\right)_x\right)v_j\xi_j \tilde{v}_j\mu_i\psi_{ji,v}\\
		&- \sum\limits_{j\neq i}\mu_j^{-1}\left(\tilde{\la}_{j,xx}+\theta_j^{\p\p}\left(\frac{w_j}{v_j}\right)^2_x+\theta_j^\p\left(\frac{w_j}{v_j}\right)_{xx}\right)v_j\xi_j \tilde{v}_j\mu_i\psi_{ji,v}\\
		&- \sum\limits_{j\neq i}\mu_j^{-1}\left(\tilde{\la}_{j,x}+\theta_j^\p\left(\frac{w_j}{v_j}\right)_x\right)v_{j,x}\xi_j\tilde{v}_j\mu_i\psi_{ji,v}\\
		&- \sum\limits_{j\neq i}\mu_j^{-1}\left(\tilde{\la}_{j,x}+\theta_j^\p\left(\frac{w_j}{v_j}\right)_x\right)v_j\xi_j^\p\left(\frac{w_j}{v_j}\right)_x\tilde{v}_j\mu_i\psi_{ji,v}\\
		&-\sum\limits_{j\neq i}\mu_j^{-1}\left(\tilde{\la}_{j,x}+\theta_j^\p\left(\frac{w_j}{v_j}\right)_x\right)v_j\xi_j\pa_x(\tilde{v}_j)\mu_i\psi_{ji,v}\\
		&- \sum\limits_{j\neq i}\mu_j^{-1}\left(\tilde{\la}_{j,x}+\theta_j^\p\left(\frac{w_j}{v_j}\right)_x\right)v_j\xi_j \tilde{v}_j[\mu_{i,x}\psi_{ji,v}+\mu_i\pa_x(\psi_{ji,v})]\\
		&=O(1)\sum\limits_{j\neq i}\left(|v_j|+|w_{j,x}|+|v_{j,x}|\right)\mathfrak{A}_j \rho_j^\e.
		\end{align*}
		Applying again Lemmas \ref{lemme6.5}, \ref{lemme6.6}, \ref{estimimpo1}, \ref{lemme11.3a}
		\begin{align*}
		\mathcal{I}_{5,3,x}&=-\sum\limits_{j\neq i}\pa_x(\mu_j^{-1})(\tilde{\la}_j-\si_j)v_{j,x}\xi_j \tilde{v}_j\mu_i\psi_{ji,v}-\sum\limits_{j\neq i}\mu_j^{-1}\left(\tilde{\la}_{j,x}+\theta_j^\p\left(\frac{w_j}{v_j}\right)_x\right)v_{j,x}\xi_j \tilde{v}_j\mu_i\psi_{ji,v}\\
		&-\sum\limits_{j\neq i}\mu_j^{-1}(\tilde{\la}_j-\si_j)v_{j,xx}\xi_j \tilde{v}_j\mu_i\psi_{ji,v}-\sum\limits_{j\neq i}\mu_j^{-1}(\tilde{\la}_j-\si_j)v_{j,x}\xi_j^\p\left(\frac{w_j}{v_j}\right)_x \tilde{v}_j\mu_i\psi_{ji,v}\\
		&-\sum\limits_{j\neq i}\mu_j^{-1}(\tilde{\la}_j-\si_j)v_{j,x}\xi_j\pa_x(\tilde{v}_j)\mu_i\psi_{ji,v}-\sum\limits_{j\neq i}\mu_j^{-1}(\tilde{\la}_j-\si_j)v_{j,x}\xi_j \tilde{v}_j[\mu_{i,x}\psi_{ji,v}+\mu_i\pa_x(\psi_{ji,v})]\\
		&=O(1)\sum\limits_{j\neq i}\left(\abs{v_j}+|v_{j,x}|\right)\mathfrak{A}_j \rho_j^\e.
		\end{align*}
		\begin{align*}
		\mathcal{I}_{5,4,x}&=-\sum\limits_{j\neq i}\pa_x(\mu_j^{-1})(\tilde{\la}_j-\si_j)v_j\xi_j^\p\left(\frac{w_j}{v_j}\right)_x \tilde{v}_j\mu_i\psi_{ji,v}\\
		&-\sum\limits_{j\neq i}\mu_j^{-1}\left(\tilde{\la}_{j,x}+\theta_j^\p\left(\frac{w_j}{v_j}\right)_x\right)v_j\xi_j^\p\left(\frac{w_j}{v_j}\right)_x \tilde{v}_j\mu_i\psi_{ji,v}\\
&-\sum\limits_{j\neq i}\mu_j^{-1}(\tilde{\la}_j-\si_j)v_{j,x}\xi_j^\p\left(\frac{w_j}{v_j}\right)_x \tilde{v}_j\mu_i\psi_{ji,v}-\sum\limits_{j\neq i}\mu_j^{-1}(\tilde{\la}_j-\si_j)v_j\xi_j^{\p\p}\left(\frac{w_j}{v_j}\right)^2_x \tilde{v}_j\mu_i\psi_{ji,v}\\
		&-\sum\limits_{j\neq i}\mu_j^{-1}(\tilde{\la}_j-\si_j)v_j\xi_j^\p\left(\frac{w_j}{v_j}\right)_{xx}\tilde{v}_j\mu_i\psi_{ji,v}-\sum\limits_{j\neq i}\mu_j^{-1}(\tilde{\la}_j-\si_j)v_j\xi_j^\p\left(\frac{w_j}{v_j}\right)_x\pa_x(\tilde{v}_j)\mu_i\psi_{ji,v}\\
		&-\sum\limits_{j\neq i}\mu_j^{-1}(\tilde{\la}_j-\si_j)v_j\xi_j^\p\left(\frac{w_j}{v_j}\right)_x \tilde{v}_j[\mu_{i,x}\psi_{ji,v}+\mu_i\pa_x(\psi_{ji,v})]\\
		&=O(1)\sum\limits_{j\neq i}\left(\abs{v_j}+|v_{j,x}|+|w_{j,x}|\right)\mathfrak{A}_j \rho_j^\e.
		\end{align*}
		Using Lemma \ref{estimimpo1}, \ref{lemme6.5}, \ref{lemme6.6}, \ref{lemme11.3a}, we get
		\begin{align*}
		\mathcal{I}_{5,5,x}&=-
		\sum\limits_{j\neq i}\pa_x(\mu_j^{-1})(\tilde{\la}_j-\si_j)v_j\xi_jv_{j,x}\pa_{v_j}\tilde{v}_j\mu_i\psi_{ji,v}\\
		&-\sum\limits_{j\neq i}\mu_j^{-1}\left(\tilde{\la}_{j,x}+\theta_j^\p\left(\frac{w_j}{v_j}\right)_x\right)v_j\xi_jv_{j,x}\pa_{v_j}\tilde{v}_j \mu_i\psi_{ji,v}\\
		&-\sum\limits_{j\neq i}\mu_j^{-1}(\tilde{\la}_j-\si_j)v^2_{j,x}\xi_j\pa_{v_j}\tilde{v}_j \mu_i\psi_{ji,v}-\sum\limits_{j\neq i}\mu_j^{-1}(\tilde{\la}_j-\si_j)v_j\xi_{j}^\p\left(\frac{w_j}{v_j}\right)_xv_{j,x}\pa_{v_j}\tilde{v}_j \mu_i\psi_{ji,v}\\
		&-\sum\limits_{j\neq i}\mu_j^{-1}(\tilde{\la}_j-\si_j)v_j\xi_jv_{j,xx}\pa_{v_j}\tilde{v}_j\mu_i\psi_{ji,v}-\sum\limits_{j\neq i}\mu_j^{-1}(\tilde{\la}_j-\si_j)v_j\xi_jv_{j,x}\pa_x(\pa_{v_j}\tilde{v}_j) \mu_i\psi_{ji,v}\\
		&-\sum\limits_{j\neq i}\mu_j^{-1}(\tilde{\la}_j-\si_j)v_j\xi_jv_{j,x}\pa_{v_j}\tilde{v}_j [\mu_{i,x}\psi_{ji,v}+\mu_i\pa_x(\psi_{ji,v})]\\
		&=O(1)\sum\limits_{j\neq i}\left(\abs{v_j}+|v_{j,x}|\right)\mathfrak{A}_j \rho_j^\e.
		\end{align*}
		Similarly we get from Lemmas \ref{lemme6.5}, \ref{lemme6.6}, \ref{estimimpo1}, \ref{lemme11.3a}
		\begin{align*}
		&\mathcal{I}_{5,6,x}=-\sum\limits_{j\neq i}\pa_x(\mu_j^{-1})(\tilde{\la}_j-\si_j)v_j\xi_j\sum\limits_{l\neq j}v_{l,x}\pa_{v_l}\tilde{v}_j\mu_i\psi_{ji,v}\\
		&-\sum\limits_{j\neq i}\mu_j^{-1}\left(\tilde{\la}_{j,x}+\theta_j^\p\left(\frac{w_j}{v_j}\right)_x\right)v_j\xi_j\sum\limits_{l\neq j}v_{l,x}\pa_{v_l}\tilde{v}_j\mu_i\psi_{ji,v}\\
		&-\sum\limits_{j\neq i}\mu_j^{-1}(\tilde{\la}_j-\si_j)v_{j,x}\xi_j\sum\limits_{l\neq j}v_{l,x}\pa_{v_l}\tilde{v}_j\mu_i\psi_{ji,v}-\sum\limits_{j\neq i}\mu_j^{-1}(\tilde{\la}_j-\si_j)v_j\xi_j\sum\limits_{l\neq j}v_{l,x}\pa_x( \pa_{v_l}\tilde{v}_j)\mu_i\psi_{ji,v}\\
		&-\sum\limits_{j\neq i}\mu_j^{-1}(\tilde{\la}_j-\si_j)v_j\xi_j^\p\left(\frac{w_j}{v_j}\right)_x\sum\limits_{l\neq j}v_{l,x}\pa_{v_l}\tilde{v}_j\mu_i\psi_{ji,v}-\sum\limits_{j\neq i}\mu_j^{-1}(\tilde{\la}_j-\si_j)v_j\xi_j\sum\limits_{l\neq j}v_{l,xx}\pa_{v_l}\tilde{v}_j\mu_i\psi_{ji,v}\\
		&-\sum\limits_{j\neq i}\mu_j^{-1}(\tilde{\la}_j-\si_j)v_j\xi_j\sum\limits_{l\neq j}v_{l,x}\pa_{v_l}\tilde{v}_j[\mu_{i,x}\psi_{ji,v}+\mu_i\pa_x(\psi_{ji,v})]\\
		&=O(1)\sum\limits_{j\neq i}\left(\abs{v_j}+|v_{j,x}|+|w_{j,x}|\right)\mathfrak{A}_j \rho_j^\e.
		\end{align*}
		We have now using Lemmas \ref{lemme6.5}, \ref{lemme6.6}, \ref{estimimpo1}, \ref{lemme11.3a}, \ref{lemme11.5}
		\begin{align*}
		\mathcal{I}_{5,7,x}&=-\sum\limits_{j\neq i}\pa_x(\mu_j^{-1})(\tilde{\la}_j-\si_j)v_j\xi_j\tilde{v}_j[\mu_{i,x}\psi_{ji,v}+\mu_i\pa_x\psi_{ji,v}]\\
		&-\sum\limits_{j\neq i}\mu_j^{-1}\left(\tilde{\la}_{j,x}+\theta_j^\p\left(\frac{w_j}{v_j}\right)_x\right)v_j\xi_j \tilde{v}_j[\mu_{i,x}\psi_{ji,v}+\mu_i\pa_x\psi_{ji,v}]\\
		&-\sum\limits_{j\neq i}\mu_j^{-1}(\tilde{\la}_j-\si_j)v_{j,x}\xi_j \tilde{v}_j[\mu_{i,x}\psi_{ji,v}+\mu_i\pa_x\psi_{ji,v}]\\
		&-\sum\limits_{j\neq i}\mu_j^{-1}(\tilde{\la}_j-\si_j)v_j\xi_j^\p\left(\frac{w_j}{v_j}\right)_x\tilde{v}_j[\mu_{i,x}\psi_{ji,v}+\mu_i\pa_x\psi_{ji,v}]\\
		&-\sum\limits_{j\neq i}\mu_j^{-1}(\tilde{\la}_j-\si_j)v_j\xi_j\tilde{v}_j[\mu_{i,xx}\psi_{ji,v}+2\mu_{i,x}\pa_x\psi_{ji,v}+\mu_{i}\pa_{xx}\psi_{ji,v}]\\
		&=O(1)\sum\limits_{j\neq i}\left(\abs{v_j}+|v_{j,x}|+|w_{j,x}|\right)\mathfrak{A}_j \rho_j^\e.
	\end{align*}
	\noi\underline{{Estimate of $\mathcal{I}_{6,x}$:}} We now estimate $\mathcal{I}_{6,l,x}$ for $1\leq l\leq 7$. Applying Lemmas \ref{lemme6.5}, \ref{lemme6.6}, \ref{lemme11.3}, \ref{lemme11.4} and using the fact that $\chi_k^\e \pa_{v_k}\eta_k=O(1)\frac{\Delta_k^\e}{|v_k|}$, $\chi_k^\e  \pa_{v_l}\eta_k=O(1)|\frac{v_l}{v_k}|\Delta_{kl}^\e $ when $l\ne k$ we get
		\begin{align*}
		\mathcal{I}_{6,1,x}&=	-\sum\limits_{k}2(v^2_{k,x}+v_kv_{k,xx})(1-\chi_k^\e\xi_k\eta_k)\langle l_i, \left[B\tilde{r}_{k,u}\tilde{r}_k+\tilde{r}_k\cdot DB\tilde{r}_k\right]\rangle \\
		&+4N\sum\limits_{k}\frac{v_k^{2N}}{\e}v_{k,x}^2(\chi')^\e_k\xi_k \eta_k\langle l_i, \left[B\tilde{r}_{k,u}\tilde{r}_k+\tilde{r}_k\cdot DB\tilde{r}_k\right]\rangle \\
		&+\sum\limits_{k}2v_kv_{k,x}\xi^\p_k\left(\frac{w_k}{v_k}\right)_x\chi_k^\e \eta_k\langle l_i, \left[B\tilde{r}_{k,u}\tilde{r}_k+\tilde{r}_k\cdot DB\tilde{r}_k\right]\rangle \\
		&+\sum\limits_{k}2v_kv_{k,x}^2\xi_k\chi_k^\e \pa_{v_k}\eta_k\langle l_i, \left[B\tilde{r}_{k,u}\tilde{r}_k+\tilde{r}_k\cdot DB\tilde{r}_k\right]\rangle \\
		&+\sum\limits_{k}\sum\limits_{l\neq k}2v_kv_{k,x}v_{l,x}\xi_k\chi_k^\e \pa_{v_l}\eta_k\langle l_i, \left[B\tilde{r}_{k,u}\tilde{r}_k+\tilde{r}_k\cdot DB\tilde{r}_k\right]\rangle \\
		&-\sum\limits_{k}2v_kv_{k,x}(1-\chi_k^\e\xi_k\eta_k)\langle u_x\cdot Dl_i, \left[B\tilde{r}_{k,u}\tilde{r}_k+\tilde{r}_k\cdot DB(u)\tilde{r}_k\right]\rangle \\
		&-\sum\limits_{k}2v_kv_{k,x}(1-\chi_k^\e\xi_k\eta_k)\langle l_i, \left[u_x\cdot DB\tilde{r}_{k,u}\tilde{r}_k+\tilde{r}_k\otimes u_x: D^2B\tilde{r}_k\right]\rangle \\
		&-\sum\limits_{k}2v_kv_{k,x}(1-\chi_k^\e\xi_k\eta_k)\langle l_i, \left[B\pa_x(\tilde{r}_{k,u}\tilde{r}_k)+\tilde{r}_{k,x}\cdot DB(u)\tilde{r}_k+\tilde{r}_k\cdot DB(u)\tilde{r}_{k,x}\right]\rangle \\
		&=O(1)\biggl(\sum_{j\ne i}(|v_j|+|v_{j,x}|)+(|v_{i,x}^2|+|v_iv_{i,x}|+|v_iv_{i,xx}|)(1-\chi_i^\e\xi_i\eta_i)
		\\
		&+\mathfrak{A}_i\rho_i^\e\big(|v_iv_{i,x}\left(\frac{w_i}{v_i}\right)_x|
		+v_{i,x}^2\Delta_i^\e\big)+\sum_k\sum_{l\ne k}|v_{k,x}||v_lv_{l,x}|
		+v_{i,x}^2(\chi')^\e_i\xi_i\eta_i \biggl).
		\end{align*}
		Similarly since $\pa_x\biggl[\langle l_i, \left[B(u)\tilde{r}_{k,u}\tilde{r}_k+\tilde{r}_k\cdot DB(u)\tilde{r}_k\right]=O(1)$ we have using Lemmas \ref{lemme6.5}, \ref{lemme6.6}, \ref{estimimpo1} 
		\begin{align*}
		\mathcal{I}_{6,2,x}&=\sum\limits_{k}\frac{2N(2N+1) v_k^{2N}}{\e}v_{k,x}^2 (\chi')_k^\e \xi_k\eta_k\langle l_i, \left[B(u)\tilde{r}_{k,u}\tilde{r}_k+\tilde{r}_k\cdot DB(u)\tilde{r}_k\right]\rangle \\
		&+\sum\limits_{k}\frac{2N v_k^{2N+1}}{\e}v_{k,xx} (\chi')_k^\e \xi_k\eta_k\langle l_i, \left[B(u)\tilde{r}_{k,u}\tilde{r}_k+\tilde{r}_k\cdot DB(u)\tilde{r}_k\right]\rangle \\
		&+\sum\limits_{k}\frac{4N^2 v_k^{4N}}{\e^2}v^2_{k,x} (\chi'')_k^\e \xi_k\eta_k\langle l_i, \left[B(u)\tilde{r}_{k,u}\tilde{r}_k+\tilde{r}_k\cdot DB(u)\tilde{r}_k\right]\rangle \\
		&+\sum\limits_{k}\frac{2N v_k^{2N+1}}{\e}v_{k,x} (\chi')_k^\e \xi'_k\left(\frac{w_k}{v_k}\right)_x\eta_k\langle l_i, \left[B(u)\tilde{r}_{k,u}\tilde{r}_k+\tilde{r}_k\cdot DB(u)\tilde{r}_k\right]\rangle \\
		&+\sum\limits_{k}\frac{2N v_k^{2N+1}}{\e}v^2_{k,x} (\chi')_k^\e \xi_k\pa_{v_k}\eta_k\langle l_i, \left[B(u)\tilde{r}_{k,u}\tilde{r}_k+\tilde{r}_k\cdot DB(u)\tilde{r}_k\right]\rangle \\
		&+\sum\limits_{k}\sum_{l\ne k}\frac{2N v_k^{2N+1}}{\e}v_{k,x} (\chi')_k^\e \xi_k\pa_{v_l}\eta_k v_{l,x}\langle l_i, \left[B(u)\tilde{r}_{k,u}\tilde{r}_k+\tilde{r}_k\cdot DB(u)\tilde{r}_k\right]\rangle \\
		&+\sum\limits_{k}\frac{2N v_k^{2N+1}}{\e}v_{k,x} (\chi')_k^\e \xi_k\eta_k\pa_x\biggl[\langle l_i, \left[B(u)\tilde{r}_{k,u}\tilde{r}_k+\tilde{r}_k\cdot DB(u)\tilde{r}_k\right]\rangle \biggl]\\
		&=O(1)\biggl(\sum_{j\ne i}(|v_{j,x}|+|v_j|)+(|v_i |(|v_{i,x}|+|v_{i,xx}|)+|v_{i,x}|(|v_{i,x}|+|w_{i,x}|))\varkappa^\ep_i\mathfrak{A}_i\rho_i^\e\\
		&+\sum_k\sum_{l\ne k}|v_{k,x}||v_lv_{l,x}|\biggl).
		\end{align*}
		Similarly we have using Lemmas \ref{lemme6.5}, \ref{lemme6.6}, \ref{lemme11.3},\ref{lemme11.4}
		\begin{align*}
		\mathcal{I}_{6,3,x}&=\sum\limits_{k}2v_kv_{k,x}\chi_k^\e\xi^\p_k\eta_k\left(\frac{w_k}{v_k}\right)_x\langle l_i, \left[B\tilde{r}_{k,u}\tilde{r}_k+\tilde{r}_k\cdot DB\tilde{r}_k\right]\rangle \\
		&+\sum\limits_{k}v_k^2v_{k,x}\frac{2N v_k^{2N-1}}{\e}(\chi')_k^\e\xi^\p_k\eta_k\left(\frac{w_k}{v_k}\right)_x\langle l_i, \left[B\tilde{r}_{k,u}\tilde{r}_k+\tilde{r}_k\cdot DB\tilde{r}_k\right]\rangle \\
		&+\sum\limits_{k}v_k^2\chi_k^\e \xi^{\p\p}_k\eta_k\left(\frac{w_k}{v_k}\right)^2_x\langle l_i, \left[B\tilde{r}_{k,u}\tilde{r}_k+\tilde{r}_k\cdot DB\tilde{r}_k\right]\rangle \\
		&+\sum\limits_{k}v_k^2v_{k,x}\chi_k^\e \xi^\p_k\pa_{v_k}\eta_k\left(\frac{w_k}{v_k}\right)_x\langle l_i, \left[B\tilde{r}_{k,u}\tilde{r}_k+\tilde{r}_k\cdot DB\tilde{r}_k\right]\rangle \\
		&+\sum\limits_{k}\sum\limits_{l\neq k}v_k^2v_{l,x}\chi_k^\e \xi^\p_k\pa_{v_l}\eta_k\left(\frac{w_k}{v_k}\right)_x\langle l_i, \left[B\tilde{r}_{k,u}\tilde{r}_k+\tilde{r}_k\cdot DB\tilde{r}_k\right]\rangle \\
		&+\sum\limits_{k}v_k^2 \chi_k^\e\xi^\p_k\eta_k\left(\frac{w_k}{v_k}\right)_{xx}\langle l_i, \left[B\tilde{r}_{k,u}\tilde{r}_k+\tilde{r}_k\cdot DB\tilde{r}_k\right]\rangle \\
		&+\sum\limits_{k}v_k^2\chi_k^\e \xi^\p_k\eta_k\left(\frac{w_k}{v_k}\right)_x
		\pa_x \biggl[\langle l_i, \left[B\tilde{r}_{k,u}\tilde{r}_k+\tilde{r}_k\cdot DB\tilde{r}_k\right]\rangle\biggl]\\
		&=O(1)\biggl(\sum_{j\ne i}(|v_{j,x}|+|v_{j}|)+\mathfrak{A}_i\rho_i^\e\big(|v_iv_{i,x}\left(\frac{w_i}{v_i}\right)_x|+|v_i^2\left(\frac{w_i}{v_i}\right)^2_x|+|v_i^2\left(\frac{w_i}{v_i}\right)_x|\\
		&+|w_{i,xx}v_i-v_{i,xx}w_i|\big)+\sum_k\sum_{l\ne k}|v_{l,x}v_l|(|v_{k,x}|+|w_{k,x}|)\biggl).
		\end{align*}
	Using again the Lemmas \ref{lemme6.5}, \ref{lemme6.6} and the fact that  $\chi_k^\e \pa_{v_k}\eta_k=O(1)\frac{\Delta_k^\e}{|v_k|}$, $\chi_k^\e  \pa_{v_l}\eta_k=O(1)|\frac{v_l}{v_k}|\Delta_{kl}^\e $ when $l\ne k$, $\chi_k^\e \pa_{v_kv_k}\eta_k=O(1)\frac{\Delta_k^\e}{v_k^2}$,  $\chi_k^\e \pa_{v_l v_k}\eta_k=O(1)|\frac{v_l}{v_k^2}|\rho_k^\e$ when $l\ne k$, $\chi_k^\e \pa_{v_l v_p}\eta_k=O(1)\Box^\ep_{kpl} |\frac{v_lv_p}{v_k^2}|\rho_k^\e$
	when $l,p\ne k$ we obtain
		\begin{align*}
		&\mathcal{I}_{6,4,x}=\sum\limits_{k}2v_kv_{k,x}\chi_k^\e \xi_k\left[v_{k,x}\pa_{v_k}\eta_k+\sum\limits_{l\neq k}v_{l,x}\pa_{v_l}\eta_k\right]\langle l_i, \left[B\tilde{r}_{k,u}\tilde{r}_k+\tilde{r}_k\cdot DB\tilde{r}_k\right]\rangle\\
		&+2N \sum\limits_{k}\frac{v_k^{2N+1}}{\e}v_{k,x} (\chi')_k^\e \xi_k\left[v_{k,x}\pa_{v_k}\eta_k+\sum\limits_{l\neq k}v_{l,x}\pa_{v_l}\eta_k\right]\langle l_i, \left[B\tilde{r}_{k,u}\tilde{r}_k+\tilde{r}_k\cdot DB\tilde{r}_k\right]\rangle\\
		&+\sum\limits_{k}v_k^2\chi_k^\e \xi_k^\p\left(\frac{w_k}{v_k}\right)_x\left[v_{k,x}\pa_{v_k}\eta_k+\sum\limits_{l\neq k}v_{l,x}\pa_{v_l}\eta_k\right]\langle l_i, \left[B\tilde{r}_{k,u}\tilde{r}_k+\tilde{r}_k\cdot DB\tilde{r}_k\right]\rangle\\
		&+\sum\limits_{k}v_k^2\chi_k^\e \xi_k\left[v_{k,xx}\pa_{v_k}\eta_k+v_{k,x}\pa_x(\pa_{v_k}\eta_k)\right]\langle l_i, \left[B\tilde{r}_{k,u}\tilde{r}_k+\tilde{r}_k\cdot DB\tilde{r}_k\right]\rangle\\
		&+\sum\limits_{k}v_k^2\chi_k^\e \xi_k\left[\sum\limits_{l\neq k}(v_{l,xx}\pa_{v_l}\eta_k+v_{l,x}\pa_x(\pa_{v_l}\eta_k))\right]\langle l_i, \left[B\tilde{r}_{k,u}\tilde{r}_k+\tilde{r}_k\cdot DB\tilde{r}_k\right]\rangle\\
		&+\sum\limits_{k}v_k^2\chi_k^\e \xi_k\left[v_{k,x}\pa_{v_k}\eta_k+\sum\limits_{l\neq k}v_{l,x}\pa_{v_l}\eta_k\right]\pa_x\biggl[\langle l_i, \left[B\tilde{r}_{k,u}\tilde{r}_k+\tilde{r}_k\cdot DB\tilde{r}_k\right]\rangle\biggl]\\
		&=O(1)\biggl(\sum_{j\ne i}(|v_j|+|v_{j,x}|+|w_{j,x}|)+\sum_k\sum_{l\ne k}(|v_{k,x}|+|w_{k,x}|+|v_k|)|v_lv_{l,x}|\\
		&+\sum_{k}\sum_{l\ne k}|v_k||v_{l,xx}v_l|
		+\sum_{k}\sum_{l\ne k}\sum_{p\ne k}|v_k| |v_{l,x}||v_{p,x}|\\
		&+\mathfrak{A}_i \Delta_i^\e\big(|v_iv_{i,x}\left(\frac{w_i}{v_i}\right)_x|
		+|v_i v_{i,x}|+v_{i,x}^2+|v_iv_{i,xx}|\big)\biggl).
		\end{align*}
		We have now using Lemmas \ref{lemme6.5}, \ref{lemme6.6}, \ref{lemme11.3}, \ref{lemme11.4}
		\begin{align*}
		\mathcal{I}_{6,5,x}&=-\sum\limits_{k}\sum\limits_{l}2v_{k,x}v_lv_k(1-\chi_k^\e\xi_k\eta_k)\langle l_i, \left[\tilde{r}_l\cdot DB\tilde{r}_{k,u}\tilde{r}_k+\tilde{r}_k\otimes \tilde{r}_l: D^2B\tilde{r}_k\right]\rangle \\
		&-\sum\limits_{k}\sum\limits_{l}v_{k}^2v_{l,x}(1-\chi_k^\e\xi_k\eta_k)\langle l_i, \left[\tilde{r}_l\cdot DB\tilde{r}_{k,u}\tilde{r}_k+\tilde{r}_k\otimes \tilde{r}_l: D^2B\tilde{r}_k\right]\rangle \\
		&+\sum\limits_{k}\sum\limits_{l}v_{k}^2v_{l}\frac{v_k^{2N-1}}{\e}v_{k,x}(\chi')_k^\e\xi_k\eta_k\langle l_i, \left[\tilde{r}_l\cdot DB\tilde{r}_{k,u}\tilde{r}_k+\tilde{r}_k\otimes \tilde{r}_l: D^2B\tilde{r}_k\right]\rangle \\
&+\sum\limits_{k}\sum\limits_{l}v_k^2v_l\chi_k^\e\xi_k^\p\left(\frac{w_k}{v_k}\right)_x\eta_k\langle l_i, \left[\tilde{r}_l\cdot DB\tilde{r}_{k,u}\tilde{r}_k+\tilde{r}_k\otimes \tilde{r}_l: D^2B\tilde{r}_k\right]\rangle \\
		&+\sum\limits_{k}\sum\limits_{l}v_k^2 v_lv_{k,x}\chi_k^\e\xi_k\pa_{v_k}\eta_k\langle l_i, \left[\tilde{r}_l\cdot DB\tilde{r}_{k,u}\tilde{r}_k+\tilde{r}_k\otimes \tilde{r}_l: D^2B\tilde{r}_k\right]\rangle \\
		&+\sum\limits_{k}\sum\limits_{l}\sum\limits_{j\neq k}v_k^2v_l v_{j,x}\chi_k^\e\xi_k\pa_{v_j}\eta_k\langle l_i, \left[\tilde{r}_l\cdot DB\tilde{r}_{k,u}\tilde{r}_k+\tilde{r}_k\otimes \tilde{r}_l: D^2B\tilde{r}_k\right]\rangle \\
		%&+\sum\limits_{k}\sum\limits_{l}\sum\limits_{j\neq k}v_k^2v_lw_{j,x}\xi_k\pa_{w_j}\eta_k\langle l_i, \left[\tilde{r}_l\cdot DB\tilde{r}_{k,u}\tilde{r}_k+\tilde{r}_k\otimes \tilde{r}_l: D^2B\tilde{r}_k\right]\rangle \\
		&-\sum\limits_{k}\sum\limits_{l}v_k^2v_l(1-\chi_k^\e\xi_k\eta_k)\langle u_x\cdot Dl_i, \left[\tilde{r}_l\cdot DB\tilde{r}_{k,u}\tilde{r}_k+\tilde{r}_k\otimes \tilde{r}_l: D^2B\tilde{r}_k\right]\rangle \\
		&-\sum\limits_{k}\sum\limits_{l}v_k^2v_l(1-\chi_k^\e\xi_k\eta_k)\langle l_i, \left[\tilde{r}_l\otimes u_x: D^2B\tilde{r}_{k,u}\tilde{r}_k+\tilde{r}_k\otimes \tilde{r}_l\otimes u_x: D^3B\tilde{r}_k\right]\rangle \\
		&-\sum\limits_{k}\sum\limits_{l}v_k^2v_l(1-\chi_k^\e\xi_k\eta_k)\langle l_i, \left[\tilde{r}_{l,x}\cdot DB\tilde{r}_{k,u}\tilde{r}_k+\tilde{r}_{k,x}\otimes \tilde{r}_l: D^2B\tilde{r}_k\right]\rangle \\
		&-\sum\limits_{k}\sum\limits_{l}v_k^2v_l(1-\chi_k^\e\xi_k\eta_k)\langle l_i, \left[\tilde{r}_l\cdot DB\pa_x(\tilde{r}_{k,u}\tilde{r}_k)+\tilde{r}_k\otimes \tilde{r}_{l,x}: D^2B\tilde{r}_k\right]\rangle \\
		&-\sum\limits_{k}\sum\limits_{l}v_k^2v_l(1-\chi_k^\e\xi_k\eta_k)\langle l_i,  \tilde{r}_k\otimes \tilde{r}_l: D^2B\tilde{r}_{k,x}\rangle \\
		&=O(1)\biggl(\sum_{j\ne i}|v_j|+\sum_{k}\sum_{l\ne k} \big(|v_k^2v_{l,x}|+(|v_k|+|w_{k,x}|+|v_{k,x}|)|v_kv_l|\big)\\
		&+\sum_k\sum_l\sum_{j\ne l}|v_kv_l||v_{j,x}v_j|
		+(|v_{i,x}|+|v_i|)v_i^2(1-\chi_i^\e\xi_i\eta_i)\\
		&+\mathfrak{A}_i\rho_i^\e\big(|v_i^3\left(\frac{w_i}{v_i}\right)_x|+
		|v_i^2v_{i,x}|\Delta_i^\e\big)+v_i^2v_{i,x}(\chi')_i^\e\xi_i\eta_i\biggl).
		\end{align*}
		Proceeding similarly, we get from Lemmas \ref{lemme6.5}, \ref{lemme6.6}, \ref{lemme11.3}, \ref{lemme11.4}
		\begin{align*}
		\mathcal{I}_{6,6,x}&=-\sum\limits_{k}\sum\limits_{l}2v_{k,x}v_lv_k(1-\chi_k^\e\xi_k\eta_k)\langle \tilde{r}_l\cdot Dl_i, \left[B \tilde{r}_{k,u}\tilde{r}_k+\tilde{r}_k\cdot DB \tilde{r}_k\right]\rangle\\
&-\sum\limits_{k}\sum\limits_{l}v_{k}^2v_{l,x}(1-\chi_k^\e\xi_k\eta_k)\langle \tilde{r}_l\cdot Dl_i, \left[B \tilde{r}_{k,u}\tilde{r}_k+\tilde{r}_k\cdot DB \tilde{r}_k\right]\rangle\\
&+2N \sum\limits_{k}\sum\limits_{l}v_lv_k^2\frac{ v_k^{2N-1}}{\e}v_{k,x}(\chi')_k^\e\xi_k\eta_k\langle \tilde{r}_l\cdot Dl_i, \left[B \tilde{r}_{k,u}\tilde{r}_k+\tilde{r}_k\cdot DB \tilde{r}_k\right]\rangle\\
		&+\sum\limits_{k}\sum\limits_{l}v_k^2v_l\chi_k^\e\xi_k^\p\left(\frac{w_k}{v_k}\right)_x\eta_k\langle \tilde{r}_l\cdot Dl_i, \left[B\tilde{r}_{k,u}\tilde{r}_k+\tilde{r}_k\cdot DB\tilde{r}_k\right]\rangle\\
		&+\sum\limits_{k}\sum\limits_{l}v_k^2v_lv_{k,x}\chi_k^\e\xi_k\pa_{v_k}\eta_k\langle \tilde{r}_l\cdot Dl_i, \left[B(u)\tilde{r}_{k,u}\tilde{r}_k+\tilde{r}_k\cdot DB(u)\tilde{r}_k\right]\rangle\\
		&+\sum\limits_{k}\sum\limits_{l}\sum\limits_{j\neq k}v_k^2v_lv_{j,x}\chi_k^\e\xi_k\pa_{v_j}\eta_k\langle \tilde{r}_l\cdot Dl_i, \left[B\tilde{r}_{k,u}\tilde{r}_k+\tilde{r}_k\cdot DB\tilde{r}_k\right]\rangle\\
		&-\sum\limits_{k}\sum\limits_{l}v_k^2v_l(1-\chi_k^\e\xi_k\eta_k)\pa_x\biggl[
		\langle \tilde{r}_l\cdot Dl_i, \left[B \tilde{r}_{k,u}\tilde{r}_k+\tilde{r}_k\cdot DB \tilde{r}_k\right]\rangle\biggl]\\
		&=O(1)\biggl(\sum_{j\ne i}|v_j|+\sum_k\sum_{l\ne k} |v_k|(|v_{k,x}v_l|+|w_{k,x}v_l|+|v_kv_{l,x}|+|v_kv_l|)\\
		&+\sum_k\sum_l\sum_{j\ne l}|v_kv_l||v_{j,x}v_j|+(|v_{i,x}|+|v_i|)v_i^2(1-\chi_i^\e\xi_i\eta_i)\\
		&+\mathfrak{A}_i\rho_i^\e\big(|v_i^3\left(\frac{w_i}{v_i}\right)_x|+|v_i^2 v_{i,x}|\Delta_i^\e\big)+v_i^2v_{i,x}(\chi')_i^\e\xi_i\eta_i\biggl).
		\end{align*}
		From the Lemmas \ref{lemme6.5}, \ref{lemme6.6}, \ref{lemme11.6}, we deduce that
		\begin{align}
		&\tilde{r}_{k,xx}=O(1)
		+O(1)\rho_k^\e\mathfrak{A}_k\bigg( |v_{k,x}\left(\frac{w_k}{v_k}\right)_x|+ \frac{v_{k,x}^2}{|v_k|}+\sum_{l\ne k}|\frac{v_lv_{k,x}v_{l,x}}{v_k}|+| \left(\frac{w_k}{v_k}\right)_x|\sum\limits_{j\neq k}|v_{j,x}v_j|\nonumber\\
		&
		+\left(\frac{w_k}{v_k}\right)_x^2|v_k|	\biggl),\label{rkxxtech1}\\
		&\pa_{xx}(\tilde{r}_{k,u}\tilde{r}_k)
		=O(1)
		+O(1)\rho_k^\e\mathfrak{A}_k\biggl( 
		  |v_k\left(\frac{w_k}{v_k}\right)^2_x|+|v_{k,x}\left(\frac{w_k}{v_k}\right)_x|+ \sum_{l\ne k}|\left(\frac{w_k}{v_k}\right)_x||v_lv_{l,x}|\nonumber\\
		&\hspace{2cm}+\frac{v_{k,x}^2}{|v_k|}+\sum_{l\ne k}|\frac{v_{l,x}v_lv_{k,x}}{v_k}|\biggl).
		\label{rkxxtech}
		\end{align}
		Applying \eqref{rkxxtech1}, \eqref{rkxxtech} and Lemmas \ref{lemme6.5}, \ref{lemme6.6}, \ref{estimimpo1}, \ref{lemme11.3}, \ref{lemme11.4}, we obtain
		\begin{align*}
		&\mathcal{I}_{6,7,x}=-2\sum\limits_{k}v_k v_{k,x}(1-\chi_k^\e\xi_k\eta_k)\langle l_i, \left[B(u)\pa_x(\tilde{r}_{k,u}\tilde{r}_{k})+\tilde{r}_k\cdot DB(u)\tilde{r}_{k,x}+\tilde{r}_{k,x}\cdot DB(u)\tilde{r}_k\right]\rangle \nonumber\\
		&+2N \sum\limits_{k}v_k^2\frac{v_k^{2N-1}}{\e}v_{k,x}(\chi')_k^\e\xi_k \eta_k\langle l_i, \left[B(u)\pa_x(\tilde{r}_{k,u}\tilde{r}_{k})+\tilde{r}_k\cdot DB(u)\tilde{r}_{k,x}+\tilde{r}_{k,x}\cdot DB(u)\tilde{r}_k\right]\rangle \\
		&+\sum\limits_{k}v_k^2\chi_k^\e\xi_k'\left(\frac{w_k}{v_k}\right)_x \eta_k\langle l_i, \left[B(u)\pa_x(\tilde{r}_{k,u}\tilde{r}_{k})+\tilde{r}_k\cdot DB(u)\tilde{r}_{k,x}+\tilde{r}_{k,x}\cdot DB(u)\tilde{r}_k\right]\rangle \\
		&+\sum\limits_{k}v_k^2\chi_k^\e\xi_k\pa_x(\eta_k)\langle l_i, \left[B(u)\pa_x(\tilde{r}_{k,u}\tilde{r}_{k})+\tilde{r}_k\cdot DB(u)\tilde{r}_{k,x}+\tilde{r}_{k,x}\cdot DB(u)\tilde{r}_k\right]\rangle \\
		&-\sum\limits_{k,l }v_l v_k^2(1-\chi_k^\e\xi_k\eta_k)\langle\tilde{r}_l\cdot D l_i, \left[B(u)\pa_x(\tilde{r}_{k,u}\tilde{r}_{k})+\tilde{r}_k\cdot DB(u)\tilde{r}_{k,x}+\tilde{r}_{k,x}\cdot DB(u)\tilde{r}_k\right]\rangle \\
		&-\sum\limits_{k}\sum_l v_l v_k^2(1-\chi_k^\e\xi_k\eta_k)\langle l_i, \left[\tilde{r}_l\cdot DB(u)\pa_x(\tilde{r}_{k,u}\tilde{r}_{k})+\tilde{r}_k\otimes \tilde{r}_l: D^2B(u)\tilde{r}_{k,x}+\tilde{r}_{k,x}\otimes \tilde{r}_l: D^2 B(u)\tilde{r}_k\right]\rangle \\
		&-\sum\limits_{k}v_k^2(1-\chi_k^\e\xi_k\eta_k)\langle l_i, \left[B(u)\pa_{xx}(\tilde{r}_{k,u}\tilde{r}_{k})+2\tilde{r}_{k,x}\cdot DB(u)\tilde{r}_{k,x}+\tilde{r}_{k}\cdot DB(u)\tilde{r}_{k,xx}+\tilde{r}_{k,xx}\cdot DB(u)\tilde{r}_k\right]\rangle \\
		&=O(1)\biggl(\sum_{j\ne i}|v_j|+\sum_k \sum_{l\ne k}|v_k|(|v_lv_{l,x}|+|v_kv_l|)+
		(|v_{i,x}|+|v_i|)|v_i|(1-{\color{black}{\chi_i^\e}}\xi_i\eta_i)\\
		&+\mathfrak{A}_i\rho_i^\e \big(|v_i^2\left(\frac{w_i}{v_i}\right)_x|+|v_i v_{i,x}|\Delta_i^\e+|v_i^3 \left(\frac{w_i}{v_i}\right)_{x}^2|\big)+v_iv_{i,x}(\chi')_i^\e \mathfrak{A}_i\rho_i^\e \biggl).
		\end{align*}
Applying the Lemmas \ref{lemme6.5}, \ref{lemme6.6}, \ref{lemme11.3}, \ref{lemme11.4} we get
\begin{align*}
		\mathcal{I}_{7,1,x}&=-\sum\limits_{k\neq j}(v_{k,xx}v_j+2v_{k,x}v_{j,x}+v_kv_{j,xx})\langle l_i,\tilde{r}_k\cdot DB\tilde{r}_j+B\tilde{r}_{k,u}\tilde{r}_j\rangle\\
		&-\sum\limits_{k\neq j}(v_{k,x}v_j+v_kv_{j,x})\langle u_x\cdot Dl_i,\tilde{r}_k\cdot DB\tilde{r}_j+B\tilde{r}_{k,u}\tilde{r}_j\rangle\\
		&-\sum\limits_{k\neq j}(v_{k,x}v_j+v_kv_{j,x})\langle l_i,\tilde{r}_k\otimes u_x: D^2B\tilde{r}_j+u_x\cdot DB\tilde{r}_{k,u}\tilde{r}_j\rangle\\
		&-\sum\limits_{k\neq j}(v_{k,x}v_j+v_kv_{j,x})\langle l_i,\tilde{r}_k\cdot DB\tilde{r}_{j,x}+\tilde{r}_{k,x}\cdot DB\tilde{r}_j+B\pa_x(\tilde{r}_{k,u}\tilde{r}_j)\rangle\\
		&=O(1)\sum_{k\ne j}(|v_{k,xx}v_j|+|v_{k,x}v_{j,x}|+|v_{k,x}v_j|).
		&\\
		\mathcal{I}_{7,2,x}&=-\sum\limits_{k\neq j}\sum\limits_{l}[v_{k,x}v_jv_l+v_kv_{j,x}v_{l}+v_kv_jv_{l,x}]\langle \tilde{r}_l\cdot Dl_i,\tilde{r}_k\cdot DB\tilde{r}_j+B\tilde{r}_{k,u}\tilde{r}_j\rangle\\
		&-\sum\limits_{k\neq j}\sum\limits_{l}v_kv_jv_l\langle \tilde{r}_{l,x}\cdot Dl_i,\tilde{r}_k\cdot DB\tilde{r}_j+B\tilde{r}_{k,u}\tilde{r}_j\rangle\\
		&-\sum\limits_{k\neq j}\sum\limits_{l}v_kv_jv_l\langle \tilde{r}_l\otimes u_x: D^2l_i,\tilde{r}_k\cdot DB\tilde{r}_j+B\tilde{r}_{k,u}\tilde{r}_j\rangle\\
		&-\sum\limits_{k\neq j}\sum\limits_{l}v_kv_jv_l\langle \tilde{r}_l\cdot Dl_i,\tilde{r}_k\otimes u_x:D^2B\tilde{r}_j+u_x\cdot DB\tilde{r}_{k,u}\tilde{r}_j\rangle\\
		&-\sum\limits_{k\neq j}\sum\limits_{l}v_kv_jv_l\langle \tilde{r}_l\cdot Dl_i,\tilde{r}_k\cdot DB\tilde{r}_{j,x}+\tilde{r}_{k,x}\cdot DB\tilde{r}_j+B\pa_x(\tilde{r}_{k,u}\tilde{r}_j)\rangle\\
		&=O(1)\sum_{k\ne j}(|v_{k,x}v_j|+|v_kv_j|).
		&\\
		\mathcal{I}_{7,3,x}&=-\sum\limits_{k\neq j}\sum\limits_{l}[v_{k,x}v_jv_l+v_kv_{j,x}v_{l}+v_kv_jv_{l,x}]\langle l_i,\tilde{r}_k\otimes\tilde{r}_l: D^2B\tilde{r}_j+\tilde{r}_l\cdot DB\tilde{r}_{k,u}\tilde{r}_j\rangle\\
		&-\sum\limits_{k\neq j}\sum\limits_{l}v_kv_jv_l\langle u_x\cdot Dl_i,\tilde{r}_k\otimes\tilde{r}_l: D^2B\tilde{r}_j+\tilde{r}_l\cdot DB\tilde{r}_{k,u}\tilde{r}_j\rangle\\
		&-\sum\limits_{k\neq j}\sum\limits_{l}v_kv_jv_l\langle l_i,\tilde{r}_k\otimes\tilde{r}_l\otimes u_x: D^3B\tilde{r}_j+\tilde{r}_l\otimes u_x D^2B\tilde{r}_{k,u}\tilde{r}_j\rangle\\
		&-\sum\limits_{k\neq j}\sum\limits_{l}v_kv_jv_l\langle l_i,\tilde{r}_{k,x}\otimes\tilde{r}_l: D^2B\tilde{r}_j+\tilde{r}_{l,x}\cdot DB\tilde{r}_{k,u}\tilde{r}_j\rangle\\
		&-\sum\limits_{k\neq j}\sum\limits_{l}v_kv_jv_l\langle l_i,\tilde{r}_k\otimes\tilde{r}_{l,x}: D^2B\tilde{r}_j+\tilde{r}_k\otimes\tilde{r}_l: D^2B\tilde{r}_{j,x}+\tilde{r}_l\cdot DB\pa_x(\tilde{r}_{k,u}\tilde{r}_j)\rangle\\
		&=O(1)\sum_{k\ne j}(|v_{k,x}v_j|+|v_kv_j|).
		\end{align*}
		From Lemma \ref{lemme11.6}, \ref{lemme6.5}, \ref{lemme6.6} it yields
	\begin{align}
		&v_kv_j \pa_{xx}(\tilde{r}_{k,u}\tilde{r}_j)
		=O(1)|v_kv_j|+O(1)\biggl(
		  |v_{k,x}v_j|+|w_{k,x}v_j|+|v_kv_{j,x}|+|v_kw_{j,x}|\biggl). \label{rkxxtech1bis5}
		  	\end{align}
		We have now using \eqref{rkxxtech1}, \eqref{rkxxtech1bis5} and Lemmas \ref{lemme11.3}, \ref{lemme11.4}, \ref{lemme11.6}
		\begin{align*}
		\mathcal{I}_{7,4,x}&=-\sum\limits_{k\neq j}(v_{k,x}v_j+v_kv_{j,x})\langle l_i,\tilde{r}_{k,x}\cdot DB\tilde{r}_j+\tilde{r}_k\cdot DB\tilde{r}_{j,x}+B(u)\pa_x(\tilde{r}_{k,u}\tilde{r}_j)\rangle\nonumber\\
		&-\sum\limits_{k\neq j}\sum_l v_lv_kv_j\langle \tilde{r}_l\cdot Dl_i(u),\tilde{r}_{k,x}\cdot DB\tilde{r}_j+\tilde{r}_k\cdot DB\tilde{r}_{j,x}+B(u)\pa_x(\tilde{r}_{k,u}\tilde{r}_j)\rangle\nonumber\\
		&-\sum\limits_{k\neq j}\sum_l v_kv_j v_l\langle l_i,\tilde{r}_l\otimes \tilde{r}_{k,x}:D^2B\tilde{r}_j+\tilde{r}_l\otimes \tilde{r}_k:D^2B\tilde{r}_{j,x}+\tilde{r}_l\cdot D B(u)\pa_x(\tilde{r}_{k,u}\tilde{r}_j)\rangle\nonumber\\
		&-\sum\limits_{k\neq j}v_kv_j\langle l_i,\tilde{r}_{k,xx}\cdot DB\tilde{r}_j+2\tilde{r}_{k,x}\cdot DB\tilde{r}_{j,x}+
		\tilde{r}_k\cdot DB\tilde{r}_{j,xx}+B(u)\pa_{xx}(\tilde{r}_{k,u}\tilde{r}_j)\rangle\nonumber\\
		&=O(1)\sum_{k\ne j}(|v_{k,x}v_j|+|w_{k,x}v_j|+|v_kv_j|).
		\end{align*}
		Combining \eqref{9.14prime}, Lemmas \ref{lemme6.5}, \ref{lemme6.6} and all the previous estimates on $I_{k,j,x}$, we obtain finally that
		\begin{align*}
		&(\mu_i v_{i,x}-(\tilde{\la}_i-\la_i^*)v_i-w_i)_{xx}\nonumber\\
		&=O(1)\biggl(
		\sum_{k\ne j}(|v_{k,xx}v_j|+|v_{k,x}v_{j,x}|+|v_{k,x}v_j|+|v_kv_j|+|w_{k,x}v_j|)\\
		&+
		\sum_{i\ne j}(|v_jv_{j,xxx}|+\abs{w_{j,xxx}v_j}+\abs{v_{j,xxx}w_j})
		+\sum_{j\ne i, k\ne j}|v_{k,xxx}v_k v_j|\mathfrak{A}_j\rho_j^\e
		\\
		&+\sum_{i\ne j}(|v_j|+|w_j|+|v_{j,x}|+|w_{j,x}|+|v_{j,xx}|+|w_{j,xx}|+|v_j|\left(\frac{w_j}{v_j}\right)_{x}^2 \mathfrak{A}_j\rho_j^\e )\\
		&+(|v_{i,x}^2|+|v_iv_{i,x}|+|v_iv_{i,xx}|+|v_i^2|)(1-{\color{black}{\chi_i^\e}}\xi_i\eta_i)+\mathfrak{A}_i \Delta_i^\e\big(
		|v_i v_{i,x}|+v_{i,x}^2+|v_iv_{i,xx}|\big)
		\\
		&+\mathfrak{A}_i\rho_i^\e\big(|v_iv_{i,x}\left(\frac{w_i}{v_i}\right)_x|
		+|v_i^2\left(\frac{w_i}{v_i}\right)^2_x|+|v_i^2\left(\frac{w_i}{v_i}\right)_x|+|w_{i,xx}v_i-v_{i,xx}w_i|
		\big)
		\\
		&+(|v_i |(|v_{i,x}|+|v_{i,xx}|)+|v_{i,x}|(|v_{i,x}|+|w_{i,x}|))\varkappa^\ep_i\mathfrak{A}_i\rho_i^\e\biggl).
\end{align*}
We have in particular obtained the desired estimate \eqref{9.14primefin}. Let us estimate now the right hand side of \eqref{9.14primefin}, first using Lemmas \ref{lemme6.5}, \ref{lemme6.6}, \eqref{6.54bis}, \eqref{6.45} and the fact that $0\leq \mathbbm{1}_{\{v_i^{2N}\geq 2\e\}}(1-\eta_i)\leq \sum_{j\ne i}\mathbbm{1}_{\{v_j^2\geq\frac{3}{4}|v_i|\}}$ we observe that
\begin{align}
&(|v_{i,x}^2|+|v_iv_{i,x}|+|v_iv_{i,xx}|+|v_i^2|)(1-\chi_i^\e\xi_i\eta_i)\nonumber\\
&=O(1)
(|v_{i,x}^2|+|v_iv_{i,x}|+|v_iv_{i,xx}|+|v_i^2|)[(1-\chi_i^\e\eta_i)\mathbbm{1}_{\{|\frac{w_i}{v_i}|\leq\frac{\delta_1}{2}\}}(\mathbbm{1}_{\{v_i^{2N}\geq 2\e\}}
+ \mathbbm{1}_{\{v_i^{2N}\leq 2\e\}})+\mathbbm{1}_{\{|\frac{w_i}{v_i}|\geq\frac{\delta_1}{2}\}}]\nonumber\\
&=O(1)
(|v_iv_{i,x}|+|v_iv_{i,xx}|+|v_i^2|)[(1-\eta_i)\mathbbm{1}_{\{|\frac{w_i}{v_i}|\leq\frac{\delta_1}{2}\}}\mathbbm{1}_{\{v_i^{2N}\geq 2\e\}}
+ \mathbbm{1}_{\{v_i^{2N}\leq 2\e\}}]\nonumber\\
&+O(1)
(v_{i,x}^2+|v_iv_{i,x}|+|v_iv_{i,xx}|+|v_i^2|)\mathbbm{1}_{\{|\frac{w_i}{v_i}|\geq\frac{\delta_1}{2}\}}+O(1)\sum_j(\La_j^1+\La_j^4)\nonumber\\
&=O(1)
(|v_iv_{i,x}|+|v_iv_{i,xx}|+|v_i^2|)\mathbbm{1}_{\{v_i^{2N}\leq 2\e\}}
+O(1) |v_{i,x} v_{i,xx}|\mathbbm{1}_{\{|\frac{w_i}{v_i}|\geq\frac{\delta_1}{2}\}}\nonumber\\
&+O(1)\sum_j(\La_j^1+\La_j^4+\La_j^6).\label{suytech1}
\end{align}
Similarly using Lemmas \ref{lemme6.5}, \ref{lemme6.6}, \eqref{6.54bis} and \eqref{ngtech5} we deduce that
\begin{align}
&\mathfrak{A}_i \Delta_i^\e\big(
		|v_i v_{i,x}|+v_{i,x}^2+|v_iv_{i,xx}|\big)+\mathfrak{A}_i\rho_i^\e\big(|v_iv_{i,x}\left(\frac{w_i}{v_i}\right)_x|
		+|v_i^2\left(\frac{w_i}{v_i}\right)^2_x|+|v_i^2\left(\frac{w_i}{v_i}\right)_x| \nonumber\\
		&+|w_{i,xx}v_i-v_{i,xx}w_i|
		\big)=O(1)\sum_j(\La_j^1+\delta_0^2\La_j^3+\La_j^4+\La_j^5)+O(1)\mathfrak{A}_i \Delta_i^\e\big(
		|v_i v_{i,x}|+|v_iv_{i,xx}|\big)\nonumber\\
		&=O(1)\sum_j(\La_j^1+\delta_0^2\La_j^3+\La_j^4+\La_j^5).
		\label{suytech2}
		\end{align}
		Combining now \eqref{9.14primefin}, \eqref{suytech1} and \eqref{suytech2}  we obtain
		\begin{align*}
		&(\mu_i v_{i,x}-(\tilde{\la}_i-\la_i^*)v_i-w_i)_{xx}\nonumber\\
		&=O(1)\biggl(
		\sum_j(\La_j^1+\La_j^2+\delta_0^2\La_j^3+\La_j^4+\La_j^5+\La_j^6)+
		\sum_{i\ne j}(|v_jv_{j,xxx}|+\abs{w_{j,xxx}v_j}+\abs{v_{j,xxx}w_j})
		\\
		&+\sum_{i\ne j}(|v_j|+|w_j|+|v_{j,x}|+|w_{j,x}|+|v_{j,xx}|+|w_{j,xx}|+|v_j|\left(\frac{w_j}{v_j}\right)_{x}^2 \mathfrak{A}_j\rho_j^\e )\\
		&+O(1)
(|v_iv_{i,x}|+|v_iv_{i,xx}|+|v_i^2|)\mathbbm{1}_{\{v_i^{2N}\leq 2\e\}}
+O(1) |v_{i,x} v_{i,xx}|\mathbbm{1}_{\{|\frac{w_i}{v_i}|\geq\frac{\delta_1}{2}\}}\nonumber\\
		&+(|v_i |(|v_{i,x}|+|v_{i,xx}|)+|v_{i,x}|(|v_{i,x}|+|w_{i,x}|))\varkappa^\ep_i\mathfrak{A}_i\rho_i^\e\biggl).
\end{align*}
	It concludes the proof of the Lemma.
	\end{proof}
\begin{lemma}\label{lemme9.9}
	For $1\leq i\leq n$ the following estimates are true.
	\begin{equation}
		\hat{r}_{i,x},\,r_{i,x}^\dagger,\,r^\ddagger_{i,x},\,r^\#_{i,x},\,r^\clubsuit_{i,x}, \,r^\diamondsuit_{i,x}=O(1).
	\end{equation}
	
\end{lemma}
\begin{proof} We recall that $r^\ddagger_{i}, r_i^\#, r_i^\clubsuit, r_i^\diamondsuit, \hat{r}_i$ and $r_i^\dagger$ are defined in \eqref{defAo} and \eqref{defAo1}. 
	From direct computation, we get using Lemmas \ref{estimimpo1}, \ref{lemme6.5}, \ref{lemme6.6}, \ref{lemme11.3} and the fact that $\tilde{r}_{i,\sig}=O(1)\xi_i\bar{v}_i$
	\begin{align*}
		\hat{r}_{i,x}&=\tilde{r}_{i,x}+\left(v_{i,x}\xi_i\pa_{v_i}\bar{v}_i+v_i\xi_i^\p\left(\frac{w_i}{v_i}\right)_x\pa_{v_i}\bar{v}_i+v_i\xi_i\pa_x(\pa_{v_i}\bar{v}_i)\right) \tilde{r}_{i,v}\\
		&-\left(\left(\frac{w_i}{v_i}\right)_x\xi_i^\p\bar{v}_i+\frac{w_i}{v_i}\left(\frac{w_i}{v_i}\right)_x\xi_i^{\p\p}\bar{v}_i+\frac{w_i}{v_i}\xi_i'\pa_x(\bar{v}_i)\right) \tilde{r}_{i,v}+\left(v_i\xi_i\pa_{v_i}\bar{v}_i-\frac{w_i}{v_i}\xi_i^\p\bar{v}_i\right)\pa_x(\tilde{r}_{i,v})\\
		&+\left(\frac{w_i}{v_i}\right)_x\theta^\p_i\tilde{r}_{i,\si}+\frac{w_i}{v_i}\left(\frac{w_i}{v_i}\right)_x\theta^{\p\p}_i\tilde{r}_{i,\si}+\frac{w_i}{v_i} \theta^\p_i\pa_x(\tilde{r}_{i,\si})\\
		&=O(1),\\
		r_{i,x}^\dagger&=\xi_i^{\p\p}\left(\frac{w_i}{v_i}\right)_x\bar{v}_i\tilde{r}_{i,v}+\xi_i^\p\pa_x(\bar{v}_i)\tilde{r}_{i,v}+\xi_i^\p\bar{v}_i\pa_x(\tilde{r}_{i,v})-\theta_i^{\p\p}\left(\frac{w_i}{v_i}\right)_x\tilde{r}_{i,\si}-\theta^\p_i\pa_x(\tilde{r}_{i,\si})\\
		&=O(1).
	\end{align*}
	Similarly we have	
	\begin{align*}
		r^\ddagger_{i,x}&=\tilde{r}_{i,x}+\left(\frac{w_i}{v_i}\right)_x\xi_i^\p\bar{v}_i\tilde{r}_{i,v}+\frac{w_i}{v_i}\xi_i^{\p\p}\left(\frac{w_i}{v_i}\right)_x\bar{v}_i\tilde{r}_{i,v}+\frac{w_i}{v_i}\xi_i^\p\pa_x(\bar{v}_i)\tilde{r}_{i,v}+\frac{w_i}{v_i}\xi_i^\p\bar{v}_i\pa_x(\tilde{r}_{i,v})\\
		&-\left(\frac{w_i}{v_i}\right)_x\theta_i^\p\tilde{r}_{i,\si}-\frac{w_i}{v_i}\theta_i^{\p\p}\left(\frac{w_i}{v_i}\right)_x\tilde{r}_{i,\si}-\frac{w_i}{v_i}\theta_i^\p\pa_x(\tilde{r}_{i,\si})\\
		&=O(1),\\
		r^\#_{i,x}&=\biggl(w_{i,x}\pa_{v_i}\bar{v}_i\xi_i+w_i\pa_x(\pa_{v_i}\bar{v}_i)\xi_i+w_i\pa_{v_i}\bar{v}_i \xi_i'\left(\frac{w_i}{v_i}\right)_x-\xi^{\p\p}_i\left(\frac{w_i}{v_i}\right)_x\bar{v}_i\frac{w_i^2}{v_i^2}-\xi^{\p}_i\pa_x(\bar{v}_i)\frac{w_i^2}{v_i^2}
		\\
		&-2\xi^{\p}_i\bar{v}_i\frac{w_i}{v_i}\left(\frac{w_i}{v_i}\right)_x\biggl)\tilde{r}_{i,v}+\left(w_i\pa_{v_i}\bar{v}_i\xi_i-\xi^\p_i\bar{v}_i\frac{w_i^2}{v_i^2}\right)\pa_x(\tilde{r}_{i,v})+2\frac{w_i}{v_i}\left(\frac{w_i}{v_i}\right)_x\theta_i^\p\tilde{r}_{i,\si}\\
		&+\frac{w^2_i}{v^2_i}\theta_i^{\p\p}\left(\frac{w_i}{v_i}\right)_x\tilde{r}_{i,\si}+\frac{w^2_i}{v^2_i}\theta_i^{\p} \pa_x(\tilde{r}_{i,\si})\\
		&=O(1).
	\end{align*}
	Since $r^\clubsuit_i=-\la_i^*\hat{r}_i$ and $r^\diamondsuit_i=-\la_i^*r_i^\dagger$ we obtain in a direct way that $\pa_x(r^\clubsuit_i)=O(1)$ and $\pa_x(r^\diamondsuit_i)=O(1)$.
\end{proof}

%------------------------------------------------------------------------------------------------------------------------------------------------
%------------------------------------------------------------------------------------------------------------------------------------------------
%------------------------------------------------------------------------------------------------------------------------------------------------
%------------------------------------------------------------------------------------------------------------------------------------------------

\section{Estimates for $\phi_i,\phi_{i,x}$ and $\psi_i, \psi_{i,x}$}\label{sectionphii}
From \eqref{supercle0} and  \eqref{supercle1} we recall that for any $i\in\{1,\cdots,n\}$
\begin{equation}
\begin{aligned}
&\phi_{i},\psi_{i}=O(1)(|\mathcal{J}_{1}|+|\mathcal{J}_{2}|),\\
&\phi_{i,x},\psi_{i,x}=O(1)(|\mathcal{J}_{1,x}|+|\mathcal{J}_{2,x}|+|\mathcal{J}_{1}|+|\mathcal{J}_{2}|).
\end{aligned}
\label{supercle2}
\end{equation}
It implies that it remains to estimate the terms $\mathcal{J}_{1}$, $\mathcal{J}_{2}$, $\mathcal{J}_{1,x}$ and $\mathcal{J}_{2,x}$.
%From \eqref{eqn-remainder-1} we recall
%\begin{align}
%	&\sum\limits_{i}\left(v_{i,t}+(\tilde{\la}_iv_i)_x-(\mu_iv_{i,x})_x-\sum\limits_{j\neq i}(\mu_i-\mu_j)b_{ij}(w_{j,x}-(w_j/v_j)v_{j,x})_x\right)\hat{r}_i \nonumber\\
%	&+\sum\limits_{i}\left(w_{i,t}+(\tilde{\la}_iw_i)_x-(\mu_iw_{i,x})_x-\sum\limits_{j\neq i}(\mu_i-\mu_j)\hat{b}_{ji}(w_{j,x}-(w_j/v_j)v_{j,x})_x\right)r^\dagger_i \nonumber\\
%	&=\sum\limits_{l=1}^{29}\sum\limits_{i}\al_{i}^{1,l}+\sum\limits_{l=1}^{6}\sum\limits_{i}\al_{i}^{5,1,l}+\sum\limits_{l=1}^{N_4}\sum\limits_{i}\al_{i}^{4,l}+\sum\limits_{l=2}^{N_5}\sum\limits_{i}\al_{i}^{5,l}+\sum\limits_{l=1}^{N_6}\sum\limits_{i}\al^{6,l}_i\nonumber\\
%	&+\sum\limits_{l=2}^{N_7}\sum\limits_{i}\al^{7,l}_i+\sum\limits_{k=8}^{10}\sum\limits_{l=1}^{N_k}\sum\limits_{i}\al^{k,l}_i+\sum\limits_{i}\sum\limits_{l=1}^{4}\al_i^{7,1,l}+\sum\limits_{i}\al_{i}^{11}.
%\end{align}
\subsection{Second order estimates for $J_i$ s}\label{section-phi-i-x}
\begin{lemma}
	For $1\leq i\leq n$ and $1\leq l\leq 27$ it follows,
	\begin{equation}
	\begin{aligned}
		&\widetilde{\al}_i^{1,l}=O(1)\sum_j(\La_j^1+\La_j^2+\delta_0^2\La_j^3+\La_j^4+\La_j^5+\La_j^6+\La_j^{6,1})+R_\e^{1,i,l},\\
		& \widetilde{\al}_{i,x}^{1,l}=O(1)\sum_j(\La_j^1+\La_j^2+\delta_0^2\La_j^3+\La_j^4+\La_j^5+\La_j^6+\La_j^{6,1})+\widetilde{R}_\e^{1,i,l},
		\end{aligned}
		\label{12.2}
	\end{equation}
	with:
	\begin{equation}
	\int_{\hat{t}}^T\int_{\R}(|R_\e^{1,i,l}|+|\widetilde{R}_\e^{1,i,l}|)dx ds=O(1)\delta_0^2,
	\label{12.2bis}
	\end{equation}
	for $\e>0$ small enough in terms of $T-\hat{t}$ and $\delta_0$.
\end{lemma}
\begin{proof}
	We are now going to compute each $\sum_i\widetilde{\al}_{i,x}^{1,l}$ with $l\in\{1,\cdots,27\}$ and give an estimate for each $\sum_i \widetilde{\al}_{i,x}^{1,l}$ and $\sum_i \widetilde{\al}_{i}^{1,l}$. We start with $\widetilde{\al}_{i}^{1,1}$. We observe that
	$$\widetilde{\al}_{i}^{1,1}=\sum\limits_{j\neq i}\tilde{r}_{i,u}\tilde{r}_j\left[(\mu_iv_{i,x}-\tilde{\la}_iv_i)v_j-(w_j-\la_j^*v_j)v_{i}\right]=O(1)\sum\limits_{j=1}^{n}\Lambda_j^1.$$
	We have now applying Lemmas \ref{lemme11.4}, \ref{lemme6.5}, \ref{lemme6.6} since $(\tilde{r}_{i,u}\tilde{r}_j)_x=O(1)$
	\begin{align*}
		\widetilde{\al}_{i,x}^{1,1}&=\sum\limits_{j\neq i}\tilde{r}_{i,u}\tilde{r}_j \left[(\mu_iv_{i,x}-\tilde{\la}_iv_i)v_j-(w_j-\la_j^*v_j)v_{i}\right]_x \\
		&+\sum\limits_{j\neq i}(\tilde{r}_{i,u}\tilde{r}_j)_x\left[(\mu_iv_{i,x}-\tilde{\la}_iv_i)v_j-(w_j-\la_j^*v_j)v_{i}\right]=O(1)\La_i^1.
		\end{align*}
		Applying \eqref{ngtech4} we deduce that
	\begin{align*}
	\widetilde{\al}_i^{1,2}&=\tilde{r}_{i,u}\tilde{r}_i\left[(\mu_iv_{i,x}-(\tilde{\la}_i-\la_i^*)v_i-w_i)v_i\right] \nonumber\\
	&=O(1)\sum_j(\La_j^1+\La_j^4+\La_j^6)%+O(1)\sum\limits_{k=1}^{n}v_k^2(1-\eta_k \chi_k^\ep)\mathbbm{1}_{\{|\frac{w_k}{v_k}|\leq\frac{\delta_1}{2}\}}
+O(1)\sum\limits_{j}v_j\xi_j (\chi')^\e_j\hat{v}_j v_i+O(1)\sum\limits_{k}^{n} \eta_k v_i v_k^2(1- \chi_k^\ep).
	\end{align*}
	Using now the Lemmas \ref{lemme6.5}, \ref{lemme6.6} we obtain
	\begin{align*}
	&\int_{\hat{t}}^T\left(\big|\sum\limits_{j}v_j\xi_j (\chi')^\e_j\hat{v}_j v_i\big|+\big|\sum\limits_{k}^{n} \eta_k v_i v_k^2(1- \chi_k^\ep)\big| \right) dx ds=O(1)(T-\hat{t})\delta_0\e^{\frac{1}{N}}.
\end{align*}
Now, we compute
	\begin{align*}
	\al_{i,x}^{1,2}&= \tilde{r}_{i,u}\tilde{r}_i\left[\mu_iv_{i,x}-(\tilde{\la}_i-\la_i^*)v_i-w_i\right]_x v_i+v_{i,x} \tilde{r}_{i,u}\tilde{r}_i\left[\mu_iv_{i,x}-(\tilde{\la}_i-\la_i^*)v_i-w_i\right] \\
	&+ (\tilde{r}_{i,u}\tilde{r}_i)_x\left[(\mu_iv_{i,x}-(\tilde{\la}_i-\la_i^*)v_i-w_i)v_i\right].	
	\end{align*}
	Using now \eqref{ngtech4} and Lemma \ref{lemme11.4}, we deduce that
	\begin{align}
	&v_{i,x} \tilde{r}_{i,u}\tilde{r}_i\left[\mu_iv_{i,x}-(\tilde{\la}_i-\la_i^*)v_i-w_i\right]+ (\tilde{r}_{i,u}\tilde{r}_i)_x\left[(\mu_iv_{i,x}-(\tilde{\la}_i-\la_i^*)v_i-w_i)v_i\right]\nonumber\\
	&=O(1)\sum_j(\La_j^1+\La_j^4+\La_j^6)%+O(1)\sum\limits_{k=1}^{n}v_k^2(1-\eta_k \chi_k^\ep)\mathbbm{1}_{\{|\frac{w_k}{v_k}|\leq\frac{\delta_1}{2}\}}
+O(1)\sum\limits_{j}(|v_{i,x}|+|w_{i,x}|+\sum_l|v_l|)| v_j\xi_j (\chi')^\e_j\hat{v}_j|\nonumber\\
&+O(1)\sum\limits_{k=1}^{n} \eta_k (|v_{i,x}|+|w_{i,x}|+\sum_l|v_l|) v_k^2(1- \chi_k^\ep).\label{ali12a}
	\end{align}
Applying now the Lemmas \ref{lemme6.5}, \ref{lemme6.6} and \eqref{estimate-v-i-xx-1aa}, it yields
\begin{align}
&\tilde{r}_{i,u}\tilde{r}_i\left[\mu_iv_{i,x}-(\tilde{\la}_i-\la_i^*)v_i-w_i\right]_x v_i
=O(1)\sum_k(\La_k^1+\La_k^4+\La_k^5+\La_k^6)\nonumber\\
&+O(1)(|v_{i,x}v_i^2|+|v_i^4|)\eta_i(1-\chi_i^\ep)%\mathbbm{1}_{\{|\frac{w_i}{v_i}|\leq\frac{\delta_1}{2}\}}
		%+O(1)|v_i^2\xi_i(\chi')^\e_i\hat{v}_i|(|v_{i,x}|+|v_i^2|)\nonumber\\
%&
%+O(1)|v_{i,x} v_i\xi_i (\chi')^\e_i\hat{v}_i|(|v_{i,x}|+|v_i^2|)
{\color{black}{+O(1) \varkappa^\ep_i\rho_i^\e\mathfrak{A}_iv_i^2v_{i,x}}}.
		\label{am12.6}
\end{align}
Finally combining \eqref{ali12a}, \eqref{am12.6} we have obtained
\begin{align}
&\widetilde{\al}_{i,x}^{1,2}=
O(1)\sum_j(\La_j^1+\La_j^4+\La_j^5+\La_j^6)%O(1)\sum\limits_{k=1}^{n}v_k^2(1-\eta_k \chi_k^\ep)\mathbbm{1}_{\{|\frac{w_k}{v_k}|\leq\frac{\delta_1}{2}\}}\nonumber\\
+O(1)\sum\limits_{k=1}^{n} \eta_k (|v_{i,x}|+|w_{i,x}|+\sum_l|v_l|) v_k^2(1- \chi_k^\ep)\nonumber\\
&+O(1)\sum\limits_{j}(|v_{i,x}|+|w_{i,x}|+\sum_l|v_l|)| v_j\xi_j (\chi')^\e_j\hat{v}_j|
%&+O(1)(|v_{i,x}v_i^2|+|v_i^4|)\eta_i(1-\chi_i^\ep)\mathbbm{1}_{\{|\frac{w_i}{v_i}|\leq\frac{\delta_1}{2}\}}
		%+O(1)|v_i^2\xi_i(\chi')^\e_j\hat{v}_i|(|v_{i,x}|+|v_i^2|)\nonumber\\
%&+O(1)|v_{i,x} v_i\xi_i (\chi')^\e_j\hat{v}_i|(|v_{i,x}|+|v_i^2|)
{\color{black}{+O(1) \varkappa^\ep_i\rho_i^\e\mathfrak{A}_iv_i^2v_{i,x}}}.
%
%&+O(1)(|v_{i,x}v_i^2|+|v_i^4|)\eta_i(1-\chi_i^\ep)\mathbbm{1}_{\{|\frac{w_i}{v_i}|\leq\frac{\delta_1}{2}\}}
%		\nonumber\\
%&+O(1)\sum\limits_{j}|v_i v_j\xi_j (\chi')^\e_j\hat{v}_j|(|v_{i,x}|+|v_i^2|)+O(1)\sum\limits_{j}|v_{i,x} v_j\xi_j (\chi')^\e_j\hat{v}_j|(|v_{i,x}|+|v_i^2|)\nonumber\\
%		&{\color{black}{+O(1)\sum_k \varkappa^\ep_k\rho_k^\e\mathfrak{A}_k|v_i v_k|}}.
\end{align}
We remark that
\begin{align}
&\widetilde{\al}_{i,x}^{1,2}=
O(1)\sum_j(\La_j^1+\La_j^4+\La_j^5+\La_j^6)+\widetilde{R}_{\e}^{1,i,2},
\end{align}
with $\int_{\hat{t}}^T\int_{\R}|\widetilde{R}_\e^{1,i,2}|dx ds=O(1)\delta_0^2,$
	for $\e>0$ small enough in terms of $T-\hat{t}$ and $\delta_0$.
%\ref{estimate-v-i-xx-1a}, the Lemma \ref{lemme6.8} and Lemmas \ref{lemme6.5}, \ref{lemme6.6} it yields
	%\begin{equation*}
	%\al_{i,x}^{1,2}=O(1)\sum_j(\Lambda_j^1+\Lambda_j^4+\Lambda_j^5)+O(1)v_i^2v_{i,x}(1-\xi_i\eta_i)+O(1)\La_i^6.
	%	\end{equation*}
	%	By proceeding as previously we deduce that
	%	\begin{equation*}
	%\al_{i,x}^{1,2}=O(1)\sum_j(\Lambda_j^1+\Lambda_j^4+\Lambda_j^5)+O(1)\La_i^6.
	%	\end{equation*}
Using now the fact that $\mbox{supp}\xi$ is included in $\{x,\theta(x)=1\}$ we observe that $\widetilde{\alpha}_{1,3}=\tilde{r}_{i,v}(2\xi_i\pa_{v_i}\bar{v}_i+v_i\xi_i\pa_{v_iv_i}\bar{v}_i)\left[(\mu_iv_{i,x}-(\tilde{\la}_i-\la_i^*)v_i-w_i)v_{i,x}\right]$. 
By using similar argument as previously and the Lemmas \ref{lemme6.8},  \ref{estimimpo1} we get
\begin{equation}
	 \widetilde{\al}_{i}^{1,3}=O(1)\sum\limits_{j=1}^{n}(\Lambda_j^1+\Lambda_j^4+\Lambda_j^6)+R_\e^{1,i,3},
\end{equation}
with $\int_{\hat{t}}^T\int_{\R}|R_\e^{1,i,3}|dx ds=O(1)\delta_0^2,$
	for $\e>0$ small enough in terms of $T-\hat{t}$ and $\delta_0$.
	Furthermore we have
\begin{align*}
	&\widetilde{\al}_{i,x}^{1,3}=  \pa_x(\tilde{r}_{i,v})(2\xi_i\pa_{v_i}\bar{v}_i+v_i\xi_i\pa_{v_iv_i}\bar{v}_i)\left[(\mu_iv_{i,x}-(\tilde{\la}_i-\la_i^*)v_i-w_i)v_{i,x}\right]\\
	&+ \tilde{r}_{i,v}\big(\xi_i'\left(\frac{w_i}{v_i}\right)_x (2\pa_{v_i}\bar{v}_i+v_i\pa_{v_iv_i}\bar{v}_i)
	+2\xi_i \pa_x( \pa_{v_i}\bar{v}_i)+v_{i,x}\xi_i\pa_{v_iv_i}\bar{v}_i+v_i\xi_i \pa_x(\pa_{v_iv_i}\bar{v}_i)\big)\\
	&\times\left[(\mu_iv_{i,x}-(\tilde{\la}_i-\la_i^*)v_i-w_i)v_{i,x}\right]\\
	&+ \tilde{r}_{i,v}(2\xi_i\pa_{v_i}\bar{v}_i+v_i\xi_i\pa_{v_iv_i}\bar{v}_i)\left[(\mu_iv_{i,x}-(\tilde{\la}_i-\la_i^*)v_i-w_i)v_{i,x}\right]_x.
	%&=\mathcal{O}(1)\La_i^1,\\
	\end{align*}
	Using Lemmas \ref{lemme11.3}, \ref{lemme6.5}, \ref{lemme6.6}  and \eqref{6.54}, \eqref{ngtech4}, \eqref{ngtech5} we deduce that
	\begin{align*}
	& \pa_x(\tilde{r}_{i,v})(2\xi_i\pa_{v_i}\bar{v}_i+v_i\xi_i\pa_{v_iv_i}\bar{v}_i)\left[(\mu_iv_{i,x}-(\tilde{\la}_i-\la_i^*)v_i-w_i)v_{i,x}\right]\\
	&=O(1)
	|\left(\frac{w_i}{v_i}\right)_x| \mathfrak{A}_i\rho_i^\e\left[(\mu_iv_{i,x}-(\tilde{\la}_i-\la_i^*)v_i-w_i)v_{i,x}\right]\big)\\
	&+O(1)\sum_j(\La_j^1+\La_j^4+\La_j^6)%+O(1)\sum\limits_{k=1}^{n}v_k^2(1-\eta_k \chi_k^\ep)\mathbbm{1}_{\{|\frac{w_k}{v_k}|\leq\frac{\delta_1}{2}\}}
+O(1)v_{i,x}v_i\xi_i (\chi')^\e_i\hat{v}_i+O(1) \eta_i  v_{i,x}v_i^2(1- \chi_i^\ep)\\
	&=O(1)\biggl(
	|\left(\frac{w_i}{v_i}\right)_x| \mathfrak{A}_i\rho_i^\e (O(1)v_i+O(1)\sum\limits_{j=1}^{n}v^2_j\left(\frac{w_j}{v_j}\right)_x \mathfrak{A}_j\rho^\e_j %+O(1)\sum\limits_{j}v_j\xi_j (\chi')^\e_j\hat{v}_j
	)\left[(\mu_iv_{i,x}-(\tilde{\la}_i-\la_i^*)v_i-w_i)\right]\biggl)\\
	&+O(1)\sum_j(\La_j^1+\La_j^4+\La_j^6)%+O(1)\sum\limits_{k=1}^{n}v_k^2(1-\eta_k \chi_k^\ep)\mathbbm{1}_{\{|\frac{w_k}{v_k}|\leq\frac{\delta_1}{2}\}}
	+O(1)|v_{i,x}v_i\xi_i (\chi')^\e_i\hat{v}_i|+O(1)| \eta_i v_{i,x} v_i^2(1- \chi_i^\ep)|\\
	&=O(1)\sum\limits_{j=1}^{n}(\Lambda_j^1+\Lambda_j^4+\Lambda_j^6)+O(1)\delta_0^2\sum_{j}
	\La_j^3+O(1)(|v_{i,x}|+|w_{i,x}|)[v_i\xi_i (\chi')^\e_i\hat{v}_i|\\
	&+O(1)\eta_i (|v_{i,x}|+|w_{i,x}|) v_i^2(1- \chi_i^\ep).
	%&+O(1)\sum_i
	%|\left(\frac{w_i}{v_i}\right)_x| \mathfrak{A}_i\rho_i^\e \sum\limits_{j}v_j\xi_j (\chi')^\e_j\hat{v}_j \left[(\mu_iv_{i,x}-(\tilde{\la}_i-\la_i^*)v_i-w_i)\right]\\
	%&=O(1)\sum\limits_{j=1}^{n}(\Lambda_j^1+\Lambda_j^4+\Lambda_j^6)+O(1)\delta_0^2\sum_{j}
	%\La_j^3+O(1)\sum\limits_{j}(|v_{i,x}|+|w_{i,x}|)[v_j\xi_j (\chi')^\e_j\hat{v}_j|\\
	%&+O(1)\sum\limits_{k=1}^{n} \eta_k (|v_{i,x}|+|w_{i,x}|) v_k^2(1- \chi_k^\ep)\\
	%&+O(1)\sum_i
	%|\left(\frac{w_i}{v_i}\right)_x| \mathfrak{A}_i\rho_i^\e \sum\limits_{j}v_j\xi_j (\chi')^\e_j\hat{v}_j  \left[ v_i +O(1)\sum\limits_{j=1}^{n}v^2_j\left(\frac{w_j}{v_j}\right)_x \mathfrak{A}_j\rho^\e_j +O(1)\sum\limits_{j}v_j\xi_j (\chi')^\e_j\hat{v}_j\right]\\
	%&=O(1)\sum\limits_{j=1}^{n}(\Lambda_j^1+\Lambda_j^4+\Lambda_j^6)+O(1)\delta_0^2\sum_{j}
	%\La_j^3+O(1)\sum\limits_{j}(|v_{i,x}|+|w_{i,x}|)[v_j\xi_j (\chi')^\e_j\hat{v}_j|\\
	%&+O(1)\sum\limits_{k=1}^{n} \eta_k (|v_{i,x}|+|w_{i,x}|) v_k^2(1- \chi_k^\ep)+O(1)\sum_i
	%|v_i\left(\frac{w_i}{v_i}\right)_x| \mathfrak{A}_i\rho_i^\e \sum\limits_{j}v_j\xi_j (\chi')^\e_j\hat{v}_j  \\
	%&+O(1)\sum_i
	%|v_i\left(\frac{w_i}{v_i}\right)_x| \mathfrak{A}_i\rho_i^\e (\sum\limits_{j}\xi_j (\chi')^\e_j\hat{v}_j )^2.%\sum\limits_{j}\xi_j (\chi')^\e_j\hat{v}_j.
	\end{align*}
	It gives in particular
	\begin{align*}
	&\pa_x(\tilde{r}_{i,v})(2\xi_i\pa_{v_i}\bar{v}_i+v_i\xi_i\pa_{v_iv_i}\bar{v}_i)\left[(\mu_iv_{i,x}-(\tilde{\la}_i-\la_i^*)v_i-w_i)v_{i,x}\right]\\
	&=O(1)\sum\limits_{j=1}^{n}(\Lambda_j^1+\Lambda_j^4+\Lambda_j^6)+O(1)\delta_0^2\sum_{j}
	\La_j^3+\widetilde{R}_\e^{1,i,3,1},
	\end{align*}
	with $\int_{\hat{t}}^T\int_{\R}|\widetilde{R}_\e^{1,i,3,1}|dx ds=O(1)\delta_0^2,$
	for $\e>0$ small enough.
	Proceeding as for $\widetilde{\alpha}_{i,x}^{1,2}$ and using the fact that $\tilde{r}_{i,v}(2\xi_i\pa_{v_i}\bar{v}_i+v_i\xi_i\pa_{v_iv_i}\bar{v}_i)=O(1)\mathfrak{A}_i\rho_i^\e$ from the Lemma \ref{estimimpo1}, we get
	% $\alpha_{i}^{1,2}$ we have since $\xi_i v_{i,x}=O(1)(v_i+\sum_{j\ne i}(|v_j|+|w_j|))$ from Lemma \ref{lemme6.8}
	\begin{align*}
	&\tilde{r}_{i,v}(2\xi_i\pa_{v_i}\bar{v}_i+v_i\xi_i\pa_{v_iv_i}\bar{v}_i)\left[(\mu_iv_{i,x}-(\tilde{\la}_i-\la_i^*)v_i-w_i)v_{i,x}\right]_x\\
	&\hspace{4cm}=
O(1)\sum_j(\La_j^1+\La_j^4+\La_j^5+\La_j^6)+\widetilde{R}_{\e}^{1,i,3,2},
\end{align*}
with $\int_{\hat{t}}^T\int_{\R}|\widetilde{R}_\e^{1,i,3,2}|dx ds=O(1)\delta_0^2,$
	for $\e>0$ small enough. Following the study of the term $\pa_x(\tilde{r}_{i,v})(2\xi_i\pa_{v_i}\bar{v}_i+v_i\xi_i\pa_{v_iv_i}\bar{v}_i)\left[(\mu_iv_{i,x}-(\tilde{\la}_i-\la_i^*)v_i-w_i)v_{i,x}\right]$ we can prove due to the fact that $\tilde{r}_{i,v} \xi_i'\left(\frac{w_i}{v_i}\right)_x (2\pa_{v_i}\bar{v}_i+v_i\pa_{v_iv_i}\bar{v}_i)=O(1))\mathfrak{A}_i\rho_i^\e\left(\frac{w_i}{v_i}\right)_x $  that
	\begin{align*}
	&\tilde{r}_{i,v} \xi_i'\left(\frac{w_i}{v_i}\right)_x (2\pa_{v_i}\bar{v}_i+v_i\pa_{v_iv_i}\bar{v}_i)
	\left[(\mu_iv_{i,x}-(\tilde{\la}_i-\la_i^*)v_i-w_i)v_{i,x}\right]\\
	&=O(1)\sum\limits_{j=1}^{n}(\Lambda_j^1+\Lambda_j^4+\Lambda_j^6)+O(1)\delta_0^2\sum_{j}
	\La_j^3+\widetilde{R}_\e^{1,i,3,3},
	\end{align*}
with $\int_{\hat{t}}^T\int_{\R}|\widetilde{R}_\e^{1,i,3,3}|dx ds=O(1)\delta_0^2,$
	for $\e>0$ small enough. Let us study now the term $ \tilde{r}_{i,v}\big(%\xi_i'\left(\frac{w_i}{v_i}\right)_x (2\pa_{v_i}\bar{v}_i+v_i\pa_{v_iv_i}\bar{v}_i)
	%+2\xi_i \pa_x( \pa_{v_i}\bar{v}_i)+
	v_{i,x}\xi_i\pa_{v_iv_i}\bar{v}_i+v_i\xi_i \pa_x(\pa_{v_iv_i}\bar{v}_i)\big)\left[(\mu_iv_{i,x}-(\tilde{\la}_i-\la_i^*)v_i-w_i)v_{i,x}\right]$. %We can proceed as previously for the term  $ \tilde{r}_{i,v}\big(\xi_i'\left(\frac{w_i}{v_i}\right)_x (2\pa_{v_i}\bar{v}_i+v_i\pa_{v_iv_i}\bar{v}_i)
	%\big)\left[(\mu_iv_{i,x}-(\tilde{\la}_i-\la_i^*)v_i-w_i)v_{i,x}\right]$. 
	We have using Lemma \ref{estimimpo1}, \ref{lemme6.5}, \ref{lemme6.6} and \eqref{ngtech5} 
	\begin{align*}
	\xi_i\pa_{v_iv_i}\bar{v}_i
	v_{i,x}^2&=O(1)\frac{\mathfrak{A}_i\rho_i^\e}{|v_i|} \biggl[v_i^2+\sum\limits_{j\ne i}v^4_j\left(\frac{w_j}{v_j}\right)_x^2  \mathfrak{A}_j\rho^\e_j  \biggl]\\
	&=O(1)\mathfrak{A}_i\rho_i^\e.
	\end{align*}
	Using now again \eqref{ngtech4} allows to prove that
	\begin{align*}
	& \tilde{r}_{i,v}v_{i,x}\xi_i\pa_{v_iv_i}\bar{v}_i\left[(\mu_iv_{i,x}-(\tilde{\la}_i-\la_i^*)v_i-w_i)v_{i,x}\right]=O(1)\sum_j(\La_j^1+\La_j^4+\La_j^6)%+O(1)\sum\limits_{k=1}^{n}v_k^2(1-\eta_k \chi_k^\ep)\mathbbm{1}_{\{|\frac{w_k}{v_k}|\leq\frac{\delta_1}{2}\}}
\\
&+O(1)\sum\limits_{j}v_j\xi_j (\chi')^\e_j\hat{v}_j+O(1)\sum\limits_{k=1}^{n} \eta_k v_k^2(1- \chi_k^\ep),
	\end{align*}
	which is the required result.
We have from \eqref{ngtech5} $\mathfrak{A}_i\rho_i^\e v_{i,x}=O(1)\mathfrak{A}_i\rho_i^\e\sqrt{|v_i|}$ it implies then using Lemma \ref{estimimpo1} that
\begin{align*}
&\xi_i\pa_x(\pa_{v_i}\bar{v}_i)v_{i,x}=\mathfrak{A}_i\rho_i^\e(\frac{v^2_{i,x}}{v_i}+|v_{i,x}|\sum_{j\ne i}\frac{|v_{j,x}|}{\sqrt{|v_i|}})=O(1))\mathfrak{A}_i\rho_i^\e.
\end{align*}
We deduce then again that
\begin{align*}
& \tilde{r}_{i,v}2\xi_i \pa_x( \pa_{v_i}\bar{v}_i)\left[(\mu_iv_{i,x}-(\tilde{\la}_i-\la_i^*)v_i-w_i)v_{i,x}\right]=O(1)\sum_j(\La_j^1+\La_j^4+\La_j^6)\\%+O(1)\sum\limits_{k=1}^{n}v_k^2(1-\eta_k \chi_k^\ep)\mathbbm{1}_{\{|\frac{w_k}{v_k}|\leq\frac{\delta_1}{2}\}}
&+O(1)\sum\limits_{j}v_j\xi_j (\chi')^\e_j\hat{v}_j+O(1)\sum\limits_{k=1}^{n} \eta_k v_k^2(1- \chi_k^\ep).
	\end{align*}
Using again Lemma \ref{estimimpo1}  and similar arguments enable to prove that
\begin{align*}
& \tilde{r}_{i,v}v_i\xi_i \pa_x(\pa_{v_iv_i}\bar{v}_i)\left[(\mu_iv_{i,x}-(\tilde{\la}_i-\la_i^*)v_i-w_i)v_{i,x}\right]=O(1)\sum_j(\La_j^1+\La_j^4+\La_j^6)\\%+O(1)\sum\limits_{k=1}^{n}v_k^2(1-\eta_k \chi_k^\ep)\mathbbm{1}_{\{|\frac{w_k}{v_k}|\leq\frac{\delta_1}{2}\}}
&+O(1)\sum\limits_{j}v_j\xi_j (\chi')^\e_j\hat{v}_j+O(1)\sum\limits_{k=1}^{n} \eta_k v_k^2(1- \chi_k^\ep).
	\end{align*}
	Combining all the previous estimates ensures that \eqref{12.2} is satisfied for $\widetilde{\al}_{i,x}^{1,3}$.
We observe easily using Lemma \ref{estimimpo1}, \ref{lemme6.5}, \ref{lemme6.6} and \eqref{ngtech5} that
\begin{align*}
\widetilde{\al}_{i}^{1,4}=\tilde{r}_{i,v}\mu_i v_{i,x}\xi^\p_i\left(w_{i,x}-\frac{w_i}{v_i}v_{i,x}\right) \pa_{v_i}\bar{v}_i
=&O(1)\sum_j \La_j^4.
\end{align*}
Next direct computation gives
	\begin{align*}
	\widetilde{\al}_{i,x}^{1,4}&=\pa_x(\tilde{r}_{i,v})\mu_i \xi_i^\p v_{i,x}\left(w_{i,x}-\frac{w_i}{v_i}v_{i,x}\right) \pa_{v_i}\bar{v}_i+\tilde{r}_{i,v}\mu_{i,x} \xi_i^\p v_{i,x}\left(w_{i,x}-\frac{w_i}{v_i}v_{i,x}\right) \pa_{v_i}\bar{v}_i\\
	&+\tilde{r}_{i,v}\mu_i \xi_i^\p  v_{i,xx}\left(w_{i,x}-\frac{w_i}{v_i}v_{i,x}\right) \pa_{v_i}\bar{v}_i+\tilde{r}_{i,v}\mu_i \xi_i^\p v_{i,x}\left(w_{i,x}-\frac{w_i}{v_i}v_{i,x}\right)_x \pa_{v_i}\bar{v}_i\\
	&+\tilde{r}_{i,v}\mu_i \xi_i^\p v_{i,x}\left(w_{i,x}-\frac{w_i}{v_i}v_{i,x}\right) \pa_x(\pa_{v_i}\bar{v}_i)+\tilde{r}_{i,v}\mu_i \xi_i^{\p\p} v_i \left(\frac{w_i}{v_i}\right)_x^2 v_{i,x} \pa_{v_i}\bar{v}_i.
	\end{align*}
	Applying Lemmas \ref{estimimpo1}, \ref{lemme11.3}, \ref{lemme6.6}, \ref{lemme6.5} and \eqref{ngtech5} we obtain
	\begin{align*}
	&\xi_i'v_{i,x}\left(w_{i,x}-\frac{w_i}{v_i}v_{i,x}\right) \pa_{v_i}\bar{v}_i\big(\pa_x(\tilde{r}_{i,v})\mu_i +\tilde{r}_{i,v}\mu_{i,x} \big)=O(1)(\sum_j\La_j^4+\delta_0^2\La_i^3)
	\end{align*}
	Let us deal with the term  $\tilde{r}_{i,v}\mu_i \xi_i^\p v_{i,xx}\left(w_{i,x}-\frac{w_i}{v_i}v_{i,x}\right) \pa_{v_i}\bar{v}_i$, applying Lemmas \ref{estimimpo1}, \ref{lemme9.7} we deduce that
	\begin{align}
		&\tilde{r}_{i,v}\mu_i \xi_i^\p v_{i,xx}\left(w_{i,x}-\frac{w_i}{v_i}v_{i,x}\right) \pa_{v_i}\bar{v}_i=O(1)(|v_{i,x}|+|w_{i,x}|+|v_i|) \xi_i^\p \rho_i^\e\left(w_{i,x}-\frac{w_i}{v_i}v_{i,x}\right)\nonumber\\
		&+O(1)\sum_k(\La_k^1+\La_k^4+\La_k^5+\La_6^k)+O(1)|v_iv_{i,x}||1-\chi_i^\e|\eta_i\xi_i^\p \rho_i^\e\left(w_{i,x}-\frac{w_i}{v_i}v_{i,x}\right)
		\nonumber\\
		&+O(1)|v_i|^3(1-{\color{black}{\chi_i^\e}})\eta_i\xi_i^\p \rho_i^\e\left(w_{i,x}-\frac{w_i}{v_i}v_{i,x}\right)
		{\color{black}{+O(1)\sum_k \varkappa^\ep_k\rho_k^\e\mathfrak{A}_k|v_{k,x}v_k|}}\xi_i^\p \rho_i^\e\left(w_{i,x}-\frac{w_i}{v_i}v_{i,x}\right)\nonumber\\
		&=O(1)\sum_k(\La_k^1+\La_k^4+\La_k^5+\La_k^6+\La_k^{6,1}).%+O(1)|v_i|^3(1-{\color{black}{\chi_i^\e}}\xi_i\eta_i)\xi_i^\p \rho_i^\e\left(w_{i,x}-\frac{w_i}{v_i}v_{i,x}\right)\nonumber\\
		%&=O(1)\sum_k(\La_k^1+\La_k^4+\La_k^5+\La_k^6+\La_k^{6,1})\nonumber\\
		%&+O(1)|v_i|^2(1-{\color{black}{\chi_i^\e}}\xi_i\eta_i)\xi_i^\p \rho_i^\e\left(w_{i,x}-\frac{w_i}{v_i}v_{i,x}\right)\big(\sum\limits_{j}v_j\xi_j (\chi')^\e_j\hat{v}_j+ \sum\limits_{k=1}^{n} \eta_k v_k^2(1- \chi_k^\ep)\big)\nonumber\\
		%&{\color{black}{+O(1)\sum_k \varkappa^\ep_k\rho_k^\e\mathfrak{A}_k|v_k|}}\xi_i^\p \rho_i^\e\left(w_{i,x}-\frac{w_i}{v_i}v_{i,x}\right)\\
	\end{align}
 We have now using again Lemmas \ref{estimimpo1}, \ref{lemme6.6} and \ref{lemme6.5} and \eqref{ngtech5}
	\begin{align*}
	&\tilde{r}_{i,v}\mu_i \xi_i^\p v_{i,x}\left(w_{i,x}-\frac{w_i}{v_i}v_{i,x}\right)_x \pa_{v_i}\bar{v}_i+\tilde{r}_{i,v}\mu_i \xi_i^{\p\p}v_i \left(\frac{w_i}{v_i}\right)_x^2 v_{i,x} \pa_{v_i}\bar{v}_i=O(1)(\sum_j\La_j^1+\delta_0^2\La_j^3+\La_j^4)+\La_i^5.
	\end{align*}
	We proceed in a similar way for the term
	$\tilde{r}_{i,v}\mu_i \xi_i^\p v_{i,x}\left(w_{i,x}-\frac{w_i}{v_i}v_{i,x}\right) \pa_x(\pa_{v_i}\bar{v}_i)$. Let us consider now $\widetilde{\al}_{i}^{1,5}$, we observe easily again using \eqref{ngtech4} , Lemma \ref{estimimpo1} and the fact that $\mbox{supp}\xi$ is included in $\{x,\;\theta(x)=x\}$  that  
	\begin{align*}
	&\widetilde{\al}_{i}^{1,5}=\tilde{r}_{i,v}v_i\xi_i^\p\pa_{v_i}\bar{v}_i(\mu_iv_{i,x}-{\color{black}{(\tilde{\la}_i-\la_i^*)v_i-w_i}})\left(\frac{w_i}{v_i}\right)_x\nonumber\\
	&=O(1)\sum_j(\La_j^1+\La_j^4+\La_j^6)%+O(1)\sum\limits_{k=1}^{n}v_k^2(1-\eta_k \chi_k^\ep)\mathbbm{1}_{\{|\frac{w_k}{v_k}|\leq\frac{\delta_1}{2}\}}
+O(1)\sum\limits_{j}v_j\xi_j (\chi')^\e_j\hat{v}_j v_i\left(\frac{w_i}{v_i}\right)_x+O(1)\sum\limits_{k=1}^{n} \eta_k v_i \left(\frac{w_i}{v_i}\right)_x v_k^2(1- \chi_k^\ep)\nonumber\\
&=O(1)\sum_j(\La_j^1+\La_j^4+\La_j^6).
	\end{align*}
	 We have now
	\begin{align*}
	\widetilde{\al}_{i,x}^{1,5}&= (\pa_x(\tilde{r}_{i,v})v_i+\tilde{r}_{i,v}v_{i,x})    \xi_i^\p\pa_{v_i}\bar{v}_i \mathcal{H}_i \left(\frac{w_i}{v_i}\right)_x+ \tilde{r}_{i,v}v_i\xi_i^{\p\p}\left(\frac{w_i}{v_i}\right)_x^2\pa_{v_i}\bar{v}_i\mathcal{H}_i\\
	&+ \tilde{r}_{i,v}v_i\xi_i^\p\pa_x(\pa_{v_i}\bar{v}_i)\mathcal{H}_i\left(\frac{w_i}{v_i}\right)_x+\tilde{r}_{i,v}v_i\xi_i^\p\pa_{v_i}\bar{v}_i(\mu_iv_{i,x}-{\color{black}{(\tilde{\la}_i-\la_i^*)v_i-w_i}})_x\left(\frac{w_i}{v_i}\right)_x\\
	&+ \tilde{r}_{i,v}v_i\xi_i^\p\pa_{v_i}\bar{v}_i \mathcal{H}_i \left(\frac{w_i}{v_i}\right)_{xx}.
	 \end{align*}
	We can proceed as before for all the terms, we are however give some details for the terms $\tilde{r}_{i,v}v_i\xi_i^\p\pa_x(\pa_{v_i}\bar{v}_i)\mathcal{H}_i\left(\frac{w_i}{v_i}\right)_x$ and $\tilde{r}_{i,v}v_i\xi_i^\p\pa_{v_i}\bar{v}_i \mathcal{H}_i \left(\frac{w_i}{v_i}\right)_{xx}$. 
	From Lemma \ref{estimimpo1} and \eqref{6.54}, \eqref{ngtech5}, \eqref{ngtech4}  we deduce that
	\begin{align*}
	&\tilde{r}_{i,v}v_i\xi_i^\p\pa_x(\pa_{v_i}\bar{v}_i)\mathcal{H}_i\left(\frac{w_i}{v_i}\right)_x=O(1)\mathfrak{A}_i\rho_i^\e(|v_{i,x}|+\sum_{k\ne i}|v_kv_{k,x}|)|\mathcal{H}_i\left(\frac{w_i}{v_i}\right)_x|\nonumber\\
	&=O(1)\mathfrak{A}_i\rho_i^\e(O(1)v_i^2+O(1)\sum\limits_{j\ne i}v^4_j\left(\frac{w_j}{v_j}\right)^2_x \mathfrak{A}_j\rho^\e_j)|\left(\frac{w_i}{v_i}\right)_x|\nonumber\\
	&+O(1)\mathfrak{A}_i\rho_i^\e \sum_{k\ne i}|v_kv_{k,x}| |\left(\frac{w_i}{v_i}\right)_x|(|v_i|+\sum\limits_{j\ne i}v^2_j\left(\frac{w_j}{v_j}\right)_x\mathfrak{A}_j\rho^\e_j)\nonumber\\
	&=O(1)\sum_j(\La_j^1+\La_j^4).%+O(1)\sum\limits_{k=1}^{n}v_k^2(1-\eta_k \chi_k^\ep)\mathbbm{1}_{\{|\frac{w_k}{v_k}|\leq\frac{\delta_1}{2}\}}
%+O(1)\sum\limits_{j}v_j\xi_j (\chi')^\e_j\hat{v}_j+O(1)\sum\limits_{k=1}^{n} \eta_k v_k^2(1- \chi_k^\ep).
	\end{align*}
	From \eqref{identity:wi-vi-xx}, \eqref{ngtech5} and Lemmas \ref{lemme6.5},  \ref{lemme6.6}, \ref{lemme6.8} we have using the fact that $w_{i,xx}v_j^2\left(\frac{w_j}{v_j}\right)_x=\La_i^1$ when $j\ne i$
	 \begin{align*}
	  \tilde{r}_{i,v}v_i\xi_i^\p\pa_{v_i}\bar{v}_i\mathcal{H}_i \left(\frac{w_i}{v_i}\right)_{xx}&= \tilde{r}_{i,v}\xi_i^\p\pa_{v_i}\bar{v}_i\mathcal{H}_i \biggl((w_{i,xx}-\frac{w_i}{v_i}v_{i,xx})-2v_{i,x}\left(\frac{w_i}{v_i}\right)_{x}\biggl)\\
	 &=O(1)(\La_i^5+\sum_j( \La_j^1+\La_j^4)).
	 \end{align*}
	 Since $\widetilde{\al}_{i}^{1,6}=\tilde{r}_{i,v}\mu_i\xi^{\p\p}_i v_i\bar{v}_i\left(\frac{w_i}{v_i}\right)_x^2$, it implies that  $\widetilde{\al}_{i}^{1,6}=O(1)\delta_0^2\La_i^3$. Furthermore we have
	 \begin{align*}
	 \widetilde{\al}_{i,x}^{1,6}&=\pa_x(\tilde{r}_{i,v})\mu_i \xi^{\p\p}_i v_i\bar{v}_i\left(\frac{w_i}{v_i}\right)_x^2+ \tilde{r}_{i,v}\mu_{i,x}\xi^{\p\p}_i v_i\bar{v}_i\left(\frac{w_i}{v_i}\right)_x^2+\tilde{r}_{i,v}\mu_i\xi^{\p\p\p}_i v_i\bar{v}_i\left(\frac{w_i}{v_i}\right)_x^3\\
	 &+ \tilde{r}_{i,v}\mu_i\xi^{\p\p}_i (v_{i,x}\bar{v}_i+v_i\pa_x(\bar{v}_i))\left(\frac{w_i}{v_i}\right)_x^2
	 +2\tilde{r}_{i,v}\mu_i\xi^{\p\p}_i v_i\bar{v}_i\left(\frac{w_i}{v_i}\right)_{xx}\left(\frac{w_i}{v_i}\right)_x.
	  %&=\mathcal{O}(1)(\La_i^1+\La_i^3+\La_i^5).
	\end{align*}
	We can proceed as previously excepted for the term $2\tilde{r}_{i,v}\mu_i\xi^{\p\p}_i v_i\bar{v}_i\left(\frac{w_i}{v_i}\right)_{xx}\left(\frac{w_i}{v_i}\right)_x$, we have then from \eqref{identity:wi-vi-xx} and Lemmas \ref{lemme6.5}, \ref{lemme6.6}
	\begin{align*}
	&2\tilde{r}_{i,v}\mu_i\xi^{\p\p}_i v_i\bar{v}_i\left(\frac{w_i}{v_i}\right)_{xx}\left(\frac{w_i}{v_i}\right)_x=2\tilde{r}_{i,v}\mu_i\xi^{\p\p}_i \bar{v}_i\left(\frac{w_i}{v_i}\right)_x\big(
	 (w_{i,xx}-\frac{w_i}{v_i}v_{i,xx})-2v_{i,x}\left(\frac{w_i}{v_i}\right)_{x}\big)\\
	 &=O(1) |v_i|(1+\left(\frac{w_i}{v_i}\right)_x^2)(w_{i,xx}-\frac{w_i}{v_i}v_{i,xx})+O(1)\delta_0^2\La_i^3,\\
	 &=O(1)\delta_0^2\La_i^3+O(1)\La_i^5.
	\end{align*}
	Next we have in a direct way since $\mbox{supp}\xi$ is included in $\{x,\;\theta'(x)=1\}$, 
	\begin{align*}
	&\widetilde{\al}_{i}^{1,7}= -v_i^2\left(\frac{w_i}{v_i}\right)_x\pa_{v_i}\bar{v}_i\tilde{r}_{i,v}\xi_i% =v_i^2\left(\frac{w_i}{v_i}\right)_x\pa_{v_i}\bar{v}_i\tilde{r}_{i,v}[(\lambda_i^*-\theta_i)\xi'_i-\xi_i ]
	=O(1)\La_i^4 .
	\end{align*}
	 Furthermore it follows
	\begin{align*}
	 & \widetilde{\al}_{i,x}^{1,7}=  -\big(2 v_i v_{i,x}\left(\frac{w_i}{v_i}\right)_x+v_i^2\left(\frac{w_i}{v_i}\right)_{xx}\big)\pa_{v_i}\bar{v}_i\tilde{r}_{i,v}\xi_i-v_i^2\left(\frac{w_i}{v_i}\right)_{x}\pa_x(\pa_{v_i}\bar{v}_i)\tilde{r}_{i,v}\xi_i\\
	 &-v_i^2\left(\frac{w_i}{v_i}\right)_{x}\pa_{v_i}\bar{v}_i\pa_x(\tilde{r}_{i,v})\xi_i-v_i^2\left(\frac{w_i}{v_i}\right)^2_{x}\pa_{v_i}\bar{v}_i\tilde{r}_{i,v}\xi'_i.
	  \end{align*}
	  This term can be treated as previously. In a direct way using Lemma \ref{estimimpo1} we get $ \widetilde{\al}_{i}^{1,8} =O(1)\La_i^1$.
	% \begin{equation}
	% \begin{aligned}
	%&\left(\frac{w_i}{v_i}\right)_{xx}=\frac{1}{v_i}\left(w_{i,xx}-\frac{w_i}{v_i}v_{i,xx}\right)-\frac{2v_{i,x}}{v_i}\left(\frac{w_i}{v_i}\right)_{x}.\\
	%&\left(\frac{w_i}{v_i}\right)_{xxx}=-\frac{v_{i,x}}{v_i^2}\left(w_{i,xx}-\frac{w_i}{v_i}v_{i,xx}\right)
	%+\frac{1}{v_i}\left(w_{i,xxx}-\frac{w_i}{v_i}v_{i,xxx}\right)-\frac{1}{v_i}\left(\frac{w_i}{v_i}\right)_x v_{i,xx}\\
	%&-\frac{2v_{i,x}}{v_i}\left(\frac{w_i}{v_i}\right)_{xx}-2(\frac{v_{i,xx}}{v_i}-\frac{v_{i,x}^2}{v_i^2})\left(\frac{w_i}{v_i}\right)_{x}.
	%\end{aligned}
	%\label{identity:wi-vi-xx}
%\end{equation}
%ON SEMBLE AVOIR BESOIN D'ENERGIE POUR $w_{i,x}$!!!!
%$+\sum\limits_{i}\sum\limits_{j\neq i}\tilde{r}_{i,v}(\xi_i\pa_{v_j}\bar{v}_i)\left[(\mu_iv_{i,x}-\tilde{\la}_iv_i)v_{j,x}\right]+\sum\limits_{i}\sum\limits_{j\neq i}\tilde{r}_{i,v}(\xi_i\pa_{w_j}\bar{v}_i)\left[(\mu_iv_{i,x}-\tilde{\la}_iv_i)w_{j,x}\right]\\
%$
Similarly from  Lemmas \ref{estimimpo1}, \ref{lemme11.3} it follows
\begin{align*}
	 & \widetilde{\al}_{i,x}^{1,8}=\sum\limits_{j\neq i}\pa_x(\tilde{r}_{i,v})(\xi_i\pa_{v_j}\bar{v}_i)(\mu_iv_{i,x}-\tilde{\la}_iv_i)v_{j,x}+\sum\limits_{j\neq i}\tilde{r}_{i,v}\xi_i^\p\left(\frac{w_i}{v_i}\right)_x\pa_{v_j}\bar{v}_i(\mu_iv_{i,x}-\tilde{\la}_iv_i)v_{j,x}\\
	  &+\sum\limits_{j\neq i}\tilde{r}_{i,v}\xi_i\pa_x(\pa_{v_j}\bar{v}_i)(\mu_iv_{i,x}-\tilde{\la}_iv_i)v_{j,x}+\sum\limits_{j\neq i}\tilde{r}_{i,v}(\xi_i\pa_{v_j}\bar{v}_i)(\mu_iv_{i,x}-\tilde{\la}_iv_i)_xv_{j,x}\\
	  &+\sum\limits_{j\neq i}\tilde{r}_{i,v}(\xi_i\pa_{v_j}\bar{v}_i)(\mu_iv_{i,x}-\tilde{\la}_iv_i)v_{j,xx}\\
	  &=O(1)\sum\limits_{j\neq i}\rho_i^\e \mathfrak{A}_i |\left(\frac{w_i}{v_i}\right)_x| |v_j v_{i,x}v_{j,x}|+O(1)\sum_i \La_i^1
	  +\sum\limits_{j\neq i}\tilde{r}_{i,v}\xi_i\pa_x(\pa_{v_j}\bar{v}_i)(\mu_iv_{i,x}-\tilde{\la}_iv_i)v_{j,x}.
	  \end{align*}
	  Using again the Lemmas \ref{lemme6.5},\ref{lemme6.6} and \eqref{ngtech5} we have
	  \begin{align*}
	  \sum\limits_{j\neq i}\rho_i^\e\mathfrak{A}_i|\left(\frac{w_i}{v_i}\right)_x| |v_jv_{i,x}v_{j,x}|&=O(1)\La_i^1+\sum\limits_{j\neq i,k\ne i}\rho_i^\e\mathfrak{A}_i |\left(\frac{w_i}{v_i}\right)_x v_k^2 \left(\frac{w_k}{v_k}\right)_x|  \mathfrak{A}_k\rho^\e_k |v_jv_{j,x}|\\
	  &=O(1)\sum_k\La_k^1.
	  \end{align*}
	  Similarly we have using Lemma \ref{estimimpo1} and \eqref{ngtech5}
	  \begin{align*}
	  \sum\limits_{j\neq i}\tilde{r}_{i,v}\xi_i\pa_x(\pa_{v_j}\bar{v}_i)(\mu_iv_{i,x}-\tilde{\la}_iv_i)v_{j,x}&=O(1)\La_i^1+O(1)\sum_{i\ne j,l} \rho_i^\e\mathfrak{A}_i |\frac{v_j v_{j,x}}{v_i}|v_l^4\left(\frac{w_l}{v_l}\right)_x^2  \mathfrak{A}_l\rho^\e_l \\
	  &=O(1)\La_i^1+O(1)\delta_0^2\sum_l\La_l^3.
	  \end{align*}
	  We observe now using again Lemma \ref{estimimpo1}
	  \begin{align*}
	  	\widetilde{\alpha}_{i}^{1,9}&=\sum\limits_{j\neq i}\tilde{r}_{i,v}v_i\xi_i\pa_{v_iv_j}\bar{v}_i\left[(\mu_iv_{i,x}-(\tilde{\la}_i-\la_i^*)v_i-w_i)v_{j,x}\right]=O(1)\La_i^1,\\
	   % \alpha_{i}^{1,11}&= \sum\limits_{i}\sum\limits_{j\neq i}\tilde{r}_{i,v}v_i\xi_i\pa_{v_iw_j}\bar{v}_i\left[(\mu_iv_{i,x}-(\tilde{\la}_i-\la_i^*)v_i-w_i)w_{j,x}\right]=O(1)\La_i^1.
	  \end{align*}
	  We have in addition	  
	  \begin{align}
	  \widetilde{\al}_{i,x}^{1,9}&=\sum\limits_{j\neq i}\pa_x(\tilde{r}_{i,v})v_i\xi_i\pa_{v_iv_j}\bar{v}_i\mathcal{H}_i v_{j,x}+\sum\limits_{j\neq i}\tilde{r}_{i,v}\left[v_{i,x}\xi_i+v_i\xi_i^\p\left(\frac{w_i}{v_i}\right)_x\right]\pa_{v_iv_j}\bar{v}_i \mathcal{H}_i  v_{j,x}\nonumber\\
	  &+\sum\limits_{j\neq i}\tilde{r}_{i,v}v_i\xi_i\pa_x(\pa_{v_iv_j}\bar{v}_i)\mathcal{H}_i v_{j,x}+\sum\limits_{j\neq i}\tilde{r}_{i,v}v_i\xi_i\pa_{v_iv_j}\bar{v}_i \mathcal{H}_i  v_{j,xx}
	  \nonumber\\
	  &+\sum\limits_{j\neq i}\tilde{r}_{i,v}v_i\xi_i\pa_{v_iv_j}\bar{v}_i(\mu_iv_{i,x}-(\tilde{\la}_i-\la_i^*)v_i-w_i)_xv_{j,x}.\label{al110}
	  \end{align}
	 Again from Lemmas  \ref{estimimpo1}, \ref{lemme11.3}, \ref{lemme6.5}, \ref{lemme6.6} and \eqref{ngtech5} we get:
	  \begin{align*}
	  &\sum\limits_{j\neq i}\pa_x(\tilde{r}_{i,v})v_i\xi_i\pa_{v_iv_j}\bar{v}_i\mathcal{H}_iv_{j,x}=O(1)\La_i^1,\\
	  &\sum\limits_{j\neq i}\tilde{r}_{i,v}v_{i,x}\xi_i\pa_{v_iv_j}\bar{v}_i\mathcal{H}_iv_{j,x}=O(1)\sum\limits_{j\neq i}\mathfrak{A}_i\rho_i^\e |\frac{v_{j}v_{j,x}}{v_i}|(v_i^2+\sum_{l\ne i}v_l^4\left(\frac{w_l}{v_l}\right)_x^2\rho_k^\e\mathfrak{A}_k)\\
	  &=O(1)\La_i^1+O(1)\delta_0^2\sum_l \La_l^3,\\
	  &\sum\limits_{j\neq i}\tilde{r}_{i,v}v_i\xi_i\pa_x(\pa_{v_iv_j}\bar{v}_i)\mathcal{H}_i v_{j,x}\\
	  &=O(1)\sum_{j\ne i}\mathfrak{A}_i\rho_i^\e(|\frac{v_jv_{i,x}v_{j,x}}{v_i}|+\sum_{k\ne i}|v_{j,x}v_{k,x}|) (|v_i|+\sum_{l\ne i}v_l^2\left(\frac{w_l}{v_l}\right)_x)=O(1)\La_i^1+O(1)\sum_l \La_l^4.
	  \end{align*}
	 We proceed similarly for the other terms in \eqref{al110}. We have now
	   \begin{align*}
	  &\widetilde{\al}_i^{1,10}=\sum\limits_{j\neq i}\tilde{r}_{i,v}\mu_i\xi_i^\p\pa_{v_j}\bar{v}_iv_i\left(\frac{w_i}{v_i}\right)_xv_{j,x}=O(1)\La_i^1.
	%&+\sum\limits_{i}\sum\limits_{j\neq i}\tilde{r}_{i,v}v_i\xi_i\pa_{v_iv_j}\bar{v}_i\left[(\mu_iv_{i,x}-(\tilde{\la}_i-\la_i^*)v_i-\theta_i v_i)v_{j,x}\right]\nonumber\\
	%&+\sum\limits_{i}\sum\limits_{j\neq i}\tilde{r}_{i,v}v_i\xi_i\pa_{v_iw_j}\bar{v}_i\left[(\mu_iv_{i,x}-(\tilde{\la}_i-\la_i^*)v_i-\theta_i v_i)w_{j,x}\right]\nonumber\\
	%& +\sum\limits_{i}\sum\limits_{j\neq i}\tilde{r}_{i,v}\mu_i\xi_i^\p\pa_{v_j}\bar{v}_iv_i\left(\frac{w_i}{v_i}\right)_xv_{j,x}+\sum\limits_{i}\sum\limits_{j\neq i}\tilde{r}_{i,v}\mu_i\xi_i^\p\pa_{w_j}\bar{v}_iv_i\left(\frac{w_i}{v_i}\right)_xw_{j,x}\nonumber\\
	 \end{align*} 
We observe easily using again Lemma \ref{estimimpo1} and \eqref{ngtech5} that
	 %$+\sum\limits_{i}\sum\limits_{j\neq i}\tilde{r}_{i,v}\mu_i\xi_i^\p\pa_{v_j}\bar{v}_iv_i\left(\frac{w_i}{v_i}\right)_xv_{j,x}+\sum\limits_{i}\sum\limits_{j\neq i}\tilde{r}_{i,v}\mu_i\xi_i^\p\pa_{w_j}\bar{v}_iv_i\left(\frac{w_i}{v_i}\right)_xw_{j,x}$
	  %\begin{align*}
	  %&=\mathcal{O}(1)(\La_i^1+\La_i^3+\La_i^5),\\
	  %\al_{i,x}^{1,10}&=\sum\limits_{j\neq i}\pa_x(\tilde{r}_{i,v})v_i\xi_i\pa_{v_iw_j}\bar{v}_i(\mu_iv_{i,x}-(\tilde{\la}_i-\la_i^*)v_i-w_i)w_{j,x}\\
	  %&+\sum\limits_{j\neq i}\tilde{r}_{i,v}\left[v_{i,x}\xi_i+v_i\xi_i^\p\left(\frac{w_i}{v_i}\right)_x\right]\pa_{v_iw_j}\bar{v}_i(\mu_iv_{i,x}-(\tilde{\la}_i-\la_i^*)v_i-w_i)w_{j,x}\\
	  %&+\sum\limits_{j\neq i}\tilde{r}_{i,v}v_i\xi_i\pa_x(\pa_{v_iw_j}\bar{v}_i)(\mu_iv_{i,x}-(\tilde{\la}_i-\la_i^*)v_i-w_i)w_{j,x}\\
	  %&+\sum\limits_{j\neq i}\tilde{r}_{i,v}v_i\xi_i\pa_{v_iw_j}\bar{v}_i(\mu_iv_{i,x}-(\tilde{\la}_i-\la_i^*)v_i-w_i)_xw_{j,x}\\
	  %&+\sum\limits_{j\neq i}\tilde{r}_{i,v}v_i\xi_i\pa_{v_iw_j}\bar{v}_i(\mu_iv_{i,x}-(\tilde{\la}_i-\la_i^*)v_i-w_i)w_{j,xx}\\
	  %&=\mathcal{O}(1)(\La_i^1+\La_i^3+\La_i^5),\\
	  \begin{align*}
	  \widetilde{\al}_{i,x}^{1,10}&=\sum\limits_{j\neq i}\pa_x(\tilde{r}_{i,v})\mu_i\xi_i^\p\pa_{v_j}\bar{v}_iv_i\left(\frac{w_i}{v_i}\right)_xv_{j,x}+\sum\limits_{j\neq i}\tilde{r}_{i,v}\mu_{i,x}\xi_i^\p\pa_{v_j}\bar{v}_iv_i\left(\frac{w_i}{v_i}\right)_xv_{j,x}\\
	  &+\sum\limits_{j\neq i}\tilde{r}_{i,v}\mu_i\xi_i^{\p\p}\left(\frac{w_i}{v_i}\right)_x^2 \pa_{v_j}\bar{v}_iv_iv_{j,x}+\sum\limits_{j\neq i}\tilde{r}_{i,v}\mu_i\xi_i^\p\pa_x(\pa_{v_j}\bar{v}_i)v_i\left(\frac{w_i}{v_i}\right)_xv_{j,x}\\
	  &+\sum\limits_{j\neq i}\tilde{r}_{i,v}\mu_i\xi_i^\p\pa_{v_j}\bar{v}_iv_{i,x}\left(\frac{w_i}{v_i}\right)_xv_{j,x}+\sum\limits_{j\neq i}\tilde{r}_{i,v}\mu_i\xi_i^\p\pa_{v_j}\bar{v}_iv_i\left(\frac{w_i}{v_i}\right)_{xx}v_{j,x}\\
	  &+\sum\limits_{j\neq i}\tilde{r}_{i,v}\mu_i\xi_i^\p\pa_{v_j}\bar{v}_iv_i\left(\frac{w_i}{v_i}\right)_xv_{j,xx}\\
	  &=O(1)(\La_i^1+\delta_0^2\La_i^3).
	  \end{align*}
	 % In the same way we obtain
	  %\begin{align*}
	  %\al_{i,x}^{1,13}=O(1)(\La_i^1+\delta_0^2\La_i^3).
	  %\end{align*}
	  We have now again using \eqref{ngtech4} and Lemmas \ref{estimimpo1}, \ref{lemme6.5}, \ref{lemme6.6}
	  \begin{align*}
	  \widetilde{\al}_i^{1,11}&=\tilde{r}_{i,uv}\tilde{r}_i v_i\xi_i\pa_{v_i}\bar{v}_i\left[(\mu_iv_{i,x}-(\tilde{\la}_i-\la_i^*)v_i- w_i)v_i\right]\nonumber\\
	  &=O(1)\sum_j(\La_j^1+\La_j^4+\La_j^6)%+O(1)\sum\limits_{k=1}^{n}v_k^2(1-\eta_k \chi_k^\ep)\mathbbm{1}_{\{|\frac{w_k}{v_k}|\leq\frac{\delta_1}{2}\}}
+O(1)v_i^3\xi_j (\chi')^\e_i\hat{v}_i+O(1) \eta_i v_i^4(1- \chi_i^\ep).\nonumber
	  \end{align*}
	  It satisfies in particular \eqref{12.2}.  We have now using again \eqref{ngtech4}, \eqref{estimate-v-i-xx-1aa} and Lemmas  \ref{estimimpo1}, \ref{lemme6.5}, \ref{lemme6.6}, \ref{lemme9.7}
	 \begin{align*}
	  &\widetilde{\al}_{i,x}^{1,11}=\pa_x(\tilde{r}_{i,uv}\tilde{r}_i)v_i\xi_i\pa_{v_i}\bar{v}_i\mathcal{H}_iv_i+\tilde{r}_{i,uv}\tilde{r}_i \left(\big(2v_{i,x}v_i\xi_i+\xi_i^\p v^2_i\left(\frac{w_i}{v_i}\right)_x\big)\pa_{v_i}\bar{v}_i+v_i^2\xi_i\pa_x(\pa_{v_i}\bar{v}_i)\right)\mathcal{H}_i\\
	  &+\tilde{r}_{i,uv}\tilde{r}_i v^2_i\xi_i\pa_{v_i}\bar{v}_i(\mu_iv_{i,x}-(\tilde{\la}_i-\la_i^*)v_i-w_i)_x\\
	  &=O(1)\big(\sum_j\La_j^1+\La_j^4+\La_j^5+\La_j^6\big)+O(1)v_i\xi_i (\chi')^\e_i\hat{v}_i+O(1)\eta_i v_i^2(1- \chi_i^\ep)\\
	  &%+
		%O(1)\big(|v_i^3v_{i,x}|\eta_i(1-\chi^\e_i)+|v_i|^5\eta_i (1-\chi^\e_i)\big)%
		%+O(1)\sum_k|v_kv_{k,x}||1-\chi_k^\e\xi_k\eta_k|
		%+O(1)\sum\limits_{k}|v_k|^3(1-{\color{black}{\chi_k^\e}}\xi_k\eta_k)
		{\color{black}{+O(1)\varkappa^\ep_i\rho_i^\e\mathfrak{A}_i|v_i^3|}}.
	  \end{align*}
	  This estimates verifies in particular \eqref{12.2}. We have now using Lemma \ref{lemme11.4}
	  \begin{align*}
	   \widetilde{ \al}_{i,x}^{1,12}&=\tilde{r}_{i,uv}\tilde{r}_i \mu_i\xi_i^\p\bar{v}_i(w_{i,x}v_i-w_iv_{i,x}) =O(1)\La_i^4\\
	 \widetilde{ \al}_{i,x}^{1,12}&=\pa_x(\tilde{r}_{i,uv}\tilde{r}_i) \mu_i\xi_i^\p\bar{v}_i(w_{i,x}v_i-w_iv_{i,x})\\
	  &+\tilde{r}_{i,uv}\tilde{r}_i \left(\mu_{i,x}\xi_i^\p\bar{v}_i+\mu_i\xi_i^{\p\p}\left(\frac{w_i}{v_i}\right)_x\bar{v}_i+\mu_i\xi_i^\p\pa_x\bar{v}_i\right)(w_{i,x}v_i-w_iv_{i,x})\\
	  &+\tilde{r}_{i,uv}\tilde{r}_i \mu_i\xi_i^\p\bar{v}_i(w_{i,xx}v_i-w_iv_{i,xx})\\
	  &=O(1)(\La_i^4+\La_i^5).
	  \end{align*}
	  Next we observe that using Lemmas \ref{estimimpo1}, \ref{lemme6.5}, \ref{lemme6.6}, \ref{lemme11.4}
	  \begin{align*}
	  & \widetilde{\al}_{i}^{1,13}=\sum\limits_{j\neq i}\tilde{r}_{i,uv}\tilde{r}_j(\xi_i\pa_{v_i}\bar{v}_i)v_jv_i\mathcal{H}_i
	 % (\mu_iv_{i,x}-(\tilde{\la}_i-\la_i^*)v_i-w_i)
	 %  \sum\limits_{j\neq i}\tilde{r}_{i,uv}\tilde{r}_j(\xi_i\pa_{v_i}\bar{v}_i)v_jv_i(\mu_iv_{i,x}-(\tilde{\la}_i-\la_i^*)v_i-w_i)
	 =O(1)\La_i^1.\\
	  &\widetilde{\al}_{i,x}^{1,13}=\sum\limits_{j\neq i}\pa_x(\tilde{r}_{i,uv}\tilde{r}_j)(\xi_i\pa_{v_i}\bar{v}_i)v_jv_i \mathcal{H}_i%(\mu_iv_{i,x}-(\tilde{\la}_i-\la_i^*)v_i-w_i)
	  +\sum\limits_{j\neq i}\tilde{r}_{i,uv}\tilde{r}_j\left(\xi_i^\p\left(\frac{w_i}{v_i}\right)_xv_jv_i+\xi_iv_{j,x}v_i+\xi_iv_jv_{i,x}\right)(\pa_{v_i}\bar{v}_i)\mathcal{H}_i\\%(\mu_iv_{i,x}-(\tilde{\la}_i-\la_i^*)v_i-w_i)\\
	  %&
	  &+\sum\limits_{j\neq i}\tilde{r}_{i,uv}\tilde{r}_j\xi_i\pa_x(\pa_{v_i}\bar{v}_i)v_jv_i\mathcal{H}_i%(\mu_iv_{i,x}-(\tilde{\la}_i-\la_i^*)v_i-w_i)\\
	  +\sum\limits_{j\neq i}\tilde{r}_{i,uv}\tilde{r}_j(\xi_i\pa_{v_i}\bar{v}_i)v_jv_i(\mu_iv_{i,x}-(\tilde{\la}_i-\la_i^*)v_i-w_i)_x\\
	  &=O(1)\sum_j\La_j^1.
	  \end{align*}
	%  \begin{align*}
	  %&+\sum\limits_{i}\sum\limits_{j\neq i}\tilde{r}_{i,uv}\tilde{r}_j(w_{i,x}-(w_i/v_i)v_{i,x})\xi_i^\p\bar{v}_iv_j\mu_i\nonumber\\
	%&+\sum\limits_{i}\tilde{r}_{i,vv}(\xi_i\pa_{v_i}\bar{v}_i)^2v_iv_{i,x}(\mu_iv_{i,x}-(\tilde{\la}_i-\la_i^*)v_i-\theta_i v_i)\nonumber\\
	%\end{align*}
	  Newt we have using Lemmas \ref{estimimpo1}, \ref{lemme11.3}, \ref{lemme11.4}, \ref{lemme6.5}, \ref{lemme6.6} and \eqref{ngtech4}, \eqref{estimate-v-i-xx-1aa}
	  \begin{align*}
	  &\widetilde{\al}_{i,x}^{1,14}=\sum\limits_{j\neq i}\tilde{r}_{i,uv}\tilde{r}_j(w_{i,x}-(w_i/v_i)v_{i,x})\xi_i^\p\bar{v}_iv_j\mu_i=O(1)\La_i^1,\\
	 &\widetilde{ \al}_{i,x}^{1,14}=\sum\limits_{j\neq i}\pa_x(\tilde{r}_{i,uv}\tilde{r}_j)\left(w_{i,x}-\frac{w_i}{v_i}v_{i,x}\right)\xi_i^\p\bar{v}_iv_j \mu_i+\sum\limits_{j\neq i}\tilde{r}_{i,uv}\tilde{r}_j\left(w_{i,x}-\frac{w_i}{v_i}v_{i,x}\right)_x\xi_i^\p\bar{v}_iv_j\mu_i\\
	  &+\sum\limits_{j\neq i}\tilde{r}_{i,uv}\tilde{r}_j\left(w_{i,x}-\frac{w_i}{v_i}v_{i,x}\right)\xi_i^{\p\p}\left(\frac{w_i}{v_i}\right)_x\bar{v}_iv_j \mu_i\\
	  &+\sum\limits_{j\neq i}\tilde{r}_{i,uv}\tilde{r}_j\left(w_{i,x}-\frac{w_i}{v_i}v_{i,x}\right)\xi_i^\p[\pa_x\bar{v}_iv_j+\bar{v}_iv_{j,x}]\mu_i+\sum\limits_{j\neq i}\tilde{r}_{i,uv}\tilde{r}_j(w_{i,x}-(w_i/v_i)v_{i,x})\xi_i^\p\bar{v}_iv_j\mu_{i,x}\\
	  &=O(1)\La_i^1,\\%+\La_i^4+\La_i^5),\\
	  &\widetilde{\al}_{i}^{1,15}= \tilde{r}_{i,vv}(\xi_i\pa_{v_i}\bar{v}_i)^2v_iv_{i,x} \mathcal{H}_i\\%(\mu_iv_{i,x}-(\tilde{\la}_i-\la_i^*)v_i-w_i)
	 &\;\;\;\; \;\;\;=O(1)\sum_j(\La_j^1+\La_j^4+\La_j^6)%+O(1)\sum\limits_{k=1}^{n}v_k^2(1-\eta_k \chi_k^\ep)\mathbbm{1}_{\{|\frac{w_k}{v_k}|\leq\frac{\delta_1}{2}\}}
+O(1)v_i^2v_{i,x}\xi_i (\chi')^\e_i\hat{v}_i+O(1) \eta_i v_i^3v_{i,x}(1- \chi_i^\ep).\\
	  	  &\widetilde{\al}_{i,x}^{1,15}=\pa_x(\tilde{r}_{i,vv})(\xi_i\pa_{v_i}\bar{v}_i)^2v_iv_{i,x}\mathcal{H}_i%(\mu_iv_{i,x}-(\tilde{\la}_i-\la_i^*)v_i-w_i)\\
	  %&
	  +\tilde{r}_{i,vv}2(\xi_i\pa_{v_i}\bar{v}_i)\left(\xi_i^\p\left(\frac{w_i}{v_i}\right)_x\pa_{v_i}\bar{v}_i+\xi_i\pa_x(\pa_{v_i}\bar{v}_i)\right)v_iv_{i,x}\mathcal{H}_i\\%(\mu_iv_{i,x}-(\tilde{\la}_i-\la_i^*)v_i-w_i)\\
	  &+\tilde{r}_{i,vv}(\xi_i\pa_{v_i}\bar{v}_i)^2(v_iv_{i,xx}+v_{i,x}^2)\mathcal{H}_i%(\mu_iv_{i,x}-(\tilde{\la}_i-\la_i^*)v_i-w_i)\\
	 % &
	 +\tilde{r}_{i,vv}(\xi_i\pa_{v_i}\bar{v}_i)^2v_iv_{i,x}(\mathcal{H}_i)_x\\%(\mu_iv_{i,x}-(\tilde{\la}_i-\la_i^*)v_i-w_i)_x\\
	  &=O(1)\sum_j(\La_j^1+\La_j^4+\La_j^5+\La_j^6)%+O(1)\sum\limits_{k=1}^{n}v_k^2(1-\eta_k \chi_k^\ep)\mathbbm{1}_{\{|\frac{w_k}{v_k}|\leq\frac{\delta_1}{2}\}}
+O(1)|v_i\xi_i (\chi')^\e_i\hat{v}_i|+O(1) |\eta_i v_i^2(1- \chi_i^\ep)|
		%+O(1)\sum_k|v_kv_{k,x}||1-\chi_k^\e\xi_k\eta_k|
		%+O(1)\sum\limits_{k}|v_k|^3(1-{\color{black}{\chi_k^\e}}\xi_k\eta_k)
		{\color{black}{+O(1) |\varkappa^\ep_i\rho_i^\e\mathfrak{A}_i|v_i^2v_{i,x}|}}|.
	  \end{align*}
These estimates satisfy in particular \eqref{12.2}. Similarly we have using Lemmas \ref{estimimpo1}, \ref{lemme6.5}, \ref{lemme6.6},
	  \ref{lemme11.3} and \eqref{ngtech5}
	  \begin{align*}
	  &\widetilde{\al}_{i}^{1,16}= \tilde{r}_{i,vv}(\xi_i\pa_{v_i}\bar{v}_i)\xi_i^\p\bar{v}_i(w_{i,x}-(w_i/v_i)v_{i,x})v_{i,x}\mu_i=O(1)\La_i^4,\\
	  & \widetilde{\al}_{i,x}^{1,16}=[\pa_x(\tilde{r}_{i,vv})\xi_i+\tilde{r}_{i,vv}\xi_i^\p\left(\frac{w_i}{v_i}\right)_x]
	  \pa_{v_i}\bar{v}_i\xi_i^\p\bar{v}_i v_i \left(\frac{w_i}{v_i}\right)_xv_{i,x}\mu_i\\
	  &+\tilde{r}_{i,vv}\xi_i[\pa_x(\pa_{v_i}\bar{v}_i)\xi_i^\p+\xi_i^{\p\p}\left(\frac{w_i}{v_i}\right)_x\pa_{v_i}\bar{v}_i]
	  \bar{v}_iv_i \left(\frac{w_i}{v_i}\right)_xv_{i,x}\mu_i\\
	  &+\tilde{r}_{i,vv}(\xi_i\pa_{v_i}\bar{v}_i)\xi_i^\p\pa_x(\bar{v}_i)v_i  \left(\frac{w_i}{v_i}\right)_xv_{i,x}\mu_i+\tilde{r}_{i,vv}(\xi_i\pa_{v_i}\bar{v}_i)\xi_i^\p\bar{v}_i(w_{i,x}-(w_i/v_i)v_{i,x})_xv_{i,x}\mu_i\\
	  &+\tilde{r}_{i,vv}(\xi_i\pa_{v_i}\bar{v}_i)\xi_i^\p\bar{v}_i v_i\left(\frac{w_i}{v_i}\right)_x[v_{i,xx}\mu_i+v_{i,x}\mu_{i,x}]\\
	  &=O(1)(\delta_0^2\La_i^3+\La_i^4+\La_i^5+\La_i^1).
	  \end{align*}
	  %\begin{align*}
	 % &+\sum\limits_{i}\tilde{r}_{i,vv}\xi_i^\p\bar{v}_i\left(\frac{w_i}{v_i}\right)_x(v_i\xi_i\pa_{v_i}\bar{v}_i)\left[(\mu_iv_{i,x}-(\tilde{\la}_i-\la_i^*)v_i-\theta_i v_i)\right]\nonumber\\
	%&+\sum\limits_{i}\tilde{r}_{i,vv}(\xi_i^\p)^2\mu_i\bar{v}_i^{2}\left(\frac{w_i}{v_i}\right)_x(w_{i,x}-(w_i/v_i)v_{i,x})\nonumber\\
	%&+\sum\limits_{i}\sum\limits_{j\neq i}\tilde{r}_{i,vv}\xi_i^2\pa_{v_j}\bar{v}_i\pa_{v_i}\bar{v}_iv_i(\mu_iv_{i,x}-(\tilde{\la}_i-\la_i^*)v_i-\theta_i v_i)v_{j,x}\nonumber\\
	%\end{align*}
	  Again we have using Lemmas \ref{estimimpo1}, \ref{lemme6.5}, \ref{lemme6.6}, \ref{lemme6.8},
	  \ref{lemme11.3}  and \eqref{ngtech5}, \eqref{identity:wi-vi-xx}, \eqref{estimate-v-i-xx-1aa}
	  \begin{align*}
	   &\widetilde{\al}_{i}^{1,17}=\tilde{r}_{i,vv}\xi_i^\p\bar{v}_i\left(\frac{w_i}{v_i}\right)_x(v_i\xi_i\pa_{v_i}\bar{v}_i)\left[(\mu_iv_{i,x}-(\tilde{\la}_i-\la_i^*)v_i-w_i)\right]=O(1)\La_i^4\\
	 &\widetilde{ \al}_{i,x}^{1,17}=[\pa_x(\tilde{r}_{i,vv})\xi_i^\p+\xi_i^{\p\p}\left(\frac{w_i}{v_i}\right)_x]
	  \bar{v}_i\left(\frac{w_i}{v_i}\right)_x(v_i\xi_i\pa_{v_i}\bar{v}_i)\mathcal{H}_i%\left[(\mu_iv_{i,x}-(\tilde{\la}_i-\la_i^*)v_i-w_i)\right]\\
	  +\tilde{r}_{i,vv}\xi_i^\p\pa_x(\bar{v}_i)\left(\frac{w_i}{v_i}\right)_x(v_i\xi_i\pa_{v_i}\bar{v}_i)\mathcal{H}_i\\%\left[(\mu_iv_{i,x}-(\tilde{\la}_i-\la_i^*)v_i-w_i)\right]\\
	  &+\tilde{r}_{i,vv}\xi_i^\p\bar{v}_i\left(\frac{w_i}{v_i}\right)_{xx}(v_i\xi_i\pa_{v_i}\bar{v}_i)\mathcal{H}_i%\left[(\mu_iv_{i,x}-(\tilde{\la}_i-\la_i^*)v_i-w_i)\right]\\
	  +\tilde{r}_{i,vv}\xi_i^\p\bar{v}_i\left(\frac{w_i}{v_i}\right)_x\left(v_{i,x}\xi_i+\xi_i^\p\left(\frac{w_i}{v_i}\right)_x\right)\pa_{v_i}\bar{v}_i\mathcal{H}_i\\%\left[(\mu_iv_{i,x}-(\tilde{\la}_i-\la_i^*)v_i-w_i)\right]\\
	  &+\tilde{r}_{i,vv}\xi_i^\p\bar{v}_i\left(\frac{w_i}{v_i}\right)_xv_i\xi_i\pa_x(\pa_{v_i}\bar{v}_i)\mathcal{H}_i%\left[(\mu_iv_{i,x}-(\tilde{\la}_i-\la_i^*)v_i-w_i)\right]\\
	  +\tilde{r}_{i,vv}\xi_i^\p\bar{v}_i\left(\frac{w_i}{v_i}\right)_x\mathcal{H}_{i,x}\\%(v_i\xi_i\pa_{v_i}\bar{v}_i)\left[(\mu_iv_{i,x}-(\tilde{\la}_i-\la_i^*)v_i-w_i)\right]_x\\
	  &=O(1)\big(\delta_0^2\La_i^3+\sum_j(\La_j^1+\La_j^4+\La_j^5+\La_j^6)
		+|\rho_i^\e \mathfrak{A}_iv_i \left(\frac{w_i}{v_i}\right)_x|(|v_iv_{i,x}|\eta_i(1-\chi^\e_i)\\
		&+|v_i|^3\eta_i (1-\chi^\e_i)\big)
		%+O(1)\sum_k|v_kv_{k,x}||1-\chi_k^\e\xi_k\eta_k|
		%+O(1)\sum\limits_{k}|v_k|^3(1-{\color{black}{\chi_k^\e}}\xi_k\eta_k)
		{\color{black}{+ \varkappa^\ep_i|v_i|}}).
\end{align*}
This last estimate satisfies in particular \eqref{12.2}.
Similarly, we can obtain using Lemmas \ref{estimimpo1}, \ref{lemme6.5}, \ref{lemme6.6}, \ref{lemme11.3}
\begin{align*}
& \widetilde{\al}_{i}^{1,18}=\tilde{r}_{i,vv}(\xi_i^\p)^2\mu_i \bar{v}_i^{2}v_i\left(\frac{w_i}{v_i}\right)^2_x=O(1)\La_i^4,\\
	&\widetilde{\al}_{i,x}^{1,18}=\pa_x(\tilde{r}_{i,vv})(\xi_i^\p)^2\mu_i\bar{v}_i^2v_i\left(\frac{w_i}{v_i}\right)^2_x+\tilde{r}_{i,vv}2\xi_i^\p\xi_i^{\p\p}\mu_i\bar{v}_i^2v_i\left(\frac{w_i}{v_i}\right)^3_x\\
	&+\tilde{r}_{i,vv}(\xi_i^\p)^2\mu_{i,x}\bar{v}_i^2v_i\left(\frac{w_i}{v_i}\right)^2_x+\tilde{r}_{i,vv}(\xi_i^\p)^2\mu_i2\pa_x(\bar{v}_i)\bar{v}_i v_i\left(\frac{w_i}{v_i}\right)^2_x\\
	&+\tilde{r}_{i,vv}(\xi_i^\p)^2\mu_i\bar{v}_i^2v_{i,x}\left(\frac{w_i}{v_i}\right)^2_x+2\tilde{r}_{i,vv}(\xi_i^\p)^2\mu_i\bar{v}_i^2v_{i}\left(\frac{w_i}{v_i}\right)_{xx}\left(\frac{w_i}{v_i}\right)_{x}\\
		&=O(1)(\La_i^4+\delta_0^2\La_i^3+\La_i^5).
	\end{align*}
	Next using Lemmas \ref{estimimpo1}, \ref{lemme6.5}, \ref{lemme6.6}, \ref{lemme11.3}  we have
	\begin{align*}
	&\widetilde{\al}_{i}^{1,19}=\sum\limits_{j\neq i}\tilde{r}_{i,vv}\xi_i^2\pa_{v_j}\bar{v}_i\pa_{v_i}\bar{v}_iv_i(\mu_iv_{i,x}-(\tilde{\la}_i-\la_i^*)v_i-w_i)v_{j,x}
	%\sum\limits_{j\neq i}\tilde{r}_{i,vv}\xi_i^2\pa_{v_j}\bar{v}_i\pa_{v_i}\bar{v}_iv_i(\mu_iv_{i,x}-(\tilde{\la}_i-\la_i^*)v_i-\theta_i v_i)v_{j,x}
	=O(1)\La_i^1,\\
	&\widetilde{\al}_{i,x}^{1,19}=\sum\limits_{j\neq i}\pa_x(\tilde{r}_{i,vv})\xi_i^2\pa_{v_j}\bar{v}_i\pa_{v_i}\bar{v}_iv_i\mathcal{H}_i%(\mu_iv_{i,x}-(\tilde{\la}_i-\la_i^*)v_i-w_i)
	v_{j,x}+\sum\limits_{j\neq i}\tilde{r}_{i,vv}2\xi_i\xi_i^\p\left(\frac{w_i}{v_i}\right)_x\pa_{v_j}\bar{v}_i\pa_{v_i}\bar{v}_iv_i \mathcal{H}_i%(\mu_iv_{i,x}-(\tilde{\la}_i-\la_i^*)v_i-w_i)
	v_{j,x}\\
	&+\sum\limits_{j\neq i}\tilde{r}_{i,vv}\xi_i^2\pa_x(\pa_{v_j}\bar{v}_i)\pa_{v_i}\bar{v}_iv_i \mathcal{H}_i%(\mu_iv_{i,x}-(\tilde{\la}_i-\la_i^*)v_i-w_i)
	v_{j,x}+\sum\limits_{j\neq i}\tilde{r}_{i,vv}\xi_i^2\pa_{v_j}\bar{v}_i\pa_x(\pa_{v_i}\bar{v}_i)v_i \mathcal{H}_i%(\mu_iv_{i,x}-(\tilde{\la}_i-\la_i^*)v_i-w_i)
	v_{j,x}\\
	&+\sum\limits_{j\neq i}\tilde{r}_{i,vv}\xi_i^2\pa_{v_j}\bar{v}_i\pa_{v_i}\bar{v}_iv_i \mathcal{H}_{i,x}%(\mu_iv_{i,x}-(\tilde{\la}_i-\la_i^*)v_i-w_i)_x
	v_{j,x}+\sum\limits_{j\neq i}\tilde{r}_{i,vv}\xi_i^2\pa_{v_j}\bar{v}_i\pa_{v_i}\bar{v}_iv_i \mathcal{H}_i %(\mu_iv_{i,x}-(\tilde{\la}_i-\la_i^*)v_i-w_i)
	v_{j,xx}\\
	&+\sum\limits_{j\neq i}\tilde{r}_{i,vv}\xi_i^2\pa_{v_j}\bar{v}_i\pa_{v_i}\bar{v}_iv_{i,x}\mathcal{H}_i %(\mu_iv_{i,x}-(\tilde{\la}_i-\la_i^*)v_i-w_i)
	v_{j,x}=O(1)\La_i^1.
	\end{align*}
	%\begin{align*}
	%&+\sum\limits_{i}\sum\limits_{j\neq i}\tilde{r}_{i,vv}\mu_i\xi_i\xi^\p_i\pa_{v_j}\bar{v}_i\bar{v}_i\left(w_{i,x}-\frac{w_i}{v_i}v_{i,x}\right)v_{j,x}\nonumber\\
	%&+\sum\limits_{i}\sum\limits_{j\neq i}\tilde{r}_{i,vv}\mu_i\xi_i\xi^\p_i\pa_{w_j}\bar{v}_i\bar{v}_i\left(w_{i,x}-\frac{w_i}{v_i}v_{i,x}\right)w_{j,x}\nonumber\\
	%&-\sum\limits_{i}\tilde{r}_{i,v\si}v_i\xi_i\pa_{v_i}\bar{v}_i(\mu_iv_{i,x}-(\tilde{\la}_i-\la_i^*)v_i-\theta_i v_i)\theta_i^\p \left(\frac{w_i}{v_i}\right)_x\nonumber\\
	%&-{\color{black}{2}}\sum\limits_{i}\tilde{r}_{i,v\si}\xi_i^\p\bar{v}_i\mu_i(w_{i,x}-(w_i/v_i)v_{i,x}) \theta_i^\p \left(\frac{w_i}{v_i}\right)_x-\sum\limits_{i}\sum\limits_{k}\tilde{r}_{i,v\si}\xi_i\pa_{v_k}\bar{v}_iv_{k,x}\theta^\p_i\mu_i\left[w_{i,x}-\frac{w_i}{v_i}v_{i,x}\right]\nonumber\\
	%\end{align*}
	 We have now using Lemmas \ref{estimimpo1}, \ref{lemme6.5}, \ref{lemme6.6}, \ref{lemme11.3}
	\begin{align*}
	&\widetilde{\al}_{i}^{1,20}=\sum\limits_{j\neq i}\tilde{r}_{i,vv}\mu_i\xi_i\xi^\p_i\pa_{v_j}\bar{v}_i\bar{v}_iv_i\left(\frac{w_i}{v_i}\right)_xv_{j,x}=O(1)\La_i^4,\\
	&\widetilde{\al}_{i,x}^{1,20}=\sum\limits_{j\neq i}(\pa_x(\tilde{r}_{i,vv})\mu_i+\tilde{r}_{i,vv}\mu_{i,x})\xi_i\xi^\p_i\pa_{v_j}\bar{v}_i\bar{v}_i v_i\left(\frac{w_i}{v_i}\right)_xv_{j,x}\\
	&+\sum\limits_{j\neq i}\tilde{r}_{i,vv}\mu_ix\big[(\xi^\p_i)^2+\xi_i\xi^{\p\p}_i\big]\pa_{v_j}\bar{v}_i\bar{v}_i v_i\left(\frac{w_i}{v_i}\right)^2_xv_{j,x}\\
	&+\sum\limits_{j\neq i}\tilde{r}_{i,vv}\mu_i\xi_i\xi^\p_i[\pa_x(\pa_{v_j}\bar{v}_i)\bar{v}_i
	+\pa_{v_j}\bar{v}_i\pa_x(\bar{v}_i)]
	v_i\left(\frac{w_i}{v_i}\right)_x v_{j,x}\\
	&+\sum\limits_{j\neq i}\tilde{r}_{i,vv}\mu_i\xi_i\xi^\p_i\pa_{v_j}\bar{v}_i\bar{v}_i\left(w_{i,x}-\frac{w_i}{v_i}v_{i,x}\right)_xv_{j,x}+\sum\limits_{j\neq i}\tilde{r}_{i,vv}\mu_i\xi_i\xi^\p_i\pa_{v_j}\bar{v}_i\bar{v}_i v_i\left(\frac{w_i}{v_i}\right)_xv_{j,xx}\\
	&=O(1)(\La_i^1+\La_i^4+\La_i^5).
	\end{align*}
	 Newt we have using Lemmas \ref{estimimpo1}, \ref{lemme6.5}, \ref{lemme6.6}, \ref{lemme11.3}, \ref{lemme9.7} and \eqref{ngtech5}, \eqref{estimate-v-i-xx-1aa}
	 	\begin{align*}
	&\widetilde{\al}_{i}^{1,21}=-\tilde{r}_{i,v\si}v_i\xi_i\pa_{v_i}\bar{v}_i(\mu_iv_{i,x}-(\tilde{\la}_i-\la_i^*)v_i-w_i)\theta_i^\p \left(\frac{w_i}{v_i}\right)_x=O(1)\sum_j \La_j^4,\\
	&\widetilde{\al}_{i,x}^{1,21}=-\pa_x(\tilde{r}_{i,v\si})v_i\xi_i\pa_{v_i}\bar{v}_i\mathcal{H}_i
	%(\mu_iv_{i,x}-(\tilde{\la}_i-\la_i^*)v_i-w_i)
	\theta_i^\p \left(\frac{w_i}{v_i}\right)_x-\tilde{r}_{i,v\si}\left(v_{i,x}\xi_i+v_i\xi_i^\p\left(\frac{w_i}{v_i}\right)_x\right)\pa_{v_i}\bar{v}_i\mathcal{H}_i %(\mu_iv_{i,x}-(\tilde{\la}_i-\la_i^*)v_i-w_i)
	\theta_i^\p \left(\frac{w_i}{v_i}\right)_x\\
	&-\tilde{r}_{i,v\si}v_i\xi_i\pa_x(\pa_{v_i}\bar{v}_i)\mathcal{H}_i%(\mu_iv_{i,x}-(\tilde{\la}_i-\la_i^*)v_i-w_i)
	\theta_i^\p \left(\frac{w_i}{v_i}\right)_x-\tilde{r}_{i,v\si}v_i\xi_i\pa_{v_i}\bar{v}_i\mathcal{H}_{i,x}%(\mu_iv_{i,x}-(\tilde{\la}_i-\la_i^*)v_i-w_i)_x
	\theta_i^\p \left(\frac{w_i}{v_i}\right)_x\\
	&-\tilde{r}_{i,v\si}v_i\xi_i\pa_{v_i}\bar{v}_i\mathcal{H}_i%(\mu_iv_{i,x}-(\tilde{\la}_i-\la_i^*)v_i-w_i)
	\theta_i^{\p\p} \left(\frac{w_i}{v_i}\right)^2_x-\tilde{r}_{i,v\si}v_i\xi_i\pa_{v_i}\bar{v}_i\mathcal{H}_i%(\mu_iv_{i,x}-(\tilde{\la}_i-\la_i^*)v_i-w_i)
	\theta_i^\p \left(\frac{w_i}{v_i}\right)_{xx}\\
	&=O(1)\sum_j(\La_j^1+\La_j^4+\delta_0^2\La_j^3+\La_j^5+\La_j^6)+
		O(1)|\rho_i^\e v_i\left(\frac{w_i}{v_i}\right)_x| \mathfrak{A}_i \big(|v_iv_{i,x}|\eta_i(1-\chi^\e_i)\\
		&+|v_i|^3\eta_i (1-\chi^\e_i)%+O(1)\sum_k|v_kv_{k,x}||1-\chi_k^\e\xi_k\eta_k|
		%+O(1)\sum\limits_{k}|v_k|^3(1-{\color{black}{\chi_k^\e}}\xi_k\eta_k)
		{\color{black}{+O(1) \varkappa^\ep_i\rho_i^\e\mathfrak{A}_i|v_i|}}\big).
\end{align*}
This last estimate satisfies \eqref{12.2}.
	In the same way, we deduce from Lemmas \ref{estimimpo1}, \ref{lemme6.5}, \ref{lemme6.6}, \ref{lemme11.3} and the fact that $\left(\frac{w_i}{v_i}\right)_x\leq 1+\left(\frac{w_i}{v_i}\right)_x^2$
	\begin{align*}
	&\widetilde{\al}_{i}^{1,22}=-{\color{black}{2}}\tilde{r}_{i,v\si}\xi_i^\p\bar{v}_i\mu_i v_i \theta_i^\p \left(\frac{w_i}{v_i}\right)^2_x%=-{\color{black}{2}}\sum\limits_{i}\tilde{r}_{i,v\si}\xi_i^\p\bar{v}_i\mu_i(w_{i,x}-(w_i/v_i)v_{i,x}) \theta_i^\p \left(\frac{w_i}{v_i}\right)_x
	=O(1)\La_i^3\delta_0^2,\\
	&\widetilde{\al}_{i,x}^{1,22}=-\pa_x(\tilde{r}_{i,v\si})\xi_i^\p\bar{v}_i\mu_i v_i\theta_i^\p \left(\frac{w_i}{v_i}\right)^2_x-\tilde{r}_{i,v\si}\xi_i^{\p\p}\bar{v}_i\mu_iv_i\theta_i^\p \left(\frac{w_i}{v_i}\right)^3_x-\tilde{r}_{i,v\si}\xi_i^\p\pa_x(\bar{v}_i)\mu_i v_i\theta_i^\p \left(\frac{w_i}{v_i}\right)^2_x\\
	&-\tilde{r}_{i,v\si}\xi_i^\p\bar{v}_i\mu_{i,x}v_i\theta_i^\p \left(\frac{w_i}{v_i}\right)^2_x-\tilde{r}_{i,v\si}\xi_i^\p\bar{v}_i\mu_i v_{i,x}\theta_i^\p \left(\frac{w_i}{v_i}\right)^2_x-\tilde{r}_{i,v\si}\xi_i^\p\bar{v}_i\mu_{i}v_i\theta_i'' \left(\frac{w_i}{v_i}\right)^3_x\\
	&-2\tilde{r}_{i,v\si}\xi_i^\p\bar{v}_i\mu_{i}v_i\theta_i'\left(\frac{w_i}{v_i}\right)_x\left(\frac{w_i}{v_i}\right)_{xx},\\
	&=O(1)(\delta_0^2\La_i^3+\La_i^5).
	\end{align*}
	We get now using Lemmas \ref{estimimpo1}, \ref{lemme6.5}, \ref{lemme6.6}, \ref{lemme11.3} and \eqref{ngtech5}
	\begin{align*}
	\widetilde{\al}_{i}^{1,23}&=-\sum\limits_{k}\tilde{r}_{i,v\si}\xi_i\pa_{v_k}\bar{v}_iv_{k,x}\theta^\p_i\mu_iv_i\left(\frac{w_i}{v_i}\right)_x=O(1)(\La_i^1+\sum_j\La_j^4),\\
	\widetilde{\al}_{i,x}^{1,23}&=-\sum_k\pa_x(\tilde{r}_{i,v\si})\xi_i\pa_{v_k}\bar{v}_iv_{k,x}\theta^\p_i\mu_iv_i\left(\frac{w_i}{v_i}\right)_x-\sum_k\tilde{r}_{i,v\si}\xi_i'\pa_{v_k}\bar{v}_iv_{k,x}\theta^\p_i\mu_iv_i\left(\frac{w_i}{v_i}\right)_x^2\\
	&-\sum_k\tilde{r}_{i,v\si}\xi_i\pa_x(\pa_{v_k}\bar{v}_i)v_{k,x}\theta^\p_i\mu_iv_i\left(\frac{w_i}{v_i}\right)_x-\sum_k\tilde{r}_{i,v\si}\xi_i\pa_{v_k}\bar{v}_iv_{k,xx}\theta^\p_i\mu_iv_i\left(\frac{w_i}{v_i}\right)_x\\
	&-\sum_k\tilde{r}_{i,v\si}\xi_i\pa_{v_k}\bar{v}_iv_{k,x}\theta_i''\mu_iv_i\left(\frac{w_i}{v_i}\right)_x^2-\sum_k\tilde{r}_{i,v\si}\xi_i\pa_{v_k}\bar{v}_iv_{k,x}\theta_i'\mu_{i,x}v_i\left(\frac{w_i}{v_i}\right)_x\\
	&-\sum_k\tilde{r}_{i,v\si}\xi_i\pa_{v_k}\bar{v}_iv_{k,x}\theta_i'\mu_{i}v_{i,x}\left(\frac{w_i}{v_i}\right)_x-\sum_k\tilde{r}_{i,v\si}\xi_i\pa_{v_k}\bar{v}_iv_{k,x}\theta_i'\mu_{i}v_{i}\left(\frac{w_i}{v_i}\right)_{xx},\\
	&=O(1)\sum_j (\delta_0^2\La_j^3+\La_j^1+\La_j^4+\Lambda_j^5+\Lambda_j^6)+O(1)\mathfrak{A}_i\rho_i^\e|v_i\left(\frac{w_i}{v_i}\right)_x|
		\big(|v_iv_{i,x}|\eta_i(1-\chi^\e_i)\\
		&+|v_i|^3\eta_i (1-\chi^\e_i){\color{black}{+ \varkappa^\ep_i\rho_i^\e\mathfrak{A}_i|v_i|}}\big).
	\end{align*}
	This last estimate verifies \eqref{12.2}.
We are going to give some details for the most difficult term $-\sum_k\tilde{r}_{i,v\si}\xi_i\pa_{v_k}\bar{v}_iv_{k,xx}\theta^\p_i\mu_iv_i\left(\frac{w_i}{v_i}\right)_x$ when $k=i$. We have then from Lemmas \ref{estimimpo1}, \ref{lemme6.5}, \ref{lemme6.6} and \eqref{ngtech5}, \eqref{estimate-v-i-xx-1aabis}
\begin{align*}
&-\tilde{r}_{i,v\si}\xi_i\pa_{v_i}\bar{v}_iv_{i,xx}\theta^\p_i\mu_iv_i\left(\frac{w_i}{v_i}\right)_x=O(1)\sum_k(\La_k^1+\La_k^4+\La_k^5+\La_k^6)\\
&+O(1)\mathfrak{A}_i\rho_i^\e |v_i\left(\frac{w_i}{v_i}\right)_x|\big(|v_i|+|v_{i,x}|+|v_i\left(\frac{w_i}{v_i}\right)_x|\big)\\
&+O(1)\mathfrak{A}_i\rho_i^\e|v_i\left(\frac{w_i}{v_i}\right)_x|
		\big(|v_iv_{i,x}|\eta_i(1-\chi^\e_i)+|v_i|^3\eta_i (1-\chi^\e_i){\color{black}{+ \varkappa^\ep_i\rho_i^\e\mathfrak{A}_i|v_i|}}\big),\\
		&=O(1)\sum_k(\La_k^1+\delta_0^2\La_k^3+\La_k^4+\La_k^5+\La_k^6)\\
%&+O(1)\mathfrak{A}_i\rho_i^\e |v_i\left(\frac{w_i}{v_i}\right)_x|\big(|v_i|+|v_{i,x}|+|v_i\left(\frac{w_i}{v_i}\right)_x|\big)\\
&+O(1)\mathfrak{A}_i\rho_i^\e|v_i\left(\frac{w_i}{v_i}\right)_x|
		\big(|v_iv_{i,x}|\eta_i(1-\chi^\e_i)+|v_i|^3\eta_i (1-\chi^\e_i){\color{black}{+ \varkappa^\ep_i\rho_i^\e\mathfrak{A}_i|v_i|}}\big).
\end{align*}
	Next we have using Lemmas \ref{estimimpo1}, \ref{lemme11.4}
	\begin{align*}
	\widetilde{\al}_{i}^{1,24}&=
	\tilde{r}_{i,u\si}\tilde{r}_i \theta^\p_i \mu_i(w_iv_{i,x}-w_{i,x}v_i)=O(1)\La_i^4,\\
	\widetilde{\al}_{i,x}^{1,24}&=\pa_x(\tilde{r}_{i,u\si}\tilde{r}_i) \theta^\p_i \mu_i(w_iv_{i,x}-w_{i,x}v_i)+\tilde{r}_{i,u\si}\tilde{r}_i \theta^{\p\p}_i\left(\frac{w_i}{v_i}\right)_x \mu_i(w_iv_{i,x}-w_{i,x}v_i)\\
	&+\tilde{r}_{i,u\si}\tilde{r}_i \theta^\p_i \mu_{i,x}(w_iv_{i,x}-w_{i,x}v_i)+\tilde{r}_{i,u\si}\tilde{r}_i \theta^\p_i \mu_i(w_iv_{i,xx}-w_{i,xx}v_i)\\
	&=O(1)(\La_i^4+\delta_0^2\La_i^3+\La_i^5).
	\end{align*}
	Using the fact that $\tilde{r}_{i,u\si}=O(1)\xi_i\bar{v}_i$ we obtain using Lemmas \ref{estimimpo1}, \ref{lemme11.4}
	\begin{align*}
	\widetilde{\al}_{i}^{1,25}&=-\sum\limits_{j\neq i}\tilde{r}_{i,u\si}\tilde{r}_j\mu_i \theta^\p_iv_j\left[w_{i,x}-(w_i/v_i)v_{i,x}\right] =O(1)\La_i^1,\\
	\widetilde{\al}_{i,x}^{1,25}&=-\sum\limits_{j\neq i}\pa_x(\tilde{r}_{i,u\si}\tilde{r}_j)\mu_i \theta^\p_iv_j\left(w_{i,x}-\frac{w_i}{v_i}v_{i,x}\right)-\sum\limits_{j\neq i}\tilde{r}_{i,u\si}\tilde{r}_j\mu_{i,x} \theta^\p_iv_j\left(w_{i,x}-\frac{w_i}{v_i}v_{i,x}\right)\\
	&-\sum\limits_{j\neq i}\tilde{r}_{i,u\si}\tilde{r}_j\mu_i \theta^{\p\p}_i\left(\frac{w_i}{v_i}\right)_xv_j\left(w_{i,x}-\frac{w_i}{v_i}v_{i,x}\right)-\sum\limits_{j\neq i}\tilde{r}_{i,u\si}\tilde{r}_j\mu_i \theta^\p_iv_{j,x}\left(w_{i,x}-\frac{w_i}{v_i}v_{i,x}\right)\\
	&-\sum\limits_{j\neq i}\tilde{r}_{i,u\si}\tilde{r}_j\mu_i \theta^\p_iv_j\left(w_{i,x}-\frac{w_i}{v_i}v_{i,x}\right)_x\\
	&=O(1)(\La_i^1+\delta_0^2\La_i^3+\La_i^5).
	\end{align*}
	Next since $\tilde{r}_{i,\si},\tilde{r}_{i,\si \si} =O(1)|\xi_i v_i|$ we have from Lemmas \ref{estimimpo1}, \ref{lemme6.5}, \ref{lemme6.6}, \ref{lemme6.8}, \ref{lemme11.3}
	\begin{align*}
	&\widetilde{\al}_{i}^{1,26}=\tilde{r}_{i,\si\si} \mu_i(\theta_i^\p)^2v_i \left(\frac{w_i}{v_i}\right)^2_x=\delta_0^2\La_i^3,\\
	&\widetilde{\al}_{i,x}^{1,26}=(\pa_x(\tilde{r}_{i,\si\si}) \mu_i+\tilde{r}_{i,\si\si} \mu_{i,x})(\theta_i^\p)^2v_i \left(\frac{w_i}{v_i}\right)^2_x+2\tilde{r}_{i,\si\si} \mu_i\theta_i^\p\theta_i^{\p\p}v_i \left(\frac{w_i}{v_i}\right)^3_x\\
	&+\tilde{r}_{i,\si\si} \mu_i(\theta_i^\p)^2v_{i,x} \left(\frac{w_i}{v_i}\right)^2_x+2\tilde{r}_{i,\si\si} \mu_i(\theta_i^\p)^2v_{i} \left(\frac{w_i}{v_i}\right)_x\left(\frac{w_i}{v_i}\right)_{xx}\\%\\
	%&+\tilde{r}_{i,\si\si} \mu_i(\theta_i^\p)^2\left(\frac{w_i}{v_i}\right)_x\left(w_{i,x}-\frac{w_i}{v_i}v_{i,x}\right)_x\\
	&=O(1)(\delta_0^2\La_i^3+\La_i^5).\\[2mm]
	&\widetilde{\al}_{i}^{1,27}=-\tilde{r}_{i,\si}\mu_i\theta_i^{\p\p}v_i\left(\frac{w_i}{v_i}\right)^2_x=\delta_0^2\La_i^3,\\
%	\end{align*}
%Next since $\tilde{r}_{i,\si}=O(1)|v_i|$ we have from Lemmas \ref{lemme6.5} 
%	$\al_{i}^{1,34}=-\tilde{r}_{i,\si}\mu_i\theta_i^{\p\p}v_i\left(\frac{w_i}{v_i}\right)^2_x=\delta_0^2\La_i^3.$
%	\begin{align*}
	&\widetilde{\al}_{i,x}^{1,27}=-(\pa_x(\tilde{r}_{i,\si})\mu_i+\tilde{r}_{i,\si}\mu_{i,x})\theta_i^{\p\p}v_i\left(\frac{w_i}{v_i}\right)^2_x
	-\tilde{r}_{i,\si}\mu_i\theta_i^{\p\p\p}v_i\left(\frac{w_i}{v_i}\right)^3_x\\
	&-\tilde{r}_{i,\si}\mu_i\theta_i^{\p\p}v_{i,x}\left(\frac{w_i}{v_i}\right)^2_x-2\tilde{r}_{i,\si}\mu_i\theta_i^{\p\p}v_{i}\left(\frac{w_i}{v_i}\right)_x\left(\frac{w_i}{v_i}\right)_{xx},\\
	&=O(1)(\delta_0^2\La_i^3+\La_i^5).
\end{align*}
\end{proof}

\begin{lemma}
\label{lemmeal4al5}
	For $1\leq l\leq 8$ the following estimates hold 
	\begin{equation}
	\begin{aligned}
		&\al_{i}^{4,l}=O(1)\sum_j (\La_j^1+\delta_0^2\La_j^3+\La_j^4+\La_j^5+\La_j^6+\La_j^{6,1}),\\
		&\al_{i,x}^{4,l}=O(1)\sum_j (\La_j^1+\La_j^2+\delta_0^2\La_j^3+\La_j^4+\La_j^5+\La_j^6+\La_j^{6,1}).
		\end{aligned}
	\end{equation}
\end{lemma}
\begin{proof}
We observe that $\al_{i}^{4,1}=\sum\limits_{j\neq i}B(u)\tilde{r}_{i,v}\pa_{v_j}\bar{v}_i\left[v_{j,xx}v_i\xi_i+v_{j,x}v_{i,x}\xi_i+v_{j,x}v_i\xi^\p_i\left(\frac{w_i}{v_i}\right)_x\right]=
O(1)\La_i^1$ and we have
using again Lemmas \ref{lemme6.5}, \ref{lemme6.6}, \ref{estimimpo1}, \ref{lemme11.3} and \eqref{ngtech5}
	%It can be checked that $ \al_{i}^{4,l}=\mathcal{O}(1)(\La_i^1+\La_i^3)$. By using Lemma \ref{} and \eqref{} we can now estimate
	\begin{align*}
		\al_{i,x}^{4,1}&=\sum\limits_{j\neq i}\sum\limits_{k}v_k\tilde{r}_k\cdot DB(u)\tilde{r}_{i,v}\pa_{v_j}\bar{v}_i\left[v_{j,xx}v_i\xi_i+v_{j,x}v_{i,x}\xi_i+v_{j,x}v_i\xi^\p_i\left(\frac{w_i}{v_i}\right)_x\right]\\
		&+\sum\limits_{j\neq i}B(u)\pa_x(\tilde{r}_{i,v})\pa_{v_j}\bar{v}_i\left[v_{j,xx}v_i\xi_i+v_{j,x}v_{i,x}\xi_i+v_{j,x}v_i\xi^\p_i\left(\frac{w_i}{v_i}\right)_x\right]\\
		&+\sum\limits_{j\neq i}B\tilde{r}_{i,v}\pa_x(\pa_{v_j}\bar{v}_i)\left[v_{j,xx}v_i\xi_i+v_{j,x}v_{i,x}\xi_i+v_{j,x}v_i\xi^\p_i\left(\frac{w_i}{v_i}\right)_x\right]\\
		&+\sum\limits_{j\neq i}B\tilde{r}_{i,v}\pa_{v_j}\bar{v}_i\left[v_{j,xxx}v_i\xi_i+v_{j,xx}v_{i,x}\xi_i+v_{j,xx}v_i\xi_i^\p\left(\frac{w_i}{v_i}\right)_x\right]\\
		&+\sum\limits_{j\neq i}B\tilde{r}_{i,v}\pa_{v_j}\bar{v}_i\left[ v_{j,xx}v_{i,x}\xi_i+v_{j,x}v_{i,xx}\xi_i+v_{j,x}v_{i,x}\xi_i^\p\left(\frac{w_i}{v_i}\right)_x\right]\\
		&+\sum\limits_{j\neq i}B\tilde{r}_{i,v}\pa_{v_j}\bar{v}_i\left[v_{j,xx}v_i\xi^\p_i\left(\frac{w_i}{v_i}\right)_x+v_{j,x}v_{i,x}\xi^\p_i\left(\frac{w_i}{v_i}\right)_x\right]\\
		&+\sum\limits_{j\neq i}B\tilde{r}_{i,v}\pa_{v_j}\bar{v}_i\left[v_{j,x}v_i\xi^{\p\p}_i\left(\frac{w_i}{v_i}\right)^2_x+v_{j,x}v_i\xi^\p_i\left(\frac{w_i}{v_i}\right)_{xx}\right]\\
        &=O(1)(\La_i^1+\sum_j\delta_0^2\La_j^3)+O(1)\La_i^2.
        \end{align*}
        We have now $ \al_{i}^{4,2}=
        \sum\limits_{j\neq i}\sum\limits_{k}B(u)\tilde{r}_{i,v}v_i\pa_{v_jv_k}\bar{v}_i\big[\xi_iv_{j,x}v_{k,x}\big]=O(1)\La_i^1$ and we get using Lemmas  \ref{estimimpo1}, \ref{lemme6.5}, \ref{lemme6.6},  \ref{lemme11.3} and \eqref{ngtech5}
        \begin{align*}
        \al_{i,x}^{4,2}&=\sum\limits_{j\neq i}\sum_k [u_x\cdot DB\tilde{r}_{i,v}+B\pa_x(\tilde{r}_{i,v})]v_i\pa_{v_jv_k}\bar{v}_i\big[\xi_iv_{j,x}v_{k,x}\big]\\
        &+\sum\limits_{j\neq i} \sum_k B\tilde{r}_{i,v}v_{i,x}\pa_{v_jv_k}\bar{v}_i\big[\xi_iv_{j,x}v_{k,x}\big]+\sum\limits_{j\neq i} \sum_kB\tilde{r}_{i,v}v_{i}\pa_x(\pa_{v_jv_k}\bar{v}_i)\big[\xi_iv_{j,x}v_{k,x}\big]\\
        &+\sum\limits_{j\neq i}\sum_k B\tilde{r}_{i,v}v_{i}\pa_{v_jv_k}\bar{v}_i\left[\xi_i^\p\left(\frac{w_i}{v_i}\right)_xv_{j,x}v_{k,x}+\xi_iv_{j,xx}v_{k,x}+\xi_iv_{j,x}v_{k,xx}\right]\\
        &=O(1)\sum_j(\La_j^1+\delta_0^2\La_j^3+\La_j^4+\La_j^6).
        \end{align*}
        Let us explain how to deal with the most delicate term $\sum\limits_{j\neq i} \sum_kB\tilde{r}_{i,v}v_{i}\pa_x(\pa_{v_jv_k}\bar{v}_i)\big[\xi_iv_{j,x}v_{k,x}\big]$, we have using Lemmas \ref{estimimpo1}, \ref{lemme6.5}, \ref{lemme6.6} and \eqref{ngtech5}, \eqref{6.54}
        \begin{align*}
        &\sum\limits_{j\neq i} \sum_kB\tilde{r}_{i,v}v_{i}\pa_x(\pa_{v_jv_k}\bar{v}_i)\big[\xi_iv_{j,x}v_{k,x}\big]=\sum\limits_{j\neq i} \sum_{k\ne i} \sum_{l\ne i} B\tilde{r}_{i,v}v_{i}\pa_{v_l v_jv_k}\bar{v}_iv_{l,x}\big[\xi_iv_{j,x}v_{k,x}\big]\\
        &\sum\limits_{j\neq i} \sum_{k\ne i} B\tilde{r}_{i,v}v_{i}\pa_{v_i v_jv_k}\bar{v}_iv_{i,x}\big[\xi_iv_{j,x}v_{k,x}\big]+\sum\limits_{j\neq i} \sum_{l\ne i} B\tilde{r}_{i,v}v_{i}\pa_{v_l v_jv_i}\bar{v}_iv_{l,x}\big[\xi_iv_{j,x}v_{i,x}\big]\\
        &+\sum\limits_{j\neq i}  B\tilde{r}_{i,v}v_{i}\pa_{v_i v_jv_i}\bar{v}_iv_{i,x}\big[\xi_iv_{j,x}v_{i,x}\big]\\
        &=O(1)\sum\limits_{j\neq i} \sum_{k\ne i} \sum_{l\ne i}\mathfrak{A}_i\rho_i^\e \frac{ |v_l v_jv_k| }{|v_i|}|v_{l,x}v_{j,x}v_{k,x}|+O(1)\sum\limits_{j\neq i} \sum_{k\ne i}\mathfrak{A}_i\rho_i^\e |v_{i,x}v_{j,x}v_{k,x}|\\
        &+O(1)\sum\limits_{j\neq i} \sum_{l\ne i}\mathfrak{A}_i\rho_i^\e  |v_{l,x}v_{j,x}v_{i,x}|
        +O(1)\sum\limits_{j\neq i} \mathfrak{A}_i\rho_i^\e \frac{|v_j|}{|v_i|}|v_{i,x}^2v_{j,x}|\\
        &=O(1)\sum\limits_{j\neq i} \sum_{k\ne i} \sum_{l\ne i}\mathfrak{A}_i\rho_i^\e  | v_k| |v_{l,x}v_{j,x}v_{k,x}|+O(1)\La_i^1
        +O(1)\sum\limits_{j\neq i} \mathfrak{A}_i\rho_i^\e \frac{|v_j|}{|v_i|}|v_{i,x}^2v_{j,x}|\\
        &=O(1)\sum_j \La_j^1
        +
        O(1)\sum\limits_{j\neq i} \mathfrak{A}_i\rho_i^\e (|v_j||v_{i}v_{j,x}|+\sum_{l\ne i} \frac{|v_j v_{j,x}|}{|v_i|} v_l^4\left(\frac{w_l}{v_l}\right)_x^2 \mathfrak{A}_k \rho_k^\e)\\
        &+O(1) \sum\limits_{j\neq i} \mathfrak{A}_i\rho_i^\e  | v_j| |v_{j,x}^3|\\
        &=O(1)\sum_j (\La_j^1+\delta_0^2\La_j^3)+O(1) \sum\limits_{j\neq i} \mathfrak{A}_i\rho_i^\e  | v_j| \xi_j |v_{j,x}^3|+O(1) \sum\limits_{j\neq i} \mathfrak{A}_i\rho_i^\e  | v_j| |v_{j,x}^3|(1-\xi_j)\\
        &=O(1)\sum_j (\La_j^1+\delta_0^2\La_j^3+\La_j^6)+O(1) \sum\limits_{j\neq i} \mathfrak{A}_i\rho_i^\e  | v_j| \xi_j |v_{j,x}^2|(|v_j|+|v_i|    +\sum\limits_{l=1}^{n}|v^2_l\left(\frac{w_l}{v_l}\right)_x|\mathfrak{A}_l\rho^\e_l)  \\
         &=O(1)\sum_j (\La_j^1+\delta_0^2\La_j^3+\La_j^4+\La_j^6).
        \end{align*}
        We have now $ \al_{i}^{4,3}=
         \sum\limits_{j\neq i}B(u)\tilde{r}_{i,uv}\tilde{r}_i\pa_{v_j}\bar{v}_i\big[v_{j,x}v^2_i\xi_i\big]=O(1)\La_i^1$ and using Lemmas \ref{estimimpo1}, \ref{lemme6.5}, \ref{lemme6.6}, \ref{lemme11.4}
        \begin{align*}
         \al_{i,x}^{4,3}&=\sum\limits_{j\neq i}[u_x\cdot DB\tilde{r}_{i,uv}\tilde{r}_i+B\pa_x(\tilde{r}_{i,uv}\tilde{r}_i)]\pa_{v_j}\bar{v}_i\left[v_{j,x}v^2_i\xi_i\right]\\
         &+\sum\limits_{j\neq i}B\tilde{r}_{i,uv}\tilde{r}_i\pa_{v_j}\bar{v}_i\left[v_{j,xx}v^2_i\xi_i+2v_{i,x}v_{j,x}v_i\xi_i+v_{j,x}v^2_i\xi_i^\p\left(\frac{w_i}{v_i}\right)_x\right]\\
        &+\sum\limits_{j\neq i}B\tilde{r}_{i,uv}\tilde{r}_i\pa_x(\pa_{v_j}\bar{v}_i)\left[v_{j,x}v^2_i\xi_i\right]=O(1)(\La_i^1+\La_i^4).
         \end{align*}
         Next we have $ \al_{i}^{4,4}
         = \sum\limits_{j\neq i}\sum\limits_{k\neq i}B(u)\tilde{r}_{i,uv}\tilde{r}_k v_i\pa_{v_j}\bar{v}_i\big[v_{j,x}v_k\xi_i\big]=O(1)\La_i^1$ and as previously we obtain using Lemmas \ref{estimimpo1}, \ref{lemme6.5}, \ref{lemme6.6}, \ref{lemme11.4}
         \begin{align*}
         \al_{i,x}^{4,4}&=\sum\limits_{j\neq i}\sum\limits_{k\neq i}[u_x\cdot DB\tilde{r}_{i,uv}\tilde{r}_k+B\pa_x(\tilde{r}_{i,uv}\tilde{r}_k)] v_i\pa_{v_j}\bar{v}_i\left[v_{j,x}v_k\xi_i\right]\\
         &+\sum\limits_{j\neq i}\sum\limits_{k\neq i} B\tilde{r}_{i,uv}\tilde{r}_k\left[ v_{i,x}\pa_{v_j}\bar{v}_i+v_i\pa_x(\pa_{v_j}\bar{v}_i)\right]\left[v_{j,x}v_k\xi_i\right]\\
         &+\sum\limits_{j\neq i}\sum\limits_{k\neq i} B \tilde{r}_{i,uv}\tilde{r}_k v_i\pa_{v_j}\bar{v}_i\left[v_{j,xx}v_k\xi_i+v_{j,x}v_{k,x}\xi_i+v_{j,x}v_k\xi_i^\p\left(\frac{w_i}{v_i}\right)_x\right]\\
         &=O(1)\La_i^1.
         \end{align*}
 We observe now that $\al_{i}^{4,5}=\sum\limits_{j\neq i}\sum\limits_{k}B(u)\tilde{r}_{i,vv}\pa_{v_j}\bar{v}_i\pa_{v_k}\bar{v}_i\big[v_{k,x}v_{j,x}v_i\xi_i^2\big]%=\sum\limits_{j\neq i}\sum\limits_{k}B(u)\tilde{r}_{i,vv}\pa_{v_j}\bar{v}_i\pa_{v_k}\bar{v}_i\big[v_{k,x}v_{j,x}v_i\xi_i^2\big]
 =O(1)\La_i^1$ and using Lemmas \ref{estimimpo1},  \ref{lemme11.3}
         \begin{align*}
         %\al_{i,x}^{4,6}&=\sum\limits_{j\neq i}[u_x\cdot DB\tilde{r}_{i,vv}+B\pa_x(\tilde{r}_{i,vv})]\pa_{v_j}\bar{v}_i\pa_{v_i}\bar{v}_i\left[v_{i,x}v_{j,x}v_i\xi_i^2\right]\\
         %&+\sum\limits_{j\neq i}B\tilde{r}_{i,vv}\pa_x(\pa_{v_j}\bar{v}_i)\pa_{v_i}\bar{v}_i\left[v_{i,x}v_{j,x}v_i\xi_i^2\right]+\sum\limits_{j\neq i}B\tilde{r}_{i,vv}\pa_{v_j}\bar{v}_i\pa_x(\pa_{v_i}\bar{v}_i)\left[v_{i,x}v_{j,x}v_i\xi_i^2\right]\\
         %&+\sum\limits_{j\neq i}B\tilde{r}_{i,vv}\pa_{v_j}\bar{v}_i\pa_{v_i}\bar{v}_i\left[(v_{i,xx}v_{j,x}v_i+v_{i,x}v_{j,xx}v_i+v^2_{i,x}v_{j,x})\xi_i^2+2v_{i,x}v_{j,x}v_i\xi_i^\p\xi_i\left(\frac{w_i}{v_i}\right)_x\right]\\
         %&=\mathcal{O}(1)(\La_i^1+\La_i^3+\La_i^5),\\
        \al_{i,x}^{4,5}&= \sum\limits_{j\neq i}\sum\limits_{k\neq i}[u_x\cdot DB\tilde{r}_{i,vv}+B\pa_x(\tilde{r}_{i,vv})]\pa_{v_j}\bar{v}_i\pa_{v_k}\bar{v}_i\left[v_{k,x}v_{j,x}v_i\xi_i^2\right]\\
        &+ \sum\limits_{j\neq i}\sum\limits_{k\neq i}B\tilde{r}_{i,vv}[\pa_x(\pa_{v_j}\bar{v}_i)\pa_{v_k}\bar{v}_i+\pa_{v_j}\bar{v}_i\pa_x(\pa_{v_k}\bar{v}_i)]\left[v_{k,x}v_{j,x}v_i\xi_i^2\right]\\
        &+ \sum\limits_{j\neq i}\sum\limits_{k\neq i}B\tilde{r}_{i,vv}\pa_{v_j}\bar{v}_i\pa_{v_k}\bar{v}_i\left[(v_{k,xx}v_{j,x}v_i+v_{k,x}v_{j,xx}v_i+v_{k,x}v_{i,x}v_{j,x})\xi_i^2+2v_{k,x}v_{j,x}v_i\xi_i^\p\xi_i\left(\frac{w_i}{v_i}\right)_x\right]\\
         &=O(1)\La_i^1.
         \end{align*}
          We have now $\al_{i}^{4,6}=\sum\limits_{j\neq i}B(u)\tilde{r}_{i,vv}\pa_{v_j}\bar{v}_i\left[v_{j,x}v_i\xi_i\xi_i^\p\left(\frac{w_i}{v_i}\right)_x\bar{v}_i\right]
         =O(1)\La_i^1$ and using Lemmas \ref{estimimpo1}, \ref{lemme6.5}, \ref{lemme6.6}, \ref{lemme11.3}
         \begin{align*}
        %\al_{i,x}^{4,9}&= \sum\limits_{j\neq i}\sum\limits_{k\neq i}[u_x\cdot DB\tilde{r}_{i,vv}+B\pa_x(\tilde{r}_{i,vv})]\pa_{v_j}\bar{v}_i\pa_{w_k}\bar{v}_i\left[w_{k,x}v_{j,x}v_i\xi_i^2\right]\\
        %&+ \sum\limits_{j\neq i}\sum\limits_{k\neq i}B\tilde{r}_{i,vv}[\pa_x(\pa_{v_j}\bar{v}_i)\pa_{w_k}\bar{v}_i+\pa_{v_j}\bar{v}_i\pa_x(\pa_{w_k}\bar{v}_i)]\left[w_{k,x}v_{j,x}v_i\xi_i^2\right]\\
        %&+ \sum\limits_{j\neq i}\sum\limits_{k\neq i}B\tilde{r}_{i,vv}\pa_{v_j}\bar{v}_i\pa_{v_k}\bar{v}_i\left[(w_{k,xx}v_{j,x}v_i+w_{k,x}v_{j,xx}v_i+w_{k,x}v_{i,x}v_{j,x})\xi_i^2+2w_{k,x}v_{j,x}v_i\xi_i^\p\xi_i\left(\frac{w_i}{v_i}\right)_x\right]\\
         %&=\mathcal{O}(1)(\La_i^1+\La_i^3+\La_i^5),\\
       &  \al_{i,x}^{4,6}=\sum\limits_{j\neq i}[u_x\cdot D B\tilde{r}_{i,vv}+B\pa_x(\tilde{r}_{i,vv})]\pa_{v_j}\bar{v}_i\left[v_{j,x}v_i\xi_i\xi_i^\p\left(\frac{w_i}{v_i}\right)_x\bar{v}_i\right]\\
         &+\sum\limits_{j\neq i}B\tilde{r}_{i,vv}\pa_x(\pa_{v_j}\bar{v}_i)\left[v_{j,x}v_i\xi_i\xi_i^\p\left(\frac{w_i}{v_i}\right)_x\bar{v}_i\right]+\sum\limits_{j\neq i}B\tilde{r}_{i,vv}\pa_{v_j}\bar{v}_i\left[v_{j,xx}v_i+v_{j,x}v_{i,x}\right]\xi_i\xi_i^\p\left(\frac{w_i}{v_i}\right)_x\bar{v}_i\\
         &+\sum\limits_{j\neq i}B\tilde{r}_{i,vv}\pa_{v_j}\bar{v}_iv_{j,x}v_i\left[\xi_i\xi_i^{\p\p}+(\xi^\p_i)^2\right]\left(\frac{w_i}{v_i}\right)^2_x\bar{v}_i+\sum\limits_{j\neq i}B\tilde{r}_{i,vv}\pa_{v_j}\bar{v}_i\left[v_{j,x}v_i\xi_i\xi_i^\p\left(\frac{w_i}{v_i}\right)_{xx}\bar{v}_i\right]\\
         &+\sum\limits_{j\neq i}B\tilde{r}_{i,vv}\pa_{v_j}\bar{v}_i\left[v_{j,x}v_i\xi_i\xi_i^\p\left(\frac{w_i}{v_i}\right)_x\pa_x(\bar{v}_i)\right]\\
         &=O(1)(\La_i^1+\delta_0^2\La_i^3).
         \end{align*}
         In a similar manner we have $ \al_{i}^{4,7}=-\sum\limits_{j\neq i} B(u)\tilde{r}_{i,v\si}\pa_{v_j}\bar{v}_i\left[v_{j,x}v_i\xi_i\left(\frac{w_i}{v_i}\right)_x\theta^\p_i\right]%=-\sum\limits_{j\neq i} B(u)\tilde{r}_{i,v\si}\pa_{v_j}\bar{v}_i\left[v_{j,x}v_i\xi_i\left(\frac{w_i}{v_i}\right)_x\theta^\p_i\right]
         =O(1)\La_i^1$ and using Lemmas  \ref{estimimpo1}, \ref{lemme6.5}, \ref{lemme6.6}, \ref{lemme11.3} and \eqref{ngtech5}
         \begin{align*}
      & \al_{i,x}^{4,7}= -\sum\limits_{j\neq i}[u_x\cdot B\tilde{r}_{i,v\si}+B\pa_x(\tilde{r}_{i,v\si})]\pa_{v_j}\bar{v}_i\left[v_{j,x}v_i\xi_i\left(\frac{w_i}{v_i}\right)_x\theta^\p_i\right]\\
         &-\sum\limits_{j\neq i} B(u)\tilde{r}_{i,v\si}\pa_x(\pa_{v_j}\bar{v}_i)\left[v_{j,x}v_i\xi_i\left(\frac{w_i}{v_i}\right)_x\theta^\p_i\right]-\sum\limits_{j\neq i} B(u)\tilde{r}_{i,v\si}\pa_{v_j}\bar{v}_i\left[(v_{j,xx}v_i+v_{j,x}v_{i,x})\xi_i\left(\frac{w_i}{v_i}\right)_x \theta^\p_i\right]\\
         &-\sum\limits_{j\neq i} B(u)\tilde{r}_{i,v\si}\pa_{v_j}\bar{v}_i\left[v_{j,x}v_i(\xi_i^\p+\theta^{\p\p}_i)\left(\frac{w_i}{v_i}\right)^2_x\right]-\sum\limits_{j\neq i} B(u)\tilde{r}_{i,v\si}\pa_{v_j}\bar{v}_iv_{j,x}v_i\xi_i\left(\frac{w_i}{v_i}\right)_{xx}\theta^\p_i\\
         &=O(1)\La_i^1+O(1)\sum_k \delta_0^2\La_k^3.
         \end{align*}
         We have finally $ \al_{i}^{4,8}=\sum\limits_{j\neq i}v_{j,x}v_i\xi_i\pa_{v_j}\bar{v}_iu_x\cdot DB(u)\tilde{r}_{i,v}=O(1)\La_i^1$ and using again Lemmas \ref{estimimpo1}, \ref{lemme11.3}
         \begin{align*}
         \al_{i,x}^{4,8}&=  \sum\limits_{j\neq i}[v_{j,xx}v_i+v_{j,x}v_{i,x}]\xi_i\pa_{v_j}\bar{v}_iu_x\cdot DB\tilde{r}_{i,v}+  \sum\limits_{j\neq i}v_{j,x}v_i\xi_i^\p\left(\frac{w_i}{v_{i}}\right)_x\pa_{v_j}\bar{v}_iu_x\cdot DB\tilde{r}_{i,v}\\
         &+  \sum\limits_{j\neq i}v_{j,x}v_i\xi_i\pa_x(\pa_{v_j}\bar{v}_i)u_x\cdot DB\tilde{r}_{i,v}+  \sum\limits_{j\neq i}v_{j,x}v_i\xi_i\pa_{v_j}\bar{v}_iu_{xx}\cdot DB\tilde{r}_{i,v}\\
         &+\sum\limits_{j\neq i}v_{j,x}v_i\xi_i\pa_{v_j}\bar{v}_i(u_x\otimes u_x):D^2B\tilde{r}_{i,v}+  \sum\limits_{j\neq i}v_{j,x}v_i\xi_i\pa_{v_j}\bar{v}_iu_x\cdot DB\pa_x(\tilde{r}_{i,v})\\
         &=O(1)\La_i^1.
	\end{align*}
This completes the estimates of $\al_{i,x}^{4,l}$ for $1\leq  l\leq 8$.
\end{proof}

\begin{lemma}
\label{lemmeal6}
	For $1\leq l\leq 10$ it follows
	\begin{align}
			&\al_{i}^{5,l}=O(1)\sum_j(\La_j^1+\La_j^2+\delta_0^2\La_j^3+\La_j^4+\La_j^5+\La_j^6+\La_j^{6,1})+R_\e^{5,i,l},\nonumber\\
			&\al_{i,x}^{5,l}=O(1)\sum_j(\La_j^1+\La_j^2+\delta_0^2\La_j^3+\La_j^4+\La_j^5+\La_j^6+\La_j^{6,1})+\widetilde{R}_\e^{5,i,l},\label{12.11i}
	\end{align}
	with:
	\begin{equation}
	\int_{\hat{t}}^T\int_{\R}(|R_\e^{5,i,l}|+|\widetilde{R}_\e^{5,i,l}|)dx ds=O(1)\delta_0^2,
	\label{12.11}
	\end{equation}
	for $\e>0$ small enough in terms of $T-\hat{t}$ and $\delta_0$.
\end{lemma}
\begin{proof}
Denote $\Xi_i=v_{i,x}-\mu_i^{-1}(\tilde{\la}_i-\la_i^*+\theta_i)v_i$, we observe that
\begin{align*}
	&J_{5,x}
	=\sum\limits_{i}\Xi_{i,x}(B-\mu_iI_n)[\tilde{r}_{i}+v_i\xi_i\pa_{v_i}\bar{v}_i\tilde{r}_{i,v}]+\sum\limits_{i}\Xi_i v_i\tilde{r}_i\cdot D(B-\mu_iI_n)[\tilde{r}_{i}+v_i\xi_i\pa_{v_i}\bar{v}_i\tilde{r}_{i,v}]\\
	&+\sum\limits_{i}\sum\limits_{j\neq i}\Xi_i v_j\tilde{r}_j\cdot D(B-\mu_iI_n)[\tilde{r}_{i}+v_i\xi_i\pa_{v_i}\bar{v}_i\tilde{r}_{i,v}]+\sum\limits_{i}\Xi_i v_{i}(B-\mu_iI_n)[\tilde{r}_{i,u}\tilde{r}_i+v_i\xi_i\pa_{v_i}\bar{v}_i\tilde{r}_{i,uv}\tilde{r}_{i}]\\
	&+\sum\limits_{i}\sum\limits_{j\neq i}\Xi_i v_{j}(B-\mu_iI_n)[\tilde{r}_{i,u}\tilde{r}_j+v_i\xi_i\pa_{v_i}\bar{v}_i\tilde{r}_{i,uv}\tilde{r}_{j}]\\
	&+\sum\limits_{i}\Xi_i v_{i,x}\xi_i(B-\mu_iI_n)[(2\pa_{v_i}\bar{v}_i+v_i\pa_{v_iv_i}\bar{v}_i)\tilde{r}_{i,v}+\xi_iv_i(\pa_{v_i}\bar{v}_i)^2\tilde{r}_{i,vv}]\\
	&+\sum\limits_{i}\sum\limits_{j\neq i}\Xi_i v_{j,x}\xi_i(B-\mu_iI_n)[(\pa_{v_j}\bar{v}_i+v_i\pa_{v_iv_j}\bar{v}_i)\tilde{r}_{i,v}+\xi_iv_i\pa_{v_i}\bar{v}_i\pa_{v_j}\bar{v}_i\tilde{r}_{i,vv}]\\
	%&+\sum\limits_{i}\sum\limits_{j\neq i}(v_{i,x}-\mu_i^{-1}(\tilde{\la}_i-\la_i^*+\theta_i)v_i)w_{j,x}\xi_i(B-\mu_iI_n)[(\pa_{w_j}\bar{v}_i+v_i\pa_{v_iw_j}\bar{v}_i)\tilde{r}_{i,v}+\xi_iv_i\pa_{v_i}\bar{v}_i\pa_{w_j}\bar{v}_i\tilde{r}_{i,vv}]\\
	&+\sum\limits_{i}\Xi_i \xi^\p_i\left(\frac{w_i}{v_i}\right)_x(B-\mu_iI_n)[(v_i\pa_{v_i}\bar{v}_i+\bar{v}_i)\tilde{r}_{i,v}+\xi_iv_i\pa_{v_i}\bar{v}_i\bar{v}_i\tilde{r}_{i,vv}]\\
	&-\sum\limits_{i}\Xi_i \theta_i^\p\left(\frac{w_i}{v_i}\right)_x(B-\mu_iI_n)[\tilde{r}_{i,\si}+\xi_iv_i\pa_{v_i}\bar{v}_i\tilde{r}_{i,v\si}]\\
	&=:\sum\limits_{l=1}^{9}\sum\limits_{i}\al_{i}^{5,l}.
\end{align*}
	We first observe that $\al_{i}^{5,1}=\Xi_{i,x}(B-\mu_iI_n)[\tilde{r}_{i}+v_i\xi_i\pa_{v_i}\bar{v}_i\tilde{r}_{i,v}]=\Xi_{i,x}\sum_{j\ne i}(\mu_j(u)-\mu_i(u))(\psi_{i,j}+v_i\xi_i\pa_{v_i}\bar{v}_i \psi_{i,j, v})r_{j}(u)$, due to the fact that $\psi_{i,j}=O(1)\xi_i\bar{v}_i$ we deduce that since $\theta(s)=s$ on $\mbox{supp}\xi$ that
 $$\al_{i}^{5,1}=(\mu_i^{-1}\mathcal{H}_i)_x\sum_{j\ne i}(\mu_j(u)-\mu_i(u))(\psi_{i,j}+v_i\xi_i\pa_{v_i}\bar{v}_i \psi_{ij, v})r_{j}(u).$$	
	Using now the Lemmas \ref{estimimpo1}, \ref{lemme6.5}, \ref{lemme6.6} and \eqref{ngtech4}, \eqref{estimate-v-i-xx-1aa}
	%\ref{lemme9.7}, \ref{lemme6.8} 
	and the fact that  $\psi_{i,j}=O(1)\xi_i\bar{v}_i$ we deduce that
	 \begin{align*}
	 & \al_{i}^{5,1}=O(1)\sum_j(\La_j^1+\La_j^4+\La_j^5+\La_j^6)%+O(1)\sum\limits_{k=1}^{n}v_k^2(1-\eta_k \chi_k^\ep)\mathbbm{1}_{\{|\frac{w_k}{v_k}|\leq\frac{\delta_1}{2}\}}
+O(1)v_i^2\xi_i (\chi')^\e_i\hat{v}_i+O(1) \eta_i v_i^3(1- \chi_i^\ep)\\
&+O(1) |v_i^2v_{i,x}|\eta_i(1-\chi^\e_i)
		{\color{black}{+O(1) \varkappa^\ep_i\rho_i^\e\mathfrak{A}_i|v_i^2v_{i,x}|}}.
\end{align*}
We observe that $\al_i^{5,1}$ satisfies \eqref{12.11}.
Let us estimate now $\al_{i,x}^{5,1}$
	\begin{align}
		\al_{i,x}^{5,1}&=(\mu_i^{-1}\mathcal{H}_i)_{xx}\sum\limits_{j\neq i}(\mu_j-\mu_i)[\psi_{ij}+v_i\xi_i\pa_{v_i}\bar{v}_i\psi_{ij,v}]r_j\nonumber\\
		&+(\mu_i^{-1}\mathcal{H}_i)_{x}\sum\limits_{j\neq i}(\mu_j-\mu_i)_x[\psi_{ij}+v_i\xi_i\pa_{v_i}\bar{v}_i\psi_{ij,v}]r_j\nonumber\\
		&+(\mu_i^{-1}\mathcal{H}_i)_{x}\sum\limits_{j\neq i}(\mu_j-\mu_i)[\pa_x(\psi_{ij})+v_{i,x}\xi_i\pa_{v_i}\bar{v}_i\psi_{ij,v}+v_i\xi_i^\p\left(\frac{w_i}{v_i}\right)_x\pa_{v_i}\bar{v}_i\psi_{ij,v}]r_j\nonumber\\
		&+(\mu_i^{-1}\mathcal{H}_i)_{x}\sum\limits_{j\neq i}(\mu_j-\mu_i)[v_i\xi_i\pa_x(\pa_{v_i}\bar{v}_i)\psi_{ij,v}+v_i\xi_i\pa_{v_i}\bar{v}_i\pa_x(\psi_{ij,v})]r_j\nonumber\\
		&+(\mu_i^{-1}\mathcal{H}_i)_{x}\sum\limits_{j\neq i}(\mu_j-\mu_i)[\psi_{ij}+v_i\xi_i\pa_{v_i}\bar{v}_i\psi_{ij,v}]u_x\cdot Dr_j.\label{ali1x}
		%&=\mathcal{O}(1)(\La_i^1+\La_i^3+\La_i^5)
		\end{align}
	First we are going to study the term $(\mu_i^{-1}\mathcal{H}_i)_{xx}\sum\limits_{j\neq i}(\mu_j-\mu_i)[\psi_{ij}+v_i\xi_i\pa_{v_i}\bar{v}_i\psi_{ij,v}]r_j$. We observe that%We then have
	%\begin{align*}
	%&\Xi_{i1,xx}=\pa_{xx}\mu_i^{-1}(\mu_i v_{i,x}-(\tilde{\la}_i-\la_i^*)v_i-w_i)+2\pa_x(\mu_i^{-1})
	%(\mu_i v_{i,x}-(\tilde{\la}_i-\la_i^*)v_i-w_i)_x\\
	%&+\mu_{i}^{-1}(\mu_i v_{i,x}-(\tilde{\la}_i-\la_i^*)v_i-w_i)_{xx}
	%\end{align*}
%We deduce that
	\begin{align*}
	&(\mu_i^{-1}\mathcal{H}_i)_{xx}\sum\limits_{j\neq i}(\mu_j-\mu_i)[\psi_{ij}+v_i\xi_i\pa_{v_i}\bar{v}_i\psi_{ij,v}]r_j=\sum_{i=1}^3 A_i\\
	&=(\pa_{xx}\mu_i^{-1}\mathcal{H}_i+2\pa_x(\mu_i^{-1})\mathcal{H}_{i,x}+\mu_{i}^{-1}\mathcal{H}_{i,xx})\sum\limits_{j\neq i}(\mu_j-\mu_i)[\psi_{ij}+v_i\xi_i\pa_{v_i}\bar{v}_i\psi_{ij,v}]r_j.
	%&+2\pa_x(\mu_i^{-1})+2\pa_x(\mu_i^{-1})\mathcal{H}_{i,x}+\mu_{i}^{-1}\mathcal{H}_{i,xx})
	%\mathcal{H}_{i,x}\sum\limits_{j\neq i}(\mu_j-\mu_i)[\psi_{ij}+v_i\xi_i\pa_{v_i}\bar{v}_i\psi_{ij,v}]r_j\\
	%&+\mu_{i}^{-1}\mathcal{H}_{i,xx}\sum\limits_{j\neq i}(\mu_j-\mu_i)[\psi_{ij}+v_i\xi_i\pa_{v_i}\bar{v}_i\psi_{ij,v}]r_j
	\end{align*}
	Let us start by studying $A_3$, from the Lemmas \ref{lemme6.5}, \ref{lemme6.6}, \eqref{6.45}, \eqref{9.15primefin} and the fact that $\psi_{ij}=O(1)\xi_i\bar{v}_i$ we deduce that
	\begin{align}
	A_3&=O(1)\biggl(%IMP\sum_k\sum_{l\ne k}(|v_{k,x}|+|w_{k,x}|+|v_k|)(|v_lv_{l,x}|+|v_kv_l|)+\sum_{k}\sum_{l\ne k}|v_k|(|v_{l,xx}v_l|+|v_kv_{l,x}|)
		%\\
		%&IMP+\sum_{k}\sum_{l\ne k}\sum_{p\ne k}|v_k| |v_{l,x}||v_{p,x}|
		%+\sum_k\sum_l\sum_{j\ne l}|v_kv_l||v_{j,x}v_j|\\
		%&+
		\sum_j(\La_j^1+\La_j^2+\delta_0^2\La_j^3+\La_j^4+\La_j^5+\La_j^6)+
		%&+\sum_k\sum_{l\ne k} |v_k|(|v_{k,x}v_l|+|w_{k,x}v_l|+|v_kv_{l,x}|+|v_kv_l|)+\sum_k\sum_l\sum_{j\ne l}|v_kv_l||v_{j,x}v_j|\\
		%&+\sum_k \sum_{l\ne k}|v_k|(|v_lv_{l,x}|+|v_kv_l|)
		%\\[5mm]
		\nonumber\\
		&+O(1)
(|v_iv_{i,x}|+|v_iv_{i,xx}|+|v_i^2|)|v_i|\mathbbm{1}_{\{v_i^{2N}\leq 2\e\}}
\nonumber\\
		%&+\mathfrak{A}_i\rho_i^\e\big(|v_iv_{i,x}\left(\frac{w_i}{v_i}\right)_x|+|v_i^2\left(\frac{w_i}{v_i}\right)^2_x|+|v_i^2\left(\frac{w_i}{v_i}\right)_x|\\
		%&+|w_{i,xx}v_i-v_{i,xx}w_i|\big)\\
		%&+\mathfrak{A}_i \Delta_i^\e\big(
		%|v_i v_{i,x}|+v_{i,x}^2+|v_iv_{i,xx}|\big)\\
		%&
		%+(|v_{i,x}|+|v_i|)v_i^2(1-{\color{black}{\chi_i^\e}}\xi_i\eta_i)\\
		%&+\mathfrak{A}_i\rho_i^\e\big(|v_i^3\left(\frac{w_i}{v_i}\right)_x|+
		%|v_i^2v_{i,x}|\Delta_i^\e\big)\\
		%&+(|v_{i,x}|+|v_i|)v_i^2(1-{\color{black}{\chi_i^\e}}\xi_i\eta_i)\\
		%&+\mathfrak{A}_i\rho_i^\e\big(|v_i^3\left(\frac{w_i}{v_i}\right)_x|+|v_i^2 v_{i,x}|\Delta_i^\e\big)+v_i^2v_{i,x}{\color{black}{(\chi')_i^\e}}\xi_i\eta_i\\
		%&+v_iv_{i,x}(\chi')_i^\e \mathfrak{A}_i\rho_i^\e
		%{\color{black}{+v_{i,x}^2(\chi')^\e_i\xi_i\eta_i }}+v_i^2v_{i,x}{\color{black}{(\chi')_i^\e}}\xi_i\eta_i
		%+(|v_i|+|v_{i,x}|)\varkappa^\ep_i\mathfrak{A}_i\rho_i^\e\\
		&+|v_i|(|v_i |(|v_{i,x}|+|v_{i,xx}|)+|v_{i,x}|(|v_{i,x}|+|w_{i,x}|))\varkappa^\ep_i\mathfrak{A}_i\rho_i^\e\biggl).
		\label{I3tech2}
	\end{align}
	We observe that $A_3$ satisfies \eqref{12.11}. Similarly from \eqref{ngtech4}, \eqref{estimate-v-i-xx-1aa} and the fact that $\pa_{xx}(\mu_i^{-1}), \pa_{x}(\mu_i^{-1})=O(1)$, $\psi_{i,j}=O(1)\xi_i\bar{v}_i$  we obtain that
	\begin{align}
	&A_1=O(1)\sum_j(\La_j^1+\La_j^4+\La_j^6)%+O(1)\sum\limits_{k=1}^{n}v_k^2(1-\eta_k \chi_k^\ep)\mathbbm{1}_{\{|\frac{w_k}{v_k}|\leq\frac{\delta_1}{2}\}}
+O(1)v_i^2\xi_i (\chi')^\e_i\hat{v}_i+O(1) \eta_i v_i^3(1- \chi_i^\ep)
,\nonumber\\
	&A_2=O(1)\sum_k(\La_k^1+\La_k^4+\La_k^5+\La_k^6)+
		O(1)\big(|v_i^2v_{i,x}|\eta_i(1-\chi^\e_i)+|v_i|^4\eta_i (1-\chi^\e_i)\big)
		%+O(1)\sum_k|v_kv_{k,x}||1-\chi_k^\e\xi_k\eta_k|
		%+O(1)\sum\limits_{k}|v_k|^3(1-{\color{black}{\chi_k^\e}}\xi_k\eta_k)
		\nonumber\\
		&{\color{black}{+O(1) \varkappa^\ep_i\rho_i^\e\mathfrak{A}_i|v_i^2 v_{i,x}|}}.
	%O(1)\sum_k(\La_k^1+\La_k^4+\La_k^5+\La_k^6).
	\label{I3tech3}
	\end{align}
	From 
	 \eqref{I3tech2} and \eqref{I3tech3} we deduce that
	\begin{align}
	&(\mu_i^{-1}\mathcal{H}_i)_{xx}\sum\limits_{j\neq i}(\mu_j-\mu_i)[\psi_{ij}+v_i\xi_i\pa_{v_i}\bar{v}_i\psi_{ij,v}]r_j&=
	O(1)\sum_k(\La_k^1+\La_k^2+\delta_0^2\La_k^3+\La_k^4+\La_k^5+\La_k^6)\nonumber\\
	&+\widetilde{R}_\e^{5,i,1,1},
	\label{I3tech4}
	\end{align}
	and  $\widetilde{R}_\e^{5,i,1,1}$ satisfies \eqref{12.11}.
Using now	  \eqref{ali1x}, \eqref{I3tech4}, Lemmas \ref{lemme6.5}, \ref{lemme6.6},\ref{estimimpo1}, \ref{lemme11.3a}, \eqref{ngtech4}, \eqref{estimate-v-i-xx-1aa}  and the fact that $\sum_{j\ne i}v_jv_{j,x}(\mu_i^{-1}\mathcal{H}_i)_x=O(1)\La_i^1$ we obtain
\begin{align}
&\al_{i,x}^{5,1}=O(1)\sum_k(\La_k^1+\La_k^2+\delta_0^2\La_k^3+\La_k^4+\La_k^5+\La_k^6)+\widetilde{R}_\e^{5,i,1},
\end{align}
and  $\widetilde{R}_\e^{5,i,1}$ satisfies \eqref{12.11}. Let us study now $\al_{i}^{5,2}=\Xi_i v_i\tilde{r}_i\cdot D(B-\mu_iI_n)[\tilde{r}_{i}+v_i\xi_i\pa_{v_i}\bar{v}_i\tilde{r}_{i,v}]$, first we observe that
		\begin{equation}
		\Xi_i= \mu_i^{-1}\mathcal{H}_i+\mu_i^{-1}(w_i-\theta_i v_i).
		\label{chichi1}
		\end{equation}
Now due to the choice of $\theta$ (see \eqref{choixtheta}) and \eqref{6.54} we obtain for $\delta_0$ sufficiently small in terms of $\delta_1$ that 
\begin{equation}
\begin{aligned}
|w_i-\theta_i v_i|&=O(1) |w_i| \mathbbm{1}_{\{|\frac{w_i}{v_i}|\geq\delta_1\}}\\
&=O(1)\mathbbm{1}_{\{|\frac{w_i}{v_i}|\geq\delta_1\}}(|v_{i,x}|+\sum_j(\La_j^1+\La_j^4)+\sum_{k\ne i} |v_k|^2(1-\chi^\e_k\xi_k\eta_k)+\sum_{k\ne i} v_k\xi_k(\chi')^\e_k\hat{v_k}).
\end{aligned}
\label{chicle}
\end{equation}
		 We deduce using \eqref{ngtech4}, \eqref{6.45}, \eqref{chichi1}, \eqref{chicle} and Lemmas \ref{estimimpo1}, \ref{lemme6.5}, \ref{lemme6.6} that
		$$\al_{i}^{5,2}=O(1)\sum_j(\La_j^1+\La_j^4+\La_j^6)+O(1)v_i^2\xi_j (\chi')^\e_i\hat{v}_i
		+O(1) \eta_i v_i^3(1- \chi_i^\ep),
		$$
		which satisfies \eqref{12.11}. We have now 
		\begin{equation}
|(w_i-\theta_i v_i)_x|\leq (|w_{i,x}|+|\theta'_i\left(\frac{w_i}{v_i}\right)_x v_i|+|\theta_i v_{i,x}|) \mathbbm{1}_{\{|\frac{w_i}{v_i}|\geq\frac{\delta_1}{2}\}}.
\label{chicle1}
\end{equation}
		We obtain now for $\al_{i,x}^{5,2}$ by using Lemmas \ref{estimimpo1}, \ref{lemme6.5}, \ref{lemme6.6}, \ref{lemme11.3} and \eqref{ngtech4}, \eqref{6.45}, \eqref{estimate-v-i-xx-1aa}, \eqref{chichi1}, \eqref{chicle}
		\begin{align*}
		\al_{i,x}^{5,2}&=\big((\mu_i^{-1}\mathcal{H}_i)_x+(\mu_i^{-1}(w_i-\theta_i v_i))_x\big)v_i\tilde{r}_i\cdot D(B-\mu_iI_n)[\tilde{r}_{i}+v_i\xi_i\pa_{v_i}\bar{v}_i\tilde{r}_{i,v}]\\
		&+\Xi_iv_{i,x}\tilde{r}_i\cdot D(B-\mu_iI_n)[\tilde{r}_{i}+v_i\xi_i\pa_{v_i}\bar{v}_i\tilde{r}_{i,v}]+\Xi_iv_i\tilde{r}_{i,x}\cdot D(B-\mu_iI_n)[\tilde{r}_{i}+v_i\xi_i\pa_{v_i}\bar{v}_i\tilde{r}_{i,v}]\\
		&+\Xi_iv_iu_x\otimes \tilde{r}_i: D^2(B-\mu_iI_n)[\tilde{r}_{i}+v_i\xi_i\pa_{v_i}\bar{v}_i\tilde{r}_{i,v}]\\
		&+\Xi_iv_i\tilde{r}_i\cdot D(B-\mu_iI_n)[\tilde{r}_{i,x}+v_{i,x}\xi_i\pa_{v_i}\bar{v}_i\tilde{r}_{i,v}+v_i\xi_i^\p\left(\frac{w_i}{v_i}\right)_x\pa_{v_i}\bar{v}_i\tilde{r}_{i,v}]\\
		&+\Xi_iv_i\tilde{r}_i\cdot D(B-\mu_iI_n)[v_i\xi_i\pa_x(\pa_{v_i}\bar{v}_i)\tilde{r}_{i,v}+v_i\xi_i\pa_{v_i}\bar{v}_i\pa_x(\tilde{r}_{i,v})]\\
		&=O(1)\sum_j(\La_j^1+\La_j^4+\La_j^5+\La_j^6+\La_j^{6,1})+\widetilde{R}_\e^{5,i,2},
		\end{align*}
		with $\widetilde{R}_\e^{5,i,2}$ satisfying \eqref{12.11}.
		Next we have  $\al_{i}^{5,3}=\sum\limits_{j\neq i}\Xi_i v_j\tilde{r}_j\cdot D(B-\mu_iI_n)[\tilde{r}_{i}+v_i\xi_i\pa_{v_i}\bar{v}_i\tilde{r}_{i,v}]=
		O(1)\sum_j\La_j^1.$ Furthermore we obtain using Lemmas \ref{estimimpo1}, \ref{lemme6.5}, \ref{lemme6.6}, \ref{lemme11.3}
		\begin{align*}
		\al_{i,x}^{5,3}&=\sum\limits_{j\neq i}[\big((\mu_i^{-1}\mathcal{H}_i)_x+(\mu_i^{-1}(w_i-\theta_i v_i))_x\big)v_j+\Xi_i v_{j,x}]\tilde{r}_j\cdot D(B-\mu_iI_n)[\tilde{r}_{i}+v_i\xi_i\pa_{v_i}\bar{v}_i\tilde{r}_{i,v}]\\
		&+\sum\limits_{j\neq i}\Xi_{i}v_j[\tilde{r}_{j,x}\cdot D(B-\mu_iI_n)+u_x\otimes \tilde{r}_j: D^2(B-\mu_iI_n)][\tilde{r}_{i}+v_i\xi_i\pa_{v_i}\bar{v}_i\tilde{r}_{i,v}]\\
		&+\sum\limits_{j\neq i}\Xi_{i}v_j\tilde{r}_j\cdot D(B-\mu_iI_n)[\tilde{r}_{i,x}+v_{i,x}\xi_i\pa_{v_i}\bar{v}_i\tilde{r}_{i,v}+v_i\xi_i^\p\left(\frac{w_i}{v_i}\right)_x\pa_{v_i}\bar{v}_i\tilde{r}_{i,v}]\\
		&+\sum\limits_{j\neq i}\Xi_{i}v_j\tilde{r}_j\cdot D(B-\mu_iI_n)[v_i\xi_i\pa_x(\pa_{v_i}\bar{v}_i)\tilde{r}_{i,v}+v_i\xi_i\pa_{v_i}\bar{v}_i\pa_x(\tilde{r}_{i,v})]\\
		&=O(1)\sum \La_j^1.
		\end{align*}
		Next we have using \eqref{6.45}, \eqref{ngtech4}, \eqref{chichi1}, \eqref{chicle}  $\al_{i}^{5,4}=\Xi_i v_{i}(B-\mu_iI_n)[\tilde{r}_{i,u}\tilde{r}_i+v_i\xi_i\pa_{v_i}\bar{v}_i\tilde{r}_{i,uv}\tilde{r}_{i}]=O(1)\sum_j(\La_j^1+\La_j^4+\La_j^6)+R_\e^{5,i,4}$ satisfying \eqref{12.11}.  Using again Lemmas \ref{estimimpo1}, \ref{lemme6.5}, \ref{lemme6.6}, \ref{lemme11.3}, \ref{lemme11.4} and \eqref{ngtech4}, \eqref{6.45}, \eqref{estimate-v-i-xx-1aa}, \eqref{chichi1}, \eqref{chicle} we get
		\begin{align*}
		\al_{i,x}^{5,4}&=[\big( (\mu_i^{-1}\mathcal{H}_i)_{x}+(\mu_i^{-1}(w_i-\theta_i v_i))_x\big)v_i+\Xi_iv_{i,x}](B-\mu_iI_n)[\tilde{r}_{i,u}\tilde{r}_i+v_i\xi_i\pa_{v_i}\bar{v}_i\tilde{r}_{i,uv}\tilde{r}_{i}]\\
		&+\Xi_iv_{i}u_x\cdot D(B-\mu_iI_n)[\tilde{r}_{i,u}\tilde{r}_i+v_i\xi_i\pa_{v_i}\bar{v}_i\tilde{r}_{i,uv}\tilde{r}_{i}]\\
		&+\Xi_i v_{i}(B-\mu_iI_n)[\pa_x(\tilde{r}_{i,u}\tilde{r}_i)+v_{i,x}\xi_i\pa_{v_i}\bar{v}_i\tilde{r}_{i,uv}\tilde{r}_{i}+v_i\xi_i^\p\left(\frac{w_i}{v_i}\right)_x\pa_{v_i}\bar{v}_i\tilde{r}_{i,uv}\tilde{r}_{i}]\\
		&+\Xi_i v_{i}(B-\mu_iI_n)[v_i\xi_i\pa_x(\pa_{v_i}\bar{v}_i)\tilde{r}_{i,uv}\tilde{r}_{i}+v_i\xi_i\pa_{v_i}\bar{v}_i\pa_x(\tilde{r}_{i,uv}\tilde{r}_{i})]\\
		&=O(1)\sum_j(\La_j^1+\La_j^4+\La_j^5+\La_j^6+\La_j^{6,1})+\widetilde{R}_\e^{5,i,4},
		\end{align*}
		with $\widetilde{R}_\e^{5,i,4}$ satisfying \eqref{12.11}.
		Newt we have $\al_{i}^{5,5}=\sum\limits_{j\neq i}\Xi_iv_{j}(B-\mu_iI_n)[\tilde{r}_{i,u}\tilde{r}_j+v_i\xi_i\pa_{v_i}\bar{v}_i\tilde{r}_{i,uv}\tilde{r}_{j}]=O(1)\sum_{j}\La_j^1.$ Furthermore we have from Lemma \ref{lemme11.4}
		\begin{align*}
		\al_{i,x}^{5,5}&=\sum\limits_{j\neq i}[\Xi_{i,x}v_{j}+\Xi_iv_{j,x}](B-\mu_iI_n)[\tilde{r}_{i,u}\tilde{r}_j+v_i\xi_i\pa_{v_i}\bar{v}_i\tilde{r}_{i,uv}\tilde{r}_{j}]\\
		&+\sum\limits_{j\neq i}\Xi_i v_j u_x\cdot D(B-\mu_iI_n)[\tilde{r}_{i,u}\tilde{r}_j+v_i\xi_i\pa_{v_i}\bar{v}_i\tilde{r}_{i,uv}\tilde{r}_{j}]\\
		&+\sum\limits_{j\neq i}\Xi_iv_j(B-\mu_iI_n)[\pa_x(\tilde{r}_{i,u}\tilde{r}_j)+v_{i,x}\xi_i\pa_{v_i}\bar{v}_i\tilde{r}_{i,uv}\tilde{r}_{j}+v_i\xi_i^\p\left(\frac{w_i}{v_i}\right)_x\pa_{v_i}\bar{v}_i\tilde{r}_{i,uv}\tilde{r}_{j}]\\
		&+\sum\limits_{j\neq i}\Xi_iv_j(B-\mu_iI_n)[ v_i\xi_i\pa_x(\pa_{v_i}\bar{v}_i)\tilde{r}_{i,uv}\tilde{r}_{j}+v_i\xi_i\pa_{v_i}\bar{v}_i\pa_x(\tilde{r}_{i,uv}\tilde{r}_{j})]\\
		&=O(1)\sum_j\La_j^1.
		\end{align*}
		We have now using the fact that $\mbox{supp}\xi$ is included in $\{s,\theta(s)=s\}$  $\al_{i}^{5,6}=\mu_i^{-1}\mathcal{H}_i v_{i,x}\xi_i(B-\mu_iI_n)[(2\pa_{v_i}\bar{v}_i+v_i\pa_{v_iv_i}\bar{v}_i)\tilde{r}_{i,v}+\xi_iv_i(\pa_{v_i}\bar{v}_i)^2\tilde{r}_{i,vv}]
		$. From the Lemmas \ref{estimimpo1}, \ref{lemme6.5}, \ref{lemme6.6} and \eqref{ngtech4}we deduce that
		$\al_{i}^{5,6}=O(1)\sum_j(\La_j^1+\La_j^4+\La_j^6)+R_\e^{5,i,6}$ satisfying \eqref{12.11}. In addition we get using again  Lemmas  \ref{estimimpo1}, \ref{lemme6.5}, \ref{lemme6.6}, \ref{lemme11.3} and \eqref{ngtech4}, \eqref{ngtech5}, \eqref{estimate-v-i-xx-1aa}
		\begin{align*}
		\al_{i,x}^{5,6}&=[(\mu_i^{-1}\mathcal{H}_i)_xv_{i,x}\xi_i+\mu_i^{-1}\mathcal{H}_iv_{i,xx}\xi_i](B-\mu_iI_n)[(2\pa_{v_i}\bar{v}_i+v_i\pa_{v_iv_i}\bar{v}_i)\tilde{r}_{i,v}+\xi_iv_i(\pa_{v_i}\bar{v}_i)^2\tilde{r}_{i,vv}]\\
		&+\mu_i^{-1}\mathcal{H}_iv_{i,x}\xi^\p_i\left(\frac{w_i}{v_i}\right)_x(B-\mu_iI_n)[(2\pa_{v_i}\bar{v}_i+v_i\pa_{v_iv_i}\bar{v}_i)\tilde{r}_{i,v}+\xi_iv_i(\pa_{v_i}\bar{v}_i)^2\tilde{r}_{i,vv}]\\
		&+\mu_i^{-1}\mathcal{H}_iv_{i,x}\xi_i u_x\cdot D(B-\mu_iI_n)[(2\pa_{v_i}\bar{v}_i+v_i\pa_{v_iv_i}\bar{v}_i)\tilde{r}_{i,v}+\xi_i v_i(\pa_{v_i}\bar{v}_i)^2\tilde{r}_{i,vv}]\\
		&+\mu_i^{-1}\mathcal{H}_iv_{i,x}\xi_i(B-\mu_iI_n)[(2\pa_x(\pa_{v_i}\bar{v}_i)+v_{i,x}\pa_{v_iv_i}\bar{v}_i+v_i\pa_x(\pa_{v_iv_i}\bar{v}_i))\tilde{r}_{i,v}]\\
		&+\mu_i^{-1}\mathcal{H}_iv_{i,x}\xi_i(B-\mu_iI_n)[(2\pa_{v_i}\bar{v}_i+v_i\pa_{v_iv_i}\bar{v}_i)\pa_x(\tilde{r}_{i,v})+(\xi_i^\p \left(\frac{w_i}{v_i}\right)_x +\xi_iv_{i,x})(\pa_{v_i}\bar{v}_i)^2\tilde{r}_{i,vv}]\\
		&+\mu_i^{-1}\mathcal{H}_iv_{i,x}\xi_i(B-\mu_iI_n)[2 \xi_iv_i(\pa_{v_i}\bar{v}_i)\pa_x(\pa_{v_i}\bar{v}_i)\tilde{r}_{i,vv}+ \xi_i v_i(\pa_{v_i}\bar{v}_i)^2\pa_x(\tilde{r}_{i,vv})]\\
		&=O(1)\sum(\La_j^1+\La_j^4+\La_j^5+\La_j^6+\delta_0^2\La_j^3)+\widetilde{R}_\e^{5,i,6},
		\end{align*}
		with $\widetilde{R}_\e^{5,i,6}$ satisfying \eqref{12.11}.
		We mention that in the previous term, we have used the fact that $\mathfrak{A}_i\rho_i^\e \frac{|v_{i,x}|^2}{|v_i|}=O(1)$ due to \eqref{ngtech5}. We have now $\al_{i}^{5,7}=\sum\limits_{j\neq i}\mu_i^{-1}\mathcal{H}_i v_{j,x}\xi_i(B-\mu_iI_n)[(\pa_{v_j}\bar{v}_i+v_i\pa_{v_iv_j}\bar{v}_i)\tilde{r}_{i,v}+ \xi_iv_i\pa_{v_i}\bar{v}_i\pa_{v_j}\bar{v}_i\tilde{r}_{i,vv}]=
		O(1)\sum_j\La_j^1$. Furthermore we have
		from  Lemmas  \ref{estimimpo1}, \ref{lemme6.5}, \ref{lemme6.6}, \ref{lemme11.3} and \eqref{ngtech5}% \ref{lemme9.7}, \ref{lemma:derivative-r-k},  \ref{lemme10.5}, \ref{lemme6.8}
		\begin{align*}
		\al_{i,x}^{5,7}
		&=\sum\limits_{j\neq i}[(\mu_i^{-1}\mathcal{H}_i)_xv_{j,x}+\mu_i^{-1}\mathcal{H}_iv_{j,xx}]\xi_i(B-\mu_iI_n)[(\pa_{v_j}\bar{v}_i+v_i\pa_{v_iv_j}\bar{v}_i)\tilde{r}_{i,v}+ \xi_iv_i\pa_{v_i}\bar{v}_i\pa_{v_j}\bar{v}_i\tilde{r}_{i,vv}]\\
		&+\sum\limits_{j\neq i}\mu_i^{-1}\mathcal{H}_i v_{j,x}\xi_i^\p\left(\frac{w_i}{v_i}\right)_x(B-\mu_iI_n)[(\pa_{v_j}\bar{v}_i+v_i\pa_{v_iv_j}\bar{v}_i)\tilde{r}_{i,v}+\xi_iv_i\pa_{v_i}\bar{v}_i\pa_{v_j}\bar{v}_i\tilde{r}_{i,vv}]\\
		&+\sum\limits_{j\neq i}\mu_i^{-1}\mathcal{H}_iv_{j,x}\xi_iu_x\cdot D(B-\mu_iI_n)[(\pa_{v_j}\bar{v}_i+v_i\pa_{v_iv_j}\bar{v}_i)\tilde{r}_{i,v}+\xi_i v_i\pa_{v_i}\bar{v}_i\pa_{v_j}\bar{v}_i\tilde{r}_{i,vv}]\\
		&+\sum\limits_{j\neq i}\mu_i^{-1}\mathcal{H}_iv_{j,x}\xi_i(B-\mu_iI_n)[(\pa_x(\pa_{v_j}\bar{v}_i)+v_{i,x}\pa_{v_iv_j}\bar{v}_i+v_i\pa_x(\pa_{v_iv_j}\bar{v}_i))\tilde{r}_{i,v}]\\
		&+\sum\limits_{j\neq i}\mu_i^{-1}\mathcal{H}_iv_{j,x}\xi_i(B-\mu_iI_n)[(\pa_{v_j}\bar{v}_i+v_i\pa_{v_iv_j}\bar{v}_i)\pa_x(\tilde{r}_{i,v})+(\xi_i^\p\left(\frac{w_i}{v_i}\right)_x v_i+\xi_i v_{i,x})\pa_{v_i}\bar{v}_i\pa_{v_j}\bar{v}_i\tilde{r}_{i,vv}]\\
		&+\sum\limits_{j\neq i}\mu_i^{-1}\mathcal{H}_iv_{j,x}\xi_i(B-\mu_iI_n)[ \xi_iv_i(\pa_x(\pa_{v_i}\bar{v}_i)\pa_{v_j}\bar{v}_i+\pa_{v_i}\bar{v}_i\pa_x(\pa_{v_j}\bar{v}_i))\tilde{r}_{i,vv}+\xi_iv_i\pa_{v_i}\bar{v}_i\pa_{v_j}\bar{v}_i\pa_x(\tilde{r}_{i,vv})]\\
		&=O(1)\sum_j\La_j^1.%+\delta_0^2\La_j^3).
		\end{align*}
	%	We have now $\al_{i}^{6,8}=\sum\limits_{j\neq i}\Xi_{i,1} w_{j,x}\xi_i(B-\mu_iI_n)[(\pa_{w_j}\bar{v}_i+v_i\pa_{v_iw_j}\bar{v}_i)\tilde{r}_{i,v}+\xi_iv_i\pa_{v_i}\bar{v}_i\pa_{w_j}\bar{v}_i\tilde{r}_{i,vv}]=O(1)\sum_j\La_j^1$ and similarly 
	%	\begin{align*}
	%\al_{i,x}^{6,8}&=\sum\limits_{j\neq i}[\Xi_{i1,x}w_{j,x}+\Xi_{i,1}w_{j,xx}]\xi_i(B-\mu_iI_n)[(\pa_{w_j}\bar{v}_i+v_i\pa_{v_iw_j}\bar{v}_i)\tilde{r}_{i,v}+\xi_iv_i\pa_{v_i}\bar{v}_i\pa_{w_j}\bar{v}_i\tilde{r}_{i,vv}]\\
	%	&+\sum\limits_{j\neq i}\Xi_{i,1}w_{j,x}\xi_i^\p\left(\frac{w_i}{v_i}\right)_x(B-\mu_iI_n)[(\pa_{w_j}\bar{v}_i+v_i\pa_{v_iw_j}\bar{v}_i)\tilde{r}_{i,v}+ \xi_iv_i\pa_{v_i}\bar{v}_i\pa_{w_j}\bar{v}_i\tilde{r}_{i,vv}]\\
	%	&+\sum\limits_{j\neq i}\Xi_{i,1}w_{j,x}\xi_iu_x\cdot D(B-\mu_iI_n)[(\pa_{w_j}\bar{v}_i+v_i\pa_{v_iw_j}\bar{v}_i)\tilde{r}_{i,v}+\xi_i v_i\pa_{v_i}\bar{v}_i\pa_{w_j}\bar{v}_i\tilde{r}_{i,vv}]\\
	%	&+\sum\limits_{j\neq i}\Xi_{i,1}w_{j,x}\xi_i(B-\mu_iI_n)[(\pa_x(\pa_{w_j}\bar{v}_i)+v_{i,x}\pa_{v_iw_j}\bar{v}_i+v_i\pa_x(\pa_{v_iw_j}\bar{v}_i))\tilde{r}_{i,v}]\\
	%	&+\sum\limits_{j\neq i}\Xi_{i,1}w_{j,x}\xi_i(B-\mu_iI_n)[(\pa_{w_j}\bar{v}_i+v_i\pa_{v_iw_j}\bar{v}_i)\pa_x(\tilde{r}_{i,v})+ (\xi_i^\p\left(\frac{w_i}{v_i}\right)_x v_i+\xi_i v_{i,x})\pa_{v_i}\bar{v}_i\pa_{w_j}\bar{v}_i\tilde{r}_{i,vv}]\\
	%	&+\sum\limits_{j\neq i}\Xi_{i,1}w_{j,x}\xi_i(B-\mu_iI_n)[ \xi_iv_i(\pa_x(\pa_{v_i}\bar{v}_i)\pa_{w_j}\bar{v}_i+\pa_{v_i}\bar{v}_i\pa_x(\pa_{w_j}\bar{v}_i))\tilde{r}_{i,vv}+\xi_iv_i\pa_{v_i}\bar{v}_i\pa_{w_j}\bar{v}_i\pa_x(\tilde{r}_{i,vv})]\\
	%	&=O(1)\sum_j(\La_j^1+\delta_0^2\La_j^3).
	%	\end{align*}
Next we have using Lemma \ref{estimimpo1} and \eqref{ngtech5} $\al_{i}^{5,8}=\mu_i^{-1}\mathcal{H}_i\xi^\p_i\left(\frac{w_i}{v_i}\right)_x(B-\mu_iI_n)[(v_i\pa_{v_i}\bar{v}_i+\bar{v}_i)\tilde{r}_{i,v}+\xi_iv_i\pa_{v_i}\bar{v}_i\bar{v}_i\tilde{r}_{i,vv}]=O(1)\sum_j(\La_j^1+\La_j^4).$ From  Lemmas  \ref{estimimpo1}, \ref{lemme6.5}, \ref{lemme6.6}, \ref{lemme11.3} and \eqref{identity:wi-vi-xx}, \eqref{ngtech5}, \eqref{estimate-v-i-xx-1aa}% \ref{lemme9.7}, \ref{lemma:derivative-r-k},  \ref{lemme10.5}, \ref{lemme6.8} 
we deduce
		\begin{align*}
		\al_{i,x}^{5,8}&=(\mu_i^{-1}\mathcal{H}_i)_x\xi^\p_i\left(\frac{w_i}{v_i}\right)_x(B-\mu_iI_n)[(v_i\pa_{v_i}\bar{v}_i+\bar{v}_i)\tilde{r}_{i,v}+\xi_iv_i\pa_{v_i}\bar{v}_i\bar{v}_i\tilde{r}_{i,vv}]\\
		&+\mu_i^{-1}\mathcal{H}_i\xi^{\p\p}_i\left(\frac{w_i}{v_i}\right)^2_x(B-\mu_iI_n)[(v_i\pa_{v_i}\bar{v}_i+\bar{v}_i)\tilde{r}_{i,v}+\xi_iv_i\pa_{v_i}\bar{v}_i\bar{v}_i\tilde{r}_{i,vv}]\\
		&+\mu_i^{-1}\mathcal{H}_i\xi^\p_i\left(\frac{w_i}{v_i}\right)_{xx}(B-\mu_iI_n)[(v_i\pa_{v_i}\bar{v}_i+\bar{v}_i)\tilde{r}_{i,v}+\xi_iv_i\pa_{v_i}\bar{v}_i\bar{v}_i\tilde{r}_{i,vv}]\\
		&+\mu_i^{-1}\mathcal{H}_i\xi^\p_i\left(\frac{w_i}{v_i}\right)_xu_x\cdot D(B-\mu_iI_n)[(v_i\pa_{v_i}\bar{v}_i+\bar{v}_i)\tilde{r}_{i,v}+\xi_iv_i\pa_{v_i}\bar{v}_i\bar{v}_i\tilde{r}_{i,vv}]\\
		&+\mu_i^{-1}\mathcal{H}_i\xi^\p_i\left(\frac{w_i}{v_i}\right)_x(B-\mu_iI_n)[(v_{i,x}\pa_{v_i}\bar{v}_i+v_i\pa_x(\pa_{v_i}\bar{v}_i)+\pa_x\bar{v}_i)\tilde{r}_{i,v}]\\
		&+\mu_i^{-1}\mathcal{H}_i\xi^\p_i\left(\frac{w_i}{v_i}\right)_x(B-\mu_iI_n)[(v_i\pa_{v_i}\bar{v}_i+\bar{v}_i)\pa_x(\tilde{r}_{i,v})]\\
		&+\mu_i^{-1}\mathcal{H}_i\xi^\p_i\left(\frac{w_i}{v_i}\right)_x(B-\mu_iI_n)[\xi^\p_iv_i\left(\frac{w_i}{v_i}\right)_x\pa_{v_i}\bar{v}_i\bar{v}_i+\xi_iv_{i,x}\pa_{v_i}\bar{v}_i\bar{v}_i+\xi_iv_i\pa_x(\pa_{v_i}\bar{v}_i)\bar{v}_i]\tilde{r}_{i,vv}\\
		&+\mu_i^{-1}\mathcal{H}_i\xi^\p_i\left(\frac{w_i}{v_i}\right)_x(B-\mu_iI_n)[\xi_iv_i\pa_{v_i}\bar{v}_i\pa_x(\bar{v}_i)\tilde{r}_{i,vv}+\xi_iv_i\pa_{v_i}\bar{v}_i\bar{v}_i\pa_x(\tilde{r}_{i,vv})]\\
		&=O(1)\sum_j (\La_j^1+\La_j^4+\La_j^5+\La_j^6+\delta_0^2\La_j^3).
		\end{align*}
		We have in particular used the fact that $\xi'_i\mu_i^{-1}\mathcal{H}_i \rho_i^\e v_i \left(\frac{w_i}{v_i}\right)_{xx}=O(1)\sum(\La_j^1+\La_j^4+\La_j^5)$ from \eqref{ngtech5}.
		Finally  using  the fact that $\tilde{r}_{i,\si}=O(1)\xi_i \bar{v}_i$ we have $\al_{i}^{5,9}= -
		\mu_i^{-1}\mathcal{H}_i\theta_i^\p\left(\frac{w_i}{v_i}\right)_x(B-\mu_iI_n)[\tilde{r}_{i,\si}+\xi_iv_i\pa_{v_i}\bar{v}_i\tilde{r}_{i,v\si}]$. Applying Lemma \ref{estimimpo1} and \eqref{ngtech5} we deduce that $\al_{i}^{5,9}= O(1)\sum_j(\La_j^1+\La_j^4)$. In addition we get as previously using in particular Lemma \ref{lemme11.3}
		\begin{align*}
		\al_{i,x}^{5,9}&=-(\mu_i^{-1}\mathcal{H}_i)_x\theta_i^\p\left(\frac{w_i}{v_i}\right)_x(B-\mu_iI_n)[\tilde{r}_{i,\si}+\xi_iv_i\pa_{v_i}\bar{v}_i\tilde{r}_{i,v\si}]\\
		&-\mu_i^{-1}\mathcal{H}_i[\theta_i^{\p\p}\left(\frac{w_i}{v_i}\right)^2_x+\theta_i^\p\left(\frac{w_i}{v_i}\right)_{xx}]
		(B-\mu_iI_n)[\tilde{r}_{i,\si}+\xi_iv_i\pa_{v_i}\bar{v}_i\tilde{r}_{i,v\si}]\\
		&-\mu_i^{-1}\mathcal{H}_i\theta_i^\p\left(\frac{w_i}{v_i}\right)_xu_x\cdot D(B-\mu_iI_n)[\tilde{r}_{i,\si}+\xi_iv_i\pa_{v_i}\bar{v}_i\tilde{r}_{i,v\si}]\\
		&-\mu_i^{-1}\mathcal{H}_i\theta_i^\p\left(\frac{w_i}{v_i}\right)_x(B-\mu_iI_n)[\pa_x(\tilde{r}_{i,\si})+\xi_i^\p\left(\frac{w_i}{v_i}\right)_xv_i\pa_{v_i}\bar{v}_i\tilde{r}_{i,v\si}+\xi_iv_{i,x}\pa_{v_i}\bar{v}_i\tilde{r}_{i,v\si}]\\
		&-\mu_i^{-1}\mathcal{H}_i\theta_i^\p\left(\frac{w_i}{v_i}\right)_x(B-\mu_iI_n)[\xi_iv_i\pa_x(\pa_{v_i}\bar{v}_i)\tilde{r}_{i,v\si}+\xi_iv_i\pa_{v_i}\bar{v}_i\pa_x(\tilde{r}_{i,v\si})]\\
		&=O(1)\sum_j (\La_j^1+\La_j^4+\La_j^5+\La_j^6+\delta_0^2\La_j^3).
	\end{align*}

\end{proof}

\begin{lemma}
\label{lemmeal7}
	For $2\leq l\leq 12$ it follows
	\begin{align}
			&\al_{i}^{6,l}=O(1)\sum_j(\La_j^1+\La_j^2+\delta_0^2\La_j^3+\La_j^4+\La_j^5+\La_j^6+\La_j^{6,1}),%+R_\e^{6,i,l},
			\nonumber\\
			&\al_{i,x}^{6,l}=O(1)\sum_j(\La_j^1+\La_j^2+\delta_0^2\La_j^3+\La_j^4+\La_j^5+\La_j^6+\La_j^{6,1}).%+\widetilde{R}_\e^{6,i,l},
	\end{align}
	%\begin{align}
	%		&\al_{i}^{5,l}=O(1)\sum_j(\La_j^1+\La_j^2+\delta_0^2\La_j^3+\La_j^4+\La_j^5+\La_j^6+\La_j^{6,1})+R_\e^{5,i,l},\nonumber\\
	%		&\al_{i,x}^{5,l}=O(1)\sum_j(\La_j^1+\La_j^2+\delta_0^2\La_j^3+\La_j^4+\La_j^5+\La_j^6+\La_j^{6,1})+\widetilde{R}_\e^{5,i,l},\label{12.11}
	%\end{align}
	%with:
	%$$\int_{\hat{t}}^T\int_{\R}(|R_\e^{6,i,l}|+|\widetilde{R}_\e^{6,i,l}|)dx ds=O(1)\delta_0^2,$$
	%for $\e>0$ small enough in terms of $T-\hat{t}$ and $\delta_0$.
\end{lemma}
\begin{proof}
	We observe that $\al_{i,x}^{6,2}=v_i\left(\frac{w_i}{v_i}\right)_x^2(B-\mu_iI_n)[\xi^{\p\p}_i\bar{v}_i\tilde{r}_{i,v}-\theta^{\p\p}_i\tilde{r}_{i,\si}]=\delta_0^2\La_i^3$ due to Lemma \ref{lemme6.5} and the fact that $\tilde{r}_{i,\si}=O(1)\xi_i\bar{v}_i$. In addition we have using Lemmas  \ref{estimimpo1}, \ref{lemme6.5}, \ref{lemme6.6}, \ref{lemme11.3} and the fact that $\tilde{r}_{i,\sig}=O(1)\xi_i v_i$
	\begin{align*}
		\al_{i,x}^{6,2}&=[v_{i,x} \left(\frac{w_i}{v_i}\right)_x^2+2 v_i \left(\frac{w_i}{v_i}\right)_x \left(\frac{w_i}{v_i}\right)_{xx}]
		(B-\mu_iI_n)[\xi^{\p\p}_i\bar{v}_i\tilde{r}_{i,v}-\theta^{\p\p}_i\tilde{r}_{i,\si}]\\%+2 v_i \left(\frac{w_i}{v_i}\right)_x \left(\frac{w_i}{v_i}\right)_{xx}(B-\mu_iI_n)[\xi^{\p\p}_i\bar{v}_i\tilde{r}_{i,v}-\theta^{\p\p}_i\tilde{r}_{i,\si}]\\
		&+ v_i\left(\frac{w_i}{v_i}\right)_x^2 u_x\cdot D(B-\mu_iI_n)[\xi^{\p\p}_i\bar{v}_i\tilde{r}_{i,v}-\theta^{\p\p}_i\tilde{r}_{i,\si}]+v_i\left(\frac{w_i}{v_i}\right)_x^3(B-\mu_iI_n)[\xi^{\p\p\p}_i\bar{v}_i\tilde{r}_{i,v}-\theta^{\p\p\p}_i\tilde{r}_{i,\si}]\\
		&+v_i\left(\frac{w_i}{v_i}\right)_x^2(B-\mu_iI_n)[\xi^{\p\p}_i\pa_x(\bar{v}_i)\tilde{r}_{i,v}+\xi^{\p\p}_i\bar{v}_i\pa_x(\tilde{r}_{i,v})-\theta^{\p\p}_i\pa_x(\tilde{r}_{i,\si})]\\
		&=O(1)(\delta_0^2\La_i^3+\La_i^5).
		\end{align*}
		Newt we have $\al_{i}^{6,3}=(w_{i,x}v_i-w_iv_{i,x})\tilde{r}_i\cdot D(B-\mu_iI_n)[\xi_i^\p\bar{v}_i\tilde{r}_{i,v}-\theta^\p_i\tilde{r}_{i,\si}]=O(1)\La_i^4$ and similarly from the Lemma \ref{lemme11.3}
		\begin{align*}
		\al_{i,x}^{6,3}&=\big[(w_{i,xx}v_i-w_iv_{i,xx})\tilde{r}_i\cdot D(B-\mu_iI_n)+(w_{i,x}v_i-w_iv_{i,x})\tilde{r}_{i,x}\cdot D(B-\mu_iI_n)\big][\xi_i^\p\bar{v}_i\tilde{r}_{i,v}-\theta^\p_i\tilde{r}_{i,\si}]\\
		&+(w_{i,x}v_i-w_iv_{i,x})\tilde{r}_i\otimes u_x:D^2(B-\mu_iI_n)
		[\xi_i^\p\bar{v}_i\tilde{r}_{i,v}-\theta^\p_i\tilde{r}_{i,\si}]\\
		&+(w_{i,x}v_i-w_iv_{i,x})\left(\frac{w_i}{v_i}\right)_x\tilde{r}_i\cdot D(B-\mu_iI_n)[\xi_i^{\p\p}\bar{v}_i\tilde{r}_{i,v}-\theta^{\p\p}_i\tilde{r}_{i,\si}]\\
		&+(w_{i,x}v_i-w_iv_{i,x})\tilde{r}_i\cdot D(B-\mu_iI_n)[\xi_i^\p\pa_x(\bar{v}_i)\tilde{r}_{i,v}+\xi_i^\p\bar{v}_i\pa_x(\tilde{r}_{i,v})-\theta^\p_i\pa_x(\tilde{r}_{i,\si})]\\
        &=O(1)(\La_i^4+\La_i^5).
        \end{align*}
        We have now $ \al_{i}^{6,4}=\sum\limits_{j\neq i}(w_{i,x}-(w_i/v_i)v_{i,x})v_j\tilde{r}_j\cdot D(B-\mu_iI_n)[\xi_i^\p\bar{v}_i\tilde{r}_{i,v}-\theta^\p_i\tilde{r}_{i,\si}]=O(1)\La_i^1$ and using the fact that 
        $\tilde{r}_{i,\si}=O(1)\xi_i \bar{v_i}$ we get from the Lemma \ref{lemme11.3}
        \begin{align*}
        \al_{i,x}^{6,4}&=\sum\limits_{j\neq i}(w_{i,x}-(w_i/v_i)v_{i,x})_xv_j\tilde{r}_j\cdot D(B-\mu_iI_n)[\xi_i^\p\bar{v}_i\tilde{r}_{i,v}-\theta^\p_i\tilde{r}_{i,\si}]\\
        &+\sum\limits_{j\neq i}(w_{i,x}-(w_i/v_i)v_{i,x})v_{j,x}\tilde{r}_j\cdot D(B-\mu_iI_n)[\xi_i^\p\bar{v}_i\tilde{r}_{i,v}-\theta^\p_i\tilde{r}_{i,\si}]\\
        &+\sum\limits_{j\neq i}(w_{i,x}-(w_i/v_i)v_{i,x})v_j(\tilde{r}_j\otimes u_x): D^2(B-\mu_iI_n)[\xi_i^\p\bar{v}_i\tilde{r}_{i,v}-\theta^\p_i\tilde{r}_{i,\si}]\\
        &+\sum\limits_{j\neq i}(w_{i,x}-(w_i/v_i)v_{i,x})v_j\left(\frac{w_i}{v_i}\right)_x\tilde{r}_j\cdot D(B-\mu_iI_n)[\xi_i^{\p\p}\bar{v}_i\tilde{r}_{i,v}-\theta^{\p\p}_i\tilde{r}_{i,\si}]\\
        &+\sum\limits_{j\neq i}(w_{i,x}-(w_i/v_i)v_{i,x})v_j\tilde{r}_j\cdot D(B-\mu_iI_n)[\xi_i^\p\pa_x(\bar{v}_i)\tilde{r}_{i,v}+\xi_i^\p\bar{v}_i\pa_x(\tilde{r}_{i,v})-\theta^\p_i\pa_x(\tilde{r}_{i,\si})]\\
        &=O(1)\La_i^1.
        \end{align*}
       We have now $\al_{i}^{6,5}= (w_{i,x}v_i-w_iv_{i,x}) (B-\mu_iI_n)[\xi_i^\p\bar{v}_i\tilde{r}_{i,uv}\tilde{r}_i-\theta^\p_i\tilde{r}_{i,u\si}\tilde{r}_i]=(w_{i,x}v_i-w_iv_{i,x}) (B-\mu_iI_n)[\xi_i^\p\bar{v}_i\tilde{r}_{i,uv}\tilde{r}_i-\theta^\p_i\tilde{r}_{i,u\si}\tilde{r}_i]=O(1)\La_i^4$ and using Lemmas  \ref{estimimpo1}, \ref{lemme11.4}
        \begin{align*}
        \al_{i,x}^{6,5}&=(w_{i,xx}v_i-w_iv_{i,xx}) (B-\mu_iI_n)[\xi_i^\p\bar{v}_i\tilde{r}_{i,uv}\tilde{r}_i-\theta^\p_i\tilde{r}_{i,u\si}\tilde{r}_i]\\
        &+(w_{i,x}v_i-w_iv_{i,x}) u_x\cdot D(B-\mu_iI_n)[\xi_i^\p\bar{v}_i\tilde{r}_{i,uv}\tilde{r}_i-\theta^\p_i\tilde{r}_{i,u\si}\tilde{r}_i]\\
        &+(w_{i,x}v_i-w_iv_{i,x})\left(\frac{w_i}{v_i}\right)_x (B-\mu_iI_n)[\xi_i^{\p\p}\bar{v}_i\tilde{r}_{i,uv}\tilde{r}_i-\theta^{\p\p}_i\tilde{r}_{i,u\si}\tilde{r}_i]\\
        &+(w_{i,x}v_i-w_iv_{i,x}) (B-\mu_iI_n)[\xi_i^\p\pa_x(\bar{v}_i)\tilde{r}_{i,uv}\tilde{r}_i+\xi_i^\p\bar{v}_i\pa_x(\tilde{r}_{i,uv}\tilde{r}_i)-\theta^\p_i\pa_x(\tilde{r}_{i,u\si}\tilde{r}_i)]\\
        &=O(1)(\La_i^4+\delta_0^2\La_i^3+\La_i^5).
        \end{align*}
        In a similar manner we get $ \al_{i}^{6,6}=\sum\limits_{j\neq i}(w_{i,x}-(w_i/v_i)v_{i,x})v_j(B-\mu_iI_n)[\xi_i^\p\bar{v}_i\tilde{r}_{i,uv}\tilde{r}_j-\theta^\p_i\tilde{r}_{i,u\si}\tilde{r}_j]=O(1)\La_i^1$ and since $\tilde{r}_{i,u\si}=\xi_i \bar{v}_i$ we have using Lemma \ref{lemme11.4}
        \begin{align*}
        \al_{i,x}^{6,6}&=\sum\limits_{j\neq i}(w_{i,x}-(w_i/v_i)v_{i,x})_xv_j(B-\mu_iI_n)[\xi_i^\p\bar{v}_i\tilde{r}_{i,uv}\tilde{r}_j-\theta^\p_i\tilde{r}_{i,u\si}\tilde{r}_j]\\
        &+\sum\limits_{j\neq i}(w_{i,x}-(w_i/v_i)v_{i,x})v_{j,x}(B-\mu_iI_n)[\xi_i^\p\bar{v}_i\tilde{r}_{i,uv}\tilde{r}_j-\theta^\p_i\tilde{r}_{i,u\si}\tilde{r}_j]\\
        &+\sum\limits_{j\neq i}(w_{i,x}-(w_i/v_i)v_{i,x})v_ju_x\cdot D(B-\mu_iI_n)[\xi_i^\p\bar{v}_i\tilde{r}_{i,uv}\tilde{r}_j-\theta^\p_i\tilde{r}_{i,u\si}\tilde{r}_j]\\
        &+\sum\limits_{j\neq i}(w_{i,x}-(w_i/v_i)v_{i,x})\left(\frac{w_i}{v_i}\right)_xv_j(B-\mu_iI_n)[\xi_i^{\p\p}\bar{v}_i\tilde{r}_{i,uv}\tilde{r}_j-\theta^{\p\p}_i\tilde{r}_{i,u\si}\tilde{r}_j]\\
        &+\sum\limits_{j\neq i}(w_{i,x}-(w_i/v_i)v_{i,x})v_j(B-\mu_iI_n)[\xi_i^\p\pa_x(\bar{v}_i)\tilde{r}_{i,uv}\tilde{r}_j+\xi_i^\p\bar{v}_i\pa_x(\tilde{r}_{i,uv}\tilde{r}_j)-\theta^\p_i\pa_x(\tilde{r}_{i,u\si}\tilde{r}_j)]\\
        &=O(1)\La_i^1.
        \end{align*}
    We have now using \eqref{ngtech5} $\al_{i}^{6,7}=(w_{i,x}-(w_i/v_i)v_{i,x})v_{i,x}\xi_i\pa_{v_i}\bar{v}_i(B-\mu_iI_n)[\xi_i^\p\bar{v}_i\tilde{r}_{i,vv}-\theta^\p_i\tilde{r}_{i,v\si}]=
    O(1)(\La_i^1+\La_i^4).$
    In addition we have using Lemmas  \ref{estimimpo1}, \ref{lemme6.5},  \ref{lemme6.6}, \ref{lemme11.3} and \eqref{ngtech5}
        \begin{align*}
        \al_{i,x}^{6,7}&=(w_{i,x}-(w_i/v_i)v_{i,x})_xv_{i,x}\xi_i\pa_{v_i}\bar{v}_i(B-\mu_iI_n)[\xi_i^\p\bar{v}_i\tilde{r}_{i,vv}-\theta^\p_i\tilde{r}_{i,v\si}]\\
        &+(w_{i,x}-(w_i/v_i)v_{i,x})v_{i,xx}\xi_i\pa_{v_i}\bar{v}_i(B-\mu_iI_n)[\xi_i^\p\bar{v}_i\tilde{r}_{i,vv}-\theta^\p_i\tilde{r}_{i,v\si}]\\
        &+(w_{i,x}-(w_i/v_i)v_{i,x})v_{i,x}\xi_i^\p\left(\frac{w_i}{v_i}\right)_x\pa_{v_i}\bar{v}_i(B-\mu_iI_n)[\xi_i^\p\bar{v}_i\tilde{r}_{i,vv}-\theta^\p_i\tilde{r}_{i,v\si}]\\
        &+(w_{i,x}-(w_i/v_i)v_{i,x})v_{i,x}\xi_i\pa_x(\pa_{v_i}\bar{v}_i)(B-\mu_iI_n)[\xi_i^\p\bar{v}_i\tilde{r}_{i,vv}-\theta^\p_i\tilde{r}_{i,v\si}]\\
        &+(w_{i,x}-(w_i/v_i)v_{i,x})v_{i,x}\xi_i\pa_{v_i}\bar{v}_iu_x\cdot D(B-\mu_iI_n)[\xi_i^\p\bar{v}_i\tilde{r}_{i,vv}-\theta^\p_i\tilde{r}_{i,v\si}]\\
        &+(w_{i,x}-(w_i/v_i)v_{i,x})\left(\frac{w_i}{v_i}\right)_xv_{i,x}\xi_i\pa_{v_i}\bar{v}_i(B-\mu_iI_n)[\xi_i^{\p\p}\bar{v}_i\tilde{r}_{i,vv}-\theta^{\p\p}_i\tilde{r}_{i,v\si}]\\
        &+(w_{i,x}-(w_i/v_i)v_{i,x})v_{i,x}\xi_i\pa_{v_i}\bar{v}_i(B-\mu_iI_n)[\xi_i^\p\pa_x(\bar{v}_i)\tilde{r}_{i,vv}+\xi_i^\p\bar{v}_i\pa_x(\tilde{r}_{i,vv})-\theta^\p_i\pa_x(\tilde{r}_{i,v\si})]\\
        &=O(1)\sum_j (\La_j^1+\delta_0^2\La_j^3+\La_j^4+\La_j^5+\La_j^6).
        \end{align*}
        We explain here how to deal with the term $(w_{i,x}-(w_i/v_i)v_{i,x})v_{i,xx}\xi_i\pa_{v_i}\bar{v}_i(B-\mu_iI_n)[\xi_i^\p\bar{v}_i\tilde{r}_{i,vv}-\theta^\p_i\tilde{r}_{i,v\si}]$. From \eqref{estimate-v-i-xx-1aabis}, we have
       %  \begin{align}
	%		&\mu_i v_{i,xx}=-\mu_{i,x}v_{i,x}+\tilde{\la}_{i,x}v_i+(\tilde{\la}_i-\la_i^*)v_{i}+w_{i,x} +O(1)\sum\limits_{j\neq i}\left(\mu_j v_{j,x}-(\tilde{\la}_j-\la_j^*)v_j-w_j\right)_x v_j\xi_j
	%+O(1)\sum_{j\ne i}\sum_k (|v_{j,x}v_j v_k|+|v_j^2 v_k|)\\ 
	%\nonumber\\
	%&+O(1)\sum_{j\ne i}(|v_j|+|v_{j,x}|+|w_{j,x}|)(|v_j|+|w_{j,x}|+|v_{j,x}|+\sum_{k} |v_k||v_j|)\nonumber\\
	%&+O(1)\sum_k( \La_k^1+\La_k^4+\La_k^5)
		%&+\mathcal{O}(1)\left(\sum\limits_{j\neq i}|v_j|^2\right)\\
	%	+O(1)\sum_k|v_kv_{k,x}||1-\xi_k\eta_k|.\nonumber\\
	%	\end{align}
		\begin{align*}
		&v_{i,xx}=O(1)(|v_{i,x}|+|w_{i,x}|+|v_i|)+O(1)\sum\limits_{j\neq i}\left(\mu_j v_{j,x}-(\tilde{\la}_j-\la_j^*)v_j-w_j\right)_x \bar{v}_j\xi_j	\nonumber\\
		&+O(1)\sum_{j\ne i}\Big(|v_j|+|v_{j,x}|+|w_{j,x}|\Big)\left(|v_j|+|w_{j,x}|+|v_{j,x}|+\sum_{k} |v_k||v_j|\right)\nonumber\\
		&+O(1)\sum_k(\La_k^1+\La_k^4+\La_k^5+\La_k^6)+
		O(1)\big(\sum_k|v_kv_{k,x}|\eta_k(1-\chi^\e_k)+\sum\limits_{k}|v_k|^3\eta_k (1-\chi^\e_k)\big)%+O(1)\sum_k|v_kv_{k,x}||1-\chi_k^\e\xi_k\eta_k|
		%+O(1)\sum\limits_{k}|v_k|^3(1-{\color{black}{\chi_k^\e}}\xi_k\eta_k)
		\nonumber\\
		&{\color{black}{+O(1)\sum_k \varkappa^\ep_k\rho_k^\e\mathfrak{A}_k|v_k v_{k,x}|}}.
	%\label{estimate-v-i-xx-1aabis}
	\end{align*}
	We deduce then using Lemmas \ref{lemme6.5}, \ref{lemme6.6} and \eqref{ngtech5} that	
	\begin{align*}
	&(w_{i,x}-(w_i/v_i)v_{i,x})v_{i,xx}\xi_i\pa_{v_i}\bar{v}_i(B-\mu_iI_n)[\xi_i^\p\bar{v}_i\tilde{r}_{i,vv}-\theta^\p_i\tilde{r}_{i,v\si}]\\
	&=O(1)|v_i \left(\frac{w_i}{v_i}\right)_x w_{i,x}\mathfrak{A}_i\rho_i^\e|+O(1)|v_i \left(\frac{w_i}{v_i}\right)_x v_{i,x}\mathfrak{A}_i\rho_i^\e|+O(1)\sum_k( \La_k^1+\La_k^4+\La_k^5+\La_k^6)\\
	&=O(1)|v_i^2 \left(\frac{w_i}{v_i}\right)_x^2\mathfrak{A}_i\rho_i^\e|+O(1)|v_i \left(\frac{w_i}{v_i}\right)_x v_{i,x}\mathfrak{A}_i\rho_i^\e|+O(1)\sum_k( \La_k^1+\La_k^4+\La_k^5+\La_k^6)\\%+O(1)v_i	\left(\frac{w_i}{v_i}\right)_xv_{i,x}\xi_i\alpha_i\\
	&=O(1)\sum_k( \La_k^1+\delta_0^2\La_k^3+\La_k^4+\La_k^5+\La_k^6).	%&+\mathcal{O}(1)\left(\sum\limits_{j\neq i}|v_j|^2\right)\\
	\end{align*}
	Similarly we have $ \al_{i}^{6,8}=\sum\limits_{j\neq i}(w_{i,x}-(w_i/v_i)v_{i,x})v_{j,x}\xi_i\pa_{v_j}\bar{v}_i\theta^\p_i(B-\mu_iI_n)[\xi_i^\p\bar{v}_i\tilde{r}_{i,vv}-\theta^\p_i\tilde{r}_{i,v\si}]=O(1)\La_i^1$ and using again Lemmas \ref{estimimpo1}, \ref{lemme6.5}, \ref{lemme6.6}, \ref{lemme11.3} and \eqref{ngtech5}
        \begin{align*}
        \al_{i,x}^{6,8}&=\sum\limits_{j\neq i}(w_{i,x}-(w_i/v_i)v_{i,x})_xv_{j,x}\xi_i\pa_{v_j}\bar{v}_i\theta^\p_i(B-\mu_iI_n)[\xi_i^\p\bar{v}_i\tilde{r}_{i,vv}-\theta^\p_i\tilde{r}_{i,v\si}]\\
        &+\sum\limits_{j\neq i}(w_{i,x}-(w_i/v_i)v_{i,x})v_{j,xx}\xi_i\pa_{v_j}\bar{v}_i\theta^\p_i(B-\mu_iI_n)[\xi_i^\p\bar{v}_i\tilde{r}_{i,vv}-\theta^\p_i\tilde{r}_{i,v\si}]\\
        &+\sum\limits_{j\neq i}(w_{i,x}-(w_i/v_i)v_{i,x})v_{j,x}\xi_i^\p\left(\frac{w_i}{v_i}\right)_x\pa_{v_j}\bar{v}_i\theta^\p_i(B-\mu_iI_n)[\xi_i^\p\bar{v}_i\tilde{r}_{i,vv}-\theta^\p_i\tilde{r}_{i,v\si}]\\
        &+\sum\limits_{j\neq i}(w_{i,x}-(w_i/v_i)v_{i,x})v_{j,x}\xi_i\pa_x(\pa_{v_j}\bar{v}_i)\theta^\p_i(B-\mu_iI_n)[\xi_i^\p\bar{v}_i\tilde{r}_{i,vv}-\theta^\p_i\tilde{r}_{i,v\si}]\\
        &+\sum\limits_{j\neq i}(w_{i,x}-(w_i/v_i)v_{i,x})v_{j,x}\xi_i\pa_{v_j}\bar{v}_i\theta^{\p\p}_i\left(\frac{w_i}{v_i}\right)_x(B-\mu_iI_n)[\xi_i^\p\bar{v}_i\tilde{r}_{i,vv}-\theta^\p_i\tilde{r}_{i,v\si}]\\
        &+\sum\limits_{j\neq i}(w_{i,x}-(w_i/v_i)v_{i,x})v_{j,x}\xi_i\pa_{v_j}\bar{v}_i\theta^\p_iu_x\cdot D(B-\mu_iI_n)[\xi_i^\p\bar{v}_i\tilde{r}_{i,vv}-\theta^\p_i\tilde{r}_{i,v\si}]\\
        &+\sum\limits_{j\neq i}(w_{i,x}-(w_i/v_i)v_{i,x})v_{j,x}\xi_i\pa_{v_j}\bar{v}_i\theta^\p_i\left(\frac{w_i}{v_i}\right)_x(B-\mu_iI_n)[\xi_i^{\p\p}\bar{v}_i\tilde{r}_{i,vv}-\theta^{\p\p}_i\tilde{r}_{i,v\si}]\\
        &+\sum\limits_{j\neq i}(w_{i,x}-(w_i/v_i)v_{i,x})v_{j,x}\xi_i\pa_{v_j}\bar{v}_i\theta^\p_i(B-\mu_iI_n)[\xi_i^\p\pa_x(\bar{v}_i)\tilde{r}_{i,vv}+\xi_i^\p\bar{v}_i\pa_x(\tilde{r}_{i,vv})-\theta^\p_i\pa_x(\tilde{r}_{i,v\si})]\\
        &=O(1)(\La_i^1+\delta_0^2\La_i^3).
        \end{align*}
        Next we have from Lemma \ref{lemme6.5} $\al_{i}^{6,9}=
         v_i\left(\frac{w_i}{v_i}\right)^2_x\xi_i^\p\bar{v}_i(B-\mu_iI_n)[\xi_i^\p\bar{v}_i\tilde{r}_{i,vv}-\theta^\p_i\tilde{r}_{i,v\si}]=\delta_0^2\La_i^3$ and using Lemmas \ref{estimimpo1}, \ref{lemme6.5}, \ref{lemme6.6}, \ref{lemme11.3} and \eqref{identity:wi-vi-xx} it yields
        \begin{align*}
        \al_{i,x}^{6,9}&=\big[v_{i,x}\left(\frac{w_i}{v_i}\right)^2_x+2v_i  \left(\frac{w_i}{v_i}\right)_{x}\left(\frac{w_i}{v_i}\right)_{xx}\big]\xi_i^\p\bar{v}_i(B-\mu_iI_n)[\xi_i^\p\bar{v}_i\tilde{r}_{i,vv}-\theta^\p_i\tilde{r}_{i,v\si}]\\
        %&+2v_i  \left(\frac{w_i}{v_i}\right)_{x}\left(\frac{w_i}{v_i}\right)_{xx}\xi_i^\p\bar{v}_i(B-\mu_iI_n)[\xi_i^\p\bar{v}_i\tilde{r}_{i,vv}-\theta^\p_i\tilde{r}_{i,v\si}]\\
        &+\big[v_i \left(\frac{w_i}{v_i}\right)^3_x\xi_i^{\p\p}\bar{v}_i+v_i\left(\frac{w_i}{v_i}\right)^2_x\xi_i^\p\pa_x(\bar{v}_i)\big]
        (B-\mu_iI_n)[\xi_i^\p\bar{v}_i\tilde{r}_{i,vv}-\theta^\p_i\tilde{r}_{i,v\si}]\\
       % &+v_i\left(\frac{w_i}{v_i}\right)^2_x\xi_i^\p\pa_x(\bar{v}_i)(B-\mu_iI_n)[\xi_i^\p\bar{v}_i\tilde{r}_{i,vv}-\theta^\p_i\tilde{r}_{i,v\si}]\\
        &+v_i\left(\frac{w_i}{v_i}\right)^2_x\xi_i^\p\bar{v}_i u_x\cdot D(B-\mu_iI_n)[\xi_i^\p\bar{v}_i\tilde{r}_{i,vv}-\theta^\p_i\tilde{r}_{i,v\si}]\\
        &+v_i\left(\frac{w_i}{v_i}\right)^3_x\xi_i^\p\bar{v}_i(B-\mu_iI_n)[\xi_i^{\p\p}\bar{v}_i\tilde{r}_{i,vv}-\theta^{\p\p}_i\tilde{r}_{i,v\si}]\\
        &+v_i\left(\frac{w_i}{v_i}\right)^2_x\xi_i^\p\bar{v}_i(B-\mu_iI_n)[\xi_i^\p\pa_x(\bar{v}_i)\tilde{r}_{i,vv}+\xi_i^\p\bar{v}_i\pa_x(\tilde{r}_{i,vv})-\theta^\p_i\pa_x(\tilde{r}_{i,v\si})]\\
        &=O(1)\sum_j(\delta_0^2\La_j^3+\La_5^j).
        \end{align*}
        Next we have $\al_{i}^{6,10}=-v_i\theta^\p_i\left(\frac{w_i}{v_i}\right)^2_x(B-\mu_iI_n)[\xi_i^\p\bar{v}_i\tilde{r}_{i,v\si}-\theta^\p_i\tilde{r}_{i,\si\si}]=O(1)\delta_0^2\La_i^3$ since $\tilde{r}_{i,\si\si}=O(1)\xi_i \bar{v}_i$, in addition we get as previously
        \begin{align*} 
        \al_{i,x}^{6,10}=O(1)\sum_j(\delta_0^2\La_j^3+\La_j^5).
        \end{align*}
        We have now using Lemma \ref{estimimpo1} and \eqref{ngtech5} $\al_{i}^{6,11}=(w_{i,x}-(w_i/v_i)v_{i,x})v_{i,x}(B-\mu_iI_n)\xi_i^\p\pa_{v_i}\bar{v}_i\tilde{r}_{i,v}=O(1)(\La_i^1+\la_i^4).$ In addition using Lemmas \ref{estimimpo1}, \ref{lemme6.5}, \ref{lemme6.6}, \ref{lemme11.3} and \eqref{ngtech5}, \eqref{estimate-v-i-xx-1aabis} we obtain
        \begin{align*}
        \al_{i,x}^{6,11}&=(w_{i,x}-(w_i/v_i)v_{i,x})_xv_{i,x}(B-\mu_iI_n)\xi_i^\p\pa_{v_i}\bar{v}_i\tilde{r}_{i,v}\\
        &+(w_{i,x}-(w_i/v_i)v_{i,x})\big[v_{i,xx}(B-\mu_iI_n)+v_{i,x}u_x\cdot D(B-\mu_iI_n)\big]\xi_i^\p\pa_{v_i}\bar{v}_i\tilde{r}_{i,v}\\
        %&+(w_{i,x}-(w_i/v_i)v_{i,x})v_{i,x}u_x\cdot D(B-\mu_iI_n)\xi_i^\p\pa_{v_i}\bar{v}_i\tilde{r}_{i,v}\\
        &+(w_{i,x}-(w_i/v_i)v_{i,x})v_{i,x}(B-\mu_iI_n)\big[\left(\frac{w_i}{v_i}\right)_x \xi_i^{\p\p}\pa_{v_i}\bar{v}_i\tilde{r}_{i,v}+\xi_i^\p\pa_x(\pa_{v_i}\bar{v}_i)\tilde{r}_{i,v}\big]\\
        %&+(w_{i,x}-(w_i/v_i)v_{i,x})v_{i,x}(B-\mu_iI_n)\xi_i^\p\pa_x(\pa_{v_i}\bar{v}_i)\tilde{r}_{i,v}\\
        &+(w_{i,x}-(w_i/v_i)v_{i,x})v_{i,x}(B-\mu_iI_n)\xi_i^\p\pa_{v_i}\bar{v}_i\pa_x(\tilde{r}_{i,v})\\
        &=O(1)\sum_k( \La_k^1+\delta_0^2\La_k^3+\La_k^4+\La_k^5+\La_k^6).
        \end{align*}
        We have now $\al_{i}^{6,12}=\sum\limits_{j\neq i}(w_{i,x}-(w_i/v_i)v_{i,x})v_{j,x}(B-\mu_iI_n)\xi_i^\p\pa_{v_j}\bar{v}_i\tilde{r}_{i,v}=O(1)\La_i^1$ and using Lemmas \ref{estimimpo1}, \ref{lemme6.5}, \ref{lemme6.6}, \ref{lemme11.3} and \eqref{ngtech5}
        \begin{align*}   
        \al_{i,x}^{6,12}&=\sum\limits_{j\neq i}(w_{i,x}-(w_i/v_i)v_{i,x})_xv_{j,x}(B-\mu_iI_n)\xi_i^\p\pa_{v_j}\bar{v}_i\tilde{r}_{i,v}\\
        &+\sum\limits_{j\neq i}(w_{i,x}-(w_i/v_i)v_{i,x})\big[v_{j,xx}(B-\mu_iI_n)+v_{j,x}u_x\cdot (B-\mu_iI_n)\big]\xi_i^\p\pa_{v_j}\bar{v}_i\tilde{r}_{i,v}\\
       % &+\sum\limits_{j\neq i}(w_{i,x}-(w_i/v_i)v_{i,x})v_{j,x}u_x\cdot (B-\mu_iI_n)\xi_i^\p\pa_{v_j}\bar{v}_i\tilde{r}_{i,v}\\
        &+\sum\limits_{j\neq i}(w_{i,x}-(w_i/v_i)v_{i,x})v_{j,x}(B-\mu_iI_n)\big[\xi_i^{\p\p}\left(\frac{w_i}{v_i}\right)_x\pa_{v_j}\bar{v}_i\tilde{r}_{i,v}+\xi_i^\p\pa_x(\pa_{v_j}\bar{v}_i)\tilde{r}_{i,v}\big]\\
       % &+\sum\limits_{j\neq i}(w_{i,x}-(w_i/v_i)v_{i,x})v_{j,x}(B-\mu_iI_n)\xi_i^\p\pa_x(\pa_{v_j}\bar{v}_i)\tilde{r}_{i,v}\\
        &+\sum\limits_{j\neq i}(w_{i,x}-(w_i/v_i)v_{i,x})v_{j,x}(B-\mu_iI_n)\xi_i^\p\pa_{v_j}\bar{v}_i\pa_x(\tilde{r}_{i,v})\\
        &=O(1)\sum_j(\La_j^1+\delta_0^2\La_j^3).
        \end{align*}
       % We proceed similarly with $ \al_{i}^{7,14}=\sum \limits_{j\neq i}(w_{i,x}-(w_i/v_i)v_{i,x})w_{j,x}(B-\mu_iI_n)\xi_i^\p\pa_{w_j}\bar{v}_i\tilde{r}_{i,v}=O(1)\La_i^1$.
  
\end{proof}
\begin{lemma}
\label{lemmeal8}
	For $1\leq l\leq 17$ we have 
	\begin{align}
			&\al_{i}^{7,l}=O(1)\sum_j(\La_j^1+\La_j^2+\delta_0^2\La_j^3+\La_j^4+\La_j^5+\La_j^6+\La_j^{6,1})+R_\e^{7,i,l},\\
			&\al_{i,x}^{7,l}=O(1)\sum_j(\La_j^1+\La_j^2+\delta_0^2\La_j^3+\La_j^4+\La_j^5+\La_j^6+\La_j^{6,1})+\widetilde{R}_\e^{7,i,l},
	\end{align}
	%\begin{align}
	%		&\al_{i}^{5,l}=O(1)\sum_j(\La_j^1+\La_j^2+\delta_0^2\La_j^3+\La_j^4+\La_j^5+\La_j^6+\La_j^{6,1})+R_\e^{5,i,l},\nonumber\\
	%		&\al_{i,x}^{5,l}=O(1)\sum_j(\La_j^1+\La_j^2+\delta_0^2\La_j^3+\La_j^4+\La_j^5+\La_j^6+\La_j^{6,1})+\widetilde{R}_\e^{5,i,l},\label{12.11}
	%\end{align}
	with:
	\begin{equation}
	\int_{\hat{t}}^T\int_{\R}(|R_\e^{7,i,l}|+|\widetilde{R}_\e^{7,i,l}|)dx ds=O(1)\delta_0^2,
	\label{12.11bis}
	\end{equation}
	for $\e>0$ small enough in terms of $T-\hat{t}$ and $\delta_0$.
\end{lemma}

\begin{proof}
First we have from \eqref{form-of-tilde-v-i}
\begin{align}
\tilde{v}_i=\pa_{v_i}\bar{v}_i v_i-\bar{v}_i=O(1)(\chi')^\e_i\hat{v}_i+O(1) \sum_{l\ne i}\Delta^\e_{i,l}v_l^2.\label{6.49bis}
\end{align}
We deduce that $\al_{i}^{7,1}=\sum\limits_{k}v_k\tilde{r}_k\cdot (D_u\mu_i^{-1})(\tilde{\la}_i-\si_i)v_i\xi_i[(\pa_{v_i}\bar{v}_i)v_i-\bar{v}_i]B\tilde{r}_{i,v}=O(1)\La_i^1+O(1)(\chi')^\e_iv_i^2\hat{v}_i=O(1)\La_i^1+R_\e^{7,i,1}$ with $R_\e^{7,i,1}$ satisfying \eqref{12.11bis}. Next we have using Lemmas  \ref{estimimpo1}, \ref{lemme6.5},  \ref{lemme6.6}, \ref{lemme11.3}, \eqref{6.49bis} and the fact that $(\pa_{v_i}\bar{v}_iv_i-\bar{v}_i)_x=O(1)\left(\Delta^\ep_i+\varkappa^\ep_i\rho^\ep_i\right)v_{i,x}+O(1)\rho_i^\e\sum_{j\ne i} v_jv_{j,x}$ (it is a direct consequence of Lemma \ref{estimimpo1}) we deduce that
	\begin{align*}
		\al_{i,x}^{7,1}&=\sum\limits_{k}\big[v_{k,x}\tilde{r}_{k}+v_k\tilde{r}_{k,x}\big]\cdot (D_u\mu_i^{-1})(\tilde{\la}_i-\si_i)v_i\xi_i[(\pa_{v_i}\bar{v}_i)v_i-\bar{v}_i]B\tilde{r}_{i,v}\\
		%&+\sum\limits_{k}v_{k,x}\tilde{r}_{k}\cdot (D_u\mu_i^{-1})(\tilde{\la}_i-\si_i)v_i\xi_i[(\pa_{v_i}\bar{v}_i)v_i-\bar{v}_i]B\tilde{r}_{i,v}\\
		&+\sum\limits_{k}v_k\tilde{r}_k\otimes u_x:D^2\mu_i^{-1}(\tilde{\la}_i-\si_i)v_i\xi_i[(\pa_{v_i}\bar{v}_i)v_i-\bar{v}_i]B\tilde{r}_{i,v}\\
		&+\sum\limits_{k}v_k\tilde{r}_k\cdot (D_u\mu_i^{-1})\left(\tilde{\la}_{i,x}+\theta_i^\p\left(\frac{w_i}{v_i}\right)_x\right) v_i\xi_i[(\pa_{v_i}\bar{v}_i)v_i-\bar{v}_i]B\tilde{r}_{i,v}\\
		&+\sum\limits_{k}v_k\tilde{r}_k\cdot (D_u\mu_i^{-1})(\tilde{\la}_i-\si_i)(v_{i,x}\xi_i+\xi_i'v_i\left(\frac{w_i}{v_i}\right)_x)(\pa_{v_i}\bar{v}_iv_i-\bar{v}_i)B\tilde{r}_{i,v}\\
		&+\sum\limits_{k}v_k\tilde{r}_k\cdot (D_u\mu_i^{-1})(\tilde{\la}_i-\si_i)v_i\xi_i(\pa_{v_i}\bar{v}_iv_i-\bar{v}_i)_xB\tilde{r}_{i,v}\\
		&+\sum\limits_{k}v_k \tilde{r}_k\cdot (D_u\mu_i^{-1})(\tilde{\la}_i-\si_i)v_i\xi_i[(\pa_{v_i}\bar{v}_i)v_i-\bar{v}_i][u_x\cdot DB\tilde{r}_{i,v}+B\pa_x\tilde{r}_{i,v}]\\
		&=O(1)(\La_i^1+\La_i^4)+\widetilde{R}_\e^{7,i,1},
		\end{align*}
		with $\widetilde{R}_\e^{7,i,1}$ satisfying \eqref{12.11bis}.
From Lemma \ref{estimimpo1} we have  $\al_{i}^{7,2}=\mu_i^{-1}(\tilde{\la}_i-\si_i)\xi^\p_i\left(\frac{w_i}{v_i}\right)_xv_i[(\pa_{v_i}\bar{v}_i)v_i-\bar{v}_i]B\tilde{r}_{i,v}=O(1)\La_i^4$, furthermore using Lemmas \ref{estimimpo1}, \ref{lemme6.5}, \ref{lemme6.6}, \ref{lemme11.3}, \eqref{identity:wi-vi-xx}, \eqref{6.49bis}, \eqref{ngtech5}
		\begin{align*}
		\al_{i,x}^{7,2}&=\big[u_x\cdot D\mu_i^{-1}(\tilde{\la}_i-\si_i)+\mu_i^{-1}\left(\tilde{\la}_{i,x}+\theta_i^\p\left(\frac{w_i}{v_i}\right)_x\right) \big]\xi^\p_i\left(\frac{w_i}{v_i}\right)_x v_i[(\pa_{v_i}\bar{v}_i)v_i-\bar{v}_i]B\tilde{r}_{i,v}\\
		%&+\mu_i^{-1}\left(\tilde{\la}_{i,x}+\theta_i^\p\left(\frac{w_i}{v_i}\right)_x\right) \xi^\p_i\left(\frac{w_i}{v_i}\right)_xv_i[(\pa_{v_i}\bar{v}_i)v_i-\bar{v}_i]B\tilde{r}_{i,v}\\
		&+\mu_i^{-1}(\tilde{\la}_i-\si_i)\big[\xi^{\p\p}_i\left(\frac{w_i}{v_i}\right)^2_xv_i+\xi^\p_i\left(\frac{w_i}{v_i}\right)_{xx}v_i\big]
		[(\pa_{v_i}\bar{v}_i)v_i-\bar{v}_i]B\tilde{r}_{i,v}\\
		%&+\mu_i^{-1}(\tilde{\la}_i-\si_i)\xi^\p_i\left(\frac{w_i}{v_i}\right)_{xx}v_i[(\pa_{v_i}\bar{v}_i)v_i-\bar{v}_i]B\tilde{r}_{i,v}\\
		&+\mu_i^{-1}(\tilde{\la}_i-\si_i)\xi^\p_i\left(\frac{w_i}{v_i}\right)_x\big[v_{i,x}[(\pa_{v_i}\bar{v}_i)v_i-\bar{v}_i]B\tilde{r}_{i,v}+v_i\pa_x((\pa_{v_i}\bar{v}_i)v_i-\bar{v}_i)B\tilde{r}_{i,v}\big]\\
		%&+\mu_i^{-1}(\tilde{\la}_i-\si_i)\xi^\p_i\left(\frac{w_i}{v_i}\right)_xv_i\pa_x((\pa_{v_i}\bar{v}_i)v_i-\bar{v}_i)B\tilde{r}_{i,v}\\
		&+\mu_i^{-1}(\tilde{\la}_i-\si_i)\xi^\p_i\left(\frac{w_i}{v_i}\right)_xv_i[(\pa_{v_i}\bar{v}_i)v_i-\bar{v}_i][u_x\cdot DB\tilde{r}_{i,v}+B\pa_x\tilde{r}_{i,v}]\\
		&=O(1)(\La_i^1+\delta_0^2\La_i^3+\La_i^4+\La_i^5).
		\end{align*}
		Next we have from \eqref{6.49bis} and Lemmas \ref{lemme6.5} $\al_{i}^{7,3}=\mu_i^{-1}(\tilde{\la}_i-\si_i)\xi_i v_{i,x}[(\pa_{v_i}\bar{v}_i)v_i-\bar{v}_i]B\tilde{r}_{i,v}=O(1)\La_i^1+O(1)(\chi')^\e_i\hat{v}_i v_{i,x}=O(1)\La_i^1+R_\e^{7,i,3}$ with $R_\e^{7,i,3}$ satisfying \eqref{12.11bis}. Furthermore we get from Lemmas \ref{estimimpo1}, \ref{lemme6.5}, \ref{lemme6.6} and \eqref{ngtech5}, \eqref{6.49bis}
		\begin{align*}
		\al_{i,x}^{7,3}&=\big[u_x\cdot D \mu_i^{-1}(\tilde{\la}_i-\si_i)+\mu_i^{-1}(\tilde{\la}_{i,x}+\theta'_i\left(\frac{w_i}{v_i}\right)_x)\big]
		\xi_i v_{i,x}[(\pa_{v_i}\bar{v}_i)v_i-\bar{v}_i]B\tilde{r}_{i,v}	\\
		%&+\mu_i^{-1}(\tilde{\la}_{i,x}+\theta'_i\left(\frac{w_i}{v_i}\right)_x)\xi_i v_{i,x}[(\pa_{v_i}\bar{v}_i)v_i-\bar{v}_i]B\tilde{r}_{i,v}\\
		&+\mu_i^{-1}(\tilde{\la}_i-\si_i)\big[\xi_i' \left(\frac{w_i}{v_i}\right)_x v_{i,x}+\xi_i v_{i,xx}\big]
		[(\pa_{v_i}\bar{v}_i)v_i-\bar{v}_i]B\tilde{r}_{i,v}\\
		%&+\mu_i^{-1}(\tilde{\la}_i-\si_i)\xi_i v_{i,xx}[(\pa_{v_i}\bar{v}_i)v_i-\bar{v}_i]B\tilde{r}_{i,v}\\
		&+\mu_i^{-1}(\tilde{\la}_i-\si_i)\xi_i v_{i,x}[(\pa_{v_i}\bar{v}_i)v_i-\bar{v}_i]_xB\tilde{r}_{i,v}\\
		&+\mu_i^{-1}(\tilde{\la}_i-\si_i)\xi_i v_{i,x}[(\pa_{v_i}\bar{v}_i)v_i-\bar{v}_i][u_x\cdot DB\tilde{r}_{i,v}+B\pa_x(\tilde{r}_{i,v})]\\
		&=O(1)\sum_j(\La_j^1+\La_j^4)+\widetilde{R}_\e^{7,i,3},
		%&+\mu_i^{-1}(\tilde{\la}_i-\si_i)\xi_i v_{i,x}[(\pa_{v_i}\bar{v}_i)v_i-\bar{v}_i]B\pa_x(\tilde{r}_{i,v})
		\end{align*}
		with $\widetilde{R}_\e^{7,i,3}$ satisfying \eqref{12.11bis}.
		In the previous estimate, let us deal with the delicate term $\mu_i^{-1}(\tilde{\la}_i-\si_i)\xi_i v_{i,x}[(\pa_{v_i}\bar{v}_i)v_i-\bar{v}_i]_xB\tilde{r}_{i,v}$, from Lemmas \ref{estimimpo1}, \ref{lemme6.5} and \eqref{ngtech5} we get
		\begin{align*}
		&\mu_i^{-1}(\tilde{\la}_i-\si_i)\xi_i v_{i,x}[(\pa_{v_i}\bar{v}_i)v_i-\bar{v}_i]_xB\tilde{r}_{i,v}=O(1)\mathfrak{A}_i\left(\Delta^\ep_i+\varkappa^\ep_i\rho^\ep_i\right)v^2_{i,x}+O(1)\rho_i^\e\sum_{j\ne i}v_{i,x} v_jv_{j,x}\\
		&=O(1)\sum_j(\La_j^1+\La_j^4)+O(1)\mathfrak{A}_i\left(\Delta^\ep_i+\varkappa^\ep_i\rho^\ep_i\right)v_{i,x}v_i\\
		&=O(1)\sum_j(\La_j^1+\La_j^4)+O(1)\mathfrak{A}_i\varkappa^\ep_i\rho^\ep_iv_{i,x}v_i=O(1)\sum_j(\La_j^1+\La_j^4)+\widetilde{R}_\e^{7,i,3,1}
\end{align*}
		 with $\widetilde{R}_\e^{7,i,3,1}$ satisfying \eqref{12.11bis}.
		%we observe first that:
		%\begin{align}
		%&v_i\pa_{v_iv_i}\bar{v}_i=\pa_{v_i}[\pa_{v_i}\bar{v}_iv_i-\bar{v}_i]=O(1)\alpha_i \sum_{k\ne i}(|\eta'(\frac{v_k^2+w_k^2}{v_i})|+|\eta''(\frac{v_k^2+w_k^2}{v_i})|)\label{superimpl}
		%\end{align}
		%It implies in particular using Lemma \ref{estimimpo1} we deduce that
		%\begin{align}
		%&\pa_x[\pa_{v_i}\bar{v}_iv_i-\bar{v}_i]=O(1)\alpha_i \sum_{k\ne i}(|\eta'(\frac{v_k^2+w_k^2}{v_i})|+|\eta''(\frac{v_k^2+w_k^2}{v_i})|) v_{i,x}+O(1)\sum_{l\ne i}(v_l v_{l,x}+w_l w_{l,x}).
		%\label{ultratech}
		%\end{align}
	%Using the Lemma \ref{lemme6.8}, we get
	%\begin{align*}
	%	&\xi_i v_{i,x}\alpha_i \sum_{k\ne i}(|\eta'(\frac{v_k^2+w_k^2}{v_i})|+|\eta''(\frac{v_k^2+w_k^2}{v_i})|)v_{i,x}\\
	%	&=O(1)v_{i}\alpha_i \sum_{k\ne i}(|\eta'(\frac{v_k^2+w_k^2}{v_i})|+|\eta''(\frac{v_k^2+w_k^2}{v_i})|)v_{i,x}+O(1)\sum_{j\ne i}\La_j^4\\
	%	&=O(1)\sum_{k\ne i}(v_k^2+w_k^2)v_{i,x}+O(1)\sum_{j\ne i}\La_j^4\\
	%	&=O(1)\sum_j(\La_j^1+\La_j^4).
	%	\end{align*}
	%	We obtain now easily that
	%	\begin{align*}
	%	&\mu_i^{-1}(\tilde{\la}_i-\si_i)\xi_i v_{i,x}[(\pa_{v_i}\bar{v}_i)v_i-\bar{v}_i]_xB\tilde{r}_{i,v}=\sum_j(\La_j^1+\La_j^4).
	%	\end{align*}
		Next we have from \eqref{6.49bis} $\al_{i}^{7,4}=\sum\limits_{k}\mu_i^{-1}\tilde{r}_k\cdot D \tilde{\la}_iv_kv_i \xi_i[(\pa_{v_i}\bar{v}_i)v_i-\bar{v}_i]B\tilde{r}_{i,v}=O(1)\La_i^1+R_\e^{7,i,4}$ with $R_\e^{7,i,4}$ satisfying \eqref{12.11bis}. 
		%$\al_{i}^{8,4}=\sum\limits_{k}\mu_i^{-1}\tilde{r}_k\cdot D \tilde{\la}_iv_kv_i \xi_i[(\pa_{v_i}\bar{v}_i)v_i-\bar{v}_i]B\tilde{r}_{i,v}$
		As previously we have using Lemmas  \ref{estimimpo1}, \ref{lemme6.5}, \ref{lemme6.6} and \eqref{6.49bis}
	\begin{align*} 
		\al_{i,x}^{7,4}&=\sum\limits_{k}\big[u_x\cdot D\mu_i^{-1}\tilde{r}_k+\mu_i^{-1}\tilde{r}_{k,x}\big]\cdot D\tilde{\la}_iv_kv_i \xi_i[(\pa_{v_i}\bar{v}_i)v_i-\bar{v}_i]B\tilde{r}_{i,v}\\
		%&+\sum\limits_{k}\mu_i^{-1}\tilde{r}_{k,x}\cdot D\tilde{\la}_iv_kv_i\xi_i[(\pa_{v_i}\bar{v}_i)v_i-\bar{v}_i]B\tilde{r}_{i,v}\\
		&+\sum\limits_{k}\mu_i^{-1}\tilde{r}_k\otimes u_x: D^2\tilde{\la}_iv_kv_i\xi_i[(\pa_{v_i}\bar{v}_i)v_i-\bar{v}_i]B\tilde{r}_{i,v}\\
		&+\sum\limits_{k}\mu_i^{-1}\tilde{r}_k\cdot D\tilde{\la}_{i,v}\pa_{v_i}\bar{v}_i \xi_i v_{i,x}v_kv_i\xi_i[(\pa_{v_i}\bar{v}_i)v_i-\bar{v}_i]B\tilde{r}_{i,v}\\
		&+\sum\limits_{k}\sum\limits_{j\neq i}\mu_i^{-1}\tilde{r}_k\cdot D\tilde{\la}_{i,v}\xi_i\pa_{v_j}\bar{v}_i v_{j,x}v_kv_i\xi_i[(\pa_{v_i}\bar{v}_i)v_i-\bar{v}_i]B\tilde{r}_{i,v}\\
		%&+\sum\limits_{k}\sum\limits_{j\neq i}\mu_i^{-1}\tilde{r}_k\cdot D\tilde{\la}_{i,v}\xi_i\pa_{w_j}\bar{v}_i w_{j,x}v_kv_i\xi_i[(\pa_{v_i}\bar{v}_i)v_i-\bar{v}_i]B\tilde{r}_{i,v}\\
		&+\sum\limits_{k}\left(\frac{w_i}{v_i}\right)_x\mu_i^{-1}\tilde{r}_k\cdot \big[D\tilde{\la}_{i,v}\xi_i'\bar{v}_i -D\tilde{\la}_{i,\sig}\theta_i'\big]
		v_kv_i\xi_i[(\pa_{v_i}\bar{v}_i)v_i-\bar{v}_i]B\tilde{r}_{i,v}\\
		&+\sum\limits_{k}\mu_i^{-1}\tilde{r}_k\cdot D\tilde{\la}_i\big[v_{k,x}v_i \xi_i+v_kv_{i,x}\xi_i+v_kv_i\xi_i^\p\left(\frac{w_i}{v_i}\right)_x\big][(\pa_{v_i}\bar{v}_i)v_i-\bar{v}_i]B\tilde{r}_{i,v}\\
		%&+\sum\limits_{k}\mu_i^{-1}\tilde{r}_k\cdot D\tilde{\la}_iv_kv_i\xi_i^\p\left(\frac{w_i}{v_i}\right)_x[(\pa_{v_i}\bar{v}_i)v_i-\bar{v}_i]B\tilde{r}_{i,v}\\
		&+\sum\limits_{k}\mu_i^{-1}\tilde{r}_k\cdot D\tilde{\la}_iv_kv_i\xi_i\pa_x((\pa_{v_i}\bar{v}_i)v_i-\bar{v}_i)B\tilde{r}_{i,v}\\
		&+\sum\limits_{k}\mu_i^{-1}\tilde{r}_k\cdot D\tilde{\la}_iv_kv_i\xi_i[(\pa_{v_i}\bar{v}_i)v_i-\bar{v}_i][u_x\cdot DB\tilde{r}_{i,v}+B\pa_x\tilde{r}_{i,v}]\\
		&=O(1)\La_i^1+\widetilde{R}_\e^{7,i,4},
		\end{align*}
		with $\widetilde{R}_\e^{7,i,4}$ satisfying \eqref{12.11bis}.
		We have now from \eqref{6.49bis} $\al_{i}^{7,5}=\mu_i^{-1}\tilde{\la}_{i,v}\pa_{v_i}\bar{v}_iv_{i,x}v_i\xi_i^2[(\pa_{v_i}\bar{v}_i)v_i-\bar{v}_i]B\tilde{r}_{i,v}=O(1)\La_i^1+R_\e^{7,i,5}$ with $R_\e^{7,i,5}$ satisfying \eqref{12.11bis}. We get using again Lemmas \ref{estimimpo1}, \ref{lemme6.5}, \ref{lemme6.6}, \ref{lemme11.3} and \eqref{6.49bis}
				\begin{align*} 
	   \al_{i,x}^{7,5}&=\big[u_x\cdot D\mu_i^{-1}\tilde{\la}_{i,v}++\mu_i^{-1}\pa_x(\tilde{\la}_{i,v})\big]\pa_{v_i}\bar{v}_iv_{i,x}v_i\xi_i^2[(\pa_{v_i}\bar{v}_i)v_i-\bar{v}_i]B\tilde{r}_{i,v}\\
	   %&+\mu_i^{-1}\pa_x(\tilde{\la}_{i,v})\pa_{v_i}\bar{v}_iv_{i,x}v_i\xi_i^2[(\pa_{v_i}\bar{v}_i)v_i-\bar{v}_i]B\tilde{r}_{i,v}\\
	   &+\mu_i^{-1}\tilde{\la}_{i,v}\big[\pa_x(\pa_{v_i}\bar{v}_i)v_{i,x}+\pa_{v_i}\bar{v}_iv_{i,xx}\big]v_i\xi_i^2[(\pa_{v_i}\bar{v}_i)v_i-\bar{v}_i]B\tilde{r}_{i,v}\\
	  % &+\mu_i^{-1}\tilde{\la}_{i,v}\pa_{v_i}\bar{v}_iv_{i,xx}v_i\xi_i^2[(\pa_{v_i}\bar{v}_i)v_i-\bar{v}_i]B\tilde{r}_{i,v}\\
	   &+\mu_i^{-1}\tilde{\la}_{i,v}\pa_{v_i}\bar{v}_i v_{i,x}\big[v_{i,x}\xi_i^2+2v_i\xi_i \xi_i^\p\left(\frac{w_i}{v_i}\right)_x\big]
	   [(\pa_{v_i}\bar{v}_i)v_i-\bar{v}_i]B\tilde{r}_{i,v}\\
	  % &+2\mu_i^{-1}\tilde{\la}_{i,v}\pa_{v_i}\bar{v}_iv_{i,x}v_i\xi_i \xi_i^\p\left(\frac{w_i}{v_i}\right)_x[(\pa_{v_i}\bar{v}_i)v_i-\bar{v}_i]B\tilde{r}_{i,v}\\
	   &+\mu_i^{-1}\tilde{\la}_{i,v}\pa_{v_i}\bar{v}_iv_{i,x}v_i\xi_i^2\pa_x((\pa_{v_i}\bar{v}_i)v_i-\bar{v}_i)B\tilde{r}_{i,v}\\
	   &+\mu_i^{-1}\tilde{\la}_{i,v}\pa_{v_i}\bar{v}_iv_{i,x}v_i\xi_i^2[(\pa_{v_i}\bar{v}_i)v_i-\bar{v}_i][u_x\cdot DB\tilde{r}_{i,v}+B\pa_x\tilde{r}_{i,v}]\\
	   &=O(1)\La_i^1+\widetilde{R}_\e^{7,i,5}
	   \end{align*}
	   with $\widetilde{R}_\e^{7,i,5}$ satisfying \eqref{12.11bis}.
	   We have now $ \al_{i}^{7,6}=\sum\limits_{j\neq i}\mu_i^{-1}\tilde{\la}_{i,v}\pa_{v_j}\bar{v}_iv_{j,x}v_i \xi_i^2[(\pa_{v_i}\bar{v}_i)v_i-\bar{v}_i]B\tilde{r}_{i,v}=O(1)\La_i^1$ and using Lemmas  \ref{estimimpo1}, \ref{lemme6.5}, \ref{lemme6.6}, \ref{lemme11.3} 
	   \begin{align*}  
	   \al_{i,x}^{7,6}&=\sum\limits_{j\neq i}\big[u_x\cdot D\mu_i^{-1}\tilde{\la}_{i,v}+\mu_i^{-1}\pa_x(\tilde{\la}_{i,v})\big]
	   \pa_{v_j}\bar{v}_iv_{j,x}v_i\xi_i^2[(\pa_{v_i}\bar{v}_i)v_i-\bar{v}_i]B\tilde{r}_{i,v}\\
	  % &+\sum\limits_{j\neq i}\mu_i^{-1}\pa_x(\tilde{\la}_{i,v})\pa_{v_j}\bar{v}_iv_{j,x}v_i\xi_i^2[(\pa_{v_i}\bar{v}_i)v_i-\bar{v}_i]B\tilde{r}_{i,v}\\
	   &+\sum\limits_{j\neq i}\mu_i^{-1}\tilde{\la}_{i,v}\big[\pa_x(\pa_{v_j}\bar{v}_i)v_{j,x}+
	   \pa_{v_j}\bar{v}_iv_{j,xx}\big]v_i\xi_i^2[(\pa_{v_i}\bar{v}_i)v_i-\bar{v}_i]B\tilde{r}_{i,v}\\
	  % &+\sum\limits_{j\neq i}\mu_i^{-1}\tilde{\la}_{i,v}\pa_{v_j}\bar{v}_iv_{j,xx}v_i\xi_i^2[(\pa_{v_i}\bar{v}_i)v_i-\bar{v}_i]B\tilde{r}_{i,v}\\
	   &+\sum\limits_{j\neq i}\mu_i^{-1}\tilde{\la}_{i,v}\pa_{v_j}\bar{v}_iv_{j,x}\big[v_{i,x}\xi_i^2
	   +2v_i\xi_i^\p\xi_i \left(\frac{w_i}{v_i}\right)_x\big]
	   [(\pa_{v_i}\bar{v}_i)v_i-\bar{v}_i]B\tilde{r}_{i,v}\\
	   %&+2\sum\limits_{j\neq i}\mu_i^{-1}\tilde{\la}_{i,v}\pa_{v_j}\bar{v}_iv_{j,x}v_i\xi_i^\p\xi_i \left(\frac{w_i}{v_i}\right)_x[(\pa_{v_i}\bar{v}_i)v_i-\bar{v}_i]B\tilde{r}_{i,v}\\
	   &+\sum\limits_{j\neq i}\mu_i^{-1}\tilde{\la}_{i,v}\pa_{v_j}\bar{v}_iv_{j,x}v_i\xi_i^2\pa_x((\pa_{v_i}\bar{v}_i)v_i-\bar{v}_i)B\tilde{r}_{i,v}\\
	   &+\sum\limits_{j\neq i}\mu_i^{-1}\tilde{\la}_{i,v}\pa_{v_j}\bar{v}_iv_{j,x}v_i\xi_i^2[(\pa_{v_i}\bar{v}_i)v_i-\bar{v}_i][u_x\cdot DB\tilde{r}_{i,v}+B\pa_x\tilde{r}_{i,v}]\\
	   &=O(1)\La_i^1.
	   \end{align*}
 Next we deal with $\al_{i}^{7,7}=\mu_i^{-1}\tilde{\la}_{i,v}\xi_i^\p\left(\frac{w_i}{v_i}\right)_x\bar{v}_i v_i\xi_i[(\pa_{v_i}\bar{v}_i)v_i-\bar{v}_i]B\tilde{r}_{i,v}=O(1)\La_i^4$ and we get from  Lemmas  \ref{estimimpo1}, \ref{lemme6.5}, \ref{lemme6.6}, \ref{lemme11.3} 
\begin{align*}
\al_{i,x}^{7,7}&=\big[u_x\cdot D\mu_i^{-1}\tilde{\la}_{i,v}+\mu_i^{-1}\pa_x(\tilde{\la}_{i,v})\big]
\xi_i^\p\left(\frac{w_i}{v_i}\right)_x\bar{v}_i v_i\xi_i[(\pa_{v_i}\bar{v}_i)v_i-\bar{v}_i]B\tilde{r}_{i,v}\\
%&+\mu_i^{-1}\pa_x(\tilde{\la}_{i,v})\xi_i^\p\left(\frac{w_i}{v_i}\right)_x\bar{v}_i v_i\xi_i[(\pa_{v_i}\bar{v}_i)v_i-\bar{v}_i]B\tilde{r}_{i,v}\\
&+\mu_i^{-1}\tilde{\la}_{i,v}\big[\xi_i''\left(\frac{w_i}{v_i}\right)^2_x+\xi_i^\p\left(\frac{w_i}{v_i}\right)_{xx}\big]
\bar{v}_i v_i\xi_i[(\pa_{v_i}\bar{v}_i)v_i-\bar{v}_i]B\tilde{r}_{i,v}\\
%&+\mu_i^{-1}\tilde{\la}_{i,v}\xi_i^\p\left(\frac{w_i}{v_i}\right)_{xx}\bar{v}_i v_i\xi_i[(\pa_{v_i}\bar{v}_i)v_i-\bar{v}_i]B\tilde{r}_{i,v}\\
&+\mu_i^{-1}\tilde{\la}_{i,v}\xi_i^\p\left(\frac{w_i}{v_i}\right)_x\big[
\pa_x(\bar{v}_i) v_i\xi_i+\bar{v}_i v_{i,x}\xi_i+\bar{v}_i v_{i}\xi'_i\left(\frac{w_i}{v_i}\right)\big]
[(\pa_{v_i}\bar{v}_i)v_i-\bar{v}_i]B\tilde{r}_{i,v}\\
%&+\mu_i^{-1}\tilde{\la}_{i,v}\xi_i^\p\left(\frac{w_i}{v_i}\right)_x\bar{v}_i v_{i,x}\xi_i[(\pa_{v_i}\bar{v}_i)v_i-\bar{v}_i]B\tilde{r}_{i,v}\\
%&+\mu_i^{-1}\tilde{\la}_{i,v}(\xi_i^\p)^2\left(\frac{w_i}{v_i}\right)^2_x\bar{v}_i v_i[(\pa_{v_i}\bar{v}_i)v_i-\bar{v}_i]B\tilde{r}_{i,v}\\
&+\mu_i^{-1}\tilde{\la}_{i,v}\xi_i^\p\left(\frac{w_i}{v_i}\right)_x\bar{v}_i v_i\xi_i\pa_x[(\pa_{v_i}\bar{v}_i)v_i-\bar{v}_i]B\tilde{r}_{i,v}\\
&+\mu_i^{-1}\tilde{\la}_{i,v}\xi_i^\p\left(\frac{w_i}{v_i}\right)_x\bar{v}_i v_i\xi_i[(\pa_{v_i}\bar{v}_i)v_i-\bar{v}_i][u_x\cdot DB\tilde{r}_{i,v}+B\pa_x(\tilde{r}_{i,v})]\\
&=O(1)(\delta_0^2\La_i^3+\La_i^4+\La_i^5).
\end{align*}
We have now  $\al_{i}^{7,8}=-\mu_i^{-1}\tilde{\la}_{i,\si}\theta_i^\p\left(\frac{w_i}{v_i}\right)_xv_i\xi_i[(\pa_{v_i}\bar{v}_i)v_i-\bar{v}_i]B\tilde{r}_{i,v}=O(1)\La_i^4$, using Lemmas \ref{tildelamb}, \ref{estimimpo1}, \ref{lemme6.5}, \ref{lemme6.6}, \ref{lemme11.3} , we deduce now
	  % \begin{align*}
	   %\al_{i,x}^{8,6}&=-\sum\limits_{j\neq i}u_x\cdot D\mu_i^{-1}\tilde{\la}_{i,v}\pa_{w_j}\bar{v}_iw_{j,x}v_i\xi_i[(\pa_{v_i}\bar{v}_i)v_i-\bar{v}_i]B\tilde{r}_{i,v}\\
	   %&-\sum\limits_{j\neq i}\mu_i^{-1}\pa_x(\tilde{\la}_{i,v})\pa_{w_j}\bar{v}_iw_{j,x}v_i\xi_i[(\pa_{v_i}\bar{v}_i)v_i-\bar{v}_i]B\tilde{r}_{i,v}\\
	   %&-\sum\limits_{j\neq i}\mu_i^{-1}\tilde{\la}_{i,v}\pa_x(\pa_{w_j}\bar{v}_i)w_{j,x}v_i\xi_i[(\pa_{v_i}\bar{v}_i)v_i-\bar{v}_i]B\tilde{r}_{i,v}\\
	   %&-\sum\limits_{j\neq i}\mu_i^{-1}\tilde{\la}_{i,v}\pa_{w_j}\bar{v}_iw_{j,xx}v_i\xi_i[(\pa_{v_i}\bar{v}_i)v_i-\bar{v}_i]B\tilde{r}_{i,v}\\
	   %&-\sum\limits_{j\neq i}\mu_i^{-1}\tilde{\la}_{i,v}\pa_{w_j}\bar{v}_iw_{j,x}v_{i,x}\xi_i[(\pa_{v_i}\bar{v}_i)v_i-\bar{v}_i]B\tilde{r}_{i,v}\\
	   %&-\sum\limits_{j\neq i}\mu_i^{-1}\tilde{\la}_{i,v}\pa_{w_j}\bar{v}_iw_{j,x}v_i\xi_i^\p\left(\frac{w_i}{v_i}\right)_x[(\pa_{v_i}\bar{v}_i)v_i-\bar{v}_i]B\tilde{r}_{i,v}\\
	   %&-\sum\limits_{j\neq i}\mu_i^{-1}\tilde{\la}_{i,v}\pa_{w_j}\bar{v}_iw_{j,x}v_i\xi_i\pa_x((\pa_{v_i}\bar{v}_i)v_i-\bar{v}_i)B\tilde{r}_{i,v}\\
	   %&-\sum\limits_{j\neq i}\mu_i^{-1}\tilde{\la}_{i,v}\pa_{=w_j}\bar{v}_iw_{j,x}v_i\xi_i[(\pa_{v_i}\bar{v}_i)v_i-\bar{v}_i][u_x\cdot DB\tilde{r}_{i,v}+B\pa_x\tilde{r}_{i,v}]\\
	   %&=\mathcal{O}(1)(\La_i^1+\La_i^3+\La_i^5),\\  
	   \begin{align*}
	   \al_{i,x}^{7,8}&=-\big[u_x\cdot D\mu_i^{-1}\tilde{\la}_{i,\si}+\mu_i^{-1}\pa_x(\tilde{\la}_{i,\si})\big]
	   \theta_i^\p\left(\frac{w_i}{v_i}\right)_xv_i\xi_i[(\pa_{v_i}\bar{v}_i)v_i-\bar{v}_i]B\tilde{r}_{i,v}\\
	  % &-\mu_i^{-1}\pa_x(\tilde{\la}_{i,\si})\theta_i^\p\left(\frac{w_i}{v_i}\right)_xv_i\xi_i[(\pa_{v_i}\bar{v}_i)v_i-\bar{v}_i]B\tilde{r}_{i,v} \\
	   &-\mu_i^{-1}\tilde{\la}_{i,\si}\big[\theta_i^{\p\p}\left(\frac{w_i}{v_i}\right)^2_x
	   +\theta_i^\p\left(\frac{w_i}{v_i}\right)_{xx}\big]
	   v_i\xi_i[(\pa_{v_i}\bar{v}_i)v_i-\bar{v}_i]B\tilde{r}_{i,v}\\
	   %&-\mu_i^{-1}\tilde{\la}_{i,\si}\theta_i^\p\left(\frac{w_i}{v_i}\right)_{xx}v_i\xi_i[(\pa_{v_i}\bar{v}_i)v_i-\bar{v}_i]B\tilde{r}_{i,v}\\
	   &-\mu_i^{-1}\tilde{\la}_{i,\si}\theta_i^\p\left(\frac{w_i}{v_i}\right)_x\big[v_{i,x}\xi_i+v_i \xi_i'\left(\frac{w_i}{v_i}\right)_x\big]
	   [(\pa_{v_i}\bar{v}_i)v_i-\bar{v}_i]B\tilde{r}_{i,v}\\
	   %&-\mu_i^{-1}\tilde{\la}_{i,\si}\theta_i^\p\left(\frac{w_i}{v_i}\right)^2_xv_i\xi_i^\p[(\pa_{v_i}\bar{v}_i)v_i-\bar{v}_i]B\tilde{r}_{i,v}\\
	   &-\mu_i^{-1}\tilde{\la}_{i,\si}\theta_i^\p\left(\frac{w_i}{v_i}\right)_xv_i\xi_i\pa_x((\pa_{v_i}\bar{v}_i)v_i-\bar{v}_i)B\tilde{r}_{i,v}\\
	   &-\mu_i^{-1}\tilde{\la}_{i,\si}\theta_i^\p\left(\frac{w_i}{v_i}\right)_xv_i\xi_i[(\pa_{v_i}\bar{v}_i)v_i-\bar{v}_i][u_x\cdot DB\tilde{r}_{i,v}+B\pa_x\tilde{r}_{i,v}]\\
	   &=O(1)(\delta_0^2\La_i^3+\La_i^4+\La_i^5).
	   \end{align*}
	   We have now $ \al_{i}^{7,9}=\mu_i^{-1}\theta_i^\p\left(\frac{w_i}{v_i}\right)_xv_i\xi_i[(\pa_{v_i}\bar{v}_i)v_i-\bar{v}_i]B\tilde{r}_{i,v}=O(1)\La_i^4$ and using  Lemmas \ref{estimimpo1}, \ref{lemme6.5}, \ref{lemme6.6}, \ref{lemme11.3} and \eqref{ngtech5}, we obtain
	   \begin{align*}
	   \al_{i,x}^{7,9}&=\big[u_x\cdot D\mu_i^{-1}\theta_i^\p\left(\frac{w_i}{v_i}\right)_x+\mu_i^{-1}\theta_i^{\p\p}\left(\frac{w_i}{v_i}\right)^2_x\big]
	   v_i\xi_i[(\pa_{v_i}\bar{v}_i)v_i-\bar{v}_i]B\tilde{r}_{i,v}\\
	 %  &+\mu_i^{-1}\theta_i^{\p\p}\left(\frac{w_i}{v_i}\right)^2_xv_i\xi_i[(\pa_{v_i}\bar{v}_i)v_i-\bar{v}_i]B\tilde{r}_{i,v}\\
	   &+\mu_i^{-1}\theta_i^\p\big[\left(\frac{w_i}{v_i}\right)_{xx}v_i \xi_i+\left(\frac{w_i}{v_i}\right)_xv_{i,x}\xi_i+\left(\frac{w_i}{v_i}\right)^2_xv_i\xi_i^\p\big]
	   [(\pa_{v_i}\bar{v}_i)v_i-\bar{v}_i]B\tilde{r}_{i,v}\\
	   %&+\mu_i^{-1}\theta_i^\p\left(\frac{w_i}{v_i}\right)_xv_{i,x}\xi_i[(\pa_{v_i}\bar{v}_i)v_i-\bar{v}_i]B\tilde{r}_{i,v}\\
	  % &+\mu_i^{-1}\theta_i^\p\left(\frac{w_i}{v_i}\right)^2_xv_i\xi_i^\p[(\pa_{v_i}\bar{v}_i)v_i-\bar{v}_i]B\tilde{r}_{i,v}\\
	   &+\mu_i^{-1}\theta_i^\p\left(\frac{w_i}{v_i}\right)_xv_i\xi_i\pa_x((\pa_{v_i}\bar{v}_i)v_i-\bar{v}_i)B\tilde{r}_{i,v}\\
	   &+\mu_i^{-1}\theta_i^\p\left(\frac{w_i}{v_i}\right)_xv_i\xi_i[(\pa_{v_i}\bar{v}_i)v_i-\bar{v}_i][u_x\cdot DB\tilde{r}_{i,v}+B\pa_x\tilde{r}_{i,v}]\\
	   &=O(1)\sum_j(\La_i^j+\delta_0^2\La_j^3+\La_j^4+\La_j^5).
	   \end{align*}
	   We have now using Lemma \ref{estimimpo1} $ \al_{i}^{7,10}=\mu_i^{-1}(\tilde{\la}_i-\si_i)v^2_iv_{i,x}\xi_i(\pa_{v_iv_i}\bar{v}_i)B\tilde{r}_{i,v}=O(1)\La_i^1+R_\e^{7,i,10}$ with $R_\e^{7,i,10}$ satisfying \eqref{12.11bis}. We deduce using Lemmas \ref{estimimpo1}, \ref{lemme6.5},  \ref{lemme6.6}, \ref{lemme11.3} and \eqref{ngtech5}
	   \begin{align*}
	   \al_{i,x}^{7,10}&=\big[u_x\cdot D\mu_i^{-1}(\tilde{\la}_i-\si_i)+\mu_i^{-1}\left(\tilde{\la}_{i,x}+\theta_i^\p\left(\frac{w_i}{v_i}\right)_x\right)\big]
	    v^2_iv_{i,x}\xi_i(\pa_{v_iv_i}\bar{v}_i)B\tilde{r}_{i,v}\\
	   &+\mu_i^{-1}(\tilde{\la}_i-\si_i)\big[2v_iv^2_{i,x}\xi_i+v_i^2v_{i,xx}\xi_i+v^2_iv_{i,x}\xi_i^\p\left(\frac{w_i}{v_i}\right)_x\big]
	   (\pa_{v_iv_i}\bar{v}_i)B\tilde{r}_{i,v}\\
	   %&+\mu_i^{-1}(\tilde{\la}_i-\si_i)v^2_iv_{i,xx}\xi_i(\pa_{v_iv_i}\bar{v}_i)B\tilde{r}_{i,v}+\mu_i^{-1}(\tilde{\la}_i-\si_i)v^2_iv_{i,x}\xi_i^\p\left(\frac{w_i}{v_i}\right)_x(\pa_{v_iv_i}\bar{v}_i)B\tilde{r}_{i,v}\\
	   &+\mu_i^{-1}(\tilde{\la}_i-\si_i)v^2_iv_{i,x}\xi_i\big[\pa_x(\pa_{v_iv_i}\bar{v}_i)B\tilde{r}_{i,v}+(\pa_{v_iv_i}\bar{v}_i)[u_x\cdot DB\tilde{r}_{i,v}+B\pa_x\tilde{r}_{i,v}]\big]\\
	   &=O(1)\sum_j(\La_j^1+\La_j^4)+\widetilde{R}_\e^{7,i,10},
	   \end{align*}
	   with $\widetilde{R}_\e^{7,i,10}$ satisfying \eqref{12.11bis}.
	   We are going to give more details here on the most delicate term $\mu_i^{-1}(\tilde{\la}_i-\si_i)v^2_iv_{i,x}\xi_i\pa_x(\pa_{v_iv_i}\bar{v}_i)B\tilde{r}_{i,v}$. We observe that using Lemmas \ref{estimimpo1}, \ref{lemme6.5}, \ref{lemme6.6} and \eqref{ngtech5} we get
	   \begin{align*}
	   &\mu_i^{-1}(\tilde{\la}_i-\si_i)v^2_iv_{i,x}\xi_i\pa_x(\pa_{v_iv_i}\bar{v}_i)B\tilde{r}_{i,v}=O(1)
	   (\Delta^\ep_i+\varkappa^\ep_i\rho^\ep_i)v_{i,x}^2+O(1) \sum_{j\ne i} \Delta^\ep_{ij}|v_jv_{j,x}v_{i,x}|\\
	   &=O(1)\xi_i
	   (\Delta^\ep_i+\varkappa^\ep_i\rho^\ep_i)v_{i,x}^2+O(1) \sum_{j\ne i} \Delta^\ep_{ij}|v_jv_{j,x}v_{i,x}|\\
	   &=O(1)\xi_i
	   (\Delta^\ep_i+\varkappa^\ep_i\rho^\ep_i)v_{i,x}v_i+O(1)\sum_j(\La_j^1+\La_j^4)=O(1)\xi_i \varkappa^\ep_i\rho^\ep_iv_{i,x}v_i+O(1)\sum_j(\La_j^1+\La_j^4),
	   \end{align*}
	   with $\widetilde{R}_\e^{7,i,10,1}=O(1)\xi_i \varkappa^\ep_i\rho^\ep_iv_{i,x}v_i$ satisfying \eqref{12.11bis}.
	   We have now $\al_{i}^{7,11}=\sum\limits_{j\neq i}\mu_i^{-1}(\tilde{\la}_i-\si_i)v_iv_{j,x}\xi_i[(\pa_{v_iv_j}\bar{v}_i)v_i-\pa_{v_j}\bar{v}_i]B\tilde{r}_{i,v}=O(1)\La_i^1$ and from Lemmas   \ref{estimimpo1}, \ref{lemme11.3}
	   \begin{align*}
	   \al_{i,x}^{7,11}&=\sum\limits_{j\neq i}\big[u_x\cdot D\mu_i^{-1}(\tilde{\la}_i-\si_i)+\mu_i^{-1}\left(\tilde{\la}_{i,x}+\theta_i^\p\left(\frac{w_i}{v_i}\right)_x\right)\big]
	   v_iv_{j,x}\xi_i[(\pa_{v_iv_j}\bar{v}_i)v_i-\pa_{v_j}\bar{v}_i]B\tilde{r}_{i,v}\\
	   %&+\sum\limits_{j\neq i}\mu_i^{-1}\left(\tilde{\la}_{i,x}+\theta_i^\p\left(\frac{w_i}{v_i}\right)_x\right)v_iv_{j,x}\xi_i[(\pa_{v_iv_j}\bar{v}_i)v_i-\pa_{v_j}\bar{v}_i]B\tilde{r}_{i,v}\\
	   &+\sum\limits_{j\neq i}\mu_i^{-1}(\tilde{\la}_i-\si_i)\big[v_{i,x}v_{j,x}\xi_i+v_iv_{j,xx}\xi_i+v_iv_{j,x}\xi_i^\p\left(\frac{w_i}{v_i}\right)_x\big]
	   [(\pa_{v_iv_j}\bar{v}_i)v_i-\pa_{v_j}\bar{v}_i]B\tilde{r}_{i,v}\\
	  % &+\sum\limits_{j\neq i}\mu_i^{-1}(\tilde{\la}_i-\si_i)v_iv_{j,xx}\xi_i[(\pa_{v_iv_j}\bar{v}_i)v_i-\pa_{v_j}\bar{v}_i]B\tilde{r}_{i,v}\\
	  % &+\sum\limits_{j\neq i}\mu_i^{-1}(\tilde{\la}_i-\si_i)v_iv_{j,x}\xi_i^\p\left(\frac{w_i}{v_i}\right)_x[(\pa_{v_iv_j}\bar{v}_i)v_i-\pa_{v_j}\bar{v}_i]B\tilde{r}_{i,v}\\
	   &+\sum\limits_{j\neq i}\mu_i^{-1}(\tilde{\la}_i-\si_i)v_iv_{j,x}\xi_i\pa_x((\pa_{v_iv_j}\bar{v}_i)v_i-\pa_{v_j}\bar{v}_i)B\tilde{r}_{i,v}\\
	   &+\sum\limits_{j\neq i}\mu_i^{-1}(\tilde{\la}_i-\si_i)v_iv_{j,x}\xi_i[(\pa_{v_iv_j}\bar{v}_i)v_i-\pa_{v_j}\bar{v}_i][u_x\cdot DB\tilde{r}_{i,v}+B\pa_x\tilde{r}_{i,v}]=O(1)\La_i^1.
	   \end{align*}
	   %We proceed similarly for $ \al_{i}^{8,13}=\sum\limits_{j\neq i}\mu_i^{-1}(\tilde{\la}_i-\si_i)v_iw_{j,x}\xi_i[(\pa_{v_iw_j}\bar{v}_i)v_i-\pa_{w_j}\bar{v}_i]B\tilde{r}_{i,v}=O(1)\La_i^1$. 
	   We have now from Lemmas \ref{lemme6.5}, \ref{lemme6.6} and \eqref{6.49bis} $ \al_{i}^{7,12}= \mu_i^{-1}(\tilde{\la}_i-\si_i)v_i\xi_i[(\pa_{v_i}\bar{v}_i)v_i-\bar{v}_i]u_x\cdot DB\tilde{r}_{i,v}=O(1)\La_i^1+R_\e^{7,i,12}$ with $R_\e^{7,i,12}$ satisfying \eqref{12.11bis}. Using Lemmas\ref{estimimpo1}, \ref{lemme6.5}, \ref{lemme6.6}, \ref{lemme11.3} and \eqref{6.49bis} it yields
	  % \begin{align*}
	   %\al_{i,x}^{8,11}&=-\sum\limits_{j\neq i}u_x\cdot D\mu_i^{-1}(\tilde{\la}_i-\si_i)v_iw_{j,x}\xi_i[(\pa_{v_iw_j}\bar{v}_i)v_i-\pa_{w_j}\bar{v}_i]B\tilde{r}_{i,v}\\
	   %&-\sum\limits_{j\neq i}\mu_i^{-1}\left(\tilde{\la}_{i,x}+\theta_i^\p\left(\frac{w_i}{v_i}\right)_x\right)v_iw_{j,x}\xi_i[(\pa_{v_iw_j}\bar{v}_i)v_i-\pa_{w_j}\bar{v}_i]B\tilde{r}_{i,v}\\
	   %&-\sum\limits_{j\neq i}\mu_i^{-1}(\tilde{\la}_i-\si_i)v_{i,x}w_{j,x}\xi_i[(\pa_{v_iw_j}\bar{v}_i)v_i-\pa_{w_j}\bar{v}_i]B\tilde{r}_{i,v}\\
	   %&-\sum\limits_{j\neq i}\mu_i^{-1}(\tilde{\la}_i-\si_i)v_iw_{j,xx}\xi_i[(\pa_{v_iw_j}\bar{v}_i)v_i-\pa_{w_j}\bar{v}_i]B\tilde{r}_{i,v}\\
	   %&-\sum\limits_{j\neq i}\mu_i^{-1}(\tilde{\la}_i-\si_i)v_iw_{j,x}\xi_i^\p\left(\frac{w_i}{v_i}\right)_x[(\pa_{v_iw_j}\bar{v}_i)v_i-\pa_{w_j}\bar{v}_i]B\tilde{r}_{i,v}\\
	   %&-\sum\limits_{j\neq i}\mu_i^{-1}(\tilde{\la}_i-\si_i)v_iw_{j,x}\xi_i\pa_x((\pa_{v_iw_j}\bar{v}_i)v_i-\pa_{w_j}\bar{v}_i)B\tilde{r}_{i,v}\\
	   %&-\sum\limits_{j\neq i}\mu_i^{-1}(\tilde{\la}_i-\si_i)v_iw_{j,x}\xi_i[(\pa_{v_iw_j}\bar{v}_i)v_i-\pa_{w_j}\bar{v}_i][u_x\cdot DB\tilde{r}_{i,v}+B\pa_x\tilde{r}_{i,v}]\\
	   %&=\mathcal{O}(1)(\La_i^1+\La_i^3+\La_i^5),\\  
	   \begin{align*}
	   \al_{i,x}^{7,12}&= \big[u_x\cdot D\mu_i^{-1}(\tilde{\la}_i-\si_i)+\mu_i^{-1}\left(\tilde{\la}_{i,x}+\theta_i^\p\left(\frac{w_i}{v_i}\right)_x\right)\big]
	   v_i\xi_i[(\pa_{v_i}\bar{v}_i)v_i-\bar{v}_i]u_x\cdot DB\tilde{r}_{i,v}\\
	   %&+\mu_i^{-1}\left(\tilde{\la}_{i,x}+\theta_i^\p\left(\frac{w_i}{v_i}\right)_x\right)v_i\xi_i[(\pa_{v_i}\bar{v}_i)v_i-\bar{v}_i]u_x\cdot DB\tilde{r}_{i,v}\\
	   &+\mu_i^{-1}(\tilde{\la}_i-\si_i)\big[v_{i,x}\xi_i+v_i\xi_i^\p\left(\frac{w_i}{v_i}\right)_x\big]
	    [(\pa_{v_i}\bar{v}_i)v_i-\bar{v}_i]u_x\cdot DB\tilde{r}_{i,v}\\
	   %&+\mu_i^{-1}(\tilde{\la}_i-\si_i)v_i\xi_i^\p\left(\frac{w_i}{v_i}\right)_x[(\pa_{v_i}\bar{v}_i)v_i-\bar{v}_i]u_x\cdot DB\tilde{r}_{i,v}\\
	   &+\mu_i^{-1}(\tilde{\la}_i-\si_i)v_i\xi_i\pa_x((\pa_{v_i}\bar{v}_i)v_i-\bar{v}_i)u_x\cdot DB\tilde{r}_{i,v}\\
	   &+\mu_i^{-1}(\tilde{\la}_i-\si_i)v_i\xi_i[(\pa_{v_i}\bar{v}_i)v_i-\bar{v}_i]\big[u_{xx}\cdot DB\tilde{r}_{i,v}+(u_x\otimes u_x):D^2B\tilde{r}_{i,v}+u_x\cdot DB\pa_x(\tilde{r}_{i,v})\big] \\
	   %&+\mu_i^{-1}(\tilde{\la}_i-\si_i)v_i\xi_i[(\pa_{v_i}\bar{v}_i)v_i-\bar{v}_i](u_x\otimes u_x):D^2B\tilde{r}_{i,v}\\
	   %&+\mu_i^{-1}(\tilde{\la}_i-\si_i)v_i\xi_i[(\pa_{v_i}\bar{v}_i)v_i-\bar{v}_i]u_x\cdot DB\pa_x(\tilde{r}_{i,v})\\
	   &=O(1)(\La_i^1+\La_i^4)+\widetilde{R}_\e^{7,i,12},
	   \end{align*}
	    with $\widetilde{R}_\e^{7,i,12}$ satisfying \eqref{12.11bis}.
	   Next we have from Lemma \ref{lemme6.5} and \eqref{6.49bis} $ \al_{i}^{7,13}=\sum\limits_{k}\mu_i^{-1}(\tilde{\la}_i-\si_i)v_i\xi_i v_k[(\pa_{v_i}\bar{v}_i)v_i-\bar{v}_i]B\tilde{r}_{i,uv}\tilde{r}_k=O(1)\La_i^1+R_\e^{7,i,13}$ with $R_\e^{7,i,13}$ satisfying \eqref{12.11bis}. From the Lemmas \ref{estimimpo1}, \ref{lemme6.5}, \ref{lemme6.6}, \ref{lemme11.3}, \ref{lemme11.4} and \eqref{6.49bis} we obtain
\begin{align*} 
	   \al_{i,x}^{7,13}&=\sum\limits_{k}\big[u_x\cdot D\mu_i^{-1}(\tilde{\la}_i-\si_i)+\mu_i^{-1}\left(\tilde{\la}_{i,x}+\theta_i^\p\left(\frac{w_i}{v_i}\right)_x\right)\big]
	   v_iv_k\xi_i[(\pa_{v_i}\bar{v}_i)v_i-\bar{v}_i]B\tilde{r}_{i,uv}\tilde{r}_k\\
	   %&+\sum\limits_{k}\mu_i^{-1}\left(\tilde{\la}_{i,x}+\theta_i^\p\left(\frac{w_i}{v_i}\right)_x\right)v_iv_k\xi_i[(\pa_{v_i}\bar{v}_i)v_i-\bar{v}_i]B\tilde{r}_{i,uv}\tilde{r}_k\
	   &+\sum\limits_{k}\mu_i^{-1}(\tilde{\la}_i-\si_i)(v_{i,x}v_k \xi_i +v_iv_{k,x}\xi_i+v_iv_k\xi_i^\p\left(\frac{w_i}{v_i}\right)_x )\xi_i[(\pa_{v_i}\bar{v}_i)v_i-\bar{v}_i]B\tilde{r}_{i,uv}\tilde{r}_k\\
	   %&+\sum\limits_{k}\mu_i^{-1}(\tilde{\la}_i-\si_i)v_iv_k\xi_i^\p\left(\frac{w_i}{v_i}\right)_x[(\pa_{v_i}\bar{v}_i)v_i-\bar{v}_i]B\tilde{r}_{i,uv}\tilde{r}_k\\
	   &+\sum\limits_{k}\mu_i^{-1}(\tilde{\la}_i-\si_i)v_iv_k\xi_i\pa_x((\pa_{v_i}\bar{v}_i)v_i-\bar{v}_i)B\tilde{r}_{i,uv}\tilde{r}_k\\
	   &+\sum\limits_{k}\mu_i^{-1}(\tilde{\la}_i-\si_i)v_iv_k\xi_i[(\pa_{v_i}\bar{v}_i)v_i-\bar{v}_i][u_x\cdot DB\tilde{r}_{i,uv}\tilde{r}_k+B\pa_x(\tilde{r}_{i,uv}\tilde{r}_k)]\\
	   &=O(1)(\La_i^1+\La_i^4)+\widetilde{R}_\e^{7,i,13},
	   \end{align*}
	   with $\widetilde{R}_\e^{7,i,13}$ satisfying \eqref{12.11bis}.
	We have now using \eqref{6.49bis} $\al_{i}^{7,14}= \mu_i^{-1}(\tilde{\la}_i-\si_i)v_iv_{i,x}\xi_i^2\pa_{v_i}\bar{v}_i[(\pa_{v_i}\bar{v}_i)v_i-\bar{v}_i]B\tilde{r}_{i,vv}=O(1)\La_i^1+R_\e^{7,i,14}$ with $R_\e^{7,i,14}$ satisfying \eqref{12.11bis}. From the Lemmas \ref{estimimpo1}, \ref{lemme6.5}, \ref{lemme6.6}, \ref{lemme11.3} and \eqref{6.49bis} we get
	   \begin{align*} 
	   \al_{i,x}^{7,14}&=\big[u_x\cdot D \mu_i^{-1}(\tilde{\la}_i-\si_i)+\mu_i^{-1}\left(\tilde{\la}_{i,x}+\theta_i^\p\left(\frac{w_i}{v_i}\right)_x\right)\big]
	   v_iv_{i,x}\xi_i^2\pa_{v_i}\bar{v}_i[(\pa_{v_i}\bar{v}_i)v_i-\bar{v}_i]B\tilde{r}_{i,vv}\\
	   %&+\mu_i^{-1}\left(\tilde{\la}_{i,x}+\theta_i^\p\left(\frac{w_i}{v_i}\right)_x\right)v_iv_{i,x}\xi_i^2\pa_{v_i}\bar{v}_i[(\pa_{v_i}\bar{v}_i)v_i-\bar{v}_i]B\tilde{r}_{i,vv}\\
	   &+\mu_i^{-1}(\tilde{\la}_i-\si_i)\big[v_{i,x}^2\xi_i^2+v_iv_{i,xx}\xi_i^2+2v_iv_{i,x}\xi_i^\p\xi_i\left(\frac{w_i}{v_i}\right)_x \big]\pa_{v_i}\bar{v}_i[(\pa_{v_i}\bar{v}_i)v_i-\bar{v}_i]B\tilde{r}_{i,vv}\\
	   %&+2\mu_i^{-1}(\tilde{\la}_i-\si_i)v_iv_{i,x}\xi_i^\p\xi_i\left(\frac{w_i}{v_i}\right)_x\pa_{v_i}\bar{v}_i[(\pa_{v_i}\bar{v}_i)v_i-\bar{v}_i]B\tilde{r}_{i,vv}\\
	   &+\mu_i^{-1}(\tilde{\la}_i-\si_i)v_iv_{i,x}\xi_i^2\big[\pa_x(\pa_{v_i}\bar{v}_i)[(\pa_{v_i}\bar{v}_i)v_i-\bar{v}_i]+\pa_{v_i}\bar{v}_i\pa_x((\pa_{v_i}\bar{v}_i)v_i-\bar{v}_i)\big]
	   B\tilde{r}_{i,vv}\\
	  % &+\mu_i^{-1}(\tilde{\la}_i-\si_i)v_iv_{i,x}\xi_i^2\pa_{v_i}\bar{v}_i\pa_x((\pa_{v_i}\bar{v}_i)v_i-\bar{v}_i)B\tilde{r}_{i,vv}\\
	   &+\mu_i^{-1}(\tilde{\la}_i-\si_i)v_iv_{i,x}\xi_i^2\pa_{v_i}\bar{v}_i[(\pa_{v_i}\bar{v}_i)v_i-\bar{v}_i][u_x\cdot DB\tilde{r}_{i,vv}+B\pa_x(\tilde{r}_{i,vv})]\\
	   &=O(1)(\La_i^1+\La_i^4)+\widetilde{R}_\e^{7,i,14},
	   \end{align*}
	   with $\widetilde{R}_\e^{7,i,14}$ satisfying \eqref{12.11bis}.
	   We have now $\al_{i}^{7,15}=\sum\limits_{j\neq i}\mu_i^{-1}(\tilde{\la}_i-\si_i)v_iv_{j,x}\xi_i^2\pa_{v_j}\bar{v}_i[(\pa_{v_i}\bar{v}_i)v_i-\bar{v}_i]B\tilde{r}_{i,vv}=O(1)\La_i^1$ and
	   \begin{align*}
	   \al_{i,x}^{7,15}&=\sum\limits_{j\neq i}\big[u_x\cdot D\mu_i^{-1}(\tilde{\la}_i-\si_i)+\mu_i^{-1}\left(\tilde{\la}_{i,x}+\theta_i^\p\left(\frac{w_i}{v_i}\right)_x\right)\big]
	   v_iv_{j,x}\xi_i^2\pa_{v_j}\bar{v}_i[(\pa_{v_i}\bar{v}_i)v_i-\bar{v}_i]B\tilde{r}_{i,vv}\\
	  % &+\sum\limits_{j\neq i}\mu_i^{-1}\left(\tilde{\la}_{i,x}+\theta_i^\p\left(\frac{w_i}{v_i}\right)_x\right)v_iv_{j,x}\xi_i^2\pa_{v_j}\bar{v}_i[(\pa_{v_i}\bar{v}_i)v_i-\bar{v}_i]B\tilde{r}_{i,vv}\\
	   &+\sum\limits_{j\neq i}\mu_i^{-1}(\tilde{\la}_i-\si_i)\big[v_{i,x}v_{j,x}\xi_i^2+v_iv_{j,xx}\xi_i^2+
	   2v_iv_{j,x}\xi_i\xi_i^\p\left(\frac{w_i}{v_i}\right)_x
	  \big]\pa_{v_j}\bar{v}_i[(\pa_{v_i}\bar{v}_i)v_i-\bar{v}_i]B\tilde{r}_{i,vv}\\
	   %&+2\sum\limits_{j\neq i}\mu_i^{-1}(\tilde{\la}_i-\si_i)v_iv_{j,x}\xi_i\xi_i^\p\left(\frac{w_i}{v_i}\right)_x\pa_{v_j}\bar{v}_i[(\pa_{v_i}\bar{v}_i)v_i-\bar{v}_i]B\tilde{r}_{i,vv}\\
	   &+\sum\limits_{j\neq i}\mu_i^{-1}(\tilde{\la}_i-\si_i)v_iv_{j,x}\xi_i^2\big[\pa_x(\pa_{v_j}\bar{v}_i)[(\pa_{v_i}\bar{v}_i)v_i-\bar{v}_i]+\pa_{v_j}\bar{v}_i\pa_x((\pa_{v_i}\bar{v}_i)v_i-\bar{v}_i)\big]
	    B\tilde{r}_{i,vv}\\
	   %&+\sum\limits_{j\neq i}\mu_i^{-1}(\tilde{\la}_i-\si_i)v_iv_{j,x}\xi_i^2\pa_{v_j}\bar{v}_i\pa_x((\pa_{v_i}\bar{v}_i)v_i-\bar{v}_i)B\tilde{r}_{i,vv}\\
	   &+\sum\limits_{j\neq i}\mu_i^{-1}(\tilde{\la}_i-\si_i)v_iv_{j,x}\xi_i^2\pa_{v_j}\bar{v}_i[(\pa_{v_i}\bar{v}_i)v_i-\bar{v}_i][u_x\cdot DB\tilde{r}_{i,vv}+B\pa_x(\tilde{r}_{i,vv})]\\
	   &=O(1)(\La_i^1+\La_i^4).
	   \end{align*}
	   %We proceed similarly with $\al_{i}^{8,18}=\sum\limits_{j\neq i}\mu_i^{-1}(\tilde{\la}_i-\si_i)v_iw_{j,x}\xi_i^2\pa_{w_j}\bar{v}_i[(\pa_{v_i}\bar{v}_i)v_i-\bar{v}_i]B\tilde{r}_{i,vv}=O(1)\La_i^1$. 
	   Concerning $\al_{i}^{7,16}=\mu_i^{-1}(\tilde{\la}_i-\si_i)v_i\xi_i\xi_i^\p\bar{v}_i\left(\frac{w_i}{v_i}\right)_x[(\pa_{v_i}\bar{v}_i)v_i-\bar{v}_i]B\tilde{r}_{i,vv}=O(1)\La_i^4$ we get using  Lemmas \ref{estimimpo1}, \ref{lemme6.5}, \ref{lemme6.6}, \ref{lemme11.3} 
	   \begin{align*}
	  & \al_{i,x}^{7,16}=\big[u_x\cdot D\mu_i^{-1}(\tilde{\la}_i-\si_i)+\mu_i^{-1}\left(\tilde{\la}_{i,x}+\theta_i^\p\left(\frac{w_i}{v_i}\right)_x\right)\big]
	   v_i\xi_i\xi_i^\p\bar{v}_i\left(\frac{w_i}{v_i}\right)_x[(\pa_{v_i}\bar{v}_i)v_i-\bar{v}_i]B\tilde{r}_{i,vv}\\
	  % &+\mu_i^{-1}\left(\tilde{\la}_{i,x}+\theta_i^\p\left(\frac{w_i}{v_i}\right)_x\right)v_i\xi_i\xi_i^\p\bar{v}_i\left(\frac{w_i}{v_i}\right)_x[(\pa_{v_i}\bar{v}_i)v_i-\bar{v}_i]B\tilde{r}_{i,vv}\\
	   &+\mu_i^{-1}(\tilde{\la}_i-\si_i)\big[v_{i,x}\xi_i\xi_i^\p+v_i(\xi_i^\p)^2\left(\frac{w_i}{v_i}\right)_x
	   +v_i\xi_i\xi_i^{\p\p}\left(\frac{w_i}{v_i}\right)_x\big]
	   \bar{v}_i\left(\frac{w_i}{v_i}\right)_x[(\pa_{v_i}\bar{v}_i)v_i-\bar{v}_i]B\tilde{r}_{i,vv}\\
	   %&+\mu_i^{-1}(\tilde{\la}_i-\si_i)v_i(\xi_i^\p)^2\bar{v}_i\left(\frac{w_i}{v_i}\right)^2_x[(\pa_{v_i}\bar{v}_i)v_i-\bar{v}_i]B\tilde{r}_{i,vv}\\
	   %&+\mu_i^{-1}(\tilde{\la}_i-\si_i)v_i\xi_i\xi_i^{\p\p}\bar{v}_i\left(\frac{w_i}{v_i}\right)^2_x[(\pa_{v_i}\bar{v}_i)v_i-\bar{v}_i]B\tilde{r}_{i,vv}\\
	   &+\mu_i^{-1}(\tilde{\la}_i-\si_i)v_i\xi_i\xi_i^\p\big[\pa_x(\bar{v}_i)\left(\frac{w_i}{v_i}\right)_x
	   +\bar{v}_i\left(\frac{w_i}{v_i}\right)_{xx}\big]
	   [(\pa_{v_i}\bar{v}_i)v_i-\bar{v}_i]B\tilde{r}_{i,vv}\\
	   %&+\mu_i^{-1}(\tilde{\la}_i-\si_i)v_i\xi_i\xi_i^\p\bar{v}_i\left(\frac{w_i}{v_i}\right)_{xx}[(\pa_{v_i}\bar{v}_i)v_i-\bar{v}_i]B\tilde{r}_{i,vv}\\
	   &+\mu_i^{-1}(\tilde{\la}_i-\si_i)v_i\xi_i\xi_i^\p\bar{v}_i\left(\frac{w_i}{v_i}\right)_{x}\pa_x((\pa_{v_i}\bar{v}_i)v_i-\bar{v}_i)B\tilde{r}_{i,vv}\\
	   &+\mu_i^{-1}(\tilde{\la}_i-\si_i)v_i\xi_i\xi_i^\p\bar{v}_i\left(\frac{w_i}{v_i}\right)_{x}[(\pa_{v_i}\bar{v}_i)v_i-\bar{v}_i][u_x\cdot DB\tilde{r}_{i,vv}+B\pa_x(\tilde{r}_{i,vv})]\\
	   &=O(1)(\La_i^4+\La_i^5).
	   \end{align*}
	   To finish, we have $ \al_{i}^{7,17}=-\mu_i^{-1}(\tilde{\la}_i-\si_i)v_i\theta_i^\p\left(\frac{w_i}{v_i}\right)_x\xi_i[(\pa_{v_i}\bar{v}_i)v_i-\bar{v}_i]B\tilde{r}_{i,v\si}
	   =O(1)\La_i^4$ and using Lemma \ref{estimimpo1}, \ref{lemme6.5}, \ref{lemme6.6}, \ref{lemme11.3} and  \eqref{identity:wi-vi-xx}, \eqref{ngtech5}
	   \begin{align*}
	   \al_{i,x}^{7,17}&=-\big[u_x\cdot D\mu_i^{-1}(\tilde{\la}_i-\si_i)+\mu_i^{-1}\left(\tilde{\la}_{i,x}+\theta_i^\p\left(\frac{w_i}{v_i}\right)_x\right)\big]
	   v_i\theta_i^\p\left(\frac{w_i}{v_i}\right)_x\xi_i[(\pa_{v_i}\bar{v}_i)v_i-\bar{v}_i]B\tilde{r}_{i,v\si}\\
	   %&-\mu_i^{-1}\left(\tilde{\la}_{i,x}+\theta_i^\p\left(\frac{w_i}{v_i}\right)_x\right)v_i\theta_i^\p\left(\frac{w_i}{v_i}\right)_x\xi_i[(\pa_{v_i}\bar{v}_i)v_i-\bar{v}_i]B\tilde{r}_{i,v\si}\\
	   %&-\mu_i^{-1}(\tilde{\la}_i-\si_i)\big[v_{i,x}\theta_i^\p+v_i\theta_i^{\p\p}\left(\frac{w_i}{v_i}\right)^2_x\big]\left(\frac{w_i}{v_i}\right)_x\xi_i[(\pa_{v_i}\bar{v}_i)v_i-\bar{v}_i]B\tilde{r}_{i,v\si}\\
	   &-\mu_i^{-1}(\tilde{\la}_i-\si_i)v_i\theta_i^{\p}\big[\left(\frac{w_i}{v_i}\right)_{xx}\xi_i+\xi_i'\left(\frac{w_i}{v_i}\right)_{x}^2\big]
	   [(\pa_{v_i}\bar{v}_i)v_i-\bar{v}_i]B\tilde{r}_{i,v\si}\\
	   &-\mu_i^{-1}(\tilde{\la}_i-\si_i)v_i\theta_i^\p\left(\frac{w_i}{v_i}\right)_{xx}\xi_i[(\pa_{v_i}\bar{v}_i)v_i-\bar{v}_i]B\tilde{r}_{i,v\si}\\
	  % &-\mu_i^{-1}(\tilde{\la}_i-\si_i)v_i\theta_i^\p\left(\frac{w_i}{v_i}\right)^2_x\xi^\p_i[(\pa_{v_i}\bar{v}_i)v_i-\bar{v}_i]B\tilde{r}_{i,v\si}\\
	   &-\mu_i^{-1}(\tilde{\la}_i-\si_i)v_i\theta_i^\p\left(\frac{w_i}{v_i}\right)_x\xi_i\pa_x((\pa_{v_i}\bar{v}_i)v_i-\bar{v}_i)B\tilde{r}_{i,v\si}\\
	   &-\mu_i^{-1}(\tilde{\la}_i-\si_i)v_i\theta_i^\p\left(\frac{w_i}{v_i}\right)_x\xi_i[(\pa_{v_i}\bar{v}_i)v_i-\bar{v}_i][u_x\cdot DB\tilde{r}_{i,vv}+B\pa_x(\tilde{r}_{i,vv})]\\
	   &=O(1)\sum(\La_j^1+\delta_0^2\La_j^3+\La_j^4+\La_j^5).
	\end{align*}

\end{proof}

\begin{lemma}\label{lemmeal9.10}
	We have for $1\leq l\leq 4$, $1\leq l'\leq 2$
	\begin{align}
			&\al_{i}^{8,l}=O(1)\sum_j(\La_j^1+\La_j^2+\delta_0^2\La_j^3+\La_j^4+\La_j^5+\La_j^6+\La_j^{6,1})+R_\e^{8,i,l},\\
			&\al_{i}^{9,l'}=O(1)\sum_j(\La_j^1+\La_j^2+\delta_0^2\La_j^3+\La_j^4+\La_j^5+\La_j^6+\La_j^{6,1})+R_\e^{9,i,l'}\\
			&\al_{i,x}^{8,l}=O(1)\sum_j(\La_j^1+\La_j^2+\delta_0^2\La_j^3+\La_j^4+\La_j^5+\La_j^6+\La_j^{6,1})+\widetilde{R}_\e^{8,i,l},\\
			&\al_{i,x}^{9,l'}=O(1)\sum_j(\La_j^1+\La_j^2+\delta_0^2\La_j^3+\La_j^4+\La_j^5+\La_j^6+\La_j^{6,1})+\widetilde{R}_\e^{9,i,l'},
	\end{align}
	with:
	\begin{equation}
	\int_{\hat{t}}^T\int_{\R}(|R_\e^{8,i,l}|+|\widetilde{R}_\e^{8,i,l}|+|R_\e^{9,i,l'}|+|\widetilde{R}_\e^{9,i,l'}|)dx ds=O(1)\delta_0^2,
	\label{12.11bis1}
	\end{equation}
	for $\e>0$ small enough in terms of $T-\hat{t}$ and $\delta_0$.
	%and
	%\begin{equation}
	%	\al_{i}^{9,l}=\mathcal{O}(1)(\La_i^1+\La_i^3)\mbox{ and }	\al_{i,x}^{9,l}=\mathcal{O}(1)(\La_i^1+\La_i^3+\La_i^5)\mbox{ for }1\leq l\leq 5.
	%\end{equation}
	%and 
	%\begin{equation}
	%	\al_{i}^{10,l}=\mathcal{O}(1)(\La_i^1+\La_i^3)\mbox{ and }	\al_{i,x}^{10,l}=\mathcal{O}(1)(\La_i^1+\La_i^3+\La_i^5)\mbox{ for }l=1,2.
	%\end{equation}
\end{lemma}
\begin{proof}
A direct computation and \eqref{6.49bis} ensures that
\begin{align}
2v_i-\xi_i\bar{v}_i-\xi_iv_i\pa_{v_i}\bar{v}_i&=2v_i(1-\chi_i^\e\xi_i\eta_i)+\xi_i(\bar{v}_i-v_i\pa_{v_i}\bar{v}_i)\nonumber\\
&=2v_i(1-\chi_i^\e\xi_i\eta_i)+O(1)\xi_i (\chi')^\e_i\hat{v}_i+O(1)\xi_i  \sum_{l\ne i}\Delta^\e_{i,l}v_l^2.
%+\xi_i\sum_{j\ne i}(v_j^2+w_j^2)\eta'(\frac{v_j^2+w_j^2}{v_i})\prod_{l\ne j,i}\eta(\frac{v_l^2+w_l^2}{v_i})
\label{superimpll}
\end{align}
First we have using \eqref{6.45} and \eqref{superimpll}
\begin{align*}
\al_{i}^{8,1}&=(2v_i-\xi_i\bar{v}_i-\xi_iv_i\pa_{v_i}\bar{v}_i)v_{i,x}\left[B(u)\tilde{r}_{i,u}\tilde{r}_i+\tilde{r}_i\cdot DB(u)\tilde{r}_i\right]\\
%&=(2v_i-\xi_i\bar{v}_i-\xi_iv_i\pa_{v_i}\bar{v}_i)v_{i,x}\left[B(u)\tilde{r}_{i,u}\tilde{r}_i+\tilde{r}_i\cdot DB(u)\tilde{r}_i\right]\\
&=O(1)\mathbbm{1}_{\{v_i^{2N}\leq\e\}}v_iv_{i,x}+
O(1)v_i(1-\eta_i)v_{i,x}\mathbbm{1}_{\{|\frac{w_i}{v_i}|\leq\frac{\delta_1}{2}\}}+O(1)v_iv_{i,x}\mathbbm{1}_{\{|\frac{w_i}{v_i}|\geq\frac{\delta_1}{2}\}}+O(1)\La_i^1\\
&=O(1)\sum_j(\La_j^1+\La_j^4+\La_j^6)+O(1)\mathbbm{1}_{\{v_i^{2N}\leq\e\}}v_iv_{i,x}=O(1)\sum_j(\La_j^1+\La_j^4+\La_j^6)+R_\e^{8,i,1},
\end{align*}
with $R_\e^{8,i,1}$ satisfying \eqref{12.11bis1}.
We have now using% the fact that $\pa_{v_i}\bar{v}_i=(\chi')_i^\e\frac{v_i^{2N}}{\e}\eta_i+\chi_i^\e\eta_i+\chi_i^\e v_i\pa_{v_i}\eta_i$ and 
Lemma \ref{estimimpo1}
\begin{align}
&2-2\xi_i\pa_{v_i}\bar{v}_i-\xi_iv_i\pa_{v_iv_i}\bar{v}_i=2(1-\chi_i^\e \xi_i\eta_i)+O(1)\xi_i \left(\Delta^\ep_i+\varkappa^\ep_i\rho^\ep_i\right)%+\xi_i \alpha_iO(1)\sum_{k\ne i}(|\eta'(\frac{v_{k}^2+w_k^2}{v_i})|+|\eta''(\frac{v_{k}^2+w_k^2}{v_i})|).
\label{superimplll}
\end{align}
Using \eqref{superimpll}, \eqref{superimplll}, \eqref{ngtech5}  and Lemmas \ref{estimimpo1}, \ref{lemme6.5}, \ref{lemme6.6}, \ref{lemme11.3}, \ref{lemme11.4}   we get
	\begin{align*}
		\al_{i,x}^{8,1}&=(2-2\xi_i\pa_{v_i}\bar{v}_i-\xi_iv_i\pa_{v_iv_i}\bar{v}_i)v^2_{i,x}\left[B\tilde{r}_{i,u}\tilde{r}_i+\tilde{r}_i\cdot DB\tilde{r}_i\right]\\
		&-\sum\limits_{j\neq i}\xi_i(\pa_{v_j}\bar{v}_i+v_i\pa_{v_iv_j}\bar{v}_i)v_{j,x}v_{i,x}\left[B\tilde{r}_{i,u}\tilde{r}_i+\tilde{r}_i\cdot DB\tilde{r}_i\right]\\
		%&-\sum\limits_{j\neq i}\xi_i(\pa_{w_j}\bar{v}_i+v_i\pa_{v_iw_j}\bar{v}_i)w_{j,x}v_{i,x}\left[B(u)\tilde{r}_{i,u}\tilde{r}_i+\tilde{r}_i\cdot DB\tilde{r}_i\right]\\
		&-\xi_i^\p\left(\frac{w_i}{v_i}\right)_x(\bar{v}_i+v_i\pa_{v_i}\bar{v}_i)v_{i,x}\left[B\tilde{r}_{i,u}\tilde{r}_i+\tilde{r}_i\cdot DB\tilde{r}_i\right]\\
		&+(2v_i-\xi_i\bar{v}_i-\xi_iv_i\pa_{v_i}\bar{v}_i)v_{i,xx}\left[B\tilde{r}_{i,u}\tilde{r}_i+\tilde{r}_i\cdot DB\tilde{r}_i\right]\\
		&+(2v_i-\xi_i\bar{v}_i-\xi_iv_i\pa_{v_i}\bar{v}_i)v_{i,x}\left[u_x\cdot DB\tilde{r}_{i,u}\tilde{r}_i+u_x\otimes \tilde{r}_i:D^2B\tilde{r}_i\right]\\
		&+(2v_i-\xi_i\bar{v}_i-\xi_iv_i\pa_{v_i}\bar{v}_i)v_{i,x}\left[B\pa_x(\tilde{r}_{i,u}\tilde{r}_i)+\tilde{r}_{i,x}\cdot DB\tilde{r}_i+\tilde{r}_i\cdot DB\tilde{r}_{i,x}\right]\\
		&=O(1)\sum (\La_j^1+\La_j^4+\La_j^5+\La_j^6+\La_j^{6,1})+\widetilde{R}_\e^{8,i,1},
		\end{align*}
		with $\widetilde{R}_\e^{8,i,1}$ satisfying \eqref{12.11bis1}.
		Let us give some details on how to treat the terms $(2-2\xi_i\pa_{v_i}\bar{v}_i-\xi_iv_i\pa_{v_iv_i}\bar{v}_i)v^2_{i,x}\left[B\tilde{r}_{i,u}\tilde{r}_i+\tilde{r}_i\cdot DB\tilde{r}_i\right]$ and $(2v_i-\xi_i\bar{v}_i-\xi_iv_i\pa_{v_i}\bar{v}_i)v_{i,xx}\left[B\tilde{r}_{i,u}\tilde{r}_i+\tilde{r}_i\cdot DB\tilde{r}_i\right]$. First we have from Lemmas \ref{lemme6.5}, \ref{lemme6.6} and \eqref{superimplll}, \eqref{ngtech5}
		\begin{align*}
		&(2-2\xi_i\pa_{v_i}\bar{v}_i-\xi_iv_i\pa_{v_iv_i}\bar{v}_i)v^2_{i,x}\left[B\tilde{r}_{i,u}\tilde{r}_i+\tilde{r}_i\cdot DB\tilde{r}_i\right]=2(1-\chi_i^\e \xi_i\eta_i)v_{i,x}^2+O(1)\xi_i v_{i,x}^2\left(\Delta^\ep_i+\varkappa^\ep_i\rho^\ep_i\right)\\
		&=2(1-\chi_i^\e \eta_i)v_{i,x}^2\mathbbm{1}_{\{|\frac{w_i}{v_i}|\leq\frac{\delta_1}{2}\}}+2(1-\chi_i^\e \eta_i)v_{i,x}^2\mathbbm{1}_{\{|\frac{w_i}{v_i}|\geq\frac{\delta_1}{2}\}}+
		O(1)\xi_i v_i v_{i,x}\left(\Delta^\ep_i+\varkappa^\ep_i\rho^\ep_i\right)+O(1)\sum_j\La_j^4\\
		&=2(1-\chi_i^\e \eta_i)v_{i,x} v_i\mathbbm{1}_{\{|\frac{w_i}{v_i}|\leq\frac{\delta_1}{2}\}}
		+O(1)\sum_j(\La_j^4+\La_j^6)+
		O(1)\xi_i v_i v_{i,x}\varkappa^\ep_i\rho^\ep_i\\
		&=O(1)\mathbbm{1}_{\{v_i^{2N}
		\leq 2\e\}}v_iv_{i,x}+O(1)\sum_j(\La_j^1+\La_j^4+\La_j^6)+
		O(1)\xi_i v_i v_{i,x}\varkappa^\ep_i\rho^\ep_i\\
		&=O(1)\sum_j(\La_j^1+\La_j^4+\La_j^6)+\widetilde{R}_\e^{8,i,1,1}
		\end{align*}
		with $\widetilde{R}_\e^{8,i,1,1}$ satisfying \eqref{12.11bis1}.
		 From \eqref{estimate-v-i-xx-1aabis}, we have
		 \begin{align*}
		&v_{i,xx}=O(1)(|v_{i,x}|+|w_{i,x}|+|v_i|)+O(1)\sum\limits_{j\neq i}\left(\mu_j v_{j,x}-(\tilde{\la}_j-\la_j^*)v_j-w_j\right)_x \bar{v}_j\xi_j	\nonumber\\
		&+O(1)\sum_{j\ne i}\Big(|v_j|+|v_{j,x}|+|w_{j,x}|\Big)\left(|v_j|+|w_{j,x}|+|v_{j,x}|+\sum_{k} |v_k||v_j|\right)\nonumber\\
		&+O(1)\sum_k(\La_k^1+\La_k^4+\La_k^5+\La_k^6)+
		O(1)\big(\sum_k|v_kv_{k,x}|\eta_k(1-\chi^\e_k)+\sum\limits_{k}|v_k|^3\eta_k (1-\chi^\e_k)\big)%+O(1)\sum_k|v_kv_{k,x}||1-\chi_k^\e\xi_k\eta_k|
		\nonumber\\
		&{\color{black}{+O(1)\sum_k \varkappa^\ep_k\rho_k^\e\mathfrak{A}_k|v_k v_{k,x}|}}.
	\end{align*}
% \begin{align}
			%&\mu_i v_{i,xx}=-\mu_{i,x}v_{i,x}+\tilde{\la}_{i,x}v_i+(\tilde{\la}_i-\la_i^*)v_{i}+w_{i,x} +O(1)\sum\limits_{j\neq i}\left(\mu_j v_{j,x}-(\tilde{\la}_j-\la_j^*)v_j-w_j\right)_x v_j\xi_j
	%+O(1)\sum_{j\ne i}\sum_k (|v_{j,x}v_j v_k|+|v_j^2 v_k|)\\ 
	%\nonumber\\
	%&+O(1)\sum_{j\ne i}(|v_j|+|v_{j,x}|+|w_{j,x}|)(|v_j|+|w_{j,x}|+|v_{j,x}|+\sum_{k} |v_k||v_j|)\nonumber\\
	%&+O(1)\sum_k( \La_k^1+\La_k^4+\La_k^5)
		%&+\mathcal{O}(1)\left(\sum\limits_{j\neq i}|v_j|^2\right)\\
	%	+O(1)\sum_k|v_kv_{k,x}||1-\xi_k\eta_k|.\nonumber\\
	%	\end{align}
		It implies from Lemmas \ref{estimimpo1}, \ref{lemme6.5}, \ref{lemme6.6} that
		\begin{align*}
		&(2v_i-\xi_i\bar{v}_i-\xi_iv_i\pa_{v_i}\bar{v}_i)v_{i,xx}\left[B\tilde{r}_{i,u}\tilde{r}_i+\tilde{r}_i\cdot DB\tilde{r}_i\right]=O(1)\sum_k(\La_k^1+\La_k^4+\La_k^5+\La_k^6)\\
		&+O(1)(2v_i-\xi_i\bar{v}_i-\xi_iv_i\pa_{v_i}\bar{v}_i)(|v_{i,x}|+|v_i|+|w_{i,x}|)
		+
		O(1)\big(|v_i^2v_{i,x}|\eta_i(1-\chi^\e_i)+|v_i|^4\eta_i (1-\chi^\e_i)\big)%+O(1)\sum_k|v_kv_{k,x}||1-\chi_k^\e\xi_k\eta_k|
		\nonumber\\
		&{\color{black}{+O(1) \varkappa^\ep_i\rho_i^\e \mathfrak{A}_i|v_i^2 v_{i,x}|}}.
		\end{align*}
		Using now \eqref{superimpll} and \eqref{6.45}, we deduce that
		\begin{align*}
		&(2v_i-\xi_i\bar{v}_i-\xi_iv_i\pa_{v_i}\bar{v}_i)v_{i,xx}\left[B\tilde{r}_{i,u}\tilde{r}_i+\tilde{r}_i\cdot DB\tilde{r}_i\right]=O(1)\sum_k(\La_k^1+\La_k^4+\La_k^5+\La_k^6)+ \widetilde{R}_\e^{8,i,l,2}
		\\
		&+O(1)v_i(1-\chi_i^\e\xi_i\eta_i)(|v_{i,x}|+|v_i|+|w_{i,x}|)\\
		&=O(1)\sum_k(\La_k^1+\La_k^4+\La_k^5+\La_k^6)+ \widetilde{R}_\e^{8,i,1,2}
		+O(1)v_i(1-\chi_i^\e\xi_i\eta_i)(|v_{i,x}|+|v_i|+|w_{i,x}|)\mathbbm{1}_{\{v_i^{2N}\leq 2\e\}}\\
		&+O(1)v_i(1-\eta_i)(|v_{i,x}|+|v_i|+|w_{i,x}|)+O(1)\mathbbm{1}_{\{|\frac{w_i}{v_i}|\geq\frac{\delta_1}{2}\}}v_i(|v_{i,x}|+|v_i|+|w_{i,x}|)\\
		%&+O(1)(2v_i-\xi_i\bar{v}_i-\xi_iv_i\pa_{v_i}\bar{v}_i)(|v_{i,x}|+|v_i|+|w_{i,x}|)\\
		&=O(1)\sum_k( \La_k^1+\La_k^4+\La_k^5+\La_k^6+\La_k^{6,1})+ \widetilde{R}_\e^{8,i,1,3},
		\end{align*}
		with $\widetilde{R}_\e^{8,i,1,2}$, $ \widetilde{R}_\e^{8,i,1,3}$ satisfying \eqref{12.11bis1}.
		We have now $\al_{i}^{8,2}=-\xi^\p_iv_i\bar{v}_i\left(\frac{w_i}{v_i}\right)_x\left[B(u)\tilde{r}_{i,u}\tilde{r}_i+\tilde{r}_i\cdot DB(u)\tilde{r}_i\right]=O(1)\La_i^4$ and using Lemmas \ref{estimimpo1}, \ref{lemme6.5}, \ref{lemme6.6}, \ref{lemme11.3}, \ref{lemme11.4} and \eqref{ngtech5}, \eqref{identity:wi-vi-xx}
		\begin{align*}
		\al_{i,x}^{8,2}&=\big[-\xi^{\p\p}_i \left(\frac{w_i}{v_i}\right)_x v_i\bar{v}_i+\xi^\p_iv_{i,x}\bar{v}_i+\xi^\p_iv_i\pa_x(\bar{v}_i)\big]
		\left(\frac{w_i}{v_i}\right)_x\left[B(u)\tilde{r}_{i,u}\tilde{r}_i+\tilde{r}_i\cdot DB(u)\tilde{r}_i\right]\\
		%&-\xi^\p_iv_{i,x}\bar{v}_i\left(\frac{w_i}{v_i}\right)_x\left[B(u)\tilde{r}_{i,u}\tilde{r}_i+\tilde{r}_i\cdot DB(u)\tilde{r}_i\right]\\
		%&-\xi^\p_iv_i\pa_x(\bar{v}_i)\left(\frac{w_i}{v_i}\right)_x\left[B(u)\tilde{r}_{i,u}\tilde{r}_i+\tilde{r}_i\cdot DB(u)\tilde{r}_i\right]\\
		&-\xi^\p_iv_i\bar{v}_i\left(\frac{w_i}{v_i}\right)_{xx}\left[B(u)\tilde{r}_{i,u}\tilde{r}_i+\tilde{r}_i\cdot DB(u)\tilde{r}_i\right]\\
		&-\xi^\p_iv_i\bar{v}_i\left(\frac{w_i}{v_i}\right)_x\left[u_x\cdot DB(u)\tilde{r}_{i,u}\tilde{r}_i+u_x\otimes \tilde{r}_i:D^2B(u)\tilde{r}_i\right]\\
		&-\xi^\p_iv_i\bar{v}_i\left(\frac{w_i}{v_i}\right)_x\left[B(u)\pa_x(\tilde{r}_{i,u}\tilde{r}_i)+\tilde{r}_{i,x}\cdot DB(u)\tilde{r}_i+\tilde{r}_i\cdot DB(u)\tilde{r}_{i,x}\right]\\
		&=O(1)\sum_j(\La_j^1+\delta_0^2\La_j^3+\La_j^4+\La_j^5).
		\end{align*}
		We have now $\al_{i}^{8,3}=-\sum\limits_{j\neq i}\xi_iv_iv_{j,x}\pa_{v_j}\bar{v}_i\left[B(u)\tilde{r}_{i,u}\tilde{r}_i+\tilde{r}_i\cdot DB(u)\tilde{r}_i\right]=O(1)\La_i^1$ and from Lemmas \ref{estimimpo1}, \ref{lemme6.5}, \ref{lemme6.6}, \ref{lemme11.3}, \ref{lemme11.4}
		\begin{align*}
		\al_{i,x}^{8,3}&=-\sum\limits_{j\neq i}\big[\xi_i^\p \left(\frac{w_i}{v_i}\right)_xv_iv_{j,x}+\xi_iv_{i,x}v_{j,x}+\xi_iv_iv_{j,xx}\big]
		\pa_{v_j}\bar{v}_i\left[B(u)\tilde{r}_{i,u}\tilde{r}_i+\tilde{r}_i\cdot DB(u)\tilde{r}_i\right]\\
		%&-\sum\limits_{j\neq i}\xi_iv_{i,x}v_{j,x}\pa_{v_j}\bar{v}_i\left[B(u)\tilde{r}_{i,u}\tilde{r}_i+\tilde{r}_i\cdot DB(u)\tilde{r}_i\right]\\
		%&-\sum\limits_{j\neq i}\xi_iv_iv_{j,xx}\pa_{v_j}\bar{v}_i\left[B(u)\tilde{r}_{i,u}\tilde{r}_i+\tilde{r}_i\cdot DB(u)\tilde{r}_i\right]\\
		&-\sum\limits_{j\neq i}\xi_iv_iv_{j,x}\pa_x(\pa_{v_j}\bar{v}_i)\left[B(u)\tilde{r}_{i,u}\tilde{r}_i+\tilde{r}_i\cdot DB(u)\tilde{r}_i\right]\\
		&-\sum\limits_{j\neq i}\xi_iv_iv_{j,x}\pa_{v_j}\bar{v}_i\left[u_x\cdot DB(u)\tilde{r}_{i,u}\tilde{r}_i+u_x\otimes \tilde{r}_i: D^2B(u)\tilde{r}_i\right]\\
		&-\sum\limits_{j\neq i}\xi_iv_iv_{j,x}\pa_{v_j}\bar{v}_i\left[B(u)\pa_x(\tilde{r}_{i,u}\tilde{r}_i)+\tilde{r}_{i,x}\cdot DB(u)\tilde{r}_i+\tilde{r}_i\cdot DB(u)\tilde{r}_{i,x}\right]\\
		&=O(1)\La_i^1.
		\end{align*}
		%We can now deal with $\al_{i}^{9,4}=-\sum\limits_{j\neq i}\xi_iv_iw_{j,x}\pa_{w_j}\bar{v}_i\left[B(u)\tilde{r}_{i,u}\tilde{r}_i+\tilde{r}_i\cdot DB(u)\tilde{r}_i\right]=O(1)\La_i^1$ in a similar way.
		Finally we have  using Lemmas \ref{lemme11.3}, \ref{lemme11.4} and \eqref{6.45}
		\begin{align*}
		\al_{i}^{8,4}=&
		(1-\chi_i^\e \eta_i\xi_i)v^2_i\left[u_x\cdot DB\tilde{r}_{i,u}\tilde{r}_i+B\pa_x(\tilde{r}_{i,u}\tilde{r}_i)+u_x\otimes \tilde{r}_i:D^2B(u)\tilde{r}_i\right]\\
		&+
		(1-\chi_i^\e\eta_i\xi_i)v^2_i\left[\tilde{r}_i\cdot DB(u)\tilde{r}_{i,x}+\tilde{r}_{i,x}\cdot DB(u)\tilde{r}_{i}\right]\\
		=&O(1)v_i^2\mathbbm{1}_{\{v_i^{2N}\leq 2\e\}}+O(1)(1-\eta_i)v_i^2+O(1)v_i^2\mathbbm{1}_{\{|\frac{w_i}{v_i}|\geq\frac{\delta_1}{2}\}}\\
		=&O(1)\sum (\La_j^1+\La_j^4+\La_j^6)+R_\e^{8,i,4},
		\end{align*} 
		with $R_\e^{8,i,4}$ satisfying \eqref{12.11bis}.
		We recall now from Lemmas \ref{lemme11.6}, \ref{lemme6.5}, \ref{lemme6.6} that
		\begin{align}
		&\pa_{xx}(\tilde{r}_{i,u}\tilde{r}_i)
		=O(1)
		+O(1)\rho_i^\e\mathfrak{A}_i\biggl(
		|v_i\left(\frac{w_i}{v_i}\right)_{xx}|
		+  |v_i\left(\frac{w_i}{v_i}\right)^2_x|+|v_{i,x}\left(\frac{w_i}{v_i}\right)_x|\nonumber\\
		&+ \sum_{k\ne i}|\left(\frac{w_i}{v_i}\right)_x||v_kv_{k,x}|+\frac{v_{i,x}^2}{|v_i|}+\sum_{k\ne i}|\frac{v_{k,x}v_kv_{i,x}}{v_i}|\biggl).\label{suyptech1}
	\end{align}
		\begin{align}
		&\tilde{r}_{i,xx}=O(1)+O(1)\rho_i^\e\mathfrak{A}_i\bigg( |v_{i,x}\left(\frac{w_i}{v_i}\right)_x|+ \frac{v_{i,x}^2}{|v_i|}+\sum_{l\ne i}|\frac{v_lv_{i,x}v_{l,x}}{v_i}|+| \left(\frac{w_i}{v_i}\right)_x|\sum\limits_{j\neq i}|v_{j,x}v_j|\nonumber\\
		&+\left(\frac{w_i}{v_i}\right)_x^2|v_i|	+| \left(\frac{w_i}{v_i}\right)_{xx} v_i|\biggl).	\label{suyptech2}
		\end{align}
		Using \eqref{suyptech1}, \eqref{suyptech1}, \eqref{6.45} and  Lemmas \ref{estimimpo1}, \ref{lemme6.5}, \ref{lemme6.6}, \ref{lemme11.3}, \ref{lemme11.4}
		\begin{align*}
		\al_{i,x}^{8,4}=&-{\color{black}{\chi_i^\e}}\pa_{v_i}\eta_i\xi_iv_{i,x}v^2_i\left[u_x\cdot DB\tilde{r}_{i,u}\tilde{r}_i+B\pa_x(\tilde{r}_{i,u}\tilde{r}_i)+u_x\otimes \tilde{r}_i:D^2B(u)\tilde{r}_i\right]\\
		&-{\color{black}{\chi_i^\e}}\pa_{v_i}\eta_i\xi_iv_{i,x}v^2_i\left[\tilde{r}_i\cdot DB(u)\tilde{r}_{i,x}+\tilde{r}_{i,x}\cdot DB(u)\tilde{r}_{i}\right]\\
		&-\sum\limits_{j\neq i}{\color{black}{\chi_i^\e}}\pa_{v_j}\eta_i\xi_iv_{j,x}v^2_i\left[u_x\cdot DB\tilde{r}_{i,u}\tilde{r}_i+B\pa_x(\tilde{r}_{i,u}\tilde{r}_i)+u_x\otimes \tilde{r}_i:D^2B(u)\tilde{r}_i\right]\\
		&-\sum\limits_{j\neq i}{\color{black}{\chi_i^\e}}\pa_{v_j}\eta_i\xi_iv_{j,x}v^2_i\left[\tilde{r}_i\cdot DB(u)\tilde{r}_{i,x}+\tilde{r}_{i,x}\cdot DB(u)\tilde{r}_{i}\right]\\
		%&-\sum\limits_{j\neq i}\pa_{w_j}\eta_i\xi_iw_{j,x}v^2_i\left[u_x\cdot DB\tilde{r}_{i,u}\tilde{r}_i+B\pa_x(\tilde{r}_{i,u}\tilde{r}_i)+u_x\otimes \tilde{r}_i:D^2B(u)\tilde{r}_i\right]\\
		%&-\sum\limits_{j\neq i}\pa_{w_j}\eta_i\xi_iw_{j,x}v^2_i\left[\tilde{r}_i\cdot DB(u)\tilde{r}_{i,x}+\tilde{r}_{i,x}\cdot DB(u)\tilde{r}_{i}\right]\\
		&-{\color{black}{\chi_i^\e}}\eta_i\xi_i^\p\left(\frac{w_i}{v_i}\right)_xv_i^2\left[u_x\cdot DB\tilde{r}_{i,u}\tilde{r}_i+B\pa_x(\tilde{r}_{i,u}\tilde{r}_i)+u_x\otimes \tilde{r}_i:D^2B(u)\tilde{r}_i\right]\\
		&-{\color{black}{\chi_i^\e}}\eta_i\xi_i^\p\left(\frac{w_i}{v_i}\right)_xv_i^2\left[\tilde{r}_i\cdot DB(u)\tilde{r}_{i,x}+\tilde{r}_{i,x}\cdot DB(u)\tilde{r}_{i}\right]\\
		&{\color{black}{-2N (\chi')_i^\e\frac{v_i^{2N+1}}{\e}\xi_i\eta_i \left[u_x\cdot DB\tilde{r}_{i,u}\tilde{r}_i+B\pa_x(\tilde{r}_{i,u}\tilde{r}_i)+u_x\otimes \tilde{r}_i:D^2B(u)\tilde{r}_i\right]}}\\
		&{\color{black}{-2N (\chi')_i^\e\frac{v_i^{2N+1}}{\e}\xi_i\eta_i \left[\tilde{r}_i\cdot DB(u)\tilde{r}_{i,x}+\tilde{r}_{i,x}\cdot DB(u)\tilde{r}_{i}\right]}}\\
		&+2(1-{\color{black}{\chi_i^\e}}\eta_i\xi_i)v_iv_{i,x}\left[u_x\cdot DB\tilde{r}_{i,u}\tilde{r}_i+B\pa_x(\tilde{r}_{i,u}\tilde{r}_i)+u_x\otimes \tilde{r}_i:D^2B(u)\tilde{r}_i\right]\\
		&+2(1-{\color{black}{\chi_i^\e}}\eta_i\xi_i)v_iv_{i,x}\left[\tilde{r}_i\cdot DB(u)\tilde{r}_{i,x}+\tilde{r}_{i,x}\cdot DB(u)\tilde{r}_{i}\right]\\
		&+(1-{\color{black}{\chi_i^\e}}\eta_i\xi_i)v_i^2\left[u_{xx}\cdot DB\tilde{r}_{i,u}\tilde{r}_i+u_x\otimes u_x: D^2B\tilde{r}_{i,u}\tilde{r}_i+2u_x\cdot DB\pa_x(\tilde{r}_{i,u}\tilde{r}_i)\right]\\
		&+(1-{\color{black}{\chi_i^\e}}\eta_i\xi_i)v_i^2\left[B\pa_{xx}(\tilde{r}_{i,u}\tilde{r}_i)+u_x\otimes u_x\otimes \tilde{r}_i:D^3B(u)\tilde{r}_i+u_{xx}\otimes \tilde{r}_i: D^2B(u)\tilde{r}_i\right]\\
		&+2(1-{\color{black}{\chi_i^\e}}\eta_i\xi_i)v_i^2\left[u_{x}\otimes \tilde{r}_{i,x}: D^2B(u)\tilde{r}_i+u_{x}\otimes \tilde{r}_{i}: D^2B(u)\tilde{r}_{i,x}\right]\\
		&+(1-{\color{black}{\chi_i^\e}}\eta_i\xi_i)v_i^2\left[\tilde{r}_{i,xx}\cdot DB(u)\tilde{r}_{i}+2\tilde{r}_{i,x}\cdot DB(u)\tilde{r}_{i,x}+\tilde{r}_i\cdot DB(u)\tilde{r}_{i,xx}\right]\\
		&=O(1)\sum_j(\La_j^1+\delta_0^2\La_j^3+\La_j^4+\La_j^5+\La_j^6)+\widetilde{R}_\e^{8,i,4},
	\end{align*}
	with $\widetilde{R}_\e^{8,i,4},$ satisfying \eqref{12.11bis}.
	Finally, for the term $\al_i^{9,1}$ and $\al_i^{9,2}$ we have  using Lemmas  \ref{estimimpo1}, \ref{lemme6.5}, \ref{lemme6.6}, \ref{lemme11.3}, \ref{lemme11.4}, \ref{lemme11.6}%\ref{lemma:derivative-r-k},  \ref{lemme10.5}. Furthermore we get
	%$\al_{i}^{10,1}=\sum\limits_{j\neq i}(v_{i,x}v_j+v_iv_{j,x})\left[\tilde{r}_i\cdot DB(u)\tilde{r}_j+B(u)\tilde{r}_{i,u}\tilde{r}_j\right]=O(1)\La_i^1$ and $\al_i^{10,2}=\sum\limits_{j\neq i}v_iv_j\left[\tilde{r}_i\cdot DB(u)\tilde{r}_j+B(u)\tilde{r}_{i,u}\tilde{r}_j\right]_x=O(1)\La_i^1$ (FAIRE $\alpha_{i,x}^{10,2}$)
	\begin{align*}
	\al_{i,x}^{9,1}&=\sum\limits_{j\neq i}(v_{i,x}v_j+v_iv_{j,x})\left[\tilde{r}_i\cdot DB(u)\tilde{r}_j+B(u)\tilde{r}_{i,u}\tilde{r}_j\right]=O(1)\La_i^1\\
		\al_{i,x}^{9,1}&=\sum\limits_{j\neq i}(v_{i,xx}v_j+2v_{i,x}v_{j,x}+v_iv_{j,xx})\left[\tilde{r}_i\cdot DB(u)\tilde{r}_j+B(u)\tilde{r}_{i,u}\tilde{r}_j\right]\\
		&+\sum\limits_{j\neq i}(v_{i,x}v_j+v_iv_{j,x})\left[\tilde{r}_i\cdot DB(u)\tilde{r}_j+B(u)\tilde{r}_{i,u}\tilde{r}_j\right]_{x}=O(1)\La_i^1,\\
		\al_{i}^{9,2}&=\sum\limits_{j\neq i}v_iv_j\left[\tilde{r}_i\cdot DB(u)\tilde{r}_j+B(u)\tilde{r}_{i,u}\tilde{r}_j\right]_x=O(1)\La_i^1,\\
		\al_{i,x}^{9,2}&= \sum\limits_{j\neq i}(v_{i,x}v_j+v_{i}v_{j,x})\left[\tilde{r}_i\cdot DB(u)\tilde{r}_j+B(u)\tilde{r}_{i,u}\tilde{r}_j\right]_x\\
		&+\sum\limits_{j\neq i}v_iv_j\left[\tilde{r}_i\cdot DB(u)\tilde{r}_j+B(u)\tilde{r}_{i,u}\tilde{r}_j\right]_{xx}\\
		&=O(1)\La_i^1.
	\end{align*}
This completes the proof of Lemma \ref{lemmeal9.10}.
\end{proof}

\begin{lemma}\label{lemma:estimate-7-1a}
	For $l=1,2,3$ it follows,
	\begin{align}
			&\al^{6,1,l}_{i}=O(1)\sum_j(\La_j^1+\La_j^5)\mbox{ and }\al^{6,1,l}_{i,x}=O(1)\sum_j(\La_j^1+\La_j^2+\delta_0^2\La_j^3+\La_j^4+\La_j^5),\\
	&\al_i^{10}=O(1)\La_i^1\mbox{ and }\al^{10}_{i,x}=O(1)\sum_j \La_j^1,\\
	&\al_i^{11}=O(1)\sum_j(\La_j^1+\La_j^5)\mbox{ and }\al^{11}_{i,x}=O(1)(\La_j^1+\La_j^3+\La_j^4+\La_j^5).
	\end{align}
\end{lemma}
\begin{proof}
	First we observe that% from \eqref{def:b-ij} and \eqref{identity:wi-vi-xx} it follows
	\begin{align}
		&b_{ij}(w_{j,x}-(w_j/v_j)v_{j,x})_{xx}\nonumber\\
		&=b_{ij}\left(w_{j,xxx}-\frac{w_j}{v_j}v_{j,xxx}\right)-2b_{ij}v_{j,xx}\left(\frac{w_j}{v_j}\right)_{x}-b_{ij}\left(\frac{w_j}{v_j}\right)_{xx}v_{j,x}.
		\label{1superimpl}
		%&=\mathcal{O}(1)(\La_j^1+\La_j^3+\La_j^4+\La_j^5).
	\end{align}
%Then we can estimate
%\begin{align*}
%	\al^{7,1,1}_i&=-\left(\sum\limits_{j\neq i}(\mu_i-\mu_j)b_{ij}(w_{j,x}-(w_j/v_j)v_{j,x})_x\right)\sum\limits_{k\neq i}\psi_{ik} r_k(u),\\
	%\al^{7,1,2}_i&=-\left(\sum\limits_{j\neq i}(\mu_i-\mu_j)b_{ij}(w_{j,x}-(w_j/v_j)v_{j,x})_x\right)\sum\limits_{k\neq i}(v_i\xi_i\pa_{v_i}\bar{v}_i-(w_i/v_i)\xi_i^\p\bar{v}_i)\psi_{ik,v} r_k,\\
%	\al^{7,1,3}_i&=-\left(\sum\limits_{j\neq i}(\mu_i-\mu_j)b_{ij}(w_{j,x}-(w_j/v_j)v_{j,x})_x\right)\sum\limits_{k\neq i}(w_i/v_i)\theta^\p_i\psi_{ik,\si} r_k.
%\end{align*}
From \eqref{def:b-ij}  and the fact that $\psi_{ji,\sig}=O(1)\xi_j\bar{v}_j$ we have $b_{ij}=O(1)\bar{v}_j \mathfrak{A}_j$ it implies using Lemmas \ref{lemme6.5}, \ref{lemme6.6} and the fact that $\psi_{ik}=O(1)\mathfrak{A}_i \bar{v}_i$
 \begin{align*}
 &\al^{6,1,1}_{i}=-\left(\sum\limits_{j\neq i}(\mu_i-\mu_j)b_{ij}(w_{j,x}-(w_j/v_j)v_{j,x})_x\right)\sum\limits_{k\neq i}\psi_{ik} r_k(u)=O(1)\sum_j(\La_j^1+\La_j^5). 
 \end{align*}
 In addition we get using \eqref{1superimpl}, \eqref{identity:wi-vi-xx}, \eqref{ngtech5} and Lemmas \ref{estimimpo1}, \ref{lemme6.5}, \ref{lemme6.6}, \ref{lemme11.3a}, \ref{lemme9.6}
 	\begin{align*}
\al^{6,1,1}_{i,x}&=-\left(\sum\limits_{j\neq i}\big[u_x\cdot D(\mu_i-\mu_j)b_{ij}+(\mu_i-\mu_j)b_{ij,x}\big]
(w_{j,x}-(w_j/v_j)v_{j,x})_x\right)\sum\limits_{k\neq i}\psi_{ik} r_k(u)\\
		%&-\left(\sum\limits_{j\neq i}(\mu_i-\mu_j)b_{ij,x}(w_{j,x}-(w_j/v_j)v_{j,x})_x\right)\sum\limits_{k\neq i}\psi_{ik} r_k(u)\\
	    &-\left(\sum\limits_{j\neq i}(\mu_i-\mu_j)b_{ij}(w_{j,x}-(w_j/v_j)v_{j,x})_{xx}\right)\sum\limits_{k\neq i}\psi_{ik} r_k(u)\\
		&-\left(\sum\limits_{j\neq i}(\mu_i-\mu_j)b_{ij}(w_{j,x}-(w_j/v_j)v_{j,x})_x\right)\sum\limits_{k\neq i}\big[\pa_x(\psi_{ik}) r_k(u)+\psi_{ik} u_x\cdot Dr_k(u)\big]\\
		%&-\left(\sum\limits_{j\neq i}(\mu_i-\mu_j)b_{ij}(w_{j,x}-(w_j/v_j)v_{j,x})_x\right)\sum\limits_{k\neq i}\psi_{ik} u_x\cdot Dr_k(u),\\
		&=O(1)\sum_j(\La_j^1+\La_j^2+\delta_0^2\La_j^3+\La_j^4+\La_j^5).
		\end{align*}
		Next we have 
		\begin{align*}
		\al^{6,1,2}_i&=-\left(\sum\limits_{j\neq i}(\mu_i-\mu_j)b_{ij}(w_{j,x}-(w_j/v_j)v_{j,x})_x\right)\sum\limits_{k\neq i}(v_i\xi_i\pa_{v_i}\bar{v}_i-(w_i/v_i)\xi_i^\p\bar{v}_i)\psi_{ik,v} r_k\\
		%&-\left(\sum\limits_{j\neq i}(\mu_i-\mu_j)b_{ij}(w_{j,x}-(w_j/v_j)v_{j,x})_x\right)\sum\limits_{k\neq i}(v_i\xi_i\pa_{v_i}\bar{v}_i-(w_i/v_i)\xi_i^\p\bar{v}_i)\psi_{ik,v} r_k\\
		&=O(1)\sum_j(\La_j^1+\La_j^5).
		\end{align*}
		Furthermore using \eqref{1superimpl}, \eqref{identity:wi-vi-xx}, \eqref{ngtech5} and Lemmas \ref{estimimpo1}, \ref{lemme6.5}, \ref{lemme6.6}, \ref{lemme11.3a}, \ref{lemme9.6} we obtain
		\begin{align*}
		\al^{6,1,2}_{i,x}&=-\left(\sum\limits_{j\neq i}u_x\cdot D(\mu_i-\mu_j)b_{ij}
		(w_{j,x}-(w_j/v_j)v_{j,x})_x\right)\sum\limits_{k\neq i}(v_i\xi_i\pa_{v_i}\bar{v}_i-(w_i/v_i)\xi_i^\p\bar{v}_i)\psi_{ik,v} r_k\\
		&-\left(\sum\limits_{j\neq i}(\mu_i-\mu_j)b_{ij,x}(w_{j,x}-(w_j/v_j)v_{j,x})_x\right)\sum\limits_{k\neq i}(v_i\xi_i\pa_{v_i}\bar{v}_i-(w_i/v_i)\xi_i^\p\bar{v}_i)\psi_{ik,v} r_k\\
		%&-\left(\sum\limits_{j\neq i}(\mu_i-\mu_j)b_{ij}(w_{j,x}-(w_j/v_j)v_{j,x})_x\right)\sum\limits_{k\neq i}(v_i\xi_i\pa_{v_i}\bar{v}_i-(w_i/v_i)\xi_i^\p\bar{v}_i)\psi_{ik,v} r_k\\
		&-\left(\sum\limits_{j\neq i}(\mu_i-\mu_j)b_{ij}(w_{j,x}-(w_j/v_j)v_{j,x})_{xx}\right)\sum\limits_{k\neq i}(v_i\xi_i\pa_{v_i}\bar{v}_i-(w_i/v_i)\xi_i^\p\bar{v}_i)\psi_{ik,v} r_k\\
		&-\left(\sum\limits_{j\neq i}(\mu_i-\mu_j)b_{ij}(w_{j,x}-(w_j/v_j)v_{j,x})_x\right)\sum\limits_{k\neq i}\left(v_{i,x}\xi_i\pa_{v_i}\bar{v}_i-\left(\frac{w_i}{v_i}\right)_x\xi_i^\p\bar{v}_i\right)\psi_{ik,v} r_k\\
		&-\left(\sum\limits_{j\neq i}(\mu_i-\mu_j)b_{ij}(w_{j,x}-(w_j/v_j)v_{j,x})_x\right)\sum\limits_{k\neq i}\left(\frac{w_i}{v_i}\right)_x(v_i\xi_i^\p\pa_{v_i}\bar{v}_i-(w_i/v_i)\xi_i^{\p\p}\bar{v}_i)\psi_{ik,v} r_k\\
		&-\left(\sum\limits_{j\neq i}(\mu_i-\mu_j)b_{ij}(w_{j,x}-(w_j/v_j)v_{j,x})_x\right)\sum\limits_{k\neq i}(v_i\xi_i\pa_x(\pa_{v_i}\bar{v}_i)-(w_i/v_i)\xi_i^\p\pa_x(\bar{v}_i))\psi_{ik,v} r_k\\
		&-\left(\sum\limits_{j\neq i}(\mu_i-\mu_j)b_{ij}(w_{j,x}-(w_j/v_j)v_{j,x})_x\right)\sum\limits_{k\neq i}(v_i\xi_i\pa_{v_i}\bar{v}_i-(w_i/v_i)\xi_i^\p\bar{v}_i)\pa_x(\psi_{ik,v}) r_k\\
		&-\left(\sum\limits_{j\neq i}(\mu_i-\mu_j)b_{ij}(w_{j,x}-(w_j/v_j)v_{j,x})_x\right)\sum\limits_{k\neq i}(v_i\xi_i\pa_{v_i}\bar{v}_i-(w_i/v_i)\xi_i^\p\bar{v}_i) \psi_{ik,v}u_x\cdot Dr_k\\
		&=O(1)\sum_j(\La_j^1+\La_j^2+\delta_0^2\La_j^3+\La_j^4+\La_j^5).
		\end{align*}
		Next we have $\al^{6,1,3}_i=-\left(\sum\limits_{j\neq i}(\mu_i-\mu_j)b_{ij}(w_{j,x}-(w_j/v_j)v_{j,x})_x\right)\sum\limits_{k\neq i}(w_i/v_i)\theta^\p_i\psi_{ik,\si} r_kO(1)\sum_j (\La_j^1+\La_j^5)$ because $\psi_{ik,\si}=O(1)\mathfrak{A}_k\bar{v}_k$. As previously we get
		\begin{align*}
		\al^{6,1,3}_{i,x}&=-\left(\sum\limits_{j\neq i}u_x\cdot D(\mu_i-\mu_j)b_{ij}(w_{j,x}-(w_j/v_j)v_{j,x})_x\right)\sum\limits_{k\neq i}(w_i/v_i)\theta^\p_i\psi_{ik,\si} r_k\\
		&-\left(\sum\limits_{j\neq i}(\mu_i-\mu_j)b_{ij,x}(w_{j,x}-(w_j/v_j)v_{j,x})_x\right)\sum\limits_{k\neq i}(w_i/v_i)\theta^\p_i\psi_{ik,\si} r_k\\
		&-\left(\sum\limits_{j\neq i}(\mu_i-\mu_j)b_{ij}(w_{j,x}-(w_j/v_j)v_{j,x})_{xx}\right)\sum\limits_{k\neq i}(w_i/v_i)\theta^\p_i\psi_{ik,\si} r_k\\
		&-\left(\sum\limits_{j\neq i}(\mu_i-\mu_j)b_{ij}(w_{j,x}-(w_j/v_j)v_{j,x})_x\right)\sum\limits_{k\neq i}\left(\frac{w_i}{v_i}\right)_x\theta^\p_i\psi_{ik,\si} r_k\\
		&-\left(\sum\limits_{j\neq i}(\mu_i-\mu_j)b_{ij}(w_{j,x}-(w_j/v_j)v_{j,x})_x\right)\sum\limits_{k\neq i}(w_i/v_i)\theta^{\p\p}_i\left(\frac{w_i}{v_i}\right)_x\psi_{ik,\si} r_k\\
		&-\left(\sum\limits_{j\neq i}(\mu_i-\mu_j)b_{ij}(w_{j,x}-(w_j/v_j)v_{j,x})_x\right)\sum\limits_{k\neq i}(w_i/v_i)\theta^\p_i\pa_x(\psi_{ik,\si}) r_k\\
		&-\left(\sum\limits_{j\neq i}(\mu_i-\mu_j)b_{ij}(w_{j,x}-(w_j/v_j)v_{j,x})_x\right)\sum\limits_{k\neq i}(w_i/v_i)\theta^\p_i(\psi_{ik,\si}) u_x\cdot Dr_k\\
		&=O(1)\sum_j(\La_j^1+\La_j^2+\delta_0^2\La_j^3+\La_j^4+\La_j^5).
	\end{align*}
	We have now using the Lemmas \ref{estimimpo1} and  \ref{lemme6.5} $\alpha_{i}^{10}=\sum\limits_{j\neq i}v_{i,t}v_j\xi_j\pa_{v_i}\bar{v_j}\tilde{r}_{j,v}	=O(1)\La_i^1$. Furthermore we get using Lemmas \ref{estimimpo1}, \ref{lemme6.5}, \ref{lemme6.6}, \ref{lemme11.3}
\begin{align*}
	\alpha_{i,x}^{10}&=\sum_{j\ne i}(v_{i,tx}v_j+v_{i,t}v_{j,x})\xi_j\pa_{v_i}\bar{v}_j\tilde{r}_{j,v}+\sum_{j\ne i}v_{i,t}v_j\xi_j'\left(\frac{w_j}{v_j}\right)_x \pa_{v_i}\bar{v}_j\tilde{r}_{j,v}\\
	&+\sum\limits_{j\ne i}v_{i,t}v_j\xi_j\big[\pa_x(\pa_{v_i}\bar{v}_j)\tilde{r}_{j,v}+\pa_{v_i}\bar{v}_j\pa_x(\tilde{r}_{j,v})\big]
	%&+\sum\limits_{j\ne i}v_{i,t}v_j\xi_j(1-\xi_i)\pa_{v_i}\bar{v}_j\pa_x(\tilde{r}_{j,v})\\
	=O(1)\sum_k \La_k^1.
\end{align*}
%We proceed similarly with $\alpha_{i}^{11,2}=\sum\limits_{j\neq i} w_{i,t} \xi_j(1-\xi_i)v_j\pa_{w_i}\bar{v}_j\tilde{r}_{j,v}=O(1)\La_i^1$. 
We recall here that
	$\hat{b}_{ij}=\psi_{ji} +\left[(\la_i^*-\la_j^*)+\frac{w_j}{v_j}\right]b_{ij}$ from  \eqref{def:hat-b-ij}.
Let us deal now with $\al_i^{11}$ which satisfies using the fact that $\tilde{r}_{i,\sig}=O(1)\xi_i\bar{v}_i$ and $\psi_{ji}, \psi_{ji,\sig}=O(1)\xi_j\bar{v}_j$ 
	\begin{align*}
	&\al_i^{11}=-\left(\sum\limits_{j\neq i}(\mu_i-\mu_j)\hat{b}_{ij}(w_{j,x}-(w_j/v_j)v_{j,x})_x\right)(\xi_i^\p\bar{v}_i\tilde{r}_{i,v}-\theta^\p_i\tilde{r}_{i,\si})=O(1)\sum_k(\La_k^1+\La_k^5).
\end{align*}
We get now from \eqref{1superimpl}, \eqref{ngtech5}, Lemmas \ref{estimimpo1}, \ref{lemme6.5}, \ref{lemme6.6}, \ref{lemme11.3},  \ref{lemme9.6} 
\begin{align*}
	\al_{i,x}^{11}&=-\left(\sum\limits_{j\neq i}u_x\cdot D(\mu_i-\mu_j)\hat{b}_{ij}(w_{j,x}-(w_j/v_j)v_{j,x})_x\right)(\xi_i^\p\bar{v}_i\tilde{r}_{i,v}-\theta^\p_i\tilde{r}_{i,\si})\\
	&-\left(\sum\limits_{j\neq i}(\mu_i-\mu_j)\hat{b}_{ij,x}(w_{j,x}-(w_j/v_j)v_{j,x})_x\right)(\xi_i^\p\bar{v}_i\tilde{r}_{i,v}-\theta^\p_i\tilde{r}_{i,\si})\\
	&-\left(\sum\limits_{j\neq i}(\mu_i-\mu_j)\hat{b}_{ij}(w_{j,x}-(w_j/v_j)v_{j,x})_{xx}\right)(\xi_i^\p\bar{v}_i\tilde{r}_{i,v}-\theta^\p_i\tilde{r}_{i,\si})\\
	&-\left(\sum\limits_{j\neq i}(\mu_i-\mu_j)\hat{b}_{ij}(w_{j,x}-(w_j/v_j)v_{j,x})_x\right)(\xi_i''\bar{v}_i\tilde{r}_{i,v}-\theta_i''\tilde{r}_{i,\si})\left(\frac{w_i}{v_i}\right)_x\\
	&-\left(\sum\limits_{j\neq i}(\mu_i-\mu_j)\hat{b}_{ij}(w_{j,x}-(w_j/v_j)v_{j,x})_x\right)(\xi_i^\p\pa_x(\bar{v}_i)\tilde{r}_{i,v}+\xi_i^\p\bar{v}_i\pa_x(\tilde{r}_{i,v})-\theta^\p_i\pa_x(\tilde{r}_{i,\si}))\\
	&=\sum_k(\La_k^1+\La_k^2+\delta_0^2\La_k^3+\La_k^4+\La_k^5).
\end{align*}
This completes the proof of Lemma \ref{lemma:estimate-7-1a}.
\end{proof}

%--------------------------------------------------------------------------------------------------------------

\subsection{Second order estimates for $K_i$'s}\label{section-psi-i-x}
\begin{lemma}
	For $1\leq i\leq n$ and $1\leq l\leq 33$ it follows,
	\begin{equation}
	\begin{aligned}
		&\B_i^{1,l}=O(1)\sum_j(\La_j^1+\La_j^2+\delta_0^2\La_j^3+\La_j^4+\La_j^5+\La_j^6+\La_j^{6,1})+R_{\e,1}^{1,i,l}\\
		& \B_{i,x}^{1,l}=O(1)\sum_j(\La_j^1+\La_j^2+\delta_0^2\La_j^3+\La_j^4+\La_j^5+\La_j^6+\La_j^{6,1})+\widetilde{R}_{\e,1}^{1,i,l},
		\end{aligned}
	\end{equation}
	with:
	\begin{equation}
	\int_{\hat{t}}^T\int_{\R}(|R_{\e,1}^{1,i,l}|+|\widetilde{R}_{\e,1}^{1,i,l}|)dx ds=O(1)\delta_0^2,
	\label{12.11bis2}
	\end{equation}
	for $\e>0$ small enough in terms of $T-\hat{t}$ and $\delta_0$.
\end{lemma}
\begin{proof}
\begin{align*}
	&K_{1,x}+K_{2,x}+K_{5,x}-\sum\limits_{i}(w_i-\la_i^*v_i)w_i\tilde{r}_{i,u}\tilde{r}_i-\sum\limits_{i}\sum\limits_{j\neq i}(w_j-\la_j^*v_j)w_i\tilde{r}_{i,u}\tilde{r}_j\\
	&=:\sum\limits_{i}(\mu_iw_{i,x}-\tilde{\la}_iw_i)_x\left[\tilde{r}_{i}+\frac{w_i}{v_i}\xi_i^\p\bar{v}_i\tilde{r}_{i,v}-\frac{w_i}{v_i}\theta_i^\p\tilde{r}_{i,\si}\right]\\
	&+\sum\limits_{i}(\mu_iv_{i,x}-\tilde{\la}_iv_i)_x\left[\left(w_i\xi_i\pa_{v_i}\bar{v}_i-\xi^\p_i\bar{v}_i\frac{w_i^2}{v_i^2}\right)\tilde{r}_{i,v}+\frac{w^2_i}{v^2_i}\theta_i^\p\tilde{r}_{i,\si}\right]+\sum\limits_{l=1}^{{\color{black}{33}}}\sum\limits_{i}\widetilde{\beta}_{i}^{1,l}.
\end{align*}
First we have $\widetilde{\B}_{i}^{1,1}=\sum\limits_{j\neq i}[(\mu_iw_{i,x}-\tilde{\la}_iw_i)v_j-(w_j-\la_j^*v_j)w_i]\tilde{r}_{i,u}\tilde{r}_j=O(1)\La_i^1$ and we get from Lemmas  \ref{lemme6.5}, \ref{lemme6.6}, \ref{lemme11.3}, \ref{lemme11.4}%s \ref{lemme6.5}, \ref{lemme6.6}, \ref{estimimpo1}, \ref{lemme10.5}
\begin{align*}
	\widetilde{\B}_{i,x}^{1,1}&= \sum\limits_{j\neq i}[(\mu_iw_{i,xx}+\mu_{i,x}w_{i,x}-\tilde{\la}_{i,x}w_{i}-\tilde{\la}_iw_{i,x})v_j-(w_{j,x}-\la_j^*v_{j,x})w_i]\tilde{r}_{i,u}\tilde{r}_j\\
	&+ \sum\limits_{j\neq i}[(\mu_iw_{i,x}-\tilde{\la}_iw_i)v_{j,x}-(w_j-\la_j^*v_j)w_{i,x}]\tilde{r}_{i,u}\tilde{r}_j\\
	&+ \sum\limits_{j\neq i}[(\mu_iw_{i,x}-\tilde{\la}_iw_i)v_j-(w_j-\la_j^*v_j)w_i]\pa_x(\tilde{r}_{i,u}\tilde{r}_j)=O(1)\La_i^1.
	%[\tilde{r}_{i,uu}(u_x\otimes \tilde{r}_j)+(\xi_{i}\bar{v}_i)_x\tilde{r}_{i,uv}\tilde{r}_j+\theta_i^\p\left(\frac{w_i}{v_i}\right)_x\tilde{r}_{u\si}\tilde{r}_j]\\
	%&+ \sum\limits_{j\neq i}[(\mu_iw_{i,x}-\tilde{\la}_iw_i)v_j-(w_j-\la_j^*v_j)w_i]\tilde{r}_{i,u}\left[\tilde{r}_{j,u}u_x+(\xi_j\bar{v}_j)_x\tilde{r}_{j,v}-\theta_{j}^\p\left(\frac{w_j}{v_j}\right)_x\tilde{r}_{j,\si}\right].
\end{align*}
Next it yields applying Lemmas \ref{lemme6.5}, \ref{lemme6.6} and \eqref{ngtech4} 
\begin{align*}
&\widetilde{\B}^{1,2}_{i}= [\mu_i(w_{i,x}v_i-w_iv_{i,x})+%(\mu_iv_{i,x}-(\tilde{\la}_i-\la_i^*)v_i-w_i)w_i
\mathcal{H}_i w_i]\tilde{r}_{i,u}\tilde{r}_i\\
&=O(1)\sum_j(\La_j^1+\La_j^4+\La_j^6)+O(1)v_i^2\xi_i (\chi')^\e_i\hat{v}_i+O(1)\eta_i w_i v_i^2(1- \chi_i^\ep)\\
&=O(1)\sum_j(\La_j^1+\La_j^4+\La_j^6)+R_{\e,1}^{1,i,2},
\end{align*}
 with $R_{\e,1}^{1,i,2}$ satisfying \eqref{12.11bis2}.
Using Lemmas \ref{lemme6.5}, \ref{lemme6.6}, \ref{lemme11.4} and \eqref{ngtech4}, \eqref{estimate-v-i-xx-1aa} we obtain
\begin{align*}
	\widetilde{\B}^{1,2}_{i,x}&=[\mu_{i,x}(w_{i,x}v_i-w_iv_{i,x})+\mu_{i}(w_{i,xx}v_i-w_iv_{i,xx})]\tilde{r}_{i,u}\tilde{r}_i\\
	&+[\mathcal{H}_{i,x}w_{i}+\mathcal{H}_i w_{i,x}]\tilde{r}_{i,u}\tilde{r}_i+[\mu_i(w_{i,x}v_i-w_iv_{i,x})+\mathcal{H}_iw_i]\pa_x(\tilde{r}_{i,u}\tilde{r}_i)\\
	&=O(1)\sum_j(\La_j^1+\La_j^4+\La_j^5+\La_j^6)+\widetilde{R}_{\e,1}^{1,i,2},
	%[\tilde{r}_{i,uu}(\tilde{r}_i\otimes u_x)+\tilde{r}_{i,u}[\tilde{r}_{i,u}u_x]]\\
	%&+[\mu_i(w_{i,x}v_i-w_iv_{i,x})+(\mu_iv_{i,x}-(\tilde{\la}_i-\la_i^*)v_i-w_i)w_i](\xi_i\bar{v}_i)_x[\tilde{r}_{i,uv}\tilde{r}_i+\tilde{r}_{i,u}\tilde{r}_{i,v}]\\
	%&-[\mu_i(w_{i,x}v_i-w_iv_{i,x})+(\mu_iv_{i,x}-(\tilde{\la}_i-\la_i^*)v_i-w_i)w_i]\theta_i^\p\left(\frac{w_i}{v_i}\right)_x[\tilde{r}_{i,u\si}\tilde{r}_i+\tilde{r}_{i,u}\tilde{r}_{i,\si}].
\end{align*}
 with $\widetilde{R}_{\e,1}^{1,i,2}$ satisfying \eqref{12.11bis2}.
Using now the fact that $\mbox{supp}\xi$ is included in $\{x,\theta(x)=x)$ we get using \eqref{ngtech4}
\begin{align*}
\widetilde{\B}^{1,3}_{i}=&%&(\mu_iw_{i,x}-(\tilde{\la}_i-\la_i^*+\theta_i)w_i)v_{i,x}\xi_i(\pa_{v_i}\bar{v}_i)\tilde{r}_{i,v}\\
(\mu_iw_{i,x}-(\tilde{\la}_i-\la_i^*+\frac{w_i}{v_i})w_i)v_{i,x}\xi_i(\pa_{v_i}\bar{v}_i)\tilde{r}_{i,v}\\
=&w_{i,x}\mathcal{H}_i\xi_i(\pa_{v_i}\bar{v}_i)\tilde{r}_{i,v}+(\tilde{\la}_i-\la_i^*)(w_{i,x}v_i-w_iv_{i,x})\xi_i(\pa_{v_i}\bar{v}_i)\tilde{r}_{i,v}+v_iw_i\left(\frac{w_i}{v_i}\right)_x\xi_i(\pa_{v_i}\bar{v}_i)\tilde{r}_{i,v}\\
=&\sum_j(\La_j^1+\La_j^4+\La_j^6)+R_{\e,1}^{1,i,3},
\end{align*}
with $R_{\e,1}^{1,i,3}$ satisfying \eqref{12.11bis2}.
Furthermore applying Lemmas  \ref{estimimpo1}, \ref{lemme6.5}, \ref{lemme6.6}, \ref{lemme11.3} and \eqref{ngtech4}, \eqref{ngtech5}, \eqref{estimate-v-i-xx-1aa}  we get
\begin{align*}
	&\widetilde{\B}^{1,3}_{i,x}=\big[w_{i,xx}\mathcal{H}_i\xi_i+w_{i,x} \mathcal{H}_{i,x}\xi_i+w_{i,x}\mathcal{H}_i\xi_i'\left(\frac{w_i}{v_i}\right)_x\big]\pa_{v_i}\bar{v}_i\tilde{r}_{i,v}\\
	&+w_{i,x}\mathcal{H}_i\xi_i\big[\pa_x(\pa_{v_i}\bar{v}_i)\big)\tilde{r}_{i,v}+\pa_{v_i}\bar{v}_i\pa_x(\tilde{r}_{i,v})\big]	+\tilde{\la}_{i,x}(w_{i,x}v_i-w_iv_{i,x})\xi_i(\pa_{v_i}\bar{v}_i)\tilde{r}_{i,v}\\
	&+(\tilde{\la}_i-\la_i^*)(w_{i,xx}v_i-w_iv_{i,xx})\xi_i(\pa_{v_i}\bar{v}_i)\tilde{r}_{i,v}+(\tilde{\la}_i-\la_i^*)(w_{i,x}v_i-w_iv_{i,x})\xi_i'\left(\frac{w_i}{v_i}\right)_x(\pa_{v_i}\bar{v}_i)\tilde{r}_{i,v}
	\\
	&+(\tilde{\la}_i-\la_i^*)(w_{i,x}v_i-w_iv_{i,x})\xi_i\big[\pa_x(\pa_{v_i}\bar{v}_i)\tilde{r}_{i,v}+(\pa_{v_i}\bar{v}_i)\pa_x(\tilde{r}_{i,v})\big]\\
	&+\big[v_{i,x}w_i \left(\frac{w_i}{v_i}\right)_x+v_i w_{i,x}\left(\frac{w_i}{v_i}\right)_x+v_iw_i\left(\frac{w_i}{v_i}\right)_{xx}\big]\xi_i(\pa_{v_i}\bar{v}_i)\tilde{r}_{i,v}\\
	&+v_iw_i\left(\frac{w_i}{v_i}\right)_x\big[\xi_i' \left(\frac{w_i}{v_i}\right)_x(\pa_{v_i}\bar{v}_i)\tilde{r}_{i,v}+\xi_i\pa_x(\pa_{v_i}\bar{v}_i)\tilde{r}_{i,v}+\xi_i(\pa_{v_i}\bar{v}_i)\pa_x(\tilde{r}_{i,v})\big]\\
	%&+v_iw_i\left(\frac{w_i}{v_i}\right)_x\xi_i\pa_x(\pa_{v_i}\bar{v}_i)\tilde{r}_{i,v}+v_iw_i\left(\frac{w_i}{v_i}\right)_x\xi_i(\pa_{v_i}\bar{v}_i)\pa_x(\tilde{r}_{i,v})\\
	&=O(1)\sum_j(\La_j^1+\delta_0^2\La_j^3+\La_j^4+\La_j^5+\La_j^6)+\widetilde{R}_{\e,1}^{1,i,3},
		%(\mu_iw_{i,x}-(\tilde{\la}_i-\la_i^*+\theta_i)w_i)_xv_{i,x}\xi_i(\pa_{v_i}\bar{v}_i)\tilde{r}_{i,v}\\
	%&+(\mu_iw_{i,x}-(\tilde{\la}_i-\la_i^*+\theta_i)w_i)v_{i,xx}\xi_i(\pa_{v_i}\bar{v}_i)\tilde{r}_{i,v}\\
	%&+(\mu_iw_{i,x}-(\tilde{\la}_i-\la_i^*+\theta_i)w_i)v_{i,x}\xi_i^\p\left(\frac{w_i}{v_i}\right)_x(\pa_{v_i}\bar{v}_i)\tilde{r}_{i,v}\\
	%&+(\mu_iw_{i,x}-(\tilde{\la}_i-\la_i^*+\theta_i)w_i)v_{i,x}\xi_i(\pa_{v_i}\bar{v}_i)_x\tilde{r}_{i,v}\\
	%&+(\mu_iw_{i,x}-(\tilde{\la}_i-\la_i^*+\theta_i)w_i)v_{i,x}\xi_i(\pa_{v_i}\bar{v}_i)[\tilde{r}_{i,uv}u_x+(\xi\bar{v}_i)_x\tilde{r}_{i,vv}-\theta_i^\p\left(\frac{w_i}{v_i}\right)_x\tilde{r}_{i,v\si}].
\end{align*}
with $\widetilde{R}_{\e,1}^{1,i,3}$ satisfying \eqref{12.11bis2}.
We give here some details on the terms $w_{i,x}\mathcal{H}_i\xi_i'\left(\frac{w_i}{v_i}\right)_x(\pa_{v_i}\bar{v}_i)\tilde{r}_{i,v}$ and $w_{i,x}\mathcal{H}_i\xi_i\pa_x(\pa_{v_i}\bar{v}_i)\tilde{r}_{i,v}$. Applying Lemma \ref{estimimpo1} and \eqref{ngtech5}, we have
\begin{align*}
&w_{i,x}\mathcal{H}_i\xi_i'\left(\frac{w_i}{v_i}\right)_x(\pa_{v_i}\bar{v}_i)\tilde{r}_{i,v}=\mathcal{H}_i\big(\xi_i'v_i\left(\frac{w_i}{v_i}\right)^{2}_x(\pa_{v_i}\bar{v}_i)+\frac{w_i}{v_i}v_{i,x}\xi_i'\left(\frac{w_i}{v_i}\right)_x(\pa_{v_i}\bar{v}_i)\big)\tilde{r}_{i,v}\\
&=O(1)\sum_j(\La_j^1+\delta_0^2\La_j^3+\La_j^4).
\end{align*}
Similarly we have using Lemmas \ref{estimimpo1}, \ref{lemme6.5}, \ref{lemme6.6} and \eqref{ngtech5}, \eqref{ngtech4}
\begin{align*}
&w_{i,x}\mathcal{H}_i\xi_i\pa_x(\pa_{v_i}\bar{v}_i)\tilde{r}_{i,v}=O(1)\mathfrak{A}_i\rho_i^\e w_{i,x}\mathcal{H}_i\big(\frac{|v_{i,x}|}{|v_i|}+\sum_{j\ne i}\frac{|v_jv_{j,x}|}{|v_i|}\big)\\
&=O(1)\mathfrak{A}_i\rho_i^\e w_{i,x}\mathcal{H}_i(1+\sum_{k\ne i}\frac{v_k^2}{|v_i|}|\left(\frac{w_k}{v_k}\right)_x|\mathfrak{A}_k\rho_k^\e\big)\\
&\hspace{2cm}+O(1)\mathfrak{A}_i\rho_i^\e w_{i,x}\big(\sum_{j\ne i}|v_jv_{j,x}|+\sum_{j\ne i}\frac{|v_jv_{j,x}|}{|v_i|}\sum_{l\ne i}v_l^2|\left(\frac{w_k}{v_k}\right)_x|\mathfrak{A}_l\rho_l^\e\big)
\\
&=O(1)\sum_{j}(\La_j^1+\La_j^4+\La_j^6)+\widetilde{R}_{\e,1}^{1,i,3,1},
\end{align*}
with $\widetilde{R}_{\e,1}^{1,i,3,1}$ satisfying \eqref{12.11bis2}.
We have now $\widetilde{\B}^{1,4}_{i}=\sum\limits_{j\neq i}(\mu_iw_{i,x}-\tilde{\la}_iw_i )v_{j,x}\xi_i(\pa_{v_j}\bar{v}_i)\tilde{r}_{i,v}=O(1)\La_i^1$ and from Lemmas \ref{estimimpo1}, \ref{lemme6.5}, \ref{lemme6.6}, \ref{lemme11.3}, \eqref{ngtech5}
\begin{align*}
	\widetilde{\B}^{1,4}_{i,x}&=\sum\limits_{j\neq i}(\mu_iw_{i,x}-\tilde{\la}_iw_i)_xv_{j,x}\xi_i(\pa_{v_j}\bar{v}_i)\tilde{r}_{i,v}+\sum\limits_{j\neq i}(\mu_iw_{i,x}-\tilde{\la}_iw_i)v_{j,xx}\xi_i(\pa_{v_j}\bar{v}_i)\tilde{r}_{i,v}\\
	&+\sum\limits_{j\neq i}(\mu_iw_{i,x}-\tilde{\la}_iw_i)v_{j,x}\big[\xi^\p_i\left(\frac{w_i}{v_i}\right)_x(\pa_{v_j}\bar{v}_i)\tilde{r}_{i,v}+\xi_i(\pa_{v_j}\bar{v}_i)_x\tilde{r}_{i,v}+\xi_i(\pa_{v_j}\bar{v}_i)\pa_x(\tilde{r}_{i,v})\big]\\
	%&+\sum\limits_{j\neq i}(\mu_iw_{i,x}-\tilde{\la}_i w_i)v_{j,x}\xi_i(\pa_{v_j}\bar{v}_i)_x\tilde{r}_{i,v}+\sum\limits_{j\neq i}(\mu_iw_{i,x}-\tilde{\la}_iw_i)v_{j,x}\xi_i(\pa_{v_j}\bar{v}_i)\pa_x(\tilde{r}_{i,v})\\
	&=O(1)\sum_j (\La_i^1+\delta_0^2\La_j^3).
\end{align*}
%It yields now $\B^{1,5}_{i}=\sum\limits_{j\neq i}(\mu_iw_{i,x}-(\tilde{\la}_i-\la_i^*+\theta_i)w_i)w_{j,x}\xi_i(\pa_{w_j}\bar{v}_i)\tilde{r}_{i,v}=O(1)\La_i^1$ and using Lemmas \ref{lemme6.5}, \ref{lemme6.6} and \ref{estimimpo1} we obtain
%\begin{align*}
%	\B^{1,5}_{i,x}&=\sum\limits_{j\neq i}(\mu_iw_{i,x}-(\tilde{\la}_i-\la_i^*+\theta_i)w_i)_x w_{j,x}\xi_i(\pa_{w_j}\bar{v}_i)\tilde{r}_{i,v}\\
%	&+\sum\limits_{j\neq i}(\mu_iw_{i,x}-(\tilde{\la}_i-\la_i^*+\theta_i)w_i)w_{j,xx}\xi_i(\pa_{w_j}\bar{v}_i)\tilde{r}_{i,v}\\
%	&+\sum\limits_{j\neq i}(\mu_iw_{i,x}-(\tilde{\la}_i-\la_i^*+\theta_i)w_i)w_{j,x}\xi_i^\p\left(\frac{w_i}{v_i}\right)_x(\pa_{w_j}\bar{v}_i)\tilde{r}_{i,v}\\
%	&+\sum\limits_{j\neq i}(\mu_iw_{i,x}-(\tilde{\la}_i-\la_i^*+\theta_i)w_i)w_{j,x}\xi_i(\pa_{w_j}\bar{v}_i)_x\tilde{r}_{i,v}\\
%	&+\sum\limits_{j\neq i}(\mu_iw_{i,x}-(\tilde{\la}_i-\la_i^*+\theta_i)w_i)w_{j,x}\xi_i(\pa_{w_j}\bar{v}_i)\pa_x(\tilde{r}_{i,v})\\
%	&=O(1)(\La_i^1+\delta_0^2\La_i^3).
%\end{align*}
We have now $\widetilde{\B}^{1,5}_{i}=
2\mu_iv_i\xi^\p_i\bar{v}_i\left(\frac{w_i}{v_i}\right)^2_x\tilde{r}_{i,v}=O(1)\delta_0^2\La_i^3$ and from Lemmas \ref{lemme6.5}, \ref{lemme6.6}
\begin{align*}
	&\widetilde{\B}^{1,5}_{i,x}=2\mu_{i,x}v_i\xi^\p_i\bar{v}_i\left(\frac{w_i}{v_i}\right)^2_x\tilde{r}_{i,v}+2\mu_iv_{i,x}\xi^\p_i\bar{v}_i\left(\frac{w_i}{v_i}\right)^2_x\tilde{r}_{i,v}+2\mu_iv_i\xi_i''\bar{v}_i\left(\frac{w_i}{v_i}\right)^3_x\tilde{r}_{i,v}\\
	&+2\mu_iv_i\xi^\p_i\pa_x(\bar{v}_i)\left(\frac{w_i}{v_i}\right)^2_x\tilde{r}_{i,v}+4\mu_iv_i\xi^\p_i\bar{v}_i\left(\frac{w_i}{v_i}\right)_{xx}\left(\frac{w_i}{v_i}\right)_{x}\tilde{r}_{i,v}
	+2\mu_iv_i\xi^\p_i\bar{v}_i\left(\frac{w_i}{v_i}\right)^2_x\pa_x(\tilde{r}_{i,v})\\
	&=O(1)(\delta_0^2\La_i^3+\La_i^5).
\end{align*}
We can treated in a similar way $\widetilde{\B}^{1,6}_{i}=\mu_i v_i\xi^{\p\p}_i\bar{v}_i\frac{w_i}{v_i}\left(\frac{w_i}{v_i}\right)^2_x\tilde{r}_{i,v}=O(1)\delta_0^2\La_i^3$. We have now using Lemma \ref{estimimpo1} and \eqref{ngtech5} $\widetilde{\B}^{1,7}_{i}=
\mu_i v_{i}\left(\frac{w_i}{v_i}\right)_xv_{i,x}\xi^\p_i(\pa_{v_i}\bar{v}_i)\frac{w_i}{v_i}\tilde{r}_{i,v}=O(1)(\La_i^1+\La_i^4)$ and from Lemma \ref{estimimpo1},  \ref{lemme6.5}, \ref{lemme6.6}, \ref{lemme11.3} and \eqref{ngtech5}, \eqref{estimate-v-i-xx-1aabis} 
\begin{align*}
	\widetilde{\B}^{1,7}_{i,x}&=\big[\mu_{i,x} v_{i}\left(\frac{w_i}{v_i}\right)_x+
	\mu_{i} v_{i,x}\left(\frac{w_i}{v_i}\right)_x+\mu_{i} v_{i}\left(\frac{w_i}{v_i}\right)_{xx}\big]
	v_{i,x}
	\xi^\p_i(\pa_{v_i}\bar{v}_i)\frac{w_i}{v_i}\tilde{r}_{i,v}	\\
	&
	+\mu_i v_{i}\left(\frac{w_i}{v_i}\right)_x \big[v_{i,xx}\xi^\p_i(\pa_{v_i}\bar{v}_i)+v_{i,x}\xi^{\p\p}_i\left(\frac{w_i}{v_i}\right)_x(\pa_{v_i}\bar{v}_i)+v_{i,x}\xi^\p_i(\pa_{v_i}\bar{v}_i)_x\big]
	\frac{w_i}{v_i}\tilde{r}_{i,v}\\
	%&+\mu_i v_{i}\left(\frac{w_i}{v_i}\right)^2_xv_{i,x}\xi^{\p\p}_i(\pa_{v_i}\bar{v}_i)\frac{w_i}{v_i}\tilde{r}_{i,v}+\mu_i v_{i}\left(\frac{w_i}{v_i}\right)_xv_{i,x}\xi^\p_i(\pa_{v_i}\bar{v}_i)_x\frac{w_i}{v_i}\tilde{r}_{i,v}\\
	&+\mu_i v_{i}\left(\frac{w_i}{v_i}\right)_x^2 v_{i,x}\xi^\p_i(\pa_{v_i}\bar{v}_i)\tilde{r}_{i,v}+\mu_iv_{i}\left(\frac{w_i}{v_i}\right)_xv_{i,x}\xi^\p_i(\pa_{v_i}\bar{v}_i)\frac{w_i}{v_i}\pa_x(\tilde{r}_{i,v})\\
	&=O(1)\sum_j(\La_j^1+\delta_0^2\La_j^3+\La_j^4+\La_j^5+\La_j^{6}+\La_j^{6,1}).
\end{align*}
We have now $\widetilde{\B}^{1,8}_{i}=\sum\limits_{j\neq i}\mu_iv_i\left(\frac{w_i}{v_i}\right)_x v_{j,x}\xi_i^\p(\pa_{v_j}\bar{v}_i)\frac{w_i}{v_i}\tilde{r}_{i,v}=O(1)\La_i^1$ and using Lemmas \ref{estimimpo1}, \ref{lemme6.5}, \ref{lemme6.6}, \ref{lemme11.3} and \eqref{ngtech5}
\begin{align*}
	\widetilde{\B}^{1,8}_{i,x}&=\sum\limits_{j\neq i}\big[\mu_{i,x}v_i\left(\frac{w_i}{v_i}\right)_x+\mu_iv_{i,x}\left(\frac{w_i}{v_i}\right)_x+\mu_iv_{i}\left(\frac{w_i}{v_i}\right)_{xx}\big]
	v_{j,x}\xi'_i(\pa_{v_j}\bar{v}_i)\frac{w_i}{v_i}\tilde{r}_{i,v}\\
	%&+\sum\limits_{j\neq i}\mu_iv_{i,x}\left(\frac{w_i}{v_i}\right)_xv_{j,x}\xi'_i(\pa_{v_j}\bar{v}_i)\frac{w_i}{v_i}\tilde{r}_{i,v}+\sum\limits_{j\neq i}\mu_iv_{i}\left(\frac{w_i}{v_i}\right)_{xx}v_{j,x}\xi_i(\pa_{v_j}\bar{v}_i)\frac{w_i}{v_i}\tilde{r}_{i,v}\\
	&+\sum\limits_{j\neq i}\mu_i v_i\left(\frac{w_i}{v_i}\right)_x\big[v_{j,xx}\xi'_i(\pa_{v_j}\bar{v}_i)
	+ v_{j,x}\xi_i''\left(\frac{w_i}{v_i}\right)_x(\pa_{v_j}\bar{v}_i)+v_{j,x}\xi_i'(\pa_{v_j}\bar{v}_i)_x\big]
	\frac{w_i}{v_i}\tilde{r}_{i,v}\\%+\sum\limits_{j\neq i}\mu_iv_i\left(\frac{w_i}{v_i}\right)^2_xv_{j,x}\xi_i''(\pa_{v_j}\bar{v}_i)\frac{w_i}{v_i}\tilde{r}_{i,v}\\
	%&+\sum\limits_{j\neq i}\mu_i v_i\left(\frac{w_i}{v_i}\right)_xv_{j,x}\xi_i'(\pa_{v_j}\bar{v}_i)_x\frac{w_i}{v_i}\tilde{r}_{i,v}
	&+\sum\limits_{j\neq i}\mu_iv_i\left(\frac{w_i}{v_i}\right)^2_x v_{j,x}\xi_i'(\pa_{v_j}\bar{v}_i)\tilde{r}_{i,v}+\sum\limits_{j\neq i}\mu_iv_i\left(\frac{w_i}{v_i}\right)_xv_{j,x}\xi'_i(\pa_{v_j}\bar{v}_i)\frac{w_i}{v_i}\pa_x(\tilde{r}_{i,v})\\
	&=O(1)(\La_i^1+\delta_0^2\La_i^3).
\end{align*}
We observe now that using Lemma \ref{estimimpo1} and \eqref{ngtech4} we have $\widetilde{\B}^{1,9}_{i}=
\mathcal{H}_i\xi_i w_iv_{i,x}\pa_{v_iv_i}\bar{v}_i\tilde{r}_{i,v}=O(1)\sum_j(\La_j^1+\La_j^4+\La_j^6)+R_{\e,1}^{1,i,9}$ with $R_{\e,1}^{1,i,9}$ satisfying \eqref{12.11bis2}. Furthermore we get employing Lemmas  \ref{estimimpo1}, \ref{lemme6.5}, \ref{lemme6.6}, \eqref{ngtech4}, \eqref{ngtech5},\eqref{estimate-v-i-xx-1aa} and the fact that $w_{i,x}=v_i\left(\frac{w_i}{v_i}\right)_x+\frac{w_i}{v_i}v_{i,x}$
\begin{align*}
\widetilde{\B}^{1,9}_{i,x}&=\big[\mathcal{H}_{i,x}\xi_i w_iv_{i,x}+\mathcal{H}_i\xi_i'\left(\frac{w_i}{v_i}\right)_x w_iv_{i,x}+\mathcal{H}_i\xi_i w_{i,x}v_{i,x}+\mathcal{H}_i\xi_i w_iv_{i,xx}\big]
\pa_{v_iv_i}\bar{v}_i\tilde{r}_{i,v}\\
%&+(\mu_i v_{i,x}-v_i(\tilde{\lambda}_i-\lambda_i^*)-w_i)\xi_i'\left(\frac{w_i}{v_i}\right)_x w_iv_{i,x}\pa_{v_iv_i}\bar{v}_i\tilde{r}_{i,v}\\
%&+(\mu_i v_{i,x}-v_i(\tilde{\lambda}_i-\lambda_i^*)-w_i)\xi_i w_{i,x}v_{i,x}\pa_{v_iv_i}\bar{v}_i\tilde{r}_{i,v}\\
%&+(\mu_i v_{i,x}-v_i(\tilde{\lambda}_i-\lambda_i^*)-w_i)\xi_i w_iv_{i,xx}\pa_{v_iv_i}\bar{v}_i\tilde{r}_{i,v}\\
&+\mathcal{H}_i\xi_i w_iv_{i,x}\big[\pa_x(\pa_{v_iv_i}\bar{v}_i)\tilde{r}_{i,v}+\pa_{v_iv_i}\bar{v}_i\pa_x(\tilde{r}_{i,v})\big]\\
%&+(\mu_i v_{i,x}-v_i(\tilde{\lambda}_i-\lambda_i^*)-w_i)\xi_i w_iv_{i,x}\pa_{v_iv_i}\bar{v}_i\pa_x(\tilde{r}_{i,v})\\
&=O(1)\sum_j(\La_j^1+\La_j^4+\La_j^5+\La_j^6)+\widetilde{R}_{\e,1}^{1,i,9},
\end{align*}
with $\widetilde{R}_{\e,1}^{1,i,9}$ satisfying \eqref{12.11bis2}.
We have now from Lemma \ref{estimimpo1} $\widetilde{\B}^{1,10}_{i}=\sum_{i\ne j}\mathcal{H}_i \xi_iw_iv_{j,x}\pa_{v_jv_i}\bar{v}_i\tilde{r}_{i,v}=O(1)\La_i^1$ and we get from Lemmas  \ref{estimimpo1}, \ref{lemme6.5}, \ref{lemme6.6} and \eqref{ngtech5}
\begin{align*}
&\widetilde{\B}^{1,10}_{i,x}=\sum_{i\ne j}\big[\mathcal{H}_{i,x}\xi_iw_iv_{j,x}+\mathcal{H}_i\xi'_i \left(\frac{w_i}{v_i}\right)_x w_iv_{j,x}+\mathcal{H}_i \xi_i(w_{i,x}v_{j,x}+w_i v_{j,xx})\big]
\pa_{v_jv_i}\bar{v}_i\tilde{r}_{i,v}\\
%&+\sum_{i\ne j}(\mu_i v_{i,x}-v_i(\tilde{\lambda}_i-\lambda_i^*)-w_i) \xi'_i \left(\frac{w_i}{v_i}\right)_x w_iv_{j,x}\pa_{v_jv_i}\bar{v}_i\tilde{r}_{i,v}\\
%&+\sum_{i\ne j}(\mu_i v_{i,x}-v_i(\tilde{\lambda}_i-\lambda_i^*)-w_i) \xi_i(w_{i,x}v_{j,x}+w_i v_{j,xx})\pa_{v_jv_i}\bar{v}_i\tilde{r}_{i,v}\\
&+\sum_{i\ne j}\mathcal{H}_i \xi_iw_iv_{j,x}\big(\pa_x(\pa_{v_jv_i})\bar{v}_i\tilde{r}_{i,v}+\pa_{v_jv_i}\bar{v}_i\pa_x(\tilde{r}_{i,v})\big)=O(1)\sum_j\La_j^1.
\end{align*}
We have now from Lemmas \ref{lemme6.5}, \ref{lemme6.6} and \eqref{ngtech4}
$\widetilde{\B}^{1,11}_{i}=\mathcal{H}_i\xi_iw_{i,x}\pa_{v_i}\bar{v}_i\tilde{r}_{i,v}=O(1)\sum_j(\La_j^1+\La_j^4+\La_j^6)+R_{\e,1}^{1,i,11}$ with $R_{\e,1}^{1,i,11}$ satisfying \eqref{12.11bis2}.
 In addition we get from Lemmas \ref{estimimpo1}, \ref{lemme6.5}, \ref{lemme6.6}, \ref{lemme11.3}, \eqref{ngtech5},\eqref{estimate-v-i-xx-1aa}  and the fact that $w_{i,x}=v_i\left(\frac{w_i}{v_i}\right)_x+\frac{w_i}{v_i}v_{i,x}$
\begin{align*}
&\widetilde{\B}^{1,11}_{i,x}=\big[\mathcal{H}_{i,x} \xi_iw_{i,x}\pa_{v_i}\bar{v}_i+\mathcal{H}_i \xi_i' \left(\frac{w_i}{v_i}\right)_x w_{i,x}\pa_{v_i}\bar{v}_i+\mathcal{H}_i \xi_iw_{i,xx}\pa_{v_i}\bar{v}_i+\mathcal{H}_i \xi_iw_{i,x}\pa_x(\pa_{v_i}\bar{v}_i )\big]
\tilde{r}_{i,v}\\
%&+( \mu_i v_{i,x}-v_i(\tilde{\lambda}_i-\lambda_i^*)-w_i) \xi_i' \left(\frac{w_i}{v_i}\right)_x w_{i,x}\pa_{v_i}\bar{v}_i\tilde{r}_{i,v}\\
%&+( \mu_i v_{i,x}-v_i(\tilde{\lambda}_i-\lambda_i^*)-w_i) \xi_i\big(w_{i,xx}\pa_{v_i}\bar{v}_i+w_{i,x}\pa_x(\pa_{v_i}\bar{v}_i)\big) \tilde{r}_{i,v}\\
&+\mathcal{H}_i\xi_iw_{i,x}\pa_{v_i}\bar{v}_i\pa_x(\tilde{r}_{i,v})=O(1)\sum_j(\La_j^1+\delta_0^2\La_j^3+\La_j^4+\La_j^5+\La_j^6)+\widetilde{R}_{\e,1}^{1,i,11},
\end{align*}
with $\widetilde{R}_{\e,1}^{1,i,11}$ satisfying \eqref{12.11bis2}.
We have now  from Lemma \ref{estimimpo1} and \eqref{ngtech5} 
$$\widetilde{\B}^{1,12}_{i}=\mathcal{H}_i w_i\xi_i^\p\left(\frac{w_i}{v_i}\right)_x\pa_{v_i}\bar{v}_i\tilde{r}_{i,v}=O(1)\sum_j(\La_j^1+\La_j^4).$$ As previously we get
\begin{align*}
&\widetilde{\B}^{1,12}_{i,x}=\big[\mathcal{H}_{i,x}w_i\xi_i'\left(\frac{w_i}{v_i}\right)_x+\mathcal{H}_i w_{i,x}\xi_i'\left(\frac{w_i}{v_i}\right)_x+\mathcal{H}_i w_i\xi_i''\left(\frac{w_i}{v_i}\right)^2_x+\mathcal{H}_iw_i\xi_i'\left(\frac{w_i}{v_i}\right)_{xx}\big]
\pa_{v_i}\bar{v}_i\tilde{r}_{i,v}\\
%&+(\mu_i v_{i,x}-v_i(\tilde{\lambda}_i-\lambda_i^*)-w_i) w_{i,x}\xi_i'\left(\frac{w_i}{v_i}\right)_x\pa_{v_i}\bar{v}_i\tilde{r}_{i,v}\\
%&+(\mu_i v_{i,x}-v_i(\tilde{\lambda}_i-\lambda_i^*)-w_i) w_i\xi_i''\left(\frac{w_i}{v_i}\right)^2_x\pa_{v_i}\bar{v}_i\tilde{r}_{i,v}\\
%&+(\mu_i v_{i,x}-v_i(\tilde{\lambda}_i-\lambda_i^*)-w_i) w_i\xi_i'\left(\frac{w_i}{v_i}\right)_{xx}\pa_{v_i}\bar{v}_i\tilde{r}_{i,v}\\
&+\mathcal{H}_i w_i\xi_i'\left(\frac{w_i}{v_i}\right)_x\big[\pa_x(\pa_{v_i}\bar{v}_i)\tilde{r}_{i,v}+\pa_{v_i}\bar{v}_i\pa_x(\tilde{r}_{i,v})\big]=O(1)\sum_j(\La_j^1+\delta_0^2\La_j^3+\La_j^4+\La_j^5+\La_j^6).
%&+(\mu_i v_{i,x}-v_i(\tilde{\lambda}_i-\lambda_i^*)-w_i) w_i\xi_i'\left(\frac{w_i}{v_i}\right)_x\pa_{v_i}\bar{v}_i\pa_x(\tilde{r}_{i,v})\\
%&=O(1)\sum_j(\La_j^1+\delta_0^2\La_j^3+\La_j^4+\La_j^5+\La_j^6).
\end{align*}
Similarly we have $\widetilde{\B}^{1,13}_{i}=
-\theta'_i\pa_{v_i}\bar{v}_iv_iw_i\left(\frac{w_i}{v_i}\right)_x\xi_i \tilde{r}_{i,v}=O(1)\La_i^4$ and using the fact that $w_{i,x}=v_i\left(\frac{w_i}{v_i}\right)_x+\frac{w_i}{v_i}v_{i,x}$
\begin{align*}
\widetilde{\B}^{1,13}_{i,x}&=-\big[\theta''_i\left(\frac{w_i}{v_i}\right)_x \pa_{v_i}\bar{v}_iv_iw_i+\theta'_i\big(\pa_x(\pa_{v_i}\bar{v}_i)v_i w_i+\pa_{v_i}\bar{v}_iv_{i,x}w_i+\pa_{v_i}\bar{v}_iv_{i}w_{i,x}\big)\big]
\left(\frac{w_i}{v_i}\right)_x\xi_i \tilde{r}_{i,v}\\%-\theta'_i\big(\pa_x(\pa_{v_i}\bar{v}_i)v_i+\pa_{v_i}\bar{v}_iv_{i,x}\big)w_i\left(\frac{w_i}{v_i}\right)_x\xi_i \tilde{r}_{i,v}\\
%&-\theta'_i\pa_{v_i}\bar{v}_iv_iw_{i,x}\left(\frac{w_i}{v_i}\right)_x\xi_i \tilde{r}_{i,v}
&-\theta'_i\pa_{v_i}\bar{v}_iv_iw_i\big[\left(\frac{w_i}{v_i}\right)_{xx}\xi_i \tilde{r}_{i,v}+\left(\frac{w_i}{v_i}\right)^2_x\xi_i' \tilde{r}_{i,v}+\left(\frac{w_i}{v_i}\right)_x\xi_i\pa_x( \tilde{r}_{i,v})\big]\\
%&-\theta'_i\pa_{v_i}\bar{v}_iv_iw_i\left(\frac{w_i}{v_i}\right)^2_x\xi_i' \tilde{r}_{i,v}-\theta'_i\pa_{v_i}\bar{v}_iv_iw_i\left(\frac{w_i}{v_i}\right)_x\xi_i\pa_x( \tilde{r}_{i,v})
&=O(1)\sum_j(\La_j^1+\delta_0^2\La_j^3+\La_j^4+\La_j^5).
\end{align*}
In the same manner we have $\widetilde{\B}^{1,14}_{i}
=-2\mu_iv_i\theta_i^\p\left(\frac{w_i}{v_i}\right)^2_x\tilde{r}_{i,\si}=O(1)\delta_0^2\La_i^3$ from Lemma \ref{lemme6.5} and the fact that $\tilde{r}_{i,\si}=O(1)\xi_i \bar{v}_i$. Applying Lemmas \ref{lemme6.5}, \ref{lemme11.3}, we deduce that
%\begin{align*}
%	\B^{1,9}_{i,x}&=\sum\limits_{j\neq i}\mu_{i,x}(w_{i,x}-(w_i/v_i)v_{i,x})\zeta_{j,x}\xi_i(\pa_{\zeta_j}\bar{v}_i)\frac{w_i}{v_i}\tilde{r}_{i,v}\\
%	&+\sum\limits_{j\neq i}\mu_i(w_{i,xx}-(w_i/v_i)v_{i,xx}-v_{i,x}(w_i/v_i)_x)\zeta_{j,x}\xi_i(\pa_{\zeta_j}\bar{v}_i)\frac{w_i}{v_i}\tilde{r}_{i,v}\\
%	&+\sum\limits_{j\neq i}\mu_i(w_{i,x}-(w_i/v_i)v_{i,x})\zeta_{j,xx}\xi_i(\pa_{\zeta_j}\bar{v}_i)\frac{w_i}{v_i}\tilde{r}_{i,v}\\
%	&+\sum\limits_{j\neq i}\mu_i(w_{i,x}-(w_i/v_i)v_{i,x})\zeta_{j,x}\xi_i^\p\left(\frac{w_i}{v_i}\right)_x(\pa_{\zeta_j}\bar{v}_i)\frac{w_i}{v_i}\tilde{r}_{i,v}\\
%	&+\sum\limits_{j\neq i}\mu_i(w_{i,x}-(w_i/v_i)v_{i,x})\zeta_{j,x}\xi_i(\pa_{\zeta_j}\bar{v}_i)_x\frac{w_i}{v_i}\tilde{r}_{i,v}\\
%	&+\sum\limits_{j\neq i}\mu_i(w_{i,x}-(w_i/v_i)v_{i,x})\zeta_{j,x}\xi_i(\pa_{\zeta_j}\bar{v}_i)\left(\frac{w_i}{v_i}\right)_x\tilde{r}_{i,v}\\
%	&+\sum\limits_{j\neq i}\mu_i(w_{i,x}-(w_i/v_i)v_{i,x})\zeta_{j,x}\xi_i(\pa_{\zeta_j}\bar{v}_i)\frac{w_i}{v_i}[\tilde{r}_{i,uv}u_x+(\xi\bar{v}_i)_x\tilde{r}_{i,vv}-\theta_i^\p\left(\frac{w_i}{v_i}\right)_x\tilde{r}_{i,v\si}].
%\end{align*}
\begin{align*}
	\widetilde{\B}^{1,14}_{i,x}&=-2\mu_{i,x}v_i\theta_i^\p\left(\frac{w_i}{v_i}\right)^2_x\tilde{r}_{i,\si}-2\mu_i v_{i,x}\theta_i^\p\left(\frac{w_i}{v_i}\right)^2_x\tilde{r}_{i,\si}-2\mu_iv_i\theta_i^{\p\p}\left(\frac{w_i}{v_i}\right)^3_x\tilde{r}_{i,\si}\\
	&-2\mu_iv_i\theta_i^\p\left(\frac{w_i}{v_i}\right)_{xx}\left(\frac{w_i}{v_i}\right)_{x}\tilde{r}_{i,\si}-2\mu_iv_i\theta_i^\p\left(\frac{w_i}{v_i}\right)^2_x\pa_x(\tilde{r}_{i,\si})=O(1)(\delta_0^2\La_i^3+\La_i^5).
\end{align*}
Similarly we obtain $\widetilde{\B}^{1,15}_{i}=-\mu_i w_i\left(\frac{w_i}{v_i}\right)^2_x\theta_i^{\p\p} \tilde{r}_{i,\si}=\delta_0^2\La_i^3$ and
\begin{align*}
	&\widetilde{\B}^{1,15}_{i,x}=-\mu_{i,x}w_i\left(\frac{w_i}{v_i}\right)^2_x\theta_i^{\p\p} \tilde{r}_{i,\si}-\mu_iw_{i,x}\left(\frac{w_i}{v_i}\right)^2_x\theta_i^{\p\p} \tilde{r}_{i,\si}-2\mu_i w_i\left(\frac{w_i}{v_i}\right)_{xx}\left(\frac{w_i}{v_i}\right)_{x}\theta_i^{\p\p} \tilde{r}_{i,\si}\\
	&-\mu_i w_i\left(\frac{w_i}{v_i}\right)^3_x\theta_i''' \tilde{r}_{i,\si}-\mu_i w_i\left(\frac{w_i}{v_i}\right)^2_x\theta_i^{\p\p}\pa_x( \tilde{r}_{i,\si})=O(1)(\La_i^3+\La_i^5).
\end{align*}
We get from \eqref{ngtech4} $\widetilde{\B}^{1,16}_{i}=\xi_i w_i v_i\pa_{v_i}\bar{v}_i \mathcal{H}_i \tilde{r}_{i,uv}\tilde{r}_i=O(1)\sum_j(\La_j^1+\La_j^4+\La_j^6)+R_{\e,1}^{1,i,16}$ with $R_{\e,1}^{1,i,11}$ satisfying \eqref{12.11bis2}. We have from Lemmas \ref{estimimpo1}, \ref{lemme6.5}, \ref{lemme6.6}, \ref{lemme11.4} and \eqref{ngtech4}, \eqref{estimate-v-i-xx-1aa} 
\begin{align*}
\widetilde{\B}^{1,16}_{i,x}&=\big[\xi_i'\left(\frac{w_i}{v_i}\right)_x w_i v_i\pa_{v_i}\bar{v}_i +\xi_i (w_{i,x} v_i+w_i v_{i,x})\pa_{v_i}\bar{v}_i +\xi_i w_i v_i\pa_x(\pa_{v_i}\bar{v}_i )\big]
\mathcal{H}_i \tilde{r}_{i,uv}\tilde{r}_i\\
%&+\xi_i (w_{i,x} v_i+w_i v_{i,x})\pa_{v_i}\bar{v}_i (\mu_i v_{i,x}-v_i(\tilde{\lambda}_i-\lambda_i^*)- w_i) \tilde{r}_{i,uv}\tilde{r}_i\\
%&+\xi_i w_i v_i\pa_x(\pa_{v_i}\bar{v}_i )(\mu_i v_{i,x}-v_i(\tilde{\lambda}_i-\lambda_i^*)- w_i) \tilde{r}_{i,uv}\tilde{r}_i\\
&+\xi_i w_i v_i\pa_{v_i}\bar{v}_i \big[\mathcal{H}_{i,x} \tilde{r}_{i,uv}\tilde{r}_i+\mathcal{H}_i \pa_x( \tilde{r}_{i,uv}\tilde{r}_i)\big]
%&+\xi_i w_i v_i\pa_{v_i}\bar{v}_i (\mu_i v_{i,x}-v_i(\tilde{\lambda}_i-\lambda_i^*)- w_i)\pa_x( \tilde{r}_{i,uv}\tilde{r}_i)\\
=O(1)\sum_j(\La_j^1+\La_j^4+\La_j^5+\La_j^6)+\widetilde{R}_{\e,1}^{1,i,16},
\end{align*}
with $\widetilde{R}_{\e,1}^{1,i,16}$ satisfying \eqref{12.11bis2}.
We have now  $\widetilde{\B}^{1,17}_{i}=\sum_{i\ne j}\xi_i w_i v_j\pa_{v_i}\bar{v}_i \mathcal{H}_i \tilde{r}_{i,uv}\tilde{r}_j=O(1)\La_i^1$ and from Lemmas \ref{estimimpo1}, \ref{lemme6.5}, \ref{lemme6.6}, \ref{lemme11.4}
\begin{align*}
\widetilde{\B}^{1,17}_{i,x}&=\sum_{i\ne j}\big[\xi_i' \left(\frac{w_i}{v_i}\right)_x w_i v_j\pa_{v_i}\bar{v}_i 
+\xi_i (w_{i,x}v_j+w_iv_{j,x})\pa_{v_i}\bar{v}_i +\xi_i w_i v_j\pa_x(\pa_{v_i}\bar{v}_i )\big]
\mathcal{H}_i \tilde{r}_{i,uv}\tilde{r}_j\\
%&+\sum_{i\ne j}\xi_i (w_{i,x}v_j+w_iv_{j,x})\pa_{v_i}\bar{v}_i (\mu_i v_{i,x}-v_i(\tilde{\lambda}_i-\lambda_i^*)-w_i) \tilde{r}_{i,uv}\tilde{r}_j\\
%&+\sum_{i\ne j}\xi_i w_i v_j\pa_x(\pa_{v_i}\bar{v}_i )(\mu_i v_{i,x}-v_i(\tilde{\lambda}_i-\lambda_i^*)-w_i) \tilde{r}_{i,uv}\tilde{r}_j\\
&+\sum_{i\ne j}\xi_i w_i v_j\pa_{v_i}\bar{v}_i \big[\mathcal{H}_{i,x} \tilde{r}_{i,uv}\tilde{r}_j+\mathcal{H}_i \pa_x( \tilde{r}_{i,uv}\tilde{r}_j)\big]%\\
%&+\sum_{i\ne j}\xi_i w_i v_j\pa_{v_i}\bar{v}_i (\mu_i v_{i,x}-v_i(\tilde{\lambda}_i-\lambda_i^*)-w_i)\pa_x( \tilde{r}_{i,uv}\tilde{r}_j)\\
%&
=O(1)\La_i^1.
\end{align*}
We obtain now $\widetilde{\B}^{1,18}_{i}=
\mu_i v_iw_i\left(\frac{w_i}{v_i}\right)_x \xi_i^\p\bar{v}_i\tilde{r}_{i,uv}\tilde{r}_i=O(1)\La_i^4$ and from Lemmas \ref{estimimpo1}, \ref{lemme6.5}, \ref{lemme6.6}, \ref{lemme11.4}
\begin{align*}
	\widetilde{\B}^{1,18}_{i,x}&=\big[\mu_{i,x}v_iw_i\left(\frac{w_i}{v_i}\right)_x+\mu_i (v_{i,x}w_i+v_i w_{i,x})\left(\frac{w_i}{v_i}\right)_x+\mu_i v_iw_i\left(\frac{w_i}{v_i}\right)_{xx} \big]
	\xi_i^\p\bar{v}_i\tilde{r}_{i,uv}\tilde{r}_i%+\mu_i (v_{i,x}w_i+v_i w_{i,x})\left(\frac{w_i}{v_i}\right)_x\xi_i^\p\bar{v}_i\tilde{r}_{i,uv}\tilde{r}_i\\
	%&+\mu_i v_iw_i\left(\frac{w_i}{v_i}\right)_{xx} \xi_i^\p\bar{v}_i\tilde{r}_{i,uv}\tilde{r}_i
	\\
	&+\mu_i v_iw_i\left(\frac{w_i}{v_i}\right)_x\big[ \xi_i''\left(\frac{w_i}{v_i}\right)_x \bar{v}_i\tilde{r}_{i,uv}\tilde{r}_i+\xi_i^\p\pa_x(\bar{v}_i)\tilde{r}_{i,uv}\tilde{r}_i+\xi_i^\p\bar{v}_i\pa_x(\tilde{r}_{i,uv}\tilde{r}_i)\big]\\
	%&+\mu_i v_iw_i\left(\frac{w_i}{v_i}\right)_x\xi_i^\p\pa_x(\bar{v}_i)\tilde{r}_{i,uv}\tilde{r}_i+\mu_{i,x}v_iw_i\left(\frac{w_i}{v_i}\right)_x\xi_i^\p\bar{v}_i\pa_x(\tilde{r}_{i,uv}\tilde{r}_i)\\
	&=O(1) (\delta_0^2\La_i^3+\La_i^4+\La_i^5).
	\end{align*}
	Next we have $\widetilde{\B}^{1,19}_{i}=\sum\limits_{j\neq i}\mu_iv_jw_i\left(\frac{w_i}{v_i}\right)_x \xi_i^\p\bar{v}_i\tilde{r}_{i,uv}\tilde{r}_i=O(1)\La_i^1$ that we can treat similarly. We have now  from \eqref{ngtech4} and Lemmas \ref{estimimpo1}, \ref{lemme6.5}, \ref{lemme6.6}, \ref{lemme11.3}
	$\widetilde{\B}^{1,20}_{i}= {\color{black}{\xi_i^2}}w_iv_{i,x}(\pa_{v_i}\bar{v}_i)^2 \mathcal{H}_i\tilde{r}_{i,vv}=O(1)\sum_j(\La_j^1+\La_j^4+\La_j^6)+R_{\e,1}^{1,i,20}$ with $R_{\e,1}^{1,i,20}$ satisfying \eqref{12.11bis2}. Using again Lemmas \ref{estimimpo1}, \ref{lemme6.5}, \ref{lemme6.6}, \ref{lemme11.3} and \eqref{ngtech4}, \eqref{estimate-v-i-xx-1aa} 
	\begin{align*}
	\widetilde{\B}^{1,20}_{i,x}&=\big[2\xi_i \xi_i' \left(\frac{w_i}{v_i}\right)_x w_iv_{i,x}(\pa_{v_i}\bar{v}_i)^2+\xi_i^2 (w_{i,x}v_{i,x}+w_{i}v_{i,xx})(\pa_{v_i}\bar{v}_i)^2+2 \xi_i^2 w_iv_{i,x} (\pa_{v_i}\bar{v}_i)_x \pa_{v_i}\bar{v}_i
	\big]\mathcal{H}_i\tilde{r}_{i,vv}\\
	%&+\xi_i(w_{i,x}v_{i,x}+w_{i}v_{i,xx})(\pa_{v_i}\bar{v}_i)^2 (\mu_iv_{i,x}-v_i(\tilde{\lambda}_i-\lambda_i^*)-w_i)\tilde{r}_{i,vv}\\
	%&+\xi_iw_iv_{i,x}2\pa_x (\pa_{v_i}\bar{v}_i)\pa_{v_i}\bar{v}_i (\mu_iv_{i,x}-v_i(\tilde{\lambda}_i-\lambda_i^*)-w_i)\tilde{r}_{i,vv}\\
	%&+\xi_iw_iv_{i,x}(\pa_{v_i}\bar{v}_i)^2 (\mu_iv_{i,x}-v_i(\tilde{\lambda}_i-\lambda_i^*)-w_i)_x\tilde{r}_{i,vv}\\
	&+\xi_i^2 w_iv_{i,x}(\pa_{v_i}\bar{v}_i)^2 \big[\mathcal{H}_{i,x} \tilde{r}_{i,vv}+\mathcal{H}_i\pa_x(\tilde{r}_{i,vv})\big]=O(1)\sum_j(\La_j^1+\La_j^4+\La_j^5+\La_j^6)+\widetilde{R}_{\e,1}^{1,i,20},
	\end{align*}
	 with $\widetilde{R}_{\e,1}^{1,i,20}$ satisfying \eqref{12.11bis2}. 
	We have now  $\widetilde{\B}^{1,21}_{i}=\sum_{i\ne j}{\color{black}{\xi_i^2}}v_{j,x}w_i \pa_{v_i}\bar{v}_i\pa_{v_j}\bar{v}_i\mathcal{H}_i\tilde{r}_{i,vv}=O(1)\La_i^1$ and from Lemmas \ref{estimimpo1}, \ref{lemme6.5}, \ref{lemme6.6}, \ref{lemme11.3}
	\begin{align*}
	&\widetilde{\B}^{1,21}_{i,x}=\sum_{i\ne j}\big[2 \xi_i'\xi_i\left(\frac{w_i}{v_i}\right)_x v_{j,x}w_i +\xi_i^2(v_{j,xx}w_i+v_{j,x}w_{i,x})\big]\pa_{v_i}\bar{v}_i\pa_{v_j}\bar{v}_i\mathcal{H}_i\tilde{r}_{i,vv}\\
	%&+\sum_{i\ne j}\xi_i(w_{j,xx}w_i+w_{j,x}w_{i,x}) \pa_{v_i}\bar{v}_i\pa_{w_j}\bar{v}_i(\mu_iv_{i,x}-v_i(\tilde{\lambda}_i-\lambda_i^*)-w_i)\tilde{r}_{i,vv}\\
	&+\sum_{i\ne j}\xi_i^2 v_{j,x}w_i\big[( \pa_x(\pa_{v_i}\bar{v}_i)\pa_{v_j}\bar{v}_i+\pa_{v_i}\bar{v}_i\pa_x(\pa_{v_j}\bar{v}_i))\mathcal{H}_i\tilde{r}_{i,vv}+\pa_{v_i}\bar{v}_i\pa_{v_j}\bar{v}_i(\mathcal{H}_{i,x}\tilde{r}_{i,vv}+\mathcal{H}_i\pa_x(\tilde{r}_{i,vv})\big]\\
	%\\
	%&+\sum_{i\ne j}\xi_i w_{j,x}w_i \pa_{v_i}\bar{v}_i\pa_{w_j}\bar{v}_i(\mu_iv_{i,x}-v_i(\tilde{\lambda}_i-\lambda_i^*)-w_i)_x\tilde{r}_{i,vv}\\
	%&+\sum_{i\ne j}\xi_i w_{j,x}w_i \pa_{v_i}\bar{v}_i\pa_{w_j}\bar{v}_i(\mu_iv_{i,x}-v_i(\tilde{\lambda}_i-\lambda_i^*)-w_i)\pa_x(\tilde{r}_{i,vv})\\
	&=O(1)\sum_j\La_j^1.
	\end{align*}
	We have now $\widetilde{\B}^{1,22}_{i}=\xi_i' \bar{v}_i w_i \xi_i \pa_{v_i}\bar{v}_i\left(\frac{w_i}{v_i}\right)_x\mathcal{H}_i\tilde{r}_{i,vv}%=
	%\xi_i' \bar{v}_i w_i \xi_i \pa_{v_i}\bar{v}_i\left(\frac{w_i}{v_i}\right)_x(\mu_iv_{i,x}-v_i(\tilde{\lambda}_i-\lambda_i^*)-w_i)\tilde{r}_{i,vv}
	=O(1)\La_i^4$ and we get using Lemmas \ref{estimimpo1}, \ref{lemme6.5}, \ref{lemme6.6}, \ref{lemme11.3} and \eqref{ngtech5}
\begin{align*}
	&\widetilde{\B}^{1,22}_{i,x}=\big[\xi_i''\left(\frac{w_i}{v_i}\right)_x \bar{v}_i w_i \xi_i +\xi_i' (\pa_x(\bar{v}_i )w_i+\bar{v}_i w_{i,x}) \xi_i +\xi_i'^2 \bar{v}_i w_i \left(\frac{w_i}{v_i}\right)_{x} \big]
	\pa_{v_i}\bar{v}_i \left(\frac{w_i}{v_i}\right)_x \mathcal{H}_i\tilde{r}_{i,vv}
	\\
	%&+\xi_i' (\pa_x(\bar{v}_i )w_i+\bar{v}_i w_{i,x}) \xi_i \pa_{v_i}\bar{v}_i\left(\frac{w_i}{v_i}\right)_x(\mu_iv_{i,x}-v_i(\tilde{\lambda}_i-\lambda_i^*)-w_i)\tilde{r}_{i,vv}\\
	%&+\xi_i'^2 \bar{v}_i w_i \left(\frac{w_i}{v_i}\right)_{x}^2 \pa_{v_i}\bar{v}_i(\mu_iv_{i,x}-v_i(\tilde{\lambda}_i-\lambda_i^*)-w_i)\tilde{r}_{i,vv}\\
	&+\xi_i' \bar{v}_i w_i \xi_i\big[(\pa_x( \pa_{v_i}\bar{v}_i) \left(\frac{w_i}{v_i}\right)_x+ \pa_{v_i}\bar{v}_i\left(\frac{w_i}{v_i}\right)_{xx})\mathcal{H}_i \tilde{r}_{i,vv}+ \pa_{v_i}\bar{v}_i\left(\frac{w_i}{v_i}\right)_{x}(\mathcal{H}_{i,x} \tilde{r}_{i,vv}+\mathcal{H}_i\pa_x( \tilde{r}_{i,vv}))\big]\\
	%&+\xi_i' \bar{v}_i w_i \xi_i \pa_{v_i}\bar{v}_i\left(\frac{w_i}{v_i}\right)_{xx}(\mu_iv_{i,x}-v_i(\tilde{\lambda}_i-\lambda_i^*)-w_i)\tilde{r}_{i,vv}\\
	%&+\xi_i' \bar{v}_i w_i \xi_i \pa_{v_i}\bar{v}_i\left(\frac{w_i}{v_i}\right)_x(\mu_iv_{i,x}-v_i(\tilde{\lambda}_i-\lambda_i^*)-w_i)_x\tilde{r}_{i,vv}\\
	%&+\xi_i' \bar{v}_i w_i \xi_i \pa_{v_i}\bar{v}_i\left(\frac{w_i}{v_i}\right)_x(\mu_iv_{i,x}-v_i(\tilde{\lambda}_i-\lambda_i^*)-w_i)\pa_x(\tilde{r}_{i,vv})\\
	&=O(1)\sum_j(\La_j^1+\delta_0^2\La_j^3+\La_j^4+\La_j^5).
\end{align*}	
We get now $\widetilde{\B}^{1,23}_{i}= \mu_i w_i v_{i,x}\left(\frac{w_i}{v_i}\right)_x \xi_i^\p\xi_i\pa_{v_i}\bar{v}_i\bar{v}_i\tilde{r}_{i,vv}=O(1)\La_i^4$ and using Lemma \ref{lemme11.3} and the fact that $\xi'_i=0$ on $\{| \frac{w_i}{v_i}|\leq\frac{\delta_1}{2}\}$
%\begin{align*}
%	\B^{1,22}_{i,x}&=\sum\limits_{j\neq i}\mu_{i,x}(w_{i,x}-(w_i/v_i)v_{i,x})v_j\frac{w_i}{v_i}\xi_i^\p\bar{v}_i\tilde{r}_{i,uv}\tilde{r}_j\\
%	&+\sum\limits_{j\neq i}\mu_i(w_{i,xx}-(w_i/v_i)v_{i,xx}-v_{i,x}(w_i/v_i)_x)v_j\frac{w_i}{v_i}\xi_i^\p\bar{v}_i\tilde{r}_{i,uv}\tilde{r}_j\\
%	&+\sum\limits_{j\neq i}\mu_i(w_{i,x}-(w_i/v_i)v_{i,x})v_{j,x}\frac{w_i}{v_i}\xi_i^\p\bar{v}_i\tilde{r}_{i,uv}\tilde{r}_j\\
%	&+\sum\limits_{j\neq i}\mu_i(w_{i,x}-(w_i/v_i)v_{i,x})v_j\left(\frac{w_i}{v_i}\right)_x\xi_i^\p\bar{v}_i\tilde{r}_{i,uv}\tilde{r}_j\\
%	&+\sum\limits_{j\neq i}\mu_i(w_{i,x}-(w_i/v_i)v_{i,x})v_j\frac{w_i}{v_i}\xi_i^{\p\p}\left(\frac{w_i}{v_i}\right)_x\bar{v}_i\tilde{r}_{i,uv}\tilde{r}_j\\
%	&+\sum\limits_{j\neq i}\mu_i(w_{i,x}-(w_i/v_i)v_{i,x})v_j\frac{w_i}{v_i}\xi_i^\p\bar{v}_{i,x}\tilde{r}_{i,uv}\tilde{r}_j\\
%	&+\sum\limits_{j\neq i}\mu_i(w_{i,x}-(w_i/v_i)v_{i,x})v_j\frac{w_i}{v_i}\xi_i^\p \bar{v}_i
%	\left[\tilde{r}_{i,uv}(\tilde{r}_j\otimes u_x)+(\xi_i\bar{v}_i)_x\tilde{r}_{i,uvv}\tilde{r}_{i}-\theta_i^\p\left(\frac{w_i}{v_i}\right)_x\tilde{r}_{i,uv\si}\tilde{r}_i\right]\\
%	&+\sum\limits_{j\neq i}\mu_i(w_{i,x}-(w_i/v_i)v_{i,x})v_j\frac{w_i}{v_i}\xi_i^\p\bar{v}_i\left[\tilde{r}_{i,uv}(\tilde{r}_{j,u} u_x)+(\xi_j\bar{v}_j)_x\tilde{r}_{i,uv}\tilde{r}_{j,v}-\theta_j^\p\left(\frac{w_j}{v_j}\right)_x\tilde{r}_{i,uv}\tilde{r}_{j,\si}\right].
%\end{align*}
\begin{align*}
	&\widetilde{\B}^{1,23}_{i,x}=\big[(\mu_{i,x} w_i +\mu_i w_{i,x})v_{i,x}\left(\frac{w_i}{v_i}\right)_x 
	+\mu_i w_i( v_{i,xx}\left(\frac{w_i}{v_i}\right)_x+ v_{i,x}\left(\frac{w_i}{v_i}\right)_{xx})\big]
	\xi_i^\p\xi_i\pa_{v_i}\bar{v}_i\bar{v}_i\tilde{r}_{i,vv}\\
	%+\mu_i w_{i,x} v_{i,x}\left(\frac{w_i}{v_i}\right)_x \xi_i^\p\xi_i\pa_{v_i}\bar{v}_i\bar{v}_i\tilde{r}_{i,vv}\\
	%&+\mu_i w_i v_{i,xx}\left(\frac{w_i}{v_i}\right)_x \xi_i^\p\xi_i\pa_{v_i}\bar{v}_i\bar{v}_i\tilde{r}_{i,vv}
	%+\mu_i w_i v_{i,x}\left(\frac{w_i}{v_i}\right)_{xx} \xi_i^\p\xi_i\pa_{v_i}\bar{v}_i\bar{v}_i\tilde{r}_{i,vv}\\
	&+\mu_i w_i v_{i,x}\left(\frac{w_i}{v_i}\right)_x\big[ (\xi_i''\xi_i+(\xi_i')^2)\left(\frac{w_i}{v_i}\right)_x
	\pa_{v_i}\bar{v}_i\bar{v}_i\tilde{r}_{i,vv}+ \xi_i^\p\xi_i(\pa_x(\pa_{v_i}\bar{v}_i)\bar{v}_i+\pa_{v_i}\bar{v}_i\pa_x(\bar{v}_i))\tilde{r}_{i,vv} \big]\\
	%&
	%+\mu_i w_i v_{i,x}\left(\frac{w_i}{v_i}\right)^2_x (\xi_i')^2\pa_{v_i}\bar{v}_i\bar{v}_i\tilde{r}_{i,vv}\\
	%&+\mu_i w_i v_{i,x}\left(\frac{w_i}{v_i}\right)_x \xi_i^\p\xi_i\pa_x(\pa_{v_i}\bar{v}_i)\bar{v}_i\tilde{r}_{i,vv}+\mu_i w_i v_{i,x}\left(\frac{w_i}{v_i}\right)_x \xi_i^\p\xi_i\pa_{v_i}\bar{v}_i\pa_x(\bar{v}_i)\tilde{r}_{i,vv}\\
	&+\mu_i w_i v_{i,x}\left(\frac{w_i}{v_i}\right)_x \xi_i^\p\xi_i\pa_{v_i}\bar{v}_i\bar{v}_i\pa_x(\tilde{r}_{i,vv})=O(1)(\La_i^1+\delta_0^2\La_i^3+\La_i^4+\La_i^5+\La_i^6+\La_i^{6,1}).
	%\mu_{i,x}(w_{i,x}-(w_i/v_i)v_{i,x})v_{i,x}\frac{w_i}{v_i}\xi_i^\p\xi_i\pa_{v_i}\bar{v}_i\bar{v}_i\tilde{r}_{i,vv}\\
	%&+\mu_i(w_{i,xx}-(w_i/v_i)v_{i,xx}-v_{i,x}(w_i/v_i)_x)v_{i,x}\frac{w_i}{v_i}\xi_i^\p\xi_i\pa_{v_i}\bar{v}_i\bar{v}_i\tilde{r}_{i,vv}\\
	%&+\mu_i(w_{i,x}-(w_i/v_i)v_{i,x})v_{i,xx}\frac{w_i}{v_i}\xi_i^\p\xi_i\pa_{v_i}\bar{v}_i\bar{v}_i\tilde{r}_{i,vv}\\
	%&+\mu_i(w_{i,x}-(w_i/v_i)v_{i,x})v_{i,x}\left(\frac{w_i}{v_i}\right)_x\xi_i^\p\xi_i\pa_{v_i}\bar{v}_i\bar{v}_i\tilde{r}_{i,vv}\\
	%&+\mu_i(w_{i,x}-(w_i/v_i)v_{i,x})v_{i,x}\frac{w_i}{v_i}\xi_i^{\p\p}\left(\frac{w_i}{v_i}\right)_x\xi_i\pa_{v_i}\bar{v}_i\bar{v}_i\tilde{r}_{i,vv}\\
	%&+\mu_i(w_{i,x}-(w_i/v_i)v_{i,x})v_{i,x}\frac{w_i}{v_i}(\xi_i^\p)^2\left(\frac{w_i}{v_i}\right)_x\pa_{v_i}\bar{v}_i\bar{v}_i\tilde{r}_{i,vv}\\
	%&+\mu_i(w_{i,x}-(w_i/v_i)v_{i,x})v_{i,x}\frac{w_i}{v_i}\xi_i^\p\xi_i(\pa_{v_i}\bar{v}_i)_x\bar{v}_i\tilde{r}_{i,vv}\\
	%&+\mu_i(w_{i,x}-(w_i/v_i)v_{i,x})v_{i,x}\frac{w_i}{v_i}\xi_i^\p\xi_i\pa_{v_i}\bar{v}_i\bar{v}_{i,x}\tilde{r}_{i,vv}\\
	%&+\mu_i(w_{i,x}-(w_i/v_i)v_{i,x})v_{i,x}\frac{w_i}{v_i}\xi_i^\p\xi_i\pa_{v_i}\bar{v}_i\bar{v}_i\left[\tilde{r}_{i,uvv} u_x+(\xi_i\bar{v}_i)_x\tilde{r}_{i,vvv}\tilde{r}_{i}-\theta_i^\p\left(\frac{w_i}{v_i}\right)_x\tilde{r}_{i,vv\si}\tilde{r}_i\right].
\end{align*}
We have now $\widetilde{\B}^{1,24}_{i}= \sum\limits_{j\neq i}\mu_i w_i\left(\frac{w_i}{v_i}\right)_x v_{j,x}\xi_i^\p\xi_i\bar{v}_i\pa_{v_j}\bar{v}_i\tilde{r}_{i,vv}=
O(1)\La_i^1$ and from Lemma \ref{lemme11.3}
\begin{align*}
	&\widetilde{\B}^{1,24}_{i,x}=\sum\limits_{j\neq i}\big[\mu_{i,x} w_i\left(\frac{w_i}{v_i}\right)_x+\mu_i w_{i,x}\left(\frac{w_i}{v_i}\right)_x+\mu_i w_i\left(\frac{w_i}{v_i}\right)_{xx}  \big]
	v_{j,x} \xi_i^\p\xi_i\bar{v}_i\pa_{v_j}\bar{v}_i\tilde{r}_{i,vv}\\%+\sum_{i\ne j}\mu_i w_{i,x}\left(\frac{w_i}{v_i}\right)_xv_{j,x}  \xi_i^\p\xi_i\bar{v}_i\pa_{w_j}\bar{v}_i\tilde{r}_{i,vv}\\
	&+\sum\limits_{j\neq i}\mu_i w_i\left(\frac{w_i}{v_i}\right)_x \big[(v_{j,xx}\bar{v}_i +v_{j,x}  \pa_x(\bar{v}_i)  )\xi_i^\p\xi_i+v_{j,x}( \xi_i'' \xi_i +(\xi_i')^2)\left(\frac{w_i}{v_i}\right)_x \bar{v}_i \big]\pa_{v_j}  \bar{v}_i\tilde{r}_{i,vv}\\
	%&+\sum\limits_{j\neq i}\mu_i w_i\left(\frac{w_i}{v_i}\right)^2_x v_{j,x} \xi_i'' \xi_i\bar{v}_i\pa_{w_j}\bar{v}_i\tilde{r}_{i,vv}+\sum\limits_{j\neq i}\mu_i w_i\left(\frac{w_i}{v_i}\right)^2_x v_{j,x}  (\xi_i^\p)^2\bar{v}_i\pa_{w_j}\bar{v}_i\tilde{r}_{i,vv}\\
	&%+\sum\limits_{j\neq i}\mu_i w_i\left(\frac{w_i}{v_i}\right)_x v_{j,x} \xi_i^\p\xi_i\pa_x(\bar{v}_i)\pa_{w_j}\bar{v}_i\tilde{r}_{i,vv}
	+\sum\limits_{j\neq i}\mu_i w_i\left(\frac{w_i}{v_i}\right)_x v_{j,x} \xi_i^\p\xi_i\bar{v}_i\big[\pa_x(\pa_{v_j}\bar{v}_i)\tilde{r}_{i,vv}+\bar{v}_i\pa_{v_j}\bar{v}_i\pa_x(\tilde{r}_{i,vv})\big] 
	%&+\sum\limits_{j\neq i}\mu_i w_i\left(\frac{w_i}{v_i}\right)_x v_{j,x} \xi_i^\p\xi_i\bar{v}_i\pa_{w_j}\bar{v}_i\pa_x(\tilde{r}_{i,vv})\\
	%&+\sum\limits_{j\neq i}\mu_i w_i\left(\frac{w_i}{v_i}\right)_x \xi_i^\p\xi_i\bar{v}_i\pa_{w_j}\bar{v}_i\pa_x(\tilde{r}_{i,vv})\\
	=O(1)(\La_i^1+\delta_0^2\La_i^3).
\end{align*}
We have now
$\widetilde{\B}^{1,25}_{i}=\mu_i w_i \left(\frac{w_i}{v_i}\right)_x^2 (\bar{v}_i\xi_i^\p)^2\tilde{r}_{i,vv}=O(1)\delta_0^2 \La_i^3$ and 
%\begin{align*}
%	\B^{1,16}_{i,x}&=\sum\limits_{j\neq i}\mu_{i,x}(w_{i,x}-(w_i/v_i)v_{i,x})\zeta_{j,x}\frac{w_i}{v_i}\xi_i^\p\xi_i\bar{v}_i\pa_{\zeta_j}\bar{v}_i\tilde{r}_{i,vv}\\
%	&+\sum\limits_{j\neq i}\mu_i(w_{i,xx}-(w_i/v_i)v_{i,xx}-v_{i,x}(w_i/v_i)_x)\zeta_{j,x}\frac{w_i}{v_i}\xi_i^\p\xi_i\bar{v}_i\pa_{\zeta_j}\bar{v}_i\tilde{r}_{i,vv}\\
%	&+\sum\limits_{j\neq i}\mu_i(w_{i,x}-(w_i/v_i)v_{i,x})\zeta_{j,xx}\frac{w_i}{v_i}\xi_i^\p\xi_i\bar{v}_i\pa_{\zeta_j}\bar{v}_i\tilde{r}_{i,vv}\\
%	&+\sum\limits_{j\neq i}\mu_i(w_{i,x}-(w_i/v_i)v_{i,x})\zeta_{j,x}\left(\frac{w_i}{v_i}\right)_x\xi_i^\p\xi_i\bar{v}_i\pa_{\zeta_j}\bar{v}_i\tilde{r}_{i,vv}\\
%	&+\sum\limits_{j\neq i}\mu_i(w_{i,x}-(w_i/v_i)v_{i,x})\zeta_{j,x}\frac{w_i}{v_i}\xi_i^{\p\p}\left(\frac{w_i}{v_i}\right)_x\xi_i\bar{v}_i\pa_{\zeta_j}\bar{v}_i\tilde{r}_{i,vv}\\
%	&+\sum\limits_{j\neq i}\mu_i(w_{i,x}-(w_i/v_i)v_{i,x})\zeta_{j,x}\frac{w_i}{v_i}(\xi_i^\p)^2\left(\frac{w_i}{v_i}\right)_x\bar{v}_i\pa_{\zeta_j}\bar{v}_i\tilde{r}_{i,vv}\\
%	&+\sum\limits_{j\neq i}\mu_i(w_{i,x}-(w_i/v_i)v_{i,x})\zeta_{j,x}\frac{w_i}{v_i}\xi_i^\p\xi_i\bar{v}_{i,x}\pa_{\zeta_j}\bar{v}_i\tilde{r}_{i,vv}\\
%	&+\sum\limits_{j\neq i}\mu_i(w_{i,x}-(w_i/v_i)v_{i,x})\zeta_{j,x}\frac{w_i}{v_i}\xi_i^\p\xi_i\bar{v}_i(\pa_{\zeta_j}\bar{v}_i)_x\tilde{r}_{i,vv}\\
%	&+\sum\limits_{j\neq i}\mu_i(w_{i,x}-(w_i/v_i)v_{i,x})\zeta_{j,x}\frac{w_i}{v_i}\xi_i^\p\xi_i\bar{v}_i\pa_{\zeta_j}\bar{v}_i\left[\tilde{r}_{i,uvv} u_x+(\xi_i\bar{v}_i)_x\tilde{r}_{i,vvv}\tilde{r}_{i}-\theta_i^\p\left(\frac{w_i}{v_i}\right)_x\tilde{r}_{i,vv\si}\tilde{r}_i\right].
%\end{align*}
\begin{align*}
	\widetilde{\B}^{1,25}_{i,x}&=\big[(\mu_{i,x} w_i +\mu_i w_{i,x})\left(\frac{w_i}{v_i}\right)_x^2 +
	2\mu_i w_i \left(\frac{w_i}{v_i}\right)_x\left(\frac{w_i}{v_i}\right)_{xx}\big]
	(\bar{v}_i\xi_i^\p)^2\tilde{r}_{i,vv}\\
	%+\mu_i w_{i,x} \left(\frac{w_i}{v_i}\right)_x^2 (\bar{v}_i\xi_i^\p)^2\tilde{r}_{i,vv}\\
	&%+2\mu_i w_i \left(\frac{w_i}{v_i}\right)_x\left(\frac{w_i}{v_i}\right)_{xx} (\bar{v}_i\xi_i^\p)^2\tilde{r}_{i,vv}
	+\mu_i w_i \left(\frac{w_i}{v_i}\right)_x^2[2 (\bar{v}_i)^2 \xi_i^\p\xi_i'' \left(\frac{w_i}{v_i}\right)_x \tilde{r}_{i,vv}+2\bar{v}_i \pa_x(\bar{v}_i) (\xi_i^\p)^2\tilde{r}_{i,vv}+(\bar{v}_i\xi_i^\p)^2\pa_x(\tilde{r}_{i,vv})\big] \\
	%&+2\mu_i w_i \left(\frac{w_i}{v_i}\right)_x^2 \bar{v}_i \pa_x(\bar{v}_i) (\xi_i^\p)^2\tilde{r}_{i,vv}
	%+\mu_i w_i \left(\frac{w_i}{v_i}\right)_x^2 (\bar{v}_i\xi_i^\p)^2\pa_x(\tilde{r}_{i,vv})\\
	&=O(1)(\delta_0^2\La_i^3+\La_i^5).
	%\mu_{i,x}(w_{i,x}-(w_i/v_i)v_{i,x})\frac{w_i}{v_i}(\bar{v}_i\xi_i^\p)^2\left(\frac{w_i}{v_i}\right)_x\tilde{r}_{i,vv}\\
	%&+\mu_i(w_{i,xx}-(w_i/v_i)v_{i,xx}-v_{i,x}(w_i/v_i)_x)\frac{w_i}{v_i}(\bar{v}_i\xi_i^\p)^2\left(\frac{w_i}{v_i}\right)_x\tilde{r}_{i,vv}\\
	%&+\mu_i(w_{i,x}-(w_i/v_i)v_{i,x})\left(\frac{w_i}{v_i}\right)_x(\bar{v}_i\xi_i^\p)^2\left(\frac{w_i}{v_i}\right)_x\tilde{r}_{i,vv}\\
	%&+2\mu_i(w_{i,x}-(w_i/v_i)v_{i,x})\frac{w_i}{v_i}(\bar{v}_{i}\xi_i^\p)(\bar{v}_i\xi_i^\p)_x\left(\frac{w_i}{v_i}\right)_x\tilde{r}_{i,vv}\\
	%&+\mu_i(w_{i,x}-(w_i/v_i)v_{i,x})\frac{w_i}{v_i}(\bar{v}_i\xi_i^\p)^2\left(\frac{w_i}{v_i}\right)_{xx}\tilde{r}_{i,vv}\\
	%&+\mu_i(w_{i,x}-(w_i/v_i)v_{i,x})\frac{w_i}{v_i}(\bar{v}_i\xi_i^\p)^2\left(\frac{w_i}{v_i}\right)_x\left[\tilde{r}_{i,uvv} u_x+(\xi_i\bar{v}_i)_x\tilde{r}_{i,vvv}\tilde{r}_{i}-\theta_i^\p\left(\frac{w_i}{v_i}\right)_x\tilde{r}_{i,vv\si}\tilde{r}_i\right].
\end{align*}
It yields now from Lemma \ref{lemme6.5} $\widetilde{\B}^{1,26}_{i}=-\mu_iw_i  \left(\frac{w_i}{v_i}\right)^2_x (\bar{v}_i\xi_i^\p)\theta_i^\p\tilde{r}_{i,v\si}=O(1)\delta_0^2\La_i^3$ and using Lemmas \ref{lemme6.5}, \ref{lemme6.6}, \ref{lemme11.3}
\begin{align*}
	\widetilde{\B}^{1,26}_{i,x}&=-\big[(\mu_{i,x}w_i+\mu_iw_{i,x})  \left(\frac{w_i}{v_i}\right)^2_x
	+2\mu_iw_i  \left(\frac{w_i}{v_i}\right)_{xx}  \left(\frac{w_i}{v_i}\right)_x\big]
 (\bar{v}_i\xi_i^\p)\theta_i^\p\tilde{r}_{i,v\si}\\%-\mu_iw_{i,x}  \left(\frac{w_i}{v_i}\right)^2_x (\bar{v}_i\xi_i^\p)\theta_i^\p\tilde{r}_{i,v\si}\\
	&%-2\mu_iw_i  \left(\frac{w_i}{v_i}\right)_{xx}  \left(\frac{w_i}{v_i}\right)_x(\bar{v}_i\xi_i^\p)\theta_i^\p\tilde{r}_{i,v\si}
	-\mu_iw_i  \left(\frac{w_i}{v_i}\right)^2_x\big[ \bar{v}_i(\xi_i''\theta_i^\p+\xi_i'\theta_i'')\left(\frac{w_i}{v_i}\right)_x\tilde{r}_{i,v\si}+ (\pa_x(\bar{v}_i)\tilde{r}_{i,v\si}+\bar{v}_i \pa_x(\tilde{r}_{i,v\si}))\xi_i^\p\theta_i^\p\big]\\
	%&-\mu_iw_i  \left(\frac{w_i}{v_i}\right)^2_x \pa_x(\bar{v}_i)\xi_i^\p\theta_i^\p\tilde{r}_{i,v\si}
	%-\mu_iw_i  \left(\frac{w_i}{v_i}\right)^3_x (\bar{v}_i\xi_i^\p)\theta_i''\tilde{r}_{i,v\si}\\
	%&-\mu_iw_i  \left(\frac{w_i}{v_i}\right)^2_x (\bar{v}_i\xi_i^\p)\theta_i^\p\pa_x(\tilde{r}_{i,v\si})\\
	&=O(1)(\delta_0^2\La_i^3+\La_i^5)
		%-\mu_{i,x}(w_{i,x}-(w_i/v_i)v_{i,x})\frac{w_i}{v_i}(\bar{v}_i\xi_i^\p)\theta_i^\p\left(\frac{w_i}{v_i}\right)_x\tilde{r}_{i,v\si}\\
	%&-\mu_i(w_{i,xx}-(w_i/v_i)v_{i,xx}-v_{i,x}(w_i/v_i)_x)\frac{w_i}{v_i}(\bar{v}_i\xi_i^\p)\theta_i^\p\left(\frac{w_i}{v_i}\right)_x\tilde{r}_{i,v\si}\\
	%&-\mu_i(w_{i,x}-(w_i/v_i)v_{i,x})\left(\frac{w_i}{v_i}\right)_x(\bar{v}_i\xi_i^\p)\theta_i^\p\left(\frac{w_i}{v_i}\right)_x\tilde{r}_{i,v\si}\\
	%&-\mu_i(w_{i,x}-(w_i/v_i)v_{i,x})\frac{w_i}{v_i}(\bar{v}_i\xi_i^\p)_x\theta_i^\p\left(\frac{w_i}{v_i}\right)_x\tilde{r}_{i,v\si}\\
	%&-\mu_i(w_{i,x}-(w_i/v_i)v_{i,x})\frac{w_i}{v_i}(\bar{v}_i\xi_i^\p)\theta_i^{\p\p}\left(\frac{w_i}{v_i}\right)^2_x\tilde{r}_{i,v\si}\\
	%&-\mu_i(w_{i,x}-(w_i/v_i)v_{i,x})\frac{w_i}{v_i}(\bar{v}_i\xi_i^\p)\theta_i^\p\left(\frac{w_i}{v_i}\right)_{xx}\tilde{r}_{i,v\si}\\
	%&-\mu_i(w_{i,x}-(w_i/v_i)v_{i,x})\frac{w_i}{v_i}(\bar{v}_i\xi_i^\p)\theta_i^\p\left(\frac{w_i}{v_i}\right)_x[\tilde{r}_{i,uv\si}u_x+(\xi\bar{v}_i)_x\tilde{r}_{i,vv\si}-\theta_i^\p\left(\frac{w_i}{v_i}\right)_x\tilde{r}_{i,v\si\si}].
\end{align*}
We have in particular used the fact that
\begin{align*}
-2\mu_iw_i  \left(\frac{w_i}{v_i}\right)_{xx}  \left(\frac{w_i}{v_i}\right)_x(\bar{v}_i\xi_i^\p)\theta_i^\p\tilde{r}_{i,v\si}&=O(1)|w_{i,xx}v_i-v_{i,xx}w_i|(1+\left(\frac{w_i}{v_i}\right)_x^2)\xi_i'\\
&+O(1) |v_{i,x}| |v_i \left(\frac{w_i}{v_i}\right)_x^2|\xi_i'=O(1)(\delta_0^2\La_i^3+\La_i^5).
\end{align*}
We have now from Lemma \ref{estimimpo1}, \ref{lemme6.5} and \eqref{ngtech5} $\widetilde{\B}^{1,27}_{i}=
-\mu_iw_i v_{i,x}\left(\frac{w_i}{v_i}\right)_x \xi_i \theta_i^\p \pa_{v_i}\bar{v}_i\tilde{r}_{i,v\si}=O(1)\sum_j\La_j^4$ and using again Lemmas \ref{estimimpo1}, \ref{lemme6.5}, \ref{lemme6.6}, \ref{lemme11.3} and \eqref{ngtech5}, \eqref{estimate-v-i-xx-1aabis}
\begin{align*}
	\widetilde{\B}^{1,27}_{i,x}&=-\big[\mu_{i,x}w_i v_{i,x} \left(\frac{w_i}{v_i}\right)_x+  \mu_i(w_{i,x} v_{i,x}+w_i v_{i,xx})\left(\frac{w_i}{v_i}\right)_x+\mu_iw_i v_{i,x}\left(\frac{w_i}{v_i}\right)_{xx}\big]
	 \xi_i \theta_i^\p \pa_{v_i}\bar{v}_i\tilde{r}_{i,v\si}\\
	 &-\mu_iw_i v_{i,x}\left(\frac{w_i}{v_i}\right)_x\big[ \left(\frac{w_i}{v_i}\right)_x(\xi_i' \theta_i^\p +\xi_i\theta_i'')\pa_{v_i}\bar{v}_i\tilde{r}_{i,v\si}+ \xi_i \theta_i^\p (\pa_x(\pa_{v_i}\bar{v}_i)\tilde{r}_{i,v\si}+ \pa_{v_i}\bar{v}_i\pa_x(\tilde{r}_{i,v\si}))\big]
	 \\
	%&-\mu_iw_i v_{i,x}\left(\frac{w_i}{v_i}\right)_x\big[ \left(\frac{w_i}{v_i}\right)_x(\xi_i' \theta_i^\p +\xi_i\theta_i'')\pa_{v_i}\bar{v}_i\tilde{r}_{i,v\si}\\
	%&--\mu_iw_i v_{i,x}\left(\frac{w_i}{v_i}\right)^2_x \xi_i \theta_i'' \pa_{v_i}\bar{v}_i\tilde{r}_{i,v\si}
	%-\mu_iw_i v_{i,x}\left(\frac{w_i}{v_i}\right)_x \xi_i \theta_i^\p \pa_x(\pa_{v_i}\bar{v}_i)\tilde{r}_{i,v\si}\\
	%&-\mu_iw_i v_{i,x}\left(\frac{w_i}{v_i}\right)_x \xi_i \theta_i^\p \pa_{v_i}\bar{v}_i\pa_x(\tilde{r}_{i,v\si})\\
	&=O(1)\sum_j(\La_j^1+\delta_0^2\La_j^3+\La_j^4+\La_j^5+\La_j^6).
	%-\mu_{i,x}(w_{i,x}-(w_i/v_i)v_{i,x}) \xi_i\frac{w_i}{v_i}v_{i,x}\theta_i^\p \pa_{v_i}\bar{v}_i\tilde{r}_{i,v\si}\\
	%&-\mu_i(w_{i,xx}-(w_i/v_i)v_{i,xx}-v_{i,x}(w_i/v_i)_x) \xi_i\frac{w_i}{v_i}v_{i,x}\theta_i^\p \pa_{v_i}\bar{v}_i\tilde{r}_{i,v\si}\\
	%&-\mu_i(w_{i,x}-(w_i/v_i)v_{i,x}) \xi_i^\p\left(\frac{w_i}{v_i}\right)_x\frac{w_i}{v_i}v_{i,x}\theta_i^\p \pa_{v_i}\bar{v}_i\tilde{r}_{i,v\si}\\
	%&-\mu_i(w_{i,x}-(w_i/v_i)v_{i,x}i) \xi_i\left(\frac{w_i}{v_i}\right)_xv_{i,x}\theta_i^\p \pa_{v_i}\bar{v}_i\tilde{r}_{i,v\si}\\
	%&-\mu_i(w_{i,x}-(w_i/v_i)v_{i,x}i) \xi_i\frac{w_i}{v_i}v_{i,xx}\theta_i^\p \pa_{v_i}\bar{v}_i\tilde{r}_{i,v\si}\\
	%&-\mu_i(w_{i,x}-(w_i/v_i)v_{i,x}i) \xi_i\frac{w_i}{v_i}v_{i,x}\theta_i^{\p\p}\left(\frac{w_i}{v_i}\right)_x \pa_{v_i}\bar{v}_i\tilde{r}_{i,v\si}\\
	%&-\mu_i(w_{i,x}-(w_i/v_i)v_{i,x}i) \xi_i\frac{w_i}{v_i}v_{i,x}\theta_i^\p (\pa_{v_i}\bar{v}_i)_x\tilde{r}_{i,v\si}\\
	%&-\mu_i(w_{i,x}-(w_i/v_i)v_{i,x}i) \xi_i\frac{w_i}{v_i}v_{i,x}\theta_i^\p \pa_{v_i}\bar{v}_i[\tilde{r}_{i,uv\si}u_x+(\xi\bar{v}_i)_x\tilde{r}_{i,vv\si}-\theta_i^\p\left(\frac{w_i}{v_i}\right)_x\tilde{r}_{i,v\si\si}].
\end{align*}
We wish to give more details on the term $\mu_iw_i v_{i,xx}\left(\frac{w_i}{v_i}\right)_x \xi_i \theta_i^\p \pa_{v_i}\bar{v}_i\tilde{r}_{i,v\si}$, applying Lemmas \ref{estimimpo1}, \ref{lemme6.5}, \ref{lemme6.6}, \ref{lemme11.3}, \eqref{ngtech5}, \eqref{estimate-v-i-xx-1aabis} and the fact that $w_{i,x}=v_i\left(\frac{w_i}{v_i}\right)_x+\frac{w_i}{v_i}v_{i,x}$ we obtain
\begin{align*}
&\mu_iw_i v_{i,xx}\left(\frac{w_i}{v_i}\right)_x \xi_i \theta_i^\p \pa_{v_i}\bar{v}_i\tilde{r}_{i,v\si}=O(1)
(|v_{i,x}|+|w_{i,x}|+|v_i|)|v_i \left(\frac{w_i}{v_i}\right)_x|\mathfrak{A}_i\rho_i^\e\\
&+O(1)\sum_k(\La_k^1+\La_k^4+\La_k^5+\La_k^6)=O(1)\sum_k(\La_k^1+\delta_0^2\La_j^3+\La_k^4+\La_k^5+\La_k^6).
\end{align*}
Next we have  $\widetilde{\B}^{1,28}_{i}=
-\sum\limits_{j\neq i}\mu_i w_i \left(\frac{w_i}{v_i}\right)_x \xi_iv_{j,x}\theta_i^\p \pa_{v_j}\bar{v}_i\tilde{r}_{i,v\si}=O(1)\La_i^1$ and from Lemmas \ref{estimimpo1}, \ref{lemme6.5}, \ref{lemme6.6}, \ref{lemme11.3} and \eqref{ngtech5} we get
\begin{align*}
	&\widetilde{\B}^{1,28}_{i,x}=-\sum\limits_{j\neq i}\big[(\mu_{i,x}w_i+\mu_iw_{i,x})\left(\frac{w_i}{v_i}\right)_x \xi_i+\mu_i w_i \left(\frac{w_i}{v_i}\right)_{xx} \xi_i+\mu_i w_i \left(\frac{w_i}{v_i}\right)_x^2 \xi_i'\big]v_{j,x}\theta_i^\p \pa_{v_j}\bar{v}_i\tilde{r}_{i,v\si}\\
	%(w_{i,x}-(w_i/v_i)v_{i,x})\xi_i \frac{w_i}{v_i}v_{j,x}\theta_i^\p \pa_{v_j}\bar{v}_i\tilde{r}_{i,v\si}\\
	%&-\sum\limits_{j\neq i}\mu_i(w_{i,xx}-(w_i/v_i)v_{i,xx}-v_{i,x}(w_i/v_i)_x)\xi_i \frac{w_i}{v_i}v_{j,x}\theta_i^\p \pa_{v_j}\bar{v}_i\tilde{r}_{i,v\si}\\
	%&-\sum\limits_{j\neq i}\mu_i(w_{i,x}-(w_i/v_i)v_{i,x})\xi_i^\p\left(\frac{w_i}{v_i}\right)_x \frac{w_i}{v_i}v_{j,x}\theta_i^\p \pa_{v_j}\bar{v}_i\tilde{r}_{i,v\si}\\
	&-\sum\limits_{j\neq i}\mu_i w_i \left(\frac{w_i}{v_i}\right)_x \xi_i\big[(v_{j,xx}\pa_{v_j}\bar{v}_i+v_{j,x}\pa_x(\pa_{v_j}\bar{v}_i))
	 \theta_i^\p\tilde{r}_{i,v\si}+v_{j,x} \pa_{v_j}\bar{v}_i(\theta_i'' \left(\frac{w_i}{v_i}\right)_x \tilde{r}_{i,v\si}+\theta'_i\pa_x( \tilde{r}_{i,v\si}))\big]\\
	%\mu_i(w_{i,x}-(w_i/v_i)v_{i,x})\xi_i \left(\frac{w_i}{v_i}\right)_xv_{j,x}\theta_i^\p \pa_{v_j}\bar{v}_i\tilde{r}_{i,v\si}\\
	%&-\sum\limits_{j\neq i}\mu_i(w_{i,x}-(w_i/v_i)v_{i,x})\xi_i \frac{w_i}{v_i}v_{j,xx}\theta_i^\p \pa_{v_j}\bar{v}_i\tilde{r}_{i,v\si}\\
	%&-\sum\limits_{j\neq i}\mu_i(w_{i,x}-(w_i/v_i)v_{i,x})\xi_i \frac{w_i}{v_i}v_{j,x}\theta_i^{\p\p}\left(\frac{w_i}{v_i}\right)_x \pa_{v_j}\bar{v}_i\tilde{r}_{i,v\si}\\
	%&-\sum\limits_{j\neq i}\mu_i(w_{i,x}-(w_i/v_i)v_{i,x})\xi_i \frac{w_i}{v_i}v_{j,x}\theta_i^\p (\pa_{v_j}\bar{v}_i)_x\tilde{r}_{i,v\si}\\
	%&-\sum\limits_{j\neq i}\mu_i(w_{i,x}-(w_i/v_i)v_{i,x})\xi_i \frac{w_i}{v_i}v_{j,x}\theta_i^\p \pa_{v_j}\bar{v}_i\pa_x(\tilde{r}_{i,v\si})\\
	&=O(1)(\La_i^1+\delta_0^2\La_i^3).
\end{align*} 
We have now from Lemma \ref{lemme6.5} $\widetilde{\B}^{1,29}_{i}=%-\mu_i(w_{i,x}-(w_i/v_i)v_{i,x}) \xi_i^\p\frac{w_i}{v_i}\left(\frac{w_i}{v_i}\right)_x\theta_i^\p \bar{v}_i\tilde{r}_{i,v\si}=
-\mu_i \xi_i^\p w_i\left(\frac{w_i}{v_i}\right)^2_x\theta_i^\p\bar{v}_i\tilde{r}_{i,v\si}=O(1)\delta_0^2\La_i^3$
and we get from Lemmas \ref{lemme6.5}, \ref{lemme6.6}, \ref{lemme11.3}
\begin{align*}
	\widetilde{\B}^{1,29}_{i,x}=&-\big[(\mu_{i,x}w_i+\mu_iw_{i,x})\xi_i^\p\theta'_i +(\xi''_i\theta'_i+\xi'_i\theta''_i)\left(\frac{w_i}{v_i}\right)_x\mu_iw_i\big]
	\bar{v}_i\left(\frac{w_i}{v_i}\right)_x^2 \tilde{r}_{i,v\si}\\
	&-\mu_i w_i\xi'_i\theta'_i\big[(\pa_x(\bar{v}_i )\tilde{r}_{i,v\si}+\bar{v}_i \pa_x (\tilde{r}_{i,v\si}))\left(\frac{w_i}{v_i}\right)_x^2+2\left(\frac{w_i}{v_i}\right)_x \left(\frac{w_i}{v_i}\right)_{xx}\bar{v}_i\tilde{r}_{i,v\si}\big]\\
	%(w_{i,xx}-(w_i/v_i)v_{i,xx}-v_{i,x}(w_i/v_i)_x) \xi_i^\p\bar{v}_i\frac{w_i}{v_i}\left(\frac{w_i}{v_i}\right)_x\theta_i^\p \tilde{r}_{i,v\si}\\
	%&-\mu_i(w_{i,x}-(w_i/v_i)v_{i,x}) \xi_i^{\p\p}\bar{v}_i\frac{w_i}{v_i}\left(\frac{w_i}{v_i}\right)^2_x\theta_i^\p \tilde{r}_{i,v\si}\\
	%&-\mu_i(w_{i,x}-(w_i/v_i)v_{i,x}) \xi_i^\p\bar{v}_{i,x}\frac{w_i}{v_i}\left(\frac{w_i}{v_i}\right)_x\theta_i^\p \tilde{r}_{i,v\si}\\
	%&-\mu_i(w_{i,x}-(w_i/v_i)v_{i,x}) \xi_i^\p\bar{v}_i\left(\frac{w_i}{v_i}\right)^2_x\theta_i^\p\tilde{r}_{i,v\si}\\
	%&-\mu_i(w_{i,x}-(w_i/v_i)v_{i,x}) \xi_i^\p\bar{v}_i\frac{w_i}{v_i}\left(\frac{w_i}{v_i}\right)_{xx}\theta_i^\p \tilde{r}_{i,v\si}\\
	%&-\mu_i(w_{i,x}-(w_i/v_i)v_{i,x}) \xi_i^\p\bar{v}_i\frac{w_i}{v_i}\left(\frac{w_i}{v_i}\right)^2_x\theta_i^{\p\p} \tilde{r}_{i,v\si}\\
	%&-\mu_i(w_{i,x}-(w_i/v_i)v_{i,x}) \xi_i^\p\bar{v}_i\frac{w_i}{v_i}\left(\frac{w_i}{v_i}\right)_x\theta_i^\p \pa_x (\tilde{r}_{i,v\si})\\
	=&O(1)(\delta_0^2\La_i^3+\La_i^5).
\end{align*}
We have now from Lemmas \ref{estimimpo1}, \ref{lemme6.5}, \ref{lemme6.6} and \eqref{ngtech5} $\widetilde{\B}^{1,30}_{i}=-w_i\xi_i\pa_{v_i}\bar{v}_i\theta'_i\left(\frac{w_i}{v_i}\right)_x\mathcal{H}_i\tilde{r}_{i,v\si}%=\\-w_i\xi_i\pa_{v_i}\bar{v}_i\theta'_i\left(\frac{w_i}{v_i}\right)_x(\mu_iv_{i,x}-v_i(\tilde{\lambda}_i-\lambda_i^*)-w_i)\tilde{r}_{i,v\si}
=\sum_j\La_j^4$. We obtain from Lemmas \ref{estimimpo1}, \ref{lemme6.5}, \ref{lemme6.6}, \ref{lemme11.3}, \eqref{ngtech5}, \eqref{estimate-v-i-xx-1aa} and using the fact that $w_{i,x}=v_i\left(\frac{w_i}{v_i}\right)_x+\frac{w_i}{v_i}v_{i,x}$
\begin{align*}
	\widetilde{\B}^{1,30}_{i,x}&=-\big[(w_{i,x}\pa_{v_i}\bar{v}_i+w_i\pa_x(\pa_{v_i}\bar{v}_i))
	\xi_i\theta'_i+w_i\pa_{v_i}\bar{v}_i(\xi'_i\theta'_i+\xi_i\theta''_i)\left(\frac{w_i}{v_i}\right)_x
	\big]
	\left(\frac{w_i}{v_i}\right)_x\mathcal{H}_i\tilde{r}_{i,v\si}\\
	&-w_i\xi_i \pa_{v_i}\bar{v}_i\theta'_i\big[\left(\frac{w_i}{v_i}\right)_{xx}\mathcal{H}_i\tilde{r}_{i,v\si}
	+\left(\frac{w_i}{v_i}\right)_{x}(\mathcal{H}_{i,x}\tilde{r}_{i,v\si}
+\mathcal{H}_i\pa_x(\tilde{r}_{i,v\si}
))\big]\\
	%&-w_i\xi_i\pa_x(\pa_{v_i}\bar{v}_i)\theta'_i\left(\frac{w_i}{v_i}\right)_x(\mu_iv_{i,x}-v_i(\tilde{\lambda}_i-\lambda_i^*)-w_i)\tilde{r}_{i,v\si}\\
	%&-w_i\xi_i\pa_{v_i}\bar{v}_i\theta''_i\left(\frac{w_i}{v_i}\right)^2_x(\mu_iv_{i,x}-v_i(\tilde{\lambda}_i-\lambda_i^*)-w_i)\tilde{r}_{i,v\si}\\
	%&-w_i\xi_i\pa_{v_i}\bar{v}_i\theta'_i\left(\frac{w_i}{v_i}\right)_{xx}(\mu_iv_{i,x}-v_i(\tilde{\lambda}_i-\lambda_i^*)-w_i)\tilde{r}_{i,v\si}\\
	%&-w_i\xi_i\pa_{v_i}\bar{v}_i\theta'_i\left(\frac{w_i}{v_i}\right)_x(\mu_iv_{i,x}-v_i(\tilde{\lambda}_i-\lambda_i^*)-w_i)_x\tilde{r}_{i,v\si}\\
	%&-w_i\xi_i\pa_{v_i}\bar{v}_i\theta'_i\left(\frac{w_i}{v_i}\right)_x(\mu_iv_{i,x}-v_i(\tilde{\lambda}_i-\lambda_i^*)-w_i)\pa_x(\tilde{r}_{i,v\si})\\
	&=\sum_j(\La_j^1+\delta_0^2\La_j^3+\La_j^4+\La_j^5+\La_j^6).
	\end{align*}
	We have now $\widetilde{\B}^{1,31}_{i}=
-\mu_i v_i w_i \left(\frac{w_i}{v_i}\right)_x \theta_i^\p \tilde{r}_{i,u\si}\tilde{r}_i=O(1)\La_i^4$ and using the fact that $\tilde{r}_{i,u\si}=O(1)\xi_i\bar{v}_i$ and Lemmas \ref{lemme6.5}, \ref{lemme6.6}, \ref{lemme11.4} we have
\begin{align*}
	\widetilde{\B}^{1,31}_{i,x}&=-\big[(\mu_{i,x} v_i w_i +\mu_iv_{i,x}w_i+\mu_i v_iw_{i,x})\left(\frac{w_i}{v_i}\right)_x+\mu_iv_i w_i \left(\frac{w_i}{v_i}\right)_{xx}\big] \theta_i^\p \tilde{r}_{i,u\si}\tilde{r}_i\\
	&-\mu_i v_i w_i \left(\frac{w_i}{v_i}\right)_x \big[\theta_i''\left(\frac{w_i}{v_i}\right)_x \tilde{r}_{i,u\si}\tilde{r}_i+\theta'_i\pa_x( \tilde{r}_{i,u\si}\tilde{r}_i)\big]
	% \theta_i^\p \tilde{r}_{i,u\si}\tilde{r}_i
	%-\mu_i( 2v_{i,x} w_i+v_i^2 \left(\frac{w_i}{v_i}\right)_x) \left(\frac{w_i}{v_i}\right)_x \theta_i^\p \tilde{r}_{i,u\si}\tilde{r}_i\\
	%&-\mu_i v_i w_i \left(\frac{w_i}{v_i}\right)_{xx} \theta_i^\p \tilde{r}_{i,u\si}\tilde{r}_i-\mu_i v_i w_i \left(\frac{w_i}{v_i}\right)^2_x \theta_i'' \tilde{r}_{i,u\si}\tilde{r}_i\\
	%&-\mu_i v_i w_i \left(\frac{w_i}{v_i}\right)_x \theta_i^\p\pa_x( \tilde{r}_{i,u\si}\tilde{r}_i)\\
	=O(1)(\La_i^4+\La_i^5).%\sum_j(\La_j^1+\delta_0^2\La_j^3+\La_j^4+\La_j^5).
	%\mu_{i,x}(w_{i,x}-(w_i/v_i)v_{i,x}) w_i\theta_i^\p \tilde{r}_{i,u\si}\tilde{r}_i\\
	%&-\mu_i(w_{i,xx}-(w_i/v_i)v_{i,xx}-v_{i,x}(w_i/v_i)_x) w_i\theta_i^\p \tilde{r}_{i,u\si}\tilde{r}_i\\
	%&-\mu_i(w_{i,x}-(w_i/v_i)v_{i,x}) w_{i,x}\theta_i^\p \tilde{r}_{i,u\si}\tilde{r}_i\\
	%&-\mu_i(w_{i,x}-(w_i/v_i)v_{i,x}) w_i\theta_i^{\p\p}\left(\frac{w_i}{v_i}\right)_x \tilde{r}_{i,u\si}\tilde{r}_i\\
	%&-\mu_i(w_{i,x}-(w_i/v_i)v_{i,x}) w_i\theta_i^\p [\tilde{r}_{i,uu\si}(\tilde{r}_i\otimes u_x)+(\xi_i\bar{v}_i)_x\tilde{r}_{i,uv\si}\tilde{r}_{i}-\theta_i^\p\left(\frac{w_i}{v_i}\tilde{r}_{i,u\si\si}\tilde{r}_i\right)]\\
	%&-\mu_i(w_{i,x}-(w_i/v_i)v_{i,x}) w_i\theta_i^\p [\tilde{r}_{i,u\si}(\tilde{r}_{i,u} u_x)+(\xi_i\bar{v}_i)_x\tilde{r}_{i,u\si}\tilde{r}_{i,v}-\theta_i^\p\left(\frac{w_i}{v_i}\tilde{r}_{i,u\si}\tilde{r}_{i,\si}\right)].
\end{align*}
Similarly we have $\widetilde{\B}^{1,32}_{i}=-\sum\limits_{j\neq i}\mu_i w_i\left(\frac{w_i}{v_i}\right)_xv_j\theta_i^\p \tilde{r}_{i,u\si}\tilde{r}_j=O(1)\La_i^1$ and using the fact that $ \tilde{r}_{i,u\si}=O(1)\xi_i \bar{v}_i$ and Lemma \ref{lemme11.4} we get
\begin{align*}
	\widetilde{\B}^{1,32}_{i,x}&=-\sum\limits_{j\neq i}\big[(\mu_{i,x} w_iv_j+\mu_i w_{i,x}v_j+\mu_i w_i v_{j,x})\left(\frac{w_i}{v_i}\right)_x+\mu_i w_i v_j\left(\frac{w_i}{v_i}\right)_{xx}\big]
	\theta_i^\p \tilde{r}_{i,u\si}\tilde{r}_j\\
	&-\sum\limits_{j\neq i}\mu_i w_{i}\left(\frac{w_i}{v_i}\right)_xv_j\big[\theta_i'' \left(\frac{w_i}{v_i}\right)_x\tilde{r}_{i,u\si}\tilde{r}_j+\theta'_i\pa_x(\tilde{r}_{i,u\si}\tilde{r}_j)\big]
	%&-\sum\limits_{j\neq i}\mu_i w_i\left(\frac{w_i}{v_i}\right)_{xx}v_j\theta_i^\p \tilde{r}_{i,u\si}\tilde{r}_j-\sum\limits_{j\neq i}\mu_i w_i\left(\frac{w_i}{v_i}\right)_xv_{j,x}\theta_i^\p \tilde{r}_{i,u\si}\tilde{r}_j\\
	%&-\sum\limits_{j\neq i}\mu_i w_i\left(\frac{w_i}{v_i}\right)^2_xv_j\theta_i'' \tilde{r}_{i,u\si}\tilde{r}_j
	%-\sum\limits_{j\neq i}\mu_i w_i\left(\frac{w_i}{v_i}\right)_xv_j\theta_i^\p \pa_x(\tilde{r}_{i,u\si}\tilde{r}_j)\\
	=O(1)\La_i^1.
	%-\sum\limits_{j\neq i}\mu_{i,x}(w_{i,x}-(w_i/v_i)v_{i,x}) (w_i/v_i)v_j\theta_i^\p \tilde{r}_{i,u\si}\tilde{r}_j\\
	%&-\sum\limits_{j\neq i}\mu_i(w_{i,xx}-(w_i/v_i)v_{i,xx}-v_{i,x}(w_i/v_i)_x) (w_i/v_i)v_j\theta_i^\p \tilde{r}_{i,u\si}\tilde{r}_j\\
	%&-\sum\limits_{j\neq i}\mu_i(w_{i,x}-(w_i/v_i)v_{i,x}) (w_i/v_i)_xv_j\theta_i^\p \tilde{r}_{i,u\si}\tilde{r}_j\\
	%&-\sum\limits_{j\neq i}\mu_i(w_{i,x}-(w_i/v_i)v_{i,x}) (w_i/v_i)v_{j,x}\theta_i^\p \tilde{r}_{i,u\si}\tilde{r}_j\\
	%&-\sum\limits_{j\neq i}\mu_i(w_{i,x}-(w_i/v_i)v_{i,x}) (w_i/v_i)v_j\theta_i^{\p\p}\left(\frac{w_i}{v_i}\right)_x \tilde{r}_{i,u\si}\tilde{r}_j\\
	%&-\sum\limits_{j\neq i}\mu_i(w_{i,x}-(w_i/v_i)v_{i,x}) (w_i/v_i)v_j\theta_i^\p [\tilde{r}_{i,uu\si}(\tilde{r}_j\otimes u_x)+(\xi_i\bar{v}_i)_x\tilde{r}_{i,uv\si}\tilde{r}_{j}-\theta_i^\p\left(\frac{w_i}{v_i}\tilde{r}_{i,u\si\si}\tilde{r}_j\right)]\\
	%&-\sum\limits_{j\neq i}\mu_i(w_{i,x}-(w_i/v_i)v_{i,x}) (w_i/v_i)v_j\theta_i^\p [\tilde{r}_{i,u\si}(\tilde{r}_{j,u} u_x)+(\xi_j\bar{v}_j)_x\tilde{r}_{i,u\si}\tilde{r}_{j,v}-\theta_j^\p\left(\frac{w_j}{v_j}\tilde{r}_{i,u\si}\tilde{r}_{j,\si}\right)].
\end{align*}
	Finally we have since $\tilde{r}_{i,\si\si}=O(1)\xi_i \bar{v}_i$, $\widetilde{\B}^{1,33}_{i}=\mu_i w_i\left(\frac{w_i}{v_i}\right)^2_x(\theta_i^\p)^2 \tilde{r}_{i,\si\si}=O(1)\delta_0^2\La_i^3$ and from Lemma \ref{lemme11.3} we obtain again since  $\tilde{r}_{i,\si\si}=O(1)\xi_i \bar{v}_i$
\begin{align*}
	\widetilde{\B}^{1,33}_{i,x}&=\big[(\mu_{i,x}w_i+\mu_i w_{i,x})\left(\frac{w_i}{v_i}\right)^2_x+2\mu_i w_i\left(\frac{w_i}{v_i}\right)_x\left(\frac{w_i}{v_i}\right)_{xx}\big](\theta_i^\p)^2 \tilde{r}_{i,\si\si}\\
	&+\mu_i w_i\left(\frac{w_i}{v_i}\right)^2_x\big[2\theta'_i\theta''_i\left(\frac{w_i}{v_i}\right)_x\tilde{r}_{i,\si\si}+(\theta_i^\p)^2\pa_x( \tilde{r}_{i,\si\si})\big]
	%(w_{i,x}-(w_i/v_i)v_{i,x}) \frac{w_i}{v_i}\left(\frac{w_i}{v_i}\right)_x(\theta_i^\p)^2 \tilde{r}_{i,\si\si}\\
	%&+\mu_i(w_{i,xx}-(w_i/v_i)v_{i,xx}-v_{i,x}(w_i/v_i)_x) \frac{w_i}{v_i}\left(\frac{w_i}{v_i}\right)_x(\theta_i^\p)^2 \tilde{r}_{i,\si\si}\\
	%&+\mu_i(w_{i,x}-(w_i/v_i)v_{i,x}) \left(\frac{w_i}{v_i}\right)^2_x(\theta_i^\p)^2 \tilde{r}_{i,\si\si}\\
	%&+\mu_i(w_{i,x}-(w_i/v_i)v_{i,x}) \frac{w_i}{v_i}\left(\frac{w_i}{v_i}\right)_{xx}(\theta_i^\p)^2 \tilde{r}_{i,\si\si}\\
	%&+2\mu_i(w_{i,x}-(w_i/v_i)v_{i,x}) \frac{w_i}{v_i}\left(\frac{w_i}{v_i}\right)^2_x(\theta_i^\p)\theta_i^{\p\p} \tilde{r}_{i,\si\si}\\
	%&+\mu_i(w_{i,x}-(w_i/v_i)v_{i,x}) \frac{w_i}{v_i}\left(\frac{w_i}{v_i}\right)_x(\theta_i^\p)^2
	%\pa_x( \tilde{r}_{i,\si\si})\\
	=O(1)(\delta_0^2\La_i^3+\La_i^5).
\end{align*}
\end{proof}
\begin{lemma}
We have
	\begin{equation}
	\begin{aligned}
			&K_{7,x}=O(1)\sum_j(\La_j^1+\La_j^2+\delta_0^2\La_j^3+\La_j^4+\La_j^5+\La_j^6+\La_j^{6,1})\\
			&K_{7,xx} =O(1)\sum_j(\La_j^1+\La_j^2+\delta_0^2\La_j^3+\La_j^4+\La_j^5+\La_j^6+\La_j^{6,1}).
						\end{aligned}
	\end{equation}
\label{K7xK8x}
\end{lemma}
\begin{proof}
We recall now that
\begin{align}
	K_{7,x}&= \sum\limits_{i\ne j}\left(\frac{w_i}{v_i}\right)_x\xi_iv_{j,x}v_i\pa_{v_j}\bar{v}_iB(u)\tilde{r}_{i,v}+
	\sum\limits_{l=1}^{{\color{black}{8}}}\sum\limits_{i}\left(\frac{w_i}{v_i}-\la_i^*\right)\al_{i}^{4,l}. \label{aK7x}
%	K_{8,x}&= \sum\limits_{i\ne j}\left(\frac{w_i}{v_i}\right)_x\xi_iw_{j,x}v_i\pa_{w_j}\bar{v}_iB\tilde{r}_{i,v}+\sum\limits_{l=1}^{10}\sum\limits_{i}\left(\frac{w_i}{v_i}-\la_i^*\right)\al_{i}^{5,l}.
\end{align}
%\begin{align}
%	K_{7,x}&= \sum\limits_{i\ne j}\left(\frac{w_i}{v_i}\right)_x\xi_iv_{j,x}v_i\pa_{v_j}\bar{v}_iB(u)\tilde{r}_{i,v}+
%	\sum\limits_{l=1}^{10}\sum\limits_{i}\left(\frac{w_i}{v_i}-\la_i^*\right)\al_{i}^{4,l}.\label{aK7x}
%\end{align}
Since each term $\al_{i}^{4,l}$  are a multiple of $\xi_i$ or $\xi'_i$ we deduce easily from \eqref{aK7x},  and Lemmas \ref{lemmeal4al5}, \ref{estimimpo1} that
\begin{align*}
&K_{7,x}=O(1)\sum_j(\La_j^1+\La_j^2+\delta_0^2\La_j^3+\La_j^4+\La_j^5+\La_j^6+\La_j^{6,1}).
\end{align*}
Let us study now $K_{7,xx}$ with
\begin{align*}
	K_{7,xx}&= \sum\limits_{i\ne j}\big[\left(\frac{w_i}{v_i}\right)_x\xi_iv_{j,x}v_i\pa_{v_j}\bar{v}_iB(u)\tilde{r}_{i,v}\big]_x +
	\sum\limits_{l=1}^{{\color{black}{8}}}\sum\limits_{i}\left(\frac{w_i}{v_i}\right)_x\al_{i}^{4,l}+\sum\limits_{l=1}^{{\color{black}{8}}}\sum\limits_{i}\left(\frac{w_i}{v_i}-\la_i^*\right)\al_{i,x}^{4,l}
	\end{align*}
Again due to the Lemma \ref{lemmeal4al5}, we have simply to study the two first terms of the previous equality. First, we have from Lemmas\ref{estimimpo1}, \ref{lemme6.5}, \ref{lemme6.6}, \ref{lemme11.3} and \eqref{ngtech5}
\begin{align*}
&\sum_{j\ne i}\big[\left(\frac{w_i}{v_i}\right)_x\xi_iv_{j,x}v_i\pa_{v_j}\bar{v}_iB(u)\tilde{r}_{i,v}\big]_x=\sum_{j\ne i}\big(\left(\frac{w_i}{v_i}\right)_{xx}\xi_i+\left(\frac{w_i}{v_i}\right)_x^2\xi_i'\big)v_{j,x}v_i\pa_{v_j}\bar{v}_iB(u)\tilde{r}_{i,v}\\
&+\sum_{j\ne i}\left(\frac{w_i}{v_i}\right)_x\xi_i \big[(v_{j,xx}v_i+v_{j,x}v_{i,x})\pa_{v_j}\bar{v}_iB(u)\tilde{r}_{i,v}+v_{j,x}v_i(\pa_x(\pa_{v_j}\bar{v}_i)B(u)+\pa_{v_j}\bar{v}_i u_x\cdot DB(u) )\tilde{r}_{i,v}\big]\\
%\sum_{j\ne i}\left(\frac{w_i}{v_i}\right)_x\xi_iv_{j,x}v_i\pa_x(\pa_{v_j}\bar{v}_i) B(u)\tilde{r}_{i,v}\\
&%+\sum_{j\ne i}\left(\frac{w_i}{v_i}\right)_x\xi_iv_{j,x}v_i\pa_{v_j}\bar{v}_iu_x\cdot DB(u)\tilde{r}_{i,v}
+\sum_{j\ne i}\left(\frac{w_i}{v_i}\right)_x\xi_iv_{j,x}v_i\pa_{v_j}\bar{v}_iB(u)\pa_x(\tilde{r}_{i,v})=O(1)\sum_j (\La_j^1+\delta_0^2\La_j^3).
\end{align*}
Let us deal now with the terms $
	\sum\limits_{l=1}^{8}\sum\limits_{i}\left(\frac{w_i}{v_i}\right)_x\al_{i}^{4,l}$. First we have from Lemmas \ref{estimimpo1}, \ref{lemme6.5}, \ref{lemme6.6} and \eqref{ngtech5}
	\begin{align*}
	\left(\frac{w_i}{v_i}\right)_x\al_{i}^{4,1}&=\sum\limits_{j\neq i}\left(\frac{w_i}{v_i}\right)_x B(u)\tilde{r}_{i,v}\pa_{v_j}\bar{v}_i\left[v_{j,xx}v_i\xi_i+v_{j,x}v_{i,x}\xi_i+v_{j,x}v_i\xi^\p_i\left(\frac{w_i}{v_i}\right)_x\right]\\
	&=O(1)(\La_i^1+\delta_0^2\La_i^3).
	\end{align*}
	Similarly we have
	%\begin{align*}
	%\left(\frac{w_i}{v_i}\right)_x\al_{i}^{4,2}&=\sum\limits_{j\neq i}\sum\limits_{k}\left(\frac{w_i}{v_i}\right)_x B(u)\tilde{r}_{i,v}v_i\pa_{v_jv_k}\bar{v}_i\big[\xi_iv_{j,x}v_{k,x}\big]=O(1)\La_i^1,\\
	%\left(\frac{w_i}{v_i}\right)_x\al_{i}^{4,3}&=\sum\limits_{j\neq i}\sum\limits_{k\neq i}\left(\frac{w_i}{v_i}\right)_xB(u)\tilde{r}_{i,v}v_i\pa_{v_jw_k}\bar{v}_i\big[\xi_iv_{j,x}w_{k,x}\big]=O(1)\La_i^1.
	%\end{align*}
	%We have now
	\begin{align*}
	\left(\frac{w_i}{v_i}\right)_x\al_{i}^{4,2}&=\sum\limits_{j\neq i}\sum\limits_{k}\left(\frac{w_i}{v_i}\right)_x B(u)\tilde{r}_{i,v}v_i\pa_{v_jv_k}\bar{v}_i\big[\xi_iv_{j,x}v_{k,x}\big]=O(1)\La_i^1,\\
	\left(\frac{w_i}{v_i}\right)_x\al_{i}^{4,3}&=\sum\limits_{j\neq i}\left(\frac{w_i}{v_i}\right)_x B(u)\tilde{r}_{i,uv}\tilde{r}_i\pa_{v_j}\bar{v}_i\big[v_{j,x}v^2_i\xi_i\big]=O(1)\La_i^1,\\
	\left(\frac{w_i}{v_i}\right)_x\al_{i}^{4,4}&=\sum\limits_{j\neq i}\sum\limits_{k\neq i}\left(\frac{w_i}{v_i}\right)_x B(u)\tilde{r}_{i,uv}\tilde{r}_k v_i\pa_{v_j}\bar{v}_i\big[v_{j,x}v_k\xi_i\big]=O(1)\La_i^1,\\
	\left(\frac{w_i}{v_i}\right)_x\al_{i}^{4,5}&=\sum\limits_{j\neq i}\sum\limits_{k}\left(\frac{w_i}{v_i}\right)_x B(u)\tilde{r}_{i,vv}\pa_{v_j}\bar{v}_i\pa_{v_k}\bar{v}_i\big[v_{k,x}v_{j,x}v_i\xi_i^2\big]=O(1)\La_i^1,\\
	%\left(\frac{w_i}{v_i}\right)_x\al_{i}^{4,7}&=\sum\limits_{j\neq i}\sum\limits_{k\neq i}\left(\frac{w_i}{v_i}\right)_x B(u)\tilde{r}_{i,vv}\pa_{v_j}\bar{v}_i\pa_{w_k}\bar{v}_i\big[w_{k,x}v_{j,x}v_i\xi_i^2\big]=O(1)\La_i^1,\\
	\left(\frac{w_i}{v_i}\right)_x\al_{i}^{4,6}&=\sum\limits_{j\neq i}B(u)\tilde{r}_{i,vv}\pa_{v_j}\bar{v}_i\left[v_{j,x}v_i\xi_i\xi_i^\p\left(\frac{w_i}{v_i}\right)^2_x\bar{v}_i\right]=O(1)\delta_0^2\La_i^3,\\
	\left(\frac{w_i}{v_i}\right)_x\al_{i}^{4,7}&=-\sum\limits_{j\neq i} B(u)\tilde{r}_{i,v\si}\pa_{v_j}\bar{v}_i\left[v_{j,x}v_i\xi_i\left(\frac{w_i}{v_i}\right)^2_x\theta^\p_i\right]=O(1)\delta_0^2\La_i^3,\\
	\left(\frac{w_i}{v_i}\right)_x\al_{i}^{4,8}&=\sum\limits_{j\neq i}\left(\frac{w_i}{v_i}\right)_x v_{j,x}v_i\xi_i\pa_{v_j}\bar{v}_iu_x\cdot DB(u)\tilde{r}_{i,v}=O(1)\La_i^1.
	\end{align*}
	It concludes the proof of the Lemma.
\end{proof}
\begin{lemma}\label{lemma:estimate-7-1}
	For $l=2,\cdots,10$, $l'=1,\cdots,10$ it follows,
	\begin{equation}
	\begin{aligned}
			&\B^{8,l}_{i}, \B^{9,l'}_i=O(1)\sum_j(\La_j^1+\La_j^2+\delta_0^2\La_j^3+\La_j^4+\La_j^5+\La_j^6+\La_j^{6,1})+R_{\e,1}^{8,9,i,l},\\
			&\B^{8,l}_{i,x},\B^{9,l'}_{i,x} =O(1)\sum_j(\La_j^1+\La_j^2+\delta_0^2\La_j^3+\La_j^4+\La_j^5+\La_j^6+\La_j^{6,1})+\widetilde{R}_{\e,1}^{8,9,i,l}
						\end{aligned}
	\end{equation}
	with:
	\begin{equation}
	\int_{\hat{t}}^T\int_{\R}(|R_{\e,1}^{8,9,i,l}|+|\widetilde{R}_{\e,1}^{8,9,i,l}|)dx ds=O(1)\delta_0^2,
	\label{12.11bis3}
	\end{equation}
	for $\e>0$ small enough in terms of $T-\hat{t}$ and $\delta_0$.
\end{lemma}
\begin{proof}

First we have $\B^{8,2}_{i}
=(w_{i,x}-\mu_i^{-1}(\tilde{\la}_i-\la^*_i+\theta_i)w_i)v_i\tilde{r}_i\cdot D(B-\mu_iI_n)\tilde{r}_{i}$ that we can rewrite as follows
\begin{align}
&\B^{8,2}_{i}=\big[v_i^2\left(\frac{w_i}{v_i}\right)_x+w_i\mu_i^{-1}\mathcal{H}_i %(\mu_i v_{i,x}-(\tilde{\la}_i-\la_i^*)v_i-w_i)
+\mu_i^{-1}w_i(w_i-\theta_i v_i)\big]\tilde{r}_i\cdot D(B-\mu_iI_n)\tilde{r}_{i}. \label{B92}
\end{align}
Using Lemmas \ref{lemme6.5}, \ref{lemme6.6}, \eqref{6.45} and the fact that $|w_i-\theta_i v_i|=O(1)|w_i|\mathbbm{1}_{\{|\frac{w_i}{v_i}|\geq \delta_1\}}$ we deduce that
\begin{align*}
w_i(w_i-\theta_i v_i)&=O(1)(|v_{i,x}|^2_{\{|\frac{w_i}{v_i}|\geq \delta_1\}}+\sum_j(\La_j^1+\La_j^4))\\
&=O(1)\sum_j(\La_j^1+\La_j^4+\La_j^6).
\end{align*}
From Lemmas \ref{lemme6.5}, \ref{lemme6.6} and \eqref{B92}, \eqref{ngtech4} we get
\begin{align*}
&\B^{8,2}_{i}=O(1)\sum_j(\La_j^1+\La_j^4+\La_j^6)+R_{\e,1}^{8,9,i,2}
\end{align*}
with $R_{\e,1}^{8,9,i,2}$ satisfying \eqref{12.11bis3}.
We have now using Lemmas \ref{lemme6.5}, \ref{lemme6.6}, \ref{lemme11.3}, \eqref{6.45}, \eqref{ngtech4}, \eqref{estimate-v-i-xx-1aa} and the fact that $w_{i,x}-\theta'_i \left(\frac{w_i}{v_i}\right)_x v_i-\theta_i v_{i,x}=0$ on $\{|\frac{w_i}{v_i}|\leq\frac{\delta_1}{2}\}$
\begin{align*}
	&\B^{8,2}_{i,x}=\big[(w_{i,xx}v_i-v_{i,xx}w_i)+(w_{i,x}\mu_i^{-1}+w_i(\mu_i^{-1})_x)\mathcal{H}_i+\mu_i^{-1}w_i\mathcal{H}_{i,x}\big]\tilde{r}_i\cdot D(B-\mu_iI_n)\tilde{r}_{i}\\
	%&+w_i\mu_i^{-1}(\mu_i v_{i,x}-(\tilde{\la}_i-\la_i^*)v_i-w_i)_x\tilde{r}_i\cdot D(B-\mu_iI_n)\tilde{r}_{i}\\
	&+\big[((\mu_i^{-1})_xw_i+\mu_i^{-1}w_{i,x})(w_i-\theta_i v_i)+\mu_i^{-1}w_i(w_{i,x}-\theta'_i \left(\frac{w_i}{v_i}\right)_x v_i-\theta_i v_{i,x})\big]
	\tilde{r}_i\cdot D(B-\mu_iI_n)\tilde{r}_{i}\\
	%&+\mu_i^{-1}w_i(w_{i,x}-\theta'_i \left(\frac{w_i}{v_i}\right)_x v_i-\theta_i v_{i,x})\big]\tilde{r}_i\cdot D(B-\mu_iI_n)\tilde{r}_{i}\\
	&+\big[v_i^2\left(\frac{w_i}{v_i}\right)_x+w_i\mu_i^{-1}\mathcal{H}_i+\mu_i^{-1}w_i(w_i-\theta_i v_i)\big]\big[\tilde{r}_{i,x}\cdot D(B-\mu_iI_n)\tilde{r}_{i}+\tilde{r}_i\cdot D(B-\mu_iI_n)\tilde{r}_{i,x}\big]\\
	&+\big[v_i^2\left(\frac{w_i}{v_i}\right)_x+w_i\mu_i^{-1}\mathcal{H}_i+\mu_i^{-1}w_i(w_i-\theta_i v_i)\big]\tilde{r}_i\otimes u_x : D^2(B-\mu_iI_n)\tilde{r}_{i}\\
	%&+\big[v_i^2\left(\frac{w_i}{v_i}\right)_x+w_i\mu_i^{-1}\mathcal{H}_i+\mu_i^{-1}w_i(w_i-\theta_i v_i)\big]\tilde{r}_i\cdot D(B-\mu_iI_n)\tilde{r}_{i,x}\\
	&=O(1)\sum_j(\La_j^1+\La_j^4+\La_j^5+\La_j^6+\La_j^{6,1})+\widetilde{R}_{\e,1}^{8,9,i,2},
	%(w_{i,xx}-(\mu_i^{-1})_x(\tilde{\la}_i-\la^*_i+\theta_i)w_i-\mu_i^{-1}(\tilde{\la}_{i,x}+\theta_i^\p\left(\frac{w_i}{v_i}\right)_x)w_i)v_i\tilde{r}_i\cdot D(B-\mu_iI_n)\tilde{r}_{i}\\
	%&+[(w_{i,x}-\mu_i^{-1}(\tilde{\la}_i-\la^*_i+\theta_i)w_i)v_{i,x}-\mu_i^{-1}(\tilde{\la}_i-\la^*_i+\theta_i)w_{i,x}v_i]\tilde{r}_i\cdot D(B-\mu_iI_n)\tilde{r}_{i}\\
	%&+(w_{i,x}-\mu_i^{-1}(\tilde{\la}_i-\la^*_i+\theta_i)w_i)v_{i} [\tilde{r}_{i,u}u_x+(\xi_i\bar{v}_i)_x\tilde{r}_{i,v}-\theta_i^\p\left(\frac{w_i}{v_i}\right)_x\tilde{r}_{i,\si}]\cdot D(B-\mu_iI_n)\tilde{r}_{i}\\
	%&+(w_{i,x}-\mu_i^{-1}(\tilde{\la}_i-\la^*_i+\theta_i)w_i)v_{i} (\tilde{r}_{i}\otimes u_x): D^2(B-\mu_iI_n)\tilde{r}_{i}\\
	%&+(w_{i,x}-\mu_i^{-1}(\tilde{\la}_i-\la^*_i+\theta_i)w_i)v_{i} \tilde{r}_{i}\cdot D(B-\mu_iI_n)[\tilde{r}_{i,u}u_x+(\xi_i\bar{v}_i)_x\tilde{r}_{i,v}-\theta_i^\p\left(\frac{w_i}{v_i}\right)_x\tilde{r}_{i,\si}].
\end{align*}
with $\widetilde{R}_{\e,1}^{8,9,i,2}$ satisfying \eqref{12.11bis3}.
We have now $\B^{8,3}_{i}=\sum\limits_{j\neq i}(w_{i,x}-\mu_i^{-1}(\tilde{\la}_i-\la^*_i+\theta_i)w_i)v_j\tilde{r}_j\cdot D(B-\mu_iI_n)\tilde{r}_{i}=O(1)\La_i^1$. From Lemma  \ref{lemme11.3} we obtain
\begin{align*}
	&\B^{8,3}_{i,x}=\big[w_{i,xx}-(\mu_i^{-1})_x(\tilde{\la}_i-\la^*_i+\theta_i)w_i-\mu_i^{-1}(\tilde{\la}_{i,x}+\theta_i^\p\left(\frac{w_i}{v_i}\right)_x)w_i\big]v_j\tilde{r}_j\cdot D(B-\mu_iI_n)\tilde{r}_{i}\\
	&-\sum\limits_{j\neq i}\mu_i^{-1}(\tilde{\la}_i-\la^*_i+\theta_i)w_{i,x} v_j\tilde{r}_j\cdot D(B-\mu_iI_n)\tilde{r}_{i}\\%+\sum\limits_{j\neq i}(w_{i,x}-\mu_i^{-1}(\tilde{\la}_i-\la^*_i+\theta_i)w_i)v_{j,x}\tilde{r}_j\cdot D(B-\mu_iI_n)\tilde{r}_{i}\\
	&+\sum\limits_{j\neq i}(w_{i,x}-\mu_i^{-1}(\tilde{\la}_i-\la^*_i+\theta_i)w_i)\big[v_{j,x}\tilde{r}_j+v_{j}\tilde{r}_{j,x}\big]\cdot D(B-\mu_iI_n)\tilde{r}_{i}\\
	&+\sum\limits_{j\neq i}(w_{i,x}-\mu_i^{-1}(\tilde{\la}_i-\la^*_i+\theta_i)w_i)\big[v_{j}(\tilde{r}_j\otimes u_x):D^2(B-\mu_iI_n)\tilde{r}_{i}+v_{j}\tilde{r}_j\cdot D(B-\mu_iI_n)\tilde{r}_{i,x}\big]\\
	%&+\sum\limits_{j\neq i}(w_{i,x}-\mu_i^{-1}(\tilde{\la}_i-\la^*_i+\theta_i)w_i)v_{j}\tilde{r}_j\cdot D(B-\mu_iI_n)\pa_x(\tilde{r}_i)\\
	&=O(1)\La_i^1.
\end{align*}
Applying the Lemma \ref{lemme11.4} we can treat $\B^{8,4}_{i}=(w_{i,x}-\mu_i^{-1}(\tilde{\la}_i-\la_i^*+\theta_i)w_i)v_{i}(B-\mu_iI_n)\tilde{r}_{i,u}\tilde{r}_i=O(1)\sum_j(\La_j^1+\La_j^4+\La_j^6)+R_{\e,1}^{8,9,i,4}$ as for  $\B_{i}^{8,2}$ with $R_{\e,1}^{8,9,i,4}$ satisfying \eqref{12.11bis3}.
 We have now $\B^{8,5}_{i}=\sum\limits_{j\neq i}(w_{i,x}-\mu_i^{-1}(\tilde{\la}_i-\la_i^*+\theta_i)w_i)v_{j}(B-\mu_iI_n)\tilde{r}_{i,u}\tilde{r}_j=O(1)\La_i^1$ which can be treatet as for $\B^{8,3}_{i}$ using Lemma \ref{lemme11.4}.
We have now $\B^{8,6}_{i}=(w_{i,x}-\mu_i^{-1}(\tilde{\la}_i-\la_i^*+\theta_i)w_i)v_{i,x}\xi_i(B-\mu_iI_n)\pa_{v_i}\bar{v}_i\tilde{r}_{i,v}$, as in \eqref{B92} we can rewrite this terms as follows using the fact that $\xi_i(w_i-\theta_i v_i)=0$
\begin{align*}
&\B^{8,6}_{i}=\big[v_i\left(\frac{w_i}{v_i}\right)_x+\frac{w_i}{v_i}\mu_i^{-1}\mathcal{H}_i\big]v_{i,x}\xi_i(B-\mu_iI_n)\pa_{v_i}\bar{v}_i\tilde{r}_{i,v}
\end{align*}
Using  Lemma \ref{estimimpo1}, \ref{lemme6.5}, \ref{lemme6.6} and \eqref{ngtech4}, \eqref{ngtech5} we deduce that
\begin{align*}
&\B^{8,6}_{i}=O(1)\sum_j(\La_j^1+\La_j^4+\La_j^6)+R_{\e,1}^{8,9,i,6},
\end{align*}
 with $R_{\e,1}^{8,9,i,6}$ satisfying \eqref{12.11bis3}.
We get now using again Lemmas \ref{estimimpo1}, \ref{lemme6.5}, \ref{lemme6.6}, \ref{lemme11.3} \eqref{ngtech4}, \eqref{ngtech5},  \eqref{estimate-v-i-xx-1aa}, \eqref{estimate-v-i-xx-1aabis} and the fact that $w_{i,x}=v_i\left(\frac{w_i}{v_i}\right)_x+\frac{w_i}{v_i}v_{i,x}$
\begin{align*}
	&\B^{8,6}_{i,x}=\big[v_{i,x}\left(\frac{w_i}{v_i}\right)_x+v_{i}\left(\frac{w_i}{v_i}\right)_{xx}
	+(\left(\frac{w_i}{v_i}\right)_x\mu_i^{-1}+\frac{w_i}{v_i}(\mu_i^{-1})_x)\mathcal{H}_i
	\big]v_{i,x}\xi_i(B-\mu_iI_n)\pa_{v_i}\bar{v}_i\tilde{r}_{i,v}\\
	%&+(\left(\frac{w_i}{v_i}\right)_x\mu_i^{-1}+\frac{w_i}{v_i}(\mu_i^{-1})_x)(\mu_i v_{i,x}-(\tilde{\la}_i-\la_i^*)v_i-w_i)v_{i,x}\xi_i(B-\mu_iI_n)\pa_{v_i}\bar{v}_i\tilde{r}_{i,v}\\
	&+\frac{w_i}{v_i}\mu_i^{-1}\mathcal{H}_{i,x}v_{i,x}\xi_i(B-\mu_iI_n)\pa_{v_i}\bar{v}_i\tilde{r}_{i,v}
	%&+((\mu_i^{-1})_x\frac{w_i}{v_i}+\mu_i^{-1}\left(\frac{w_i}{v_i}\right)_x)(w_i-\theta_i v_i)\big]v_{i,x}\xi_i(B-\mu_iI_n)\pa_{v_i}\bar{v}_i\tilde{r}_{i,v}\\
	%&+\mu_i^{-1}\frac{w_i}{v_i}(w_{i,x}-\theta_i'\left(\frac{w_i}{v_i}\right)_x v_i-\theta_i v_{i,x})v_{i,x}\xi_i(B-\mu_iI_n)\pa_{v_i}\bar{v}_i\tilde{r}_{i,v}\\
	+\big[v_i\left(\frac{w_i}{v_i}\right)_x+\frac{w_i}{v_i}\mu_i^{-1}\mathcal{H}_i\big]v_{i,xx}\xi_i(B-\mu_iI_n)\pa_{v_i}\bar{v}_i\tilde{r}_{i,v}\\
	&+\big[v_i\left(\frac{w_i}{v_i}\right)_x+\frac{w_i}{v_i}\mu_i^{-1}\mathcal{H}_i\big]v_{i,x}\big[\xi_i'\left(\frac{w_i}{v_i}\right)_x (B-\mu_iI_n)+\xi_i u_x\cdot D(B-\mu_iI_n)\big]
	\pa_{v_i}\bar{v}_i\tilde{r}_{i,v}\\
	&+\big[v_i\left(\frac{w_i}{v_i}\right)_x+\frac{w_i}{v_i}\mu_i^{-1}\mathcal{H}_i\big]v_{i,x}\xi_i(B-\mu_iI_n)\big[\pa_x(\pa_{v_i}\bar{v}_i)\tilde{r}_{i,v}+\pa_{v_i}\bar{v}_i\pa_x(\tilde{r}_{i,v})\big]\\
	%\xi_i u_x\cdot D(B-\mu_iI_n)\pa_{v_i}\bar{v}_i\tilde{r}_{i,v}\\
	%&+\big[v_i\left(\frac{w_i}{v_i}\right)_x+\frac{w_i}{v_i}\mu_i^{-1}(\mu_i v_{i,x}-(\tilde{\la}_i-\la_i^*)v_i-w_i)\big]v_{i,x}\xi_i(B-\mu_iI_n)\pa_x(\pa_{v_i}\bar{v}_i)\tilde{r}_{i,v}\\
	%&+\big[v_i\left(\frac{w_i}{v_i}\right)_x+\frac{w_i}{v_i}\mu_i^{-1}(\mu_i v_{i,x}-(\tilde{\la}_i-\la_i^*)v_i-w_i)\big]v_{i,x}\xi_i(B-\mu_iI_n)\pa_{v_i}\bar{v}_i\pa_x(\tilde{r}_{i,v})\\
	&=O(1)\sum_j(\La_j^1+\delta_0^2\La_j^3+\La_j^4+\La_j^5+\La_j^6)+\widetilde{R}_{\e,1}^{8,9,i,6},
\end{align*}
 with $\widetilde{R}_{\e,1}^{8,9,i,6}$ satisfying \eqref{12.11bis3}.
We have now $\B^{8,7}_{i}=\sum\limits_{j\neq i}(w_{i,x}-\mu_i^{-1}(\tilde{\la}_i-\la_i^*+\theta_i)w_i)v_{j,x}\xi_i(B-\mu_iI_n)(\pa_{v_j}\bar{v}_i)\tilde{r}_{i,v}=O(1)\La_i^1$ and from the Lemma \ref{estimimpo1}, \ref{lemme11.3}, \eqref{ngtech5} and the fact that $w_{i,x}=v_i\left(\frac{w_i}{v_i}\right)_x+\frac{w_i}{v_i}v_{i,x}$
\begin{align*}
	&\B^{8,7}_{i,x}=\sum\limits_{j\neq i}(w_{i,x}-\mu_i^{-1}(\tilde{\la}_i-\la_i^*+\theta_i)w_i)_xv_{j,x}\xi_i(B-\mu_iI_n)(\pa_{v_j}\bar{v}_i)\tilde{r}_{i,v}\\
	&+\sum\limits_{j\neq i}(w_{i,x}-\mu_i^{-1}(\tilde{\la}_i-\la_i^*+\theta_i)w_i)\big[v_{j,xx}\xi_i+v_{j,x}\xi_i^\p\left(\frac{w_i}{v_i}\right)_x\big]
	(B-\mu_iI_n)(\pa_{v_j}\bar{v}_i)\tilde{r}_{i,v}\\
	%&+\sum\limits_{j\neq i}(w_{i,x}-\mu_i^{-1}(\tilde{\la}_i-\la_i^*+\theta_i)w_i)v_{j,x}\xi_i^\p\left(\frac{w_i}{v_i}\right)_x(B-\mu_iI_n)(\pa_{v_j}\bar{v}_i)\tilde{r}_{i,v}\\
	&+\sum\limits_{j\neq i}(w_{i,x}-\mu_i^{-1}(\tilde{\la}_i-\la_i^*+\theta_i)w_i)v_{j,x}\xi_i\big[u_x\cdot D(B-\mu_iI_n)(\pa_{v_j}\bar{v}_i)+(B-\mu_iI_n)(\pa_{v_j}\bar{v}_i)_x)\big]\tilde{r}_{i,v}\\
	%&+\sum\limits_{j\neq i}(w_{i,x}-\mu_i^{-1}(\tilde{\la}_i-\la_i^*+\theta_i)w_i)v_{j,x}\xi_i(B-\mu_iI_n)(\pa_{v_j}\bar{v}_i)_x\tilde{r}_{i,v}\\
	&+\sum\limits_{j\neq i}(w_{i,x}-\mu_i^{-1}(\tilde{\la}_i-\la_i^*+\theta_i)w_i)v_{j,x}\xi_i(B-\mu_iI_n)(\pa_{v_j}\bar{v}_i)\pa_x(\tilde{r}_{i,v})
	=O(1)(\La_i^1+\sum_j \delta_0^2\La_j^3).
\end{align*}
We have now as in \eqref{B92} and using Lemmas \ref{estimimpo1}, \ref{lemme6.5}, \ref{lemme6.6}, \eqref{ngtech5}
\begin{align*}
&\B^{8,8}_{i}=(w_{i,x}-\mu_i^{-1}(\tilde{\la}_i-\la_i^*+\theta_i)w_i)\bar{v}_i\xi^\p_i\left(\frac{w_i}{v_i}\right)_x(B-\mu_iI_n)\tilde{r}_{i,v}\\
&=\big[v_i\bar{v_i}\left(\frac{w_i}{v_i}\right)_x+\frac{w_i \bar{v}_i}{v_i}\mu_i^{-1}\mathcal{H}_i\big]\xi^\p_i\left(\frac{w_i}{v_i}\right)_x(B-\mu_iI_n)\tilde{r}_{i,v}=O(1)(\delta_0^2\La_i^3+\sum_j\La_j^4).
\end{align*}
We get now using again Lemmas \ref{estimimpo1}, \ref{lemme6.5}, \ref{lemme6.6}, \ref{lemme11.3} and \eqref{ngtech5}, \eqref{paraestim1}
\begin{align*}
	&\B^{8,8}_{i,x}=\big[\big(v_{i,x}\bar{v_i}+v_i \pa_x(\bar{v}_i)\big)\left(\frac{w_i}{v_i}\right)_x+
	v_i\bar{v_i}\left(\frac{w_i}{v_i}\right)_{xx}\big]\xi^\p_i\left(\frac{w_i}{v_i}\right)_x(B-\mu_iI_n)\tilde{r}_{i,v}\\
	&+\big[(\left(\frac{w_i }{v_i}\right)_x \bar{v}_i+\frac{w_i}{v_i}\pa_x(\bar{v}_i))\mu_i^{-1}\mathcal{H}_i+\frac{w_i}{v_i}\bar{v}_i((\mu_i^{-1})_x\mathcal{H}_i+\mu_i^{-1}\mathcal{H}_{i,x})\big]
	\xi^\p_i\left(\frac{w_i}{v_i}\right)_x(B-\mu_iI_n)\tilde{r}_{i,v}\\
	%&+\frac{w_i\bar{v}_i}{v_i}(\mu_{i}^{-1})_x(\mu_i v_{i,x}-(\tilde{\la}_i-\la_i^*)v_i-w_i)\xi^\p_i\left(\frac{w_i}{v_i}\right)_x(B-\mu_iI_n)\tilde{r}_{i,v}\\
	%&+\frac{w_i \bar{v}_i}{v_i}\mu_i^{-1}(\mu_i v_{i,x}-(\tilde{\la}_i-\la_i^*)v_i-w_i)_x \xi^\p_i\left(\frac{w_i}{v_i}\right)_x(B-\mu_iI_n)\tilde{r}_{i,v}\\
	&+\big[v_i\bar{v_i}\left(\frac{w_i}{v_i}\right)_x+\frac{w_i \bar{v}_i}{v_i}\mu_i^{-1}\mathcal{H}_i\big]\big[\xi_i''\left(\frac{w_i}{v_i}\right)^2_x+\xi^\p_i\left(\frac{w_i}{v_i}\right)_{xx}\big]
	(B-\mu_iI_n)\tilde{r}_{i,v}\\
	%&+\big[v_i\bar{v_i}\left(\frac{w_i}{v_i}\right)_x+\frac{w_i \bar{v}_i}{v_i}\mu_i^{-1}(\mu_i v_{i,x}-(\tilde{\la}_i-\la_i^*)v_i-w_i)\big]\xi^\p_i\left(\frac{w_i}{v_i}\right)_{xx}(B-\mu_iI_n)\tilde{r}_{i,v}\\
	&+\big[v_i\bar{v_i}\left(\frac{w_i}{v_i}\right)_x+\frac{w_i \bar{v}_i}{v_i}\mu_i^{-1}\mathcal{H}_i\big]\xi^\p_i\left(\frac{w_i}{v_i}\right)_x \big[u_x\cdot D(B-\mu_iI_n)\tilde{r}_{i,v}+(B-\mu_iI_n)\pa_x(\tilde{r}_{i,v})\big]\\
	%&+\big[v_i\bar{v_i}\left(\frac{w_i}{v_i}\right)_x+\frac{w_i \bar{v}_i}{v_i}\mu_i^{-1}(\mu_i v_{i,x}-(\tilde{\la}_i-\la_i^*)v_i-w_i)\big]\xi^\p_i\left(\frac{w_i}{v_i}\right)_x(B-\mu_iI_n)\pa_x(\tilde{r}_{i,v})\\
	&=\sum_j(\La_j^1+\delta_0^2\La_j^3+\La_j^4+\La_j^5).
	\end{align*}
	Using \eqref{ngtech5}, the fact that $(B-\mu_iI_n)\tilde{r}_{i,\si}=\sum_{k\ne i}(\mu_k-\mu_i)\psi_{ik,\sig}r_k(u)=O(1)\xi_i\bar{v}_i$ because $\psi_{ik,\sig}=O(1)\xi_i\bar{v}_i$ and $\theta(x)=x$ on the support of $\xi$, we deduce as in \eqref{B92} that
	\begin{align*}
	\B^{8,9}_{i}&=-(w_{i,x}-\mu_i^{-1}(\tilde{\la}_i-\la_i^*+\theta_i)w_i)\theta_i^\p\left(\frac{w_i}{v_i}\right)_x(B-\mu_iI_n)\tilde{r}_{i,\si}\\
	%&-(w_{i,x}-\mu_i^{-1}(\tilde{\la}_i-\la_i^*+\theta_i)w_i)\theta_i^\p\left(\frac{w_i}{v_i}\right)_x(B-\mu_iI_n)\tilde{r}_{i,\si}\\
	&=-\big[v_i\left(\frac{w_i}{v_i}\right)_x+\frac{w_i}{v_i}\mu_i^{-1}\mathcal{H}_i\big]\left(\frac{w_i}{v_i}\right)_x \sum_{k\ne i}(\mu_k-\mu_i)\psi_{ik,\sig}r_k(u),\\
	&=\delta_0^2\La_i^3+O(1)\sum_j\La_j^4.
	\end{align*}
	We have now using Lemmas  \ref{estimimpo1}, \ref{lemme6.5}, \ref{lemme6.6}, \ref{lemme11.3}, \eqref{ngtech5}, \eqref{paraestim1}
 and the fact that $\tilde{r}_{i,\sig}=O(1)\xi_i\bar{v}_i$, 
\begin{align*}
	&\B^{8,9}_{i,x}=-\big[v_{i,x}\left(\frac{w_i}{v_i}\right)_x+v_i\left(\frac{w_i}{v_i}\right)_{xx}
	\big]\left(\frac{w_i}{v_i}\right)_x(B-\mu_iI_n)\tilde{r}_{i,\si}\\
	&-\big[(\left(\frac{w_i}{v_i}\right)_x\mu_i^{-1}+\frac{w_i}{v_i}(\mu_i^{-1})_x)\mathcal{H}_i+\frac{w_i}{v_i}\mu_i^{-1}\mathcal{H}_{i,x}]\left(\frac{w_i}{v_i}\right)_x(B-\mu_iI_n)\tilde{r}_{i,\si}\\
	%-\big[(\left(\frac{w_i}{v_i}\right)_x\mu_i^{-1}+\frac{w_i}{v_i}(\mu_i^{-1})_x)(\mu_i v_{i,x}-(\tilde{\la}_i-\la_i^*)v_i-w_i)\big]\theta_i^\p\left(\frac{w_i}{v_i}\right)_x(B-\mu_iI_n)\tilde{r}_{i,\si}\\
	%&-\frac{w_i}{v_i}\mu_i^{-1}(\mu_i v_{i,x}-(\tilde{\la}_i-\la_i^*)v_i-w_i)_x\theta_i^\p\left(\frac{w_i}{v_i}\right)_x(B-\mu_iI_n)\tilde{r}_{i,\si}\\
	&-\big[v_i\left(\frac{w_i}{v_i}\right)_x+\frac{w_i}{v_i}\mu_i^{-1}\mathcal{H}_i\big]\big[\left(\frac{w_i}{v_i}\right)_{xx}(B-\mu_iI_n)\tilde{r}_{i,\si}+\left(\frac{w_i}{v_i}\right)_{x} u_x\cdot D(B-\mu_iI_n)\tilde{r}_{i,\si}\big]\\
	&-\big[v_i\left(\frac{w_i}{v_i}\right)_x+\frac{w_i}{v_i}\mu_i^{-1}\mathcal{H}_i\big]\theta_i^\p\left(\frac{w_i}{v_i}\right)_{x}(B-\mu_iI_n)\pa_x(\tilde{r}_{i,\si})
%\theta_i''\left(\frac{w_i}{v_i}\right)_x^2 (B-\mu_iI_n)\tilde{r}_{i,\si}\\
	%&-\big[v_i\left(\frac{w_i}{v_i}\right)_x+\frac{w_i}{v_i}\mu_i^{-1}(\mu_i v_{i,x}-(\tilde{\la}_i-\la_i^*)v_i-w_i)\big]\theta_i^\p\left(\frac{w_i}{v_i}\right)_{xx}(B-\mu_iI_n)\tilde{r}_{i,\si}\\
	%&-\big[v_i\left(\frac{w_i}{v_i}\right)_x+\frac{w_i}{v_i}\mu_i^{-1}(\mu_i v_{i,x}-(\tilde{\la}_i-\la_i^*)v_i-w_i)\big]\theta_i^\p\left(\frac{w_i}{v_i}\right)_{x}u_x\cdot D(B-\mu_iI_n)\tilde{r}_{i,\si}\\
	%&-\big[v_i\left(\frac{w_i}{v_i}\right)_x+\frac{w_i}{v_i}\mu_i^{-1}(\mu_i v_{i,x}-(\tilde{\la}_i-\la_i^*)v_i-w_i)\big]\theta_i^\p\left(\frac{w_i}{v_i}\right)_{x}(B-\mu_iI_n)\pa_x(\tilde{r}_{i,\si})\\
	=\sum_j(\La_j^1+\delta_0^2\La_j^3+\La_j^4+\La_j^5).
	\end{align*}
	%-(w_{i,x}-\mu_i^{-1}(\tilde{\la}_i-\la_i^*+\theta_i)w_i)_x\theta_i^\p\left(\frac{w_i}{v_i}\right)_x(B-\mu_iI_n)\tilde{r}_{i,\si}\\
	%&-(w_{i,x}-\mu_i^{-1}(\tilde{\la}_i-\la_i^*+\theta_i)w_i)\theta_i^{\p\p}\left(\frac{w_i}{v_i}\right)^2_x(B-\mu_iI_n)\tilde{r}_{i,\si}\\
	%&-(w_{i,x}-\mu_i^{-1}(\tilde{\la}_i-\la_i^*+\theta_i)w_i)\theta_i^\p\left(\frac{w_i}{v_i}\right)_{xx}(B-\mu_iI_n)\tilde{r}_{i,\si}\\
	%&-(w_{i,x}-\mu_i^{-1}(\tilde{\la}_i-\la_i^*+\theta_i)w_i)\theta_i^\p\left(\frac{w_i}{v_i}\right)_xu_x\cdot D(B-\mu_iI_n)\tilde{r}_{i,\si}\\
	%&-(w_{i,x}-\mu_i^{-1}(\tilde{\la}_i-\la_i^*+\theta_i)w_i)\theta_i^\p\left(\frac{w_i}{v_i}\right)_x(B-\mu_iI_n)(\tilde{r}_{i,\si})_x. %\left[\tilde{r}_{i,u\si}u_x+(\xi_i\bar{v}_i)_x\tilde{r}_{i,v\si}-\theta_i^\p\left(\frac{w_i}{v_i}\right)_x\tilde{r}_{i,\si\si}\right].
	From \eqref{estimate-v-i-xx-1aa} we get now $\B^{9,1}_{i}=
\mu_i^{-1}\mathcal{H}_{i,x}(B-\mu_iI_n)[w_i\xi_i\pa_{v_i}\bar{v}_i\tilde{r}_{i,v}]=\sum_j(\La_j^1+\La_j^4+\La_j^5+\La_j^6)+R_{\e,1}^{8,9,i,1}$ with $R_{\e,1}^{8,9,i,1}$ satisfying \eqref{12.11bis3}. Combining Lemmas \ref{lemme6.5}, \ref{lemme6.6}, \ref{lemme11.3}  and \eqref{ngtech5}, \eqref{estimate-v-i-xx-1aa}, \eqref{9.15primefin} we obtain 
\begin{align*}
	&\B^{9,1}_{i,x}=\big[((\mu_i^{-1})_x \mathcal{H}_{i,x}+\mu_i^{-1}\mathcal{H}_{i,xx})(B-\mu_iI_n)+\mu_i^{-1}\mathcal{H}_{i,x} u_x\cdot D(B-\mu_iI_n) \big][w_i\xi_i\pa_{v_i}\bar{v}_i\tilde{r}_{i,v}]\\
	%&\mu_i^{-1}(\mu_i v_{i,x}-(\tilde{\la}_i-\la_i^*)v_i-w_i)_{xx}(B-\mu_iI_n)[w_i\xi_i\pa_{v_i}\bar{v}_i\tilde{r}_{i,v}]\\
	%&+\mu_i^{-1}(\mu_i v_{i,x}-(\tilde{\la}_i-\la_i^*)v_i-w_i)_xu_x\cdot D(B-\mu_iI_n)[w_i\xi_i\pa_{v_i}\bar{v}_i\tilde{r}_{i,v}]\\
	&+\mu_i^{-1}\mathcal{H}_{i,x}(B-\mu_iI_n)\big[(w_{i,x}\xi_i+w_i \left(\frac{w_i}{v_i}\right)_x\xi_i')
	\pa_{v_i}\bar{v}_i\tilde{r}_{i,v}+w_i\xi_i((\pa_{v_i}\bar{v}_i)_x\tilde{r}_{i,v}+\pa_{v_i}\bar{v}_i \pa_x(\tilde{r}_{i,v})) ]\\
	%&+\mu_i^{-1}(\mu_iv_{i,x}-(\tilde{\la}_i-\la_i^*)v_i-w_i)_x\left(\frac{w_i}{v_i}\right)_x(B-\mu_iI_n)[w_i\xi^\p_i\pa_{v_i}\bar{v}_i\tilde{r}_{i,v}]\\
	%&+\mu_i^{-1}(\mu_iv_{i,x}-(\tilde{\la}_i-\la_i^*)v_i-w_i)_x(B-\mu_iI_n)[w_i\xi_i(\pa_{v_i}\bar{v}_i)_x\tilde{r}_{i,v}]\\
	%&+\mu_i^{-1}(\mu_i v_{i,x}-\mu_i^{-1}(\tilde{\la}_i-\la_i^*)v_i-w_i)_x[w_i\xi_i\pa_{v_i}\bar{v}_i](B-\mu_iI_n)\pa_x(\tilde{r}_{i,v})\\
	&=\sum_j(\La_j^1+\La^2_j+\delta_0^2\La_j^3+\La_j^4+\La_j^5+\La_j^6)+\widetilde{R}_{\e,1}^{8,9,i,1},
\end{align*}
with $\widetilde{R}_{\e,1}^{8,9,i,1}$ satisfying \eqref{12.11bis3}.
We here just give few details on how to deal with the term $\mu_i^{-1}\mathcal{H}_{i,xx}(B-\mu_iI_n)[w_i\xi_i\pa_{v_i}\bar{v}_i\tilde{r}_{i,v}]$. Applying Lemmas \ref{lemme6.5}, \ref{lemme6.6} and \eqref{6.45}, \eqref{9.15primefin} we deduce that
\begin{align*}
		&\mu_i^{-1}\mathcal{H}_{i,xx}(B-\mu_iI_n)[w_i\xi_i\pa_{v_i}\bar{v}_i\tilde{r}_{i,v}]
		=O(1)\biggl(%IMP\sum_k\sum_{l\ne k}(|v_{k,x}|+|w_{k,x}|+|v_k|)(|v_lv_{l,x}|+|v_kv_l|)+\sum_{k}\sum_{l\ne k}|v_k|(|v_{l,xx}v_l|+|v_kv_{l,x}|)
		%\\
		%&IMP+\sum_{k}\sum_{l\ne k}\sum_{p\ne k}|v_k| |v_{l,x}||v_{p,x}|
		%+\sum_k\sum_l\sum_{j\ne l}|v_kv_l||v_{j,x}v_j|\\
		%&+
		\sum_j(\La_j^1+\La_j^2+\delta_0^2\La_j^3+\La_j^4+\La_j^5+\La_j^6)\\
		&+
|w_i|(|v_iv_{i,x}|+|v_iv_{i,xx}|+|v_i^2|)\mathbbm{1}_{\{v_i^{2N}\leq 2\e\}}
+ |w_i v_{i,x} v_{i,xx}|\mathbbm{1}_{\{|\frac{w_i}{v_i}|\geq\frac{\delta_1}{2}\}}\nonumber\\
		%&+\mathfrak{A}_i\rho_i^\e\big(|v_iv_{i,x}\left(\frac{w_i}{v_i}\right)_x|+|v_i^2\left(\frac{w_i}{v_i}\right)^2_x|+|v_i^2\left(\frac{w_i}{v_i}\right)_x|\\
		%&+|w_{i,xx}v_i-v_{i,xx}w_i|\big)\\
		%&+\mathfrak{A}_i \Delta_i^\e\big(
		%|v_i v_{i,x}|+v_{i,x}^2+|v_iv_{i,xx}|\big)\\
		%&
		%+(|v_{i,x}|+|v_i|)v_i^2(1-{\color{black}{\chi_i^\e}}\xi_i\eta_i)\\
		%&+\mathfrak{A}_i\rho_i^\e\big(|v_i^3\left(\frac{w_i}{v_i}\right)_x|+
		%|v_i^2v_{i,x}|\Delta_i^\e\big)\\
		%&+(|v_{i,x}|+|v_i|)v_i^2(1-{\color{black}{\chi_i^\e}}\xi_i\eta_i)\\
		%&+\mathfrak{A}_i\rho_i^\e\big(|v_i^3\left(\frac{w_i}{v_i}\right)_x|+|v_i^2 v_{i,x}|\Delta_i^\e\big)+v_i^2v_{i,x}{\color{black}{(\chi')_i^\e}}\xi_i\eta_i\\
		%&+v_iv_{i,x}(\chi')_i^\e \mathfrak{A}_i\rho_i^\e
		%{\color{black}{+v_{i,x}^2(\chi')^\e_i\xi_i\eta_i }}+v_i^2v_{i,x}{\color{black}{(\chi')_i^\e}}\xi_i\eta_i
		%+(|v_i|+|v_{i,x}|)\varkappa^\ep_i\mathfrak{A}_i\rho_i^\e\\
		&+(|w_i v_i |(|v_{i,x}|+|v_{i,xx}|)+|w_iv_{i,x}|(|v_{i,x}|+|w_{i,x}|))\varkappa^\ep_i\mathfrak{A}_i\rho_i^\e\biggl)\\
		&=\sum_j(\La_j^1+\La_j^2+\delta_0^2\La_j^3+\La_j^4+\La_j^5+\La_j^6)+\widetilde{R}_{\e,1}^{8,9,i,1,1},
\end{align*}
with $\widetilde{R}_{\e,1}^{8,9,i,1,1}=O(1)|w_i|(|v_iv_{i,x}|+|v_iv_{i,xx}|+|v_i^2|)\mathbbm{1}_{\{v_i^{2N}\leq 2\e\}}+O(1)(|w_i v_i |(|v_{i,x}|+|v_{i,xx}|)+|w_iv_{i,x}|(|v_{i,x}|+|w_{i,x}|))\varkappa^\ep_i\mathfrak{A}_i\rho_i^\e$ satisfying \eqref{12.11bis3}.

We have now from \eqref{ngtech4} $\B^{9,2}_{i}=\mu_i^{-1}\mathcal{H}_i v_i\tilde{r}_i\cdot D(B-\mu_iI_n)[w_i\xi_i\pa_{v_i}\bar{v}_i\tilde{r}_{i,v}]=O(1)\sum_j(\La_j^1+\La_j^4+\La_j^6)+R_{\e,1}^{8,9,i,2}$ with $R_{\e,1}^{8,9,i,2}$ satisfying \eqref{12.11bis3}. Using Lemmas \ref{lemme6.8}, \ref{lemme6.5}, \ref{lemme6.5}, \ref{lemme11.3} and \eqref{ngtech4}, \eqref{estimate-v-i-xx-1aa} it yields
\begin{align*}
	&\B^{9,2}_{i,x}=\big[((\mu_i^{-1})_x\mathcal{H}_i+\mu_i^{-1}\mathcal{H}_{i,x})v_i\tilde{r}_i
	+\mu_i^{-1}\mathcal{H}_i( v_{i,x}\tilde{r}_i+v_i \pa_x(\tilde{r}_i))\big]
	\cdot D(B-\mu_iI_n)[w_i\xi_i\pa_{v_i}\bar{v}_i\tilde{r}_{i,v}]\\
	%&+\mu_i^{-1}(\mu_i v_{i,x}-(\tilde{\la}_i-\la^*_i)v_i-w_i)_xv_i\tilde{r}_i\cdot D(B-\mu_iI_n)[w_i\xi_i\pa_{v_i}\bar{v}_i\tilde{r}_{i,v}]\\
	%&+\mu_i^{-1}(\mu_i v_{i,x}-(\tilde{\la}_i-\la^*_i)v_i-w_i)v_{i,x}\tilde{r}_i\cdot D(B-\mu_iI_n)[w_i\xi_i\pa_{v_i}\bar{v}_i\tilde{r}_{i,v}]\\
	%&+\mu_i^{-1}(\mu_i v_{i,x}-(\tilde{\la}_i-\la^*_i)v_i-w_i)v_i \pa_x(\tilde{r}_i)\cdot D(B-\mu_iI_n)[w_i\xi_i\pa_{v_i}\bar{v}_i\tilde{r}_{i,v}]\\
	&+\mu_i^{-1}\mathcal{H}_i v_i\big[(\tilde{r}_i\otimes u_x):D^2(B-\mu_iI_n)w_i+ \tilde{r}_i\cdot D(B-\mu_iI_n)w_{i,x}\big]
	[\xi_i\pa_{v_i}\bar{v}_i\tilde{r}_{i,v}]\\
	%&+\mu_i^{-1}(\mu_i v_{i,x}-(\tilde{\la}_i-\la^*_i)v_i-w_i)v_i \tilde{r}_i\cdot D(B-\mu_iI_n)[w_{i,x}\xi_i\pa_{v_i}\bar{v}_i\tilde{r}_{i,v}]\\
	&+\mu_i^{-1}\mathcal{H}_i v_i\tilde{r}_i\cdot D(B-\mu_iI_n)w_i\big[\left(\frac{w_i}{v_i}\right)_x\xi_i^\p\pa_{v_i}\bar{v}_i\tilde{r}_{i,v}+\xi_i((\pa_{v_i}\bar{v}_i)_x\tilde{r}_{i,v}+\pa_{v_i}\bar{v}_i\pa_x(\tilde{r}_{i,v})) \big]\\
	%&+\mu_i^{-1}(\mu_i v_{i,x}-(\tilde{\la}_i-\la^*_i)v_i-w_i)v_i\tilde{r}_i\cdot D(B-\mu_iI_n)[w_i\xi_i(\pa_{v_i}\bar{v}_i)_x\tilde{r}_{i,v}+xi_i((\pa_{v_i}\bar{v}_i)_x\tilde{r}_{i,v}+\pa_{v_i}\bar{v}_i\pa_x(\tilde{r}_{i,v})\big]\\
	%&+\mu_i^{-1}(\mu_i v_{i,x}-(\tilde{\la}_i-\la^*_i)v_i-w_i)v_i[w_i\xi_i\pa_{v_i}\bar{v}_i]\tilde{r}_i\cdot D(B-\mu_iI_n) \pa_x(\tilde{r}_{i,v})\\
	&=O(1)\sum_j(\La_j^1+\La_j^4+\La_j^5+\La_j^6)+\widetilde{R}_{\e,1}^{8,9,i,2},
\end{align*}
with $\widetilde{R}_{\e,1}^{8,9,i,2}$ satisfying \eqref{12.11bis3}. We have now $\B^{9,3}_{i}=
\sum\limits_{j\neq i}\mu_i^{-1}\mathcal{H}_iv_j\tilde{r}_j\cdot D(B-\mu_iI_n)[w_i\xi_i\pa_{v_i}\bar{v}_i\tilde{r}_{i,v}]=O(1)\La_i^1$ and from Lemmas \ref{estimimpo1}, \ref{lemme6.5}, \ref{lemme6.6}, \ref{lemme11.3}
\begin{align*}
	&\B^{9,3}_{i,x}=\sum\limits_{j\neq i}\big[((\mu_i^{-1})_{x}\mathcal{H}_iv_j+\mu_i^{-1} \mathcal{H}_{i,x}v_j+\mu_i^{-1}\mathcal{H}_iv_{j,x}) \tilde{r}_j\cdot D(B-\mu_iI_n)[w_i\xi_i\pa_{v_i}\bar{v}_i\tilde{r}_{i,v}]\\
	&+\sum\limits_{j\neq i}\mu_i^{-1}\mathcal{H}_iv_j\big[
\tilde{r}_{j,x}\cdot D(B-\mu_iI_n)+(\tilde{r}_j\otimes u_x):D^2(B-\mu_iI_n)\big] [w_i\xi_i\pa_{v_i}\bar{v}_i\tilde{r}_{i,v}]\\
	%&+\sum\limits_{j\neq i}(v_{i,x}-\mu_i^{-1}(\tilde{\la}_i-\la^*_i+\theta_i)v_i)v_j \tilde{r}_{j,x}\cdot D(B-\mu_iI_n)[w_i\xi_i\pa_{v_i}\bar{v}_i\tilde{r}_{i,v}]\\
	%&+\sum\limits_{j\neq i}(v_{i,x}-\mu_i^{-1}(\tilde{\la}_i-\la^*_i+\theta_i)v_i)v_j(\tilde{r}_j\otimes u_x):D^2(B-\mu_iI_n)[w_i\xi_i\pa_{v_i}\bar{v}_i\tilde{r}_{i,v}]\\
	&+\sum\limits_{j\neq i}\mu_i^{-1}\mathcal{H}_iv_j \tilde{r}_j\cdot D(B-\mu_iI_n)\big[(w_{i,x}\xi_i+w_i \left(\frac{w_i}{v_i}\right)_x\xi_i')
	\pa_{v_i}\bar{v}_i\tilde{r}_{i,v}+w_i\xi_i(\pa_{v_i}\bar{v}_i)_x\tilde{r}_{i,v}
	\big]\\
	%&+\sum\limits_{j\neq i}(v_{i,x}-\mu_i^{-1}(\tilde{\la}_i-\la^*_i+\theta_i)v_i)v_j\left(\frac{w_i}{v_i}\right)_x\tilde{r}_j\cdot D(B-\mu_iI_n)[w_i\xi_i^\p\pa_{v_i}\bar{v}_i\tilde{r}_{i,v}]\\
	%&+\sum\limits_{j\neq i}(v_{i,x}-\mu_i^{-1}(\tilde{\la}_i-\la^*_i+\theta_i)v_i)v_j \tilde{r}_j\cdot D(B-\mu_iI_n)[w_i\xi_i(\pa_{v_i}\bar{v}_i)_x\tilde{r}_{i,v}]\\
	&+\sum\limits_{j\neq i}\mu_i^{-1}\mathcal{H}_i v_j \tilde{r}_j\cdot D(B-\mu_iI_n)[w_i\xi_i\pa_{v_i}\bar{v}_i](\tilde{r}_{i,v})_x
	=O(1)\La_i^1.
	%\left[\tilde{r}_{i,uv}u_x+(\xi_i\bar{v}_i)_x\tilde{r}_{i,vv}-\theta_i^\p\left(\frac{w_i}{v_i}\right)_x\tilde{r}_{i,v\si}\right].
\end{align*}
From \eqref{ngtech4} and Lemmas \ref{lemme6.5}, \ref{lemme6.6} we get $\B^{9,4}_{i}=\mu_i^{-1}\mathcal{H}_i v_{i}(B-\mu_iI_n)[w_i\xi_i\pa_{v_i}\bar{v}_i\tilde{r}_{i,uv}\tilde{r}_{i}]=O(1)\sum_j(\La_j^1+\La_j^4+\La_j^6)+R_{\e,1}^{8,9,i,4}$ with $R_{\e,1}^{8,9,i,4}$ satisfying \eqref{12.11bis3}. Combining Lemmas \ref{estimimpo1}, \ref{lemme6.5}, \ref{lemme6.6}, \ref{lemme11.4} and \eqref{ngtech4}, \eqref{estimate-v-i-xx-1aa}, we obtain 
\begin{align*}
	&\B^{9,4}_{i,x}=\big[((\mu_i^{-1})_x\mathcal{H}_i+\mu_i^{-1}\mathcal{H}_{i,x})v_{i}+\mu_i^{-1}\mathcal{H}_iv_{i,x}\big](B-\mu_iI_n)
	[w_i\xi_i\pa_{v_i}\bar{v}_i\tilde{r}_{i,uv}\tilde{r}_{i}]\\
	%&+\mu_i^{-1}(\mu_i v_{i,x}-(\tilde{\la}_i-\la_i^*)v_i-w_i)_xv_{i}(B-\mu_iI_n)[w_i\xi_i\pa_{v_i}\bar{v}_i\tilde{r}_{i,uv}\tilde{r}_{i}]\\
	%&+\mu_i^{-1}( \mu_iv_{i,x}-(\tilde{\la}_i-\la_i^*)v_i-w_i)v_{i,x}(B-\mu_iI_n)[w_i\xi_i\pa_{v_i}\bar{v}_i\tilde{r}_{i,uv}\tilde{r}_{i}]\\
	&+\mu_i^{-1}\mathcal{H}_iv_{i}\big[u_x\cdot D(B-\mu_iI_n)w_i \xi_i+(B-\mu_iI_n)w_{i,x}\xi_i +(B-\mu_iI_n)w_i\xi_i^\p\left(\frac{w_i}{v_i}\right)_x    \big]  [\pa_{v_i}\bar{v}_i\tilde{r}_{i,uv}\tilde{r}_{i}]\\
	&+\mu_i^{-1}\mathcal{H}_i v_{i}(B-\mu_iI_n)w_{i}\xi_i\big[(\pa_{v_i}\bar{v}_i)_x\tilde{r}_{i,uv}\tilde{r}_{i}+\pa_{v_i}\bar{v}_i\pa_x(\tilde{r}_{i,uv}\tilde{r}_i)\big]\\
	%&+\mu_i^{-1}(\mu_i v_{i,x}-(\tilde{\la}_i-\la_i^*)v_i-w_i)v_{i}\left(\frac{w_i}{v_i}\right)_x(B-\mu_iI_n)[w_i\xi_i^\p\pa_{v_i}\bar{v}_i\tilde{r}_{i,uv}\tilde{r}_{i}]\\
	%&+\mu_i^{-1}(\mu_i v_{i,x}-(\tilde{\la}_i-\la_i^*)v_i-w_i)v_{i}(B-\mu_iI_n)[w_i\xi_i(\pa_{v_i}\bar{v}_i)_x\tilde{r}_{i,uv}\tilde{r}_{i}]\\
	%&+\mu_i^{-1}(\mu_i v_{i,x}-(\tilde{\la}_i-\la_i^*)v_i-w_i)v_{i}[w_i\xi_i\pa_{v_i}\bar{v}_i](B-\mu_iI_n)\pa_x(\tilde{r}_{i,uv}\tilde{r}_i)\\
	%&+\mu_i^{-1}(\mu_i v_{i,x}-(\tilde{\la}_i-\la_i^*)v_i-w_i)v_{i}[w_i\xi_i\pa_{v_i}\bar{v}_i](B-\mu_iI_n)\tilde{r}_{i,u}\tilde{r}_{i,x}\\
	&=O(1)\sum_j(\La_j^1+\La_j^4+\La_j^5+\La_j^6)+\widetilde{R}_{\e,1}^{8,9,i,4},
\end{align*}
with $\widetilde{R}_{\e,1}^{8,9,i,4}$ satisfying \eqref{12.11bis3}.
Next we have  $\B^{9,5}_{i}=\sum\limits_{j\neq i}\mu_i^{-1}\mathcal{H}_iv_{j}(B-\mu_iI_n)[w_i\xi_i\pa_{v_i}\bar{v}_i\tilde{r}_{i,uv}\tilde{r}_{j}]=O(1)\La_i^1$ and from Lemma \ref{estimimpo1}, \ref{lemme6.5}, \ref{lemme6.6}, \ref{lemme11.4} it yields with $\widetilde{\mathcal{H}_i}=\mu_i^{-1}\mathcal{H}_i$
\begin{align*}
	&\B^{9,5}_{i,x}=\sum\limits_{j\neq i}\big[((\mu_i^{-1})_x\mathcal{H}_i+\mu_i^{-1}\mathcal{H}_{i,x})v_{j}+\mu_i^{-1}\mathcal{H}_iv_{j,x}\big]
	(B-\mu_iI_n)[w_i\xi_i\pa_{v_i}\bar{v}_i\tilde{r}_{i,uv}\tilde{r}_{j}]\\
	&+\sum\limits_{j\neq i}\widetilde{\mathcal{H}_i} v_{j}\big[\big(u_x\cdot D(B-\mu_iI_n)w_i +(B-\mu_iI_n)w_{i,x}\big)\xi_i +(B-\mu_iI_n)w_i\xi_i^\p\left(\frac{w_i}{v_i}\right)_x    \big]  [\pa_{v_i}\bar{v}_i\tilde{r}_{i,uv}\tilde{r}_{j}]\\
	&+\sum\limits_{j\neq i}\mu_i^{-1}\mathcal{H}_i v_{j}(B-\mu_iI_n)w_{i}\xi_i\big[(\pa_{v_i}\bar{v}_i)_x\tilde{r}_{i,uv}\tilde{r}_{j}+\pa_{v_i}\bar{v}_i\pa_x(\tilde{r}_{i,uv}\tilde{r}_j)\big]
%+\sum\limits_{j\neq i}(v_{i,x}-\mu_i^{-1}(\tilde{\la}_i-\la_i^*+\theta_i)v_i)v_{j,x}(B-\mu_iI_n)[w_i\xi_i\pa_{v_i}\bar{v}_i\tilde{r}_{i,uv}\tilde{r}_{j}]\\
	%&+\sum\limits_{j\neq i}(v_{i,x}-\mu_i^{-1}(\tilde{\la}_i-\la_i^*+\theta_i)v_i)v_{j}u_x\cdot D(B-\mu_iI_n)[w_i\xi_i\pa_{v_i}\bar{v}_i\tilde{r}_{i,uv}\tilde{r}_{j}]\\
	%&+\sum\limits_{j\neq i}(v_{i,x}-\mu_i^{-1}(\tilde{\la}_i-\la_i^*+\theta_i)v_i)v_{j}(B-\mu_iI_n)[w_{i,x}\xi_i\pa_{v_i}\bar{v}_i\tilde{r}_{i,uv}\tilde{r}_{j}]\\
	%&+\sum\limits_{j\neq i}(v_{i,x}-\mu_i^{-1}(\tilde{\la}_i-\la_i^*+\theta_i)v_i)v_{j}\left(\frac{w_i}{v_i}\right)_x(B-\mu_iI_n)[w_i\xi_i^\p\pa_{v_i}\bar{v}_i\tilde{r}_{i,uv}\tilde{r}_{j}]\\
	%&+\sum\limits_{j\neq i}(v_{i,x}-\mu_i^{-1}(\tilde{\la}_i-\la_i^*+\theta_i)v_i)v_{j}(B-\mu_iI_n)[w_i\xi_i\pa_x(\pa_{v_i}\bar{v}_i)\tilde{r}_{i,uv}\tilde{r}_{j}]\\
	%&+\sum\limits_{j\neq i}(v_{i,x}-\mu_i^{-1}(\tilde{\la}_i-\la_i^*+\theta_i)v_i)v_{j}[w_i\xi_i\pa_{v_i}\bar{v}_i]	(B-\mu_iI_n)\pa_x(\tilde{r}_{i,uv}\tilde{r}_j)\\
	=O(1)\La_i^1.
	%&+\sum\limits_{j\neq i}(v_{i,x}-\mu_i^{-1}(\tilde{\la}_i-\la_i^*+\theta_i)v_i)v_{j}[w_i\xi_i\pa_{v_i}\bar{v}_i](B-\mu_iI_n)\tilde{r}_{i,uv}\tilde{r}_{j,x}.
\end{align*}
We deduce from \eqref{ngtech4} that $\B^{9,6}_{i}=
\mu_i^{-1}\mathcal{H}_i v_{i,x}\xi_i(B-\mu_iI_n)[(w_i\pa_{v_iv_i}\bar{v}_i)\tilde{r}_{i,v}+w_i(\pa_{v_i}\bar{v}_i)^2\tilde{r}_{i,vv}]=O(1)\sum_j(\La_j^1+\La_j^4+\La_j^6)+R_{\e,1}^{8,9,i,6}$ with $R_{\e,1}^{8,9,i,6}$ satisfying \eqref{12.11bis3}. We get now from Lemmas \ref{estimimpo1}, \ref{lemme6.5}, \ref{lemme6.6}, \ref{lemme11.3} and \eqref{ngtech4}, \eqref{ngtech5},  \eqref{estimate-v-i-xx-1aa}
\begin{align*}
	&\B^{9,6}_{i,x}=\big[((\mu_i^{-1})_x\mathcal{H}_i+\mu_i^{-1}\mathcal{H}_{i,x})v_{i,x}+\mu_i^{-1}\mathcal{H}_iv_{i,xx}\big] \xi_i(B-\mu_iI_n)[(w_i\pa_{v_iv_i}\bar{v}_i)\tilde{r}_{i,v}+w_i(\pa_{v_i}\bar{v}_i)^2\tilde{r}_{i,vv}]\\
	&+\mu_i^{-1}\mathcal{H}_iv_{i,x}\big[u_x\cdot D(B-\mu_iI_n)w_i+(B-\mu_iI_n)w_{i,x}  \big]  \xi_i[\pa_{v_iv_i}\bar{v}_i\tilde{r}_{i,v}+(\pa_{v_i}\bar{v}_i)^2\tilde{r}_{i,vv}]\\
	&+ \mu_i^{-1}\mathcal{H}_i v_{i,x}\xi_i w_i(B-\mu_iI_n)\big[\pa_x (\pa_{v_iv_i}\bar{v}_i)\tilde{r}_{i,v}+\pa_{v_iv_i}\bar{v}_i\pa_x(\tilde{r}_{i,v})+2\pa_x (\pa_{v_i}\bar{v}_i)(\pa_{v_i}\bar{v}_i)\tilde{r}_{i,vv}]\\
	&
	+\mu_i^{-1}\mathcal{H}_i v_{i,x}w_i\xi_i'  \left(\frac{w_i}{v_i}\right)_x(B-\mu_iI_n)[(\pa_{v_iv_i}\bar{v}_i)\tilde{r}_{i,v}+(\pa_{v_i}\bar{v}_i)^2\tilde{r}_{i,vv}]\\
	&+\mu_i^{-1}\mathcal{H}_i v_{i,x} \xi_i w_i(B-\mu_iI_n)(\pa_{v_i}\bar{v}_i)^2\pa_x(\tilde{r}_{i,vv})
	=O(1)\sum_j(\La_j^1+\La_j^4+\La_j^5+\La_j^6)+\widetilde{R}_{\e,1}^{8,9,i,6},
\end{align*}
with $\widetilde{R}_{\e,1}^{8,9,i,6}$ satisfying \eqref{12.11bis3}.
We give some details here for the term $\mu_i^{-1}\mathcal{H}_iv_{i,x}\xi_i(B-\mu_iI_n)(w_{i,x}\pa_{v_iv_i}\bar{v}_i)\tilde{r}_{i,v}$, using Lemmas \ref{estimimpo1}, \ref{lemme6.5}, \ref{lemme6.6} and \eqref{ngtech4}, \eqref{ngtech5} we get
\begin{align*}
&\mu_i^{-1}\mathcal{H}_iv_{i,x}\xi_i(B-\mu_iI_n)(w_{i,x}\pa_{v_iv_i}\bar{v}_i)\tilde{r}_{i,v}
=\mu_i^{-1}\mathcal{H}_iv_{i,x}\left(\frac{w_i}{v_i}\right)_x\xi_i(B-\mu_iI_n)(v_i\pa_{v_iv_i}\bar{v}_i)\tilde{r}_{i,v}\\
&+\mu_i^{-1}\mathcal{H}_i\frac{w_i}{v_i}v_{i,x}^2(\pa_{v_iv_i}\bar{v}_i) \xi_i(B-\mu_iI_n)\tilde{r}_{i,v}\\
&=O(1)\big(v_i^2+\sum_{j\ne i}v_j^4\left(\frac{w_j}{v_j}\right)_x^2\mathfrak{A}_j\rho_j^\e )\big[\left(\frac{w_i}{v_i}\right)_x+\frac{\mathcal{H}_i}{v_i}\big]\mathfrak{A}_i\rho_i^\e
=\sum_j(\La_j^1+\La_j^4+\La_j^6)+\widetilde{R}_{\e,1}^{8,9,i,6,1},
\end{align*}
with $\widetilde{R}_{\e,1}^{8,9,i,6,1}$ satisfying \eqref{12.11bis3}.
We have now $\B^{9,7}_{i}=\sum\limits_{j\neq i}\mu_i^{-1}\mathcal{H}_i v_{j,x}\xi_iw_i(B-\mu_iI_n)[(\pa_{v_iv_j}\bar{v}_i)\tilde{r}_{i,v}+\pa_{v_i}\bar{v}_i\pa_{v_j}\bar{v}_i\tilde{r}_{i,vv}]=O(1)\La_i^1$ and from Lemmas \ref{estimimpo1}, \ref{lemme6.5}, \ref{lemme6.6}, \ref{lemme11.3} and \eqref{ngtech5}
\begin{align*}
	&\B^{9,7}_{i,x}=\sum\limits_{j\neq i}\big[((\mu_i^{-1})_x\mathcal{H}_i+\mu_i^{-1}\mathcal{H}_{i,x})v_{j,x}+\mu_i^{-1}\mathcal{H}_i v_{j,xx}\big]\xi_i w_i(B-\mu_iI_n)[\pa_{v_iv_j}\bar{v}_i\tilde{r}_{i,v}+\pa_{v_i}\bar{v}_i\pa_{v_j}\bar{v}_i\tilde{r}_{i,vv}]\\
	%&+\sum\limits_{j\neq i}(v_{i,x}-\mu_i^{-1}(\tilde{\la}_i-\la_i^*+\theta_i)v_i)v_{j,xx}\xi_i(B-\mu_iI_n)[(w_i\pa_{v_iv_j}\bar{v}_i)\tilde{r}_{i,v}+w_i\pa_{v_i}\bar{v}_i\pa_{v_j}\bar{v}_i\tilde{r}_{i,vv}]\\
	&+\sum\limits_{j\neq i}\mu_i^{-1}\mathcal{H}_i  v_{j,x}\big[\xi_i^\p\left(\frac{w_i}{v_i}\right)_xw_i+\xi_iw_{i,x}\big](B-\mu_iI_n)[\pa_{v_iv_j}\bar{v}_i\tilde{r}_{i,v}+\pa_{v_i}\bar{v}_i\pa_{v_j}\bar{v}_i\tilde{r}_{i,vv}]\\
	&+\sum\limits_{j\neq i}\mu_i^{-1}\mathcal{H}_i v_{j,x}\xi_i w_i u_x\cdot D(B-\mu_iI_n)[\pa_{v_iv_j}\bar{v}_i\tilde{r}_{i,v}+\pa_{v_i}\bar{v}_i\pa_{v_j}\bar{v}_i\tilde{r}_{i,vv}]\\
	&+\sum\limits_{j\neq i}\mu_i^{-1}\mathcal{H}_i v_{j,x}\xi_i w_i(B-\mu_iI_n)\big[(\pa_{v_iv_j}\bar{v}_i)_x\tilde{r}_{i,v}+\pa_{v_iv_j}\bar{v}_i(\tilde{r}_{i,v})_x+(\pa_{v_i}\bar{v}_i\pa_{v_j}\bar{v}_i)_x\tilde{r}_{i,vv}\big]\\
	&+\sum\limits_{j\neq i}\mu_i^{-1}\mathcal{H}_i v_{j,x}\xi_i w_i(B-\mu_iI_n)\pa_{v_i}\bar{v}_i\pa_{v_j}\bar{v}_i(\tilde{r}_{i,vv})_x
	%[(w_{i,x}\pa_{v_iv_j}\bar{v}_i)\tilde{r}_{i,v}+w_{i,x}\pa_{v_i}\bar{v}_i\pa_{v_j}\bar{v}_i\tilde{r}_{i,vv}]\\
	%&+\sum\limits_{j\neq i}(v_{i,x}-\mu_i^{-1}(\tilde{\la}_i-\la_i^*+\theta_i)v_i)v_{j,x}\xi_i(B-\mu_iI_n)[(w_i(\pa_{v_iv_j}\bar{v}_i))_x\tilde{r}_{i,v}+w_i(\pa_{v_i}\bar{v}_i\pa_{v_j}\bar{v}_i)_x\tilde{r}_{i,vv}]\\
	%&+\sum\limits_{j\neq i}(v_{i,x}-\mu_i^{-1}(\tilde{\la}_i-\la_i^*+\theta_i)v_i)v_{j,x}\xi_i(B-\mu_iI_n)[(w_i\pa_{v_iv_j}\bar{v}_i)(\tilde{r}_{i,v})_x+w_i\pa_{v_i}\bar{v}_i\pa_{v_j}\bar{v}_i(\tilde{r}_{i,vv})_x]\\
	=O(1)\La_i^1.
\end{align*}
 We have now from Lemmas \ref{lemme6.5}, \ref{lemme6.6} and \eqref{ngtech5} $\B^{9,8}_{i}=\mu_i^{-1}
 \mathcal{H}_i \xi^\p_i w_i\left(\frac{w_i}{v_i}\right)_x(B-\mu_iI_n)[\pa_{v_i}\bar{v}_i\tilde{r}_{i,v}+\xi_i\pa_{v_i}\bar{v}_i\bar{v}_i\tilde{r}_{i,vv}]=O(1)\sum_j\La_j^4$ and using Lemmas
  \ref{estimimpo1}, \ref{lemme6.5}, \ref{lemme6.6}, \ref{lemme11.3} and \eqref{ngtech5},   \eqref{estimate-v-i-xx-1aa}
\begin{align*}
	&\B^{9,8}_{i,x}=\big[((\mu_i^{-1})_x\mathcal{H}_i+\mu_i^{-1}\mathcal{H}_{i,x})
	\xi^\p_i+\widetilde{\mathcal{H}_i}\xi_i''\left(\frac{w_i}{v_i}\right)_x  \big]
	w_i\left(\frac{w_i}{v_i}\right)_x(B-\mu_iI_n)[\pa_{v_i}\bar{v}_i\tilde{r}_{i,v}+\xi_i\pa_{v_i}\bar{v}_i\bar{v}_i\tilde{r}_{i,vv}]\\
	&+\mu_i^{-1}\mathcal{H}_{i}\xi^\p_i\big[ (w_{i,x}\left(\frac{w_i}{v_i}\right)_x+w_i\left(\frac{w_i}{v_i}\right)_{xx})\big](B-\mu_iI_n) [\pa_{v_i}\bar{v}_i\tilde{r}_{i,v}+\xi_i\pa_{v_i}\bar{v}_i\bar{v}_i\tilde{r}_{i,vv}]\\
	%&+\mu_i^{-1}(\mu_i v_{i,x}-\mu_i^{-1}(\tilde{\la}_i-\la_i^*)v_i-w_i)_x\xi^\p_i\left(\frac{w_i}{v_i}\right)_x(B-\mu_iI_n)[(w_i\pa_{v_i}\bar{v}_i)\tilde{r}_{i,v}+\xi_iw_i\pa_{v_i}\bar{v}_i\bar{v}_i\tilde{r}_{i,vv}]\\
	%&+\mu_i^{-1}(\mu_i v_{i,x}-\mu_i^{-1}(\tilde{\la}_i-\la_i^*)v_i-w_i)\xi^{\p\p}_i\left(\frac{w_i}{v_i}\right)^2_x(B-\mu_iI_n)[(w_i\pa_{v_i}\bar{v}_i)\tilde{r}_{i,v}+\xi_iw_i\pa_{v_i}\bar{v}_i\bar{v}_i\tilde{r}_{i,vv}]\\
	%&+\mu_i^{-1}(\mu_i v_{i,x}-\mu_i^{-1}(\tilde{\la}_i-\la_i^*)v_i-w_i)\xi^\p_i\left(\frac{w_i}{v_i}\right)_{xx}(B-\mu_iI_n)[(w_i\pa_{v_i}\bar{v}_i)\tilde{r}_{i,v}+\xi_iw_i\pa_{v_i}\bar{v}_i\bar{v}_i\tilde{r}_{i,vv}]\\
	&+\mu_i^{-1}\mathcal{H}_i\xi^\p_iw_i\left(\frac{w_i}{v_i}\right)_{x}u_x\cdot D(B-\mu_iI_n)[\pa_{v_i}\bar{v}_i\tilde{r}_{i,v}+\xi_i\pa_{v_i}\bar{v}_i\bar{v}_i\tilde{r}_{i,vv}]\\
	&+\mu_i^{-1}\mathcal{H}_i\xi^\p_iw_i\left(\frac{w_i}{v_i}\right)_x(B-\mu_iI_n)[\pa_x (\pa_{v_i}\bar{v}_i)\tilde{r}_{i,v}+\pa_{v_i}\bar{v}_i(\tilde{r}_{i,v})_x+\xi_i((\pa_{v_i}\bar{v}_i)_x \bar{v}_i+\pa_{v_i}\bar{v}_i   (\bar{v}_i)_x) \tilde{r}_{i,vv}
	\big]\\
	&+\mu_i^{-1}\mathcal{H}_i\xi^\p_i w_i\left(\frac{w_i}{v_i}\right)_x(B-\mu_iI_n)\big[(\xi_i' \left(\frac{w_i}{v_i}\right)_x\tilde{r}_{i,vv}+  \xi_i(\tilde{r}_{i,vv})_x)  \pa_{\bar{v}_i}\bar{v}_i\bar{v}_i\big]\\
	%&+\mu_i^{-1}(\mu_i v_{i,x}-\mu_i^{-1}(\tilde{\la}_i-\la_i^*)v_i-w_i)\xi^\p_i\left(\frac{w_i}{v_i}\right)_x(B-\mu_iI_n)[(w_i\pa_{v_i}\bar{v}_i)(\tilde{r}_{i,v})_x+\xi_iw_i\pa_{v_i}\bar{v}_i\bar{v}_i(\tilde{r}_{i,vv})_x]\\
	&=O(1)\sum_j(\La_j^1+\delta_0^2\La_j^3+\La_j^4+\La_j^5+\La_j^6).
\end{align*}
We have now using \eqref{ngtech5} and the fact that $\theta'(x)=1$ on $\mbox{supp}\xi$ $\B^{9,9}_{i}=
-\mu_i^{-1}\mathcal{H}_i\left(\frac{w_i}{v_i}\right)_x(B-\mu_iI_n)[\xi_iw_i\pa_{v_i}\bar{v}_i\tilde{r}_{i,v\si}]=O(1)\sum_j\La_j^4$. Furthermore using  Lemmas \ref{estimimpo1}, \ref{lemme6.5}, \ref{lemme6.6}, \ref{lemme11.3} and \eqref{ngtech5},   \eqref{estimate-v-i-xx-1aa}
\begin{align*}
	&\B^{9,9}_{i,x}=-\big[((\mu_i^{-1})_x\mathcal{H}_i+\mu_i^{-1}\mathcal{H}_{i,x})
	\left(\frac{w_i}{v_i}\right)_x+\widetilde{\mathcal{H}_i}\left(\frac{w_i}{v_i}\right)_{xx}  \big](B-\mu_iI_n)[\xi_iw_i\pa_{v_i}\bar{v}_i\tilde{r}_{i,v\si}]\\
	%&-\mu_i^{-1}(v_{i,x}-\mu_i^{-1}(\tilde{\la}_i-\la_i^*)v_i-w_i)_x\theta_i^\p\left(\frac{w_i}{v_i}\right)_x(B-\mu_iI_n)[\xi_iw_i\pa_{v_i}\bar{v}_i\tilde{r}_{i,v\si}]\\
	%&-\mu_i^{-1}(\mu_i v_{i,x}-(\tilde{\la}_i-\la_i^*)v_i-w_i)\theta_i^{\p\p}\left(\frac{w_i}{v_i}\right)^2_x(B-\mu_iI_n)[\xi_iw_i\pa_{v_i}\bar{v}_i\tilde{r}_{i,v\si}]\\
	%&-\mu_i^{-1}(\mu_i v_{i,x}-\mu_i^{-1}(\tilde{\la}_i-\la_i^*)v_i-w_i)\theta_i^\p\left(\frac{w_i}{v_i}\right)_{xx}(B-\mu_iI_n)[\xi_iw_i\pa_{v_i}\bar{v}_i\tilde{r}_{i,v\si}]\\
	&-\mu_i^{-1}\mathcal{H}_i \left(\frac{w_i}{v_i}\right)_x\big[u_x\cdot D(B-\mu_iI_n)\xi_i
	+(B-\mu_iI_n)\xi^\p_i\left(\frac{w_i}{v_i}\right)_x\big]
	[w_i\pa_{v_i}\bar{v}_i\tilde{r}_{i,v\si}]\\
	&-\mu_i^{-1}\mathcal{H}_i\left(\frac{w_i}{v_i}\right)_x(B-\mu_iI_n)\xi_i\big[w_{i,x}\pa_{v_i}\bar{v}_i\tilde{r}_{i,v\si}+w_i((\pa_{v_i}\bar{v}_i)_x\tilde{r}_{i,v\si}+\pa_{v_i}\bar{v}_i(\tilde{r}_{i,v\si})_x)\big]\\
	%&-\mu_i^{-1}(v_{i,x}-\mu_i^{-1}(\tilde{\la}_i-\la_i^*)v_i-w_i)\theta_i^\p\left(\frac{w_i}{v_i}\right)^2_x(B-\mu_iI_n)[\xi^\p_iw_i\pa_{v_i}\bar{v}_i\tilde{r}_{i,v\si}]\\
	%&-\mu_i^{-1}(\mu_i v_{i,x}-\mu_i^{-1}(\tilde{\la}_i-\la_i^*)v_i-w_i)\theta_i^\p\left(\frac{w_i}{v_i}\right)_x(B-\mu_iI_n)[\xi_iw_{i,x}\pa_{v_i}\bar{v}_i\tilde{r}_{i,v\si}]\\
	%&-\mu_i^{-1}(\mu_i v_{i,x}-\mu_i^{-1}(\tilde{\la}_i-\la_i^*)v_i-w_i)\theta_i^\p\left(\frac{w_i}{v_i}\right)_x(B-\mu_iI_n)[\xi_iw_i(\pa_{v_i}\bar{v}_i)_x\tilde{r}_{i,v\si}]\\
	%&-\mu_i^{-1}(\mu_i v_{i,x}-\mu_i^{-1}(\tilde{\la}_i-\la_i^*)v_i-w_i)\theta_i^\p\left(\frac{w_i}{v_i}\right)_x[\xi_iw_i\pa_{v_i}\bar{v}_i](B-\mu_iI_n)(\tilde{r}_{i,v\si})_x\\
	&=O(1)\sum_j(\La_j^1+\delta_0^2\La_j^3+\La_j^4+\La_j^5+\La_j^6).
\end{align*}
To finish, it yields from Lemma \ref{lemme6.5} and \eqref{ngtech4} $\B^{9,10}_{i}=
\mu_i^{-1}\mathcal{H}_i(B-\mu_iI_n)[\xi_iw_{i,x}\pa_{v_i}\bar{v}_i\tilde{r}_{i,v}]=O(1)\sum_j(\La_j^1+\La_j^4+\La_j^6)+R_{\e,1}^{8,9,i,10}$ with $R_{\e,1}^{8,9,i,10}$ satisfying \eqref{12.11bis3}. We get finally from Lemmas \ref{estimimpo1}, \ref{lemme6.5}, \ref{lemme6.6}, \ref{lemme11.3}, \eqref{ngtech4}, \eqref{ngtech5}, \eqref{estimate-v-i-xx-1aa} and the fact that $w_{i,x}=v_i\left(\frac{w_i}{v_i}\right)_x+v_{i,x}$
\begin{align*}
	&\B^{9,10}_{i,x}=\big[((\mu_i^{-1})_x\mathcal{H}_i+\mu_i^{-1}\mathcal{H}_{i,x})
	(B-\mu_iI_n)+\widetilde{\mathcal{H}_i}u_x\cdot D(B-\mu_iI_n)   \big] [\xi_iw_{i,x}\pa_{v_i}\bar{v}_i\tilde{r}_{i,v}]\\
	&+\mu_i^{-1}\mathcal{H}_i  (B-\mu_iI_n)\big[ (\xi_i'\left(\frac{w_i}{v_i}\right)_x 
	w_{i,x}+\xi_i w_{i,xx})\pa_{v_i}\bar{v}_i\tilde{r}_{i,v}
	+\xi_i w_{i,x}( \pa_x(\pa_{v_i}\bar{v}_i)\tilde{r}_{i,v}+\pa_{v_i}\bar{v}_i\pa_x(\tilde{r}_{i,v}))\big] \\
	%&+\mu_i^{-1}(\mu_i v_{i,x}-(\tilde{\la}_i-\la_i^*)v_i-w_i)u_x\cdot D(B-\mu_iI_n)[\xi_iw_{i,x}\pa_{v_i}\bar{v}_i\tilde{r}_{i,v}]\\
	%&+\mu_i^{-1}(\mu_i v_{i,x}-(\tilde{\la}_i-\la_i^*)v_i-w_i)\left(\frac{w_i}{v_i}\right)_x (B-\mu_iI_n)[\xi_i'w_{i,x}\pa_{v_i}\bar{v}_i\tilde{r}_{i,v}]\\
	%&+\mu_i^{-1}(\mu_i v_{i,x}-(\tilde{\la}_i-\la_i^*)v_i-w_i)(B-\mu_iI_n)[\xi_iw_{i,xx}\pa_{v_i}\bar{v}_i\tilde{r}_{i,v}]\\
	%&+\mu_i^{-1}(\mu_i v_{i,x}-(\tilde{\la}_i-\la_i^*)v_i-w_i)(B-\mu_iI_n)[\xi_iw_{i,x}\pa_x(\pa_{v_i}\bar{v}_i)\tilde{r}_{i,v}]\\
	%&+\mu_i^{-1}(\mu_i v_{i,x}-(\tilde{\la}_i-\la_i^*)v_i-w_i)(B-\mu_iI_n)[\xi_iw_{i,x}\pa_{v_i}\bar{v}_i\pa_x(\tilde{r}_{i,v})]\\
		&=O(1)\sum_j(\La_j^1+\delta_0^2\La_j^3+\La_j^4+\La_j^5+\La_j^6)+\widetilde{R}_{\e,1}^{8,9,i,10},
\end{align*}
with $\widetilde{R}_{\e,1}^{8,9,i,10}$ satisfying \eqref{12.11bis3}.
\end{proof}

\begin{lemma}%\label{dernierK}
We have
	\begin{align}
			&K_{10,x}+\sum_{k=2}^{12}\sum_i (\frac{w_i}{v_i}-\la_i^*)\al_{i}^{6,k}+\sum_{j=13}^{15}K_{j,x}\nonumber\\
			&\hspace{3cm}=O(1)\sum_j(\La_j^1+\La_j^2+\delta_0^2\La_j^3+\La_j^4+\La_j^5+\La_j^6+\La_j^{6,1})+R_{\e,1}^{10},\label{K11x}\\
			&K_{10,xx}+\sum_{k=2}^{12}\sum_i \big[(\frac{w_i}{v_i}-\la_i^*)\al_{i}^{6,k}\big]_x +\sum_{j=13}^{15}K_{j,xx}\nonumber \\
			&\hspace{3cm}=O(1)\sum_j(\La_j^1+\La_j^2+\delta_0^2\La_j^3+\La_j^4+\La_j^5+\La_j^6+\La_j^{6,1})+\widetilde{R}_{\e,1}^{5,10}\label{K11xx}
						\end{align}
with:
	\begin{equation}
	\int_{\hat{t}}^T\int_{\R}(|R_{\e,1}^{10}|+|\widetilde{R}_{\e,1}^{10}|)dx ds=O(1)\delta_0^2,
	\label{12.11bis4}
	\end{equation}
	for $\e>0$ small enough in terms of $T-\hat{t}$ and $\delta_0$.
\end{lemma}
\begin{proof}
First since
\begin{align*}
	K_{10,x}&=-\sum\limits_{k=1}^{9}\sum\limits_{i}\la_i^*\alpha_{i}^{5,k},
	%-\sum\limits_{k=1}^{10}\sum\limits_{i}\la_i^*\alpha_{i}^{6,k},
	\end{align*}
	we deduce the result for $K_{10,x}, K_{10,xx}$ simply by applying the Lemma \ref{lemmeal6}.  Let us deal now with the term $\sum_{k=2}^{12}\sum_i (\frac{w_i}{v_i}-\la_i^*)\al_{i}^{6,k}$, since each term $\al_{i}^{6,k}$ is a multiple of $\xi_i, \xi'_i, \theta'_i$ or $\theta''_i$ we deduce easily \eqref{K11x} by applying the Lemma \ref{lemmeal7}. We want now to estimate $\sum_{k=2}^{12}\sum_i \big[(\frac{w_i}{v_i}-\la_i^*)\al_{i}^{6,k}\big]_x$ and again from the Lemma \ref{lemmeal7}, it suffices to consider the terms 
	$\sum_{k=2}^{12}\sum_i \left(\frac{w_i}{v_i}\right)_x\al_{i}^{6,k}$.
	 It yields from the Lemmas \ref{estimimpo1}, \ref{lemme6.5}, \ref{lemme6.6} and since $\tilde{r}_{i,\si}, \tilde{r}_{i,\si\si},\tilde{r}_{i,u\si}=O(1)\xi_i\bar{v}_i$
\begin{align*}
\left(\frac{w_i}{v_i}\right)_x\al_{i}^{6,2}&=v_i \left(\frac{w_i}{v_i}\right)^3_x(B-\mu_iI_n)[\xi^{\p\p}_i\bar{v}_i\tilde{r}_{i,v}-\theta^{\p\p}_i\tilde{r}_{i,\si}]=O(1)\delta_0^2\La_i^3,\\
\left(\frac{w_i}{v_i}\right)_x\al_{i}^{6,3}&=\left(\frac{w_i}{v_i}\right)_x(w_{i,x}v_i-w_iv_{i,x})\tilde{r}_i\cdot D(B-\mu_iI_n)[\xi_i^\p\bar{v}_i\tilde{r}_{i,v}-\theta^\p_i\tilde{r}_{i,\si}]=O(1)\La_i^4,\\
\left(\frac{w_i}{v_i}\right)_x\al_{i}^{6,4}&=\sum\limits_{j\neq i}v_i \left(\frac{w_i}{v_i}\right)^2_x v_j\tilde{r}_j\cdot D(B-\mu_iI_n)[\xi_i^\p\bar{v}_i\tilde{r}_{i,v}-\theta^\p_i\tilde{r}_{i,\si}]=O(1)\delta_0^2\La_i^3,\\
\left(\frac{w_i}{v_i}\right)_x\al_{i}^{6,5}&=\left(\frac{w_i}{v_i}\right)_x (w_{i,x}v_i-w_iv_{i,x}) (B-\mu_iI_n)[\xi_i^\p\bar{v}_i\tilde{r}_{i,uv}\tilde{r}_i-\theta^\p_i\tilde{r}_{i,u\si}\tilde{r}_i]=O(1)\La_i^4\\
\left(\frac{w_i}{v_i}\right)_x\al_{i}^{6,6}&=\sum\limits_{j\neq i}v_i \left(\frac{w_i}{v_i}\right)_x^2 v_j(B-\mu_iI_n)[\xi_i^\p\bar{v}_i\tilde{r}_{i,uv}\tilde{r}_j-\theta^\p_i\tilde{r}_{i,u\si}\tilde{r}_j]=O(1)\delta_0^2\La_i^3,\\
\left(\frac{w_i}{v_i}\right)_x\al_{i}^{6,7}&=v_i \left(\frac{w_i}{v_i}\right)_x^2 v_{i,x}\xi_i\pa_{v_i}\bar{v}_i(B-\mu_iI_n)[\xi_i^\p\bar{v}_i\tilde{r}_{i,vv}-\theta^\p_i\tilde{r}_{i,v\si}]=O(1)\delta_0^2\La_i^3,\\
\left(\frac{w_i}{v_i}\right)_x\al_{i}^{6,8}&=\sum\limits_{j\neq i}v_i \left(\frac{w_i}{v_i}\right)_x^2 v_{j,x}\xi_i\pa_{v_j}\bar{v}_i\theta^\p_i(B-\mu_iI_n)[\xi_i^\p\bar{v}_i\tilde{r}_{i,vv}-\theta^\p_i\tilde{r}_{i,v\si}]=O(1)\delta_0^2\La_i^3,\\
%\left(\frac{w_i}{v_i}\right)_x\al_{i}^{7,9}&=\sum\limits_{j\neq i}v_i \left(\frac{w_i}{v_i}\right)_x^2 w_{j,x}\xi_i\pa_{w_j}\bar{v}_i\theta^\p_i(B-\mu_iI_n)[\xi_i^\p\bar{v}_i\tilde{r}_{i,vv}-\theta^\p_i\tilde{r}_{i,v\si}]=O(1)\delta_0^2\La_i^3,\\
\left(\frac{w_i}{v_i}\right)_x\al_{i}^{6,9}&=v_i \left(\frac{w_i}{v_i}\right)^3_x\xi_i^\p\bar{v}_i(B-\mu_iI_n)[\xi_i^\p\bar{v}_i\tilde{r}_{i,vv}-\theta^\p_i\tilde{r}_{i,v\si}]=O(1)\delta_0^2\La_i^3,\\
\left(\frac{w_i}{v_i}\right)_x\al_{i}^{6,10}&=-\theta^\p_iv_i \left(\frac{w_i}{v_i}\right)_x^3(B-\mu_iI_n)[\xi_i^\p\bar{v}_i\tilde{r}_{i,v\si}-\theta^\p_i\tilde{r}_{i,\si\si}]=O(1)\delta_0^2\La_i^3,\\
\left(\frac{w_i}{v_i}\right)_x\al_{i}^{6,11}&=v_i \left(\frac{w_i}{v_i}\right)_x^2 v_{i,x}(B-\mu_iI_n)\xi_i^\p\pa_{v_i}\bar{v}_i\tilde{r}_{i,v}=O(1)\delta_0^2\La_i^3,\\
\left(\frac{w_i}{v_i}\right)_x\al_{i}^{6,12}&=\sum\limits_{j\neq i} v_i \left(\frac{w_i}{v_i}\right)_x^2 v_{j,x}(B-\mu_iI_n)\xi_i^\p\pa_{v_j}\bar{v}_i\tilde{r}_{i,v}=O(1)\delta_0^2\La_i^3.
%\left(\frac{w_i}{v_i}\right)_x\al_{i}^{7,14}&=\sum\limits_{j\neq i}v_i \left(\frac{w_i}{v_i}\right)_x^2 w_{j,x}(B-\mu_iI_n)\xi_i^\p\pa_{w_j}\bar{v}_i\tilde{r}_{i,v}=O(1)\delta_0^2\La_i^3.
\end{align*}
We are now going to deal with the terms $K_{13,x}, K_{14,x}$ and $K_{15,x}$. We recall first that we have
\begin{align*}
K_{13,x}&=\sum\limits_{i}\mu_i^{-1}(\tilde{\la}_i-\si_i)\left(\frac{w_i}{v_i}\right)_xv_i\xi_i[(\pa_{v_i}\bar{v}_i)v_i-\bar{v}_i]B\tilde{r}_{i,v}+\sum\limits_{k=1}^{17}\sum\limits_{i}\left(\frac{w_i}{v_i}-\la_i^*\right)\alpha_{i}^{7,k}.
	\end{align*}
First we have from Lemmas \ref{estimimpo1}, \ref{lemme6.5}, \ref{lemme6.6}, \ref{lemme11.3} and \eqref{ngtech5}
\begin{align*}
&\mu_i^{-1}(\tilde{\la}_i-\si_i)\left(\frac{w_i}{v_i}\right)_xv_i\xi_i[(\pa_{v_i}\bar{v}_i)v_i-\bar{v}_i]B\tilde{r}_{i,v}=O(1)\La_i^4,\\
&\big[\mu_i^{-1}(\tilde{\la}_i-\si_i)\left(\frac{w_i}{v_i}\right)_x v_i\xi_i[\pa_{v_i}\bar{v}_i\,v_i-\bar{v}_i]B\tilde{r}_{i,v}\big]_x=(\mu_i^{-1})_x(\tilde{\la}_i-\si_i)\left(\frac{w_i}{v_i}\right)_xv_i\xi_i[\pa_{v_i}\bar{v}_i\,v_i-\bar{v}_i]B\tilde{r}_{i,v}\\
&+\mu_i^{-1}\big[(\tilde{\la}_{i,x}+\theta'_i\left(\frac{w_i}{v_i}\right)_x)\left(\frac{w_i}{v_i}\right)_x+(\tilde{\la}_i-\si_i)\left(\frac{w_i}{v_i}\right)_{xx}\big]
v_i\xi_i[(\pa_{v_i}\bar{v}_i)v_i-\bar{v}_i]B\tilde{r}_{i,v}\\
%&+\mu_i^{-1}(\tilde{\la}_i-\si_i)\left(\frac{w_i}{v_i}\right)_{xx}v_i\xi_i[(\pa_{v_i}\bar{v}_i)v_i-\bar{v}_i]B\tilde{r}_{i,v}\\
&+\mu_i^{-1}(\tilde{\la}_i-\si_i)\left(\frac{w_i}{v_i}\right)_x\big[v_{i,x}\xi_i+v_i\xi_i' \left(\frac{w_i}{v_i}\right)_x\big]
[(\pa_{v_i}\bar{v}_i)v_i-\bar{v}_i]B\tilde{r}_{i,v}\\
%&+\mu_i^{-1}(\tilde{\la}_i-\si_i)\left(\frac{w_i}{v_i}\right)_x^2 v_i\xi_i'[(\pa_{v_i}\bar{v}_i)v_i-\bar{v}_i]B\tilde{r}_{i,v}\\
&+\mu_i^{-1}(\tilde{\la}_i-\si_i)\left(\frac{w_i}{v_i}\right)_xv_i\xi_i\big[\pa_x(\pa_{v_i}\bar{v}_i \,v_i-\bar{v}_i)B\tilde{r}_{i,v}+(\pa_{v_i}\bar{v}_i\,v_i-\bar{v}_i)u_x\cdot DB(u)\tilde{r}_{i,v}\big]\\
%&+\mu_i^{-1}(\tilde{\la}_i-\si_i)\left(\frac{w_i}{v_i}\right)_xv_i\xi_i\big[(\pa_{v_i}\bar{v}_i)v_i-\bar{v}_i]u_x\cdot DB(u)\tilde{r}_{i,v}\\
&+\mu_i^{-1}(\tilde{\la}_i-\si_i)\left(\frac{w_i}{v_i}\right)_xv_i\xi_i( \pa_{v_i}\bar{v}_i\,v_i-\bar{v}_i)B\pa_x(\tilde{r}_{i,v})=O(1)\sum_j(\La_j^1+\delta_0^2\La_j^3+\La_j^4+\La_j^5).
\end{align*}
Now since each term $\al_{i}^{7,l}$ is a multiple of $\xi_i$ or $\xi'_i$ we deduce easily from Lemma \ref{lemmeal8}, that $K_{14,x}$ satisfies the estimate \eqref{K11x}. As previously we want now to verify that $K_{14,xx}$ satisfies \eqref{K11xx}. To do this, we need to estimates the terms $\sum\limits_{k=1}^{17}\sum\limits_{i}\left(\frac{w_i}{v_i}\right)_x\alpha_{i}^{7,k}$. We have then using Lemmas \ref{estimimpo1}, \ref{lemme6.5}, \ref{lemme6.6} and \eqref{ngtech5}
\begin{align*}
\left(\frac{w_i}{v_i}\right)_x\al_{i}^{7,1}&=\sum\limits_{k}\left(\frac{w_i}{v_i}\right)_x v_k\tilde{r}_k\cdot (D_u\mu_i^{-1})(\tilde{\la}_i-\si_i)v_i\xi_i[(\pa_{v_i}\bar{v}_i)v_i-\bar{v}_i]B\tilde{r}_{i,v}=O(1)\La_i^4,\\
\left(\frac{w_i}{v_i}\right)_x\al_{i}^{7,2}&=\mu_i^{-1}(\tilde{\la}_i-\si_i)\xi^\p_i\left(\frac{w_i}{v_i}\right)_x^2v_i[(\pa_{v_i}\bar{v}_i)v_i-\bar{v}_i]B\tilde{r}_{i,v}=O(1)\delta_0^2\La_i^3,\\
\left(\frac{w_i}{v_i}\right)_x\al_{i}^{7,3}&=\mu_i^{-1}(\tilde{\la}_i-\si_i)\xi_i \left(\frac{w_i}{v_i}\right)_x v_{i,x}[(\pa_{v_i}\bar{v}_i)v_i-\bar{v}_i]B\tilde{r}_{i,v}=O(1)(\La_i^1+\La_i^4),\\
\left(\frac{w_i}{v_i}\right)_x\al_{i}^{7,4}&=\sum\limits_{k}\mu_i^{-1}\left(\frac{w_i}{v_i}\right)_x\tilde{r}_k\cdot D \tilde{\la}_iv_kv_i \xi_i[(\pa_{v_i}\bar{v}_i)v_i-\bar{v}_i]B\tilde{r}_{i,v}=O(1)\La_i^4,\\
\left(\frac{w_i}{v_i}\right)_x\al_{i}^{7,5}&=\left(\frac{w_i}{v_i}\right)_x\mu_i^{-1}\tilde{\la}_{i,v}\pa_{v_i}\bar{v}_iv_{i,x}v_i\xi_i^2[(\pa_{v_i}\bar{v}_i)v_i-\bar{v}_i]B\tilde{r}_{i,v}=O(1)\La_i^4,\\
\left(\frac{w_i}{v_i}\right)_x\al_{i}^{7,6}&=\sum\limits_{j\neq i}\left(\frac{w_i}{v_i}\right)_x\mu_i^{-1}\tilde{\la}_{i,v}\pa_{v_j}\bar{v}_iv_{j,x}v_i \xi_i^2[(\pa_{v_i}\bar{v}_i)v_i-\bar{v}_i]B\tilde{r}_{i,v}=O(1)\La_i^4,\\
%\left(\frac{w_i}{v_i}\right)_x\al_{i}^{8,7}&=\sum\limits_{j\neq i}\left(\frac{w_i}{v_i}\right)_x\mu_i^{-1}\tilde{\la}_{i,v}\pa_{w_j}\bar{v}_iw_{j,x}v_i\xi_i^2[(\pa_{v_i}\bar{v}_i)v_i-\bar{v}_i]B\tilde{r}_{i,v}=O(1)\La_i^4,\\
\left(\frac{w_i}{v_i}\right)_x\al_{i}^{7,7}&=\mu_i^{-1}\tilde{\la}_{i,v}\xi_i^\p\left(\frac{w_i}{v_i}\right)_x^2\bar{v}_i v_i\xi_i[(\pa_{v_i}\bar{v}_i)v_i-\bar{v}_i]B\tilde{r}_{i,v}=O(1)\La_i^4,\\
\left(\frac{w_i}{v_i}\right)_x\al_{i}^{7,8}&=-\mu_i^{-1}\tilde{\la}_{i,\si}\theta_i^\p\left(\frac{w_i}{v_i}\right)_x^2 v_i\xi_i[(\pa_{v_i}\bar{v}_i)v_i-\bar{v}_i]B\tilde{r}_{i,v}=O(1)\delta_0^2\La_i^3,\\
\left(\frac{w_i}{v_i}\right)_x\al_{i}^{7,9}&=\mu_i^{-1}\theta_i^\p\left(\frac{w_i}{v_i}\right)_x^2 v_i\xi_i[(\pa_{v_i}\bar{v}_i)v_i-\bar{v}_i]B\tilde{r}_{i,v}=O(1)\delta_0^2\La_i^3,\\
\left(\frac{w_i}{v_i}\right)_x\al_{i}^{7,10}&=\left(\frac{w_i}{v_i}\right)_x\mu_i^{-1}(\tilde{\la}_i-\si_i)v^2_iv_{i,x}\xi_i(\pa_{v_iv_i}\bar{v}_i)B\tilde{r}_{i,v}=O(1)\sum_j\La_j^4,\\
\left(\frac{w_i}{v_i}\right)_x\al_{i}^{7,11}&=\sum\limits_{j\neq i}\left(\frac{w_i}{v_i}\right)_x\mu_i^{-1}(\tilde{\la}_i-\si_i)v_iv_{j,x}\xi_i[(\pa_{v_iv_j}\bar{v}_i)v_i-\pa_{v_j}\bar{v}_i]B\tilde{r}_{i,v}=O(1)\La_i^1,\\
%\left(\frac{w_i}{v_i}\right)_x\al_{i}^{8,13}&=\sum\limits_{j\neq i}\left(\frac{w_i}{v_i}\right)_x\mu_i^{-1}(\tilde{\la}_i-\si_i)v_iw_{j,x}\xi_i[(\pa_{v_iw_j}\bar{v}_i)v_i-\pa_{w_j}\bar{v}_i]B\tilde{r}_{i,v}=O(1)\La_i^1,\\
\left(\frac{w_i}{v_i}\right)_x\al_{i}^{7,12}&=\left(\frac{w_i}{v_i}\right)_x\mu_i^{-1}(\tilde{\la}_i-\si_i)v_i\xi_i[(\pa_{v_i}\bar{v}_i)v_i-\bar{v}_i]u_x\cdot DB\tilde{r}_{i,v}=O(1)\La_i^4,\\
\left(\frac{w_i}{v_i}\right)_x\al_{i}^{7,13}&=\sum\limits_{k}\left(\frac{w_i}{v_i}\right)_x\mu_i^{-1}(\tilde{\la}_i-\si_i)v_i\xi_i v_k[(\pa_{v_i}\bar{v}_i)v_i-\bar{v}_i]B\tilde{r}_{i,uv}\tilde{r}_k=O(1)\La_i^4,\\
\left(\frac{w_i}{v_i}\right)_x\al_{i}^{7,14}&= \left(\frac{w_i}{v_i}\right)_x\mu_i^{-1}(\tilde{\la}_i-\si_i)v_iv_{i,x}\xi_i^2\pa_{v_i}\bar{v}_i[(\pa_{v_i}\bar{v}_i)v_i-\bar{v}_i]B\tilde{r}_{i,vv}=O(1)\La_i^4,\\
\left(\frac{w_i}{v_i}\right)_x\al_{i}^{7,15}&=\sum\limits_{j\neq i}\left(\frac{w_i}{v_i}\right)_x\mu_i^{-1}(\tilde{\la}_i-\si_i)v_iv_{j,x}\xi_i^2\pa_{v_j}\bar{v}_i[(\pa_{v_i}\bar{v}_i)v_i-\bar{v}_i]B\tilde{r}_{i,vv}=O(1)\La_i^4,\\
%\left(\frac{w_i}{v_i}\right)_x\al_{i}^{8,18}&=\sum\limits_{j\neq i}\left(\frac{w_i}{v_i}\right)_x\mu_i^{-1}(\tilde{\la}_i-\si_i)v_iw_{j,x}\xi_i^2\pa_{w_j}\bar{v}_i[(\pa_{v_i}\bar{v}_i)v_i-\bar{v}_i]B\tilde{r}_{i,vv}=O(1)\La_i^4,\\
\left(\frac{w_i}{v_i}\right)_x\al_{i}^{7,16}&=\mu_i^{-1}(\tilde{\la}_i-\si_i)v_i\xi_i\xi_i^\p\bar{v}_i\left(\frac{w_i}{v_i}\right)^2_x[(\pa_{v_i}\bar{v}_i)v_i-\bar{v}_i]B\tilde{r}_{i,vv}=O(1)\delta_0^2\La_i^3,\\
\left(\frac{w_i}{v_i}\right)_x\al_{i}^{7,17}&=-\mu_i^{-1}(\tilde{\la}_i-\si_i)v_i\theta_i^\p\left(\frac{w_i}{v_i}\right)_x^2 \xi_i[(\pa_{v_i}\bar{v}_i)v_i-\bar{v}_i]B\tilde{r}_{i,v\si}=O(1)\delta_0^2\La_i^3.
\end{align*}
Similarly we consider now $K_{14,x}$ and we recall that
\begin{align*}
K_{14,x}&=\sum\limits_{i}v_i^2\left(\frac{w_i}{v_i}\right)_x(1-\xi_i\chi_i^\e \eta_i)\left[B(u)\tilde{r}_{i,u}\tilde{r}_i+\tilde{r}_i\cdot DB(u)\tilde{r}_i\right]+\sum\limits_{k=1}^{4}\sum\limits_{i}\left(\frac{w_i}{v_i}-\la_i^*\right)\al_{i}^{8,k}.
\end{align*}
Applying Lemmas \ref{estimimpo1}, \ref{lemme6.5}, \ref{lemme6.6}, \ref{lemme11.3}, \ref{lemme11.4} and \eqref{ngtech5} we get%and the fact that $\pa_{v_i}\eta_i=\frac{O(1)}{|v_i|}\alpha_i$, $\pa_{v_k}\eta_i=O(1)\frac{|v_k|}{|v_i|}\alpha_i$, $\pa_{w_k}\eta_i=O(1)\frac{|w_k|}{|v_i|}\alpha_i$
\begin{align*}
&v_i^2\left(\frac{w_i}{v_i}\right)_x(1-{\color{black}{\chi_i^\e}}\xi_i\eta_i)\left[B(u)\tilde{r}_{i,u}\tilde{r}_i+\tilde{r}_i\cdot DB(u)\tilde{r}_i\right]=O(1)\La_i^4\\
&\big[v_i^2\left(\frac{w_i}{v_i}\right)_x(1-{\color{black}{\chi_i^\e}}\xi_i\eta_i)\left[B(u)\tilde{r}_{i,u}\tilde{r}_i+\tilde{r}_i\cdot DB(u)\tilde{r}_i\right]\big]_x\\
&=(w_{i,xx}v_i-v_{i,xx}w_i)(1-{\color{black}{\chi_i^\e}}\xi_i \eta_i)\left[B(u)\tilde{r}_{i,u}\tilde{r}_i+\tilde{r}_i\cdot DB(u)\tilde{r}_i\right]\\
&-v_i^2\left(\frac{w_i}{v_i}\right)_x\big[\xi'_i\left(\frac{w_i}{v_i}\right)_x\chi_i^\e \eta_i+2N\frac{v_i^{2N-1}}{\e}(\chi')_i^\e v_{i,x} \xi_i\eta_i+\xi_i\chi_i^\e\pa_x(\eta_i)\big]\left[B(u)\tilde{r}_{i,u}\tilde{r}_i+\tilde{r}_i\cdot DB(u)\tilde{r}_i\right]\\
&+v_i^2\left(\frac{w_i}{v_i}\right)_x(1-{\color{black}{\chi_i^\e}}\xi_i\eta_i)\big[u_x\cdot DB(u)\tilde{r}_{i,u}\tilde{r}_i+B(u)\pa_x(\tilde{r}_{i,u}\tilde{r}_i)+\tilde{r}_i\otimes u_x: D^2B(u)\tilde{r}_i\big]\\
&+v_i^2\left(\frac{w_i}{v_i}\right)_x(1-{\color{black}{\chi_i^\e}}\xi_i\eta_i)\big[\pa_x(\tilde{r}_i)\cdot DB(u)\tilde{r}_i+
\tilde{r}_i\cdot DB(u)\pa_x(\tilde{r}_i)\big]
%-(w_{i,x}v_i-v_{i,x}w_i)\xi'_i\left(\frac{w_i}{v_i}\right)_x\eta_i\left[B(u)\tilde{r}_{i,u}\tilde{r}_i+\tilde{r}_i\cdot DB(u)\tilde{r}_i\right]\\
%&-(w_{i,x}v_i-v_{i,x}w_i)\xi_i\pa_x(\eta_i)\left[B(u)\tilde{r}_{i,u}\tilde{r}_i+\tilde{r}_i\cdot DB(u)\tilde{r}_i\right]\\
%&=\left(\frac{w_i}{v_i}\right)_{xx}(v_i-\xi_i\bar{v}_i)v_i\left[B(u)\tilde{r}_{i,u}\tilde{r}_i+\tilde{r}_i\cdot DB(u)\tilde{r}_i\right]\\
%&+\left(\frac{w_i}{v_i}\right)_x\pa_x (v_i-\xi_i\bar{v}_i)v_i\left[B(u)\tilde{r}_{i,u}\tilde{r}_i+\tilde{r}_i\cdot DB(u)\tilde{r}_i\right]\\
%&+\left(\frac{w_i}{v_i}\right)_x(v_i-\xi_i\bar{v}_i)v_{i,x}\left[B(u)\tilde{r}_{i,u}\tilde{r}_i+\tilde{r}_i\cdot DB(u)\tilde{r}_i\right]\\
%&+\left(\frac{w_i}{v_i}\right)_x v_i^2(1-\xi_i\eta_i)\left[u_x\cdot DB(u)\tilde{r}_{i,u}\tilde{r}_i+\tilde{r}_i\otimes u_x: D^2B(u)\tilde{r}_i\right]\\
%&+\left(\frac{w_i}{v_i}\right)_xv_i^2(1-\xi_i\eta_i)\left[B(u)\pa_x(\tilde{r}_{i,u}\tilde{r}_i)+\pa_x(\tilde{r}_i)\cdot DB(u)\tilde{r}_i+
%\tilde{r}_i\cdot DB(u)\pa_x(\tilde{r}_i)\right]\\
=O(1)\sum_j(\delta_0^2\La_j^3+\La_j^4+\La_j^5).
\end{align*}
We wish now to estimates the terms in $\sum\limits_{k=1}^{4}\sum\limits_{i}\left(\frac{w_i}{v_i}-\la_i^*\right)\al_{i}^{8,k}$, from the Lemma \ref{lemmeal9.10} we have directly the result for $\sum\limits_{k=2}^{3}\sum\limits_{i}\left(\frac{w_i}{v_i}-\la_i^*\right)\al_{i}^{8,k}$. We have now simply to estimate $\frac{w_i}{v_i}\al_{i}^{8,1}$ and  $\frac{w_i}{v_i}\al_{i}^{8,4}$. Using the formula \eqref{superimpll}, we get
\begin{align*}
&\frac{w_i}{v_i}\al_{i}^{8,1}=O(1)v_{i,x}\biggl(2w_i(1-\chi_i^\e\xi_i\eta_i)+O(1)\xi_i \frac{w_i}{v_i} (\chi')^\e_i\hat{v}_i +O(1)\xi_i \frac{w_i}{v_i} \sum_{l\ne i}\Delta^\e_{i,l}v_l^2\biggl).
%\big(2w_i(1-\xi_i\eta_i)+\xi_i\frac{w_i}{v_i}\sum_{j\ne i}(v_j^2+w_j^2)\eta'(\frac{v_j^2+w_j^2}{v_i})\prod_{l\ne j,i}\eta(\frac{v_l^2+w_l^2}{v_i})\big)v_{i,x}\\
%&\hspace{5cm}\times\left[B(u)\tilde{r}_{i,u}\tilde{r}_i+\tilde{r}_i\cdot DB(u)\tilde{r}_i\right].
\end{align*}
We deduce then using \eqref{6.45} and Lemma \ref{lemme6.5} that
\begin{align*}
\frac{w_i}{v_i}\al_{i}^{8,1}&=O(1)\biggl( |w_i v_{i,x}|(1-\eta_i)\mathbbm{1}_{\{|\frac{w_i}{v_i}|\leq \frac{\delta_1}{2}\}}+|w_i v_{i,x}|\mathbbm{1}_{\{|\frac{w_i}{v_i}|\geq \frac{\delta_1}{2}\}}+O(1)\mathbbm{1}_{\{v_i^{2N}\leq 2\e, \,|\frac{w_i}{v_i}|\leq \frac{\delta_1}{2}\}}|w_i v_{i,x}|\\
&+\xi_i \frac{w_i}{v_i} (\chi')^\e_i\hat{v}_i v_{i,x}+\La_i^1\biggl)=O(1)(\sum_j (\La_j^1+\La_j^4+\La_j^6)+R_{\e,1}^{10},
\end{align*}
with $R_{\e,1}^{10}$ satisfying \eqref{12.11bis4}.
Similarly we have using Lemmas \ref{lemme6.5}, \ref{lemme6.6}, \ref{lemme11.3}, \ref{lemme11.4} and \eqref{6.45}
\begin{align*}
\frac{w_i}{v_i}\al_{i}^{8,4}&=(1-\eta_i\chi_i^\e\xi_i)w_i v_i\left[u_x\cdot DB(u)\tilde{r}_{i,u}\tilde{r}_i+B(u)\pa_x(\tilde{r}_{i,u}\tilde{r}_i)
+\tilde{r}_{i,x}\cdot DB(u)\tilde{r}_i\right]\\
&+(1-\eta_i\chi_i^\e \xi_i)w_i v_i \left[\tilde{r}_{i}\otimes u_x:D^2 B(u)\tilde{r}_i+\tilde{r}_{i}\cdot DB(u)\tilde{r}_{i,x}\right]\\
&=O(1)\sum_j(\La_j^1+\La_j^4+\La_j^6)+R_{\e,1}^{10},
\end{align*}
with $R_{\e,1}^{10}$ satisfying \eqref{12.11bis4}.
It remains now to deal with the terms $\sum\limits_{k=1}^{4}\sum\limits_{i}\big[\left(\frac{w_i}{v_i}-\la_i^*\right)\al_{i}^{8,k}\big]_x$. From the Lemma \ref{lemmeal9.10}, it suffices to study the following terms $\sum_{l=2}^3 \left(\frac{w_i}{v_i}\right)_x\al_{i}^{8,l}$, $\big(\frac{w_i}{v_i}\al_i^{8,1}\big)_x$ and 
$\big(\frac{w_i}{v_i}\al_i^{8,4}\big)_x$.
 First we have
\begin{align*}
&\left(\frac{w_i}{v_i}\right)_x\al_{i}^{8,2}=-\xi^\p_iv_i\bar{v}_i\left(\frac{w_i}{v_i}\right)_x^2\left[B(u)\tilde{r}_{i,u}\tilde{r}_i+\tilde{r}_i\cdot DB(u)\tilde{r}_i\right]=O(1)\delta_0^2\La_i^3,\\
&\left(\frac{w_i}{v_i}\right)_x\al_{i}^{8,3}=-\sum\limits_{j\neq i}\xi_i \left(\frac{w_i}{v_i}\right)_x v_iv_{j,x}\pa_{v_j}\bar{v}_i\left[B(u)\tilde{r}_{i,u}\tilde{r}_i+\tilde{r}_i\cdot DB(u)\tilde{r}_i\right]=O(1)\La_i^1.
\end{align*}
Now a direct computation gives using Lemma \ref{estimimpo1}
\begin{equation}\label{csuperimplll}
\widetilde{E}_i=2-\xi_i\chi_i^\e \eta_i-\xi_i\pa_{v_i}\bar{v}_i=2(1-\xi_i\chi_i^\e\eta_i)+\xi_i O(1)\left(\Delta^\ep_i+\varkappa^\ep_i\rho^\ep_i\right).
\end{equation}
Using \eqref{csuperimplll}, \eqref{6.45}, \eqref{6.54bis}, \eqref{ngtech5}, \eqref{estimate-v-i-xx-1aabis}  and Lemmas \ref{estimimpo1}, \ref{lemme6.5},  \ref{lemme6.6}, \ref{lemme11.3}, \ref{lemme11.4} we get
	\begin{align*}
		&\big(\frac{w_i}{v_i}\al_{i}^{8,1}\big)_x=\big[\widetilde{E}_i%(2-\xi_i\chi_i^\e \eta_i-\xi_i\pa_{v_i}\bar{v}_i)
		(w_{i,x}v_{i,x}+w_i v_{i,xx})-\xi_i\pa_x(\pa_{v_i}\bar{v}_i)w_iv_{i,x}\big]     \left[B\tilde{r}_{i,u}\tilde{r}_i+\tilde{r}_i\cdot DB\tilde{r}_i\right]\\
		&-\big[\xi_i^\p\left(\frac{w_i}{v_i}\right)_x(\chi_i^\e \eta_i+\pa_{v_i}\bar{v}_i)+\xi_i(2N\frac{ v_i^{2N-1}}{\e}(\chi')_i^\e v_{i,x}\eta_i+ \chi_i^\e\pa_x(\eta_i))\big] w_i v_{i,x}\left[B\tilde{r}_{i,u}\tilde{r}_i+\tilde{r}_i\cdot DB\tilde{r}_i\right]\\
		&+\widetilde{E}_iw_i v_{i,x}  \left[u_x\cdot DB\tilde{r}_{i,u}\tilde{r}_i+u_x\otimes \tilde{r}_i:D^2B\tilde{r}_i+B\pa_x(\tilde{r}_{i,u}\tilde{r}_i)+\tilde{r}_{i,x}\cdot DB\tilde{r}_i+\tilde{r}_i\cdot DB\tilde{r}_{i,x}\right]\\
		&=O(1)\sum (\La_j^1+\La_j^4+\La_j^5+\La_j^6+\La_j^{6,1})+\widetilde{R}_{\e,1}^{10},
		\end{align*}
	with $\widetilde{R}_{\e,1}^{10}$ satisfying \eqref{12.11bis4}. Let us deal here with the most delicate term $\widetilde{E}_i w_i v_{i,xx}   \left[B\tilde{r}_{i,u}\tilde{r}_i+\tilde{r}_i\cdot DB\tilde{r}_i\right]$, using Lemmas \ref{lemme6.5}, \ref{lemme6.6} and \eqref{csuperimplll}, \eqref{6.45}, \eqref{estimate-v-i-xx-1aabis}  we obtain
	\begin{align*}
	&\widetilde{E}_i w_i v_{i,xx}   \left[B\tilde{r}_{i,u}\tilde{r}_i+\tilde{r}_i\cdot DB\tilde{r}_i\right]=O(1)\biggl(|\widetilde{E}_i w_i |(|v_{i,x} |+|w_{i,x}|+|v_i|)+\sum_k(\La_k^1+\La_k^4+\La_k^5+\La_k^6)+\widetilde{R}_{\e,1}^{10,1}\biggl)\\
	&=O(1)\biggl((\mathbbm{1}_{\{|\frac{w_i}{v_i}|\geq\frac{\delta_1}{2}\}}+\mathbbm{1}_{\{|\frac{w_i}{v_i}|\leq\frac{\delta_1}{2}\}}(1-\eta_i)+\mathbbm{1}_{\{|\frac{w_i}{v_i}|\leq\frac{\delta_1}{2}, \, v_i^{2N}\leq 2\e\}})
	|w_i |(|v_{i,x} |+|w_{i,x}|+|v_i|)\\
	&+\sum_k(\La_k^1+\La_k^4+\La_k^5+\La_k^6)+\widetilde{R}_{\e,1}^{10,1}\biggl)=O(1)\sum_k(\La_k^1+\La_k^4+\La_k^5+\La_k^6+\La_k^{6,1})+\widetilde{R}_{\e,1}^{10,2},
		\end{align*}
	with $\widetilde{R}_{\e,1}^{10,1}, \widetilde{R}_{\e,1}^{10,2}$ satisfying \eqref{12.11bis4}
		%Let us give some details on how to treat the term $(2v_i-\xi_i\bar{v}_i-\xi_iv_i\pa_{v_i}\bar{v}_i)v_{i,xx}\left[B\tilde{r}_{i,u}\tilde{r}_i+\tilde{r}_i\cdot DB\tilde{r}_i\right]$
		 %From the Lemma \ref{lemme9.7}, we have
         %\begin{align}
	%		&\mu_i v_{i,xx}=-\mu_{i,x}v_{i,x}+\tilde{\la}_{i,x}v_i+(\tilde{\la}_i-\la_i^*)v_{i}+w_{i,x} +O(1)\sum\limits_{j\neq i}\left(\mu_j v_{j,x}-(\tilde{\la}_j-\la_j^*)v_j-w_j\right)_x v_j\xi_j
	%+O(1)\sum_{j\ne i}\sum_k (|v_{j,x}v_j v_k|+|v_j^2 v_k|)\\ 
	%\nonumber\\
	%&+O(1)\sum_{j\ne i}(|v_j|+|v_{j,x}|+|w_{j,x}|)(|v_j|+|w_{j,x}|+|v_{j,x}|+\sum_{k} |v_k||v_j|)\nonumber\\
	%&+O(1)\sum_k( \La_k^1+\La_k^4+\La_k^5)
		%&+\mathcal{O}(1)\left(\sum\limits_{j\neq i}|v_j|^2\right)\\
	%	+O(1)\sum_k|v_kv_{k,x}||1-\xi_k\eta_k|.\nonumber\\
	%	\end{align}
Next we have using \eqref{suyptech1}, \eqref{suyptech2}
, \eqref{6.45},  Lemmas \ref{estimimpo1}, \ref{lemme6.5}, \ref{lemme6.6}, \ref{lemme11.3}, \ref{lemme11.4}
\begin{align*}
&\big(\frac{w_i}{v_i}\al_i^{8,4}\big)_x=-\big[(\pa_x(\eta_i)\xi_i+\xi_i' \left(\frac{w_i}{v_i}\right)_x \eta_i)\chi_i^\e\big]
v_iw_i \left[\tilde{r}_{i}\otimes u_x:D^2 B(u)\tilde{r}_i+\tilde{r}_{i}\cdot DB(u)\tilde{r}_{i,x}\right]\\
&-\big[(\pa_x(\eta_i)\xi_i+\xi_i' \left(\frac{w_i}{v_i}\right)_x \eta_i)\chi_i^\e\big]
v_iw_i \left[u_x\cdot DB(u)\tilde{r}_{i,u}\tilde{r}_i+B(u)\pa_x(\tilde{r}_{i,u}\tilde{r}_i)
+\tilde{r}_{i,x}\cdot DB(u)\tilde{r}_i\right]\\
&-2N\frac{v_i^{2N}}{\e}v_{i,x}(\chi')_i^\e\xi_i\eta_i w_i \left[\tilde{r}_{i}\otimes u_x:D^2 B(u)\tilde{r}_i+\tilde{r}_{i}\cdot DB(u)\tilde{r}_{i,x}\right]\\
&-2N\frac{v_i^{2N}}{\e}v_{i,x}(\chi')_i^\e\xi_i\eta_i w_i \left[u_x\cdot DB(u)\tilde{r}_{i,u}\tilde{r}_i+B(u)\pa_x(\tilde{r}_{i,u}\tilde{r}_i)
+\tilde{r}_{i,x}\cdot DB(u)\tilde{r}_i\right]\\
%&-\big[(\pa_x(\eta_i)\xi_i+\xi_i' \left(\frac{w_i}{v_i}\right)_x \eta_i)\xhi_i^\e\big]
%v_iw_i \left[\tilde{r}_{i}\otimes u_x:D^2 B(u)\tilde{r}_i+\tilde{r}_{i}\cdot DB(u)\tilde{r}_{i,x}\right]\\
%&-\xi_i' \left(\frac{w_i}{v_i}\right)_x \eta_iv_iw_i \left[u_x\cdot DB(u)\tilde{r}_{i,u}\tilde{r}_i+B(u)\pa_x(\tilde{r}_{i,u}\tilde{r}_i)
%+\tilde{r}_{i,x}\cdot DB(u)\tilde{r}_i\right]\\
%&-\xi_i' \left(\frac{w_i}{v_i}\right)_x \eta_iv_iw_i \left[\tilde{r}_{i}\otimes u_x:D^2 B(u)\tilde{r}_i+\tilde{r}_{i}\cdot DB(u)\tilde{r}_{i,x}\right]\\
&+(1-\xi_i \chi_i^\e \eta_i)(v_{i,x}w_i+w_{i,x}v_i)\left[u_x\cdot DB(u)\tilde{r}_{i,u}\tilde{r}_i+B(u)\pa_x(\tilde{r}_{i,u}\tilde{r}_i)
+\tilde{r}_{i,x}\cdot DB(u)\tilde{r}_i\right]\\
&+(1-\xi_i  \chi_i^\e \eta_i)(v_{i,x}w_i+w_{i,x}v_i)\left[\tilde{r}_{i}\otimes u_x:D^2 B(u)\tilde{r}_i+\tilde{r}_{i}\cdot DB(u)\tilde{r}_{i,x}\right]\\
&+(1-\xi_i  \chi_i^\e  \eta_i)v_i w_i\left[u_{xx}\cdot DB(u)\tilde{r}_{i,u}\tilde{r}_i+u_{x}\otimes u_x:D^2 B(u)\tilde{r}_{i,u}\tilde{r}_i+2 u_x\cdot DB(u)\pa_x(\tilde{r}_{i,u}\tilde{r}_i)\right]\\
&+(1-\xi_i  \chi_i^\e \eta_i)v_i w_i\left[B(u)\pa_{xx}(\tilde{r}_{i,u}\tilde{r}_i)+\tilde{r}_{i,xx}\cdot DB(u)\tilde{r}_i+\tilde{r}_{i,x}\otimes u_x: D^2 B(u)\tilde{r}_i\right]\\
&+(1-\xi_i \chi_i^\e  \eta_i)v_i w_i\left[\tilde{r}_{i x}\otimes u_x:D^2 B(u)\tilde{r}_i+\tilde{r}_{i}\otimes u_{xx}:D^2 B(u)\tilde{r}_i+\tilde{r}_{i}\otimes u_x\otimes u_x:D^3 B(u)\tilde{r}_i\right]\\
&+(1-\xi_i  \chi_i^\e \eta_i)v_i w_i\left[2
\tilde{r}_{i}\otimes u_x: D^2 B(u)\tilde{r}_{i,x}+\tilde{r}_{i}\cdot DB(u)\tilde{r}_{i,xx}+2\tilde{r}_{i,x}\cdot DB(u)\tilde{r}_{i,x}\right]\\
&=O(1)\sum_j(\La_j^1+\delta_0^2\La_j^3+\La_j^4+\La_j^5+\La_j^6+\La_j^{6,1})+\widetilde{R}_{\e,1}^{10},
\end{align*}
with $\widetilde{R}_{\e,1}^{10}$ satisfying \eqref{12.11bis4}
It remains now to estimate $K_{15,x}$ and $K_{15,xx}$. First using Lemmas \ref{lemme11.3}, \ref{lemme11.4} we get
\begin{align*}
	K_{15,x}&=\sum\limits_{i}\sum\limits_{j\neq i}((w_{i,x}-\la_i^*v_{i,x})v_j+(w_{i}-\la_i^*v_{i})v_{j,x})\left[\tilde{r}_i\cdot DB(u)\tilde{r}_j+B(u)\tilde{r}_{i,u}\tilde{r}_j\right]\\
	&+\sum\limits_{i}\sum\limits_{j\neq i}(w_{i}-\la_i^*v_{i})v_j\left[\tilde{r}_{i,x}\cdot DB(u)\tilde{r}_j+\tilde{r}_{i}\otimes u_x: D^2 B(u)\tilde{r}_j+\tilde{r}_{i}\cdot DB(u)\tilde{r}_{j,x}\right]\\
	&+\sum\limits_{i}\sum\limits_{j\neq i}(w_{i}-\la_i^*v_{i})v_j\left[
	u_x\cdot DB(u)\tilde{r}_{i,u}\tilde{r}_j+ B(u)\pa_x(\tilde{r}_{i,u}\tilde{r}_j)\right]=O(1)\La_i^1.
	%&=:\sum\limits_{l=1}^{2}\sum\limits_{i}\B^{16,l}_i.
\end{align*}
From the Lemmas \ref{lemme6.5}, \ref{lemme6.6}, \ref{lemme11.6}, we observe that
\begin{align}
		&\pa_{xx}(\tilde{r}_{i,u}\tilde{r}_j)
		=O(1)%\biggl(\sum\limits_{k=1}|v_{k,x}|+\sum_{k,l} |v_k v_l|+\sum_k\rho_k^\e \mathfrak{A}_k|v_k\bar{v}_k\left(\frac{w_k}{v_k}\right)_x|
		%\biggl)\\
		+O(1)\rho_i^\e\mathfrak{A}_i\biggl(%+\sum_k \sum_{l\ne i}|v_k|
		|v_i\left(\frac{w_i}{v_i}\right)_{xx}|
		+  |v_i\left(\frac{w_i}{v_i}\right)^2_x|+|v_{i,x}\left(\frac{w_i}{v_i}\right)_x|\nonumber\\
		&+ \sum_{k\ne i}|\left(\frac{w_i}{v_i}\right)_x||v_kv_{k,x}|+|\left(\frac{w_i}{v_i}\right)_x v_i\left(\frac{w_j}{v_j}\right)_x|\rho_j^{\e}\mathfrak{A}_j+\frac{v_{i,x}^2}{|v_i|}+\sum_{k\ne i}|\frac{v_{k,x}v_kv_{i,x}}{v_i}|\biggl)\nonumber\\
		&+O(1)\rho_j^\e\mathfrak{A}_j \biggl(|v_j \left(\frac{w_j}{v_j}\right)_{xx}|+  |v_j\left(\frac{w_j}{v_j}\right)^2_x|+|v_{j,x}\left(\frac{w_j}{v_j}\right)_x|\nonumber\\
		&+|\left(\frac{w_j}{v_j}\right)_x|\sum_{k\ne j}|v_kv_{k,x}|+\abs{\left(\frac{w_j}{v_j}\right)_x v_j\left(\frac{w_i}{v_i}\right)_x}\rho_i^\e\mathfrak{A}_i+\frac{v_{j,x}^2}{|v_j|}+\sum_{k\ne j}|\frac{v_{k,x}v_kv_{j,x}}{v_j}|
		\biggl).
		%&+\sum\limits_{k\neq i}[\left( v_{k,x}\pa_{v_k}\bar{v}_i\right)\xi_i\tilde{r}_{i,uv}\tilde{r}_j]_x+\sum\limits_{k\neq j}[\left( v_{k,x}\pa_{v_k}\bar{v}_j\right)\xi_j\tilde{r}_{i,u}\tilde{r}_{j,v}]_x.
		\label{increditechg}
	\end{align}
Similarly we have  using \eqref{suyptech2}, \eqref{increditechg}, Lemmas  \ref{lemme11.3}, \ref{lemme11.4}
\begin{align*}
	&K_{15,xx}=\sum\limits_{i}\sum\limits_{j\neq i}((w_{i,xx}-\la_i^*v_{i,xx})v_j+2(w_{i,x}-\la_i^*v_{i,x})v_{j,x})\left[\tilde{r}_i\cdot DB(u)\tilde{r}_j+B(u)\tilde{r}_{i,u}\tilde{r}_j\right]\\
	&+2\sum\limits_{i}\sum\limits_{j\neq i}E_{i,j}^1 \left[\tilde{r}_{i,x}\cdot DB(u)\tilde{r}_j+\tilde{r}_{i}\otimes u_x: D^2 B(u)\tilde{r}_j+\tilde{r}_{i}\cdot DB(u)\tilde{r}_{j,x}\right]\\
	&+2\sum\limits_{i}\sum\limits_{j\neq i}E_{i,j}^1\left[
	u_x\cdot DB(u)\tilde{r}_{i,u}\tilde{r}_j+ B(u)\pa_x(\tilde{r}_{i,u}\tilde{r}_j)\right]\\
	%&+\sum\limits_{i}\sum\limits_{j\neq i}\left((w_{i,x}-\la_i^*v_{i,x})v_j+(w_{i}-\la_i^*v_{i})v_{j,x}\right)\left[\tilde{r}_{i,x}\cdot DB(u)\tilde{r}_j+\tilde{r}_{i}\otimes u_x: D^2 B(u)\tilde{r}_j+\tilde{r}_{i}\cdot DB(u)\tilde{r}_{j,x}\right]\\
	%&+\sum\limits_{i}\sum\limits_{j\neq i}\left((w_{i,x}-\la_i^*v_{i,x})v_j+(w_{i}-\la_i^*v_{i})v_{j,x}\right)\left[
	%u_x\cdot DB(u)\tilde{r}_{i,u}\tilde{r}_j+ B(u)\pa_x(\tilde{r}_{i,u}\tilde{r}_j)\right]\\
	&+\sum\limits_{i}\sum\limits_{j\neq i}E_{i,j}^2\left[\tilde{r}_{i,xx}\cdot DB(u)\tilde{r}_j+
	\tilde{r}_{i,x}\otimes u_x: D^2 B(u)\tilde{r}_j+2\tilde{r}_{i,x}\cdot DB(u)\tilde{r}_{j,x}\right]\\
	&+\sum\limits_{i}\sum\limits_{j\neq i}E_{i,j}^2\left[
	\tilde{r}_{i,x}\otimes u_x: D^2 B(u)\tilde{r}_j+\tilde{r}_{i}\otimes u_{xx}: D^2 B(u)\tilde{r}_j+\tilde{r}_{i}\otimes u_x\otimes u_x : D^3 B(u)\tilde{r}_j\right]\\
	&+\sum\limits_{i}\sum\limits_{j\neq i}E_{i,j}^2\left[2 \tilde{r}_{i}\otimes u_x: D^2 B(u)\tilde{r}_{j,x}+
	\tilde{r}_{i}\cdot DB(u)\tilde{r}_{j,xx}+u_{xx}\cdot DB(u)\tilde{r}_{i,u}\tilde{r}_j
	\right]\\
	&+\sum\limits_{i}\sum\limits_{j\neq i}E_{i,j}^2 \left[
	u_{x}\otimes u_x: D^2B(u)\tilde{r}_{i,u}\tilde{r}_j +2u_x\cdot DB(u)\pa_x(\tilde{r}_{i,u}\tilde{r}_j)+ B(u)\pa_{xx}(\tilde{r}_{i,u}\tilde{r}_j)\right]\\
	&+\sum\limits_{i}\sum\limits_{j\neq i}
	(w_{i}-\la_i^*v_{i})v_{j,xx}\left[\tilde{r}_i\cdot DB(u)\tilde{r}_j+B(u)\tilde{r}_{i,u}\tilde{r}_j\right]
	%&+\sum\limits_{i}\sum\limits_{j\neq i}E_{i,j}^2 \left[ 2u_x\cdot DB(u)\pa_x(\tilde{r}_{i,u}\tilde{r}_j)+ B(u)\pa_{xx}(\tilde{r}_{i,u}\tilde{r}_j)\right]\\
	=O(1)\sum_k \La_k^1.%+\delta_0^2\La_k^3+\La_k^5).
	%&=:\sum\limits_{l=1}^{2}\sum\limits_{i}\B^{16,l}_i.
\end{align*}
with $E_{i,j}^1=(w_{i,x}-\la_i^*v_{i,x})v_j+(w_{i}-\la_i^*v_{i})v_{j,x}$ and $E_{i,j}^2=(w_{i}-\la_i^*v_{i})v_j$.
\end{proof}
\begin{lemma}%\label{dernierK}
We have
	\begin{align}
			&\sum\limits_{l=1}^{5}\sum\limits_{i}\B^{16,1,l}_i+\sum\limits_{k=1}^{2}\sum\limits_{i}\B^{8,1,k} \nonumber\\
			&\hspace{3cm}=O(1)\sum_j(\La_j^1+\La_j^2+\delta_0^2\La_j^3+\La_j^4+\La_j^5+\La_j^6+\La_j^{6,1})+R_{\e,1}^{8,16},\label{K11xa}\\
			&\sum\limits_{l=1}^{5}\sum\limits_{i}\B^{16,1,l}_{i,x}+\sum\limits_{k=1}^{2}\sum\limits_{i}\B_x^{8,1,k}\nonumber \\
			&\hspace{3cm}=O(1)\sum_j(\La_j^1+\La_j^2+\delta_0^2\La_j^3+\La_j^4+\La_j^5+\La_j^6+\La_j^{6,1})+\widetilde{R}_{\e,1}^{8,16}\label{K11xxa}
						\end{align}
with:
	\begin{equation}
	\int_{\hat{t}}^T\int_{\R}(|R_{\e,1}^{8,16}|+|\widetilde{R}_{\e,1}^{8,16}|)dx ds=O(1)\delta_0^2,
	\label{12.11bis5}
	\end{equation}
	for $\e>0$ small enough in terms of $T-\hat{t}$ and $\delta_0$.
\end{lemma}
\begin{proof}
First we observe from Lemmas \ref{lemme6.5}, \ref{lemme6.6}, \eqref{ngtech4}, \eqref{ngtech5},  \eqref{estimate-v-i-xx-1aa}  and since $\psi_{i,j}=O(1)\xi_i\bar{v}_i$ that
\begin{align*}
%\B^{8,1,1}& =\theta_i^{\p\p}v_i\left(\frac{w_i}{v_i}\right)_x^2\sum\limits_{j\neq i}(\mu_j-\mu_i)\psi_{ij} r_{j}=O(1)\delta_0^2\La_i^3,\\
\B^{8,1,1}& =\theta_i[(\mu_i^{-1})_x\mathcal{H}_i+\mu_i^{-1}\mathcal{H}_{i,x}]\sum\limits_{j\neq i}(\mu_j-\mu_i)\psi_{ij}  r_{j}=O(1)\sum_j(\La_j^1+\La_j^4+\La_j^5+\La_j^6)+R_{\e,1}^{8,16,1},\\
\B^{8,1,2}& =\mu_i^{-1}\left(\frac{w_i}{v_i}\right)_x\mathcal{H}_i\sum\limits_{j\neq i}(\mu_j-\mu_i)\psi_{ij}  r_{j} =O(1)\sum_j\La_j^4.
\end{align*}
with $R_{\e,1}^{8,16,1}$ satisfying \eqref{12.11bis5}.
%In the same way, we have using Lemma \ref{lemma:derivative-r-k}
%\begin{align*}
%\B^{9,1,1}_x& =\theta_i'''v_i\left(\frac{w_i}{v_i}\right)_x^3(B-\mu_i Id)\tilde{r}_i+\theta_i'' v_{i,x}\left(\frac{w_i}{v_i}\right)_x^2 (B-\mu_i Id)\tilde{r}_i\\
%&+2\theta_i^{\p\p}v_i\left(\frac{w_i}{v_i}\right)_{xx}\left(\frac{w_i}{v_i}\right)_x(B-\mu_i Id)\tilde{r}_i
%+\theta_i^{\p\p}v_i\left(\frac{w_i}{v_i}\right)_x^2 u_x\cdot D(B-\mu_i Id)\tilde{r}_i\\
%&+\theta_i^{\p\p}v_i\left(\frac{w_i}{v_i}\right)_x^2 (B-\mu_i Id)\pa_x(\tilde{r}_i)\\
%&=O(1)(\delta_0^2\La_i^3+\La_i^5).
%\B^{9,1,2}& =\theta_i(v_{i,x}-\mu_i^{-1}(\tilde{\la}_i-\la_i^*)v_i-w_i)_x\sum\limits_{j\neq i}(\mu_j-\mu_i)\psi_{ij}  r_{j}=O(1)\sum_j(\La_j^1+\La_j^4+\La_j^5+\La_j^6),\\
%\B^{9,1,3}& =\mu_i^{-1}\theta_i^\p\left(\frac{w_i}{v_i}\right)_x(\mu_i v_{i,x}-(\tilde{\la}_i-\la_i^*)v_i-w_i)\sum\limits_{j\neq i}(\mu_j-\mu_i)\psi_{ij}  r_{j} =O(1)\sum_j(\La_j^1+\La_j^4+\La_j^6).
%\end{align*}
We have now using Lemmas \ref{lemme6.5}, \ref{lemme6.6}, \ref{lemme11.3a}, \eqref{6.45}, \eqref{ngtech4}, \eqref{ngtech5}, \eqref{estimate-v-i-xx-1aa}, \eqref{9.15primefin} and the fact that $(B-\mu_i Id)\tilde{r}_i=O(1)\xi_i\bar{v}_i$ (CAS DIFFICILE)
\begin{align*}
&\B^{8,1,1}_x =\big[\theta_i'\left(\frac{w_i}{v_i}\right)_x[\mu_{i,x}^{-1}\mathcal{H}_i+\mu_i^{-1}\mathcal{H}_{i,x}]+\theta_i [\mu_{i,xx}^{-1}\mathcal{H}_i+2\mu_{i,x}^{-1}\mathcal{H}_{i,x}+ \mu_i^{-1}\mathcal{H}_{i,xx}]\big]
(B-\mu_i Id)\tilde{r}_i\\
%&+\theta_i(v_{i,x}-\mu_i^{-1}(\tilde{\la}_i-\la_i^*)v_i-\mu_i^{-1}w_i)_{xx}(B-\mu_i Id)\tilde{r}_i\\
&+\theta_i[(\mu_i^{-1})_x\mathcal{H}_i+\mu_i^{-1}\mathcal{H}_{i,x}] \sum\limits_{j\neq i}\big[u_x\cdot D(\mu_j-\mu_i)\psi_{ij}  r_{j}+(\mu_j-\mu_i)[\pa_x(\psi_{ij} ) r_{j} +\psi_{ij} \sum_k  v_k r_{j,u}\tilde{r}_k]\big].  \\
%&+\theta_i(v_{i,x}-\mu_i^{-1}(\tilde{\la}_i-\la_i^*)v_i-\mu_i^{-1}w_i)_{x} \sum\limits_{j\neq i}(\mu_j-\mu_i)\pa_x(\psi_{ij} ) r_{j} \\
%&+\theta_i(v_{i,x}-\mu_i^{-1}(\tilde{\la}_i-\la_i^*)v_i-\mu_i^{-1}w_i)_{x} \sum\limits_{j\neq i}\sum_k v_k(\mu_j-\mu_i)\psi_{ij}  r_{j,u}\tilde{r}_k \\
%u_x\cdot D(B-\mu_i Id)\tilde{r}_i\\
%&+\theta_i(v_{i,x}-\mu_i^{-1}(\tilde{\la}_i-\la_i^*)v_i-\mu_i^{-1}w_i)_{x}(B-\mu_i Id)\pa_x(\tilde{r}_i)\\
&=O(1)\sum_j(\La_j^1+\La_2^j+\delta_0^2\La_j^3+\La_j^4+\La_j^5+\La_j^6)+R_{\e,1}^{8,16,1},
\end{align*}
with $\widetilde{R}_{\e,1}^{8,16,2}$ satisfying \eqref{12.11bis5}.
%We mention that we have treated the second term $\theta_i(v_{i,x}-\mu_i^{-1}(\tilde{\la}_i-\la_i^*)v_i-\mu_i^{-1}w_i)_{xx}(B-\mu_i Id)\tilde{r}_i$ in the same way than the second term in $\beta^{10,1}_{i,x}$ and we have in particular
%\begin{align*}
%&\theta_i(v_{i,x}-\mu_i^{-1}(\tilde{\la}_i-\la_i^*)v_i-\mu_i^{-1}w_i)_{xx}(B-\mu_i Id)\tilde{r}_i=O(1)\sum_k(\La_k^1+\La_k^2+\delta_0^2\La_k^3+\La_k^4+\La_k^5+\La_k^6).
%\end{align*}
In a similar way we have from Lemmas \ref{lemme6.5}, \ref{lemme6.6}, \ref{lemme11.3a} and \eqref{ngtech5}, \eqref{estimate-v-i-xx-1aa}
\begin{align*}
\B^{8,1,2}_x&=\big[\left(\frac{w_i}{v_i}\right)_x[(\mu_i^{-1})_x\mathcal{H}_i+\mu_i^{-1}\mathcal{H}_{i,x}]+\left(\frac{w_i}{v_i}\right)_{xx}\mu_i^{-1}\mathcal{H}_i\big]
\sum\limits_{j\neq i}(\mu_j-\mu_i)\psi_{ij}  r_{j}\\
%\theta_i''\left(\frac{w_i}{v_i}\right)_x^2\mu_i^{-1}(\mu_i v_{i,x}-(\tilde{\la}_i-\la_i^*)v_i-w_i) \sum\limits_{j\neq i}(\mu_j-\mu_i)\psi_{ij}  r_{j}\\
%&+\theta_i'\left(\frac{w_i}{v_i}\right)_{xx}\mu_i^{-1}(\mu_i v_{i,x}-(\tilde{\la}_i-\la_i^*)v_i-w_i)\sum\limits_{j\neq i}(\mu_j-\mu_i)\psi_{ij}  r_{j}\\
%&+\theta_i'\left(\frac{w_i}{v_i}\right)_{x}(\mu_i^{-1})_x (\mu_i v_{i,x}-(\tilde{\la}_i-\la_i^*)v_i-w_i)\sum\limits_{j\neq i}(\mu_j-\mu_i)\psi_{ij}  r_{j}\\
%&+\theta_i'\left(\frac{w_i}{v_i}\right)_{x}\mu_i^{-1} (\mu_i v_{i,x}-(\tilde{\la}_i-\la_i^*)v_i-w_i)_x\sum\limits_{j\neq i}(\mu_j-\mu_i)\psi_{ij}  r_{j}\\
&+\left(\frac{w_i}{v_i}\right)_x \mu_i^{-1}\mathcal{H}_i \sum\limits_{j\neq i}\big[u_x\cdot D(\mu_j-\mu_i)\psi_{ij}  r_{j}+(\mu_j-\mu_i)(\pa_x(\psi_{ij})  r_{j}+\sum_k v_k \psi_{ij}  r_{j u}\tilde{r}_k)\big]\\
%+\theta_i'\left(\frac{w_i}{v_i}\right)_{x}\mu_i^{-1} (\mu_i v_{i,x}-(\tilde{\la}_i-\la_i^*)v_i-w_i)\sum\limits_{j\neq i}u_x\cdot D(\mu_j-\mu_i)\psi_{ij}  r_{j}\\
%&+\theta_i'\left(\frac{w_i}{v_i}\right)_{x}\mu_i^{-1} (\mu_i v_{i,x}-(\tilde{\la}_i-\la_i^*)v_i-w_i)\sum\limits_{j\neq i}(\mu_j-\mu_i)\pa_x(\psi_{ij})  r_{j}\\
%&+\theta_i'\left(\frac{w_i}{v_i}\right)_{x}\mu_i^{-1} (\mu_i v_{i,x}-(\tilde{\la}_i-\la_i^*)v_i-w_i)\sum\limits_{j\neq i}\sum_k v_k (\mu_j-\mu_i)\psi_{ij}  r_{j u}\tilde{r}_k\\
&=O(1)\sum_j(\La_j^1+ \delta_0^2\La_j^3+\La_j^4+\La_j^5+\La_j^6).
\end{align*}
We recall now from \eqref{def:b-ij} and \eqref{def:hat-b-ij} that 
\begin{align*}
b_{ij}&= \xi_j^\p\bar{v}_j\psi_{ji,v} -\psi_{ji,\si},\,
\hat{b}_{ij}=\psi_{ji}+\left[(\la_i^*-\la_j^*)+\frac{w_j}{v_j}\right]\big[\xi_j^\p\bar{v}_j\psi_{ji,v}-\psi_{ji,\si}\big].
\end{align*}
We have in particular using Lemmas \ref{estimimpo1}, \ref{lemme6.5}, \ref{lemme6.6} and the fact that  $b_{ij}, \hat{b}_{ij}=\bar{v}_j\mathfrak{A}_j$, $\psi_{i,k}, \tilde{r}_{i,\sig}=O(1)\bar{v}_i\mathfrak{A}_i$
\begin{align*}
	\B^{16,1,1}&=-\sum\limits_{j\neq i}\left(w_{j,x}-\frac{w_j}{v_j}v_{j,x}\right)_x(\mu_i-\mu_j)\hat{b}_{ij}    \sum\limits_{k\neq i}\psi_{ik}r_k=O(1)\sum_j(\La_j^1+\La_j^5),\\
	\B^{16,1,2}&=-\sum\limits_{j\neq i}\left(w_{j,x}-\frac{w_j}{v_j}v_{j,x}\right)_x(\mu_i-\mu_j)\hat{b}_{ij}  (\frac{w_i}{v_i}-\la_i^*)\xi_i^\p\bar{v}_i\tilde{r}_{i,v}=O(1)\sum_j(\La_j^1+\La_j^5),\\
	\B^{16,1,3}&=\sum\limits_{j\neq i}\left(w_{j,x}-\frac{w_j}{v_j}v_{j,x}\right)_x(\mu_i-\mu_j)\hat{b}_{ij}  (\frac{w_i}{v_i}-\la_i^*)\theta_i^\p\tilde{r}_{i,\si}=O(1)\sum_j(\La_j^1+\La_j^5),\\
	\B^{16,1,4}&=\sum\limits_{j\neq i}\left(w_{j,x}-\frac{w_j}{v_j}v_{j,x}\right)_x(\mu_i-\mu_j){b}_{ij}   \la_i^* \sum\limits_{k\neq i}\psi_{ik}r_k=O(1)\sum_j(\La_j^1+\La_j^5),\\
	\B^{16,1,5}&=\sum\limits_{j\neq i}\left(w_{j,x}-\frac{w_j}{v_j}v_{j,x}\right)_x(\mu_i-\mu_j){b}_{ij}   (\la_i^*-\frac{w_i}{v_i}) \Big((v_i\xi_i\pa_{v_i}\bar{v}_i-\frac{w_i}{v_i}\xi_i^\p\bar{v}_i) \tilde{r}_{i,v}\\
	&+\frac{w_i}{v_i}\theta^\p_i\tilde{r}_{i,\si}\Big)=O(1)\sum_j(\La_j^1+\La_j^5).
\end{align*}
%\begin{equation}\label{def:b-ij}
%	b_{ij}:= \xi_j^\p\bar{v}_j\psi_{ji,v} -\psi_{ji,\si}.
%\end{equation}
We have now
\begin{align}
&\left(w_{j,x}-\frac{w_j}{v_j}v_{j,x}\right)_{xx}=w_{j,xxx}-\frac{w_j}{v_j}v_{j,xxx}-2\left(\frac{w_j}{v_j}\right)_x v_{j,xx}-v_{j,x} \left(\frac{w_j}{v_j}\right)_{xx}.
\label{calculimp} 
\end{align}
We have in addition from Lemmas \ref{estimimpo1}, \ref{lemme6.5}, \ref{lemme6.6}, \ref{lemme11.3a}, \ref{lemme9.6},  \eqref{calculimp}, \eqref{ngtech5} and the fact that $\hat{b}_{ij}=\bar{v}_j\mathfrak{A}_j$, $\psi_{ik}=O(1)\xi_i\bar{v}_i$
\begin{align*}
	\B^{16,1,1}_x&=-\sum\limits_{j\neq i}\big[\left(w_{j,x}-\frac{w_j}{v_j}v_{j,x}\right)_{xx}(\mu_i-\mu_j)
	+\left(w_{j,x}-\frac{w_j}{v_j}v_{j,x}\right)_{x} u_x\cdot D(\mu_i-\mu_j)\big]
	\hat{b}_{ij}    \sum\limits_{k\neq i}\psi_{ik}r_k\\
	%& -\sum\limits_{j\neq i}\left(w_{j,x}-\frac{w_j}{v_j}v_{j,x}\right)_{x}u_x\cdot D(\mu_i-\mu_j)\hat{b}_{ij}    \sum\limits_{k\neq i}\psi_{ik}r_k\\
	&-\sum\limits_{j\neq i}\left(w_{j,x}-\frac{w_j}{v_j}v_{j,x}\right)_{x}(\mu_i-\mu_j)\big[\hat{b}_{ij,x}    \sum\limits_{k\neq i}\psi_{ik}r_k+\hat{b}_{ij}    \sum\limits_{k\neq i}(\pa_x(\psi_{ik})r_k+\psi_{ik}u_x\cdot Dr_k(u))\big]\\
	%&-\sum\limits_{j\neq i}\left(w_{j,x}-\frac{w_j}{v_j}v_{j,x}\right)_{x}(\mu_i-\mu_j)\hat{b}_{ij}    \sum\limits_{k\neq i}\pa_x(\psi_{ik})r_k\\
	%&-\sum\limits_{j\neq i}\left(w_{j,x}-\frac{w_j}{v_j}v_{j,x}\right)_{x}(\mu_i-\mu_j)\hat{b}_{ij}    \sum\limits_{k\neq i}\psi_{ik}u_x\cdot Dr_k(u),\\
	&=\sum_j (\La_j^1+\La_j^2+\La_4^j+\La_j^5).
	%\B^{17,1,2}&=-\sum\limits_{j\neq i}\left(w_{j,x}-\frac{w_j}{v_j}v_{j,x}\right)_x(\mu_i-\mu_j)\hat{b}_{ij}  (\frac{w_i}{v_i}-\la_i^*)\xi_i^\p\bar{v}_i\tilde{r}_{i,v}=O(1)\sum_j(\La_j^1+\La_j^5),\\
	%\B^{17,1,3}&=\sum\limits_{j\neq i}\left(w_{j,x}-\frac{w_j}{v_j}v_{j,x}\right)_x(\mu_i-\mu_j)\hat{b}_{ij}  (\frac{w_i}{v_i}-\la_i^*)\theta_i^\p\tilde{r}_{i,\si}=O(1)\sum_j(\La_j^1+\La_j^5),\\
	%\B^{17,1,4}&=\sum\limits_{j\neq i}\left(w_{j,x}-\frac{w_j}{v_j}v_{j,x}\right)_x(\mu_i-\mu_j){b}_{ij}   \la_i^* \sum\limits_{k\neq i}\psi_{ik}r_k=O(1)\sum_j(\La_j^1+\La_j^5),\\
	%\B^{17,1,5}&=\sum\limits_{j\neq i}\left(w_{j,x}-\frac{w_j}{v_j}v_{j,x}\right)_x(\mu_i-\mu_j){b}_{ij}   (\la_i^*-\frac{w_i}{v_i}) \Big((v_i\xi_i\pa_{v_i}\bar{v}_i-\frac{w_i}{v_i}\xi_i^\p\bar{v}_i) \tilde{r}_{i,v}\\
	%&+\frac{w_i}{v_i}\theta^\p_i\tilde{r}_{i,\si}\Big)=O(1)\sum_j(\La_j^1+\La_j^5).
\end{align*}
Applying Lemmas \ref{estimimpo1}, \ref{lemme6.5}, \ref{lemme6.6}, \ref{lemme11.3}, \ref{lemme9.6}, \eqref{calculimp}, \eqref{ngtech5} and the fact that $\hat{b}_{ij}=\bar{v}_j\mathfrak{A}_j$, we obtain
\begin{align*}
	\B^{16,1,2}_x &=-\sum\limits_{j\neq i}\left(w_{j,x}-\frac{w_j}{v_j}v_{j,x}\right)_{xx}(\mu_i-\mu_j)	\hat{b}_{ij}  (\frac{w_i}{v_i}-\la_i^*)\xi_i^\p\bar{v}_i\tilde{r}_{i,v}\\
	&-\sum\limits_{j\neq i}\left(w_{j,x}-\frac{w_j}{v_j}v_{j,x}\right)_x\big[u_x\cdot D(\mu_i-\mu_j)\hat{b}_{ij} + (\mu_i-\mu_j)\hat{b}_{ij,x}  \big]  (\frac{w_i}{v_i}-\la_i^*)\xi_i^\p\bar{v}_i\tilde{r}_{i,v}\\
	%&-\sum\limits_{j\neq i}\left(w_{j,x}-\frac{w_j}{v_j}v_{j,x}\right)_x(\mu_i-\mu_j)\hat{b}_{ij,x}  (\frac{w_i}{v_i}-\la_i^*)\xi_i^\p\bar{v}_i\tilde{r}_{i,v}\\
	&-\sum\limits_{j\neq i}\left(w_{j,x}-\frac{w_j}{v_j}v_{j,x}\right)_x(\mu_i-\mu_j)\hat{b}_{ij}  \big[(\left(\frac{w_i}{v_i}\right)_x \xi_i^\p+ (\frac{w_i}{v_i}-\la_i^*) \xi_i''\left(\frac{w_i}{v_i}\right)_x\big]
	\bar{v}_i\tilde{r}_{i,v}\\
	%&-\sum\limits_{j\neq i}\left(w_{j,x}-\frac{w_j}{v_j}v_{j,x}\right)_x(\mu_i-\mu_j)\hat{b}_{ij}  (\frac{w_i}{v_i}-\la_i^*)\xi_i''\left(\frac{w_i}{v_i}\right)_x\bar{v}_i\tilde{r}_{i,v}\\
	&-\sum\limits_{j\neq i}\left(w_{j,x}-\frac{w_j}{v_j}v_{j,x}\right)_x(\mu_i-\mu_j)\hat{b}_{ij}  (\frac{w_i}{v_i}-\la_i^*)\xi_i^\p\big[\pa_x(\bar{v}_i)\tilde{r}_{i,v}+\bar{v}_i\pa_x(\tilde{r}_{i,v})\big]\\
	%&-\sum\limits_{j\neq i}\left(w_{j,x}-\frac{w_j}{v_j}v_{j,x}\right)_x(\mu_i-\mu_j)\hat{b}_{ij}  (\frac{w_i}{v_i}-\la_i^*)\xi_i^\p\bar{v}_i\pa_x(\tilde{r}_{i,v}),\\
	&=O(1)\sum_j(\La_j^1+\La_j^2+\La_j^4+\La_j^5).
\end{align*}
We proceed similarly for $\B^{16,1,3}_x$ and we have using the fact that $\tilde{r}_{i,\sig}, \psi_{i,k}=O(1)\xi_i\bar{v}_i$ and Lemmas \ref{lemme11.3}, \ref{lemme11.3a}
\begin{align*}
	\B^{16,1,3}_x, \B^{16,1,4}_x&%=\sum\limits_{j\neq i}\left(w_{j,x}-\frac{w_j}{v_j}v_{j,x}\right)_{xx}(\mu_i-\mu_j)\hat{b}_{ij}  (\frac{w_i}{v_i}-\la_i^*)\theta_i^\p\tilde{r}_{i,\si}\\
	%&+\sum\limits_{j\neq i}\left(w_{j,x}-\frac{w_j}{v_j}v_{j,x}\right)_{x}u_x\cdot D(\mu_i-\mu_j)\hat{b}_{ij}  (\frac{w_i}{v_i}-\la_i^*)\theta_i^\p\tilde{r}_{i,\si}\\
	%&+\sum\limits_{j\neq i}\left(w_{j,x}-\frac{w_j}{v_j}v_{j,x}\right)_{x}(\mu_i-\mu_j)\hat{b}_{ij,x}  (\frac{w_i}{v_i}-\la_i^*)\theta_i^\p\tilde{r}_{i,\si}\\
	%&+\sum\limits_{j\neq i}\left(w_{j,x}-\frac{w_j}{v_j}v_{j,x}\right)_{x}(\mu_i-\mu_j)\hat{b}_{ij}  \left(\frac{w_i}{v_i}\right)_x\theta_i^\p\tilde{r}_{i,\si}\\
	%&+\sum\limits_{j\neq i}\left(w_{j,x}-\frac{w_j}{v_j}v_{j,x}\right)_{x}(\mu_i-\mu_j)\hat{b}_{ij}  (\frac{w_i}{v_i}-\la_i^*)\theta_i''\left(\frac{w_i}{v_i}\right)_x\tilde{r}_{i,\si}\\
	%&+\sum\limits_{j\neq i}\left(w_{j,x}-\frac{w_j}{v_j}v_{j,x}\right)_{x}(\mu_i-\mu_j)\hat{b}_{ij}  (\frac{w_i}{v_i}-\la_i^*)\theta_i^\p\pa_x(\tilde{r}_{i,\si})\\
	=O(1)\sum_j(\La_j^1+\La_j^2+\La_j^4+\La_j^5).
	\end{align*}
	%In the same way, we get FINIR
	%\begin{align*}
	%\B^{17,1,4}_x&=\sum\limits_{j\neq i}\left(w_{j,x}-\frac{w_j}{v_j}v_{j,x}\right)_{xx}(\mu_i-\mu_j){b}_{ij}   \la_i^* \sum\limits_{k\neq i}\psi_{ik}r_k\\
	%&+\sum\limits_{j\neq i}\left(w_{j,x}-\frac{w_j}{v_j}v_{j,x}\right)_xu_x\cdot D(\mu_i-\mu_j){b}_{ij}   \la_i^* \sum\limits_{k\neq i}\psi_{ik}r_k\\
	%&+\sum\limits_{j\neq i}\left(w_{j,x}-\frac{w_j}{v_j}v_{j,x}\right)_x(\mu_i-\mu_j){b}_{ij,x}   \la_i^* \sum\limits_{k\neq i}\psi_{ik}r_k\\
	%&+\sum\limits_{j\neq i}\left(w_{j,x}-\frac{w_j}{v_j}v_{j,x}\right)_x(\mu_i-\mu_j){b}_{ij}   \la_i^* \sum\limits_{k\neq i}\psi_{ik,x}r_k\\
	%&+\sum\limits_{j\neq i}\left(w_{j,x}-\frac{w_j}{v_j}v_{j,x}\right)_x(\mu_i-\mu_j){b}_{ij}   \la_i^* \sum\limits_{k\neq i}\psi_{ik}u_x\cdot Dr_k(u)\\
	%&=O(1)\sum_j(\La_j^1+\La_j^2+\delta_0^2\La_j^3+\La_j^5).
%\end{align*}
From the Lemmas \ref{estimimpo1}, \ref{lemme6.5}, \ref{lemme6.6}, \ref{lemme11.3}, \ref{lemme9.6}, \eqref{calculimp},\eqref{ngtech5} and the fact that  $\widetilde{r}_{i,\sig},\widetilde{{\cal{M}}}_i=(v_i\xi_i\pa_{v_i}\bar{v}_i-\frac{w_i}{v_i}\xi_i^\p\bar{v}_i) \tilde{r}_{i,v}+\frac{w_i}{v_i}\theta^\p_i\tilde{r}_{i,\si}=O(1)\mathfrak{A}_i \bar{v}_i$, $b_{ij}=O(1)\mathfrak{A}_j\bar{v}_j$
\begin{align*}
	&\B^{16,1,5}_x=\sum\limits_{j\neq i}\left(w_{j,x}-\frac{w_j}{v_j}v_{j,x}\right)_x(\mu_i-\mu_j)\big[{b}_{ij,x}   (\la_i^*-\frac{w_i}{v_i})+{b}_{ij}  \left(\frac{w_i}{v_i}\right)_x\big]\widetilde{{\cal{M}}}_i\\
	% \Big((v_i\xi_i\pa_{v_i}\bar{v}_i-\frac{w_i}{v_i}\xi_i^\p\bar{v}_i) \tilde{r}_{i,v}+\frac{w_i}{v_i}\theta^\p_i\tilde{r}_{i,\si}\Big)\\
	%&+\sum\limits_{j\neq i}\left(w_{j,x}-\frac{w_j}{v_j}v_{j,x}\right)_{xx}(\mu_i-\mu_j){b}_{ij}   (\la_i^*-\frac{w_i}{v_i}) \Big((v_i\xi_i\pa_{v_i}\bar{v}_i-\frac{w_i}{v_i}\xi_i^\p\bar{v}_i) \tilde{r}_{i,v}+\frac{w_i}{v_i}\theta^\p_i\tilde{r}_{i,\si}\Big)\\
	%&+\sum\limits_{j\neq i}\left(w_{j,x}-\frac{w_j}{v_j}v_{j,x}\right)_xu_x\cdot D(\mu_i-\mu_j){b}_{ij}   (\la_i^*-\frac{w_i}{v_i}) \widetilde{{\cal{M}}}_i\\%\Big((v_i\xi_i\pa_{v_i}\bar{v}_i-\frac{w_i}{v_i}\xi_i^\p\bar{v}_i) \tilde{r}_{i,v}+\frac{w_i}{v_i}\theta^\p_i\tilde{r}_{i,\si}\Big)\\
	&+\sum\limits_{j\neq i}\big[\left(w_{j,x}-\frac{w_j}{v_j}v_{j,x}\right)_{xx}(\mu_i-\mu_j)
	+\left(w_{j,x}-\frac{w_j}{v_j}v_{j,x}\right)_xu_x\cdot D(\mu_i-\mu_j)\big]
	{b}_{ij}   (\la_i^*-\frac{w_i}{v_i}) \widetilde{{\cal{M}}}_i\\
	%+\sum\limits_{j\neq i}\left(w_{j,x}-\frac{w_j}{v_j}v_{j,x}\right)_x(\mu_i-\mu_j)\big[{b}_{ij,x}   (\la_i^*-\frac{w_i}{v_i})+{b}_{ij}  \left(\frac{w_i}{v_i}\right)_x\big]
	%\widetilde{{\cal{M}}}_i\\% \Big((v_i\xi_i\pa_{v_i}\bar{v}_i-\frac{w_i}{v_i}\xi_i^\p\bar{v}_i) \tilde{r}_{i,v}+\frac{w_i}{v_i}\theta^\p_i\tilde{r}_{i,\si}\Big)\\
	%&-\sum\limits_{j\neq i}\left(w_{j,x}-\frac{w_j}{v_j}v_{j,x}\right)_x(\mu_i-\mu_j){b}_{ij}  \left(\frac{w_i}{v_i}\right)_x \widetilde{{\cal{M}}}_i\\% \Big((v_i\xi_i\pa_{v_i}\bar{v}_i-\frac{w_i}{v_i}\xi_i^\p\bar{v}_i) \tilde{r}_{i,v}+\frac{w_i}{v_i}\theta^\p_i\tilde{r}_{i,\si}\Big)\\
	&+\sum\limits_{j\neq i}\left(w_{j,x}-\frac{w_j}{v_j}v_{j,x}\right)_x(\mu_i-\mu_j){b}_{ij}   (\la_i^*-\frac{w_i}{v_i}) \big[(v_{i,x}\xi_i +v_i\xi_i'\left(\frac{w_i}{v_i}\right)_x)
	\pa_{v_i}\bar{v}_i+v_i\xi_i\pa_{x}  (\pa_{v_i}\bar{v}_i)\big]\tilde{r}_{i,v}\\
	& -\sum\limits_{j\neq i}\left(w_{j,x}-\frac{w_j}{v_j}v_{j,x}\right)_{x}(\mu_i-\mu_j){b}_{ij}   (\la_i^*-\frac{w_i}{v_i}) \big[( \xi_i^\p+\frac{w_i}{v_i}\xi_i'') \left(\frac{w_i}{v_i}\right)_x \bar{v}_i+\frac{w_i}{v_i} \xi_i' \pa_x(\bar{v}_i)\big] \tilde{r}_{i,v}\\
	&+\sum\limits_{j\neq i}\left(w_{j,x}-\frac{w_j}{v_j}v_{j,x}\right)_{x}(\mu_i-\mu_j){b}_{ij}   (\la_i^*-\frac{w_i}{v_i})\left(\frac{w_i}{v_i}\right)_x  (\theta^\p_i+\frac{w_i}{v_i}\theta''_i)\tilde{r}_{i,\si} \\
	&+\sum\limits_{j\neq i}\left(w_{j,x}-\frac{w_j}{v_j}v_{j,x}\right)_{x}(\mu_i-\mu_j){b}_{ij}   (\la_i^*-\frac{w_i}{v_i}) \big[(v_i\xi_i\pa_{v_i}\bar{v}_i-\frac{w_i}{v_i}\xi_i^\p\bar{v}_i)\pa_x( \tilde{r}_{i,v})+\frac{w_i}{v_i}\theta^\p_i\pa_x(\tilde{r}_{i,\si})\big]\\
	&=O(1)\sum_j(\La_j^1+\La_j^2+\La_j^4+\La_j^5).
\end{align*}
\end{proof}
We conclude now this section by estimating the last terms appearing on the right hand side of \eqref{J2}. First we have from Lemma \ref{lemme11.3}
\begin{align*}
&\sum\limits_{j\neq i} (w_i-\la_i^*v_i)v_j\left[\tilde{r}_j\cdot DA\tilde{r}_i-\tilde{r}_i\cdot DA\tilde{r}_j\right]=O(1)\La_i^1,\\
&\sum\limits_{j\neq i} \big((w_i-\la_i^*v_i)v_j\left[\tilde{r}_j\cdot DA\tilde{r}_i-\tilde{r}_i\cdot DA\tilde{r}_j\right]\big)_x=\sum\limits_{j\neq i} (w_{i,x}-\la_i^*v_{i,x})v_j\left[\tilde{r}_j\cdot DA\tilde{r}_i-\tilde{r}_i\cdot DA\tilde{r}_j\right]\\
&+\sum\limits_{j\neq i} (w_i-\la_i^*v_i)v_{j,x}\left[\tilde{r}_j\cdot DA\tilde{r}_i-\tilde{r}_i\cdot DA\tilde{r}_j\right]+\sum\limits_{j\neq i} (w_i-\la_i^*v_i)v_j\left[\tilde{r}_{j,x}\cdot DA\tilde{r}_i-\tilde{r}_{i,x}\cdot DA\tilde{r}_j\right]\\
&+\sum\limits_{j\neq i} (w_i-\la_i^*v_i)v_j\left[\tilde{r}_j\otimes u_x: D^2A\tilde{r}_i-\tilde{r}_i\otimes u_x:D^2A\tilde{r}_j+\tilde{r}_j\cdot DA\tilde{r}_{i,x}-\tilde{r}_i\cdot DA\tilde{r}_{j,x}\right]\\
&=O(1)\La_i^1.
\end{align*}
Similarly we have from Lemmas \ref{estimimpo1}, \ref{lemme6.5}, \ref{lemme6.6}, \ref{lemme11.3}
\begin{align*}
&-\sum\limits_{j\neq i}v_{i,t}(w_j-\la_j^*v_j)\xi_j(\pa_{v_i}\bar{v}_j)\tilde{r}_{j,v}=O(1)\La_i^1,\\
&-\sum\limits_{j\neq i}\big(v_{i,t}(w_j-\la_j^*v_j)\xi_j(\pa_{v_i}\bar{v}_j)\tilde{r}_{j,v}\big)_x=O(1)\sum_k\La_k^1.
\end{align*}
We have now to finish using Lemma \ref{lemme11.3} and the fact that $\tilde{r}_{i,\sig}=O(1)\xi_i\bar{v}_i$
\begin{align*}
&v_i\left(\frac{w_i}{v_i}\right)_x^2(B-\mu_iI_n)[\xi_i^\p\bar{v}_i\tilde{r}_{i,v}-\theta^\p_i\tilde{r}_{i,\si}]=O(1)\delta_0^2 \La_i^3\\
&\Big(v_i\left(\frac{w_i}{v_i}\right)_x^2 (B-\mu_iI_n)[\xi_i^\p\bar{v}_i\tilde{r}_{i,v}-\theta^\p_i\tilde{r}_{i,\si}]\Big)_x=v_{i,x}\left(\frac{w_i}{v_i}\right)_x^2(B-\mu_iI_n)[\xi_i^\p\bar{v}_i\tilde{r}_{i,v}-\theta^\p_i\tilde{r}_{i,\si}]\\
&+v_{i}\left(\frac{w_i}{v_i}\right)_x\big[ 2\left(\frac{w_i}{v_i}\right)_{xx}(B-\mu_iI_n)+\left(\frac{w_i}{v_i}\right)_x\ u_x\cdot D(B-\mu_iI_n)\big][\xi_i^\p\bar{v}_i\tilde{r}_{i,v}-\theta^\p_i\tilde{r}_{i,\si}]\\
&+v_{i}\left(\frac{w_i}{v_i}\right)_x^2 (B-\mu_iI_n)\big[(\xi_i'' \left(\frac{w_i}{v_i}\right)_x\bar{v}_i+\xi_i' \pa_x(\bar{v}_i))\tilde{r}_{i,v}-\theta''_i\left(\frac{w_i}{v_i}\right)_x\tilde{r}_{i,\si}\big]\\
&+v_{i}\left(\frac{w_i}{v_i}\right)_x^2 (B-\mu_iI_n)\big[\xi_i^\p\bar{v}_i\pa_x(\tilde{r}_{i,v})-\theta^\p_i\pa_x(\tilde{r}_{i,\si})\big]=O(1)(\delta_0\La_i^3+\La_i^5).
\end{align*}

\section{Study of the terms $\Phi_i$ and $\Psi_i$}
\label{Phii}
\subsection{Estimate of $\Phi_i$}
We are now interested in studying the term $\Phi_i$ in \eqref{eqfonda1}. We mention that we have studied in the previous section the behavior of $\phi_{i,x}$, it remains then to study each $\mathcal{T}_l$ since $\Phi_i=\mu_i \phi_{i,x}+\sum_{l=1}^9 \mathcal{T}_l$. 
\begin{lemma}
\label{Phiilemme}
We have\\
	\begin{equation}
			\Phi_i=O(1)\sum_j(\La_j^1+\La_j^2+\delta_0^2\La_j^3+\La_j^4+\La_j^5+\La_j^6+\La_j^{6,1}+\La_j^7+\La_j^8)+R_{\e,2},\label{K11xab}
			\end{equation}
with:
	\begin{equation}
	\int_{\hat{t}}^T\int_{\R}|R_{\e,2}| dx ds=O(1)\delta_0^2,
	\label{12.11bis6}
	\end{equation}
	for $\e>0$ small enough in terms of $T-\hat{t}$ and $\delta_0$.
\end{lemma}
\begin{proof}

Let us start with the term $\mathcal{T}_1$, we observe from \eqref{T1} that
\begin{align*}
&\mathcal{T}_1=\mathcal{T}_{1,1}+\mathcal{T}_{1,2}+\mathcal{T}_{1,3},
\end{align*}
with
\begin{align*}
&\mathcal{T}_{1,1}=\sum\limits_{j\neq i}{\color{black}{a_{ij}}}\mu_{j,xx}v_j\left(\frac{w_j}{v_j}\right)_x+\sum\limits_{j\neq i}{\color{black}{a_{ij}}}(2\mu_{j,x}-\tilde{\la}_{j})\left(w_{j,xx}-\frac{w_j}{v_j}v_{j,xx}\right)-\sum\limits_{j\neq i}2{\color{black}{a_{ij}}}\tilde{\la}_{j,x}v_j\left(\frac{w_j}{v_j}\right)_x,\\
&-\sum\limits_{j\neq i}{\color{black}{a_{ij}}} \frac{v_{j,x}}{v_j}\left(\psi_j-\mathcal{F}_j\right)+\sum\limits_{j\neq i}{\color{black}{a_{ij}}} \frac{w_jv_{j,x}}{v^2_j}\left(\phi_j-\mathcal{E}_j\right) +\sum\limits_{j\neq i}{\color{black}{a_{ij}}}\big(\psi_{j,x}-\frac{w_j}{v_j}\phi_{j,x}\big),\nonumber\\
&\mathcal{T}_{1,2}={\color{black}{-}}\mu_i\sum\limits_{j\neq i}(\mu_i-\mu_j)b_{ij}\frac{v_{j,x}}{v_j}\left(w_{j,xx}-\frac{w_j}{v_j}v_{j,xx}\right){\color{black}{+}}\mu_i\sum\limits_{j\neq i}2(\mu_i-\mu_j)b_{ij}\frac{v^2_{j,x}}{v_j}\left(\frac{w_j}{v_j}\right)_x\nonumber\\
	&{\color{black}{+}}\mu_i\sum\limits_{j\neq i}(\mu_i-\mu_j)_xb_{ij}\left(w_{j,x}-\frac{w_j}{v_j}v_{j,x}\right)_{x}
	-\sum\limits_{j\neq i}{\color{black}{a_{ij}}}(\mu_{j,x}-\tilde{\la}_j)v_{j,x}\left(\frac{w_j}{v_j}\right)_x\nonumber\\
	&-\sum\limits_{j\neq i}{\color{black}{a_{ij}}}\frac{v_{j,x}}{v_j}(\mu_j-\mu_i)\left(w_{j,xx}-\frac{w_j}{v_j}v_{j,xx}\right)-\sum\limits_{j\neq i}\frac{2 {\color{black}{a_{ij}}} \mu_iv_{j,x}^2}{v_j}\left(\frac{w_j}{v_j}\right)_x,
	\nonumber\\
	%\mu_i\sum\limits_{j\neq i}2(\mu_i-\mu_j)b_{ij}\left(\frac{w_j}{v_j}\right)_xv_{j,xx}\nonumber\\
	%&+\mu_i\sum\limits_{j\neq i}(\mu_i-\mu_j)b_{ij}\frac{v_{j,x}}{v_j}\left(w_{j,xx}-\frac{w_j}{v_j}v_{j,xx}\right)-\mu_i\sum\limits_{j\neq i}2(\mu_i-\mu_j)b_{ij}\frac{v^2_{j,x}}{v_j}\left(\frac{w_j}{v_j}\right)_x\nonumber\\
	%&-\mu_i\sum\limits_{j\neq i}(\mu_i-\mu_j)b_{ij,x}\left(w_{j,x}-\frac{w_j}{v_j}v_{j,x}\right)_{x}-\mu_i\sum\limits_{j\neq i}(\mu_i-\mu_j)_xb_{ij}\left(w_{j,x}-\frac{w_j}{v_j}v_{j,x}\right)_{x}\\
	%&+\sum\limits_{j\neq i}\mu_{j,xx}v_j\left(\frac{w_j}{v_j}\right)_x+\sum\limits_{j\neq i}(2\mu_{j,x}-\tilde{\la}_{j})\left(w_{j,xx}-\frac{w_j}{v_j}v_{j,xx}\right)-\sum\limits_{j\neq i}2\tilde{\la}_{j,x}v_j\left(\frac{w_j}{v_j}\right)_x\\
	%&+\sum\limits_{j\neq i}\big(\psi_{j,x}-\frac{w_j}{v_j}\phi_{j,x}+\frac{w_j}{v_j}\mathcal{E}_{j,x}-\mathcal{F}_{j,x}\big)-\sum\limits_{j\neq i}(\mu_{j,x}-\tilde{\la}_j)v_{j,x}\left(\frac{w_j}{v_j}\right)_x\\
	%&-\sum\limits_{j\neq i}\frac{v_{j,x}}{v_j}(\mu_j-\mu_i)\left(w_{j,xx}-\frac{w_j}{v_j}v_{j,xx}\right)-\sum\limits_{j\neq i}\frac{v_{j,x}}{v_j}\left(\psi_j-\mathcal{F}_j\right)+\sum\limits_{j\neq i}\frac{w_jv_{j,x}}{v^2_j}\left(\phi_j-\mathcal{E}_j\right)\\
	%&-\sum\limits_{j\neq i}\frac{2\mu_iv_{j,x}^2}{v_j}\left(\frac{w_j}{v_j}\right)_x+\sum\limits_{j\neq i}2\mu_i\left(\frac{w_j}{v_j}\right)_x v_{j,xx}.
&\mathcal{T}_{1,3}=\mu_i\sum\limits_{j\neq i}(\mu_i-\mu_j)b_{ij,x}\left(w_{j,x}-\frac{w_j}{v_j}v_{j,x}\right)_{x}+\sum\limits_{j\neq i}2{\color{black}{a_{ij}}}\mu_i\left(\frac{w_j}{v_j}\right)_x v_{j,xx}
+\sum\limits_{j\neq i}{\color{black}{a_{ij}}}\big(\frac{w_j}{v_j}\mathcal{E}_{j,x}-\mathcal{F}_{j,x}\big).
\end{align*}
We are going now to estimate each term $\mathcal{T}_{1,i}$ with $i\in\{1,\cdots 3\}$. First of all we are going to estimate $\mathcal{T}_{1,1}$, using Lemmas \ref{lemme6.5}, \ref{lemme6.6}, \ref{lemme11.3} and the fact that $a_{ij}=\mu_i b_{ij}=O(1)\xi_j\bar{v}_j$ we deduce that:
\begin{align}
&\mathcal{T}_{1,1}=O(1)\sum_j(\La_j^4+\La_j^5)+O(1)\sum_j\big(|\mathcal{E}_j|+|\mathcal{F}_{j}|+|\psi_{j}|+|\psi_{j,x}|+|\phi_{j}|+|\phi_{j,x}|\big).\label{T11a}
\end{align}
Using \eqref{ngtech5} and the fact that $b_{i,j}, \hat{b}_{i,j}=O(1)\mathfrak{A}_j \bar{v}_j$
we deduce that
\begin{equation}
\begin{aligned}
	&\mathcal{E}_j=-\sum\limits_{k\neq j}(\mu_j-\mu_k)b_{jk}\left(w_{k,x}-\frac{w_k}{v_k}v_{k,x}\right)_x=\sum_k(\La_k^4+\La_k^5),\\
	& \mathcal{F}_j=-\sum\limits_{k\neq j}(\mu_j-\mu_k)\hat{b}_{jk}\left(w_{k,x}-\frac{w_k}{v_k}v_{k,x}\right)_x=\sum_k(\La_k^4+\La_k^5).
	\end{aligned}
	\label{defEjFja}
\end{equation}
Combining \eqref{T11a}, \eqref{defEjFja}, the estimates of the previous section on $\psi_{j}, \psi_{j,x}, \phi_{j,}, \phi_{j,x}$ we deduce that $\mathcal{T}_{1,1}$ satisfies the estimate \eqref{K11xab}. We have now using Lemmas \ref{lemme6.5}, \ref{lemme6.6}, \eqref{ngtech5} and the fact that $b_{ij},a_{ij}=O(1)\mathfrak{A}_j \bar{v}_j$
\begin{align*}
&\mathcal{T}_{1,2}=O(1)\sum_j(\La_j^1+\La_j^4+\La_j^5).
\end{align*}
Let us estimate now the term $\mathcal{T}_{1,3}$, first we have applying Lemmas \ref{lemme6.5}, \ref{lemme6.6}, \ref{lemme9.6} and \eqref{ngtech5}
\begin{align*}
&\mu_i\sum\limits_{j\neq i}(\mu_i-\mu_j)b_{ij,x}\left(w_{j,x}-\frac{w_j}{v_j}v_{j,x}\right)_{x}=
O(1)\biggl(\sum_j(\La_j^1+\delta_0^2\La_j^3+\La_j^4+\La_j^5)\\
&+\rho_j^\e \mathfrak{A}_j\big[
|v_{j,x}(w_{j,xx}-\frac{w_j}{v_j}v_{j,xx})|+|\left(\frac{w_j}{v_j}\right)_x(v_jw_{j,xx}-w_jv_{j,xx})|+|v_{j,x}^2\left(\frac{w_j}{v_j}\right)_x|\\
&+|v_{j,x}\left(\frac{w_j}{v_j}\right)_x\sum_{k\ne j}v_kv_{k,x}|\big]\biggl)=O(1)\sum_j(\La_j^1+\delta_0^2\La_j^3+\La_j^4+\La_j^5).
\end{align*}
Combining now Lemmas \ref{lemme6.5}, \ref{lemme6.6}, \eqref{ngtech5}, \eqref{estimate-v-i-xx-1aabis} and the fact that $a_{ij}=O(1)\mathfrak{A}_j \bar{v}_j$, we deduce that
\begin{align*}
&\sum\limits_{j\neq i}2{\color{black}{a_{ij}}}\mu_i\left(\frac{w_j}{v_j}\right)_x v_{j,xx}=O(1)\sum_{j}(\La_j^1+\La_j^4+\La_j^5+\La_j^6)+O(1)\rho_j^\e\mathfrak{A}_j|v_j\left(\frac{w_j}{v_j}\right)_x(v_{j,x}+w_{j,x})|\\
&=O(1)\sum_{j}(\La_j^1+\delta_0^2\La_j^3+\La_j^4+\La_j^5+\La_j^6).
\end{align*}
Let us deal now with the last term of $\mathcal{T}_{1,3}$. First we observe that proceeding as for the first term of $\mathcal{T}_{1,3}$, using Lemmas \ref{lemme6.5}, \ref{lemme6.6}, \eqref{calculimp}, \eqref{ngtech5} and the fact that $b_{jk}=O(1)\mathfrak{A}_k\bar{v}_k$, $a_{ij}=O(1)\mathfrak{A}_j\bar{v}_j$ we obtain
\begin{align*}
&{\color{black}{a_{ij}}}\frac{w_j}{v_j} \mathcal{E}_{j,x}=-{\color{black}{a_{ij}}}\frac{w_j}{v_j} \sum\limits_{k\neq j}\big[(\mu_j-\mu_k)_x b_{jk}+(\mu_j-\mu_k) b_{jk,x}\big] \left(w_{k,x}-\frac{w_k}{v_k}v_{k,x}\right)_x\\
&-{\color{black}{a_{ij}}}\frac{w_j}{v_j} \sum\limits_{k\neq j}(\mu_j-\mu_k) b_{jk}\left(w_{k,x}-\frac{w_k}{v_k}v_{k,x}\right)_{xx}=O(1)\sum_j(\La_j^1+\La_j^2+\delta_0^2\La_j^3+\La_j^4+\La_j^5).
\end{align*}
In a similar way, we can deal with the term $-\sum\limits_{j\neq i}{\color{black}{a_{ij}}}\mathcal{F}_{j,x}$, it ensures by collecting the previous estimate that $\mathcal{T}_{1,1}$ satisfies \eqref{K11xab}. Let us deal now with $\mathcal{T}_2$ which is defined in \eqref{defglobT} and we start with treating the term $\sum\limits_{j\neq i}(a_{ij,t}-\mu_i a_{ij,xx})\left(w_{j,x}-\frac{w_j}{v_j}v_{j,x}\right)$. 
First applying Lemmas \ref{lemme6.5}, \ref{lemme6.6} and \eqref{estimate-v-i-xx-1aabis}, \eqref{aijt} we get
\begin{align*}
&\sum\limits_{j\neq i}a_{ij,t}\left(w_{j,x}-\frac{w_j}{v_j}v_{j,x}\right)= O(1)\sum_j(\La_j^1+\delta_0^2\La_j^3+\La_j^4)\\
&+O(1)\sum_{j\ne i} \mathfrak{A}_j\rho_j^\e |v_j \left(\frac{w_j}{v_j}\right)_x|\big( |w_{j,xx}-v_{j,xx}\frac{w_j}{v_j}|+|\tilde{\psi}_j|+|\tilde{\phi}_j|+|v_{j,xx}|+|v_{j,x}|\big)\\
%&+O(1)\sum_{j\ne i}| v_j \left(\frac{w_j}{v_j}\right)_x| \tilde{\al}_j \sum_{k\ne j}(|v_{k,t}v_k|+|w_{k,t}w_k|)\\
&=\sum_j(\La_j^1+\delta_0^2\La_j^3+\La_j^4+\La_j^5)+O(1)\sum_j(|\tilde{\psi}_j|+|\tilde{\phi}_j|)+\sum_j  |v_j \left(\frac{w_j}{v_j}\right)_x v_{j,xx}|\mathfrak{A}_j\rho_j^\e\\
&=\sum_j(\La_j^1+\delta_0^2\La_j^3+\La_j^4+\La_j^5+\La_j^6)+O(1)\sum_j(|\tilde{\psi}_j|+|\tilde{\phi}_j|).
\end{align*}
We deduce now using \eqref{defEjFja} and the results of the previous section on $\psi_j$, $\phi_j$ that the term $\sum\limits_{j\neq i}a_{ij,t}\left(w_{j,x}-\frac{w_j}{v_j}v_{j,x}\right)$ satisfies \eqref{K11xab}. We consider the term $-\sum\limits_{j\neq i}\mu_i a_{ij,xx}\left(w_{j,x}-\frac{w_j}{v_j}v_{j,x}\right)$. Applying now the Lemmas \ref{lemme6.5}, \ref{lemme6.6}, \eqref{ngtech5}, \eqref{aijxx}, \eqref{estimate-v-i-xx-1aabis} and the fact that $w_{j,x}=v_j\left(\frac{w_j}{v_j}\right)_x+v_{j,x}$, we obtain
\begin{align*}
&-\sum\limits_{j\neq i}\mu_i a_{ij,xx}\left(w_{j,x}-\frac{w_j}{v_j}v_{j,x}\right)=O(1)\biggl(\sum_j(\La_j^1+\delta_0^2\La_j^3+\La_j^4)+\rho_j^\e\mathfrak{A}_j
\big(| v_j^2\left(\frac{w_j}{v_j}\right)_{xx} \left(\frac{w_j}{v_j}\right)_x|\nonumber\\
	&+|v_{j,x}v_j \left(\frac{w_j}{v_j}\right)_x|+|v_{j,xx}v_j\left(\frac{w_j}{v_j}\right)_x|+|\left(\frac{w_j}{v_j}\right)_x|\sum\limits_{l\neq j}\abs{v_l v_{l,x}v_{j,x}}+v_{j,x}^2|\left(\frac{w_j}{v_j}\right)_x | \big)\biggl)\nonumber\\
	&=O(1)\sum_j(\La_j^1+\delta_0^2\La_j^3+\La_j^4+\La_j^5+\La_j^6).
	%&+v_j\left(\frac{w_j}{v_j}\right)_{x}^2\sum_{k\ne j}|w_kw_{k,x}|
	%+v_j\left(\frac{w_j}{v_j}\right)_x
	%\sum_{k\ne j}(|v_k v_{k,x}|+|w_{k,x}w_{k}|)+v_j \left(\frac{w_j}{v_j}\right)_x \sum_{l,k\ne l}\abs{w_{l,x}w_{k,x}}\nonumber\\
	%&+|v_j\left(\frac{w_j}{v_j}\right)_x^2| \sum_k| v_k|+v_j \left(\frac{w_j}{v_j}\right)_x\sum\limits_{k}\abs{v_k w_{j,x}}+v_j \left(\frac{w_j}{v_j}\right)_x\sum\limits_{k}\sum\limits_{l\neq j}\abs{v_k}(\abs{v_{l,x}}+\abs{w_{l,x}})\\
	%&+v_j \left(\frac{w_j}{v_j}\right)_x \sum\limits_{k,l}\abs{v_{k,x}v_{l,x}}+v_j\left(\frac{w_j}{v_j}\right)_x\sum\limits_{k}\sum\limits_{l\neq j}\abs{v_{k,x}w_{l,x}}+v_j \left(\frac{w_j}{v_j}\right)_x\sum\limits_{k\neq j,l\neq j}\abs{w_{k,x}w_{l,x}} \biggl).
\end{align*}
By proceeding as for the first term of $\mathcal{T}_{1,3}$, we observe that the last term of $\mathcal{T}_2$ satisfies
\eqref{K11xab} which achieves the estimates for $\mathcal{T}_2$.  $\mathcal{T}_3$ satisfies directly 
\eqref{K11xab} by using \eqref{defEjFja} and the estimate of the previous section on $\phi_i$. In a similar manner, by using the Lemma \ref{lemme9.6}, \eqref{ngtech5} and the fact that $a_{ij}=O(1)\mathfrak{A}_j\bar{v}_j$ we observe that $\mathcal{T}_4$ satisfies
\eqref{K11xab}. Let us deal now with $\mathcal{T}_5$, from Lemmas \ref{lemme6.5}, \ref{lemme6.6}, \eqref{T5n}, \eqref{Ai}, \eqref{6.45} and the fact that $a_{ij}=O(1)\mathfrak{A}_j\bar{v}_j$ we observe that
\begin{align*}
&\mathcal{T}_5=O(1)\big(\sum_j (\La_j^1+\La_j^7+\La_j^8)+{\cal A}_i v_{i,x}\big)=O(1)\big(\sum_j (\La_j^1+\La_j^7+\La_j^8)+v_{i,x}v_i^2(1-\xi_i\chi_i^\e\eta_i)\big)\\
&=O(1)\big(\sum_j (\La_j^1+\La_j^7+\La_j^8)+(v_{i,x}v_i^2(1-\eta_i)\mathbbm{1}_{\{|\frac{w_i}{v_i}|\leq\frac{\delta_1}{2}\}}+v_{i,x}v_i^2 \mathbbm{1}_{\{|\frac{w_i}{v_i}|\geq\frac{\delta_1}{2}\}}+
v_{i,x}v_i^2\mathbbm{1}_{\{v_i^{2N}\leq 2\e\}})\big)\\
&=O(1)\sum_j(\La_j^1+\La_j^4+\La_j^7+\La_j^8)+R_{\e,2}^{5},
\end{align*}
with $R_{\e,2}^5$ satisfying \eqref{12.11bis6} for $\e$ small enough. . Proceeding in a similar way, we deduce from Lemmas \ref{lemme6.5}, \ref{lemme6.6}, \eqref{T6n}, \eqref{6.45} and the fact that $a_{ij}=O(1)\mathfrak{A}_j\bar{v}_j$ that
\begin{align*}
&\mathcal{T}_6=O(1)\big(\sum_j (\La_j^1+\La_j^5+\La_j^7)+v_i^3(1-\xi_i\chi_i^\e\eta_i)\big)
%&=O(1)\big(\sum_j (\La_j^1+\La_j^7+\La_j^8)+(v_{i,x}v_i^2(1-\eta_i)\mathbbm{1}_{\{|\frac{w_i}{v_i}|\leq\frac{\delta_1}{2}\}}+v_{i,x}v_i^2 \mathbbm{1}_{\{|\frac{w_i}{v_i}|\geq\frac{\delta_1}{2}\}}+
%v_{i,x}v_i^2\mathbbm{1}_{\{v_i^{2N}\leq 2\e\}})\big)\\
=O(1)\sum_j(\La_j^1+\La_j^4+\La_j^5+\La_j^7)+R_{\e,2}^{6},
\end{align*}
with $R_{\e,2}^6$ satisfying \eqref{12.11bis6} for $\e$ small enough. Using again Lemmas \ref{lemme6.5}, \ref{lemme6.6}, \ref{lemme9.6}, \eqref{T7n}, \eqref{defEjFja}, the fact that $a_{ij}=O(1)\mathfrak{A}_j\bar{v}_j$ and the estimates on $\phi_i$ of the previous section we get
\begin{align*}
	\mathcal{T}_7&=O(1)\sum_j(\La_j^1+\La_j^2+\delta_0^2\La_j^3+\La_j^4+\La_j^5+\La_j^6+\La_j^{6,1})+O(1)\La_i^8+R_{\e,2}^7,
\end{align*}
with  $R_{\e,2}^7$ satisfying \eqref{12.11bis6} for $\e$ small enough. Similarly applying  Lemmas \ref{lemme6.5}, \ref{lemme6.6}, \ref{lemme9.6}, \eqref{T8n}, \eqref{defEjFja}  and the estimates on $\phi_j$ of the previous section, it yields
\begin{align*}
	\mathcal{T}_8&=O(1)\sum_j(\La_j^1+\La_j^2+\delta_0^2\La_j^3+\La_j^4+\La_j^5+\La_j^6+\La_j^{6,1}+\La_j^7)+R_{\e,2}^8,
\end{align*}
with  $R_{\e,2}^8$ satisfying \eqref{12.11bis6} for $\e$ small enough. Let us deal now with the last term $\mathcal{T}_9$, using Lemmas \ref{lemme.laisig}, \ref{lemme6.5}, \ref{lemme6.6}, \eqref{T9n} ,\eqref{defEjFja}, the fact that $a_{ij},\hat{a}_{ij}=O(1)\mathfrak{A}_j\bar{v}_j$ and the estimates on $\phi_i, \psi_i$ of the previous section, it gives
\begin{align*}
	\mathcal{T}_{9}&=O(1)\sum_j(\La_j^1+\La_j^2+\delta_0^2\La_j^3+\La_j^4+\La_j^5+\La_j^6+\La_j^{6,1})+O(1)\La_i^8+R_{\e,2}^9,
	%(\tilde{\la}_{i,v}\xi_i^\p\bar{v}_i-\tilde{\lambda}_{i,\sig}\theta^\p_i)\left[\left(\frac{w_i}{v_i}\right)_x(z_i-\la_i^*v_{i})-\left(\frac{w_i}{v_i}\right)_t v_i\right]\nonumber\\
	%&=(\tilde{\la}_{i,v}\xi_i^\p\bar{v}_i-\tilde{\lambda}_{i,\sig}\theta^\p_i)\frac{1}{v_i}\left[\frac{w_i}{v_i}(z_{i,x}v_i-z_iv_{i,x})+v_{i,x}\hat{z}_i-v_i\hat{z}_{i,x}+(\tilde{\la}_i-\la_i^*)(v_{i,x}w_i-v_iw_{i,x})\right]\nonumber\\
	%&+(\tilde{\la}_{i,v}\xi_i^\p\bar{v}_i-\tilde{\lambda}_{i,\sig}\theta^\p_i)\left[-\psi_i+\mathcal{F}_i+\frac{w_i}{v_i}(\phi_i-\mathcal{E}_i)+\sum\limits_{j\neq i}a_{ij}\frac{w_{i,x}}{v_i}\left(w_{j,x}-\frac{w_j}{v_j}v_{j,x}\right)\right]\nonumber\\
	%&-(\tilde{\la}_{i,v}\xi_i^\p\bar{v}_i-\tilde{\lambda}_{i,\sig}\theta^\p_i)\sum\limits_{j\neq i}\hat{a}_{ij}\frac{v_{i,x}}{v_i}\left(w_{j,x}-\frac{w_j}{v_j}v_{j,x}\right)\nonumber\\
	%&+ (\tilde{\la}_{i,v}\xi_i^\p\bar{v}_i-\tilde{\lambda}_{i,\sig}\theta^\p_i)\left[\sum\limits_{j\neq i}\hat{a}_{ij}\left(w_{j,x}-\frac{w_j}{v_j}v_{j,x}\right)_x-\frac{w_i}{v_i}\sum\limits_{j\neq i}a_{ij}\left(w_{j,x}-\frac{w_j}{v_j}v_{j,x}\right)_x\right]\nonumber\\
	%&+ (\tilde{\la}_{i,v}\xi_i^\p\bar{v}_i-\tilde{\lambda}_{i,\sig}\theta^\p_i)\left[\sum\limits_{j\neq i}\hat{a}_{ij,x}\left(w_{j,x}-\frac{w_j}{v_j}v_{j,x}\right)-\frac{w_i}{v_i}\sum\limits_{j\neq i}a_{ij,x}\left(w_{j,x}-\frac{w_j}{v_j}v_{j,x}\right)\right].\label{T9n}
\end{align*}
with  $R_{\e,2}^9$ satisfying \eqref{12.11bis6} for $\e$ small enough. It achieves the proof of Lemma \ref{Phiilemme}.
%	,\\
%	\mathcal{T}_9&=
%	\sum\limits_{j\neq i}\tilde{\la}_{i,v}(w_{j,x}z_i-\hat{z}_{j,x}v_{i})\xi_i\pa_{w_j}\bar{v}_i-\sum\limits_{j\neq i}\tilde{\la}_{i,v}v_{i}\xi_i(\pa_{w_j}\bar{v}_i)[\psi_j-\mathcal{F}_j]\\
%	&+\sum\limits_{j\neq i}\tilde{\la}_{i,v}v_{i}\xi_i(\pa_{w_j}\bar{v}_i)\sum\limits_{k\neq j}\hat{a}_{jk}\left(w_{k,xx}-\frac{w_k}{v_k}v_{k,xx}-\left(\frac{w_k}{v_k}\right)_xv_{k,x}\right)\\
%	&+\sum\limits_{j\neq i}\tilde{\la}_{i,v}v_{i}\xi_i(\pa_{w_j}\bar{v}_i)\sum\limits_{k\neq j}\left(w_{k,x}-\frac{w_k}{v_k}v_{k,x}\right)\hat{a}_{jk,x}-\sum\limits_{j\neq i}\tilde{\la}_{i,v}\xi_i\pa_{w_j}\bar{v}_i(\la_i^*-\la_j^*)v_iw_{j,x}.
\end{proof}
\subsection{Estimate of $\Psi_i$}
We are now interested in studying the term $\Psi_i$ in \eqref{eqfonda2}. We mention that we have studied in the previous section the behavior of $\psi_{i,x}$, it remains then to study each $\widehat{\mathcal{T}}_l$ since $\Phi_i=\mu_i \psi_{i,x}+\sum\limits_{l=1}^{9}\widehat{\mathcal{T}}_l$. 
\begin{lemma}
\label{Psiilemme}
We have\\
	\begin{equation}
			\Psi_i=O(1)\sum_j(\La_j^1+\La_j^2+\delta_0^2\La_j^3+\La_j^4+\La_j^5+\La_j^6+\La_j^{6,1}+\La_j^7+\La_j^8)+R_{\e,3},\label{K11xabb}
			\end{equation}
with:
	\begin{equation}
	\int_{\hat{t}}^T\int_{\R}|R_{\e,3}| dx ds=O(1)\delta_0^2,
	\label{12.11bis7}
	\end{equation}
	for $\e>0$ small enough in terms of $T-\hat{t}$ and $\delta_0$.
\end{lemma}
\begin{proof}

Let us start with the term $\widehat{\mathcal{T}}_1$, we observe that
\begin{align*}
\widehat{\mathcal{T}}_1=\widehat{\mathcal{T}}_{1,1}+\widehat{\mathcal{T}}_{1,2}+\widehat{\mathcal{T}}_{1,3},
\end{align*}
with
\begin{align*}
&\widehat{\mathcal{T}}_{1,1}=\sum\limits_{j\neq i}\hat{a}_{ij}\left[\mu_{j,xx}v_j\left(\frac{w_j}{v_j}\right)_x+(2\mu_{j,x}-\tilde{\la}_{j})\left(w_{j,xx}-\frac{w_j}{v_j}v_{j,xx}\right)-2\tilde{\la}_{j,x}v_j\left(\frac{w_j}{v_j}\right)_x\right]\nonumber\\
&+\sum\limits_{j\neq i}\hat{a}_{ij}\left[-\frac{v_{j,x}}{v_j}\left(\psi_j-\mathcal{F}_j\right)+\frac{w_jv_{j,x}}{v^2_j}\left(\phi_j-\mathcal{E}_j\right)+\psi_{j,x}-\frac{w_j}{v_j}\phi_{j,x} \right]\nonumber\\
&\widehat{\mathcal{T}}_{1,2}%&=\sum\limits_{j\neq i}(\mu_j-\mu_i)(\hat{a}_{ij}+\mu_i\hat{b}_{ij})\left(w_{j,xxx}-\frac{w_j}{v_j}v_{j,xxx}\right)\\
	=\mu_i\sum\limits_{j\neq i}2(\mu_i-\mu_j)\hat{b}_{ij}\frac{v^2_{j,x}}{v_j}\left(\frac{w_j}{v_j}\right)_x{\color{black}{+}}\mu_i\sum\limits_{j\neq i}(\mu_i-\mu_j)_x\hat{b}_{ij}\left(w_{j,x}-\frac{w_j}{v_j}v_{j,x}\right)_{x}\nonumber\\
	&+\sum\limits_{j\neq i}\hat{a}_{ij}\left[-(\mu_{j,x}-\tilde{\la}_j)v_{j,x}\left(\frac{w_j}{v_j}\right)_x\right]{\color{black}{-}}\mu_i\sum\limits_{j\neq i}(\mu_i-\mu_j)\hat{b}_{ij}\frac{v_{j,x}}{v_j}\left(w_{j,xx}-\frac{w_j}{v_j}v_{j,xx}\right)\nonumber\\
	&+\sum\limits_{j\neq i}\hat{a}_{ij}\left[-\frac{v_{j,x}}{v_j}(\mu_j-\mu_i)\left(w_{j,xx}-\frac{w_j}{v_j}v_{j,xx}\right)\right]+\sum\limits_{j\neq i}\hat{a}_{ij}\left[-\frac{2\mu_iv_{j,x}^2}{v_j}\left(\frac{w_j}{v_j}\right)_x\right]\nonumber\\
&\widehat{\mathcal{T}}_{1,3}=\mu_i\sum\limits_{j\neq i}(\mu_i-\mu_j)\hat{b}_{ij,x}\left(w_{j,x}-\frac{w_j}{v_j}v_{j,x}\right)_{x}{\color{black}{-}}\mu_i\sum\limits_{j\neq i}2(\mu_i-\mu_j)\hat{b}_{ij}\left(\frac{w_j}{v_j}\right)_xv_{j,xx}
\\
&+2\sum\limits_{j\neq i}\hat{a}_{ij}\mu_i\left(\frac{w_j}{v_j}\right)_x v_{j,xx}+\sum\limits_{j\neq i}\hat{a}_{ij}(\frac{w_j}{v_j}\mathcal{E}_{j,x}-\mathcal{F}_{j,x})
\end{align*}
Using \eqref{defEjFja}, the estimates of the previous section on $\psi_j,\phi_j,\psi_{j,x}, \phi_{j,x}$ and the fact that $\hat{a}_{ij}=O(1)\mathfrak{A}_j\bar{v}_j$ we deduce that
\begin{align*}
&\widehat{\mathcal{T}}_{1,1}=O(1)\sum_j(\La_j^1+\La_j^2+\delta_0^2\La_j^3+\La_j^4+\La_j^5+\La_j^6+\La_j^{6,1})+R_{\e,3}^{1,1},
\end{align*}
with $R_{\e,3}^{1,1}$ satisfying \eqref{12.11bis6} for $\e$ small enough. Using now Lemmas \ref{estimimpo1}, \ref{lemme6.5}, \ref{lemme6.6}, \eqref{ngtech5} and the fact that  $\hat{a}_{ij},\hat{b}_{ij}=O(1)\mathfrak{A}_j\bar{v}_j$ , it yields
\begin{align*}
&\widehat{\mathcal{T}}_{1,2}=O(1)\sum_j(\La_j^1+\La_j^4+\La_j^5+\La_j^6).
\end{align*}
By proceeding as for $\mathcal{T}_{1,3}$ we obtain that \ref{lemme9.6}
$$\widehat{\mathcal{T}}_{1,3}=O(1)\sum_j(\La_j^1+\La_j^2+\delta_0^2\La_j^3+\La_j^4+\La_j^5+\La_j^6).$$
From the estimates on $\psi_i$ of the previous section and \eqref{defEjFja}, we deduce directly that $\widehat{\mathcal{T}}_2$ satisfies \eqref{K11xabb}. Proceeding in a similar manner as for $\mathcal{T}_2$, we get
\begin{align*}
&\widehat{\mathcal{T}}_3=O(1)\sum_j(\La_j^1+\La_j^2+\delta_0^2\La_j^3+\La_j^4+\La_j^5+\La_j^6+\La_j^{6,1})+R_{\e,3}^{3},
\end{align*}
with $R_{\e,3}^{3}$ satisfying \eqref{12.11bis6} for $\e$ small enough. Using again Lemmas \ref{estimimpo1}, \ref{lemme6.5}, \ref{lemme6.6}, \ref{lemme9.6} and \eqref{ngtech5}, we can show that $\widehat{\mathcal{T}}_4$ satisfies \eqref{K11xabb}. From \eqref{That5}, \eqref{Ai}, \eqref{6.45} and Lemmas \ref{lemme6.5}, \ref{lemme6.6} we obtain now
\begin{align*}
&\widehat{\mathcal{T}}_5=O(1)\sum_j(\La_j^1+\La_j^4+\La_j^6+\La_j^7+\La_j^{8})+R_{\e,3}^{5},
\end{align*}
with  $R_{\e,3}^{5}$ satisfying \eqref{12.11bis6} for $\e$ small enough. Similarly from \eqref{That6} we deduce that
\begin{align*}
&\widehat{\mathcal{T}}_6=O(1)\sum_j(\La_j^1+\La_j^4+\La_j^6+\La_j^7)+R_{\e,3}^{6},
\end{align*}
with  $R_{\e,3}^{6}$ satisfying \eqref{12.11bis6} for $\e$ small enough. Using now \eqref{That7}, the estimates of the previous section on $\phi_i$ we deduce that
\begin{align*}
&\widehat{\mathcal{T}}_7=O(1)\sum_j(\La_j^1+\La_j^2+\delta_0^2\La_j^3+\La_j^4+\La_j^5+\La_j^6+\La_j^{6,1} +\La_j^8)+R_{\e,3}^{7},
\end{align*}
with  $R_{\e,3}^{7}$ satisfying \eqref{12.11bis6} for $\e$ small enough. Using now Lemmas \ref{estimimpo1}, \ref{lemme6.5}, \ref{lemme6.6}, \ref{lemme9.6}, \eqref{That8}, \eqref{defEjFja}, the estimates of the previous section on $\phi_j$ we get
\begin{align*}
&\widehat{\mathcal{T}}_8=O(1)\sum_j(\La_j^1+\La_j^2+\delta_0^2\La_j^3+\La_j^4+\La_j^5+\La_j^6+\La_j^{6,1} +\La_j^7)+R_{\e,3}^{8},
\end{align*}
with  $R_{\e,3}^{8}$ satisfying \eqref{12.11bis6} for $\e$ small enough. To finish, from Lemmas \eqref{lemme.laisig}, \ref{estimimpo1},  \ref{lemme6.5}, \ref{lemme6.6}, \ref{lemme9.6}, \eqref{That9}, \eqref{defEjFja}, the estimates of the previous section on $\phi_i,\psi_i$ we obtain
\begin{align*}
&\widehat{\mathcal{T}}_9=O(1)\sum_j(\La_j^1+\La_j^2+\delta_0^2\La_j^3+\La_j^4+\La_j^5+\La_j^6+\La_j^{6,1} +\La_j^8)+R_{\e,3}^{9},
\end{align*}
with  $R_{\e,3}^{9}$ satisfying \eqref{12.11bis6} for $\e$ small enough. It achieves the proof of  Lemma \ref{Psiilemme}.
\end{proof}

\appendix

\section{Proof of Lemma \ref{lemma:transformation}}
\begin{proof} First we observe that for any $i\in\{1,\cdots, n\}$, $X_i$ is a $C^1$ diffeomorphism from $\R$ to $\R$ since $X_i$ has a positive derivate. It enables in particular to ensure that the application $\mathcal{T}$ is invertible.
	By using the definition of $\mathcal{T}(f)_i$ it follows that
	\begin{equation*}
		\int\limits_{\R}\abs{\mathcal{T}(f)_i(y)}^p\,dy=	\int\limits_{\R}\abs{f_i(X_i(y))}^p\,dy=	\int\limits_{\R}\abs{f_i(x)}^p\,\sqrt{d_i(u(X_i^{-1}(x)))}\,dx.
	\end{equation*}
	Therefore, the first inequality in \eqref{inequality-T-1} follows. Similarly, we get the second inequality in \eqref{inequality-T-1}. Now, to see \eqref{inequality-T-2}, we first calculate the following derivatives.
	\begin{equation*}
		\frac{d\mathcal{T}(f)_i}{dy}(y)=f'_i(x)\frac{1}{\sqrt{d_i(u(y))}}\mbox{ where }x=X_i(y).
	\end{equation*}
	Hence, we have
	\begin{align*}
		\int\limits_{\R}\abs{\frac{d\mathcal{T}(f)_i}{dy}(y)}^p\,dy&=\int\limits_{\R}\abs{\frac{1}{d_i(u(y))}}^\frac{p}{2}\abs{f_i'(x)}^p\,dy=	\int\limits_{\R}\abs{\frac{1}{d_i(u(X_i^{-1}(x)))}}^\frac{p-1}{2}\abs{f'_i(x)}^p\,dx.
	\end{align*}
	Therefore, the first inequality in \eqref{inequality-T-2} follows and similarly, one can achieve the second inequality in \eqref{inequality-T-2}. Now, we can also derive
	\begin{equation}\label{lemma-1-cal-1}
		\frac{d^2\mathcal{T}(f)_i}{dy^2}(y)= \frac{1}{d_i(u(y))}\frac{d^2f_i}{dx^2}(x)-\frac{1}{ 2 d_i(u(y))^{\frac{3}{2}}} u_y(y)\cdot Dd_i(u(y))\frac{df_i}{dx}(x)\mbox{ where }x=X_i(y).
	\end{equation}
	By using the inequality $(a+b)^p\leq 2^{p-1}(a^p+b^p)$ for $p\geq 1, a>0,b>0$, we can obtain
	\begin{align*}
		&\int\limits_{\R}\abs{\frac{d^2\mathcal{T}(f)_i}{dy^2}(y)}^p\,dy\leq 2^{p-1}\int\limits_{\R}\abs{\frac{1}{d_i(u(y))}}^p\abs{\frac{d^2f_i}{dx^2}(x)}^p\,dy\\
		&\hspace{2cm}+2^{p-1}\int\limits_{\R}\abs{\frac{1}{ 2 d_i(u(y))^{\frac{3}{2}}} u_y(y)\cdot Dd_i(u(y))\frac{df_i}{dx}(x)}^p\,dy\\
		&\leq 2^{p-1}\int\limits_{\R}[\frac{1}{d_i(u(X_i^{-1}(x)))}]^\frac{2p-1}{2}\abs{\frac{d^2f_i}{dx^2}(x)}^p\,dx\\
		&\hspace{2cm}+\frac{1}{2} \norm{u_x}^p_{L^\f}\norm{Dd_i(u)}^p_{L^\f}\int\limits_{\R}[\frac{1}{d_i(u(X_i^{-1}(x)))}]^{\frac{3p-1}{2}}\abs{\frac{df_i}{dx}(x)}^p\,dx.
	\end{align*}
	Hence, the inequality \eqref{inequality-T-3a} follows using the fact that for $p\geq 1$, $a,b >0$ we have $(a+b)^{\frac{1}{p}}\leq a^{\frac{1}{p}}+b^{\frac{1}{p}}$. Again from \eqref{lemma-1-cal-1} we have
	\begin{align*}
		&\int\limits_{\R}\abs{\frac{d^2f_i}{dx^2}(x)}^p\,dx\leq 2^{p-1}\int\limits_{\R}\abs{d_i(u(y))}^p \abs{\frac{d^2\mathcal{T}(f)_i}{dy^2}(y)}^p\,dx\\
		&\hspace{2cm}+\frac{1}{2}\int\limits_{\R}\abs{\frac{1}{\sqrt{d_i(u(y))}}}^p\abs{u_y(y)\cdot Dd_i(u(y))\frac{df_i}{dx}(x)}^p\,dx\\
		&\leq 2^{p-1}\int\limits_{\R} \abs{d_i(u(y))}^{\frac{2p-1}{2}} \abs{\frac{d^2\mathcal{T}(f)_i}{dy^2}(y)}^p\,dy+\frac{1}{2}c_0^{-\frac{p}{2}}\norm{u_x}^p_{L^\f}\norm{Dd_i(u)}^p_{L^\f}\int\limits_{\R}\abs{\frac{df_i}{dx}(x)}^p\,dx.
	\end{align*}
	Therefore, by using \eqref{inequality-T-2}, we obtain the inequality \eqref{inequality-T-3b}. Deriving \eqref{lemma-1-cal-1} again in $y$, we obtain with $x=X_i(y)$
	\begin{equation}\label{lemma-1-cal-1h}
	\begin{aligned}
		&\frac{d^3\mathcal{T}(f)_i}{dy^3}(y)=-\frac{3}{2} \frac{1}{d_i(u(y))^2} u_y(y)\cdot D d_i(u(y))\frac{d^2f_i}{dx^2}(x)+\frac{1}{d_i(u(y))^{\frac{3}{2}}}\frac{d^3 f_i}{dx^3}(x)\\
		&+\frac{3}{4}\frac{1}{d_i(u(y))^{\frac{5}{2}}} (u_y(y)\cdot Dd_i(u(y)))^2 \frac{d f_i}{dx}(x)
		-\frac{1}{ 2 d_i(u(y))^{\frac{3}{2}}} u_{yy}(y)\cdot Dd_i(u(y))\frac{df_i}{dx}(x)
		\\
		&-\frac{1}{ 2 d_i(u(y))^{\frac{3}{2}}} u_y(y)\otimes u_y(y):D^2 d_i(u(y))\frac{df_i}{dx}(x).
		\end{aligned}
	\end{equation}
	It implies again by convexity that
	\begin{align*}
		&\int\limits_{\R}\abs{\frac{d^3\mathcal{T}(f)_i}{dy^3}(y)}^p\,dy\leq 5^{p-1}\big(\frac{3}{2}\big)^p m_1^p \norm{u_y}^p_{L^\infty}\int\limits_{\R}\abs{\frac{1}{d_i(u(y))}}^{2p}\abs{\frac{d^2f_i}{dx^2}(x)}^p\,dy\\
		&+5^{p-1} \int\limits_{\R}   \frac{1}{d_i(u(y))^{\frac{3p}{2}}}\abs{\frac{d^3 f_i}{dx^3}(x)}^p dy
		+5^{p-1}\big(\frac{3}{4}\big)^p \norm{u_y}_{L^\infty}^{2p}m_1^{2p}\int_{\R} \frac{1}{d_i(u(y))^{\frac{5p}{2}}}\abs{ \frac{d f_i}{dx}(x) }^p dy  \\
		&+\frac{5^{p-1}}{2^p}\norm{u_{yy}}_{L^\infty}^p m_1^p \int_{\R} \frac{1}{ d_i(u(y))^{\frac{3p}{2}}} \abs{\frac{df_i}{dx}(x)}^p dy+\frac{5^{p-1}}{2^p}\norm{u_{y}}_{L^\infty}^{2p} m_1^p \int_{\R} \frac{1}{ d_i(u(y))^{\frac{3p}{2}}} \abs{\frac{df_i}{dx}(x)}^p dy\\
		&\leq 5^{p-1}\big(\frac{3}{2}\big)^p m_1^p \norm{u_x}^p_{L^\infty} c_o^{-\frac{4p-1}{2}}\|\frac{d^2f_i}{dx^2}\|_{L^p}^p+5^{p-1}(c_0)^{-\frac{3p-1}{2}}\|\frac{d^3 f_i}{dx^3}\|_{L^p}^p
		\\
		&+5^{p-1}\big(\frac{3}{4}\big)^p \norm{u_x}_{L^\infty}^{2p}m_1^{2p}(c_0)^{-\frac{5p-1}{2}}\norm{\frac{d f_i}{dx}}_{L^p}^p  +\frac{5^{p-1}}{2^p}\norm{u_{xx}}_{L^\infty}^p m_1^p (c_0)^{-\frac{3p-1}{2}} \norm{\frac{df_i}{dx}}_{L^p}^p \\
		&+ \frac{5^{p-1}}{2^p}\norm{u_{x}}_{L^\infty}^{2p} m_1^p(c_o)^{-\frac{3p-1}{2}} \norm{\frac{df_i}{dx}}_{L^p}^p.
	\end{align*}
	We deduce now directly the estimate \eqref{4.7}.
\end{proof}

%--------------------------------------------------------------------------------------------------------------

\section{Proof Lemma \ref{lemme6.6}}\label{appendix:higher-regularity}
We prove the Lemma \ref{lemme6.6}. An analogous version can also be found in \cite{BB-vv-lim-ann-math}. We prove the result by an explicit computation of the terms.
\begin{proof}[Proof Lemma \ref{lemme6.6}:]
	By taking derivative of the decomposition \eqref{eqn-v-i} with respect to $x$ we obtain
	\begin{align}
		u_{xx}&=\sum\limits_{i}v_{i,x}\tilde{r}_i+\sum\limits_{i,j}v_iv_j\tilde{r}_{i,u}\tilde{r}_j+\sum\limits_{i}v_{i,x}v_i\pa_{v_i}\bar{v}_i\xi_i\tilde{r}_{i,v}+\sum\limits_{i}\xi_i^\p\bar{v}_i\left(w_{i,x}-\frac{w_i}{v_i}v_{i,x}\right)\tilde{r}_{i,v}\nonumber\\
		&-\sum\limits_{i}\theta_i^\p\left(w_{i,x}-\frac{w_i}{v_i}v_{i,x}\right)\tilde{r}_{i,\si}+\sum\limits_i\sum\limits_{j\neq i}v_{j,x}v_i(\pa_{v_j}\bar{v}_i)\xi_i\tilde{r}_{i,v}. \label{eqn:u-xx}
	\end{align}
	Then, taking one more derivative with respect to $x$, it gives %we have
	\begin{align}
		&u_{xxx}
		=\sum\limits_{i}v_{i,xx}\tilde{r}_i+\sum\limits_{i}v_{i,x}\tilde{r}_{i,x}\nonumber\\
		&+\sum\limits_{i,j}(v_{i,x}v_j+v_iv_{j,x})\tilde{r}_{i,u}\tilde{r}_j+\sum\limits_{i,j}v_iv_j\pa_x( \tilde{r}_{i,u}\tilde{r}_j)\nonumber\\
		&+\sum\limits_{i}[v_{i,xx}v_i+v_{i,x}^2]\pa_{v_i}\bar{v}_i\xi_i\tilde{r}_{i,v}+\sum\limits_{i}v_{i,x}v_i\pa_x(\pa_{v_i}\bar{v}_i)\xi_i\tilde{r}_{i,v}\nonumber\\
		&+\sum\limits_{i}v_{i,x}v_i\pa_{v_i}\bar{v}_i\xi_i^\p\left(\frac{w_i}{v_i}\right)_x\tilde{r}_{i,v}+\sum\limits_{i}v_{i,x}v_i\pa_{v_i}\bar{v}_i\xi_i\pa_x(\tilde{r}_{i,v})\nonumber\\
		&+\sum\limits_{i}\xi_i^{\p\p}\left(\frac{w_i}{v_i}\right)_x\bar{v}_i\left(w_{i,x}-\frac{w_i}{v_i}v_{i,x}\right)\tilde{r}_{i,v}+\sum\limits_{i}\xi_i^\p\bar{v}_{i,x}\left(w_{i,x}-\frac{w_i}{v_i}v_{i,x}\right)\tilde{r}_{i,v}\nonumber\\
		&+\sum\limits_{i}\xi_i^\p\bar{v}_i\left(w_{i,xx}-\frac{w_i}{v_i}v_{i,xx}\right)\tilde{r}_{i,v}-\sum\limits_{i}\xi_i^\p\bar{v}_i\left(\frac{w_i}{v_i}\right)_xv_{i,x}\tilde{r}_{i,v}\nonumber\\
		&+\sum\limits_{i}\xi_i^\p\bar{v}_i\left(w_{i,x}-\frac{w_i}{v_i}v_{i,x}\right)\pa_x(\tilde{r}_{i,v})\nonumber\\
		&-\sum\limits_{i}\theta_i^{\p\p}\left(\frac{w_i}{v_i}\right)_x\left(w_{i,x}-\frac{w_i}{v_i}v_{i,x}\right)\tilde{r}_{i,\si}-\sum\limits_{i}\theta_i^\p\left(w_{i,xx}-\frac{w_i}{v_i}v_{i,xx}\right)\tilde{r}_{i,\si}\nonumber\\
		&+\sum\limits_{i}\theta_i^\p\left( \frac{w_i}{v_i}\right)_xv_{i,x}\tilde{r}_{i,\si}-\sum\limits_{i}\theta_i^\p\left(w_{i,x}-\frac{w_i}{v_i}v_{i,x}\right)\pa_x(\tilde{r}_{i,\si})\nonumber\\
		&+\sum\limits_i\sum\limits_{j\neq i}v_{j,xx}v_i(\pa_{v_j}\bar{v}_i)\xi_i\tilde{r}_{i,v}\nonumber\\
		&+\sum\limits_i\sum\limits_{j\neq i}[v_{j,x}v_{i,x}\pa_{v_j}\bar{v}_i+v_{j,x}v_i\pa_x(\pa_{v_j}\bar{v}_i)]\xi_i\tilde{r}_{i,v}\nonumber\\
		&+\sum\limits_i\sum\limits_{j\neq i}v_{j,x}v_i\pa_{v_j}\bar{v}_i\xi_i^\p\left(\frac{w_i}{v_i}\right)_x\tilde{r}_{i,v}+\sum\limits_i\sum\limits_{j\neq i}v_{j,x}v_i\pa_{v_j}\bar{v}_i\xi_i\pa_x(\tilde{r}_{i,v})=\sum_{l=1}^{12}\mathcal{Z}_l.\label{estimuxxx}
	\end{align}
	Note that $\mathcal{Z}_{l}$ for $l=1,3,6,8,10$ has terms involving second order derivatives of $v_{j}$ and $w_{j}$. Therefore, we wish to write
	\begin{align}
		&\mathcal{Z}_{1}+\mathcal{Z}_{3}+\mathcal{Z}_{6}+\mathcal{Z}_{8}+\mathcal{Z}_{10}\nonumber\\
		&=\sum\limits_{i}v_{i,xx}\tilde{r}_i+\sum\limits_{i}v_{i,xx}v_i\pa_{v_i}\bar{v}_i\xi_i\tilde{r}_{i,v}+\sum\limits_{i}\xi_i^\p\bar{v}_i\left(w_{i,xx}-\frac{w_i}{v_i}v_{i,xx}\right)\tilde{r}_{i,v}\nonumber\\
		&-\sum\limits_{i}\theta_i^\p\left(w_{i,xx}-\frac{w_i}{v_i}v_{i,xx}\right)\tilde{r}_{i,\si}+\sum\limits_i\sum\limits_{j\neq i}v_{j,xx}v_i(\pa_{v_j}\bar{v}_i)\xi_i\tilde{r}_{i,v}+\sum\limits_{i}v_{i,x}\tilde{r}_{i,x}+\sum\limits_{i} v_{i,x}^2\pa_{v_i}\bar{v}_i\xi_i\tilde{r}_{i,v}\nonumber\\
		&+\sum\limits_{i}v_{i,x}v_i\pa_x(\pa_{v_i}\bar{v}_i)\xi_i\tilde{r}_{i,v}-\sum\limits_{i}\xi_i^\p\bar{v}_i\left(\frac{w_i}{v_i}\right)_xv_{i,x}\tilde{r}_{i,v}-\sum\limits_{i}\theta_i^{\p\p}\left(\frac{w_i}{v_i}\right)_x\left(w_{i,x}-\frac{w_i}{v_i}v_{i,x}\right)\tilde{r}_{i,\si}\nonumber\\
		&=\sum\limits_{i}v_{i,xx}\bigg[\tilde{r}_i+v_i\pa_{v_i}\bar{v}_i\xi_i\tilde{r}_{i,v}-\xi_i^\p\bar{v}_i \frac{w_i}{v_i}\tilde{r}_{i,v}+ \theta_i^\p\frac{w_i}{v_i} \tilde{r}_{i,\si}+ \sum\limits_{j\neq i} v_j\pa_{v_i}\bar{v}_j\xi_j\tilde{r}_{j,v}\bigg]\nonumber\\
		&+\sum\limits_{i}w_{i,xx}\bigg[\xi_i^\p\bar{v}_i \tilde{r}_{i,v}-\theta_i^\p \tilde{r}_{i,\si}\bigg]+\sum\limits_{i}v_{i,x}\tilde{r}_{i,x}+\sum\limits_{i} v_{i,x}^2\pa_{v_i}\bar{v}_i\xi_i\tilde{r}_{i,v}+\sum\limits_{i}v_{i,x}v_i\pa_x(\pa_{v_i}\bar{v}_i)\xi_i\tilde{r}_{i,v}\nonumber\\
		&-\sum\limits_{i}\xi_i^\p\bar{v}_i\left(\frac{w_i}{v_i}\right)_xv_{i,x}\tilde{r}_{i,v}-\sum\limits_{i}\theta_i^{\p\p}\left(\frac{w_i}{v_i}\right)_x\left(w_{i,x}-\frac{w_i}{v_i}v_{i,x}\right)\tilde{r}_{i,\si}.\label{collectZi}
	\end{align}
	From Lemmas \ref{estimimpo1}, \ref{lemme6.5}, \ref{lemme11.3} and the fact that $\tilde{r}_{i,\sig}=O(1)\xi_i\bar{v}_i$ it follows that
	\begin{align}
		&\sum\limits_{i}v_{i,x}\tilde{r}_{i,x}+\sum\limits_{i} v_{i,x}^2\pa_{v_i}\bar{v}_i\xi_i\tilde{r}_{i,v}+\sum\limits_{i}v_{i,x}v_i\pa_x(\pa_{v_i}\bar{v}_i)\xi_i\tilde{r}_{i,v}-\sum\limits_{i}\xi_i^\p\bar{v}_i\left(\frac{w_i}{v_i}\right)_xv_{i,x}\tilde{r}_{i,v}\nonumber\\
		&-\sum\limits_{i}\theta_i^{\p\p}\left(\frac{w_i}{v_i}\right)_x\left(w_{i,x}-\frac{w_i}{v_i}v_{i,x}\right)\tilde{r}_{i,\si}=O(1)\de_0^2\sum\limits_{i}(\abs{w_{i,x}}+\abs{v_{i,x}}).\label{aestimZi}
	\end{align}
	Similarly using again Lemmas \ref{estimimpo1}, \ref{lemme6.5}, \ref{lemme11.3}, \ref{lemme11.4} and the fact that $\tilde{r}_{i,\sig}=O(1)\xi_i\bar{v}_i$ it gives for $l=2,4,5,7,9,11,12$,
	\begin{align}
		\mathcal{Z}_{l}=O(1)\de_0^2\sum\limits_{i}(\abs{w_{i,x}}+\abs{v_{i,x}}+\abs{w_i}+\abs{v_i}).\label{bestimZi}
	\end{align}
	Consider $r_i^{VV},r_i^{VW}$ defined as in \eqref{defriVVs}, it yields from \eqref{estimuxxx},
	\eqref{collectZi}, \eqref{aestimZi} and \eqref{bestimZi}
	%\begin{align*}
	%	r_i^{VV}&:= \tilde{r}_i+v_i\pa_{v_i}\bar{v}_i\xi_i\tilde{r}_{i,v}-\xi_i^\p\bar{v}_i \frac{w_i}{v_i}\tilde{r}_{i,v}+ \theta_i^\p\frac{w_i}{v_i} \tilde{r}_{i,\si}+ \sum\limits_{j\neq i} v_j\pa_{v_i}\bar{v}_j\xi_j\tilde{r}_{j,v},\\
	%	r^{VW}_{i}&:=\xi_i^\p\bar{v}_i \tilde{r}_{i,v}-\theta_i^\p \tilde{r}_{i,\si}.
	%\end{align*}
	%Hence, it yields,
	\begin{equation}\label{eqn:u-xxx}
		u_{xxx}=\sum\limits_{i}(v_{i,xx}r_i^{VV}+w_{i,xx}r_i^{  VW})+O(1)\de_0^2\sum\limits_{i}(\abs{w_{i,x}}+\abs{v_{i,x}}+\abs{w_i}+\abs{v_i}).
	\end{equation}
	Similarly, by taking derivative of \eqref{eqn-w-i} with respect to $x$ it follows
	\begin{align}
		u_{tx}&=\sum\limits_{i}(w_{i,x}-\la_i^*v_{i,x})\tilde{r}_i+\sum\limits_{i,j}(w_i-\la_i^*v_i)v_j\tilde{r}_{i,u}\tilde{r}_j \nonumber\\
		&+\sum\limits_{i}v_{i,x}(w_i-\la_i^*v_i)\pa_{v_i}\bar{v}_i\xi_i\tilde{r}_{i,v}+\sum\limits_{i}\xi_i^\p\bar{v}_i\left(w_{i,x}-\frac{w_i}{v_i}v_{i,x}\right)\left(\frac{w_i}{v_i}-\la_i^*\right)\tilde{r}_{i,v}\nonumber\\
		&-\sum\limits_{i}\theta_i^\p\left(w_{i,x}-\frac{w_i}{v_i}v_{i,x}\right)\left(\frac{w_i}{v_i}-\la_i^*\right)\tilde{r}_{i,\si} +\sum\limits_i\sum\limits_{j\neq i} v_{j,x}(w_i-\la_i^*v_i)(\pa_{v_j}\bar{v}_i)\xi_i\tilde{r}_{i,v}.\label{eqn:u-tx}
	\end{align}
	Then we have
	\begin{align}
		&u_{txx}=\sum\limits_{i}(w_{i,xx}-\la_i^*v_{i,xx})\tilde{r}_i+\sum\limits_{i}(w_{i,x}-\la_i^*v_{i,x})\tilde{r}_{i,x}\nonumber\\
		&+\sum\limits_{i}((w_{i,x}-\la_i^*v_{i,x})v_j+(w_{i}-\la_i^*v_{i})v_{j,x})\tilde{r}_{i,u}\tilde{r}_j+\sum\limits_{i,j}(w_{i}-\la_i^*v_{i})v_j\pa_x(\tilde{r}_{i,u}\tilde{r}_j)\nonumber\\
		&+\sum\limits_{i}[v_{i,xx}(w_{i}-\la_i^*v_{i})+v_{i,x}(w_{i,x}-\la_i^*v_{i,x})]\pa_{v_i}\bar{v}_i\xi_i\tilde{r}_{i,v}+\sum\limits_{i}v_{i,x}(w_{i}-\la_i^*v_{i})\pa_x(\pa_{v_i}\bar{v}_i)\xi_i\tilde{r}_{i,v}\nonumber\\
		&+\sum\limits_{i}v_{i,x}(w_{i}-\la_i^*v_{i})\pa_{v_i}\bar{v}_i\xi_i^\p\left(\frac{w_i}{v_i}\right)_x\tilde{r}_{i,v}+\sum\limits_{i}v_{i,x}(w_{i}-\la_i^*v_{i})\pa_{v_i}\bar{v}_i\xi_i\pa_x(\tilde{r}_{i,v})\nonumber\\
		&+\sum\limits_{i}\big[\xi_i^{\p\p}\left(\frac{w_i}{v_i}\right)_x\bar{v}_i+\xi_i^\p\bar{v}_{i,x}\big]
		\left(w_{i,x}-\frac{w_i}{v_i}v_{i,x}\right)\left(\frac{w_i}{v_i}-\la_i^*\right)\tilde{r}_{i,v}\nonumber\\%+\sum\limits_{i}\xi_i^\p\bar{v}_{i,x}\left(w_{i,x}-\frac{w_i}{v_i}v_{i,x}\right)\left(\frac{w_i}{v_i}-\la_i^*\right)\tilde{r}_{i,v}\\
		&+\sum\limits_{i}\xi_i^\p\bar{v}_i\left(w_{i,xx}-\frac{w_i}{v_i}v_{i,xx}\right)\left(\frac{w_i}{v_i}-\la_i^*\right)\tilde{r}_{i,v}-\sum\limits_{i}\xi_i^\p\bar{v}_i\left(\frac{w_i}{v_i}\right)_xv_{i,x}\left(\frac{w_i}{v_i}-\la_i^*\right)\tilde{r}_{i,v}\nonumber\\
		&+\sum\limits_{i}\xi_i^\p\bar{v}_i\left(w_{i,x}-\frac{w_i}{v_i}v_{i,x}\right)\big[\left(\frac{w_i}{v_i}-\la_i^*\right)\pa_x(\tilde{r}_{i,v})+\left(\frac{w_i}{v_i}\right)_x\tilde{r}_{i,v}\big]\nonumber\\
		%&+\sum\limits_{i}\xi_i^\p\bar{v}_i\left(w_{i,x}-\frac{w_i}{v_i}v_{i,x}\right)\left(\frac{w_i}{v_i}\right)_x\tilde{r}_{i,v}\\
		&-\sum\limits_{i}\big[\theta_i^{\p\p}\left(\frac{w_i}{v_i}\right)_x\left(w_{i,x}-\frac{w_i}{v_i}v_{i,x}\right)+\theta_i^\p\left(w_{i,xx}-\frac{w_i}{v_i}v_{i,xx}\right)\big]
		\left(\frac{w_i}{v_i}-\la_i^*\right)\tilde{r}_{i,\si}\nonumber\\%-\sum\limits_{i}\theta_i^\p\left(w_{i,xx}-\frac{w_i}{v_i}v_{i,xx}\right)\left(\frac{w_i}{v_i}-\la_i^*\right)\tilde{r}_{i,\si}\\
		&+\sum\limits_{i}\big[\theta_i^\p\left( \frac{w_i}{v_i}\right)_x\big(v_{i,x}\left(\frac{w_i}{v_i}-\la_i^*\right)-v_i\left(\frac{w_i}{v_i}\right)_x\big)
		\tilde{r}_{i,\si}-\theta'_iv_i \left(\frac{w_i}{v_i}\right)_x \left(\frac{w_i}{v_i}-\la_i^*\right)\pa_x(\tilde{r}_{i,\si})\big]\nonumber\\
%-\sum\limits_{i}\theta_i^\p\left(w_{i,x}-\frac{w_i}{v_i}v_{i,x}\right)\left(\frac{w_i}{v_i}\right)_x\tilde{r}_{i,\si}-\sum\limits_{i}\theta_i^\p\left(w_{i,x}-\frac{w_i}{v_i}v_{i,x}\right)\left(\frac{w_i}{v_i}-\la_i^*\right)\pa_x(\tilde{r}_{i,\si})\\
		&+\sum\limits_i\sum\limits_{j\neq i}v_{j,xx} (w_i-\la_i^*v_i)\pa_{v_j}\bar{v}_i \left(\frac{w_i}{v_i}-\la_i^*\right) \xi_i\tilde{r}_{i,v}\nonumber\\
		&+\sum\limits_i\sum\limits_{j\neq i}\big[v_{j,x} (w_{i,x}-\la_i^*v_{i,x})\pa_{v_j}\bar{v}_i+ v_{j,x}(w_i-\la_i^*v_i)\pa_x(\pa_{v_j}\bar{v}_i)\big]
		\left(\frac{w_i}{v_i}-\la_i^*\right)\xi_i \tilde{r}_{i,v}\nonumber\\
		%&+\sum\limits_i\sum\limits_{j\neq i} v_{j,x}(w_i-\la_i^*v_i)\pa_x(\pa_{v_j}\bar{v}_i)\left(\frac{w_i}{v_i}-\la_i^*\right)\xi_i\tilde{r}_{i,v}\\
		&+\sum\limits_i\sum\limits_{j\neq i} v_{j,x} (w_i-\la_i^*v_i)\pa_{v_j}\bar{v}_i \big[\left(\frac{w_i}{v_i}\right)_x \xi_i+\left(\frac{w_i}{v_i}-\la_i^*\right)\xi'_i \left(\frac{w_i}{v_i}\right)_x\big]
		\tilde{r}_{i,v}\nonumber\\
		&+\sum\limits_i\sum\limits_{j\neq i} v_{j,x} (w_i-\la_i^*v_i)\pa_{v_j}\bar{v}_i \left(\frac{w_i}{v_i}-\la_i^*\right)\xi_i\pa_x(\tilde{r}_{i,v})=\sum_{l=1}^{13}\mathcal{X}_l.\label{estimutxx}
		%&+\sum\limits_i\sum\limits_{j\neq i} v_{j,x}(w_i-\la_i^*v_i)\pa_{v_j}\bar{v}_i \left(\frac{w_i}{v_i}-\la_i^*\right)\xi'_i \left(\frac{w_i}{v_i}\right)_x\tilde{r}_{i,v}.
	\end{align}
	Observe that $\mathcal{X}_{l}$ for $l=1,3,6,8,10$ has terms involving second order derivatives of $v_{j}$ and $w_{j}$. Therefore, one can have
	\begin{align}
		&\mathcal{X}_{1}+\mathcal{X}_{3}+\mathcal{X}_{6}+\mathcal{X}_{8}+\mathcal{X}_{10}\nonumber\\
		&=\sum\limits_{i}(w_{i,xx}-\la_i^*v_{i,xx})\tilde{r}_i+\sum\limits_{i}v_{i,xx}(w_i-\la_i^*v_i)\pa_{v_i}\bar{v}_i\xi_i\tilde{r}_{i,v}\nonumber\\
		&+\sum\limits_{i}\xi_i^\p\bar{v}_i\left(w_{i,xx}-\frac{w_i}{v_i}v_{i,xx}\right)\left(\frac{w_i}{v_i}-\la_i^*\right)\tilde{r}_{i,v}-\sum\limits_{i}\theta_i^\p\left(w_{i,xx}-\frac{w_i}{v_i}v_{i,xx}\right)\left(\frac{w_i}{v_i}-\la_i^*\right)\tilde{r}_{i,\si}\nonumber\\
		&+\sum\limits_i\sum\limits_{j\neq i} v_{j,xx}(w_i-\la_i^*v_i)\pa_{v_j}\bar{v}_i \xi_i\left(\frac{w_i}{v_i}-\la_i^*\right)\tilde{r}_{i,v}+\sum\limits_{i}(w_{i,x}-\la_i^*v_{i,x})\tilde{r}_{i,x}\nonumber\\
		&+\sum\limits_{i} v_{i,x}(w_{i,x}-\la_i^*v_{i,x})\pa_{v_i}\bar{v}_i\xi_i\tilde{r}_{i,v}+\sum\limits_{i}v_{i,x}(w_{i}-\la_i^*v_{i})\pa_x(\pa_{v_i}\bar{v}_i)\xi_i\tilde{r}_{i,v}\nonumber\\
		&-\sum\limits_{i}\xi_i^\p\bar{v}_i\left(\frac{w_i}{v_i}\right)_xv_{i,x}\left(\frac{w_i}{v_i}-\la_i^*\right)\tilde{r}_{i,v}-\sum\limits_{i}\theta_i^{\p\p}\left(\frac{w_i}{v_i}\right)_x\left(w_{i,x}-\frac{w_i}{v_i}v_{i,x}\right) \left(\frac{w_i}{v_i}-\la_i^*\right)\tilde{r}_{i,\si}\nonumber\\
		&=\sum\limits_{i}w_{i,xx}\bigg[\tilde{r}_i+\xi_i^\p\bar{v}_i\left(\frac{w_i}{v_i}-\la_i^*\right) \tilde{r}_{i,v}-\theta_i^\p\left(\frac{w_i}{v_i}-\la_i^*\right) \tilde{r}_{i,\si}\bigg]\nonumber\\
		&+\sum\limits_{i}v_{i,xx}\bigg[-\la_i^*\tilde{r}_i+(w_i-\la_i^*v_i)\pa_{v_i}\bar{v}_i\xi_i\tilde{r}_{i,v}-\frac{w_i}{v_i}\left(\frac{w_i}{v_i}-\la_i^*\right) [\xi_i^\p\bar{v}_i\tilde{r}_{i,v}-\theta_i^\p\tilde{r}_{i,\si}]\nonumber\\
		&\hspace{7cm}+\sum\limits_{j\neq i} (w_j-\la_j^*v_j)\pa_{v_i}\bar{v}_j \xi_j\left(\frac{w_j}{v_j}-\la_j^*\right)\tilde{r}_{j,v}\bigg]\nonumber\\
		&+\sum\limits_{i}(w_{i,x}-\la_i^*v_{i,x})\tilde{r}_{i,x}+\sum\limits_{i} v_{i,x}(w_{i,x}-\la_i^*v_{i,x})\pa_{v_i}\bar{v}_i\xi_i\tilde{r}_{i,v}\nonumber\\
		&+\sum\limits_{i}v_{i,x}(w_{i}-\la_i^*v_{i})\pa_x(\pa_{v_i}\bar{v}_i)\xi_i\tilde{r}_{i,v}-\sum\limits_{i}\xi_i^\p\bar{v}_i\left(\frac{w_i}{v_i}\right)_xv_{i,x}\left(\frac{w_i}{v_i}-\la_i^*\right)\tilde{r}_{i,v}\nonumber\\
		&-\sum\limits_{i}\theta_i^{\p\p}\left(\frac{w_i}{v_i}\right)_x\left(w_{i,x}-\frac{w_i}{v_i}v_{i,x}\right) \left(\frac{w_i}{v_i}-\la_i^*\right)\tilde{r}_{i,\si}.\label{collectXi}
	\end{align}
	As before, from Lemmas \ref{estimimpo1}, \ref{lemme6.5}, \ref{lemme11.3} and the fact that $\tilde{r}_{i,\sig}=O(1)\mathfrak{A}_i \bar{v}_i$ it follows that
	\begin{align}
		&\sum\limits_{i}\big[(w_{i,x}-\la_i^*v_{i,x})\tilde{r}_{i,x}+v_{i,x}(w_{i,x}-\la_i^*v_{i,x})\pa_{v_i}\bar{v}_i\xi_i\tilde{r}_{i,v}+v_{i,x}(w_{i}-\la_i^*v_{i})\pa_x(\pa_{v_i}\bar{v}_i)\xi_i\tilde{r}_{i,v}\big]\nonumber\\
		&-\sum\limits_{i}\big[\xi_i^\p\bar{v}_i\left(\frac{w_i}{v_i}\right)_xv_{i,x}\left(\frac{w_i}{v_i}-\la_i^*\right)\tilde{r}_{i,v}+\theta_i^{\p\p}\left(\frac{w_i}{v_i}\right)_x\left(w_{i,x}-\frac{w_i}{v_i}v_{i,x}\right) \left(\frac{w_i}{v_i}-\la_i^*\right)\tilde{r}_{i,\si}\big]\nonumber\\
		&=O(1)\de_0^2\sum\limits_{i}(\abs{w_{i,x}}+\abs{v_{i,x}}).\label{aestimXi}
	\end{align}
	Similarly, it gives for $l=2,4,5,7,9,11,12,13$ using Lemmas \ref{estimimpo1}, \ref{lemme6.5}, \ref{lemme11.3}, \ref{lemme11.4}  and the fact that $\tilde{r}_{i,\sig}=\mathfrak{A}_i\bar{v}_i$,
	\begin{align}
		\mathcal{X}_{l}= O(1)\de_0^2\sum\limits_{i}(\abs{w_{i,x}}+\abs{v_{i,x}}+\abs{w_i}+\abs{v_i}).
		\label{bestimXi}
	\end{align}
	%Consider $r_i^{WW},r_i^{WV}$ defined as follows
	%\begin{align*}
	%	r_i^{WW}&:=\tilde{r}_i+\xi_i^\p\bar{v}_i\left(\frac{w_i}{v_i}-\la_i^*\right) \tilde{r}_{i,v}-\theta_i^\p\left(\frac{w_i}{v_i}-\la_i^*\right) \tilde{r}_{i,\si},\\
	%	r^{WV}_{i}&:=-\la_i^*\tilde{r}_i+(w_i-\la_i^*v_i)\pa_{v_i}\bar{v}_i\xi_i\tilde{r}_{i,v}-\frac{w_i}{v_i}\left(\frac{w_i}{v_i}-\la_i^*\right) [\xi_i^\p\bar{v}_i\tilde{r}_{i,v}-\theta_i^\p\tilde{r}_{i,\si}]\\
	%	&+\sum\limits_{j\neq i} (w_j-\la_j^*v_j)\pa_{v_i}\bar{v}_j \xi_j\left(\frac{w_j}{v_j}-\la_j^*\right)\tilde{r}_{j,v}.
	%\end{align*}
	Consider $r_i^{WW},r_i^{WV}$ defined as in \eqref{defriVVs}, it yields from \eqref{estimutxx},
	\eqref{collectXi}, \eqref{aestimXi} and \eqref{bestimXi}
	Hence, it yields,
	\begin{equation}\label{eqn:u-txx}
		u_{txx}=\sum\limits_{i}w_{i,xx}r_i^{WW}+v_{i,xx}r_i^{WV}+O(1)\de_0^2\sum\limits_{i}(\abs{w_{i,x}}+\abs{v_{i,x}}+\abs{w_i}+\abs{v_i}).
	\end{equation}
	From \eqref{eqn:u-xxx} and \eqref{eqn:u-txx} we have
	\begin{equation}
		\begin{pmatrix}
			u_{xxx}\\
			u_{txx}
		\end{pmatrix}=\begin{pmatrix}
			\mathcal{Q}^{VV}&	\mathcal{Q}^{VW}\\
			\mathcal{Q}^{WV}&	\mathcal{Q}^{WW}
		\end{pmatrix}\begin{pmatrix}
			v_{1,xx}\\
			\vdots\\
			v_{n,xx}\\
			w_{1,xx}\\
			\vdots\\
			w_{n,xx}
		\end{pmatrix}+O(1)\de_0^2\sum\limits_{i}(\abs{w_{i,x}}+\abs{v_{i,x}}+\abs{w_i}+\abs{v_i}),
	\end{equation}
	where the $n\times n$ matrices $\mathcal{Q}^{VV}, \mathcal{Q}^{VW},\mathcal{Q}^{WW}$ and $\mathcal{Q}^{WV}$ are defined as follows
	\begin{align*}
		&\mathcal{Q}^{VV}:=\left[r_1^{VV}\cdots r_n^{VV}\right],\quad\quad\hspace{0,5cm} \mathcal{Q}^{VW}:=\left[r_1^{VW}\cdots r_n^{VW}\right],\\
		&\mathcal{Q}^{WW}:=\left[r_1^{WW}\cdots r_n^{WW}\right],\quad\quad \mathcal{Q}^{WV}	:=\left[r_1^{WV}\cdots r_n^{WV}\right].
	\end{align*}
	Since $(\widetilde{r}_i)_{1\leq i \leq n}$ is a basis of $\R^n$, we observe that $\mathcal{Q}^{VV}$ and $\mathcal{Q}^{WW}$ are invertible matrices due to \eqref{defriVVs}, Lemma \ref{lemme6.5} and the fact that $\widetilde{r}_{i,\sig}=O(1)\mathfrak{A}_i\bar{v}_i$. 
	Furthermore similarly we remark that  $\mathcal{Q}^{VW}=O(1)\de_0$. It implies that the matrix
	$\begin{pmatrix}
			\mathcal{Q}^{VV}&	\mathcal{Q}^{VW}\\
			\mathcal{Q}^{WV}&	\mathcal{Q}^{WW}
		\end{pmatrix}$ is invertible and in addition we have
		$\begin{pmatrix}
			\mathcal{Q}^{VV}&	\mathcal{Q}^{VW}\\
			\mathcal{Q}^{WV}&	\mathcal{Q}^{WW}
		\end{pmatrix}^{-1}=O(1).$ We deduce then that for any $i\in\{1,\cdot,n\}$ we have using  the equation satisfies by $u_{t,xx}$ and for $t\in[\hat{t},T]$
		\begin{align*}
		&v_{i,xx}(t,\cdot), w_{i,xx}(t,\cdot)=O(1)(|u_{xxx}|+|u_{t,xx}|)(t,\cdot)+O(1)\de_0^2\sum\limits_{i}(\abs{w_{i,x}}+\abs{v_{i,x}}+\abs{w_i}+\abs{v_i})(t,\cdot),\\
		&=O(1)\big(|u_x| |u_{xx}|+|u_x|^3+|u_{xxx}|+|u_x|  |u_{xxx}|+|u_x|^2 |u_{xx}|+|u_{xx}|^2+|u_x|^4+|u_{xxxx}|\big)(t,\cdot)\\
		&+O(1)\de_0^2\sum\limits_{i}(\abs{w_{i,x}}+\abs{v_{i,x}}+\abs{w_i}+\abs{v_i})(t,\cdot).
\end{align*}
Using now Lemma \ref{lemme6.5} and Corollary \ref{coro2.2} we deduce the estimates in \eqref{estimate-v-w-xx}.\\
We wish now prove the estimates in \eqref{estimate-v-w-tx}. Taking derivative of \eqref{eqn:u-xx} with respect to $t$, we have
	\begin{align*}
		&u_{xxt}=\sum\limits_{i}\big[v_{i,tx}\tilde{r}_i+v_{i,x}\tilde{r}_{i,t}\big]+\sum\limits_{i,j}\big[(v_{i,t}v_j+v_iv_{j,t})\tilde{r}_{i,u}\tilde{r}_j+v_iv_j\pa_t(\tilde{r}_{i,u}\tilde{r}_j)\big]\\
		&+\sum\limits_{i}[v_{i,tx}v_i+v_{i,x}v_{i,t}]\pa_{v_i}\bar{v}_i\xi_i\tilde{r}_{i,v}+\sum\limits_{i}v_{i,x}v_i\big[\pa_t(\pa_{v_i}\bar{v}_i)\xi_i+\pa_{v_i}\bar{v}_i\xi_i^\p\left(\frac{w_i}{v_i}\right)_t\big]\tilde{r}_{i,v}
	\\
		&+\sum\limits_{i}v_{i,x}v_i\pa_{v_i}\bar{v}_i\xi_i\pa_t(\tilde{r}_{i,v})+\sum\limits_{i}\big[\xi_i^{\p\p}\left(\frac{w_i}{v_i}\right)_t\bar{v}_i+\xi_i^\p\bar{v}_{i,t}\big]
		\left(w_{i,x}-\frac{w_i}{v_i}v_{i,x}\right)\tilde{r}_{i,v}\\%+\sum\limits_{i}\xi_i^\p\bar{v}_{i,t}\left(w_{i,x}-\frac{w_i}{v_i}v_{i,x}\right)\tilde{r}_{i,v}\\
		&+\sum\limits_{i}\xi_i^\p\bar{v}_i\left(w_{i,tx}-\frac{w_i}{v_i}v_{i,tx}\right)\tilde{r}_{i,v}-\sum\limits_{i}\xi_i^\p\bar{v}_i\left(\frac{w_i}{v_i}\right)_tv_{i,x}\tilde{r}_{i,v}-\sum\limits_{i}\xi_i^\p\bar{v}_iv_i \left(\frac{w_i}{v_i}\right)_x \pa_t(\tilde{r}_{i,v})\\
		&-\sum\limits_{i}\theta_i^{\p\p}\left(\frac{w_i}{v_i}\right)_t\left(w_{i,x}-\frac{w_i}{v_i}v_{i,x}\right)\tilde{r}_{i,\si}-\sum\limits_{i}\theta_i^\p\left(w_{i,tx}-\frac{w_i}{v_i}v_{i,tx}\right)\tilde{r}_{i,\si}\\
		&+\sum\limits_{i}\theta_i^\p\left( \frac{w_i}{v_i}\right)_tv_{i,x}\tilde{r}_{i,\si}-\sum\limits_{i}\theta_i^\p\left(w_{i,x}-\frac{w_i}{v_i}v_{i,x}\right)\pa_t(\tilde{r}_{i,\si})+\sum\limits_i\sum\limits_{j\neq i} v_{j,tx}v_i\pa_{v_j}\bar{v}_i \xi_i\tilde{r}_{i,v}\\
		&+\sum\limits_i\sum\limits_{j\neq i}  v_{j,x}\big[\big(v_{i,t}\pa_{v_j}\bar{v}_i +v_i\pa_t(\pa_{v_j}\bar{v}_i)\big)\xi_i\tilde{r}_{i,v}+v_i\pa_{v_j}\bar{v}_i\big( \xi_i^\p\left(\frac{w_i}{v_i}\right)_t\tilde{r}_{i,v}+\xi_i\pa_t(\tilde{r}_{i,v})\big)\big].
		%\\
		%&+\sum\limits_i\sum\limits_{j\neq i}  v_{j,x}v_i\pa_t(\pa_{v_j}\bar{v}_i) \xi_i\tilde{r}_{i,v}\\
		%&+\sum\limits_i\sum\limits_{j\neq i}  v_{j,x}v_i\pa_{v_j}\bar{v}_i \xi_i^\p\left(\frac{w_i}{v_i}\right)_t\tilde{r}_{i,v}\\
		%&+\sum\limits_i\sum\limits_{j\neq i} v_{j,x}v_i(\pa_{v_j}\bar{v}_i)\xi_i\pa_t(\tilde{r}_{i,v}).
	\end{align*}
	We rearrange the terms to get that
	\begin{align*}
		u_{xxt}&=\sum\limits_{i}v_{i,tx}\tilde{r}_i+\sum\limits_{i} v_{i,tx}v_i \pa_{v_i}\bar{v}_i\xi_i\tilde{r}_{i,v}+\sum\limits_{i}\xi_i^\p\bar{v}_i\left(w_{i,tx}-\frac{w_i}{v_i}v_{i,tx}\right)\tilde{r}_{i,v}\\
		&-\sum\limits_{i}\theta_i^\p\left(w_{i,tx}-\frac{w_i}{v_i}v_{i,tx}\right)\tilde{r}_{i,\si}+\sum\limits_i\sum\limits_{j\neq i}\left[v_{j,tx}v_i\pa_{v_j}\bar{v}_i\right]\xi_i\tilde{r}_{i,v}+\widetilde{\mathcal{Z}},
	\end{align*}
	where $\widetilde{\mathcal{Z}}$ is defined as follows
	\begin{align*}
		&\widetilde{\mathcal{Z}}=\sum\limits_{i}v_{i,x}\tilde{r}_{i,t}+\sum\limits_{i,j}(v_{i,t}v_j+v_iv_{j,t})\tilde{r}_{i,u}\tilde{r}_j+\sum\limits_{i,j}v_iv_j\pa_t(\tilde{r}_{i,u}\tilde{r}_j)+\sum\limits_{i} v_{i,x}v_{i,t} \pa_{v_i}\bar{v}_i\xi_i\tilde{r}_{i,v}\\
		&+\sum\limits_{i}v_{i,x}v_i\pa_t(\pa_{v_i}\bar{v}_i)\xi_i\tilde{r}_{i,v}+\sum\limits_{i}v_{i,x}v_i\pa_{v_i}\bar{v}_i\xi_i^\p\left(\frac{w_i}{v_i}\right)_t\tilde{r}_{i,v}+\sum\limits_{i}v_{i,x}v_i\pa_{v_i}\bar{v}_i\xi_i\pa_t(\tilde{r}_{i,v})\\
		&+\sum\limits_{i}\xi_i^{\p\p}\left(\frac{w_i}{v_i}\right)_t\bar{v}_i\left(w_{i,x}-\frac{w_i}{v_i}v_{i,x}\right)\tilde{r}_{i,v}+\sum\limits_{i}\xi_i^\p\bar{v}_{i,t}\left(w_{i,x}-\frac{w_i}{v_i}v_{i,x}\right)\tilde{r}_{i,v}\\
		&-\sum\limits_{i}\xi_i^\p\bar{v}_i\left(\frac{w_i}{v_i}\right)_tv_{i,x}\tilde{r}_{i,v}-\sum\limits_{i}\xi_i^\p\bar{v}_i\left(w_{i,x}-\frac{w_i}{v_i}v_{i,x}\right)\pa_t(\tilde{r}_{i,v})\\
		&-\sum\limits_{i}\theta_i^{\p\p}\left(\frac{w_i}{v_i}\right)_t\left(w_{i,x}-\frac{w_i}{v_i}v_{i,x}\right)\tilde{r}_{i,\si}+\sum\limits_{i}\big[\theta_i^\p\left( \frac{w_i}{v_i}\right)_tv_{i,x}\tilde{r}_{i,\si}-\theta_i^\p v_i \left(\frac{w_i}{v_i}\right)_x \pa_t(\tilde{r}_{i,\si})\big]\\
				&+\sum\limits_i\sum\limits_{j\neq i}  v_{j,x}\big[\big(v_{i,t}\pa_{v_j}\bar{v}_i +v_i\pa_t(\pa_{v_j}\bar{v}_i)\big)\xi_i\tilde{r}_{i,v}+v_i\pa_{v_j}\bar{v}_i\big( \xi_i^\p\left(\frac{w_i}{v_i}\right)_t\tilde{r}_{i,v}+\xi_i\pa_t(\tilde{r}_{i,v})\big)\big].
	\end{align*}
	It yields from \eqref{defriVVs}
	\begin{equation}\label{eqn:u-xxt}
		u_{xxt}=\sum\limits_{i}\left[v_{i,tx}r^{VV}_i+w_{i,tx}r^{VW}_i\right]+\widetilde{\mathcal{Z}}.
	\end{equation}
	We further observe from Lemmas \ref{estimimpo1}, \ref{lemme6.5}, \ref{lemme11.3}, \ref{lemme11.4} and using the fact that $\tilde{r}_{i,\sig}=O(1)\mathfrak{A}_i\bar{v}_i$
	\begin{equation}
		\widetilde{\mathcal{Z}}=O(1)\de_0^2\sum\limits_{i}(\abs{w_{i,t}}+\abs{v_{i,t}}+\abs{w_{i,x}}+\abs{v_{i,x}}+\abs{w_i}+\abs{v_i}).
		\label{B15}
	\end{equation}
	By taking derivative in time of \eqref{eqn:u-tx} we obtain
	\begin{align*}
		&u_{txt}=\sum\limits_{i}(w_{i,tx}-\la_i^*v_{i,tx})\tilde{r}_i+\sum\limits_{i}(w_{i,x}-\la_i^*v_{i,x})\tilde{r}_{i,t}\\
		&+\sum\limits_{i,j}\big[\big((w_{i,t}-\la_i^*v_{i,t})v_j+(w_{i}-\la_i^*v_{i})v_{j,t}\big)\tilde{r}_{i,u}\tilde{r}_j+(w_{i}-\la_i^*v_{i})v_j\pa_t(\tilde{r}_{i,u}\tilde{r}_j)\big]\\
		&+\sum\limits_{i}\big[[v_{i,tx}(w_{i}-\la_i^*v_{i})+v_{i,x}(w_{i,t}-\la_i^*v_{i,t})]\pa_{v_i}\bar{v}_i+v_{i,x}(w_{i}-\la_i^*v_{i})\pa_t(\pa_{v_i}\bar{v}_i)\big]\xi_i\tilde{r}_{i,v}\\%+\sum\limits_{i}v_{i,x}(w_{i}-\la_i^*v_{i})\pa_t(\pa_{v_i}\bar{v}_i)\xi_i\tilde{r}_{i,v}\\
		&+\sum\limits_{i}v_{i,x}(w_{i}-\la_i^*v_{i})\pa_{v_i}\bar{v}_i\xi_i^\p\left(\frac{w_i}{v_i}\right)_t\tilde{r}_{i,v}+\sum\limits_{i}v_{i,x}(w_{i}-\la_i^*v_{i})\pa_{v_i}\bar{v}_i\xi_i\pa_t(\tilde{r}_{i,v})\\
		&+\sum\limits_{i}\big[\xi_i^{\p\p}\left(\frac{w_i}{v_i}\right)_t\bar{v}_i+\xi_i^\p\bar{v}_{i,t}\big]
		\left(w_{i,x}-\frac{w_i}{v_i}v_{i,x}\right)\left(\frac{w_i}{v_i}-\la_i^*\right)\tilde{r}_{i,v}\\
		%&+\sum\limits_{i}\xi_i^\p\bar{v}_{i,t}\left(w_{i,x}-\frac{w_i}{v_i}v_{i,x}\right)\left(\frac{w_i}{v_i}-\la_i^*\right)\tilde{r}_{i,v}\\
		&+\sum\limits_{i}\xi_i^\p\bar{v}_i\left(w_{i,tx}-\frac{w_i}{v_i}v_{i,tx}\right)\left(\frac{w_i}{v_i}-\la_i^*\right)\tilde{r}_{i,v}-\sum\limits_{i}\xi_i^\p\bar{v}_i\left(\frac{w_i}{v_i}\right)_tv_{i,x}\left(\frac{w_i}{v_i}-\la_i^*\right)\tilde{r}_{i,v}\\
		&+\sum\limits_{i}\xi_i^\p\bar{v}_i\left(w_{i,x}-\frac{w_i}{v_i}v_{i,x}\right)\big[\left(\frac{w_i}{v_i}-\la_i^*\right)\pa_t(\tilde{r}_{i,v})+\left(\frac{w_i}{v_i}\right)_t\tilde{r}_{i,v}\big]\\
		%&+\sum\limits_{i}\xi_i^\p\bar{v}_i\left(w_{i,x}-\frac{w_i}{v_i}v_{i,x}\right)\left(\frac{w_i}{v_i}\right)_t\tilde{r}_{i,v}\\
		&-\sum\limits_{i}\theta_i^{\p\p}\left(\frac{w_i}{v_i}\right)_t v_i\left(\frac{w_i}{v_i}\right)_x \left(\frac{w_i}{v_i}-\la_i^*\right)\tilde{r}_{i,\si}-\sum\limits_{i}\theta_i^\p\left(w_{i,tx}-\frac{w_i}{v_i}v_{i,tx}\right)\left(\frac{w_i}{v_i}-\la_i^*\right)\tilde{r}_{i,\si}\\
		&+\sum\limits_{i}\theta_i^\p\left( \frac{w_i}{v_i}\right)_tv_{i,x}\left(\frac{w_i}{v_i}-\la_i^*\right)\tilde{r}_{i,\si}-\sum\limits_{i}\theta_i^\p v_i \left(\frac{w_i}{v_i}\right)_x\left(\frac{w_i}{v_i}\right)_t\tilde{r}_{i,\si}\\
		&-\sum\limits_{i}\theta_i^\p\left(w_{i,x}-\frac{w_i}{v_i}v_{i,x}\right)\left(\frac{w_i}{v_i}-\la_i^*\right)\pa_t(\tilde{r}_{i,\si})+\sum\limits_i\sum\limits_{j\neq i}\left[ v_{j,tx}(w_i-\la_i^*v_i)\pa_{v_j}\bar{v}_i\right]\xi_i\tilde{r}_{i,v}\\
		&+\sum\limits_i\sum\limits_{j\neq i}v_{j,x} \big[ (w_{i,t}-\la_i^*v_{i,t})\pa_{v_j}\bar{v}_i+(w_i-\la_i^*v_i)\pa_t(\pa_{v_j}\bar{v}_i)
		\big]\xi_i\tilde{r}_{i,v}\\
		%&+\sum\limits_i\sum\limits_{j\neq i}\left[ v_{j,x}(w_i-\la_i^*v_i)\pa_t(\pa_{v_j}\bar{v}_i)\right]\xi_i\tilde{r}_{i,v}\\
		&+\sum\limits_i\sum\limits_{j\neq i}\left[ v_{j,x}(w_i-\la_i^*v_i)\pa_{v_j}\bar{v}_i\right]\big[\xi_i\pa_t(\tilde{r}_{i,v})+\xi'_i\left(\frac{w_i}{v_i}\right)_t\tilde{r}_{i,v}\big].
		%&+\sum\limits_i\sum\limits_{j\neq i}\left[ v_{j,x} (w_i-\la_i^*v_i)\pa_{v_j}\bar{v}_i\right]\left(\frac{w_i}{v_i}\right)_t\xi_i\tilde{r}_{i,v}\\
		%&+\sum\limits_i\sum\limits_{j\neq i}\left[ v_{j,x}(w_i-\la_i^*v_i)\pa_{v_j}\bar{v}_i\right]\xi'_i\left(\frac{w_i}{v_i}\right)_t\tilde{r}_{i,v}.
	\end{align*}
	By rearranging the terms we get
	\begin{align*}
		&u_{txt}=\sum\limits_{i}(w_{i,tx}-\la_i^*v_{i,tx})\tilde{r}_i+\sum\limits_{i}v_{i,tx}(w_{i}-\la_i^*v_{i}) \pa_{v_i}\bar{v}_i\xi_i\tilde{r}_{i,v}\\
		&+\sum\limits_{i}\xi_i^\p\bar{v}_i\left(w_{i,tx}-\frac{w_i}{v_i}v_{i,tx}\right)\left(\frac{w_i}{v_i}-\la_i^*\right)\tilde{r}_{i,v}-\sum\limits_{i}\theta_i^\p\left(w_{i,tx}-\frac{w_i}{v_i}v_{i,tx}\right)\left(\frac{w_i}{v_i}-\la_i^*\right)\tilde{r}_{i,\si}\\
		&+\sum\limits_i\sum\limits_{j\neq i}v_{j,tx}(w_i-\la_i^*v_i)\pa_{v_j}\bar{v}_i \xi_i\tilde{r}_{i,v}+	\widetilde{\mathcal{X}},
	\end{align*}
	%where $	\widetilde{\mathcal{X}}$ is defined as follows
	%\begin{align*}
	%	\widetilde{\mathcal{X}}&=\sum\limits_{i}(w_{i,tx}-\la_i^*v_{i,tx})\tilde{r}_i+\sum\limits_{i}v_{i,tx}(w_{i}-\la_i^*v_{i}) \pa_{v_i}\bar{v}_i\xi_i\tilde{r}_{i,v}\\
	%	&+\sum\limits_{i}\xi_i^\p\bar{v}_i\left(w_{i,tx}-\frac{w_i}{v_i}v_{i,tx}\right)\left(\frac{w_i}{v_i}-\la_i^*\right)\tilde{r}_{i,v}-\sum\limits_{i}\theta_i^\p\left(w_{i,tx}-\frac{w_i}{v_i}v_{i,tx}\right)\left(\frac{w_i}{v_i}-\la_i^*\right)\tilde{r}_{i,\si}\\
	%	&+\sum\limits_i\sum\limits_{j\neq i}v_{j,tx}(w_i-\la_i^*v_i)\pa_{v_j}\bar{v}_i \xi_i\tilde{r}_{i,v}+	\widetilde{\mathcal{X}},
	%\end{align*}
	where $	\widetilde{\mathcal{X}}$ is defined as follows
	\begin{align*}
		&\widetilde{\mathcal{X}}=
		\sum_{ij}\big[\big((w_{i,t}-\la_i^*v_{i,t})v_j+(w_{i}-\la_i^*v_{i})v_{j,t}\big)\tilde{r}_{i,u}\tilde{r}_j+(w_{i}-\la_i^*v_{i})v_j\pa_t(\tilde{r}_{i,u}\tilde{r}_j)\big]\\
		&+\sum\limits_{i}(w_{i,x}-\la_i^*v_{i,x})\tilde{r}_{i,t}+\sum\limits_{i}\big[v_{i,x}(w_{i,t}-\la_i^*v_{i,t})\pa_{v_i}\bar{v}_i+v_{i,x}(w_{i}-\la_i^*v_{i})\pa_t(\pa_{v_i}\bar{v}_i)\big]\xi_i\tilde{r}_{i,v}\\%+\sum\limits_{i}v_{i,x}(w_{i}-\la_i^*v_{i})\pa_t(\pa_{v_i}\bar{v}_i)\xi_i\tilde{r}_{i,v}\\
		&+\sum\limits_{i}v_{i,x}(w_{i}-\la_i^*v_{i})\pa_{v_i}\bar{v}_i\xi_i^\p\left(\frac{w_i}{v_i}\right)_t\tilde{r}_{i,v}+\sum\limits_{i}v_{i,x}(w_{i}-\la_i^*v_{i})\pa_{v_i}\bar{v}_i\xi_i\pa_t(\tilde{r}_{i,v})\\
		&+\sum\limits_{i}\big[\big(\xi_i^{\p\p}\left(\frac{w_i}{v_i}\right)_t\bar{v}_i+\xi_i^\p\bar{v}_{i,t}\big)v_i \left(\frac{w_i}{v_i}\right)_x-\xi_i^\p\bar{v}_i\left(\frac{w_i}{v_i}\right)_tv_{i,x}\big]
		\left(\frac{w_i}{v_i}-\la_i^*\right)\tilde{r}_{i,v}\\
		%&+\sum\limits_{i}\xi_i^\p\bar{v}_{i,t}\left(w_{i,x}-\frac{w_i}{v_i}v_{i,x}\right)\left(\frac{w_i}{v_i}-\la_i^*\right)\tilde{r}_{i,v}\\
		%-\sum\limits_{i}\xi_i^\p\bar{v}_i\left(\frac{w_i}{v_i}\right)_tv_{i,x}\left(\frac{w_i}{v_i}-\la_i^*\right)\tilde{r}_{i,v}\\
		&+\sum\limits_{i}\xi_i^\p\bar{v}_i\left(w_{i,x}-\frac{w_i}{v_i}v_{i,x}\right)\big[\left(\frac{w_i}{v_i}-\la_i^*\right)\pa_t(\tilde{r}_{i,v})+\left(\frac{w_i}{v_i}\right)_t\tilde{r}_{i,v}\big]\\
		%&+\sum\limits_{i}\xi_i^\p\bar{v}_i\left(w_{i,x}-\frac{w_i}{v_i}v_{i,x}\right)\left(\frac{w_i}{v_i}\right)_t\tilde{r}_{i,v}\\
		&-\sum\limits_{i}\theta_i^{\p\p}\left(\frac{w_i}{v_i}\right)_t v_i\left(\frac{w_i}{v_i}\right)_x \left(\frac{w_i}{v_i}-\la_i^*\right)\tilde{r}_{i,\si}+\sum\limits_{i}\theta_i^\p\left( \frac{w_i}{v_i}\right)_tv_{i,x}\left(\frac{w_i}{v_i}-\la_i^*\right)\tilde{r}_{i,\si}\\
		&-\sum\limits_{i}\theta_i^\p v_i \left(\frac{w_i}{v_i}\right)_x\left(\frac{w_i}{v_i}\right)_t\tilde{r}_{i,\si}-\sum\limits_{i}\theta_i^\p\left(w_{i,x}-\frac{w_i}{v_i}v_{i,x}\right)\left(\frac{w_i}{v_i}-\la_i^*\right)\pa_t(\tilde{r}_{i,\si})\\
		&+\sum\limits_i\sum\limits_{j\neq i}v_{j,x} \big[ (w_{i,t}-\la_i^*v_{i,t})\pa_{v_j}\bar{v}_i+(w_i-\la_i^*v_i)\pa_t(\pa_{v_j}\bar{v}_i)
		\big]\xi_i\tilde{r}_{i,v}\\
		%&+\sum\limits_i\sum\limits_{j\neq i}\left[ v_{j,x}(w_i-\la_i^*v_i)\pa_t(\pa_{v_j}\bar{v}_i)\right]\xi_i\tilde{r}_{i,v}\\
		&+\sum\limits_i\sum\limits_{j\neq i}\left[ v_{j,x}(w_i-\la_i^*v_i)\pa_{v_j}\bar{v}_i\right]\big[\xi_i\pa_t(\tilde{r}_{i,v})+\xi'_i\left(\frac{w_i}{v_i}\right)_t\tilde{r}_{i,v}\big].
		%&+\sum\limits_i\sum\limits_{j\neq i}\left[ v_{j,x} (w_i-\la_i^*v_i)\pa_{v_j}\bar{v}_i\right]\left(\frac{w_i}{v_i}\right)_t\xi_i\tilde{r}_{i,v}\\
		%&+\sum\limits_i\sum\limits_{j\neq i}\left[ v_{j,x}(w_i-\la_i^*v_i)\pa_{v_j}\bar{v}_i\right]\xi'_i\left(\frac{w_i}{v_i}\right)_t\tilde{r}_{i,v}.
	\end{align*}
	Note that from Lemmas \ref{estimimpo1}, \ref{lemme6.5}, \ref{lemme11.3}, \ref{lemme11.4} and the fact that $\tilde{r}_{i,\sig}=O(1)\mathfrak{A}_i\bar{v}_i$
	\begin{equation}
		\widetilde{\mathcal{X}}=O(1)\de_0^2 \sum\limits_{i}(\abs{w_{i,t}}+\abs{v_{i,t}}+\abs{w_{i,x}}+\abs{v_{i,x}}+\abs{w_i}+\abs{v_i}).\label{B16bis}
	\end{equation}
	Furthermore, it follows
	\begin{equation}\label{eqn:u-txt}
		u_{txt}=\sum\limits_{i}\left[w_{i,tx}r^{WW}_{i}+v_{i,tx}r^{WV}_i\right]+\widetilde{\mathcal{X}}.
	\end{equation}
	Combining \eqref{eqn:u-xxt}, \eqref{B15},  \eqref{eqn:u-txt} and \eqref{B16bis} it follows that for $t\in[\hat{t},T]$
	\begin{equation*}
		\begin{pmatrix}
			u_{xxt}\\
			u_{txt}
		\end{pmatrix}=\begin{pmatrix}
			\mathcal{Q}^{VV}&	\mathcal{Q}^{VW}\\
			\mathcal{Q}^{WV}&	\mathcal{Q}^{WW}
		\end{pmatrix}\begin{pmatrix}
			v_{1,tx}\\
			\vdots\\
			v_{n,tx}\\
			w_{1,tx}\\
			\vdots\\
			w_{n,tx}
		\end{pmatrix}+O(1)\de_0^2 \sum\limits_{i}(\abs{w_{i,t}}+\abs{v_{i,t}}+\abs{w_{i,x}}+\abs{v_{i,x}}+\abs{w_i}+\abs{v_i})%+\begin{pmatrix}
			%\widetilde{\mathcal{Z}}\\
			%\widetilde{\mathcal{X}}
		%\end{pmatrix}.
	\end{equation*}
By a similar argument as before, we deduce then that for any $i\in\{1,\cdot,n\}$ we have using  Corollary \ref{coro2.2} and the equations satisfy by $u_{t,xx}$, $u_{ttx}$ that for $t\in[\hat{t},T]$
		\begin{align*}
		&v_{i,xt}(t,\cdot), w_{i,xt}(t,\cdot)=O(1)(|u_{txx}|+|u_{ttx}|)(t,\cdot)\\
		&\hspace{5cm}+\de_0^2\sum\limits_{i}(|w_{i,t}|+|v_{i,t}|+\abs{w_{i,x}}+\abs{v_{i,x}}+\abs{w_i}+\abs{v_i})(t,\cdot),\\
		&=O(1)\big(|u_{tx}||u_x|+|u_t||u_x|^2+|u_t| |u_{xx}|+|u_{t,xx}|+|u_{t,xx}| |u_x|+|u_{t,x}| |u_x|^2+|u_{tx}||u_{xx}|\\
		&+|u_{xx}||u_x||u_t|+|u_x|^3|u_t|+|u_t||u_{xxx}|+|u_{xxxt}|\big)(t,\cdot)\\
		&+\de_0^2\sum\limits_{i}(|w_{i,t}|+|v_{i,t}|+\abs{w_{i,x}}+\abs{v_{i,x}}+\abs{w_i}+\abs{v_i})(t,\cdot),\\
		&=O(1)\big(|u_{tx}||u_x|+|u_t||u_x|^2+|u_t| |u_{xx}|+|u_{t,xx}|+|u_{tx}||u_{xx}|+|u_t||u_{xxx}|+|u_{xxxt}|\big)(t,\cdot)\\
		&+\de_0^2\sum\limits_{i}(|w_{i,t}|+|v_{i,t}|+\abs{w_{i,x}}+\abs{v_{i,x}}+\abs{w_i}+\abs{v_i})(t,\cdot),\\
&=O(1)\big(|u_{tx}||u_x|+|u_t||u_x|^2+|u_t| |u_{xx}|+|u_{tx}||u_{xx}|+|u_t||u_{xxx}|+|u_{xxxx}|\big)(t,\cdot)\\
&+O(1)\big(|u_x| |u_{xx}|+|u_x|^3+|u_{xxx}|+|u_x|  |u_{xxx}|+|u_x|^2 |u_{xx}|+|u_{xx}|^2+|u_x|^4+|u_{xxxx}|\big)(t,\cdot)\\
&+\de_0^2\sum\limits_{i}(|w_{i,t}|+|v_{i,t}|+\abs{w_{i,x}}+\abs{v_{i,x}}+\abs{w_i}+\abs{v_i})(t,\cdot).
\end{align*}
Applying now Lemma \ref{lemme6.5} and Corollary \ref{coro2.2}, we obtain the desired estimate \eqref{estimate-v-w-tx}. This completes the proof of Lemma \ref{lemme6.6}.
\end{proof}

%--------------------------------------------------------------------------------------------------------------

\section{Proof of Lemma \ref{lemma:transversal-1} and \ref{lemma:transversal-2}}\label{sec:proof-of-transversal}
Here we first give a proof for Lemma \ref{lemma:transversal-1}. It can also be found in \cite{HJ-temple-class} when $\mu,\mu^\#$ are non-constant. For constant $\mu,\mu^\#$, the result was first established in \cite{BB-temple,BB-vv-lim-ann-math}. For the sake of completeness we give a sketch of the proof.
\begin{proof}[Proof of Lemma \ref{lemma:transversal-1}:]
	Set $c_1:=\norm{\mu,\mu^\#}_{L^\f}$. Let $z,z^\#$ be the solution to \eqref{eqn-z-1}, \eqref{eqn-z-2} with $\varphi=\varphi^\#=0$. Consider
	\begin{equation}
		Q(z,z^\#):=\int\limits_{\R}\int\limits_{\R} K(x-y)\abs{z(x)}\abs{z^\#(y)}\,dxdy,
	\end{equation}
	where $K$ is defined as follows
	\begin{equation}
		K(s):=\left\{\begin{array}{rl}
			1/c&\mbox{ if }s\geq 0,\\
			1/c e^{\frac{cs}{2c_1}}&\mbox{ if }s<0.
		\end{array}\right.
	\end{equation}
	Now, we can calculate
	\begin{align*}
		\frac{d}{dt}Q(z(t),z^\#(t))&=\int\limits_{\R}\int\limits_{\R} K(x-y)[sgn(z(x))z_t(x)\abs{z^\#(y)}+sgn(z^\#(y))z^\#_t(y)\abs{z(x)}]\,dxdy\\
		&=\int\int K(x-y)\bigg[sgn(z(x))((\mu(x) z_x(x))_x-(\la z(x))_x)\abs{z^\#(y)}\\
		&\quad\quad+sgn(z^\#(y))((\mu^\#(y) z^\#_y(y))_y-(\la^\# z^\#(y))_y)\abs{z(x)}\bigg]\,dxdy\\
		&=\int\limits_{\R}\int\limits_{\R} K^\p(x-y)\bigg\{\la\abs{z(x)}\abs{z^\#(y)}-\la^\#\abs{z(x)}\abs{z^\#(y)}\bigg\}\,dxdy\\
		&+\int\limits_{\R}\int\limits_{\R} K^\p(x-y)\bigg\{\mu_x\abs{z(x)}\abs{z^\#(y)}-\mu_y^\#\abs{z(x)}\abs{z^\#(y)}\bigg\}\,dxdy\\
		&+\int\limits_{\R}\int\limits_{\R} K^{\p\p}(x-y)\bigg\{\mu(x)\abs{z(x)}\abs{z^\#(y)}+\mu^\#(y)\abs{z(x)}\abs{z^\#(y)}\bigg\}\,dxdy\\
		&\leq -\int\limits_{\R}\int\limits_{\R} (cK^\p(x-y)-2c_1K^{\p\p}(x-y))\abs{z(x)}\abs{z^\#(y)}\,dxdy\\
		&\leq -\int\limits_{\R} \abs{z(x)}\abs{z^\#(x)}\,dx.
	\end{align*}
	Hence, we get
	\begin{equation}
		\int_{0}^T\int\limits_{\R} \abs{z(t,x)}\abs{z^\#(t,x)}\,dxdt\leq Q(z(0),z^\#(0))\leq \frac{1}{c}\norm{z(0)}_{L^1}\norm{z^\#(0)}_{L^1}.
	\end{equation}
	Now for $\varphi,\varphi^\#\neq 0$, the proof follows exactly same way as in \cite{BB-vv-lim-ann-math}. This completes the proof of Lemma \ref{lemma:transversal-1}.
\end{proof}
Next, we prove Lemma \ref{lemma:transversal-2}. Proof is based on a similar argument as in \cite[Lemma 7.3]{BB-vv-lim-ann-math}. For sake of completeness, we give a sketch of the proof below.
\begin{proof}[Proof of Lemma \ref{lemma:transversal-2}:]
	We define
	\begin{equation}
		\mathcal{I}_1(T):=\sup\limits_{(\tau,\xi)\in[0,T]\times R}\int\limits_{0}^{T-\tau}\int\limits_{\R}\abs{z_x(t,x)z^\#(t+\tau,x+\xi)}\,dxdt\leq (C^*\de_0^2)^2T.
	\end{equation}
	Let us denote the solution of $z_t=(\mu z_x)_x$ as $G^\mu$. Then $G^\mu$ satisfies the following estimate
	\begin{equation}
		\abs{G^\mu_x(t,x;s,y)}\leq C_1e^{-q_*\frac{(x-y)^2}{t-s}},\mbox{ for some }C_1,q_*>0.
	\end{equation}
	Such $G^\mu$ can be constructed by method of parametrix (see \cite[Theorem 11, Chapter 1]{Friedman}). Therefore, we have 
	\begin{equation}\label{estimate-G-mu}
		\norm{G^\mu_x(1,x;0,x-\cdot)}_{L^1(\R)},\, \norm{G^\mu_x(t,x;t-\cdot,x-\cdot)}_{L^1(0,1;L^1(\R))}\leq C_1^*,
	\end{equation} 
	for some $C_1^*>0$. We can write
	\begin{align*}
		z_x(t,x)&=\int\limits_{\R} G^\mu_x(1,x;0,y)z(t-1,y)\,dy+\int\limits_{0}^1\int\limits_{\R} G^\mu_x(1,x;s,y)[\varphi-(\la z)_x](t-1+s,y)\,dyds.
	\end{align*}
	By using \eqref{estimate-G-mu} we get
	\begin{align*}
		&\int\limits_{1}^{T-\tau}\int\limits_{\R} \abs{z_x(t,x)z^\#(t+\tau,x+\xi)}\,dxdt\\
		&\leq \int\limits_{1}^{T-\tau}\int\limits_{\R}\int\limits_{\R} \abs{G^\mu_x(1,x;0,x-y) z(t-1,x-y)z^\#(t+\ta,x+\xi)}\,dydxdt\\
		&+\int\limits_{1}^{T-\tau}\int\limits_{\R}\int\limits_{0}^{1}\int\limits_{\R} \norm{\la_x}_{L^\f}\abs{G^\mu_x(1,x;s,x-y) z(t-1+s,x-y)z^\#(t+\ta,x+\xi)}\,dydsdxdt\\
		&+\int\limits_{1}^{T-\tau}\int\limits_{\R}\int\limits_{0}^{1}\int\limits_{\R} \norm{\la}_{L^\f}\abs{G^\mu_x(1,x;s,x-y) z_x(t-1+s,x-y)z^\#(t+\ta,x+\xi)}\,dydsdxdt\\
		&+\int\limits_{1}^{T-\tau}\int\limits_{\R}\int\limits_{0}^{1}\int\limits_{\R} \abs{G^\mu_x(1,x;s,x-y) \varphi(t-1+s,x-y)z^\#(t+\ta,x+\xi)}\,dydsdxdt\\
		&\leq \left(\int\limits_{\R}\abs{G_x^\mu(1,x)}\,dx+\sup\limits_{x}\int\limits_{0}^{1}\int\limits_{\R} \norm{\la_x}_{L^\f}\abs{G^\mu_x(1,x;s,x-y)}\,dyds\right)\cdot\\
		&\quad\quad\quad\quad\quad\quad\quad\quad\quad\quad\quad\quad\sup\limits_{s,y,\tau,\xi}\Big(\int\limits_{1}^{T-\tau}\int\limits_{\R}\abs{ z(t-1+s,x-y)z^\#(t+\ta,x+\xi)}\,dxdt\Big)\\
		&+\left(\sup\limits_{x}\int\limits_{0}^{1}\int\limits_{\R} \norm{\la}_{L^\f}\abs{G^\mu_x(1,x;s,x-y)}\,dyds\right)\cdot\\
		&\quad\quad\quad\quad\quad\quad\quad\quad\quad\quad\quad\quad\sup\limits_{s,y,\tau,\xi}\Big(\int\limits_{1}^{T-\tau}\int\limits_{\R}\abs{ z_x(t-1+s,x-y)z^\#(t+\ta,x+\xi)}\,dxdt\Big)\\
		&+\norm{z^\#}_{L^\f([0,T]\times\R)}\left[\sup\limits_{x}\int\limits_{0}^{1}\int\limits_{\R} \abs{G^\mu_x(1,x;s,x-y)}\,dyds\right]\cdot \int\limits_{0}^{T}\int\limits_{\R}\abs{\varphi(t,x)}\,dxdt\\
		&\leq C^*_1(1+\norm{\la_x}_{L^\f})\frac{4\de_0^2}{c}+C_1^*\norm{\la}_{L^\f}\mathcal{I}_1(T)+C_1\de_0^3.
	\end{align*}
	We note that $\mathcal{I}_1(1)\leq C^2_1\de_0^2$. Since $C_1^*\norm{\la}_{L^\f}\leq C_1^*\norm{\la_x}_{L^\f}\leq C_1^*\de_0<1/2$, we obtain
	\begin{equation}
		\mathcal{I}_1(T)\leq C^*_1(1+\norm{\la_x}_{L^\f})\frac{4\de_0^2}{c}+C^2_1\de_0^2+C_1\de_0^3+\frac{1}{2}\mathcal{I}_1(T).
	\end{equation}
	Hence, we get \eqref{est:transversal-2}. This completes the proof of Lemma \ref{lemma:transversal-2}.
\end{proof}

%----------------------------------------------------------------------------
%---------------------start--of--new--section-------------------------------
%----------------------------------------------------------------------------
%---------------------start--of--new--section-------------------------------
%----------------------------------------------------------------------------

\section{Proof of Proposition \ref{prop:parabolic}}\label{appendix:proof-prarabolic}
Here, we prove Proposition \ref{prop:parabolic}. For constant $B(u)$, similar estimates have been established in \cite{BB-vv-lim-ann-math}. For non-constant $B$, proof of Proposition \ref{prop:parabolic} for $k\leq2$ can be found in \cite{HJ-temple-class}. For higher order estimates the proof follows in a similar manner. For the sake of completeness we give a complete proof with explicit calculations.
\begin{proof}[Proof of Proposition \ref{prop:parabolic}:] We can take derivative w.r.t $x$ on both sides of \eqref{eqn-main} to get
	\begin{equation}
		(u_x)_t+A(u)u_{xx}=B(u)u_{xxx}+(u_x\bullet B(u)u_x)_x-u_x\bullet A(u)u_x+{\color{black}{u_x\bullet B(u)u_{xx}}}.
		\label{D1}
	\end{equation}
	We would like to make a change of variable $v=P(u)u_x$. To this end, multiplying \eqref{D1} on the left by $P(u)$ we get
	\begin{align*}
		&(Pu_x)_t+PA(u)P^{-1}(P(u)u_{x})_x=PBP^{-1}(Pu_x)_{xx}-PBP^{-1}(u_x\bullet P(u)u_x)_x\\
		&-PBP^{-1}[u_x\bullet Pu_{xx}]+u_t\bullet P(u)u_x+PA(u)P^{-1}(u_x\bullet P(u)u_{x})+P(u_x\bullet B(u)u_x)_x\\
		&-Pu_x\bullet A(u)u_x+{\color{black}{P u_x\bullet B(u)u_{xx}}}
	\end{align*}
	Then, if $A_1(u)=P(u)A(u)P(u)^{-1}$ and $B_1(u)=\mbox{diag}(\mu_1(u),\cdots,\mu_n(u))=P(u)B(u)P(u)^{-1}$, we have
	\begin{align*}
		&v_t+A_1(u)v_x=B_1(u)v_{xx}-B_1((P^{-1}(u)v)\bullet PP^{-1}v)_x-B_1[(P^{-1}v)\bullet P(P^{-1}v)_x]\\
		&+u_t\bullet P(u)P^{-1}v+A_1(u)(P^{-1}(u)v\bullet P(u)P^{-1}v)+P(P^{-1}v\bullet B(u)P^{-1}v)_x\\
		&-P(P^{-1}v)\bullet A(u)P^{-1}v+{\color{black}{P (P^{-1}v)\bullet B(u)(P^{-1}v)_x}}.
	\end{align*}
	Hence,
	\begin{align*}
		&v_t+A_1(u)v_x=B_1(u)v_{xx}-B_1((P^{-1}v \bullet P^{-1}(u)v)\bullet PP^{-1}v)-B_1((P^{-1}(u)v_x)\bullet PP^{-1}v)\\
		&-B_1((P^{-1}(u)v) \otimes (P^{-1}(u)v): D^2PP^{-1}v)-2B_1[(P^{-1}v)\bullet P(P^{-1}v)\bullet P^{-1}v]\\
		&-2B_1[(P^{-1}v)\bullet PP^{-1}v_x]+u_t\bullet P(u)P^{-1}v+A_1(u)(P^{-1}(u)v\bullet P(u)P^{-1}v)\\
		&+P(P^{-1}v)\bullet P^{-1}v\bullet B(u)P^{-1}v+PP^{-1}v_x\bullet B(u)P^{-1}v\\
		&+P(P^{-1}v)\otimes(P^{-1}v):D^2B(u)P^{-1}v+{\color{black}{2}}P(P^{-1}v)\bullet B(u)(P^{-1}v)\bullet P^{-1}v\\
		&+{\color{black}{2}}P(P^{-1}v)\bullet B(u)P^{-1}v_x-P(P^{-1}v)\bullet A(u)P^{-1}v.
	\end{align*}
	We observe that $u_t=BP^{-1}v_x+B[(P^{-1}(u)v)\bullet P^{-1}(u)v]-AP^{-1}v+(P^{-1}(u)v)\bullet B(u)P^{-1}(u)v$. Then we get
	\begin{equation}		\label{princi}
		v_t+A_1^*v_x=B_1(u)v_{xx}+\mathcal{R},
	\end{equation}
	where $\mathcal{R}=\mathcal{R}(u,v,v_x)$ is defined as follows
	\begin{align}
		\mathcal{R}(u,v,v_x)&:=[A_1^*-A_1(u)]v_x-B_1(((P^{-1}(u)v)\bullet P^{-1}(u)v)\bullet PP^{-1}v)\nonumber\\
		&-B_1((P^{-1}(u)v_x)\bullet PP^{-1}v)-B_1((P^{-1}(u)v)\otimes (P^{-1}(u)v): D^2PP^{-1}v)\nonumber\\
		&-2B_1[(P^{-1}v)\bullet P(P^{-1}v)\bullet P^{-1}v]-2B_1[(P^{-1}v)\bullet PP^{-1}v_x]\nonumber\\
		&+[BP^{-1}v_x+B[(P^{-1}(u)v)\bullet P^{-1}(u)v]-AP^{-1}v]\bullet P(u)P^{-1}v\nonumber\\
		&+[(P^{-1}(u)v)\bullet B(u)P^{-1}(u)v]\bullet P(u)P^{-1}v+A_1(u)(P^{-1}(u)v\bullet P(u)P^{-1}v)\nonumber\\
		&+P(P^{-1}v)\bullet P^{-1}v\bullet B(u)P^{-1}v+PP^{-1}v_x\bullet B(u)P^{-1}v\nonumber\\
		&+P(P^{-1}v)\otimes(P^{-1}v):D^2B(u)P^{-1}v+{\color{black}{2}}P(P^{-1}v)\bullet B(u)(P^{-1}v)\bullet P^{-1}v\nonumber\\
		&+{\color{black}{2}}P(P^{-1}v)\bullet B(u)P^{-1}v_x-P(P^{-1}v)\bullet A(u)P^{-1}v.\label{def:remainder-R}
	\end{align}
	As in \cite{HJ-temple-class}, we apply a change of variable $v\mapsto \tilde{v}$ such that $v_i(t,x)=\tilde{v}_i(t,X_i(t,x))$ where $ (X_i)_x(t,x)=\frac{1}{\sqrt{\mu_i(u(t,x))}}$ for any $i\in\{1,\cdots,n\}$. Note that for any $t>0$,
	\begin{equation}\label{def:T-t-f}
		\mathcal{T}_t(f)_i(x)=f_i(X_i(t,x))\mbox{ where }X_i(t,x)=\int\limits_{0}^{x}\frac{1}{\sqrt{\mu_i(u(t,z))}}\,dz,
	\end{equation}
	it implies that $v= \mathcal{T}_t \tilde{v}$. We observe that for any $(t,x)\in]0,+\infty[\times\R$
	\begin{align*}
		&v_{t}(t,x)=\mathcal{T}_t(\tilde{v}_t{{\color{black}(t,\cdot)}})(x)+(\tilde{v}_{i,x}(t,X_i(t,x))X_{i,t})_{1\leq i\leq n},\\
		&v_x(t,x)=B_1^{-1/2}(u(t,x))\mathcal{T}_t({\color{black}{\tilde{v}_{x}(t,\cdot)}})(x),\\
		&v_{xx}(t,x)=B_1^{-1}(u(t,x))\mathcal{T}_t({\color{black}{\tilde{v}_{xx}(t,\cdot)}})(x)+(P^{-1}(u(t,x))v(t,x))\bullet B_1(u(t,x))^{-1/2}\mathcal{T}_t(\tilde{v}_{x}{\color{black}{(t,\cdot)}})(x).
	\end{align*}
	Then we have from \eqref{princi}
	\begin{align}
		& \mathcal{T}_t(\tilde{v}_t(t,\cdot))+(\tilde{v}_{i,x}(t,X_i(t,x))X_{i,t})_{1\leq i\leq n}+A_1^*B_1^{-1/2}(u)\mathcal{T}_t(\tilde{v}_{x}(t,\cdot))-\mathcal{T}_t(\tilde{v}_{xx}(t,\cdot))\\
		&-
		B_1(u) (P^{-1}(u)v)\bullet B_1(u)^{-1/2}\mathcal{T}_t(\tilde{v}_{x}(t,\cdot))=\mathcal{R}.
	\end{align}
	From Lemma \ref{lemma:transformation}, we know that for any $t>0$, $\mathcal{T}_t$ is invertible we get then
	\begin{align}
		&\tilde{v}_t+A_2^*\tilde{v}_x-\tilde{v}_{xx}=\mathcal{T}_t^{-1}\big([A_2^*-A^*_1B_1^{-1/2}(u)]\mathcal{T}_t(\tilde{v}_x(t,\cdot))\nonumber\\
		&+B_1(u)(P^{-1}(u) v)\bullet B_1^{-1/2}(u)\mathcal{T}_t(\tilde{v}_{x}(t,\cdot))\big)-\mathcal{T}_t^{-1} (M(X_t)\mathcal{T}_t( \tilde{v}_x(t,\cdot)))+ \mathcal{T}_t^{-1}(\mathcal{R}),\label{eqn-tilde-v}
	\end{align}
	where $A_2^*=A_1^*B_1^{-1/2}(u^*)$ and $M(X_t)=\mbox{diag}(X_{1,t},\cdots,X_{n,t}).$ We first prove \eqref{estimate:parabolic-1} for smooth initial data.  We argue by contradiction. To this end, first we assume that the conclusion \eqref{estimate:parabolic-1} does not hold for $k=1$. Due to the assumption of smoothness of initial data, by continuity \eqref{estimate:parabolic-1} is satisfied in a small time interval $[0,\alpha]$, we can assume again by continuity that there exists a time $t^*<\hat{t}$ such that \eqref{estimate:parabolic-1} holds for $t\in[0,t^*]$ with $k=1$ and there is equality  at time $t=t^*$. We compute 
	\begin{align}
		X_{i,t}(t,x)&=-\int\limits_{0}^{x}\frac{u_t(t,z)\cdot D\mu_i(u(t,z))}{2(\mu_i(u(t,z)))^{3/2}}\,dz\nonumber\\
		&=-\int\limits_{0}^{x}\frac{(B(u)u_x)_x(t,z)\cdot D\mu_i(u(t,z))-A(u)u_x(t,z)\cdot D\mu_i(u(t,z))}{2(\mu_i(u))^{3/2}(t,z)}\,dz.\nonumber
	\end{align}
	Using integration by parts we can have
	\begin{align*}
		&X_{i,t}(t,x)
		=-\int\limits_{0}^{x}\left[-\frac{(B(u)u_x)\otimes u_x: D^2\mu_i(u)}{2(\mu_i(u))^{3/2}}-\frac{3((B(u)u_x)\cdot D\mu_i(u))(u_x\cdot D\mu_i(u))}{4(\mu_i(u))^{5/2}}\right](t,z)\,dz\\
		&-\frac{B(u)u_x\cdot D\mu_i(u)}{2(\mu_i(u))^{3/2}}(t,x)+\frac{B(u)u_x\cdot D\mu_i(u)}{2(\mu_i(u))^{3/2}}(t,0)+\int\limits_{0}^{x}\frac{A(u)u_x\cdot D\mu_i(u)}{2(\mu_i(u))^{3/2}}(t,z)\,dz.
	\end{align*}
	Then we get for $t\in[0,t^*]$ using the fact that $\|u_x(t,\cdot)\|_{L^1}\leq \delta_0$, $\|u_x\|_{L^\infty}\leq\kappa\|v\|_{L^\infty}$ and the smallness assumption on $\delta_0$
	\begin{equation}	\label{eqtech}
		\norm{X_{i,t}(t,\cdot)}_{L^\f(\R)}\leq 8\kappa^9[\norm{v}_{L^\f}+\de_0].
	\end{equation}
	Furthermore, we observe from Lemma \ref{lemma:transformation} since $\tilde{v}(t,\cdot)={\cal T}_t^{-1} v(t,\cdot)$ and the smallness assumption on $\delta_0$ that
	\begin{align}
		&\norm{\tilde{v}(t,\cdot)}_{L^1}\leq \kappa_1\norm{v(t,\cdot)}_{L^1}\leq \kappa^2\norm{u_x(t,\cdot)}_{L^1},\label{D9}\\
		&\norm{\tilde{v}_x(t,\cdot)}_{L^1}\leq \norm{v_x(t,\cdot)}_{L^1}\leq 2\kappa\norm{u_{xx}(t,\cdot)}_{L^1},\label{D10}\\
		&\norm{u_{xx}(t,\cdot)}_{L^1}\leq \kappa_P^2\delta_0 \norm{u_{xx}(t,\cdot)}_{L^1}+\kappa_P \norm{v_x(t,\cdot)}_{L^1}\leq 2\kappa\norm{v_x(t,\cdot)}_{L^1}\leq 2\kappa\norm{\tilde{v}_x(t,\cdot)}_{L^1}.\label{D11}
	\end{align}
	From \eqref{eqn-tilde-v} we can write for $0<t\leq t^*$
	\begin{align*}
		&\tilde{v}_{x}(t,\cdot)=G_x(t/2)\star \tilde{v}(t/2)+\int\limits_{t/2}^{t}G_x(t-s)\star \Big\{\mathcal{T}_s^{-1}\big([A_2^*-A^*_1B_1^{-1/2}(u)]\mathcal{T}_s(\tilde{v}_x(s,\cdot))\\
		&+B_1(u)(P^{-1}(u) v)\bullet B_1^{-1/2}(u)\mathcal{T}_s(\tilde{v}_{x}(s,\cdot))\big)-\mathcal{T}_s^{-1} (M(X_s)\mathcal{T}_s( \tilde{v}_x(s,\cdot)))+ \mathcal{T}_s^{-1}(\mathcal{R})(s,\cdot)\Big\}\,ds.
	\end{align*}
	Using \eqref{def:kappa}, \eqref{D9}, Lemma \ref{lemma:transformation} we have
	\begin{align}
		&\norm{\tilde{v}_{x}(t)}_{L^1}\leq \norm{G_x(t/2)}_{L^1}\norm{\tilde{v}(t/2)}_{L^1}+\int\limits_{t/2}^{t}\norm{G_x(t-s)}_{L^1}\norm{\mathcal{T}_s^{-1}([A_2^*-A^*_1B_1^{-1/2}(u)]\mathcal{T}_s (\tilde{v}_x))}_{L^1}\,ds\nonumber\\
		&+\int\limits_{t/2}^{t}\norm{G_x(t-s)}_{L^1}\norm{\mathcal{T}_s^{-1}(B_1(P^{-1}{v})\bullet B_1^{-1/2}\mathcal{T}_s(\tilde{v}_{x}))-\mathcal{T}_s^{-1} (M(X_s)\mathcal{T}_s( \tilde{v}_x))+ \mathcal{T}_s^{-1}(\mathcal{R})}_{L^1}(s,\cdot)\,ds\nonumber\\
		&\leq \frac{\sqrt{2}\kappa^3\de_0}{\sqrt{t}}+\int\limits_{t/2}^{t}\frac{\kappa\kappa_1}{\sqrt{t-s}}\norm{[A_2^*-A^*_1B_1^{-1/2}(u)]\mathcal{T}_s(\tilde{v}_x)}_{L^1}(s,\cdot)\,ds\nonumber\\
		&+\int\limits_{t/2}^{t}\frac{\kappa\kappa_1}{\sqrt{t-s}}\norm{ B_1(P^{-1}{v})\bullet B_1^{-1/2}\mathcal{T}_s(\tilde{v}_{x})-M(X_s)\mathcal{T}_s( \tilde{v}_x)+\mathcal{R}}_{L^1}(s,\cdot)\,ds.\label{estimucle1}	
	\end{align}
	We note that using again the Lemma  \ref{lemma:transformation} and \eqref{eqtech}, \eqref{D10}
	\begin{align}
		&\norm{[A_2^*-A^*_1B_1^{-1/2}(u)]\mathcal{T}_s(\tilde{v}_x(s,\cdot))}_{L^1}
		\leq 2 \kappa^3\delta_0\norm{u_{xx}(s,\cdot)}_{L^1},\nonumber\\
		&\norm{B_1(P^{-1}{v})\bullet B_1^{-1/2}\mathcal{T}_s(\tilde{v}_{x}(s,\cdot))}_{L^1}
		\leq 2\kappa^6\norm{u_{xx}(s,\cdot)}^2_{L^1},\nonumber\\
		&\norm{-M(X_s)\mathcal{T}_s( \tilde{v}_x(s,\cdot))}_{L^1}
		\leq 16\kappa^{12}\norm{u_{xx}(s,\cdot)}_{L^1}[\norm{u_{xx}(s,\cdot)}_{L^1}+\de_0]
		.\label{estimucle2}
	\end{align}
	From \eqref{absolucle1} we get
	\begin{align}
		\norm{\mathcal{R}}_{L^1}&\leq  20\cdot 2^8\kappa^8\Big(\norm{u^*-u}_{L^\f}\norm{v_x}_{L^1}+\norm{v}_{L^\f}\norm{v}_{L^1}+\norm{v}_{L^\f}^2\norm{v}_{L^1}\nonumber\\
		&\hspace{2.4 cm}+\norm{v}_{L^\f}\norm{v_x}_{L^1}\Big)\nonumber\\
		&\leq 120\cdot 2^{10}\kappa^{11}\left[\norm{u_{xx}}_{L^1}^2+\de_0\norm{u_{xx}}_{L^1}\right].\label{estimucle3}
	\end{align}
	Therefore, using \eqref{estimucle1},  \eqref{estimucle2}, \eqref{estimucle3},  \eqref{estimate:parabolic-1} we get for any $t\in[0,t^*]$ and the fact that  $\int\limits_{t/2}^{t}\frac{1}{\sqrt{t-s}}\frac{1}{s}ds\leq \frac{2\sqrt{2}}{\sqrt{t}}$% and  $\int\limits_{t/2}^{t}\frac{1}{\sqrt{t-s}}\frac{1}{\sqrt{s}}ds\leq 2$
	\begin{align}
		&\norm{\tilde{v}_{x}(t)}_{L^1}
		\leq \frac{\sqrt{2}\kappa^3 \de_0}{\sqrt{t}}+150\cdot2^{11}\kappa^{20}\int\limits_{t/2}^{t}\frac{1}{\sqrt{t-s}}\left[\frac{\de_0^2}{s}+\frac{\de_0^2}{\sqrt{s}}\right]\,ds\nonumber\\
		&\leq \frac{\sqrt{2}\kappa^3\de_0}{\sqrt{t}}+\frac{150\cdot2^{13}\kappa^{20}\delta_0^2}{\sqrt{t}}+150\cdot2^{13}\kappa^{20}\delta_0^2\nonumber\\
		&= \frac{\kappa^3\de_0}{\sqrt{t}}\big(\sqrt{2}+150\cdot 2^{13}\kappa^{17}\delta_0+150\cdot 2^{13}\kappa^{17}\delta_0 \sqrt{t}\big)<\frac{(2+\sqrt{2} )\kappa^3 \de_0}{2\sqrt{t}},
		%&= \frac{\sqrt{2}\kappa^3\de_0}{\sqrt{t}}+\left(150\cdot2^{15}\kappa^{10}\delta_0+150\cdot2^{15}\kappa^{10}\delta_0\sqrt{t}\right)\frac{\kappa\de_0}{4\sqrt{t}}\\
		%&<\frac{(2-\sqrt{\de_0})\kappa\de_0}{\sqrt{t}},
		\label{estimucle4}
	\end{align}
	since $150\cdot 2^{13}\kappa^{17}\delta_0+150\cdot 2^{13}\kappa^{17}\delta_0 \sqrt{t}
	\leq \frac{150}{\kappa^{12}2^{33}}+\frac{150}{2^{32}\kappa^{12}}<\frac{2-\sqrt{2}}{2}$. Recall from \eqref{D11} that $\norm{u_{xx}}_{L^1}\leq \kappa\norm{\tilde{v}_x}_{L^1}+\kappa^2\de_0\norm{u_{xx}}_{L^1}$, hence for $t=t^*$ using \eqref{estimucle4}, \eqref{estimate:parabolic-1} it follows
	\begin{align}
		&\norm{u_{xx}(t^*,\cdot)}_{L^1}\leq \kappa\norm{\tilde{v}_x(t^*,\cdot)}_{L^1}+\frac{2 \kappa^6\de_0^2}{\sqrt{t^*}}\leq \frac{\kappa^4 \de_0}{\sqrt{t}}(\frac{2+\sqrt{2}}{2}+2\kappa^2 \delta_0)<\frac{2\kappa^4\de_0}{\sqrt{t^*}},
	\end{align}
	since $2\kappa^2 \delta_0<\frac{2-\sqrt{2}}{2}$.
	This contradicts with the assumption that equality holds in \eqref{estimate:parabolic-1} with $k=1$ at $t=t^*$. Hence, by density argument we conclude that \eqref{estimate:parabolic-1} holds for BV initial data as well.\\
	Next we show the estimate \eqref{estimate:parabolic-1} with $k=2$. We proceed again by contradiction, we can assume again by continuity provided that we consider smooth initial data that there exists a time $t^*<\hat{t}$ such that \eqref{estimate:parabolic-1} holds for $t\in[0,t^*]$ with $k=2$ and there is equality  at time $t=t^*$.
	Again, from \eqref{eqn-tilde-v} we get
	\begin{align*}
		&\tilde{v}_{xx}=G_x(t/2)\star \tilde{v}_{x}(t/2)+\int\limits_{t/2}^{t}G_x(t-s)\star \Big\{ (\mathcal{T}_s^{-1}([A_2^*-A^*_1B_1^{-1/2}(u)]\mathcal{T}_s(\tilde{v}_x)))_x\\
		&+\big(\mathcal{T}_s^{-1}(B_1(u)(P^{-1}(u) v)\bullet B_1^{-1/2}(u)\mathcal{T}_s(\tilde{v}_{x}))\big)_x-\big(\mathcal{T}_s^{-1} (M(X_s)\mathcal{T}_s( \tilde{v}_x))\big)_x+ \big(\mathcal{T}_s^{-1}(\mathcal{R})\big)_x\Big\}(s,\cdot)\,ds.
	\end{align*}
	We obtain from Lemma \ref{lemma:transformation}
	\begin{align}
		\norm{\tilde{v}_{xx}(t)}_{L^1}
		&\leq \norm{G_x(t/2)}_{L^1}\norm{ \tilde{v}_x(t/2)}_{L^1}+\int\limits_{t/2}^{t}\norm{G_x(t-s)}_{L^1} \Big\{\norm{([A_2^*-A^*_1B_1^{-1/2}(u)]\mathcal{T}_s(\tilde{v}_{x}))_x}_{L^1}\nonumber\\
		&+\norm{[(B_1(u) u_x\bullet B_1^{-1/2})(\mathcal{T}_s (\tilde{v}_{x}))]_x}_{L^1}+\norm{\big( M(X_s)\mathcal{T}_s( \tilde{v}_x)\big)_x}_{L^1} +\norm{(\mathcal{R})_x}_{L^1}\Big\}(s,\cdot)\,ds.
		\label{estimulcle1}
	\end{align}
	By using Lemma \ref{lemma:transformation} we have again
	\begin{align}
		\norm{\tilde{v}_{xx}}_{L^1}&\leq \kappa\norm{v_{xx}}_{L^1}+\kappa^3 \norm{v}_{L^\f}\norm{\tilde{v}_x}_{L^1}\leq \kappa\norm{v_{xx}}_{L^1}+\kappa^3 \norm{v}_{L^\f}\norm{v_x}_{L^1}\nonumber\\
		&\leq \kappa^2(\norm{u_{xxx}}_{L^1}+3\norm{u_x}_{L^\f}\norm{u_{xx}}_{L^1}+\norm{u_x}_{L^\f}^2\norm{u_x}_{L^1})\nonumber\\
		&+\kappa^5 \norm{u_x}_{L^\f}(\norm{u_{xx}}_{L^1}+\norm{u_{x}}_{L^\f}\norm{u_x}_{L^1})\nonumber\\
		&\leq \kappa^2\norm{u_{xxx}}_{L^1}+6\kappa^5\norm{u_{xx}}_{L^1}^2.\label{tilde-v-xx-and-u-xxxp}
	\end{align}
	Similarly, we have from Lemma \ref{lemma:transformation},  \eqref{estimate:parabolic-1} with $k=1$ for $t\in]0,t^*]$
	\begin{equation}\label{u-xxx-tilde-v-xx}
		\norm{u_{xxx}(t,\cdot)}_{L^1}\leq \kappa^2\norm{\tilde{v}_{xx}(t,\cdot)}_{L^1}+6\kappa^5\norm{u_{xx}(t,\cdot)}_{L^1}^2\leq \kappa^2\norm{\tilde{v}_{xx}(t,\cdot)}_{L^1}+\frac{24\kappa^{13}\de_0^2}{t}.
	\end{equation}
	We note that
	\begin{align}
		X_{i,tx}=\frac{-1}{2(\mu_i(u))^{3/2}}\left[(u_x\bullet B(u)u_x+B(u)u_{xx})\cdot D\mu_i(u)-(A(u)u_x)\cdot D\mu_i(u)\right].
		\label{D20}
	\end{align}
	which implies,
	\begin{equation}
		\norm{X_{i,tx}}_{L^\f}\leq 4\kappa^9(\norm{u_{xxx}}_{L^1}+\norm{u_{xx}}_{L^1}+\norm{u_{xx}}_{L^1}^2).
		\label{D21}
	\end{equation}
	Now we have using Lemma  \ref{lemma:transformation},\eqref{tilde-v-xx-and-u-xxxp} \eqref{D10}, \eqref{eqtech} for $s\in[0,t^*]$
	\begin{align}
		&\norm{\pa_x([A_2^*-A^*_1B_1^{-1/2}(u)]\mathcal{T}_s(\tilde{v}_x))}_{L^1}
		\leq \kappa\delta_0\|\tilde{v}_{xx}\|_{L^1}+\kappa^2\|u_{xx}\|_{L^1}\|\tilde{v}_x\|_{L^1}\nonumber\\
		&\leq  \kappa^3\delta_0\norm{u_{xxx}}_{L^1}+12\kappa^6\norm{u_{xx}}_{L^1}^2,\nonumber\\
		&\norm{\pa_x(B_1(P^{-1}{v})\bullet B_1^{-1/2}\mathcal{T}_s(\tilde{v}_{x}))(s,\cdot)}_{L^1}
		\leq 10 \kappa^8\norm{u_{xx}}_{L^1}\norm{u_{xxx}}_{L^1}+60\kappa^{11}\norm{u_{xx}}_{L^1}^3,\nonumber\\
		&\norm{\pa_x(M(X_s)\mathcal{T}_s( \tilde{v}_x)(s,\cdot))}_{L^1}
		\leq 10^2\kappa^{13}(\norm{u_{xx}}_{L^1}^2+\norm{u_{xx}}_{L^1}\norm{u_{xxx}}_{L^1})\nonumber\\
		&\hspace{6cm}+10^2\kappa^{14}(\norm{u_{xx}}_{L^1}^3+\de_0\norm{u_{xxx}}_{L^1}).\label{estimulcle2}
	\end{align}
	From \eqref{absolucle1}, we have
	\begin{align}
		&\norm{\pa_x(\mathcal{R})}_{L^1}\leq  20\cdot 2^{10}\kappa^9\Big(\norm{u^*-u}_{L^\f}\norm{v_{xx}}_{L^1}+\norm{v}_{L^\f}\norm{v_{x}}_{L^1}+\norm{v}_{L^\f}^3\norm{v}_{L^1}\nonumber\\
		&+{\color{black}{\norm{v}_{L^\f}^2\norm{v}_{L^1}}}+\norm{v}_{L^\f}^2\norm{v_{x}}_{L^1}+\norm{v_x}_{L^\f}\norm{v_x}_{L^1}+\norm{v}_{L^\f}\norm{v_{xx}}_{L^1}\Big)\nonumber\\
		&\leq  20\cdot 2^{10}\kappa^{14}\left(\de_0\norm{u_{xxx}}_{L^1}+{\color{black}{8}}\de_0\norm{u_{xx}}_{L^1}^2+4\norm{u_{xx}}_{L^1}^2+16\norm{u_{xx}}_{L^1}^3\right)\nonumber\\
		&+  20\cdot 2^{12}\kappa^{12}\left(\norm{u_{xx}}_{L^1}\norm{u_{xxx}}_{L^1}+4\norm{u_{xx}}_{L^1}^3\right)\nonumber\\
		&\leq  20\cdot 2^{18}\kappa^{14}\left((\norm{u_{xx}}_{L^1}+\de_0)\norm{u_{xxx}}_{L^1}+\norm{u_{xx}}_{L^1}^2+\norm{u_{xx}}_{L^1}^3\right).
		\label{estimulcle3}
	\end{align}
	From \eqref{estimucle4} we know that for  $t\in]0,t^*]$ well have
	%\eqref{D10} and \eqref{estimate:parabolic-1} with $k=1$, we deduce that for $t\in]0,t^*]$ we have
	$\norm{\tilde{v}_x(t,\cdot)}_{L^1}\leq \frac{(2+\sqrt{2} )\kappa^3 \de_0}{2\sqrt{t}}$. 
	Combining \eqref{estimulcle1}, \eqref{estimulcle2}, \eqref{estimulcle3} and \eqref{estimate:parabolic-1} with $k=1,2$, it follows that for any $t\in]0,t^*]$ we have
	\begin{align}
		&\norm{\tilde{v}_{xx}(t)}_{L^1}
		\leq \frac{4\kappa^4 \de_0}{t}+20\cdot2^{26}\kappa^{26}\int\limits_{t/2}^{t}\frac{1}{\sqrt{t-s}}\left[\frac{\de_0^2}{s^{3/2}}+\frac{\de_0^2}{s}\right]\,ds\nonumber\\
		&\leq \frac{4\kappa^4\de_0}{t}+\frac{20 \cdot2^{29}\kappa^{26}\delta_0^2}{t}+\frac{20\cdot2^{29}\kappa^{26}\delta_0^2}{\sqrt{t}}\nonumber\\
		&= \frac{4\kappa^4\de_0}{t}\big(1+5\cdot 2^{29}\kappa^{22}\delta_0+5\cdot 2^{29}\kappa^{22}\delta_0\sqrt{t}\big)
		%+\left(20 \cdot2^{21}\kappa^{24}\delta_0+20 \cdot2^{22}\kappa^{24}\delta_0\sqrt{t}\right)\frac{\de_0}{2t}\\
		\leq  \frac{4(2+\sqrt{2})\kappa^4\de_0}{2 t},\label{superimpogh4}
	\end{align}
	since $5\cdot 2^{29}\kappa^{22}\delta_0+5\cdot 2^{29}\kappa^{22}\delta_0\sqrt{t}<\frac{1}{\sqrt{2}}$. Therefore, from \eqref{u-xxx-tilde-v-xx} 
	\begin{equation}
		\norm{u_{xxx}(t^*,\cdot)}_{L^1}\leq \kappa^2\norm{\tilde{v}_{xx}(t^*,\cdot)}_{L^1}+\frac{24\kappa^{13}\de_0^2}{t^*}
		<\frac{8 \kappa^6\de_0}{t^*},
	\end{equation}
	since $24\kappa^7\delta_0<4-\frac{4}{\sqrt{2}}$. It gives a contradiction and it implies the estimate
	\eqref{estimate:parabolic-1} for $k=2$.\\ 
	Next we show the estimate \eqref{estimate:parabolic-1} with $k=3$ by proceeding again by contradiction, with $t^*$ defined as previously for $k=3$.  Again, from \eqref{eqn-tilde-v} we obtain
	\begin{align*}
		&\tilde{v}_{xxx}=G_x(t/2)\star \tilde{v}_{xx}(t/2)+\int\limits_{t/2}^{t}G_x(t-s)\star \Big\{ (\mathcal{T}_s^{-1}([A_2^*-A^*_1B_1^{-1/2}(u)]\mathcal{T}_s(\tilde{v}_x)))_{xx}+\\
		&\big(\mathcal{T}_s^{-1}(B_1(u)(P^{-1}(u) v)\bullet B_1^{-1/2}(u)\mathcal{T}_s(\tilde{v}_{x}))\big)_{xx}-\big(\mathcal{T}_s^{-1} (M(X_s)\mathcal{T}_s( \tilde{v}_x))\big)_{xx}+ \big(\mathcal{T}_s^{-1}(\mathcal{R})\big)_{xx}\Big\}(s,\cdot)\,ds.
	\end{align*}
	Now, from Lemma \ref{lemma:transformation} it gives
	\begin{align}
		&\norm{\tilde{v}_{xxx}(t)}_{L^1}
		\leq \norm{G_x(t/2)}_{L^1}\norm{ \tilde{v}_{xx}(t/2)}_{L^1}+\kappa\int\limits_{t/2}^{t}\norm{G_x(t-s)}_{L^1} \Big\{\norm{([A_2^*-A^*_1B_1^{-1/2}(u)]\tilde{v}_{x})_{xx}}_{L^1}\nonumber\\
		&+\norm{([B_1(u_x\bullet B_1^{-1/2})]\mathcal{T}_s(\tilde{v}_{x}))_{xx}}_{L^1}+\norm{\big( M(X_s)\mathcal{T}_s( \tilde{v}_x)\big)_{xx}}_{L^1} +\norm{(\mathcal{R})_{xx}}_{L^1}\nonumber\\
		&+\kappa\big(\norm{u_{xx}}_{L^1}\norm{([A_2^*-A^*_1B_1^{-1/2}(u)]\tilde{v}_{x})_{x}}_{L^1}+\norm{u_{xx}}_{L^1}\norm{([B_1(u_x\bullet B_1^{-1/2})]\mathcal{T}_s(\tilde{v}_{x}))_{x}}_{L^1}\nonumber\\
		&+\norm{u_{xx}}_{L^1}\norm{\big( M(X_s)\mathcal{T}_s( \tilde{v}_x)\big)_{x}}_{L^1} +\norm{u_{xx}}_{L^1}\norm{(\mathcal{R})_{x}}_{L^1}\big)\Big\}(s,\cdot)\,ds.\label{tilde-v-xxx}
	\end{align}
	Furthermore, we get using the fact $v_x(t,x)=B_1^{-1/2}(u)\mathcal{T}_t(v(t,\cdot))(x)$
	\begin{align}
		v_{xxx}(t,x)&=B_1^{-3/2}(u(t,x))\mathcal{T}_t\tilde{v}_{xxx}(t,x)+u_x\cdot DB_1^{-1}(u(t,x))\mathcal{T}_t\tilde{v}_{xx}(t,x)\nonumber\\
		&+(u_{xx}(t,x))\bullet B_1(u(t,x))^{-1/2}\mathcal{T}_t(\tilde{v}_{x})(t,x)\nonumber\\
		&+(u_{x}(t,x)\otimes u_x(t,x)):D^2B_1(u(t,x))^{-1/2}\mathcal{T}_t(\tilde{v}_{x})(t,x)\nonumber\\
		&+u_{x}(t,x)\bullet B_1(u(t,x))^{-1/2}[B_1(u(t,x))^{-1/2}\mathcal{T}_t(\tilde{v}_{xx})(t,x)].\label{vxxx}
	\end{align}	
	Therefore from Lemma \ref{lemma:transformation} we have
	\begin{align}
		&\norm{\tilde{v}_{xxx}}_{L^1}\leq \kappa^3\norm{v_{xxx}}_{L^1}+4\kappa^6\norm{u_{xx}}_{L^1}\norm{\tilde{v}_{xx}}_{L^1}\label{D26}\\
		&\norm{v_{xxx}}_{L^1}\leq \kappa^4\norm{\tilde{v}_{xxx}}_{L^1}+4\kappa^3 \norm{u_{xx}}_{L^1}\norm{\tilde{v}_{xx}}_{L^1}\label{D27}\\
		&\norm{v_{xxx}}_{L^1}\leq \kappa\norm{u_{xxxx}}_{L^1}+7\kappa\norm{u_{xx}}_{L^1}\norm{u_{xxx}}_{L^1}+7\kappa\|u_{xx}\|_{L^1}^3,\label{D28}
	\end{align}
	It yields using \eqref{D27}, \eqref{tilde-v-xx-and-u-xxxp} and the fact that $\|v_{xx}\|_{L^1}\leq 4\kappa\|u_{xx}\|_{L^1}^2+\kappa\|u_{xxx}\|_{L^1}$
	\begin{align}
		&\norm{u_{xxxx}}_{L^1}\leq \kappa\big(\|v_{xxx}\|_{L^1}+6\|u_{xx}\|_{L^1}\|v_{xx}\|_{L^1}+3\|u_{xx}\|_{L^1}^2\|v_x\|_{L^1}+\delta_0\kappa \|u_{xx}\|_{L^1}^3\nonumber\\
		&+3\delta_0\kappa\|u_{xx}\|_{L^1}^3+\kappa\|u_{xxx}\|_{L^1}\|u_{xx}\|_{L^1}\big) 
		\nonumber\\
		&\leq \kappa\norm{v_{xxx}}_{L^1}+7\kappa^2\norm{u_{xx}}_{L^1}\norm{u_{xxx}}_{L^1}+34\kappa^2\norm{u_{xx}}_{L^1}^3 \nonumber\\
		%&\leq \kappa^4\norm{\tilde{v}_{xxx}}_{L^1}+2\de_0\kappa^5\norm{u_{xx}}_{L^1}\norm{\tilde{v}_{xx}}_{L^1}.      
		&\leq \kappa^5 \norm{\tilde{v}_{xxx}}_{L^1}+11\kappa^6\norm{u_{xx}}_{L^1}\norm{u_{xxx}}_{L^1}+58\kappa^9\norm{u_{xx}}_{L^1}^3.
		\label{u-xxxx-tilde-v-xxx}
	\end{align}
	Combining \eqref{D26}, \eqref{D28}, \eqref{tilde-v-xx-and-u-xxxp} we get
	\begin{align}
		\norm{\tilde{v}_{xxx}}_{L^1}&\leq \kappa^4\norm{u_{xxxx}}_{L^1}+11 \kappa^8\norm{u_{xxx}}_{L^1}\norm{u_{xx}}_{L^1}+31\kappa^{11}\norm{u_{xx}}_{L^1}^3.
		\label{superimpol46}%
		%&+8\kappa^7\de_0^2\norm{u_{xx}}_{L^1}^2\\
		%&\leq \kappa^3\norm{u_{xxxx}}_{L^1}+\left(10^2\kappa^7\cdot 4\cdot 2+8\kappa^7\cdot 2^3\de_0^2+8\kappa^7\cdot 2\de_0^2\sqrt{t}\right)\frac{\de_0^2}{t^{3/2}}\\
		%&\leq \kappa^3\norm{u_{xxxx}}_{L^1}+\frac{10^3\kappa^7\de_0^2}{t^{3/2}}.
	\end{align}
	Moreover, by using \eqref{D20} and the computations as in Appendix \ref{appendix:R-derivatives}, we get
	\begin{align}
		\norm{X_{i,txx}}_{L^1}
		&\leq \kappa^8\big(4\norm{u_{xx}}_{L^1}^2+\norm{u_{xxx}}_{L^1}+3\norm{u_{xx}}_{L^1}\big).\label{Xitechutile}
		%&\leq 2^{9}\kappa^{14}\left(\frac{\de_0}{s}+\frac{\de_0}{\sqrt{s}}\right).
	\end{align}
	Now using Lemma \ref{lemma:transformation},  the estimates \eqref{estimate:parabolic-1} with $k=1,2$,
	\eqref{superimpol46}, \eqref{tilde-v-xx-and-u-xxxp}, \eqref{D10}, \eqref{Xitechutile}, \eqref{eqtech}, \eqref{D21} we have for $s\in[0,t^*]$
	\begin{align}
		&\norm{\pa_{xx}([A_2^*-A^*_1B_1^{-1/2}(u)]\mathcal{T}_s(\tilde{v}_x))(s,\cdot)}_{L^1}
		\leq\big( \kappa^2\delta_0\norm{\tilde{v}_{xxx}}_{L^1}+3\kappa\norm{u_{xx}}_{L^1}\norm{\tilde{v}_{xx}}_{L^1},\nonumber\\
		&+\kappa^2\norm{u_{xxx}}_{L^1}\norm{\tilde{v}_x}_{L^1}+\kappa^2\norm{u_{xx}}_{L^1}^2\norm{\tilde{v}_x}_{L^1}\big)(s,\cdot),\nonumber\\
		&\leq \big(\kappa^6\de_0\norm{u_{xxxx}}_{L^1}+
		6\kappa^{10}\norm{u_{xx}}_{L^1}\norm{u_{xxx}}_{L^1}+21\kappa^6\norm{u_{xx}}_{L^1}^3)\big)(s,\cdot)\nonumber\\
		&\leq \kappa^6\de_0\norm{u_{xxxx}}_{L^1}(s,\cdot)+\frac{3 \cdot 2^5\kappa^{18}\de_0^2}{s^{3/2}},\nonumber\\[1,5mm]
		%\\
		&\norm{\pa_{xx}(B_1(u_x)\bullet B_1^{-1/2}\mathcal{T}_s(\tilde{v}_{x}))(s,\cdot)}_{L^1}
		\leq\big(  \kappa^3(4\norm{u_{xx}}^3_{L^1}+6\norm{u_{xxx}}_{L^1}\norm{u_{xx}}_{L^1}\nonumber\\
		&+\norm{u_{xxxx}}_{L^1}
		)\norm{\tilde{v}_x}_{L^1}+\kappa^2(\norm{u_{xxx}}_{L^1}+5\kappa^3\norm{u_{xx}}_{L^1}^2)\norm{\tilde{v}_{xx}}_{L^1}+\kappa^3\norm{u_{xx}}_{L^1}\norm{\tilde{v}_{xxx}}_{L^1}\big)(s,\cdot),\nonumber\\
		&\leq \big(70\kappa^{14}\norm{u_{xx}}_{L^1}^4+28\kappa^{11}\|u_{xxx}\|_{L^1}\|u_{xx}\|_{L^1}^2+2\kappa^7\|u_{xxxx}\|_{L^1}\|u_{xx}\|_{L^1}+\kappa^4\|u_{xxx}\|_{L^1}^2\big)(s,\cdot),\nonumber\\
		&\leq 4\kappa^{11}\|u_{xxxx}(s,\cdot)\|_{L^1}\frac{\delta_0}{\sqrt{s}}+\frac{2^{13}\kappa^{20}\delta_0^2}{s^2}\nonumber\\[1,5mm]
		%&\leq 8\kappa^9\norm{u_{xxxx}}_{L^1}+\frac{2^{13}\kappa^{13}\de_0^2}{s^2},\\
		%\\
		&\norm{\pa_{xx}(M(X_s)\mathcal{T}_s( \tilde{v}_x)(s,\cdot))}_{L^1}
		\leq\big( \norm{X_{sxx}}_{L^1}\norm{\tilde{v}_{xx}}_{L^1}+2\norm{X_{sx}}_{L^\f}\norm{\tilde{v}_{xx}}_{L^1}\nonumber\\
		&+\norm{X_{s}}_{L^\infty}(\kappa \norm{\tilde{v}_{xxx}}_{L^1}+\|u_{xx}\|_{L^1}\|\tilde{v}_{xx}\|_{L^1})\big(s,\cdot)\nonumber\\
		&\leq 8\kappa^{13}(\kappa\|u_{xx}(s,\cdot)\|_{L^1}+\delta_0)\|u_{xxxx}(s,\cdot)\|_{L^1}+
		2^{14}\frac{\kappa^{23}\delta_0^2}{s^{3/2}}+2^{14}\frac{\kappa^{27}\delta_0^2}{s^2}.\label{estimysupercru}
		%\left(\frac{\de_0^2}{s^2}+\frac{\de_0}{s^{3/2}}\right)+2^4\kappa^{12}\left(\frac{\de_0}{s}+\frac{\de_0}{\sqrt{s}}\right)\norm{u_{xxxx}}_{L^1}.
	\end{align}
	By similar arguments, we have using \eqref{absolucle1}  and the estimates \eqref{estimate:parabolic-1} with $k=1,2$ for $s\in[0,t^*]$
	\begin{align}
		&\norm{	\pa_{xx}\mathcal{R}(u,v,v_x)}_{L^1}\leq  10^2\cdot 2^{12}\kappa^{10}\Big(\norm{u^*-u}_{L^\f}\norm{v_{xxx}}_{L^1}+\norm{v}^2_{L^\f}\norm{v_{x}}_{L^1}+\norm{v}_{L^\f}\norm{v_{xx}}_{L^1}\nonumber\\
		&+{\color{black}{\norm{v}_{L^\infty}^3\norm{v}_{L^1}+\norm{v_x}_{L^\infty}\norm{v_x}_{L^1}}}+\norm{v}_{L^\f}^4\norm{v}_{L^1}+\norm{v}_{L^\f}^3\norm{v_x}_{L^1}+\norm{v}_{L^\f}\norm{v_{x}}_{L^\f}\norm{v_x}_{L^1}\nonumber\\
		&\hspace*{2cm}+\norm{v}^2_{L^\f}\norm{v_{xx}}_{L^1}+\norm{v_x}_{L^\f}\norm{v_{xx}}_{L^1}+\norm{v}_{L^\f}\norm{v_{xxx}}_{L^1}\Big),\nonumber\\
		&\leq  10^2\cdot 2^{12}\kappa^{10}\big(\kappa(\delta_0+\kappa\|u_{xx}(s,\cdot)\|_{L^1})\|u_{xxxx}(s,\cdot)\|_{L^1}+
		2^{16}\frac{\kappa^{26}\delta_0^2}{s^{3/2}}+2^{15}\frac{\kappa^{28}\delta_0^2}{s^2}\big).
		\label{estimysupercru1}
	\end{align}
	From \eqref{superimpogh4}, we have
	%\eqref{tilde-v-xx-and-u-xxx}, the estimates \eqref{estimate:parabolic-1} with $k=1,2$ and the smallness assumption on $\delta_0$, we note that for $t\in[0,t^*]$
	\begin{equation}
		\norm{\tilde{v}_{xx}(t,\cdot)}_{L^1}\leq %\frac{9 \kappa^6\de_0}{t}.
		\frac{4(2+\sqrt{2})\kappa^4\de_0}{2 t}
		\label{initialfg}
	\end{equation} 
	Therefore, from \eqref{tilde-v-xxx}, \eqref{estimysupercru}, \eqref{estimysupercru1}, \eqref{estimulcle2}, \eqref{estimulcle3}  and the estimates \eqref{estimate:parabolic-1} with $k=1,2,3$ (when $t\in]0,t^*]$) it follows that for $t\in]0,t^*]$ we have
	\begin{align}
		\norm{\tilde{v}_{xxx}(t,\cdot)}_{L^1}
		&\leq \frac{4(2+\sqrt{2})\kappa^5\de_0}{ t^{3/2}}
		%\frac{9\sqrt{2}\kappa^7\de_0}{t^{3/2}}
		+2^{32}\kappa^{30}\int\limits_{t/2}^{t}\frac{1}{\sqrt{t-s}}\left[\frac{\de_0^2}{s^{3/2}}+\frac{\de_0^2}{s^2}\right]\,ds\nonumber\\
		&\leq  \frac{2^4\kappa^5\de_0}{t^{3/2}}+\left(2^{34}\kappa^{25}\sqrt{2}\delta_0+2^{34}\kappa^{25}\delta_0\sqrt{t}\right)\frac{\de_0\kappa^5}{t^{3/2}}\leq \frac{(2^{5}-1)\kappa^5\de_0}{t^{3/2}},\label{superumesti}
	\end{align}
	since $\left(2^{34}\kappa^{25}\sqrt{2}\delta_0+2^{35}\kappa^{23}\delta_0\sqrt{t}\right)<1$ due to the assumptions on $\hat{t}$ and $\delta_0$. From \eqref{u-xxxx-tilde-v-xxx}, \eqref{superumesti}  and the estimates \eqref{estimate:parabolic-1} with $k=1,2$ (when $t\in]0,t^*]$) we have for $t\in]0,t^*]$
	\begin{align}
		&\norm{u_{xxxx}(t,\cdot)}\leq\kappa^5 \norm{\tilde{v}_{xxx}(t,\cdot)}_{L^1}+11 \cdot 2^4\kappa^{16}\frac{\delta_0^2}{t^{3/2}}+58\kappa^{21}2^3 \frac{\delta_0^3}{t^{3/2}}%+12\delta_0^{3}\frac{\kappa^{10}}{t}
		\nonumber\\
		&\leq \frac{(2^{5}-1)\kappa^{10}\de_0}{t^{3/2}}+11 \cdot 2^4\kappa^{16}\frac{\delta_0^2}{t^{3/2}}+58\kappa^{21}2^3 \frac{\delta_0^3}{t^{3/2}}%+12\delta_0^{3}\frac{\kappa^{10}}{t}
		\leq \frac{(2^{5}-1/2)\kappa^{10}\de_0}{t^{3/2}},\nonumber%\|u_{xx}\|_{L^1}^2.
		%\kappa^4\norm{\tilde{v}_{xxx}}_{L^1}+\frac{2^6\kappa^{16}\de_0^2}{t^{3/2}}.
	\end{align}
	since $11 \cdot 2^4\kappa^{6}\delta_0+58\kappa^{11}2^3\delta_0^2 <\frac{1}{2}$.
	Subsequently  we obtain
	\begin{equation}
		\norm{u_{xxxx}(t^*,\cdot)}_{L^1}<\frac{2^5\kappa^{10}\de_0}{t^{3/2}},
	\end{equation}
	what is a contradiction. It implies that we have obtained the estimate  \eqref{estimate:parabolic-1} for $k=3$.
	Finally, we show the estimate \eqref{estimate:parabolic-1} with $k=4$ by proceeding again by contradiction with $t^*$ defined in the same way.  Again, from \eqref{eqn-tilde-v} we get
	\begin{align*}
		\tilde{v}_{xxxx}&=G_x(t/2)\star \tilde{v}_{xxx}(t/2)+\int\limits_{t/2}^{t}G_x(t-s)\star \Big\{ (\mathcal{T}_s^{-1}([A_2^*-A^*_1B_1^{-1/2}(u)]\mathcal{T}_s(\tilde{v}_x)))_{xxx}\\
		&\hspace{4.5cm}+\big(\mathcal{T}_s^{-1}(B_1(u)(P^{-1}(u) v)\bullet B_1^{-1/2}(u)\mathcal{T}_s(\tilde{v}_{x}))\big)_{xxx}\\
		&\hspace{4.5cm}-\big(\mathcal{T}_s^{-1} (M(X_s)\mathcal{T}_s( \tilde{v}_x))\big)_{xxx}+ \big(\mathcal{T}_s^{-1}(\mathcal{R})\big)_{xxx}\Big\}(s,\cdot)\,ds.
	\end{align*}
	Then it follows from Lemma \ref{lemma:transformation},
	\begin{align}
		&\norm{\tilde{v}_{xxxx}(t)}_{L^1}\leq \norm{G_x(t/2)}_{L^1}\norm{ \tilde{v}_{xxx}(t/2)}_{L^1}+\kappa^2 \int\limits_{t/2}^{t}\norm{G_x(t-s)}_{L^1}\\
		& \Big\{\norm{([A_2^*-A^*_1B_1^{-1/2}(u)]\mathcal{T}_s(\tilde{v}_{x}))_{xxx}}_{L^1}+\norm{([B_1(u_x\bullet B_1^{-1/2})]\mathcal{T}_s(\tilde{v}_{x}))_{xxx}}_{L^1}\nonumber\\
		&+\norm{\big( M(X_s)\mathcal{T}_s( \tilde{v}_x)\big)_{xxx}}_{L^1} +\norm{(\mathcal{R})_{xxx}}_{L^1}+\kappa \norm{u_{x}}_{L^\f}\Big(\norm{([A_2^*-A^*_1B_1^{-1/2}(u)]\tilde{v}_{x})_{xx}}_{L^1}\nonumber\\
		&+\norm{([B_1(u_x\bullet B_1^{-1/2})]\mathcal{T}_s(\tilde{v}_{x}))_{xx}}_{L^1}+\norm{\big( M(X_s)\mathcal{T}_s( \tilde{v}_x)\big)_{xx}}_{L^1} +\norm{(\mathcal{R})_{xx}}_{L^1}\Big)\nonumber\\
		&+(3\kappa^3 \norm{u_{x}}^2_{L^\f}+\norm{u_{xx}}_{L^\f})\Big(\norm{([A_2^*-A^*_1B_1^{-1/2}(u)]\tilde{v}_{x})_{x}}_{L^1}+\norm{([B_1(u_x\bullet B_1^{-1/2})]\mathcal{T}_s(\tilde{v}_{x}))_{x}}_{L^1}\nonumber\\
		&+\norm{\big( M(X_s)\mathcal{T}_s( \tilde{v}_x)\big)_{x}}_{L^1} +\norm{(\mathcal{R})_{x}}_{L^1}\Big)\Big\}(s,\cdot)\,ds.\label{tilde-v-xxxx}
	\end{align}
	We recall now from \eqref{estimucle4}, \eqref{superimpogh4}, \eqref{superumesti} that for $s\in]0,t^*]$ we have
	\begin{align}
		&\|\tilde{v}_x(s,\cdot)\|_{L^1}\leq \frac{(2+\sqrt{2} )\kappa^3 \de_0}{2\sqrt{s}}, \;\|\tilde{v}_{xx}(s,\cdot)\|_{L^1}\leq \frac{4(2+\sqrt{2})\kappa^4\de_0}{2 s},\;\|\tilde{v}_{xxx}(s,\cdot)\|_{L^1}\leq \frac{(2^{5}-1)\kappa^5\de_0}{s^{3/2}}
		\label{superimpml1}
	\end{align}
	Now, by deriving \eqref{vxxx} we obtain
	\begin{align}
		v_{xxxx}(t,x)&=B_1^{-2}(u(t,x))\mathcal{T}_t\tilde{v}_{xxxx}(t,x)+u_x\cdot DB_1^{-3/2}(u(t,x))\mathcal{T}_t\tilde{v}_{xxx}(t,x)\nonumber\\
		&+u_{xx}\cdot DB_1^{-1}(u(t,x))\mathcal{T}_t\tilde{v}_{xx}(t,x)+(u_x\otimes u_{x}):D^2B_1^{-1}(u(t,x))\mathcal{T}_t\tilde{v}_{xx}(t,x)\nonumber\\
		&+u_{x}\cdot DB_1^{-1}(u(t,x))B_1^{-1/2}\mathcal{T}_t\tilde{v}_{xxx}(t,x)+(u_{xxx}(t,x))\bullet B_1(u(t,x))^{-1/2}\mathcal{T}_t(\tilde{v}_{x})(t,x)\nonumber\\
		&+(u_{xx}\otimes u_x):D^2B_1^{-1/2}\mathcal{T}_t(\tilde{v}_{x})(t,x)+(u_{xx}(t,x))\bullet B_1(u(t,x))^{-1/2}B_1^{-1/2}\mathcal{T}_t(\tilde{v}_{xx})(t,x)\nonumber\\
		&+2(u_{xx}\otimes u_x):D^2B_1^{-1/2}\mathcal{T}_t(\tilde{v}_{x})(t,x)+(u_{x}\otimes u_x\otimes u_x):D^3B_1^{-1/2}\mathcal{T}_t(\tilde{v}_{x})(t,x)\nonumber\\
		&+(u_{x}\otimes u_x):D^2B_1^{-1/2}B_1^{-1/2}\mathcal{T}_t(\tilde{v}_{xx})(t,x)+(u_x\otimes u_{x}): D^2 B_1^{-1/2}[B_1^{-1/2}\mathcal{T}_t(\tilde{v}_{xx})(t,x)]\nonumber\\
		&+u_{xx}\bullet B_1^{-1/2}[B_1^{-1/2}\mathcal{T}_t(\tilde{v}_{xx})(t,x)]+u_{x}\bullet B_1^{-1/2}[u_x\cdot DB_1^{-1/2}\mathcal{T}_t(\tilde{v}_{xx})(t,x)]\nonumber\\
		&+u_{x}(t,x)\bullet B_1(u(t,x))^{-1/2}[B_1^{-1}\mathcal{T}_t(\tilde{v}_{xxx})(t,x)].
		\label{superitych}
	\end{align}
	Then we have applying Lemma  \ref{lemma:transformation}, \eqref{superitych}, \eqref{superimpml1}, \eqref{estimate:parabolic-1} with $k=1,2,3$ for $s\in]0,t^*]$ and using the smallness of $\delta_0$
	\begin{align}
		&\norm{\tilde{v}_{xxxx}(s,\cdot)}_{L^1}\leq \kappa^3\big( \norm{v_{xxxx}}_{L^1}+
		4\kappa^3\|u_{xx}\|^2_{L^1}\|\tilde{v}_{xx}\|_{L^1}+3\kappa^3\|u_{xx}\|_{L^1}\|\tilde{v}_{xxx}\|_{L^1}\nonumber\\
		&+3\kappa^3 \|u_{xxx}\|_{L^1}\|\tilde{v}_{xx}\|_{L^1}+\kappa^2\|\tilde{v}_x\|_{L^1} (\|u_{xx}\|_{L^1}^3 +\|u_{xxxx}\|_{L^1}+3\|u_{xxx}\|_{L^1}\|u_{xx}\|_{L^1})\big)(s,\cdot)\nonumber\\
		&\leq  \kappa^3\norm{v_{xxxx}(s,\cdot)}_{L^1}+\frac{2^5\kappa^{18}\delta_0^2}{s^2}.
		\label{superimpot11}
		%5\kappa^5\norm{u_{xx}}_{L^1}\norm{\tilde{v}_{xxx}}_{L^1}+6\kappa^5 \norm{u_{xxx}}_{L^1}\norm{\tilde{v}_{xx}}_{L^1}\\
		%&+4\kappa^5 \norm{u_{xx}}^2_{L^1}\norm{\tilde{v}_{xx}}_{L^1}.
	\end{align}	
	By a similar argument, it follows using also \eqref{tilde-v-xx-and-u-xxxp}, \eqref{D10}, \eqref{superimpol46} and \eqref{estimate:parabolic-1} with $k=1,2,3$ for $s\in]0,t^*]$
	\begin{align}
		\norm{v_{xxxx}(s,\cdot)}_{L^1}&\leq\big( \kappa^3\norm{\tilde{v}_{xxxx}}_{L^1}+
		4\kappa^3\|u_{xx}\|^2_{L^1}\|\tilde{v}_{xx}\|_{L^1}+3\kappa^3\|u_{xx}\|_{L^1}\|\tilde{v}_{xxx}\|_{L^1}\nonumber\\
		&+3\kappa^3 \|u_{xxx}\|_{L^1}\|\tilde{v}_{xx}\|_{L^1}
		+\kappa^2\|\tilde{v}_x\|_{L^1} (\|u_{xx}\|_{L^1}^3 +\|u_{xxxx}\|_{L^1}+3\|u_{xxx}\|_{L^1}\|u_{xx}\|_{L^1})\big)(s,\cdot),\nonumber\\
		&\leq \kappa^3\norm{\tilde{v}_{xxxx}(s,\cdot)}_{L^1}+2^{27}\kappa^{24}\frac{\delta_0^2}{s^2},\nonumber\\
		\label{vsuperit1}
		% \kappa^3 \norm{\tilde{v}_{xxxx}}_{L^1}+5\kappa^5\norm{u_{xx}}_{L^1}\norm{\tilde{v}_{xxx}}_{L^1}+6\kappa^5 \norm{u_{xxx}}_{L^1}\norm{\tilde{v}_{xx}}_{L^1}\\
		%&+4\kappa^5 \norm{u_{xx}}^2_{L^1}\norm{\tilde{v}_{xx}}_{L^1}.
	\end{align}	
	%Recall that
	%\begin{align*}
	%	\norm{\tilde{v}_{x}}_{L^1}\leq \frac{2\kappa\de_0}{t^{1/2}},\quad\norm{\tilde{v}_{xx}}_{L^1}\leq \frac{2^3\kappa^2\de_0}{t},\quad\norm{\tilde{v}_{xxx}}_{L^1}\leq \frac{2^5\kappa^3\de_0}{t^{3/2}}.
	%\end{align*}
	%Consequently, it follows,
	%\begin{align*}
	%	\norm{v_{xxxx}(s,\cdot)}_{L^1}&\leq \kappa^4 \norm{\tilde{v}_{xxxx}}_{L^1}+5\kappa^5\norm{u_{xx}}_{L^1}\norm{\tilde{v}_{xxx}}_{L^1}+6\kappa^5 \norm{u_{xxx}}_{L^1}\norm{\tilde{v}_{xx}}_{L^1}\\
	%	&+4\kappa^5 \norm{u_{xx}}^2_{L^1}\norm{\tilde{v}_{xx}}_{L^1}.
	%\end{align*}	
	Note that applying \eqref{vsuperit1} and \eqref{estimate:parabolic-1} with $k=1,2,3,4$ for $s\in]0,t^*]$
	, we have
	\begin{align}
		&\norm{v_{xxxx}(s,\cdot)}_{L^1}\leq \kappa\big( \norm{u_{xxxxx}}_{L^1}+11\|u_{xx}\|_{L^1}^4+10\|u_{xxx}\|_{L^1}^2+25\|u_{xx}\|_{L^1}^2\|u_{xxx}\|_{L^1}\nonumber\\
		&+5\|u_{xx}\|_{L^1}\|u_{xxxx}\|_{L^1}\big)(s,\cdot)\leq \frac{2^9\kappa^{11}\delta_0}{s^2}.\nonumber\\
		%2^6\kappa\big(\norm{u_{xx}}_{L^1}\norm{u_{xxxx}}_{L^1}+\norm{u_{xx}}^2_{L^1}\norm{u_{xxx}}_{L^1}\big)\\
		%&+2^6\kappa\big( \norm{u_{xxx}}^2_{L^1}+\norm{u_{xx}}_{L^1}^4\big)\\
		%&\leq \frac{2^9\kappa^{16}\de_0}{s^2}+\frac{2^{12}\kappa ^{13}\de_0^2}{s^2},\\
		&\norm{u_{xxxxx}(s,\cdot)}_{L^1}\leq \kappa \|v_{xxxx}(s,\cdot)\|_{L^1}+\kappa^2\big(11\|u_{xx}\|_{L^1}^4+10\|u_{xxx}\|_{L^1}^2+25\|u_{xx}\|_{L^1}^2\|u_{xxx}\|_{L^1}\nonumber\\
		&+5\|u_{xx}\|_{L^1}\|u_{xxxx}\|_{L^1}\big)(s,\cdot)\nonumber\\
		&\leq \kappa^4 \norm{\tilde{v}_{xxxx}(s,\cdot)}_{L^1}+2^{28}\kappa^{25}\frac{\delta_0^2}{s^2}.
		%\leq \kappa\norm{v_{xxxx}}_{L^1}+2^6\kappa\big(\norm{u_{xx}}_{L^1}\norm{u_{xxxx}}_{L^1}+\norm{u_{xx}}^2_{L^1}\norm{u_{xxx}}_{L^1}\big)\\
		%&+2^6\kappa\big( \norm{u_{xxx}}^2_{L^1}+\norm{u_{xx}}_{L^1}^4\big)\\
		%&\leq \kappa\norm{v_{xxxx}}_{L^1}+\frac{2^{12}\kappa ^{13}\de_0^2}{s^2}\\
		\label{supercruclea1}
	\end{align}
	Combining now \eqref{superimpot11} and \eqref{supercruclea1}, we have
	\begin{align}
		&\norm{\tilde{v}_{xxxx}(s,\cdot)}_{L^1}\leq  \kappa^3\norm{v_{xxxx}(s,\cdot)}_{L^1}+\frac{2^5\kappa^{18}\delta_0^2}{s^2}\leq  \frac{2^{10}\kappa^{14}\delta_0}{s^2}.
		\label{D46}
	\end{align}
	%We also have
	%\begin{equation*}
	%	\norm{\tilde{v}_{xxxx}(s,\cdot)}_{L^1}\leq 			\frac{16\kappa^{18}\de_0}{s^2}+\frac{2^{17}\kappa ^{18}\de_0^2}{s^2}.
	%\end{equation*}
	Deriving two times \eqref{D20}, we deduce that
	\begin{align}
		\norm{X_{i,txxx}(s,\cdot)}_{L^1}&\leq 2^3 \kappa^5\big(\norm{u_{xxx}}_{L^1}
		+2^{9}\kappa^3\norm{u_{xxx}}_{L^1}\norm{u_{xx}}_{L^1}+\kappa^{6}2^{9}\norm{u_{xx}}^3_{L^1}\nonumber\\
		&\hspace{7cm}+\kappa^3 2^{8} \norm{u_{xx}}^2_{L^1}\big)(s,\cdot).
		\label{vmtechg1}
		%&\leq \frac{2^{8}\kappa^{18}\de_0}{s^{3/2}}+\frac{2^{7}\kappa^{14}\de_0}{s}.
	\end{align}
	%\begin{equation*}
	%	\norm{X_{i,sx}}_{L^1}\leq 2^4\kappa^{6}\left(\de_0+\frac{\de_0}{\sqrt{s}}\right),
	%	\norm{X_{i,sxx}}_{L^1}\leq 2^{9}\kappa^{16}\left(\frac{\de_0}{s}+\frac{\de_0}{\sqrt{s}}\right).
	%\end{equation*}
	%&= \frac{\sqrt{2}\kappa^3\de_0}{\sqrt{t}}+\left(150\cdot2^{15}\kappa^{10}\delta_0+150\cdot2^{15}\kappa^{10}\delta_0\sqrt{t}\right)\frac{\kappa\de_0}{4\sqrt{t}}\\
	%&<\frac{(2-\sqrt{\de_0})\kappa\de_0}{\sqrt{t}},
	Next applying Lemma \ref{lemma:transformation}, \eqref{vmtechg1}, \eqref{eqtech}, \eqref{superimpml1}, \eqref{D21}, \eqref{Xitechutile}, \eqref{D46}
	and \eqref{estimate:parabolic-1} with $k=1,2,3$, we obtain for $s\in]0,t^*]$
	\begin{align}
		&\norm{\pa_{xxx}([A_2^*-A^*_1B_1^{-1/2}(u)]\mathcal{T}_s(\tilde{v}_x))(s,\cdot)}_{L^1}
		\leq \kappa^3\delta_0\norm{\tilde{v}_{xxxx}(s,\cdot)}_{L^1}+16\kappa^4\big(\norm{u_{xx}}_{L^1}\norm{\tilde{v}_{xxx}}_{L^1}\nonumber\\
		&+\norm{u_{xx}}^2_{L^1}\norm{\tilde{v}_{xx}}_{L^1}+\norm{u_{xxx}}_{L^1}\norm{\tilde{v}_{xx}}_{L^1}+\norm{u_{xx}}_{L^1}^3\norm{\tilde{v}_{x}}_{L^1}\big)(s,\cdot)\leq 	\frac{2^{18}\kappa ^{23}\de_0^2}{s^2},\nonumber\\
		%\\
		&\norm{\pa_{xxx}(B_1(u_x)\bullet B_1^{-1/2}\mathcal{T}_s(\tilde{v}_{x}))(s,\cdot)}_{L^1}
		\leq  2^{10}\kappa^8(\norm{u_{xxxx}}_{L^1}+\norm{u_{xx}}^3_{L^1}\nonumber\\
		&+\norm{u_{xxx}}_{L^1}\norm{u_{xx}}_{L^1})(s,\cdot)\norm{\tilde{v}_{xx}(s,\cdot)}_{L^1}+2^5\kappa^8(\norm{u_{xxx}}_{L^1}+\norm{u_{xx}}_{L^1}^2)(s,\cdot)\norm{\tilde{v}_{xxx}(s,\cdot)}_{L^1}\nonumber\\
		&+\kappa^8\norm{u_{xx}(s,\cdot)}_{L^1}\norm{\tilde{v}_{xxxx}(s,\cdot)}_{L^1}\leq \frac{2^{12}\kappa^{26}\de_0^2}{s^{5/2}},\nonumber\\
		%\\
		&\norm{\pa_{xxx}(M(X_s)\mathcal{T}_s( \tilde{v}_x)(s,\cdot))}_{L^1}
		\leq \kappa^2\norm{X_{s}}_{L^\f}\norm{\tilde{v}_{xxxx}}_{L^1}(s,\cdot)+
		\norm{\tilde{v}_{xxx}}_{L^1}\big(\frac{3}{2}\kappa^4\|u_x\|_{L^\infty}\|X_s\|_{L^\infty}\nonumber\\
		&+3\kappa\|X_{sx}\|_{L^\infty}+3\kappa\|X_{sxx}\|_{L^1}\big)(s,\cdot)+
		\norm{\tilde{v}_{xx}(s,\cdot)}_{L^1}\big(3\kappa^3\|X_{sx}\|_{L^\infty}\|u_x\|_{L^\infty}+\|X_{sxxx}\|_{L^1}\big)(s,\cdot)\nonumber\\
		&\leq 8\kappa^{11}[\kappa \norm{u_{xx}}_{L^1}+\de_0]\norm{\tilde{v}_{xxxx}}_{L^1}(s,\cdot)+
		%&2^4\kappa\norm{X_{sxxx}}_{L^1}\norm{\tilde{v}_{xx}}_{L^1}+2^4\kappa^2\norm{X_{sxx}}_{L^\f}\norm{\tilde{v}_{xx}}_{L^1}\\
		%&+2^4\kappa^2\norm{X_{sx}}_{L^\f}\norm{u_{xx}}_{L^1}\norm{\tilde{v}_{xx}}_{L^1}\\
		%&+2^4\kappa^2\norm{X_{sx}}_{L^\f}\norm{u_{xx}}^2_{L^1}\norm{\tilde{v}_{x}}_{L^1}\\
		%&+2^4\kappa^2\norm{X_{s}}_{L^\f}\norm{u_{xx}}_{L^1}\norm{\tilde{v}_{xx}}_{L^1}\\
		%&+2^4\kappa^2\norm{X_{s}}_{L^\f}\norm{u_{xxx}}_{L^1}\norm{\tilde{v}_{xx}}_{L^1}\\
		%&+2^4\kappa^2\norm{X_{s}}_{L^\f}\norm{\tilde{v}_{xxx}}_{L^1}+2^4\kappa^2\norm{X_{s}}_{L^\f}\norm{\tilde{v}_{xxxx}}_{L^1}\\
		2^{25}\kappa^{38}\left(\frac{\de_0^2}{s^{5/2}}+\frac{\de_0^2}{s^2}\right).\label{D48}
	\end{align}
	By similar arguments, we have using \eqref{absolucle1}, \eqref{supercruclea1}, \eqref{D27}, \eqref{D10}, \eqref{superimpml1}  and the estimates \eqref{estimate:parabolic-1} with $k=1,2,3,4$ for $s\in[0,t^*]$
	\begin{align}
		&\norm{	\pa_{xxx}\mathcal{R}(u,v,v_x)}_{L^1}\leq  10^3\cdot 2^{14}\kappa^{12}\Big(\norm{u^*-u}_{L^\f}\norm{v_{xxxx}}_{L^1}+\norm{v}^3_{L^\f}\norm{v_{x}}_{L^1}
		\nonumber\\
		&+\norm{v}_{L^\f}\norm{v_x}_{L^\f}\norm{v_{x}}_{L^1}+\norm{v}_{L^\f}\norm{v_{xxx}}_{L^1}+\norm{v}_{L^\f}^5\norm{v}_{L^1}+\norm{v}_{L^\f}^4\norm{v_x}_{L^1}\nonumber\\
		&+\norm{v}_{L^\f}^2\norm{v_x}_{L^\f}\norm{v_x}_{L^1}+\norm{v}_{L^\f}^3\norm{v_{xx}}_{L^1}+\norm{v}_{L^\f}\norm{v_x}_{L^\f}\norm{v_{xx}}_{L^1}+\norm{v_x}_{L^\f}^2\norm{v_{x}}_{L^1}\nonumber\\
		&+\norm{v_x}_{L^\f}\norm{v_{xxx}}_{L^1}+\norm{v}_{L^\f}\norm{v_{xxxx}}_{L^1}
		{\color{black}{+\norm{v_{xxx}}_{L^1}\norm{v}_{L^\infty}^2+\norm{v}_{L^\infty}^4\norm{v}_{L^1}
		}}\nonumber\\
		&{\color{black}{+\norm{v_{xx}}_{L^1}\norm{v}_{L^\infty}^2+\norm{v_{xx}}_{L^1}\norm{v_x}_{L^\infty}+\norm{v_{xx}}_{L^1}\norm{v_{xx}}_{L^\infty}     }}\Big)(s,\cdot)\leq\frac{ 2^{35}\kappa^{36}\delta_0^2}{s^2}+\frac{ 2^{35}\kappa^{38}\delta_0^2}{s^{5/2}}.
		\label{D49}
	\end{align}
	%We observe now from \eqref{superimpol46} and \eqref{estimate:parabolic-1} that for $s\in]0,t^*]$
	%we have using the smallness assumption on $\delta_0$:
	%\begin{align}
	%	\norm{\tilde{v}_{xxx}(s,\cdot)}_{L^1}&\leq \big(\kappa^4\norm{u_{xxxx}}_{L^1}+11 \kappa^8\norm{u_{xxx}}_{L^1}\norm{u_{xx}}_{L^1}+31\kappa^{11}\norm{u_{xx}}_{L^1}^3\big)(s,\cdot)\nonumber\\
	%	&\leq \frac{\delta_0}{s^{3/2}}\big(2^5\kappa^{16}+11\cdot 2^4\kappa^{20}\delta_0+31\cdot 2^3 \kappa^{23}\delta_0^2\big)\leq  \frac{2^6\kappa^{16}\delta_0}{s^{3/2}}.
	%
	%&+8\kappa^7\de_0^2\norm{u_{xx}}_{L^1}^2\\
	%&\leq \kappa^3\norm{u_{xxxx}}_{L^1}+\left(10^2\kappa^7\cdot 4\cdot 2+8\kappa^7\cdot 2^3\de_0^2+8\kappa^7\cdot 2\de_0^2\sqrt{t}\right)\frac{\de_0^2}{t^{3/2}}\\
	%&\leq \kappa^3\norm{u_{xxxx}}_{L^1}+\frac{10^3\kappa^7\de_0^2}{t^{3/2}}.
	%	\label{msuperimpoi}
	%\end{align}
	Hence, it follows from \eqref{tilde-v-xxxx}, \eqref{superimpml1}, \eqref{D48}, \eqref{D49}, \eqref{estimysupercru}, \eqref{estimysupercru1}, \eqref{estimulcle2}, \eqref{estimulcle3} for $t\in]0,t^*]$
	\begin{align}
		\norm{\tilde{v}_{xxxx}(t)}_{L^1}
		&\leq  \frac{2^6\sqrt{2}\kappa^{6}\delta_0}{t^2}+2^{42}\kappa^{40}\int\limits_{t/2}^{t}\frac{1}{\sqrt{t-s}}\left[\frac{\de_0^2}{s^{5/2}}+\frac{\de_0^2}{s^2}\right]\,ds\nonumber\\
		&\leq \frac{2^7\kappa^{6}\de_0}{t^2}+2^{44}\kappa^{34}\left(\delta_0+\delta_0\sqrt{t}\right)\frac{\kappa^6\de_0}{t^{2}}< \frac{(2^8-1)\kappa^{6}\de_0}{t^{2}},\label{msuperimpoi1}
	\end{align}
	since $2^{44}\kappa^{34}\left(\delta_0+\delta_0\sqrt{t}\right)<1$. From \eqref{supercruclea1} we deduce that
	\begin{equation}
		\norm{u_{xxxxx}}_{L^1}<\frac{\kappa^{10}\de_0}{(t^*)^{2}}( (2^8-1)+2^{28}\kappa^{15}\delta_0)\leq \frac{\kappa^{10}2^8 \de_0}{(t^*)^{2}},
	\end{equation}
	since $2^{28}\kappa^{15}\delta_0\leq 1$. It gives a contradiction on $t^*$ and it implies in particular  \eqref{estimate:parabolic-1} for $k=4$.\\
	\\
	Next we prove the bound for $\norm{u_{xxxxx}(t,\cdot)}_{L^\f}$ with $t\in]0,\hat{t}]$. Similar to the previous cases, it is enough to prove for smooth initial data. Then there exists $t^*>0$ such that bound \eqref{estimate:parabolic-2} holds for $t\in(0,t^*]$. Then we claim that $t^*\geq\hat{t}$. Assume by absurd that $t^*<\hat{t}$. 
	First we recall from \eqref{superimpml1} and \eqref{msuperimpoi1} that we have for $s\in]0,t^*]$
	\begin{align}
		\norm{\tilde{v}_{x}(s,\cdot)}_{L^\f}\leq \frac{4(2+\sqrt{2})\kappa^4\de_0}{2 s},\quad\norm{\tilde{v}_{xx}(s,\cdot)}_{L^\f}\leq \frac{2^{5}\kappa^5\de_0}{s^{3/2}},\quad\norm{\tilde{v}_{xxx}(s,\cdot)}_{L^\f}\leq  \frac{2^8\kappa^{6}\de_0}{s^{2}}.
		\label{D52}
	\end{align}
	In addition using  \eqref{superitych} we observe that 
	\begin{align}
		&\norm{\tilde{v}_{xxxx}(s,\cdot)}_{L^\f}\leq \big(\kappa^3 \norm{v_{xxxx}}_{L^\f}+3\kappa^5\norm{u_{x}}_{L^\f}\norm{\tilde{v}_{xxx}}_{L^\f}+4\kappa^5(\norm{u_{x}}^2_{L^\f}+\norm{u_{xx}}_{L^\f})\norm{\tilde{v}_{xx}}_{L^\f}\nonumber\\
		&+\kappa^5(\norm{u_{x}}^3_{L^\f}+3\norm{u_{x}}_{L^\f}\norm{u_{xx}}_{L^\f}+\norm{u_{xxx}}_{L^\f})\norm{\tilde{v}_{x}}_{L^\f}\big)(s,\cdot),\nonumber\\
		&\norm{v_{xxxx}(s,\cdot)}_{L^\f}\leq\big( \kappa^3 \norm{\tilde{v}_{xxxx}}_{L^\f}+3\kappa^5\norm{u_{x}}_{L^\f}\norm{\tilde{v}_{xxx}}_{L^\f}+\kappa^5(3\norm{u_{x}}^2_{L^\f}+4\norm{u_{xx}}_{L^\f})\norm{\tilde{v}_{xx}}_{L^\f}\nonumber\\
		&+\kappa^5(\norm{u_{x}}^3_{L^\f}+3\norm{u_{x}}_{L^\f}\norm{u_{xx}}_{L^\f}+\norm{u_{xxx}}_{L^\f})\norm{\tilde{v}_{x}}_{L^\f}\big)(s,\cdot).\label{Dr56}
	\end{align}
	From the definition of $v$ and using \eqref{estimate:parabolic-2},  the estimates \eqref{estimate:parabolic-1} with $k=1,2,3,4$  we obtain
	\begin{align}
		&\norm{v_{xxxx}(s,\cdot)}_{L^\f}\leq \big(\kappa\norm{u_{xxxxx}}_{L^\f}+2^6\kappa\big(\norm{u_{x}}_{L^\f}\norm{u_{xxxx}}_{L^\f}+\norm{u_{x}}^2_{L^\f}\norm{u_{xxx}}_{L^\f}\big)\nonumber\\
		&+2^6\kappa\big( \norm{u_{xx}}_{L^\f}\norm{u_{xxx}}_{L^\f}+\norm{u_{xx}}^2_{L^\f}\norm{u_{x}}_{L^\f}+\norm{u_{xx}}_{L^\f}\norm{u_{x}}^3_{L^\f}+\norm{u_x}^5_{L^\f}\big)\big)(s,\cdot)\nonumber\\
		&\leq\big( \kappa\norm{u_{xxxxx}}_{L^\f}+2^6\kappa\big(\norm{u_{xx}}_{L^1}\norm{u_{xxxxx}}_{L^1}+\norm{u_{xx}}^2_{L^1}\norm{u_{xxxx}}_{L^1}\big)\nonumber\\
		&+2^6\kappa\big( \norm{u_{xxx}}_{L^1}\norm{u_{xxxx}}_{L^1}+\norm{u_{xxx}}^2_{L^1}\norm{u_{xx}}_{L^1}+\norm{u_{xxx}}_{L^1}\norm{u_{xx}}^3_{L^1}+\norm{u_{xx}}^5_{L^1}\big)\big)(s,\cdot)\nonumber\\
		&\leq \frac{2^{16}\kappa^{18}\de_0}{s^{5/2}}.
		\label{Dr57}
	\end{align}
	Similarly from the definition of $v$, using  \eqref{D52}, \eqref{Dr56} and the estimates \eqref{estimate:parabolic-1} with $k=1,2,3,4$  we obtain
	\begin{align}
		&\norm{u_{xxxxx}(s,\cdot)}_{L^\f}\leq\big( \kappa \norm{v_{xxxx}}_{L^\f}+2^6\kappa\big(\norm{u_{xx}}_{L^1}\norm{u_{xxxxx}}_{L^\f}+\norm{u_{xx}}^2_{L^1}\norm{u_{xxxx}}_{L^1}\big)\nonumber\\
		&+2^6\kappa\big( \norm{u_{xxx}}_{L^1}\norm{u_{xxxx}}_{L^1}+\norm{u_{xxx}}^2_{L^1}\norm{u_{xx}}_{L^1}+\norm{u_{xxx}}_{L^1}\norm{u_{xx}}^3_{L^1}+\norm{u_{xx}}^5_{L^1}\big)\big)(s,\cdot)\nonumber\\
		&\leq \kappa\norm{v_{xxxx}(s,\cdot)}_{L^\f}+\frac{2^{16}\kappa ^{16}\de_0^2}{s^{5/2}}\leq \kappa^4\norm{\tilde{v}_{xxxx}(s,\cdot)}_{L^\f}+\frac{2^{24}\kappa ^{22}\de_0^2}{s^{5/2}}.
		\label{sD55}
	\end{align}
	Then, we get from \eqref{Dr56}, \eqref{Dr57}, \eqref{D52} for $s\in]0,t^*]$
	\begin{align}
		& \norm{\tilde{v}_{xxxx}(s,\cdot)}_{L^\f}\leq \frac{2^{25}\kappa^{24}\de_0}{s^{5/2}}.\label{sD57}
		%\norm{v_{xxxx}(s,\cdot)}_{L^\f}&\leq \kappa^3 \norm{\tilde{v}_{xxxx}(s,\cdot)}_{L^\f}+\frac{2^{20}\kappa^{16}\de_0^2}{s^{5/2}},\\
		%\norm{u_{xxxxx}(s,\cdot)}_{L^\f}&\leq \kappa^4\norm{\tilde{v}_{xxxx}(s,\cdot)}_{L^\f}+\frac{2^{21}\kappa^{16}\de_0^2}{s^{5/2}}.
	\end{align}
	By a similar argument as before we have using Lemma \ref{lemma:transformation},
	
	\begin{align}
		&\norm{\tilde{v}_{xxxx}(t)}_{L^\f}
		\leq \norm{G_x(t/2)}_{L^1}\norm{ \tilde{v}_{xxx}(t/2)}_{L^\f}+5\kappa^3 \int\limits_{t/2}^{t}\norm{G_x(t-s)}_{L^1} \nonumber\\
		&\Big\{\norm{([A_2^*-A^*_1B_1^{-1/2}(u)]\mathcal{T}_s(\tilde{v}_{x}))_{xxx}}_{L^\f}+\norm{([B_1(u_x\bullet B_1^{-1/2})]\mathcal{T}_s(\tilde{v}_{x}))_{xxx}}_{L^\f}\nonumber\\
		&+\norm{\big( M(X_s)\mathcal{T}_s( \tilde{v}_x)\big)_{xxx}}_{L^\f} +\norm{(\mathcal{R})_{xxx}}_{L^\f}+\frac{3}{2}\kappa \norm{u_{x}}_{L^\f}\Big(\norm{([A_2^*-A^*_1B_1^{-1/2}(u)]\tilde{v}_{x})_{xx}}_{L^\f}\nonumber\\
		&+\norm{([B_1(u_x\bullet B_1^{-1/2})]\mathcal{T}_s(\tilde{v}_{x}))_{xx}}_{L^\f}+\norm{\big( M(X_s)\mathcal{T}_s( \tilde{v}_x)\big)_{xx}}_{L^\f} +\norm{(\mathcal{R})_{xx}}_{L^\f}\Big)\nonumber\\
		&+(\frac{11}{4}\norm{u_{x}}^2_{L^\f}+\frac{1}{2}\norm{u_{xx}}_{L^\f})\Big(\norm{([A_2^*-A^*_1B_1^{-1/2}(u)]\tilde{v}_{x})_{x}}_{L^\f}+\norm{(\mathcal{R})_{x}}_{L^\f}\nonumber\\
		&+\norm{([B_1(u_x\bullet B_1^{-1/2})]\mathcal{T}_s(\tilde{v}_{x}))_{x}}_{L^\f}+
		\norm{\big( M(X_s)\mathcal{T}_s( \tilde{v}_x)\big)_{x}}_{L^\f} \Big)\Big\}(s,\cdot)\,ds.\label{Duhamel1}
	\end{align}
	Next proceeding as previously, by computations and using \eqref{D52}, \eqref{sD57}, the estimates \eqref{estimate:parabolic-1} 
	\begin{align}
		&\norm{\pa_{xxx}([A_2^*-A^*_1B_1^{-1/2}(u)]\mathcal{T}_s(\tilde{v}_x))(s,\cdot)}_{L^\f}
		\leq \kappa^5\delta_0\norm{\tilde{v}_{xxxx}(s,\cdot)}_{L^\f}+16\kappa^4\norm{u_{x}}_{L^\f}\norm{\tilde{v}_{xxx}}_{L^\f}\nonumber\\
		&+16\kappa^4(\norm{u_{xx}}_{L^\f}+\norm{u_{x}}^2_{L^\f})\norm{\tilde{v}_{xx}}_{L^\f}+\kappa^4\norm{u_{x}}^3_{L^\f}\norm{\tilde{v}_{x}}_{L^\f}\leq 	\frac{2^{25}\kappa^{29}\de_0^2}{s^{5/2}}+\frac{2^{20}\kappa ^{25}\de_0^2}{s^{5/2}},\nonumber\\
		%\\
		&\norm{\pa_{xxx}(B_1(u_x)\bullet B_1^{-1/2}\mathcal{T}_s(\tilde{v}_{x}))(s,\cdot)}_{L^\f}
		\leq  2^{10}\kappa^8(\norm{u_{x}}^3_{L^\f}+\norm{u_{xx}}_{L^\f}\norm{u_{x}}_{L^\f}\nonumber\\
		&+\norm{u_{xxx}}_{L^\f})\norm{\tilde{v}_{xx}}_{L^\f}+2^5\kappa^8(\norm{u_{xxx}}_{L^\f}\norm{u_x}_{L^\f}+\norm{u_{xx}}_{L^\f}^2+\norm{u_{x}}_{L^\f}^4)\norm{\tilde{v}_{x}}_{L^\f}\nonumber\\
		&+2^5\kappa^8(\norm{u_{xxxx}}_{L^\f}+\norm{u_{x}}_{L^\f}^2\norm{u_{xx}}_{L^\f})\norm{\tilde{v}_{x}}_{L^\f}+2^5\kappa^8(\norm{u_{xx}}_{L^\f}+\norm{u_{x}}_{L^\f}^2)\norm{\tilde{v}_{xxx}}_{L^\f}\nonumber\\
		&+\kappa^8\norm{u_{x}}_{L^\f}\norm{\tilde{v}_{xxxx}}_{L^\f}\leq \frac{2^{35}\kappa^{36}\de_0^2}{s^{3}}\nonumber\\
		%\\
		&\norm{\pa_{xxx}(M(X_s)\mathcal{T}_s( \tilde{v}_x)(s,\cdot))}_{L^\f}
		\leq 2^4\kappa\norm{X_{sxxx}}_{L^\f}\norm{\tilde{v}_{x}}_{L^\f}+2^4\kappa^2\norm{X_{sxx}}_{L^\f}\norm{\tilde{v}_{xx}}_{L^\f}\nonumber\\
		&+2^4\kappa^2\norm{X_{sx}}_{L^\f}\norm{u_{x}}_{L^\f}\norm{\tilde{v}_{xx}}_{L^\f}+2^4\kappa^2\norm{X_{s}}_{L^\f}\norm{u_{x}}_{L^\f}\norm{\tilde{v}_{xxx}}_{L^\f}\nonumber\\
		&+2^4\kappa^2\norm{X_{s}}_{L^\f}\norm{u_{x}}^2_{L^\f}\norm{\tilde{v}_{xx}}_{L^\f}+2^4\kappa^2\norm{X_{s}}_{L^\f}\norm{u_{xx}}_{L^\f}\norm{\tilde{v}_{xx}}_{L^\f}\nonumber\\
		&+2^4\kappa^2\norm{X_{s}}_{L^\f}\norm{\tilde{v}_{xxxx}}_{L^\f}\leq 2^{24}\kappa^{27}\left(\frac{\de_0^2}{s^{3}}+\frac{\de_0^2}{s^{5/2}}\right).\label{scruR58}
	\end{align}
	As previously using \eqref{absolucle1},\eqref{Dr57}, \eqref{supercruclea1}, \eqref{D27}, \eqref{superimpml1},  the estimates \eqref{estimate:parabolic-1} we obtain for $s\in]0,t^*]$
	\begin{align}
		&\norm{	\pa_{xxx}\mathcal{R}(u,v,v_x)(s,\cdot)}_{L^\infty}\leq  10^3\cdot 2^{14}\kappa^{12}\Big(\norm{u^*-u}_{L^\f}\norm{v_{xxxx}}_{L^\infty}+\norm{v}^3_{L^\f}\norm{v_{x}}_{L^\infty}
		\nonumber\\
		&+\norm{v}_{L^\f}\norm{v_x}_{L^\f}^2+\norm{v}_{L^\f}\norm{v_{xxx}}_{L^\infty}+\norm{v}_{L^\f}^6+\norm{v}_{L^\f}^4\norm{v_x}_{L^\infty}\nonumber\\
		&+\norm{v}_{L^\f}^2\norm{v_x}_{L^\f}^2+\norm{v}_{L^\f}^3\norm{v_{xx}}_{L^\infty}+\norm{v}_{L^\f}\norm{v_x}_{L^\f}\norm{v_{xx}}_{L^\infty}+\norm{v_x}_{L^\f}^3\nonumber\\
		&+\norm{v_x}_{L^\f}\norm{v_{xxx}}_{L^\infty}+\norm{v}_{L^\f}\norm{v_{xxxx}}_{L^\infty}
		{\color{black}{+\norm{v_{xxx}}_{L^\infty}\norm{v}_{L^\infty}^2+\norm{v}_{L^\infty}^5}}\nonumber\\
		&{\color{black}{+\norm{v_{xx}}_{L^\infty}\norm{v}_{L^\infty}^2+\norm{v_{xx}}_{L^\infty}\norm{v_x}_{L^\infty}+\norm{v_{xx}}_{L^\infty}\norm{v_{xx}}_{L^\infty}     }}\Big)(s,\cdot)\nonumber\\
		&\leq \frac{2^{36}\kappa^{35}\de_0^2}{s^{5/2}}+ \frac{2^{42}\kappa^{38}\de_0^2}{s^{3}}.
		\label{absolucle1yt}
	\end{align}
	Hence, it follows from \eqref{D52}, \eqref{Duhamel1}, \eqref{scruR58}, \eqref{absolucle1yt},  the estimates \eqref{estimate:parabolic-1},\ eqref{D48}, \eqref{D49}, \eqref{estimysupercru}, \eqref{estimysupercru1}, \eqref{estimulcle2}, \eqref{estimulcle3} 
	that for $t\in]0,t^*]$
	\begin{align}
		\norm{\tilde{v}_{xxxx}(t)}_{L^\f}
		&\leq \frac{2^{8}\sqrt{2}\kappa^{7}\de_0}{t^{5/2}}+2^{47}\kappa^{38}\int\limits_{t/2}^{t}\frac{1}{\sqrt{t-s}}\left[\frac{\de_0^2}{s^{3}}+\frac{\de_0^2}{s^{5/2}}\right]\,ds\nonumber\\
		&\leq \frac{2^{8}\sqrt{2}\kappa^{7}\de_0}{t^{5/2}}+2^{49}\kappa^{38}\left(\delta_0+\delta_0\sqrt{t}\right)\frac{\de_0}{t^{{5/2}}}< \frac{(2^{11}-1)\kappa^{7}\de_0}{t^{5/2}},\label{sD56}
	\end{align}
	since $2^{49}\kappa^{38}\left(\delta_0+\delta_0\sqrt{t}\right)<1$. We deduce now from \eqref{sD55}, \eqref{sD56} that
	\begin{equation}
		\norm{u_{xxxxx}(t^*,\cdot)}_{L^\f}\leq \frac{(2^{11}-1)\kappa^{11}\de_0}{(t^*)^{5/2}}+\frac{2^{24}\kappa ^{22}\de_0^2}{(t^*)^{5/2}}
		<\frac{2^{11}\kappa^{11}\de_0}{(t^*)^{5/2}},
	\end{equation}
	since $2^{24}\kappa ^{22}\de_0<1$. We have then obtained a contradiction on $t^*$ and
	this completes the proof of Proposition \ref{prop:parabolic}.
\end{proof}

\subsection{Estimates of $\pa_x^k\mathcal{R}$}\label{appendix:R-derivatives}
Consider a function $g_1: \bar{B}(u^*,\bar{\delta})\times (\R^n)^{k}\rightarrow\R^n$ be a $C^4$ function with $k\in\{2,3\}$ such that $\forall u\in \bar{B}(u^*,\bar{\delta})$, $g_1(u,\cdot)$ is $k$ linear. We observe in particular that this condition implies $\frac{\delta^2 g_1}{\delta^2 z_i}=0$ for $i\in\{1,\cdots,k\}$.
%\begin{equation}\label{condition-g-1}
%\forall u\in B(u^*,\bar{\delta}),\, \mbox{$g(u,\cdot)$ is i linear.}
%g_1(\cdot,0,\cdot,\cdot)=0,g_1(\cdot,\cdot,0,\cdot)=0,g_1(\cdot,\cdot,\cdot,0)=0\mbox{ and }\frac{\pa^2 g_1}{\pa z_i\pa z_i}=0\mbox{ for }i\geq1.
%\end{equation} 
Denote $C^*_1$ as 
\begin{align}
	&C^*_1:=\max\biggl(\sup_{\{u\in  \bar{B}(u^*,\bar{\delta}),\|z_1\|=\cdots=\|z_k\|=1\}} \|g_1(u,z_1,\cdots,z_k)\|,\cdots,\nonumber\\
	&\hspace{1cm} \sup_{\{u\in  \bar{B}(u^*,\bar{\delta}),\|v_1\|=\cdots=\|v_3\|=\|z_1\|=\cdots=\|z_k\|=1\}} \|\frac{\pa^3 g_1}{\pa^3 u}(u,z_1,\cdots,z_k)(v_1,v_2,v_3)\|\biggl).\label{defC1*}
\end{align}
Let $\hat{g}_1(t,x)=g_1(z_0(t,x),z_1(t,x),z_2(t,x),z_3(t,x))$. Then we have for $k=3$ using the fact that $\frac{\delta^2 g_1}{\delta^2 z_i}=0$ for $i\in\{1,\cdots,3\}$
\begin{align}
	\hat{g}_{1,x}&=\sum\limits_{i=0}^{3}z_{i,x}\frac{\pa g_1}{\pa z_i},\\
	\hat{g}_{1,xx}&=(z_{0,x}\otimes z_{0,x})\frac{\pa^2 g_1}{\pa z_0\pa z_0}+2\sum\limits_{i=0}^{3}\sum\limits_{j=i+1}^{3}(z_{i,x}\otimes z_{j,x})\frac{\pa^2 g_1}{\pa z_i\pa z_j}+\sum\limits_{i=0}^{3}z_{i,xx}\frac{\pa g_1}{\pa z_i},\\
	\hat{g}_{1,xxx}&=(z_{0,x}\otimes z_{0,x}\otimes z_{0,x})\frac{\pa^3 g_1}{\pa^3  z_0}+3\sum\limits_{i=1}^3  (z_{0,x}\otimes z_{0,x}\otimes  z_{i,x})\frac{\pa^3 g_1}{\pa^2 z_0 \pa z_i}\nonumber\\
	&+\sum\limits_{i=0}^{3}\sum\limits_{j\ne i, l\ne i,j}(z_{i,x}\otimes z_{j,x}\otimes  z_{l,x})\frac{\pa^3 g_1}{\pa z_i\pa z_j\pa z_l}+3(z_{0,x}\otimes z_{0,xx})\frac{\pa^2 g_1}{\pa^2 z_0}\nonumber\\
	&+3\sum\limits_{i=0}^{3}\sum\limits_{j\ne i}z_{i,xx}\otimes z_{j,x}\frac{\pa^2 g_1}{\pa z_i\pa z_j}+\sum\limits_{i=0}^{3}z_{i,xxx}\frac{\pa g_1}{\pa z_i}.
	%BON\hat{g}_{1,xxx}&=\sum\limits_{i=0}^{3}(z_{i,x}\otimes z_{i,x}\otimes z_{i,x})\frac{\pa^3 g_1}{\pa z_i\pa z_i\pa z_i}+3\sum\limits_{i\ne j} (z_{i,x}\otimes z_{i,x}\otimes  z_{j,x})\frac{\pa^3 g_1}{\pa z_i\pa z_i\pa z_j}\nonumber\\
	%&+\sum\limits_{i=0}^{3}\sum\limits_{j\ne i, l\ne i,j}^{3}(z_{i,x}\otimes z_{j,x}\otimes  z_{l,x})\frac{\pa^3 g_1}{\pa z_i\pa z_j\pa z_l}+3\sum\limits_{i=0}^{3}(z_{i,x}\otimes z_{i,xx})\frac{\pa^2 g_1}{\pa z_i\pa z_i}\nonumber\\
	%&+3\sum\limits_{i=0}^{3}\sum\limits_{j\ne i}^{3}z_{i,xx}\otimes z_{j,x}\frac{\pa^2 g_1}{\pa z_i\pa z_j}+\sum\limits_{i=0}^{3}z_{i,xxx}\frac{\pa g_1}{\pa z_i}.
	%\hat{g}_{1,xxx}&=\sum\limits_{i=0}^{3}(z_{i,x}\otimes z_{i,x}\otimes z_{i,x})\frac{\pa^3 g_1}{\pa z_i\pa z_i\pa z_i}+3\sum\limits_{i=0}^{3}\sum\limits_{j=i+1}^{3}(z_{i,x}\otimes z_{i,x}\otimes  z_{j,x})\frac{\pa^3 g_1}{\pa z_i\pa z_i\pa z_j}\nonumber\\
	%&+3\sum\limits_{i=0}^{3}\sum\limits_{j=i+1}^{3}(z_{i,x}\otimes z_{j,x}\otimes   z_{j,x})\frac{\pa^3 g_1}{\pa z_i\pa z_j\pa z_j}+\sum\limits_{i=0}^{3}\sum\limits_{j=i+1}^{3}\sum\limits_{k=j+1}^{3}(z_{i,x}\otimes z_{j,x}\otimes  z_{k,x})\frac{\pa^3 g_1}{\pa z_i\pa z_j\pa z_k}\nonumber\\
	%&+3\sum\limits_{i=0}^{3}(z_{i,x}\otimes z_{i,xx})\frac{\pa^2 g_1}{\pa z_i\pa z_i}+3\sum\limits_{i=0}^{3}\sum\limits_{j=i+1}^{3}(z_{i,xx}\otimes z_{j,x}+z_{i,x}\otimes z_{j,xx})\frac{\pa^2 g_1}{\pa z_i\pa z_j}\nonumber\\
	%&+\sum\limits_{i=0}^{3}z_{i,xxx}\frac{\pa g_1}{\pa z_i}.
\end{align}
From the assumptions on $g_1$ we have with $i,j,k\in\{1,\cdots,3\}$
\begin{align*}
	&\|g_1(z)\|,\|\frac{\pa g_1}{\pa z_0}(z)\|,\|\frac{\pa^2 g_1}{\pa z_0^2}(z)\|,\|\frac{\pa^3 g_1}{\pa z_0^3}(z)\|\leq C_1^*\|z_1\|\|z_2\| \|z_3\|,\\
	&\|\frac{\pa g_1}{\pa z_i}(z)\|,\|\frac{\pa^2 g_1}{\pa z_0\pa z_i}(z)\|,\|\frac{\pa^3 g_1}{\pa z_0^2\pa z_i}(z)\|\leq C_1^*\|z_j\|\|z_k\|,\mbox{ for }j,k\neq i,\\
	&\|\frac{\pa^2 g_1}{\pa z_i\pa z_j}(z)\|,\|\frac{\pa^3 g_1}{\pa z_0\pa z_i\pa z_j}(z)\|\leq C_1^*\|z_k\|,\mbox{ for }k\neq i,j, i\ne j.
\end{align*}
Then we have the following estimates
\begin{align}
	&\norm{\hat{g}_{1,x}}_{L^\f}\leq C_1^*\left(\norm{z_{0,x}}_{L^\f}\prod\limits_{i=1}^{3}\norm{z_{i}}_{L^\f}+\sum\limits_{i=1}^{3}\norm{z_{i,x}}_{L^\f}\prod\limits_{1\leq j\leq 3,\,j\neq i}\norm{z_{j}}_{L^\f}\right),\nonumber\\
	&\norm{\hat{g}_{1,x}}_{L^1}\leq C_1^*\left(\norm{z_{0,x}}_{L^1}\prod\limits_{i=1}^{3}\norm{z_{i}}_{L^\f}+\sum\limits_{i=1}^{3}\norm{z_{i,x}}_{L^1}\prod\limits_{1\leq j\leq 3,\,j\neq i}\norm{z_{j}}_{L^\f}\right),\nonumber\\
	&\norm{\hat{g}_{1,xx}}_{L^\f}\leq C_1^*\Big((\norm{z_{0,x}}^2_{L^\f}+\norm{z_{0,xx}}_{L^\f})\prod\limits_{i=1}^{3}\norm{z_{i}}_{L^\f}\nonumber\\
	&+\sum\limits_{i=1}^{3}({\color{black}{2}}\norm{z_{0,x}}_{L^\f}\norm{z_{i,x}}_{L^\f}+\norm{z_{i,xx}}_{L^\f})\prod\limits_{1\leq j\leq 3,j\neq i}\norm{z_{j}}_{L^\f}\nonumber\\
	&\hspace{3cm}+{\color{black}{2}}\sum\limits_{i=1}^{3}\sum\limits_{j=i+1}^{3}\norm{z_{i,x}}_{L^\f}\norm{z_{j,x}}_{L^\f}\prod\limits_{1\leq k\leq 3,k\neq i,j}\norm{z_{k}}_{L^\f}\Big),\nonumber\\
	&\norm{\hat{g}_{1,xx}}_{L^1}\leq C_1^*\Big((\norm{z_{0,x}}_{L^1}\norm{z_{0,x}}_{L^\f}+\norm{z_{0,xx}}_{L^1})\prod\limits_{i=1}^{3}\norm{z_{i}}_{L^\f}\nonumber\\
	&+\sum\limits_{i=1}^{3}({\color{black}{2}}\norm{z_{0,x}}_{L^\f}\norm{z_{i,x}}_{L^1}+\norm{z_{i,xx}}_{L^1})\prod\limits_{1\leq j\leq 3,j\neq i}\norm{z_{j}}_{L^\f}\nonumber\\
	&\hspace{3cm}+{\color{black}{2}}\sum\limits_{i=1}^{3}\sum\limits_{j=i+1}^{3}\norm{z_{i,x}}_{L^1}\norm{z_{j,x}}_{L^\f}\prod\limits_{1\leq k\leq3, k\neq i,j}\norm{z_{k}}_{L^\f}\Big),\nonumber\\
	&\norm{\hat{g}_{1,xxx}}_{L^\f}\leq C_1^*\Big((\norm{z_{0,x}}^3_{L^\f}+{\color{black}{3}}\norm{z_{0,x}}_{L^\f}\norm{z_{0,xx}}_{L^\f}+\norm{z_{0,xxx}}_{L^\f})\prod\limits_{i=1}^{3}\norm{z_{i}}_{L^\f}\nonumber\\
	&+\sum\limits_{i=1}^{3}({\color{black}{3}}\norm{z_{0,x}}^2_{L^\f}\norm{z_{i,x}}_{L^\f}+\norm{z_{0,xx}}_{L^\f}\norm{z_{i,x}}_{L^\f})\prod\limits_{1\leq j\leq 3,j\neq i}\norm{z_{j}}_{L^\f}\nonumber\\
	&+\sum\limits_{i=1}^{3}(\norm{z_{i,xxx}}_{L^\f}+\norm{z_{0,x}}_{L^\f}\norm{z_{i,xx}}_{L^\f})\prod\limits_{1\leq j\leq 3,j\neq i}\norm{z_{j}}_{L^\f}\nonumber\\
	&+\sum\limits_{i=1}^{3}\sum\limits_{j=1, i\ne j}^{3}(\norm{z_{i,xx}}_{L^\f}\norm{z_{j,x}}_{L^\f}+{\color{black}{3\norm{z_{0,x}}_{L^\f}\norm{z_{i,x}}_{L^\f}\norm{z_{j,x}}_{L^\f}}} )\prod\limits_{1\leq k\leq 3, k\neq i,j}\norm{z_{k}}_{L^\f}\nonumber\\
	&%+\sum\limits_{i=1}^{3}\sum\limits_{j=i+1}^{3}\norm{z_{0,x}}_{L^\f}\norm{z_{i,x}}_{L^\f}\norm{z_{j,x}}_{L^\f}\prod\limits_{k\neq i,j}\norm{z_{k}}_{L^\f}
	+6\prod\limits_{k=1}^{3}\norm{z_{k,x}}_{L^\f}\Big),\nonumber\\
	&\norm{\hat{g}_{1,xxx}}_{L^1}\leq C_1^*\Big((\norm{z_{0,x}}^2_{L^\f}\norm{z_{0,x}}_{L^1}+{\color{black}{3}}\norm{z_{0,x}}_{L^\f}\norm{z_{0,xx}}_{L^1}+\norm{z_{0,xxx}}_{L^1})\prod\limits_{i=1}^{3}\norm{z_{i}}_{L^\f}\nonumber\\
	&+\sum\limits_{i=1}^{3}({\color{black}{3}}\norm{z_{0,x}}^2_{L^\f}\norm{z_{i,x}}_{L^1}+\norm{z_{0,xx}}_{L^\f}\norm{z_{i,x}}_{L^1})\prod\limits_{1\leq j\leq 3,j\neq i}\norm{z_{j}}_{L^\f}\nonumber\\
	&+\sum\limits_{i=1}^{3}(\norm{z_{i,xxx}}_{L^1}+\norm{z_{0,x}}_{L^\f}\norm{z_{i,xx}}_{L^1})\prod\limits_{1\leq j\leq 3,j\neq i}\norm{z_{j}}_{L^\f}\nonumber\\
	&+\sum\limits_{i=1}^{3}\sum\limits_{j=1, i\ne j}^{3}(\norm{z_{i,xx}}_{L^1}\norm{z_{j,x}}_{L^\f}+{\color{black}{3\norm{z_{0,x}}_{L^\f}\norm{z_{i,x}}_{L^\f}\norm{z_{j,x}}_{L^1}}} )\prod\limits_{1\leq k\leq 3, k\neq i,j}\norm{z_{k}}_{L^\f}\nonumber\\
	%+\sum\limits_{i=1}^{3}\sum\limits_{j=i+1}^{3}(\norm{z_{i,xx}}_{L^1}\norm{z_{j,x}}_{L^\f}+\norm{z_{i,x}}_{L^\f}\norm{z_{j,xx}}_{L^1})\prod\limits_{k\neq i,j}\norm{z_{k}}_{L^\f}\\
	&%+\sum\limits_{i=1}^{3}\sum\limits_{j=i+1}^{3}\norm{z_{0,x}}_{L^\f}\norm{z_{i,x}}_{L^\f}\norm{z_{j,x}}_{L^1}\prod\limits_{k\neq i,j}\norm{z_{k}}_{L^\f}
	+{\color{black}{6}}\norm{z_{1,x}}_{L^1}\prod\limits_{k=2}^{3}\norm{z_{k,x}}_{L^\f}\Big).\label{estimk3}
\end{align}
We proceed similarly when $k=2$ and we have
%Consider a function $g_2:\R^n\times \R^n\times \R^n\rr\R^n$ be a $C^4$ function such that
%\begin{equation}\label{condition-g-2}
%	g_2(\cdot,0,\cdot)=0,g_2(\cdot,\cdot,0)=0,\,\mbox{ and }\frac{\pa^2 g_2}{\pa z_i\pa z_i}=0\mbox{ for }i\geq1.
%\end{equation} 
%Denote $C^*_2$ as $C^*_2:=\norm{g_2}_{C^3(B(0,1))}$. Let $\hat{g}_2(t,x)=g_2(z_0(t,x),z_1(t,x),z_2(t,x))$. Then have
\begin{align*}
	&\hat{g}_{1,x}=\sum\limits_{i=0}^{2}z_{i,x}\frac{\pa g_1}{\pa z_i},\\
	&\hat{g}_{1,xx}=z_{0,x}\otimes z_{0,x}\frac{\pa^2 g_1}{\pa^2 z_0}+\sum\limits_{i=0}^{2}\sum\limits_{i\ne j}(z_{i,x}\otimes z_{j,x})\frac{\pa^2 g_1}{\pa z_i\pa z_j}+\sum\limits_{i=0}^{2}z_{i,xx}\frac{\pa g_1}{\pa z_i},\\
	&\hat{g}_{1,xxx}=3 z_{0,x}\otimes z_{0,xx}\frac{\pa^2 g_1}{\pa^2 z_0}+z_{0,x}\otimes z_{0,x}\otimes z_{0,x}\frac{\pa^3 g_1}{\pa^3 z_0}+3\sum_{j=1}^2 z_{0,x}\otimes z_{0,x}\otimes z_{j,x}\frac{\pa^3 g_1}{\pa^2 z_0\pa z_j}\\
	&+3\sum_{i=0}^2\sum_{j\ne i} z_{i,xx}\otimes z_{j,x}\frac{\pa^2 g_1}{\pa z_i\pa z_j}+6 (z_{0,x}\otimes z_{1,x}\otimes  z_{2,x})\frac{\pa^3 g_1}{\pa z_0\pa z_1\pa z_2}+\sum\limits_{i=0}^{2}z_{i,xxx}\frac{\pa g_1}{\pa z_i}.
	%=\sum\limits_{i=0}^{2}(z_{i,x}\otimes z_{i,x}\otimes z_{i,x})\frac{\pa^3 g_2}{\pa z_i\pa z_i\pa z_i}+3\sum\limits_{i=0}^{2}\sum\limits_{j=i+1}^{2}(z_{i,x}\otimes z_{i,x}\otimes  z_{j,x})\frac{\pa^3 g_2}{\pa z_i\pa z_i\pa z_j}\nonumber\\
	%&+3\sum\limits_{i=0}^{2}\sum\limits_{j=i+1}^{2}(z_{i,x}\otimes z_{j,x}\otimes   z_{j,x})\frac{\pa^3 g_2}{\pa z_i\pa z_j\pa z_j}+(z_{0,x}\otimes z_{1,x}\otimes  z_{2,x})\frac{\pa^3 g_2}{\pa z_0\pa z_1\pa z_2}\nonumber\\
	%&+3\sum\limits_{i=0}^{2}(z_{i,x}\otimes z_{i,xx})\frac{\pa^2 g_2}{\pa z_i\pa z_i}+3\sum\limits_{i=0}^{2}\sum\limits_{j=i+1}^{2}(z_{i,xx}\otimes z_{j,x}+z_{i,x}\otimes z_{j,xx})\frac{\pa^2 g_2}{\pa z_i\pa z_j}\nonumber\\
	%&+\sum\limits_{i=0}^{2}z_{i,xxx}\frac{\pa g_2}{\pa z_i}.
\end{align*}
From the condition of $k$ linearity on $g_1$ we have for $i,j \in\{1,2\}$
\begin{align*}
	\|g_1(z)\|,\|\frac{\pa g_1}{\pa z_0}(z)\|,\|\frac{\pa^2 g_1}{\pa z_0^2}(z)\|,\|\frac{\pa^3 g_1}{\pa z_0^3}(z)\|&\leq C_1^*\|z_1\|\|z_2\|,\\
	\|\frac{\pa g_1}{\pa z_i}(z)\|,\|\frac{\pa^2 g_1}{\pa z_0\pa z_i}(z)\|,\|\frac{\pa^3 g_1}{\pa z_0^2\pa z_i}(z)\|&\leq C_1^*\|z_j\|, \mbox{ for }j\neq i,\\
	\|\frac{\pa^2 g_1}{\pa z_i\pa z_j}(z)\|,\|\frac{\pa^3 g_1}{\pa z_0\pa z_i\pa z_j}(z)\|&\leq C_1^*,\mbox{ for }i\ne j.
\end{align*}
Then we have the following estimates
\begin{align}
	&\norm{\hat{g}_{1,x}}_{L^\f}\leq C_1^* \left(\norm{z_{0,x}}_{L^\f}\prod\limits_{i=1}^{2}\norm{z_{i}}_{L^\f}+\norm{z_{1,x}}_{L^\f}\norm{z_{2}}_{L^\f}+\norm{z_{1}}_{L^\f}\norm{z_{2,x}}_{L^\f}\right),\nonumber\\
	&\norm{\hat{g}_{1,x}}_{L^1}\leq C_1^*\left(\norm{z_{0,x}}_{L^1}\prod\limits_{i=1}^{2}\norm{z_{i}}_{L^\f}+\norm{z_{1,x}}_{L^1}\norm{z_{2}}_{L^\f}+\norm{z_{1}}_{L^\f}\norm{z_{2,x}}_{L^1}\right),\nonumber\\
	&\norm{\hat{g}_{1,xx}}_{L^\f}\leq C_1^* \Big((\norm{z_{0,x}}^2_{L^\f}+\norm{z_{0,xx}}_{L^\f})\prod\limits_{i=1}^{2}\norm{z_{i}}_{L^\f}+{\color{black}{2}}\norm{z_{1,x}}_{L^\f}\norm{z_{2,x}}_{L^\f}\nonumber\\
	&\hspace{1cm}+\sum\limits_{i=1}^{2}({\color{black}{2}}\norm{z_{0,x}}_{L^\f}\norm{z_{i,x}}_{L^\f}+\norm{z_{i,xx}}_{L^\f})\prod\limits_{1\leq j\leq 2,j\neq i}\norm{z_{j}}_{L^\f}\Big),\nonumber\\
	&\norm{\hat{g}_{1,xx}}_{L^1}\leq C_1^* \Big((\norm{z_{0,x}}_{L^\f}\norm{z_{0,x}}_{L^1}+\norm{z_{0,xx}}_{L^1})\prod\limits_{i=1}^{2}\norm{z_{i}}_{L^\f}+{\color{black}{2}}\norm{z_{1,x}}_{L^1}\norm{z_{2,x}}_{L^\f}\nonumber\\
	&\hspace{1cm}+\sum\limits_{i=1}^{2}({\color{black}{2}}\norm{z_{0,x}}_{L^\f}\norm{z_{i,x}}_{L^1}+\norm{z_{i,xx}}_{L^1})\prod\limits_{1\leq j\leq 2,j\neq i}\norm{z_{j}}_{L^\f}\Big),\nonumber\\
	&\norm{\hat{g}_{2,xxx}}_{L^\f}\leq C_1^*\Big((\norm{z_{0,x}}^3_{L^\f}+\norm{z_{0,x}}_{L^\f}\norm{z_{0,xx}}_{L^\f}+\norm{z_{0,xxx}}_{L^\f})\prod\limits_{i=1}^{2}\norm{z_{i}}_{L^\f}\nonumber\\
	&+\sum\limits_{i=1}^{2}({\color{black}{3}}\norm{z_{0,x}}^2_{L^\f}\norm{z_{i,x}}_{L^\f}+{\color{black}{3}}\norm{z_{0,xx}}_{L^\f}\norm{z_{i,x}}_{L^\f})\prod\limits_{1\leq j\leq 2,j\neq i}\norm{z_{j}}_{L^\f}\nonumber\\
	&+\sum\limits_{i=1}^{2}(\norm{z_{i,xxx}}_{L^\f}+{\color{black}{3}}\norm{z_{0,x}}_{L^\f}\norm{z_{i,xx}}_{L^\f})\prod\limits_{1\leq j\leq 2,j\neq i}\norm{z_{j}}_{L^\f}\nonumber\\
	&+{\color{black}{3}}\norm{z_{1,xx}}_{L^\f}\norm{z_{2,x}}_{L^\f}+{\color{black}{3}}\norm{z_{1,x}}_{L^\f}\norm{z_{2,xx}}_{L^\f}+{\color{black}{6}}\norm{z_{0,x}}_{L^\f}\norm{z_{1,x}}_{L^\f}\norm{z_{2,x}}_{L^\f}\Big),\nonumber\\
	&\norm{\hat{g}_{1,xxx}}_{L^1}\leq C_1^*\Big((\norm{z_{0,x}}^2_{L^\f}\norm{z_{0,x}}_{L^1}+\norm{z_{0,x}}_{L^\f}\norm{z_{0,xx}}_{L^1}+\norm{z_{0,xxx}}_{L^1})\prod\limits_{i=1}^{2}\norm{z_{i}}_{L^\f}\nonumber\\
	&+\sum\limits_{i=1}^{2}({\color{black}{3}}\norm{z_{0,x}}^2_{L^\f}\norm{z_{i,x}}_{L^1}+{\color{black}{3}}\norm{z_{0,xx}}_{L^\f}\norm{z_{i,x}}_{L^1})\prod\limits_{1\leq j\leq 2,j\neq i}\norm{z_{j}}_{L^\f}\nonumber\\
	&+\sum\limits_{i=1}^{2}(\norm{z_{i,xxx}}_{L^1}+{\color{black}{3}}\norm{z_{0,x}}_{L^\f}\norm{z_{i,xx}}_{L^1})\prod\limits_{1\leq j\leq 2,j\neq i}\norm{z_{j}}_{L^\f}\nonumber\\
	&+{\color{black}{3}}\norm{z_{1,xx}}_{L^1}\norm{z_{2,x}}_{L^\f}+{\color{black}{3}}\norm{z_{1,x}}_{L^\f}\norm{z_{2,xx}}_{L^1}+{\color{black}{6}}\norm{z_{0,x}}_{L^\f}\norm{z_{1,x}}_{L^\f}\norm{z_{2,x}}_{L^1}\Big).\label{estimk2}
\end{align}
From \eqref{def:remainder-R} recall that
\begin{align}
	\mathcal{R}(u,v,v_x)&:=[A_1^*-A_1(u)]v_x-B_1(((P^{-1}(u)v)\bullet P^{-1}(u)v)\bullet PP^{-1}v)\nonumber\\
	&-B_1((P^{-1}(u)v_x)\bullet PP^{-1}v)-B_1((P^{-1}(u)v)\otimes (P^{-1}(u)v)\bullet D^2PP^{-1}v)\nonumber\\
	&-2B_1[(P^{-1}v)\bullet P(P^{-1}v)\bullet P^{-1}v]-2B_1[(P^{-1}v)\bullet PP^{-1}v_x]\nonumber\\
	&+[BP^{-1}v_x+B[(P^{-1}(u)v)\bullet P^{-1}(u)v]-AP^{-1}v]\bullet P(u)P^{-1}v\nonumber\\
	&+[(P^{-1}(u)v)\bullet B(u)P^{-1}(u)v]\bullet P(u)P^{-1}v+A_1(u)(P^{-1}(u)v\bullet P(u)P^{-1}v)\nonumber\\
	&+P(P^{-1}v)\bullet P^{-1}v\bullet B(u)P^{-1}v+PP^{-1}v_x\bullet B(u)P^{-1}v\nonumber\\
	&+P(P^{-1}v)\otimes(P^{-1}v):D^2B(u)P^{-1}v+{\color{black}{2}}P(P^{-1}v)\bullet B(u)(P^{-1}v)\bullet P^{-1}v\nonumber\\
	&+{\color{black}{2}}P(P^{-1}v)\bullet B(u)P^{-1}v_x-P(P^{-1}v)\bullet A(u)P^{-1}v.
	\label{D34}
\end{align}
We observe now that if we consider $g_1(u,z_1,z_2,z_3)= B_1(((P^{-1}(u)z_1)\bullet P^{-1}(u)z_2)\bullet PP^{-1}z_3)$ on $B(u^*,\bar{\delta})\times(\R^n)^3$, then it follows from the definition of $C_1^*$ in \eqref{defC1*} that by convenience we note $\norm{g_1}_{C^3(\mathcal{W})}$ with $\mathcal{W}=B(u^*,\bar{\delta})\times (B(0,1))^3$
\begin{equation}
	\norm{g_1}_{C^3(\mathcal{W})}\leq 2^8\kappa^8.\label{estimR1}
\end{equation} 
Similarly, we have
\begin{align}
	\norm{(u,z_1,z_2)\mapsto B_1(u)((P^{-1}(u)z_2)\bullet P(u)P^{-1}(u)z_1)}_{C^3(\mathcal{W})}&\leq 2^7\kappa^7,\nonumber\\
	\norm{(u,z_1,z_2,z_3)\mapsto B_1(u)((P^{-1}(u)z_1)\otimes (P^{-1}(u)z_2): D^2P(u)P^{-1}(u)z_3)}_{C^3(\mathcal{W})}&\leq 2^8\kappa^8,\nonumber\\
	\norm{(u,z_1,z_2,z_3)\mapsto B_1(u)[(P^{-1}(u)z_1)\bullet P(u)(P^{-1}(u)z_2)\bullet P^{-1}(u)z_3]}_{C^3(\mathcal{W})}&\leq 2^8\kappa^8,\nonumber\\
	\norm{(u,z_1,z_2)\mapsto B_1(u)[(P^{-1}(u)z_1)\bullet P(u)P^{-1}(u)z_2]}_{C^3(\mathcal{W})}&\leq 2^7\kappa^7,\nonumber\\
	\norm{(u,z_1,z_2)\mapsto [B(u)P^{-1}(u)z_2]\bullet P(u)P^{-1}z_1}_{C^3(\mathcal{W})}&\leq 2^7\kappa^7,\nonumber\\
	\norm{(u,z_1,z_2,z_3)\mapsto [B(u)((P^{-1}(u)z_1)\bullet P^{-1}(u)z_2)]\bullet P(u)P^{-1}z_3}_{C^3(\mathcal{W})}&\leq 2^8\kappa^8,\nonumber\\
	\norm{(u,z_1,z_2)\mapsto [A(u)P^{-1}(u)z_1]\bullet P(u)P^{-1}(u)z_2}_{C^3(\mathcal{W})}&\leq 2^7\kappa^7,\nonumber\\
	\norm{(u,z_1,z_2)\mapsto [(P^{-1}(u)z_1)\bullet B(u)P^{-1}(u)z_2]\bullet P(u)P^{-1}z_3}_{C^3(\mathcal{W})}&\leq 2^8\kappa^8,\nonumber\\
	\norm{(u,z_1,z_2)\mapsto A_1(u)(P^{-1}(u)z_1\bullet P(u)P^{-1}z_2)}_{C^3(\mathcal{W})}&\leq 2^7\kappa^7,\nonumber\\
	\norm{(u,z_1,z_2,z_3)\mapsto P(u)(P^{-1}(u)z_1)\bullet P^{-1}(u)z_2\bullet B(u)P^{-1}(u)z_3}_{C^3(\mathcal{W})}&\leq 2^8\kappa^8,\nonumber\\
	\norm{(u,z_1,z_2)\mapsto P(u)P^{-1}(u)z_2\bullet B(u)P^{-1}(u)z_1}_{C^3(\mathcal{W})}&\leq 2^7\kappa^7,\nonumber\\
	\norm{(u,z_1,z_2,z_3)\mapsto P(u)(P^{-1}(u)z_1)\otimes(P^{-1}(u)z_2):D^2B(u)P^{-1}(u)z_3}_{C^3(\mathcal{W})}&\leq 2^8\kappa^8,\nonumber\\
	\norm{(u,z_1,z_2,z_3)\mapsto P(u)(P^{-1}(u)z_1)\bullet B(u)(P^{-1}(u)z_2)\bullet P^{-1}(u)z_3}_{C^3(\mathcal{W})}&\leq 2^8\kappa^8,\nonumber\\
	\norm{(u,z_1,z_2)\mapsto P(u)(P^{-1}(u)z_1)\bullet B(u)P^{-1}(u)z_2}_{C^3(\mathcal{W})}&\leq 2^7\kappa^7,\nonumber\\
	\norm{(u,z_1,z_2)\mapsto P(u)(P^{-1}(u)z_1)\bullet A(u)P^{-1}(u)z_2}_{C^3(\mathcal{W})}&\leq 2^7\kappa^7.\label{estimg1}
\end{align}
We wish now to apply the estimates \eqref{estimk3}, \eqref{estimk2} for dealing with the term $\mathcal{R}(u,v,v_x)$ with $z_{1}=z_2=v$, $z_{1}=z_2=z_3=v$ or $z_1=v$, $z_2=v_x$.
We will use in the sequel the fact that 
\begin{align}
	&\norm{u_x}_{L^\f}\leq \kappa\norm{v}_{L^\f}, \,\norm{u_x}_{L^1}\leq \kappa\norm{v}_{L^1},\,\norm{u_{xx}}_{L^\infty}\leq \kappa(\norm{v_x}_{L^\infty}+\kappa \norm{v}_{L^\infty}^2),\nonumber\\
	&\norm{u_{xx}}_{L^1}\leq \kappa(\norm{v_x}_{L^1}+\kappa \norm{v}_{L^\infty}\norm{v}_{L^1}),\nonumber\\
	&\norm{u_{xxx}}_{L^\infty}\leq \kappa(\norm{v_{xx}}_{L^\infty}+\kappa \norm{v}_{L^\infty}\norm{v_x}_{L^\infty}+\kappa^2\norm{v}_{L^\infty}^3 ),\nonumber\\
	&\norm{u_{xxx}}_{L^1}\leq \kappa(\norm{v_{xx}}_{L^1}+\kappa \norm{v}_{L^\infty}\norm{v_x}_{L^1}+\kappa^2\norm{v}_{L^\infty}^2\norm{v}_{L^1} ).
	\label{vbtech}
\end{align}
Combining now \eqref{vbtech}, \eqref{D34}, \eqref{estimR1}, \eqref{estimk3}, \eqref{estimk2}, \eqref{estimg1} then we have
\begin{align}
	&\norm{	\mathcal{R}(u,v,v_x)}_{L^\f}\leq  20\cdot 2^8\kappa^8\left(\norm{u^*-u}_{L^\f}\norm{v_x}_{L^\f}+\norm{v}_{L^\f}^2+\norm{v}_{L^\f}^3+\norm{v}_{L^\f}\norm{v_x}_{L^\f}\right),\nonumber\\
	&\norm{	\mathcal{R}(u,v,v_x)}_{L^1}\leq  20\cdot 2^8\kappa^8\Big(\norm{u^*-u}_{L^\f}\norm{v_x}_{L^1}+\norm{v}_{L^\f}\norm{v}_{L^1}+\norm{v}_{L^\f}^2\norm{v}_{L^1}\nonumber\\
	&\hspace{2.4 cm}+\norm{v}_{L^\f}\norm{v_x}_{L^1}\Big),\nonumber\\
	%\\
	%\\
	%\\
	&\norm{	\pa_x\mathcal{R}(u,v,v_x)}_{L^\f}\leq  20\cdot 2^{10}\kappa^9\Big(\norm{u^*-u}_{L^\f}\norm{v_{xx}}_{L^\f}+\norm{v}_{L^\f}\norm{v_{x}}_{L^\f}+\norm{v}_{L^\f}^4+ {\color{black}{\norm{v}_{L^\f}^3}}\nonumber\\
	&\hspace*{2cm}+\norm{v}_{L^\f}^2\norm{v_{x}}_{L^\f}+\norm{v_x}_{L^\f}^2+\norm{v}_{L^\f}\norm{v_{xx}}_{L^\f}\Big),\nonumber\\
	%\\
	%\\
	%\\
	&\norm{	\pa_x\mathcal{R}(u,v,v_x)}_{L^1}\leq  20\cdot 2^{10}\kappa^9\Big(\norm{u^*-u}_{L^\f}\norm{v_{xx}}_{L^1}+\norm{v}_{L^\f}\norm{v_{x}}_{L^1}\nonumber\\
	&+\norm{v}_{L^\f}^3\norm{v}_{L^1}+{\color{black}{\norm{v}_{L^\f}^2}\norm{v}_{L^1}}+\norm{v}_{L^\f}^2\norm{v_{x}}_{L^1}+\norm{v_x}_{L^\f}\norm{v_x}_{L^1}+\norm{v}_{L^\f}\norm{v_{xx}}_{L^1}\Big),\nonumber\\
	%\\
	%\\
	%\\
	&\norm{	\pa_{xx}\mathcal{R}(u,v,v_x)}_{L^\f}\leq  10^2 \cdot 2^{12}\kappa^{10}\Big(\norm{u^*-u}_{L^\f}\norm{v_{xxx}}_{L^\f}+\norm{v}^2_{L^\f}\norm{v_{x}}_{L^\f}+\norm{v}_{L^\f}\norm{v_{xx}}_{L^\f}\nonumber\\
	&+{\color{black}{\norm{v}_{L^\infty}^4+\norm{v_x}_{L^\infty}^2}}+\norm{v}_{L^\f}^5+\norm{v}_{L^\f}^3\norm{v_x}_{L^\f}+\norm{v}_{L^\f}\norm{v_{x}}^2_{L^\f}+\norm{v}^2_{L^\f}\norm{v_{xx}}_{L^\f}\nonumber\\
	&\hspace*{2cm}+\norm{v_x}_{L^\f}\norm{v_{xx}}_{L^\f}+\norm{v}_{L^\f}\norm{v_{xxx}}_{L^\f}\Big),\nonumber\\
	%\\
	%\\
	%\\
	&\norm{	\pa_{xx}\mathcal{R}(u,v,v_x)}_{L^1}\leq  10^2\cdot 2^{12}\kappa^{10}\Big(\norm{u^*-u}_{L^\f}\norm{v_{xxx}}_{L^1}+\norm{v}^2_{L^\f}\norm{v_{x}}_{L^1}+\norm{v}_{L^\f}\norm{v_{xx}}_{L^1}\nonumber\\
	&+{\color{black}{\norm{v}_{L^\infty}^3\norm{v}_{L^1}+\norm{v_x}_{L^\infty}\norm{v_x}_{L^1}}}+\norm{v}_{L^\f}^4\norm{v}_{L^1}+\norm{v}_{L^\f}^3\norm{v_x}_{L^1}+\norm{v}_{L^\f}\norm{v_{x}}_{L^\f}\norm{v_x}_{L^1}\nonumber\\
	&\hspace*{2cm}+\norm{v}^2_{L^\f}\norm{v_{xx}}_{L^1}+\norm{v_x}_{L^\f}\norm{v_{xx}}_{L^1}+\norm{v}_{L^\f}\norm{v_{xxx}}_{L^1}\Big),\nonumber\\
	%\\
	%\\
	%\\
	&\norm{	\pa_{xxx}\mathcal{R}(u,v,v_x)}_{L^1}\leq  10^3\cdot 2^{14}\kappa^{12}\Big(\norm{u^*-u}_{L^\f}\norm{v_{xxxx}}_{L^1}+\norm{v}^3_{L^\f}\norm{v_{x}}_{L^1}
	\nonumber\\
	&+\norm{v}_{L^\f}\norm{v_x}_{L^\f}\norm{v_{x}}_{L^1}+\norm{v}_{L^\f}\norm{v_{xxx}}_{L^1}+\norm{v}_{L^\f}^5\norm{v}_{L^1}+\norm{v}_{L^\f}^4\norm{v_x}_{L^1}\nonumber\\
	&+\norm{v}_{L^\f}^2\norm{v_x}_{L^\f}\norm{v_x}_{L^1}+\norm{v}_{L^\f}^3\norm{v_{xx}}_{L^1}+\norm{v}_{L^\f}\norm{v_x}_{L^\f}\norm{v_{xx}}_{L^1}+\norm{v_x}_{L^\f}^2\norm{v_{x}}_{L^1}\nonumber\\
	&+\norm{v_x}_{L^\f}\norm{v_{xxx}}_{L^1}+\norm{v}_{L^\f}\norm{v_{xxxx}}_{L^1}
	{\color{black}{+\norm{v_{xxx}}_{L^1}\norm{v}_{L^\infty}^2+\norm{v}_{L^\infty}^4\norm{v}_{L^1}
	}}\nonumber\\
	&{\color{black}{+\norm{v_{xx}}_{L^1}\norm{v}_{L^\infty}^2+\norm{v_{xx}}_{L^1}\norm{v_x}_{L^\infty}+\norm{v_{xx}}_{L^1}\norm{v_{xx}}_{L^\infty}     }}\Big).\nonumber\\
	&\norm{	\pa_{xxx}\mathcal{R}(u,v,v_x)}_{L^\infty}\leq  10^3\cdot 2^{14}\kappa^{12}\Big(\norm{u^*-u}_{L^\f}\norm{v_{xxxx}}_{L^\infty}+\norm{v}^3_{L^\f}\norm{v_{x}}_{L^\infty}
	\nonumber\\
	&+\norm{v}_{L^\f}\norm{v_x}_{L^\f}^2+\norm{v}_{L^\f}\norm{v_{xxx}}_{L^\infty}+\norm{v}_{L^\f}^6+\norm{v}_{L^\f}^4\norm{v_x}_{L^\infty}\nonumber\\
	&+\norm{v}_{L^\f}^2\norm{v_x}_{L^\f}^2+\norm{v}_{L^\f}^3\norm{v_{xx}}_{L^\infty}+\norm{v}_{L^\f}\norm{v_x}_{L^\f}\norm{v_{xx}}_{L^\infty}+\norm{v_x}_{L^\f}^3\nonumber\\
	&+\norm{v_x}_{L^\f}\norm{v_{xxx}}_{L^\infty}+\norm{v}_{L^\f}\norm{v_{xxxx}}_{L^\infty}
	{\color{black}{+\norm{v_{xxx}}_{L^\infty}\norm{v}_{L^\infty}^2+\norm{v}_{L^\infty}^5}}\nonumber\\
	&{\color{black}{+\norm{v_{xx}}_{L^\infty}\norm{v}_{L^\infty}^2+\norm{v_{xx}}_{L^\infty}\norm{v_x}_{L^\infty}+\norm{v_{xx}}_{L^\infty}\norm{v_{xx}}_{L^\infty}     }}\Big).\label{absolucle1}
\end{align}

\section{Detailed justification for \eqref{equaJ1x} and \eqref{equaK1x}}
\label{sectionE}
\subsection{Calculation for $J_{1,x}+J_{2,x}+J_{3,x}$ in  \eqref{equaJ1x}}
We wish here to justify the tedious computation of  \eqref{equaJ1x}
\begin{align*}
	&\al_i^{1,4}+\al_i^{2,8}+\al_i^{3,1}+\al_i^{3,3}=\sum\limits_{i}\tilde{r}_{i,v}(2\xi_i\pa_{v_i}\bar{v}_i+v_i\xi_i\pa_{v_iv_i}\bar{v}_i)\left[(\mu_iv_{i,x}-(\tilde{\la}_i-\la_i^*)v_i-\theta_i v_i)v_{i,x}\right]\nonumber\\
	&+\sum\limits_{i}\tilde{r}_{i,v}\mu_i v_{i,x}\xi^\p_i\left(w_{i,x}-\frac{w_i}{v_i}v_{i,x}\right) \pa_{v_i}\bar{v}_i.\\
	&\al_i^{1,5}+\al_i^{2,10}+\al_i^{3,2}+\al_i^{3,11}=\sum\limits_{i}\tilde{r}_{i,v}v_i\xi_i^\p\pa_{v_i}\bar{v}_i(\mu_iv_{i,x}-{\color{black}{(\tilde{\la}_i-\la_i^*)v_i-\theta_i v_i}})\left(\frac{w_i}{v_i}\right)_x\nonumber\\
	&+\sum\limits_{i}\tilde{r}_{i,v}\mu_i\xi^{\p\p}_i\left(w_{i,x}-\frac{w_i}{v_i}v_{i,x}\right) \bar{v}_i\left(\frac{w_i}{v_i}\right)_x-\sum_i v_i^2\left(\frac{w_i}{v_i}\right)_x\pa_{v_i}\bar{v}_i\tilde{r}_{i,v}\theta'_i\xi_i  . \nonumber\\[2mm]
	&\al_i^{1,6}+\al_i^{2,9}+\al_i^{3,4}=\sum\limits_{i}\sum\limits_{j\neq i}\tilde{r}_{i,v}(\xi_i\pa_{v_j}\bar{v}_i)\left[(\mu_iv_{i,x}-\tilde{\la}_iv_i)v_{j,x}\right]\\
	&+\sum\limits_{i}\sum\limits_{j\neq i}\tilde{r}_{i,v}v_i\xi_i\pa_{v_iv_j}\bar{v}_i\left[(\mu_iv_{i,x}-(\tilde{\la}_i-\la_i^*)v_i-\theta_i v_i)v_{j,x}\right] +\sum\limits_{i}\sum\limits_{j\neq i}\tilde{r}_{i,v}\mu_i\xi_i^\p\pa_{v_j}\bar{v}_iv_i\left(\frac{w_i}{v_i}\right)_xv_{j,x}.\nonumber\\[2mm]
	&\al_i^{1,7}+\al_i^{2,2}+\al_i^{3,5}=\sum\limits_{i}\tilde{r}_{i,uv}\tilde{r}_i v_i\xi_i\pa_{v_i}\bar{v}_i\left[(\mu_iv_{i,x}-(\tilde{\la}_i-\la_i^*)v_i- \theta_i v_i)v_i\right]\nonumber\\
	&+\sum\limits_{i}\tilde{r}_{i,uv}\tilde{r}_i \mu_i\xi_i^\p\bar{v}_i(w_{i,x}v_i-w_iv_{i,x}).\\[2mm]
	&\al_i^{1,8}+\al_i^{2,3}+\al_i^{3,6}=\sum\limits_{i}\sum\limits_{j\neq i}\tilde{r}_{i,uv}\tilde{r}_j(\xi_i\pa_{v_i}\bar{v}_i)v_jv_i(\mu_iv_{i,x}-(\tilde{\la}_i-\la_i^*)v_i-\theta_i v_i)\nonumber\\
	&+\sum\limits_{i}\sum\limits_{j\neq i}\tilde{r}_{i,uv}\tilde{r}_j(w_{i,x}-(w_i/v_i)v_{i,x})\xi_i^\p\bar{v}_iv_j\mu_i.\nonumber\\[2mm]
	&\al_i^{1,10}+\al_i^{2,4}+\al_i^{3,7}=\sum\limits_{i}\tilde{r}_{i,vv}(\xi_i\pa_{v_i}\bar{v}_i)^2v_iv_{i,x}(\mu_iv_{i,x}-(\tilde{\la}_i-\la_i^*)v_i-\theta_i v_i)\nonumber\\
	&+\sum\limits_{i}\tilde{r}_{i,vv}(\xi_i\pa_{v_i}\bar{v}_i)\xi_i^\p\bar{v}_i(w_{i,x}-(w_i/v_i)v_{i,x})v_{i,x}\mu_i.\nonumber\\[2mm]
	&\al_i^{1,11}+\al_i^{2,6}+\al_i^{3,9}=\sum\limits_{i}\tilde{r}_{i,vv}\xi_i^\p\bar{v}_i\left(\frac{w_i}{v_i}\right)_x(v_i\xi_i\pa_{v_i}\bar{v}_i)\left[(\mu_iv_{i,x}-(\tilde{\la}_i-\la_i^*)v_i-\theta_i v_i)\right]\nonumber\\
	&+\sum\limits_{i}\tilde{r}_{i,vv}(\xi_i^\p)^2\mu_i\bar{v}_i^{2}\left(\frac{w_i}{v_i}\right)_x(w_{i,x}-(w_i/v_i)v_{i,x}).\nonumber\\[2mm]
	&\al_i^{1,12}+\al_i^{2,5}+\al_i^{3,8}=\sum\limits_{i}\sum\limits_{j\neq i}\tilde{r}_{i,vv}\xi_i^2\pa_{v_j}\bar{v}_i\pa_{v_i}\bar{v}_iv_i(\mu_iv_{i,x}-(\tilde{\la}_i-\la_i^*)v_i-\theta_i v_i)v_{j,x}\nonumber\\
	%&+\sum\limits_{i}\sum\limits_{j\neq i}\tilde{r}_{i,vv}\xi_i^2\pa_{w_j}\bar{v}_i\pa_{v_i}\bar{v}_iv_i(\mu_iv_{i,x}-(\tilde{\la}_i-\la_i^*)v_i-\theta_i v_i)w_{j,x} \nonumber\\
	&+\sum\limits_{i}\sum\limits_{j\neq i}\tilde{r}_{i,vv}\mu_i\xi_i\xi^\p_i\pa_{v_j}\bar{v}_i\bar{v}_i\left(w_{i,x}-\frac{w_i}{v_i}v_{i,x}\right)v_{j,x}.\nonumber\\[2mm]
	&\al_i^{1,13}+\al_i{1,18}+\al_i^{2,7}+\al_i^{2,15}+\al_i^{3,10}=-\sum\limits_{i}\tilde{r}_{i,v\si}v_i\xi_i\pa_{v_i}\bar{v}_i(\mu_iv_{i,x}-(\tilde{\la}_i-\la_i^*)v_i-\theta_i v_i)\theta_i^\p \left(\frac{w_i}{v_i}\right)_x\nonumber\\
	&-{\color{black}{2}}\sum\limits_{i}\tilde{r}_{i,v\si}\xi_i^\p\bar{v}_i\mu_i(w_{i,x}-(w_i/v_i)v_{i,x}) \theta_i^\p \left(\frac{w_i}{v_i}\right)_x\\[2mm]
	&\al_i^{1,16}+\al_i^{1,17}+\al_i^{2,13}+\al_i^{2,14}=-\sum\limits_{i}\sum\limits_{k}\tilde{r}_{i,v\si}\xi_i\pa_{v_k}\bar{v}_iv_{k,x}\theta^\p_i\mu_i\left[w_{i,x}-\frac{w_i}{v_i}v_{i,x}\right].\nonumber\\[2mm]
	&\al_i^{1,14}+\al_i^{2,11}=\sum\limits_{i}\tilde{r}_{i,u\si}\tilde{r}_i \theta^\p_i \mu_i(w_iv_{i,x}-w_{i,x}v_i).\\[2mm]
	&\al_i^{1,15}+\al_i^{2,12}=-\sum\limits_{i}\sum\limits_{j\neq i}\tilde{r}_{i,u\si}\tilde{r}_j\mu_i \theta^\p_iv_j\left[w_{i,x}-(w_i/v_i)v_{i,x}\right].\nonumber\\
	&\al_i^{1,19}+\al_i^{2,16}=\sum\limits_{i}\tilde{r}_{i,\si\si} \mu_i(\theta_i^\p)^2\left(\frac{w_i}{v_i}\right)_x\left[w_{i,x}-(w_i/v_i)v_{i,x}\right].\nonumber\\
	&\al_i^{1,7}+\al_i^{1,20}+\al_i^{2,17}=-\sum\limits_{i}\tilde{r}_{i,\si}\mu_i\theta_i^{\p\p}\left(\frac{w_i}{v_i}\right)_x\left[w_{i,x}-(w_i/v_i)v_{i,x}\right].
\end{align*}

\subsection{Calculation for $K_{1,x}+K_{2,x}+K_{5,x}$ in \eqref{equaK1x}}
Furthermore similarly we justify the computation of \eqref{equaK1x} by observing that
\begin{align*}
	&\beta_i^{1,3}-\sum_i\sum_{j\ne i}(w_j-\la_j^*v_j)w_i\tilde{r}_{i,u}\tilde{r}_j=\sum\limits_{i}\sum\limits_{j\neq i}[(\mu_iw_{i,x}-\tilde{\la}_iw_i)v_j-(w_j-\la_j^*v_j)w_i]\tilde{r}_{i,u}\tilde{r}_j.\\
	&\beta_i^{1,2}-\sum_i(w_i-\la_i^*v_jiw_i\tilde{r}_{i,u}\tilde{r}_i=
	\sum\limits_{i}[\mu_i(w_{i,x}v_i-w_iv_{i,x})+(\mu_iv_{i,x}-(\tilde{\la}_i-\la_i^*)v_i-w_i)w_i]\tilde{r}_{i,u}\tilde{r}_i.\\[2mm]
	&\beta_i^{1,4}+\beta_i^{2,2}+\beta_i^{2,3}+\beta_i^{5,2}+\beta_i^{5,4}=\sum\limits_{i}(\mu_iw_{i,x}-(\tilde{\la}_i-\la_i^*+\theta_i)w_i)v_{i,x}\xi_i(\pa_{v_i}\bar{v}_i)\tilde{r}_{i,v}\\
	&+\sum_i(\mu_i v_{i,x}-v_i(\tilde{\lambda}_i-\lambda_i^*)-\theta_i v_i)\xi_i w_iv_{i,x}\pa_{v_iv_i}\bar{v}_i\tilde{r}_{i,v}+\sum_{i}(\mu_i v_{i,x}-v_i(\tilde{\lambda}_i-\lambda_i^*)-\theta_i v_i) \xi_iw_{i,x}\pa_{v_i}\bar{v}_i\tilde{r}_{i,v}.\\[2mm]
	&\beta_i^{1,6}+\beta_i^{1,7}+\beta_i^{2,12}+\beta_i^{2,15}=2\sum\limits_{i}\mu_i(w_{i,x}-(w_i/v_i)v_{i,x})\xi^\p_i\bar{v}_i\left(\frac{w_i}{v_i}\right)_x\tilde{r}_{i,v}\\
	&+\sum\limits_{i}\mu_i(w_{i,x}-(w_i/v_i)v_{i,x})\xi^{\p\p}_i\bar{v}_i\frac{w_i}{v_i}\left(\frac{w_i}{v_i}\right)_x\tilde{r}_{i,v}.\\[2mm]
	&\beta_i^{1,8}+\beta_i^{2,11}+\beta_i^{2,13}+\beta_i^{5,3}=\sum\limits_{i}\mu_i(w_{i,x}-(w_i/v_i)v_{i,x})v_{i,x}\xi^\p_i(\pa_{v_i}\bar{v}_i)\frac{w_i}{v_i}\tilde{r}_{i,v}\\
	&+\sum_{i}(\mu_i v_{i,x}-v_i(\tilde{\lambda}_i-\lambda_i^*)-\theta_i v_i) w_i\xi_i^\p\left(\frac{w_i}{v_i}\right)_x\pa_{v_i}\bar{v}_i\tilde{r}_{i,v}.\\[2mm]
	&\beta_i^{1,5}+\beta_i^{1,9}+\beta_i^{2,4}+\beta_i^{2,14}+\beta_i^{5,5}=\sum\limits_{i}\sum\limits_{j\neq i}(\mu_iw_{i,x}-\tilde{\la}_iw_i )v_{j,x}\xi_i(\pa_{v_j}\bar{v}_i)\tilde{r}_{i,v}\\
	&+\sum\limits_{i}\sum\limits_{j\neq i}\mu_i(w_{i,x}-(w_i/v_i)v_{i,x})v_{j,x}\xi_i^\p(\pa_{v_j}\bar{v}_i)\frac{w_i}{v_i}\tilde{r}_{i,v}+\sum_{i\ne j}(\mu_i v_{i,x}-v_i(\tilde{\lambda}_i-\lambda_i^*)-\theta_i v_i) \xi_iw_iv_{j,x}\pa_{v_jv_i}\bar{v}_i\tilde{r}_{i,v}.\\[2mm]
	&\beta_i^{1,10}+\beta_i^{1,11}+\beta_i^{2,22}+\beta_i^{2,23}=-2\sum\limits_{i}\mu_i(w_{i,x}-(w_i/v_i)v_{i,x})\theta_i^\p\left(\frac{w_i}{v_i}\right)_x\tilde{r}_{i,\si}\\
	&-\sum\limits_{i}\mu_i(w_{i,x}-(w_i/v_i)v_{i,x})\frac{w_i}{v_i}\left(\frac{w_i}{v_i}\right)_x\theta_i^{\p\p} \tilde{r}_{i,\si}.\\[2mm]
	&\beta_i^{1,12}+\beta_i^{2,16}+\beta_i^{2,5}+\beta_i^{5,6}=\sum_{i}\xi_i w_i v_i\pa_{v_i}\bar{v}_i (\mu_i v_{i,x}-v_i(\tilde{\lambda}_i-\lambda_i^*)-\theta_i v_i) \tilde{r}_{i,uv}\tilde{r}_i\\
	&+\sum\limits_{i}\mu_i(w_{i,x}-(w_i/v_i)v_{i,x})w_i \xi_i^\p\bar{v}_i\tilde{r}_{i,uv}\tilde{r}_i.\\[2mm]
	&\beta_i^{1,13}+\beta_i^{2,6}+\beta_i^{2,17}+\beta_i^{5,7}=\sum_{i\ne j}\xi_i w_i v_j\pa_{v_i}\bar{v}_i (\mu_i v_{i,x}-v_i(\tilde{\lambda}_i-\lambda_i^*)-\theta_i v_i) \tilde{r}_{i,uv}\tilde{r}_j\\
	&+\sum\limits_{i}\sum\limits_{j\neq i}\mu_i(w_{i,x}-(w_i/v_i)v_{i,x})v_j\frac{w_i}{v_i}\xi_i^\p\bar{v}_i\tilde{r}_{i,uv}\tilde{r}_i.\\[2mm]
	&\beta_i^{1,14}+\beta_i^{1,16}+\beta_i^{2,18}+\beta_i^{2,20}=\sum\limits_{i}\mu_i(w_{i,x}-(w_i/v_i)v_{i,x})v_{i,x}\frac{w_i}{v_i}\xi_i^\p\xi_i\pa_{v_i}\bar{v}_i\bar{v}_i\tilde{r}_{i,vv}\\
	&+\sum\limits_{i}\mu_i(w_{i,x}-(w_i/v_i)v_{i,x})\frac{w_i}{v_i}(\bar{v}_i\xi_i^\p)^2\left(\frac{w_i}{v_i}\right)_x\tilde{r}_{i,vv}.\\[2mm]
	&\beta_i^{2,7}+\beta_i^{2,10}+\beta_i^{5,8}+\beta_i^{5,9}=\sum_i {\color{black}{\xi_i^2}}w_iv_{i,x}(\pa_{v_i}\bar{v}_i)^2 (\mu_iv_{i,x}-v_i(\tilde{\lambda}_i-\lambda_i^*)-\theta_i v_i)\tilde{r}_{i,vv}\\
	&+\sum_{i}\xi_i' \bar{v}_i w_i \xi_i \pa_{v_i}\bar{v}_i\left(\frac{w_i}{v_i}\right)_x(\mu_iv_{i,x}-v_i(\tilde{\lambda}_i-\lambda_i^*)-\theta_i v_i)\tilde{r}_{i,vv}.\\[2mm]
	&\beta_i^{1,15}+\beta_i^{2,8}+\beta_i^{2,19}+\beta_i^{5,10}=\sum_{i\ne j}{\color{black}{\xi_i^2}}v_{j,x}w_i \pa_{v_i}\bar{v}_i\pa_{v_j}\bar{v}_i(\mu_iv_{i,x}-v_i(\tilde{\lambda}_i-\lambda_i^*)-\theta_i v_i)\tilde{r}_{i,vv}\\
	&+\sum\limits_{i}\sum\limits_{j\neq i}\mu_i(w_{i,x}-(w_i/v_i)v_{i,x})v_{j,x}\frac{w_i}{v_i}\xi_i^\p\xi_i\bar{v}_i\pa_{v_j}\bar{v}_i\tilde{r}_{i,vv}.\\[2mm]
	&\beta_i^{1,17}+\beta_i^{20}+\beta_i^{22}+\beta_i^{2,21}+\beta_i^{2,26}+\beta_i^{2,28}=-\sum\limits_{i}\mu_i(w_{i,x}-(w_i/v_i)v_{i,x})\frac{w_i}{v_i}(\bar{v}_i\xi_i^\p)\theta_i^\p\left(\frac{w_i}{v_i}\right)_x\tilde{r}_{i,v\si} \\
	&-\sum\limits_{i}\mu_i(w_{i,x}-(w_i/v_i)v_{i,x}) \xi_i\frac{w_i}{v_i}v_{i,x}\theta_i^\p \pa_{v_i}\bar{v}_i\tilde{r}_{i,v\si}-\sum\limits_{i}\mu_i(w_{i,x}-(w_i/v_i)v_{i,x}) \xi_i^\p\frac{w_i}{v_i}\left(\frac{w_i}{v_i}\right)_x\theta_i^\p \bar{v}_i\tilde{r}_{i,v\si}.\\[2mm]
	&\beta_i^{2,9}+\beta_i^{5,11}=-\sum_i w_i\xi_i\pa_{v_i}\bar{v}_i\theta'_i\left(\frac{w_i}{v_i}\right)_x(\mu_iv_{i,x}-v_i(\tilde{\lambda}_i-\lambda_i^*)-\theta_i v_i)\tilde{r}_{i,v\si}.\\[2mm]
	&\beta_i^{1,21}+\beta_i^{2,27}=-\sum\limits_{i}\sum\limits_{j\neq i}\mu_i(w_{i,x}-(w_i/v_i)v_{i,x})\xi_i \frac{w_i}{v_i}v_{j,x}\theta_i^\p \pa_{v_j}\bar{v}_i\tilde{r}_{i,v\si}.\\[2mm]
	&\beta_i^{1,18}+\beta_i^{1,19}+\beta_i^{2,24}+\beta_i^{2,25}=-\sum\limits_{i}\mu_i(w_{i,x}-(w_i/v_i)v_{i,x}) w_i\theta_i^\p \tilde{r}_{i,u\si}\tilde{r}_i\\
	&-\sum\limits_{i}\sum\limits_{j\neq i}\mu_i(w_{i,x}-(w_i/v_i)v_{i,x}) (w_i/v_i)v_j\theta_i^\p \tilde{r}_{i,u\si}\tilde{r}_j.\\[2mm]
	&\beta_i^{1,23}+\beta_i^{2,29}=\sum\limits_{i}\mu_i(w_{i,x}-(w_i/v_i)v_{i,x}) \frac{w_i}{v_i}\left(\frac{w_i}{v_i}\right)_x(\theta_i^\p)^2 \tilde{r}_{i,\si\si}.
\end{align*}

\bigskip

\noi\textbf{Acknowledgements:} BH would like to thank the ANR HEAD project ANR-24-CE40-3260.

AJ would like to thank Anusandhan National Research Foundation (ANRF) for supporting his position as a Ramanujan fellow at HRI, Prayagraj, India, under project file no. RJF/2025/000642. The initial version of this article has been made when AJ was a post-doctoral fellow at IMATI-CNR, Pavia, Italy and he would like to thank the project PRIN 2022YXWSLR “BOUNDARY ANALYSIS FOR DISPERSIVE AND VISCOUS FLUIDS” - DIT.PN012.008 for supporting his postdoctoral position at IMATI-CNR.
%AJ would like to thank the project PRIN 2022YXWSLR “BOUNDARY ANALYSIS FOR DISPERSIVE AND VISCOUS FLUIDS” - DIT.PN012.008 for supporting his postdoctoral position at IMATI-CNR, Pavia, Italy.

\end{document}